\documentclass[11pt]{book}  
\usepackage{amssymb}
\usepackage{amsfonts}
\usepackage{amsmath}
\usepackage{geometry}
\usepackage{bm}
\numberwithin{equation}{section}
\usepackage[toc,page]{appendix}
\usepackage{rotating}
\usepackage{graphicx,pdflscape}
\usepackage{qtree}
\newtheorem{theorem}{Theorem}[section]
\newtheorem{definition}[theorem]{Definition}

\newtheorem{remark}[theorem]{Remark}

\newenvironment{proof}[1][Proof]{\begin{trivlist}
\item[\hskip \labelsep {\bfseries #1}]}{\end{trivlist}}
\newcommand{\qed}{\nobreak \ifvmode \relax \else
      \ifdim\lastskip<1.5em \hskip-\lastskip
      \hskip1.5em plus0em minus0.5em \fi \nobreak
      \vrule height0.75em width0.5em depth0.25em\fi}

\raggedright                            

\title{\bf Complete polynomials using 3-term and reversible 3-term recurrence formulas (3TRF and R3TRF)}    
\author{\textsc{Yoon-Seok Choun} \footnote{The first series is `Special functions and three-term recurrence formula (3TRF).' It is available as arXiv.  The second series is `Special functions and reversible three-term recurrence formula (R3TRF).' It is also available as arXiv.
 This book is a third series.} \footnote{ychoun@gradcenter.cuny.edu; ychoun@gmail.com}\\
  Baruch College, The City University of New York,\\
  Natural Science Department, A506,\\
  17 Lexington Avenue,\\
  New York, NY 10010
} 
\date{\today}                           

\begin{document}
\frontmatter                            
\maketitle 
\tableofcontents                        
\mainmatter 
\chapter*{Preface}
The recurrence relation of coefficients starts to appear by putting a function $ y(x)=\sum_{n=0}^{\infty }c_n x^{n+\lambda } $, where $\lambda $ is an indicial root, into a linear ordinary differential equation (ODE). There can be between 2-term and infinity-term in the recursive relation. From the past until now, we have built the power series in the closed form for the 2-term recurrence relation in a linear ODE.
The recurrence relation for the 2-term is given by
\begin{equation}
c_{n+1}= A_n c_n \nonumber
\end{equation}
Power series solutions for the 2-term recurrence relation in a linear ODE have an infinite series and a polynomial.

Currently, the power series in the closed form and the integral representation for more than 3-term recurrence relation are unknown. 
In general, the recurrence relation for the 3-term is written by
\begin{equation}
c_{n+1}=A_n \;c_n +B_n \;c_{n-1} \hspace{1cm};n\geq 1
\nonumber
\end{equation}
where
\begin{equation}
c_1= A_0 \;c_0
\nonumber
\end{equation}
where $c_1= A_0 \;c_0$ and $\lambda $ is an indicial root. On the above, $A_n$ and $B_n$ are themselves polynomials of degree $m$: for the second-order ODEs, a numerator and a denominator of $A_n$ are usually equal or less than polynomials of degrees 2. Likewise, a numerator and a denominator of $B_n$ are also equal or less than polynomials of degrees 2

With my definition, power series solutions of the 3-term recurrence relation in a linear ODE has 4 types which are (1) an infinite series, (2) a polynomial which makes $B_n$ term terminated, referred as `a polynomial of type 1,' (3) a polynomial which makes $A_n$ term terminated, defined as `a polynomial of type 2,' and (4) a polynomial which makes $A_n$ and $B_n$ terms terminated, designated as `a polynomial of type 3' or `complete polynomial.' 

Heun, Grand Confluent Hypergeoemtric (GCH), Mathieu, Lam\'{e}, Confluent Heun (CH) and Double Confluent Heun (DCH) differential equations have the 3-term recurrence relation between successive coefficients.

In the first series ``Special functions and three term recurrence formula (3TRF)'',  I show how to obtain power series solutions of the above first four, namely Heun, GCH, Mathieu, Lam\'{e} equations for an infinite series and a polynomial of type 1. The method of proof for an infinite series and a polynomial of type 1 in the 3-term recurrence relation is called as three term recurrence formula (3TRF). And integral forms and generating functions of the above 4 equations are constructed analytically.

In the second series ``Special functions and reversible three-term recurrence formula (R3TRF)'',  I show how to obtain (1) power series solutions, (2) integral solutions and (3) generating functions of the above 5 equations (Heun, GCH, Mathieu, Lam\'{e} and CH equations) for an infinite series and a polynomial of type 2. The method of proof for an infinite series and a polynomial of type 2 in the 3-term recurrence relation is called as reversible three term recurrence formula (R3TRF).

In this series, natural numbers $\mathbb{N}_{0}$ means $\{ 0,1,2,3,\cdots \}$. Pochhammer symbol $(x)_n$ is used to represent the rising factorial: $(x)_n = \frac{\Gamma (x+n)}{\Gamma (x)}$. The complete polynomial has two different types which are (1) the first species complete polynomial and (2) the second species complete polynomial. The former is applicable if there are only one eigenvalue in $B_n$ term and an eigenvalue in $A_n$ term. And the latter is applicable if there are two eigenvalues in $B_n$ term and an eigenvalue in $A_n$ term.  
By applying 3TRF and R3TRF, I generalize the 3-term recurrence relation in the above 5 equations (Heun, GCH, Lam\'{e}, CH and DCH equations) for complete polynomials of two types in the form of power series expansions.

\section*{Structure of book}

\textbf{Chapter 1.} ``Complete polynomial using the three-term recurrence formula and its applications''--by applying 3TRF, I generalize the 3-term recurrence relation for complete polynomials of the first and second species in a linear ODE.
\vspace{1mm}

\textbf{Chapter 2.} ``Complete polynomial using the reversible three-term recurrence formula and its applications''--by applying R3TRF, I generalize the 3-term recurrence relation for the first and second species complete polynomials in a linear ODE.
\vspace{1mm}



\textbf{Chapter 3.} ``Complete polynomials of Heun equation using three-term recurrence formula''--By applying $A_n$ and $B_n$ terms in a recurrence relation of Heun equation into general summation formulas of the first and second species complete polynomials using 3TRF, power series solutions of Heun equation for their polynomials of type 3 are consturcted.
\vspace{1mm}

\textbf{Chapter 4.} ``Complete polynomials of Heun equation using reversible three-term recurrence formula''--By applying $A_n$ and $B_n$ terms in a recurrence relation of Heun equation into general summation formulas of the first and second species complete polynomials using R3TRF, power series solutions of Heun equation for their polynomials of type 3 are consturcted.
\vspace{1mm}
 
\textbf{Chapter 5.} ``Complete polynomials of Confluent Heun equation about the regular singular point at zero''--By applying $A_n$ and $B_n$ terms in a recurrence relation of Confluent Heun equation (CHE) around $x=0$ into general summation formulas of the first species complete polynomials using 3TRF and R3TRF, power series solutions of the CHE for their polynomials of type 3 are consturcted. 
\vspace{1mm}

\textbf{Chapter 6.} ``Complete polynomials of Confluent Heun equation about the irregular singular point at infinity''--By applying $A_n$ and $B_n$ terms in a recurrence relation of the CHE around $x=\infty $ into general summation formulas of the first and second species complete polynomials using 3TRF and R3TRF, power series solutions of the CHE for their polynomials of type 3 are consturcted. 
\vspace{1mm}

\textbf{Chapter 7.} ``Complete polynomials of Grand Confluent Hypergeometric equation about the regular singular point at zero''--By applying $A_n$ and $B_n$ terms in a recurrence relation of Grand Confluent Hypergeometric equation (GCH) around $x=0$ into general summation formulas of the first species complete polynomials using 3TRF and R3TRF, power series solutions of the GCH for their polynomials of type 3 are consturcted. 
\vspace{1mm}

\textbf{Chapter 8.} ``Complete polynomials of Grand Confluent Hypergeometric equation about the irregular singular point at infinity''--By applying $A_n$ and $B_n$ terms in a recurrence relation of the GCH around $x=\infty $ into general summation formulas of the first and second species complete polynomials using 3TRF and R3TRF, power series solutions of the GCH for their polynomials of type 3 are consturcted. 
\vspace{1mm}

\textbf{Chapter 9.} ``Complete polynomials of Lam\'{e} equation''--By applying $A_n$ and $B_n$ terms in recurrence relations of Lam\'{e} equations in the algebraic and Weierstrass's forms into general summation formulas of the first species complete polynomials using 3TRF and R3TRF, power series solutions of Lam\'{e} equations in two different forms for their polynomials of type 3 are constructed. 
\vspace{1mm}

\textbf{Chapter 10.} ``Double Confluent Heun functions using reversible three term recurrence formula''-- (1) apply R3TRF and analyze power series expansions in closed forms of Double Confluent Heun equation (DCHE) with irregular singularities at the origin and infinity and their integral forms for a polynomial of type 2. (2) apply R3TRF and construct generating functions for Double confluent Heun polynomials of type 2 around $x=0$ and $x=\infty $. 
\vspace{1mm}

\textbf{Chapter 11.} ``Complete polynomials of Double Confluent Heun equation''--By applying $A_n$ and $B_n$ terms in recurrence relations of the DCHE around $x=0$ and $x=\infty $ into general summation formulas of the first species complete polynomials using 3TRF and R3TRF, power series solutions of the DCHE for their polynomials of type 3 are consturcted. 

\section*{Acknowledgements}
I thank Bogdan Nicolescu.  The endless discussions I had with him were of great joy.  
\chapter{Complete polynomial using the three-term recurrence formula and its applications}
\chaptermark{Complete polynomial using 3TRF} 
In the first series\cite{Choun2012}, by substituting a power series with unknown coefficients into a linear ordinary differential equation (ODE), I generalize the three term recurrence relation for an infinite series and a polynomial which makes $B_n$ term terminated.
 
On chapter 1 in the second series \cite{Choun2013}, by putting a Taylor series with unknown coefficients into a linear ODE, I generalize the three term recurrence relation in a backward for an infinite series and a polynomial which makes $A_n$ term terminated.

In this chapter, by using Frobenius method, I generalize the three term recurrence relation for a polynomial which makes $A_n$ and $B_n$ terms terminated at the same time in a linear ODE.  

\section{Introduction}
The linear (ordinary) differential equations have very long history, over 350 years: Newton and Leibniz are considered as frontiers of differential equations. 
In general, by using Frobenius method, we attempt to seek a power series solution to certain differential equations especially linear homogeneous ODEs. The $n^{th}$ term recurrence relation (there can be between two and infinity term of coefficients in the recursive relation) for coefficients starts to appear since we substitute a power series with unknown coefficients into the ODE. 
The power series of the 2-term recursive relation in an ODE and its definite or contour integral forms have been constructed analytically. Currently, due to its complex mathematical calculations, the analytic solutions of a power series for the recursion relation having the recurrence relation among three different consecutive coefficients in a linear ODE have been obscured including its integral representation. \cite{Erde1955,Hobs1931,Whit1952}   
For instance, Heun, Lam\'{e} and Mathieu equations are some of well-known differential equations: the power series in their equations consist of three different coefficients in the recurrence relation. 

Until the 20th century, we have built the power series in the closed form of the two term recurrence relation in a linear ODE to describe physical phenomenons in linear systems. However, since 21th century, more precisely when modern physics (quantum mechanic, QCD, supersymmetric field theories, string theories, general relativity, etc) are born in the world, we do not have the 2-term recurrence relation of the power series in a linear ODE any more. More than three term recurrence relation of the power series might be required in most of difficult physics problems. 
 
In general, by using Frobenius method, the 2-term recursive relation in any linear ODEs is taken by
\begin{equation}
c_{n+1}=A_n \;c_n \hspace{1cm};n\geq 0
\nonumber
\end{equation}
And the 3-term recurrence relation in any linear ODEs is given by
\begin{equation}
c_{n+1}=A_n \;c_n +B_n \;c_{n-1} \hspace{1cm};n\geq 1
\nonumber
\end{equation}
where
\begin{equation}
c_1= A_0 \;c_0
\nonumber
\end{equation}
On the above, $A_n$ and $B_n$ are themselves polynomials of degree $m$: for the second-order ODEs, a numerator and a denominator of $A_n$ are usually equal or less than polynomials of degrees 2. Likewise, a numerator and a denominator of $B_n$ are also equal or less than polynomials of degrees 2.   
 
As we all know, the power series of the two term recursive relation in a linear ODE have two types which are an infinite series and a polynomial. In contrast, with my definition, the power series in the three term recurrence relation of a linear ODE has an infinite series and three types of polynomials: (1) a polynomial which makes $B_n$ term terminated; $A_n$ term is not terminated, (2) a polynomial which makes $A_n$ term terminated; $B_n$ term is not terminated, (3) a polynomial which makes $A_n$ and $B_n$ terms terminated.\cite{Choun2012,Choun2013}\footnote{If $A_n$ and $B_n$ terms are not terminated, it turns to be an infinite series.}

In the first series\cite{Choun2012}, by substituting a Taylor series with unknown coefficients into a linear ODE, I show how to obtain the Frobenius solutions of the three term recurrence relation for an infinite series and a polynomial of type 1: The sequence $c_n$ is composed of combinations of $A_n$ and $B_n$ terms. I observe the term of sequence $c_n$ which includes zero term of $A_n's$, one term of $A_n's$, two terms of $A_n's$, three terms of $A_n's$, etc to construct the general expressions of the power series for an infinite series and a polynomial of type 1. The method of proof for an infinite series and a polynomial of type 1 in the 3-term recurrence relation is called as three term recurrence formula (3TRF).

In the second series\cite{Choun2013}, by substituting a power series with unknown coefficients into a linear ODE, I generalize the three term recurrence relation in a backward for an infinite series and a polynomial of type 2: I observe the term of sequence $c_n$ which includes zero term of $B_n's$, one term of $B_n's$, two terms of $B_n's$, three terms of $B_n's$, etc to construct the general expressions of the power series for an infinite series and a polynomial of type 2. The method of proof for an infinite series and a polynomial of type 2 in the 3-term recurrence relation is designated as reversible three term recurrence formula (R3TRF).

The general power series expansion for an infinite series using 3TRF is equivalent to a Taylor series expression for an infinite series using R3TRF. 
In the former, $A_n$ is the leading term of sub-power series in a power series form of a function $y(x)$. In the latter, $B_n$ is the leading term of sub-power series in a Taylor series expansion of a function $y(x)$. 

In the mathematical definition, a polynomial of type 3 is referred as the spectral polynomial in general. In the three term recursive relation of a linear ODE, I categorize a type 3 polynomial as a complete polynomial.
In this chapter I generalize three term recurrence relation in the form of power series expansion for a complete polynomial by applying 3TRF. 

\section{The general 3-term recursion relation of Lam\'{e} equation}
Lam\'{e} equation (ellipsoidal harmonic equation), is a second-order linear ODE which has four regular singular points, was introduced by Gabriel Lam\'{e} (1837)\cite{Lame1837}.
This equation is derived from the solution of Laplace equation in elliptic coordinates in the technique of separation of variables. 
In contrast, Lam\'{e} Wave equation, is also a second-order linear ODE which has four regular singular points, arises when the method of separation of variables is applied to the Helmholtz equation in ellipsoidal coordinates.\cite{Erde1955,Arsc1964} 
The various scholars just leave the general solution of Lam\'{e} equation, having the 3-term recursive relation in a power series, as solutions of recurrences because of its complicated computations.

Let us look at the Lam\'{e} differential equation in Weierstrass's form,
\begin{equation}
\frac{d^2{y}}{d{z}^2} = \{ \alpha (\alpha +1)\rho^2\;sn^2(z,\rho )-h\} y(z)\label{eq:1}
\end{equation}
where $\rho$, $\alpha $  and $h$ are real parameters such that $0<\rho <1$ and $\alpha \geq -\frac{1}{2}$.
The algebraic form, which is obtained from (\ref{eq:1}) by setting $\xi= sn^2(z,\rho) $ as an independent variable is
\begin{equation}
\frac{d^2{y}}{d{\xi }^2} + \frac{1}{2}\left(\frac{1}{\xi } +\frac{1}{\xi -1} + \frac{1}{\xi -\rho ^{-2}}\right) \frac{d{y}}{d{\xi }} +  \frac{-\alpha (\alpha +1) \xi +h\rho ^{-2}}{4 \xi (\xi -1)(\xi -\rho ^{-2})} y(\xi ) = 0\label{eq:2}
\end{equation}
This is a Fuchsian ordinary differential equation with four regular singular singularities at $\xi=0, 1, \rho ^{-2}, \infty $; the exponents at the first three are all $\{0, 1/2\}$, and those at $\infty$ are $\{-\alpha/2,(\alpha +1)/2\}$. 
 
we assume the solution takes the form 
\begin{equation}
y(\xi) = \sum_{n=0}^{\infty } c_n \xi^{n+\lambda }
\label{eq:3}
\end{equation}
where $\lambda $ is an indicial root. 
Substituting (\ref{eq:3}) into (\ref{eq:2}) gives for the coefficients $c_n$ the recurrence relations
\begin{equation}
c_{n+1}=A_n \;c_n +B_n \;c_{n-1} \hspace{1cm};n\geq 1
\label{eq:4}
\end{equation}
where,
\begin{subequations}
\begin{equation}
A_n = \frac{ (1+\rho ^2)(n+\lambda )^2 -\frac{1}{4}h}{(n+\lambda +1)\left( n+\lambda +\frac{1}{2}\right)}
\label{eq:5a}
\end{equation}
\begin{equation}
B_n = -\rho ^2 \frac{\left( n+\lambda -\frac{1}{2}+\frac{\alpha }{2}\right) \left( n+\lambda -1-\frac{\alpha }{2}\right)}{(n+\lambda +1)\left( n+\lambda +\frac{1}{2}\right)}
\label{eq:5b}
\end{equation}
\begin{equation}
c_1= A_0 \;c_0 \label{eq:5c}
\end{equation}
\end{subequations}
We have two indicial roots which are $\lambda =0$ and $\frac{1}{2}$.

For the Lam\'{e} polynomial of type 1 in Ref.\cite{Chou2012g}, I treat $\alpha $ as a fixed value and  $h$ as a free variable: (1) I define $\alpha = 2(2\alpha_j +j)$ or  $-2(2\alpha_j +j)-1$ where $j, \alpha_j \in \mathbb{N}_{0}$ for the Lam\'{e} polynomial of type 1 of the first kind. (2) I define $\alpha = 2(2\alpha_j +j)+1$ or $-2(2\alpha_j +j+1)$ for the Lam\'{e} polynomial of type 1 of the second kind.\footnote{If we take $\alpha \geq -\frac{1}{2}$, $\alpha =  -2(2\alpha_j +j)-1 $ and $-2(2\alpha_j +j+1)$  are not available any more. In this chapter, I consider $\alpha $ as arbitrary.} 

For the Lam\'{e} polynomial of type 2 in the chapter 9 of Ref.\cite{Choun2013}, I treat $h$ as a fixed value and $\alpha$ as a free variable: (1) I define $h= 4(1+\rho ^2)(h_j+2j )^2 $ where $j,h_j \in \mathbb{N}_{0}$ for the Lam\'{e} polynomial of type 2 of the first kind. (2) I define $h= 4(1+\rho ^2)\left( h_j+2j +1/2\right)^2$ for the Lam\'{e} polynomial of type 2 of the second kind.

For the Lam\'{e} polynomial of type 3 (complete polynomial), let us suppose that $\alpha $ is $-2\left( j+1/2+\lambda \right)$ or $2(j+\lambda )$ where $j=0,1,2,\cdots$, and also an accessory parameter $h$ has been so chosen that the coefficients $c_{j+1}=0$: If $\lambda =0$ for the Lam\'{e} polynomial of the first kind, $\alpha  =-1,-3,-5,\cdots$ or $0,2,4,\cdots$. As $\lambda =\frac{1}{2}$ for the Lam\'{e} polynomial of the second kind, $\alpha  =-2,-4,-6,\cdots$ or $1,3,5,\cdots$.

Moreover, (\ref{eq:4}) with $n=j+1$ then gives successively $c_{j+1}=c_{j+2}=c_{j+3}=\cdots=0$ and the solution becomes a polynomial in $\xi $ of degree $j$. Now the condition $c_j=0$ is definitely an polynomial equation of degree $j+1$ for the determination of $h$.
In general, various authors have used matrix forms to determine the characteristic value of an accessory parameter in a 3-term recursive relation of a linear ODE.

In this chapter, instead of using the matrix representation, I generalize the algebraic equation for the determination of the accessory parameter in the form of partial sums of the sequences $\{A_n\}$ and $\{B_n\}$ using 3TRF for computational practice. I also generalize the 3-term recurrence relation in a linear ODE for a complete polynomial in the form of a power series expansion in the closed form.
\newpage
\section{Complete polynomial using 3TRF}
A complete polynomial (a polynomial of type 3) has two different types:

(1) If there is only one eigenvalue in $B_n$ term of 3-term recurrence relation, there are multi-valued roots of an eigenvalue (spectral parameter) in $A_n$ term.

(2) If there are two eigenvalues in $B_n$ term of 3-term recurrence relation, an eigenvalue in $A_n$ term has only one valued root. 

I designate the former as the first species complete polynomial and the latter as the second species complete polynomial. 
For instance, in general solutions of Heun's equation around $x=0$ \cite{Heun1889,Ronv1995}, if an exponent parameter $\alpha $ (or $\beta $) and an accessory parameter $q$ are fixed constants, the first species complete polynomial must to be utilized. And as parameters $\alpha $, $\beta $ and $q$ are fixed constants, the second species complete polynomial must to be applied. 

I denominate the mathematical formulas of a polynomial which makes $A_n$ and $B_n$ terms terminated, by allowing $A_n$ as the leading term in each sub-power series of the general power series, as ``complete polynomials using 3-term recurrence formula (3TRF)'' that will be expressed below.
\subsection{The first species complete polynomial using 3TRF} 
Assume that
\begin{equation}
c_1= A_0 \;c_0
\label{eq:6}
\end{equation}
(\ref{eq:6}) is a necessary boundary condition. The 3-term recurrence relation in all differential equations follows (\ref{eq:6}).
\begin{equation}
\prod _{n=a_i}^{a_i-1} B_n =1 \hspace{1cm} \mathrm{where}\;a_i=0,1,2,\cdots 
\label{eq:7a}
\end{equation}
(\ref{eq:7a}) is also a necessary condition. Every differential equations take satisfied with (\ref{eq:7a}). 
\newpage
My definition of $B_{i,j,k,l}$ refer to $B_{i}B_{j}B_{k}B_{l}$. Also, $A_{i,j,k,l}$ refer to $A_{i}A_{j}A_{k}A_{l}$. For $n=0,1,2,3,\cdots $, (\ref{eq:4}) gives
\begin{equation}
\begin{tabular}{  l  }
  \vspace{2 mm}
  $c_0$ \\
  \vspace{2 mm}
  $c_1 = A_0 c_0 $ \\
  \vspace{2 mm}
  $c_2 = (B_1+ A_{0,1}) c_0  $ \\
  \vspace{2 mm}
  $c_3 = ( A_0 B_2 +A_2 B_1 +A_{0,1,2}) c_0  $\\
  \vspace{2 mm}
  $c_4 = ( B_{1,3}+ A_{0,1} B_3 + A_{0,3} B_2 + A_{2,3} B_1 + A_{0,1,2,3}) c_0  $ \\
  \vspace{2 mm}
  $c_5 = ( A_0 B_{2,4}+ A_2 B_{1,4}+ A_4 B_{1,3} + A_{0,1,2}B_4 + A_{0,1,4} B_3  + A_{0,3,4} B_2 + A_{2,3,4} B_1 + A_{0,1,2,3,4}) c_0  $  \\     
  \vspace{2 mm}
  $c_6 = (B_{1,3,5} + A_{0,1} B_{3,5} + A_{0,3} B_{2,5} + A_{0,5} B_{2,4}+ A_{2,3} B_{1,5} + A_{2,5} B_{1,4} + A_{4,5} B_{1,3} $ \\
  \vspace{2 mm}
  \hspace{0.8 cm} $ + A_{0,1,2,3} B_5 + A_{0,1,2,5}B_4+ A_{0,1,4,5} B_3  + A_{0,3,4,5} B_2 + A_{2,3,4,5} B_1 + A_{0,1,2,3,4,5}) c_0 $ \\     
  \vspace{2 mm}
  $c_7 = ( A_0 B_{2,4,6} + A_2 B_{1,4,6}+ A_4 B_{1,3,6}+ A_6 B_{1,3,5} + A_{0,1,2} B_{4,6} + A_{0,1,4} B_{3,6} + A_{0,1,6} B_{3,5} + A_{0,3,4} B_{2,6}  $ \\
  \vspace{2 mm}
  \hspace{0.8 cm} $+ A_{0,3,6} B_{2,5} + A_{0,5,6} B_{2,4} + A_{2,3,4} B_{1,6} + A_{2,3,6} B_{1,5} + A_{2,5,6} B_{1,4} + A_{4,5,6} B_{1,3}  $\\
  \vspace{2 mm}                 
  \hspace{0.8 cm} $ + A_{0,1,2,3,4}B_6 + A_{2,3,4,5,6} B_1+ A_{0,3,4,5,6} B_2 + A_{0,1,4,5,6} B_3 + A_{0,1,2,5,6}B_4 +A_{0,1,2,3,6} B_5+A_{0,1,2,3,4,5,6}) c_0 $ \\
  \vspace{2 mm}
  $c_8 = ( B_{1,3,5,7}+ A_{0,1}B_{3,5,7} + A_{0,3}B_{2,5,7} + A_{0,5}B_{2,4,7} + A_{0,7} B_{2,4,6} + A_{2,3}B_{1,5,7} + A_{2,5}B_{1,4,7} + A_{2,7} B_{1,4,6}  $ \\
  \vspace{2 mm}
  \hspace{0.8 cm} $  + A_{4,5}B_{1,3,7}+ A_{4,7} B_{1,3,6} + A_{6,7} B_{1,3,5}  + A_{0,1,2,3} B_{5,7} + A_{0,1,2,5}B_{4,7} + A_{0,1,2,7} B_{4,6}  + A_{0,1,4,5}B_{3,7}  $\\
  \vspace{2 mm}
  \hspace{0.8 cm} $+ A_{0,1,4,7} B_{3,6} + A_{0,1,6,7} B_{3,5} + A_{0,3,4,5} B_{2,7} + A_{0,3,4,7} B_{2,6} + A_{0,3,6,7} B_{2,5} + A_{0,5,6,7} B_{2,4} + A_{2,3,4,5} B_{1,7}  $\\
  \vspace{2 mm}                 
  \hspace{0.8 cm} $ + A_{2,3,4,7} B_{1,6} + A_{2,3,6,7} B_{1,5} + A_{2,5,6,7} B_{1,4}+ A_{4,5,6,7} B_{1,3}+ A_{0,1,2,3,4,5} B_{7}+ A_{0,1,2,3,4,7}B_6 $\\
   \vspace{2 mm}
   \hspace{0.8 cm} $  + A_{0,1,2,3,6,7} B_5  + A_{0,1,2,5,6,7}B_4 + A_{0,1,4,5,6,7} B_3 + A_{0,3,4,5,6,7} B_2 +A_{2,3,4,5,6,7} B_1 +A_{0,1,2,3,4,5,6,7}) c_0$\\ 
\hspace{2 mm}\large{\vdots} \hspace{5cm}\large{\vdots}\\ 
\end{tabular}
\label{eq:7}
\end{equation}
In (\ref{eq:7}) the number of individual sequence $c_n$ follows Fibonacci sequence: 1,1,2,3,5,8,13,21,34,55,$\cdots$.
The sequence $c_n$ consists of combinations $A_n$ and $B_n$ in (\ref{eq:7}). 

The coefficient $c_n$ with even subscripts consists of even terms of $A_n's$. The sequence $c_n$ with odd subscripts consists of odd terms of $A_n's$. And classify sequences $c_n$ to its even and odd parts in (\ref{eq:7}). 
\begin{equation}
\begin{tabular}{  l  l }
  \vspace{2 mm}
  \large{$c_0 = \tilde{A}^0$} &\hspace{1cm} \large{$c_1= \tilde{A}^1$}  \\
  \vspace{2 mm}
  \large{$c_2 = \tilde{A}^0+\tilde{A}^2$} &\hspace{1cm} \large{$c_3 = \tilde{A}^1+\tilde{A}^3$} \\
  \vspace{2 mm}
  \large{$c_4 = \tilde{A}^0+\tilde{A}^2+\tilde{A}^4$} &\hspace{1cm}  \large{$c_5 = \tilde{A}^1+\tilde{A}^3+\tilde{A}^5$}\\
  \vspace{2 mm}
  \large{$c_6 = \tilde{A}^0+\tilde{A}^2+\tilde{A}^4+\tilde{A}^6$} &\hspace{1cm}  \large{$c_7 = \tilde{A}^1+\tilde{A}^3+\tilde{A}^5+\tilde{A}^7$}\\
  \vspace{2 mm}
  \large{$c_8 = \tilde{A}^0+\tilde{A}^2+\tilde{A}^4+\tilde{A}^6+\tilde{A}^8$} &\hspace{1cm} \large{$c_9 = \tilde{A}^1+\tilde{A}^3+\tilde{A}^5+\tilde{A}^7+\tilde{A}^9$} \\
 \hspace{2cm} \large{\vdots} & \hspace{3cm}\large{\vdots} \\
\end{tabular}
\label{eq:8}
\end{equation}
In the above, $\tilde{A}^{\tau } $= $\tau$ terms of $A_n's$ where $\tau =0,1,2,\cdots$.

When a function $y(x)$, analytic at $x=0$, is expanded in a power series, we write
\begin{equation}
y(x)= \sum_{n=0}^{\infty } c_n x^{n+\lambda }= \sum_{\tau =0}^{\infty } y_{\tau}(x) = y_0(x)+ y_1(x)+y_2(x)+ \cdot \label{eq:9}
\end{equation}
where
\begin{equation}
y_{\tau}(x)= \sum_{l=0}^{\infty } c_l^{\tau} x^{l+\lambda }\label{eq:10}
\end{equation}
$\lambda $ is the indicial root. $y_{\tau}(x)$ is sub-power series that have sequences $c_n$ including $\tau$ term of $A_n$'s in (\ref{eq:7}). For example $y_0(x)$ has sequences $c_n$ including zero term of $A_n$'s in (\ref{eq:7}), $y_1(x)$ has sequences $c_n$ including one term of $A_n$'s in (\ref{eq:7}), $y_2(x)$ has sequences $c_n$ including two term of $A_n$'s in (\ref{eq:7}), etc.

First observe the term inside parentheses of sequences $c_n$ which does not include any $A_n$'s in (\ref{eq:7}) and (\ref{eq:8}): $c_n$ with even subscripts ($c_0$, $c_2$, $c_4$,$\cdots$). 

\begin{equation}
\begin{tabular}{  l  }
  \vspace{2 mm}
  $c_0$ \\
  \vspace{2 mm}
  $c_2 = B_1 c_0  $ \\
  \vspace{2 mm}
  $c_4 = B_{1,3} c_0  $ \\
  \vspace{2 mm}
  $c_6 = B_{1,3,5}c_0 $ \\
  \vspace{2 mm}
  $c_8 = B_{1,3,5,7}c_0 $\\
  \hspace{2 mm}
  \large{\vdots}\hspace{1cm}\large{\vdots} \\ 
\end{tabular}
\label{eq:11}
\end{equation}
(\ref{eq:11}) gives the indicial equation
\begin{equation}
c_{2n}= c_0  \prod _{i_0=0}^{n-1}B_{2i_0+1} \hspace{1cm} \mathrm{where}\;n=0,1,2,\cdots 
\label{eq:12}
\end{equation}
Substitute (\ref{eq:12}) in (\ref{eq:10}) putting $\tau = 0$. 
\begin{eqnarray}
y_0^m(x) &=& c_0 \sum_{n=0}^{m} \left\{ \prod _{i_0=0}^{n-1}B_{2i_0+1} \right\} x^{2n+\lambda }\label{eq:13a}\\
&=& c_0 \sum_{i_0=0}^{m} \left\{ \prod _{i_1=0}^{i_0-1}B_{2i_1+1} \right\} x^{2i_0+\lambda } \label{eq:13b}
\end{eqnarray}
(\ref{eq:13a}) and (\ref{eq:13b}) are the sub-power series that has sequences $c_0, c_2, c_4, \cdots, c_{2m}$ including zero term of $A_n$'s where $m\in \mathbb{N}_0$. If $m\rightarrow \infty $, the sub-power series $y_0^m(x)$ turns to be an infinite series $y_0(x)$.

Observe the terms inside parentheses of sequence $c_n$ which include one term of $A_n$'s in (\ref{eq:7}) and (\ref{eq:8}): $c_n$ with odd subscripts ($c_1$, $c_3$, $c_5$,$\cdots$). 

\begin{equation}
\begin{tabular}{  l  }
  \vspace{2 mm}
  $c_1= A_0 c_0$ \\
  \vspace{2 mm}
  $c_3 = \left\{ A_0 \Big( \frac{B_2}{1}\Big)1+ A_2\Big(\frac{B_2}{B_2}\Big) B_1\right\} c_0  $ \\
  \vspace{2 mm}
  $c_5 = \left\{ A_0 \Big( \frac{B_{2,4}}{1}\Big)1+ A_2\Big(\frac{B_{2,4}}{B_2}\Big) B_1+ A_4\Big(\frac{B_{2,4}}{B_{2,4}}\Big) B_{1,3}\right\} c_0  $ \\
  \vspace{2 mm}
  $c_7 = \left\{ A_0 \Big( \frac{B_{2,4,6}}{1}\Big)1+ A_2\Big(\frac{B_{2,4,6}}{B_2}\Big) B_1+ A_4\Big(\frac{B_{2,4,6}}{B_{2,4}}\Big) B_{1,3}+  A_6\Big(\frac{B_{2,4,6}}{B_{2,4,6}}\Big) B_{1,3,5}\right\} c_0  $ \\
  \vspace{2 mm}
  $c_9 = \left\{ A_0 \Big( \frac{B_{2,4,6,8}}{1}\Big)1+ A_2\Big(\frac{B_{2,4,6,8}}{B_2}\Big) B_1+ A_4\Big(\frac{B_{2,4,6,8}}{B_{2,4}}\Big) B_{1,3}+  A_6\Big(\frac{B_{2,4,6,8}}{B_{2,4,6}}\Big) B_{1,3,5} \right.$\\
  \vspace{2 mm}
  \hspace{0.8 cm} $+\left.  A_8\Big(\frac{B_{2,4,6,8}}{B_{2,4,6,8}}\Big) B_{1,3,5,7}\right\} c_0 $\\
  $c_{11} =\left\{ A_0 \Big( \frac{B_{2,4,6,8,10}}{1}\Big)1+ A_2\Big(\frac{B_{2,4,6,8,10}}{B_2}\Big) B_1+ A_4\Big(\frac{B_{2,4,6,8,10}}{B_{2,4}}\Big) B_{1,3}+  A_6\Big(\frac{B_{2,4,6,8,10}}{B_{2,4,6}}\Big) B_{1,3,5} \right.$\\
   \vspace{2 mm}
   \hspace{0.8 cm} $+\left.  A_8\Big(\frac{B_{2,4,6,8,10}}{B_{2,4,6,8}}\Big) B_{1,3,5,7}+  A_{10}\Big(\frac{B_{2,4,6,8,10}}{B_{2,4,6,8,10}}\Big)B_{1,3,5,7,9} \right\} c_0 $\\
  \hspace{2 mm}
  \large{\vdots}\hspace{3cm}\large{\vdots} \\ 
\end{tabular}
\label{eq:14}
\end{equation}
(\ref{eq:14}) gives the indicial equation
\begin{equation}
c_{2n+1}= c_0  \sum_{i_0=0}^{n} \left\{ A_{2i_0} \prod _{i_1=0}^{i_0-1}B_{2i_1+1} \prod _{i_2=i_0}^{n-1}B_{2i_2+2} \right\} 
\label{eq:15}
\end{equation}
Substitute (\ref{eq:15}) in (\ref{eq:10}) putting $\tau = 1$. 
\begin{eqnarray}
y_1^m(x)&=& c_0 \sum_{n=0}^{m}\left\{ \sum_{i_0=0}^{n} \left\{ A_{2i_0} \prod _{i_1=0}^{i_0-1}B_{2i_1+1} \prod _{i_2=i_0}^{n-1}B_{2i_2+2} \right\} \right\} x^{2n+1+\lambda } \label{eq:16a}\\ 
&=& c_0 \sum_{i_0=0}^{m}\left\{ A_{2i_0} \prod _{i_1=0}^{i_0-1}B_{2i_1+1}  \sum_{i_2=i_0}^{m} \left\{ \prod _{i_3=i_0}^{i_2-1}B_{2i_3+2} \right\}\right\} x^{2i_2+1+\lambda } \label{eq:16b}
\end{eqnarray}
(\ref{eq:16a}) and (\ref{eq:16b}) are the sub-power series that has sequences $c_1, c_3, c_5, \cdots, c_{2m+1}$ including one term of $A_n$'s. If $m\rightarrow \infty $, the sub-power series $y_1^m(x)$ turns to be an infinite series $y_1(x)$.

Observe the terms inside parentheses of sequence $c_n$ which include two terms of $A_n$'s in  (\ref{eq:7}) and (\ref{eq:8}): $c_n$ with even subscripts ($c_2$, $c_4$, $c_6$,$\cdots$). 

\begin{equation}
\begin{tabular}{  l  }
  \vspace{2 mm}
  $c_2= A_{0,1} c_0$ \\
  \vspace{2 mm}
  $c_4 = \Bigg\{ A_0 \left\{ A_1 \left( \frac{1}{1}\right) \left(\frac{B_{1,3}}{B_1}\right)1+ A_3 \left(\frac{B_2}{1}\right) \left( \frac{B_{1,3}}{B_{1,3}}\right)1\right\} + A_2 \left\{ A_3 \left( \frac{B_2}{B_2}\right) \left( \frac{B_{1,3}}{B_{1,3}}\right) B_1 \right\} \Bigg\} c_0  $ \\
  \vspace{2 mm}
  $c_6 = \Bigg\{ A_0 \left\{ A_1 \left( \frac{1}{1}\right) \left(\frac{B_{1,3,5}}{B_1}\right)1+ A_3 \left(\frac{B_2}{1}\right) \left( \frac{B_{1,3,5}}{B_{1,3}}\right)1+ A_5 \left(\frac{B_{2,4}}{1}\right) \left( \frac{B_{1,3,5}}{B_{1,3,5}}\right)1\right\}$\\
  \vspace{2 mm}
  \hspace{0.8 cm} $+ A_2 \left\{ A_3 \left( \frac{B_2}{B_2}\right) \left( \frac{B_{1,3,5}}{B_{1,3}}\right) B_1 + A_5 \left( \frac{B_{2,4}}{B_2}\right) \left( \frac{B_{1,3,5}}{B_{1,3,5}}\right) B_1\right\} + A_4 \left\{ A_5 \left( \frac{B_{2,4}}{B_{2,4}}\right) \left( \frac{B_{1,3,5}}{B_{1,3,5}}\right) B_{1,3} \right\} \Bigg\} c_0 $ \\
\vspace{2 mm}
    $c_8 = \Bigg\{ A_0 \left\{ A_1 \left( \frac{1}{1}\right) \left(\frac{B_{1,3,5,7}}{B_1}\right)1+ A_3 \left(\frac{B_2}{1}\right) \left( \frac{B_{1,3,5,7}}{B_{1,3}}\right)1+ A_5 \left(\frac{B_{2,4}}{1}\right) \left( \frac{B_{1,3,5,7}}{B_{1,3,5}}\right)1 +  A_7 \left(\frac{B_{2,4,6}}{1}\right) \left( \frac{B_{1,3,5,7}}{B_{1,3,5,7}}\right)1  \right\}$\\
\vspace{2 mm}
  \hspace{0.8 cm} $+ A_2 \left\{ A_3 \left( \frac{B_2}{B_2}\right) \left( \frac{B_{1,3,5,7}}{B_{1,3}}\right) B_1 + A_5 \left( \frac{B_{2,4}}{B_2}\right) \left( \frac{B_{1,3,5,7}}{B_{1,3,5}}\right) B_1 + A_7 \left( \frac{B_{2,4,6}}{B_2}\right) \left( \frac{B_{1,3,5,7}}{B_{1,3,5,7}}\right) B_1\right\}$\\
\vspace{2 mm}
  \hspace{0.8 cm} $+ A_4 \left\{ A_5 \left( \frac{B_{2,4}}{B_{2,4}}\right) \left( \frac{B_{1,3,5,7}}{B_{1,3,5}}\right) B_{1,3} +  A_7 \left( \frac{B_{2,4,6}}{B_{2,4}}\right) \left( \frac{B_{1,3,5,7}}{B_{1,3,5,7}}\right) B_{1,3}\right\}$\\
\vspace{2 mm}
  \hspace{0.8 cm} $+ A_6 \left\{ A_7 \Big( \frac{B_{2,4,6}}{B_{2,4,6}}\Big) \Big( \frac{B_{1,3,5,7}}{B_{1,3,5,7}}\Big) B_{1,3,5} \right\} \Bigg\} c_0$\\
  \hspace{2 mm}
  \large{\vdots}\hspace{5cm}\large{\vdots} \\  
\end{tabular}
\label{eq:17}
\end{equation}
\newpage
(\ref{eq:17}) gives the indicial equation
\begin{eqnarray}
 c_{2n+2} &=& c_0\sum_{i_0=0}^{n} \left\{ A_{2i_0}\sum_{i_1=i_0}^{n} \left\{ A_{2i_1+1} \prod _{i_2=0}^{i_0-1}B_{2i_2+1} \prod _{i_3=i_0}^{i_1-1}B_{2i_3+2}\prod _{i_4=i_1}^{n-1}B_{2i_4+3} \right\}\right\}  
\label{eq:18}\\
&=& c_0\sum_{i_0=0}^{n} \left\{ A_{2i_0} \prod _{i_1=0}^{i_0-1}B_{2i_1+1} \sum_{i_2=i_0}^{n} \left\{ A_{2i_2+1}  \prod _{i_3=i_0}^{i_2-1}B_{2i_3+2} \prod _{i_4=i_2}^{n-1}B_{2i_4+3} \right\}\right\} \label{eq:x18}   
\end{eqnarray}
Substitute (\ref{eq:x18}) in (\ref{eq:10}) putting $\tau = 2$. 
\begin{eqnarray}
 y_2^m(x) &=& c_0 \sum_{n=0}^{m}\left\{ \sum_{i_0=0}^{n} \left\{ A_{2i_0} \prod _{i_1=0}^{i_0-1}B_{2i_1+1} \sum_{i_2=i_0}^{n} \left\{ A_{2i_2+1}  \prod _{i_3=i_0}^{i_2-1}B_{2i_3+2} \prod _{i_4=i_2}^{n-1}B_{2i_4+3} \right\}\right\} \right\} x^{2n+2+\lambda } \hspace{2cm} \label{eq:19a}\\
&=& c_0 \sum_{i_0=0}^{m}\left\{ A_{2i_0} \prod _{i_1=0}^{i_0-1}B_{2i_1+1}  \sum_{i_2=i_0}^{m} \left\{ A_{2i_2+1} \prod _{i_3=i_0}^{i_2-1}B_{2i_3+2} \sum_{i_4=i_2}^{m}\left\{ \prod _{i_5=i_2}^{i_4-1}B_{2i_5+3}\right\}  \right\}\right\} x^{2i_4+2+\lambda } \label{eq:19b}
\end{eqnarray}
(\ref{eq:19a}) and (\ref{eq:19b}) are the sub-power series that has sequences $c_2, c_4, c_6, \cdots, c_{2m+2}$ including two terms of $A_n$'s. If $m\rightarrow \infty $, the sub-power series $y_2^m(x)$ turns to be an infinite series $y_2(x)$.

Observe the terms inside parentheses of sequence $c_n$ which include three terms of $A_n$'s in (\ref{eq:7}) and (\ref{eq:8}): $c_n$ with odd subscripts ($c_3$, $c_5$, $c_7$,$\cdots$). 
\begin{equation}
\begin{tabular}{  l  }
  \vspace{2 mm}
  $c_3= A_{0,1,2} \;c_0$ \\
  \vspace{2 mm}
  $c_5 = \Bigg\{ A_0 \bigg\{ A_1 \Big[ A_2 \cdot 1 \left( \frac{1}{1}\right) \left( \frac{1}{1}\right) \left(\frac{B_{4}}{1}\right)+ A_4 \cdot 1 \left( \frac{1}{1}\right) \left(\frac{B_3}{1}\right) \left( \frac{B_{4}}{B_{4}}\right)\Big] $\\
 \vspace{2 mm}
  \hspace{0.8 cm} $+ A_3 A_4 \cdot 1 \left( \frac{B_2}{1}\right) \left( \frac{B_{3}}{B_{3}}\right) \left( \frac{B_{4}}{B_{4}}\right) \bigg\} + A_2 \bigg\{ A_3 \Big[ A_4 B_1 \left( \frac{B_{2}}{B_{2}}\right)\left( \frac{B_{3}}{B_{3}}\right)\left( \frac{B_{4}}{B_{4}}\right)\Big] \bigg\} \Bigg\} c_0  $ \\
  \vspace{2 mm}
  $c_7 = \Bigg\{ A_0 \bigg\{ A_1 \Big[ A_2 \cdot 1 \left( \frac{1}{1}\right) \left( \frac{1}{1}\right) \left(\frac{B_{4,6}}{1}\right)+ A_4 \cdot 1 \left( \frac{1}{1}\right) \left(\frac{B_3}{1}\right) \left( \frac{B_{4,6}}{B_{4}}\right)+ A_6 \cdot 1 \left( \frac{1}{1}\right) \left(\frac{B_{3,5}}{1}\right) \left( \frac{B_{4,6}}{B_{{4,6}}}\right)\Big]$\\
 \vspace{2 mm}
  \hspace{0.8 cm} $+ A_3 \Big[ A_4 \cdot 1 \left( \frac{B_2}{1}\right) \left( \frac{B_{3}}{B_{3}}\right) \left( \frac{B_{4,6}}{B_{4}}\right) + A_6 \cdot 1 \left( \frac{B_2}{1}\right) \left( \frac{B_{3,5}}{B_{3}}\right) \left( \frac{B_{4,6}}{B_{4,6}}\right)\Big]+ A_5 \Big[ A_6\cdot 1 \left( \frac{B_{2,4}}{1}\right) \left( \frac{B_{3,5}}{B_{3,5}}\right) \left( \frac{B_{4,6}}{B_{4,6}}\right) \Big] \bigg\} $\\
 \vspace{2 mm}
  \hspace{0.8 cm} $+ A_2 \bigg\{ A_3 \Big[ A_4 B_1 \left( \frac{B_{2}}{B_{2}}\right)\left( \frac{B_{3}}{B_{3}}\right)\left( \frac{B_{4,6}}{B_{4}}\right)+ A_6 B_1 \left( \frac{B_{2}}{B_{2}}\right)\left( \frac{B_{3,5}}{B_{3}}\right)\left( \frac{B_{4,6}}{B_{4,6}}\right)\Big] $\\
\vspace{2 mm}
  \hspace{0.8 cm}$ + A_5 \Big[ A_6 B_1 \left( \frac{B_{2,4}}{B_{2}}\right) \left( \frac{B_{3,5}}{B_{3,5}}\right) \left( \frac{B_{4,6}}{B_{4,6}}\right) \Big] \bigg\}+ A_4 \bigg\{ A_5 \Big[ A_6 B_{1,3} \left( \frac{B_{2,4}}{B_{2,4}}\right) \left( \frac{B_{3,5}}{B_{3,5}}\right) \left( \frac{B_{4,6}}{B_{4,6}}\right) \Big]\bigg\} \Bigg\} c_0 $ \\
 \vspace{2 mm}
   $c_9 = \Bigg\{ A_0 \bigg\{ A_1 \Big[ A_2 \cdot 1 \left( \frac{1}{1}\right) \left( \frac{1}{1}\right) \left(\frac{B_{4,6,8}}{1}\right)+ A_4 \cdot 1 \left( \frac{1}{1}\right) \left(\frac{B_3}{1}\right) \left( \frac{B_{4,6,8}}{B_{4}}\right)+ A_6 \cdot 1 \left( \frac{1}{1}\right) \left(\frac{B_{3,5}}{1}\right) \left( \frac{B_{4,6,8}}{B_{{4,6}}}\right) $\\ 
 \vspace{2 mm}
\hspace{0.8 cm} $+ A_8 \cdot 1 \left( \frac{1}{1}\right) \left(\frac{B_{3,5,7}}{1}\right) \left( \frac{B_{4,6,8}}{B_{{4,6,8}}}\right)\Big]$\\
\vspace{2 mm}
\hspace{0.8 cm} $+ A_3 \Big[ A_4 \cdot 1 \left( \frac{B_2}{1}\right) \left( \frac{B_{3}}{B_{3}}\right) \left( \frac{B_{4,6,8}}{B_{4}}\right) + A_6 \cdot 1 \left( \frac{B_2}{1}\right) \left( \frac{B_{3,5}}{B_{3}}\right) \left( \frac{B_{4,6,8}}{B_{4,6}}\right)+ A_8 \cdot 1 \left( \frac{B_2}{1}\right) \left( \frac{B_{3,5,7}}{B_{3}}\right) \left( \frac{B_{4,6,8}}{B_{4,6,8}}\right) \Big]$\\
  \vspace{2 mm}
  \hspace{0.8 cm} $+ A_5 \Big[ A_6\cdot 1 \left( \frac{B_{2,4}}{1}\right) \left( \frac{B_{3,5}}{B_{3,5}}\right) \left( \frac{B_{4,6,8}}{B_{4,6}}\right)+ A_8\cdot 1 \left( \frac{B_{2,4}}{1}\right) \left( \frac{B_{3,5,7}}{B_{3,5}}\right) \left( \frac{B_{4,6,8}}{B_{4,6,8}}\right)\Big] $\\
  \vspace{2 mm}
  \hspace{0.8 cm} $+ A_7\Big[ A_8\cdot 1 \left( \frac{B_{2,4,6}}{1}\right) \left( \frac{B_{3,5,7}}{B_{3,5,7}}\right) \left( \frac{B_{4,6,8}}{B_{4,6,8}}\right) \Big] \bigg\} $\\
 \vspace{2 mm}
  \hspace{0.8 cm} $+ A_2 \bigg\{ A_3 \Big[ A_4 B_1 \left( \frac{B_{2}}{B_{2}}\right)\left( \frac{B_{3}}{B_{3}}\right)\left( \frac{B_{4,6,8}}{B_{4}}\right)+ A_6 B_1 \left( \frac{B_{2}}{B_{2}}\right)\left( \frac{B_{3,5}}{B_{3}}\right)\left( \frac{B_{4,6,8}}{B_{4,6}}\right) + A_8 B_1 \left( \frac{B_{2}}{B_{2}}\right)\left( \frac{B_{3,5,7}}{B_{3}}\right)\left( \frac{B_{4,6,8}}{B_{4,6,8}}\right)\Big]$\\
  \vspace{2 mm}
  \hspace{0.8 cm} $+ A_5 \Big[ A_6 B_1 \left( \frac{B_{2,4}}{B_{2}}\right) \left( \frac{B_{3,5}}{B_{3,5}}\right) \left( \frac{B_{4,6,8}}{B_{4,6}}\right) + A_8 B_1 \left( \frac{B_{2,4}}{B_{2}}\right) \left( \frac{B_{3,5,7}}{B_{3,5}}\right) \left( \frac{B_{4,6,8}}{B_{4,6,8}}\right)\Big] $\\
\vspace{2 mm}
  \hspace{0.8 cm} $+ A_7 \Big[ A_8 B_1 \left( \frac{B_{2,4,6}}{B_{2}}\right) \left( \frac{B_{3,5,7}}{B_{3,5,7}}\right) \left( \frac{B_{4,6,8}}{B_{4,6,8}}\right) \Big] \bigg\}$\\
\vspace{2 mm}
  \hspace{0.8 cm} $+ A_4 \bigg\{ A_5 \Big[ A_6 B_{1,3} \left( \frac{B_{2,4}}{B_{2,4}}\right) \left( \frac{B_{3,5}}{B_{3,5}}\right) \left( \frac{B_{4,6,8}}{B_{4,6}}\right) + A_8 B_{1,3} \left( \frac{B_{2,4}}{B_{2,4}}\right) \left( \frac{B_{3,5,7}}{B_{3,5}}\right) \left( \frac{B_{4,6,8}}{B_{4,6,8}}\right)\Big] $\\
\vspace{2 mm}
 \hspace{0.8 cm} $+ A_7 \Big[ A_8 B_{1,3} \left( \frac{B_{2,4,6}}{B_{2,4}}\right) \left( \frac{B_{3,5,7}}{B_{3,5,7}}\right) \left( \frac{B_{4,6,8}}{B_{4,6,8}}\right) \Big] \bigg\}+ A_6 \bigg\{ A_7 \Big[ A_8 B_{1,3,5}  \left( \frac{B_{2,4,6}}{B_{2,4,6}}\right) \left( \frac{B_{3,5,7}}{B_{3,5,7}}\right) \left( \frac{B_{4,6,8}}{B_{4,6,8}}\right)\Big] \bigg\} \Bigg\} c_0 $ \\
  \hspace{2 mm}
  \large{\vdots}\hspace{6cm}\large{\vdots} \\ 
\end{tabular}
\label{eq:20}
\end{equation}
(\ref{eq:20}) gives the indicial equation
\begin{eqnarray}
 c_{2n+3} &=& c_0 \sum_{i_0=0}^{n} \left\{ A_{2i_0}\sum_{i_1=i_0}^{n} \left\{ A_{2i_1+1}\sum_{i_2=i_1}^{n} \left\{ A_{2i_2+2} \prod _{i_3=0}^{i_0-1}B_{2i_3+1} \right.\right.\right.\nonumber\\
&&\times \left.\left.\left. \prod _{i_4=i_0}^{i_1-1}B_{2i_4+2} \prod _{i_5=i_1}^{i_2-1}B_{2i_5+3}\prod _{i_6=i_2}^{n-1}B_{2i_6+4} \right\} \right\} \right\} \label{eq:21}\\
&=& c_0 \sum_{i_0=0}^{n} \left\{ A_{2i_0}\prod _{i_1=0}^{i_0-1}B_{2i_1+1} \sum_{i_2=i_0}^{n} \left\{ A_{2i_2+1} \prod _{i_3=i_0}^{i_2-1}B_{2i_3+2} \right.\right. \nonumber\\
&&\times \left.\left. \sum_{i_4=i_2}^{n} \left\{ A_{2i_4 +2} \prod _{i_5=i_2}^{i_4-1}B_{2i_5+3}\prod _{i_6=i_4}^{n-1}B_{2i_6+4} \right\} \right\} \right\} \label{eq:x21} 
\end{eqnarray}
Substitute (\ref{eq:x21}) in (\ref{eq:10}) putting $\tau = 3$.
\begin{eqnarray}
y_3^m(x) &=& c_0 \sum_{n=0}^{m}\left\{ \sum_{i_0=0}^{n} \left\{ A_{2i_0}\prod _{i_1=0}^{i_0-1}B_{2i_1+1} \sum_{i_2=i_0}^{n} \left\{ A_{2i_2+1} \prod _{i_3=i_0}^{i_2-1}B_{2i_3+2}\right.\right.\right. \nonumber\\
&&\times \left.\left.\left. \sum_{i_4=i_2}^{n} \left\{ A_{2i_4 +2} \prod _{i_5=i_2}^{i_4-1}B_{2i_5+3}\prod _{i_6=i_4}^{n-1}B_{2i_6+4} \right\} \right\} \right\}\right\} x^{2n+3+\lambda } \label{eq:22a}\\
&=&  c_0 \sum_{i_0=0}^{m}\left\{ A_{2i_0} \prod _{i_1=0}^{i_0-1}B_{2i_1+1}  \sum_{i_2=i_0}^{m} \left\{ A_{2i_2+1} \prod _{i_3=i_0}^{i_2-1}B_{2i_3+2}\right.\right.\nonumber\\
&&\times   \left.\left.\sum_{i_4=i_2}^{m}\left\{ A_{2i_4+2}\prod _{i_5=i_2}^{i_4-1}B_{2i_5+3} \sum_{i_6=i_4}^{m} \left\{\prod _{i_7=i_4}^{i_6-1}B_{2i_7+4}  \right\}\right\} \right\}\right\} x^{2i_6+3+\lambda }\label{eq:22b}
\end{eqnarray}
(\ref{eq:22a}) and (\ref{eq:22b}) are the sub-power series that has sequences $c_3, c_5, c_7, \cdots, c_{2m+3}$ including three terms of $A_n$'s. If $m\rightarrow \infty $, the sub-power series $y_3^m(x)$ turns to be an infinite series $y_3(x)$.

By repeating this process for all higher terms of $A$'s, we obtain every indicial equations for the sequence $c_{2n+\tau }$ and $y_{\tau }^m(x)$ terms where $\tau  \geq 4$. 

According (\ref{eq:12}), (\ref{eq:15}), (\ref{eq:x18}), (\ref{eq:x21}) and every $c_{2n+\tau }$ where $\tau  \geq 4$, the general expression of coefficients $c_{2n+\tau }$ for a fixed $n,\tau \in \mathbb{N}_{0}$ is taken by
\begin{eqnarray}
\tilde{c}(0,n) \;=\;\;c_{2n} &=& c_0\prod _{i_0=0}^{n-1}B_{2i_0+1} \label{eq:23a}\\
\tilde{c}(1,n) = c_{2n+1} &=& c_0\sum_{i_0=0}^{n} \left\{ A_{2i_0} \prod _{i_1=0}^{i_0-1}B_{2i_1+1} \prod _{i_2=i_0}^{n-1}B_{2i_2+2} \right\} \label{eq:23b}\\
\tilde{c}(\tau ,n) = c_{2n+\tau } &=& c_0\sum_{i_0=0}^{n} \left\{A_{2i_0}\prod _{i_1=0}^{i_0-1} B_{2i_1+1} 
\prod _{k=1}^{\tau -1} \left( \sum_{i_{2k}= i_{2(k-1)}}^{n} A_{2i_{2k}+k}\prod _{i_{2k+1}=i_{2(k-1)}}^{i_{2k}-1}B_{2i_{2k+1}+(k+1)}\right) \right. \nonumber\\
&& \times \left. \prod _{i_{2\tau}=i_{2(\tau -1)}}^{n-1} B_{2i_{2\tau }+(\tau +1)} \right\}\;\; \mathrm{where}\;\tau \geq 2
\hspace{2cm}\label{eq:23c}  
\end{eqnarray}
According (\ref{eq:13b}), (\ref{eq:16b}), (\ref{eq:19b}), (\ref{eq:22b}) and every $y_{\tau }^m(x)$ where $\tau  \geq 4$, the general expression of sub-power series $y_{\tau }^m(x)$ for a fixed $\tau \in \mathbb{N}_{0}$ is given by
\begin{eqnarray}
y_0^m(x) &=& c_0 \sum_{i_0=0}^{m} \left\{ \prod _{i_1=0}^{i_0-1}B_{2i_1+1} \right\} x^{2i_0+\lambda } \label{eq:24a}\\
y_1^m(x) &=& c_0 \sum_{i_0=0}^{m}\left\{ A_{2i_0} \prod _{i_1=0}^{i_0-1}B_{2i_1+1}  \sum_{i_2=i_0}^{m} \left\{ \prod _{i_3=i_0}^{i_2-1}B_{2i_3+2} \right\}\right\} x^{2i_2+1+\lambda } \label{eq:24b}\\
y_{\tau }^m(x) &=& c_0 \sum_{i_0=0}^{m} \left\{A_{2i_0}\prod _{i_1=0}^{i_0-1} B_{2i_1+1} 
\prod _{k=1}^{\tau -1} \left( \sum_{i_{2k}= i_{2(k-1)}}^{m} A_{2i_{2k}+k}\prod _{i_{2k+1}=i_{2(k-1)}}^{i_{2k}-1}B_{2i_{2k+1}+(k+1)}\right) \right. \nonumber\\
&& \times \left. \sum_{i_{2\tau} = i_{2(\tau -1)}}^{m} \left( \prod _{i_{2\tau +1}=i_{2(\tau -1)}}^{i_{2\tau}-1} B_{2i_{2\tau +1}+(\tau +1)} \right) \right\} x^{2i_{2\tau}+\tau+\lambda }\;\;\mathrm{where}\;\tau \geq 2
\label{eq:24c} 
\end{eqnarray}
\begin{definition}
For the first species complete polynomial,  we need a condition which is given by
\begin{equation}
B_{j+1}= c_{j+1}=0\hspace{1cm}\mathrm{where}\;j=0,1,2,\cdots   
 \label{eq:25}
\end{equation}
\end{definition}
(\ref{eq:25}) gives successively $c_{j+2}=c_{j+3}=c_{j+4}=\cdots=0$. And $c_{j+1}=0$ is defined by a polynomial equation of degree $j+1$ for the determination of an accessory parameter in $A_n$ term.
\begin{theorem}
The general expression of a function $y(x)$ for the first species complete polynomial using 3-term recurrence formula and its algebraic equation for the determination of an accessory parameter in $A_n$ term are given by
\begin{enumerate}
\item As $B_1=0$,
\begin{equation}
0 =\bar{c}(1,0) \nonumber 
\end{equation}
\begin{equation}
y(x) = y_{0}^{0}(x) \nonumber 
\end{equation}
\item As $B_{2N+2}=0$ where $N \in \mathbb{N}_{0}$,
\begin{equation}
0  = \sum_{r=0}^{N+1}\bar{c}\left( 2r, N+1-r\right) \nonumber 
\end{equation}
\begin{equation}
y(x)= \sum_{r=0}^{N} y_{2r}^{N-r}(x)+ \sum_{r=0}^{N} y_{2r+1}^{N-r}(x)  \nonumber 
\end{equation}
\item As $B_{2N+3}=0$ where $N \in \mathbb{N}_{0}$,
\begin{equation}
0  = \sum_{r=0}^{N+1}\bar{c}\left( 2r+1, N+1-r\right) \nonumber 
\end{equation}
\begin{equation}
y(x)= \sum_{r=0}^{N+1} y_{2r}^{N+1-r}(x)+ \sum_{r=0}^{N} y_{2r+1}^{N-r}(x)  \nonumber 
\end{equation}
In the above,
\begin{eqnarray}
\bar{c}(0,n)  &=& \prod _{i_0=0}^{n-1}B_{2i_0+1} \nonumber \\
\bar{c}(1,n) &=&  \sum_{i_0=0}^{n} \left\{ A_{2i_0} \prod _{i_1=0}^{i_0-1}B_{2i_1+1} \prod _{i_2=i_0}^{n-1}B_{2i_2+2} \right\} \nonumber \\
\bar{c}(\tau ,n) &=& \sum_{i_0=0}^{n} \left\{A_{2i_0}\prod _{i_1=0}^{i_0-1} B_{2i_1+1} 
\prod _{k=1}^{\tau -1} \left( \sum_{i_{2k}= i_{2(k-1)}}^{n} A_{2i_{2k}+k}\prod _{i_{2k+1}=i_{2(k-1)}}^{i_{2k}-1}B_{2i_{2k+1}+(k+1)}\right) \right. \nonumber\\
&&\times \left. \prod _{i_{2\tau}=i_{2(\tau -1)}}^{n-1} B_{2i_{2\tau }+(\tau +1)} \right\} 
\nonumber   
\end{eqnarray}
and
\begin{eqnarray}
y_0^m(x) &=& c_0 x^{\lambda } \sum_{i_0=0}^{m} \left\{ \prod _{i_1=0}^{i_0-1}B_{2i_1+1} \right\} x^{2i_0 } \nonumber \\
y_1^m(x) &=& c_0 x^{\lambda } \sum_{i_0=0}^{m}\left\{ A_{2i_0} \prod _{i_1=0}^{i_0-1}B_{2i_1+1}  \sum_{i_2=i_0}^{m} \left\{ \prod _{i_3=i_0}^{i_2-1}B_{2i_3+2} \right\}\right\} x^{2i_2+1 } \nonumber \\
y_{\tau }^m(x) &=& c_0 x^{\lambda } \sum_{i_0=0}^{m} \left\{A_{2i_0}\prod _{i_1=0}^{i_0-1} B_{2i_1+1} 
\prod _{k=1}^{\tau -1} \left( \sum_{i_{2k}= i_{2(k-1)}}^{m} A_{2i_{2k}+k}\prod _{i_{2k+1}=i_{2(k-1)}}^{i_{2k}-1}B_{2i_{2k+1}+(k+1)}\right) \right. \nonumber\\
&& \times \left. \sum_{i_{2\tau} = i_{2(\tau -1)}}^{m} \left( \prod _{i_{2\tau +1}=i_{2(\tau -1)}}^{i_{2\tau}-1} B_{2i_{2\tau +1}+(\tau +1)} \right) \right\} x^{2i_{2\tau}+\tau }\hspace{1cm}\mathrm{where}\;\tau \geq 2
\nonumber  
\end{eqnarray}
\end{enumerate}
\end{theorem}
\begin{proof}
For instance, if $B_1= c_1=0$ in (\ref{eq:25}), then gives $c_{2}=c_{3}=c_{4}=\cdots=0$. According to (\ref{eq:9}), its power series is given by
\begin{equation}
y(x)= \sum_{n=0}^{0} c_n x^{n+\lambda } = c_0 x^{\lambda } = y_{0}^{0}(x)\label{eq:26}
\end{equation}
The sub-power series $y_{0}^{0}(x)$ in (\ref{eq:26}) is obtain by putting $m=0$ in (\ref{eq:24a}).
And a polynomial equation of degree 1 for the determination of an accessory parameter in $A_n$ term is taken by
\begin{equation}
0 = c_1 = \tilde{c}(1,0) = c_0 A_0  \label{eq:27}
\end{equation}
A coefficient $\tilde{c}(1,0)$ in (\ref{eq:27}) is obtained by putting $n=0$ in (\ref{eq:23b}).

If $B_2= c_2=0$ in (\ref{eq:25}), then gives $c_{3}=c_{4}=c_{5}=\cdots=0$. According to (\ref{eq:9}), its power series is given by
\begin{equation}
y(x)= \sum_{n=0}^{1} c_n x^{n+\lambda } = (c_0 +c_1 x) x^{\lambda } \label{eq:28}
\end{equation}
First observe sequences $c_0$ and $c_1$ in (\ref{eq:7}) and (\ref{eq:8}). 
The sub-power series that has a sequence $c_0$ is given by  putting $m=0$ in (\ref{eq:24a}). The sub-power series that has a sequence $c_1$ is given by  putting $m=0$ in (\ref{eq:24b}). Take the new (\ref{eq:24a}) and (\ref{eq:24b}) into (\ref{eq:28}).
\begin{equation}
y(x)= y_{0}^{0}(x)+ y_{1}^{0}(x) = c_0x^{\lambda } \left( 1+A_0 x\right)\label{eq:29}
\end{equation}
A sequence $c_2$ consists of zero and two terms of $A_n's$ in (\ref{eq:7}) and (\ref{eq:8}). 
Putting $n=1$ in (\ref{eq:23a}), a coefficient $c_2$ for zero term of $A_n's$ is denoted by $\tilde{c}(0,1)$. 
Taking $\tau=2$ and $n=0$  in (\ref{eq:23c}), a coefficient $c_2$ for two terms of $A_n's$ is denoted by $\tilde{c}(2,0)$.
Since the sum of $\tilde{c}(0,1)$ and $\tilde{c}(2,0)$ is equivalent to zero, we obtain a polynomial equation of degree 2 for the determination of an accessory parameter in $A_n$ term which is given by
\begin{equation}
0 = c_2 = \tilde{c}(0,1) + \tilde{c}(2,0)  = c_0 \left( A_{0,1} +B_1 \right)  \label{eq:30}
\end{equation} 

If $B_3= c_3=0$ in (\ref{eq:25}), then gives $c_{4}=c_{5}=c_{6}=\cdots=0$. According to (\ref{eq:9}), its power series is given by
\begin{equation}
y(x)= \sum_{n=0}^{2} c_n x^{n+\lambda } = (c_0 +c_1 x+c_2 x^2) x^{\lambda } \label{eq:31}
\end{equation}
Observe sequences $c_0$-$c_2$ in (\ref{eq:7}) and (\ref{eq:8}). 
The sub-power series, having sequences $c_0$ and $c_2$ including zero term of $A_n's$, is given by  putting $m=1$ in (\ref{eq:24a}) denoted by $y_0^1(x)$.
The sub-power series, having a sequence $c_1$ including one term of $A_n's$, is given by  putting $m=0$ in (\ref{eq:24b}) denoted by $y_1^0(x)$.
The sub-power series, having a sequence $c_2$ including two terms of $A_n's$, is given by  putting $\tau=2$ and $m=0$ in (\ref{eq:24c}) denoted by $y_2^0(x)$.
Taking $y_0^1(x)$, $y_1^0(x)$ and $y_2^0(x)$ into (\ref{eq:31}),
\begin{equation}
y(x)= y_{0}^{1}(x) + y_{2}^{0}(x)+ y_{1}^{0}(x)= c_0x^{\lambda } \left( 1+A_0 x + \left( A_{0,1} +B_1\right) x^2\right) \label{eq:32}
\end{equation} 
A sequence $c_3$ consists of one and three terms of $A_n's$ in (\ref{eq:7}) and (\ref{eq:8}).
Putting $n=1$ in (\ref{eq:23b}), a coefficient $c_3$ for one term of $A_n's$ is denoted by $\tilde{c}(1,1)$.
Taking $\tau=3$ and $n=0$  in (\ref{eq:23c}), a coefficient $c_3$ for three terms of $A_n's$ is denoted by $\tilde{c}(3,0)$.
Since the sum of $\tilde{c}(1,1)$ and $\tilde{c}(3,0)$ is equivalent to zero, we obtain a polynomial equation of degree 3 for the determination of an accessory parameter in $A_n$ term which is given by
\begin{equation}
0 = c_3 = \tilde{c}(1,1) + \tilde{c}(3,0)  = c_0 \left( A_{0,1,2} +B_1 A_2 + B_2 A_0 \right)  \label{eq:33}
\end{equation} 

If $B_4= c_4=0$ in (\ref{eq:25}), then gives $c_{5}=c_{6}=c_{7}=\cdots=0$. According to (\ref{eq:9}), its power series is given by
\begin{equation}
y(x)= \sum_{n=0}^{3} c_n x^{n+\lambda } = (c_0 +c_1 x+c_2 x^2+c_3 x^3) x^{\lambda } \label{eq:34}
\end{equation}
Observe sequences $c_0$-$c_3$ in (\ref{eq:7}) and (\ref{eq:8}). 
The sub-power series, having sequences $c_0$ and $c_2$ including zero term of $A_n's$, is given by  putting $m=1$ in (\ref{eq:24a}) denoted by $y_0^1(x)$.
The sub-power series, having sequences $c_1$ and $c_3$ including one term of $A_n's$, is given by  putting $m=1$ in (\ref{eq:24b}) denoted by $y_1^1(x)$.
The sub-power series, having a sequence $c_2$ including two terms of $A_n's$, is given by  putting $\tau=2$ and $m=0$ in (\ref{eq:24c}) denoted by $y_2^0(x)$.
The sub-power series, having a sequence $c_3$ including three terms of $A_n's$, is given by  putting $\tau=3$ and $m=0$ in (\ref{eq:24c}) denoted by $y_3^0(x)$.
Taking $y_0^1(x)$, $y_1^1(x)$, $y_2^0(x)$ and $y_3^0(x)$ into (\ref{eq:34}),
\begin{equation}
y(x)= y_{0}^{1}(x)+y_{2}^{0}(x)+ y_{1}^{1}(x) + y_{3}^{0}(x) \label{eq:35}
\end{equation}
A sequence $c_4$ consists of zero, two and four terms of $A_n's$ in (\ref{eq:7}) and (\ref{eq:8}).
Putting $n=2$ in (\ref{eq:23a}), a coefficient $c_4$ for zero term of $A_n's$ is denoted by $\tilde{c}(0,2)$.
Taking $\tau=2$ and $n=1$  in (\ref{eq:23c}), a coefficient $c_4$ for two terms of $A_n's$ is denoted by $\tilde{c}(2 ,1)$.
Taking $\tau=4$ and $n=0$  in (\ref{eq:23c}), a coefficient $c_4$ for four terms of $A_n's$ is denoted by $\tilde{c}(4 ,0)$.
Since the sum of $\tilde{c}(0,2)$, $\tilde{c}(2,1)$ and $\tilde{c}(4,0)$ is equivalent to zero, we obtain a polynomial equation of degree 4 for the determination of an accessory parameter in $A_n$ term which is given by
\begin{equation}
0 = c_4 = \tilde{c}(0,2) + \tilde{c}(2 ,1) + \tilde{c}(4 ,0) \label{eq:36}
\end{equation}

If $B_5= c_5=0$ in (\ref{eq:25}), then gives $c_{6}=c_{7}=c_{8}=\cdots=0$. According to (\ref{eq:9}), its power series is given by
\begin{equation}
y(x)= \sum_{n=0}^{4} c_n x^{n+\lambda } = (c_0 +c_1 x+c_2 x^2+c_3 x^3+c_4 x^4) x^{\lambda } \label{eq:37}
\end{equation}
Observe sequences $c_0$-$c_4$ in (\ref{eq:7}) and (\ref{eq:8}). 
The sub-power series, having sequences $c_0$, $c_2$ and $c_4$ including zero term of $A_n's$, is given by  putting $m=2$ in (\ref{eq:24a}) denoted by $y_0^2(x)$.
The sub-power series, having sequences $c_1$ and $c_3$ including one term of $A_n's$, is given by  putting $m=1$ in (\ref{eq:24b}) denoted by $y_1^1(x)$.
The sub-power series, having sequences $c_2$ and $c_4$ including two terms of $A_n's$, is given by  putting $\tau=2$ and $m=1$ in (\ref{eq:24c}) denoted by $y_2^1(x)$.
The sub-power series, having a sequence $c_3$ including three terms of $A_n's$, is given by  putting $\tau=3$ and $m=0$ in (\ref{eq:24c}) denoted by $y_3^0(x)$.
The sub-power series, having a sequence $c_4$ including four terms of $A_n's$, is given by  putting $\tau=4$ and $m=0$ in (\ref{eq:24c}) denoted by $y_4^0(x)$.
Taking $y_0^2(x)$, $y_1^1(x)$, $y_2^1(x)$, $y_3^0(x)$ and $y_4^0(x)$ into (\ref{eq:37}),
\begin{equation}
y(x)= y_0^2(x)+y_2^1(x)+y_4^0(x)+y_1^1(x)+y_3^0(x) \label{eq:38}
\end{equation}
A sequence $c_5$ consists of one, three and five terms of $A_n's$ in (\ref{eq:7}) and (\ref{eq:8}).
Putting $n=2$ in (\ref{eq:23b}), a coefficient $c_5$ for one term of $A_n's$ is denoted by $\tilde{c}(1,2)$.
Taking $\tau=3$ and $n=1$  in (\ref{eq:23c}), a coefficient $c_5$ for three terms of $A_n's$ is denoted by $\tilde{c}(3,1)$.
Taking $\tau=5$ and $n=0$  in (\ref{eq:23c}), a coefficient $c_5$ for five terms of $A_n's$ is denoted by $\tilde{c}(5,0)$.
Since the sum of $\tilde{c}(1,2)$, $\tilde{c}(3,1)$ and $\tilde{c}(5,0)$ is equivalent to zero, we obtain a polynomial equation of degree 5 for the determination of an accessory parameter in $A_n$ term which is given by
\begin{equation}
0 = c_5 = \tilde{c}(1,2) + \tilde{c}(3,1) + \tilde{c}(5,0) \label{eq:39}
\end{equation}

If $B_6= c_6=0$ in (\ref{eq:25}), then gives $c_{7}=c_{8}=c_{9}=\cdots=0$. According to (\ref{eq:9}), its power series is given by
\begin{equation}
y(x)= \sum_{n=0}^{5} c_n x^{n+\lambda } = (c_0 +c_1 x+c_2 x^2+c_3 x^3+c_4 x^4+c_5 x^5) x^{\lambda } \label{eq:40}
\end{equation}
Observe sequences $c_0$-$c_5$ in (\ref{eq:7}) and (\ref{eq:8}). 
The sub-power series, having sequences $c_0$, $c_2$ and $c_4$ including zero term of $A_n's$, is given by  putting $m=2$ in (\ref{eq:24a}) denoted by $y_0^2(x)$.
The sub-power series, having sequences $c_1$, $c_3$ and $c_5$ including one term of $A_n's$, is given by  putting $m=2$ in (\ref{eq:24b}) denoted by $y_1^2(x)$.
The sub-power series, having sequences $c_2$ and $c_4$ including two terms of $A_n's$, is given by  putting $\tau=2$ and $m=1$ in (\ref{eq:24c}) denoted by $y_2^1(x)$.
The sub-power series, having sequences $c_3$ and $c_5$ including three terms of $A_n's$, is given by  putting $\tau=3$ and $m=1$ in (\ref{eq:24c}) denoted by $y_3^1(x)$.
The sub-power series, having a sequence $c_4$ including four terms of $A_n's$, is given by  putting $\tau=4$ and $m=0$ in (\ref{eq:24c}) denoted by $y_4^0(x)$.
The sub-power series, having a sequence $c_5$ including five terms of $A_n's$, is given by  putting $\tau=5$ and $m=0$ in (\ref{eq:24c}) denoted by $y_5^0(x)$.
Taking $y_0^2(x)$, $y_1^2(x)$, $y_2^1(x)$, $y_3^1(x)$, $y_4^0(x)$ and $y_5^0(x)$ into (\ref{eq:40}),
\begin{equation}
y(x)= y_0^2(x)+y_2^1(x)+y_4^0(x)+y_1^2(x)+y_3^1(x)+y_5^0(x)= \sum_{r=0}^{2}y_{2r}^{2-r}(x) +\sum_{r=0}^{2}y_{2r+1}^{2-r}(x)\label{eq:41}
\end{equation}
A sequence $c_6$ consists of zero, two, four and six terms of $A_n's$ in (\ref{eq:7}) and (\ref{eq:8}).
Putting $n=3$ in (\ref{eq:23a}), a coefficient $c_6$ for zero term of $A_n's$ is denoted by $\tilde{c}(0,3)$.
Taking $\tau=2$ and $n=2$  in (\ref{eq:23c}), a coefficient $c_6$ for two terms of $A_n's$ is denoted by $\tilde{c}(2,2)$.
Taking $\tau=4$ and $n=1$  in (\ref{eq:23c}), a coefficient $c_6$ for four terms of $A_n's$ is denoted by $\tilde{c}(4,1)$.
Taking $\tau=6$ and $n=0$  in (\ref{eq:23c}), a coefficient $c_6$ for six terms of $A_n's$ is denoted by $\tilde{c}(6,0)$.
Since the sum of $\tilde{c}(0,3)$, $\tilde{c}(2,2)$, $\tilde{c}(4,1)$ and $\tilde{c}(6,0)$ is equivalent to zero, we obtain a polynomial equation of degree 6 for the determination of an accessory parameter in $A_n$ term which is given by
\begin{equation}
0 = c_6 = \tilde{c}(0,3) + \tilde{c}(2,2) + \tilde{c}(4,1) + \tilde{c}(6,0) = \sum_{r=0}^{3}\tilde{c}(2r,3-r)\label{eq:42}
\end{equation}

If $B_7= c_7=0$ in (\ref{eq:25}), then gives $c_{8}=c_{9}=c_{10}=\cdots=0$. According to (\ref{eq:9}), its power series is given by
\begin{eqnarray}
y(x)&=& \sum_{n=0}^{6} c_n x^{n+\lambda } = y_0^3(x)+y_2^2(x)+y_4^1(x)+y_6^0(x)+y_1^2(x)+y_3^1(x)+y_5^0(x)\nonumber\\
&=& \sum_{r=0}^{3}y_{2r}^{3-r}(x) +\sum_{r=0}^{2}y_{2r+1}^{2-r}(x)\label{eq:43}
\end{eqnarray}
A polynomial equation of degree 7 for the determination of an accessory parameter in $A_n$ term is given by
\begin{equation}
0 = c_7 = \tilde{c}(1,3) + \tilde{c}(3,2) + \tilde{c}(5,1) + \tilde{c}(7,0) = \sum_{r=0}^{3}\tilde{c}(2r+1,3-r)\label{eq:44}
\end{equation}
By repeating this process for $B_{j+1}= c_{j+1}=0$ where $j\geq 7$, we obtain the first species complete polynomial of degree $j$ and an polynomial equation of the $j+1^{th}$ order in an accessory parameter.
\qed
\end{proof}
\subsection{The second species complete polynomial using 3TRF}
There is another way to construct a polynomial which makes $A_n$ and $B_n$ terms terminated at the same time. 

For instance, according to Karl Heun(1889) \cite{Heun1889,Ronv1995}, Heun's equation is a second-order linear ODE of the form 
\begin{equation}
\frac{d^2{y}}{d{x}^2} + \left(\frac{\gamma }{x} +\frac{\delta }{x-1} + \frac{\epsilon }{x-a}\right) \frac{d{y}}{d{x}} +  \frac{\alpha \beta x-q}{x(x-1)(x-a)} y = 0 \label{eq:50}
\end{equation}
With the condition $\epsilon = \alpha +\beta -\gamma -\delta +1$. The parameters play different roles: $a \ne 0 $ is the singularity parameter, $\alpha $, $\beta $, $\gamma $, $\delta $, $\epsilon $ are exponent parameters, $q$ is the accessory parameter. Also, $\alpha $ and $\beta $ are identical to each other. It has four regular singular points which are 0, 1, $a$ and $\infty $ with exponents $\{ 0, 1-\gamma \}$, $\{ 0, 1-\delta \}$, $\{ 0, 1-\epsilon \}$ and $\{ \alpha, \beta \}$.

we assume the solution takes the form
\begin{equation}
y(x) = \sum_{n=0}^{\infty } c_n x^{n+\lambda } \label{eq:51}
\end{equation}
Substituting (\ref{eq:51}) into (\ref{eq:50}) gives for the coefficients $c_n$ the recurrence relations
\begin{equation}
c_{n+1}=A_n \;c_n +B_n \;c_{n-1} \hspace{1cm};n\geq 1
\label{eq:52}
\end{equation}
where,
\begin{subequations}
\begin{equation}
A_n =  \frac{(n+\lambda )(n+\alpha +\beta -\delta +\lambda +a(n+\delta +\gamma -1+\lambda ))+q}{a(n+1+\lambda )(n+\gamma +\lambda )} 
\label{eq:53a}
\end{equation}
\begin{equation}
B_n = - \frac{(n-1+\lambda +\alpha )(n-1+\lambda +\beta )}{a(n+1+\lambda )(n+\gamma +\lambda )}
\label{eq:53b}
\end{equation}
\begin{equation}
c_1= A_0 \;c_0 \label{eq:53c}
\end{equation}
\end{subequations}
We have two indicial roots which are $\lambda =0$ and $1-\gamma$.

For a fixed $j \in \mathbb{N}_{0}$, let us suppose that 
\begin{equation}
\begin{cases} \alpha = -j-\lambda  \cr
\beta  = -j+1-\lambda \cr
q  = -(j+\lambda )(-\delta +1-j-\lambda +a(\gamma +\delta -1+j+\lambda )) 
\end{cases}\label{eq:54}
\end{equation} 
Plug (\ref{eq:54}) into (\ref{eq:53a})-(\ref{eq:53c}).
\begin{equation}
\begin{cases} A_n =  \frac{1+a}{a}\frac{(n-j)\left\{ n+\frac{1}{1+a}\left[ -\delta +1-j+a(\gamma +\delta -1+j+2\lambda )\right]\right\}}{(n+1+\lambda )(n+\gamma +\lambda )}  \cr
B_n = -\frac{1}{a} \frac{(n-j-1 )(n-j)}{(n+1+\lambda )(n+\gamma +\lambda )} \cr
c_1= A_0 \;c_0
\end{cases}\label{eq:55}
\end{equation}
Moreover, $B_j=B_{j+1}=A_j=0$ for a fixed $j \in \mathbb{N}_{0}$ in (\ref{eq:55}), then (\ref{eq:52}) gives successively 
$ c_{j+1}=c_{j+2}=c_{j+3}=\cdots=0 $. As I mentioned before, I categorize this type polynomial as the second species complete polynomial. With the quantum mechanical point of view, we are allowed to use this type of polynomials as long as there are two different eigenenergies in $B_n$ term and an eigenenergy in $A_n$ term in a 3-term recursive relation of a second-order ODE.
For the first species complete polynomial, there are $j+1$ values of a spectral parameter (eigenvalue or eigenenergy) of a numerator in $A_n$ term. In contrast, the second species complete polynomial only has one value of a parameter of a numerator in $A_n$ term for each $j$.
\begin{definition}
For the second species complete polynomial, we need a condition which is defined by
\begin{equation}
B_{j}=B_{j+1}= A_{j}=0\hspace{1cm}\mathrm{where}\;j \in \mathbb{N}_{0}    
 \label{eq:56}
\end{equation}
\end{definition}
\begin{theorem}
The general expression of a function $y(x)$ for the second species complete polynomial using 3-term recurrence formula is given by
\begin{enumerate}
\item As $B_1=A_0=0$,
\begin{equation}
y(x) = y_{0}^{0}(x) \nonumber 
\end{equation}
\item As $B_{2N+1}=B_{2N+2}=A_{2N+1}=0$  where $N \in \mathbb{N}_{0}$,
\begin{equation}
y(x)= \sum_{r=0}^{N} y_{2r}^{N-r}(x)+ \sum_{r=0}^{N} y_{2r+1}^{N-r}(x)  \nonumber 
\end{equation}
\item As $B_{2N+2}=B_{2N+3}=A_{2N+2}=0$  where $N \in \mathbb{N}_{0}$,
\begin{equation}
y(x)= \sum_{r=0}^{N+1} y_{2r}^{N+1-r}(x)+ \sum_{r=0}^{N} y_{2r+1}^{N-r}(x)  \nonumber 
\end{equation}
In the above,
\begin{eqnarray}
y_0^m(x) &=& c_0 x^{\lambda } \sum_{i_0=0}^{m} \left\{ \prod _{i_1=0}^{i_0-1}B_{2i_1+1} \right\} x^{2i_0 } \nonumber \\
y_1^m(x) &=& c_0 x^{\lambda } \sum_{i_0=0}^{m}\left\{ A_{2i_0} \prod _{i_1=0}^{i_0-1}B_{2i_1+1}  \sum_{i_2=i_0}^{m} \left\{ \prod _{i_3=i_0}^{i_2-1}B_{2i_3+2} \right\}\right\} x^{2i_2+1 } \nonumber \\
y_{\tau }^m(x) &=& c_0 x^{\lambda } \sum_{i_0=0}^{m} \left\{A_{2i_0}\prod _{i_1=0}^{i_0-1} B_{2i_1+1} 
\prod _{k=1}^{\tau -1} \left( \sum_{i_{2k}= i_{2(k-1)}}^{m} A_{2i_{2k}+k}\prod _{i_{2k+1}=i_{2(k-1)}}^{i_{2k}-1}B_{2i_{2k+1}+(k+1)}\right) \right. \nonumber\\
&& \times \left. \sum_{i_{2\tau} = i_{2(\tau -1)}}^{m} \left( \prod _{i_{2\tau +1}=i_{2(\tau -1)}}^{i_{2\tau}-1} B_{2i_{2\tau +1}+(\tau +1)} \right) \right\} x^{2i_{2\tau}+\tau }\hspace{1cm}\mathrm{where}\;\tau \geq 2
\nonumber  
\end{eqnarray}
\end{enumerate}
\end{theorem}
\begin{proof}
First if $j=0$ in (\ref{eq:56}),  $B_{1}=A_0=0$ and then (\ref{eq:52}) gives successively 
$ c_{1}=c_{2}=c_{3}=\cdots=0 $. According to (\ref{eq:9}), its power series is given by
\begin{equation}
y(x)= \sum_{n=0}^{0} c_n x^{n+\lambda } = c_0 x^{\lambda } = y_{0}^{0}(x)\label{eq:57}
\end{equation}
The sub-power series $y_{0}^{0}(x)$ in (\ref{eq:57}) is obtain by putting $m=0$ in (\ref{eq:24a}).
 
If $j=1$ in (\ref{eq:56}), $B_{1}=B_{2}=A_1=0$ and then (\ref{eq:52}) gives successively 
$ c_{2}=c_{3}=c_{4}=\cdots=0 $. According to (\ref{eq:9}), its power series is given by
\begin{equation}
y(x)= \sum_{n=0}^{1} c_n x^{n+\lambda } = (c_0 +c_1 x) x^{\lambda } \label{eq:58}
\end{equation}
Observe sequences $c_0$ and $c_1$ in (\ref{eq:7}) and (\ref{eq:8}). 
The sub-power series that has a sequence $c_0$ is given by  putting $m=0$ in (\ref{eq:24a}). The sub-power series that has a sequence $c_1$ is given by  putting $m=0$ in (\ref{eq:24b}). Take the new (\ref{eq:24a}) and (\ref{eq:24b}) into (\ref{eq:58}).
\begin{equation}
y(x)= y_{0}^{0}(x)+ y_{1}^{0}(x) = c_0x^{\lambda } \left( 1+A_0 x\right)\label{eq:59}
\end{equation}

If $j=2$ in (\ref{eq:56}), $B_{2}=B_{3}=A_2=0$ and then (\ref{eq:52}) gives successively 
$ c_{3}=c_{4}=c_{5}=\cdots=0 $. According to (\ref{eq:9}), its power series is given by
\begin{equation}
y(x)= \sum_{n=0}^{2} c_n x^{n+\lambda } = (c_0 +c_1 x+c_2 x^2) x^{\lambda } \label{eq:60}
\end{equation}
Observe sequences $c_0$-$c_2$ in (\ref{eq:7}) and (\ref{eq:8}). 
The sub-power series, having sequences $c_0$ and $c_2$ including zero term of $A_n's$, is given by  putting $m=1$ in (\ref{eq:24a}) denoted by $y_0^1(x)$.
The sub-power series, having a sequence $c_1$ including one term of $A_n's$, is given by  putting $m=0$ in (\ref{eq:24b}) denoted by $y_1^0(x)$.
The sub-power series, having a sequence $c_2$ including two terms of $A_n's$, is given by  putting $\tau=2$ and $m=0$ in (\ref{eq:24c}) denoted by $y_2^0(x)$.
Taking $y_0^1(x)$, $y_1^0(x)$ and $y_2^0(x)$ into (\ref{eq:60}),
\begin{equation}
y(x)= y_{0}^{1}(x) + y_{2}^{0}(x)+ y_{1}^{0}(x)= c_0x^{\lambda } \left( 1+A_0 x + \left( A_{0,1} +B_1\right) x^2\right) \label{eq:61}
\end{equation}

If $j=3$ in (\ref{eq:56}), $B_{3}=B_{4}=A_3=0$ and then (\ref{eq:52}) gives successively 
$ c_{4}=c_{5}=c_{6}=\cdots=0 $. According to (\ref{eq:9}), its power series is given by
\begin{equation}
y(x)= \sum_{n=0}^{3} c_n x^{n+\lambda } = (c_0 +c_1 x+c_2 x^2+c_3 x^3) x^{\lambda } \label{eq:62}
\end{equation}
Observe sequences $c_0$-$c_3$ in (\ref{eq:7}) and (\ref{eq:8}). 
The sub-power series, having sequences $c_0$ and $c_2$ including zero term of $A_n's$, is given by  putting $m=1$ in (\ref{eq:24a}) denoted by $y_0^1(x)$.
The sub-power series, having sequences $c_1$ and $c_3$ including one term of $A_n's$, is given by  putting $m=1$ in (\ref{eq:24b}) denoted by $y_1^1(x)$.
The sub-power series, having a sequence $c_2$ including two terms of $A_n's$, is given by  putting $\tau=2$ and $m=0$ in (\ref{eq:24c}) denoted by $y_2^0(x)$.
The sub-power series, having a sequence $c_3$ including three terms of $A_n's$, is given by  putting $\tau=3$ and $m=0$ in (\ref{eq:24c}) denoted by $y_3^0(x)$.
Taking $y_0^1(x)$, $y_1^1(x)$, $y_2^0(x)$ and $y_3^0(x)$ into (\ref{eq:62}),
\begin{equation}
y(x)= y_{0}^{1}(x)+y_{2}^{0}(x)+ y_{1}^{1}(x) + y_{3}^{0}(x) \label{eq:63}
\end{equation}

If $j=4$ in (\ref{eq:56}), $B_{4}=B_{5}=A_4=0$ and then (\ref{eq:52}) gives successively 
$ c_{5}=c_{6}=c_{7}=\cdots=0 $. According to (\ref{eq:9}), its power series is given by
\begin{equation}
y(x)= \sum_{n=0}^{4} c_n x^{n+\lambda } = (c_0 +c_1 x+c_2 x^2+c_3 x^3+c_4 x^4) x^{\lambda } \label{eq:64}
\end{equation}
Observe sequences $c_0$-$c_4$ in (\ref{eq:7}) and (\ref{eq:8}). 
The sub-power series, having sequences $c_0$, $c_2$ and $c_4$ including zero term of $A_n's$, is given by  putting $m=2$ in (\ref{eq:24a}) denoted by $y_0^2(x)$.
The sub-power series, having sequences $c_1$ and $c_3$ including one term of $A_n's$, is given by  putting $m=1$ in (\ref{eq:24b}) denoted by $y_1^1(x)$.
The sub-power series, having sequences $c_2$ and $c_4$ including two terms of $A_n's$, is given by  putting $\tau=2$ and $m=1$ in (\ref{eq:24c}) denoted by $y_2^1(x)$.
The sub-power series, having a sequence $c_3$ including three terms of $A_n's$, is given by  putting $\tau=3$ and $m=0$ in (\ref{eq:24c}) denoted by $y_3^0(x)$.
The sub-power series, having a sequence $c_4$ including four terms of $A_n's$, is given by  putting $\tau=4$ and $m=0$ in (\ref{eq:24c}) denoted by $y_4^0(x)$.
Taking $y_0^2(x)$, $y_1^1(x)$, $y_2^1(x)$, $y_3^0(x)$ and $y_4^0(x)$ into (\ref{eq:64}),
\begin{equation}
y(x)= y_0^2(x)+y_2^1(x)+y_4^0(x)+y_1^1(x)+y_3^0(x) \label{eq:65}
\end{equation}

If $j=5$ in (\ref{eq:56}), $B_{5}=B_{6}=A_5=0$ and then (\ref{eq:52}) gives successively 
$ c_{6}=c_{7}=c_{8}=\cdots=0 $. According to (\ref{eq:9}), its power series is given by
\begin{eqnarray}
y(x)&=& \sum_{n=0}^{5} c_n x^{n+\lambda } = y_0^2(x)+y_2^1(x)+y_4^0(x)+y_1^2(x)+y_3^1(x)+y_5^0(x)\nonumber\\
&=& \sum_{r=0}^{2}y_{2r}^{2-r}(x) +\sum_{r=0}^{2}y_{2r+1}^{2-r}(x) \label{eq:66}
\end{eqnarray}

If $j=6$ in (\ref{eq:56}), $B_{6}=B_{7}=A_6=0$ and then (\ref{eq:52}) gives successively 
$ c_{7}=c_{8}=c_{9}=\cdots=0 $. According to (\ref{eq:9}), its power series is given by 
\begin{eqnarray}
y(x)&=& \sum_{n=0}^{6} c_n x^{n+\lambda } = y_0^3(x)+y_2^2(x)+y_4^1(x)+y_6^0(x)+y_1^2(x)+y_3^1(x)+y_5^0(x)\nonumber\\
&=& \sum_{r=0}^{3}y_{2r}^{3-r}(x) +\sum_{r=0}^{2}y_{2r+1}^{2-r}(x)\label{eq:67}
\end{eqnarray}
By repeating this process for $B_{j}=B_{j+1}= A_{j}=0$ where $j\geq 7$, we obtain the second species complete polynomial of degree $j$.
\qed
\end{proof}
\section{Summary}
The Method of Frobenius tells us that there are two types of power series in the 2-term recurrence relation of a linear ODE which are an infinite series and a polynomial. 
In contrast, the power series in the 3-term recursive relation have an infinite series and three types of polynomials: (1) a polynomial which makes $B_n$ term terminated, (2) a polynomial which makes $A_n$ term terminated, (3) a polynomial which makes $A_n$ and $B_n$ terms terminated.

I designate a type 3 polynomial as a complete polynomial in the three term recursive relation of a linear ODE.
With my definition, complete polynomials has two different types which are the first and second species complete polynomials. The former is applicable if there are only one eigenvalue in $B_n$ term and an eigenvalue in $A_n$ term. And the latter is applicable if there are two eigenvalues in $B_n$ term and an eigenvalue in $A_n$. 

In this chapter, by substituting a power series with unknown coefficients into a linear ODE, I construct the general expression for the first and second complete polynomials by allowing $A_n$ as the leading term in each sub-power series of the general power series: I observe the term of sequence $c_n$ which includes zero term of $A_n's$, one term of $A_n's$, two terms of $A_n's$, three terms of $A_n's$, etc. 
The mathematical structures of the first and second complete polynomials are equivalent to each other as you see theorem 1.3.2 and 1.3.4. The only difference between two type polynomials is that the first complete polynomial has multi-valued roots of an eigenvalue (spectral parameter) in $A_n$ term, but the second complete polynomial has only one valued root of its eigenvalue in $A_n$ term.

\addcontentsline{toc}{section}{Bibliography}
\bibliographystyle{model1a-num-names}
\bibliography{<your-bib-database>}

%
\chapter{Complete polynomial using the reversible three-term recurrence formula and its applications}
\chaptermark{Complete polynomial using R3TRF}  
On chapter 1, by putting a power series with unknown coefficients into a linear ordinary differential equation (ODE), I generalize the three term recurrence relation for a polynomial which makes $A_n$ and $B_n$ terms terminated. 

In this chapter, by applying the Frobenius method, I generalize the three term recurrence relation in a backward for a polynomial which makes $A_n$ and $B_n$ terms terminated. 

The general expression for the former is designated as complete polynomials using 3-term recurrence formula (3TRF). And the latter is denominated as complete polynomials using reversible 3-term recurrence formula (R3TRF). 

\section{Introduction}
In general, many linear (ordinary) differential equations can not be solved explicitly in terms of well-known simpler functions. 
Instead, we try to looking for power series solutions by putting a power series with unknown coefficients into ODEs. The recursive relation of coefficients starts to appear in a power series; there can be between two term and infinity term in the recurrence relation which includes all parameters in ODEs we require. 

Power series solutions for the 2-term recursive relation in ODEs have been obtained with simple familiar functions such as hypergeometric, Kummer, Legendre, Bessel functions, etc. The definite or contour integral representations for those power series expansions are constructed analytically. In contrast, there are no such analytic solutions in the form of power series for the 3-term recurrence relation because of its complex mathematical computations.\cite{1Hortacsu:2011rr,1Whit1914,1Erde1955,1Hobs1931,1Whit1952} 
For instance, some of linear ODEs which contain three different coefficients in the recursive relation are Heun, Confluent Heun, Double Confluent Heun, Biconfluent Heun, Lam\'{e} and Mathieu equations.\footnote{Lam\'{e} Wave and Triconfluent Heun equations are composed of the 4-term recurrence relation in a form of its power series. The power series forms of these mathematical structures will be discussed in the future series.} 

In general, by using the method of Frobenius, the 3-term recursive relation in any linear ODEs is given by
\begin{equation}
c_{n+1}=A_n \;c_n +B_n \;c_{n-1} \hspace{1cm};n\geq 1
\nonumber
\end{equation}
where
\begin{equation}
c_1= A_0 \;c_0
\nonumber
\end{equation}
By putting subscript $n=0,1,2,\cdots$ in succession in the recurrence relation, we can compute the remaining coefficients in terms of $c_0$ recursively.
When a coefficient $c_0$ is factored out of parentheses for $c_n$ in succession, the sequence $c_n$ inside parentheses consists of combinations $A_n$ and $B_n$ terms. Of course, we are able to construct a general  summation formula for $c_n's$ with $c_0$.

In the first and second series\cite{1Choun2012,1Choun2013}, I apply 3TRF and R3TRF to the power series in closed forms of 5 different ODEs (Heun, Confluent Heun, Biconfluent Heun, Mathieu and Lam\'{e} equations) for an infinite series, polynomials of type 1 and 2.
The general power series expansions of the above 5 ODEs for infinite series using 3TRF and R3TRF are equivalent to each other. In the analytic solutions using 3TRF, $A_n$ is the leading term of each sub-power series in power series forms of 5 different functions in the above ODEs. In the analytic solutions using R3TRF, $B_n$ is the leading term of each sub-power series in power series expansions of them.  
 
In the mathematical definition, various authors defines a polynomial of type 3 as the spectral polynomial. With my definition, I categorize a type 3 polynomial as a complete polynomial in the 3-term recurrence relation of a ODE.  
And the complete polynomial has two different types which are the first and second species complete polynomials.
For the first species complete polynomial, we need a condition which is $B_{j+1}= c_{j+1}=0$ where $j=0,1,2,\cdots$. This condition gives successively $c_{j+2}=c_{j+3}=c_{j+4}=\cdots=0$. There are $j+1$ values of an eigenvalue (spectral parameter) in $A_n$ term and one value of an eigenvalue in $B_n$ term for each $j$.
For the second species complete polynomial, we need a condition which is $B_{j}=B_{j+1}= A_{j}=0$. It gives successively $c_{j+1}=c_{j+2}=c_{j+3}=\cdots=0$. There are two eigenvalues in $B_n$ term and only one eigenvalue in $A_n$ term for each $j$.
 
On chapter 1, by putting a power series with unknown coefficients into a linear ODE, I generalize the 3-term recurrence relation for complete polynomials of the first and second species including all higher terms of $A_n's$ using 3TRF: I observe the term parenthesis of sequences $c_0, c_2, c_4,\cdots, c_{2m}$ where $m\in \mathbb{N}_{0}$ which includes zero term of $A_n's$, sequences $c_1, c_3, c_5,\cdots, c_{2m+1}$ which includes one term of $A_n's$, sequences $c_2, c_4, c_6,\cdots, c_{2m+2}$ which includes two terms of $A_n's$, sequences $c_3, c_5, c_7,\cdots, c_{2m+3}$ which includes three terms of $A_n's$, etc.  
For the first species complete polynomial, I generalize the algebraic equation for the determination of the accessory parameter in the form of partial sums of the sequences $\{A_n\}$ and $\{B_n\}$ using 3TRF, by letting $A_n$ as a leading term in its polynomial equation, for computational practice.

In this chapter, I construct the general summation formulas of the 3-term recurrence relation in a linear ODE for two different complete polynomials including all higher terms of $B_n's$ using R3TRF: I observe the term parenthesis of sequences $c_0, c_1, c_2,\cdots, c_{m}$ which includes zero term of $B_n's$, sequences $c_2, c_3, c_4,\cdots, c_{m+2}$ which includes one term of $B_n's$, sequences $c_4, c_5, c_6,\cdots, c_{m+4}$ which includes two terms of $B_n's$, sequences $c_6, c_7, c_8,\cdots, c_{m+6}$ which includes three terms of $B_n's$, etc.   
For the first species complete polynomial, I generalize the polynomial equation for the determination of the accessory parameter in the form of partial sums of the sequences $\{A_n\}$ and $\{B_n\}$ using R3TRF, by letting $B_n$ as a leading term in its algebraic equation.

\section{Complete polynomial using R3TRF}
As I mentioned before, a complete polynomial (a polynomial of type 3) has two different types:

(1) If there is only one eigenvalue in $B_n$ term of the 3-term recurrence relation, there are multi-valued roots of an eigenvalue (spectral parameter) in $A_n$ term.

(2) If there are two eigenvalues in $B_n$ term of the 3-term recurrence relation, an eigenvalue in $A_n$ term has only one valued root. 

The former is referred as the first species complete polynomial and the latter is designated as the second species complete polynomial.
  
For instance, Heun function is considered as the mother of all all well-known special functions such as: Spheroidal Wave, Lame, Mathieu, and hypergeometric $_2F_1$, $_1F_1$ and $_0F_1$ functions. The coefficients in a power series of Heun's equation have a relation between three different coefficients. 
Heun's equation has four different types of a confluent form such as Confluent Heun, Doubly-Confluent Heun, Biconfluent Heun and Triconfluent Heun equations. 
According to Karl Heun(1889) \cite{1Heun1889,1Ronv1995}, Heun's equation is a second-order linear ODE of the form 
\begin{equation}
\frac{d^2{y}}{d{x}^2} + \left(\frac{\gamma }{x} +\frac{\delta }{x-1} + \frac{\epsilon }{x-a}\right) \frac{d{y}}{d{x}} +  \frac{\alpha \beta x-q}{x(x-1)(x-a)} y = 0 \label{eq:1001}
\end{equation}
With the condition $\epsilon = \alpha +\beta -\gamma -\delta +1$. The parameters play different roles: $a \ne 0 $ is the singularity parameter, $\alpha $, $\beta $, $\gamma $, $\delta $, $\epsilon $ are exponent parameters, $q$ is the accessory parameter. Also, $\alpha $ and $\beta $ are identical to each other. It has four regular singular points which are 0, 1, $a$ and $\infty $ with exponents $\{ 0, 1-\gamma \}$, $\{ 0, 1-\delta \}$, $\{ 0, 1-\epsilon \}$ and $\{ \alpha, \beta \}$.

we assume the solution takes the form
\begin{equation}
y(x) = \sum_{n=0}^{\infty } c_n x^{n+\lambda } \label{eq:1002}
\end{equation}
Substituting (\ref{eq:1002}) into (\ref{eq:1001}) gives for the coefficients $c_n$ the recurrence relations
\begin{equation}
c_{n+1}=A_n \;c_n +B_n \;c_{n-1} \hspace{1cm};n\geq 1
\label{eq:1003}
\end{equation}
where,
\begin{subequations}
\begin{equation}
A_n =  \frac{(n+\lambda )(n+\alpha +\beta -\delta +\lambda +a(n+\delta +\gamma -1+\lambda ))+q}{a(n+1+\lambda )(n+\gamma +\lambda )} 
\label{eq:1004a}
\end{equation}
\begin{equation}
B_n = - \frac{(n-1+\lambda +\alpha )(n-1+\lambda +\beta )}{a(n+1+\lambda )(n+\gamma +\lambda )}
\label{eq:1004b}
\end{equation}
\begin{equation}
c_1= A_0 \;c_0 \label{eq:1004c}
\end{equation}
\end{subequations}
We have two indicial roots which are $\lambda =0$ and $1-\gamma$.

If an exponent parameter $\alpha $ (or $\beta $) and an accessory (spectral) parameter $q$ are fixed constants, the first species complete polynomial must to be applied. As parameters $\alpha $, $\beta $ and $q$ are fixed constants, the second species complete polynomial should to be utilized. 

I designate the general summation formulas of a polynomial which makes $A_n$ and $B_n$ terms terminated, by allowing $B_n$ as the leading term in each sub-power series of a power series solution, as ``complete polynomials using reversible 3-term recurrence formula (R3TRF)'' that will be explained below.
\subsection{The first species complete polynomial using R3TRF} 

My definition of $B_{i,j,k,l}$ refer to $B_{i}B_{j}B_{k}B_{l}$. Also, $A_{i,j,k,l}$ refer to $A_{i}A_{j}A_{k}A_{l}$. For $n=0,1,2,3,\cdots $, (\ref{eq:1003}) gives
\begin{equation}
\begin{tabular}{  l  }
  \vspace{2 mm}
  $c_0$ \\
  \vspace{2 mm}
  $c_1 = A_0 c_0 $ \\
  \vspace{2 mm}
  $c_2 = (A_{0,1}+B_1) c_0 $ \\
  \vspace{2 mm}
  $c_3 = (A_{0,1,2}+ B_1 A_2+ B_2 A_0) c_0 $\\
  \vspace{2 mm}
  $c_4 = (A_{0,1,2,3}+B_1 A_{2,3} + B_2 A_{0,3} + B_3 A_{0,1} + B_{1,3}) c_0 $ \\
  \vspace{2 mm}
  $c_5 = (A_{0,1,2,3,4}+ B_1 A_{2,3,4}+ B_2 A_{0,3,4} + B_3 A_{0,1,4} + B_4 A_{0,1,2} $  \\
  \vspace{2 mm}
  \hspace{0.8 cm} $ + B_{1,3} A_4 + B_{1,4} A_2 + B_{2,4} A_0 ) c_0 $ \\     
  \vspace{2 mm}
  $c_6 = (A_{0,1,2,3,4,5}+ B_1 A_{2,3,4,5}+ B_2 A_{0,3,4,5}+ B_3 A_{0,1,4,5} + B_4 A_{0,1,2,5}+ B_5 A_{0,1,2,3} $ \\
  \vspace{2 mm}
  \hspace{0.8 cm} $+ B_{1,3} A_{4,5} + B_{1,4} A_{2,5} + B_{2,4} A_{0,5} + B_{1,5} A_{2,3} + B_{2,5} A_{0,3} + B_{3,5} A_{0,1} + B_{1,3,5}) c_0 $ \\     
  \vspace{2 mm}
  $c_7 = (A_{0,1,2,3,4,5,6}+ B_1 A_{2,3,4,5,6}+ B_2 A_{0,3,4,5,6}+ B_3 A_{0,1,4,5,6} + B_4 A_{0,1,2,5,6}+B_5 A_{0,1,2,3,6} $ \\
  \vspace{2 mm}
  \hspace{0.8 cm} $+ B_6 A_{0,1,2,3,4}+ B_{1,4} A_{2,5,6} + B_{1,3} A_{4,5,6}+ B_{2,4} A_{0,5,6}+ B_{1,5} A_{2,3,6} +B_{2,5} A_{0,3,6} $\\
  \vspace{2 mm}
  \hspace{0.8 cm} $+ B_{3,5} A_{0,1,6} + B_{1,6} A_{2,3,4}+ B_{2,6} A_{0,3,4} + B_{3,6} A_{0,1,4} + B_{4,6} A_{0,1,2}$\\
  \vspace{2 mm}                 
  \hspace{0.8 cm} $+ B_{1,3,6} A_4 + B_{1,3,5} A_6 + B_{1,4,6} A_2 + B_{2,4,6} A_0 ) c_0 $ \\
  \vspace{2 mm}
  $c_8 = (A_{0,1,2,3,4,5,6,7}+B_1 A_{2,3,4,5,6,7} +B_2 A_{0,3,4,5,6,7} + B_3 A_{0,1,4,5,6,7} + B_4 A_{0,1,2,5,6,7}+ B_5 A_{0,1,2,3,6,7} $ \\
  \vspace{2 mm}
  \hspace{0.8 cm} $+B_6 A_{0,1,2,3,4,7}+ B_{7} A_{0,1,2,3,4,5} + B_{1,4} A_{2,5,6,7} + B_{2,4} A_{0,5,6,7}+ B_{1,5} A_{2,3,6,7}+ B_{1,3} A_{4,5,6,7}  $\\
  \vspace{2 mm}
  \hspace{0.8 cm} $+ B_{2,5} A_{0,3,6,7} + B_{3,5} A_{0,1,6,7} + B_{1,6} A_{2,3,4,7} + B_{2,6} A_{0,3,4,7} + B_{3,6} A_{0,1,4,7} + B_{4,6} A_{0,1,2,7} $\\
  \vspace{2 mm}                 
  \hspace{0.8 cm} $+ B_{1,7} A_{2,3,4,5} +B_{3,7} A_{0,1,4,5}+ B_{2,7} A_{0,3,4,5}+ B_{4,7} A_{0,1,2,5}+ B_{5,7} A_{0,1,2,3} $\\
   \vspace{2 mm}
   \hspace{0.8 cm} $+ B_{1,4,6} A_{2,7} + B_{2,4,6} A_{0,7} +B_{1,3,7} A_{4,5}+ B_{1,4,7} A_{2,5}+ B_{1,3,6} A_{4,7} +B_{2,4,7} A_{0,5}$ \\
   \vspace{2 mm}
   \hspace{0.8 cm} $+ B_{1,3,5} A_{6,7}+ B_{1,5,7} A_{2,3}+ B_{2,5,7} A_{0,3}+ B_{3,5,7} A_{0,1} + B_{1,3,5,7}) c_0$\\ 
\hspace{2 mm}\large{\vdots} \hspace{5cm}\large{\vdots}\\ 
\end{tabular}
\label{eq:1007}
\end{equation}
In (\ref{eq:1007}) the number of individual sequence $c_n$ follows Fibonacci sequence: 1,1,2,3,5,8,13,21,34,55,$\cdots$.
The sequence $c_n$ consists of combinations $A_n$ and $B_n$ in (\ref{eq:1007}).\footnote{(\ref{eq:1007}) is equivalent to an expansion of the 3-term recursive relation for $n=0,1,2,\cdots$ on chapter 1. The former is that $B_n$ term is the leading one of each term, composed of $A_n's$ and $B_n's$, inside parenthesis of sequences $c_n's$. The latter is that $A_n$ term is the leading one of each term inside parenthesis of sequences $c_n's$.}

The coefficient $c_n$ with even subscripts consists of even terms of $B_n's$. The sequence $c_n$ with odd subscripts consists of odd terms of $B_n's$. And classify sequences $c_n$ to its even and odd parts in (\ref{eq:1007}). 
\begin{equation}
\begin{tabular}{  l  l }
  \vspace{2 mm}
  \large{$c_0 = \tilde{B}^0$} &\hspace{1cm} \large{$c_1= \tilde{B}^0$}  \\
  \vspace{2 mm}
  \large{$c_2 = \tilde{B}^0+\tilde{B}^1$} &\hspace{1cm} \large{$c_3 = \tilde{B}^0+\tilde{B}^1$} \\
  \vspace{2 mm}
  \large{$c_4 = \tilde{B}^0+\tilde{B}^1+\tilde{B}^2$} &\hspace{1cm}  \large{$c_5 = \tilde{B}^0+\tilde{B}^1+\tilde{B}^2$}\\
  \vspace{2 mm}
  \large{$c_6 = \tilde{B}^0+\tilde{B}^1+\tilde{B}^2+\tilde{B}^3$} &\hspace{1cm}  \large{$c_7 = \tilde{B}^0+\tilde{B}^1+\tilde{B}^2+\tilde{B}^3$}\\
  \vspace{2 mm}
  \large{$c_8 = \tilde{B}^0+\tilde{B}^1+\tilde{B}^2+\tilde{B}^3+\tilde{B}^4$} &\hspace{1cm} \large{$c_9 = \tilde{B}^0+\tilde{B}^1+\tilde{B}^2+\tilde{B}^3+\tilde{B}^4$} \\
 \hspace{2cm} \large{\vdots} & \hspace{3cm}\large{\vdots} \\
\end{tabular}
\label{eq:1008}
\end{equation}
In the above, $\tilde{B}^{\tau } $= $\tau$ terms of $B_n's$ where $\tau =0,1,2,\cdots$.

When a function $y(x)$, analytic at $x=0$, is expanded in a power series, we write
\begin{equation}
y(x)= \sum_{n=0}^{\infty } c_n x^{n+\lambda }= \sum_{\tau =0}^{\infty } y_{\tau}(x) = y_0(x)+ y_1(x)+y_2(x)+ \cdot \label{eq:1009}
\end{equation}
where
\begin{equation}
y_{\tau}(x)= \sum_{l=0}^{\infty } c_l^{\tau} x^{l+\lambda }\label{eq:10010}
\end{equation}
$\lambda $ is the indicial root. $y_{\tau}(x)$ is sub-power series that has sequences $c_n$ including $\tau$ term of $B_n$'s in (\ref{eq:1007}). For example $y_0(x)$ has sequences $c_n$ including zero term of $B_n$'s in (\ref{eq:1007}), $y_1(x)$ has sequences $c_n$ including one term of $B_n$'s in (\ref{eq:1007}), $y_2(x)$ has sequences $c_n$ including two term of $B_n$'s in (\ref{eq:1007}), etc.

First observe the term inside parentheses of sequences $c_n$ which does not include any $B_n$'s in (\ref{eq:1007}) and (\ref{eq:1008}): $c_n$ with every subscripts ($c_0$, $c_1$, $c_2$,$\cdots$). 

\begin{equation}
\begin{tabular}{  l  }
  \vspace{2 mm}
  $c_0$ \\
  \vspace{2 mm}
  $c_1 = A_0 c_0  $ \\
  \vspace{2 mm}
  $c_2 = A_{0,1} c_0  $ \\
  \vspace{2 mm}
  $c_3 = A_{0,1,2}c_0 $ \\
  \vspace{2 mm}
  $c_4 = A_{0,1,2,3}c_0 $\\
  \hspace{2 mm}
  \large{\vdots}\hspace{1cm}\large{\vdots} \\ 
\end{tabular}
\label{eq:10011}
\end{equation}
(\ref{eq:10011}) gives the indicial equation
\begin{equation}
c_{n}= c_0 \prod _{i_0=0}^{n-1}A_{i_0} \hspace{1cm} \mathrm{where}\;n=0,1,2,\cdots 
\label{eq:10012}
\end{equation}
Substitute (\ref{eq:10012}) in (\ref{eq:10010}) putting $\tau = 0$. 
\begin{eqnarray}
y_0^m(x) &=& c_0 \sum_{n=0}^{m} \left\{ \prod _{i_0=0}^{n-1}A_{i_0} \right\} x^{n+\lambda }\label{eq:10013a}\\
&=& c_0 \sum_{i_0=0}^{m} \left\{ \prod _{i_1=0}^{i_0-1}A_{i_1} \right\} x^{i_0+\lambda } \label{eq:10013b}
\end{eqnarray}
(\ref{eq:10013a}) and (\ref{eq:10013b}) are the sub-power series that has sequences $c_0, c_1, c_2, \cdots, c_{m}$ including zero term of $B_n$'s where $m=0,1,2,\cdots$. If $m\rightarrow \infty $, the sub-power series $y_0^m(x)$ turns to be an infinite series $y_0(x)$.

Observe the terms inside parentheses of sequence $c_n$ which include one term of $B_n$'s in (\ref{eq:1007}) and (\ref{eq:1008}): $c_n$ with every index except $c_0$ and $c_1$ ($c_2$, $c_3$, $c_4$,$\cdots$). 

\begin{equation}
\begin{tabular}{  l  }
  \vspace{2 mm}
  $c_2= B_1 c_0$ \\
  \vspace{2 mm}
  $c_3 = \bigg\{ B_1 \cdot 1\cdot \Big( \frac{A_2}{1}\Big) + B_2 A_0\Big(\frac{A_2}{A_2}\Big) \bigg\} c_0  $ \\
  \vspace{2 mm}
  $c_4 = \bigg\{ B_1 \cdot 1\cdot \Big( \frac{A_{2,3}}{1}\Big) + B_2 A_0 \Big(\frac{A_{2,3}}{A_2}\Big) + B_3 A_{0,1} \Big(\frac{A_{2,3}}{A_{2,3}}\Big) \bigg\} c_0  $ \\
  \vspace{2 mm}
  $c_5 = \bigg\{ B_1 \cdot 1\cdot \Big( \frac{A_{2,3,4}}{1}\Big)+ B_2 A_0 \Big(\frac{A_{2,3,4}}{A_2}\Big) + B_3 A_{0,1}\Big(\frac{A_{2,3,4}}{A_{2,3}}\Big) +  B_4 A_{0,1,2}\Big(\frac{A_{2,3,4}}{A_{2,3,4}}\Big) \bigg\} c_0  $ \\
  \vspace{2 mm}
  $c_6 = \bigg\{ B_1 \cdot 1\cdot \Big( \frac{A_{2,3,4,5}}{1}\Big)+ B_2 A_0 \Big(\frac{A_{2,3,4,5}}{A_2}\Big) + B_3 A_{0,1} \Big(\frac{A_{2,3,4,5}}{A_{2,3}}\Big) +  B_4 A_{0,1,2}\Big(\frac{A_{2,3,4,5}}{A_{2,3,4}}\Big)  $\\
  \vspace{2 mm}
  \hspace{0.8 cm} $+  B_5 A_{0,1,2,3}\Big(\frac{A_{2,3,4,5}}{A_{2,3,4,5}}\Big) \bigg\} c_0 $\\
  $c_7 =\bigg\{ B_1 \cdot 1\cdot  \Big( \frac{A_{2,3,4,5,6}}{1}\Big)+ B_2 A_0 \Big(\frac{A_{2,3,4,5,6}}{A_2}\Big) + B_3 A_{0,1}\Big(\frac{A_{2,3,4,5,6}}{A_{2,3}}\Big) +  B_4 A_{0,1,2}\Big(\frac{A_{2,3,4,5,6}}{A_{2,3,4}}\Big) $\\
   \vspace{2 mm}
   \hspace{0.8 cm} $+  B_5 A_{0,1,2,3} \Big(\frac{A_{2,3,4,5,6}}{A_{2,3,4,5}}\Big) +  B_6 A_{0,1,2,3,4}\Big(\frac{A_{2,3,4,5,6}}{A_{2,3,4,5,6}}\Big) \bigg\} c_0 $\\
  \hspace{2 mm}
  \large{\vdots}\hspace{3cm}\large{\vdots} \\ 
\end{tabular}
\label{eq:10014}
\end{equation}
(\ref{eq:10014}) gives the indicial equation
\begin{equation}
c_{n+2}= c_0 \sum_{i_0=0}^{n} \left\{ B_{i_0+1} \prod _{i_1=0}^{i_0-1}A_{i_1} \prod _{i_2=i_0}^{n-1}A_{i_2+2} \right\}  
\label{eq:10015}
\end{equation}
Substitute (\ref{eq:10015}) in (\ref{eq:10010}) putting $\tau = 1$. 
\begin{eqnarray}
y_1^m(x)&=& c_0 \sum_{n=0}^{m}\left\{ \sum_{i_0=0}^{n} \left\{ B_{i_0+1} \prod _{i_1=0}^{i_0-1}A_{i_1} \prod _{i_2=i_0}^{n-1}A_{i_2+2} \right\} \right\} x^{n+2+\lambda } \label{eq:10016a}\\ 
&=& c_0 \sum_{i_0=0}^{m}\left\{ B_{i_0+1} \prod _{i_1=0}^{i_0-1}A_{i_1} \sum_{i_2=i_0}^{m} \left\{ \prod _{i_3=i_0}^{i_2-1}A_{i_3+2} \right\}\right\} x^{i_2+2+\lambda } \label{eq:10016b}
\end{eqnarray}
(\ref{eq:10016a}) and (\ref{eq:10016b}) are the sub-power series that has sequences $c_2, c_3, c_4, \cdots, c_{m+2}$ including one term of $B_n$'s. If $m\rightarrow \infty $, the sub-power series $y_1^m(x)$ turns to be an infinite series $y_1(x)$.

Observe the terms inside parentheses of sequence $c_n$ which include two terms of $B_n$'s in  (\ref{eq:1007}) and (\ref{eq:1008}): $c_n$ with every index except $c_0$--$c_3$ ($c_4$, $c_5$, $c_6$,$\cdots$). 

\begin{equation}
\begin{tabular}{  l  }
  \vspace{2 mm}
  $c_4= B_{1,3} c_0$ \\
  \vspace{2 mm}
  $c_5 = \Bigg\{ B_1 \cdot 1\cdot \bigg\{ B_3 \left( \frac{1}{1}\right) \left(\frac{A_4}{1}\right) + B_4 \left(\frac{A_2}{1}\right) \left( \frac{A_4}{A_4}\right)\bigg\} + B_2 A_0\bigg\{ B_4 \left( \frac{A_2}{A_2}\right) \left( \frac{A_4}{A_4}\right) \bigg\} \Bigg\} c_0  $ \\
  \vspace{2 mm}
  $c_6 = \Bigg\{ B_1 \cdot 1\cdot \bigg\{ B_3 \left( \frac{1}{1}\right) \left(\frac{A_{4,5}}{1}\right)+ B_4 \left(\frac{A_2}{1}\right) \left( \frac{A_{4,5}}{A_4}\right)+ B_5 \left(\frac{A_{2,3}}{1}\right) \left( \frac{A_{4,5}}{A_{4,5}}\right)\bigg\}$\\
  \vspace{2 mm}
  \hspace{0.8 cm} $+ B_2 A_0\bigg\{ B_4 \left( \frac{A_2}{A_2}\right) \left( \frac{A_{4,5}}{A_4}\right) + B_5 \left( \frac{A_{2,3}}{A_2}\right) \left( \frac{A_{4,5}}{A_{4,5}}\right) \bigg\} + B_3  A_{0,1} \bigg\{ B_5 \left( \frac{A_{2,3}}{A_{2,3}}\right) \left( \frac{A_{4,5}}{A_{4,5}}\right)\bigg\} \Bigg\} c_0 $ \\
\vspace{2 mm}
    $c_7 = \Bigg\{ B_1\cdot  1 \cdot \bigg\{ B_3 \left( \frac{1}{1}\right) \left(\frac{A_{4,5,6}}{1}\right) + B_4 \left(\frac{A_2}{1}\right) \left( \frac{A_{4,5,6}}{A_4}\right) + B_5 \left(\frac{A_{2,3}}{1}\right) \left( \frac{A_{4,5,6}}{A_{4,5}}\right) +  B_6 \left(\frac{A_{2,3,4}}{1}\right) \left( \frac{A_{4,5,6}}{A_{4,5,6}}\right) \bigg\}$\\
\vspace{2 mm}
  \hspace{0.8 cm} $+ B_2 A_0\bigg\{ B_4 \left( \frac{A_2}{A_2}\right) \left( \frac{A_{4,5,6}}{A_4}\right)  + B_5 \left( \frac{A_{2,3}}{A_2}\right) \left( \frac{A_{4,5,6}}{A_{4,5}}\right) + B_6 \left( \frac{A_{2,3,4}}{A_2}\right) \left( \frac{A_{4,5,6}}{A_{4,5,6}}\right) \bigg\}$\\
\vspace{2 mm}
  \hspace{0.8 cm} $+ B_3 A_{0,1} \bigg\{ B_5 \left( \frac{A_{2,3}}{A_{2,3}}\right) \left( \frac{A_{4,5,6}}{A_{4,5}}\right) +  B_6 \left( \frac{A_{2,3,4}}{A_{2,3}}\right) \left( \frac{A_{4,5,6}}{A_{4,5,6}}\right) \bigg\}$\\
\vspace{2 mm}
  \hspace{0.8 cm} $+ B_4 A_{0,1,2}\bigg\{ B_6 \Big( \frac{A_{2,3,4}}{A_{2,3,4}}\Big) \Big( \frac{A_{4,5,6}}{A_{4,5,6}}\Big)  \bigg\} \Bigg\} c_0$\\
  \hspace{2 mm}
  \large{\vdots}\hspace{5cm}\large{\vdots} \\  
\end{tabular}
\label{eq:10017}
\end{equation}
(\ref{eq:10017}) gives the indicial equation
\begin{eqnarray}
 c_{n+4} &=& c_0 \sum_{i_0=0}^{n} \left\{ B_{i_0+1}\sum_{i_1=i_0}^{n} \left\{ B_{i_1+3}{\displaystyle \prod _{i_2=0}^{i_0-1}A_{i_2} \prod _{i_3=i_0}^{i_1-1}A_{i_3+2}\prod _{i_4=i_1}^{n-1}A_{i_4+4}}\right\}\right\}  
\label{eq:10018}\\
 &=&  c_0 \sum_{i_0=0}^{n} \left\{ B_{i_0+1}\prod _{i_1=0}^{i_0-1}A_{i_1} \sum_{i_2=i_0}^{n} \left\{ B_{i_2+3}  \prod _{i_3=i_0}^{i_2-1}A_{i_3+2}\prod _{i_4=i_2}^{n-1}A_{i_4+4} \right\}\right\} \label{eq:x10018}
\end{eqnarray}
Substitute (\ref{eq:x10018}) in (\ref{eq:10010}) putting $\tau = 2$. 
\begin{eqnarray}
 y_2^m(x) &=&  c_0 \sum_{n=0}^{m}\left\{ \sum_{i_0=0}^{n} \left\{ B_{i_0+1}\prod _{i_1=0}^{i_0-1}A_{i_1} \sum_{i_2=i_0}^{n} \left\{ B_{i_2+3}  \prod _{i_3=i_0}^{i_2-1}A_{i_3+2}\prod _{i_4=i_2}^{n-1}A_{i_4+4} \right\}\right\} \right\} x^{n+4+\lambda }  \label{eq:10019a}\\
&=& c_0 \sum_{i_0=0}^{m}\left\{ B_{i_0+1} \prod _{i_1=0}^{i_0-1}A_{i_1}  \sum_{i_2=i_0}^{m} \left\{ B_{i_2+3} \prod _{i_3=i_0}^{i_2-1}A_{i_3+2} \sum_{i_4=i_2}^{m}\left\{ \prod _{i_5=i_2}^{i_4-1}A_{i_5+4}\right\}  \right\}\right\} x^{i_4+4+\lambda } \hspace{2cm} \label{eq:10019b}
\end{eqnarray}
(\ref{eq:10019a}) and (\ref{eq:10019b}) are the sub-power series that has sequences $c_4, c_5, c_6, \cdots, c_{m+4}$ including two terms of $B_n$'s. If $m\rightarrow \infty $, the sub-power series $y_2^m(x)$ turns to be an infinite series $y_2(x)$.

Observe the terms inside parentheses of sequence $c_n$ which include three terms of $B_n$'s in (\ref{eq:1007}) and (\ref{eq:1008}): $c_n$ with every index except $c_0$--$c_5$ ($c_6$, $c_7$, $c_8$,$\cdots$). 
\begin{equation}
\begin{tabular}{  l  }
  \vspace{2 mm}
  $c_6= B_{1,3,5} \;c_0$ \\
  \vspace{2 mm}
  $c_7 = \Bigg\{ B_1 \bigg\{ B_3 \cdot 1\cdot \Big[ B_5 \left( \frac{1}{1}\right) \left( \frac{1}{1}\right) \left(\frac{A_{6}}{1}\right)+ B_6 \left( \frac{1}{1}\right) \left(\frac{A_4}{1}\right) \left( \frac{A_{6}}{A_{6}}\right)\Big] $\\
 \vspace{2 mm}
  \hspace{0.8 cm} $+ B_4 \cdot 1 \cdot \Big[ B_6 \left( \frac{A_2}{1}\right) \left( \frac{A_{4}}{A_{4}}\right) \left( \frac{A_{6}}{A_{6}}\right) \Big] \bigg\} + B_2 \bigg\{ B_4 A_0 \Big[ B_6 \left( \frac{A_{2}}{A_{2}}\right)\left( \frac{A_{4}}{A_{4}}\right)\left( \frac{A_{6}}{A_{6}}\right)\Big] \bigg\} \Bigg\} c_0  $ \\
  \vspace{2 mm}
  $c_8 = \Bigg\{ B_1 \bigg\{ B_3 \cdot 1\cdot  \Big[ B_5 \left( \frac{1}{1}\right) \left( \frac{1}{1}\right) \left(\frac{A_{6,7}}{1}\right)+ B_6 \left( \frac{1}{1}\right) \left(\frac{A_4}{1}\right) \left( \frac{A_{6,7}}{A_{6}}\right)+ B_7 \left( \frac{1}{1}\right) \left(\frac{A_{4,5}}{1}\right) \left( \frac{A_{6,7}}{A_{{6,7}}}\right)\Big]$\\
 \vspace{2 mm}
  \hspace{0.8 cm} $+ B_4 \cdot 1 \cdot \Big[ B_6 \left( \frac{A_2}{1}\right) \left( \frac{A_{4}}{A_{4}}\right) \left( \frac{A_{6,7}}{A_{6}}\right) + B_7 \left( \frac{A_2}{1}\right) \left( \frac{A_{4,5}}{A_{4}}\right) \left( \frac{A_{6,7}}{A_{6,7}}\right)\Big]+ B_5 \cdot 1 \cdot \Big[ B_7 \left( \frac{A_{2,3}}{1}\right) \left( \frac{A_{4,5}}{A_{4,5}}\right) \left( \frac{A_{6,7}}{A_{6,7}}\right) \Big] \bigg\} $\\
 \vspace{2 mm}
  \hspace{0.8 cm} $+ B_2 \bigg\{ B_4  A_0\Big[ B_6 \left( \frac{A_{2}}{A_{2}}\right)\left( \frac{A_{4}}{A_{4}}\right)\left( \frac{A_{6,7}}{A_{6}}\right)+ B_7 \left( \frac{A_{2}}{A_{2}}\right)\left( \frac{A_{4,5}}{A_{4}}\right)\left( \frac{A_{6,7}}{A_{6,7}}\right)\Big] $\\
\vspace{2 mm}
  \hspace{0.8 cm}$ + B_5 A_0\Big[ B_7 \left( \frac{A_{2,3}}{A_{2}}\right) \left( \frac{A_{4,5}}{A_{4,5}}\right) \left( \frac{A_{6,7}}{A_{6,7}}\right) \Big] \bigg\}+ B_3 \bigg\{ B_5 A_{0,1} \Big[ B_7 \left( \frac{A_{2,3}}{A_{2,3}}\right) \left( \frac{A_{4,5}}{A_{4,5}}\right) \left( \frac{A_{6,7}}{A_{6,7}}\right) \Big]\bigg\} \Bigg\} c_0 $ \\
 \vspace{2 mm}
   $c_9 = \Bigg\{ B_1 \bigg\{ B_3 \cdot 1 \cdot \Big[ B_5 \left( \frac{1}{1}\right) \left( \frac{1}{1}\right) \left(\frac{A_{6,7,8}}{1}\right)+ B_6 \left( \frac{1}{1}\right) \left(\frac{A_4}{1}\right) \left( \frac{A_{6,7,8}}{A_{6}}\right)+ B_7 \left( \frac{1}{1}\right) \left(\frac{A_{4,5}}{1}\right) \left( \frac{A_{6,7,8}}{A_{{6,7}}}\right) $\\ 
\vspace{2 mm}
  \hspace{0.8 cm} $+ B_8 \left( \frac{1}{1}\right) \left(\frac{A_{4,5,6}}{1}\right) \left( \frac{A_{6,7,8}}{A_{{6,7,8}}}\right)\Big]$\\
\vspace{2 mm}
\hspace{0.8 cm} $+ B_4 \cdot 1\cdot  \Big[ B_6 \left( \frac{A_2}{1}\right) \left( \frac{A_{4}}{A_{4}}\right) \left( \frac{A_{6,7,8}}{A_{6}}\right) + B_7 \left( \frac{A_2}{1}\right) \left( \frac{A_{4,5}}{A_{4}}\right) \left( \frac{A_{6,7,8}}{A_{6,7}}\right)+ B_8 \left( \frac{A_2}{1}\right) \left( \frac{A_{4,5,6}}{A_{4}}\right) \left( \frac{A_{6,7,8}}{A_{6,7,8}}\right) \Big]$\\
  \vspace{2 mm}
  \hspace{0.8 cm} $+ B_5 \cdot 1 \cdot \Big[ B_7 \left( \frac{A_{2,3}}{1}\right) \left( \frac{A_{4,5}}{A_{4,5}}\right) \left( \frac{A_{6,7,8}}{A_{6,7}}\right)+ B_8 \left( \frac{A_{2,3}}{1}\right) \left( \frac{A_{4,5,6}}{A_{4,5}}\right) \left( \frac{A_{6,7,8}}{A_{6,7,8}}\right)\Big] $\\
\vspace{2 mm}
  \hspace{0.8 cm} $+ B_6 \cdot 1 \cdot \Big[ B_8 \left( \frac{A_{2,3,4}}{1}\right) \left( \frac{A_{4,5,6}}{A_{4,5,6}}\right) \left( \frac{A_{6,7,8}}{A_{6,7,8}}\right) \Big] \bigg\} $\\
 \vspace{2 mm}
  \hspace{0.8 cm} $+ B_2 \bigg\{ B_4 A_0 \Big[ B_6  \left( \frac{A_{2}}{A_{2}}\right)\left( \frac{A_{4}}{A_{4}}\right)\left( \frac{A_{6,7,8}}{A_{6}}\right)+ B_7 \left( \frac{A_{2}}{A_{2}}\right)\left( \frac{A_{4,5}}{A_{4}}\right)\left( \frac{A_{6,7,8}}{A_{6,7}}\right) + B_8 \left( \frac{A_{2}}{A_{2}}\right)\left( \frac{A_{4,5,6}}{A_{4}}\right)\left( \frac{A_{6,7,8}}{A_{6,7,8}}\right)\Big]$\\
  \vspace{2 mm}
  \hspace{0.8 cm} $+ B_5 A_0 \Big[ B_7 \left( \frac{A_{2,3}}{A_{2}}\right) \left( \frac{A_{4,5}}{A_{4,5}}\right) \left( \frac{A_{6,7,8}}{A_{6,7}}\right) + B_8 \left( \frac{A_{2,3}}{A_{2}}\right) \left( \frac{A_{4,5,6}}{A_{4,5}}\right) \left( \frac{A_{6,7,8}}{A_{6,7,8}}\right)\Big] $\\
  \vspace{2 mm}
  \hspace{0.8 cm} $+ B_6 A_0 \Big[ B_8 \left( \frac{A_{2,3,4}}{A_{2}}\right) \left( \frac{A_{4,5,6}}{A_{4,5,6}}\right) \left( \frac{A_{6,7,8}}{A_{6,7,8}}\right) \Big] \bigg\}$\\
\end{tabular}
\end{equation}  
\begin{equation}
\begin{tabular}{  l  }  
\vspace{2 mm}
  \hspace{0.8 cm} $+ B_3 \bigg\{ B_5  A_{0,1}\Big[ B_7 \left( \frac{A_{2,3}}{A_{2,3}}\right) \left( \frac{A_{4,5}}{A_{4,5}}\right) \left( \frac{A_{6,7,8}}{A_{6,7}}\right) + B_8 \left( \frac{A_{2,3}}{A_{2,3}}\right) \left( \frac{A_{4,5,6}}{A_{4,5}}\right) \left( \frac{A_{6,7,8}}{A_{6,7,8}}\right)\Big] $\\
\vspace{2 mm}
 \hspace{0.8 cm} $+ B_6 A_{0,1}\Big[ B_8 \left( \frac{A_{2,3,4}}{A_{2,3}}\right) \left( \frac{A_{4,5,6}}{A_{4,5,6}}\right) \left( \frac{A_{6,7,8}}{A_{6,7,8}}\right) \Big] \bigg\}+ B_4 \bigg\{ B_6 A_{0,1,2} \Big[ B_8 \left( \frac{A_{2,3,4}}{A_{2,3,4}}\right) \left( \frac{A_{4,5,6}}{A_{4,5,6}}\right) \left( \frac{A_{6,7,8}}{A_{6,7,8}}\right)\Big] \bigg\} \Bigg\} c_0 $ \\
  \hspace{2 mm}
  \large{\vdots}\hspace{6cm}\large{\vdots} \\  
\end{tabular}
\label{eq:10020}
\end{equation}
(\ref{eq:10020}) gives the indicial equation
\begin{eqnarray}
 c_{n+6} &=& c_0 \sum_{i_0=0}^{n} \left\{ B_{i_0+1}\sum_{i_1=i_0}^{n} \left\{ B_{i_1+3}\sum_{i_2=i_1}^{n} \left\{ B_{i_2+5} \prod _{i_3=0}^{i_0-1}A_{i_3} \right.\right.\right.\nonumber\\
&&\times \left.\left.\left. \prod _{i_4=i_0}^{i_1-1}A_{i_4+2} \prod _{i_5=i_1}^{i_2-1}A_{i_5+4}\prod _{i_6=i_2}^{n-1}A_{i_6+6} \right\} \right\} \right\} \label{eq:10021}\\
 &=& c_0 \sum_{i_0=0}^{n} \left\{ B_{i_0+1}\prod _{i_1=0}^{i_0-1}A_{i_1} \sum_{i_2=i_0}^{n} \left\{ B_{i_2+3}\prod _{i_3=i_0}^{i_2-1}A_{i_3+2} \right.\right. \nonumber\\
&&\times \left.\left. \sum_{i_4=i_2}^{n} \left\{ B_{i_4+5} \prod _{i_5=i_2}^{i_4-1}A_{i_5+4} \prod _{i_6=i_4}^{n-1}A_{i_6+6} \right\} \right\} \right\} \label{eq:x10021}
\end{eqnarray}
Substitute (\ref{eq:x10021}) in (\ref{eq:10010}) putting $\tau = 3$.
\begin{eqnarray}
 y_3^m(x) &=& c_0 \sum_{n=0}^{m}\left\{\sum_{i_0=0}^{n} \left\{ B_{i_0+1}\prod _{i_1=0}^{i_0-1}A_{i_1} \sum_{i_2=i_0}^{n} \left\{ B_{i_2+3}\prod _{i_3=i_0}^{i_2-1}A_{i_3+2} \right.\right.\right. \nonumber\\
&&\times \left.\left.\left. \sum_{i_4=i_2}^{n} \left\{ B_{i_4+5} \prod _{i_5=i_2}^{i_4-1}A_{i_5+4} \prod _{i_6=i_4}^{n-1}A_{i_6+6} \right\} \right\} \right\}\right\} x^{n+6+\lambda } \label{eq:10022a}\\
&=&  c_0 \sum_{i_0=0}^{m}\left\{ B_{i_0+1} \prod _{i_1=0}^{i_0-1}A_{i_1} \sum_{i_2=i_0}^{m} \left\{ B_{i_2+3} \prod _{i_3=i_0}^{i_2-1}A_{i_3+2}\right.\right.\nonumber\\
&&\times   \left.\left.\sum_{i_4=i_2}^{m}\left\{ B_{i_4+5}\prod _{i_5=i_2}^{i_4-1}A_{i_5+4} \sum_{i_6=i_4}^{m} \left\{\prod _{i_7=i_4}^{i_6-1}A_{i_7+6}  \right\}\right\} \right\}\right\} x^{i_6+6+\lambda }\label{eq:10022b}
\end{eqnarray}
(\ref{eq:10022a}) and (\ref{eq:10022b}) are the sub-power series that has sequences $c_6, c_7, c_8, \cdots, c_{m+6}$ including three terms of $B_n$'s. If $m\rightarrow \infty $, the sub-power series $y_3^m(x)$ turns to be an infinite series $y_3(x)$.

By repeating this process for all higher terms of $B$'s, we obtain every indicial equations for the sequence $c_{n+2\tau }$ and the sub-power series $y_{\tau }^m(x)$ where $\tau  \geq 4$. 

According (\ref{eq:10012}), (\ref{eq:10015}), (\ref{eq:x10018}), (\ref{eq:x10021}) and every $c_{n+2\tau }$ where $\tau  \geq 4$, the general expression of coefficients $c_{n+2\tau }$ for a fixed $n,\tau \in \mathbb{N}_{0}$ is taken by
\begin{eqnarray}
\tilde{c}(0,n) \;=\;\;c_{n} &=& c_0 \prod _{i_0=0}^{n-1}A_{i_0} \label{eq:10023a}\\
\tilde{c}(1,n) = c_{n+2} &=& c_0 \sum_{i_0=0}^{n} \left\{ B_{i_0+1} \prod _{i_1=0}^{i_0-1}A_{i_1} \prod _{i_2=i_0}^{n-1}A_{i_2+2} \right\} \label{eq:10023b}\\
\tilde{c}(\tau ,n) = c_{n+2\tau } &=& c_0 \sum_{i_0=0}^{n} \left\{B_{i_0+1}\prod _{i_1=0}^{i_0-1} A_{i_1} 
\prod _{k=1}^{\tau -1} \left( \sum_{i_{2k}= i_{2(k-1)}}^{n} B_{i_{2k}+(2k+1)}\prod _{i_{2k+1}=i_{2(k-1)}}^{i_{2k}-1}A_{i_{2k+1}+2k}\right) \right. \nonumber\\
&& \times \left. \prod _{i_{2\tau} = i_{2(\tau -1)}}^{n-1} A_{i_{2\tau }+ 2\tau} \right\}\;\; \mathrm{where}\;\tau \geq 2
\hspace{2cm} \label{eq:10023c}  
\end{eqnarray}
According (\ref{eq:10013b}), (\ref{eq:10016b}), (\ref{eq:10019b}), (\ref{eq:10022b}) and every $y_{\tau }^m(x)$ where $\tau  \geq 4$, the general expression of the sub-power series $y_{\tau }^m(x)$ for a fixed $\tau \in \mathbb{N}_{0}$ is given by
\begin{eqnarray}
y_0^m(x) &=& c_0 \sum_{i_0=0}^{m} \left\{ \prod _{i_1=0}^{i_0-1}A_{i_1} \right\} x^{i_0+\lambda} \label{eq:10024a}\\
y_1^m(x) &=& c_0 \sum_{i_0=0}^{m}\left\{ B_{i_0+1} \prod _{i_1=0}^{i_0-1}A_{i_1}  \sum_{i_2=i_0}^{m} \left\{ \prod _{i_3=i_0}^{i_2-1}A_{i_3+2} \right\}\right\} x^{i_2+2+\lambda } \label{eq:10024b}\\
y_{\tau }^m(x) &=& c_0 \sum_{i_0=0}^{m} \left\{B_{i_0+1}\prod _{i_1=0}^{i_0-1} A_{i_1} 
\prod _{k=1}^{\tau -1} \left( \sum_{i_{2k}= i_{2(k-1)}}^{m} B_{i_{2k}+(2k+1)}\prod _{i_{2k+1}=i_{2(k-1)}}^{i_{2k}-1}A_{i_{2k+1}+2k}\right) \right. \nonumber\\
&& \times \left. \sum_{i_{2\tau} = i_{2(\tau -1)}}^{m} \left( \prod _{i_{2\tau +1}=i_{2(\tau -1)}}^{i_{2\tau}-1} A_{i_{2\tau +1}+ 2\tau} \right) \right\} x^{i_{2\tau}+2\tau +\lambda}\;\;\mathrm{where}\;\tau \geq 2
\label{eq:10024c} 
\end{eqnarray}
\begin{definition}
For the first species complete polynomial, we need a condition which is given by
\begin{equation}
B_{j+1}= c_{j+1}=0\hspace{1cm}\mathrm{where}\;j=0,1,2,\cdots   
 \label{eq:10025}
\end{equation}
\end{definition}
(\ref{eq:10025}) gives successively $c_{j+2}=c_{j+3}=c_{j+4}=\cdots=0$. And $c_{j+1}=0$ is defined by a polynomial equation of degree $j+1$ for the determination of an accessory parameter in $A_n$ term.
\begin{theorem}
The general expression of a function $y(x)$ for the first species complete polynomial using reversible 3-term recurrence formula and its algebraic equation for the determination of an accessory parameter in $A_n$ term are given by
\begin{enumerate}
\item As $B_1=0$,
\begin{equation}
0 =\bar{c}(0,1) \nonumber 
\end{equation}
\begin{equation}
y(x) = y_{0}^{0}(x) \nonumber 
\end{equation}
\item As $B_2=0$, 
\begin{equation}
0 = \bar{c}(0,2)+\bar{c}(1,0) \nonumber 
\end{equation}
\begin{equation}
y(x)= y_{0}^{1}(x)  \nonumber 
\end{equation}
\item As $B_{2N+3}=0$ where $N \in \mathbb{N}_{0}$,
\begin{equation}
0  = \sum_{r=0}^{N+1}\bar{c}\left( r, 2(N-r)+3\right) \nonumber 
\end{equation}
\begin{equation}
y(x)= \sum_{r=0}^{N+1} y_{r}^{2(N+1-r)}(x) \nonumber 
\end{equation}
\item As $B_{2N+4}=0$ where $N \in \mathbb{N}_{0}$,
\begin{equation}
0  = \sum_{r=0}^{N+2}\bar{c}\left( r, 2(N+2-r)\right) \nonumber 
\end{equation}
\begin{equation}
y(x)= \sum_{r=0}^{N+1} y_{r}^{2(N-r)+3}(x) \nonumber 
\end{equation}
In the above,
\begin{eqnarray}
\bar{c}(0,n) &=& \prod _{i_0=0}^{n-1}A_{i_0} \nonumber\\
\bar{c}(1,n) &=& \sum_{i_0=0}^{n} \left\{ B_{i_0+1} \prod _{i_1=0}^{i_0-1}A_{i_1} \prod _{i_2=i_0}^{n-1}A_{i_2+2} \right\} \nonumber\\
\bar{c}(\tau ,n) &=& \sum_{i_0=0}^{n} \left\{B_{i_0+1}\prod _{i_1=0}^{i_0-1} A_{i_1} 
\prod _{k=1}^{\tau -1} \left( \sum_{i_{2k}= i_{2(k-1)}}^{n} B_{i_{2k}+(2k+1)}\prod _{i_{2k+1}=i_{2(k-1)}}^{i_{2k}-1}A_{i_{2k+1}+2k}\right)\right.\nonumber\\
&&\times \left.\prod _{i_{2\tau} = i_{2(\tau -1)}}^{n-1} A_{i_{2\tau }+ 2\tau} \right\} 
\nonumber  
\end{eqnarray}
and
\begin{eqnarray}
y_0^m(x) &=& c_0 x^{\lambda} \sum_{i_0=0}^{m} \left\{ \prod _{i_1=0}^{i_0-1}A_{i_1} \right\} x^{i_0 } \nonumber\\
y_1^m(x) &=& c_0 x^{\lambda} \sum_{i_0=0}^{m}\left\{ B_{i_0+1} \prod _{i_1=0}^{i_0-1}A_{i_1}  \sum_{i_2=i_0}^{m} \left\{ \prod _{i_3=i_0}^{i_2-1}A_{i_3+2} \right\}\right\} x^{i_2+2 } \nonumber\\
y_{\tau }^m(x) &=& c_0 x^{\lambda} \sum_{i_0=0}^{m} \left\{B_{i_0+1}\prod _{i_1=0}^{i_0-1} A_{i_1} 
\prod _{k=1}^{\tau -1} \left( \sum_{i_{2k}= i_{2(k-1)}}^{m} B_{i_{2k}+(2k+1)}\prod _{i_{2k+1}=i_{2(k-1)}}^{i_{2k}-1}A_{i_{2k+1}+2k}\right) \right. \nonumber\\
&& \times \left. \sum_{i_{2\tau} = i_{2(\tau -1)}}^{m} \left( \prod _{i_{2\tau +1}=i_{2(\tau -1)}}^{i_{2\tau}-1} A_{i_{2\tau +1}+ 2\tau} \right) \right\} x^{i_{2\tau}+2\tau }\hspace{1cm}\mathrm{where}\;\tau \geq 2
\nonumber 
\end{eqnarray}
\end{enumerate}
\end{theorem}
\begin{proof}
For instance, if $B_1= c_1=0$ in (\ref{eq:10025}), then gives $c_{2}=c_{3}=c_{4}=\cdots=0$. According to (\ref{eq:1009}), its power series is given by
\begin{equation}
y(x)= \sum_{n=0}^{0} c_n x^{n+\lambda } = c_0 x^{\lambda } = y_{0}^{0}(x)\label{eq:10026}
\end{equation}
The sub-power series $y_{0}^{0}(x)$ in (\ref{eq:10026}) is obtain by putting $m=0$ in (\ref{eq:10024a}).
And a polynomial equation of degree 1 for the determination of an accessory parameter in $A_n$ term is taken by
\begin{equation}
0 = c_1 = \tilde{c}(0,1) = c_0 A_0  \label{eq:10027}
\end{equation}
A coefficient $\tilde{c}(0,1)$ in (\ref{eq:10027}) is obtained by putting $n=1$ in (\ref{eq:10023a}).

If $B_2= c_2=0$ in (\ref{eq:10025}), then gives $c_{3}=c_{4}=c_{5}=\cdots=0$. According to (\ref{eq:1009}), its power series is given by
\begin{equation}
y(x)= \sum_{n=0}^{1} c_n x^{n+\lambda } = (c_0 +c_1 x) x^{\lambda } \label{eq:10028}
\end{equation}
First observe sequences $c_0$ and $c_1$ in (\ref{eq:1007}) and (\ref{eq:1008}). 
The sub-power series that has sequences $c_0$ and $c_1$ is given by  putting $m=1$ in (\ref{eq:10024a}). Take the new (\ref{eq:10024a}) into (\ref{eq:10028}).
\begin{equation}
y(x)= y_{0}^{1}(x) = c_0 \left( 1+A_0 x\right)x^{\lambda } \label{eq:10029}
\end{equation}
A sequence $c_2$ consists of zero and one terms of $B_n's$ in (\ref{eq:1007}) and (\ref{eq:1008}). 
Putting $n=2$ in (\ref{eq:10023a}), a coefficient $c_2$ for zero term of $B_n's$ is denoted by $\tilde{c}(0,2)$. 
Taking  $n=0$ in (\ref{eq:10023b}), a coefficient $c_2$ for one term of $B_n's$ is denoted by $\tilde{c}(1,0)$.
Since the sum of $\tilde{c}(0,2)$ and $\tilde{c}(1,0)$ is equivalent to zero, we obtain a polynomial equation of degree 2 for the determination of an accessory parameter in $A_n$ term which is given by
\begin{equation}
0 = c_2 = \tilde{c}(0,2) + \tilde{c}(1,0)  = c_0 \left( A_{0,1} +B_1 \right)  \label{eq:10030}
\end{equation} 

If $B_3= c_3=0$ in (\ref{eq:10025}), then gives $c_{4}=c_{5}=c_{6}=\cdots=0$. According to (\ref{eq:1009}), its power series is given by
\begin{equation}
y(x)= \sum_{n=0}^{2} c_n x^{n+\lambda } = (c_0 +c_1 x+c_2 x^2) x^{\lambda } \label{eq:10031}
\end{equation}
Observe sequences $c_0$--$c_2$ in (\ref{eq:1007}) and (\ref{eq:1008}). 
The sub-power series, having sequences $c_0$--$c_2$ including zero term of $B_n's$, is given by  putting $m=2$ in (\ref{eq:10024a}) denoted by $y_0^2(x)$.
The sub-power series, having a sequence $c_2$ including one term of $B_n's$, is given by  putting $m=0$ in (\ref{eq:10024b}) denoted by $y_1^0(x)$.
Taking $y_0^2(x)$ and $y_1^0(x)$ into (\ref{eq:10031}),
\begin{equation}
y(x)= y_{0}^{2}(x) + y_{1}^{0}(x) = c_0x^{\lambda } \left( 1+A_0 x + \left( A_{0,1} +B_1\right) x^2\right) \label{eq:10032}
\end{equation} 
A sequence $c_3$ consists of zero and one terms of $B_n's$ in (\ref{eq:1007}) and (\ref{eq:1008}).
Putting $n=3$ in (\ref{eq:10023a}), a coefficient $c_3$ for zero term of $B_n's$ is denoted by $\tilde{c}(0,3)$.
Taking $n=1$  in (\ref{eq:10023b}), a coefficient $c_3$ for one term of $B_n's$ is denoted by $\tilde{c}(1,1)$.
Since the sum of $\tilde{c}(0,3)$ and $\tilde{c}(1,1)$ is equivalent to zero, we obtain a polynomial equation of degree 3 for the determination of an accessory parameter in $A_n$ term which is given by
\begin{equation}
0 = c_3 = \tilde{c}(0,3) + \tilde{c}(1,1)  = c_0 \left( A_{0,1,2} +B_1 A_2 + B_2 A_0 \right)  \label{eq:10033}
\end{equation} 

If $B_4= c_4=0$ in (\ref{eq:10025}), then gives $c_{5}=c_{6}=c_{7}=\cdots=0$. According to (\ref{eq:1009}), its power series is given by
\begin{equation}
y(x)= \sum_{n=0}^{3} c_n x^{n+\lambda } = (c_0 +c_1 x+c_2 x^2+c_3 x^3) x^{\lambda } \label{eq:10034}
\end{equation}
Observe sequences $c_0$--$c_3$ in (\ref{eq:1007}) and (\ref{eq:1008}). 
The sub-power series, having sequences $c_0$--$c_3$ including zero term of $B_n's$, is given by  putting $m=3$ in (\ref{eq:10024a}) denoted by $y_0^3(x)$.
The sub-power series, having sequences $c_2$ and $c_3$ including one term of $B_n's$, is given by  putting $m=1$ in (\ref{eq:10024b}) denoted by $y_1^1(x)$.
Taking $y_0^3(x)$ and $y_1^1(x)$ into (\ref{eq:10034}),
\begin{equation}
y(x)= y_{0}^{3}(x)+y_{1}^{1}(x) \label{eq:10035}
\end{equation}
A sequence $c_4$ consists of zero, one and two terms of $B_n's$ in (\ref{eq:1007}) and (\ref{eq:1008}).
Putting $n=4$ in (\ref{eq:10023a}), a coefficient $c_4$ for zero term of $B_n's$ is denoted by $\tilde{c}(0,4)$.
Taking $n=2$  in (\ref{eq:10023b}), a coefficient $c_4$ for one term of $B_n's$ is denoted by $\tilde{c}(1 ,2)$.
Taking $\tau=2$ and $n=0$  in (\ref{eq:10023c}), a coefficient $c_4$ for two terms of $B_n's$ is denoted by $\tilde{c}(2 ,0)$.
Since the sum of $\tilde{c}(0,4)$, $\tilde{c}(1,2)$ and $\tilde{c}(2,0)$ is equivalent to zero, we obtain a polynomial equation of degree 4 for the determination of an accessory parameter in $A_n$ term which is given by
\begin{equation}
0 = c_4 = \tilde{c}(0,4) + \tilde{c}(1 ,2) + \tilde{c}(2 ,0) \label{eq:10036}
\end{equation}

If $B_5= c_5=0$ in (\ref{eq:10025}), then gives $c_{6}=c_{7}=c_{8}=\cdots=0$. According to (\ref{eq:1009}), its power series is given by
\begin{equation}
y(x)= \sum_{n=0}^{4} c_n x^{n+\lambda } = (c_0 +c_1 x+c_2 x^2+c_3 x^3+c_4 x^4) x^{\lambda } \label{eq:10037}
\end{equation}
Observe sequences $c_0$--$c_4$ in (\ref{eq:1007}) and (\ref{eq:1008}). 
The sub-power series, having sequences  $c_0$--$c_4$ including zero term of $B_n's$, is given by  putting $m=4$ in (\ref{eq:10024a}) denoted by $y_0^4(x)$.
The sub-power series, having sequences  $c_2$--$c_4$ including one term of $B_n's$, is given by  putting $m=2$ in (\ref{eq:10024b}) denoted by $y_1^2(x)$.
The sub-power series, having sequences $c_4$ including two terms of $B_n's$, is given by  putting $\tau=2$ and $m=0$ in (\ref{eq:10024c}) denoted by $y_2^0(x)$.
Taking $y_0^4(x)$, $y_1^2(x)$ and $y_2^0(x)$ into (\ref{eq:10037}),
\begin{equation}
y(x)= y_0^4(x)+y_1^2(x)+y_2^0(x) \label{eq:10038}
\end{equation}
A sequence $c_5$ consists of zero, one and two terms of $B_n's$ in (\ref{eq:1007}) and (\ref{eq:1008}).
Putting $n=5$ in (\ref{eq:10023a}), a coefficient $c_5$ for zero term of $B_n's$ is denoted by $\tilde{c}(0,5)$.
Taking $n=3$  in (\ref{eq:10023b}), a coefficient $c_5$ for one term of $B_n's$ is denoted by $\tilde{c}(1,3)$.
Taking $\tau=2$ and $n=1$  in (\ref{eq:10023c}), a coefficient $c_5$ for two terms of $B_n's$ is denoted by $\tilde{c}(2,1)$.
Since the sum of $\tilde{c}(0,5)$, $\tilde{c}(1,3)$ and $\tilde{c}(2,1)$ is equivalent to zero, we obtain a polynomial equation of degree 5 for the determination of an accessory parameter in $A_n$ term which is given by
\begin{equation}
0 = c_5 = \tilde{c}(0,5) + \tilde{c}(1,3) + \tilde{c}(2,1) \label{eq:10039}
\end{equation}

If $B_6= c_6=0$ in (\ref{eq:10025}), then gives $c_{7}=c_{8}=c_{9}=\cdots=0$. According to (\ref{eq:1009}), its power series is given by
\begin{equation}
y(x)= \sum_{n=0}^{5} c_n x^{n+\lambda } = (c_0 +c_1 x+c_2 x^2+c_3 x^3+c_4 x^4+c_5 x^5) x^{\lambda } \label{eq:10040}
\end{equation}
Observe sequences $c_0$--$c_5$ in (\ref{eq:1007}) and (\ref{eq:1008}). 
The sub-power series, having sequences $c_0$--$c_5$ including zero term of $B_n's$, is given by  putting $m=5$ in (\ref{eq:10024a}) denoted by $y_0^5(x)$.
The sub-power series, having sequences $c_2$--$c_5$ including one term of $B_n's$, is given by  putting $m=3$ in (\ref{eq:10024b}) denoted by $y_1^3(x)$.
The sub-power series, having sequences $c_4$--$c_5$ including two terms of $B_n's$, is given by  putting $\tau=2$ and $m=1$ in (\ref{eq:10024c}) denoted by $y_2^1(x)$.
Taking $y_0^5(x)$, $y_1^3(x)$ and $y_2^1(x)$ into (\ref{eq:10040}),
\begin{equation}
y(x)= y_0^5(x)+y_1^3(x)+y_2^1(x) = \sum_{r=0}^{2}y_{r}^{5-2r}(x)\label{eq:10041}
\end{equation}
A sequence $c_6$ consists of zero, one, two and three terms of $B_n's$ in (\ref{eq:1007}) and (\ref{eq:1008}).
Putting $n=6$ in (\ref{eq:10023a}), a coefficient $c_6$ for zero term of $B_n's$ is denoted by $\tilde{c}(0,6)$.
Taking $n=4$  in (\ref{eq:10023b}), a coefficient $c_6$ for one term of $B_n's$ is denoted by $\tilde{c}(1,4)$.
Taking $\tau=2$ and $n=2$  in (\ref{eq:10023c}), a coefficient $c_6$ for two terms of $B_n's$ is denoted by $\tilde{c}(2,2)$.
Taking $\tau=3$ and $n=0$  in (\ref{eq:10023c}), a coefficient $c_6$ for three terms of $B_n's$ is denoted by $\tilde{c}(3,0)$.
Since the sum of $\tilde{c}(0,6)$, $\tilde{c}(1,4)$, $\tilde{c}(2,2)$ and $\tilde{c}(3,0)$ is equivalent to zero, we obtain a polynomial equation of degree 6 for the determination of an accessory parameter in $A_n$ term which is given by
\begin{equation}
0 = c_6 = \tilde{c}(0,6) + \tilde{c}(1,4) + \tilde{c}(2,2) + \tilde{c}(3,0) = \sum_{r=0}^{3}\tilde{c}(r,6-2r)\label{eq:10042}
\end{equation}

If $B_7= c_7=0$ in (\ref{eq:10025}), then gives $c_{8}=c_{9}=c_{10}=\cdots=0$. According to (\ref{eq:1009}), its power series is given by
\begin{eqnarray}
y(x)&=& \sum_{n=0}^{6} c_n x^{n+\lambda } = y_0^6(x)+y_1^4(x)+y_2^2(x)+y_3^0(x) \nonumber\\
&=& \sum_{r=0}^{3}y_{r}^{6-2r}(x) \label{eq:10043}
\end{eqnarray}
A polynomial equation of degree 7 for the determination of an accessory parameter in $A_n$ term is given by
\begin{equation}
0 = c_7 = \tilde{c}(0,7) + \tilde{c}(1,5) + \tilde{c}(2,3) + \tilde{c}(3,1) = \sum_{r=0}^{3}\tilde{c}(r,7-2r)\label{eq:10044}
\end{equation}
By repeating this process for $B_{j+1}= c_{j+1}=0$ where $j\geq 7$, we obtain the first species complete polynomial of degree $j$ and an polynomial equation of the $j+1^{th}$ order in an accessory parameter.
\qed
\end{proof}
\subsection{The second species complete polynomial using R3TRF}
As I mention on chapter 1, one of examples for the second complete polynomial is Heun's differential equation. 
For a fixed $j \in \mathbb{N}_{0}$, we need a condition such as 
\begin{equation}
\begin{cases} \alpha = -j-\lambda  \cr
\beta  = -j+1-\lambda \cr
q  = -(j+\lambda )(-\delta +1-j-\lambda +a(\gamma +\delta -1+j+\lambda )) 
\end{cases}\label{eq:10045}
\end{equation} 
Plug (\ref{eq:10045}) into (\ref{eq:1004a})-(\ref{eq:1004c}).
\begin{equation}
\begin{cases} A_n =  \frac{1+a}{a}\frac{(n-j)\left\{ n+\frac{1}{1+a}\left[ -\delta +1-j+a(\gamma +\delta -1+j+2\lambda )\right]\right\}}{(n+1+\lambda )(n+\gamma +\lambda )}  \cr
B_n = -\frac{1}{a} \frac{(n-j-1 )(n-j)}{(n+1+\lambda )(n+\gamma +\lambda )} \cr
c_1= A_0 \;c_0
\end{cases}\label{eq:10046}
\end{equation}
$B_j=B_{j+1}=A_j=0$ for a fixed $j \in \mathbb{N}_{0}$ is held by (\ref{eq:10046}), then (\ref{eq:1003}) gives successively 
$ c_{j+1}=c_{j+2}=c_{j+3}=\cdots=0 $. 
For the first species complete polynomial, there are $j+1$ values of a spectral parameter $q$ in $A_n$ term. In contrast, the second species complete polynomial only has one value of $q$ for each $j$.
\begin{definition}
For the second species complete polynomial, we need a condition which is defined by
\begin{equation}
B_{j}=B_{j+1}= A_{j}=0\hspace{1cm}\mathrm{where}\;j \in \mathbb{N}_{0}    
 \label{eq:10047}
\end{equation}
\end{definition}
\begin{theorem}
The general expression of a function $y(x)$ for the second species complete polynomial using reversible 3-term recurrence formula is given by
\begin{enumerate}
\item As $B_1=A_0=0$,
\begin{equation}
y(x) = y_{0}^{0}(x) \nonumber 
\end{equation}
\item As $B_1=B_2=A_1=0$, 
\begin{equation}
y(x)= y_{0}^{1}(x) \nonumber 
\end{equation}
\item As $B_{2N+2}=B_{2N+3}=A_{2N+2}=0$ where $N \in \mathbb{N}_{0}$,
\begin{equation}
y(x)= \sum_{r=0}^{N+1} y_{r}^{2(N+1-r)}(x) \nonumber 
\end{equation}
\item As $B_{2N+3}=B_{2N+4}=A_{2N+3}=0$ where $N \in \mathbb{N}_{0}$,
\begin{equation}
y(x)= \sum_{r=0}^{N+1} y_{r}^{2(N-r)+3}(x) \nonumber 
\end{equation}
In the above,
\begin{eqnarray}
y_0^m(x) &=& c_0 x^{\lambda} \sum_{i_0=0}^{m} \left\{ \prod _{i_1=0}^{i_0-1}A_{i_1} \right\} x^{i_0 } \nonumber\\
y_1^m(x) &=& c_0 x^{\lambda} \sum_{i_0=0}^{m}\left\{ B_{i_0+1} \prod _{i_1=0}^{i_0-1}A_{i_1}  \sum_{i_2=i_0}^{m} \left\{ \prod _{i_3=i_0}^{i_2-1}A_{i_3+2} \right\}\right\} x^{i_2+2 } \nonumber\\
y_{\tau }^m(x) &=& c_0 x^{\lambda} \sum_{i_0=0}^{m} \left\{B_{i_0+1}\prod _{i_1=0}^{i_0-1} A_{i_1} 
\prod _{k=1}^{\tau -1} \left( \sum_{i_{2k}= i_{2(k-1)}}^{m} B_{i_{2k}+(2k+1)}\prod _{i_{2k+1}=i_{2(k-1)}}^{i_{2k}-1}A_{i_{2k+1}+2k}\right) \right. \nonumber\\
&& \times \left. \sum_{i_{2\tau} = i_{2(\tau -1)}}^{m} \left( \prod _{i_{2\tau +1}=i_{2(\tau -1)}}^{i_{2\tau}-1} A_{i_{2\tau +1}+ 2\tau} \right) \right\} x^{i_{2\tau}+2\tau }\hspace{1cm}\mathrm{where}\;\tau \geq 2
\nonumber 
\end{eqnarray} 
\end{enumerate}
\end{theorem}
\begin{proof}
First if $j=0$ in (\ref{eq:10047}),  $B_{1}=A_0=0$ and then (\ref{eq:1003}) gives successively 
$ c_{1}=c_{2}=c_{3}=\cdots=0 $. According to (\ref{eq:1009}), its power series is given by
\begin{equation}
y(x)= \sum_{n=0}^{0} c_n x^{n+\lambda } = c_0 x^{\lambda } = y_{0}^{0}(x)\label{eq:10048}
\end{equation}
The sub-power series $y_{0}^{0}(x)$ in (\ref{eq:10048}) is obtain by putting $m=0$ in (\ref{eq:10024a}).

If $j=1$ in (\ref{eq:10047}), $B_{1}=B_{2}=A_1=0$ and then (\ref{eq:1003}) gives successively 
$ c_{2}=c_{3}=c_{4}=\cdots=0 $. According to (\ref{eq:1009}), its power series is given by
\begin{equation}
y(x)= \sum_{n=0}^{1} c_n x^{n+\lambda } = (c_0 +c_1 x) x^{\lambda } \label{eq:10049}
\end{equation}
Observe sequences $c_0$ and $c_1$ in (\ref{eq:1007}) and (\ref{eq:1008}). 
The sub-power series that has sequences $c_0$ and $c_1$ is given by  putting $m=1$ in (\ref{eq:10024a}). Take the new (\ref{eq:10024a}) into (\ref{eq:10049}).
\begin{equation}
y(x)= y_{0}^{1}(x) = c_0 \left( 1+A_0 x\right)x^{\lambda } \label{eq:10050}
\end{equation}

If $j=2$ in (\ref{eq:10047}), $B_{2}=B_{3}=A_2=0$ and then (\ref{eq:1003}) gives successively 
$ c_{3}=c_{4}=c_{5}=\cdots=0 $. According to (\ref{eq:1009}), its power series is given by
\begin{equation}
y(x)= \sum_{n=0}^{2} c_n x^{n+\lambda } = (c_0 +c_1 x+c_2 x^2) x^{\lambda } \label{eq:10051}
\end{equation}
Observe sequences $c_0$--$c_2$ in (\ref{eq:1007}) and (\ref{eq:1008}). 
The sub-power series, having sequences $c_0$--$c_2$ including zero term of $B_n's$, is given by  putting $m=2$ in (\ref{eq:10024a}) denoted by $y_0^2(x)$.
The sub-power series, having a sequence $c_2$ including one term of $B_n's$, is given by  putting $m=0$ in (\ref{eq:10024b}) denoted by $y_1^0(x)$.
Taking $y_0^2(x)$ and $y_1^0(x)$ into (\ref{eq:10051}),
\begin{equation}
y(x)= y_{0}^{2}(x) + y_{1}^{0}(x) = c_0x^{\lambda } \left( 1+A_0 x + \left( A_{0,1} +B_1\right) x^2\right) \label{eq:10052}
\end{equation} 

If $j=3$ in (\ref{eq:10047}), $B_{3}=B_{4}=A_3=0$ and then (\ref{eq:1003}) gives successively 
$ c_{4}=c_{5}=c_{6}=\cdots=0 $. According to (\ref{eq:1009}), its power series is given by
\begin{equation}
y(x)= \sum_{n=0}^{3} c_n x^{n+\lambda } = (c_0 +c_1 x+c_2 x^2+c_3 x^3) x^{\lambda } \label{eq:10053}
\end{equation}
Observe sequences $c_0$--$c_3$ in (\ref{eq:1007}) and (\ref{eq:1008}). 
The sub-power series, having sequences $c_0$--$c_3$ including zero term of $B_n's$, is given by  putting $m=3$ in (\ref{eq:10024a}) denoted by $y_0^3(x)$.
The sub-power series, having sequences $c_2$ and $c_3$ including one term of $B_n's$, is given by  putting $m=1$ in (\ref{eq:10024b}) denoted by $y_1^1(x)$.
Taking $y_0^3(x)$ and $y_1^1(x)$ into (\ref{eq:10053}),
\begin{equation}
y(x)= y_{0}^{3}(x)+y_{1}^{1}(x) \label{eq:10054}
\end{equation}

If $j=4$ in (\ref{eq:10047}), $B_{4}=B_{5}=A_4=0$ and then (\ref{eq:1003}) gives successively 
$ c_{5}=c_{6}=c_{7}=\cdots=0 $. According to (\ref{eq:1009}), its power series is given by
\begin{equation}
y(x)= \sum_{n=0}^{4} c_n x^{n+\lambda } = (c_0 +c_1 x+c_2 x^2+c_3 x^3+c_4 x^4) x^{\lambda } \label{eq:10055}
\end{equation}
Observe sequences $c_0$--$c_4$ in (\ref{eq:1007}) and (\ref{eq:1008}). 
The sub-power series, having sequences  $c_0$--$c_4$ including zero term of $B_n's$, is given by  putting $m=4$ in (\ref{eq:10024a}) denoted by $y_0^4(x)$.
The sub-power series, having sequences  $c_2$--$c_4$ including one term of $B_n's$, is given by  putting $m=2$ in (\ref{eq:10024b}) denoted by $y_1^2(x)$.
The sub-power series, having sequences $c_4$ including two terms of $B_n's$, is given by  putting $\tau=2$ and $m=0$ in (\ref{eq:10024c}) denoted by $y_2^0(x)$.
Taking $y_0^4(x)$, $y_1^2(x)$ and $y_2^0(x)$ into (\ref{eq:10055}),
\begin{equation}
y(x)= y_0^4(x)+y_1^2(x)+y_2^0(x) \label{eq:10056}
\end{equation}

If $j=5$ in (\ref{eq:10047}), $B_{5}=B_{6}=A_5=0$ and then (\ref{eq:1003}) gives successively 
$ c_{6}=c_{7}=c_{8}=\cdots=0 $. According to (\ref{eq:1009}), its power series is given by
\begin{equation}
y(x) = \sum_{n=0}^{5} c_n x^{n+\lambda } = y_0^5(x)+y_1^3(x)+y_2^1(x)= \sum_{r=0}^{2}y_{r}^{5-2r}(x) \label{eq:10057}
\end{equation}

If $j=6$ in (\ref{eq:10047}), $B_{6}=B_{7}=A_6=0$ and then (\ref{eq:1003}) gives successively 
$ c_{7}=c_{8}=c_{9}=\cdots=0 $. According to (\ref{eq:1009}), its power series is given by 
\begin{equation}
y(x) = \sum_{n=0}^{6} c_n x^{n+\lambda } = y_0^6(x)+y_1^4(x)+y_2^2(x)+y_3^0(x) = \sum_{r=0}^{3}y_{r}^{6-2r}(x) \label{eq:10058}
\end{equation}
By repeating this process for $B_{j}=B_{j+1}= A_{j}=0$ where $j\geq 7$, we obtain the second species complete polynomial of degree $j$.
\qed
\end{proof}
\section{\label{sec:level6}Summary}

There are two types of power series of 2-term recursive relation in a ODE which are an infinite series and a polynomial. 
With my definition, the power series of the 3-term recurrence relation in a ODE has $2^{3-1}$ types. Complete polynomials, one of polynomials in 3-recursion relation, has two different types which are the first and second species complete polynomials. The former is applicable since a parameter of a numerator in $B_n$ term and a (spectral) parameter of a numerator in $A_n$ term are fixed constants: its polynomial has multi-valued roots of an eigenvalue (spectral parameter) in $A_n$ term. The latter is utilized as two parameters of a numerator in $B_n$ term and a parameter of a numerator in $A_n$ term are fixed constants: its polynomial has only one valued root of its eigenvalue in $A_n$ term.

In this chapter, by substituting a power series with unknown coefficients into a linear ODE, I construct general expressions for the first species complete polynomial with a condition $B_{j+1}= c_{j+1}=0 $ where $j \in \mathbb{N}_{0}$ and second species complete polynomial with $B_{j}=B_{j+1}= A_{j}=0$ by allowing $B_n$ as the leading term in each sub-power series of the general power series: I observe the term of sequence $c_n$ which includes zero term of $B_n's$, one term of $B_n's$, two terms of $B_n's$, three terms of $B_n's$, etc. 

The general summation formulas of the first \& second species complete polynomials on chapter 1 and this chapter are equivalent to each other. The former is that $A_n$ is the leading term in each sub-power series of the general summation formulas. And the latter is that $B_n$ is the leading term in each sub-power series of the general power series solutions.

The general power series solutions of complete polynomials using 3TRF in chapter 1 consist of the sum of two sub-power series. And the Frobenius solutions of complete polynomials using R3TRF in this paper is only composed of one sub-formal series.
Because of this simple mathematical structure, the latter is more efficient to be applied into special functions even if numerical values of these two type polynomials are equivalent to each other.
\addcontentsline{toc}{section}{Bibliography}
\bibliographystyle{model1a-num-names}
\bibliography{<your-bib-database>}

%
\chapter{Complete polynomials of Heun equation using three-term recurrence formula}
\chaptermark{Complete polynomials of Heun equation using 3TRF} 

For an infinite series and a polynomial which makes $B_n$ term terminated, power series expansions and integral representations of Heun equation are constructed analytically by applying three term recurrence formula (3TRF). \cite{4Chou2012c,4Chou2012d}

In chapter 2 of Ref.\cite{4Choun2013}, I apply reversible three term recurrence formula (R3TRF) to power series expansions in closed forms and  integral forms of Heun equation for an infinite series and a polynomial which makes $A_n$ term terminated.

In this chapter I apply mathematical formulas of complete polynomials using 3TRF to the Frobenius solutions in closed forms of Heun equation for a polynomial which makes $A_n$ and $B_n$ terms terminated in which the coefficients are given explicitly.

%
Nine examples of 192 local solutions of the Heun equation (Maier, 2007) are provided in the appendix.  For each example, I construct power series expansions of Heun equation for complete polynomials using 3TRF.
\section{Introduction}
The Heun function is considered as the mother of all well-known special functions such as: Spheroidal Wave, Lame, Mathieu, hypergeometric type functions, etc. The power series of Heun's equation has a 3-term recursive relation between successive coefficients. According to Karl Heun \cite{4Heun1889,4Ronv1995}, Heun equation is a second-order linear ODE which has four regular singular points. Heun's equation has four different types of a confluent form such as Confluent Heun, Doubly-Confluent Heun, Biconfluent Heun and Triconfluent Heun equations.
Even if Heun's equation was discovered in 1889, it has not utilized in mathematics and physics areas until recently  because of its three different successive coefficients in its formal series. Power series solutions for the 3-term recurrence relation in a linear ODE are really hard to be derived analytically in which the coefficients are given in an explicit manner, and even their numerical calculation proved. \cite{4Arsc1983} However, since around 1990,  Heun's equation has started to arise again in fields of modern physics such as quantum theories, general relativity, string theory, etc.  Its equation has appeared in various areas such as Kerr-de Sitter black holes problems \cite{4Teuk1973,4Leav1985,4Bati2006,4Bati2007,4Bati2010,4Take2006}, Calogero-Moser-Sutherland systems 
\cite{4Take2003,4Take2004,4Take2004a,4Take2005}, quantum inozemtsev model \cite{4Inoz1989}, etc. 

Heun ordinary differential equation is given by \cite{4Heun1889,4Ronv1995,4Slavy2000,4Erde1955}
\begin{equation}
\frac{d^2{y}}{d{x}^2} + \left(\frac{\gamma }{x} +\frac{\delta }{x-1} + \frac{\epsilon }{x-a}\right) \frac{d{y}}{d{x}} +  \frac{\alpha \beta x-q}{x(x-1)(x-a)} y = 0 \label{eq:4001}
\end{equation}
With the condition $\epsilon = \alpha +\beta -\gamma -\delta +1$. The parameters play different roles: $a \ne 0 $ is the singularity parameter, $\alpha $, $\beta $, $\gamma $, $\delta $, $\epsilon $ are exponent parameters, $q$ is the accessory parameter. Also, $\alpha $ and $\beta $ are identical to each other. The total number of free parameters is six. It has four regular singular points which are 0, 1, $a$ and $\infty $ with exponents $\{ 0, 1-\gamma \}$, $\{ 0, 1-\delta \}$, $\{ 0, 1-\epsilon \}$ and $\{ \alpha, \beta \}$. 

We assume the solution takes the form 
\begin{equation}
y(x)= \sum_{n=0}^{\infty } c_n x^{n+\lambda } \label{eq:4002}
\end{equation}
where $\lambda $ is an indicial root. 
Substituting (\ref{eq:4002}) into (\ref{eq:4001}) gives for the coefficients $c_n$ the recurrence relations
\begin{equation}
c_{n+1}=A_n \;c_n +B_n \;c_{n-1} \hspace{1cm};n\geq 1 \label{eq:4003}
\end{equation}
where
\begin{subequations}
\begin{eqnarray}
A_n &=& \frac{(n+\lambda )(n-1+\gamma +\epsilon +\lambda + a(n-1+\gamma +\lambda +\delta ))+q}{a(n+1+\lambda )(n+\gamma +\lambda )}\nonumber\\
&=& \frac{(n+\lambda )(n+\alpha +\beta -\delta +\lambda +a(n+\delta +\gamma -1+\lambda ))+q}{a(n+1+\lambda )(n+\gamma +\lambda )} \label{eq:4004a}
\end{eqnarray}
\begin{equation}
B_n = -\frac{(n-1+\lambda )(n+\gamma +\delta +\epsilon -2+\lambda )+\alpha \beta }{a(n+1+\lambda )(n+\gamma +\lambda )}=- \frac{(n-1+\lambda +\alpha )(n-1+\lambda +\beta )}{a(n+1+\lambda )(n+\gamma +\lambda )} \label{eq:4004b}
\end{equation}
\begin{equation}
c_1= A_0 \;c_0 \label{eq:4004c}
\end{equation}
\end{subequations}
We have two indicial roots which are $\lambda = 0$ and $ 1-\gamma $

\section{Power series}
With my definition, power series solutions in the three term recurrence relation of a linear ODE have an infinite series and three types of polynomials: (1) a polynomial which makes $B_n$ term terminated; $A_n$ term is not terminated, (2) a polynomial which makes $A_n$ term terminated; $B_n$ term is not terminated, (3) a polynomial which makes $A_n$ and $B_n$ terms terminated at the same time.

In Ref.\cite{4Chou2012c,4Chou2012d} I show how to obtain power series expansions in closed forms of Heun equation around $x=0$ and its combined definite and contour integrals for an infinite series and a polynomial of type 1 by applying 3TRF. 
The sequence $c_n$ combines into combinations of $A_n$ and $B_n$ terms in (\ref{eq:4003}). This is done by letting $A_n$ in the sequence $c_n$ is the leading term in the analytic function $y(x)$: I observe the term of sequence $c_n$ which includes zero term of $A_n's$, one term of $A_n's$, two terms of $A_n's$, three terms of $A_n's$, etc. 
For a polynomial of type 1, I treat $\gamma $, $\delta $ and $q$ as free variables and fixed values of $\alpha $ and/or $\beta $.

In chapter 2 of Ref.\cite{4Choun2013} I show how to describe the Frobenius solutions in closed forms of Heun equation around $x=0$ and its integral representations for an infinite series and a polynomial of type 2 by applying R3TRF. 
This is done by letting $B_n$ in the sequence $c_n$ is the leading term in the analytic function $y(x)$: I observe the term of sequence $c_n$ which includes zero term of $B_n's$, one term of $B_n's$, two terms of $B_n's$, three terms of $B_n's$, etc. 
For a polynomial of type 2, I treat $\alpha $, $\beta $, $\gamma $ and $\delta $ as free variables and a fixed value of $q$. 

In general, Heun (spectral) polynomial has been defined as a type 3 polynomial where $A_n$ and $B_n$ terms terminated.
According to Shapiro \textit {et al.}, \cite{4Shap2010,4Shap2011,4Shap2012}
 Heun's equation has the form
\begin{equation}
\left\{ Q(z)\frac{d^2}{d{z}^2} +  P(z)\frac{d}{d{z}} +  V(z) \right\} S(z)= 0 \nonumber
\end{equation}
where $Q(z)$ is a cubic complex polynomial, $P(z)$ is a polynomial of degree at most 2 and $V(z)$ is at most linear.  
They investigate spectral polynomials of the Heun equation in the case of real roots of these polynomials
and asymptotic root distribution when complex roots are present. The Heun (spectral) polynomials has a polynomial solution $y(x)$ for a given degree $n$. They study the union of the roots of the latter set of $V(z)$'s when $n\rightarrow \infty$. 
 
Mart\'{i}nez-Finkelshtein and Rakhmanov \cite{4Mart2011} investigate the asymptotic zero distribution of Heine-Stieltjes polynomials with complex polynomial coefficients --Heun polynomials is the special case of Stieltjes polynomials. 
Kalnins \textit {et al.} \cite{4Kaln1989,4Kaln1990} describe Heun polynomials in terms of Jacobi polynomials. 
This is done by utilizing group theory and its connection with the method of separation of variables applied to the Laplace-Beltrami eigenvalue equation on the $n$-sphere.

With an accessory parameter belonging to a hyperelliptic spectral curve, Smirnov \cite{4Smir2001} derives an integral representation of finite-gap solutions of Heun's equation (at a close connection with Treibich-Verdier equation) with nonnegative exponent parameters from the Liouville formula. He shows that power series solutions of Heun's equation are analyzed in terms of Heun polynomials since an accessory parameter is equal to branching points of the hyperelliptic curve with any integer exponent parameters.

 Heun polynomial comes from the Heun differential equation around $x=0$ in (\ref{eq:4001}) that has a fixed value of $\alpha $ or $\beta $, just as it has a fixed value of $q$. In the 3-term recurrence relation in a linear ODE, a polynomial of type 3 I categorize as complete polynomials. 
Complete polynomials can be divided into two different types which are the first species complete polynomial and the second species complete polynomial. 

If a parameter of a numerator in $B_n$ term and a (spectral) parameter of a numerator in $A_n$ term are fixed constants, we should apply the first species complete polynomial. And if two parameters of a numerator in $B_n$ term and a parameter of a numerator in $A_n$ term are fixed constants, the second species complete polynomial is able to be utilized. 
The former has multi-valued roots of an eigenvalue (spectral parameter) in $A_n$ term, but the latter has only one valued root of each eigenvalue in $A_n$ term.
For the first species complete polynomial of Heun's equation around  $x=0$, I treat an exponent parameter $\alpha $ (or $\beta $) and an accessory parameter $q$ are fixed constants. And as parameters $\alpha $, $\beta $ and $q$ are fixed constants, the second species complete polynomial must to be applied.  

In chapter 1 I construct the mathematical formulas of complete polynomials for the first and second species, by allowing $A_n$ as the leading term in each sub-power series of the general power series $y(x)$, as ``complete polynomials using 3-term recurrence formula (3TRF)''
In this chapter, by applying complete polynomials using 3TRF, I construct the power series expansion in closed forms of Heun equation around $x=0$ for a polynomial which makes $A_n$ and $B_n$ terms terminated.

\subsection{The first species complete polynomial using 3TRF}
For the first species complete polynomial, we need a condition which is given by
\begin{equation}
B_{j+1}= c_{j+1}=0\hspace{1cm}\mathrm{where}\;j\in \mathbb{N}_{0}  
 \label{eq:4005}
\end{equation}
(\ref{eq:4005}) gives successively $c_{j+2}=c_{j+3}=c_{j+4}=\cdots=0$. And $c_{j+1}=0$ is defined by a polynomial equation of degree $j+1$ for the determination of an accessory parameter in $A_n$ term. 
\begin{theorem}
In chapter 1, the general expression of a function $y(x)$ for the first species complete polynomial using 3-term recurrence formula and its algebraic equation for the determination of an accessory parameter in $A_n$ term are given by
\begin{enumerate} 
\item As $B_1=0$,
\begin{equation}
0 =\bar{c}(1,0) \label{eq:4006a}
\end{equation}
\begin{equation}
y(x) = y_{0}^{0}(x) \label{eq:4006b}
\end{equation}
\item As $B_{2N+2}=0$ where $N \in \mathbb{N}_{0}$,
\begin{equation}
0  = \sum_{r=0}^{N+1}\bar{c}\left( 2r, N+1-r\right) \label{eq:4007a}
\end{equation}
\begin{equation}
y(x)= \sum_{r=0}^{N} y_{2r}^{N-r}(x)+ \sum_{r=0}^{N} y_{2r+1}^{N-r}(x)  \label{eq:4007b}
\end{equation}
\item As $B_{2N+3}=0$ where $N \in \mathbb{N}_{0}$,
\begin{equation}
0  = \sum_{r=0}^{N+1}\bar{c}\left( 2r+1, N+1-r\right) \label{eq:4008a}
\end{equation}
\begin{equation}
y(x)= \sum_{r=0}^{N+1} y_{2r}^{N+1-r}(x)+ \sum_{r=0}^{N} y_{2r+1}^{N-r}(x)  \label{eq:4008b}
\end{equation}
In the above,
\begin{eqnarray}
\bar{c}(0,n)  &=& \prod _{i_0=0}^{n-1}B_{2i_0+1} \label{eq:4009a}\\
\bar{c}(1,n) &=&  \sum_{i_0=0}^{n} \left\{ A_{2i_0} \prod _{i_1=0}^{i_0-1}B_{2i_1+1} \prod _{i_2=i_0}^{n-1}B_{2i_2+2} \right\} 
\label{eq:4009b}\\
\bar{c}(\tau ,n) &=& \sum_{i_0=0}^{n} \left\{A_{2i_0}\prod _{i_1=0}^{i_0-1} B_{2i_1+1} 
\prod _{k=1}^{\tau -1} \left( \sum_{i_{2k}= i_{2(k-1)}}^{n} A_{2i_{2k}+k}\prod _{i_{2k+1}=i_{2(k-1)}}^{i_{2k}-1}B_{2i_{2k+1}+(k+1)}\right) \right. \nonumber\\
&&\times \left. \prod _{i_{2\tau}=i_{2(\tau -1)}}^{n-1} B_{2i_{2\tau }+(\tau +1)} \right\} 
\hspace{1cm}\label{eq:4009c} 
\end{eqnarray}
and
\begin{eqnarray}
y_0^m(x) &=& c_0 x^{\lambda } \sum_{i_0=0}^{m} \left\{ \prod _{i_1=0}^{i_0-1}B_{2i_1+1} \right\} x^{2i_0 } \label{eq:40010a}\\
y_1^m(x) &=& c_0 x^{\lambda } \sum_{i_0=0}^{m}\left\{ A_{2i_0} \prod _{i_1=0}^{i_0-1}B_{2i_1+1}  \sum_{i_2=i_0}^{m} \left\{ \prod _{i_3=i_0}^{i_2-1}B_{2i_3+2} \right\}\right\} x^{2i_2+1 } \label{eq:40010b}\\
y_{\tau }^m(x) &=& c_0 x^{\lambda } \sum_{i_0=0}^{m} \left\{A_{2i_0}\prod _{i_1=0}^{i_0-1} B_{2i_1+1} 
\prod _{k=1}^{\tau -1} \left( \sum_{i_{2k}= i_{2(k-1)}}^{m} A_{2i_{2k}+k}\prod _{i_{2k+1}=i_{2(k-1)}}^{i_{2k}-1}B_{2i_{2k+1}+(k+1)}\right) \right. \nonumber\\
&& \times \left. \sum_{i_{2\tau} = i_{2(\tau -1)}}^{m} \left( \prod _{i_{2\tau +1}=i_{2(\tau -1)}}^{i_{2\tau}-1} B_{2i_{2\tau +1}+(\tau +1)} \right) \right\} x^{2i_{2\tau}+\tau }\hspace{1cm}\mathrm{where}\;\tau \geq 2
\label{eq:40010c} 
\end{eqnarray}
\end{enumerate}
\end{theorem}
Put $n= j+1$ in (\ref{eq:4004b}) and use the condition $B_{j+1}=0$ for $\alpha $.  
\begin{equation}
\alpha = -j-\lambda 
\label{eq:40011}
\end{equation}
Take (\ref{eq:40011}) into (\ref{eq:4004a}) and (\ref{eq:4004b}).
\begin{subequations}
\begin{equation}
A_n = \frac{(n+\lambda )(n-j +\beta -\delta +a(n+\delta +\gamma -1+\lambda ))+q}{a(n+1+\lambda )(n+\gamma +\lambda )} \label{eq:40012a}
\end{equation}
\begin{equation}
B_n = - \frac{(n-1-j)(n-1+\lambda +\beta )}{a(n+1+\lambda )(n+\gamma +\lambda )} \label{eq:40012b}
\end{equation}
\end{subequations}
Now the condition $c_{j+1}=0$ is clearly an algebraic equation in $q$ of degree $j+1$ and thus has $j+1$ zeros denoted them by $q_j^m$ eigenvalues where $m = 0,1,2, \cdots, j$. They can be arranged in the following order: $q_j^0 < q_j^1 < q_j^2 < \cdots < q_j^j$.
 
Substitute (\ref{eq:40012a}) and (\ref{eq:40012b}) into (\ref{eq:4009a})--(\ref{eq:40010c}).

As $B_{1}= c_{1}=0$, take the new (\ref{eq:4009b}) into (\ref{eq:4006a}) putting $j=0$. Substitute the new (\ref{eq:40010a}) into (\ref{eq:4006b}) putting $j=0$. 

As $B_{2N+2}= c_{2N+2}=0$, take the new (\ref{eq:4009a})--(\ref{eq:4009c}) into (\ref{eq:4007a}) putting $j=2N+1$. Substitute the new 
(\ref{eq:40010a})--(\ref{eq:40010c}) into (\ref{eq:4007b}) putting $j=2N+1$ and $q=q_{2N+1}^m$.

As $B_{2N+3}= c_{2N+3}=0$, take the new (\ref{eq:4009a})--(\ref{eq:4009c}) into (\ref{eq:4008a}) putting $j=2N+2$. Substitute the new 
(\ref{eq:40010a})--(\ref{eq:40010c}) into (\ref{eq:4008b}) putting $j=2N+2$ and $q=q_{2N+2}^m$.

After the replacement process, the general expression of power series of Heun equation about $x=0$ for the first species complete polynomial using 3-term recurrence formula and its algebraic equation for the determination of an accessory parameter $q$ are given by
\begin{enumerate} 
\item As $\alpha =-\lambda $,

An algebraic equation of degree 1 for the determination of $q$ is given by
\begin{equation}
0= \bar{c}(1,0;0,q)= q + \lambda (\beta -\delta +a(\gamma+ \delta -1+\lambda )) \label{eq:40013a}
\end{equation}
The eigenvalue of $q$ is written by $q_0^0$. Its eigenfunction is given by
\begin{equation}
y(x) = y_0^0\left( 0,q_0^0;x\right)= c_0 x^{\lambda } \label{eq:40013b}
\end{equation}
\item As $\alpha =-2N-1-\lambda $ where $N \in \mathbb{N}_{0}$,

An algebraic equation of degree $2N+2$ for the determination of $q$ is given by
\begin{equation}
0 = \sum_{r=0}^{N+1}\bar{c}\left( 2r, N+1-r; 2N+1,q\right)  \label{eq:40014a}
\end{equation}
The eigenvalue of $q$ is written by $q_{2N+1}^m$ where $m = 0,1,2,\cdots,2N+1 $; $q_{2N+1}^0 < q_{2N+1}^1 < \cdots < q_{2N+1}^{2N+1}$. Its eigenfunction is given by 
\begin{equation} 
y(x) = \sum_{r=0}^{N} y_{2r}^{N-r}\left( 2N+1,q_{2N+1}^m;x\right)+ \sum_{r=0}^{N} y_{2r+1}^{N-r}\left( 2N+1,q_{2N+1}^m;x\right)
\label{eq:40014b} 
\end{equation}
\item As $\alpha =-2N-2-\lambda $ where $N \in \mathbb{N}_{0}$,

An algebraic equation of degree $2N+3$ for the determination of $q$ is given by
\begin{equation}  
0 = \sum_{r=0}^{N+1}\bar{c}\left( 2r+1, N+1-r; 2N+2,q\right) \label{eq:40015a}
\end{equation}
The eigenvalue of $q$ is written by $q_{2N+2}^m$ where $m = 0,1,2,\cdots,2N+2 $; $q_{2N+2}^0 < q_{2N+2}^1 < \cdots < q_{2N+2}^{2N+2}$. Its eigenfunction is given by
\begin{equation} 
y(x) =  \sum_{r=0}^{N+1} y_{2r}^{N+1-r}\left( 2N+2,q_{2N+2}^m;x\right) + \sum_{r=0}^{N} y_{2r+1}^{N-r}\left( 2N+2,q_{2N+2}^m;x\right) \label{eq:40015b}
\end{equation}
In the above,
\begin{eqnarray}
\bar{c}(0,n;j,q)  &=& \frac{\left( -\frac{j}{2}\right)_{n} \left( \frac{\beta }{2}+ \frac{\lambda }{2}\right)_{n}}{\left( 1+\frac{\lambda }{2}\right)_{n} \left( \frac{\gamma }{2}+ \frac{1}{2}+ \frac{\lambda }{2}\right)_{n}} \left(-\frac{1}{a} \right)^{n}\label{eq:40016a}\\
\bar{c}(1,n;j,q) &=& \left( \frac{1+a}{a}\right) \sum_{i_0=0}^{n}\frac{\left( i_0+\frac{\lambda }{2}\right)\left( i_0 +\Gamma _0^{(F)} \right) +\frac{q}{4(1+a)}}{\left( i_0+\frac{1}{2}+\frac{\lambda }{2}\right) \left( i_0+\frac{\gamma }{2}+\frac{\lambda }{2}\right)} \frac{\left( -\frac{j}{2}\right)_{i_0} \left( \frac{\beta }{2}+ \frac{\lambda }{2}\right)_{i_0}}{\left( 1+\frac{\lambda }{2}\right)_{i_0} \left( \frac{\gamma }{2}+ \frac{1}{2}+ \frac{\lambda }{2}\right)_{i_0}} \nonumber\\
&&\times  \frac{\left( -\frac{j}{2}+\frac{1}{2} \right)_{n} \left( \frac{\beta }{2}+ \frac{1 }{2}+ \frac{\lambda }{2}\right)_{n}\left( \frac{3}{2}+\frac{\lambda }{2}\right)_{i_0} \left( \frac{\gamma }{2}+1+ \frac{\lambda }{2}\right)_{i_0}}{\left( -\frac{j}{2}+\frac{1}{2}\right)_{i_0} \left( \frac{\beta }{2}+ \frac{1}{2}+ \frac{\lambda }{2}\right)_{i_0}\left( \frac{3}{2}+\frac{\lambda }{2}\right)_{n} \left( \frac{\gamma }{2}+1+ \frac{\lambda }{2}\right)_{n}} \left(-\frac{1}{a} \right)^{n }  
\label{eq:40016b}
\end{eqnarray}
\begin{eqnarray}
\bar{c}(\tau ,n;j,q) &=& \left( \frac{1+a}{a}\right)^{\tau} \sum_{i_0=0}^{n}\frac{\left( i_0+\frac{\lambda }{2}\right)\left( i_0 +\Gamma _0^{(F)} \right) +\frac{q}{4(1+a)}}{\left( i_0+\frac{1}{2}+\frac{\lambda }{2}\right) \left( i_0+\frac{\gamma }{2}+\frac{\lambda }{2}\right)} \frac{\left( -\frac{j}{2}\right)_{i_0} \left( \frac{\beta }{2}+ \frac{\lambda }{2}\right)_{i_0}}{\left( 1+\frac{\lambda }{2}\right)_{i_0} \left( \frac{\gamma }{2}+ \frac{1}{2}+ \frac{\lambda }{2}\right)_{i_0}} \nonumber\\
&&\times \prod_{k=1}^{\tau -1} \left( \sum_{i_k = i_{k-1}}^{n} \frac{\left( i_k+ \frac{k}{2}+\frac{\lambda }{2}\right)\left(i_k +\Gamma _k^{(F)} \right)+\frac{q}{4(1+a)}}{\left( i_k+\frac{k}{2}+\frac{1}{2}+\frac{\lambda }{2}\right) \left( i_k+\frac{k}{2}+\frac{\gamma }{2}+\frac{\lambda }{2}\right)} \right. \nonumber\\
&&\times \left. \frac{\left( -\frac{j}{2}+\frac{k}{2}\right)_{i_k} \left( \frac{\beta }{2}+ \frac{k}{2}+ \frac{\lambda }{2}\right)_{i_k}\left( 1+ \frac{k}{2}+\frac{\lambda }{2}\right)_{i_{k-1}} \left( \frac{\gamma }{2}+ \frac{1}{2}+ \frac{k}{2}+ \frac{\lambda }{2}\right)_{i_{k-1}}}{\left( -\frac{j}{2}+\frac{k}{2}\right)_{i_{k-1}} \left( \frac{\beta }{2}+ \frac{k}{2}+ \frac{\lambda }{2}\right)_{i_{k-1}}\left( 1+ \frac{k}{2}+\frac{\lambda }{2}\right)_{i_k} \left( \frac{\gamma }{2}+ \frac{1}{2}+ \frac{k}{2}+ \frac{\lambda }{2}\right)_{i_k}} \right) \nonumber\\
&&\times \frac{\left( -\frac{j}{2}+\frac{\tau }{2}\right)_{n} \left( \frac{\beta }{2}+ \frac{\tau }{2}+ \frac{\lambda }{2}\right)_{n}\left( 1+\frac{\tau }{2}+\frac{\lambda }{2}\right)_{i_{\tau -1}} \left( \frac{\gamma }{2}+ \frac{1}{2}+\frac{\tau }{2}+ \frac{\lambda }{2}\right)_{i_{\tau -1}}}{\left( -\frac{j}{2}+\frac{\tau }{2}\right)_{i_{\tau -1}} \left( \frac{\beta }{2}+\frac{\tau }{2}+ \frac{\lambda }{2}\right)_{i_{\tau -1}}\left( 1+\frac{\tau }{2}+\frac{\lambda }{2}\right)_{n} \left( \frac{\gamma }{2}+ \frac{1}{2}+\frac{\tau }{2}+ \frac{\lambda }{2}\right)_{n}} \left(-\frac{1}{a} \right)^{n } \hspace{1.5cm}\label{eq:40016c} 
\end{eqnarray}
\begin{eqnarray}
y_0^m(j,q;x) &=& c_0 x^{\lambda }  \sum_{i_0=0}^{m} \frac{\left( -\frac{j}{2}\right)_{i_0} \left( \frac{\beta }{2}+ \frac{\lambda }{2}\right)_{i_0}}{\left( 1+\frac{\lambda }{2}\right)_{i_0} \left( \frac{\gamma }{2}+ \frac{1}{2}+ \frac{\lambda }{2}\right)_{i_0}} z^{i_0} \label{eq:40017a}\\
y_1^m(j,q;x) &=& c_0 x^{\lambda } \left\{\sum_{i_0=0}^{m} \frac{\left( i_0+\frac{\lambda }{2}\right)\left(i_0 +\Gamma _0^{(F)} \right)+\frac{q}{4(1+a)}}{\left( i_0+\frac{1}{2}+\frac{\lambda }{2}\right) \left( i_0+\frac{\gamma }{2}+\frac{\lambda }{2}\right)} \frac{\left( -\frac{j}{2}\right)_{i_0} \left( \frac{\beta }{2}+ \frac{\lambda }{2}\right)_{i_0}}{\left( 1+\frac{\lambda }{2}\right)_{i_0} \left( \frac{\gamma }{2}+ \frac{1}{2}+ \frac{\lambda }{2}\right)_{i_0}} \right. \nonumber\\
&&\times  \left. \sum_{i_1 = i_0}^{m} \frac{\left( -\frac{j}{2}+\frac{1}{2} \right)_{i_1} \left( \frac{\beta }{2}+ \frac{1}{2}+ \frac{\lambda }{2}\right)_{i_1}\left( \frac{3}{2}+\frac{\lambda }{2}\right)_{i_0} \left( \frac{\gamma }{2}+1+ \frac{\lambda }{2}\right)_{i_0}}{\left(-\frac{j}{2}+\frac{1}{2} \right)_{i_0} \left( \frac{\beta }{2}+ \frac{1}{2}+ \frac{\lambda }{2}\right)_{i_0}\left( \frac{3}{2}+\frac{\lambda }{2}\right)_{i_1} \left( \frac{\gamma }{2} +1+ \frac{\lambda }{2}\right)_{i_1}} z^{i_1}\right\} \eta \hspace{1.5cm}
\label{eq:40017b}\\
y_{\tau }^m(j,q;x) &=& c_0 x^{\lambda } \left\{ \sum_{i_0=0}^{m} \frac{\left( i_0+\frac{\lambda }{2}\right)\left(i_0 +\Gamma _0^{(F)} \right)+\frac{q}{4(1+a)}}{\left( i_0+\frac{1}{2}+\frac{\lambda }{2}\right) \left( i_0+\frac{\gamma }{2}+\frac{\lambda }{2}\right)} \frac{\left( -\frac{j}{2}\right)_{i_0} \left( \frac{\beta }{2}+ \frac{\lambda }{2}\right)_{i_0}}{\left( 1+\frac{\lambda }{2}\right)_{i_0} \left( \frac{\gamma }{2}+ \frac{1}{2}+ \frac{\lambda }{2}\right)_{i_0}} \right.\nonumber\\
&&\times \prod_{k=1}^{\tau -1} \left( \sum_{i_k = i_{k-1}}^{m} \frac{\left( i_k+ \frac{k}{2}+\frac{\lambda }{2}\right)\left(i_k +\Gamma _k^{(F)} \right)+\frac{q}{4(1+a)}}{\left( i_k+\frac{k}{2}+\frac{1}{2}+\frac{\lambda }{2}\right) \left( i_k+\frac{k}{2}+\frac{\gamma }{2}+\frac{\lambda }{2}\right)} \right. \nonumber\\
&&\times \left. \frac{\left( -\frac{j}{2}+\frac{k}{2}\right)_{i_k} \left( \frac{\beta }{2}+ \frac{k}{2}+ \frac{\lambda }{2}\right)_{i_k}\left( 1+ \frac{k}{2}+\frac{\lambda }{2}\right)_{i_{k-1}} \left( \frac{\gamma }{2}+ \frac{1}{2}+ \frac{k}{2}+ \frac{\lambda }{2}\right)_{i_{k-1}}}{\left( -\frac{j}{2}+\frac{k}{2}\right)_{i_{k-1}} \left( \frac{\beta }{2}+ \frac{k}{2}+ \frac{\lambda }{2}\right)_{i_{k-1}}\left( 1+ \frac{k}{2}+\frac{\lambda }{2}\right)_{i_k} \left( \frac{\gamma }{2}+ \frac{1}{2}+ \frac{k}{2}+ \frac{\lambda }{2}\right)_{i_k}} \right) \label{eq:40017c} \\
&&\times  \left. \sum_{i_{\tau } = i_{\tau -1}}^{m} \frac{\left( -\frac{j}{2}+\frac{\tau }{2}\right)_{i_{\tau }} \left( \frac{\beta }{2}+ \frac{\tau }{2}+ \frac{\lambda }{2}\right)_{i_{\tau }}\left( 1+\frac{\tau }{2}+\frac{\lambda }{2}\right)_{i_{\tau -1}} \left( \frac{\gamma }{2}+ \frac{1}{2}+\frac{\tau }{2}+ \frac{\lambda }{2}\right)_{i_{\tau -1}}}{\left( -\frac{j}{2}+\frac{\tau }{2}\right)_{i_{\tau -1}} \left( \frac{\beta }{2}+ \frac{\tau }{2}+ \frac{\lambda }{2}\right)_{i_{\tau -1}}\left( 1+\frac{\tau }{2}+\frac{\lambda }{2}\right)_{i_{\tau }} \left( \frac{\gamma }{2}+ \frac{1}{2}+\frac{\tau }{2}+ \frac{\lambda }{2}\right)_{i_{\tau }}} z^{i_{\tau }}\right\} \eta ^{\tau } \nonumber
\end{eqnarray}
where
\begin{equation}
\begin{cases} \tau \geq 2 \cr
z = -\frac{1}{a}x^2 \cr
\eta = \frac{(1+a)}{a} x \cr
\Gamma _0^{(F)} = \frac{1}{2(1+a)}\left( \beta -\delta -j +a\left( \gamma +\delta -1+\lambda \right)\right) \cr
\Gamma _k^{(F)} = \frac{1}{2(1+a)}\left( \beta -\delta -j+k +a\left( \gamma +\delta -1+k+\lambda \right)\right)
\end{cases}\nonumber
\end{equation}
\end{enumerate}
Put $c_0$= 1 as $\lambda =0$ for the first kind of independent solutions of Heun equation and $\displaystyle{ c_0= \left( \frac{1+a}{a}\right)^{1-\gamma }}$ as $\lambda = 1-\gamma $ for the second one in (\ref{eq:40013a})--(\ref{eq:40017c}). 
\begin{remark}
The power series expansion of Heun equation of the first kind for the first species complete polynomial using 3TRF about $x=0$ is given by
\begin{enumerate} 
\item As $\alpha =0$ and $q=q_0^0=0$,

The eigenfunction is given by
\begin{equation}
y(x) = H_pF_{0,0} \left( a, q=q_0^0=0; \alpha =0, \beta, \gamma, \delta ; \eta = \frac{(1+a)}{a} x ; z= -\frac{1}{a} x^2 \right) =1 \label{eq:40018}
\end{equation}
\item As $\alpha =-2N-1$ where $N \in \mathbb{N}_{0}$,

An algebraic equation of degree $2N+2$ for the determination of $q$ is given by
\begin{equation}
0 = \sum_{r=0}^{N+1}\bar{c}\left( 2r, N+1-r; 2N+1,q\right)\label{eq:40019a}
\end{equation}
The eigenvalue of $q$ is written by $q_{2N+1}^m$ where $m = 0,1,2,\cdots,2N+1 $; $q_{2N+1}^0 < q_{2N+1}^1 < \cdots < q_{2N+1}^{2N+1}$. Its eigenfunction is given by
\begin{eqnarray} 
y(x) &=& H_pF_{2N+1,m} \left( a, q=q_{2N+1}^m; \alpha =-2N-1, \beta, \gamma, \delta ; \eta = \frac{(1+a)}{a} x ; z= -\frac{1}{a} x^2 \right)\nonumber\\
&=& \sum_{r=0}^{N} y_{2r}^{N-r}\left( 2N+1,q_{2N+1}^m;x\right)+ \sum_{r=0}^{N} y_{2r+1}^{N-r}\left( 2N+1,q_{2N+1}^m;x\right)  
\label{eq:40019b}
\end{eqnarray}
\item As $\alpha =-2N-2$ where $N \in \mathbb{N}_{0}$,

An algebraic equation of degree $2N+3$ for the determination of $q$ is given by
\begin{eqnarray}
0  = \sum_{r=0}^{N+1}\bar{c}\left( 2r+1, N+1-r; 2N+2,q\right)\label{eq:40020a}
\end{eqnarray}
The eigenvalue of $q$ is written by $q_{2N+2}^m$ where $m = 0,1,2,\cdots,2N+2 $; $q_{2N+2}^0 < q_{2N+2}^1 < \cdots < q_{2N+2}^{2N+2}$. Its eigenfunction is given by
\begin{eqnarray} 
y(x) &=& H_pF_{2N+2,m} \left( a, q=q_{2N+2}^m; \alpha =-2N-2, \beta, \gamma, \delta ; \eta = \frac{(1+a)}{a} x ; z= -\frac{1}{a} x^2 \right)\nonumber\\
&=& \sum_{r=0}^{N+1} y_{2r}^{N+1-r}\left( 2N+2,q_{2N+2}^m;x\right) + \sum_{r=0}^{N} y_{2r+1}^{N-r}\left( 2N+2,q_{2N+2}^m;x\right) 
\label{eq:40020b}
\end{eqnarray}
In the above,
\begin{eqnarray}
\bar{c}(0,n;j,q)  &=& \frac{\left( -\frac{j}{2}\right)_{n} \left( \frac{\beta }{2} \right)_{n}}{\left( 1 \right)_{n} \left( \frac{\gamma }{2}+ \frac{1}{2} \right)_{n}} \left(-\frac{1}{a} \right)^{n}\label{eq:40021a}\\
\bar{c}(1,n;j,q) &=& \left( \frac{1+a}{a}\right) \sum_{i_0=0}^{n}\frac{ i_0 \left(i_0 +\Gamma _0^{(F)} \right)+\frac{q}{4(1+a)}}{\left( i_0+\frac{1}{2} \right) \left( i_0+\frac{\gamma }{2} \right)} \frac{\left( -\frac{j}{2}\right)_{i_0} \left( \frac{\beta }{2} \right)_{i_0}}{\left( 1 \right)_{i_0} \left( \frac{\gamma }{2}+ \frac{1}{2} \right)_{i_0}}\nonumber\\
&&\times  \frac{\left( -\frac{j}{2}+\frac{1}{2} \right)_{n} \left( \frac{\beta }{2}+ \frac{1}{2} \right)_{n}\left( \frac{3}{2} \right)_{i_0} \left( \frac{\gamma }{2}+1 \right)_{i_0}}{\left( -\frac{j}{2}+\frac{1}{2}\right)_{i_0} \left( \frac{\beta }{2}+ \frac{1}{2} \right)_{i_0}\left( \frac{3}{2} \right)_{n} \left( \frac{\gamma }{2}+1 \right)_{n}} \left(-\frac{1}{a} \right)^{n}  
\label{eq:40021b}\\
\bar{c}(\tau ,n;j,q) &=& \left( \frac{1+a}{a}\right)^{\tau} \sum_{i_0=0}^{n}\frac{ i_0 \left( i_0 +\Gamma _0^{(F)} \right)+\frac{q}{4(1+a)}}{\left( i_0+\frac{1}{2} \right) \left( i_0+\frac{\gamma }{2} \right)}  \frac{\left( -\frac{j}{2}\right)_{i_0} \left( \frac{\beta }{2} \right)_{i_0}}{\left( 1 \right)_{i_0} \left( \frac{\gamma }{2}+ \frac{1}{2} \right)_{i_0}} \nonumber\\
&&\times \prod_{k=1}^{\tau -1} \left( \sum_{i_k = i_{k-1}}^{n} \frac{\left( i_k+ \frac{k}{2} \right)\left(i_k +\Gamma _k^{(F)} \right)+\frac{q}{4(1+a)}}{\left( i_k+\frac{k}{2}+\frac{1}{2} \right) \left( i_k+\frac{k}{2}+\frac{\gamma }{2} \right)} \right. \nonumber\\
&&\times \left. \frac{\left( -\frac{j}{2}+\frac{k}{2}\right)_{i_k} \left( \frac{\beta }{2}+ \frac{k}{2} \right)_{i_k}\left( \frac{k}{2}+1 \right)_{i_{k-1}} \left( \frac{\gamma }{2}+ \frac{1}{2}+ \frac{k}{2} \right)_{i_{k-1}}}{\left( -\frac{j}{2}+\frac{k}{2}\right)_{i_{k-1}} \left( \frac{\beta }{2}+ \frac{k}{2} \right)_{i_{k-1}}\left( \frac{k}{2} +1\right)_{i_k} \left( \frac{\gamma }{2}+ \frac{1}{2}+ \frac{k}{2} \right)_{i_k}} \right) \nonumber\\
&&\times \frac{\left( -\frac{j}{2}+\frac{\tau }{2}\right)_{n} \left( \frac{\beta }{2}+ \frac{\tau }{2} \right)_{n}\left( \frac{\tau }{2} +1\right)_{i_{\tau -1}} \left( \frac{\gamma }{2}+ \frac{1}{2}+\frac{\tau }{2} \right)_{i_{\tau -1}}}{\left( -\frac{j}{2}+\frac{\tau }{2}\right)_{i_{\tau -1}} \left( \frac{\beta }{2}+\frac{\tau }{2} \right)_{i_{\tau -1}}\left( \frac{\tau }{2}+1 \right)_{n} \left( \frac{\gamma }{2}+ \frac{1}{2}+\frac{\tau }{2} \right)_{n}} \left(-\frac{1}{a} \right)^{n } \hspace{1.5cm} \label{eq:40021c} 
\end{eqnarray}
\begin{eqnarray}
y_0^m(j,q;x) &=& \sum_{i_0=0}^{m} \frac{\left( -\frac{j}{2}\right)_{i_0} \left( \frac{\beta }{2} \right)_{i_0}}{\left( 1 \right)_{i_0} \left( \frac{\gamma }{2}+ \frac{1}{2} \right)_{i_0}} z^{i_0} \label{eq:40022a}\\
y_1^m(j,q;x) &=& \left\{\sum_{i_0=0}^{m} \frac{ i_0 \left(i_0 +\Gamma _0^{(F)} \right)+\frac{q}{4(1+a)}}{\left( i_0+\frac{1}{2} \right) \left( i_0+\frac{\gamma }{2} \right)} \frac{\left( -\frac{j}{2}\right)_{i_0} \left( \frac{\beta }{2} \right)_{i_0}}{\left( 1 \right)_{i_0} \left( \frac{\gamma }{2}+ \frac{1}{2} \right)_{i_0}} \right. \nonumber\\
&&\times \left. \sum_{i_1 = i_0}^{m} \frac{\left( -\frac{j}{2}+\frac{1}{2} \right)_{i_1} \left( \frac{\beta }{2}+ \frac{1}{2} \right)_{i_1}\left( \frac{3}{2} \right)_{i_0} \left( \frac{\gamma }{2}+1 \right)_{i_0}}{\left(-\frac{j}{2}+\frac{1}{2} \right)_{i_0} \left( \frac{\beta }{2}+ \frac{1}{2} \right)_{i_0}\left( \frac{3}{2} \right)_{i_1} \left( \frac{\gamma }{2} +1 \right)_{i_1}} z^{i_1}\right\} \eta 
\label{eq:40022b}\\
y_{\tau }^m(j,q;x) &=& \left\{ \sum_{i_0=0}^{m} \frac{ i_0 \left(i_0 +\Gamma _0^{(F)} \right)+\frac{q}{4(1+a)}}{\left( i_0+\frac{1}{2} \right) \left( i_0+\frac{\gamma }{2} \right)} \frac{\left( -\frac{j}{2}\right)_{i_0} \left( \frac{\beta }{2} \right)_{i_0}}{\left( 1 \right)_{i_0} \left( \frac{\gamma }{2}+ \frac{1}{2} \right)_{i_0}} \right.\nonumber\\
&&\times \prod_{k=1}^{\tau -1} \left( \sum_{i_k = i_{k-1}}^{m} \frac{\left( i_k+ \frac{k}{2} \right)\left(i_k +\Gamma _k^{(F)} \right)+\frac{q}{4(1+a)}}{\left( i_k+\frac{k}{2}+\frac{1}{2} \right) \left( i_k+\frac{k}{2}+\frac{\gamma }{2} \right)} \right. \nonumber\\
&&\times \left. \frac{\left( -\frac{j}{2}+\frac{k}{2}\right)_{i_k} \left( \frac{\beta }{2}+ \frac{k}{2} \right)_{i_k}\left( \frac{k}{2}+1 \right)_{i_{k-1}} \left( \frac{\gamma }{2}+ \frac{1}{2}+ \frac{k}{2} \right)_{i_{k-1}}}{\left( -\frac{j}{2}+\frac{k}{2}\right)_{i_{k-1}} \left( \frac{\beta }{2}+ \frac{k}{2} \right)_{i_{k-1}}\left( \frac{k}{2} +1\right)_{i_k} \left( \frac{\gamma }{2}+ \frac{1}{2}+ \frac{k}{2} \right)_{i_k}} \right) \nonumber\\
&&\times \left. \sum_{i_{\tau } = i_{\tau -1}}^{m} \frac{\left( -\frac{j}{2}+\frac{\tau }{2}\right)_{i_{\tau }} \left( \frac{\beta }{2}+ \frac{\tau }{2} \right)_{i_{\tau }}\left( \frac{\tau }{2}+1 \right)_{i_{\tau -1}} \left( \frac{\gamma }{2}+ \frac{1}{2}+\frac{\tau }{2} \right)_{i_{\tau -1}}}{\left( -\frac{j}{2}+\frac{\tau }{2}\right)_{i_{\tau -1}} \left( \frac{\beta }{2}+ \frac{\tau }{2} \right)_{i_{\tau -1}}\left( \frac{\tau }{2}+1 \right)_{i_{\tau }} \left( \frac{\gamma }{2}+ \frac{1}{2}+\frac{\tau }{2} \right)_{i_{\tau }}} z^{i_{\tau }}\right\} \eta ^{\tau } \hspace{1.5cm} \label{eq:40022c}
\end{eqnarray}
where
\begin{equation}
\begin{cases} \tau \geq 2 \cr
\Gamma _0^{(F)} = \frac{1}{2(1+a)}\left( \beta -\delta -j +a\left( \gamma +\delta -1 \right)\right) \cr
\Gamma _k^{(F)} = \frac{1}{2(1+a)}\left( \beta -\delta -j+k +a\left( \gamma +\delta -1+k \right)\right)
\end{cases}\nonumber
\end{equation}
\end{enumerate}
\end{remark}
\begin{remark}
The power series expansion of Heun equation of the second kind for the first species complete polynomial using 3TRF about $x=0$ is given by
\begin{enumerate} 
\item As $\alpha = \gamma -1$ and $q=q_0^0= (\gamma -1)(\beta -\delta +a\delta )$,

The eigenfunction is given by
\begin{eqnarray}
y(x) &=& H_pS_{0,0} \left( a, q=q_0^0=(\gamma -1)(\beta -\delta +a\delta ); \alpha =\gamma -1, \beta, \gamma, \delta ; \eta = \frac{(1+a)}{a} x ; z= -\frac{1}{a} x^2 \right) \nonumber\\
&=& \eta ^{1-\gamma } \label{eq:40023}
\end{eqnarray}
\item As $\alpha =\gamma -2N-2 $ where $N \in \mathbb{N}_{0}$, 

An algebraic equation of degree $2N+2$ for the determination of $q$ is given by
\begin{equation}
0 = \sum_{r=0}^{N+1}\bar{c}\left( 2r, N+1-r; 2N+1,q\right) \label{eq:40024a}
\end{equation}
The eigenvalue of $q$ is written by $q_{2N+1}^m$ where $m = 0,1,2,\cdots,2N+1 $; $q_{2N+1}^0 < q_{2N+1}^1 < \cdots < q_{2N+1}^{2N+1}$. Its eigenfunction is given by
\begin{eqnarray} 
y(x) &=& H_pS_{2N+1,m} \left( a, q=q_{2N+1}^m; \alpha =\gamma -2N-2, \beta, \gamma, \delta ; \eta = \frac{(1+a)}{a} x ; z= -\frac{1}{a} x^2 \right) \nonumber\\
&=&  \sum_{r=0}^{N} y_{2r}^{N-r}\left( 2N+1,q_{2N+1}^m;x\right)+ \sum_{r=0}^{N} y_{2r+1}^{N-r}\left( 2N+1,q_{2N+1}^m;x\right)
\label{eq:40024b}
\end{eqnarray}
\item As $\alpha =\gamma -2N-3 $ where  $N \in \mathbb{N}_{0}$,

An algebraic equation of degree $2N+3$ for the determination of $q$ is given by
\begin{equation}
0 = \sum_{r=0}^{N+1}\bar{c}\left( 2r+1, N+1-r; 2N+2,q\right)  \label{eq:40025a}
\end{equation}
The eigenvalue of $q$ is written by $q_{2N+2}^m$ where $m = 0,1,2,\cdots,2N+2 $; $q_{2N+2}^0 < q_{2N+2}^1 < \cdots < q_{2N+2}^{2N+2}$. Its eigenfunction is given by
\begin{eqnarray} 
y(x) &=& H_pS_{2N+2,m} \left( a, q=q_{2N+2}^m; \alpha =\gamma -2N-3, \beta, \gamma, \delta ; \eta = \frac{(1+a)}{a} x ; z= -\frac{1}{a} x^2 \right) \nonumber\\
&=& \sum_{r=0}^{N+1} y_{2r}^{N+1-r}\left( 2N+2,q_{2N+2}^m;x\right) + \sum_{r=0}^{N} y_{2r+1}^{N-r}\left( 2N+2,q_{2N+2}^m;x\right) 
\label{eq:40025b}
\end{eqnarray}
In the above,
\begin{eqnarray}
\bar{c}(0,n;j,q)  &=& \frac{\left( -\frac{j}{2}\right)_{n} \left( \frac{\beta }{2}-\frac{\gamma }{2}+ \frac{1}{2}\right)_{n}}{\left( \frac{3}{2}-\frac{\gamma }{2} \right)_{n} \left( 1\right)_{n}} \left(-\frac{1}{a} \right)^{n}\label{eq:40026a}\\
\bar{c}(1,n;j,q) &=& \left( \frac{1+a}{a}\right) \sum_{i_0=0}^{n}\frac{\left( i_0+\frac{1}{2}-\frac{\gamma }{2}\right)\left( i_0 +\Gamma _0^{(F)} \right) +\frac{q}{4(1+a)}}{\left( i_0+1-\frac{\gamma }{2}\right) \left( i_0 +\frac{1}{2} \right)} \frac{\left( -\frac{j}{2}\right)_{i_0} \left( \frac{\beta }{2}-\frac{\gamma}{2}+ \frac{1}{2} \right)_{i_0}}{\left( \frac{3}{2}-\frac{\gamma }{2}\right)_{i_0} \left( 1\right)_{i_0}}\nonumber\\
&&\times  \frac{\left( -\frac{j}{2}+\frac{1}{2} \right)_{n} \left( \frac{\beta }{2} -\frac{\gamma }{2}+1\right)_{n}\left( 2-\frac{\gamma }{2}\right)_{i_0} \left( \frac{3}{2}\right)_{i_0}}{\left( -\frac{j}{2}+\frac{1}{2}\right)_{i_0} \left( \frac{\beta }{2} -\frac{\gamma }{2}+1\right)_{i_0}\left( 2-\frac{\gamma }{2}\right)_{n} \left( \frac{3}{2} \right)_{n}} \left(-\frac{1}{a} \right)^{n }  
\label{eq:40026b}\\
\bar{c}(\tau ,n;j,q) &=& \left( \frac{1+a}{a}\right)^{\tau} \sum_{i_0=0}^{n}\frac{\left( i_0+\frac{1}{2}-\frac{\gamma }{2}\right)\left( i_0 +\Gamma _0^{(F)} \right) +\frac{q}{4(1+a)}}{\left( i_0+1-\frac{\gamma }{2}\right) \left( i_0 +\frac{1}{2} \right)} \frac{\left( -\frac{j}{2}\right)_{i_0} \left( \frac{\beta }{2}-\frac{\gamma}{2}+ \frac{1}{2} \right)_{i_0}}{\left( \frac{3}{2}-\frac{\gamma }{2}\right)_{i_0} \left( 1\right)_{i_0}} \nonumber\\
&&\times \prod_{k=1}^{\tau -1} \left( \sum_{i_k = i_{k-1}}^{n} \frac{\left( i_k+ \frac{k}{2}+\frac{1}{2}-\frac{\gamma}{2}\right)\left( i_k +\Gamma _k^{(F)} \right) +\frac{q}{4(1+a)}}{\left( i_k+\frac{k}{2}+1-\frac{\gamma}{2}\right) \left( i_k+\frac{k}{2} +\frac{1}{2} \right)} \right. \nonumber\\
&&\times \left. \frac{\left( -\frac{j}{2}+\frac{k}{2}\right)_{i_k} \left( \frac{\beta }{2}-\frac{\gamma}{2}+ \frac{k}{2}+ \frac{1}{2} \right)_{i_k}\left( \frac{k}{2}+\frac{3}{2}-\frac{\gamma}{2}\right)_{i_{k-1}} \left( \frac{k}{2}+1 \right)_{i_{k-1}}}{\left( -\frac{j}{2}+\frac{k}{2}\right)_{i_{k-1}} \left( \frac{\beta }{2}-\frac{\gamma}{2}+ \frac{k}{2}+\frac{1}{2} \right)_{i_{k-1}}\left(  \frac{k}{2}+\frac{3}{2}-\frac{\gamma}{2}\right)_{i_k} \left( \frac{k}{2}+ 1 \right)_{i_k}} \right) \label{eq:40026c}\\
&&\times \frac{\left( -\frac{j}{2}+\frac{\tau }{2}\right)_{n} \left( \frac{\beta }{2}-\frac{\gamma}{2}+ \frac{\tau }{2}+ \frac{1}{2} \right)_{n}\left( \frac{\tau }{2}+\frac{3}{2}-\frac{\gamma }{2}\right)_{i_{\tau -1}} \left( \frac{\tau }{2}+1 \right)_{i_{\tau -1}}}{\left( -\frac{j}{2}+\frac{\tau }{2}\right)_{i_{\tau -1}} \left( \frac{\beta }{2}-\frac{\gamma}{2}+\frac{\tau }{2}+ \frac{1}{2}\right)_{i_{\tau -1}}\left( \frac{\tau }{2}+\frac{3}{2}-\frac{\gamma }{2}\right)_{n} \left( \frac{\tau }{2}+1 \right)_{n}} \left(-\frac{1}{a} \right)^{n } \nonumber 
\end{eqnarray}
\begin{eqnarray}
y_0^m(j,q;x) &=& \eta ^{1-\gamma }  \sum_{i_0=0}^{m} \frac{\left( -\frac{j}{2}\right)_{i_0} \left( \frac{\beta }{2} -\frac{\gamma}{2}+\frac{1}{2}\right)_{i_0}}{\left( \frac{3}{2}-\frac{\gamma}{2}\right)_{i_0} \left( 1\right)_{i_0}} z^{i_0} \label{eq:40027a}\\
y_1^m(j,q;x) &=& \eta ^{1-\gamma } \left\{\sum_{i_0=0}^{m} \frac{\left( i_0+\frac{1}{2}-\frac{\gamma}{2}\right)\left( i_0 \Gamma _0^{(F)} \right) +\frac{q}{4(1+a)}}{\left( i_0+1-\frac{\gamma}{2}\right) \left( i_0 +\frac{1}{2} \right)} \frac{\left( -\frac{j}{2}\right)_{i_0} \left( \frac{\beta }{2} -\frac{\gamma}{2}+\frac{1}{2}\right)_{i_0}}{\left( \frac{3}{2}-\frac{\gamma}{2}\right)_{i_0} \left( 1\right)_{i_0}} \right. \nonumber\\
&&\times \left. \sum_{i_1 = i_0}^{m} \frac{\left( -\frac{j}{2}+\frac{1}{2} \right)_{i_1} \left( \frac{\beta }{2} -\frac{\gamma}{2}+1\right)_{i_1}\left( 2-\frac{\gamma}{2}\right)_{i_0} \left( \frac{3}{2} \right)_{i_0}}{\left(-\frac{j}{2}+\frac{1}{2} \right)_{i_0} \left( \frac{\beta }{2} -\frac{\gamma}{2}+1\right)_{i_0}\left( 2-\frac{\gamma}{2}\right)_{i_1} \left( \frac{3}{2} \right)_{i_1}} z^{i_1}\right\} \eta 
\label{eq:40027b} 
\end{eqnarray}
\begin{eqnarray}
y_{\tau }^m(j,q;x) &=& \eta ^{1-\gamma } \left\{ \sum_{i_0=0}^{m} \frac{\left( i_0+\frac{1}{2}-\frac{\gamma}{2}\right)\left( i_0 +\Gamma _0^{(F)} \right) +\frac{q}{4(1+a)}}{\left( i_0+1-\frac{\gamma}{2}\right) \left( i_0 +\frac{1}{2} \right)} \frac{\left( -\frac{j}{2}\right)_{i_0} \left( \frac{\beta }{2} -\frac{\gamma}{2}+\frac{1}{2}\right)_{i_0}}{\left( \frac{3}{2}-\frac{\gamma}{2}\right)_{i_0} \left( 1\right)_{i_0}} \right.\nonumber\\
&&\times \prod_{k=1}^{\tau -1} \left( \sum_{i_k = i_{k-1}}^{m} \frac{\left( i_k+ \frac{k}{2}+\frac{1}{2}-\frac{\gamma}{2}\right)\left( i_k +\Gamma _k^{(F)} \right) +\frac{q}{4(1+a)}}{\left( i_k+\frac{k}{2}+1-\frac{\gamma}{2}\right) \left( i_k+\frac{k}{2} +\frac{1}{2} \right)} \right. \nonumber\\
&&\times \left. \frac{\left( -\frac{j}{2}+\frac{k}{2}\right)_{i_k} \left( \frac{\beta }{2}-\frac{\gamma}{2}+ \frac{k}{2}+ \frac{1}{2} \right)_{i_k}\left( \frac{k}{2}+\frac{3}{2}-\frac{\gamma}{2}\right)_{i_{k-1}} \left( \frac{k}{2}+1 \right)_{i_{k-1}}}{\left( -\frac{j}{2}+\frac{k}{2}\right)_{i_{k-1}} \left( \frac{\beta }{2}-\frac{\gamma}{2}+ \frac{k}{2}+\frac{1}{2} \right)_{i_{k-1}}\left(  \frac{k}{2}+\frac{3}{2}-\frac{\gamma}{2}\right)_{i_k} \left( \frac{k}{2}+ 1 \right)_{i_k}} \right) \label{eq:40027c}\\
&&\times \left. \sum_{i_{\tau } = i_{\tau -1}}^{m} \frac{\left( -\frac{j}{2}+\frac{\tau }{2}\right)_{i_{\tau }} \left( \frac{\beta }{2}-\frac{\gamma}{2}+ \frac{\tau }{2}+ \frac{1}{2} \right)_{i_{\tau }}\left( \frac{\tau }{2}+\frac{3}{2}-\frac{\gamma}{2}\right)_{i_{\tau -1}} \left( \frac{\tau }{2}+ 1 \right)_{i_{\tau -1}}}{\left( -\frac{j}{2}+\frac{\tau }{2}\right)_{i_{\tau -1}} \left( \frac{\beta }{2}-\frac{\gamma}{2}+ \frac{\tau }{2}+ \frac{1}{2} \right)_{i_{\tau -1}}\left( \frac{\tau }{2}+\frac{3}{2}-\frac{\gamma}{2}\right)_{i_{\tau }} \left( \frac{\tau }{2}+1 \right)_{i_{\tau }}} z^{i_{\tau }}\right\} \eta ^{\tau } \nonumber
\end{eqnarray}
\end{enumerate}
where
\begin{equation}
\begin{cases} \tau \geq 2 \cr
\Gamma _0^{(F)} = \frac{1}{2(1+a)}\left( \beta -\delta -j +a \delta \right) \cr
\Gamma _k^{(F)} = \frac{1}{2(1+a)}\left( \beta -\delta -j+k +a\left( \delta +k \right)\right)
\end{cases}\nonumber
\end{equation}
\end{remark}
\subsection{The second species complete polynomial using 3TRF}

For the second species complete polynomial, we need a condition which is defined by
\begin{equation}
B_{j}=B_{j+1}= A_{j}=0\hspace{1cm}\mathrm{where}\;j \in \mathbb{N}_{0}    
 \label{eq:40028}
\end{equation}
\begin{theorem}
In chapter 1, the general expression of a function $y(x)$ for the second species complete polynomial using 3-term recurrence formula is given by
\begin{enumerate} 
\item As $B_1=A_0=0$,
\begin{equation}
y(x) = y_{0}^{0}(x) \label{eq:40029a}
\end{equation}
\item As $B_{2N+1}=B_{2N+2}=A_{2N+1}=0$  where $N \in \mathbb{N}_{0}$,
\begin{equation}
y(x)= \sum_{r=0}^{N} y_{2r}^{N-r}(x)+ \sum_{r=0}^{N} y_{2r+1}^{N-r}(x)  \label{eq:40029b}
\end{equation}
\item As $B_{2N+2}=B_{2N+3}=A_{2N+2}=0$  where $N \in \mathbb{N}_{0}$,
\begin{equation}
y(x)= \sum_{r=0}^{N+1} y_{2r}^{N+1-r}(x)+ \sum_{r=0}^{N} y_{2r+1}^{N-r}(x)  \label{eq:40029c}
\end{equation}
In the above,
\begin{eqnarray}
y_0^m(x) &=& c_0 x^{\lambda } \sum_{i_0=0}^{m} \left\{ \prod _{i_1=0}^{i_0-1}B_{2i_1+1} \right\} x^{2i_0 } \label{eq:40030a}\\
y_1^m(x) &=& c_0 x^{\lambda } \sum_{i_0=0}^{m}\left\{ A_{2i_0} \prod _{i_1=0}^{i_0-1}B_{2i_1+1}  \sum_{i_2=i_0}^{m} \left\{ \prod _{i_3=i_0}^{i_2-1}B_{2i_3+2} \right\}\right\} x^{2i_2+1 } \label{eq:40030b}\\
y_{\tau }^m(x) &=& c_0 x^{\lambda } \sum_{i_0=0}^{m} \left\{A_{2i_0}\prod _{i_1=0}^{i_0-1} B_{2i_1+1} 
\prod _{k=1}^{\tau -1} \left( \sum_{i_{2k}= i_{2(k-1)}}^{m} A_{2i_{2k}+k}\prod _{i_{2k+1}=i_{2(k-1)}}^{i_{2k}-1}B_{2i_{2k+1}+(k+1)}\right) \right. \nonumber\\
&& \times \left. \sum_{i_{2\tau} = i_{2(\tau -1)}}^{m} \left( \prod _{i_{2\tau +1}=i_{2(\tau -1)}}^{i_{2\tau}-1} B_{2i_{2\tau +1}+(\tau +1)} \right) \right\} x^{2i_{2\tau}+\tau }\hspace{1cm}\mathrm{where}\;\tau \geq 2
\label{eq:40030c}  
\end{eqnarray}
\end{enumerate}
\end{theorem}
Put $n= j+1$ in (\ref{eq:4004b}) and use the condition $B_{j+1}=0$ for $\alpha $.  
\begin{equation}
\alpha = -j-\lambda 
\label{eq:40031}
\end{equation}
Put $n= j$ in (\ref{eq:4004b}) and use the condition $B_{j}=0$ for $\beta $.  
\begin{equation}
\beta = -j+1-\lambda 
\label{eq:40032}
\end{equation}
Substitute (\ref{eq:40031}) and (\ref{eq:40032}) into (\ref{eq:4004a}). Put $n= j$ in the new (\ref{eq:4004a}) and use the condition $A_{j}=0$ for $q$.  
\begin{equation}
q = -(j+\lambda )\left[ -\delta -j+1-\lambda +a(\gamma +\delta +j-1+\lambda )\right]
\label{eq:40033}
\end{equation}
Take (\ref{eq:40031}), (\ref{eq:40032}) and (\ref{eq:40033}) into (\ref{eq:4004a}) and (\ref{eq:4004b}).
\begin{subequations}
\begin{equation}
A_n = \frac{1+a}{a}\frac{(n-j)\left[n+\frac{1}{1+a}\left( -\delta -j+1+a\left( \gamma +\delta +j-1+2\lambda \right)\right) \right]}{(n+1+\lambda )(n+\gamma +\lambda )} \label{eq:40034a}
\end{equation}
\begin{equation}
B_n = -\frac{1}{a}\frac{(n-j)(n-j-1)}{(n+1+\lambda )(n+\gamma +\lambda )} \label{eq:40034b}
\end{equation}
\end{subequations}
Substitute (\ref{eq:40034a}) and (\ref{eq:40034b}) into (\ref{eq:40030a})--(\ref{eq:40030c}).

As $B_1=A_0=0$, substitute the new (\ref{eq:40030a}) into (\ref{eq:40029a}) putting $j=0$. 

As $B_{2N+1}=B_{2N+2}=A_{2N+1}=0$, substitute the new (\ref{eq:40030a})--(\ref{eq:40030c}) into (\ref{eq:40029b}) putting $j=2N+1$.

As $B_{2N+2}=B_{2N+3}=A_{2N+2}=0$, substitute the new (\ref{eq:40030a})--(\ref{eq:40030c}) into (\ref{eq:40029c}) putting $j=2N+2$.

After the replacement process, the general expression of power series of Heun equation about $x=0$ for the second species complete polynomial using 3-term recurrence formula is given by
\begin{enumerate} 
\item As $\alpha =-\lambda $, $\beta =1-\lambda $ and $q= -\lambda \left[ -\delta +1-\lambda +a\left( \gamma +\delta -1+\lambda \right)\right]$,
 
Its eigenfunction is given by
\begin{equation}
y(x) = y_0^0(0;x)= c_0 x^{\lambda } \label{eq:40035a}
\end{equation}
\item As $\alpha =-2N-1-\lambda $, $\beta = -2N-\lambda $ and $q= -\left(2N+1+\lambda \right) \left[ -\delta -2N-\lambda +a\left( \gamma +\delta+2N+\lambda \right)\right]$ where $N \in \mathbb{N}_{0}$,

Its eigenfunction is given by
\begin{equation}
y(x)= \sum_{r=0}^{N} y_{2r}^{N-r}\left( 2N+1;x\right) + \sum_{r=0}^{N} y_{2r+1}^{N-r}\left( 2N+1;x\right)  \label{eq:40035b}
\end{equation}
\item As $\alpha =-2N-2-\lambda $, $\beta = -2N-1-\lambda $ and $q= -\left(2N+2+\lambda \right) \left[ -\delta -2N-1-\lambda +a\left( \gamma +\delta+2N+1+\lambda \right)\right]$ where $N \in \mathbb{N}_{0}$,

Its eigenfunction is given by
\begin{equation}
y(x)= \sum_{r=0}^{N+1} y_{2r}^{N+1-r}\left( 2N+2;x\right) + \sum_{r=0}^{N} y_{2r+1}^{N-r}\left( 2N+2;x\right)  \label{eq:40035c}
\end{equation}
In the above,
\begin{eqnarray}
y_0^m(j;x) &=& c_0 x^{\lambda }  \sum_{i_0=0}^{m} \frac{\left( -\frac{j}{2}\right)_{i_0} \left( - \frac{j}{2}+\frac{1}{2}\right)_{i_0}}{\left( 1+\frac{\lambda }{2}\right)_{i_0} \left( \frac{\gamma }{2}+ \frac{1}{2}+ \frac{\lambda }{2}\right)_{i_0}} z^{i_0} \label{eq:40036a}\\
y_1^m(j;x) &=& c_0 x^{\lambda } \left\{\sum_{i_0=0}^{m} \frac{\left( i_0-\frac{j}{2}\right)\left( i_0 +\Gamma ^{(S)} \right)}{\left( i_0+\frac{1}{2}+\frac{\lambda }{2}\right) \left( i_0+\frac{\gamma }{2}+\frac{\lambda }{2}\right)} \frac{\left( -\frac{j}{2}\right)_{i_0} \left( -\frac{j}{2}+\frac{1}{2}\right)_{i_0}}{\left( 1+\frac{\lambda }{2}\right)_{i_0} \left( \frac{\gamma }{2}+ \frac{1}{2}+ \frac{\lambda }{2}\right)_{i_0}} \right. \nonumber\\
&&\times  \left. \sum_{i_1 = i_0}^{m} \frac{\left( -\frac{j}{2}+\frac{1}{2} \right)_{i_1} \left( -\frac{j}{2}+1\right)_{i_1}\left( \frac{3}{2}+\frac{\lambda }{2}\right)_{i_0} \left( \frac{\gamma }{2}+1+ \frac{\lambda }{2}\right)_{i_0}}{\left(-\frac{j}{2}+\frac{1}{2} \right)_{i_0} \left( -\frac{j}{2}+1\right)_{i_0}\left( \frac{3}{2}+\frac{\lambda }{2}\right)_{i_1} \left( \frac{\gamma }{2} +1+ \frac{\lambda }{2}\right)_{i_1}} z^{i_1}\right\} \eta \hspace{1.5cm}\label{eq:40036b}
\end{eqnarray}
\begin{eqnarray} 
y_{\tau }^m(j;x) &=& c_0 x^{\lambda } \left\{ \sum_{i_0=0}^{m} \frac{\left( i_0-\frac{j}{2}\right)\left(i_0 +\Gamma ^{(S)} \right)}{\left( i_0+\frac{1}{2}+\frac{\lambda }{2}\right) \left( i_0+\frac{\gamma }{2}+\frac{\lambda }{2}\right)} \frac{\left( -\frac{j}{2}\right)_{i_0} \left( -\frac{j}{2}+\frac{1}{2}\right)_{i_0}}{\left( 1+\frac{\lambda }{2}\right)_{i_0} \left( \frac{\gamma }{2}+ \frac{1}{2}+ \frac{\lambda }{2}\right)_{i_0}} \right.\nonumber\\
&&\times   \prod_{k=1}^{\tau -1} \left( \sum_{i_k = i_{k-1}}^{m} \frac{\left( i_k+ \frac{k}{2}-\frac{j}{2}\right)\left( i_k + \frac{k}{2}+\Gamma ^{(S)} \right)}{\left( i_k+\frac{k}{2}+\frac{1}{2}+\frac{\lambda }{2}\right) \left( i_k+\frac{k}{2}+\frac{\gamma }{2}+\frac{\lambda }{2}\right)} \right. \nonumber\\
&&\times  \left. \frac{\left( -\frac{j}{2}+\frac{k}{2}\right)_{i_k} \left( -\frac{j}{2}+ \frac{k}{2}+ \frac{1}{2}\right)_{i_k}\left( 1+ \frac{k}{2}+\frac{\lambda }{2}\right)_{i_{k-1}} \left( \frac{\gamma }{2}+ \frac{1}{2}+ \frac{k}{2}+ \frac{\lambda }{2}\right)_{i_{k-1}}}{\left( -\frac{j}{2}+\frac{k}{2}\right)_{i_{k-1}} \left( -\frac{j}{2}+ \frac{k}{2}+ \frac{1}{2}\right)_{i_{k-1}}\left( 1+ \frac{k}{2}+\frac{\lambda }{2}\right)_{i_k} \left( \frac{\gamma }{2}+ \frac{1}{2}+ \frac{k}{2}+ \frac{\lambda }{2}\right)_{i_k}} \right) \label{eq:40036c}\\
&&\times  \left. \sum_{i_{\tau } = i_{\tau -1}}^{m} \frac{\left( -\frac{j}{2}+\frac{\tau }{2}\right)_{i_{\tau }} \left( -\frac{j}{2}+ \frac{\tau }{2}+ \frac{1}{2}\right)_{i_{\tau }}\left( 1+\frac{\tau }{2}+\frac{\lambda }{2}\right)_{i_{\tau -1}} \left( \frac{\gamma }{2}+ \frac{1}{2}+\frac{\tau }{2}+ \frac{\lambda }{2}\right)_{i_{\tau -1}}}{\left( -\frac{j}{2}+\frac{\tau }{2}\right)_{i_{\tau -1}} \left(-\frac{j}{2}+ \frac{\tau }{2}+ \frac{1}{2}\right)_{i_{\tau -1}}\left( 1+\frac{\tau }{2}+\frac{\lambda }{2}\right)_{i_{\tau }} \left( \frac{\gamma }{2}+ \frac{1}{2}+\frac{\tau }{2}+ \frac{\lambda }{2}\right)_{i_{\tau }}} z^{i_{\tau }}\right\} \eta ^{\tau } \nonumber
\end{eqnarray}
where
\begin{equation}
\begin{cases} \tau \geq 2 \cr
z = -\frac{1}{a}x^2 \cr
\eta = \frac{(1+a)}{a} x \cr
\Gamma ^{(S)} = \frac{1}{2(1+a)}\left( -\delta -j+1 +a\left( \gamma +\delta -1+j +2\lambda \right)\right)  
\end{cases}\nonumber
\end{equation}
\end{enumerate}
Put $c_0$= 1 as $\lambda =0$ for the first kind of independent solutions of Heun equation and $\displaystyle{ c_0= \left( \frac{1+a}{a}\right)^{1-\gamma }}$ as $\lambda = 1-\gamma $ for the second one in (\ref{eq:40035a})--(\ref{eq:40036c}). 
\begin{remark}
The power series expansion of Heun equation of the first kind for the second species complete polynomial using 3TRF about $x=0$ is given by
\begin{enumerate} 
\item As $\alpha =0 $, $\beta =1 $ and $q=0$,
 
Its eigenfunction is given by
\begin{equation}
y(x) = H_pF_0 \left( a, q=0; \alpha =0, \beta =1, \gamma, \delta ; \eta = \frac{(1+a)}{a} x ; z= -\frac{1}{a} x^2 \right) =1 \label{eq:40037a}
\end{equation}
\item As $\alpha =-2N-1 $, $\beta = -2N $ and $q= \left(2N+1 \right) \left[ \delta +2N -a\left( \gamma +\delta+2N \right)\right]$ where $N \in \mathbb{N}_{0}$,

Its eigenfunction is given by
\begin{eqnarray} 
y(x) &=& H_pF_{2N+1} \Bigg( a, q= \left( 2N+1\right) \left[ \delta +2N -a\left( \gamma +\delta+2N \right)\right]; \alpha =-2N-1 \nonumber\\
&&, \beta =-2N, \gamma, \delta; \eta = \frac{(1+a)}{a} x ; z= -\frac{1}{a} x^2 \Bigg) \nonumber\\
&=& \sum_{r=0}^{N} y_{2r}^{N-r}\left( 2N+1;x\right) + \sum_{r=0}^{N} y_{2r+1}^{N-r}\left( 2N+1;x\right) 
\label{eq:40037b}
\end{eqnarray}
\item As $\alpha =-2N-2 $, $\beta = -2N-1 $ and $q= \left(2N+2 \right) \left[ \delta +2N+1 -a\left( \gamma +\delta +2N+1 \right)\right]$ where $N \in \mathbb{N}_{0}$,

Its eigenfunction is given by
\begin{eqnarray} 
y(x) &=& H_pF_{2N+2} \Bigg( a, q= \left( 2N+2 \right) \left[ \delta +2N+1 -a\left( \gamma +\delta +2N+1 \right)\right]; \alpha =-2N-2 \nonumber\\
&&, \beta =-2N-1, \gamma, \delta; \eta = \frac{(1+a)}{a} x ; z= -\frac{1}{a} x^2 \Bigg) \nonumber\\
&=& \sum_{r=0}^{N+1} y_{2r}^{N+1-r}\left( 2N+2;x\right) + \sum_{r=0}^{N} y_{2r+1}^{N-r}\left( 2N+2;x\right) 
\label{eq:40037c}
\end{eqnarray}
In the above,
\begin{eqnarray}
y_0^m(j;x) &=&  \sum_{i_0=0}^{m} \frac{\left( -\frac{j}{2}\right)_{i_0} \left( - \frac{j}{2}+\frac{1}{2}\right)_{i_0}}{\left( 1 \right)_{i_0} \left( \frac{\gamma }{2}+ \frac{1}{2} \right)_{i_0}} z^{i_0} \label{eq:40038a}\\
y_1^m(j;x) &=&  \left\{\sum_{i_0=0}^{m} \frac{\left( i_0-\frac{j}{2}\right)\left( i_0 +\Gamma ^{(S)} \right)}{\left( i_0+\frac{1}{2} \right) \left( i_0+\frac{\gamma }{2} \right)} \frac{\left( -\frac{j}{2}\right)_{i_0} \left( -\frac{j}{2}+\frac{1}{2}\right)_{i_0}}{\left( 1 \right)_{i_0} \left( \frac{\gamma }{2}+ \frac{1}{2} \right)_{i_0}} \right. \nonumber\\
&&\times \left. \sum_{i_1 = i_0}^{m} \frac{\left( -\frac{j}{2}+\frac{1}{2} \right)_{i_1} \left( -\frac{j}{2}+1\right)_{i_1}\left( \frac{3}{2} \right)_{i_0} \left( \frac{\gamma }{2}+1 \right)_{i_0}}{\left(-\frac{j}{2}+\frac{1}{2} \right)_{i_0} \left( -\frac{j}{2}+1\right)_{i_0}\left( \frac{3}{2} \right)_{i_1} \left( \frac{\gamma }{2} +1 \right)_{i_1}} z^{i_1}\right\} \eta \label{eq:40038b}
\end{eqnarray}
\begin{eqnarray}
y_{\tau }^m(j;x) &=&  \left\{ \sum_{i_0=0}^{m} \frac{\left( i_0-\frac{j}{2}\right)\left( i_0 +\Gamma ^{(S)} \right)}{\left( i_0+\frac{1}{2} \right) \left( i_0+\frac{\gamma }{2} \right)} \frac{\left( -\frac{j}{2}\right)_{i_0} \left( -\frac{j}{2}+\frac{1}{2}\right)_{i_0}}{\left( 1 \right)_{i_0} \left( \frac{\gamma }{2}+ \frac{1}{2} \right)_{i_0}} \right.\nonumber\\
&&\times \prod_{k=1}^{\tau -1} \left( \sum_{i_k = i_{k-1}}^{m} \frac{\left( i_k+ \frac{k}{2}-\frac{j}{2}\right)\left(i_k + \frac{k}{2}+\Gamma ^{(S)} \right)}{\left( i_k+\frac{k}{2}+\frac{1}{2} \right) \left( i_k+\frac{k}{2}+\frac{\gamma }{2} \right)} \right. \nonumber\\
&&\times \left. \frac{\left( -\frac{j}{2}+\frac{k}{2}\right)_{i_k} \left( -\frac{j}{2}+ \frac{k}{2}+ \frac{1}{2}\right)_{i_k}\left(  \frac{k}{2} +1\right)_{i_{k-1}} \left( \frac{\gamma }{2}+ \frac{1}{2}+ \frac{k}{2} \right)_{i_{k-1}}}{\left( -\frac{j}{2}+\frac{k}{2}\right)_{i_{k-1}} \left( -\frac{j}{2}+ \frac{k}{2}+ \frac{1}{2}\right)_{i_{k-1}}\left( \frac{k}{2}+1 \right)_{i_k} \left( \frac{\gamma }{2}+ \frac{1}{2}+ \frac{k}{2} \right)_{i_k}} \right) \label{eq:40038c}\\
&&\times \left. \sum_{i_{\tau } = i_{\tau -1}}^{m} \frac{\left( -\frac{j}{2}+\frac{\tau }{2}\right)_{i_{\tau }} \left( -\frac{j}{2}+ \frac{\tau }{2}+ \frac{1}{2}\right)_{i_{\tau }}\left(  \frac{\tau }{2} +1\right)_{i_{\tau -1}} \left( \frac{\gamma }{2}+ \frac{1}{2}+\frac{\tau }{2} \right)_{i_{\tau -1}}}{\left( -\frac{j}{2}+\frac{\tau }{2}\right)_{i_{\tau -1}} \left(-\frac{j}{2}+ \frac{\tau }{2}+ \frac{1}{2}\right)_{i_{\tau -1}}\left(  \frac{\tau }{2}+1 \right)_{i_{\tau }} \left( \frac{\gamma }{2}+ \frac{1}{2}+\frac{\tau }{2} \right)_{i_{\tau }}} z^{i_{\tau }}\right\} \eta ^{\tau }  \nonumber
\end{eqnarray}
where
\begin{equation}
\begin{cases} \tau \geq 2 \cr
\Gamma ^{(S)} = \frac{1}{2(1+a)}\left( -\delta -j+1 +a\left( \gamma +\delta -1+j  \right)\right)  
\end{cases}\nonumber
\end{equation}
\end{enumerate}
\end{remark} 
\begin{remark}
The power series expansion of Heun equation of the second kind for the second species complete polynomial using 3TRF about $x=0$ is given by
\begin{enumerate} 
\item As $\alpha =\gamma -1 $, $\beta =\gamma $ and $q= ( \gamma -1) \left[ \gamma- \delta +a\delta \right]$,
 
Its eigenfunction is given by
\begin{eqnarray}
y(x) &=& H_pS_0 \left( a, q= ( \gamma -1) \left[ \gamma- \delta +a\delta \right]; \alpha =\gamma -1, \beta =\gamma, \gamma, \delta ; \eta = \frac{(1+a)}{a} x ; z= -\frac{1}{a} x^2 \right) \nonumber\\
&=& \eta ^{1-\gamma }  \label{eq:40039a}
\end{eqnarray}
\item As $\alpha =\gamma -2N-2 $, $\beta = \gamma -2N-1 $ and $q= \left(\gamma-2N-2 \right) \left[\gamma -\delta -2N-1 +a\left( \delta+2N+1 \right)\right]$ where $N \in \mathbb{N}_{0}$,

Its eigenfunction is given by
\begin{eqnarray} 
y(x) &=& H_pS_{2N+1} \Bigg( a, q= \left(\gamma -2N-2 \right) \left[\gamma -\delta -2N-1 +a\left( \delta+2N+1 \right)\right]; \alpha =\gamma -2N-2  \nonumber\\
&&, \beta =\gamma -2N-1, \gamma, \delta ; \eta = \frac{(1+a)}{a} x ; z= -\frac{1}{a} x^2 \Bigg) \nonumber\\
&=& \sum_{r=0}^{N} y_{2r}^{N-r}\left( 2N+1;x\right) + \sum_{r=0}^{N} y_{2r+1}^{N-r}\left( 2N+1;x\right) 
\label{eq:40039b}
\end{eqnarray}
\item As $\alpha =\gamma -2N-3 $, $\beta = \gamma -2N-2 $ and $q= \left(\gamma-2N-3 \right) \left[\gamma -\delta -2N-2 +a\left( \delta+2N+2 \right)\right]$ where $N \in \mathbb{N}_{0}$,

Its eigenfunction is given by
\begin{eqnarray} 
y(x) &=& H_pS_{2N+2} \Bigg( a, q= \left(\gamma -2N-3 \right) \left[\gamma -\delta -2N-2 +a\left( \delta+2N+2 \right)\right]; \alpha =\gamma -2N-3  \nonumber\\
&&, \beta =\gamma -2N-2, \gamma, \delta ; \eta = \frac{(1+a)}{a} x ; z= -\frac{1}{a} x^2 \Bigg) \nonumber\\
&=& \sum_{r=0}^{N+1} y_{2r}^{N+1-r}\left( 2N+2;x\right) + \sum_{r=0}^{N} y_{2r+1}^{N-r}\left( 2N+2;x\right) 
\label{eq:40039c}
\end{eqnarray}
In the above,
\begin{eqnarray}
y_0^m(j;x) &=& \eta ^{1-\gamma } \sum_{i_0=0}^{m} \frac{\left( -\frac{j}{2}\right)_{i_0} \left( - \frac{j}{2}+\frac{1}{2}\right)_{i_0}}{\left( \frac{3}{2}-\frac{\gamma}{2}\right)_{i_0} \left( 1\right)_{i_0}} z^{i_0} \label{eq:40040a}\\
y_1^m(j;x) &=& \eta ^{1-\gamma } \left\{\sum_{i_0=0}^{m} \frac{\left( i_0-\frac{j}{2}\right)\left(i_0 +\Gamma ^{(S)} \right)}{\left( i_0+1-\frac{\gamma}{2}\right) \left( i_0 +\frac{1}{2} \right)} \frac{\left( -\frac{j}{2}\right)_{i_0} \left( -\frac{j}{2}+\frac{1}{2}\right)_{i_0}}{\left( \frac{3}{2}-\frac{\gamma}{2}\right)_{i_0} \left( 1\right)_{i_0}} \right. \nonumber\\
&&\times \left. \sum_{i_1 = i_0}^{m} \frac{\left( -\frac{j}{2}+\frac{1}{2} \right)_{i_1} \left( -\frac{j}{2}+1\right)_{i_1}\left( 2-\frac{\gamma}{2}\right)_{i_0} \left( \frac{3}{2} \right)_{i_0}}{\left(-\frac{j}{2}+\frac{1}{2} \right)_{i_0} \left( -\frac{j}{2}+1\right)_{i_0}\left( 2-\frac{\gamma}{2}\right)_{i_1} \left( \frac{3}{2} \right)_{i_1}} z^{i_1}\right\} \eta \label{eq:40040b}\\
y_{\tau }^m(j;x) &=& \eta ^{1-\gamma } \left\{ \sum_{i_0=0}^{m} \frac{\left( i_0-\frac{j}{2}\right)\left( i_0 +\Gamma ^{(S)} \right)}{\left( i_0+1-\frac{\gamma}{2}\right) \left( i_0+ \frac{1}{2} \right)} \frac{\left( -\frac{j}{2}\right)_{i_0} \left( -\frac{j}{2}+\frac{1}{2}\right)_{i_0}}{\left( \frac{3}{2}-\frac{\gamma}{2}\right)_{i_0} \left( 1\right)_{i_0}} \right.\nonumber\\
&&\times \prod_{k=1}^{\tau -1} \left( \sum_{i_k = i_{k-1}}^{m} \frac{\left( i_k+ \frac{k}{2}-\frac{j}{2}\right)\left( i_k + \frac{k}{2}+\Gamma ^{(S)} \right)}{\left( i_k+\frac{k}{2}+1-\frac{\gamma}{2}\right) \left( i_k+\frac{k}{2} +\frac{1}{2} \right)} \right. \nonumber\\
&&\times \left. \frac{\left( -\frac{j}{2}+\frac{k}{2}\right)_{i_k} \left( -\frac{j}{2}+ \frac{k}{2}+ \frac{1}{2}\right)_{i_k}\left(  \frac{k}{2}+\frac{3}{2}-\frac{\gamma}{2}\right)_{i_{k-1}} \left( \frac{k}{2}+1 \right)_{i_{k-1}}}{\left( -\frac{j}{2}+\frac{k}{2}\right)_{i_{k-1}} \left( -\frac{j}{2}+ \frac{k}{2}+ \frac{1}{2}\right)_{i_{k-1}}\left( \frac{k}{2}+\frac{3}{2}-\frac{\gamma}{2}\right)_{i_k} \left( \frac{k}{2}+1 \right)_{i_k}} \right) \label{eq:40040c}\\
&&\times \left. \sum_{i_{\tau } = i_{\tau -1}}^{m} \frac{\left( -\frac{j}{2}+\frac{\tau }{2}\right)_{i_{\tau }} \left( -\frac{j}{2}+ \frac{\tau }{2}+ \frac{1}{2}\right)_{i_{\tau }}\left( \frac{\tau }{2}+\frac{3}{2}-\frac{\gamma}{2}\right)_{i_{\tau -1}} \left( \frac{\tau }{2}+ 1 \right)_{i_{\tau -1}}}{\left( -\frac{j}{2}+\frac{\tau }{2}\right)_{i_{\tau -1}} \left(-\frac{j}{2}+ \frac{\tau }{2}+ \frac{1}{2}\right)_{i_{\tau -1}}\left( \frac{\tau }{2}+\frac{3}{2}-\frac{\gamma}{2}\right)_{i_{\tau }} \left( \frac{\tau }{2}+1 \right)_{i_{\tau }}} z^{i_{\tau }}\right\} \eta ^{\tau } \nonumber
\end{eqnarray}
where
\begin{equation}
\begin{cases} \tau \geq 2 \cr
\Gamma ^{(S)} = \frac{1}{2(1+a)}\left( -\delta -j+1 +a\left( -\gamma +\delta +j+1 \right)\right) 
\end{cases}\nonumber
\end{equation}
\end{enumerate}
\end{remark} 
It is required that $\gamma \ne 0,-1,-2,\cdots$ for the first kind of independent solutions of Heun equation about $x=0$ for the first and second species complete polynomials. Because if it does not, its solutions will be divergent. By same reasons, it is required that $\gamma \ne 2,3,4, \cdots$ for the second kind of independent solutions of Heun equation about $x=0$ for the first and second species complete polynomials.
\section{Summary}

The power series expansion of Heun equation consists of a 3-term recursive relation between consecutive coefficients. Hypergeometric equation has the 2-term recursive relations between successive coefficients: the mathematical structures of it have been built analytically including its integral representation and the generating function of it. In contrast, the form of a power series in Heun equation is unknown currently involving its integral forms and the orthogonal relation of it because its recursive relation involves three terms. Instead of reducing 3-term in the recurrence relation into 2-term in order to describe Heun functions in the form of hypergeometric type, I construct the power series expansion in closed forms of Heun equation with a regular singularity at the origin by substituting a power series with unknown coefficients into Heun equation directly.

Power Series solutions and integral representations of Heun equation for an infinite series and a polynomial of type 1 are obtained by applying 3TRF in Ref.\cite{4Chou2012c,4Chou2012d}. Also, Frobenius solutions of Heun equation including its integral forms for an infinite series and a polynomial of type 2 are derived by applying R3TRF in chapter 2 of Ref.\cite{4Choun2013}. Infinite series solutions of Heun equation by applying 3TRF and R3TRF are equivalent to each other. For infinite series by applying 3TRF, $B_n$ is the leading term in sequence $c_n$ in each sub-power series of the analytic function $y(x)$. For infinite series by applying R3TRF, $A_n$ is the leading term in each sub-power series of the analytic function $y(x)$.

A polynomial which makes $A_n$ and $B_n$ terms terminated at the same time has two different types: (1) the first species complete polynomial and (2) the second species complete polynomial. Complete polynomials using 3TRF are derived  by allowing $A_n$ as the leading term in each sub-power series of the general power series $y(x)$. In the form of a power series of Heun equation around $x=0$ for the first species complete polynomial by using 3TRF in this chapter, I treat $\beta $, $\gamma $ and $\delta $ as free variables and fixed values of $\alpha $ and $q$: the spectral parameter $q$ has multi-valued roots of an eigenvalue. I construct the algebraic equation for the determination of the spectral parameter $q$ in the form of partial sums of the sequences $\{A_n\}$ and $\{B_n\}$ using 3TRF for computational practice rather than using the matrix formalism. In the Frobenius solution of Heun equation for the second species complete polynomial by using 3TRF, I treat $\gamma $ and $\delta $ as free variables and fixed values of $\alpha $,  $\beta $ and $q$: the parameter $q$ has only one valued root of each eigenvalue. 
  
For polynomials of type 1 and 2 in Heun equation, the general power series $y(x)$ consists of infinite numbers of finite sub-power series. For a polynomial of type 1, $A_n$ in sequence $c_n$ is the leading term in each finite sub-power series. For a polynomial of type 2, $B_n$ is the leading term in each finite sub-power series. In contrast, for the first and second species complete polynomials using 3TRF in Heun equation, the general power series $y(x)$ consists of finite numbers of finite sub-power series: $A_n$ is the leading term in each finite sub-power series.   

As we observe the power series expansions in closed forms of Heun equation for the first and second species complete polynomials, the denominators and numerators in all $B_n$ terms of each finite sub-power series arise with Pochhammer symbol. The differential expressions of hypergeometric (Jacobi) polynomial derive from converting ratio of a Pochhammer symbol in numerator to another Pochhammer symbol in denominator into the contour integral form. And a generating function and an orthogonal relation of Jacobi polynomial are analyzed from its differential forms. By utilizing a similar method, we are able to obtain differential representations, generating functions and orthogonal relations of Heun equation for complete polynomials of two types because of Pochhammer symbol in $B_n$ terms. 

\begin{appendices}
\section{Power series expansion of 192 Heun functions}
A machine-generated list of 192 (isomorphic to the Coxeter group of the Coxeter diagram $D_4$) local solutions of the Heun equation was obtained by Robert S. Maier(2007) \cite{4Maie2007}. 
In this appendix, replacing coefficients in the power series expansion of Heun equation of the first kind around $x=0$ for the first and second complete polynomials using 3TRF, I derive the Frobenius solutions for the first and second complete polynomials of nine out of the 192 local solutions of Heun function in Table 2 \cite{4Maie2007}.
\addtocontents{toc}{\protect\setcounter{tocdepth}{1}}
\subsection{${\displaystyle (1-x)^{1-\delta } Hl(a, q - (\delta  - 1)\gamma a; \alpha - \delta  + 1, \beta - \delta + 1, \gamma ,2 - \delta ; x)}$ }
\subsubsection{The first species complete polynomial}
Replacing coefficients $q$, $\alpha$, $\beta$ and $\delta$ by $q - (\delta - 1)\gamma a $, $\alpha - \delta  + 1 $, $\beta - \delta + 1$ and $2 - \delta$ into (\ref{eq:40018})--(\ref{eq:40022c}). Multiply $(1-x)^{1-\delta }$ and the new (\ref{eq:40018}), (\ref{eq:40019b}) and (\ref{eq:40020b})  together.\footnote{I treat $\beta $, $\gamma$ and $\delta$ as free variables and fixed values of $\alpha $ and $q$.} 
\begin{enumerate} 
\item As $\alpha  = \delta  -1$ and $q = (\delta - 1)\gamma a+q_0^0$ where $q_0^0=0$,

The eigenfunction is given by
\begin{eqnarray}
& &(1-x)^{1-\delta } y(x)\nonumber\\
&=& (1-x)^{1-\delta } Hl\left( a, 0; 0, \beta - \delta + 1, \gamma ,2 - \delta ; x\right)\nonumber\\
&=& (1-x)^{1-\delta } \nonumber
\end{eqnarray}
\item As $\alpha = \delta -2N-2$ where $N \in \mathbb{N}_{0}$,

An algebraic equation of degree $2N+2$ for the determination of $q$ is given by
\begin{equation}
0 = \sum_{r=0}^{N+1}\bar{c}\left( 2r, N+1-r; 2N+1,\tilde{q}\right)\nonumber
\end{equation}
The eigenvalue of $q$ is written by $(\delta  - 1)\gamma a+ q_{2N+1}^m$ where $m = 0,1,2,\cdots,2N+1 $; $q_{2N+1}^0 < q_{2N+1}^1 < \cdots < q_{2N+1}^{2N+1}$. Its eigenfunction is given by
\begin{eqnarray} 
& &(1-x)^{1-\delta } y(x)\nonumber\\
&=& (1-x)^{1-\delta } Hl\left( a, q_{2N+1}^m; -2N-1, \beta - \delta + 1, \gamma ,2 - \delta ; x\right)\nonumber\\
&=& (1-x)^{1-\delta } \left\{ \sum_{r=0}^{N} y_{2r}^{N-r}\left( 2N+1,q_{2N+1}^m;x\right)+ \sum_{r=0}^{N} y_{2r+1}^{N-r}\left( 2N+1,q_{2N+1}^m;x\right) \right\}  
\nonumber
\end{eqnarray}   
\item As $\alpha =\delta-2N-3$ where $N \in \mathbb{N}_{0}$,

An algebraic equation of degree $2N+3$ for the determination of $q$ is given by
\begin{eqnarray}
0  = \sum_{r=0}^{N+1}\bar{c}\left( 2r+1, N+1-r; 2N+2,\tilde{q}\right)\nonumber
\end{eqnarray}
The eigenvalue of $q$ is written by $(\delta  - 1)\gamma a+ q_{2N+2}^m$ where $m = 0,1,2,\cdots,2N+2 $; $q_{2N+2}^0 < q_{2N+2}^1 < \cdots < q_{2N+2}^{2N+2}$. Its eigenfunction is given by
\begin{eqnarray} 
& &(1-x)^{1-\delta } y(x)\nonumber\\
&=& (1-x)^{1-\delta } Hl\left( a, q_{2N+2}^m; -2N-2, \beta - \delta + 1, \gamma ,2 - \delta ; x\right)\nonumber\\
&=& (1-x)^{1-\delta } \left\{ \sum_{r=0}^{N+1} y_{2r}^{N+1-r}\left( 2N+2,q_{2N+2}^m;x\right) + \sum_{r=0}^{N} y_{2r+1}^{N-r}\left( 2N+2,q_{2N+2}^m;x\right) \right\}   
\nonumber
\end{eqnarray}
In the above,
\begin{eqnarray}
\bar{c}(0,n;j,\tilde{q})  &=& \frac{\left( -\frac{j}{2}\right)_{n} \left( \frac{\beta}{2} -\frac{\delta}{2} +\frac{1}{2} \right)_{n}}{\left( 1 \right)_{n} \left( \frac{\gamma }{2}+ \frac{1}{2} \right)_{n}} \left(-\frac{1}{a} \right)^{n}\nonumber\\
\bar{c}(1,n;j,\tilde{q}) &=& \left( \frac{1+a}{a}\right) \sum_{i_0=0}^{n}\frac{ i_0 \left(i_0 + \Gamma _0^{(F)} \right)+\frac{\tilde{q}}{4(1+a)}}{\left( i_0+\frac{1}{2} \right) \left( i_0+\frac{\gamma }{2} \right)} \frac{\left( -\frac{j}{2}\right)_{i_0} \left( \frac{\beta }{2}-\frac{\delta}{2}+\frac{1}{2} \right)_{i_0}}{\left( 1 \right)_{i_0} \left( \frac{\gamma }{2}+ \frac{1}{2} \right)_{i_0}} \nonumber\\
&&\times  \frac{\left( -\frac{j}{2}+\frac{1}{2} \right)_{n} \left( \frac{\beta }{2}-\frac{\delta}{2}+1 \right)_{n}\left( \frac{3}{2} \right)_{i_0} \left( \frac{\gamma }{2}+1 \right)_{i_0}}{\left( -\frac{j}{2}+\frac{1}{2}\right)_{i_0} \left( \frac{\beta }{2}-\frac{\delta}{2}+1 \right)_{i_0}\left( \frac{3}{2} \right)_{n} \left( \frac{\gamma }{2}+1 \right)_{n}} \left(-\frac{1}{a} \right)^{n}  
\nonumber\\
\bar{c}(\tau ,n;j,\tilde{q}) &=& \left( \frac{1+a}{a}\right)^{\tau} \sum_{i_0=0}^{n}\frac{ i_0 \left( i_0 +\Gamma _0^{(F)} \right) +\frac{\tilde{q}}{4(1+a)}}{\left( i_0+\frac{1}{2} \right) \left( i_0+\frac{\gamma }{2} \right)} \frac{\left( -\frac{j}{2}\right)_{i_0} \left( \frac{\beta }{2}-\frac{\delta}{2}+\frac{1}{2} \right)_{i_0}}{\left( 1 \right)_{i_0} \left( \frac{\gamma }{2}+ \frac{1}{2} \right)_{i_0}} \nonumber\\
&&\times \prod_{k=1}^{\tau -1} \left( \sum_{i_k = i_{k-1}}^{n} \frac{\left( i_k+ \frac{k}{2} \right)\left( i_k +\Gamma _k^{(F)} \right) +\frac{\tilde{q}}{4(1+a)}}{\left( i_k+\frac{k}{2}+\frac{1}{2} \right) \left( i_k+\frac{k}{2}+\frac{\gamma }{2} \right)} \right. \nonumber\\
&&\times \left. \frac{\left( -\frac{j}{2}+\frac{k}{2}\right)_{i_k} \left( \frac{\beta }{2}-\frac{\delta}{2}+ \frac{k}{2}+\frac{1}{2} \right)_{i_k}\left( \frac{k}{2}+1 \right)_{i_{k-1}} \left( \frac{\gamma }{2}+ \frac{1}{2}+ \frac{k}{2} \right)_{i_{k-1}}}{\left( -\frac{j}{2}+\frac{k}{2}\right)_{i_{k-1}} \left( \frac{\beta }{2}-\frac{\delta}{2}+ \frac{k}{2}+\frac{1}{2} \right)_{i_{k-1}}\left( \frac{k}{2} +1\right)_{i_k} \left( \frac{\gamma }{2}+ \frac{1}{2}+ \frac{k}{2} \right)_{i_k}} \right) \nonumber\\
&&\times \frac{\left( -\frac{j}{2}+\frac{\tau }{2}\right)_{n} \left( \frac{\beta }{2}-\frac{\delta}{2}+ \frac{\tau }{2} +\frac{1}{2}\right)_{n}\left( \frac{\tau }{2} +1\right)_{i_{\tau -1}} \left( \frac{\gamma }{2}+ \frac{1}{2}+\frac{\tau }{2} \right)_{i_{\tau -1}}}{\left( -\frac{j}{2}+\frac{\tau }{2}\right)_{i_{\tau -1}} \left( \frac{\beta }{2}-\frac{\delta}{2}+\frac{\tau }{2} +\frac{1}{2} \right)_{i_{\tau -1}}\left( \frac{\tau }{2}+1 \right)_{n} \left( \frac{\gamma }{2}+ \frac{1}{2}+\frac{\tau }{2} \right)_{n}} \left(-\frac{1}{a} \right)^{n } \nonumber
\end{eqnarray}
\begin{eqnarray}
y_0^m(j, \tilde{q}; x) &=& \sum_{i_0=0}^{m} \frac{\left( -\frac{j}{2}\right)_{i_0} \left( \frac{\beta }{2} -\frac{\delta}{2}+\frac{1}{2}\right)_{i_0}}{\left( 1 \right)_{i_0} \left( \frac{\gamma }{2}+ \frac{1}{2} \right)_{i_0}} z^{i_0} \nonumber\\
y_1^m(j,\tilde{q}; x) &=& \left\{\sum_{i_0=0}^{m} \frac{ i_0 \left( i_0 +\Gamma _0^{(F)} \right) +\frac{\tilde{q}}{4(1+a)}}{\left( i_0+\frac{1}{2} \right) \left( i_0+\frac{\gamma }{2} \right)} \frac{\left( -\frac{j}{2}\right)_{i_0} \left( \frac{\beta }{2}-\frac{\delta}{2}+\frac{1}{2} \right)_{i_0}}{\left( 1 \right)_{i_0} \left( \frac{\gamma }{2}+ \frac{1}{2} \right)_{i_0}} \right. \nonumber\\
&&\times \left. \sum_{i_1 = i_0}^{m} \frac{\left( -\frac{j}{2}+\frac{1}{2} \right)_{i_1} \left( \frac{\beta }{2}-\frac{\delta}{2}+ 1 \right)_{i_1}\left( \frac{3}{2} \right)_{i_0} \left( \frac{\gamma }{2}+1 \right)_{i_0}}{\left(-\frac{j}{2}+\frac{1}{2} \right)_{i_0} \left( \frac{\beta }{2}-\frac{\delta}{2}+ 1 \right)_{i_0}\left( \frac{3}{2} \right)_{i_1} \left( \frac{\gamma }{2} +1 \right)_{i_1}} z^{i_1}\right\} \eta 
\nonumber\\
y_{\tau }^m(j,\tilde{q}; x) &=& \left\{ \sum_{i_0=0}^{m} \frac{ i_0 \left( i_0 +\Gamma _0^{(F)} \right) +\frac{\tilde{q}}{4(1+a)}}{\left( i_0+\frac{1}{2} \right) \left( i_0+\frac{\gamma }{2} \right)} \frac{\left( -\frac{j}{2}\right)_{i_0} \left( \frac{\beta }{2}-\frac{\delta}{2}+\frac{1}{2} \right)_{i_0}}{\left( 1 \right)_{i_0} \left( \frac{\gamma }{2}+ \frac{1}{2} \right)_{i_0}} \right.\nonumber\\
&&\times \prod_{k=1}^{\tau -1} \left( \sum_{i_k = i_{k-1}}^{m} \frac{\left( i_k+ \frac{k}{2} \right)\left( i_k +\Gamma _k^{(F)} \right) +\frac{\tilde{q}}{4(1+a)}}{\left( i_k+\frac{k}{2}+\frac{1}{2} \right) \left( i_k+\frac{k}{2}+\frac{\gamma }{2} \right)} \right. \nonumber\\
&&\times \left. \frac{\left( -\frac{j}{2}+\frac{k}{2}\right)_{i_k} \left( \frac{\beta }{2}-\frac{\delta}{2}+ \frac{k}{2}+\frac{1}{2} \right)_{i_k}\left( \frac{k}{2}+1 \right)_{i_{k-1}} \left( \frac{\gamma }{2}+ \frac{1}{2}+ \frac{k}{2} \right)_{i_{k-1}}}{\left( -\frac{j}{2}+\frac{k}{2}\right)_{i_{k-1}} \left( \frac{\beta }{2}-\frac{\delta}{2}+ \frac{k}{2}+\frac{1}{2} \right)_{i_{k-1}}\left( \frac{k}{2} +1\right)_{i_k} \left( \frac{\gamma }{2}+ \frac{1}{2}+ \frac{k}{2} \right)_{i_k}} \right) \nonumber\\
&&\times \left. \sum_{i_{\tau } = i_{\tau -1}}^{m} \frac{\left( -\frac{j}{2}+\frac{\tau }{2}\right)_{i_{\tau }} \left( \frac{\beta }{2}-\frac{\delta}{2}+ \frac{\tau }{2}+\frac{1}{2} \right)_{i_{\tau }}\left( \frac{\tau }{2}+1 \right)_{i_{\tau -1}} \left( \frac{\gamma }{2}+ \frac{1}{2}+\frac{\tau }{2} \right)_{i_{\tau -1}}}{\left( -\frac{j}{2}+\frac{\tau }{2}\right)_{i_{\tau -1}} \left( \frac{\beta }{2}-\frac{\delta}{2}+ \frac{\tau }{2}+\frac{1}{2} \right)_{i_{\tau -1}}\left( \frac{\tau }{2}+1 \right)_{i_{\tau }} \left( \frac{\gamma }{2}+ \frac{1}{2}+\frac{\tau }{2} \right)_{i_{\tau }}} z^{i_{\tau }}\right\} \eta ^{\tau } \nonumber
\end{eqnarray}
where
\begin{equation}
\begin{cases} \tau \geq 2 \cr
z = -\frac{1}{a}x^2 \cr
\eta = \frac{(1+a)}{a} x \cr
\tilde{q} = q-(\delta -1)\gamma a \cr
\Gamma _0^{(F)} = \frac{1}{2(1+a)}\left( \beta -j-1 +a\left( \gamma -\delta +1 \right)\right) \cr
\Gamma _k^{(F)} = \frac{1}{2(1+a)}\left( \beta -j-1+k +a\left( \gamma -\delta +1+k \right)\right)
\end{cases}\nonumber
\end{equation}
\end{enumerate}
\subsubsection{The second species complete polynomial}
Replacing coefficients $q$, $\alpha$, $\beta$ and $\delta$ by $q - (\delta - 1)\gamma a $, $\alpha - \delta  + 1 $, $\beta - \delta + 1$ and $2 - \delta$ into (\ref{eq:40037a})--(\ref{eq:40038c}). Multiply $(1-x)^{1-\delta }$ and the new (\ref{eq:40037a})--(\ref{eq:40037c}) together.\footnote{I treat $\gamma$ and $\delta$ as free variables and fixed values of $\alpha $, $\beta $ and $q$.}  
\begin{enumerate} 
\item As $\alpha  = \delta -1 $, $\beta   =\delta $ and $q = (\delta  - 1)\gamma a$,

Its eigenfunction is given by
\begin{eqnarray}
& &(1-x)^{1-\delta } y(x)\nonumber\\
&=& (1-x)^{1-\delta } Hl\left( a, 0; 0, 0, \gamma ,2 - \delta ; x\right)\nonumber\\
&=& (1-x)^{1-\delta }  \nonumber
\end{eqnarray}
\item As $\alpha = \delta -2N-2 $, $\beta = \delta -2N -1$ and $q = (\delta  - 1)\gamma a+ \left(2N+1 \right) \left[ -\delta +2N +2-a\left( \gamma -\delta+2N +2\right)\right]$ where $N \in \mathbb{N}_{0}$,

Its eigenfunction is given by
\begin{eqnarray} 
& &(1-x)^{1-\delta } y(x)\nonumber\\
&=& (1-x)^{1-\delta } Hl\left( a, \left( 2N+1 \right) \left[ -\delta +2N +2-a\left( \gamma -\delta+2N +2\right)\right]; -2N-1, -2N, \gamma ,2 - \delta ; x\right)\nonumber\\
&=& (1-x)^{1-\delta } \left\{ \sum_{r=0}^{N} y_{2r}^{N-r}\left( 2N+1;x\right) + \sum_{r=0}^{N} y_{2r+1}^{N-r}\left( 2N+1;x\right)   \right\} 
\nonumber
\end{eqnarray}
\item As $\alpha = \delta -2N-3 $, $\beta = \delta -2N-2 $ and $q = (\delta  - 1)\gamma a +\left(2N+2 \right) \left[ -\delta +2N+3 -a\left( \gamma -\delta +2N+3 \right)\right]$ where $N \in \mathbb{N}_{0}$,

Its eigenfunction is given by
\begin{eqnarray} 
& &(1-x)^{1-\delta } y(x)\nonumber\\
&=& (1-x)^{1-\delta } Hl\left( a, \left( 2N+2\right) \left[ -\delta +2N+3 -a\left( \gamma -\delta +2N+3 \right)\right]; -2N-2, -2N-1, \gamma \right. \nonumber\\
&&, \left. 2 - \delta ; x\right)\nonumber\\
&=& (1-x)^{1-\delta } \left\{  \sum_{r=0}^{N+1} y_{2r}^{N+1-r}\left( 2N+2;x\right) + \sum_{r=0}^{N} y_{2r+1}^{N-r}\left( 2N+2;x\right)  \right\} 
\nonumber
\end{eqnarray}
In the above,
\begin{eqnarray}
y_0^m(j;x) &=&  \sum_{i_0=0}^{m} \frac{\left( -\frac{j}{2}\right)_{i_0} \left( - \frac{j}{2}+\frac{1}{2}\right)_{i_0}}{\left( 1 \right)_{i_0} \left( \frac{\gamma }{2}+ \frac{1}{2} \right)_{i_0}} z^{i_0} \nonumber\\
y_1^m(j;x) &=&  \left\{\sum_{i_0=0}^{m} \frac{\left( i_0-\frac{j}{2}\right)\left(i_0 +\Gamma ^{(S)} \right)}{\left( i_0+\frac{1}{2} \right) \left( i_0+\frac{\gamma }{2} \right)} \frac{\left( -\frac{j}{2}\right)_{i_0} \left( -\frac{j}{2}+\frac{1}{2}\right)_{i_0}}{\left( 1 \right)_{i_0} \left( \frac{\gamma }{2}+ \frac{1}{2} \right)_{i_0}} \right. \nonumber\\
&&\times \left. \sum_{i_1 = i_0}^{m} \frac{\left( -\frac{j}{2}+\frac{1}{2} \right)_{i_1} \left( -\frac{j}{2}+1\right)_{i_1}\left( \frac{3}{2} \right)_{i_0} \left( \frac{\gamma }{2}+1 \right)_{i_0}}{\left(-\frac{j}{2}+\frac{1}{2} \right)_{i_0} \left( -\frac{j}{2}+1\right)_{i_0}\left( \frac{3}{2} \right)_{i_1} \left( \frac{\gamma }{2} +1 \right)_{i_1}} z^{i_1}\right\} \eta \nonumber
\end{eqnarray}
\begin{eqnarray}
y_{\tau }^m(j;x) &=&  \left\{ \sum_{i_0=0}^{m} \frac{\left( i_0-\frac{j}{2}\right)\left(i_0 +\Gamma ^{(S)} \right)}{\left( i_0+\frac{1}{2} \right) \left( i_0+\frac{\gamma }{2} \right)} \frac{\left( -\frac{j}{2}\right)_{i_0} \left( -\frac{j}{2}+\frac{1}{2}\right)_{i_0}}{\left( 1 \right)_{i_0} \left( \frac{\gamma }{2}+ \frac{1}{2} \right)_{i_0}} \right.\nonumber\\
&&\times \prod_{k=1}^{\tau -1} \left( \sum_{i_k = i_{k-1}}^{m} \frac{\left( i_k+ \frac{k}{2}-\frac{j}{2}\right)\left(i_k + \frac{k}{2}+\Gamma ^{(S)} \right)}{\left( i_k+\frac{k}{2}+\frac{1}{2} \right) \left( i_k+\frac{k}{2}+\frac{\gamma }{2} \right)} \right. \nonumber\\
&&\times \left. \frac{\left( -\frac{j}{2}+\frac{k}{2}\right)_{i_k} \left( -\frac{j}{2}+ \frac{k}{2}+ \frac{1}{2}\right)_{i_k}\left(  \frac{k}{2} +1\right)_{i_{k-1}} \left( \frac{\gamma }{2}+ \frac{1}{2}+ \frac{k}{2} \right)_{i_{k-1}}}{\left( -\frac{j}{2}+\frac{k}{2}\right)_{i_{k-1}} \left( -\frac{j}{2}+ \frac{k}{2}+ \frac{1}{2}\right)_{i_{k-1}}\left( \frac{k}{2}+1 \right)_{i_k} \left( \frac{\gamma }{2}+ \frac{1}{2}+ \frac{k}{2} \right)_{i_k}} \right) \nonumber\\
&&\times \left. \sum_{i_{\tau } = i_{\tau -1}}^{m} \frac{\left( -\frac{j}{2}+\frac{\tau }{2}\right)_{i_{\tau }} \left( -\frac{j}{2}+ \frac{\tau }{2}+ \frac{1}{2}\right)_{i_{\tau }}\left(  \frac{\tau }{2} +1\right)_{i_{\tau -1}} \left( \frac{\gamma }{2}+ \frac{1}{2}+\frac{\tau }{2} \right)_{i_{\tau -1}}}{\left( -\frac{j}{2}+\frac{\tau }{2}\right)_{i_{\tau -1}} \left(-\frac{j}{2}+ \frac{\tau }{2}+ \frac{1}{2}\right)_{i_{\tau -1}}\left(  \frac{\tau }{2}+1 \right)_{i_{\tau }} \left( \frac{\gamma }{2}+ \frac{1}{2}+\frac{\tau }{2} \right)_{i_{\tau }}} z^{i_{\tau }}\right\} \eta ^{\tau } \nonumber 
\end{eqnarray}
where
\begin{equation}
\begin{cases} \tau \geq 2 \cr
z = -\frac{1}{a}x^2 \cr
\eta = \frac{(1+a)}{a} x \cr
\Gamma ^{(S)} = \frac{1}{2(1+a)}\left( \delta -j-1 +a\left( \gamma -\delta +j+1 \right)\right)
\end{cases}\nonumber
\end{equation}
\end{enumerate} 
\subsection{\footnotesize ${\displaystyle x^{1-\gamma } (1-x)^{1-\delta }Hl(a, q-(\gamma +\delta -2)a -(\gamma -1)(\alpha +\beta -\gamma -\delta +1), \alpha - \gamma -\delta +2, \beta - \gamma -\delta +2, 2-\gamma, 2 - \delta ; x)}$ \normalsize}
\subsubsection{The first species complete polynomial}
Replacing coefficients $q$, $\alpha$, $\beta$, $\gamma $ and $\delta$ by $q-(\gamma +\delta -2)a-(\gamma -1)(\alpha +\beta -\gamma -\delta +1)$, $\alpha - \gamma -\delta +2$, $\beta - \gamma -\delta +2, 2-\gamma$ and $2 - \delta$ into (\ref{eq:40018})--(\ref{eq:40022c}). Multiply $x^{1-\gamma } (1-x)^{1-\delta }$ and the new (\ref{eq:40018}), (\ref{eq:40019b}) and (\ref{eq:40020b})  together.\footnote{I treat $\beta $, $\gamma$ and $\delta$ as free variables and fixed values of $\alpha $ and $q$.} 
\begin{enumerate} 
\item As $\alpha = \gamma +\delta -2 $ and $q=(\gamma +\delta -2)a+(\gamma -1)( \beta -1)+q_0^0$ where $q_0^0=0$,

The eigenfunction is given by
\begin{eqnarray}
& &x^{1-\gamma } (1-x)^{1-\delta } y(x)\nonumber\\
&=& x^{1-\gamma } (1-x)^{1-\delta } Hl(a, 0; 0, \beta - \gamma -\delta +2, 2-\gamma, 2 - \delta ; x)\nonumber\\
&=& x^{1-\gamma } (1-x)^{1-\delta } \nonumber
\end{eqnarray}
\item As $\alpha = \gamma +\delta -2N-3$ where $N \in \mathbb{N}_{0}$,

An algebraic equation of degree $2N+2$ for the determination of $q$ is given by
\begin{equation}
0 = \sum_{r=0}^{N+1}\bar{c}\left( 2r, N+1-r; 2N+1, q-(\gamma +\delta -2)a -(\gamma -1)( \beta -2N-2)\right)  \nonumber
\end{equation}
The eigenvalue of $q$ is written by $(\gamma +\delta -2)a+(\gamma -1)( \beta -2N-2)+q_{2N+1}^m$ where $m = 0,1,2,\cdots,2N+1 $; $q_{2N+1}^0 < q_{2N+1}^1 < \cdots < q_{2N+1}^{2N+1}$. Its eigenfunction is given by
\begin{eqnarray} 
 & &x^{1-\gamma } (1-x)^{1-\delta } y(x)\nonumber\\
&=& x^{1-\gamma } (1-x)^{1-\delta } Hl(a, q_{2N+1}^m; -2N-1, \beta - \gamma -\delta +2, 2-\gamma, 2 - \delta ; x)\nonumber\\
&=& x^{1-\gamma } (1-x)^{1-\delta } \left\{ \sum_{r=0}^{N} y_{2r}^{N-r}\left( 2N+1,q_{2N+1}^m;x\right)+ \sum_{r=0}^{N} y_{2r+1}^{N-r}\left( 2N+1,q_{2N+1}^m;x\right)  \right\}  
\nonumber
\end{eqnarray}
\item As $\alpha = \gamma +\delta -2N-4$ where $N \in \mathbb{N}_{0}$,

An algebraic equation of degree $2N+3$ for the determination of $q$ is given by
\begin{equation}
0 =  \sum_{r=0}^{N+1}\bar{c}\left( 2r+1, N+1-r; 2N+2, q-(\gamma +\delta -2)a -(\gamma -1)( \beta -2N-3)\right)  \nonumber
\end{equation}
The eigenvalue of $q$ is written by $(\gamma +\delta -2)a +(\gamma -1)(\beta -2N-3)+q_{2N+2}^m$ where $m = 0,1,2,\cdots,2N+2 $; $q_{2N+2}^0 < q_{2N+2}^1 < \cdots < q_{2N+2}^{2N+2}$. Its eigenfunction is given by
\begin{eqnarray} 
& &x^{1-\gamma } (1-x)^{1-\delta } y(x)\nonumber\\
&=& x^{1-\gamma } (1-x)^{1-\delta } Hl(a, q_{2N+2}^m; -2N-2, \beta - \gamma -\delta +2, 2-\gamma, 2 - \delta ; x)\nonumber\\
&=& x^{1-\gamma } (1-x)^{1-\delta } \left\{ \sum_{r=0}^{N+1} y_{2r}^{N+1-r}\left( 2N+2,q_{2N+2}^m;x\right) + \sum_{r=0}^{N} y_{2r+1}^{N-r}\left( 2N+2,q_{2N+2}^m;x\right) \right\} 
\nonumber
\end{eqnarray}
In the above,
\begin{eqnarray}
\bar{c}(0,n;j,\tilde{q})  &=& \frac{\left( -\frac{j}{2}\right)_{n} \left( \frac{\beta }{2}-\frac{\gamma }{2}-\frac{\delta }{2}+1 \right)_{n}}{\left( 1 \right)_{n} \left( -\frac{\gamma }{2}+ \frac{3}{2} \right)_{n}} \left(-\frac{1}{a} \right)^{n}\nonumber\\
\bar{c}(1,n;j,\tilde{q}) &=& \left( \frac{1+a}{a}\right) \sum_{i_0=0}^{n}\frac{ i_0 \left(i_0 +\Gamma _0^{(F)} \right)+\frac{\tilde{q}}{4(1+a)}}{\left( i_0+\frac{1}{2} \right) \left( i_0+1-\frac{ \gamma }{2} \right)}  \frac{\left( -\frac{j}{2}\right)_{i_0} \left( \frac{\beta }{2}-\frac{\gamma }{2}-\frac{\delta }{2} +1\right)_{i_0}}{\left( 1 \right)_{i_0} \left( -\frac{\gamma }{2}+ \frac{3}{2} \right)_{i_0}}\nonumber\\
&&\times  \frac{\left( -\frac{j}{2}+\frac{1}{2} \right)_{n} \left( \frac{\beta }{2}-\frac{\gamma }{2}-\frac{\delta }{2}+ \frac{3}{2} \right)_{n}\left( \frac{3}{2} \right)_{i_0} \left( -\frac{\gamma }{2}+2 \right)_{i_0}}{\left( -\frac{j}{2}+\frac{1}{2}\right)_{i_0} \left( \frac{\beta }{2}-\frac{\gamma }{2}-\frac{\delta }{2}+\frac{3}{2} \right)_{i_0}\left( \frac{3}{2} \right)_{n} \left( -\frac{\gamma }{2}+2 \right)_{n}} \left(-\frac{1}{a} \right)^{n}  
\nonumber\\
\bar{c}(\tau ,n;j,\tilde{q}) &=& \left( \frac{1+a}{a}\right)^{\tau} \sum_{i_0=0}^{n}\frac{ i_0 \left(i_0 +\Gamma _0^{(F)} \right)+\frac{\tilde{q}}{4(1+a)}}{\left( i_0+\frac{1}{2} \right) \left( i_0+1-\frac{\gamma }{2} \right)} \frac{\left( -\frac{j}{2}\right)_{i_0} \left( \frac{\beta }{2}-\frac{\gamma }{2}-\frac{\delta }{2}+1 \right)_{i_0}}{\left( 1 \right)_{i_0} \left( -\frac{\gamma }{2}+ \frac{3}{2} \right)_{i_0}} \nonumber\\
&&\times \prod_{k=1}^{\tau -1} \left( \sum_{i_k = i_{k-1}}^{n} \frac{\left( i_k+ \frac{k}{2} \right)\left(i_k +\Gamma _k^{(F)} \right)+\frac{\tilde{q}}{4(1+a)}}{\left( i_k+\frac{k}{2}+\frac{1}{2} \right) \left( i_k+\frac{k}{2}+1-\frac{\gamma }{2} \right)} \right. \nonumber\\
&&\times \left. \frac{\left( -\frac{j}{2}+\frac{k}{2}\right)_{i_k} \left( \frac{\beta }{2}-\frac{\gamma }{2}-\frac{\delta }{2}+1+ \frac{k}{2} \right)_{i_k}\left( \frac{k}{2}+1 \right)_{i_{k-1}} \left( -\frac{\gamma }{2}+ \frac{3}{2}+ \frac{k}{2} \right)_{i_{k-1}}}{\left( -\frac{j}{2}+\frac{k}{2}\right)_{i_{k-1}} \left( \frac{\beta }{2}-\frac{\gamma }{2}-\frac{\delta }{2}+1+ \frac{k}{2} \right)_{i_{k-1}}\left( \frac{k}{2} +1\right)_{i_k} \left( -\frac{\gamma }{2}+ \frac{3}{2}+ \frac{k}{2} \right)_{i_k}} \right) \nonumber\\
&&\times  \frac{\left( -\frac{j}{2}+\frac{\tau }{2}\right)_{n} \left( \frac{\beta }{2}-\frac{\gamma }{2}-\frac{\delta }{2}+1+ \frac{\tau }{2} \right)_{n}\left( \frac{\tau }{2} +1\right)_{i_{\tau -1}} \left( -\frac{\gamma }{2}+ \frac{3}{2}+\frac{\tau }{2} \right)_{i_{\tau -1}}}{\left( -\frac{j}{2}+\frac{\tau }{2}\right)_{i_{\tau -1}} \left( \frac{\beta }{2}-\frac{\gamma }{2}-\frac{\delta }{2}+1+\frac{\tau }{2} \right)_{i_{\tau -1}}\left( \frac{\tau }{2}+1 \right)_{n} \left( -\frac{\gamma }{2}+ \frac{3}{2}+\frac{\tau }{2} \right)_{n}} \left(-\frac{1}{a} \right)^{n } \nonumber
\end{eqnarray}
\begin{eqnarray}
y_0^m(j,\tilde{q};x) &=& \sum_{i_0=0}^{m} \frac{\left( -\frac{j}{2}\right)_{i_0} \left( \frac{\beta }{2}-\frac{\gamma }{2}-\frac{\delta }{2}+1 \right)_{i_0}}{\left( 1 \right)_{i_0} \left( -\frac{\gamma }{2}+ \frac{3}{2} \right)_{i_0}} z^{i_0} \nonumber\\
y_1^m(j,\tilde{q};x) &=& \left\{\sum_{i_0=0}^{m} \frac{ i_0 \left(i_0 +\Gamma _0^{(F)} \right)+\frac{\tilde{q}}{4(1+a)}}{\left( i_0+\frac{1}{2} \right) \left( i_0+1-\frac{\gamma }{2} \right)} \frac{\left( -\frac{j}{2}\right)_{i_0} \left( \frac{\beta }{2}-\frac{\gamma }{2}-\frac{\delta }{2}+1 \right)_{i_0}}{\left( 1 \right)_{i_0} \left( -\frac{\gamma }{2}+ \frac{3}{2} \right)_{i_0}} \right. \nonumber\\
&&\times \left. \sum_{i_1 = i_0}^{m} \frac{\left( -\frac{j}{2}+\frac{1}{2} \right)_{i_1} \left( \frac{\beta }{2}-\frac{\gamma }{2}-\frac{\delta }{2}+ \frac{3}{2} \right)_{i_1}\left( \frac{3}{2} \right)_{i_0} \left( -\frac{\gamma }{2}+2 \right)_{i_0}}{\left(-\frac{j}{2}+\frac{1}{2} \right)_{i_0} \left( \frac{\beta }{2}-\frac{\gamma }{2}-\frac{\delta }{2}+ \frac{3}{2} \right)_{i_0}\left( \frac{3}{2} \right)_{i_1} \left( -\frac{\gamma }{2} +2 \right)_{i_1}} z^{i_1}\right\} \eta 
\nonumber
\end{eqnarray}
\begin{eqnarray}
y_{\tau }^m(j,\tilde{q};x) &=& \left\{ \sum_{i_0=0}^{m} \frac{ i_0 \left(i_0 +\Gamma _0^{(F)} \right)+\frac{\tilde{q}}{4(1+a)}}{\left( i_0+\frac{1}{2} \right) \left( i_0+1-\frac{\gamma }{2} \right)} \frac{\left( -\frac{j}{2}\right)_{i_0} \left( \frac{\beta }{2}-\frac{\gamma }{2}-\frac{\delta }{2}+1 \right)_{i_0}}{\left( 1 \right)_{i_0} \left( -\frac{\gamma }{2}+ \frac{3}{2} \right)_{i_0}} \right.\nonumber\\
&&\times \prod_{k=1}^{\tau -1} \left( \sum_{i_k = i_{k-1}}^{m} \frac{\left( i_k+ \frac{k}{2} \right)\left(i_k +\Gamma _k^{(F)} \right)+\frac{\tilde{q}}{4(1+a)}}{\left( i_k+\frac{k}{2}+\frac{1}{2} \right) \left( i_k+\frac{k}{2}+1-\frac{\gamma }{2} \right)} \right. \nonumber\\
&&\times \left. \frac{\left( -\frac{j}{2}+\frac{k}{2}\right)_{i_k} \left( \frac{\beta }{2}-\frac{\gamma }{2}-\frac{\delta }{2}+1+ \frac{k}{2} \right)_{i_k}\left( \frac{k}{2}+1 \right)_{i_{k-1}} \left( -\frac{\gamma }{2}+ \frac{3}{2}+ \frac{k}{2} \right)_{i_{k-1}}}{\left( -\frac{j}{2}+\frac{k}{2}\right)_{i_{k-1}} \left( \frac{\beta }{2}-\frac{\gamma }{2}-\frac{\delta }{2}+1+ \frac{k}{2} \right)_{i_{k-1}}\left( \frac{k}{2} +1\right)_{i_k} \left( -\frac{\gamma }{2}+ \frac{3}{2}+ \frac{k}{2} \right)_{i_k}} \right) \nonumber\\
&&\times \left. \sum_{i_{\tau } = i_{\tau -1}}^{m} \frac{\left( -\frac{j}{2}+\frac{\tau }{2}\right)_{i_{\tau }} \left( \frac{\beta }{2}-\frac{\gamma }{2}-\frac{\delta }{2}+1+ \frac{\tau }{2} \right)_{i_{\tau }}\left( \frac{\tau }{2}+1 \right)_{i_{\tau -1}} \left( -\frac{\gamma }{2}+ \frac{3}{2}+\frac{\tau }{2} \right)_{i_{\tau -1}}}{\left( -\frac{j}{2}+\frac{\tau }{2}\right)_{i_{\tau -1}} \left( \frac{\beta }{2}-\frac{\gamma }{2}-\frac{\delta }{2}+1+ \frac{\tau }{2} \right)_{i_{\tau -1}}\left( \frac{\tau }{2}+1 \right)_{i_{\tau }} \left( -\frac{\gamma }{2}+ \frac{3}{2}+\frac{\tau }{2} \right)_{i_{\tau }}} z^{i_{\tau }}\right\} \eta ^{\tau } \nonumber
\end{eqnarray}
where
\begin{equation}
\begin{cases} \tau \geq 2 \cr
z = -\frac{1}{a}x^2 \cr
\eta = \frac{(1+a)}{a} x \cr
\tilde{q} = q-(\gamma +\delta -2)a -(\gamma -1)(\alpha +\beta -\gamma -\delta +1) \cr
\Gamma _0^{(F)} =  \frac{1}{2(1+a)}\left( \beta -\gamma -j +a\left( -\gamma -\delta +3 \right)\right) \cr
\Gamma _k^{(F)} = \frac{1}{2(1+a)}\left( \beta -\gamma -j+k +a\left( -\gamma -\delta +3+k \right)\right) 
\end{cases}\nonumber
\end{equation}
\end{enumerate}
\subsubsection{The second species complete polynomial}
Replacing coefficients $q$, $\alpha$, $\beta$, $\gamma $ and $\delta$ by $q-(\gamma +\delta -2)a-(\gamma -1)(\alpha +\beta -\gamma -\delta +1)$, $\alpha - \gamma -\delta +2$, $\beta - \gamma -\delta +2, 2-\gamma$ and $2 - \delta$ into (\ref{eq:40037a})--(\ref{eq:40038c}). Multiply $x^{1-\gamma } (1-x)^{1-\delta }$ and the new (\ref{eq:40037a})--(\ref{eq:40037c}) together.\footnote{I treat $\gamma$ and $\delta$ as free variables and fixed values of $\alpha $, $\beta $ and $q$.}  
\begin{enumerate} 
\item As $\alpha =\gamma +\delta -2 $, $\beta =\gamma +\delta -1 $ and $q=(\gamma +\delta -2)a +(\gamma -1)(\gamma +\delta -2) $,
 
Its eigenfunction is given by 
\begin{eqnarray}
& &x^{1-\gamma } (1-x)^{1-\delta } y(x)\nonumber\\
&=& x^{1-\gamma } (1-x)^{1-\delta } Hl(a, 0; 0, 1, 2-\gamma, 2 - \delta ; x)\nonumber\\
&=& x^{1-\gamma } (1-x)^{1-\delta } \nonumber
\end{eqnarray}
\item As $\alpha = \gamma +\delta -2N-3 $, $\beta = \gamma +\delta -2N-2 $ and $q=(\gamma +\delta -2)a+(\gamma -1)( \gamma +\delta -4N-4) + \left(2N+1 \right) \left[ -\delta +2N+2 -a\left( -\gamma -\delta+2N+4 \right)\right]$ where $N \in \mathbb{N}_{0}$,

Its eigenfunction is given by
\begin{eqnarray} 
 & &x^{1-\gamma } (1-x)^{1-\delta } y(x)\nonumber\\
&=& x^{1-\gamma } (1-x)^{1-\delta } Hl(a, \left( 2N+1 \right) \left[ -\delta +2N+2 -a\left( -\gamma -\delta +2N+4 \right)\right]; -2N-1, -2N, 2-\gamma, 2 - \delta ; x)\nonumber\\
&=& x^{1-\gamma } (1-x)^{1-\delta } \left\{ \sum_{r=0}^{N} y_{2r}^{N-r}\left( 2N+1;x\right) + \sum_{r=0}^{N} y_{2r+1}^{N-r}\left( 2N+1;x\right) \right\}
\nonumber
\end{eqnarray}
\item As $\alpha = \gamma +\delta -2N-4 $, $\beta = \gamma +\delta -2N-3 $ and $q=(\gamma +\delta -2)a+(\gamma -1)( \gamma +\delta -4N-6) + \left( 2N+2 \right) \left[ -\delta +2N+3 -a\left( -\gamma -\delta +2N+5 \right)\right]$ where $N \in \mathbb{N}_{0}$,

Its eigenfunction is given by
\begin{eqnarray} 
 & &x^{1-\gamma } (1-x)^{1-\delta } y(x)\nonumber\\
&=& x^{1-\gamma } (1-x)^{1-\delta } Hl(a, \left( 2N+2 \right) \left[ -\delta +2N+3 -a\left( -\gamma -\delta +2N+5 \right)\right]; -2N-2, -2N-1 \nonumber\\
&&, 2-\gamma, 2 - \delta ; x)\nonumber\\
&=& x^{1-\gamma } (1-x)^{1-\delta } \left\{ \sum_{r=0}^{N+1} y_{2r}^{N+1-r}\left( 2N+2;x\right) + \sum_{r=0}^{N} y_{2r+1}^{N-r}\left( 2N+2;x\right) \right\}
\nonumber
\end{eqnarray}
In the above,
\begin{eqnarray}
y_0^m(j;x) &=&  \sum_{i_0=0}^{m} \frac{\left( -\frac{j}{2}\right)_{i_0} \left( - \frac{j}{2}+\frac{1}{2}\right)_{i_0}}{\left( 1 \right)_{i_0} \left( -\frac{\gamma }{2}+ \frac{3}{2} \right)_{i_0}} z^{i_0} \nonumber\\
y_1^m(j;x) &=&  \left\{\sum_{i_0=0}^{m} \frac{\left( i_0-\frac{j}{2}\right)\left(i_0 +\Gamma ^{(S)} \right)}{\left( i_0+\frac{1}{2} \right) \left( i_0+1-\frac{\gamma }{2} \right)} \frac{\left( -\frac{j}{2}\right)_{i_0} \left( -\frac{j}{2}+\frac{1}{2}\right)_{i_0}}{\left( 1 \right)_{i_0} \left( -\frac{\gamma }{2}+ \frac{3}{2} \right)_{i_0}} \right. \nonumber\\
&&\times \left. \sum_{i_1 = i_0}^{m} \frac{\left( -\frac{j}{2}+\frac{1}{2} \right)_{i_1} \left( -\frac{j}{2}+1\right)_{i_1}\left( \frac{3}{2} \right)_{i_0} \left( -\frac{\gamma }{2}+2 \right)_{i_0}}{\left(-\frac{j}{2}+\frac{1}{2} \right)_{i_0} \left( -\frac{j}{2}+1\right)_{i_0}\left( \frac{3}{2} \right)_{i_1} \left( -\frac{\gamma }{2} +2 \right)_{i_1}} z^{i_1}\right\} \eta \nonumber
\end{eqnarray}
\begin{eqnarray} 
y_{\tau }^m(j;x) &=&  \left\{ \sum_{i_0=0}^{m} \frac{\left( i_0-\frac{j}{2}\right)\left(i_0 +\Gamma ^{(S)} \right)}{\left( i_0+\frac{1}{2} \right) \left( i_0+1-\frac{\gamma }{2} \right)} \frac{\left( -\frac{j}{2}\right)_{i_0} \left( -\frac{j}{2}+\frac{1}{2}\right)_{i_0}}{\left( 1 \right)_{i_0} \left( -\frac{\gamma }{2}+ \frac{3}{2} \right)_{i_0}} \right.\nonumber\\
&&\times \prod_{k=1}^{\tau -1} \left( \sum_{i_k = i_{k-1}}^{m} \frac{\left( i_k+ \frac{k}{2}-\frac{j}{2}\right)\left(i_k + \frac{k}{2}+\Gamma ^{(S)} \right)}{\left( i_k+\frac{k}{2}+\frac{1}{2} \right) \left( i_k+\frac{k}{2}+1-\frac{\gamma }{2} \right)} \right. \nonumber\\
&&\times \left. \frac{\left( -\frac{j}{2}+\frac{k}{2}\right)_{i_k} \left( -\frac{j}{2}+ \frac{k}{2}+ \frac{1}{2}\right)_{i_k}\left(  \frac{k}{2} +1\right)_{i_{k-1}} \left( -\frac{\gamma }{2}+ \frac{3}{2}+ \frac{k}{2} \right)_{i_{k-1}}}{\left( -\frac{j}{2}+\frac{k}{2}\right)_{i_{k-1}} \left( -\frac{j}{2}+ \frac{k}{2}+ \frac{1}{2}\right)_{i_{k-1}}\left( \frac{k}{2}+1 \right)_{i_k} \left( -\frac{\gamma }{2}+ \frac{3}{2}+ \frac{k}{2} \right)_{i_k}} \right) \nonumber\\
&&\times \left. \sum_{i_{\tau } = i_{\tau -1}}^{m} \frac{\left( -\frac{j}{2}+\frac{\tau }{2}\right)_{i_{\tau }} \left( -\frac{j}{2}+ \frac{\tau }{2}+ \frac{1}{2}\right)_{i_{\tau }}\left(  \frac{\tau }{2} +1\right)_{i_{\tau -1}} \left( -\frac{\gamma }{2}+ \frac{3}{2}+\frac{\tau }{2} \right)_{i_{\tau -1}}}{\left( -\frac{j}{2}+\frac{\tau }{2}\right)_{i_{\tau -1}} \left(-\frac{j}{2}+ \frac{\tau }{2}+ \frac{1}{2}\right)_{i_{\tau -1}}\left(  \frac{\tau }{2}+1 \right)_{i_{\tau }} \left( -\frac{\gamma }{2}+ \frac{3}{2}+\frac{\tau }{2} \right)_{i_{\tau }}} z^{i_{\tau }}\right\} \eta ^{\tau } \nonumber 
\end{eqnarray}
where
\begin{equation}
\begin{cases} \tau \geq 2 \cr
z = -\frac{1}{a}x^2 \cr
\eta = \frac{(1+a)}{a} x \cr
\Gamma ^{(S)} =  \frac{1}{2(1+a)}\left( \delta -j-1 +a\left( -\gamma -\delta +3+j  \right)\right)
\end{cases}\nonumber
\end{equation}
\end{enumerate}
 \subsection{${\displaystyle  Hl(1-a,-q+\alpha \beta; \alpha,\beta, \delta, \gamma; 1-x)}$} 
\subsubsection{The first species complete polynomial}
Replacing coefficients $a$, $q$, $\gamma $, $\delta$ and $x$ by $1-a$, $-q +\alpha \beta $, $\delta $, $\gamma $ and $1-x$  into (\ref{eq:40018})--(\ref{eq:40022c}).\footnote{I treat $\beta $, $\gamma$ and $\delta$ as free variables and fixed values of $\alpha $ and $q$.}
\begin{enumerate} 
\item As $\alpha =0$ and $ q = -q_0^0$ where $q_0^0=0$, 

The eigenfunction is given by
\begin{equation}
y(\xi ) = Hl(1-a,0; 0,\beta, \delta, \gamma; 1-x) =1 \nonumber
\end{equation}
\item As $\alpha =-2N-1$ where $N \in \mathbb{N}_{0}$,

An algebraic equation of degree $2N+2$ for the determination of $q$ is given by
\begin{equation}
0 = \sum_{r=0}^{N+1}\bar{c}\left( 2r, N+1-r; 2N+1, -q-(2N+1)\beta \right)\nonumber
\end{equation}
The eigenvalue of $q $ is written by $-(2N+1) \beta - q_{2N+1}^m$ where $m = 0,1,2,\cdots,2N+1 $; $q_{2N+1}^0 < q_{2N+1}^1 < \cdots < q_{2N+1}^{2N+1}$. Its eigenfunction is given by
\begin{eqnarray} 
y(\xi ) &=& Hl(1-a, q_{2N+1}^m; -2N-1,\beta, \delta, \gamma; 1-x)\nonumber\\
&=& \sum_{r=0}^{N} y_{2r}^{N-r}\left( 2N+1,q_{2N+1}^m;\xi \right)+ \sum_{r=0}^{N} y_{2r+1}^{N-r}\left( 2N+1,q_{2N+1}^m;\xi \right)  
\nonumber
\end{eqnarray}
\item As $\alpha =-2N-2$ where $N \in \mathbb{N}_{0}$,

An algebraic equation of degree $2N+3$ for the determination of $q$ is given by
\begin{eqnarray}
0  = \sum_{r=0}^{N+1}\bar{c}\left( 2r+1, N+1-r; 2N+2, -q-(2N+2)\beta\right)\nonumber
\end{eqnarray}
The eigenvalue of $q$ is written by $-(2N+2) \beta - q_{2N+2}^m$ where $m = 0,1,2,\cdots,2N+2 $; $q_{2N+2}^0 < q_{2N+2}^1 < \cdots < q_{2N+2}^{2N+2}$. Its eigenfunction is given by

\begin{eqnarray} 
y(\xi ) &=& Hl(1-a, q_{2N+2}^m; -2N-2,\beta, \delta, \gamma; 1-x)\nonumber\\
&=& \sum_{r=0}^{N+1} y_{2r}^{N+1-r}\left( 2N+2,q_{2N+2}^m;\xi \right) + \sum_{r=0}^{N} y_{2r+1}^{N-r}\left( 2N+2,q_{2N+2}^m;\xi \right) 
\nonumber
\end{eqnarray}
In the above,
\begin{eqnarray}
\bar{c}(0,n;j,\tilde{q})  &=& \frac{\left( -\frac{j}{2}\right)_{n} \left( \frac{\beta }{2} \right)_{n}}{\left( 1 \right)_{n} \left( \frac{\delta }{2}+ \frac{1}{2} \right)_{n}} \left( \frac{-1}{1-a} \right)^{n}\nonumber\\
\bar{c}(1,n;j,\tilde{q}) &=& \left( \frac{2-a}{1-a}\right) \sum_{i_0=0}^{n}\frac{ i_0 \left(i_0 +\Gamma _0^{(F)} \right)+\frac{\tilde{q}}{4(2-a)}}{\left( i_0+\frac{1}{2} \right) \left( i_0+\frac{\delta }{2} \right)} \frac{\left( -\frac{j}{2}\right)_{i_0} \left( \frac{\beta }{2} \right)_{i_0}}{\left( 1 \right)_{i_0} \left( \frac{\delta }{2}+ \frac{1}{2} \right)_{i_0}}\nonumber\\
&&\times   \frac{\left( -\frac{j}{2}+\frac{1}{2} \right)_{n} \left( \frac{\beta }{2}+ \frac{1}{2} \right)_{n}\left( \frac{3}{2} \right)_{i_0} \left( \frac{\delta }{2}+1 \right)_{i_0}}{\left( -\frac{j}{2}+\frac{1}{2}\right)_{i_0} \left( \frac{\beta }{2}+ \frac{1}{2} \right)_{i_0}\left( \frac{3}{2} \right)_{n} \left( \frac{\delta }{2}+1 \right)_{n}} \left( \frac{-1}{1-a} \right)^{n} \nonumber
\end{eqnarray}
\begin{eqnarray}
\bar{c}(\tau ,n;j,\tilde{q}) &=& \left( \frac{2-a}{1-a}\right)^{\tau} \sum_{i_0=0}^{n}\frac{ i_0 \left(i_0 +\Gamma _0^{(F)} \right)+\frac{\tilde{q}}{4(2-a)}}{\left( i_0+\frac{1}{2} \right) \left( i_0+\frac{\delta }{2} \right)} \frac{\left( -\frac{j}{2}\right)_{i_0} \left( \frac{\beta }{2} \right)_{i_0}}{\left( 1 \right)_{i_0} \left( \frac{\delta }{2}+ \frac{1}{2} \right)_{i_0}}  \nonumber\\
&&\times \prod_{k=1}^{\tau -1} \left( \sum_{i_k = i_{k-1}}^{n} \frac{\left( i_k+ \frac{k}{2} \right)\left(i_k +\Gamma _k^{(F)} \right)+\frac{\tilde{q}}{4(2-a)}}{\left( i_k+\frac{k}{2}+\frac{1}{2} \right) \left( i_k+\frac{k}{2}+\frac{\delta }{2} \right)} \right. \nonumber\\
&&\times \left. \frac{\left( -\frac{j}{2}+\frac{k}{2}\right)_{i_k} \left( \frac{\beta }{2}+ \frac{k}{2} \right)_{i_k}\left( \frac{k}{2}+1 \right)_{i_{k-1}} \left( \frac{\delta }{2}+ \frac{1}{2}+ \frac{k}{2} \right)_{i_{k-1}}}{\left( -\frac{j}{2}+\frac{k}{2}\right)_{i_{k-1}} \left( \frac{\beta }{2}+ \frac{k}{2} \right)_{i_{k-1}}\left( \frac{k}{2} +1\right)_{i_k} \left( \frac{\delta }{2}+ \frac{1}{2}+ \frac{k}{2} \right)_{i_k}} \right) \nonumber\\
&&\times \frac{\left( -\frac{j}{2}+\frac{\tau }{2}\right)_{n} \left( \frac{\beta }{2}+ \frac{\tau }{2} \right)_{n}\left( \frac{\tau }{2} +1\right)_{i_{\tau -1}} \left( \frac{\delta }{2}+ \frac{1}{2}+\frac{\tau }{2} \right)_{i_{\tau -1}}}{\left( -\frac{j}{2}+\frac{\tau }{2}\right)_{i_{\tau -1}} \left( \frac{\beta }{2}+\frac{\tau }{2} \right)_{i_{\tau -1}}\left( \frac{\tau }{2}+1 \right)_{n} \left( \frac{\delta }{2}+ \frac{1}{2}+\frac{\tau }{2} \right)_{n}} \left(\frac{-1}{1-a} \right)^{n } \nonumber
\end{eqnarray}
\begin{eqnarray}
y_0^m(j,\tilde{q};\xi) &=& \sum_{i_0=0}^{m} \frac{\left( -\frac{j}{2}\right)_{i_0} \left( \frac{\beta }{2} \right)_{i_0}}{\left( 1 \right)_{i_0} \left( \frac{\delta }{2}+ \frac{1}{2} \right)_{i_0}} z^{i_0} \nonumber\\
y_1^m(j,\tilde{q};\xi) &=& \left\{\sum_{i_0=0}^{m} \frac{ i_0 \left(i_0 +\Gamma _0^{(F)} \right)+\frac{\tilde{q}}{4(2-a)}}{\left( i_0+\frac{1}{2} \right) \left( i_0+\frac{\delta }{2} \right)} \frac{\left( -\frac{j}{2}\right)_{i_0} \left( \frac{\beta }{2} \right)_{i_0}}{\left( 1 \right)_{i_0} \left( \frac{\delta }{2}+ \frac{1}{2} \right)_{i_0}} \right. \nonumber\\
&&\times \left. \sum_{i_1 = i_0}^{m} \frac{\left( -\frac{j}{2}+\frac{1}{2} \right)_{i_1} \left( \frac{\beta }{2}+ \frac{1}{2} \right)_{i_1}\left( \frac{3}{2} \right)_{i_0} \left( \frac{\delta }{2}+1 \right)_{i_0}}{\left(-\frac{j}{2}+\frac{1}{2} \right)_{i_0} \left( \frac{\beta }{2}+ \frac{1}{2} \right)_{i_0}\left( \frac{3}{2} \right)_{i_1} \left( \frac{\delta }{2} +1 \right)_{i_1}} z^{i_1}\right\} \eta 
\nonumber\\
y_{\tau }^m(j,\tilde{q};\xi) &=& \left\{ \sum_{i_0=0}^{m} \frac{ i_0 \left(i_0 +\Gamma _0^{(F)} \right)+\frac{\tilde{q}}{4(2-a)}}{\left( i_0+\frac{1}{2} \right) \left( i_0+\frac{\delta }{2} \right)} \frac{\left( -\frac{j}{2}\right)_{i_0} \left( \frac{\beta }{2} \right)_{i_0}}{\left( 1 \right)_{i_0} \left( \frac{\delta }{2}+ \frac{1}{2} \right)_{i_0}} \right.\nonumber\\
&&\times \prod_{k=1}^{\tau -1} \left( \sum_{i_k = i_{k-1}}^{m} \frac{\left( i_k+ \frac{k}{2} \right)\left(i_k +\Gamma _k^{(F)} \right)+\frac{\tilde{q}}{4(2-a)}}{\left( i_k+\frac{k}{2}+\frac{1}{2} \right) \left( i_k+\frac{k}{2}+\frac{\delta }{2} \right)} \right. \nonumber\\
&&\times \left. \frac{\left( -\frac{j}{2}+\frac{k}{2}\right)_{i_k} \left( \frac{\beta }{2}+ \frac{k}{2} \right)_{i_k}\left( \frac{k}{2}+1 \right)_{i_{k-1}} \left( \frac{\delta }{2}+ \frac{1}{2}+ \frac{k}{2} \right)_{i_{k-1}}}{\left( -\frac{j}{2}+\frac{k}{2}\right)_{i_{k-1}} \left( \frac{\beta }{2}+ \frac{k}{2} \right)_{i_{k-1}}\left( \frac{k}{2} +1\right)_{i_k} \left( \frac{\delta }{2}+ \frac{1}{2}+ \frac{k}{2} \right)_{i_k}} \right) \nonumber\\
&&\times \left. \sum_{i_{\tau } = i_{\tau -1}}^{m} \frac{\left( -\frac{j}{2}+\frac{\tau }{2}\right)_{i_{\tau }} \left( \frac{\beta }{2}+ \frac{\tau }{2} \right)_{i_{\tau }}\left( \frac{\tau }{2}+1 \right)_{i_{\tau -1}} \left( \frac{\delta }{2}+ \frac{1}{2}+\frac{\tau }{2} \right)_{i_{\tau -1}}}{\left( -\frac{j}{2}+\frac{\tau }{2}\right)_{i_{\tau -1}} \left( \frac{\beta }{2}+ \frac{\tau }{2} \right)_{i_{\tau -1}}\left( \frac{\tau }{2}+1 \right)_{i_{\tau }} \left( \frac{\delta }{2}+ \frac{1}{2}+\frac{\tau }{2} \right)_{i_{\tau }}} z^{i_{\tau }}\right\} \eta ^{\tau } \nonumber
\end{eqnarray}
where
\begin{equation}
\begin{cases} \tau \geq 2 \cr
\xi =1-x \cr
z = \frac{-1}{1-a}\xi^2 \cr
\eta = \frac{2-a}{1-a}\xi \cr
\tilde{q} = -q+\alpha \beta  \cr
\Gamma _0^{(F)} =  \frac{1}{2(2-a)}\left( \beta -\gamma -j +(1-a)\left( \gamma +\delta -1 \right)\right) \cr
\Gamma _k^{(F)} = \frac{1}{2(2-a)}\left( \beta -\gamma -j+k +(1-a)\left( \gamma +\delta -1+k \right)\right) 
\end{cases}\nonumber 
\end{equation}
\end{enumerate}
\subsubsection{The second species complete polynomial}
Replacing coefficients $a$, $q$, $\gamma $, $\delta$ and $x$ by $1-a$, $-q +\alpha \beta $, $\delta $, $\gamma $ and $1-x$  into (\ref{eq:40037a})--(\ref{eq:40038c}).\footnote{I treat $\gamma$ and $\delta$ as free variables and fixed values of $\alpha $, $\beta $ and $q$.}
\begin{enumerate} 
\item As $\alpha =0 $, $\beta =1 $ and $ q=0$,
 
Its eigenfunction is given by
\begin{equation}
y(\xi ) = Hl(1-a,0; 0, 1, \delta, \gamma; 1-x)  =1 \nonumber
\end{equation}
\item As $\alpha =-2N-1 $, $\beta = -2N $ and $q= - \left(2N+1 \right) \left[ \gamma -(1-a)\left( \gamma +\delta+2N \right)\right]$ where $N \in \mathbb{N}_{0}$,

Its eigenfunction is given by
\begin{eqnarray} 
y(\xi ) &=& Hl(1-a,\left( 2N+1 \right) \left[ \gamma +2N -(1-a)\left( \gamma +\delta +2N \right)\right]; -2N-1,-2N, \delta, \gamma; 1-x)\nonumber\\
&=& \sum_{r=0}^{N} y_{2r}^{N-r}\left( 2N+1;\xi \right) + \sum_{r=0}^{N} y_{2r+1}^{N-r}\left( 2N+1;\xi \right) 
\nonumber
\end{eqnarray}
\item As $\alpha =-2N-2 $, $\beta = -2N-1 $ and $q=  - \left(2N+2 \right) \left[ \gamma  -(1-a)\left( \gamma +\delta +2N+1 \right)\right]$ where $N \in \mathbb{N}_{0}$,

Its eigenfunction is given by
\begin{eqnarray} 
y(\xi ) &=& Hl(1-a,\left( 2N+2 \right) \left[ \gamma +2N+1 -(1-a)\left( \gamma +\delta +2N+1 \right)\right]; -2N-2, -2N-1, \delta, \gamma; 1-x)\nonumber\\
&=& \sum_{r=0}^{N+1} y_{2r}^{N+1-r}\left( 2N+2;\xi \right) + \sum_{r=0}^{N} y_{2r+1}^{N-r}\left( 2N+2;\xi \right) 
\nonumber
\end{eqnarray}
In the above,
\begin{eqnarray}
y_0^m(j;\xi) &=&  \sum_{i_0=0}^{m} \frac{\left( -\frac{j}{2}\right)_{i_0} \left( - \frac{j}{2}+\frac{1}{2}\right)_{i_0}}{\left( 1 \right)_{i_0} \left( \frac{\delta }{2}+ \frac{1}{2} \right)_{i_0}} z^{i_0} \nonumber\\
y_1^m(j;\xi) &=&  \left\{\sum_{i_0=0}^{m} \frac{\left( i_0-\frac{j}{2}\right)\left(i_0 +\Gamma ^{(S)} \right)}{\left( i_0+\frac{1}{2} \right) \left( i_0+\frac{\delta }{2} \right)} \frac{\left( -\frac{j}{2}\right)_{i_0} \left( -\frac{j}{2}+\frac{1}{2}\right)_{i_0}}{\left( 1 \right)_{i_0} \left( \frac{\delta }{2}+ \frac{1}{2} \right)_{i_0}} \right. \nonumber\\
&&\times \left. \sum_{i_1 = i_0}^{m} \frac{\left( -\frac{j}{2}+\frac{1}{2} \right)_{i_1} \left( -\frac{j}{2}+1\right)_{i_1}\left( \frac{3}{2} \right)_{i_0} \left( \frac{\delta }{2}+1 \right)_{i_0}}{\left(-\frac{j}{2}+\frac{1}{2} \right)_{i_0} \left( -\frac{j}{2}+1\right)_{i_0}\left( \frac{3}{2} \right)_{i_1} \left( \frac{\delta }{2} +1 \right)_{i_1}} z^{i_1}\right\} \eta \nonumber\\
y_{\tau }^m(j;\xi) &=&  \left\{ \sum_{i_0=0}^{m} \frac{\left( i_0-\frac{j}{2}\right)\left(i_0 +\Gamma ^{(S)} \right)}{\left( i_0+\frac{1}{2} \right) \left( i_0+\frac{\delta }{2} \right)} \frac{\left( -\frac{j}{2}\right)_{i_0} \left( -\frac{j}{2}+\frac{1}{2}\right)_{i_0}}{\left( 1 \right)_{i_0} \left( \frac{\delta }{2}+ \frac{1}{2} \right)_{i_0}} \right.\nonumber\\
&&\times \prod_{k=1}^{\tau -1} \left( \sum_{i_k = i_{k-1}}^{m} \frac{\left( i_k+ \frac{k}{2}-\frac{j}{2}\right)\left(i_k + \frac{k}{2}+\Gamma ^{(S)} \right)}{\left( i_k+\frac{k}{2}+\frac{1}{2} \right) \left( i_k+\frac{k}{2}+\frac{\delta }{2} \right)} \right. \nonumber\\
&&\times \left. \frac{\left( -\frac{j}{2}+\frac{k}{2}\right)_{i_k} \left( -\frac{j}{2}+ \frac{k}{2}+ \frac{1}{2}\right)_{i_k}\left(  \frac{k}{2} +1\right)_{i_{k-1}} \left( \frac{\delta }{2}+ \frac{1}{2}+ \frac{k}{2} \right)_{i_{k-1}}}{\left( -\frac{j}{2}+\frac{k}{2}\right)_{i_{k-1}} \left( -\frac{j}{2}+ \frac{k}{2}+ \frac{1}{2}\right)_{i_{k-1}}\left( \frac{k}{2}+1 \right)_{i_k} \left( \frac{\delta }{2}+ \frac{1}{2}+ \frac{k}{2} \right)_{i_k}} \right) \nonumber\\
&&\times \left. \sum_{i_{\tau } = i_{\tau -1}}^{m} \frac{\left( -\frac{j}{2}+\frac{\tau }{2}\right)_{i_{\tau }} \left( -\frac{j}{2}+ \frac{\tau }{2}+ \frac{1}{2}\right)_{i_{\tau }}\left(  \frac{\tau }{2} +1\right)_{i_{\tau -1}} \left( \frac{\delta }{2}+ \frac{1}{2}+\frac{\tau }{2} \right)_{i_{\tau -1}}}{\left( -\frac{j}{2}+\frac{\tau }{2}\right)_{i_{\tau -1}} \left(-\frac{j}{2}+ \frac{\tau }{2}+ \frac{1}{2}\right)_{i_{\tau -1}}\left(  \frac{\tau }{2}+1 \right)_{i_{\tau }} \left( \frac{\delta }{2}+ \frac{1}{2}+\frac{\tau }{2} \right)_{i_{\tau }}} z^{i_{\tau }}\right\} \eta ^{\tau } \nonumber 
\end{eqnarray}
where
\begin{equation}
\begin{cases} \tau \geq 2 \cr
\xi =1-x \cr
z = \frac{-1}{1-a}\xi^2 \cr
\eta = \frac{2-a}{1-a}\xi \cr
\Gamma ^{(S)} =\frac{1}{2(2-a)}\left( -\gamma -j+1 +(1-a)\left( \gamma +\delta -1+j \right)\right) 
\end{cases}\nonumber 
\end{equation}
\end{enumerate}
\subsection{\footnotesize ${\displaystyle (1-x)^{1-\delta } Hl(1-a,-q+(\delta -1)\gamma a+(\alpha -\delta +1)(\beta -\delta +1); \alpha-\delta +1,\beta-\delta +1, 2-\delta, \gamma; 1-x)}$ \normalsize}
\subsubsection{The first species complete polynomial}
Replacing coefficients $a$, $q$, $\alpha $, $\beta $, $\gamma $, $\delta$ and $x$ by $1-a$, $-q+(\delta -1)\gamma a+(\alpha -\delta +1)(\beta -\delta +1)$, $\alpha-\delta +1 $, $\beta-\delta +1 $, $2-\delta$, $\gamma $ and $1-x$ into (\ref{eq:40018})--(\ref{eq:40022c}). Multiply $(1-x)^{1-\delta }$ and the new (\ref{eq:40018}), (\ref{eq:40019b}) and (\ref{eq:40020b})  together.\footnote{I treat $\beta $, $\gamma$ and $\delta$ as free variables and fixed values of $\alpha $ and $q$.}
\begin{enumerate} 
\item As $\alpha = \delta -1$ and $ q= (\delta -1)\gamma a -q_0^0$ where $q_0^0=0$,

The eigenfunction is given by
\begin{eqnarray}
& &(1-x)^{1-\delta } y(\xi)\nonumber\\
&=& (1-x)^{1-\delta } Hl(1-a,0; 0,\beta-\delta +1, 2-\delta, \gamma; 1-x) \nonumber\\
&=& (1-x)^{1-\delta } \nonumber
\end{eqnarray}
\item As $\alpha =\delta -2N-2$ where $N \in \mathbb{N}_{0}$,

An algebraic equation of degree $2N+2$ for the determination of $q$ is given by
\begin{equation}
0 = \sum_{r=0}^{N+1}\bar{c}\left( 2r, N+1-r; 2N+1, -q+(\delta -1)\gamma a-(2N+1)(\beta -\delta +1)\right)\nonumber
\end{equation}
The eigenvalue of $q$ is written by $(\delta -1)\gamma a-(2N+1)(\beta -\delta +1)- q_{2N+1}^m$ where $m = 0,1,2,\cdots,2N+1 $; $q_{2N+1}^0 < q_{2N+1}^1 < \cdots < q_{2N+1}^{2N+1}$. Its eigenfunction is given by
\begin{eqnarray} 
& &(1-x)^{1-\delta } y(\xi)\nonumber\\
&=& (1-x)^{1-\delta } Hl(1-a,q_{2N+1}^m; -2N-1,\beta-\delta +1, 2-\delta, \gamma; 1-x) \nonumber\\
&=& (1-x)^{1-\delta } \left\{ \sum_{r=0}^{N} y_{2r}^{N-r}\left( 2N+1,q_{2N+1}^m;\xi \right)+ \sum_{r=0}^{N} y_{2r+1}^{N-r}\left( 2N+1,q_{2N+1}^m;\xi \right) \right\} 
\nonumber
\end{eqnarray}
\item As $\alpha =\delta-2N-3$ where $N \in \mathbb{N}_{0}$,

An algebraic equation of degree $2N+3$ for the determination of $q$ is given by
\begin{eqnarray}
0  = \sum_{r=0}^{N+1}\bar{c}\left( 2r+1, N+1-r; 2N+2, -q+(\delta -1)\gamma a-(2N+2)(\beta -\delta +1)\right)\nonumber
\end{eqnarray}
The eigenvalue of $q$ is written by $(\delta -1)\gamma a-(2N+2)(\beta -\delta +1)- q_{2N+2}^m$ where $m = 0,1,2,\cdots,2N+2 $; $q_{2N+2}^0 < q_{2N+2}^1 < \cdots < q_{2N+2}^{2N+2}$. Its eigenfunction is given by
\begin{eqnarray} 
& &(1-x)^{1-\delta } y(\xi)\nonumber\\
&=& (1-x)^{1-\delta } Hl(1-a,q_{2N+1}^m; -2N-1,\beta-\delta +1, 2-\delta, \gamma; 1-x) \nonumber\\
&=& (1-x)^{1-\delta } \left\{ \sum_{r=0}^{N+1} y_{2r}^{N+1-r}\left( 2N+2,q_{2N+2}^m;\xi \right) + \sum_{r=0}^{N} y_{2r+1}^{N-r}\left( 2N+2,q_{2N+2}^m;\xi \right) \right\}
\nonumber
\end{eqnarray}
In the above,
\begin{eqnarray}
\bar{c}(0,n;j,\tilde{q})  &=& \frac{\left( -\frac{j}{2}\right)_{n} \left( \frac{\beta }{2}-\frac{\delta }{2}+\frac{1}{2} \right)_{n}}{\left( 1 \right)_{n} \left( -\frac{\delta }{2}+ \frac{3}{2} \right)_{n}} \left( \frac{-1}{1-a} \right)^{n}\nonumber \\
\bar{c}(1,n;j,\tilde{q}) &=& \left( \frac{2-a}{1-a}\right) \sum_{i_0=0}^{n}\frac{ i_0 \left(i_0 +\Gamma _0^{(F)} \right)+\frac{\tilde{q}}{4(2-a)}}{\left( i_0+\frac{1}{2} \right) \left( i_0+1-\frac{\delta }{2} \right)} \frac{\left( -\frac{j}{2}\right)_{i_0} \left( \frac{\beta }{2}-\frac{\delta }{2}+\frac{1}{2}\right)_{i_0}}{\left( 1 \right)_{i_0} \left( -\frac{\delta }{2}+ \frac{3}{2} \right)_{i_0}}\nonumber\\
&&\times  \frac{\left( -\frac{j}{2}+\frac{1}{2} \right)_{n} \left( \frac{\beta }{2}-\frac{\delta }{2}+1 \right)_{n}\left( \frac{3}{2} \right)_{i_0} \left( -\frac{\delta }{2}+2 \right)_{i_0}}{\left( -\frac{j}{2}+\frac{1}{2}\right)_{i_0} \left( \frac{\beta }{2}-\frac{\delta }{2}+1 \right)_{i_0}\left( \frac{3}{2} \right)_{n} \left( -\frac{\delta }{2}+2 \right)_{n}} \left( \frac{-1}{1-a} \right)^{n}  
\nonumber \\
\bar{c}(\tau ,n;j,\tilde{q}) &=& \left( \frac{2-a}{1-a}\right)^{\tau} \sum_{i_0=0}^{n}\frac{ i_0 \left(i_0 +\Gamma _0^{(F)}  \right)+\frac{\tilde{q}}{4(2-a)}}{\left( i_0+\frac{1}{2} \right) \left( i_0+1-\frac{\delta }{2} \right)} \frac{\left( -\frac{j}{2}\right)_{i_0} \left( \frac{\beta }{2}-\frac{\delta }{2}+\frac{1}{2} \right)_{i_0}}{\left( 1 \right)_{i_0} \left( -\frac{\delta }{2}+ \frac{3}{2} \right)_{i_0}} \nonumber\\
&&\times \prod_{k=1}^{\tau -1} \left( \sum_{i_k = i_{k-1}}^{n} \frac{\left( i_k+ \frac{k}{2} \right)\left( i_k +\Gamma _k^{(F)} \right) +\frac{\tilde{q}}{4(2-a)}}{\left( i_k+\frac{k}{2}+\frac{1}{2} \right) \left( i_k+\frac{k}{2}+1-\frac{\delta }{2} \right)} \right. \nonumber\\
&&\times \left. \frac{\left( -\frac{j}{2}+\frac{k}{2}\right)_{i_k} \left( \frac{\beta }{2}-\frac{\delta }{2}+\frac{1}{2}+ \frac{k}{2} \right)_{i_k}\left( \frac{k}{2}+1 \right)_{i_{k-1}} \left( -\frac{\delta }{2}+ \frac{3}{2}+ \frac{k}{2} \right)_{i_{k-1}}}{\left( -\frac{j}{2}+\frac{k}{2}\right)_{i_{k-1}} \left( \frac{\beta }{2}-\frac{\delta }{2}+\frac{1}{2}+ \frac{k}{2} \right)_{i_{k-1}}\left( \frac{k}{2} +1\right)_{i_k} \left( -\frac{\delta }{2}+ \frac{3}{2}+ \frac{k}{2} \right)_{i_k}} \right) \nonumber\\
&&\times \frac{\left( -\frac{j}{2}+\frac{\tau }{2}\right)_{n} \left( \frac{\beta }{2}-\frac{\delta }{2}+\frac{1}{2}+ \frac{\tau }{2} \right)_{n}\left( \frac{\tau }{2} +1\right)_{i_{\tau -1}} \left( -\frac{\delta }{2}+ \frac{3}{2}+\frac{\tau }{2} \right)_{i_{\tau -1}}}{\left( -\frac{j}{2}+\frac{\tau }{2}\right)_{i_{\tau -1}} \left( \frac{\beta }{2}-\frac{\delta }{2}+\frac{1}{2}+\frac{\tau }{2} \right)_{i_{\tau -1}}\left( \frac{\tau }{2}+1 \right)_{n} \left( -\frac{\delta }{2}+ \frac{3}{2}+\frac{\tau }{2} \right)_{n}} \left(\frac{-1}{1-a} \right)^{n } \nonumber 
\end{eqnarray}
\begin{eqnarray}
y_0^m(j,\tilde{q};\xi) &=& \sum_{i_0=0}^{m} \frac{\left( -\frac{j}{2}\right)_{i_0} \left( \frac{\beta }{2}-\frac{\delta }{2}+\frac{1}{2} \right)_{i_0}}{\left( 1 \right)_{i_0} \left( -\frac{\delta }{2}+ \frac{3}{2} \right)_{i_0}} z^{i_0} \nonumber \\
y_1^m(j,\tilde{q};\xi) &=& \left\{\sum_{i_0=0}^{m} \frac{ i_0 \left(i_0 +\Gamma _0^{(F)} \right)+\frac{\tilde{q}}{4(2-a)}}{\left( i_0+\frac{1}{2} \right) \left( i_0+1-\frac{\delta }{2} \right)} \frac{\left( -\frac{j}{2}\right)_{i_0} \left( \frac{\beta }{2}-\frac{\delta }{2}+\frac{1}{2} \right)_{i_0}}{\left( 1 \right)_{i_0} \left( -\frac{\delta }{2}+ \frac{3}{2} \right)_{i_0}} \right. \nonumber\\
&&\times \left. \sum_{i_1 = i_0}^{m} \frac{\left( -\frac{j}{2}+\frac{1}{2} \right)_{i_1} \left( \frac{\beta }{2}-\frac{\delta }{2}+1 \right)_{i_1}\left( \frac{3}{2} \right)_{i_0} \left( -\frac{\delta }{2}+2 \right)_{i_0}}{\left(-\frac{j}{2}+\frac{1}{2} \right)_{i_0} \left( \frac{\beta }{2}-\frac{\delta }{2}+ 1 \right)_{i_0}\left( \frac{3}{2} \right)_{i_1} \left( -\frac{\delta }{2} +2 \right)_{i_1}} z^{i_1}\right\} \eta 
\nonumber\\
y_{\tau }^m(j,\tilde{q};\xi) &=& \left\{ \sum_{i_0=0}^{m} \frac{ i_0 \left(i_0 +\Gamma _0^{(F)} \right)+\frac{\tilde{q}}{4(2-a)}}{\left( i_0+\frac{1}{2} \right) \left( i_0+1-\frac{\delta }{2} \right)} \frac{\left( -\frac{j}{2}\right)_{i_0} \left( \frac{\beta }{2}-\frac{\delta }{2}+\frac{1}{2} \right)_{i_0}}{\left( 1 \right)_{i_0} \left( -\frac{\delta }{2}+ \frac{3}{2} \right)_{i_0}} \right.\nonumber\\
&&\times \prod_{k=1}^{\tau -1} \left( \sum_{i_k = i_{k-1}}^{m} \frac{\left( i_k+ \frac{k}{2} \right)\left( i_k +\Gamma _k^{(F)} \right) +\frac{\tilde{q}}{4(2-a)}}{\left( i_k+\frac{k}{2}+\frac{1}{2} \right) \left( i_k+\frac{k}{2}+1-\frac{\delta }{2} \right)} \right. \nonumber\\
&&\times \left. \frac{\left( -\frac{j}{2}+\frac{k}{2}\right)_{i_k} \left( \frac{\beta }{2}-\frac{\delta }{2}+\frac{1}{2}+ \frac{k}{2} \right)_{i_k}\left( \frac{k}{2}+1 \right)_{i_{k-1}} \left( -\frac{\delta }{2}+ \frac{3}{2}+ \frac{k}{2} \right)_{i_{k-1}}}{\left( -\frac{j}{2}+\frac{k}{2}\right)_{i_{k-1}} \left( \frac{\beta }{2}-\frac{\delta }{2}+\frac{1}{2}+ \frac{k}{2} \right)_{i_{k-1}}\left( \frac{k}{2} +1\right)_{i_k} \left( -\frac{\delta }{2}+ \frac{3}{2}+ \frac{k}{2} \right)_{i_k}} \right) \nonumber\\
&&\times \left. \sum_{i_{\tau } = i_{\tau -1}}^{m} \frac{\left( -\frac{j}{2}+\frac{\tau }{2}\right)_{i_{\tau }} \left( \frac{\beta }{2}-\frac{\delta }{2}+\frac{1}{2} + \frac{\tau }{2} \right)_{i_{\tau }}\left( \frac{\tau }{2}+1 \right)_{i_{\tau -1}} \left( -\frac{\delta }{2}+ \frac{3}{2}+\frac{\tau }{2} \right)_{i_{\tau -1}}}{\left( -\frac{j}{2}+\frac{\tau }{2}\right)_{i_{\tau -1}} \left( \frac{\beta }{2}-\frac{\delta }{2}+\frac{1}{2}+ \frac{\tau }{2} \right)_{i_{\tau -1}}\left( \frac{\tau }{2}+1 \right)_{i_{\tau }} \left( -\frac{\delta }{2}+ \frac{3}{2}+\frac{\tau }{2} \right)_{i_{\tau }}} z^{i_{\tau }}\right\} \eta ^{\tau }  \nonumber 
\end{eqnarray}
where
\begin{equation}
\begin{cases} \tau \geq 2 \cr
\xi =1-x \cr
z = \frac{-1}{1-a}\xi^2 \cr
\eta = \frac{2-a}{1-a}\xi \cr
\tilde{q} = -q+(\delta -1)\gamma a+(\alpha -\delta +1)(\beta -\delta +1) \cr
\Gamma _0^{(F)} = \frac{1}{2(2-a)}\left( \beta -\gamma -\delta +1-j +(1-a)\left( \gamma -\delta + 1 \right)\right)  \cr
\Gamma _k^{(F)} =  \frac{1}{2(2-a)}\left( \beta -\gamma -\delta +1 -j+k +(1-a)\left( \gamma-\delta +1+k \right)\right)
\end{cases}\nonumber 
\end{equation}  
\end{enumerate}
\subsubsection{The second species complete polynomial}
Replacing coefficients $a$, $q$, $\alpha $, $\beta $, $\gamma $, $\delta$ and $x$ by $1-a$, $-q+(\delta -1)\gamma a+(\alpha -\delta +1)(\beta -\delta +1)$, $\alpha-\delta +1 $, $\beta-\delta +1 $, $2-\delta$, $\gamma $ and $1-x$ into (\ref{eq:40037a})--(\ref{eq:40038c}). Multiply $(1-x)^{1-\delta }$ and the new (\ref{eq:40037a})--(\ref{eq:40037c}) together.\footnote{I treat $\gamma$ and $\delta$ as free variables and fixed values of $\alpha $, $\beta $ and $q$.}
\begin{enumerate} 
\item As $\alpha =\delta -1 $, $\beta =\delta $ and $ q=(\delta -1)\gamma a  $,
 
Its eigenfunction is given by
\begin{eqnarray}
& &(1-x)^{1-\delta } y(\xi)\nonumber\\
&=& (1-x)^{1-\delta } Hl(1-a,0; 0, 1, 2-\delta, \gamma; 1-x) \nonumber\\
&=& (1-x)^{1-\delta } \nonumber
\end{eqnarray}
\item As $ \alpha = \delta-2N-2 $, $\beta = \delta-2N-1 $ and $q = (\delta -1)\gamma a- \left( 2N+1 \right) \left[  \gamma -(1-a)\left( \gamma-\delta +2N +2 \right)\right]$ where $N \in \mathbb{N}_{0}$,

Its eigenfunction is given by
\begin{eqnarray} 
& &(1-x)^{1-\delta } y(\xi)\nonumber\\
&=& (1-x)^{1-\delta } Hl(1-a,\left( 2N+1 \right) \left[ \gamma +2N -(1-a)\left( \gamma-\delta +2N +2 \right)\right]; -2N-1,-2N \nonumber\\
&&, 2-\delta, \gamma; 1-x) \nonumber\\
&=& (1-x)^{1-\delta } \left\{ \sum_{r=0}^{N} y_{2r}^{N-r}\left( 2N+1; \xi \right) + \sum_{r=0}^{N} y_{2r+1}^{N-r}\left( 2N+1;\xi \right) \right\}
\nonumber
\end{eqnarray}
\item As $\alpha =\delta -2N-3 $, $\beta = \delta-2N-2 $ and $q = (\delta -1)\gamma a-\left( 2N+2 \right) \left[ \gamma -(1-a)\left( \gamma -\delta +2N+3 \right)\right]$ where $N \in \mathbb{N}_{0}$,

Its eigenfunction is given by
\begin{eqnarray} 
& &(1-x)^{1-\delta } y(\xi)\nonumber\\
&=& (1-x)^{1-\delta } Hl(1-a,\left( 2N+2 \right) \left[ \gamma +2N+1 -(1-a)\left( \gamma -\delta +2N+3 \right)\right]; -2N-2, -2N-1\nonumber\\
&&, 2-\delta, \gamma; 1-x) \nonumber\\
&=& (1-x)^{1-\delta } \left\{ \sum_{r=0}^{N+1} y_{2r}^{N+1-r}\left( 2N+2;\xi \right) + \sum_{r=0}^{N} y_{2r+1}^{N-r}\left( 2N+2;\xi \right) \right\}
\nonumber
\end{eqnarray}
In the above,
\begin{eqnarray}
y_0^m(j;\xi ) &=&  \sum_{i_0=0}^{m} \frac{\left( -\frac{j}{2}\right)_{i_0} \left( - \frac{j}{2}+\frac{1}{2}\right)_{i_0}}{\left( 1 \right)_{i_0} \left( -\frac{\delta }{2}+ \frac{3}{2} \right)_{i_0}} z^{i_0} \nonumber\\
y_1^m(j;\xi) &=&  \left\{\sum_{i_0=0}^{m} \frac{\left( i_0-\frac{j}{2}\right)\left(i_0 +\Gamma ^{(S)} \right)}{\left( i_0+\frac{1}{2} \right) \left( i_0+1-\frac{\delta }{2} \right)} \frac{\left( -\frac{j}{2}\right)_{i_0} \left( -\frac{j}{2}+\frac{1}{2}\right)_{i_0}}{\left( 1 \right)_{i_0} \left( -\frac{\delta }{2}+ \frac{3}{2} \right)_{i_0}} \right. \nonumber\\
&&\times \left. \sum_{i_1 = i_0}^{m} \frac{\left( -\frac{j}{2}+\frac{1}{2} \right)_{i_1} \left( -\frac{j}{2}+1\right)_{i_1}\left( \frac{3}{2} \right)_{i_0} \left( -\frac{\delta }{2}+2 \right)_{i_0}}{\left(-\frac{j}{2}+\frac{1}{2} \right)_{i_0} \left( -\frac{j}{2}+1\right)_{i_0}\left( \frac{3}{2} \right)_{i_1} \left( -\frac{\delta }{2} +2 \right)_{i_1}} z^{i_1}\right\} \eta \nonumber\\ 
y_{\tau }^m(j;\xi) &=&  \left\{ \sum_{i_0=0}^{m} \frac{\left( i_0-\frac{j}{2}\right)\left(i_0 +\Gamma ^{(S)} \right)}{\left( i_0+\frac{1}{2} \right) \left( i_0+1-\frac{\delta }{2} \right)} \frac{\left( -\frac{j}{2}\right)_{i_0} \left( -\frac{j}{2}+\frac{1}{2}\right)_{i_0}}{\left( 1 \right)_{i_0} \left( -\frac{\delta }{2}+ \frac{3}{2} \right)_{i_0}} \right.\nonumber\\
&&\times \prod_{k=1}^{\tau -1} \left( \sum_{i_k = i_{k-1}}^{m} \frac{\left( i_k+ \frac{k}{2}-\frac{j}{2}\right)\left(i_k + \frac{k}{2}+\Gamma ^{(S)} \right)}{\left( i_k+\frac{k}{2}+\frac{1}{2} \right) \left( i_k+\frac{k}{2}+1-\frac{\delta }{2} \right)} \right. \nonumber\\
&&\times \left. \frac{\left( -\frac{j}{2}+\frac{k}{2}\right)_{i_k} \left( -\frac{j}{2}+ \frac{k}{2}+ \frac{1}{2}\right)_{i_k}\left(  \frac{k}{2} +1\right)_{i_{k-1}} \left( -\frac{\delta }{2}+ \frac{3}{2}+ \frac{k}{2} \right)_{i_{k-1}}}{\left( -\frac{j}{2}+\frac{k}{2}\right)_{i_{k-1}} \left( -\frac{j}{2}+ \frac{k}{2}+ \frac{1}{2}\right)_{i_{k-1}}\left( \frac{k}{2}+1 \right)_{i_k} \left( -\frac{\delta }{2}+ \frac{3}{2}+ \frac{k}{2} \right)_{i_k}} \right) \nonumber\\
&&\times \left. \sum_{i_{\tau } = i_{\tau -1}}^{m} \frac{\left( -\frac{j}{2}+\frac{\tau }{2}\right)_{i_{\tau }} \left( -\frac{j}{2}+ \frac{\tau }{2}+ \frac{1}{2}\right)_{i_{\tau }}\left(  \frac{\tau }{2} +1\right)_{i_{\tau -1}} \left( -\frac{\delta }{2}+ \frac{3}{2}+\frac{\tau }{2} \right)_{i_{\tau -1}}}{\left( -\frac{j}{2}+\frac{\tau }{2}\right)_{i_{\tau -1}} \left(-\frac{j}{2}+ \frac{\tau }{2}+ \frac{1}{2}\right)_{i_{\tau -1}}\left(  \frac{\tau }{2}+1 \right)_{i_{\tau }} \left( -\frac{\delta }{2}+ \frac{3}{2}+\frac{\tau }{2} \right)_{i_{\tau }}} z^{i_{\tau }}\right\} \eta ^{\tau } \nonumber
\end{eqnarray}
where
\begin{equation}
\begin{cases} \tau \geq 2 \cr
\xi =1-x \cr
z = \frac{-1}{1-a}\xi^2 \cr
\eta = \frac{2-a}{1-a}\xi \cr
\Gamma ^{(S)} = \frac{1}{2(2-a)}\left( -\gamma -j+1 +(1-a)\left( \gamma -\delta +1+j  \right)\right) 
\end{cases}\nonumber 
\end{equation} 
\end{enumerate}
\subsection{\footnotesize ${\displaystyle x^{-\alpha } Hl\left(\frac{1}{a},\frac{q+\alpha [(\alpha -\gamma -\delta +1)a-\beta +\delta ]}{a}; \alpha , \alpha -\gamma +1, \alpha -\beta +1,\delta ;\frac{1}{x}\right)}$\normalsize}
\subsubsection{The first species complete polynomial}
Replacing coefficients $a$, $q$, $\alpha $, $\beta $, $\gamma $ and $x$ by $\frac{1}{a}$, $\frac{q+\alpha [(\alpha -\gamma -\delta +1)a-\beta +\delta ]}{a}$, $\alpha-\gamma +1 $, $\alpha $,  $\alpha -\beta +1 $ and $\frac{1}{x}$ into (\ref{eq:40018})--(\ref{eq:40022c}). Multiply $x^{-\alpha }$ and the new (\ref{eq:40018}), (\ref{eq:40019b}) and (\ref{eq:40020b})  together.\footnote{I treat $\alpha $, $\beta$ and $\delta$ as free variables and fixed values of $\gamma $ and $q$.}
\begin{enumerate} 
\item As $\gamma =\alpha -1 $ and $ q = \alpha ((a-1)\delta +\beta )+ a q_0^0 $ where $q_0^0=0$,

The eigenfunction is given by
\begin{eqnarray}
& &x^{-\alpha } y(\xi)\nonumber\\
&=& x^{-\alpha }  Hl\left(\frac{1}{a},0; 0, \alpha , \alpha -\beta +1,\delta ;\frac{1}{x}\right) \nonumber\\
&=& x^{-\alpha } \nonumber 
\end{eqnarray}
\item As $\gamma =\alpha +2N+2 $ where $N \in \mathbb{N}_{0}$,

An algebraic equation of degree $2N+2$ for the determination of $q$ is given by
\begin{equation}
0 = \sum_{r=0}^{N+1}\bar{c}\left( 2r, N+1-r; 2N+1, \frac{q-\alpha [ (\delta +2N+1)a+\beta -\delta ]}{a}\right)\nonumber 
\end{equation}
The eigenvalue of $q$ is written by $ \alpha [ ( \delta +2N+1)a +\beta -\delta ] +a q_{2N+1}^m$ where $m = 0,1,2,\cdots,2N+1 $; $q_{2N+1}^0 < q_{2N+1}^1 < \cdots < q_{2N+1}^{2N+1}$. Its eigenfunction is given by
\begin{eqnarray} 
& &x^{-\alpha } y(\xi)\nonumber\\
&=& x^{-\alpha }  Hl\left(\frac{1}{a},q_{2N+1}^m; -2N-1, \alpha , \alpha -\beta +1,\delta ;\frac{1}{x}\right) \nonumber\\
&=& x^{-\alpha } \left\{ \sum_{r=0}^{N} y_{2r}^{N-r}\left( 2N+1,q_{2N+1}^m;\xi\right)+ \sum_{r=0}^{N} y_{2r+1}^{N-r}\left( 2N+1,q_{2N+1}^m;\xi\right) \right\} 
\nonumber 
\end{eqnarray}
\item As $\gamma =\alpha +2N+3$ where $N \in \mathbb{N}_{0}$,

An algebraic equation of degree $2N+3$ for the determination of $q$ is given by
\begin{eqnarray}
0  = \sum_{r=0}^{N+1}\bar{c}\left( 2r+1, N+1-r; 2N+2, \frac{q-\alpha [ ( \delta +2N+2)a+\beta -\delta ]}{a}\right)\nonumber 
\end{eqnarray}
The eigenvalue of $q$ is written by $\alpha [ ( \delta +2N+2)a+\beta -\delta ]+ a q_{2N+2}^m$ where $m = 0,1,2,\cdots,2N+2 $; $q_{2N+2}^0 < q_{2N+2}^1 < \cdots < q_{2N+2}^{2N+2}$. Its eigenfunction is given by
\begin{eqnarray} 
& &x^{-\alpha } y(\xi)\nonumber\\
&=& x^{-\alpha }  Hl\left(\frac{1}{a}, q_{2N+2}^m; -2N-2, \alpha , \alpha -\beta +1,\delta ;\frac{1}{x}\right) \nonumber\\
&=& x^{-\alpha } \left\{ \sum_{r=0}^{N+1} y_{2r}^{N+1-r}\left( 2N+2,q_{2N+2}^m; \xi\right) + \sum_{r=0}^{N} y_{2r+1}^{N-r}\left( 2N+2,q_{2N+2}^m;\xi \right) \right\}
\nonumber 
\end{eqnarray}
In the above,
\begin{eqnarray}
\bar{c}(0,n;j,\tilde{q})  &=& \frac{\left( -\frac{j}{2}\right)_{n} \left( \frac{\alpha }{2} \right)_{n}}{\left( 1 \right)_{n} \left( \frac{\alpha}{2} -\frac{\beta}{2} +1 \right)_{n}} \left( -a \right)^{n}\nonumber \\
\bar{c}(1,n;j,\tilde{q}) &=& \left( 1+a \right) \sum_{i_0=0}^{n}\frac{ i_0 \left(i_0 +\Gamma _0^{(F)} \right)+\frac{a\tilde{q}}{4(1+a)}}{\left( i_0+\frac{1}{2} \right) \left( i_0+\frac{\alpha}{2} -\frac{\beta}{2} +\frac{1}{2} \right)} \frac{\left( -\frac{j}{2}\right)_{i_0} \left( \frac{\alpha }{2} \right)_{i_0}}{\left( 1 \right)_{i_0} \left( \frac{\alpha}{2} -\frac{\beta}{2} +1 \right)_{i_0}}\nonumber\\
&&\times  \frac{\left( -\frac{j}{2}+\frac{1}{2} \right)_{n} \left( \frac{\alpha }{2}+ \frac{1}{2} \right)_{n}\left( \frac{3}{2} \right)_{i_0} \left( \frac{\alpha}{2} -\frac{\beta}{2} +\frac{3}{2} \right)_{i_0}}{\left( -\frac{j}{2}+\frac{1}{2}\right)_{i_0} \left( \frac{\alpha }{2}+ \frac{1}{2} \right)_{i_0}\left( \frac{3}{2} \right)_{n} \left( \frac{\alpha}{2} -\frac{\beta}{2} +\frac{3}{2}  \right)_{n}} \left( -a \right)^{n}  
\nonumber \\
\bar{c}(\tau ,n;j,\tilde{q}) &=& \left( 1+a \right)^{\tau} \sum_{i_0=0}^{n}\frac{ i_0 \left(i_0 +\Gamma _0^{(F)} \right)+\frac{a\tilde{q}}{4(1+a)}}{\left( i_0+\frac{1}{2} \right) \left( i_0+\frac{\alpha}{2} -\frac{\beta}{2} +\frac{1}{2} \right)} \frac{\left( -\frac{j}{2}\right)_{i_0} \left( \frac{\alpha }{2} \right)_{i_0}}{\left( 1 \right)_{i_0} \left( \frac{\alpha}{2} -\frac{\beta}{2} +1 \right)_{i_0}} \nonumber\\
&&\times \prod_{k=1}^{\tau -1} \left( \sum_{i_k = i_{k-1}}^{n} \frac{\left( i_k+ \frac{k}{2} \right)\left(i_k +\Gamma _k^{(F)} \right)+\frac{a\tilde{q}}{4(1+a)}}{\left( i_k+\frac{k}{2}+\frac{1}{2} \right) \left( i_k+\frac{k}{2}+\frac{\alpha}{2} -\frac{\beta}{2} +\frac{1}{2} \right)} \right. \nonumber\\
&&\times \left. \frac{\left( -\frac{j}{2}+\frac{k}{2}\right)_{i_k} \left( \frac{\alpha }{2}+ \frac{k}{2} \right)_{i_k}\left( \frac{k}{2}+1 \right)_{i_{k-1}} \left( \frac{\alpha}{2} -\frac{\beta}{2} +1+ \frac{k}{2} \right)_{i_{k-1}}}{\left( -\frac{j}{2}+\frac{k}{2}\right)_{i_{k-1}} \left( \frac{\alpha }{2}+ \frac{k}{2} \right)_{i_{k-1}}\left( \frac{k}{2} +1\right)_{i_k} \left( \frac{\alpha}{2} -\frac{\beta}{2} +1+ \frac{k}{2} \right)_{i_k}} \right) \nonumber\\
&&\times \frac{\left( -\frac{j}{2}+\frac{\tau }{2}\right)_{n} \left( \frac{\alpha }{2}+ \frac{\tau }{2} \right)_{n}\left( \frac{\tau }{2} +1\right)_{i_{\tau -1}} \left( \frac{\alpha}{2} -\frac{\beta}{2} +1+\frac{\tau }{2} \right)_{i_{\tau -1}}}{\left( -\frac{j}{2}+\frac{\tau }{2}\right)_{i_{\tau -1}} \left( \frac{\alpha }{2}+\frac{\tau }{2} \right)_{i_{\tau -1}}\left( \frac{\tau }{2}+1 \right)_{n} \left( \frac{\alpha}{2} -\frac{\beta}{2} +1+\frac{\tau }{2} \right)_{n}} \left( -a\right)^{n } \nonumber 
\end{eqnarray}
\begin{eqnarray}
y_0^m(j,\tilde{q};\xi) &=& \sum_{i_0=0}^{m} \frac{\left( -\frac{j}{2}\right)_{i_0} \left( \frac{\alpha }{2} \right)_{i_0}}{\left( 1 \right)_{i_0} \left( \frac{\alpha}{2} -\frac{\beta}{2} +1 \right)_{i_0}} z^{i_0} \nonumber\\
y_1^m(j,\tilde{q};\xi) &=& \left\{\sum_{i_0=0}^{m} \frac{ i_0 \left(i_0 +\Gamma _0^{(F)} \right)+\frac{a\tilde{q}}{4(1+a)}}{\left( i_0+\frac{1}{2} \right) \left( i_0+\frac{\alpha}{2} -\frac{\beta}{2} +\frac{1}{2} \right)} \frac{\left( -\frac{j}{2}\right)_{i_0} \left( \frac{\alpha }{2} \right)_{i_0}}{\left( 1 \right)_{i_0} \left( \frac{\alpha}{2} -\frac{\beta}{2} +1 \right)_{i_0}} \right. \nonumber\\
&&\times \left. \sum_{i_1 = i_0}^{m} \frac{\left( -\frac{j}{2}+\frac{1}{2} \right)_{i_1} \left( \frac{\alpha  }{2}+ \frac{1}{2} \right)_{i_1}\left( \frac{3}{2} \right)_{i_0} \left( \frac{\alpha}{2} -\frac{\beta}{2} +\frac{3}{2} \right)_{i_0}}{\left(-\frac{j}{2}+\frac{1}{2} \right)_{i_0} \left( \frac{\alpha }{2}+ \frac{1}{2} \right)_{i_0}\left( \frac{3}{2} \right)_{i_1} \left( \frac{\alpha}{2} -\frac{\beta}{2} +\frac{3}{2} \right)_{i_1}} z^{i_1}\right\} \eta 
\nonumber\\
y_{\tau }^m(j,\tilde{q};\xi) &=& \left\{ \sum_{i_0=0}^{m} \frac{ i_0 \left(i_0 +\Gamma _0^{(F)} \right)+\frac{a\tilde{q}}{4(1+a)}}{\left( i_0+\frac{1}{2} \right) \left( i_0+\frac{\alpha}{2} -\frac{\beta}{2} +\frac{1}{2} \right)} \frac{\left( -\frac{j}{2}\right)_{i_0} \left( \frac{\alpha }{2} \right)_{i_0}}{\left( 1 \right)_{i_0} \left(\frac{\alpha}{2} -\frac{\beta}{2} +1 \right)_{i_0}} \right.\nonumber\\
&&\times \prod_{k=1}^{\tau -1} \left( \sum_{i_k = i_{k-1}}^{m} \frac{\left( i_k+ \frac{k}{2} \right)\left(i_k +\Gamma _k^{(F)} \right)+\frac{a\tilde{q}}{4(1+a)}}{\left( i_k+\frac{k}{2}+\frac{1}{2} \right) \left( i_k+\frac{k}{2}+\frac{\alpha}{2} -\frac{\beta}{2} +\frac{1}{2} \right)} \right. \nonumber\\
&&\times \left. \frac{\left( -\frac{j}{2}+\frac{k}{2}\right)_{i_k} \left( \frac{\alpha }{2}+ \frac{k}{2} \right)_{i_k}\left( \frac{k}{2}+1 \right)_{i_{k-1}} \left( \frac{\alpha}{2} -\frac{\beta}{2} +1+ \frac{k}{2} \right)_{i_{k-1}}}{\left( -\frac{j}{2}+\frac{k}{2}\right)_{i_{k-1}} \left( \frac{\alpha }{2}+ \frac{k}{2} \right)_{i_{k-1}}\left( \frac{k}{2} +1\right)_{i_k} \left( \frac{\alpha}{2} -\frac{\beta}{2} +1+ \frac{k}{2} \right)_{i_k}} \right) \nonumber\\
&&\times \left. \sum_{i_{\tau } = i_{\tau -1}}^{m} \frac{\left( -\frac{j}{2}+\frac{\tau }{2}\right)_{i_{\tau }} \left( \frac{\alpha }{2}+ \frac{\tau }{2} \right)_{i_{\tau }}\left( \frac{\tau }{2}+1 \right)_{i_{\tau -1}} \left( \frac{\alpha}{2} -\frac{\beta}{2} +1+\frac{\tau }{2} \right)_{i_{\tau -1}}}{\left( -\frac{j}{2}+\frac{\tau }{2}\right)_{i_{\tau -1}} \left( \frac{\alpha }{2}+ \frac{\tau }{2} \right)_{i_{\tau -1}}\left( \frac{\tau }{2}+1 \right)_{i_{\tau }} \left( \frac{\alpha}{2} -\frac{\beta}{2} +1+\frac{\tau }{2} \right)_{i_{\tau }}} z^{i_{\tau }}\right\} \eta ^{\tau } \nonumber
\end{eqnarray}
where
\begin{equation}
\begin{cases} \tau \geq 2 \cr
\xi =\frac{1}{x} \cr
z = -a \xi^2 \cr
\eta = (1+a)\xi \cr
\tilde{q} = \frac{q+\alpha [(\alpha -\gamma -\delta +1)a-\beta +\delta ]}{a} \cr
\Gamma _0^{(F)} =  \frac{a}{2(1+a)}\left( \alpha -\delta -j +\frac{1}{a}\left( \alpha -\beta +\delta \right)\right) \cr
\Gamma _k^{(F)} = \frac{a}{2(1+a)}\left( \alpha -\delta -j+k +\frac{1}{a}\left( \alpha -\beta +\delta +k \right)\right) 
\end{cases}\nonumber 
\end{equation}
\end{enumerate}
\subsection{ \footnotesize ${\displaystyle \left(1-\frac{x}{a} \right)^{-\beta } Hl\left(1-a, -q+\gamma \beta; -\alpha +\gamma +\delta, \beta, \gamma, \delta; \frac{(1-a)x}{x-a} \right)}$\normalsize}
\subsubsection{The first species complete polynomial}
\paragraph{The case of $\gamma $, $q$ = fixed values and $\alpha, \beta, \delta $ = free variables }
Replacing coefficients $a$, $q$, $\alpha $ and $x$ by $1-a$, $-q+\gamma \beta $, $-\alpha+\gamma +\delta $ and $\frac{(1-a)x}{x-a}$ into (\ref{eq:40018})--(\ref{eq:40022c}). Replacing $\gamma $ by $\alpha -\delta -j$ into the new (\ref{eq:40021a})--(\ref{eq:40022c}). Multiply $\left(1-\frac{x}{a} \right)^{-\beta }$ and the new (\ref{eq:40018}), (\ref{eq:40019b}) and (\ref{eq:40020b})  together. 
\begin{enumerate} 
\item As $ \gamma =\alpha - \delta $ and $ q= (\alpha -\delta ) \beta -q_0^0$  where $q_0^0=0$,

The eigenfunction is given by
\begin{eqnarray}
 && \left(1-\frac{x}{a} \right)^{-\beta } y(\xi ) \nonumber\\
 &=& \left( 1-\frac{x}{a} \right)^{-\beta } Hl\left( 1-a, 0; 0, \beta, \gamma, \delta; \frac{(1-a)x}{x-a} \right) \nonumber\\
&=& \left(1-\frac{x}{a} \right)^{-\beta } \nonumber 
\end{eqnarray}
\item As $ \gamma =\alpha -\delta -2N-1$ where $N \in \mathbb{N}_{0}$,

An algebraic equation of degree $2N+2$ for the determination of $q$ is given by
\begin{equation}
0 = \sum_{r=0}^{N+1}\bar{c}\left( 2r, N+1-r; 2N+1, -q +\left( \alpha -\delta -2N-1\right) \beta\right)\nonumber 
\end{equation}
The eigenvalue of $q$ is written by $\left( \alpha -\delta -2N-1\right) \beta -q_{2N+1}^m$ where $m = 0,1,2,\cdots,2N+1 $; $q_{2N+1}^0 < q_{2N+1}^1 < \cdots < q_{2N+1}^{2N+1}$. Its eigenfunction is given by
\begin{eqnarray} 
&& \left(1-\frac{x}{a} \right)^{-\beta } y(\xi ) \nonumber\\
 &=& \left(1-\frac{x}{a} \right)^{-\beta } Hl\left( 1-a, q_{2N+1}^m; -2N-1, \beta, \gamma, \delta; \frac{(1-a)x}{x-a} \right) \nonumber\\
&=& \left(1-\frac{x}{a} \right)^{-\beta } \left\{ \sum_{r=0}^{N} y_{2r}^{N-r}\left( 2N+1, q_{2N+1}^m;\xi\right)+ \sum_{r=0}^{N} y_{2r+1}^{N-r}\left( 2N+1, q_{2N+1}^m;\xi\right) \right\} 
\nonumber 
\end{eqnarray}
\item As $ \gamma =\alpha -\delta -2N-2$ where $N \in \mathbb{N}_{0}$,

An algebraic equation of degree $2N+3$ for the determination of $q$ is given by
\begin{eqnarray}
0  = \sum_{r=0}^{N+1}\bar{c}\left( 2r+1, N+1-r; 2N+2, -q+\left(\alpha- \delta -2N-2\right) \beta \right)\nonumber 
\end{eqnarray}
The eigenvalue of $q$ is written by $\left(\alpha- \delta -2N-2\right) \beta -q_{2N+2}^m$ where $m = 0,1,2,\cdots,2N+2 $; $q_{2N+2}^0 < q_{2N+2}^1 < \cdots < q_{2N+2}^{2N+2}$. Its eigenfunction is given by
\begin{eqnarray} 
&& \left(1-\frac{x}{a} \right)^{-\beta } y(\xi ) \nonumber\\
 &=& \left(1-\frac{x}{a} \right)^{-\beta } Hl\left( 1-a, q_{2N+2}^m; -2N-2, \beta, \gamma, \delta; \frac{(1-a)x}{x-a} \right) \nonumber\\
&=& \left(1-\frac{x}{a} \right)^{-\beta } \left\{ \sum_{r=0}^{N+1} y_{2r}^{N+1-r}\left( 2N+2, q_{2N+2}^m;\xi\right) + \sum_{r=0}^{N} y_{2r+1}^{N-r}\left( 2N+2, q_{2N+2}^m;\xi\right) \right\}
\nonumber 
\end{eqnarray}
In the above,
\begin{eqnarray}
\bar{c}(0,n;j, \tilde{q})  &=& \frac{\left( -\frac{j}{2}\right)_{n} \left( \frac{\beta }{2} \right)_{n}}{\left( 1 \right)_{n} \left( \frac{\alpha}{2} -\frac{\delta}{2} -\frac{j}{2}+ \frac{1}{2} \right)_{n}} \left( \frac{-1}{1-a} \right)^{n}\nonumber\\
\bar{c}(1,n;j, \tilde{q}) &=& \left( \frac{2-a}{1-a}\right) \sum_{i_0=0}^{n}\frac{ i_0 \left(i_0 +\Gamma _0^{(F)} \right)+\frac{\tilde{q}}{4(2-a)}}{\left( i_0+\frac{1}{2} \right) \left( i_0+\frac{\alpha}{2} -\frac{\delta}{2} -\frac{j}{2} \right)} \frac{\left( -\frac{j}{2}\right)_{i_0} \left( \frac{\beta }{2} \right)_{i_0}}{\left( 1 \right)_{i_0} \left( \frac{\alpha}{2} -\frac{\delta}{2} -\frac{j}{2}+ \frac{1}{2} \right)_{i_0}}\nonumber\\
&&\times  \frac{\left( -\frac{j}{2}+\frac{1}{2} \right)_{n} \left( \frac{\beta }{2}+ \frac{1}{2} \right)_{n}\left( \frac{3}{2} \right)_{i_0} \left( \frac{\alpha}{2} -\frac{\delta}{2} -\frac{j}{2}+1 \right)_{i_0}}{\left( -\frac{j}{2}+\frac{1}{2}\right)_{i_0} \left( \frac{\beta }{2}+ \frac{1}{2} \right)_{i_0}\left( \frac{3}{2} \right)_{n} \left( \frac{\alpha}{2} -\frac{\delta}{2} -\frac{j}{2}+1 \right)_{n}} \left( \frac{-1}{1-a} \right)^{n}  
\nonumber\\
\bar{c}(\tau ,n;j, \tilde{q}) &=& \left( \frac{2-a}{1-a}\right)^{\tau} \sum_{i_0=0}^{n}\frac{ i_0 \left(i_0 +\Gamma _0^{(F)} \right)+\frac{\tilde{q}}{4(2-a)}}{\left( i_0+\frac{1}{2} \right) \left( i_0+\frac{\alpha}{2} -\frac{\delta}{2} -\frac{j}{2} \right)} \frac{\left( -\frac{j}{2}\right)_{i_0} \left( \frac{\beta }{2} \right)_{i_0}}{\left( 1 \right)_{i_0} \left( \frac{\alpha}{2} -\frac{\delta}{2} -\frac{j}{2}+ \frac{1}{2} \right)_{i_0}} \nonumber\\
&&\times \prod_{k=1}^{\tau -1} \left( \sum_{i_k = i_{k-1}}^{n} \frac{\left( i_k+ \frac{k}{2} \right)\left(i_k +\Gamma _k^{(F)} \right)+\frac{\tilde{q}}{4(2-a)}}{\left( i_k+\frac{k}{2}+\frac{1}{2} \right) \left( i_k+\frac{k}{2}+\frac{\alpha}{2} -\frac{\delta}{2} -\frac{j}{2} \right)} \right. \nonumber\\
&&\times \left. \frac{\left( -\frac{j}{2}+\frac{k}{2}\right)_{i_k} \left( \frac{\beta }{2}+ \frac{k}{2} \right)_{i_k}\left( \frac{k}{2}+1 \right)_{i_{k-1}} \left( \frac{\alpha}{2} -\frac{\delta}{2} -\frac{j}{2}+ \frac{1}{2}+ \frac{k}{2} \right)_{i_{k-1}}}{\left( -\frac{j}{2}+\frac{k}{2}\right)_{i_{k-1}} \left( \frac{\beta }{2}+ \frac{k}{2} \right)_{i_{k-1}}\left( \frac{k}{2} +1\right)_{i_k} \left( \frac{\alpha}{2} -\frac{\delta}{2} -\frac{j}{2}+ \frac{1}{2}+ \frac{k}{2} \right)_{i_k}} \right) \nonumber\\
&&\times \frac{\left( -\frac{j}{2}+\frac{\tau }{2}\right)_{n} \left( \frac{\beta }{2}+ \frac{\tau }{2} \right)_{n}\left( \frac{\tau }{2} +1\right)_{i_{\tau -1}} \left( \frac{\alpha}{2} -\frac{\delta}{2} -\frac{j}{2}+ \frac{1}{2}+\frac{\tau }{2} \right)_{i_{\tau -1}}}{\left( -\frac{j}{2}+\frac{\tau }{2}\right)_{i_{\tau -1}} \left( \frac{\beta }{2}+\frac{\tau }{2} \right)_{i_{\tau -1}}\left( \frac{\tau }{2}+1 \right)_{n} \left( \frac{\alpha}{2} -\frac{\delta}{2} -\frac{j}{2}+ \frac{1}{2}+\frac{\tau }{2} \right)_{n}} \left( \frac{-1}{1-a} \right)^{n } \nonumber  
\end{eqnarray}
\begin{eqnarray}
y_0^m(j, \tilde{q};\xi) &=& \sum_{i_0=0}^{m} \frac{\left( -\frac{j}{2}\right)_{i_0} \left( \frac{\beta }{2} \right)_{i_0}}{\left( 1 \right)_{i_0} \left( \frac{\alpha}{2} -\frac{\delta}{2} -\frac{j}{2}+ \frac{1}{2} \right)_{i_0}} z^{i_0} \nonumber\\
y_1^m(j, \tilde{q};\xi) &=& \left\{\sum_{i_0=0}^{m} \frac{ i_0 \left(i_0 +\Gamma _0^{(F)} \right)+\frac{\tilde{q}}{4(2-a)}}{\left( i_0+\frac{1}{2} \right) \left( i_0+\frac{\alpha}{2} -\frac{\delta}{2} -\frac{j}{2} \right)} \frac{\left( -\frac{j}{2}\right)_{i_0} \left( \frac{\beta }{2} \right)_{i_0}}{\left( 1 \right)_{i_0} \left( \frac{\alpha}{2} -\frac{\delta}{2} -\frac{j}{2}+ \frac{1}{2} \right)_{i_0}} \right. \nonumber\\
&&\times \left. \sum_{i_1 = i_0}^{m} \frac{\left( -\frac{j}{2}+\frac{1}{2} \right)_{i_1} \left( \frac{\beta }{2}+ \frac{1}{2} \right)_{i_1}\left( \frac{3}{2} \right)_{i_0} \left( \frac{\alpha}{2} -\frac{\delta}{2} -\frac{j}{2}+1 \right)_{i_0}}{\left(-\frac{j}{2}+\frac{1}{2} \right)_{i_0} \left( \frac{\beta }{2}+ \frac{1}{2} \right)_{i_0}\left( \frac{3}{2} \right)_{i_1} \left( \frac{\alpha}{2} -\frac{\delta}{2} -\frac{j}{2}+1 \right)_{i_1}} z^{i_1}\right\} \eta 
\nonumber\\
y_{\tau }^m(j, \tilde{q};\xi) &=& \left\{ \sum_{i_0=0}^{m} \frac{ i_0 \left(i_0 +\Gamma _0^{(F)} \right)+\frac{\tilde{q}}{4(2-a)}}{\left( i_0+\frac{1}{2} \right) \left( i_0+\frac{\alpha}{2} -\frac{\delta}{2} -\frac{j}{2} \right)} \frac{\left( -\frac{j}{2}\right)_{i_0} \left( \frac{\beta }{2} \right)_{i_0}}{\left( 1 \right)_{i_0} \left( \frac{\alpha}{2} -\frac{\delta}{2} -\frac{j}{2}+ \frac{1}{2} \right)_{i_0}} \right.\nonumber\\
&&\times \prod_{k=1}^{\tau -1} \left( \sum_{i_k = i_{k-1}}^{m} \frac{\left( i_k+ \frac{k}{2} \right)\left(i_k +\Gamma _k^{(F)} \right)+\frac{\tilde{q}}{4(2-a)}}{\left( i_k+\frac{k}{2}+\frac{1}{2} \right) \left( i_k+\frac{k}{2}+\frac{\alpha}{2} -\frac{\delta}{2} -\frac{j}{2} \right)} \right. \nonumber\\
&&\times \left. \frac{\left( -\frac{j}{2}+\frac{k}{2}\right)_{i_k} \left( \frac{\beta }{2}+ \frac{k}{2} \right)_{i_k}\left( \frac{k}{2}+1 \right)_{i_{k-1}} \left( \frac{\alpha}{2} -\frac{\delta}{2} -\frac{j}{2}+ \frac{1}{2}+ \frac{k}{2} \right)_{i_{k-1}}}{\left( -\frac{j}{2}+\frac{k}{2}\right)_{i_{k-1}} \left( \frac{\beta }{2}+ \frac{k}{2} \right)_{i_{k-1}}\left( \frac{k}{2} +1\right)_{i_k} \left( \frac{\alpha}{2} -\frac{\delta}{2} -\frac{j}{2}+ \frac{1}{2}+ \frac{k}{2} \right)_{i_k}} \right) \nonumber\\
&&\times \left. \sum_{i_{\tau } = i_{\tau -1}}^{m} \frac{\left( -\frac{j}{2}+\frac{\tau }{2}\right)_{i_{\tau }} \left( \frac{\beta }{2}+ \frac{\tau }{2} \right)_{i_{\tau }}\left( \frac{\tau }{2}+1 \right)_{i_{\tau -1}} \left( \frac{\alpha}{2} -\frac{\delta}{2} -\frac{j}{2}+ \frac{1}{2}+\frac{\tau }{2} \right)_{i_{\tau -1}}}{\left( -\frac{j}{2}+\frac{\tau }{2}\right)_{i_{\tau -1}} \left( \frac{\beta }{2}+ \frac{\tau }{2} \right)_{i_{\tau -1}}\left( \frac{\tau }{2}+1 \right)_{i_{\tau }} \left( \frac{\alpha}{2} -\frac{\delta}{2} -\frac{j}{2}+ \frac{1}{2}+\frac{\tau }{2} \right)_{i_{\tau }}} z^{i_{\tau }}\right\} \eta ^{\tau }  \nonumber 
\end{eqnarray}
where 
\begin{equation}
\begin{cases} \tau \geq 2 \cr
\xi = \frac{(1-a)x}{x-a} \cr
z = -\frac{1}{1-a}\xi^2 \cr
\eta = \frac{2-a}{1-a} \xi \cr
\tilde{q} =  -q+(\alpha -\delta -j)\beta \cr
\Gamma _0^{(F)} = \frac{1}{2(2-a)}\left( \beta -\delta -j +(1-a)\left( \alpha -1-j \right)\right)  \cr
\Gamma _k^{(F)} =  \frac{1}{2(2-a)}\left( \beta -\delta -j+k +(1-a)\left( \alpha -1-j+k \right)\right)
\end{cases}\nonumber 
\end{equation}
\end{enumerate}
\paragraph{The case of $\delta $, $q$ = fixed values and $\alpha, \beta, \gamma$ = free variables}
Replacing coefficients $a$, $q$, $\alpha $ and $x$ by $1-a$, $-q+\gamma \beta $, $-\alpha+\gamma +\delta $ and $\frac{(1-a)x}{x-a}$ into (\ref{eq:40018})--(\ref{eq:40022c}). Multiply $\left(1-\frac{x}{a} \right)^{-\beta }$ and the new (\ref{eq:40018}), (\ref{eq:40019b}) and (\ref{eq:40020b})  together. 
\begin{enumerate} 
\item As $ \delta =\alpha -\gamma$ and $q= \gamma \beta - q_0^0 $ where $q_0^0=0$,

The eigenfunction is given by
\begin{eqnarray}
&& \left(1-\frac{x}{a} \right)^{-\beta } y(\xi ) \nonumber\\
 &=& \left(1-\frac{x}{a} \right)^{-\beta } Hl\left( 1-a, 0; 0, \beta, \gamma, \delta; \frac{(1-a)x}{x-a} \right) \nonumber\\
&=& \left( 1-\frac{x}{a} \right)^{-\beta } \nonumber
\end{eqnarray}
\item As $ \delta =\alpha -\gamma -2N-1$ where $N \in \mathbb{N}_{0}$,

An algebraic equation of degree $2N+2$ for the determination of $q$ is given by
\begin{equation}
0 = \sum_{r=0}^{N+1}\bar{c}\left( 2r, N+1-r; 2N+1,\tilde{q}\right)\nonumber
\end{equation}
The eigenvalue of $q$ is written by $\gamma \beta -q_{2N+1}^m$ where $m = 0,1,2,\cdots,2N+1 $; $q_{2N+1}^0 < q_{2N+1}^1 < \cdots < q_{2N+1}^{2N+1}$. Its eigenfunction is given by
\begin{eqnarray} 
&& \left(1-\frac{x}{a} \right)^{-\beta } y(\xi ) \nonumber\\
 &=& \left(1-\frac{x}{a} \right)^{-\beta } Hl\left( 1-a, q_{2N+1}^m; -2N-1, \beta, \gamma, \delta; \frac{(1-a)x}{x-a} \right) \nonumber\\
&=& \left(1-\frac{x}{a} \right)^{-\beta } \left\{ \sum_{r=0}^{N} y_{2r}^{N-r}\left( 2N+1,q_{2N+1}^m;\xi\right)+ \sum_{r=0}^{N} y_{2r+1}^{N-r}\left( 2N+1,q_{2N+1}^m;\xi\right) \right\}  
\nonumber
\end{eqnarray}
\item As $ \delta =\alpha -\gamma-2N-2$ where $N \in \mathbb{N}_{0}$,

An algebraic equation of degree $2N+3$ for the determination of $q$ is given by
\begin{eqnarray}
0  = \sum_{r=0}^{N+1}\bar{c}\left( 2r+1, N+1-r; 2N+2,\tilde{q}\right)\nonumber
\end{eqnarray}
The eigenvalue of $q$ is written by $\gamma \beta -q_{2N+2}^m$ where $m = 0,1,2,\cdots,2N+2 $; $q_{2N+2}^0 < q_{2N+2}^1 < \cdots < q_{2N+2}^{2N+2}$. Its eigenfunction is given by
\begin{eqnarray} 
&& \left(1-\frac{x}{a} \right)^{-\beta } y(\xi ) \nonumber\\
 &=& \left(1-\frac{x}{a} \right)^{-\beta } Hl\left(1-a, -q+\gamma \beta; -\alpha +\gamma +\delta, \beta, \gamma, \delta; \frac{(1-a)x}{x-a} \right) \nonumber\\
&=& \left(1-\frac{x}{a} \right)^{-\beta } \left\{ \sum_{r=0}^{N+1} y_{2r}^{N+1-r}\left( 2N+2,q_{2N+2}^m;\xi\right) + \sum_{r=0}^{N} y_{2r+1}^{N-r}\left( 2N+2,q_{2N+2}^m;\xi\right) \right\} 
\nonumber
\end{eqnarray}
In the above,
\begin{eqnarray}
\bar{c}(0,n;j,\tilde{q})  &=& \frac{\left( -\frac{j}{2}\right)_{n} \left( \frac{\beta }{2} \right)_{n}}{\left( 1 \right)_{n} \left( \frac{\gamma }{2}+ \frac{1}{2} \right)_{n}} \left(\frac{-1}{1-a} \right)^{n}\nonumber\\
\bar{c}(1,n;j,\tilde{q}) &=& \left( \frac{2-a}{1-a}\right) \sum_{i_0=0}^{n}\frac{ i_0 \left(i_0 +\Gamma _0^{(F)} \right)+\frac{\tilde{q}}{4(2-a)}}{\left( i_0+\frac{1}{2} \right) \left( i_0+\frac{\gamma }{2} \right)} \frac{\left( -\frac{j}{2}\right)_{i_0} \left( \frac{\beta }{2} \right)_{i_0}}{\left( 1 \right)_{i_0} \left( \frac{\gamma }{2}+ \frac{1}{2} \right)_{i_0}}\nonumber\\
&&\times  \frac{\left( -\frac{j}{2}+\frac{1}{2} \right)_{n} \left( \frac{\beta }{2}+ \frac{1}{2} \right)_{n}\left( \frac{3}{2} \right)_{i_0} \left( \frac{\gamma }{2}+1 \right)_{i_0}}{\left( -\frac{j}{2}+\frac{1}{2}\right)_{i_0} \left( \frac{\beta }{2}+ \frac{1}{2} \right)_{i_0}\left( \frac{3}{2} \right)_{n} \left( \frac{\gamma }{2}+1 \right)_{n}} \left(\frac{-1}{1-a} \right)^{n}  
\nonumber\\
\bar{c}(\tau ,n;j,\tilde{q}) &=& \left( \frac{2-a}{1-a}\right)^{\tau} \sum_{i_0=0}^{n}\frac{ i_0 \left(i_0 +\Gamma _0^{(F)} \right)+\frac{\tilde{q}}{4(2-a)}}{\left( i_0+\frac{1}{2} \right) \left( i_0+\frac{\gamma }{2} \right)} \frac{\left( -\frac{j}{2}\right)_{i_0} \left( \frac{\beta }{2} \right)_{i_0}}{\left( 1 \right)_{i_0} \left( \frac{\gamma }{2}+ \frac{1}{2} \right)_{i_0}} \nonumber\\
&&\times \prod_{k=1}^{\tau -1} \left( \sum_{i_k = i_{k-1}}^{n} \frac{\left( i_k+ \frac{k}{2} \right)\left(i_k +\Gamma _k^{(F)} \right)+\frac{\tilde{q}}{4(2-a)}}{\left( i_k+\frac{k}{2}+\frac{1}{2} \right) \left( i_k+\frac{k}{2}+\frac{\gamma }{2} \right)} \right. \nonumber\\
&&\times  \left. \frac{\left( -\frac{j}{2}+\frac{k}{2}\right)_{i_k} \left( \frac{\beta }{2}+ \frac{k}{2} \right)_{i_k}\left( \frac{k}{2}+1 \right)_{i_{k-1}} \left( \frac{\gamma }{2}+ \frac{1}{2}+ \frac{k}{2} \right)_{i_{k-1}}}{\left( -\frac{j}{2}+\frac{k}{2}\right)_{i_{k-1}} \left( \frac{\beta }{2}+ \frac{k}{2} \right)_{i_{k-1}}\left( \frac{k}{2} +1\right)_{i_k} \left( \frac{\gamma }{2}+ \frac{1}{2}+ \frac{k}{2} \right)_{i_k}} \right) \nonumber\\
&&\times \frac{\left( -\frac{j}{2}+\frac{\tau }{2}\right)_{n} \left( \frac{\beta }{2}+ \frac{\tau }{2} \right)_{n}\left( \frac{\tau }{2} +1\right)_{i_{\tau -1}} \left( \frac{\gamma }{2}+ \frac{1}{2}+\frac{\tau }{2} \right)_{i_{\tau -1}}}{\left( -\frac{j}{2}+\frac{\tau }{2}\right)_{i_{\tau -1}} \left( \frac{\beta }{2}+\frac{\tau }{2} \right)_{i_{\tau -1}}\left( \frac{\tau }{2}+1 \right)_{n} \left( \frac{\gamma }{2}+ \frac{1}{2}+\frac{\tau }{2} \right)_{n}} \left( \frac{-1}{1-a} \right)^{n } \nonumber 
\end{eqnarray}
\begin{eqnarray}
y_0^m(j,\tilde{q};\xi) &=& \sum_{i_0=0}^{m} \frac{\left( -\frac{j}{2}\right)_{i_0} \left( \frac{\beta }{2} \right)_{i_0}}{\left( 1 \right)_{i_0} \left( \frac{\gamma }{2}+ \frac{1}{2} \right)_{i_0}} z^{i_0} \nonumber\\
y_1^m(j,\tilde{q};\xi) &=& \left\{\sum_{i_0=0}^{m} \frac{ i_0 \left(i_0 +\Gamma _0^{(F)} \right)+\frac{\tilde{q}}{4(2-a)}}{\left( i_0+\frac{1}{2} \right) \left( i_0+\frac{\gamma }{2} \right)} \frac{\left( -\frac{j}{2}\right)_{i_0} \left( \frac{\beta }{2} \right)_{i_0}}{\left( 1 \right)_{i_0} \left( \frac{\gamma }{2}+ \frac{1}{2} \right)_{i_0}} \right. \nonumber\\
&&\times \left. \sum_{i_1 = i_0}^{m} \frac{\left( -\frac{j}{2}+\frac{1}{2} \right)_{i_1} \left( \frac{\beta }{2}+ \frac{1}{2} \right)_{i_1}\left( \frac{3}{2} \right)_{i_0} \left( \frac{\gamma }{2}+1 \right)_{i_0}}{\left(-\frac{j}{2}+\frac{1}{2} \right)_{i_0} \left( \frac{\beta }{2}+ \frac{1}{2} \right)_{i_0}\left( \frac{3}{2} \right)_{i_1} \left( \frac{\gamma }{2} +1 \right)_{i_1}} z^{i_1}\right\} \eta 
\nonumber
\end{eqnarray}
\begin{eqnarray}
y_{\tau }^m(j,\tilde{q};\xi) &=& \left\{ \sum_{i_0=0}^{m} \frac{ i_0 \left(i_0 +\Gamma _0^{(F)} \right)+\frac{\tilde{q}}{4(2-a)}}{\left( i_0+\frac{1}{2} \right) \left( i_0+\frac{\gamma }{2} \right)} \frac{\left( -\frac{j}{2}\right)_{i_0} \left( \frac{\beta }{2} \right)_{i_0}}{\left( 1 \right)_{i_0} \left( \frac{\gamma }{2}+ \frac{1}{2} \right)_{i_0}} \right.\nonumber\\
&&\times \prod_{k=1}^{\tau -1} \left( \sum_{i_k = i_{k-1}}^{m} \frac{\left( i_k+ \frac{k}{2} \right)\left(i_k +\Gamma _k^{(F)} \right)+\frac{\tilde{q}}{4(2-a)}}{\left( i_k+\frac{k}{2}+\frac{1}{2} \right) \left( i_k+\frac{k}{2}+\frac{\gamma }{2} \right)} \right. \nonumber\\
&&\times  \left. \frac{\left( -\frac{j}{2}+\frac{k}{2}\right)_{i_k} \left( \frac{\beta }{2}+ \frac{k}{2} \right)_{i_k}\left( \frac{k}{2}+1 \right)_{i_{k-1}} \left( \frac{\gamma }{2}+ \frac{1}{2}+ \frac{k}{2} \right)_{i_{k-1}}}{\left( -\frac{j}{2}+\frac{k}{2}\right)_{i_{k-1}} \left( \frac{\beta }{2}+ \frac{k}{2} \right)_{i_{k-1}}\left( \frac{k}{2} +1\right)_{i_k} \left( \frac{\gamma }{2}+ \frac{1}{2}+ \frac{k}{2} \right)_{i_k}} \right) \nonumber\\
&&\times \left. \sum_{i_{\tau } = i_{\tau -1}}^{m} \frac{\left( -\frac{j}{2}+\frac{\tau }{2}\right)_{i_{\tau }} \left( \frac{\beta }{2}+ \frac{\tau }{2} \right)_{i_{\tau }}\left( \frac{\tau }{2}+1 \right)_{i_{\tau -1}} \left( \frac{\gamma }{2}+ \frac{1}{2}+\frac{\tau }{2} \right)_{i_{\tau -1}}}{\left( -\frac{j}{2}+\frac{\tau }{2}\right)_{i_{\tau -1}} \left( \frac{\beta }{2}+ \frac{\tau }{2} \right)_{i_{\tau -1}}\left( \frac{\tau }{2}+1 \right)_{i_{\tau }} \left( \frac{\gamma }{2}+ \frac{1}{2}+\frac{\tau }{2} \right)_{i_{\tau }}} z^{i_{\tau }}\right\} \eta ^{\tau }  \nonumber 
\end{eqnarray}
where 
\begin{equation}
\begin{cases} \tau \geq 2 \cr
\xi = \frac{(1-a)x}{x-a} \cr
z = -\frac{1}{1-a}\xi^2 \cr
\eta = \frac{2-a}{1-a} \xi \cr
\tilde{q} =  -q+\gamma \beta \cr
\Gamma _0^{(F)} = \frac{1}{2(2-a)}\left( -\alpha +\beta +\gamma +(1-a)\left( \alpha -1-j \right)\right) \cr
\Gamma _k^{(F)} = \frac{1}{2(2-a)}\left(  -\alpha +\beta +\gamma +k +(1-a)\left( \alpha -1-j+k \right)\right) 
\end{cases}\nonumber 
\end{equation}
\end{enumerate}
\subsection{ \footnotesize ${\displaystyle (1-x)^{1-\delta }\left(1-\frac{x}{a} \right)^{-\beta+\delta -1} Hl\left(1-a, -q+\gamma [(\delta -1)a+\beta -\delta +1]; -\alpha +\gamma +1, \beta -\delta+1, \gamma, 2-\delta; \frac{(1-a)x}{x-a} \right)}$ \normalsize}
\subsubsection{The first species complete polynomial}
\paragraph{The case of $\alpha $, $q$ = fixed values and $\gamma, \beta, \delta $ = free variables }
Replacing coefficients $a$, $q$, $\alpha $, $\beta $, $\delta $ and $x$ by $1-a$, $-q+\gamma [(\delta -1)a+\beta -\delta +1]$, $-\alpha +\gamma +1$, $\beta -\delta+1$, $2-\delta $ and $\frac{(1-a)x}{x-a}$ into (\ref{eq:40018})--(\ref{eq:40022c}). Multiply $(1-x)^{1-\delta }\left(1-\frac{x}{a} \right)^{-\beta+\delta -1}$ and the new (\ref{eq:40018}), (\ref{eq:40019b}) and (\ref{eq:40020b}) together. 
\begin{enumerate} 
\item As $ \alpha =\gamma +1 $ and $q =\gamma [(\delta -1)a+\beta -\delta +1]-q_0^0$ where $q_0^0=0$,

The eigenfunction is given by
\begin{eqnarray}
&& (1-x)^{1-\delta }\left(1-\frac{x}{a} \right)^{-\beta+\delta -1} y(\xi ) \nonumber\\
 &=& (1-x)^{1-\delta }\left(1-\frac{x}{a} \right)^{-\beta+\delta -1} Hl\bigg( 1-a, 0; 0, \beta -\delta +1, \gamma , 2-\delta; \frac{(1-a)x}{x-a} \bigg) \nonumber\\
&=& (1-x)^{1-\delta }\left(1-\frac{x}{a} \right)^{-\beta+\delta -1}  \nonumber
\end{eqnarray}
\item As $ \alpha =\gamma +1 +2N+1 $ where $N \in \mathbb{N}_{0}$,

An algebraic equation of degree $2N+2$ for the determination of $q$ is given by
\begin{equation}
0 = \sum_{r=0}^{N+1}\bar{c}\left( 2r, N+1-r; 2N+1,\tilde{q}\right)\nonumber
\end{equation}
The eigenvalue of $q$ is written by $\gamma [(\delta -1)a+\beta -\delta +1]-q_{2N+1}^m$ where $m = 0,1,2,\cdots,2N+1 $; $q_{2N+1}^0 < q_{2N+1}^1 < \cdots < q_{2N+1}^{2N+1}$. Its eigenfunction is given by
\begin{eqnarray} 
&& (1-x)^{1-\delta }\left(1-\frac{x}{a} \right)^{-\beta+\delta -1} y(\xi ) \nonumber\\
 &=& (1-x)^{1-\delta }\left(1-\frac{x}{a} \right)^{-\beta+\delta -1} Hl\bigg( 1-a, q_{2N+1}^m; -2N-1, \beta -\delta +1, \gamma , 2-\delta; \frac{(1-a)x}{x-a} \bigg) \nonumber\\
&=& (1-x)^{1-\delta }\left(1-\frac{x}{a} \right)^{-\beta+\delta -1} \left\{ \sum_{r=0}^{N} y_{2r}^{N-r}\left( 2N+1,q_{2N+1}^m;\xi\right)+ \sum_{r=0}^{N} y_{2r+1}^{N-r}\left( 2N+1,q_{2N+1}^m;\xi\right) \right\} 
\nonumber
\end{eqnarray}
\item As $\alpha = \gamma +1 +2N+2 $ where $N \in \mathbb{N}_{0}$,

An algebraic equation of degree $2N+3$ for the determination of $q$ is given by
\begin{eqnarray}
0  = \sum_{r=0}^{N+1}\bar{c}\left( 2r+1, N+1-r; 2N+2,\tilde{q}\right)\nonumber
\end{eqnarray}
The eigenvalue of $q$ is written by $\gamma [(\delta -1)a+\beta -\delta +1]-q_{2N+2}^m$ where $m = 0,1,2,\cdots,2N+2 $; $q_{2N+2}^0 < q_{2N+2}^1 < \cdots < q_{2N+2}^{2N+2}$. Its eigenfunction is given by
\begin{eqnarray} 
&& (1-x)^{1-\delta }\left(1-\frac{x}{a} \right)^{-\beta+\delta -1} y(\xi ) \nonumber\\
 &=& (1-x)^{1-\delta }\left(1-\frac{x}{a} \right)^{-\beta+\delta -1} Hl\bigg( 1-a, q_{2N+2}^m; -2N-2, \beta -\delta+1, \gamma , 2-\delta; \frac{(1-a)x}{x-a} \bigg) \nonumber\\
&=& (1-x)^{1-\delta }\left(1-\frac{x}{a} \right)^{-\beta+\delta -1} \left\{ \sum_{r=0}^{N+1} y_{2r}^{N+1-r}\left( 2N+2,q_{2N+2}^m;\xi\right) + \sum_{r=0}^{N} y_{2r+1}^{N-r}\left( 2N+2,q_{2N+2}^m;\xi\right) \right\} 
\nonumber
\end{eqnarray}
In the above,
\begin{eqnarray}
\bar{c}(0,n;j,\tilde{q})  &=& \frac{\left( -\frac{j}{2}\right)_{n} \left( \frac{\beta}{2} -\frac{\delta}{2} +\frac{1}{2} \right)_{n}}{\left( 1 \right)_{n} \left( \frac{\gamma }{2}+ \frac{1}{2} \right)_{n}} \left( \frac{-1}{1-a} \right)^{n}\nonumber\\
\bar{c}(1,n;j,\tilde{q}) &=& \left( \frac{2-a}{1-a}\right) \sum_{i_0=0}^{n}\frac{ i_0 \left(i_0 +\Gamma _0^{(F)} \right)+\frac{\tilde{q}}{4(2-a)}}{\left( i_0+\frac{1}{2} \right) \left( i_0+\frac{\gamma }{2} \right)} \frac{\left( -\frac{j}{2}\right)_{i_0} \left(  \frac{\beta}{2} -\frac{\delta}{2} +\frac{1}{2} \right)_{i_0}}{\left( 1 \right)_{i_0} \left( \frac{\gamma }{2}+ \frac{1}{2} \right)_{i_0}}\nonumber\\
&&\times  \frac{\left( -\frac{j}{2}+\frac{1}{2} \right)_{n} \left(  \frac{\beta}{2}  -\frac{\delta}{2}  +1 \right)_{n}\left( \frac{3}{2} \right)_{i_0} \left( \frac{\gamma }{2}+1 \right)_{i_0}}{\left( -\frac{j}{2}+\frac{1}{2}\right)_{i_0} \left(  \frac{\beta}{2}  -\frac{\delta}{2} +1 \right)_{i_0}\left( \frac{3}{2} \right)_{n} \left( \frac{\gamma }{2}+1 \right)_{n}} \left( \frac{-1}{1-a} \right)^{n}  
\nonumber\\
\bar{c}(\tau ,n;j,\tilde{q}) &=& \left( \frac{2-a}{1-a}\right)^{\tau} \sum_{i_0=0}^{n}\frac{ i_0 \left(i_0 +\Gamma _0^{(F)} \right)+\frac{\tilde{q}}{4(2-a)}}{\left( i_0+\frac{1}{2} \right) \left( i_0+\frac{\gamma }{2} \right)} \frac{\left( -\frac{j}{2}\right)_{i_0} \left( \frac{\beta}{2}  -\frac{\delta}{2} +\frac{1}{2} \right)_{i_0}}{\left( 1 \right)_{i_0} \left( \frac{\gamma }{2}+ \frac{1}{2} \right)_{i_0}} \nonumber\\
&&\times \prod_{k=1}^{\tau -1} \left( \sum_{i_k = i_{k-1}}^{n} \frac{\left( i_k+ \frac{k}{2} \right)\left(i_k +\Gamma _k^{(F)} \right) +\frac{\tilde{q}}{4(2-a)}}{\left( i_k+\frac{k}{2}+\frac{1}{2} \right) \left( i_k+\frac{k}{2}+\frac{\gamma }{2} \right)} \right. \nonumber\\
&&\times \left. \frac{\left( -\frac{j}{2}+\frac{k}{2}\right)_{i_k} \left( \frac{\beta}{2}  -\frac{\delta}{2} +\frac{1}{2}+ \frac{k}{2} \right)_{i_k}\left( \frac{k}{2}+1 \right)_{i_{k-1}} \left( \frac{\gamma }{2}+ \frac{1}{2}+ \frac{k}{2} \right)_{i_{k-1}}}{\left( -\frac{j}{2}+\frac{k}{2}\right)_{i_{k-1}} \left( \frac{\beta}{2} -\frac{\delta}{2} +\frac{1}{2}+ \frac{k}{2} \right)_{i_{k-1}}\left( \frac{k}{2} +1\right)_{i_k} \left( \frac{\gamma }{2}+ \frac{1}{2}+ \frac{k}{2} \right)_{i_k}} \right) \nonumber\\
&&\times \frac{\left( -\frac{j}{2}+\frac{\tau }{2}\right)_{n} \left( \frac{\beta}{2}  -\frac{\delta}{2} +\frac{1}{2}+ \frac{\tau }{2} \right)_{n}\left( \frac{\tau }{2} +1\right)_{i_{\tau -1}} \left( \frac{\gamma }{2}+ \frac{1}{2}+\frac{\tau }{2} \right)_{i_{\tau -1}}}{\left( -\frac{j}{2}+\frac{\tau }{2}\right)_{i_{\tau -1}} \left( \frac{\beta}{2}  -\frac{\delta}{2} +\frac{1}{2}+\frac{\tau }{2} \right)_{i_{\tau -1}}\left( \frac{\tau }{2}+1 \right)_{n} \left( \frac{\gamma }{2}+ \frac{1}{2}+\frac{\tau }{2} \right)_{n}} \left( \frac{-1}{1-a} \right)^{n } \nonumber 
\end{eqnarray}
\begin{eqnarray}
y_0^m(j,\tilde{q};\xi) &=& \sum_{i_0=0}^{m} \frac{\left( -\frac{j}{2}\right)_{i_0} \left( \frac{\beta}{2}  -\frac{\delta}{2} +\frac{1}{2} \right)_{i_0}}{\left( 1 \right)_{i_0} \left( \frac{\gamma }{2}+ \frac{1}{2} \right)_{i_0}} z^{i_0} \nonumber\\
y_1^m(j,\tilde{q};\xi) &=& \left\{\sum_{i_0=0}^{m} \frac{ i_0 \left(i_0 +\Gamma _0^{(F)} \right)+\frac{\tilde{q}}{4(2-a)}}{\left( i_0+\frac{1}{2} \right) \left( i_0+\frac{\gamma }{2} \right)} \frac{\left( -\frac{j}{2}\right)_{i_0} \left( \frac{\beta}{2}  -\frac{\delta}{2} +\frac{1}{2} \right)_{i_0}}{\left( 1 \right)_{i_0} \left( \frac{\gamma }{2}+ \frac{1}{2} \right)_{i_0}} \right. \nonumber\\
&&\times \left. \sum_{i_1 = i_0}^{m} \frac{\left( -\frac{j}{2}+\frac{1}{2} \right)_{i_1} \left( \frac{\beta}{2}  -\frac{\delta}{2} +1 \right)_{i_1}\left( \frac{3}{2} \right)_{i_0} \left( \frac{\gamma }{2}+1 \right)_{i_0}}{\left(-\frac{j}{2}+\frac{1}{2} \right)_{i_0} \left( \frac{\beta}{2}  -\frac{\delta}{2} +1 \right)_{i_0}\left( \frac{3}{2} \right)_{i_1} \left( \frac{\gamma }{2} +1 \right)_{i_1}} z^{i_1}\right\} \eta 
\nonumber
\end{eqnarray}
\begin{eqnarray}
y_{\tau }^m(j,\tilde{q};\xi) &=& \left\{ \sum_{i_0=0}^{m} \frac{ i_0 \left(i_0 +\Gamma _0^{(F)} \right)+\frac{\tilde{q}}{4(2-a)}}{\left( i_0+\frac{1}{2} \right) \left( i_0+\frac{\gamma }{2} \right)} \frac{\left( -\frac{j}{2}\right)_{i_0} \left( \frac{\beta}{2}  -\frac{\delta}{2} +\frac{1}{2} \right)_{i_0}}{\left( 1 \right)_{i_0} \left( \frac{\gamma }{2}+ \frac{1}{2} \right)_{i_0}} \right.\nonumber\\
&&\times \prod_{k=1}^{\tau -1} \left( \sum_{i_k = i_{k-1}}^{m} \frac{\left( i_k+ \frac{k}{2} \right)\left(i_k +\Gamma _k^{(F)} \right) +\frac{\tilde{q}}{4(2-a)}}{\left( i_k+\frac{k}{2}+\frac{1}{2} \right) \left( i_k+\frac{k}{2}+\frac{\gamma }{2} \right)} \right. \nonumber\\
&&\times \left. \frac{\left( -\frac{j}{2}+\frac{k}{2}\right)_{i_k} \left( \frac{\beta}{2}  -\frac{\delta}{2} +\frac{1}{2}+ \frac{k}{2} \right)_{i_k}\left( \frac{k}{2}+1 \right)_{i_{k-1}} \left( \frac{\gamma }{2}+ \frac{1}{2}+ \frac{k}{2} \right)_{i_{k-1}}}{\left( -\frac{j}{2}+\frac{k}{2}\right)_{i_{k-1}} \left( \frac{\beta}{2} -\frac{\delta}{2} +\frac{1}{2}+ \frac{k}{2} \right)_{i_{k-1}}\left( \frac{k}{2} +1\right)_{i_k} \left( \frac{\gamma }{2}+ \frac{1}{2}+ \frac{k}{2} \right)_{i_k}} \right) \nonumber\\
&&\times \left. \sum_{i_{\tau } = i_{\tau -1}}^{m} \frac{\left( -\frac{j}{2}+\frac{\tau }{2}\right)_{i_{\tau }} \left( \frac{\beta}{2}  -\frac{\delta}{2} +\frac{1}{2}+ \frac{\tau }{2} \right)_{i_{\tau }}\left( \frac{\tau }{2}+1 \right)_{i_{\tau -1}} \left( \frac{\gamma }{2}+ \frac{1}{2}+\frac{\tau }{2} \right)_{i_{\tau -1}}}{\left( -\frac{j}{2}+\frac{\tau }{2}\right)_{i_{\tau -1}} \left( \frac{\beta}{2} -\frac{\delta}{2} +\frac{1}{2}+ \frac{\tau }{2} \right)_{i_{\tau -1}}\left( \frac{\tau }{2}+1 \right)_{i_{\tau }} \left( \frac{\gamma }{2}+ \frac{1}{2}+\frac{\tau }{2} \right)_{i_{\tau }}} z^{i_{\tau }}\right\} \eta ^{\tau } \nonumber
\end{eqnarray}
where  
\begin{equation}
\begin{cases} 
\tau \geq 2 \cr
\xi = \frac{(1-a)x}{x-a} \cr
z = -\frac{1}{1-a}\xi^2 \cr
\eta = \frac{(2-a)}{1-a} \xi \cr
\tilde{q} = -q+\gamma [(\delta -1)a+\beta -\delta +1]\cr
\Gamma _0^{(F)} = \frac{1}{2(2-a)}\left( \beta -1-j +(1-a)\left( \gamma -\delta +1 \right)\right) \cr
\Gamma _k^{(F)} = \frac{1}{2(2-a)}\left( \beta -1-j+k +(1-a)\left( \gamma -\delta +1+k \right)\right)  
\end{cases}\nonumber 
\end{equation} 
\end{enumerate}
\paragraph{The case of $\gamma $, $q$ = fixed values and $\alpha , \beta, \delta $ = free variables }
Replacing coefficients $a$, $q$, $\alpha $, $\beta $, $\delta $ and $x$ by $1-a$, $-q+\gamma [(\delta -1)a+\beta -\delta +1]$, $-\alpha +\gamma +1$, $\beta -\delta+1$, $2-\delta $ and $\frac{(1-a)x}{x-a}$ into (\ref{eq:40018})--(\ref{eq:40022c}). Replacing $\gamma $ by $\alpha -1-j$  into the new (\ref{eq:40021a})--(\ref{eq:40022c}). Multiply $(1-x)^{1-\delta }\left(1-\frac{x}{a} \right)^{-\beta+\delta -1}$ and the new (\ref{eq:40018}), (\ref{eq:40019b}) and (\ref{eq:40020b}) together. 
\begin{enumerate} 
\item As $ \gamma =\alpha -1$ and $q=(\alpha -1) [(\delta -1)a+\beta -\delta +1]- q_0^0 $ where $q_0^0=0$,

The eigenfunction is given by
\begin{eqnarray}
&& (1-x)^{1-\delta }\left(1-\frac{x}{a} \right)^{-\beta+\delta -1} y(\xi ) \nonumber\\
 &=& (1-x)^{1-\delta }\left(1-\frac{x}{a} \right)^{-\beta+\delta -1} Hl\bigg( 1-a, 0; 0, \beta -\delta +1, \gamma , 2-\delta; \frac{(1-a)x}{x-a} \bigg) \nonumber\\
&=& (1-x)^{1-\delta }\left(1-\frac{x}{a} \right)^{-\beta+\delta -1} \nonumber 
\end{eqnarray}
\item As $ \gamma =\alpha -2N-2$ where $N \in \mathbb{N}_{0}$,

An algebraic equation of degree $2N+2$ for the determination of $q$ is given by
\begin{equation}
0 = \sum_{r=0}^{N+1}\bar{c}\left( 2r, N+1-r; 2N+1,-q+(\alpha -2N-2) [(\delta -1)a+\beta -\delta +1]\right)\nonumber 
\end{equation}
The eigenvalue of $q$ is written by $(\alpha -2N-2) [(\delta -1)a+\beta -\delta +1]-q_{2N+1}^m$ where $m = 0,1,2,\cdots,2N+1 $; $q_{2N+1}^0 < q_{2N+1}^1 < \cdots < q_{2N+1}^{2N+1}$. Its eigenfunction is given by
\begin{eqnarray} 
&& (1-x)^{1-\delta }\left(1-\frac{x}{a} \right)^{-\beta+\delta -1} y(\xi ) \nonumber\\
 &=& (1-x)^{1-\delta }\left(1-\frac{x}{a} \right)^{-\beta+\delta -1} Hl\bigg( 1-a, q_{2N+1}^m; -2N-1, \beta -\delta+1, \gamma , 2-\delta; \frac{(1-a)x}{x-a} \bigg) \nonumber\\
&=& (1-x)^{1-\delta }\left(1-\frac{x}{a} \right)^{-\beta+\delta -1} \left\{ \sum_{r=0}^{N} y_{2r}^{N-r}\left( 2N+1,q_{2N+1}^m;\xi\right)+ \sum_{r=0}^{N} y_{2r+1}^{N-r}\left( 2N+1,q_{2N+1}^m;\xi\right) \right\}  
\nonumber 
\end{eqnarray}
\item As $ \gamma =\alpha -2N-3$ where $N \in \mathbb{N}_{0}$,

An algebraic equation of degree $2N+3$ for the determination of $q$ is given by
\begin{eqnarray}
0  = \sum_{r=0}^{N+1}\bar{c}\left( 2r+1, N+1-r; 2N+2,-q+(\alpha -2N-3)[(\delta -1)a+\beta -\delta +1] \right)\nonumber 
\end{eqnarray}
The eigenvalue of $q$ is written by $(\alpha -2N-3)[(\delta -1)a+\beta -\delta +1]-q_{2N+2}^m$ where $m = 0,1,2,\cdots,2N+2 $; $q_{2N+2}^0 < q_{2N+2}^1 < \cdots < q_{2N+2}^{2N+2}$. Its eigenfunction is given by
\begin{eqnarray} 
&& (1-x)^{1-\delta }\left(1-\frac{x}{a} \right)^{-\beta+\delta -1} y(\xi ) \nonumber\\
 &=& (1-x)^{1-\delta }\left(1-\frac{x}{a} \right)^{-\beta+\delta -1} Hl\bigg( 1-a, q_{2N+2}^m; -2N-2, \beta -\delta+1, \gamma , 2-\delta; \frac{(1-a)x}{x-a} \bigg) \nonumber\\
&=& (1-x)^{1-\delta }\left(1-\frac{x}{a} \right)^{-\beta+\delta -1} \left\{ \sum_{r=0}^{N+1} y_{2r}^{N+1-r}\left( 2N+2,q_{2N+2}^m;\xi\right) + \sum_{r=0}^{N} y_{2r+1}^{N-r}\left( 2N+2,q_{2N+2}^m;\xi\right) \right\} 
\nonumber 
\end{eqnarray}
In the above,
\begin{eqnarray}
\bar{c}(0,n;j,\tilde{q})  &=& \frac{\left( -\frac{j}{2}\right)_{n} \left( \frac{\beta }{2}-\frac{\delta}{2}+\frac{1}{2} \right)_{n}}{\left( 1 \right)_{n} \left( \frac{\alpha }{2} -\frac{j}{2} \right)_{n}} \left( \frac{-1}{1-a} \right)^{n}\nonumber \\
\bar{c}(1,n;j,\tilde{q}) &=& \left( \frac{2-a}{1-a}\right) \sum_{i_0=0}^{n}\frac{ i_0 \left(i_0 +\Gamma _0^{(F)} \right)+\frac{\tilde{q}}{4(2-a)}}{\left( i_0+\frac{1}{2} \right) \left( i_0+\frac{\alpha }{2}-\frac{1}{2}-\frac{j}{2} \right)} \frac{\left( -\frac{j}{2}\right)_{i_0} \left( \frac{\beta }{2}-\frac{\delta}{2}+\frac{1}{2} \right)_{i_0}}{\left( 1 \right)_{i_0} \left( \frac{\alpha }{2} -\frac{j}{2} \right)_{i_0}}\nonumber\\
&&\times  \frac{\left( -\frac{j}{2}+\frac{1}{2} \right)_{n} \left( \frac{\beta }{2}-\frac{\delta}{2}+1 \right)_{n}\left( \frac{3}{2} \right)_{i_0} \left( \frac{\alpha }{2}+\frac{1}{2}-\frac{j}{2} \right)_{i_0}}{\left( -\frac{j}{2}+\frac{1}{2}\right)_{i_0} \left( \frac{\beta }{2}-\frac{\delta}{2}+1 \right)_{i_0}\left( \frac{3}{2} \right)_{n} \left( \frac{\alpha }{2}+\frac{1}{2}-\frac{j}{2} \right)_{n}} \left( \frac{-1}{1-a} \right)^{n}  
\nonumber \\
\bar{c}(\tau ,n;j,\tilde{q}) &=& \left( \frac{2-a}{1-a}\right)^{\tau} \sum_{i_0=0}^{n}\frac{ i_0 \left(i_0 +\Gamma _0^{(F)} \right)+\frac{\tilde{q}}{4(2-a)}}{\left( i_0+\frac{1}{2} \right) \left( i_0+\frac{\alpha }{2}-\frac{1}{2}-\frac{j}{2} \right)} \frac{\left( -\frac{j}{2}\right)_{i_0} \left( \frac{\beta }{2}-\frac{\delta}{2}+\frac{1}{2} \right)_{i_0}}{\left( 1 \right)_{i_0} \left( \frac{\alpha }{2} -\frac{j}{2} \right)_{i_0}} \nonumber\\
&&\times \prod_{k=1}^{\tau -1} \left( \sum_{i_k = i_{k-1}}^{n} \frac{\left( i_k+ \frac{k}{2} \right)\left(i_k +\Gamma _k^{(F)} \right)+\frac{\tilde{q}}{4(2-a)}}{\left( i_k+\frac{k}{2}+\frac{1}{2} \right) \left( i_k+\frac{k}{2}+\frac{\alpha }{2}-\frac{1}{2}-\frac{j}{2} \right)} \right. \nonumber\\
&&\times \left. \frac{\left( -\frac{j}{2}+\frac{k}{2}\right)_{i_k} \left( \frac{\beta }{2}-\frac{\delta}{2}+\frac{1}{2}+ \frac{k}{2} \right)_{i_k}\left( \frac{k}{2}+1 \right)_{i_{k-1}} \left( \frac{\alpha }{2} -\frac{j}{2} + \frac{k}{2} \right)_{i_{k-1}}}{\left( -\frac{j}{2}+\frac{k}{2}\right)_{i_{k-1}} \left( \frac{\beta }{2}-\frac{\delta}{2}+\frac{1}{2}+ \frac{k}{2} \right)_{i_{k-1}}\left( \frac{k}{2} +1\right)_{i_k} \left( \frac{\alpha }{2} -\frac{j}{2} + \frac{k}{2} \right)_{i_k}} \right) \nonumber\\
&&\times \frac{\left( -\frac{j}{2}+\frac{\tau }{2}\right)_{n} \left( \frac{\beta }{2}-\frac{\delta}{2}+\frac{1}{2}+ \frac{\tau }{2} \right)_{n}\left( \frac{\tau }{2} +1\right)_{i_{\tau -1}} \left( \frac{\alpha }{2} -\frac{j}{2} +\frac{\tau }{2} \right)_{i_{\tau -1}}}{\left( -\frac{j}{2}+\frac{\tau }{2}\right)_{i_{\tau -1}} \left( \frac{\beta }{2}-\frac{\delta}{2}+\frac{1}{2}+\frac{\tau }{2} \right)_{i_{\tau -1}}\left( \frac{\tau }{2}+1 \right)_{n} \left( \frac{\alpha }{2} -\frac{j}{2} +\frac{\tau }{2} \right)_{n}} \left( \frac{-1}{1-a} \right)^{n } \nonumber  
\end{eqnarray}
\begin{eqnarray}
y_0^m(j,\tilde{q};\xi) &=& \sum_{i_0=0}^{m} \frac{\left( -\frac{j}{2}\right)_{i_0} \left( \frac{\beta }{2}-\frac{\delta}{2}+\frac{1}{2} \right)_{i_0}}{\left( 1 \right)_{i_0} \left( \frac{\alpha }{2} -\frac{j}{2} \right)_{i_0}} z^{i_0} \nonumber \\
y_1^m(j,\tilde{q};\xi) &=& \left\{\sum_{i_0=0}^{m} \frac{ i_0 \left(i_0 +\Gamma _0^{(F)} \right)+\frac{\tilde{q}}{4(2-a)}}{\left( i_0+\frac{1}{2} \right) \left( i_0+\frac{\alpha }{2}-\frac{1}{2}-\frac{j}{2} \right)} \frac{\left( -\frac{j}{2}\right)_{i_0} \left( \frac{\beta }{2}-\frac{\delta}{2}+\frac{1}{2} \right)_{i_0}}{\left( 1 \right)_{i_0} \left( \frac{\alpha }{2} -\frac{j}{2} \right)_{i_0}} \right. \nonumber\\
&&\times \left. \sum_{i_1 = i_0}^{m} \frac{\left( -\frac{j}{2}+\frac{1}{2} \right)_{i_1} \left( \frac{\beta }{2}-\frac{\delta}{2}+1 \right)_{i_1}\left( \frac{3}{2} \right)_{i_0} \left( \frac{\alpha }{2}+\frac{1}{2}-\frac{j}{2} \right)_{i_0}}{\left(-\frac{j}{2}+\frac{1}{2} \right)_{i_0} \left( \frac{\beta }{2}-\frac{\delta}{2}+1 \right)_{i_0}\left( \frac{3}{2} \right)_{i_1} \left( \frac{\alpha }{2}+\frac{1}{2}-\frac{j}{2} \right)_{i_1}} z^{i_1}\right\} \eta 
\nonumber
\end{eqnarray}
\begin{eqnarray}
y_{\tau }^m(j,\tilde{q};\xi) &=& \left\{ \sum_{i_0=0}^{m} \frac{ i_0 \left(i_0 +\Gamma _0^{(F)} \right)+\frac{\tilde{q}}{4(2-a)}}{\left( i_0+\frac{1}{2} \right) \left( i_0+\frac{\alpha }{2}-\frac{1}{2}-\frac{j}{2} \right)} \frac{\left( -\frac{j}{2}\right)_{i_0} \left( \frac{\beta }{2}-\frac{\delta}{2}+\frac{1}{2} \right)_{i_0}}{\left( 1 \right)_{i_0} \left( \frac{\alpha }{2} -\frac{j}{2} \right)_{i_0}} \right.\nonumber\\
&&\times \prod_{k=1}^{\tau -1} \left( \sum_{i_k = i_{k-1}}^{m} \frac{\left( i_k+ \frac{k}{2} \right)\left(i_k +\Gamma _k^{(F)} \right)+\frac{\tilde{q}}{4(2-a)}}{\left( i_k+\frac{k}{2}+\frac{1}{2} \right) \left( i_k+\frac{k}{2}+\frac{\alpha }{2}-\frac{1}{2}-\frac{j}{2} \right)} \right. \nonumber\\
&&\times \left. \frac{\left( -\frac{j}{2}+\frac{k}{2}\right)_{i_k} \left( \frac{\beta }{2}-\frac{\delta}{2}+\frac{1}{2}+ \frac{k}{2} \right)_{i_k}\left( \frac{k}{2}+1 \right)_{i_{k-1}} \left( \frac{\alpha }{2} -\frac{j}{2} + \frac{k}{2} \right)_{i_{k-1}}}{\left( -\frac{j}{2}+\frac{k}{2}\right)_{i_{k-1}} \left( \frac{\beta }{2}-\frac{\delta}{2}+\frac{1}{2}+ \frac{k}{2} \right)_{i_{k-1}}\left( \frac{k}{2} +1\right)_{i_k} \left( \frac{\alpha }{2} -\frac{j}{2} + \frac{k}{2} \right)_{i_k}} \right) \nonumber\\
&&\times \left. \sum_{i_{\tau } = i_{\tau -1}}^{m} \frac{\left( -\frac{j}{2}+\frac{\tau }{2}\right)_{i_{\tau }} \left( \frac{\beta }{2}-\frac{\delta}{2}+\frac{1}{2}+ \frac{\tau }{2} \right)_{i_{\tau }}\left( \frac{\tau }{2}+1 \right)_{i_{\tau -1}} \left( \frac{\alpha }{2} -\frac{j}{2} +\frac{\tau }{2} \right)_{i_{\tau -1}}}{\left( -\frac{j}{2}+\frac{\tau }{2}\right)_{i_{\tau -1}} \left( \frac{\beta }{2}-\frac{\delta}{2}+\frac{1}{2}+ \frac{\tau }{2} \right)_{i_{\tau -1}}\left( \frac{\tau }{2}+1 \right)_{i_{\tau }} \left( \frac{\alpha }{2} -\frac{j}{2} +\frac{\tau }{2} \right)_{i_{\tau }}} z^{i_{\tau }}\right\} \eta ^{\tau } \nonumber
\end{eqnarray}
where  
\begin{equation}
\begin{cases} 
\tau \geq 2 \cr
\xi = \frac{(1-a)x}{x-a} \cr
z = -\frac{1}{1-a}\xi^2 \cr
\eta = \frac{(2-a)}{1-a} \xi \cr
\tilde{q} = -q+(\alpha -1-j)[(\delta -1)a+\beta -\delta +1]\cr
\Gamma _0^{(F)} = \frac{1}{2(2-a)}\left( \beta -1-j +(1-a)\left( \alpha -\delta -j\right)\right) \cr
\Gamma _k^{(F)} = \frac{1}{2(2-a)}\left( \beta -1-j+k +(1-a)\left( \alpha -\delta -j+k \right)\right) 
\end{cases}\nonumber 
\end{equation} 
\end{enumerate}
\subsection{ \footnotesize ${\displaystyle x^{-\alpha } Hl\left(\frac{a-1}{a}, \frac{-q+\alpha (\delta a+\beta -\delta )}{a}; \alpha, \alpha -\gamma +1, \delta , \alpha -\beta +1; \frac{x-1}{x} \right)}$\normalsize}
\subsubsection{The first species complete polynomial}
Replacing coefficients $a$, $q$, $\alpha $, $\beta $, $\gamma $, $\delta $ and $x$ by $\frac{a-1}{a}$, $\frac{-q+\alpha (\delta a+\beta -\delta )}{a}$, $\alpha -\gamma +1$, $\alpha $, $\delta $, $\alpha -\beta +1$ and $\frac{x-1}{x}$ into (\ref{eq:40018})--(\ref{eq:40022c}). Multiply $x^{-\alpha }$ and the new (\ref{eq:40018}), (\ref{eq:40019b}) and (\ref{eq:40020b}) together.\footnote{I treat $\alpha $, $\beta $ and $\delta$ as free variables and fixed values of $\gamma $ and $q$.}
\begin{enumerate} 
\item As $\gamma =\alpha +1 $ and $q=\alpha (\delta a+\beta -\delta )-a q_0^0 $  where $q_0^0=0$,

The eigenfunction is given by
\begin{eqnarray}
&& x^{-\alpha } y(\xi ) \nonumber\\
 &=& x^{-\alpha } Hl\left(\frac{a-1}{a}, 0; 0, \alpha, \delta , \alpha -\beta +1; \frac{x-1}{x} \right) \nonumber\\
&=& x^{-\alpha } \nonumber
\end{eqnarray}
\item As $\gamma = \alpha +2N+2 $ where $N \in \mathbb{N}_{0}$,

An algebraic equation of degree $2N+2$ for the determination of $q$ is given by
\begin{equation}
0 = \sum_{r=0}^{N+1}\bar{c}\left( 2r, N+1-r; 2N+1,\tilde{q}\right) \nonumber
\end{equation}
The eigenvalue of $q$ is written by $\alpha (\delta a+\beta -\delta )-a q_{2N+1}^m$ where $m = 0,1,2,\cdots,2N+1 $; $q_{2N+1}^0 < q_{2N+1}^1 < \cdots < q_{2N+1}^{2N+1}$. Its eigenfunction is given by
\begin{eqnarray} 
&& x^{-\alpha } y(\xi ) \nonumber\\
 &=& x^{-\alpha } Hl\left(\frac{a-1}{a}, q_{2N+1}^m; -2N-1, \alpha, \delta , \alpha -\beta +1; \frac{x-1}{x} \right) \nonumber\\
&=& x^{-\alpha } \left\{ \sum_{r=0}^{N} y_{2r}^{N-r}\left( 2N+1,q_{2N+1}^m;\xi\right)+ \sum_{r=0}^{N} y_{2r+1}^{N-r}\left( 2N+1,q_{2N+1}^m;\xi\right) \right\}
\nonumber
\end{eqnarray}
\item As $\gamma =\alpha +2N+3 $ where $N \in \mathbb{N}_{0}$,

An algebraic equation of degree $2N+3$ for the determination of $q$ is given by
\begin{eqnarray}
0  = \sum_{r=0}^{N+1}\bar{c}\left( 2r+1, N+1-r; 2N+2,\tilde{q}\right)\nonumber
\end{eqnarray}
The eigenvalue of $q$ is written by $\alpha (\delta a+\beta -\delta )-a q_{2N+2}^m$ where $m = 0,1,2,\cdots,2N+2 $; $q_{2N+2}^0 < q_{2N+2}^1 < \cdots < q_{2N+2}^{2N+2}$. Its eigenfunction is given by
\begin{eqnarray} 
&& x^{-\alpha } y(\xi ) \nonumber\\
 &=& x^{-\alpha } Hl\left(\frac{a-1}{a}, q_{2N+2}^m; -2N-2, \alpha , \delta , \alpha -\beta +1; \frac{x-1}{x} \right) \nonumber\\
&=& x^{-\alpha } \left\{ \sum_{r=0}^{N+1} y_{2r}^{N+1-r}\left( 2N+2,q_{2N+2}^m;\xi\right) + \sum_{r=0}^{N} y_{2r+1}^{N-r}\left( 2N+2,q_{2N+2}^m;\xi\right) \right\}
\nonumber
\end{eqnarray}
In the above,
\begin{eqnarray}
\bar{c}(0,n;j,\tilde{q})  &=& \frac{\left( -\frac{j}{2}\right)_{n} \left( \frac{\alpha }{2} \right)_{n}}{\left( 1 \right)_{n} \left( \frac{\delta }{2}+ \frac{1}{2} \right)_{n}} \left(-\frac{a}{a-1} \right)^{n}\nonumber\\
\bar{c}(1,n;j,\tilde{q}) &=& \left( \frac{2a-1}{a-1}\right) \sum_{i_0=0}^{n}\frac{ i_0 \left(i_0 +\Gamma _0^{(F)} \right)+\frac{a\tilde{q}}{4(2a-1)}}{\left( i_0+\frac{1}{2} \right) \left( i_0+\frac{\delta }{2} \right)} \frac{\left( -\frac{j}{2}\right)_{i_0} \left( \frac{\alpha }{2} \right)_{i_0}}{\left( 1 \right)_{i_0} \left( \frac{\delta }{2}+ \frac{1}{2} \right)_{i_0}} \nonumber\\
&&\times  \frac{\left( -\frac{j}{2}+\frac{1}{2} \right)_{n} \left( \frac{\alpha }{2}+ \frac{1}{2} \right)_{n}\left( \frac{3}{2} \right)_{i_0} \left( \frac{\delta }{2}+1 \right)_{i_0}}{\left( -\frac{j}{2}+\frac{1}{2}\right)_{i_0} \left( \frac{\alpha }{2}+ \frac{1}{2} \right)_{i_0}\left( \frac{3}{2} \right)_{n} \left( \frac{\delta }{2}+1 \right)_{n}} \left(-\frac{a}{a-1} \right)^{n}  
\nonumber\\
\bar{c}(\tau ,n;j,\tilde{q}) &=& \left( \frac{2a-1}{a-1}\right)^{\tau} \sum_{i_0=0}^{n}\frac{ i_0 \left(i_0 +\Gamma _0^{(F)} \right) +\frac{a\tilde{q}}{4(2a-1)}}{\left( i_0+\frac{1}{2} \right) \left( i_0+\frac{\delta }{2} \right)} \frac{\left( -\frac{j}{2}\right)_{i_0} \left( \frac{\alpha }{2} \right)_{i_0}}{\left( 1 \right)_{i_0} \left( \frac{\delta }{2}+ \frac{1}{2} \right)_{i_0}} \nonumber\\
&&\times \prod_{k=1}^{\tau -1} \left( \sum_{i_k = i_{k-1}}^{n} \frac{\left( i_k+ \frac{k}{2} \right)\left(i_k +\Gamma _k^{(F)} \right) +\frac{a\tilde{q}}{4(2a-1)}}{\left( i_k+\frac{k}{2}+\frac{1}{2} \right) \left( i_k+\frac{k}{2}+\frac{\delta }{2} \right)} \right. \nonumber\\
&&\times \left. \frac{\left( -\frac{j}{2}+\frac{k}{2}\right)_{i_k} \left( \frac{\alpha }{2}+ \frac{k}{2} \right)_{i_k}\left( \frac{k}{2}+1 \right)_{i_{k-1}} \left( \frac{\delta }{2}+ \frac{1}{2}+ \frac{k}{2} \right)_{i_{k-1}}}{\left( -\frac{j}{2}+\frac{k}{2}\right)_{i_{k-1}} \left( \frac{\alpha }{2}+ \frac{k}{2} \right)_{i_{k-1}}\left( \frac{k}{2} +1\right)_{i_k} \left( \frac{\delta }{2}+ \frac{1}{2}+ \frac{k}{2} \right)_{i_k}} \right) \nonumber\\
&&\times \frac{\left( -\frac{j}{2}+\frac{\tau }{2}\right)_{n} \left( \frac{\alpha }{2}+ \frac{\tau }{2} \right)_{n}\left( \frac{\tau }{2} +1\right)_{i_{\tau -1}} \left( \frac{\delta }{2}+ \frac{1}{2}+\frac{\tau }{2} \right)_{i_{\tau -1}}}{\left( -\frac{j}{2}+\frac{\tau }{2}\right)_{i_{\tau -1}} \left( \frac{\alpha }{2}+\frac{\tau }{2} \right)_{i_{\tau -1}}\left( \frac{\tau }{2}+1 \right)_{n} \left( \frac{\delta }{2}+ \frac{1}{2}+\frac{\tau }{2} \right)_{n}} \left(-\frac{a}{a-1} \right)^{n } \nonumber
\end{eqnarray}
\begin{eqnarray}
y_0^m(j,\tilde{q};\xi) &=& \sum_{i_0=0}^{m} \frac{\left( -\frac{j}{2}\right)_{i_0} \left( \frac{\alpha }{2} \right)_{i_0}}{\left( 1 \right)_{i_0} \left( \frac{\delta }{2}+ \frac{1}{2} \right)_{i_0}} z^{i_0} \nonumber\\
y_1^m(j,\tilde{q};\xi) &=& \left\{\sum_{i_0=0}^{m} \frac{ i_0 \left( i_0 +\Gamma _0^{(F)} \right) +\frac{a\tilde{q}}{4(2a-1)}}{\left( i_0+\frac{1}{2} \right) \left( i_0+\frac{\delta }{2} \right)} \frac{\left( -\frac{j}{2}\right)_{i_0} \left( \frac{\alpha }{2} \right)_{i_0}}{\left( 1 \right)_{i_0} \left( \frac{\delta }{2}+ \frac{1}{2} \right)_{i_0}} \right. \nonumber\\
&&\times \left. \sum_{i_1 = i_0}^{m} \frac{\left( -\frac{j}{2}+\frac{1}{2} \right)_{i_1} \left( \frac{\alpha }{2}+ \frac{1}{2} \right)_{i_1}\left( \frac{3}{2} \right)_{i_0} \left( \frac{\delta }{2}+1 \right)_{i_0}}{\left(-\frac{j}{2}+\frac{1}{2} \right)_{i_0} \left( \frac{\alpha }{2}+ \frac{1}{2} \right)_{i_0}\left( \frac{3}{2} \right)_{i_1} \left( \frac{\delta }{2} +1 \right)_{i_1}} z^{i_1}\right\} \eta 
\nonumber
\end{eqnarray}
\begin{eqnarray}
y_{\tau }^m(j,\tilde{q};\xi) &=& \left\{ \sum_{i_0=0}^{m} \frac{ i_0 \left(i_0 +\Gamma _0^{(F)} \right)+\frac{a\tilde{q}}{4(2a-1)}}{\left( i_0+\frac{1}{2} \right) \left( i_0+\frac{\delta }{2} \right)} \frac{\left( -\frac{j}{2}\right)_{i_0} \left( \frac{\alpha }{2} \right)_{i_0}}{\left( 1 \right)_{i_0} \left( \frac{\delta }{2}+ \frac{1}{2} \right)_{i_0}} \right.\nonumber\\
&&\times \prod_{k=1}^{\tau -1} \left( \sum_{i_k = i_{k-1}}^{m} \frac{\left( i_k+ \frac{k}{2} \right)\left(i_k +\Gamma _k^{(F)} \right) +\frac{a\tilde{q}}{4(2a-1)}}{\left( i_k+\frac{k}{2}+\frac{1}{2} \right) \left( i_k+\frac{k}{2}+\frac{\delta }{2} \right)} \right. \nonumber\\
&&\times \left. \frac{\left( -\frac{j}{2}+\frac{k}{2}\right)_{i_k} \left( \frac{\alpha }{2}+ \frac{k}{2} \right)_{i_k}\left( \frac{k}{2}+1 \right)_{i_{k-1}} \left( \frac{\delta }{2}+ \frac{1}{2}+ \frac{k}{2} \right)_{i_{k-1}}}{\left( -\frac{j}{2}+\frac{k}{2}\right)_{i_{k-1}} \left( \frac{\alpha }{2}+ \frac{k}{2} \right)_{i_{k-1}}\left( \frac{k}{2} +1\right)_{i_k} \left( \frac{\delta }{2}+ \frac{1}{2}+ \frac{k}{2} \right)_{i_k}} \right) \nonumber\\
&&\times \left. \sum_{i_{\tau } = i_{\tau -1}}^{m} \frac{\left( -\frac{j}{2}+\frac{\tau }{2}\right)_{i_{\tau }} \left( \frac{\alpha }{2}+ \frac{\tau }{2} \right)_{i_{\tau }}\left( \frac{\tau }{2}+1 \right)_{i_{\tau -1}} \left( \frac{\delta }{2}+ \frac{1}{2}+\frac{\tau }{2} \right)_{i_{\tau -1}}}{\left( -\frac{j}{2}+\frac{\tau }{2}\right)_{i_{\tau -1}} \left( \frac{\alpha }{2}+ \frac{\tau }{2} \right)_{i_{\tau -1}}\left( \frac{\tau }{2}+1 \right)_{i_{\tau }} \left( \frac{\delta }{2}+ \frac{1}{2}+\frac{\tau }{2} \right)_{i_{\tau }}} z^{i_{\tau }}\right\} \eta ^{\tau }  \nonumber 
\end{eqnarray}
where  
\begin{equation}
\begin{cases} 
\tau \geq 2 \cr
\xi= \frac{x-1}{x} \cr
z = \frac{-a}{a-1}\xi^2 \cr
\eta = \frac{2a-1}{a-1} \xi \cr
\tilde{q} =  \frac{-q+\alpha (\delta a+\beta -\delta )}{a} \cr
\Gamma _0^{(F)} = \frac{a}{2(2a-1)}\left( \beta -1 -j +\frac{a-1}{a}\left( \alpha -\beta +\delta \right)\right) \cr
\Gamma _k^{(F)} = \frac{a}{2(2a-1)}\left( \beta -1 -j+k +\frac{a-1}{a}\left( \alpha -\beta +\delta +k \right)\right) 
\end{cases}\nonumber 
\end{equation} 
\end{enumerate}
\subsection{ \footnotesize ${\displaystyle \left(\frac{x-a}{1-a} \right)^{-\alpha } Hl\left(a, q-(\beta -\delta )\alpha ; \alpha , -\beta+\gamma +\delta , \delta , \gamma; \frac{a(x-1)}{x-a} \right)}$\normalsize}
\subsubsection{The first species complete polynomial}
\paragraph{The case of $\beta $, $q$ = fixed values and $\alpha , \gamma, \delta $ = free variables }
Replacing coefficients $q$, $\alpha $, $\beta $, $\gamma $, $\delta $ and $x$ by $q-(\beta -\delta )\alpha $, $-\beta+\gamma +\delta $, $\alpha $, $\delta $,  $\gamma $ and $\frac{a(x-1)}{x-a}$ into (\ref{eq:40018})--(\ref{eq:40022c}). Multiply $\left(\frac{x-a}{1-a} \right)^{-\alpha }$ and the new (\ref{eq:40018}), (\ref{eq:40019b}) and (\ref{eq:40020b}) together. 
\begin{enumerate} 
\item As $ \beta =\gamma +\delta $ and $q =\gamma \alpha +q_0^0 $ where $q_0^0=0$,

The eigenfunction is given by
\begin{eqnarray}
&& \left(\frac{x-a}{1-a} \right)^{-\alpha } y(\xi ) \nonumber\\
 &=& \left(\frac{x-a}{1-a} \right)^{-\alpha } Hl\left( a, 0; 0, \alpha , \delta , \gamma; \frac{a(x-1)}{x-a} \right) \nonumber\\
&=& \left(\frac{x-a}{1-a} \right)^{-\alpha } \nonumber
\end{eqnarray}
\item As $ \beta =\gamma +\delta +2N+1$ where $N \in \mathbb{N}_{0}$,

An algebraic equation of degree $2N+2$ for the determination of $q$ is given by
\begin{equation}
0 = \sum_{r=0}^{N+1}\bar{c}\left( 2r, N+1-r; 2N+1, q-(\gamma +2N+1 )\alpha \right) \nonumber
\end{equation}
The eigenvalue of $q$ is written by $(\gamma +2N+1 )\alpha + q_{2N+1}^m$ where $m = 0,1,2,\cdots,2N+1 $; $q_{2N+1}^0 < q_{2N+1}^1 < \cdots < q_{2N+1}^{2N+1}$. Its eigenfunction is given by
\begin{eqnarray} 
&& \left(\frac{x-a}{1-a} \right)^{-\alpha } y(\xi ) \nonumber\\
 &=& \left(\frac{x-a}{1-a} \right)^{-\alpha } Hl\left( a, q_{2N+1}^m; -2N-1, \alpha , \delta , \gamma; \frac{a(x-1)}{x-a} \right) \nonumber\\
&=& \left(\frac{x-a}{1-a} \right)^{-\alpha } \left\{ \sum_{r=0}^{N} y_{2r}^{N-r}\left( 2N+1,q_{2N+1}^m;\xi\right)+ \sum_{r=0}^{N} y_{2r+1}^{N-r}\left( 2N+1,q_{2N+1}^m;\xi\right) \right\} 
\nonumber
\end{eqnarray}
\item As $ \beta =\gamma +\delta +2N+2$ where $N \in \mathbb{N}_{0}$,

An algebraic equation of degree $2N+3$ for the determination of $q$ is given by
\begin{eqnarray}
0  = \sum_{r=0}^{N+1}\bar{c}\left( 2r+1, N+1-r; 2N+2, q-(\gamma +2N+2 )\alpha \right) \nonumber
\end{eqnarray}
The eigenvalue of $q$ is written by $(\gamma +2N+2 )\alpha +q_{2N+2}^m$ where $m = 0,1,2,\cdots,2N+2 $; $q_{2N+2}^0 < q_{2N+2}^1 < \cdots < q_{2N+2}^{2N+2}$. Its eigenfunction is given by
\begin{eqnarray} 
&& \left(\frac{x-a}{1-a} \right)^{-\alpha } y(\xi ) \nonumber\\
 &=& \left(\frac{x-a}{1-a} \right)^{-\alpha } Hl\left( a, q_{2N+2}^m ; -2N-2, \alpha , \delta , \gamma; \frac{a(x-1)}{x-a} \right) \nonumber\\
&=& \left(\frac{x-a}{1-a} \right)^{-\alpha } \left\{ \sum_{r=0}^{N+1} y_{2r}^{N+1-r}\left( 2N+2,q_{2N+2}^m;\xi\right) + \sum_{r=0}^{N} y_{2r+1}^{N-r}\left( 2N+2,q_{2N+2}^m;\xi\right) \right\} 
\nonumber
\end{eqnarray}
In the above,
\begin{eqnarray}
\bar{c}(0,n;j,\tilde{q})  &=& \frac{\left( -\frac{j}{2}\right)_{n} \left( \frac{\alpha }{2} \right)_{n}}{\left( 1 \right)_{n} \left( \frac{\delta }{2}+ \frac{1}{2} \right)_{n}} \left(-\frac{1}{a} \right)^{n}\nonumber\\
\bar{c}(1,n;j,\tilde{q}) &=& \left( \frac{1+a}{a}\right) \sum_{i_0=0}^{n}\frac{ i_0 \left(i_0 +\Gamma _0^{(F)} \right)+\frac{\tilde{q}}{4(1+a)}}{\left( i_0+\frac{1}{2} \right) \left( i_0+\frac{\delta }{2} \right)} \frac{\left( -\frac{j}{2}\right)_{i_0} \left( \frac{\alpha }{2} \right)_{i_0}}{\left( 1 \right)_{i_0} \left( \frac{\delta }{2}+ \frac{1}{2} \right)_{i_0}}\nonumber\\
&&\times  \frac{\left( -\frac{j}{2}+\frac{1}{2} \right)_{n} \left( \frac{\alpha }{2}+ \frac{1}{2} \right)_{n}\left( \frac{3}{2} \right)_{i_0} \left( \frac{\delta }{2}+1 \right)_{i_0}}{\left( -\frac{j}{2}+\frac{1}{2}\right)_{i_0} \left( \frac{\alpha }{2}+ \frac{1}{2} \right)_{i_0}\left( \frac{3}{2} \right)_{n} \left( \frac{\delta }{2}+1 \right)_{n}} \left(-\frac{1}{a} \right)^{n}  
\nonumber\\
\bar{c}(\tau ,n;j,\tilde{q}) &=& \left( \frac{1+a}{a}\right)^{\tau} \sum_{i_0=0}^{n}\frac{ i_0 \left(i_0 +\Gamma _0^{(F)} \right)+\frac{\tilde{q}}{4(1+a)}}{\left( i_0+\frac{1}{2} \right) \left( i_0+\frac{\delta }{2} \right)} \frac{\left( -\frac{j}{2}\right)_{i_0} \left( \frac{\alpha }{2} \right)_{i_0}}{\left( 1 \right)_{i_0} \left( \frac{\delta }{2}+ \frac{1}{2} \right)_{i_0}}  \nonumber\\
&&\times \prod_{k=1}^{\tau -1} \left( \sum_{i_k = i_{k-1}}^{n} \frac{\left( i_k+ \frac{k}{2} \right)\left(i_k +\Gamma _k^{(F)} \right)+\frac{\tilde{q}}{4(1+a)}}{\left( i_k+\frac{k}{2}+\frac{1}{2} \right) \left( i_k+\frac{k}{2}+\frac{\delta }{2} \right)} \right. \nonumber\\
&&\times \left. \frac{\left( -\frac{j}{2}+\frac{k}{2}\right)_{i_k} \left( \frac{\alpha }{2}+ \frac{k}{2} \right)_{i_k}\left( \frac{k}{2}+1 \right)_{i_{k-1}} \left( \frac{\delta }{2}+ \frac{1}{2}+ \frac{k}{2} \right)_{i_{k-1}}}{\left( -\frac{j}{2}+\frac{k}{2}\right)_{i_{k-1}} \left( \frac{\alpha }{2}+ \frac{k}{2} \right)_{i_{k-1}}\left( \frac{k}{2} +1\right)_{i_k} \left( \frac{\delta }{2}+ \frac{1}{2}+ \frac{k}{2} \right)_{i_k}} \right) \nonumber\\
&&\times \frac{\left( -\frac{j}{2}+\frac{\tau }{2}\right)_{n} \left( \frac{\alpha }{2}+ \frac{\tau }{2} \right)_{n}\left( \frac{\tau }{2} +1\right)_{i_{\tau -1}} \left( \frac{\delta }{2}+ \frac{1}{2}+\frac{\tau }{2} \right)_{i_{\tau -1}}}{\left( -\frac{j}{2}+\frac{\tau }{2}\right)_{i_{\tau -1}} \left( \frac{\alpha }{2}+\frac{\tau }{2} \right)_{i_{\tau -1}}\left( \frac{\tau }{2}+1 \right)_{n} \left( \frac{\delta }{2}+ \frac{1}{2}+\frac{\tau }{2} \right)_{n}} \left(-\frac{1}{a} \right)^{n } \nonumber 
\end{eqnarray}
\begin{eqnarray}
y_0^m(j,\tilde{q};\xi) &=& \sum_{i_0=0}^{m} \frac{\left( -\frac{j}{2}\right)_{i_0} \left( \frac{\alpha }{2} \right)_{i_0}}{\left( 1 \right)_{i_0} \left( \frac{\delta }{2}+ \frac{1}{2} \right)_{i_0}} z^{i_0} \nonumber\\
y_1^m(j,\tilde{q};\xi) &=& \left\{\sum_{i_0=0}^{m} \frac{ i_0 \left(i_0 +\Gamma _0^{(F)} \right)+\frac{\tilde{q}}{4(1+a)}}{\left( i_0+\frac{1}{2} \right) \left( i_0+\frac{\delta }{2} \right)} \frac{\left( -\frac{j}{2}\right)_{i_0} \left( \frac{\alpha }{2} \right)_{i_0}}{\left( 1 \right)_{i_0} \left( \frac{\delta }{2}+ \frac{1}{2} \right)_{i_0}} \right. \nonumber\\
&&\times \left. \sum_{i_1 = i_0}^{m} \frac{\left( -\frac{j}{2}+\frac{1}{2} \right)_{i_1} \left( \frac{\alpha }{2}+ \frac{1}{2} \right)_{i_1}\left( \frac{3}{2} \right)_{i_0} \left( \frac{\delta }{2}+1 \right)_{i_0}}{\left(-\frac{j}{2}+\frac{1}{2} \right)_{i_0} \left( \frac{\alpha }{2}+ \frac{1}{2} \right)_{i_0}\left( \frac{3}{2} \right)_{i_1} \left( \frac{\delta }{2} +1 \right)_{i_1}} z^{i_1}\right\} \eta 
\nonumber
\end{eqnarray}
\begin{eqnarray}
y_{\tau }^m(j,\tilde{q};\xi) &=& \left\{ \sum_{i_0=0}^{m} \frac{ i_0 \left(i_0 +\Gamma _0^{(F)} \right)+\frac{\tilde{q}}{4(1+a)}}{\left( i_0+\frac{1}{2} \right) \left( i_0+\frac{\delta }{2} \right)} \frac{\left( -\frac{j}{2}\right)_{i_0} \left( \frac{\alpha }{2} \right)_{i_0}}{\left( 1 \right)_{i_0} \left( \frac{\delta }{2}+ \frac{1}{2} \right)_{i_0}} \right.\nonumber\\
&&\times \prod_{k=1}^{\tau -1} \left( \sum_{i_k = i_{k-1}}^{m} \frac{\left( i_k+ \frac{k}{2} \right)\left(i_k +\Gamma _k^{(F)} \right)+\frac{\tilde{q}}{4(1+a)}}{\left( i_k+\frac{k}{2}+\frac{1}{2} \right) \left( i_k+\frac{k}{2}+\frac{\delta }{2} \right)} \right. \nonumber\\
&&\times \left. \frac{\left( -\frac{j}{2}+\frac{k}{2}\right)_{i_k} \left( \frac{\alpha }{2}+ \frac{k}{2} \right)_{i_k}\left( \frac{k}{2}+1 \right)_{i_{k-1}} \left( \frac{\delta }{2}+ \frac{1}{2}+ \frac{k}{2} \right)_{i_{k-1}}}{\left( -\frac{j}{2}+\frac{k}{2}\right)_{i_{k-1}} \left( \frac{\alpha }{2}+ \frac{k}{2} \right)_{i_{k-1}}\left( \frac{k}{2} +1\right)_{i_k} \left( \frac{\delta }{2}+ \frac{1}{2}+ \frac{k}{2} \right)_{i_k}} \right) \nonumber\\
&&\times \left. \sum_{i_{\tau } = i_{\tau -1}}^{m} \frac{\left( -\frac{j}{2}+\frac{\tau }{2}\right)_{i_{\tau }} \left( \frac{\alpha }{2}+ \frac{\tau }{2} \right)_{i_{\tau }}\left( \frac{\tau }{2}+1 \right)_{i_{\tau -1}} \left( \frac{\delta }{2}+ \frac{1}{2}+\frac{\tau }{2} \right)_{i_{\tau -1}}}{\left( -\frac{j}{2}+\frac{\tau }{2}\right)_{i_{\tau -1}} \left( \frac{\alpha }{2}+ \frac{\tau }{2} \right)_{i_{\tau -1}}\left( \frac{\tau }{2}+1 \right)_{i_{\tau }} \left( \frac{\delta }{2}+ \frac{1}{2}+\frac{\tau }{2} \right)_{i_{\tau }}} z^{i_{\tau }}\right\} \eta ^{\tau } \nonumber 
\end{eqnarray}
where  
\begin{equation}
\begin{cases}
\tau \geq 2  \cr 
\xi= \frac{a(x-1)}{x-a} \cr
z = -\frac{1}{a}\xi^2 \cr
\eta = \frac{(1+a)}{a} \xi \cr
\tilde{q} =  q-(\beta -\delta )\alpha \cr
\Gamma _0^{(F)} = \frac{1}{2(1+a)}\left( \alpha -\gamma -j +a\left( \gamma +\delta -1 \right)\right) \cr
\Gamma _k^{(F)} =  \frac{1}{2(1+a)}\left( \alpha -\gamma -j+k +a\left( \gamma +\delta -1+k \right)\right)
\end{cases}\nonumber 
\end{equation}  
\end{enumerate}
\paragraph{The case of $\gamma $, $q$ = fixed values and $\alpha , \beta, \delta $ = free variables }
Replacing coefficients $q$, $\alpha $, $\beta $, $\gamma $, $\delta $ and $x$ by $q-(\beta -\delta )\alpha $, $-\beta+\gamma +\delta $, $\alpha $, $\delta $,  $\gamma $ and $\frac{a(x-1)}{x-a}$ into (\ref{eq:40018})--(\ref{eq:40022c}). Replacing $\gamma $ by $\beta -\delta -j$ into the new (\ref{eq:40021a})--(\ref{eq:40022c}). Multiply $\left(\frac{x-a}{1-a} \right)^{-\alpha }$ and the new (\ref{eq:40018}), (\ref{eq:40019b}) and (\ref{eq:40020b}) together.  
\begin{enumerate} 
\item As $ \gamma =\beta -\delta $ and $q=(\beta -\delta )\alpha+ q_0^0$ where $q_0^0=0$,

The eigenfunction is given by
\begin{eqnarray}
&& \left(\frac{x-a}{1-a} \right)^{-\alpha } y(\xi ) \nonumber\\
 &=& \left(\frac{x-a}{1-a} \right)^{-\alpha } Hl\left(a, 0 ; 0, \alpha , \delta , \gamma; \frac{a(x-1)}{x-a} \right) \nonumber\\
&=& \left(\frac{x-a}{1-a} \right)^{-\alpha } \nonumber
\end{eqnarray}
\item As $ \gamma  =\beta -\delta -2N-1$ where $N \in \mathbb{N}_{0}$,

An algebraic equation of degree $2N+2$ for the determination of $q$ is given by
\begin{equation}
0 = \sum_{r=0}^{N+1}\bar{c}\left( 2r, N+1-r; 2N+1,\tilde{q}\right)\nonumber
\end{equation}
The eigenvalue of $q$ is written by $(\beta -\delta )\alpha + q_{2N+1}^m$ where $m = 0,1,2,\cdots,2N+1 $; $q_{2N+1}^0 < q_{2N+1}^1 < \cdots < q_{2N+1}^{2N+1}$. Its eigenfunction is given by
\begin{eqnarray} 
&& \left(\frac{x-a}{1-a} \right)^{-\alpha } y(\xi ) \nonumber\\
 &=& \left(\frac{x-a}{1-a} \right)^{-\alpha } Hl\left( a, q_{2N+1}^m; -2N-1, \alpha , \delta , \gamma; \frac{a(x-1)}{x-a} \right) \nonumber\\
&=& \left(\frac{x-a}{1-a} \right)^{-\alpha } \left\{ \sum_{r=0}^{N} y_{2r}^{N-r}\left( 2N+1,q_{2N+1}^m;\xi\right)+ \sum_{r=0}^{N} y_{2r+1}^{N-r}\left( 2N+1,q_{2N+1}^m;\xi\right) \right\} 
\nonumber
\end{eqnarray}
\item As $ \gamma =\beta -\delta -2N-2$ where $N \in \mathbb{N}_{0}$,

An algebraic equation of degree $2N+3$ for the determination of $q$ is given by
\begin{eqnarray}
0  = \sum_{r=0}^{N+1}\bar{c}\left( 2r+1, N+1-r; 2N+2,\tilde{q}\right)\nonumber
\end{eqnarray}
The eigenvalue of $q$ is written by $(\beta -\delta )\alpha + q_{2N+2}^m$ where $m = 0,1,2,\cdots,2N+2 $; $q_{2N+2}^0 < q_{2N+2}^1 < \cdots < q_{2N+2}^{2N+2}$. Its eigenfunction is given by
\begin{eqnarray} 
&& \left(\frac{x-a}{1-a} \right)^{-\alpha } y(\xi ) \nonumber\\
 &=& \left(\frac{x-a}{1-a} \right)^{-\alpha } Hl\left(a, q_{2N+2}^m ; -2N-2, \alpha , \delta , \gamma; \frac{a(x-1)}{x-a} \right) \nonumber\\
&=& \left(\frac{x-a}{1-a} \right)^{-\alpha } \left\{ \sum_{r=0}^{N+1} y_{2r}^{N+1-r}\left( 2N+2,q_{2N+2}^m;\xi\right) + \sum_{r=0}^{N} y_{2r+1}^{N-r}\left( 2N+2,q_{2N+2}^m;\xi\right) \right\}
\nonumber
\end{eqnarray}
In the above,
\begin{eqnarray}
\bar{c}(0,n;j,\tilde{q})  &=& \frac{\left( -\frac{j}{2}\right)_{n} \left( \frac{\alpha }{2} \right)_{n}}{\left( 1 \right)_{n} \left( \frac{\delta }{2}+ \frac{1}{2} \right)_{n}} \left(-\frac{1}{a} \right)^{n}\nonumber\\
\bar{c}(1,n;j,\tilde{q}) &=& \left( \frac{1+a}{a}\right) \sum_{i_0=0}^{n}\frac{ i_0 \left(i_0 +\Gamma _0^{(F)} \right)+\frac{\tilde{q}}{4(1+a)}}{\left( i_0+\frac{1}{2} \right) \left( i_0+\frac{\delta }{2} \right)} \frac{\left( -\frac{j}{2}\right)_{i_0} \left( \frac{\alpha }{2} \right)_{i_0}}{\left( 1 \right)_{i_0} \left( \frac{\delta }{2}+ \frac{1}{2} \right)_{i_0}} \nonumber\\
&&\times  \frac{\left( -\frac{j}{2}+\frac{1}{2} \right)_{n} \left( \frac{\alpha }{2}+ \frac{1}{2} \right)_{n}\left( \frac{3}{2} \right)_{i_0} \left( \frac{\delta }{2}+1 \right)_{i_0}}{\left( -\frac{j}{2}+\frac{1}{2}\right)_{i_0} \left( \frac{\alpha }{2}+ \frac{1}{2} \right)_{i_0}\left( \frac{3}{2} \right)_{n} \left( \frac{\delta }{2}+1 \right)_{n}} \left(-\frac{1}{a} \right)^{n}  
\nonumber\\
\bar{c}(\tau ,n;j,\tilde{q}) &=& \left( \frac{1+a}{a}\right)^{\tau} \sum_{i_0=0}^{n}\frac{ i_0 \left(i_0 +\Gamma _0^{(F)} \right)+\frac{\tilde{q}}{4(1+a)}}{\left( i_0+\frac{1}{2} \right) \left( i_0+\frac{\delta }{2} \right)} \frac{\left( -\frac{j}{2}\right)_{i_0} \left( \frac{\alpha }{2} \right)_{i_0}}{\left( 1 \right)_{i_0} \left( \frac{\delta }{2}+ \frac{1}{2} \right)_{i_0}} \nonumber\\
&&\times \prod_{k=1}^{\tau -1} \left( \sum_{i_k = i_{k-1}}^{n} \frac{\left( i_k+ \frac{k}{2} \right)\left(i_k +\Gamma _k^{(F)} \right)+\frac{\tilde{q}}{4(1+a)}}{\left( i_k+\frac{k}{2}+\frac{1}{2} \right) \left( i_k+\frac{k}{2}+\frac{\delta }{2} \right)} \right. \nonumber\\
&&\times \left. \frac{\left( -\frac{j}{2}+\frac{k}{2}\right)_{i_k} \left( \frac{\alpha }{2}+ \frac{k}{2} \right)_{i_k}\left( \frac{k}{2}+1 \right)_{i_{k-1}} \left( \frac{\delta }{2}+ \frac{1}{2}+ \frac{k}{2} \right)_{i_{k-1}}}{\left( -\frac{j}{2}+\frac{k}{2}\right)_{i_{k-1}} \left( \frac{\alpha }{2}+ \frac{k}{2} \right)_{i_{k-1}}\left( \frac{k}{2} +1\right)_{i_k} \left( \frac{\delta }{2}+ \frac{1}{2}+ \frac{k}{2} \right)_{i_k}} \right) \nonumber\\
&&\times \frac{\left( -\frac{j}{2}+\frac{\tau }{2}\right)_{n} \left( \frac{\alpha }{2}+ \frac{\tau }{2} \right)_{n}\left( \frac{\tau }{2} +1\right)_{i_{\tau -1}} \left( \frac{\delta }{2}+ \frac{1}{2}+\frac{\tau }{2} \right)_{i_{\tau -1}}}{\left( -\frac{j}{2}+\frac{\tau }{2}\right)_{i_{\tau -1}} \left( \frac{\alpha }{2}+\frac{\tau }{2} \right)_{i_{\tau -1}}\left( \frac{\tau }{2}+1 \right)_{n} \left( \frac{\delta }{2}+ \frac{1}{2}+\frac{\tau }{2} \right)_{n}} \left(-\frac{1}{a} \right)^{n } \nonumber 
\end{eqnarray}
\begin{eqnarray}
y_0^m(j,\tilde{q};\xi) &=& \sum_{i_0=0}^{m} \frac{\left( -\frac{j}{2}\right)_{i_0} \left( \frac{\alpha }{2} \right)_{i_0}}{\left( 1 \right)_{i_0} \left( \frac{\delta }{2}+ \frac{1}{2} \right)_{i_0}} z^{i_0} \nonumber\\
y_1^m(j,\tilde{q};\xi) &=& \left\{\sum_{i_0=0}^{m} \frac{ i_0 \left( i_0 +\Gamma _0^{(F)} \right)+\frac{\tilde{q}}{4(1+a)}}{\left( i_0+\frac{1}{2} \right) \left( i_0+\frac{\delta }{2} \right)} \frac{\left( -\frac{j}{2}\right)_{i_0} \left( \frac{\alpha }{2} \right)_{i_0}}{\left( 1 \right)_{i_0} \left( \frac{\delta }{2}+ \frac{1}{2} \right)_{i_0}} \right. \nonumber\\
&&\times \left. \sum_{i_1 = i_0}^{m} \frac{\left( -\frac{j}{2}+\frac{1}{2} \right)_{i_1} \left( \frac{\alpha }{2}+ \frac{1}{2} \right)_{i_1}\left( \frac{3}{2} \right)_{i_0} \left( \frac{\delta }{2}+1 \right)_{i_0}}{\left(-\frac{j}{2}+\frac{1}{2} \right)_{i_0} \left( \frac{\alpha }{2}+ \frac{1}{2} \right)_{i_0}\left( \frac{3}{2} \right)_{i_1} \left( \frac{\delta }{2} +1 \right)_{i_1}} z^{i_1}\right\} \eta 
\nonumber
\end{eqnarray}
\begin{eqnarray}
y_{\tau }^m(j,\tilde{q};\xi) &=& \left\{ \sum_{i_0=0}^{m} \frac{ i_0 \left(i_0 +\Gamma _0^{(F)} \right)+\frac{\tilde{q}}{4(1+a)}}{\left( i_0+\frac{1}{2} \right) \left( i_0+\frac{\delta }{2} \right)} \frac{\left( -\frac{j}{2}\right)_{i_0} \left( \frac{\alpha }{2} \right)_{i_0}}{\left( 1 \right)_{i_0} \left( \frac{\delta }{2}+ \frac{1}{2} \right)_{i_0}} \right.\nonumber\\
&&\times \prod_{k=1}^{\tau -1} \left( \sum_{i_k = i_{k-1}}^{m} \frac{\left( i_k+ \frac{k}{2} \right)\left(i_k +\Gamma _k^{(F)} \right)+\frac{\tilde{q}}{4(1+a)}}{\left( i_k+\frac{k}{2}+\frac{1}{2} \right) \left( i_k+\frac{k}{2}+\frac{\delta }{2} \right)} \right. \nonumber\\
&&\times \left. \frac{\left( -\frac{j}{2}+\frac{k}{2}\right)_{i_k} \left( \frac{\alpha }{2}+ \frac{k}{2} \right)_{i_k}\left( \frac{k}{2}+1 \right)_{i_{k-1}} \left( \frac{\delta }{2}+ \frac{1}{2}+ \frac{k}{2} \right)_{i_{k-1}}}{\left( -\frac{j}{2}+\frac{k}{2}\right)_{i_{k-1}} \left( \frac{\alpha }{2}+ \frac{k}{2} \right)_{i_{k-1}}\left( \frac{k}{2} +1\right)_{i_k} \left( \frac{\delta }{2}+ \frac{1}{2}+ \frac{k}{2} \right)_{i_k}} \right) \nonumber\\
&&\times \left. \sum_{i_{\tau } = i_{\tau -1}}^{m} \frac{\left( -\frac{j}{2}+\frac{\tau }{2}\right)_{i_{\tau }} \left( \frac{\alpha }{2}+ \frac{\tau }{2} \right)_{i_{\tau }}\left( \frac{\tau }{2}+1 \right)_{i_{\tau -1}} \left( \frac{\delta }{2}+ \frac{1}{2}+\frac{\tau }{2} \right)_{i_{\tau -1}}}{\left( -\frac{j}{2}+\frac{\tau }{2}\right)_{i_{\tau -1}} \left( \frac{\alpha }{2}+ \frac{\tau }{2} \right)_{i_{\tau -1}}\left( \frac{\tau }{2}+1 \right)_{i_{\tau }} \left( \frac{\delta }{2}+ \frac{1}{2}+\frac{\tau }{2} \right)_{i_{\tau }}} z^{i_{\tau }}\right\} \eta ^{\tau } \nonumber 
\end{eqnarray}
where  
\begin{equation}
\begin{cases}
\tau \geq 2  \cr 
\xi= \frac{a(x-1)}{x-a} \cr
z = -\frac{1}{a}\xi^2 \cr
\eta = \frac{(1+a)}{a} \xi \cr
\tilde{q} =  q-(\beta -\delta )\alpha \cr
\Gamma _0^{(F)} = \frac{1}{2(1+a)}\left( \alpha -\beta +\delta +a\left( \beta -1-j \right)\right) \cr
\Gamma _k^{(F)} = \frac{1}{2(1+a)}\left( \alpha -\beta +\delta +k +a\left( \beta -1-j+k \right)\right)   
\end{cases}\nonumber 
\end{equation} 
\end{enumerate}  
\paragraph{The case of $\delta $, $q$ = fixed values and $\alpha , \beta, \gamma $ = free variables } 
Replacing coefficients $q$, $\alpha $, $\beta $, $\gamma $, $\delta $ and $x$ by $q-(\beta -\delta )\alpha $, $-\beta+\gamma +\delta $, $\alpha $, $\delta $,  $\gamma $ and $\frac{a(x-1)}{x-a}$ into (\ref{eq:40018})--(\ref{eq:40022c}). Replacing $\delta $ by $\beta -\gamma -j$ into the new (\ref{eq:40021a})--(\ref{eq:40022c}). Multiply $\left(\frac{x-a}{1-a} \right)^{-\alpha }$ and the new (\ref{eq:40018}), (\ref{eq:40019b}) and (\ref{eq:40020b}) together.
\begin{enumerate} 
\item As $ \delta =\beta -\gamma $ and $q = \gamma \alpha +q_0^0 $ where $q_0^0=0$,

The eigenfunction is given by
\begin{eqnarray}
&& \left(\frac{x-a}{1-a} \right)^{-\alpha } y(\xi ) \nonumber\\
 &=& \left(\frac{x-a}{1-a} \right)^{-\alpha } Hl\left( a, 0; 0, \alpha , \delta , \gamma; \frac{a(x-1)}{x-a} \right) \nonumber\\
&=& \left(\frac{x-a}{1-a} \right)^{-\alpha } \nonumber
\end{eqnarray}
\item As $ \delta  = \beta -\gamma -2N-1$ where $N \in \mathbb{N}_{0}$,

An algebraic equation of degree $2N+2$ for the determination of $q$ is given by
\begin{equation}
0 = \sum_{r=0}^{N+1}\bar{c}\left( 2r, N+1-r; 2N+1,q-(\gamma +2N+1)\alpha \right)\nonumber
\end{equation}
The eigenvalue of $q$ is written by $(\gamma +2N+1)\alpha + q_{2N+1}^m$ where $m = 0,1,2,\cdots,2N+1 $; $q_{2N+1}^0 < q_{2N+1}^1 < \cdots < q_{2N+1}^{2N+1}$. Its eigenfunction is given by
\begin{eqnarray} 
&& \left(\frac{x-a}{1-a} \right)^{-\alpha } y(\xi ) \nonumber\\
 &=& \left(\frac{x-a}{1-a} \right)^{-\alpha } Hl\left( a, q_{2N+1}^m; -2N-1, \alpha , \delta , \gamma; \frac{a(x-1)}{x-a} \right) \nonumber\\
&=& \left(\frac{x-a}{1-a} \right)^{-\alpha } \left\{ \sum_{r=0}^{N} y_{2r}^{N-r}\left( 2N+1,q_{2N+1}^m;\xi\right)+ \sum_{r=0}^{N} y_{2r+1}^{N-r}\left( 2N+1,q_{2N+1}^m;\xi\right) \right\} 
\nonumber
\end{eqnarray}
\item As $ \delta  =\beta -\gamma -2N-2$ where $N \in \mathbb{N}_{0}$,

An algebraic equation of degree $2N+3$ for the determination of $q$ is given by
\begin{eqnarray}
0  = \sum_{r=0}^{N+1}\bar{c}\left( 2r+1, N+1-r; 2N+2, q-(\gamma +2N+2)\alpha \right)\nonumber
\end{eqnarray}
The eigenvalue of $q$ is written by $(\gamma +2N+2)\alpha + q_{2N+2}^m$ where $m = 0,1,2,\cdots,2N+2 $; $q_{2N+2}^0 < q_{2N+2}^1 < \cdots < q_{2N+2}^{2N+2}$. Its eigenfunction is given by
\begin{eqnarray} 
&& \left(\frac{x-a}{1-a} \right)^{-\alpha } y(\xi ) \nonumber\\
 &=& \left(\frac{x-a}{1-a} \right)^{-\alpha } Hl\left( a, q_{2N+2}^m; -2N-2, \alpha , \delta , \gamma; \frac{a(x-1)}{x-a} \right) \nonumber\\
&=& \left(\frac{x-a}{1-a} \right)^{-\alpha } \left\{ \sum_{r=0}^{N+1} y_{2r}^{N+1-r}\left( 2N+2,q_{2N+2}^m;\xi\right) + \sum_{r=0}^{N} y_{2r+1}^{N-r}\left( 2N+2,q_{2N+2}^m;\xi\right) \right\}
\nonumber
\end{eqnarray}
In the above,
\begin{eqnarray}
\bar{c}(0,n;j,\tilde{q})  &=& \frac{\left( -\frac{j}{2}\right)_{n} \left( \frac{\alpha }{2} \right)_{n}}{\left( 1 \right)_{n} \left( \frac{\beta}{2} -\frac{\gamma}{2}+ \frac{1}{2} -\frac{j}{2} \right)_{n}} \left(-\frac{1}{a} \right)^{n}\nonumber\\
\bar{c}(1,n;j,\tilde{q}) &=& \left( \frac{1+a}{a}\right) \sum_{i_0=0}^{n}\frac{ i_0 \left( i_0 +\Gamma _0^{(F)} \right) +\frac{\tilde{q}}{4(1+a)}}{\left( i_0+\frac{1}{2} \right) \left( i_0+\frac{\beta}{2} -\frac{\gamma}{2} -\frac{j}{2} \right)} \frac{\left( -\frac{j}{2}\right)_{i_0} \left( \frac{\alpha }{2} \right)_{i_0}}{\left( 1 \right)_{i_0} \left( \frac{\beta}{2} -\frac{\gamma}{2}+ \frac{1}{2} -\frac{j}{2} \right)_{i_0}}\nonumber\\
&&\times  \frac{\left( -\frac{j}{2}+\frac{1}{2} \right)_{n} \left( \frac{\alpha }{2}+ \frac{1}{2} \right)_{n}\left( \frac{3}{2} \right)_{i_0} \left( \frac{\beta}{2} -\frac{\gamma}{2}+1-\frac{j}{2} \right)_{i_0}}{\left( -\frac{j}{2}+\frac{1}{2}\right)_{i_0} \left( \frac{\alpha }{2}+ \frac{1}{2} \right)_{i_0}\left( \frac{3}{2} \right)_{n} \left( \frac{\beta}{2} -\frac{\gamma}{2}+1-\frac{j}{2} \right)_{n}} \left(-\frac{1}{a} \right)^{n} \nonumber\\
\bar{c}(\tau ,n;j,\tilde{q}) &=& \left( \frac{1+a}{a}\right)^{\tau} \sum_{i_0=0}^{n}\frac{ i_0 \left( i_0 +\Gamma _0^{(F)} \right) +\frac{\tilde{q}}{4(1+a)}}{\left( i_0+\frac{1}{2} \right) \left( i_0+\frac{\beta}{2} -\frac{\gamma}{2} -\frac{j}{2} \right)} \frac{\left( -\frac{j}{2}\right)_{i_0} \left( \frac{\alpha }{2} \right)_{i_0}}{\left( 1 \right)_{i_0} \left( \frac{\beta}{2} -\frac{\gamma}{2}+ \frac{1}{2}-\frac{j}{2} \right)_{i_0}} \nonumber\\
&&\times \prod_{k=1}^{\tau -1} \left( \sum_{i_k = i_{k-1}}^{n} \frac{\left( i_k+ \frac{k}{2} \right)\left( i_k +\Gamma _k^{(F)} \right) +\frac{\tilde{q}}{4(1+a)}}{\left( i_k+\frac{k}{2}+\frac{1}{2} \right) \left( i_k+\frac{k}{2}+\frac{\beta}{2} -\frac{\gamma}{2}-\frac{j}{2} \right)} \right. \nonumber\\
&&\times \left. \frac{\left( -\frac{j}{2}+\frac{k}{2}\right)_{i_k} \left( \frac{\alpha }{2}+ \frac{k}{2} \right)_{i_k}\left( \frac{k}{2}+1 \right)_{i_{k-1}} \left( \frac{\beta}{2} -\frac{\gamma}{2}+ \frac{1}{2}+ \frac{k}{2}-\frac{j}{2} \right)_{i_{k-1}}}{\left( -\frac{j}{2}+\frac{k}{2}\right)_{i_{k-1}} \left( \frac{\alpha }{2}+ \frac{k}{2} \right)_{i_{k-1}}\left( \frac{k}{2} +1\right)_{i_k} \left( \frac{\beta}{2} -\frac{\gamma}{2}+ \frac{1}{2}+ \frac{k}{2}-\frac{j}{2} \right)_{i_k}} \right) \nonumber\\
&&\times \frac{\left( -\frac{j}{2}+\frac{\tau }{2}\right)_{n} \left( \frac{\alpha }{2}+ \frac{\tau }{2} \right)_{n}\left( \frac{\tau }{2} +1\right)_{i_{\tau -1}} \left( \frac{\beta}{2} -\frac{\gamma}{2}+ \frac{1}{2}+\frac{\tau }{2}-\frac{j}{2} \right)_{i_{\tau -1}}}{\left( -\frac{j}{2}+\frac{\tau }{2}\right)_{i_{\tau -1}} \left( \frac{\alpha }{2}+\frac{\tau }{2} \right)_{i_{\tau -1}}\left( \frac{\tau }{2}+1 \right)_{n} \left( \frac{\beta}{2} -\frac{\gamma}{2}+ \frac{1}{2}+\frac{\tau }{2} -\frac{j}{2}\right)_{n}} \left(-\frac{1}{a} \right)^{n } \nonumber 
\end{eqnarray}
\begin{eqnarray}
y_0^m(j,\tilde{q};\xi) &=& \sum_{i_0=0}^{m} \frac{\left( -\frac{j}{2}\right)_{i_0} \left( \frac{\alpha }{2} \right)_{i_0}}{\left( 1 \right)_{i_0} \left( \frac{\beta}{2} -\frac{\gamma}{2}+ \frac{1}{2}-\frac{j}{2} \right)_{i_0}} z^{i_0} \nonumber\\
y_1^m(j,\tilde{q};\xi) &=& \left\{\sum_{i_0=0}^{m} \frac{ i_0 \left(i_0 +\Gamma _0^{(F)} \right)+\frac{\tilde{q}}{4(1+a)}}{\left( i_0+\frac{1}{2} \right) \left( i_0+\frac{\beta}{2} -\frac{\gamma}{2}-\frac{j}{2} \right)} \frac{\left( -\frac{j}{2}\right)_{i_0} \left( \frac{\alpha }{2} \right)_{i_0}}{\left( 1 \right)_{i_0} \left( \frac{\beta}{2} -\frac{\gamma}{2}+ \frac{1}{2}-\frac{j}{2} \right)_{i_0}} \right. \nonumber\\
&&\times \left. \sum_{i_1 = i_0}^{m} \frac{\left( -\frac{j}{2}+\frac{1}{2} \right)_{i_1} \left( \frac{\alpha }{2}+ \frac{1}{2} \right)_{i_1}\left( \frac{3}{2} \right)_{i_0} \left( \frac{\beta}{2} -\frac{\gamma}{2}+1-\frac{j}{2} \right)_{i_0}}{\left(-\frac{j}{2}+\frac{1}{2} \right)_{i_0} \left( \frac{\alpha }{2}+ \frac{1}{2} \right)_{i_0}\left( \frac{3}{2} \right)_{i_1} \left( \frac{\beta}{2} -\frac{\gamma}{2}+1-\frac{j}{2} \right)_{i_1}} z^{i_1}\right\} \eta 
\nonumber
\end{eqnarray}
\begin{eqnarray}
y_{\tau }^m(j,\tilde{q};\xi) &=& \left\{ \sum_{i_0=0}^{m} \frac{ i_0 \left(i_0 +\Gamma _0^{(F)} \right)+\frac{\tilde{q}}{4(1+a)}}{\left( i_0+\frac{1}{2} \right) \left( i_0+\frac{\beta}{2} -\frac{\gamma}{2}-\frac{j}{2} \right)} \frac{\left( -\frac{j}{2}\right)_{i_0} \left( \frac{\alpha }{2} \right)_{i_0}}{\left( 1 \right)_{i_0} \left( \frac{\beta}{2} -\frac{\gamma}{2}+ \frac{1}{2}-\frac{j}{2} \right)_{i_0}} \right.\nonumber\\
&&\times \prod_{k=1}^{\tau -1} \left( \sum_{i_k = i_{k-1}}^{m} \frac{\left( i_k+ \frac{k}{2} \right)\left( i_k +\Gamma _k^{(F)} \right) +\frac{\tilde{q}}{4(1+a)}}{\left( i_k+\frac{k}{2}+\frac{1}{2} \right) \left( i_k+\frac{k}{2}+\frac{\beta}{2} -\frac{\gamma}{2}-\frac{j}{2} \right)} \right. \nonumber\\
&&\times \left. \frac{\left( -\frac{j}{2}+\frac{k}{2}\right)_{i_k} \left( \frac{\alpha }{2}+ \frac{k}{2} \right)_{i_k}\left( \frac{k}{2}+1 \right)_{i_{k-1}} \left( \frac{\beta}{2} -\frac{\gamma}{2}+ \frac{1}{2}+ \frac{k}{2}-\frac{j}{2} \right)_{i_{k-1}}}{\left( -\frac{j}{2}+\frac{k}{2}\right)_{i_{k-1}} \left( \frac{\alpha }{2}+ \frac{k}{2} \right)_{i_{k-1}}\left( \frac{k}{2} +1\right)_{i_k} \left( \frac{\beta}{2} -\frac{\gamma}{2}+ \frac{1}{2}+ \frac{k}{2}-\frac{j}{2} \right)_{i_k}} \right) \nonumber\\
&&\times \left. \sum_{i_{\tau } = i_{\tau -1}}^{m} \frac{\left( -\frac{j}{2}+\frac{\tau }{2}\right)_{i_{\tau }} \left( \frac{\alpha }{2}+ \frac{\tau }{2} \right)_{i_{\tau }}\left( \frac{\tau }{2}+1 \right)_{i_{\tau -1}} \left( \frac{\beta}{2} -\frac{\gamma}{2}+ \frac{1}{2}+\frac{\tau }{2}-\frac{j}{2} \right)_{i_{\tau -1}}}{\left( -\frac{j}{2}+\frac{\tau }{2}\right)_{i_{\tau -1}} \left( \frac{\alpha }{2}+ \frac{\tau }{2} \right)_{i_{\tau -1}}\left( \frac{\tau }{2}+1 \right)_{i_{\tau }} \left(\frac{\beta}{2} -\frac{\gamma}{2}+ \frac{1}{2}+\frac{\tau }{2}-\frac{j}{2} \right)_{i_{\tau }}} z^{i_{\tau }}\right\} \eta ^{\tau } \nonumber 
\end{eqnarray} 
where  
\begin{equation}
\begin{cases}
\tau \geq 2  \cr 
\xi= \frac{a(x-1)}{x-a} \cr
z = -\frac{1}{a}\xi^2 \cr
\eta = \frac{(1+a)}{a} \xi \cr
\tilde{q} =  q-(\gamma +j )\alpha \cr
\Gamma _0^{(F)} = \frac{1}{2(1+a)}\left( \alpha -\gamma -j +a\left(  \beta -1-j \right)\right) \cr
\Gamma _k^{(F)} =  \frac{1}{2(1+a)}\left( \alpha -\gamma -j+k +a\left( \beta -1-j+k \right)\right) 
\end{cases}\nonumber 
\end{equation} 
\end{enumerate}
\end{appendices} 

\addcontentsline{toc}{section}{Bibliography}
\bibliographystyle{model1a-num-names}
\bibliography{<your-bib-database>}
\chapter{Complete polynomials of Heun equation using reversible three-term recurrence formula}
\chaptermark{Complete polynomials of Heun equation using R3TRF} 
 \addtocontents{toc}{\protect\setcounter{tocdepth}{2}}
In this chapter I construct power series solutions of Heun equation for a polynomial which makes $A_n$ and $B_n$ terms terminated by applying mathematical expression of complete polynomials using reversible 3-term recurrence formula (R3TRF).

Nine examples of 192 local solutions of the Heun equation (Maier, 2007) are provided in the appendix.  For each example, I show the power series expansions of Heun equation for complete polynomials using R3TRF.

\section{Introduction}
Since we substitute a power series with unknown coefficients into a linear ordinary differential equation (ODE), a recurrence relation starts to appear. There can be between two term and multi-term in the recursion relation which includes all parameters in a ODE. A formal series solution of Hypergeoemtric equation consists of 2-term recurrence relation between successive coefficients. This equation generalizes all well-known ODEs having 2-term recursive relation in their series solutions such as Laguerre, Kummer, Legendre, Bessel, Coulomb Wave equations, etc. The Frobenius solution in a closed form of this equation including its definite or contour integral forms have been expressed analytically. 

Unfortunately, the phenomenon of the physical world is the non-linearized system. For simpler numerical computations and better apprehension of the nature, we usually linearize the physical system. Since we describe linearized natural phenomena in fields of the classical physics such as E \& M, statistical mechanic, Newtonian mechanic, quantum mechanic, thermodynamics , etc, differential equations come out of their physical system generally.  From the birth of differential equations until now, we have described the linearized natural sciences using simpler functions such as functions of hypergeometric type having a 2-term recursive relation in a power series.

However, since the modern physics such as QCD, quantum field theory, general relativity, SUSY, particle physics, etc came into   existence, we can not describe the nature with a linear ODE having a 2-term recursion relation in a formal series any more. 
When we deal with head on difficult mathematical physics problems with more complicated metircs or in higher dimensions, we are confronted with more than 3-term recursion relation between consecutive coefficients in a linear ODE. \cite{5Hortacsu:2011rr} 
 Even though a linearized ODE is the most simple form among all differential equations, its formal series solutions and integral representations are unknown.
 
According to Karl Heun \cite{5Heun1889,5Ronv1995}, Heun's differential equation is a second-order linear ODE having four regular singular points. The 4 different confluent forms, such as Confluent Heun, Doubly-Confluent Heun, Biconfluent Heun and Triconfluent Heun equations, are generated from Heun equation by a confluent process of the singular points. Its confluent process is similar as the derivation of confluent hypergeometric equation from hypergeometric equation.
 Recently Heun's equation starts to appear in the hydrogen-molecule ion \cite{5Wils1928}, in the Stark effect \cite{5Epst1926}, in black hole problems with the Kerr metric \cite{5Teuk1973,5Leav1985,5Bati2006,5Bati2007,5Bati2010,5Take2006}, in the fluid dynamics  \cite{5Cras1998,5Chug2009,5Chug2009a}, in quantum inozemtsev model \cite{5Inoz1989}, etc. 

The coefficients in a formal series of Heun's equation have a 3-term recurrence relation. The Heun's equation generalizes the most well-known ODEs having 3-term recurrence relation between successive coefficients including 2-term recursion relation such as: Spheroidal Wave, Lame, Mathieu, and hypergeometric $_2F_1$, $_1F_1$ and $_0F_1$ functions, etc. Sometimes several authors treat Heun equation as a successor of hypergeomegtric equations in $21^{th}$ century. \cite{5Fizi2012} For definite or contour integral forms of Heun equation, no such solutions have found because of its three different consecutive coefficients in a series: indeed, there are no analytic solutions in a form of its power series until now.  \cite{5Arsc1983} Instead, global properties of Heun equation is investigated by utilizing simple Fredholm integral equations. \cite{5Lamb1934,5Erde1942,5Slee1969a,5Slee1969b,5Arsc1964} 

Heun ordinary differential equation is given by \cite{5Slavy2000,5Erde1955,5Heun1889,5Ronv1995}
\begin{equation}
\frac{d^2{y}}{d{x}^2} + \left(\frac{\gamma }{x} +\frac{\delta }{x-1} + \frac{\epsilon }{x-a}\right) \frac{d{y}}{d{x}} +  \frac{\alpha \beta x-q}{x(x-1)(x-a)} y = 0 \label{eq:5001}
\end{equation}
With the condition $\epsilon = \alpha +\beta -\gamma -\delta +1$. The parameters play different roles: $a \ne 0 $ is the singularity parameter, determining the radius of convergence of a power series, $\alpha $, $\beta $, $\gamma $, $\delta $, $\epsilon $ are exponent parameters, $q$ is the accessory (spectral) parameter. Also, $\alpha $ and $\beta $ are identical to each other. The total number of free parameters is six. It has four regular singular points which are 0, 1, $a$ and $\infty $ with exponents $\{ 0, 1-\gamma \}$, $\{ 0, 1-\delta \}$, $\{ 0, 1-\epsilon \}$ and $\{ \alpha, \beta \}$. 

We assume the solution takes the form 
\begin{equation}
y(x)= \sum_{n=0}^{\infty } c_n x^{n+\lambda } \label{eq:5002}
\end{equation}
where $\lambda $ is an indicial root. 
Substituting (\ref{eq:5002}) into (\ref{eq:5001}) gives for the coefficients $c_n$ the recurrence relations
\begin{equation}
c_{n+1}=A_n \;c_n +B_n \;c_{n-1} \hspace{1cm};n\geq 1 \label{eq:5003}
\end{equation}
where
\begin{subequations}
\begin{eqnarray}
A_n &=& \frac{(n+\lambda )(n-1+\gamma +\epsilon +\lambda + a(n-1+\gamma +\lambda +\delta ))+q}{a(n+1+\lambda )(n+\gamma +\lambda )}\nonumber\\
&=& \frac{(n+\lambda )(n+\alpha +\beta -\delta +\lambda +a(n+\delta +\gamma -1+\lambda ))+q}{a(n+1+\lambda )(n+\gamma +\lambda )} \label{eq:5004a}
\end{eqnarray}
\begin{eqnarray}
B_n &=& -\frac{(n-1+\lambda )(n+\gamma +\delta +\epsilon -2+\lambda )+\alpha \beta }{a(n+1+\lambda )(n+\gamma +\lambda )} \nonumber\\
&=& - \frac{(n-1+\lambda +\alpha )(n-1+\lambda +\beta )}{a(n+1+\lambda )(n+\gamma +\lambda )} \label{eq:5004b}
\end{eqnarray}
\begin{equation}
c_1= A_0 \;c_0 \label{eq:5004c}
\end{equation}
\end{subequations}
We have two indicial roots which are $\lambda = 0$ and $ 1-\gamma $

\section{Power series}
As we all recognize, there are only 2 types of a formal series for the 2-term recursive relation between consecutive coefficients  in a linear ODE which are an infinite series and a polynomial. In contrast, there are $2^{3-1}$ possible power series  solutions of Heun equation having 3-term recurrence relation between successive coefficients. 

Table 6.1 tells us that the Frobenius solutions of Heun equation have dissimilarly an infinite series and three types of polynomials: (1) a polynomial which makes $B_n$ term terminated (2) a polynomial which makes $A_n$ term terminated and (3) a polynomial which makes $A_n$ and $B_n$ terms terminated, referred as `a complete polynomial.' 

\begin{table}[h] 
\begin{center}
\thispagestyle{plain}
\hspace*{-0.1\linewidth}\resizebox{1.2\linewidth}{!}
{
 \Tree[.{\Huge Heun's differential equation} [.{\Huge 3TRF} [.{\Huge Infinite series} ]
              [.{\Huge Polynomials} [[.{\Huge Polynomial of type 1} ]
               [.{\Huge Polynomial of type 3} [.{\Huge $ \begin{array}{lcll}  1^{\mbox{st}}\;  \mbox{species}\\ \mbox{complete} \\ \mbox{polynomial} \end{array}$} ] [.{\Huge $ \begin{array}{lcll}  2^{\mbox{nd}}\;  \mbox{species}\\ \mbox{complete} \\ \mbox{polynomial} \end{array}$} ]]]]]                         
 [.{\Huge R3TRF} [.{\Huge Infinite series} ]
   [.{\Huge Polynomials} [[.{\Huge Polynomial of type 2} ]
    [.{\Huge  Polynomial of type 3} [.{\Huge $ \begin{array}{lcll}  1^{\mbox{st}}\;  \mbox{species} \\ \mbox{complete} \\ \mbox{polynomial} \end{array}$} ] [.{\Huge $ \begin{array}{lcll}  2^{\mbox{nd}}\;  \mbox{species} \\ \mbox{complete} \\ \mbox{polynomial} \end{array}$} ]]]]]]
}
\end{center}
\caption{Power series of Heun's differential equation}
\end{table}
The sequence $c_n$ combines into combinations of $A_n$ and $B_n$ terms in (\ref{eq:5003}). 
By allowing $A_n$ in the sequence $c_n$ is the leading term of each sub-power series in a function $y(x)$ \cite{5Choun2012}, I construct the general summation formulas of the 3-term recurrence relation in a linear ODE for an infinite series and a polynomial of type 1: I observe the term of sequence $c_n$ which includes zero term of $A_n's$, one term of $A_n's$, two terms of $A_n's$, three terms of $A_n's$, etc. I designate this mathematical technique as `three term recurrence formula (3TRF).' 
I construct power series solutions in closed form of Heun equation around $x=0$ for an infinite series and a polynomial of type 1 by applying 3TRF. \cite{5Chou2012c,5Chou2012d} For a polynomial of type 1, I treat $\gamma $, $\delta $ and $q$ as free variables and fixed values of $\alpha $ and/or $\beta $.

By allowing $B_n$ in the sequence $c_n$ is the leading term of each sub-power series in a function $y(x)$ \cite{5Choun2013}, I obtain the general summation formulas of the 3-term recurrence relation in a linear ODE for an infinite series and a polynomial of type 2: I observe the term of sequence $c_n$ which includes zero term of $B_n's$, one term of $B_n's$, two terms of $B_n's$, three terms of $B_n's$, etc. I denominate this mathematical technique as `reversible three term recurrence formula (R3TRF).' 
I derive power series solutions of Heun equation around $x=0$ for an infinite series and a polynomial of type 2 by applying R3TRF in chapter 2 of Ref.\cite{5Choun2013}. For a polynomial of type 2, I treat $\alpha $, $\beta $, $\gamma $ and $\delta $ as free variables and a fixed value of $q$. Infinite series solutions using 3TRF of Heun equation around $x=0$ are equivalent to their infinite series using R3TRF. The former is that $A_n$ is the leading term in each sub-power series of Heun equation. The latter is that $B_n$ is the leading term in each sub-power series of it.
  
In general, a spectral polynomial of Heun equation has been known as a polynomial of type 3 where $A_n$ and $B_n$ terms terminated. For a formal solution of Heun equation around $x=0$ in (\ref{eq:5001}), Heun spectral polynomial has a fixed value of $\alpha $ or $\beta $, just as it has a fixed value of $q$. Heun (spectral) polynomial has been investigated by many great scholars. 
Darboux notices that Heun's equation is related with the generalization of Lam\'{e}'s equation. \cite{5Darb1882} He shows the general solutions of the generalized Lam\'{e}'s equation with half-odd-integer parameters. Treibich and Verdier prove the finite-gapness of the potential in the generalized Lam\'{e}'s equation for any integer parameters: its general solutions with integer parameters has been studied by various mathematicians. \cite{5Trei1990,5Trei1990a,5Belo1994,5Smir1989,5Smir1994}    

Smirnov proves that Frobenius solutions of Heun equation is obtained in terms of Heun polynomials since an accessory parameter is located at branching points of the hyperelliptic curve with any integer exponent parameters. \cite{5Smir2001}
He also shows an integral of finite-gap solutions of Heun equation with an accessory parameter belonging to a hyperelliptic spectral curve.

Recently, Shapiro \textit {et al.} study spectral polynomials of the Heun equation in the case of real roots of these polynomials and asymptotic root distribution when complex roots are present. \cite{5Shap2010,5Shap2011,5Shap2012}
Mart\'{i}nez-Finkelshtein and Rakhmanov investigate the asymptotic zero distribution of Heine-Stieltjes polynomials, which is the general form of Heun polynomials, with complex polynomial coefficients. \cite{5Mart2011}

Kalnins \textit {et al.} describe Heun polynomials in terms of Jacobi polynomials by using group theory and its connection with the method of separation of variables applied to the Laplace-Beltrami eigenvalue equation on the $n$-sphere. \cite{5Kaln1989,5Kaln1990}
Patera and Winternitz show a new basis in the representation theory of the rotation group $O(3)$. And this new basis is composed of products of two Lam\'{e} polynomials in terms of functions on an $O(3)$ sphere. Also, they notice that the basis are Heun polynomials in a space of functions of one complex variable. \cite{5Pate1973}  
 
Complete polynomials in the 3-term recurrence relation in a linear ODE can be classified into two different types which are (1) the first species complete polynomial where a parameter of a numerator in $B_n$ term and a (spectral) parameter of a numerator in $A_n$ term are fixed constants and (2) the second species complete polynomial where two parameters of a numerator in $B_n$ term and a parameter of a numerator in $A_n$ term are fixed constants.
For the first species complete polynomial of Heun equation around  $x=0$ in (\ref{eq:5001}), an exponent parameter $\alpha $ (or $\beta $) is a fixed constant and an accessory parameter $q$ has multi-valued roots. For the second species complete polynomial of Heun equation, parameters $\alpha $, $\beta $ and $q$ are fixed constants.


In chapter 2, I obtain the mathematical expressions of complete polynomials for the first and second species, by allowing $B_n$ as the leading term in each sub-power series of the general power series $y(x)$, as ``complete polynomials using reversible 3-term recurrence formula (R3TRF)''
In this chapter I show Frobenius solutions in compact forms, called summation notation, of Heun equation around $x=0$ for two types of polynomials which make $A_n$ and $B_n$ terms terminated by applying complete polynomials using R3TRF.

\subsection{The first species complete polynomial using R3TRF}
For the first species complete polynomial, we need a condition which is given by
\begin{equation}
B_{j+1}= c_{j+1}=0\hspace{1cm}\mathrm{where}\;j\in \mathbb{N}_{0}  
 \label{eq:5005}
\end{equation}
(\ref{eq:5005}) gives successively $c_{j+2}=c_{j+3}=c_{j+4}=\cdots=0$. And $c_{j+1}=0$ is defined by a polynomial equation of degree $j+1$ for the determination of an accessory parameter in $A_n$ term. 
\begin{theorem}
In chapter 2, the general expression of a function $y(x)$ for the first species complete polynomial using reversible 3-term recurrence formula and its algebraic equation for the determination of an accessory parameter in $A_n$ term are given by
\begin{enumerate} 
\item As $B_1=0$,
\begin{equation}
0 =\bar{c}(0,1) \label{eq:5006a}
\end{equation}
\begin{equation}
y(x) = y_{0}^{0}(x) \label{eq:5006b}
\end{equation}
\item As $B_2=0$, 
\begin{equation}
0 = \bar{c}(0,2)+\bar{c}(1,0) \label{eq:5007a}
\end{equation}
\begin{equation}
y(x)= y_{0}^{1}(x) \label{eq:5007b}
\end{equation}
\item As $B_{2N+3}=0$ where $N \in \mathbb{N}_{0}$,
\begin{equation}
0  = \sum_{r=0}^{N+1}\bar{c}\left( r, 2(N-r)+3\right) \label{eq:5008a}
\end{equation}
\begin{equation}
y(x)= \sum_{r=0}^{N+1} y_{r}^{2(N+1-r)}(x) \label{eq:5008b}
\end{equation}
\item As $B_{2N+4}=0$ where$N \in \mathbb{N}_{0}$,
\begin{equation}
0  = \sum_{r=0}^{N+2}\bar{c}\left( r, 2(N+2-r)\right) \label{eq:5009a}
\end{equation}
\begin{equation}
y(x)=  \sum_{r=0}^{N+1} y_{r}^{2(N-r)+3}(x) \label{eq:5009b}
\end{equation}
In the above,
\begin{eqnarray}
\bar{c}(0,n) &=& \prod _{i_0=0}^{n-1}A_{i_0} \label{eq:50010a}\\
\bar{c}(1,n) &=& \sum_{i_0=0}^{n} \left\{ B_{i_0+1} \prod _{i_1=0}^{i_0-1}A_{i_1} \prod _{i_2=i_0}^{n-1}A_{i_2+2} \right\} \label{eq:50010b}\\
\bar{c}(\tau ,n) &=& \sum_{i_0=0}^{n} \left\{B_{i_0+1}\prod _{i_1=0}^{i_0-1} A_{i_1} 
\prod _{k=1}^{\tau -1} \left( \sum_{i_{2k}= i_{2(k-1)}}^{n} B_{i_{2k}+(2k+1)}\prod _{i_{2k+1}=i_{2(k-1)}}^{i_{2k}-1}A_{i_{2k+1}+2k}\right)\right.\nonumber\\
&&\times \left. \prod _{i_{2\tau} = i_{2(\tau -1)}}^{n-1} A_{i_{2\tau }+ 2\tau} \right\} 
 \label{eq:50010c}
\end{eqnarray}
and
\begin{eqnarray}
y_0^m(x) &=& c_0 x^{\lambda} \sum_{i_0=0}^{m} \left\{ \prod _{i_1=0}^{i_0-1}A_{i_1} \right\} x^{i_0 } \label{eq:50011a}\\
y_1^m(x) &=& c_0 x^{\lambda} \sum_{i_0=0}^{m}\left\{ B_{i_0+1} \prod _{i_1=0}^{i_0-1}A_{i_1}  \sum_{i_2=i_0}^{m} \left\{ \prod _{i_3=i_0}^{i_2-1}A_{i_3+2} \right\}\right\} x^{i_2+2 } \label{eq:50011b}\\
y_{\tau }^m(x) &=& c_0 x^{\lambda} \sum_{i_0=0}^{m} \left\{B_{i_0+1}\prod _{i_1=0}^{i_0-1} A_{i_1} 
\prod _{k=1}^{\tau -1} \left( \sum_{i_{2k}= i_{2(k-1)}}^{m} B_{i_{2k}+(2k+1)}\prod _{i_{2k+1}=i_{2(k-1)}}^{i_{2k}-1}A_{i_{2k+1}+2k}\right) \right. \nonumber\\
&& \times \left. \sum_{i_{2\tau} = i_{2(\tau -1)}}^{m} \left( \prod _{i_{2\tau +1}=i_{2(\tau -1)}}^{i_{2\tau}-1} A_{i_{2\tau +1}+ 2\tau} \right) \right\} x^{i_{2\tau}+2\tau }\hspace{1cm}\mathrm{where}\;\tau \geq 2
\label{eq:50011c}
\end{eqnarray}
\end{enumerate}
\end{theorem}
 Put $n= j+1$ in (\ref{eq:5004b}) and use the condition $B_{j+1}=0$ for $\alpha $.  
\begin{equation}
\alpha = -j-\lambda 
\label{eq:50011}
\end{equation}
Take (\ref{eq:50011}) into (\ref{eq:5004a}) and (\ref{eq:5004b}).
\begin{subequations}
\begin{equation}
A_n = \frac{1+a}{a}\frac{\left( n+ \Delta_0^{-} \left( j,q\right) \right) \left( n+ \Delta_0^{+} \left( j,q\right)\right)}{(n+1+\lambda )(n+\gamma +\lambda )} \label{eq:50012a}
\end{equation}
\begin{equation}
B_n = -\frac{1}{a} \frac{(n-1-j)(n-1+\lambda +\beta )}{(n+1+\lambda )(n+\gamma +\lambda )} \label{eq:50012b}
\end{equation}
\end{subequations}
where
\begin{equation}
\begin{cases} \Delta_k^{\pm} \left( j,q\right) = \frac{\varphi +2(1+a)\left( \lambda +2k\right) \pm \sqrt{\varphi ^2-4(1+a)q}}{2(1+a)}  \cr
\varphi = \beta -\delta -\lambda -j+a(\gamma +\delta -1)  
\end{cases}\nonumber
\end{equation} 
Now the condition $c_{j+1}=0$ is clearly an algebraic equation in $q$ of degree $j+1$ and thus has $j+1$ zeros denoted them by $q_j^m$ eigenvalues where $m = 0,1,2, \cdots, j$. They can be arranged in the following order: $q_j^0 < q_j^1 < q_j^2 < \cdots < q_j^j$.
 
Substitute (\ref{eq:50012a}) and (\ref{eq:50012b}) into (\ref{eq:50010a})--(\ref{eq:50011c}).

As $B_{1}= c_{1}=0$, take the new (\ref{eq:50010a}) into (\ref{eq:5006a}) putting $j=0$. Substitute the new (\ref{eq:50011a}) into (\ref{eq:5006b}) putting $j=0$.

As $B_{2}= c_{2}=0$, take the new (\ref{eq:50010a}) and (\ref{eq:50010b}) into (\ref{eq:5007a}) putting $j=1$. Substitute the new (\ref{eq:50011a}) into (\ref{eq:5007b}) putting $j=1$ and $q=q_1^m$. 

As $B_{2N+3}= c_{2N+3}=0$, take the new (\ref{eq:50010a})--(\ref{eq:50010c}) into (\ref{eq:5008a}) putting $j=2N+2$. Substitute the new 
(\ref{eq:50011a})--(\ref{eq:50011c}) into (\ref{eq:5008b}) putting $j=2N+2$ and $q=q_{2N+2}^m$.

As $B_{2N+4}= c_{2N+4}=0$, take the new (\ref{eq:50010a})--(\ref{eq:50010c}) into (\ref{eq:5009a}) putting $j=2N+3$. Substitute the new 
(\ref{eq:50011a})--(\ref{eq:50011c}) into (\ref{eq:5009b}) putting $j=2N+3$ and $q=q_{2N+3}^m$.

After the replacement process, the general expression of power series of Heun equation about $x=0$ for the first species complete polynomial using reversibe 3-term recurrence formula and its algebraic equation for the determination of an accessory parameter $q$ are given by
\begin{enumerate} 
\item As $\alpha = -\lambda $,

An algebraic equation of degree 1 for the determination of $q$ is given by
\begin{equation}
0= \bar{c}(0,1;0,q)= q + \lambda (\beta -\delta +a(\gamma+ \delta -1+\lambda )) \label{eq:50013a}
\end{equation}
The eigenvalue of $q$ is written by $q_0^0$. Its eigenfunction is given by
\begin{equation}
y(x) = y_0^0\left( 0,q_0^0;x\right)= c_0 x^{\lambda } \label{eq:50013b}  
\end{equation}
\item As $\alpha =-1-\lambda $,

An algebraic equation of degree 2 for the determination of $q$ is given by
\begin{eqnarray}
0 &=& \bar{c}(0,2;1,q)+\bar{c}(1,0;1,q) \label{eq:50014a}\\
&=& a(1+\lambda )(\beta +\lambda )(\gamma +\lambda ) 
+\prod_{l=0}^{1}\Big( q+(\lambda +l)(\beta -\delta -1+l+a(\gamma +\delta +\lambda -1+l))\Big) \nonumber
\end{eqnarray}
The eigenvalue of $q$ is written by $q_1^m$ where $m = 0,1 $; $q_{1}^0 < q_{1}^1$. Its eigenfunction is given by
\begin{eqnarray}
y(x) &=& y_{0}^{1}\left( 1,q_1^m;x\right) \nonumber\\
&=& c_0 x^{\lambda } \left\{ 1+\frac{\lambda \left( \beta -\delta -1+a(\gamma +\delta +\lambda -1)\right) + q_1^m}{(1+a)(1+\lambda )(\gamma +\lambda )} \eta \right\} \label{eq:50014b}  
\end{eqnarray}
\item As $\alpha =-2N-2-\lambda $ where $N \in \mathbb{N}_{0}$,

An algebraic equation of degree $2N+3$ for the determination of $q$ is given by
\begin{equation}
0 =  \sum_{r=0}^{N+1}\bar{c}\left( r, 2(N-r)+3; 2N+2,q\right)  \label{eq:50015a}
\end{equation}
The eigenvalue of $q$ is written by $q_{2N+2}^m$ where $m = 0,1,2,\cdots,2N+2 $; $q_{2N+2}^0 < q_{2N+2}^1 < \cdots < q_{2N+2}^{2N+2}$. Its eigenfunction is given by 
\begin{equation} 
y(x) = \sum_{r=0}^{N+1} y_{r}^{2(N+1-r)}\left( 2N+2, q_{2N+2}^m; x \right)  
\label{eq:50015b} 
\end{equation}
\item As $\alpha =-2N-3-\lambda $ where $N \in \mathbb{N}_{0}$,

An algebraic equation of degree $2N+4$ for the determination of $q$ is given by
\begin{equation}  
0 =  \sum_{r=0}^{N+2}\bar{c}\left( r, 2(N+2-r); 2N+3,q\right) \label{eq:50016a}
\end{equation}
The eigenvalue of $q$ is written by $q_{2N+3}^m$ where $m = 0,1,2,\cdots,2N+3 $; $q_{2N+3}^0 < q_{2N+3}^1 < \cdots < q_{2N+3}^{2N+3}$. Its eigenfunction is given by
\begin{equation} 
y(x) =  \sum_{r=0}^{N+1} y_{r}^{2(N-r)+3} \left( 2N+3,q_{2N+3}^m;x\right) \label{eq:50016b}
\end{equation}
In the above,
\begin{eqnarray}
\bar{c}(0,n;j,q)  &=& \frac{\left( \Delta_0^{-} \left( j,q\right) \right)_{n}\left( \Delta_0^{+} \left( j,q\right) \right)_{n}}{\left( 1+ \lambda \right)_{n} \left( \gamma + \lambda \right)_{n}} \left( \frac{1+a}{a} \right)^{n}\label{eq:50017a}\\
\bar{c}(1,n;j,q) &=& \left( -\frac{1}{a}\right) \sum_{i_0=0}^{n}\frac{\left( i_0 -j\right)\left( i_0+\beta +\lambda  \right) }{\left( i_0+2+ \lambda \right) \left( i_0+1+ \gamma + \lambda \right)} \frac{ \left( \Delta_0^{-} \left( j,q\right) \right)_{i_0}\left( \Delta_0^{+} \left( j,q\right) \right)_{i_0}}{\left( 1+ \lambda \right)_{i_0} \left( \gamma + \lambda \right)_{i_0}} \nonumber\\
&\times&  \frac{ \left( \Delta_1^{-} \left( j,q\right) \right)_{n}\left( \Delta_1^{+} \left( j,q\right) \right)_{n} \left( 3 + \lambda \right)_{i_0} \left( 2+ \gamma + \lambda \right)_{i_0}}{\left( \Delta_1^{-} \left( j,q\right) \right)_{i_0}\left( \Delta_1^{+} \left( j,q\right) \right)_{i_0}\left( 3 + \lambda \right)_{n} \left( 2+ \gamma + \lambda \right)_{n}} \left(\frac{1+a}{a} \right)^{n }  
\label{eq:50017b}\\
\bar{c}(\tau ,n;j,q) &=& \left( -\frac{1}{a}\right)^{\tau} \sum_{i_0=0}^{n}\frac{\left( i_0 -j\right)\left( i_0+\beta +\lambda  \right) }{\left( i_0+2+ \lambda \right) \left( i_0+1+ \gamma + \lambda \right)} \frac{ \left( \Delta_0^{-} \left( j,q\right) \right)_{i_0}\left( \Delta_0^{+} \left( j,q\right) \right)_{i_0}}{\left( 1+ \lambda \right)_{i_0} \left( \gamma + \lambda \right)_{i_0}} \nonumber\\
&\times&  \prod_{k=1}^{\tau -1} \left( \sum_{i_k = i_{k-1}}^{n} \frac{\left( i_k+ 2k-j\right)\left( i_k +2k+\beta +\lambda \right)}{\left( i_k+2k+2+ \lambda \right) \left( i_k+2k+1+ \gamma + \lambda \right)} \right. \nonumber\\
&\times &  \left. \frac{ \left( \Delta_k^{-} \left( j,q\right) \right)_{i_k}\left( \Delta_k^{+} \left( j,q\right) \right)_{i_k} \left( 2k+1 + \lambda \right)_{i_{k-1}} \left( 2k+ \gamma + \lambda \right)_{i_{k-1}}}{\left( \Delta_k^{-} \left( j,q\right) \right)_{i_{k-1}}\left( \Delta_k^{+} \left( j,q\right) \right)_{i_{k-1}}\left( 2k+1 + \lambda \right)_{i_k} \left( 2k+ \gamma + \lambda \right)_{i_k}} \right)  \nonumber\\
&\times& \frac{ \left( \Delta_{\tau }^{-} \left( j,q\right) \right)_{n}\left( \Delta_{\tau }^{+} \left( j,q\right) \right)_{n} \left( 2\tau +1+ \lambda \right)_{i_{\tau -1}} \left( 2\tau + \gamma + \lambda \right)_{i_{\tau -1}}}{\left( \Delta_{\tau }^{-} \left( j,q\right) \right)_{i_{\tau -1}}\left( \Delta_{\tau }^{+} \left( j,q\right) \right)_{i_{\tau -1}}\left( 2\tau +1 + \lambda \right)_{n} \left( 2\tau + \gamma + \lambda \right)_{n}} \left(\frac{1+a}{a}\right)^{n } \hspace{1.5cm}\label{eq:50017c}
\end{eqnarray}
\begin{eqnarray}
y_0^m(j,q;x) &=& c_0 x^{\lambda }  \sum_{i_0=0}^{m} \frac{\left( \Delta_0^{-} \left( j,q\right) \right)_{i_0}\left( \Delta_0^{+} \left( j,q\right) \right)_{i_0}}{\left( 1+ \lambda \right)_{i_0} \left( \gamma + \lambda \right)_{i_0}} \eta ^{i_0} \label{eq:50018a}\\
y_1^m(j,q;x) &=& c_0 x^{\lambda } \left\{\sum_{i_0=0}^{m}\frac{\left( i_0 -j\right)\left( i_0+\beta +\lambda  \right) }{\left( i_0+2+ \lambda \right) \left( i_0+1+ \gamma + \lambda \right)} \frac{ \left( \Delta_0^{-} \left( j,q\right) \right)_{i_0}\left( \Delta_0^{+} \left( j,q\right) \right)_{i_0}}{\left( 1+ \lambda \right)_{i_0} \left( \gamma + \lambda \right)_{i_0}} \right. \nonumber\\
&\times& \left. \sum_{i_1 = i_0}^{m} \frac{ \left( \Delta_1^{-} \left( j,q\right) \right)_{i_1}\left( \Delta_1^{+} \left( j,q\right) \right)_{i_1} \left( 3 + \lambda \right)_{i_0} \left( 2+ \gamma + \lambda \right)_{i_0}}{\left( \Delta_1^{-} \left( j,q\right) \right)_{i_0}\left( \Delta_1^{+} \left( j,q\right) \right)_{i_0}\left( 3 + \lambda \right)_{i_1} \left( 2+ \gamma + \lambda \right)_{i_1}} \eta ^{i_1}\right\} z 
\label{eq:50018b} 
\end{eqnarray}
\begin{eqnarray}
y_{\tau }^m(j,q;x) &=& c_0 x^{\lambda } \left\{ \sum_{i_0=0}^{m} \frac{\left( i_0 -j\right)\left( i_0+\beta +\lambda  \right) }{\left( i_0+2+ \lambda \right) \left( i_0+1+ \gamma + \lambda \right)} \frac{ \left( \Delta_0^{-} \left( j,q\right) \right)_{i_0}\left( \Delta_0^{+} \left( j,q\right) \right)_{i_0}}{\left( 1+ \lambda \right)_{i_0} \left( \gamma + \lambda \right)_{i_0}} \right.\nonumber\\
&\times& \prod_{k=1}^{\tau -1} \left( \sum_{i_k = i_{k-1}}^{m} \frac{\left( i_k+ 2k-j\right)\left( i_k +2k+\beta +\lambda \right)}{\left( i_k+2k+2+ \lambda \right) \left( i_k+2k+1+ \gamma + \lambda \right)} \right. \nonumber\\
&\times &  \left. \frac{ \left( \Delta_k^{-} \left( j,q\right) \right)_{i_k}\left( \Delta_k^{+} \left( j,q\right) \right)_{i_k} \left( 2k+1 + \lambda \right)_{i_{k-1}} \left( 2k+ \gamma + \lambda \right)_{i_{k-1}}}{\left( \Delta_k^{-} \left( j,q\right) \right)_{i_{k-1}}\left( \Delta_k^{+} \left( j,q\right) \right)_{i_{k-1}}\left( 2k+1 + \lambda \right)_{i_k} \left( 2k+ \gamma + \lambda \right)_{i_k}} \right) \nonumber\\
&\times & \left. \sum_{i_{\tau } = i_{\tau -1}}^{m}  \frac{ \left( \Delta_{\tau }^{-} \left( j,q\right) \right)_{i_{\tau }}\left( \Delta_{\tau }^{+} \left( j,q\right) \right)_{i_{\tau }} \left( 2\tau +1 + \lambda \right)_{i_{\tau -1}} \left( 2\tau + \gamma + \lambda \right)_{i_{\tau -1}}}{\left( \Delta_{\tau }^{-} \left( j,q\right) \right)_{i_{\tau -1}}\left( \Delta_{\tau }^{+} \left( j,q\right) \right)_{i_{\tau -1}}\left( 2\tau +1 + \lambda \right)_{i_\tau } \left( 2\tau + \gamma + \lambda \right)_{i_{\tau }}} \eta ^{i_{\tau }}\right\} z^{\tau } \hspace{1.5cm}\label{eq:50018c} 
\end{eqnarray}
where
\begin{equation}
\begin{cases} \tau \geq 2 \cr
z = -\frac{1}{a}x^2 \cr
\eta = \frac{(1+a)}{a} x \cr
\Delta_k^{\pm} \left( j,q\right) = \frac{\varphi +2(1+a)\left( \lambda +2k\right) \pm \sqrt{\varphi ^2-4(1+a)q}}{2(1+a)}  \cr
\varphi = \beta -\delta -\lambda -j+a(\gamma +\delta -1)
\end{cases}\nonumber
\end{equation}
\end{enumerate}
Put $c_0$= 1 as $\lambda =0$ for the first kind of independent solutions of Heun equation and $\displaystyle{ c_0= \left( \frac{1+a}{a}\right)^{1-\gamma }}$ as $\lambda = 1-\gamma $ for the second one in (\ref{eq:50013a})--(\ref{eq:50018c}). 
\begin{remark}
The power series expansion of Heun equation of the first kind for the first species complete polynomial using R3TRF about $x=0$ is given by
\begin{enumerate} 
\item As $\alpha =0$ and $q=q_0^0=0$,

The eigenfunction is given by
\begin{equation}
y(x) = H_pF_{0,0}^{R} \left( a, q=q_0^0=0; \alpha =0, \beta, \gamma, \delta ; \eta = \frac{(1+a)}{a} x ; z= -\frac{1}{a} x^2 \right) =1 \label{eq:50019}
\end{equation}
\item As $\alpha =-1$,

An algebraic equation of degree 2 for the determination of $q$ is given by
\begin{equation}
0 = a\beta \gamma 
+\prod_{l=0}^{1}\Big( q+ l(\beta -\delta -1+l+a(\gamma +\delta -1+l))\Big) \label{eq:50020a}
\end{equation}
The eigenvalue of $q$ is written by $q_1^m$ where $m = 0,1 $; $q_{1}^0 < q_{1}^1$. Its eigenfunction is given by
\begin{eqnarray}
y(x) &=& H_pF_{1,m}^R \left( a, q=q_1^m; \alpha = -1, \beta, \gamma, \delta ; \eta = \frac{(1+a)}{a} x ; z= -\frac{1}{a} x^2 \right)\nonumber\\
&=&  1+\frac{ q_1^m}{(1+a)\gamma } \eta \label{eq:50020b}  
\end{eqnarray}
\item As $\alpha =-2N-2 $ where $N \in \mathbb{N}_{0}$,

An algebraic equation of degree $2N+3$ for the determination of $q$ is given by
\begin{equation}
0 =  \sum_{r=0}^{N+1}\bar{c}\left( r, 2(N-r)+3; 2N+2,q\right)  \label{eq:50021a}
\end{equation}
The eigenvalue of $q$ is written by $q_{2N+2}^m$ where $m = 0,1,2,\cdots,2N+2 $; $q_{2N+2}^0 < q_{2N+2}^1 < \cdots < q_{2N+2}^{2N+2}$. Its eigenfunction is given by 
\begin{eqnarray} 
y(x) &=& H_pF_{2N+2,m}^R \left( a, q=q_{2N+2}^m; \alpha =-2N-2, \beta, \gamma, \delta ; \eta = \frac{(1+a)}{a} x ; z= -\frac{1}{a} x^2 \right)\nonumber\\
&=& \sum_{r=0}^{N+1} y_{r}^{2(N+1-r)}\left( 2N+2, q_{2N+2}^m; x \right)  
\label{eq:50021b} 
\end{eqnarray}
\item As $\alpha =-2N-3 $ where $N \in \mathbb{N}_{0}$,

An algebraic equation of degree $2N+4$ for the determination of $q$ is given by
\begin{equation}  
0 = \sum_{r=0}^{N+2}\bar{c}\left( r, 2(N+2-r); 2N+3,q\right) \label{eq:50022a}
\end{equation}
The eigenvalue of $q$ is written by $q_{2N+3}^m$ where $m = 0,1,2,\cdots,2N+3 $; $q_{2N+3}^0 < q_{2N+3}^1 < \cdots < q_{2N+3}^{2N+3}$. Its eigenfunction is given by
\begin{eqnarray} 
y(x) &=& H_pF_{2N+3,m}^R \left( a, q=q_{2N+3}^m; \alpha =-2N-3, \beta, \gamma, \delta ; \eta = \frac{(1+a)}{a} x ; z= -\frac{1}{a} x^2 \right)\nonumber\\
&=&   \sum_{r=0}^{N+1} y_{r}^{2(N-r)+3} \left( 2N+3,q_{2N+3}^m;x\right) \label{eq:50022b}
\end{eqnarray}
In the above,
\begin{eqnarray}
\bar{c}(0,n;j,q)  &=& \frac{\left( \Delta_0^{-} \left( j,q\right) \right)_{n}\left( \Delta_0^{+} \left( j,q\right) \right)_{n}}{\left( 1 \right)_{n} \left( \gamma \right)_{n}} \left( \frac{1+a}{a} \right)^{n}\label{eq:50023a}\\
\bar{c}(1,n;j,q) &=& \left( -\frac{1}{a}\right) \sum_{i_0=0}^{n}\frac{\left( i_0 -j\right)\left( i_0+\beta \right) }{\left( i_0+2 \right) \left( i_0+1+ \gamma \right)} \frac{ \left( \Delta_0^{-} \left( j,q\right) \right)_{i_0}\left( \Delta_0^{+} \left( j,q\right) \right)_{i_0}}{\left( 1 \right)_{i_0} \left( \gamma \right)_{i_0}} \nonumber\\
&\times&  \frac{ \left( \Delta_1^{-} \left( j,q\right) \right)_{n}\left( \Delta_1^{+} \left( j,q\right) \right)_{n} \left( 3 \right)_{i_0} \left( 2+ \gamma \right)_{i_0}}{\left( \Delta_1^{-} \left( j,q\right) \right)_{i_0}\left( \Delta_1^{+} \left( j,q\right) \right)_{i_0}\left( 3 \right)_{n} \left( 2+ \gamma \right)_{n}} \left(\frac{1+a}{a} \right)^{n }  
\label{eq:50023b}\\
\bar{c}(\tau ,n;j,q) &=& \left( -\frac{1}{a}\right)^{\tau} \sum_{i_0=0}^{n}\frac{\left( i_0 -j\right)\left( i_0+\beta \right) }{\left( i_0+2 \right) \left( i_0+1+ \gamma \right)} \frac{ \left( \Delta_0^{-} \left( j,q\right) \right)_{i_0}\left( \Delta_0^{+} \left( j,q\right) \right)_{i_0}}{\left( 1 \right)_{i_0} \left( \gamma \right)_{i_0}} \nonumber\\
&\times& \prod_{k=1}^{\tau -1} \left( \sum_{i_k = i_{k-1}}^{n} \frac{\left( i_k+ 2k-j\right)\left( i_k +2k+\beta \right)}{\left( i_k+2k+2 \right) \left( i_k+2k+1+ \gamma \right)} \right. \nonumber\\
&\times &  \left. \frac{ \left( \Delta_k^{-} \left( j,q\right) \right)_{i_k}\left( \Delta_k^{+} \left( j,q\right) \right)_{i_k} \left( 2k+1 \right)_{i_{k-1}} \left( 2k+ \gamma \right)_{i_{k-1}}}{\left( \Delta_k^{-} \left( j,q\right) \right)_{i_{k-1}}\left( \Delta_k^{+} \left( j,q\right) \right)_{i_{k-1}}\left( 2k+1 \right)_{i_k} \left( 2k+ \gamma \right)_{i_k}} \right) \nonumber\\
&\times& \frac{ \left( \Delta_{\tau }^{-} \left( j,q\right) \right)_{n}\left( \Delta_{\tau }^{+} \left( j,q\right) \right)_{n} \left( 2\tau +1 \right)_{i_{\tau -1}} \left( 2\tau + \gamma \right)_{i_{\tau -1}}}{\left( \Delta_{\tau }^{-} \left( j,q\right) \right)_{i_{\tau -1}}\left( \Delta_{\tau }^{+} \left( j,q\right) \right)_{i_{\tau -1}}\left( 2\tau +1 \right)_{n} \left( 2\tau + \gamma  \right)_{n}} \left(\frac{1+a}{a}\right)^{n } \hspace{1.5cm}\label{eq:50023c} 
\end{eqnarray}
\begin{eqnarray}
y_0^m(j,q;x) &=& \sum_{i_0=0}^{m} \frac{\left( \Delta_0^{-} \left( j,q\right) \right)_{i_0}\left( \Delta_0^{+} \left( j,q\right) \right)_{i_0}}{\left( 1 \right)_{i_0} \left( \gamma \right)_{i_0}} \eta ^{i_0} \label{eq:50024a}\\
y_1^m(j,q;x) &=& \left\{\sum_{i_0=0}^{m}\frac{\left( i_0 -j\right)\left( i_0+\beta \right) }{\left( i_0+2 \right) \left( i_0+1+ \gamma \right)} \frac{ \left( \Delta_0^{-} \left( j,q\right) \right)_{i_0}\left( \Delta_0^{+} \left( j,q\right) \right)_{i_0}}{\left( 1 \right)_{i_0} \left( \gamma \right)_{i_0}} \right. \nonumber\\
&\times& \left. \sum_{i_1 = i_0}^{m} \frac{ \left( \Delta_1^{-} \left( j,q\right) \right)_{i_1}\left( \Delta_1^{+} \left( j,q\right) \right)_{i_1} \left( 3 \right)_{i_0} \left( 2+ \gamma \right)_{i_0}}{\left( \Delta_1^{-} \left( j,q\right) \right)_{i_0}\left( \Delta_1^{+} \left( j,q\right) \right)_{i_0}\left( 3 \right)_{i_1} \left( 2+ \gamma \right)_{i_1}} \eta ^{i_1}\right\} z 
\label{eq:50024b}
\end{eqnarray}
\begin{eqnarray}
y_{\tau }^m(j,q;x) &=& \left\{ \sum_{i_0=0}^{m} \frac{\left( i_0 -j\right)\left( i_0+\beta \right) }{\left( i_0+2 \right) \left( i_0+1+ \gamma \right)} \frac{ \left( \Delta_0^{-} \left( j,q\right) \right)_{i_0}\left( \Delta_0^{+} \left( j,q\right) \right)_{i_0}}{\left( 1 \right)_{i_0} \left( \gamma \right)_{i_0}} \right.\nonumber\\
&\times& \prod_{k=1}^{\tau -1} \left( \sum_{i_k = i_{k-1}}^{m} \frac{\left( i_k+ 2k-j\right)\left( i_k +2k+\beta \right)}{\left( i_k+2k+2 \right) \left( i_k+2k+1+ \gamma \right)} \right. \nonumber\\
&\times &  \left. \frac{ \left( \Delta_k^{-} \left( j,q\right) \right)_{i_k}\left( \Delta_k^{+} \left( j,q\right) \right)_{i_k} \left( 2k+1 \right)_{i_{k-1}} \left( 2k+ \gamma \right)_{i_{k-1}}}{\left( \Delta_k^{-} \left( j,q\right) \right)_{i_{k-1}}\left( \Delta_k^{+} \left( j,q\right) \right)_{i_{k-1}}\left( 2k+1 \right)_{i_k} \left( 2k+ \gamma \right)_{i_k}} \right) \nonumber\\
&\times & \left. \sum_{i_{\tau } = i_{\tau -1}}^{m}  \frac{ \left( \Delta_{\tau }^{-} \left( j,q\right) \right)_{i_{\tau }}\left( \Delta_{\tau }^{+} \left( j,q\right) \right)_{i_{\tau }} \left( 2\tau +1  \right)_{i_{\tau -1}} \left( 2\tau + \gamma \right)_{i_{\tau -1}}}{\left( \Delta_{\tau }^{-} \left( j,q\right) \right)_{i_{\tau -1}}\left( \Delta_{\tau }^{+} \left( j,q\right) \right)_{i_{\tau -1}}\left( 2\tau +1 \right)_{i_\tau } \left( 2\tau + \gamma \right)_{i_{\tau }}} \eta ^{i_{\tau }}\right\} z^{\tau } \hspace{1.5cm}\label{eq:50024c} 
\end{eqnarray}
where
\begin{equation}
\begin{cases} \tau \geq 2 \cr
\Delta_k^{\pm} \left( j,q\right) = \frac{\varphi +4(1+a)k \pm \sqrt{\varphi ^2-4(1+a)q}}{2(1+a)}  \cr
\varphi = \beta -\delta -j+a(\gamma +\delta -1)
\end{cases}\nonumber
\end{equation}
\end{enumerate}
\end{remark}
\begin{remark}
The power series expansion of Heun equation of the second kind for the first species complete polynomial using R3TRF about $x=0$ is given by
\begin{enumerate} 
\item As $\alpha = \gamma -1$ and $q=q_0^0= (\gamma -1)(\beta -\delta +a\delta )$,

The eigenfunction is given by
\begin{eqnarray}
y(x) &=& H_pS_{0,0}^R \left( a, q=q_0^0=(\gamma -1)(\beta -\delta +a\delta ); \alpha =\gamma -1, \beta, \gamma, \delta ; \eta = \frac{(1+a)}{a} x ; z= -\frac{1}{a} x^2 \right) \nonumber\\
&=& \eta ^{1-\gamma } \label{eq:50025}
\end{eqnarray}
\item As $\alpha =\gamma -2$,

An algebraic equation of degree 2 for the determination of $q$ is given by
\begin{equation}
0 = a(2-\gamma )(\beta -\gamma +1 )  
+\prod_{l=0}^{1}\Big( q+(1+l-\gamma )(\beta -\delta -1+l+a( \delta + l))\Big) \label{eq:50026a}
\end{equation}
The eigenvalue of $q$ is written by $q_1^m$ where $m = 0,1 $; $q_{1}^0 < q_{1}^1$. Its eigenfunction is given by
\begin{eqnarray}
y(x) &=& H_pS_{1,m}^R \left( a, q=q_1^m; \alpha =\gamma -2, \beta, \gamma, \delta ; \eta = \frac{(1+a)}{a} x ; z= -\frac{1}{a} x^2 \right) \nonumber\\
&=& \eta ^{1-\gamma } \left\{ 1+\frac{(1-\gamma ) \left( \beta -\delta -1+a\delta \right) + q_1^m}{(1+a)(2-\gamma )} \eta \right\} \label{eq:50026b}  
\end{eqnarray}
\item As $\alpha =\gamma -2N-3  $ where $N \in \mathbb{N}_{0}$,

An algebraic equation of degree $2N+3$ for the determination of $q$ is given by
\begin{equation}
0 = \sum_{r=0}^{N+1}\bar{c}\left( r, 2(N-r)+3; 2N+2,q\right)  \label{eq:50027a}
\end{equation}
The eigenvalue of $q$ is written by $q_{2N+2}^m$ where $m = 0,1,2,\cdots,2N+2 $; $q_{2N+2}^0 < q_{2N+2}^1 < \cdots < q_{2N+2}^{2N+2}$. Its eigenfunction is given by 
\begin{eqnarray} 
y(x) &=& H_pS_{2N+2,m}^R \left( a, q=q_{2N+2}^m; \alpha =\gamma -2N-3, \beta, \gamma, \delta ; \eta = \frac{(1+a)}{a} x ; z= -\frac{1}{a} x^2 \right) \nonumber\\
&=&  \sum_{r=0}^{N+1} y_{r}^{2(N+1-r)}\left( 2N+2, q_{2N+2}^m; x \right)  
\label{eq:50027b} 
\end{eqnarray}
\item As $\alpha =\gamma -2N-4 $ where $N \in \mathbb{N}_{0}$,

An algebraic equation of degree $2N+4$ for the determination of $q$ is given by
\begin{equation}  
0 =  \sum_{r=0}^{N+2}\bar{c}\left( r, 2(N+2-r); 2N+3,q\right) \label{eq:50028a}
\end{equation}
The eigenvalue of $q$ is written by $q_{2N+3}^m$ where $m = 0,1,2,\cdots,2N+3 $; $q_{2N+3}^0 < q_{2N+3}^1 < \cdots < q_{2N+3}^{2N+3}$. Its eigenfunction is given by
\begin{eqnarray} 
y(x) &=& H_pS_{2N+3,m}^R \left( a, q=q_{2N+3}^m; \alpha =\gamma -2N-4, \beta, \gamma, \delta ; \eta = \frac{(1+a)}{a} x ; z= -\frac{1}{a} x^2 \right) \nonumber\\
&=&   \sum_{r=0}^{N+1} y_{r}^{2(N-r)+3} \left( 2N+3,q_{2N+3}^m;x\right) \label{eq:50028b}
\end{eqnarray}
In the above,
\begin{eqnarray}
\bar{c}(0,n;j,q)  &=& \frac{\left( \Delta_0^{-} \left( j,q\right) \right)_{n}\left( \Delta_0^{+} \left( j,q\right) \right)_{n}}{\left( 2-\gamma  \right)_{n} \left( 1\right)_{n}} \left( \frac{1+a}{a} \right)^{n}\label{eq:50029a}\\
\bar{c}(1,n;j,q) &=& \left( -\frac{1}{a}\right) \sum_{i_0=0}^{n}\frac{\left( i_0 -j\right)\left( i_0+1+\beta -\gamma \right) }{\left( i_0+3-\gamma  \right) \left( i_0+2 \right)} \frac{ \left( \Delta_0^{-} \left( j,q\right) \right)_{i_0}\left( \Delta_0^{+} \left( j,q\right) \right)_{i_0}}{\left( 2-\gamma  \right)_{i_0} \left( 1\right)_{i_0}} \nonumber\\
&\times&  \frac{ \left( \Delta_1^{-} \left( j,q\right) \right)_{n}\left( \Delta_1^{+} \left( j,q\right) \right)_{n} \left( 4-\gamma  \right)_{i_0} \left( 3\right)_{i_0}}{\left( \Delta_1^{-} \left( j,q\right) \right)_{i_0}\left( \Delta_1^{+} \left( j,q\right) \right)_{i_0}\left( 4-\gamma  \right)_{n} \left( 3 \right)_{n}} \left(\frac{1+a}{a} \right)^{n}  
\label{eq:50029b}\\
\bar{c}(\tau ,n;j,q) &=& \left( -\frac{1}{a}\right)^{\tau} \sum_{i_0=0}^{n}\frac{\left( i_0 -j\right)\left( i_0+\beta -\gamma +1 \right) }{\left( i_0+3-\gamma \right) \left( i_0+2 \right)} \frac{ \left( \Delta_0^{-} \left( j,q\right) \right)_{i_0}\left( \Delta_0^{+} \left( j,q\right) \right)_{i_0}}{\left( 2-\gamma \right)_{i_0} \left( 1 \right)_{i_0}} \nonumber\\
&\times& \prod_{k=1}^{\tau -1} \left( \sum_{i_k = i_{k-1}}^{n} \frac{\left( i_k+ 2k-j\right)\left( i_k +2k+1+\beta -\gamma \right)}{\left( i_k+2k+3-\gamma  \right) \left( i_k+2k+2\right)} \right. \nonumber\\
&\times &  \left. \frac{ \left( \Delta_k^{-} \left( j,q\right) \right)_{i_k}\left( \Delta_k^{+} \left( j,q\right) \right)_{i_k} \left( 2k+2-\gamma \right)_{i_{k-1}} \left( 2k+1 \right)_{i_{k-1}}}{\left( \Delta_k^{-} \left( j,q\right) \right)_{i_{k-1}}\left( \Delta_k^{+} \left( j,q\right) \right)_{i_{k-1}}\left( 2k+2-\gamma \right)_{i_k} \left( 2k+1\right)_{i_k}} \right) \nonumber\\
&\times& \frac{ \left( \Delta_{\tau }^{-} \left( j,q\right) \right)_{n}\left( \Delta_{\tau }^{+} \left( j,q\right) \right)_{n} \left( 2\tau +2-\gamma  \right)_{i_{\tau -1}} \left( 2\tau +1 \right)_{i_{\tau -1}}}{\left( \Delta_{\tau }^{-} \left( j,q\right) \right)_{i_{\tau -1}}\left( \Delta_{\tau }^{+} \left( j,q\right) \right)_{i_{\tau -1}}\left( 2\tau +2-\gamma  \right)_{n} \left( 2\tau +1\right)_{n}} \left(\frac{1+a}{a}\right)^{n }   \hspace{1.5cm}\label{eq:50029c}
\end{eqnarray}
\begin{eqnarray}
y_0^m(j,q;x) &=& \eta ^{1-\gamma }  \sum_{i_0=0}^{m} \frac{\left( \Delta_0^{-} \left( j,q\right) \right)_{i_0}\left( \Delta_0^{+} \left( j,q\right) \right)_{i_0}}{\left( 2-\gamma  \right)_{i_0} \left( 1 \right)_{i_0}} \eta ^{i_0} \label{eq:50030a}\\
y_1^m(j,q;x) &=& \eta ^{1-\gamma } \left\{\sum_{i_0=0}^{m}\frac{\left( i_0 -j\right)\left( i_0+1+\beta -\gamma   \right) }{\left( i_0+3-\gamma \right) \left( i_0+2 \right)} \frac{ \left( \Delta_0^{-} \left( j,q\right) \right)_{i_0}\left( \Delta_0^{+} \left( j,q\right) \right)_{i_0}}{\left( 2-\gamma \right)_{i_0} \left( 1 \right)_{i_0}} \right. \nonumber\\
&\times& \left. \sum_{i_1 = i_0}^{m} \frac{ \left( \Delta_1^{-} \left( j,q\right) \right)_{i_1}\left( \Delta_1^{+} \left( j,q\right) \right)_{i_1} \left( 4-\gamma  \right)_{i_0} \left( 3 \right)_{i_0}}{\left( \Delta_1^{-} \left( j,q\right) \right)_{i_0}\left( \Delta_1^{+} \left( j,q\right) \right)_{i_0}\left( 4-\gamma \right)_{i_1} \left( 3\right)_{i_1}} \eta ^{i_1}\right\} z 
\label{eq:50030b}
\end{eqnarray}
\begin{eqnarray}
y_{\tau }^m(j,q;x) &=& \eta ^{1-\gamma } \left\{ \sum_{i_0=0}^{m} \frac{\left( i_0 -j\right)\left( i_0+1+\beta -\gamma \right) }{\left( i_0+3-\gamma \right) \left( i_0+2 \right)} \frac{ \left( \Delta_0^{-} \left( j,q\right) \right)_{i_0}\left( \Delta_0^{+} \left( j,q\right) \right)_{i_0}}{\left( 2-\gamma \right)_{i_0} \left( 1\right)_{i_0}} \right.\nonumber\\
&\times& \prod_{k=1}^{\tau -1} \left( \sum_{i_k = i_{k-1}}^{m} \frac{\left( i_k+ 2k-j\right)\left( i_k +2k+1+\beta -\gamma \right)}{\left( i_k+2k+3-\gamma  \right) \left( i_k+2k+2\right)} \right. \nonumber\\
&\times &  \left. \frac{ \left( \Delta_k^{-} \left( j,q\right) \right)_{i_k}\left( \Delta_k^{+} \left( j,q\right) \right)_{i_k} \left( 2k+2-\gamma \right)_{i_{k-1}} \left( 2k+1 \right)_{i_{k-1}}}{\left( \Delta_k^{-} \left( j,q\right) \right)_{i_{k-1}}\left( \Delta_k^{+} \left( j,q\right) \right)_{i_{k-1}}\left( 2k+2-\gamma \right)_{i_k} \left( 2k+1\right)_{i_k}} \right) \nonumber\\
&\times & \left. \sum_{i_{\tau } = i_{\tau -1}}^{m}  \frac{ \left( \Delta_{\tau }^{-} \left( j,q\right) \right)_{i_{\tau }}\left( \Delta_{\tau }^{+} \left( j,q\right) \right)_{i_{\tau }} \left( 2\tau +2-\gamma \right)_{i_{\tau -1}} \left( 2\tau +1 \right)_{i_{\tau -1}}}{\left( \Delta_{\tau }^{-} \left( j,q\right) \right)_{i_{\tau -1}}\left( \Delta_{\tau }^{+} \left( j,q\right) \right)_{i_{\tau -1}}\left( 2\tau +2-\gamma \right)_{i_\tau } \left( 2\tau +1 \right)_{i_{\tau }}} \eta ^{i_{\tau }}\right\} z^{\tau } \hspace{1.5cm}\label{eq:50030c} 
\end{eqnarray}
where
\begin{equation}
\begin{cases} \tau \geq 2 \cr
\Delta_k^{\pm} \left( j,q\right) = \frac{\varphi +2(1+a)\left( 2k+1-\gamma \right) \pm \sqrt{\varphi ^2-4(1+a)q}}{2(1+a)}  \cr
\varphi = \beta +\gamma-\delta -1-j+a(\gamma +\delta -1)
\end{cases}\nonumber
\end{equation}
\end{enumerate}
\end{remark}
\subsection{The second species complete polynomial using R3TRF}

For the second species complete polynomial, we need a condition which is defined by
\begin{equation}
B_{j}=B_{j+1}= A_{j}=0\hspace{1cm}\mathrm{where}\;j \in \mathbb{N}_{0}    
 \label{eq:50031}
\end{equation}
\begin{theorem}
In chapter 2, the general expression of a function $y(x)$ for the second species complete polynomial using 3-term recurrence formula is given by
\begin{enumerate} 
\item As $B_1=A_0=0$,
\begin{equation}
y(x) = y_{0}^{0}(x) \label{eq:50032a}
\end{equation}
\item As $B_1=B_2=A_1=0$, 
\begin{equation}
y(x)= y_{0}^{1}(x)  \label{eq:50032b}
\end{equation}
\item As $B_{2N+2}=B_{2N+3}=A_{2N+2}=0$ where $N \in \mathbb{N}_{0}$,
\begin{equation}
y(x)= \sum_{r=0}^{N+1} y_{r}^{2(N+1-r)}(x) \label{eq:50032c}
\end{equation}
\item As $B_{2N+3}=B_{2N+4}=A_{2N+3}=0$ where $N \in \mathbb{N}_{0}$,
\begin{equation}
y(x)= \sum_{r=0}^{N+1} y_{r}^{2(N-r)+3}(x) \label{eq:50032d}
\end{equation}
In the above,
\begin{eqnarray}
y_0^m(x) &=& c_0 x^{\lambda} \sum_{i_0=0}^{m} \left\{ \prod _{i_1=0}^{i_0-1}A_{i_1} \right\} x^{i_0 }\label{eq:50033a}\\
y_1^m(x) &=& c_0 x^{\lambda} \sum_{i_0=0}^{m}\left\{ B_{i_0+1} \prod _{i_1=0}^{i_0-1}A_{i_1}  \sum_{i_2=i_0}^{m} \left\{ \prod _{i_3=i_0}^{i_2-1}A_{i_3+2} \right\}\right\} x^{i_2+2 } \label{eq:50033b}\\
y_{\tau }^m(x) &=& c_0 x^{\lambda} \sum_{i_0=0}^{m} \left\{B_{i_0+1}\prod _{i_1=0}^{i_0-1} A_{i_1} 
\prod _{k=1}^{\tau -1} \left( \sum_{i_{2k}= i_{2(k-1)}}^{m} B_{i_{2k}+(2k+1)}\prod _{i_{2k+1}=i_{2(k-1)}}^{i_{2k}-1}A_{i_{2k+1}+2k}\right) \right. \nonumber\\
&& \times \left. \sum_{i_{2\tau} = i_{2(\tau -1)}}^{m} \left( \prod _{i_{2\tau +1}=i_{2(\tau -1)}}^{i_{2\tau}-1} A_{i_{2\tau +1}+ 2\tau} \right) \right\} x^{i_{2\tau}+2\tau }\hspace{1cm}\mathrm{where}\;\tau \geq 2
\label{eq:50033c}
\end{eqnarray} 
\end{enumerate}
\end{theorem}
Put $n= j+1$ in (\ref{eq:5004b}) and use the condition $B_{j+1}=0$ for $\alpha $.  
\begin{equation}
\alpha = -j-\lambda 
\label{eq:50034}
\end{equation}
Put $n= j$ in (\ref{eq:5004b}) and use the condition $B_{j}=0$ for $\beta $.  
\begin{equation}
\beta = -j+1-\lambda 
\label{eq:50035}
\end{equation}
Substitute (\ref{eq:50034}) and (\ref{eq:50035}) into (\ref{eq:5004a}). Put $n= j$ in the new (\ref{eq:5004a}) and use the condition $A_{j}=0$ for $q$.  
\begin{equation}
q = -(j+\lambda )\left[ -\delta -j+1-\lambda +a(\gamma +\delta +j-1+\lambda )\right]
\label{eq:50036}
\end{equation}
Take (\ref{eq:50034}), (\ref{eq:50035}) and (\ref{eq:50036}) into (\ref{eq:5004a}) and (\ref{eq:5004b}).
\begin{subequations}
\begin{equation}
A_n = \frac{1+a}{a}\frac{(n-j)\left(n+\Pi _0\left( j\right) \right)}{(n+1+\lambda )(n+\gamma +\lambda )} \label{eq:50037a}
\end{equation}
\begin{equation}
B_n = -\frac{1}{a}\frac{(n-j)(n-j-1)}{(n+1+\lambda )(n+\gamma +\lambda )} \label{eq:50037b}
\end{equation}
where
\begin{equation}
\Pi _k\left( j\right) = \frac{1}{1+a}\left( -\delta -j+1+a\left( \gamma +\delta +j-1+2\lambda \right)\right) +2k \label{eq:50037c}
\end{equation}
\end{subequations}
Substitute (\ref{eq:50037a}) and (\ref{eq:50037c}) into (\ref{eq:50033a})--(\ref{eq:50033c}).

As $B_1=A_0=0$, substitute the new (\ref{eq:50033a}) into (\ref{eq:50032a}) putting $j=0$. 

As $B_1=B_2=A_1=0$, substitute the new (\ref{eq:50033a}) into (\ref{eq:50032b}) putting $j=1$. 

As $B_{2N+2}=B_{2N+3}=A_{2N+2}=0$, substitute the new (\ref{eq:50033a})--(\ref{eq:50033c}) into (\ref{eq:50032c}) putting $j=2N+2$.

As $B_{2N+3}=B_{2N+4}=A_{2N+3}=0$, substitute the new (\ref{eq:50033a})--(\ref{eq:50033c}) into (\ref{eq:50032d}) putting $j=2N+3$.

After the replacement process, the general expression of power series of Heun equation about $x=0$ for the second species complete polynomial using reversible 3-term recurrence formula is given by
\begin{enumerate} 
\item As $\alpha =-\lambda $, $\beta =1-\lambda $ and $q= -\lambda \left[ -\delta +1-\lambda +a\left( \gamma +\delta -1+\lambda \right)\right]$,
 
Its eigenfunction is given by
\begin{equation}
y(x) = y_0^0(0;x)= c_0 x^{\lambda } \label{eq:50038a}
\end{equation}
\item As $\alpha =-1-\lambda $, $\beta = -\lambda $ and $q= -(1+\lambda )\left[ -\delta -\lambda +a\left( \gamma +\delta +\lambda \right)\right]$,
 
Its eigenfunction is given by
\begin{equation}
y(x) = y_0^1(1;x)= c_0 x^{\lambda } \left\{ 1+ \frac{\delta -a\left( \gamma +\delta +2\lambda \right)}{(1+a)(1+\lambda )(\gamma +\lambda )}\eta \right\} \label{eq:50038b}
\end{equation}
\item As $\alpha =-2N-2-\lambda $, $\beta = -2N-1-\lambda $ and $q= -\left(2N+2+\lambda \right) \left[ -\delta -2N-1-\lambda +a\left( \gamma +\delta+2N+1+\lambda \right)\right]$ where $N \in \mathbb{N}_{0}$,

Its eigenfunction is given by
\begin{equation}
y(x)=  \sum_{r=0}^{N+1} y_{r}^{2(N+1-r)}\left( 2N+2;x\right)  \label{eq:50038c}
\end{equation}
\item As $\alpha =-2N-3-\lambda $, $\beta = -2N-2-\lambda $ and $q= -\left(2N+3+\lambda \right) \left[ -\delta -2N-2-\lambda +a\left( \gamma +\delta+2N+2+\lambda \right)\right]$ where $N \in \mathbb{N}_{0}$,

Its eigenfunction is given by
\begin{equation}
y(x)= \sum_{r=0}^{N+1} y_{r}^{2(N-r)+3}\left( 2N+3;x\right)  \label{eq:50038d}
\end{equation}
In the above,
\begin{eqnarray}
y_0^m(j;x) &=& c_0 x^{\lambda }  \sum_{i_0=0}^{m} \frac{\left( -j\right)_{i_0} \left( \Pi _0\left( j\right)\right)_{i_0}}{\left( 1+\lambda \right)_{i_0} \left( \gamma + \lambda \right)_{i_0}} \eta ^{i_0} \label{eq:50039a}\\
y_1^m(j;x) &=& c_0 x^{\lambda } \left\{\sum_{i_0=0}^{m} \frac{\left( i_0-j\right)\left( i_0 +1-j \right)}{\left( i_0+2+\lambda \right) \left( i_0+1+\gamma +\lambda \right)} \frac{\left( -j\right)_{i_0} \left( \Pi _0\left( j\right)\right)_{i_0}}{\left( 1+\lambda \right)_{i_0} \left(\gamma + \lambda \right)_{i_0}} \right. \nonumber\\
&\times& \left. \sum_{i_1 = i_0}^{m} \frac{\left( 2-j\right)_{i_1} \left( \Pi _1\left( j\right)\right)_{i_1}\left( 3+ \lambda \right)_{i_0} \left(  \gamma +2+\lambda \right)_{i_0}}{\left( 2-j \right)_{i_0} \left( \Pi _1\left( j\right)\right)_{i_0}\left( 3+\lambda \right)_{i_1} \left(  \gamma +2+\lambda \right)_{i_1}} \eta ^{i_1}\right\} z \label{eq:50039b}\\
y_{\tau }^m(j;x) &=& c_0 x^{\lambda } \left\{ \sum_{i_0=0}^{m} \frac{\left( i_0-j\right)\left( i_0 +1-j \right)}{\left( i_0+2+\lambda \right) \left( i_0+1+\gamma +\lambda \right)} \frac{\left( -j\right)_{i_0} \left( \Pi _0\left( j\right)\right)_{i_0}}{\left( 1+\lambda \right)_{i_0} \left(\gamma + \lambda \right)_{i_0}} \right.\nonumber\\
&\times& \prod_{k=1}^{\tau -1} \left( \sum_{i_k = i_{k-1}}^{m} \frac{\left( i_k+ 2k-j\right)\left( i_k +2k+1-j \right)}{\left( i_k+2k+2+ \lambda \right) \left( i_k+2k+1+ \gamma + \lambda \right)} \right. \nonumber\\
&\times& \left. \frac{\left( 2k-j\right)_{i_k} \left( \Pi _k\left( j\right) \right)_{i_k}\left( 2k+1+\lambda \right)_{i_{k-1}} \left( 2k+\gamma +\lambda \right)_{i_{k-1}}}{\left( 2k-j\right)_{i_{k-1}} \left( \Pi _k\left( j\right) \right)_{i_{k-1}}\left( 2k+1+\lambda \right)_{i_k} \left( 2k+\gamma +\lambda \right)_{i_k}} \right) \nonumber\\
&\times& \left. \sum_{i_{\tau } = i_{\tau -1}}^{m} \frac{\left( 2\tau -j\right)_{i_{\tau }} \left( \Pi _{\tau }\left( j\right)\right)_{i_{\tau }}\left( 2\tau +1+\lambda \right)_{i_{\tau -1}} \left( 2\tau +\gamma +\lambda \right)_{i_{\tau -1}}}{\left( 2\tau -j\right)_{i_{\tau -1}} \left(\Pi _{\tau }\left( j\right)\right)_{i_{\tau -1}}\left( 2\tau +1+\lambda \right)_{i_{\tau }} \left( 2\tau +\gamma +\lambda \right)_{i_{\tau }}} \eta ^{i_{\tau }}\right\} z^{\tau } \hspace{1.5cm}\label{eq:50039c}
\end{eqnarray}
where
\begin{equation}
\begin{cases} \tau \geq 2 \cr
z = -\frac{1}{a}x^2 \cr
\eta = \frac{(1+a)}{a} x \cr
\Pi _k\left( j\right) = \frac{1}{1+a}\left( -\delta -j+1+a\left( \gamma +\delta +j-1+2\lambda \right)\right) +2k  
\end{cases}\nonumber
\end{equation}
\end{enumerate}
Put $c_0$= 1 as $\lambda =0$ for the first kind of independent solutions of Heun equation and $\displaystyle{ c_0= \left( \frac{1+a}{a}\right)^{1-\gamma }}$ as $\lambda = 1-\gamma $ for the second one in (\ref{eq:50038a})--(\ref{eq:50039c}). 
\begin{remark}
The power series expansion of Heun equation of the first kind for the second species complete polynomial using R3TRF about $x=0$ is given by
\begin{enumerate} 
\item As $\alpha =0$, $\beta =1 $ and $q= 0$,
 
Its eigenfunction is given by
\begin{equation}
y(x) = H_pF_0^R \left( a, q=0; \alpha =0, \beta =1, \gamma, \delta ; \eta = \frac{(1+a)}{a} x ; z= -\frac{1}{a} x^2 \right) =1 \label{eq:50040a}
\end{equation}
\item As $\alpha =-1 $, $\beta = 0$ and $q= \delta -a\left( \gamma +\delta \right) $,
 
Its eigenfunction is given by
\begin{eqnarray}
y(x) &=& H_pF_1^R \Bigg( a, q= \delta -a\left( \gamma +\delta \right); \alpha =-1, \beta =0, \gamma, \delta ; \eta = \frac{(1+a)}{a} x ; z= -\frac{1}{a} x^2 \Bigg) \nonumber\\ 
&=&  1+ \frac{\delta -a\left( \gamma +\delta \right)}{(1+a) \gamma }\eta \label{eq:50040b}
\end{eqnarray}
\item As $\alpha =-2N-2 $, $\beta = -2N-1 $ and $q=  \left( 2N+2 \right) \left[  \delta +2N+1 -a\left( \gamma +\delta+2N+1  \right)\right]$ where $N \in \mathbb{N}_{0}$,

Its eigenfunction is given by
\begin{eqnarray}
y(x) &=& H_pF_{2N+2}^R \Bigg( a, q= \left( 2N+2 \right) \left[  \delta +2N+1 -a\left( \gamma +\delta +2N+1  \right)\right]; \alpha =-2N-2 \nonumber\\
&&, \beta =-2N-1, \gamma, \delta; \eta = \frac{(1+a)}{a} x ; z= -\frac{1}{a} x^2 \Bigg) \nonumber\\
&=&  \sum_{r=0}^{N+1} y_{r}^{2(N+1-r)}\left( 2N+2;x\right)  \label{eq:50040c}
\end{eqnarray}
\item As $\alpha =-2N-3 $, $\beta = -2N-2 $ and $q= \left( 2N+3 \right) \left[ \delta +2N+2 -a\left( \gamma +\delta +2N+2 \right)\right]$ where $N \in \mathbb{N}_{0}$,

Its eigenfunction is given by
\begin{eqnarray}
y(x) &=& H_pF_{2N+3}^R \Bigg( a, q= \left( 2N+3 \right) \left[ \delta +2N+2 -a\left( \gamma +\delta +2N+2 \right)\right]; \alpha =-2N-3 \nonumber\\
&&, \beta =-2N-2, \gamma, \delta; \eta = \frac{(1+a)}{a} x ; z= -\frac{1}{a} x^2 \Bigg) \nonumber\\
&=& \sum_{r=0}^{N+1} y_{r}^{2(N-r)+3}\left( 2N+3;x\right)  \label{eq:50040d}
\end{eqnarray}
In the above,
\begin{eqnarray}
y_0^m(j;x) &=&  \sum_{i_0=0}^{m} \frac{\left( -j\right)_{i_0} \left( \Pi _0\left( j\right)\right)_{i_0}}{\left( 1\right)_{i_0} \left( \gamma \right)_{i_0}} \eta ^{i_0} \label{eq:50041a}\\
y_1^m(j;x) &=& \left\{\sum_{i_0=0}^{m} \frac{\left( i_0-j\right)\left( i_0 +1-j \right)}{\left( i_0+2 \right) \left( i_0+1+\gamma \right)} \frac{\left( -j\right)_{i_0} \left( \Pi _0\left( j\right)\right)_{i_0}}{\left( 1 \right)_{i_0} \left(\gamma \right)_{i_0}} \right. \nonumber\\
&\times& \left. \sum_{i_1 = i_0}^{m} \frac{\left( 2-j\right)_{i_1} \left( \Pi _1\left( j\right)\right)_{i_1}\left( 3 \right)_{i_0} \left(  \gamma +2 \right)_{i_0}}{\left( 2-j \right)_{i_0} \left( \Pi _1\left( j\right)\right)_{i_0}\left( 3 \right)_{i_1} \left(  \gamma +2 \right)_{i_1}} \eta ^{i_1}\right\} z \label{eq:50041b}\\
y_{\tau }^m(j;x) &=& \left\{ \sum_{i_0=0}^{m} \frac{\left( i_0-j\right)\left( i_0 +1-j \right)}{\left( i_0+2 \right) \left( i_0+1+\gamma  \right)} \frac{\left( -j\right)_{i_0} \left( \Pi _0\left( j\right)\right)_{i_0}}{\left( 1 \right)_{i_0} \left(\gamma \right)_{i_0}} \right.\nonumber\\
&\times& \prod_{k=1}^{\tau -1} \left( \sum_{i_k = i_{k-1}}^{m} \frac{\left( i_k+ 2k-j\right)\left( i_k +2k+1-j \right)}{\left( i_k+2k+2 \right) \left( i_k+2k+1+ \gamma \right)} \right. \nonumber\\
&\times& \left. \frac{\left( 2k-j\right)_{i_k} \left( \Pi _k\left( j\right) \right)_{i_k}\left( 2k+1 \right)_{i_{k-1}} \left( 2k+\gamma \right)_{i_{k-1}}}{\left( 2k-j\right)_{i_{k-1}} \left( \Pi _k\left( j\right) \right)_{i_{k-1}}\left( 2k+1 \right)_{i_k} \left( 2k+\gamma \right)_{i_k}} \right) \nonumber\\
&\times& \left. \sum_{i_{\tau } = i_{\tau -1}}^{m} \frac{\left( 2\tau -j\right)_{i_{\tau }} \left( \Pi _{\tau }\left( j\right)\right)_{i_{\tau }}\left( 2\tau +1 \right)_{i_{\tau -1}} \left( 2\tau +\gamma \right)_{i_{\tau -1}}}{\left( 2\tau -j\right)_{i_{\tau -1}} \left(\Pi _{\tau }\left( j\right)\right)_{i_{\tau -1}}\left( 2\tau +1 \right)_{i_{\tau }} \left( 2\tau +\gamma \right)_{i_{\tau }}} \eta ^{i_{\tau }}\right\} z^{\tau } \hspace{1.5cm} \label{eq:50041c}
\end{eqnarray}
where
\begin{equation}
\begin{cases} \tau \geq 2 \cr
\Pi _k\left( j\right) = \frac{1}{1+a}\left( -\delta -j+1+a\left( \gamma +\delta +j-1 \right)\right) +2k  
\end{cases}\nonumber
\end{equation}
\end{enumerate}
\end{remark} 
\begin{remark}
The power series expansion of Heun equation of the second kind for the second species complete polynomial using R3TRF about $x=0$ is given by
\begin{enumerate} 
\item As $\alpha =\gamma -1 $, $\beta =\gamma  $ and $q= ( \gamma -1)\left( \gamma-\delta +a \delta \right)$,
 
Its eigenfunction is given by
\begin{eqnarray}
y(x) &=& H_pS_0^R \left( a, q= ( \gamma -1)\left( \gamma-\delta +a \delta \right); \alpha =\gamma -1, \beta =\gamma, \gamma, \delta ; \eta = \frac{(1+a)}{a} x ; z= -\frac{1}{a} x^2 \right)\nonumber\\
 &=& \eta ^{1-\gamma } \label{eq:50042a}
\end{eqnarray}
\item As $\alpha =\gamma -2 $, $\beta = \gamma -1 $ and $q= ( \gamma -2)\left[ \gamma -\delta -1 +a\left( \delta +1 \right)\right]$,
 
Its eigenfunction is given by
\begin{eqnarray}
y(x) &=& H_pS_1^R \Bigg( a, q=  ( \gamma -2)\left[ \gamma -\delta -1 +a\left( \delta +1 \right)\right]; \alpha =\gamma -2, \beta =\gamma -1, \gamma, \delta \nonumber\\
&&; \eta = \frac{(1+a)}{a} x ; z= -\frac{1}{a} x^2 \Bigg) \nonumber\\ 
&=&  \eta ^{1-\gamma } \left\{ 1+ \frac{\delta +a\left( \gamma -\delta -2 \right)}{(1+a)(2-\gamma )}\eta \right\} \label{eq:50042b}
\end{eqnarray}
\item As $\alpha =\gamma -2N-3 $, $\beta = \gamma -2N-2 $ and $q= \left( \gamma -2N-3 \right) \left[ \gamma-\delta -2N-2 +a\left( \delta+2N+2 \right)\right]$ where $N \in \mathbb{N}_{0}$,

Its eigenfunction is given by
\begin{eqnarray}
y(x) &=& H_pS_{2N+2}^R \Bigg( a, q= \left( \gamma -2N-3 \right) \left[ \gamma-\delta -2N-2 +a\left( \delta+2N+2 \right)\right]; \alpha =\gamma -2N-3 \nonumber\\
&&, \beta =\gamma -2N-2, \gamma, \delta; \eta = \frac{(1+a)}{a} x ; z= -\frac{1}{a} x^2 \Bigg) \nonumber\\
&=&  \sum_{r=0}^{N+1} y_{r}^{2(N+1-r)}\left( 2N+2;x\right)  \label{eq:50042c}
\end{eqnarray}
\item As $\alpha =\gamma -2N-4 $, $\beta = \gamma -2N-3 $ and $q= \left(\gamma -2N-4\right) \left[ \gamma -\delta -2N-3  +a\left( \delta+2N+3 \right)\right]$ where $N \in \mathbb{N}_{0}$,

Its eigenfunction is given by
\begin{eqnarray}
y(x) &=& H_pS_{2N+3}^R \Bigg( a, q= \left(\gamma -2N-4\right) \left[ \gamma -\delta -2N-3  +a\left( \delta +2N+3 \right)\right]; \alpha =\gamma -2N-4 \nonumber\\
&&, \beta =\gamma -2N-3, \gamma, \delta; \eta = \frac{(1+a)}{a} x ; z= -\frac{1}{a} x^2 \Bigg) \nonumber\\
&=& \sum_{r=0}^{N+1} y_{r}^{2(N-r)+3}\left( 2N+3;x\right)  \label{eq:50042d}
\end{eqnarray}
In the above,
\begin{eqnarray}
y_0^m(j;x) &=&  \eta ^{1-\gamma }  \sum_{i_0=0}^{m} \frac{\left( -j\right)_{i_0} \left( \Pi _0\left( j\right)\right)_{i_0}}{\left( 2-\gamma \right)_{i_0} \left( 1\right)_{i_0}} \eta ^{i_0} \label{eq:50043a}\\
y_1^m(j;x) &=&  \eta ^{1-\gamma } \left\{\sum_{i_0=0}^{m} \frac{\left( i_0-j\right)\left( i_0 +1-j \right)}{\left( i_0+3-\gamma  \right) \left( i_0+2\right)} \frac{\left( -j\right)_{i_0} \left( \Pi _0\left( j\right)\right)_{i_0}}{\left( 2-\gamma \right)_{i_0} \left( 1\right)_{i_0}} \right. \nonumber\\
&\times& \left. \sum_{i_1 = i_0}^{m} \frac{\left( 2-j\right)_{i_1} \left( \Pi _1\left( j\right)\right)_{i_1}\left( 4-\gamma  \right)_{i_0} \left( 3\right)_{i_0}}{\left( 2-j \right)_{i_0} \left( \Pi _1\left( j\right)\right)_{i_0}\left( 4-\gamma  \right)_{i_1} \left( 3\right)_{i_1}} \eta ^{i_1}\right\} z \label{eq:50043b}\\
y_{\tau }^m(j;x) &=&  \eta ^{1-\gamma } \left\{ \sum_{i_0=0}^{m} \frac{\left( i_0-j\right)\left( i_0 +1-j \right)}{\left( i_0+3-\gamma \right) \left( i_0+2 \right)} \frac{\left( -j\right)_{i_0} \left( \Pi _0\left( j\right)\right)_{i_0}}{\left( 2-\gamma  \right)_{i_0} \left( 1\right)_{i_0}} \right.\nonumber\\
&\times& \prod_{k=1}^{\tau -1} \left( \sum_{i_k = i_{k-1}}^{m} \frac{\left( i_k+ 2k-j\right)\left( i_k +2k+1-j \right)}{\left( i_k+2k+3-\gamma \right) \left( i_k+2k+2 \right)} \right. \nonumber\\
&\times& \left. \frac{\left( 2k-j\right)_{i_k} \left( \Pi _k\left( j\right) \right)_{i_k}\left( 2k+2-\gamma  \right)_{i_{k-1}} \left( 2k+1 \right)_{i_{k-1}}}{\left( 2k-j\right)_{i_{k-1}} \left( \Pi _k\left( j\right) \right)_{i_{k-1}}\left( 2k+2-\gamma \right)_{i_k} \left( 2k+1 \right)_{i_k}} \right) \nonumber\\
&\times& \left. \sum_{i_{\tau } = i_{\tau -1}}^{m} \frac{\left( 2\tau -j\right)_{i_{\tau }} \left( \Pi _{\tau }\left( j\right)\right)_{i_{\tau }}\left( 2\tau +2-\gamma  \right)_{i_{\tau -1}} \left( 2\tau +1 \right)_{i_{\tau -1}}}{\left( 2\tau -j\right)_{i_{\tau -1}} \left(\Pi _{\tau }\left( j\right)\right)_{i_{\tau -1}}\left( 2\tau +2-\gamma  \right)_{i_{\tau }} \left( 2\tau +1 \right)_{i_{\tau }}} \eta ^{i_{\tau }}\right\} z^{\tau } \hspace{1.5cm}\label{eq:50043c}
\end{eqnarray}
where
\begin{equation}
\begin{cases} \tau \geq 2 \cr
\Pi _k\left( j\right) = \frac{1}{1+a}\left( -\delta -j+1+a\left( -\gamma +\delta +j+1\right)\right) +2k  
\end{cases}\nonumber
\end{equation}
\end{enumerate}
\end{remark} 
It is required that $\gamma \ne 0,-1,-2,\cdots$ for the first kind of independent solutions of Heun equation about $x=0$ for the first and second species complete polynomials by applying complete polynomials using either 3TRF or R3TRF. If $\gamma $ is zero or a negative integer, the coefficients $c_n$ in a formal series go to infinity at a certain value of $n$, and a power series solution of (\ref{eq:5002}) can not be analyzed as $\lambda =0$.
 By same reasons, it is required that $\gamma \ne 2,3,4, \cdots$ for the second independent solutions of Heun equation about $x=0$ for the first and second species complete polynomials as $\lambda =1-\gamma $.
\section{Summary} 
  
In chapter 5, I show formal series solutions in closed forms of Heun equation around $x=0$ for the first and second species complete polynomial by applying a mathematical expression of complete polynomials using 3TRF. This is done by letting $A_n$ in sequences $c_n$ is the leading term in each finite sub-power series of the general series solution $y(x)$. 
In this chapter, the first and second species complete polynomials of Heun equation around $x=0$ are constructed by applying a general summation formula of complete polynomials using R3TRF. This is done by letting $B_n$ in sequences $c_n$ is the leading term in each finite sub-power series of the general series solution $y(x)$.

These two complete polynomial solutions are identical to each other analytically. Complete polynomial solutions of Heun equation around $x=0$ by applying 3TRF are composed of the sum of two sub-power series, $y_i^j(x)$ where $i,j \in \mathbb{N}_{0}$. In contrast, complete polynomial expressions of Heun equation by applying R3TRF only consists of one sub-formal series. The latter is more applicable into any special functions analysis for the transformation to other simpler functions such as hypergeometric type functions because of its simple form of formal series solutions. 
 
In chapter 5 and this chapter, two polynomial equations of Heun equation around $x=0$ for the determination of the spectral parameter $q$ in the form of partial sums of the sequences $\{A_n\}$ and $\{B_n\}$ using 3TRF and R3TRF are equivalent to each other analytically.
For an algebraic equation for $q$ of Heun equation using 3TRF, $A_n$ is the leading term in each sequences $c_{2n+l}$ where $l \in \mathbb{N}_{0}$: in chapter 1, I observe the term inside parentheses of sequences $c_n$ which does not include any $A_n$'s for $c_{2n}$ with every subscripts ($c_0$, $c_1$, $c_2$,$\cdots$), I observe the terms inside parentheses of sequence $c_n$ which include one term of $A_n$'s for $c_{2n+1}$ with odd subscripts ($c_1$, $c_3$, $c_5$,$\cdots$), I observe the terms inside parentheses of sequence $c_n$ which include two terms of $A_n$'s for $c_{2n+2}$ with even subscripts ($c_2$, $c_4$, $c_6$,$\cdots$), I observe the terms inside parentheses of sequence $c_n$ which include three terms of $A_n$'s for $c_{2n+3}$ with odd subscripts ($c_3$, $c_5$, $c_7$,$\cdots$), etc.   

For a polynomial equation for $q$ of Heun equation using R3TRF, $B_n$ is the leading term in each sequences $c_{n+2l}$: in chapter 2, I observe the term inside parentheses of sequences $c_n$ which does not include any $B_n$'s for $c_{n}$ with every subscripts ($c_0$, $c_1$, $c_2$,$\cdots$), I observe the terms inside parentheses of sequence $c_n$ which include one term of $B_n$'s for $c_{n+2}$ with every index except $c_0$ and $c_1$ ($c_2$, $c_3$, $c_4$,$\cdots$), I observe the terms inside parentheses of sequence $c_n$ which include two terms of $B_n$'s for $c_{n+4}$ with every index except $c_0$--$c_3$ ($c_4$, $c_5$, $c_6$,$\cdots$), I observe the terms inside parentheses of sequence $c_n$ which include three terms of $B_n$'s for $c_{n+6}$ with every index except $c_0$--$c_5$ ($c_6$, $c_7$, $c_8$,$\cdots$), etc.

For the first and second complete polynomials of Heun equation around $x=0$ by applying 3TRF in chapter 5, the denominators and numerators in all $B_n$ terms of each finite sub-power series arise with Pochhammer symbol. For the first and second complete polynomials of Heun equation around $x=0$ by applying R3TRF in this chapter, the denominators and numerators in all $A_n$ terms of each finite sub-power series arise with Pochhammer symbol.
And hypergeometric type equations having the 2-term recurrence relation are composed of a ratio of a Pochhammer symbol in numerator to another Pochhammer symbol in denominator in sequence $c_n$. 

In chapter 2 of Ref.\cite{5Choun2013} and \cite{5Chou2012d}, combined definite and contour integrals for polynomials of type 1 and 2 of Heun equation around $x=0$ are constructed by applying an integral representation of a generalized hypergeometric polynomial $_5F_4$; it is derived from power series solutions in closed forms of Heun equation.
By using similar methods, integral (differential) representations of Heun equation for complete polynomials of two types will be constructed analytically in the future series including its generating function and the orthogonal relation of it. 

\begin{appendices}
\section{Power series expansion of 192 Heun functions}
A machine-generated list of 192 (isomorphic to the Coxeter group of the Coxeter diagram $D_4$) local solutions of the Heun equation was obtained by Robert S. Maier(2007) \cite{5Maie2007}. 
In this appendix, replacing coefficients in Frobenius solutions of Heun equation around $x=0$ of the first kind for the first and second complete polynomials using R3TRF, I construct power series solutions for the first and second complete polynomials of nine out of the 192 local solutions of Heun equation in Table 2 \cite{5Maie2007}
\addtocontents{toc}{\protect\setcounter{tocdepth}{1}}
\subsection{${\displaystyle (1-x)^{1-\delta } Hl(a, q - (\delta  - 1)\gamma a; \alpha - \delta  + 1, \beta - \delta + 1, \gamma ,2 - \delta ; x)}$ }
\subsubsection{The first species complete polynomial}
Replacing coefficients $q$, $\alpha$, $\beta$ and $\delta$ by $q - (\delta - 1)\gamma a $, $\alpha - \delta  + 1 $, $\beta - \delta + 1$ and $2 - \delta$ into (\ref{eq:50019})--(\ref{eq:50024c}). Multiply $(1-x)^{1-\delta }$ and the new (\ref{eq:50019}), (\ref{eq:50020b}), (\ref{eq:50021b}) and (\ref{eq:50022b})  together.\footnote{I treat $\beta $, $\gamma$ and $\delta$ as free variables and fixed values of $\alpha $ and $q$.} 
\begin{enumerate} 
\item As $\alpha  = \delta  -1$ and $q = (\delta - 1)\gamma a+q_0^0$ where $q_0^0=0$,

The eigenfunction is given by
\begin{eqnarray}
& &(1-x)^{1-\delta } y(x)\nonumber\\
&=& (1-x)^{1-\delta } Hl\left( a, 0; 0, \beta - \delta + 1, \gamma ,2 - \delta ; x\right)\nonumber\\
&=& (1-x)^{1-\delta } \nonumber
\end{eqnarray}
\item As $\alpha =\delta -2$,

An algebraic equation of degree 2 for the determination of $q$ is given by
\begin{equation}
0 = a\gamma (\beta -\delta +1)  
+\prod_{l=0}^{1}\Big( q-(\delta -1)\gamma a+ l(\beta +l+a(\gamma -\delta +1+l))\Big) \nonumber
\end{equation}
The eigenvalue of $q$ is written by $(\delta -1)\gamma a +q_1^m$ where $m = 0,1 $; $q_{1}^0 < q_{1}^1$. Its eigenfunction is given by
\begin{eqnarray}
& &(1-x)^{1-\delta } y(x)\nonumber\\
&=& (1-x)^{1-\delta } Hl\left( a, q_1^m; -1, \beta - \delta + 1, \gamma ,2 - \delta ; x\right)\nonumber\\
&=& (1-x)^{1-\delta } \left\{ 1+\frac{ q_1^m}{(1+a)\gamma } \eta\right\} \nonumber
\end{eqnarray}
 \item As $\alpha =\delta -2N-3 $ where $N \in \mathbb{N}_{0}$,

An algebraic equation of degree $2N+3$ for the determination of $q$ is given by
\begin{equation}
0 = \sum_{r=0}^{N+1}\bar{c}\left( r, 2(N-r)+3; 2N+2,\tilde{q}\right)  \nonumber
\end{equation}
The eigenvalue of $q$ is written by $(\delta  - 1)\gamma a+ q_{2N+2}^m$ where $m = 0,1,2,\cdots,2N+2 $; $q_{2N+2}^0 < q_{2N+2}^1 < \cdots < q_{2N+2}^{2N+2}$. Its eigenfunction is given by 
\begin{eqnarray} 
& &(1-x)^{1-\delta } y(x)\nonumber\\
&=& (1-x)^{1-\delta } Hl\left( a, q_{2N+2}^m; -2N-2, \beta - \delta + 1, \gamma ,2 - \delta ; x\right)\nonumber\\
&=& (1-x)^{1-\delta } \sum_{r=0}^{N+1} y_{r}^{2(N+1-r)}\left( 2N+2, q_{2N+2}^m; x \right) 
\nonumber 
\end{eqnarray}
\item As $\alpha =\delta -2N-4 $ where $N \in \mathbb{N}_{0}$,

An algebraic equation of degree $2N+4$ for the determination of $q$ is given by
\begin{equation}  
0 = \sum_{r=0}^{N+2}\bar{c}\left( r, 2(N+2-r); 2N+3,\tilde{q}\right) \nonumber
\end{equation}
The eigenvalue of $q$ is written by $(\delta  - 1)\gamma a+q_{2N+3}^m$ where $m = 0,1,2,\cdots,2N+3 $; $q_{2N+3}^0 < q_{2N+3}^1 < \cdots < q_{2N+3}^{2N+3}$. Its eigenfunction is given by
\begin{eqnarray} 
& &(1-x)^{1-\delta } y(x)\nonumber\\
&=& (1-x)^{1-\delta } Hl\left( a, q=q_{2N+3}^m; -2N-3, \beta - \delta + 1, \gamma ,2 - \delta ; x\right)\nonumber\\
&=& (1-x)^{1-\delta } \sum_{r=0}^{N+1} y_{r}^{2(N-r)+3} \left( 2N+3,q_{2N+3}^m;x\right)  \nonumber
\end{eqnarray}
In the above,
\begin{eqnarray}
\bar{c}(0,n;j,\tilde{q})  &=& \frac{\left( \Delta_0^{-} \left( j,\tilde{q}\right) \right)_{n}\left( \Delta_0^{+} \left( j,\tilde{q}\right) \right)_{n}}{\left( 1 \right)_{n} \left( \gamma \right)_{n}} \left( \frac{1+a}{a} \right)^{n}\nonumber\\
\bar{c}(1,n;j,\tilde{q}) &=& \left( -\frac{1}{a}\right) \sum_{i_0=0}^{n}\frac{\left( i_0 -j\right)\left( i_0+1+\beta -\delta \right) }{\left( i_0+2 \right) \left( i_0+1+ \gamma \right)} \frac{ \left( \Delta_0^{-} \left( j,\tilde{q}\right) \right)_{i_0}\left( \Delta_0^{+} \left( j,\tilde{q}\right) \right)_{i_0}}{\left( 1 \right)_{i_0} \left( \gamma \right)_{i_0}} \nonumber\\
&&\times  \frac{ \left( \Delta_1^{-} \left( j,\tilde{q}\right) \right)_{n}\left( \Delta_1^{+} \left( j,\tilde{q}\right) \right)_{n} \left( 3 \right)_{i_0} \left( 2+ \gamma \right)_{i_0}}{\left( \Delta_1^{-} \left( j,\tilde{q}\right) \right)_{i_0}\left( \Delta_1^{+} \left( j,\tilde{q}\right) \right)_{i_0}\left( 3 \right)_{n} \left( 2+ \gamma \right)_{n}} \left(\frac{1+a}{a} \right)^{n }  
\nonumber\\
\bar{c}(\tau ,n;j,\tilde{q}) &=& \left( -\frac{1}{a}\right)^{\tau} \sum_{i_0=0}^{n}\frac{\left( i_0 -j\right)\left( i_0+1+\beta -\delta \right) }{\left( i_0+2 \right) \left( i_0+1+ \gamma \right)} \frac{ \left( \Delta_0^{-} \left( j,\tilde{q}\right) \right)_{i_0}\left( \Delta_0^{+} \left( j,\tilde{q}\right) \right)_{i_0}}{\left( 1 \right)_{i_0} \left( \gamma \right)_{i_0}}   \nonumber\\
&&\times  \prod_{k=1}^{\tau -1} \left( \sum_{i_k = i_{k-1}}^{n} \frac{\left( i_k+ 2k-j\right)\left( i_k +2k+1+\beta -\delta \right)}{\left( i_k+2k+2 \right) \left( i_k+2k+1+ \gamma \right)} \right. \nonumber\\
&&\times  \left. \frac{ \left( \Delta_k^{-} \left( j,\tilde{q}\right) \right)_{i_k}\left( \Delta_k^{+} \left( j,\tilde{q}\right) \right)_{i_k} \left( 2k+1 \right)_{i_{k-1}} \left( 2k+ \gamma \right)_{i_{k-1}}}{\left( \Delta_k^{-} \left( j,\tilde{q}\right) \right)_{i_{k-1}}\left( \Delta_k^{+} \left( j,\tilde{q}\right) \right)_{i_{k-1}}\left( 2k+1 \right)_{i_k} \left( 2k+ \gamma \right)_{i_k}} \right) \nonumber\\
&&\times \frac{ \left( \Delta_{\tau }^{-} \left( j,\tilde{q}\right) \right)_{n}\left( \Delta_{\tau }^{+} \left( j,\tilde{q}\right) \right)_{n} \left( 2\tau +1 \right)_{i_{\tau -1}} \left( 2\tau + \gamma \right)_{i_{\tau -1}}}{\left( \Delta_{\tau }^{-} \left( j,\tilde{q}\right) \right)_{i_{\tau -1}}\left( \Delta_{\tau }^{+} \left( j,\tilde{q}\right) \right)_{i_{\tau -1}}\left( 2\tau +1 \right)_{n} \left( 2\tau + \gamma  \right)_{n}} \left(\frac{1+a}{a}\right)^{n } \nonumber 
\end{eqnarray}
\begin{eqnarray}
y_0^m(j,\tilde{q};x) &=& \sum_{i_0=0}^{m} \frac{\left( \Delta_0^{-} \left( j,\tilde{q}\right) \right)_{i_0}\left( \Delta_0^{+} \left( j,\tilde{q}\right) \right)_{i_0}}{\left( 1 \right)_{i_0} \left( \gamma \right)_{i_0}} \eta ^{i_0} \nonumber\\
y_1^m(j,\tilde{q};x) &=& \left\{\sum_{i_0=0}^{m}\frac{\left( i_0 -j\right)\left( i_0+1+\beta -\delta \right) }{\left( i_0+2 \right) \left( i_0+1+ \gamma \right)} \frac{ \left( \Delta_0^{-} \left( j,\tilde{q}\right) \right)_{i_0}\left( \Delta_0^{+} \left( j,\tilde{q}\right) \right)_{i_0}}{\left( 1 \right)_{i_0} \left( \gamma \right)_{i_0}} \right. \nonumber\\
&&\times \left. \sum_{i_1 = i_0}^{m} \frac{ \left( \Delta_1^{-} \left( j,\tilde{q}\right) \right)_{i_1}\left( \Delta_1^{+} \left( j,\tilde{q}\right) \right)_{i_1} \left( 3 \right)_{i_0} \left( 2+ \gamma \right)_{i_0}}{\left( \Delta_1^{-} \left( j,\tilde{q}\right) \right)_{i_0}\left( \Delta_1^{+} \left( j,\tilde{q}\right) \right)_{i_0}\left( 3 \right)_{i_1} \left( 2+ \gamma \right)_{i_1}} \eta ^{i_1}\right\} z 
\nonumber
\end{eqnarray}
\begin{eqnarray}
y_{\tau }^m(j,\tilde{q};x) &=& \left\{ \sum_{i_0=0}^{m} \frac{\left( i_0 -j\right)\left( i_0+1+\beta -\delta \right) }{\left( i_0+2 \right) \left( i_0+1+ \gamma \right)} \frac{ \left( \Delta_0^{-} \left( j,\tilde{q}\right) \right)_{i_0}\left( \Delta_0^{+} \left( j,\tilde{q}\right) \right)_{i_0}}{\left( 1 \right)_{i_0} \left( \gamma \right)_{i_0}} \right.\nonumber\\
&&\times \prod_{k=1}^{\tau -1} \left( \sum_{i_k = i_{k-1}}^{m} \frac{\left( i_k+ 2k-j\right)\left( i_k +2k+1+\beta -\delta \right)}{\left( i_k+2k+2 \right) \left( i_k+2k+1+ \gamma \right)} \right. \nonumber\\
&&\times  \left. \frac{ \left( \Delta_k^{-} \left( j,\tilde{q}\right) \right)_{i_k}\left( \Delta_k^{+} \left( j,\tilde{q}\right) \right)_{i_k} \left( 2k+1 \right)_{i_{k-1}} \left( 2k+ \gamma \right)_{i_{k-1}}}{\left( \Delta_k^{-} \left( j,\tilde{q}\right) \right)_{i_{k-1}}\left( \Delta_k^{+} \left( j,\tilde{q}\right) \right)_{i_{k-1}}\left( 2k+1 \right)_{i_k} \left( 2k+ \gamma \right)_{i_k}} \right) \nonumber\\
&&\times \left. \sum_{i_{\tau } = i_{\tau -1}}^{m}  \frac{ \left( \Delta_{\tau }^{-} \left( j,\tilde{q}\right) \right)_{i_{\tau }}\left( \Delta_{\tau }^{+} \left( j,\tilde{q}\right) \right)_{i_{\tau }} \left( 2\tau +1  \right)_{i_{\tau -1}} \left( 2\tau + \gamma \right)_{i_{\tau -1}}}{\left( \Delta_{\tau }^{-} \left( j,\tilde{q}\right) \right)_{i_{\tau -1}}\left( \Delta_{\tau }^{+} \left( j,\tilde{q}\right) \right)_{i_{\tau -1}}\left( 2\tau +1 \right)_{i_\tau } \left( 2\tau + \gamma \right)_{i_{\tau }}} \eta ^{i_{\tau }}\right\} z^{\tau } \nonumber
\end{eqnarray}
where
\begin{equation}
\begin{cases} \tau \geq 2 \cr
z = -\frac{1}{a}x^2 \cr
\eta = \frac{(1+a)}{a} x \cr
\tilde{q} = q-(\delta -1)\gamma a \cr
\Delta_k^{\pm} \left( j,\tilde{q}\right) = \frac{\varphi +4(1+a)k \pm \sqrt{\varphi ^2-4(1+a)\tilde{q}}}{2(1+a)}  \cr
\varphi = \beta -1-j+a(\gamma -\delta +1)
\end{cases}\nonumber
\end{equation}
\end{enumerate}
\subsubsection{The second species complete polynomial}
Replacing coefficients $q$, $\alpha$, $\beta$ and $\delta$ by $q - (\delta - 1)\gamma a $, $\alpha - \delta  + 1 $, $\beta - \delta + 1$ and $2 - \delta$ into (\ref{eq:50040a})--(\ref{eq:50041c}). Multiply $(1-x)^{1-\delta }$ and the new (\ref{eq:50040a})--(\ref{eq:50041c}) together.\footnote{I treat $\gamma$ and $\delta$ as free variables and fixed values of $\alpha $, $\beta $ and $q$.}  
\begin{enumerate} 
\item As $\alpha  = \delta -1 $, $\beta   =\delta $ and $q = (\delta  - 1)\gamma a$,
 
Its eigenfunction is given by
\begin{eqnarray}
& &(1-x)^{1-\delta } y(x)\nonumber\\
&=& (1-x)^{1-\delta } Hl\left( a, 0; 0, 1, \gamma ,2 - \delta ; x\right)\nonumber\\
&=& (1-x)^{1-\delta }  \nonumber
\end{eqnarray}
\item As $\alpha = \delta -2 $, $\beta =\delta -1$ and $q= (\delta -1)\gamma a -\delta +2-a\left( \gamma -\delta +2\right) $,
 
Its eigenfunction is given by
\begin{eqnarray}
& &(1-x)^{1-\delta } y(x)\nonumber\\
&=& (1-x)^{1-\delta } Hl\left( a, -\delta +2-a\left( \gamma -\delta +2\right); -1, 0, \gamma ,2 - \delta ; x\right)\nonumber\\
&=& (1-x)^{1-\delta } \left\{  1+ \frac{-\delta +2-a\left( \gamma -\delta +2\right)}{(1+a) \gamma }\eta 
\right\}  \nonumber
\end{eqnarray}
\item As $\alpha =\delta -2N-3 $, $\beta = \delta -2N-2 $ and $q= (\delta -1)\gamma a + \left( 2N+2 \right) \left[  -\delta +2N+3 -a\left( \gamma -\delta +2N+3  \right)\right]$ where $N \in \mathbb{N}_{0}$,

Its eigenfunction is given by
\begin{eqnarray}
& &(1-x)^{1-\delta } y(x)\nonumber\\
&=& (1-x)^{1-\delta } Hl\left( a, \left( 2N+2 \right) \left[  -\delta +2N+3 -a\left( \gamma -\delta+2N+3  \right)\right] ; -2N-2, -2N-1 \right. \nonumber\\
&&, \left. \gamma ,2 - \delta ; x\right)\nonumber\\
&=&  (1-x)^{1-\delta }  \sum_{r=0}^{N+1} y_{r}^{2(N+1-r)}\left( 2N+2;x\right)  \nonumber
\end{eqnarray}
\item As $\alpha =\delta -2N-4 $, $\beta = \delta -2N-3 $ and $q= (\delta -1)\gamma a +  \left( 2N+3 \right) \left[ -\delta +2N+4 -a\left( \gamma -\delta +2N+4 \right)\right]$ where $N \in \mathbb{N}_{0}$,

Its eigenfunction is given by
\begin{eqnarray}
& &(1-x)^{1-\delta } y(x)\nonumber\\
&=& (1-x)^{1-\delta } Hl\left( a, \left( 2N+3 \right) \left[ -\delta +2N+4 -a\left( \gamma -\delta +2N+4 \right)\right] ; -2N-3, -2N-2, \gamma ,2 - \delta ; x\right)\nonumber\\
&=&  (1-x)^{1-\delta } \sum_{r=0}^{N+1} y_{r}^{2(N-r)+3}\left( 2N+3;x\right)  \nonumber
\end{eqnarray}
In the above,
\begin{eqnarray}
y_0^m(j;x) &=&  \sum_{i_0=0}^{m} \frac{\left( -j\right)_{i_0} \left( \Pi _0\left( j\right)\right)_{i_0}}{\left( 1\right)_{i_0} \left( \gamma \right)_{i_0}} \eta ^{i_0} \nonumber\\
y_1^m(j;x) &=& \left\{\sum_{i_0=0}^{m} \frac{\left( i_0-j\right)\left( i_0 +1-j \right)}{\left( i_0+2 \right) \left( i_0+1+\gamma \right)} \frac{\left( -j\right)_{i_0} \left( \Pi _0\left( j\right)\right)_{i_0}}{\left( 1 \right)_{i_0} \left(\gamma \right)_{i_0}} \right. \nonumber\\
&&\times \left. \sum_{i_1 = i_0}^{m} \frac{\left( 2-j\right)_{i_1} \left( \Pi _1\left( j\right)\right)_{i_1}\left( 3 \right)_{i_0} \left(  \gamma +2 \right)_{i_0}}{\left( 2-j \right)_{i_0} \left( \Pi _1\left( j\right)\right)_{i_0}\left( 3 \right)_{i_1} \left(  \gamma +2 \right)_{i_1}} \eta ^{i_1}\right\} z \nonumber\\
y_{\tau }^m(j;x) &=& \left\{ \sum_{i_0=0}^{m} \frac{\left( i_0-j\right)\left( i_0 +1-j \right)}{\left( i_0+2 \right) \left( i_0+1+\gamma  \right)} \frac{\left( -j\right)_{i_0} \left( \Pi _0\left( j\right)\right)_{i_0}}{\left( 1 \right)_{i_0} \left(\gamma \right)_{i_0}} \right.\nonumber\\
&&\times \prod_{k=1}^{\tau -1} \left( \sum_{i_k = i_{k-1}}^{m} \frac{\left( i_k+ 2k-j\right)\left( i_k +2k+1-j \right)}{\left( i_k+2k+2 \right) \left( i_k+2k+1+ \gamma \right)} \right. \nonumber\\
&&\times \left. \frac{\left( 2k-j\right)_{i_k} \left( \Pi _k\left( j\right) \right)_{i_k}\left( 2k+1 \right)_{i_{k-1}} \left( 2k+\gamma \right)_{i_{k-1}}}{\left( 2k-j\right)_{i_{k-1}} \left( \Pi _k\left( j\right) \right)_{i_{k-1}}\left( 2k+1 \right)_{i_k} \left( 2k+\gamma \right)_{i_k}} \right) \nonumber\\
&&\times \left. \sum_{i_{\tau } = i_{\tau -1}}^{m} \frac{\left( 2\tau -j\right)_{i_{\tau }} \left( \Pi _{\tau }\left( j\right)\right)_{i_{\tau }}\left( 2\tau +1 \right)_{i_{\tau -1}} \left( 2\tau +\gamma \right)_{i_{\tau -1}}}{\left( 2\tau -j\right)_{i_{\tau -1}} \left(\Pi _{\tau }\left( j\right)\right)_{i_{\tau -1}}\left( 2\tau +1 \right)_{i_{\tau }} \left( 2\tau +\gamma \right)_{i_{\tau }}} \eta ^{i_{\tau }}\right\} z^{\tau } \nonumber
\end{eqnarray}
where
\begin{equation}
\begin{cases} \tau \geq 2 \cr
z = -\frac{1}{a}x^2 \cr
\eta = \frac{(1+a)}{a} x \cr
\Pi _k\left( j\right) = \frac{1}{1+a}\left(  \delta -1-j +a\left( \gamma -\delta +1+j  \right)\right) +2k  
\end{cases}\nonumber
\end{equation}
\end{enumerate}
\subsection{\footnotesize ${\displaystyle x^{1-\gamma } (1-x)^{1-\delta }Hl(a, q-(\gamma +\delta -2)a -(\gamma -1)(\alpha +\beta -\gamma -\delta +1), \alpha - \gamma -\delta +2, \beta - \gamma -\delta +2, 2-\gamma, 2 - \delta ; x)}$ \normalsize}
\subsubsection{The first species complete polynomial}
Replacing coefficients $q$, $\alpha$, $\beta$, $\gamma $ and $\delta$ by $q-(\gamma +\delta -2)a-(\gamma -1)(\alpha +\beta -\gamma -\delta +1)$, $\alpha - \gamma -\delta +2$, $\beta - \gamma -\delta +2, 2-\gamma$ and $2 - \delta$ into (\ref{eq:50019})--(\ref{eq:50024c}). Multiply $x^{1-\gamma } (1-x)^{1-\delta }$ and the new (\ref{eq:50019}), (\ref{eq:50020b}), (\ref{eq:50021b}) and (\ref{eq:50022b})  together.\footnote{I treat $\beta $, $\gamma$ and $\delta$ as free variables and fixed values of $\alpha $ and $q$.} 
\begin{enumerate} 
\item As $\alpha = \gamma +\delta -2 $ and $q=(\gamma +\delta -2)a+(\gamma -1)( \beta -1)+q_0^0$ where $q_0^0=0$,

The eigenfunction is given by
\begin{eqnarray}
& &x^{1-\gamma } (1-x)^{1-\delta } y(x)\nonumber\\
&=& x^{1-\gamma } (1-x)^{1-\delta } Hl(a, 0; 0, \beta - \gamma -\delta +2, 2-\gamma, 2 - \delta ; x)\nonumber\\
&=& x^{1-\gamma } (1-x)^{1-\delta } \nonumber
\end{eqnarray}
\item As $\alpha =\gamma +\delta -3$,

An algebraic equation of degree 2 for the determination of $q$ is given by
\begin{equation}
0 = a (\beta -\gamma -\delta +2) (2-\gamma )
+\prod_{l=0}^{1}\Big( q-(\gamma +\delta -2)a -(\gamma -1)(\beta -2)+ l(\beta -\gamma -1+l-a( \gamma +\delta -3-l))\Big) \nonumber
\end{equation}
The eigenvalue of $q$ is written by $(\gamma +\delta -2)a +(\gamma -1)(\beta -2) +q_1^m$ where $m = 0,1 $; $q_{1}^0 < q_{1}^1$. Its eigenfunction is given by
\begin{eqnarray}
& &x^{1-\gamma } (1-x)^{1-\delta } y(x)\nonumber\\
&=& x^{1-\gamma } (1-x)^{1-\delta } Hl(a, q_1^m; -1, \beta - \gamma -\delta +2, 2-\gamma, 2 - \delta ; x)\nonumber\\
&=& x^{1-\gamma } (1-x)^{1-\delta } \left\{  1+\frac{ q_1^m}{(1+a)(2-\gamma )} \eta \right\} \nonumber
\end{eqnarray}
\item As $\alpha =\gamma +\delta -2N-4 $ where $N \in \mathbb{N}_{0}$,

An algebraic equation of degree $2N+3$ for the determination of $q$ is given by
\begin{equation}
0 = \sum_{r=0}^{N+1}\bar{c}\left( r, 2(N-r)+3; 2N+2,q-(\gamma +\delta -2)a -(\gamma -1)(\beta -2N-3)\right)  \nonumber
\end{equation}
The eigenvalue of $q$ is written by $(\gamma +\delta -2)a +(\gamma -1)(\beta -2N-3) +q_{2N+2}^m$ where $m = 0,1,2,\cdots,2N+2 $; $q_{2N+2}^0 < q_{2N+2}^1 < \cdots < q_{2N+2}^{2N+2}$. Its eigenfunction is given by 
\begin{eqnarray} 
& &x^{1-\gamma } (1-x)^{1-\delta } y(x)\nonumber\\
&=& x^{1-\gamma } (1-x)^{1-\delta } Hl(a, q_{2N+2}^m; -2N-2, \beta - \gamma -\delta +2, 2-\gamma, 2 - \delta ; x)\nonumber\\
&=& x^{1-\gamma } (1-x)^{1-\delta } \sum_{r=0}^{N+1} y_{r}^{2(N+1-r)}\left( 2N+2, q_{2N+2}^m; x \right)  
\nonumber
\end{eqnarray} 
\item As $\alpha =\gamma +\delta -2N-5 $ where $N \in \mathbb{N}_{0}$,

An algebraic equation of degree $2N+4$ for the determination of $q$ is given by
\begin{equation}  
0 = \sum_{r=0}^{N+2}\bar{c}\left( r, 2(N+2-r); 2N+3,q-(\gamma +\delta -2)a -(\gamma -1)(\beta -2N-4)\right) \nonumber
\end{equation}
The eigenvalue of $q$ is written by $(\gamma +\delta -2)a +(\gamma -1)(\beta -2N-4)+ q_{2N+3}^m$ where $m = 0,1,2,\cdots,2N+3 $; $q_{2N+3}^0 < q_{2N+3}^1 < \cdots < q_{2N+3}^{2N+3}$. Its eigenfunction is given by
\begin{eqnarray} 
& &x^{1-\gamma } (1-x)^{1-\delta } y(x)\nonumber\\
&=& x^{1-\gamma } (1-x)^{1-\delta } Hl(a, q_{2N+3}^m; -2N-3, \beta - \gamma -\delta +2, 2-\gamma, 2 - \delta ; x)\nonumber\\
&=& x^{1-\gamma } (1-x)^{1-\delta } \sum_{r=0}^{N+1} y_{r}^{2(N-r)+3} \left( 2N+3,q_{2N+3}^m;x\right) \nonumber
\end{eqnarray}
In the above,
\begin{eqnarray}
\bar{c}(0,n;j,\tilde{q})  &=& \frac{\left( \Delta_0^{-} \left( j,\tilde{q}\right) \right)_{n}\left( \Delta_0^{+} \left( j,\tilde{q}\right) \right)_{n}}{\left( 1 \right)_{n} \left( 2-\gamma \right)_{n}} \left( \frac{1+a}{a} \right)^{n}\nonumber\\
\bar{c}(1,n;j,\tilde{q}) &=& \left( -\frac{1}{a}\right) \sum_{i_0=0}^{n}\frac{\left( i_0 -j\right)\left( i_0 +2+\beta -\gamma -\delta \right) }{\left( i_0+2 \right) \left( i_0+3-\gamma \right)} \frac{ \left( \Delta_0^{-} \left( j,\tilde{q}\right) \right)_{i_0}\left( \Delta_0^{+} \left( j,\tilde{q}\right) \right)_{i_0}}{\left( 1 \right)_{i_0} \left( 2- \gamma \right)_{i_0}} \nonumber\\
&&\times  \frac{ \left( \Delta_1^{-} \left( j,\tilde{q}\right) \right)_{n}\left( \Delta_1^{+} \left( j,\tilde{q}\right) \right)_{n} \left( 3 \right)_{i_0} \left( 4- \gamma \right)_{i_0}}{\left( \Delta_1^{-} \left( j,\tilde{q}\right) \right)_{i_0}\left( \Delta_1^{+} \left( j,\tilde{q}\right) \right)_{i_0}\left( 3 \right)_{n} \left( 4-\gamma \right)_{n}} \left(\frac{1+a}{a} \right)^{n }  
\nonumber\\
\bar{c}(\tau ,n;j,\tilde{q}) &=& \left( -\frac{1}{a}\right)^{\tau} \sum_{i_0=0}^{n}\frac{\left( i_0 -j\right)\left( i_0+2+\beta -\gamma -\delta \right) }{\left( i_0+2 \right) \left( i_0+3- \gamma \right)} \frac{ \left( \Delta_0^{-} \left( j,\tilde{q}\right) \right)_{i_0}\left( \Delta_0^{+} \left( j,\tilde{q}\right) \right)_{i_0}}{\left( 1 \right)_{i_0} \left( 2-\gamma \right)_{i_0}} \nonumber\\
&&\times \prod_{k=1}^{\tau -1} \left( \sum_{i_k = i_{k-1}}^{n} \frac{\left( i_k+ 2k-j\right)\left( i_k +2k+2+\beta -\gamma -\delta \right)}{\left( i_k+2k+2 \right) \left( i_k+2k+3- \gamma \right)} \right. \nonumber\\
&&\times  \left. \frac{ \left( \Delta_k^{-} \left( j,\tilde{q}\right) \right)_{i_k}\left( \Delta_k^{+} \left( j,\tilde{q}\right) \right)_{i_k} \left( 2k+1 \right)_{i_{k-1}} \left( 2k+2- \gamma \right)_{i_{k-1}}}{\left( \Delta_k^{-} \left( j,\tilde{q}\right) \right)_{i_{k-1}}\left( \Delta_k^{+} \left( j,\tilde{q}\right) \right)_{i_{k-1}}\left( 2k+1 \right)_{i_k} \left( 2k+2- \gamma \right)_{i_k}} \right) \nonumber\\
&&\times \frac{ \left( \Delta_{\tau }^{-} \left( j,\tilde{q}\right) \right)_{n}\left( \Delta_{\tau }^{+} \left( j,\tilde{q}\right) \right)_{n} \left( 2\tau +1 \right)_{i_{\tau -1}} \left( 2\tau +2- \gamma \right)_{i_{\tau -1}}}{\left( \Delta_{\tau }^{-} \left( j,\tilde{q}\right) \right)_{i_{\tau -1}}\left( \Delta_{\tau }^{+} \left( j,\tilde{q}\right) \right)_{i_{\tau -1}}\left( 2\tau +1 \right)_{n} \left( 2\tau +2- \gamma  \right)_{n}} \left(\frac{1+a}{a}\right)^{n } \nonumber
\end{eqnarray}
\begin{eqnarray}
y_0^m(j,\tilde{q};x) &=& \sum_{i_0=0}^{m} \frac{\left( \Delta_0^{-} \left( j,\tilde{q}\right) \right)_{i_0}\left( \Delta_0^{+} \left( j,\tilde{q}\right) \right)_{i_0}}{\left( 1 \right)_{i_0} \left( 2-\gamma \right)_{i_0}} \eta ^{i_0} \nonumber\\
y_1^m(j,\tilde{q};x) &=& \left\{\sum_{i_0=0}^{m}\frac{\left( i_0 -j\right)\left( i_0+2+\beta -\gamma -\delta \right) }{\left( i_0+2 \right) \left( i_0+3- \gamma \right)} \frac{ \left( \Delta_0^{-} \left( j,\tilde{q}\right) \right)_{i_0}\left( \Delta_0^{+} \left( j,\tilde{q}\right) \right)_{i_0}}{\left( 1 \right)_{i_0} \left( 2-\gamma \right)_{i_0}} \right. \nonumber\\
&&\times \left. \sum_{i_1 = i_0}^{m} \frac{ \left( \Delta_1^{-} \left( j,\tilde{q}\right) \right)_{i_1}\left( \Delta_1^{+} \left( j,\tilde{q}\right) \right)_{i_1} \left( 3 \right)_{i_0} \left( 4-\gamma \right)_{i_0}}{\left( \Delta_1^{-} \left( j,\tilde{q}\right) \right)_{i_0}\left( \Delta_1^{+} \left( j,\tilde{q}\right) \right)_{i_0}\left( 3 \right)_{i_1} \left( 4-\gamma \right)_{i_1}} \eta ^{i_1}\right\} z 
\nonumber
\end{eqnarray}
\begin{eqnarray}
y_{\tau }^m(j,\tilde{q};x) &=& \left\{ \sum_{i_0=0}^{m} \frac{\left( i_0 -j\right)\left( i_0+2+\beta -\gamma -\delta \right) }{\left( i_0+2 \right) \left( i_0+3-\gamma \right)} \frac{ \left( \Delta_0^{-} \left( j,\tilde{q}\right) \right)_{i_0}\left( \Delta_0^{+} \left( j,\tilde{q}\right) \right)_{i_0}}{\left( 1 \right)_{i_0} \left( 2-\gamma \right)_{i_0}} \right.\nonumber\\
&&\times \prod_{k=1}^{\tau -1} \left( \sum_{i_k = i_{k-1}}^{m} \frac{\left( i_k+ 2k-j\right)\left( i_k +2k+2+\beta -\gamma -\delta \right)}{\left( i_k+2k+2 \right) \left( i_k+2k+3- \gamma \right)} \right. \nonumber\\
&&\times  \left. \frac{ \left( \Delta_k^{-} \left( j,\tilde{q}\right) \right)_{i_k}\left( \Delta_k^{+} \left( j,\tilde{q}\right) \right)_{i_k} \left( 2k+1 \right)_{i_{k-1}} \left( 2k+2- \gamma \right)_{i_{k-1}}}{\left( \Delta_k^{-} \left( j,\tilde{q}\right) \right)_{i_{k-1}}\left( \Delta_k^{+} \left( j,\tilde{q}\right) \right)_{i_{k-1}}\left( 2k+1 \right)_{i_k} \left( 2k+2- \gamma \right)_{i_k}} \right) \nonumber\\
&&\times \left. \sum_{i_{\tau } = i_{\tau -1}}^{m}  \frac{ \left( \Delta_{\tau }^{-} \left( j,\tilde{q}\right) \right)_{i_{\tau }}\left( \Delta_{\tau }^{+} \left( j,\tilde{q}\right) \right)_{i_{\tau }} \left( 2\tau +1  \right)_{i_{\tau -1}} \left( 2\tau +2- \gamma \right)_{i_{\tau -1}}}{\left( \Delta_{\tau }^{-} \left( j,\tilde{q}\right) \right)_{i_{\tau -1}}\left( \Delta_{\tau }^{+} \left( j,\tilde{q}\right) \right)_{i_{\tau -1}}\left( 2\tau +1 \right)_{i_\tau } \left( 2\tau +2- \gamma \right)_{i_{\tau }}} \eta ^{i_{\tau }}\right\} z^{\tau } \nonumber
\end{eqnarray}
where
\begin{equation}
\begin{cases} \tau \geq 2 \cr
z = -\frac{1}{a}x^2 \cr
\eta = \frac{(1+a)}{a} x \cr
\tilde{q} = q-(\gamma +\delta -2)a -(\gamma -1)(\alpha +\beta -\gamma -\delta +1) \cr
\Delta_k^{\pm} \left( j,\tilde{q} \right) = \frac{\varphi +4(1+a)k \pm \sqrt{\varphi ^2-4(1+a)\tilde{q} }}{2(1+a)}  \cr
\varphi = \beta -\gamma -j-a( \gamma +\delta -3)
\end{cases}\nonumber
\end{equation}
\end{enumerate}
\subsubsection{The second species complete polynomial}
Replacing coefficients $q$, $\alpha$, $\beta$, $\gamma $ and $\delta$ by $q-(\gamma +\delta -2)a-(\gamma -1)(\alpha +\beta -\gamma -\delta +1)$, $\alpha - \gamma -\delta +2$, $\beta - \gamma -\delta +2, 2-\gamma$ and $2 - \delta$ into (\ref{eq:50040a})--(\ref{eq:50041c}). Multiply $x^{1-\gamma } (1-x)^{1-\delta }$ and the new (\ref{eq:50040a})--(\ref{eq:50041c}) together.\footnote{I treat $\gamma$ and $\delta$ as free variables and fixed values of $\alpha $, $\beta $ and $q$.}  
\begin{enumerate} 
\item As $\alpha =\gamma +\delta -2 $, $\beta =\gamma +\delta -1 $ and $q=(\gamma +\delta -2)a +(\gamma -1)(\gamma +\delta -2) $,
 
Its eigenfunction is given by
\begin{eqnarray}
& &x^{1-\gamma } (1-x)^{1-\delta } y(x)\nonumber\\
&=& x^{1-\gamma } (1-x)^{1-\delta } Hl(a, 0; 0, 1, 2-\gamma, 2 - \delta ; x)\nonumber\\
&=& x^{1-\gamma } (1-x)^{1-\delta } \nonumber
\end{eqnarray}
\item As $\alpha =\gamma +\delta -3 $, $\beta =\gamma +\delta -2 $ and $q =2(\gamma +\delta -3)a +(\gamma -1)(\gamma +\delta -4) +(2-\delta ) $,
 
Its eigenfunction is given by
\begin{eqnarray}
& &x^{1-\gamma } (1-x)^{1-\delta } y(x)\nonumber\\
&=& x^{1-\gamma } (1-x)^{1-\delta } Hl(a, (2-\delta ) +a\left( \gamma +\delta -4 \right); -1, 0, 2-\gamma, 2 - \delta ; x)\nonumber\\
&=& x^{1-\gamma } (1-x)^{1-\delta } \left\{ 1+ \frac{2-\delta +a\left( \gamma +\delta -4\right)}{(1+a) (2-\gamma )}\eta \right\} \nonumber
\end{eqnarray}
\item As $\alpha =\gamma +\delta -2N-4 $, $\beta =\gamma +\delta -2N-3 $ and $q =(\gamma +\delta -2)a +(\gamma -1)(\gamma +\delta -4N-6) + \left( 2N+2 \right) \left[ -\delta +2N+3 +a\left( \gamma +\delta -2N-5 \right)\right]$ where $N \in \mathbb{N}_{0}$,

Its eigenfunction is given by
\begin{eqnarray}
 & &x^{1-\gamma } (1-x)^{1-\delta } y(x)\nonumber\\
&=& x^{1-\gamma } (1-x)^{1-\delta } Hl(a, \left( 2N+2 \right) \left[ -\delta +2N+3 +a\left( \gamma +\delta -2N-5 \right)\right ]; -2N-2, -2N-1 \nonumber\\
&&, 2-\gamma, 2 - \delta ; x)\nonumber\\
&=& x^{1-\gamma } (1-x)^{1-\delta }  \sum_{r=0}^{N+1} y_{r}^{2(N+1-r)}\left( 2N+2;x\right)  \nonumber
\end{eqnarray}
\item As $\alpha =\gamma +\delta -2N-5 $, $\beta =\gamma +\delta -2N-4 $ and $q =(\gamma +\delta -2)a +(\gamma -1)( \gamma +\delta -4N-8) + \left( 2N+3 \right) \left[ -\delta +2N+4 +a\left( \gamma +\delta -2N-6 \right)\right]$ where $N \in \mathbb{N}_{0}$,

Its eigenfunction is given by
\begin{eqnarray}
& &x^{1-\gamma } (1-x)^{1-\delta } y(x)\nonumber\\
&=& x^{1-\gamma } (1-x)^{1-\delta } Hl(a, \left( 2N+3 \right) \left[ -\delta +2N+4 +a\left( \gamma +\delta -2N-6 \right)\right]; -2N-3, -2N-2 \nonumber\\
&&, 2-\gamma, 2 - \delta ; x)\nonumber\\
&=& x^{1-\gamma } (1-x)^{1-\delta } \sum_{r=0}^{N+1} y_{r}^{2(N-r)+3}\left( 2N+3;x\right)  \nonumber
\end{eqnarray}
In the above,
\begin{eqnarray}
y_0^m(j;x) &=&  \sum_{i_0=0}^{m} \frac{\left( -j\right)_{i_0} \left( \Pi _0\left( j\right)\right)_{i_0}}{\left( 1\right)_{i_0} \left( 2-\gamma \right)_{i_0}} \eta ^{i_0} \nonumber\\
y_1^m(j;x) &=& \left\{\sum_{i_0=0}^{m} \frac{\left( i_0-j\right)\left( i_0 +1-j \right)}{\left( i_0+2 \right) \left( i_0+3-\gamma \right)} \frac{\left( -j\right)_{i_0} \left( \Pi _0\left( j\right)\right)_{i_0}}{\left( 1 \right)_{i_0} \left( 2-\gamma \right)_{i_0}} \right. \nonumber\\
&&\times \left. \sum_{i_1 = i_0}^{m} \frac{\left( 2-j\right)_{i_1} \left( \Pi _1\left( j\right)\right)_{i_1}\left( 3 \right)_{i_0} \left(  4-\gamma \right)_{i_0}}{\left( 2-j \right)_{i_0} \left( \Pi _1\left( j\right)\right)_{i_0}\left( 3 \right)_{i_1} \left(  4-\gamma \right)_{i_1}} \eta ^{i_1}\right\} z \nonumber\\
y_{\tau }^m(j;x) &=& \left\{ \sum_{i_0=0}^{m} \frac{\left( i_0-j\right)\left( i_0 +1-j \right)}{\left( i_0+2 \right) \left( i_0+3-\gamma  \right)} \frac{\left( -j\right)_{i_0} \left( \Pi _0\left( j\right)\right)_{i_0}}{\left( 1 \right)_{i_0} \left( 2-\gamma \right)_{i_0}} \right.\nonumber\\
&&\times \prod_{k=1}^{\tau -1} \left( \sum_{i_k = i_{k-1}}^{m} \frac{\left( i_k+ 2k-j\right)\left( i_k +2k+1-j \right)}{\left( i_k+2k+2 \right) \left( i_k+2k+3- \gamma \right)} \right. \nonumber\\
&&\times \left. \frac{\left( 2k-j\right)_{i_k} \left( \Pi _k\left( j\right) \right)_{i_k}\left( 2k+1 \right)_{i_{k-1}} \left( 2k+2-\gamma \right)_{i_{k-1}}}{\left( 2k-j\right)_{i_{k-1}} \left( \Pi _k\left( j\right) \right)_{i_{k-1}}\left( 2k+1 \right)_{i_k} \left( 2k+2-\gamma \right)_{i_k}} \right) \nonumber\\
&&\times \left. \sum_{i_{\tau } = i_{\tau -1}}^{m} \frac{\left( 2\tau -j\right)_{i_{\tau }} \left( \Pi _{\tau }\left( j\right)\right)_{i_{\tau }}\left( 2\tau +1 \right)_{i_{\tau -1}} \left( 2\tau +2-\gamma \right)_{i_{\tau -1}}}{\left( 2\tau -j\right)_{i_{\tau -1}} \left(\Pi _{\tau }\left( j\right)\right)_{i_{\tau -1}}\left( 2\tau +1 \right)_{i_{\tau }} \left( 2\tau +2-\gamma \right)_{i_{\tau }}} \eta ^{i_{\tau }}\right\} z^{\tau } \nonumber
\end{eqnarray}
where
\begin{equation}
\begin{cases} \tau \geq 2 \cr
z = -\frac{1}{a}x^2 \cr
\eta = \frac{(1+a)}{a} x \cr
\Pi _k\left( j\right) = \frac{1}{1+a}\left( \delta -1-j -a\left( \gamma +\delta -3-j \right)\right) +2k  
\end{cases}\nonumber
\end{equation}
\end{enumerate}
\subsection{${\displaystyle  Hl(1-a,-q+\alpha \beta; \alpha,\beta, \delta, \gamma; 1-x)}$} 
\subsubsection{The first species complete polynomial}
Replacing coefficients $a$, $q$, $\gamma $, $\delta$ and $x$ by $1-a$, $-q +\alpha \beta $, $\delta $, $\gamma $ and $1-x$  into (\ref{eq:50019})--(\ref{eq:50024c}).\footnote{I treat $\beta $, $\gamma$ and $\delta$ as free variables and fixed values of $\alpha $ and $q$.}
\begin{enumerate} 
\item As $\alpha =0$ and $ q = -q_0^0$ where $q_0^0=0$, 

The eigenfunction is given by
\begin{equation}
y(\xi ) = Hl(1-a, 0; 0,\beta, \delta, \gamma; 1-x) =1 \nonumber
\end{equation}
\item As $\alpha =-1$,

An algebraic equation of degree 2 for the determination of $q$ is given by
\begin{equation}
0 = (1-a)\beta \gamma 
+\prod_{l=0}^{1}\Big( -q -\beta + l(\beta -\gamma -1+l+(1-a)(\gamma +\delta -1+l))\Big) \nonumber
\end{equation}
The eigenvalue of $q$ is written by $-\beta -q_1^m$ where $m = 0,1 $; $q_{1}^0 < q_{1}^1$. Its eigenfunction is given by
\begin{eqnarray}
y(\xi ) &=& Hl(1-a, q_1^m; -1,\beta, \delta, \gamma; 1-x)\nonumber\\
&=&  1+\frac{ q_1^m}{(2-a)\delta } \eta \nonumber 
\end{eqnarray}
\item As $\alpha =-2N-2 $ where $N \in \mathbb{N}_{0}$,

An algebraic equation of degree $2N+3$ for the determination of $q$ is given by
\begin{equation}
0 = \sum_{r=0}^{N+1}\bar{c}\left( r, 2(N-r)+3; 2N+2,-q-(2N+2)\beta \right)  \nonumber
\end{equation}
The eigenvalue of $q$ is written by $-(2N+2)\beta -q_{2N+2}^m$ where $m = 0,1,2,\cdots,2N+2 $; $q_{2N+2}^0 < q_{2N+2}^1 < \cdots < q_{2N+2}^{2N+2}$. Its eigenfunction is given by 
\begin{eqnarray} 
y(\xi ) &=& Hl(1-a, q_{2N+2}^m; -2N-2,\beta, \delta, \gamma; 1-x)\nonumber\\
&=& \sum_{r=0}^{N+1} y_{r}^{2(N+1-r)}\left( 2N+2, q_{2N+2}^m; \xi \right)  
\nonumber 
\end{eqnarray}
\item As $\alpha =-2N-3 $ where $N \in \mathbb{N}_{0}$,

An algebraic equation of degree $2N+4$ for the determination of $q$ is given by
\begin{equation}  
0 = \sum_{r=0}^{N+2}\bar{c}\left( r, 2(N+2-r); 2N+3,-q-(2N+3)\beta\right) \nonumber
\end{equation}
The eigenvalue of $q $ is written by $-(2N+3)\beta -q_{2N+3}^m$ where $m = 0,1,2,\cdots,2N+3 $; $q_{2N+3}^0 < q_{2N+3}^1 < \cdots < q_{2N+3}^{2N+3}$. Its eigenfunction is given by
\begin{eqnarray} 
y(\xi ) &=& Hl(1-a, q_{2N+3}^m; -2N-3,\beta, \delta, \gamma; 1-x)\nonumber\\
&=&   \sum_{r=0}^{N+1} y_{r}^{2(N-r)+3} \left( 2N+3,q_{2N+3}^m;\xi \right) \nonumber
\end{eqnarray}
In the above,
\begin{eqnarray}
\bar{c}(0,n;j,\tilde{q})  &=& \frac{\left( \Delta_0^{-} \left( j,\tilde{q}\right) \right)_{n}\left( \Delta_0^{+} \left( j,\tilde{q}\right) \right)_{n}}{\left( 1 \right)_{n} \left( \delta  \right)_{n}} \left( \frac{2-a}{1-a} \right)^{n}\nonumber\\
\bar{c}(1,n;j,\tilde{q}) &=& \left( -\frac{1}{1-a}\right) \sum_{i_0=0}^{n}\frac{\left( i_0 -j\right)\left( i_0+\beta \right) }{\left( i_0+2 \right) \left( i_0+1+ \delta \right)} \frac{ \left( \Delta_0^{-} \left( j,\tilde{q}\right) \right)_{i_0}\left( \Delta_0^{+} \left( j,\tilde{q}\right) \right)_{i_0}}{\left( 1 \right)_{i_0} \left( \delta \right)_{i_0}} \nonumber\\
&&\times  \frac{ \left( \Delta_1^{-} \left( j,\tilde{q}\right) \right)_{n}\left( \Delta_1^{+} \left( j,\tilde{q}\right) \right)_{n} \left( 3 \right)_{i_0} \left( 2+ \delta \right)_{i_0}}{\left( \Delta_1^{-} \left( j,\tilde{q}\right) \right)_{i_0}\left( \Delta_1^{+} \left( j,\tilde{q}\right) \right)_{i_0}\left( 3 \right)_{n} \left( 2+ \delta \right)_{n}} \left(\frac{2-a}{1-a} \right)^{n }  
\nonumber\\
\bar{c}(\tau ,n;j,\tilde{q}) &=& \left( -\frac{1}{1-a}\right)^{\tau} \sum_{i_0=0}^{n}\frac{\left( i_0 -j\right)\left( i_0+\beta \right) }{\left( i_0+2 \right) \left( i_0+1+ \delta \right)}\frac{ \left( \Delta_0^{-} \left( j,\tilde{q}\right) \right)_{i_0}\left( \Delta_0^{+} \left( j,\tilde{q}\right) \right)_{i_0}}{\left( 1 \right)_{i_0} \left( \delta \right)_{i_0}}  \nonumber\\
&&\times \prod_{k=1}^{\tau -1} \left( \sum_{i_k = i_{k-1}}^{n} \frac{\left( i_k+ 2k-j\right)\left( i_k +2k+\beta \right)}{\left( i_k+2k+2 \right) \left( i_k+2k+1+ \delta \right)} \right. \nonumber\\
&&\times  \left. \frac{ \left( \Delta_k^{-} \left( j,\tilde{q}\right) \right)_{i_k}\left( \Delta_k^{+} \left( j,\tilde{q}\right) \right)_{i_k} \left( 2k+1 \right)_{i_{k-1}} \left( 2k+ \delta \right)_{i_{k-1}}}{\left( \Delta_k^{-} \left( j,\tilde{q}\right) \right)_{i_{k-1}}\left( \Delta_k^{+} \left( j,\tilde{q}\right) \right)_{i_{k-1}}\left( 2k+1 \right)_{i_k} \left( 2k+ \delta \right)_{i_k}} \right)  \nonumber\\
&&\times \frac{ \left( \Delta_{\tau }^{-} \left( j,\tilde{q}\right) \right)_{n}\left( \Delta_{\tau }^{+} \left( j,\tilde{q}\right) \right)_{n} \left( 2\tau +1 \right)_{i_{\tau -1}} \left( 2\tau + \delta \right)_{i_{\tau -1}}}{\left( \Delta_{\tau }^{-} \left( j,\tilde{q}\right) \right)_{i_{\tau -1}}\left( \Delta_{\tau }^{+} \left( j,\tilde{q}\right) \right)_{i_{\tau -1}}\left( 2\tau +1 \right)_{n} \left( 2\tau + \delta  \right)_{n}} \left(\frac{2-a}{1-a}\right)^{n } \nonumber
\end{eqnarray}
\begin{eqnarray}
y_0^m(j,\tilde{q};\xi ) &=& \sum_{i_0=0}^{m} \frac{\left( \Delta_0^{-} \left( j,\tilde{q}\right) \right)_{i_0}\left( \Delta_0^{+} \left( j,\tilde{q}\right) \right)_{i_0}}{\left( 1 \right)_{i_0} \left( \delta \right)_{i_0}} \eta ^{i_0} \nonumber\\
y_1^m(j,\tilde{q};\xi ) &=& \left\{\sum_{i_0=0}^{m}\frac{\left( i_0 -j\right)\left( i_0+\beta \right) }{\left( i_0+2 \right) \left( i_0+1+ \delta \right)} \frac{ \left( \Delta_0^{-} \left( j,\tilde{q}\right) \right)_{i_0}\left( \Delta_0^{+} \left( j,\tilde{q}\right) \right)_{i_0}}{\left( 1 \right)_{i_0} \left( \delta \right)_{i_0}} \right. \nonumber\\
&&\times \left. \sum_{i_1 = i_0}^{m} \frac{ \left( \Delta_1^{-} \left( j,\tilde{q}\right) \right)_{i_1}\left( \Delta_1^{+} \left( j,\tilde{q}\right) \right)_{i_1} \left( 3 \right)_{i_0} \left( 2+ \delta \right)_{i_0}}{\left( \Delta_1^{-} \left( j,\tilde{q}\right) \right)_{i_0}\left( \Delta_1^{+} \left( j,\tilde{q}\right) \right)_{i_0}\left( 3 \right)_{i_1} \left( 2+ \delta \right)_{i_1}} \eta ^{i_1}\right\} z 
\nonumber
\end{eqnarray}
\begin{eqnarray}
y_{\tau }^m(j,\tilde{q};\xi ) &=& \left\{ \sum_{i_0=0}^{m} \frac{\left( i_0 -j\right)\left( i_0+\beta \right) }{\left( i_0+2 \right) \left( i_0+1+ \delta \right)} \frac{ \left( \Delta_0^{-} \left( j,\tilde{q}\right) \right)_{i_0}\left( \Delta_0^{+} \left( j,\tilde{q}\right) \right)_{i_0}}{\left( 1 \right)_{i_0} \left( \delta \right)_{i_0}} \right.\nonumber\\
&&\times \prod_{k=1}^{\tau -1} \left( \sum_{i_k = i_{k-1}}^{m} \frac{\left( i_k+ 2k-j\right)\left( i_k +2k+\beta \right)}{\left( i_k+2k+2 \right) \left( i_k+2k+1+ \delta \right)} \right. \nonumber\\
&&\times  \left. \frac{ \left( \Delta_k^{-} \left( j,\tilde{q}\right) \right)_{i_k}\left( \Delta_k^{+} \left( j,\tilde{q}\right) \right)_{i_k} \left( 2k+1 \right)_{i_{k-1}} \left( 2k+ \delta \right)_{i_{k-1}}}{\left( \Delta_k^{-} \left( j,\tilde{q}\right) \right)_{i_{k-1}}\left( \Delta_k^{+} \left( j,\tilde{q}\right) \right)_{i_{k-1}}\left( 2k+1 \right)_{i_k} \left( 2k+ \delta \right)_{i_k}} \right) \nonumber\\
&&\times \left. \sum_{i_{\tau } = i_{\tau -1}}^{m}  \frac{ \left( \Delta_{\tau }^{-} \left( j,\tilde{q}\right) \right)_{i_{\tau }}\left( \Delta_{\tau }^{+} \left( j,\tilde{q}\right) \right)_{i_{\tau }} \left( 2\tau +1  \right)_{i_{\tau -1}} \left( 2\tau + \delta \right)_{i_{\tau -1}}}{\left( \Delta_{\tau }^{-} \left( j,\tilde{q}\right) \right)_{i_{\tau -1}}\left( \Delta_{\tau }^{+} \left( j,\tilde{q}\right) \right)_{i_{\tau -1}}\left( 2\tau +1 \right)_{i_\tau } \left( 2\tau + \delta \right)_{i_{\tau }}} \eta ^{i_{\tau }}\right\} z^{\tau } \nonumber 
\end{eqnarray}
where
\begin{equation}
\begin{cases} \tau \geq 2 \cr
\xi =1-x \cr
z = \frac{-1}{1-a}\xi^2 \cr
\eta = \frac{2-a}{1-a}\xi \cr
\tilde{q} = -q+\alpha \beta  \cr
\Delta_k^{\pm} \left( j,\tilde{q}\right) = \frac{\varphi +4(2-a)k \pm \sqrt{\varphi ^2-4(2-a)\tilde{q}}}{2(2-a)}  \cr
\varphi = \beta -\gamma -j+(1-a)(\gamma +\delta -1)
\end{cases}\nonumber
\end{equation}
\end{enumerate}
\subsubsection{The second species complete polynomial}
Replacing coefficients $a$, $q$, $\gamma $, $\delta$ and $x$ by $1-a$, $-q +\alpha \beta $, $\delta $, $\gamma $ and $1-x$  into (\ref{eq:50040a})--(\ref{eq:50041c}).\footnote{I treat $\gamma$ and $\delta$ as free variables and fixed values of $\alpha $, $\beta $ and $q$.}
\begin{enumerate}
\item As $\alpha =0$, $\beta =1 $ and $q= 0$, 
 
Its eigenfunction is given by
\begin{equation}
y(\xi ) = Hl(1-a,0; 0, 1, \delta, \gamma; 1-x)  =1 \nonumber
\end{equation}
\item As $\alpha =-1 $, $\beta = 0$ and $q= -\gamma +(1-a)\left( \gamma +\delta \right) $,
 
Its eigenfunction is given by
\begin{eqnarray}
y(\xi ) &=& Hl(1-a,\gamma -(1-a)\left( \gamma +\delta \right); -1,0, \delta, \gamma; 1-x)\nonumber\\
&=&  1+ \frac{\gamma  -(1-a)\left( \gamma +\delta \right)}{(2-a) \delta }\eta \nonumber
\end{eqnarray}
\item As $\alpha =-2N-2 $, $\beta = -2N-1 $ and $q =(2N+1)(2N+2)-  \left( 2N+2 \right) \left[ \gamma +2N+1 -(1-a)\left( \gamma +\delta +2N+1  \right)\right]$ where $N \in \mathbb{N}_{0}$,

Its eigenfunction is given by
\begin{eqnarray}
y(\xi ) &=& Hl(1-a,\left( 2N+2 \right) \left[ \gamma +2N+1 -(1-a)\left( \gamma +\delta +2N+1 \right)\right]; -2N-2,-2N-1 \nonumber\\
&&, \delta, \gamma; 1-x)\nonumber\\
&=&  \sum_{r=0}^{N+1} y_{r}^{2(N+1-r)}\left( 2N+2; \xi\right)  \nonumber
\end{eqnarray}
\item As $\alpha =-2N-3 $, $\beta = -2N-2 $ and $q = (2N+2)(2N+3) -\left( 2N+3 \right) \left[ \gamma +2N+2 -(1-a)\left( \gamma +\delta +2N+2 \right)\right]$ where $N \in \mathbb{N}_{0}$,

Its eigenfunction is given by
\begin{eqnarray}
y(\xi ) &=& Hl(1-a,\left( 2N+3 \right) \left[ \gamma +2N+2 -(1-a)\left( \gamma +\delta +2N+2 \right)\right]; -2N-3,-2N-2, \delta, \gamma; 1-x)\nonumber\\
&=& \sum_{r=0}^{N+1} y_{r}^{2(N-r)+3}\left( 2N+3; \xi\right)  \nonumber
\end{eqnarray}
In the above,
\begin{eqnarray}
y_0^m(j;\xi) &=&  \sum_{i_0=0}^{m} \frac{\left( -j\right)_{i_0} \left( \Pi _0\left( j\right)\right)_{i_0}}{\left( 1\right)_{i_0} \left( \delta \right)_{i_0}} \eta ^{i_0} \nonumber\\
y_1^m(j;\xi) &=& \left\{\sum_{i_0=0}^{m} \frac{\left( i_0-j\right)\left( i_0 +1-j \right)}{\left( i_0+2 \right) \left( i_0+1+\delta \right)} \frac{\left( -j\right)_{i_0} \left( \Pi _0\left( j\right)\right)_{i_0}}{\left( 1 \right)_{i_0} \left(\delta \right)_{i_0}} \right. \nonumber\\
&&\times \left. \sum_{i_1 = i_0}^{m} \frac{\left( 2-j\right)_{i_1} \left( \Pi _1\left( j\right)\right)_{i_1}\left( 3 \right)_{i_0} \left(  \delta +2 \right)_{i_0}}{\left( 2-j \right)_{i_0} \left( \Pi _1\left( j\right)\right)_{i_0}\left( 3 \right)_{i_1} \left(  \delta +2 \right)_{i_1}} \eta ^{i_1}\right\} z \nonumber\\
y_{\tau }^m(j;\xi) &=& \left\{ \sum_{i_0=0}^{m} \frac{\left( i_0-j\right)\left( i_0 +1-j \right)}{\left( i_0+2 \right) \left( i_0+1+\delta  \right)} \frac{\left( -j\right)_{i_0} \left( \Pi _0\left( j\right)\right)_{i_0}}{\left( 1 \right)_{i_0} \left(\delta \right)_{i_0}} \right.\nonumber\\
&&\times \prod_{k=1}^{\tau -1} \left( \sum_{i_k = i_{k-1}}^{m} \frac{\left( i_k+ 2k-j\right)\left( i_k +2k+1-j \right)}{\left( i_k+2k+2 \right) \left( i_k+2k+1+ \delta \right)} \right. \nonumber\\
&&\times \left. \frac{\left( 2k-j\right)_{i_k} \left( \Pi _k\left( j\right) \right)_{i_k}\left( 2k+1 \right)_{i_{k-1}} \left( 2k+\delta \right)_{i_{k-1}}}{\left( 2k-j\right)_{i_{k-1}} \left( \Pi _k\left( j\right) \right)_{i_{k-1}}\left( 2k+1 \right)_{i_k} \left( 2k+\delta \right)_{i_k}} \right) \nonumber\\
&&\times \left. \sum_{i_{\tau } = i_{\tau -1}}^{m} \frac{\left( 2\tau -j\right)_{i_{\tau }} \left( \Pi _{\tau }\left( j\right)\right)_{i_{\tau }}\left( 2\tau +1 \right)_{i_{\tau -1}} \left( 2\tau +\delta \right)_{i_{\tau -1}}}{\left( 2\tau -j\right)_{i_{\tau -1}} \left(\Pi _{\tau }\left( j\right)\right)_{i_{\tau -1}}\left( 2\tau +1 \right)_{i_{\tau }} \left( 2\tau +\delta \right)_{i_{\tau }}} \eta ^{i_{\tau }}\right\} z^{\tau } \nonumber
\end{eqnarray}
where
\begin{equation}
\begin{cases} \tau \geq 2 \cr
\xi =1-x \cr
z = \frac{-1}{1-a}\xi^2 \cr
\eta = \frac{2-a}{1-a}\xi \cr
\Pi _k\left( j\right) = \frac{1}{2-a}\left( -\gamma -j+1+(1-a)\left( \gamma +\delta +j-1 \right)\right) +2k  
\end{cases}\nonumber
\end{equation}
\end{enumerate}
\subsection{\footnotesize ${\displaystyle (1-x)^{1-\delta } Hl(1-a,-q+(\delta -1)\gamma a+(\alpha -\delta +1)(\beta -\delta +1); \alpha-\delta +1,\beta-\delta +1, 2-\delta, \gamma; 1-x)}$ \normalsize}
\subsubsection{The first species complete polynomial}
Replacing coefficients $a$, $q$, $\alpha $, $\beta $, $\gamma $, $\delta$ and $x$ by $1-a$, $-q+(\delta -1)\gamma a+(\alpha -\delta +1)(\beta -\delta +1)$, $\alpha-\delta +1 $, $\beta-\delta +1 $, $2-\delta$, $\gamma $ and $1-x$ into (\ref{eq:50019})--(\ref{eq:50024c}). Multiply $(1-x)^{1-\delta }$ and the new (\ref{eq:50019}), (\ref{eq:50020b}), (\ref{eq:50021b}) and (\ref{eq:50022b})  together.\footnote{I treat $\beta $, $\gamma$ and $\delta$ as free variables and fixed values of $\alpha $ and $q$.}
\begin{enumerate} 
\item As $\alpha = \delta -1$ and $ q= (\delta -1)\gamma a -q_0^0$ where $q_0^0=0$,

The eigenfunction is given by
\begin{eqnarray}
& &(1-x)^{1-\delta } y(\xi)\nonumber\\
&=& (1-x)^{1-\delta } Hl(1-a,0; 0,\beta-\delta +1, 2-\delta, \gamma; 1-x) \nonumber\\
&=& (1-x)^{1-\delta } \nonumber
\end{eqnarray}
\item As $\alpha =\delta -2$,

An algebraic equation of degree 2 for the determination of $q$ is given by
\begin{equation}
0 = (1-a)(\beta -\delta +1)(2-\delta ) 
+\prod_{l=0}^{1}\Big( -q+(\delta -1)\gamma a -(\beta -\delta +1)+ l(\beta -\gamma -\delta +l+(1-a)(\gamma -\delta +1+l))\Big) \nonumber
\end{equation}
The eigenvalue of $q $ is written by $(\delta -1)\gamma a -(\beta -\delta +1)-q_1^m$ where $m = 0,1 $; $q_{1}^0 < q_{1}^1$. Its eigenfunction is given by
\begin{eqnarray}
& &(1-x)^{1-\delta } y(\xi)\nonumber\\
&=& (1-x)^{1-\delta } Hl(1-a,q_1^m; -1,\beta -\delta +1, 2-\delta, \gamma; 1-x) \nonumber\\
&=& (1-x)^{1-\delta } \left\{ 1+\frac{ q_1^m}{(2-a)(2-\delta ) } \eta \right\} \nonumber  
\end{eqnarray}
\item As $\alpha =\delta -2N-3 $ where $N \in \mathbb{N}_{0}$,

An algebraic equation of degree $2N+3$ for the determination of $q$ is given by
\begin{equation}
0 = \sum_{r=0}^{N+1}\bar{c}\left( r, 2(N-r)+3; 2N+2,-q+(\delta -1)\gamma a-(2N+2)(\beta -\delta +1) \right)  \nonumber
\end{equation}
The eigenvalue of $q$ is written by $(\delta -1)\gamma a-(2N+2)(\beta -\delta +1)-q_{2N+2}^m$ where $m = 0,1,2,\cdots,2N+2 $; $q_{2N+2}^0 < q_{2N+2}^1 < \cdots < q_{2N+2}^{2N+2}$. Its eigenfunction is given by 
\begin{eqnarray} 
& &(1-x)^{1-\delta } y(\xi)\nonumber\\
&=& (1-x)^{1-\delta } Hl(1-a,q_{2N+2}^m; -2N-2,\beta-\delta +1, 2-\delta, \gamma; 1-x) \nonumber\\
&=& (1-x)^{1-\delta }\sum_{r=0}^{N+1} y_{r}^{2(N+1-r)}\left( 2N+2, q_{2N+2}^m; \xi \right)  
\nonumber
\end{eqnarray}
\item As $\alpha =\delta -2N-4 $ where $N \in \mathbb{N}_{0}$,

An algebraic equation of degree $2N+4$ for the determination of $q$ is given by
\begin{equation}  
0 = \sum_{r=0}^{N+2}\bar{c}\left( r, 2(N+2-r); 2N+3,-q + (\delta -1)\gamma a -(2N+3)(\beta -\delta +1)\right) \nonumber
\end{equation}
The eigenvalue of $q$ is written by $(\delta -1)\gamma a -(2N+3)(\beta -\delta +1) -q_{2N+3}^m$ where $m = 0,1,2,\cdots,2N+3 $; $q_{2N+3}^0 < q_{2N+3}^1 < \cdots < q_{2N+3}^{2N+3}$. Its eigenfunction is given by
\begin{eqnarray} 
& &(1-x)^{1-\delta } y(\xi)\nonumber\\
&=& (1-x)^{1-\delta } Hl(1-a,q_{2N+3}^m; -2N-3,\beta-\delta +1, 2-\delta, \gamma; 1-x) \nonumber\\
&=& (1-x)^{1-\delta }  \sum_{r=0}^{N+1} y_{r}^{2(N-r)+3} \left( 2N+3,q_{2N+3}^m;\xi\right) \nonumber
\end{eqnarray}
In the above,
\begin{eqnarray}
\bar{c}(0,n;j,\tilde{q})  &=& \frac{\left( \Delta_0^{-} \left( j,\tilde{q}\right) \right)_{n}\left( \Delta_0^{+} \left( j,\tilde{q}\right) \right)_{n}}{\left( 1 \right)_{n} \left( 2-\delta  \right)_{n}} \left( \frac{2-a}{1-a} \right)^{n}\nonumber\\
\bar{c}(1,n;j,\tilde{q}) &=& \left( -\frac{1}{1-a}\right) \sum_{i_0=0}^{n}\frac{\left( i_0 -j\right)\left( i_0+1+\beta -\delta  \right) }{\left( i_0+2 \right) \left( i_0+3-\delta  \right)} \frac{ \left( \Delta_0^{-} \left( j,\tilde{q}\right) \right)_{i_0}\left( \Delta_0^{+} \left( j,\tilde{q}\right) \right)_{i_0}}{\left( 1 \right)_{i_0} \left( 2-\delta \right)_{i_0}} \nonumber\\
&&\times  \frac{ \left( \Delta_1^{-} \left( j,\tilde{q}\right) \right)_{n}\left( \Delta_1^{+} \left( j,\tilde{q}\right) \right)_{n} \left( 3 \right)_{i_0} \left( 4-\delta \right)_{i_0}}{\left( \Delta_1^{-} \left( j,\tilde{q}\right) \right)_{i_0}\left( \Delta_1^{+} \left( j,\tilde{q}\right) \right)_{i_0}\left( 3 \right)_{n} \left( 4-\delta \right)_{n}} \left(\frac{2-a}{1-a} \right)^{n }  
\nonumber\\
\bar{c}(\tau ,n;j,\tilde{q}) &=& \left( -\frac{1}{1-a}\right)^{\tau} \sum_{i_0=0}^{n}\frac{\left( i_0 -j\right)\left( i_0+1+\beta -\delta  \right) }{\left( i_0+2 \right) \left( i_0+3-\delta  \right)}  \frac{ \left( \Delta_0^{-} \left( j,\tilde{q}\right) \right)_{i_0}\left( \Delta_0^{+} \left( j,\tilde{q}\right) \right)_{i_0}}{\left( 1 \right)_{i_0} \left( 2-\delta \right)_{i_0}}  \nonumber\\
&&\times \prod_{k=1}^{\tau -1} \left( \sum_{i_k = i_{k-1}}^{n} \frac{\left( i_k+ 2k-j\right)\left( i_k +2k+1+\beta -\delta  \right)}{\left( i_k+2k+2 \right) \left( i_k+2k+3-\delta  \right)} \right. \nonumber\\
&&\times  \left. \frac{ \left( \Delta_k^{-} \left( j,\tilde{q}\right) \right)_{i_k}\left( \Delta_k^{+} \left( j,\tilde{q}\right) \right)_{i_k} \left( 2k+1 \right)_{i_{k-1}} \left( 2k+2-\delta \right)_{i_{k-1}}}{\left( \Delta_k^{-} \left( j,\tilde{q}\right) \right)_{i_{k-1}}\left( \Delta_k^{+} \left( j,\tilde{q}\right) \right)_{i_{k-1}}\left( 2k+1 \right)_{i_k} \left( 2k+2-\delta \right)_{i_k}} \right)\nonumber\\
&&\times \frac{ \left( \Delta_{\tau }^{-} \left( j,\tilde{q}\right) \right)_{n}\left( \Delta_{\tau }^{+} \left( j,\tilde{q}\right) \right)_{n} \left( 2\tau +1 \right)_{i_{\tau -1}} \left( 2\tau +2-\delta  \right)_{i_{\tau -1}}}{\left( \Delta_{\tau }^{-} \left( j,\tilde{q}\right) \right)_{i_{\tau -1}}\left( \Delta_{\tau }^{+} \left( j,\tilde{q}\right) \right)_{i_{\tau -1}}\left( 2\tau +1 \right)_{n} \left( 2\tau +2-\delta \right)_{n}} \left(\frac{2-a}{1-a}\right)^{n } \nonumber
\end{eqnarray}
\begin{eqnarray}
y_0^m(j,\tilde{q};\xi) &=& \sum_{i_0=0}^{m} \frac{\left( \Delta_0^{-} \left( j,\tilde{q}\right) \right)_{i_0}\left( \Delta_0^{+} \left( j,\tilde{q}\right) \right)_{i_0}}{\left( 1 \right)_{i_0} \left( 2-\delta  \right)_{i_0}} \eta ^{i_0} \nonumber\\
y_1^m(j,\tilde{q};\xi) &=& \left\{\sum_{i_0=0}^{m}\frac{\left( i_0 -j\right)\left( i_0+1+\beta -\delta \right) }{\left( i_0+2 \right) \left( i_0+3-\delta \right)} \frac{ \left( \Delta_0^{-} \left( j,\tilde{q}\right) \right)_{i_0}\left( \Delta_0^{+} \left( j,\tilde{q}\right) \right)_{i_0}}{\left( 1 \right)_{i_0} \left( 2-\delta  \right)_{i_0}} \right. \nonumber\\
&&\times \left. \sum_{i_1 = i_0}^{m} \frac{ \left( \Delta_1^{-} \left( j,\tilde{q}\right) \right)_{i_1}\left( \Delta_1^{+} \left( j,\tilde{q}\right) \right)_{i_1} \left( 3 \right)_{i_0} \left( 4-\delta \right)_{i_0}}{\left( \Delta_1^{-} \left( j,\tilde{q}\right) \right)_{i_0}\left( \Delta_1^{+} \left( j,\tilde{q}\right) \right)_{i_0}\left( 3 \right)_{i_1} \left( 4-\delta \right)_{i_1}} \eta ^{i_1}\right\} z 
\nonumber
\end{eqnarray}
\begin{eqnarray}
y_{\tau }^m(j,\tilde{q};\xi) &=& \left\{ \sum_{i_0=0}^{m} \frac{\left( i_0 -j\right)\left( i_0+1+\beta -\delta  \right) }{\left( i_0+2 \right) \left( i_0+3-\delta  \right)} \frac{ \left( \Delta_0^{-} \left( j,\tilde{q}\right) \right)_{i_0}\left( \Delta_0^{+} \left( j,\tilde{q}\right) \right)_{i_0}}{\left( 1 \right)_{i_0} \left( 2-\delta  \right)_{i_0}} \right.\nonumber\\
&&\times \prod_{k=1}^{\tau -1} \left( \sum_{i_k = i_{k-1}}^{m} \frac{\left( i_k+ 2k-j\right)\left( i_k +2k+1+\beta -\delta  \right)}{\left( i_k+2k+2 \right) \left( i_k+2k+3-\delta  \right)} \right. \nonumber\\
&&\times  \left. \frac{ \left( \Delta_k^{-} \left( j,\tilde{q}\right) \right)_{i_k}\left( \Delta_k^{+} \left( j,\tilde{q}\right) \right)_{i_k} \left( 2k+1 \right)_{i_{k-1}} \left( 2k+2-\delta \right)_{i_{k-1}}}{\left( \Delta_k^{-} \left( j,\tilde{q}\right) \right)_{i_{k-1}}\left( \Delta_k^{+} \left( j,\tilde{q}\right) \right)_{i_{k-1}}\left( 2k+1 \right)_{i_k} \left( 2k+2-\delta \right)_{i_k}} \right) \nonumber\\
&&\times \left. \sum_{i_{\tau } = i_{\tau -1}}^{m}  \frac{ \left( \Delta_{\tau }^{-} \left( j,\tilde{q}\right) \right)_{i_{\tau }}\left( \Delta_{\tau }^{+} \left( j,\tilde{q}\right) \right)_{i_{\tau }} \left( 2\tau +1  \right)_{i_{\tau -1}} \left( 2\tau +2-\delta \right)_{i_{\tau -1}}}{\left( \Delta_{\tau }^{-} \left( j,\tilde{q}\right) \right)_{i_{\tau -1}}\left( \Delta_{\tau }^{+} \left( j,\tilde{q}\right) \right)_{i_{\tau -1}}\left( 2\tau +1 \right)_{i_\tau } \left( 2\tau +2-\delta \right)_{i_{\tau }}} \eta ^{i_{\tau }}\right\} z^{\tau } \nonumber
\end{eqnarray}
where
\begin{equation}
\begin{cases} \tau \geq 2 \cr
\xi =1-x \cr
z = \frac{-1}{1-a}\xi^2 \cr
\eta = \frac{2-a}{1-a}\xi \cr
\tilde{q} = -q+(\delta -1)\gamma a+(\alpha -\delta +1)(\beta -\delta +1) \cr
\Delta_k^{\pm} \left( j,\tilde{q}\right) = \frac{\varphi +4(2-a)k \pm \sqrt{\varphi ^2-4(2-a)\tilde{q}}}{2(2-a)}  \cr
\varphi = \beta -\gamma -\delta +1-j+(1-a)(\gamma-\delta +1)
\end{cases}\nonumber
\end{equation}
\end{enumerate}
\subsubsection{The second species complete polynomial}
Replacing coefficients $a$, $q$, $\alpha $, $\beta $, $\gamma $, $\delta$ and $x$ by $1-a$, $-q+(\delta -1)\gamma a+(\alpha -\delta +1)(\beta -\delta +1)$, $\alpha-\delta +1 $, $\beta-\delta +1 $, $2-\delta$, $\gamma $ and $1-x$ into (\ref{eq:50040a})--(\ref{eq:50041c}). Multiply $(1-x)^{1-\delta }$ and the new (\ref{eq:50040a})--(\ref{eq:50041c}) together.\footnote{I treat $\gamma$ and $\delta$ as free variables and fixed values of $\alpha $, $\beta $ and $q$.}
\begin{enumerate} 
\item As $\alpha =\delta -1 $, $\beta =\delta $ and $ q=(\delta -1)\gamma a  $,
 
Its eigenfunction is given by
\begin{eqnarray}
& &(1-x)^{1-\delta } y(\xi)\nonumber\\
&=& (1-x)^{1-\delta } Hl(1-a,0; 0, 1, 2-\delta, \gamma; 1-x) \nonumber\\
&=& (1-x)^{1-\delta } \nonumber
\end{eqnarray}
\item As $\alpha =\delta -2 $, $\beta  =\delta -1$ and $ q =  (\delta -1)\gamma a -\gamma +(1-a)\left(  \gamma-\delta +2\right) $,
 
Its eigenfunction is given by
\begin{eqnarray}
& &(1-x)^{1-\delta } y(\xi)\nonumber\\
&=& (1-x)^{1-\delta } Hl(1-a,\gamma -(1-a)\left(  \gamma-\delta +2\right); -1,0, 2-\delta, \gamma; 1-x) \nonumber\\
&=& (1-x)^{1-\delta } \left\{  1+ \frac{\gamma  -(1-a)\left( \gamma -\delta +2\right)}{(2-a) (2-\delta ) }\eta \right\} \nonumber
\end{eqnarray}
\item As $\alpha =\delta -2N-3 $, $\beta = \delta -2N-2 $ and $q= (\delta -1)\gamma a +(2N+1)(2N+2) -\left( 2N+2 \right) \left[ \gamma +2N+1 -(1-a)\left( \gamma -\delta +2N+3\right)\right]$ where $N \in \mathbb{N}_{0}$,

Its eigenfunction is given by
\begin{eqnarray}
& &(1-x)^{1-\delta } y(\xi)\nonumber\\
&=& (1-x)^{1-\delta } Hl(1-a,\left( 2N+2 \right) \left[ \gamma +2N+1 -(1-a)\left( \gamma -\delta +2N+3\right)\right]; -2N-2, -2N-1 \nonumber\\
&&, 2-\delta, \gamma; 1-x) \nonumber\\
&=& (1-x)^{1-\delta }  \sum_{r=0}^{N+1} y_{r}^{2(N+1-r)}\left( 2N+2;\xi\right)  \nonumber
\end{eqnarray}
\item As $\alpha =\delta -2N-4 $, $\beta = \delta -2N-3 $ and $q = (\delta -1)\gamma a +(2N+2)(2N+3)- \left( 2N+3 \right) \left[ \delta  +2N+2 -(1-a)\left( \gamma -\delta +2N+4\right)\right]$ where $N \in \mathbb{N}_{0}$,

Its eigenfunction is given by
\begin{eqnarray}
& &(1-x)^{1-\delta } y(\xi)\nonumber\\
&=& (1-x)^{1-\delta } Hl(1-a,\left( 2N+3 \right) \left[ \delta  +2N+2 -(1-a)\left( \gamma -\delta +2N+4\right)\right]; -2N-3, -2N-2 \nonumber\\
&&, 2-\delta, \gamma; 1-x) \nonumber\\
&=& (1-x)^{1-\delta } \sum_{r=0}^{N+1} y_{r}^{2(N-r)+3}\left( 2N+3;\xi\right)  \nonumber
\end{eqnarray}
In the above,
\begin{eqnarray}
y_0^m(j;\xi) &=&  \sum_{i_0=0}^{m} \frac{\left( -j\right)_{i_0} \left( \Pi _0\left( j\right)\right)_{i_0}}{\left( 1\right)_{i_0} \left( 2-\delta  \right)_{i_0}} \eta ^{i_0} \nonumber\\
y_1^m(j;\xi) &=& \left\{\sum_{i_0=0}^{m} \frac{\left( i_0-j\right)\left( i_0 +1-j \right)}{\left( i_0+2 \right) \left( i_0+3-\delta  \right)} \frac{\left( -j\right)_{i_0} \left( \Pi _0\left( j\right)\right)_{i_0}}{\left( 1 \right)_{i_0} \left( 2-\delta  \right)_{i_0}} \right. \nonumber\\
&&\times \left. \sum_{i_1 = i_0}^{m} \frac{\left( 2-j\right)_{i_1} \left( \Pi _1\left( j\right)\right)_{i_1}\left( 3 \right)_{i_0} \left(  4-\delta \right)_{i_0}}{\left( 2-j \right)_{i_0} \left( \Pi _1\left( j\right)\right)_{i_0}\left( 3 \right)_{i_1} \left(  4-\delta \right)_{i_1}} \eta ^{i_1}\right\} z \nonumber\\
y_{\tau }^m(j;\xi) &=& \left\{ \sum_{i_0=0}^{m} \frac{\left( i_0-j\right)\left( i_0 +1-j \right)}{\left( i_0+2 \right) \left( i_0+3-\delta \right)} \frac{\left( -j\right)_{i_0} \left( \Pi _0\left( j\right)\right)_{i_0}}{\left( 1 \right)_{i_0} \left( 2-\delta  \right)_{i_0}} \right.\nonumber\\
&&\times \prod_{k=1}^{\tau -1} \left( \sum_{i_k = i_{k-1}}^{m} \frac{\left( i_k+ 2k-j\right)\left( i_k +2k+1-j \right)}{\left( i_k+2k+2 \right) \left( i_k+2k+3-\delta \right)} \right. \nonumber\\
&&\times \left. \frac{\left( 2k-j\right)_{i_k} \left( \Pi _k\left( j\right) \right)_{i_k}\left( 2k+1 \right)_{i_{k-1}} \left( 2k+2-\delta  \right)_{i_{k-1}}}{\left( 2k-j\right)_{i_{k-1}} \left( \Pi _k\left( j\right) \right)_{i_{k-1}}\left( 2k+1 \right)_{i_k} \left( 2k+2-\delta  \right)_{i_k}} \right) \nonumber\\
&&\times \left. \sum_{i_{\tau } = i_{\tau -1}}^{m} \frac{\left( 2\tau -j\right)_{i_{\tau }} \left( \Pi _{\tau }\left( j\right)\right)_{i_{\tau }}\left( 2\tau +1 \right)_{i_{\tau -1}} \left( 2\tau +2-\delta \right)_{i_{\tau -1}}}{\left( 2\tau -j\right)_{i_{\tau -1}} \left(\Pi _{\tau }\left( j\right)\right)_{i_{\tau -1}}\left( 2\tau +1 \right)_{i_{\tau }} \left( 2\tau +2-\delta \right)_{i_{\tau }}} \eta ^{i_{\tau }}\right\} z^{\tau } \nonumber
\end{eqnarray}
where
\begin{equation}
\begin{cases} \tau \geq 2 \cr
\xi =1-x \cr
z = \frac{-1}{1-a}\xi^2 \cr
\eta = \frac{2-a}{1-a}\xi \cr
\Pi _k\left( j\right) = \frac{1}{2-a}\left( -\gamma -j+1+(1-a)\left( \gamma -\delta +j+1 \right)\right) +2k  
\end{cases}\nonumber
\end{equation}
\end{enumerate}
\subsection{\footnotesize ${\displaystyle x^{-\alpha } Hl\left(\frac{1}{a},\frac{q+\alpha [(\alpha -\gamma -\delta +1)a-\beta +\delta ]}{a}; \alpha , \alpha -\gamma +1, \alpha -\beta +1,\delta ;\frac{1}{x}\right)}$\normalsize}
\subsubsection{The first species complete polynomial}
Replacing coefficients $a$, $q$, $\alpha $, $\beta $, $\gamma $ and $x$ by $\frac{1}{a}$, $\frac{q+\alpha [(\alpha -\gamma -\delta +1)a-\beta +\delta ]}{a}$, $\alpha-\gamma +1 $, $\alpha $,  $\alpha -\beta +1 $ and $\frac{1}{x}$ into (\ref{eq:50019})--(\ref{eq:50024c}). Multiply $x^{-\alpha }$ and the new (\ref{eq:50019}), (\ref{eq:50020b}), (\ref{eq:50021b}) and (\ref{eq:50022b})  together.\footnote{I treat $\alpha $, $\beta$ and $\delta$ as free variables and fixed values of $\gamma $ and $q$.}
\begin{enumerate} 
\item As $\gamma =\alpha +1 $ and $ q = \alpha ((a-1)\delta +\beta )+ a q_0^0 $ where $q_0^0=0$,

The eigenfunction is given by
\begin{eqnarray}
& &x^{-\alpha } y(\xi)\nonumber\\
&=& x^{-\alpha }  Hl\left(\frac{1}{a},0; 0, \alpha , \alpha -\beta +1,\delta ;\frac{1}{x}\right) \nonumber\\
&=& x^{-\alpha } \nonumber 
\end{eqnarray}
\item As $\gamma =\alpha +2$,

An algebraic equation of degree 2 for the determination of $q$ is given by
\begin{equation}
0 = \frac{1}{a}\alpha (\alpha -\beta +1)
+\prod_{l=0}^{1}\left( \frac{q-\alpha \left( ( \delta +1)a +\beta -\delta \right)}{a}  + l(\alpha -\delta -1+l+a(\alpha -\beta  +\delta +l))\right) \nonumber
\end{equation}
The eigenvalue of $ q  $ is written by $\alpha \left( ( \delta +1)a +\beta -\delta \right) +a q_1^m$ where $m = 0,1 $; $q_{1}^0 < q_{1}^1$. Its eigenfunction is given by
\begin{eqnarray}
& &x^{-\alpha } y(\xi)\nonumber\\
&=& x^{-\alpha }  Hl\left(\frac{1}{a},q_1^m; -1, \alpha , \alpha -\beta +1,\delta ;\frac{1}{x}\right) \nonumber\\
&=& x^{-\alpha }\left\{ 1+\frac{ a q_1^m}{(1+a)(\alpha -\beta +1) } \eta \right\} \nonumber 
\end{eqnarray}
\item As $\gamma =\alpha +2N+3 $ where $N \in \mathbb{N}_{0}$,

An algebraic equation of degree $2N+3$ for the determination of $q$ is given by
\begin{equation}
0 = \sum_{r=0}^{N+1}\bar{c}\left( r, 2(N-r)+3; 2N+2,\frac{q-\alpha \left( ( \delta +2N+2)a +\beta -\delta \right)}{a}\right)  \nonumber
\end{equation}
The eigenvalue of $ q $ is written by $\alpha \left( ( \delta +2N+2)a +\beta -\delta \right) + a q_{2N+2}^m$ where $m = 0,1,2,\cdots,2N+2 $; $q_{2N+2}^0 < q_{2N+2}^1 < \cdots < q_{2N+2}^{2N+2}$. Its eigenfunction is given by 
\begin{eqnarray} 
& &x^{-\alpha } y(\xi)\nonumber\\
&=& x^{-\alpha }  Hl\left(\frac{1}{a}, q_{2N+2}^m; -2N-2, \alpha , \alpha -\beta +1,\delta ;\frac{1}{x}\right) \nonumber\\
&=& x^{-\alpha } \sum_{r=0}^{N+1} y_{r}^{2(N+1-r)}\left( 2N+2, q_{2N+2}^m; \xi \right)  
\nonumber
\end{eqnarray}
\item As $\alpha -\gamma +1 =-2N-3 $ where $N \in \mathbb{N}_{0}$,

An algebraic equation of degree $2N+4$ for the determination of $q$ is given by
\begin{equation}  
0 = \sum_{r=0}^{N+2}\bar{c}\left( r, 2(N+2-r); 2N+3,\frac{q-\alpha \left( ( \delta +2N+3)a +\beta -\delta \right)}{a}\right) \nonumber
\end{equation}
The eigenvalue of $q$ is written by $\alpha \left( ( \delta +2N+3)a +\beta -\delta \right)+ a q_{2N+3}^m$ where $m = 0,1,2,\cdots,2N+3 $; $q_{2N+3}^0 < q_{2N+3}^1 < \cdots < q_{2N+3}^{2N+3}$. Its eigenfunction is given by
\begin{eqnarray} 
& &x^{-\alpha } y(\xi)\nonumber\\
&=& x^{-\alpha }  Hl\left(\frac{1}{a}, q_{2N+3}^m; -2N-3, \alpha , \alpha -\beta +1,\delta ;\frac{1}{x}\right) \nonumber\\
&=& x^{-\alpha }   \sum_{r=0}^{N+1} y_{r}^{2(N-r)+3} \left( 2N+3,q_{2N+3}^m; \xi\right) \nonumber
\end{eqnarray}
In the above,
\begin{eqnarray}
\bar{c}(0,n;j,\tilde{q})  &=& \frac{\left( \Delta_0^{-} \left( j,\tilde{q}\right) \right)_{n}\left( \Delta_0^{+} \left( j,\tilde{q}\right) \right)_{n}}{\left( 1 \right)_{n} \left( \alpha -\beta +1 \right)_{n}} \left(  1+a \right)^{n} \nonumber\\
\bar{c}(1,n;j,\tilde{q}) &=& \left( -a\right) \sum_{i_0=0}^{n}\frac{\left( i_0 -j\right)\left( i_0+\alpha  \right) }{\left( i_0+2 \right) \left( i_0+2+\alpha -\beta  \right)} \frac{ \left( \Delta_0^{-} \left( j,\tilde{q}\right) \right)_{i_0}\left( \Delta_0^{+} \left( j,\tilde{q}\right) \right)_{i_0}}{\left( 1 \right)_{i_0} \left( \alpha -\beta +1 \right)_{i_0}} \nonumber\\
&&\times \frac{ \left( \Delta_1^{-} \left( j,\tilde{q}\right) \right)_{n}\left( \Delta_1^{+} \left( j,\tilde{q}\right) \right)_{n} \left( 3 \right)_{i_0} \left( \alpha -\beta +3 \right)_{i_0}}{\left( \Delta_1^{-} \left( j,\tilde{q}\right) \right)_{i_0}\left( \Delta_1^{+} \left( j,\tilde{q}\right) \right)_{i_0}\left( 3 \right)_{n} \left( \alpha -\beta +3 \right)_{n}} \left( 1+a \right)^{n } \nonumber\\
\bar{c}(\tau ,n;j,\tilde{q}) &=& \left( -a\right)^{\tau} \sum_{i_0=0}^{n}\frac{\left( i_0 -j\right)\left( i_0+\alpha \right) }{\left( i_0+2 \right) \left( i_0+2+\alpha -\beta \right)} \frac{ \left( \Delta_0^{-} \left( j,\tilde{q}\right) \right)_{i_0}\left( \Delta_0^{+} \left( j,\tilde{q}\right) \right)_{i_0}}{\left( 1 \right)_{i_0} \left( \alpha -\beta +1 \right)_{i_0}}\nonumber\\
&&\times \prod_{k=1}^{\tau -1} \left( \sum_{i_k = i_{k-1}}^{n} \frac{\left( i_k+ 2k-j\right)\left( i_k +2k+\alpha \right)}{\left( i_k+2k+2 \right) \left( i_k+2k+2+\alpha -\beta \right)} \right. \nonumber\\
&&\times  \left. \frac{ \left( \Delta_k^{-} \left( j,\tilde{q}\right) \right)_{i_k}\left( \Delta_k^{+} \left( j,\tilde{q}\right) \right)_{i_k} \left( 2k+1 \right)_{i_{k-1}} \left( 2k+1+ \alpha -\beta \right)_{i_{k-1}}}{\left( \Delta_k^{-} \left( j,\tilde{q}\right) \right)_{i_{k-1}}\left( \Delta_k^{+} \left( j,\tilde{q}\right) \right)_{i_{k-1}}\left( 2k+1 \right)_{i_k} \left( 2k+1+ \alpha -\beta \right)_{i_k}} \right) \nonumber\\
&&\times \frac{ \left( \Delta_{\tau }^{-} \left( j,\tilde{q}\right) \right)_{n}\left( \Delta_{\tau }^{+} \left( j,\tilde{q}\right) \right)_{n} \left( 2\tau +1 \right)_{i_{\tau -1}} \left( 2\tau +1+ \alpha -\beta \right)_{i_{\tau -1}}}{\left( \Delta_{\tau }^{-} \left( j,\tilde{q}\right) \right)_{i_{\tau -1}}\left( \Delta_{\tau }^{+} \left( j,\tilde{q}\right) \right)_{i_{\tau -1}}\left( 2\tau +1 \right)_{n} \left( 2\tau +1+ \alpha -\beta   \right)_{n}} \left( 1+a \right)^{n } \nonumber
\end{eqnarray}
\begin{eqnarray}
y_0^m(j,\tilde{q};\xi) &=& \sum_{i_0=0}^{m} \frac{\left( \Delta_0^{-} \left( j,\tilde{q}\right) \right)_{i_0}\left( \Delta_0^{+} \left( j,\tilde{q}\right) \right)_{i_0}}{\left( 1 \right)_{i_0} \left( \alpha -\beta +1 \right)_{i_0}} \eta ^{i_0} \nonumber\\
y_1^m(j,\tilde{q};\xi) &=& \left\{\sum_{i_0=0}^{m}\frac{\left( i_0 -j\right)\left( i_0+\alpha  \right) }{\left( i_0+2 \right) \left( i_0+2+\alpha -\beta \right)} \frac{ \left( \Delta_0^{-} \left( j,\tilde{q}\right) \right)_{i_0}\left( \Delta_0^{+} \left( j,\tilde{q}\right) \right)_{i_0}}{\left( 1 \right)_{i_0} \left( \alpha -\beta +1 \right)_{i_0}} \right. \nonumber\\
&&\times \left. \sum_{i_1 = i_0}^{m} \frac{ \left( \Delta_1^{-} \left( j,\tilde{q}\right) \right)_{i_1}\left( \Delta_1^{+} \left( j,\tilde{q}\right) \right)_{i_1} \left( 3 \right)_{i_0} \left( \alpha -\beta +3 \right)_{i_0}}{\left( \Delta_1^{-} \left( j,\tilde{q}\right) \right)_{i_0}\left( \Delta_1^{+} \left( j,\tilde{q}\right) \right)_{i_0}\left( 3 \right)_{i_1} \left( \alpha -\beta +3 \right)_{i_1}} \eta ^{i_1}\right\} z 
\nonumber\\
y_{\tau }^m(j,\tilde{q};\xi) &=& \left\{ \sum_{i_0=0}^{m} \frac{\left( i_0 -j\right)\left( i_0+\alpha  \right) }{\left( i_0+2 \right) \left( i_0+2+\alpha -\beta \right)} \frac{ \left( \Delta_0^{-} \left( j,\tilde{q}\right) \right)_{i_0}\left( \Delta_0^{+} \left( j,\tilde{q}\right) \right)_{i_0}}{\left( 1 \right)_{i_0} \left( \alpha -\beta +1 \right)_{i_0}} \right.\nonumber\\
&&\times \prod_{k=1}^{\tau -1} \left( \sum_{i_k = i_{k-1}}^{m} \frac{\left( i_k+ 2k-j\right)\left( i_k +2k+\alpha \right)}{\left( i_k+2k+2 \right) \left( i_k+2k+2+\alpha -\beta \right)} \right. \nonumber\\
&&\times  \left. \frac{ \left( \Delta_k^{-} \left( j,\tilde{q}\right) \right)_{i_k}\left( \Delta_k^{+} \left( j,\tilde{q}\right) \right)_{i_k} \left( 2k+1 \right)_{i_{k-1}} \left( 2k+1+ \alpha -\beta \right)_{i_{k-1}}}{\left( \Delta_k^{-} \left( j,\tilde{q}\right) \right)_{i_{k-1}}\left( \Delta_k^{+} \left( j,\tilde{q}\right) \right)_{i_{k-1}}\left( 2k+1 \right)_{i_k} \left( 2k+1+ \alpha -\beta \right)_{i_k}} \right) \nonumber\\
&&\times \left. \sum_{i_{\tau } = i_{\tau -1}}^{m}  \frac{ \left( \Delta_{\tau }^{-} \left( j,\tilde{q}\right) \right)_{i_{\tau }}\left( \Delta_{\tau }^{+} \left( j,\tilde{q}\right) \right)_{i_{\tau }} \left( 2\tau +1  \right)_{i_{\tau -1}} \left( 2\tau +1+ \alpha -\beta \right)_{i_{\tau -1}}}{\left( \Delta_{\tau }^{-} \left( j,\tilde{q}\right) \right)_{i_{\tau -1}}\left( \Delta_{\tau }^{+} \left( j,\tilde{q}\right) \right)_{i_{\tau -1}}\left( 2\tau +1 \right)_{i_\tau } \left( 2\tau +1+ \alpha -\beta \right)_{i_{\tau }}} \eta ^{i_{\tau }}\right\} z^{\tau }  \nonumber
\end{eqnarray}
where
\begin{equation}
\begin{cases} \tau \geq 2 \cr
\xi =\frac{1}{x} \cr
z = -a \xi^2 \cr
\eta = (1+a)\xi \cr
\tilde{q} = \frac{q+\alpha [(\alpha -\gamma -\delta +1)a-\beta +\delta ]}{a} \cr
\Delta_k^{\pm} \left( j,\tilde{q}\right) = \frac{a\varphi +4(1+a)k \pm \sqrt{(a\varphi )^2-4a(1+a)\tilde{q}}}{2(1+a)}  \cr
\varphi = \alpha  -\delta -j+\frac{1}{a}(\alpha -\beta +\delta )
\end{cases}\nonumber
\end{equation}
\end{enumerate}
\subsection{ \footnotesize ${\displaystyle \left(1-\frac{x}{a} \right)^{-\beta } Hl\left(1-a, -q+\gamma \beta; -\alpha +\gamma +\delta, \beta, \gamma, \delta; \frac{(1-a)x}{x-a} \right)}$\normalsize}
\subsubsection{The first species complete polynomial}
\paragraph{The case of $\gamma $, $q$ = fixed values and $\alpha, \beta, \delta $ = free variables }
Replacing coefficients $a$, $q$, $\alpha $ and $x$ by $1-a$, $-q+\gamma \beta $, $-\alpha+\gamma +\delta $ and $\frac{(1-a)x}{x-a}$ into (\ref{eq:50019})--(\ref{eq:50024c}). Replacing $\gamma $ by $\alpha -\delta -j$ into the new (\ref{eq:50023a})--(\ref{eq:50024c}). Multiply $\left(1-\frac{x}{a} \right)^{-\beta }$ and the new (\ref{eq:50019}), (\ref{eq:50020b}), (\ref{eq:50021b}) and (\ref{eq:50022b}) together. 
\begin{enumerate} 
\item As $ \gamma =\alpha - \delta $ and $ q= (\alpha -\delta ) \beta -q_0^0$  where $q_0^0=0$,

The eigenfunction is given by
\begin{eqnarray}
 && \left(1-\frac{x}{a} \right)^{-\beta } y(\xi ) \nonumber\\
 &=& \left( 1-\frac{x}{a} \right)^{-\beta } Hl\left( 1-a, 0; 0, \beta, \gamma, \delta; \frac{(1-a)x}{x-a} \right) \nonumber\\
&=& \left(1-\frac{x}{a} \right)^{-\beta } \nonumber 
\end{eqnarray}
\item As $ \gamma =\alpha -\delta -1$, 

An algebraic equation of degree 2 for the determination of $q$ is given by
\begin{equation}
0 = (1-a)\beta \gamma 
+\prod_{l=0}^{1}\Big( -q +(\alpha -\delta -1)\beta + l(\beta -\delta -1+l+(1-a)(\alpha -2+l))\Big) \nonumber
\end{equation}
The eigenvalue of $q$ is written by $(\alpha -\delta -1)\beta -q_1^m$ where $m = 0,1 $; $q_{1}^0 < q_{1}^1$. Its eigenfunction is given by
\begin{eqnarray}
&& \left(1-\frac{x}{a} \right)^{-\beta } y(\xi ) \nonumber\\
 &=& \left(1-\frac{x}{a} \right)^{-\beta } Hl\left( 1-a, q_1^m; -1, \beta, \gamma, \delta; \frac{(1-a)x}{x-a} \right) \nonumber\\
&=& \left( 1-\frac{x}{a} \right)^{-\beta } \left\{ 1+\frac{ q_1^m}{(2-a)(\alpha -\delta -1) } \eta \right\} \nonumber
\end{eqnarray}
\item As $ \gamma =\alpha -\delta -2N-2 $ where $N \in \mathbb{N}_{0}$,

An algebraic equation of degree $2N+3$ for the determination of $q$ is given by
\begin{equation}
0 = \sum_{r=0}^{N+1}\bar{c}\left( r, 2(N-r)+3; 2N+2,-q +(\alpha -\delta -2N-2)\beta \right)  
\nonumber
\end{equation}
The eigenvalue of $q$ is written by $(\alpha -\delta -2N-2)\beta -q_{2N+2}^m$ where $m = 0,1,2,\cdots,2N+2 $; $q_{2N+2}^0 < q_{2N+2}^1 < \cdots < q_{2N+2}^{2N+2}$. Its eigenfunction is given by 
\begin{eqnarray} 
&& \left(1-\frac{x}{a} \right)^{-\beta } y(\xi ) \nonumber\\
 &=& \left(1-\frac{x}{a} \right)^{-\beta } Hl\left( 1-a, q_{2N+2}^m; -2N-2, \beta, \gamma, \delta; \frac{(1-a)x}{x-a} \right) \nonumber\\
&=& \left(1-\frac{x}{a} \right)^{-\beta } \sum_{r=0}^{N+1} y_{r}^{2(N+1-r)}\left( 2N+2, q_{2N+2}^m; \xi \right)  
\nonumber
\end{eqnarray}
\item As $ \gamma =\alpha -\delta -2N-3 $ where $N \in \mathbb{N}_{0}$,

An algebraic equation of degree $2N+4$ for the determination of $q$ is given by
\begin{eqnarray}  
0 = \sum_{r=0}^{N+2}\bar{c}\left( r, 2(N+2-r); 2N+3,-q +(\alpha -\delta -2N-3)\beta \right) 
\nonumber
\end{eqnarray}
The eigenvalue of $q$ is written by $ (\alpha -\delta -2N-3)\beta -q_{2N+3}^m$ where $m = 0,1,2,\cdots,2N+3 $; $q_{2N+3}^0 < q_{2N+3}^1 < \cdots < q_{2N+3}^{2N+3}$. Its eigenfunction is given by
\begin{eqnarray} 
&& \left(1-\frac{x}{a} \right)^{-\beta } y(\xi ) \nonumber\\
 &=& \left(1-\frac{x}{a} \right)^{-\beta } Hl\left( 1-a, q_{2N+3}^m; -2N-3, \beta, \gamma, \delta; \frac{(1-a)x}{x-a} \right) \nonumber\\
&=& \left(1-\frac{x}{a} \right)^{-\beta }  \sum_{r=0}^{N+1} y_{r}^{2(N-r)+3} \left( 2N+3,q_{2N+3}^m;\xi \right) \nonumber
\end{eqnarray}
In the above,
\begin{eqnarray}
\bar{c}(0,n;j,\tilde{q})  &=& \frac{\left( \Delta_0^{-} \left( j,\tilde{q}\right) \right)_{n}\left( \Delta_0^{+} \left( j,\tilde{q}\right) \right)_{n}}{\left( 1 \right)_{n} \left( \alpha -\delta -j \right)_{n}} \left( \frac{2-a}{1-a} \right)^{n}\nonumber\\
\bar{c}(1,n;j,\tilde{q}) &=& \left( -\frac{1}{1-a}\right) \sum_{i_0=0}^{n}\frac{\left( i_0 -j\right)\left( i_0+\beta \right) }{\left( i_0+2 \right) \left( i_0+1-j+ \alpha -\delta   \right)} \frac{ \left( \Delta_0^{-} \left( j,\tilde{q}\right) \right)_{i_0}\left( \Delta_0^{+} \left( j,\tilde{q}\right) \right)_{i_0}}{\left( 1 \right)_{i_0} \left( \alpha -\delta -j \right)_{i_0}} \nonumber\\
&&\times \frac{ \left( \Delta_1^{-} \left( j,\tilde{q}\right) \right)_{n}\left( \Delta_1^{+} \left( j,\tilde{q}\right) \right)_{n} \left( 3 \right)_{i_0} \left( 2-j+ \alpha -\delta \right)_{i_0}}{\left( \Delta_1^{-} \left( j,\tilde{q}\right) \right)_{i_0}\left( \Delta_1^{+} \left( j,\tilde{q}\right) \right)_{i_0}\left( 3 \right)_{n} \left( 2-j+ \alpha -\delta \right)_{n}} \left(\frac{2-a}{1-a} \right)^{n }  
\nonumber\\
\bar{c}(\tau ,n;j,\tilde{q}) &=& \left( -\frac{1}{1-a}\right)^{\tau} \sum_{i_0=0}^{n}\frac{\left( i_0 -j\right)\left( i_0+\beta \right) }{\left( i_0+2 \right) \left( i_0+1-j+ \alpha -\delta \right)} \frac{ \left( \Delta_0^{-} \left( j,\tilde{q}\right) \right)_{i_0}\left( \Delta_0^{+} \left( j,\tilde{q}\right) \right)_{i_0}}{\left( 1 \right)_{i_0} \left( \alpha -\delta -j\right)_{i_0}}\nonumber\\
&&\times \prod_{k=1}^{\tau -1} \left( \sum_{i_k = i_{k-1}}^{n} \frac{\left( i_k+ 2k-j\right)\left( i_k +2k+\beta \right)}{\left( i_k+2k+2 \right) \left( i_k+2k+1-j+ \alpha -\delta \right)} \right. \nonumber\\
&&\times  \left. \frac{ \left( \Delta_k^{-} \left( j,\tilde{q}\right) \right)_{i_k}\left( \Delta_k^{+} \left( j,\tilde{q}\right) \right)_{i_k} \left( 2k+1 \right)_{i_{k-1}} \left( 2k-j+ \alpha -\delta \right)_{i_{k-1}}}{\left( \Delta_k^{-} \left( j,\tilde{q}\right) \right)_{i_{k-1}}\left( \Delta_k^{+} \left( j,\tilde{q}\right) \right)_{i_{k-1}}\left( 2k+1 \right)_{i_k} \left( 2k-j+ \alpha -\delta \right)_{i_k}} \right) \nonumber\\
&&\times \frac{ \left( \Delta_{\tau }^{-} \left( j,\tilde{q}\right) \right)_{n}\left( \Delta_{\tau }^{+} \left( j,\tilde{q}\right) \right)_{n} \left( 2\tau +1 \right)_{i_{\tau -1}} \left( 2\tau -j+ \alpha -\delta \right)_{i_{\tau -1}}}{\left( \Delta_{\tau }^{-} \left( j,\tilde{q}\right) \right)_{i_{\tau -1}}\left( \Delta_{\tau }^{+} \left( j,\tilde{q}\right) \right)_{i_{\tau -1}}\left( 2\tau +1 \right)_{n} \left( 2\tau -j+ \alpha -\delta \right)_{n}} \left(\frac{2-a}{1-a}\right)^{n } \nonumber 
\end{eqnarray}
\begin{eqnarray}
y_0^m(j,\tilde{q};\xi) &=& \sum_{i_0=0}^{m} \frac{\left( \Delta_0^{-} \left( j,\tilde{q}\right) \right)_{i_0}\left( \Delta_0^{+} \left( j,\tilde{q}\right) \right)_{i_0}}{\left( 1 \right)_{i_0} \left( \alpha -\delta -j \right)_{i_0}} \eta ^{i_0} \nonumber\\
y_1^m(j,\tilde{q};\xi) &=& \left\{\sum_{i_0=0}^{m}\frac{\left( i_0 -j\right)\left( i_0+\beta \right) }{\left( i_0+2 \right) \left( i_0+1-j+ \alpha -\delta \right)} \frac{ \left( \Delta_0^{-} \left( j,\tilde{q}\right) \right)_{i_0}\left( \Delta_0^{+} \left( j,\tilde{q}\right) \right)_{i_0}}{\left( 1 \right)_{i_0} \left( \alpha -\delta -j \right)_{i_0}} \right. \nonumber\\
&&\times \left. \sum_{i_1 = i_0}^{m} \frac{ \left( \Delta_1^{-} \left( j,\tilde{q}\right) \right)_{i_1}\left( \Delta_1^{+} \left( j,\tilde{q}\right) \right)_{i_1} \left( 3 \right)_{i_0} \left( 2-j+\alpha -\delta \right)_{i_0}}{\left( \Delta_1^{-} \left( j,\tilde{q}\right) \right)_{i_0}\left( \Delta_1^{+} \left( j,\tilde{q}\right) \right)_{i_0}\left( 3 \right)_{i_1} \left( 2-j+\alpha -\delta \right)_{i_1}} \eta ^{i_1}\right\} z 
\nonumber\\
y_{\tau }^m(j,\tilde{q};\xi) &=& \left\{ \sum_{i_0=0}^{m} \frac{\left( i_0 -j\right)\left( i_0+\beta \right) }{\left( i_0+2 \right) \left( i_0+1-j+ \alpha -\delta  \right)} \frac{ \left( \Delta_0^{-} \left( j,\tilde{q}\right) \right)_{i_0}\left( \Delta_0^{+} \left( j,\tilde{q}\right) \right)_{i_0}}{\left( 1 \right)_{i_0} \left( \alpha -\delta -j \right)_{i_0}} \right.\nonumber\\
&&\times \prod_{k=1}^{\tau -1} \left( \sum_{i_k = i_{k-1}}^{m} \frac{\left( i_k+ 2k-j\right)\left( i_k +2k+\beta \right)}{\left( i_k+2k+2 \right) \left( i_k+2k+1-j+ \alpha -\delta \right)} \right. \nonumber\\
&&\times  \left. \frac{ \left( \Delta_k^{-} \left( j,\tilde{q}\right) \right)_{i_k}\left( \Delta_k^{+} \left( j,\tilde{q}\right) \right)_{i_k} \left( 2k+1 \right)_{i_{k-1}} \left( 2k-j+ \alpha -\delta \right)_{i_{k-1}}}{\left( \Delta_k^{-} \left( j,\tilde{q}\right) \right)_{i_{k-1}}\left( \Delta_k^{+} \left( j,\tilde{q}\right) \right)_{i_{k-1}}\left( 2k+1 \right)_{i_k} \left( 2k-j+ \alpha -\delta \right)_{i_k}} \right) \nonumber\\
&&\times \left. \sum_{i_{\tau } = i_{\tau -1}}^{m}  \frac{ \left( \Delta_{\tau }^{-} \left( j,\tilde{q}\right) \right)_{i_{\tau }}\left( \Delta_{\tau }^{+} \left( j,\tilde{q}\right) \right)_{i_{\tau }} \left( 2\tau +1  \right)_{i_{\tau -1}} \left( 2\tau -j+ \alpha -\delta \right)_{i_{\tau -1}}}{\left( \Delta_{\tau }^{-} \left( j,\tilde{q}\right) \right)_{i_{\tau -1}}\left( \Delta_{\tau }^{+} \left( j,\tilde{q}\right) \right)_{i_{\tau -1}}\left( 2\tau +1 \right)_{i_\tau } \left( 2\tau -j+ \alpha -\delta \right)_{i_{\tau }}} \eta ^{i_{\tau }}\right\} z^{\tau }  \nonumber 
\end{eqnarray}
where
\begin{equation}
\begin{cases} \tau \geq 2 \cr
\xi = \frac{(1-a)x}{x-a} \cr
z = -\frac{1}{1-a}\xi^2 \cr
\eta = \frac{2-a}{1-a} \xi \cr
\tilde{q} =  -q+(\alpha -\delta -j)\beta \cr
\Delta_k^{\pm} \left( j,\tilde{q}\right) = \frac{\varphi +4(2-a)k \pm \sqrt{\varphi ^2-4(2-a)\tilde{q}}}{2(2-a)}  \cr
\varphi = \beta -\delta -j+a(\alpha -1-j)
\end{cases}\nonumber
\end{equation}
\end{enumerate}
\paragraph{The case of $\delta $, $q$ = fixed values and $\alpha, \beta, \gamma$ = free variables}
Replacing coefficients $a$, $q$, $\alpha $ and $x$ by $1-a$, $-q+\gamma \beta $, $-\alpha+\gamma +\delta $ and $\frac{(1-a)x}{x-a}$ into (\ref{eq:50019})--(\ref{eq:50024c}) with $\delta \rightarrow \alpha -\gamma -j$ in $\varphi $. Multiply $\left(1-\frac{x}{a} \right)^{-\beta }$ and the new (\ref{eq:50019}), (\ref{eq:50020b}), (\ref{eq:50021b}) and (\ref{eq:50022b})  together. 
\begin{enumerate} 
\item As $ \delta =\alpha -\gamma$ and $q= \gamma \beta - q_0^0 $ where $q_0^0=0$,

The eigenfunction is given by
\begin{eqnarray}
&& \left(1-\frac{x}{a} \right)^{-\beta } y(\xi ) \nonumber\\
 &=& \left(1-\frac{x}{a} \right)^{-\beta } Hl\left( 1-a, 0; 0, \beta, \gamma, \delta; \frac{(1-a)x}{x-a} \right) \nonumber\\
&=& \left( 1-\frac{x}{a} \right)^{-\beta } \nonumber
\end{eqnarray}
\item As $ \delta = \alpha -\gamma -1$,

An algebraic equation of degree 2 for the determination of $q$ is given by
\begin{equation}
0 = (1-a)\beta \gamma 
+\prod_{l=0}^{1}\Big( \tilde{q} + l(-\alpha +\beta +\gamma +l+(1-a)( \alpha -2+l))\Big) \nonumber
\end{equation}
The eigenvalue of $q$ is written by $\gamma \beta -q_1^m$ where $m = 0,1 $; $q_{1}^0 < q_{1}^1$. Its eigenfunction is given by
\begin{eqnarray}
&& \left(1-\frac{x}{a} \right)^{-\beta } y(\xi ) \nonumber\\
 &=& \left(1-\frac{x}{a} \right)^{-\beta } Hl\left( 1-a, q_1^m; -1, \beta, \gamma, \delta; \frac{(1-a)x}{x-a} \right) \nonumber\\
&=& \left( 1-\frac{x}{a} \right)^{-\beta } \left\{  1+\frac{ q_1^m}{(2-a)\gamma } \eta \right\} \nonumber
\end{eqnarray}
\item As $ \delta =\alpha -\gamma -2N-2 $ where $N \in \mathbb{N}_{0}$,

An algebraic equation of degree $2N+3$ for the determination of $q$ is given by
\begin{equation}
0 =  \sum_{r=0}^{N+1}\bar{c}\left( r, 2(N-r)+3; 2N+2,\tilde{q}\right)  \nonumber
\end{equation}
The eigenvalue of $q $ is written by $ \gamma \beta -q_{2N+2}^m$ where $m = 0,1,2,\cdots,2N+2 $; $q_{2N+2}^0 < q_{2N+2}^1 < \cdots < q_{2N+2}^{2N+2}$. Its eigenfunction is given by 
\begin{eqnarray} 
&& \left(1-\frac{x}{a} \right)^{-\beta } y(\xi ) \nonumber\\
 &=& \left(1-\frac{x}{a} \right)^{-\beta } Hl\left( 1-a, q_{2N+2}^m; -2N-2, \beta, \gamma, \delta; \frac{(1-a)x}{x-a} \right) \nonumber\\
&=& \left(1-\frac{x}{a} \right)^{-\beta } \sum_{r=0}^{N+1} y_{r}^{2(N+1-r)}\left( 2N+2, q_{2N+2}^m; \xi \right)  
\nonumber 
\end{eqnarray}
\item As $ \delta =\alpha -\gamma -2N-3 $ where $N \in \mathbb{N}_{0}$,

An algebraic equation of degree $2N+4$ for the determination of $q$ is given by
\begin{equation}  
0 = \sum_{r=0}^{N+2}\bar{c}\left( r, 2(N+2-r); 2N+3,\tilde{q}\right) \nonumber
\end{equation}
The eigenvalue of $q$ is written by $\gamma \beta -q_{2N+3}^m$ where $m = 0,1,2,\cdots,2N+3 $; $q_{2N+3}^0 < q_{2N+3}^1 < \cdots < q_{2N+3}^{2N+3}$. Its eigenfunction is given by
\begin{eqnarray} 
&& \left(1-\frac{x}{a} \right)^{-\beta } y(\xi ) \nonumber\\
 &=& \left(1-\frac{x}{a} \right)^{-\beta } Hl\left( 1-a, q_{2N+3}^m; -2N-3, \beta, \gamma, \delta; \frac{(1-a)x}{x-a} \right) \nonumber\\
&=& \left(1-\frac{x}{a} \right)^{-\beta }  \sum_{r=0}^{N+1} y_{r}^{2(N-r)+3} \left( 2N+3,q_{2N+3}^m; \xi\right) \nonumber
\end{eqnarray}
In the above,
\begin{eqnarray}
\bar{c}(0,n;j,\tilde{q})  &=& \frac{\left( \Delta_0^{-} \left( j,\tilde{q}\right) \right)_{n}\left( \Delta_0^{+} \left( j,\tilde{q}\right) \right)_{n}}{\left( 1 \right)_{n} \left( \gamma \right)_{n}} \left( \frac{2-a}{1-a} \right)^{n}\nonumber\\
\bar{c}(1,n;j,\tilde{q}) &=& \left( -\frac{1}{1-a}\right) \sum_{i_0=0}^{n}\frac{\left( i_0 -j\right)\left( i_0+\beta \right) }{\left( i_0+2 \right) \left( i_0+1+ \gamma \right)} \frac{ \left( \Delta_0^{-} \left( j,\tilde{q}\right) \right)_{i_0}\left( \Delta_0^{+} \left( j,\tilde{q}\right) \right)_{i_0}}{\left( 1 \right)_{i_0} \left( \gamma \right)_{i_0}} \nonumber\\
&&\times  \frac{ \left( \Delta_1^{-} \left( j,\tilde{q}\right) \right)_{n}\left( \Delta_1^{+} \left( j,\tilde{q}\right) \right)_{n} \left( 3 \right)_{i_0} \left( 2+ \gamma \right)_{i_0}}{\left( \Delta_1^{-} \left( j,\tilde{q}\right) \right)_{i_0}\left( \Delta_1^{+} \left( j,\tilde{q}\right) \right)_{i_0}\left( 3 \right)_{n} \left( 2+ \gamma \right)_{n}} \left(\frac{2-a}{1-a} \right)^{n }  
\nonumber\\
\bar{c}(\tau ,n;j,\tilde{q}) &=& \left( -\frac{1}{1-a}\right)^{\tau} \sum_{i_0=0}^{n}\frac{\left( i_0 -j\right)\left( i_0+\beta \right) }{\left( i_0+2 \right) \left( i_0+1+ \gamma \right)} \frac{ \left( \Delta_0^{-} \left( j,\tilde{q}\right) \right)_{i_0}\left( \Delta_0^{+} \left( j,\tilde{q}\right) \right)_{i_0}}{\left( 1 \right)_{i_0} \left( \gamma \right)_{i_0}}  \nonumber\\
&&\times \prod_{k=1}^{\tau -1} \left( \sum_{i_k = i_{k-1}}^{n} \frac{\left( i_k+ 2k-j\right)\left( i_k +2k+\beta \right)}{\left( i_k+2k+2 \right) \left( i_k+2k+1+ \gamma \right)} \right. \nonumber\\
&&\times  \left. \frac{ \left( \Delta_k^{-} \left( j,\tilde{q}\right) \right)_{i_k}\left( \Delta_k^{+} \left( j,\tilde{q}\right) \right)_{i_k} \left( 2k+1 \right)_{i_{k-1}} \left( 2k+ \gamma \right)_{i_{k-1}}}{\left( \Delta_k^{-} \left( j,\tilde{q}\right) \right)_{i_{k-1}}\left( \Delta_k^{+} \left( j,\tilde{q}\right) \right)_{i_{k-1}}\left( 2k+1 \right)_{i_k} \left( 2k+ \gamma \right)_{i_k}} \right) \nonumber\\
&&\times \frac{ \left( \Delta_{\tau }^{-} \left( j,\tilde{q}\right) \right)_{n}\left( \Delta_{\tau }^{+} \left( j,\tilde{q}\right) \right)_{n} \left( 2\tau +1 \right)_{i_{\tau -1}} \left( 2\tau + \gamma \right)_{i_{\tau -1}}}{\left( \Delta_{\tau }^{-} \left( j,\tilde{q}\right) \right)_{i_{\tau -1}}\left( \Delta_{\tau }^{+} \left( j,\tilde{q}\right) \right)_{i_{\tau -1}}\left( 2\tau +1 \right)_{n} \left( 2\tau + \gamma  \right)_{n}} \left(\frac{2-a}{1-a}\right)^{n } \nonumber 
\end{eqnarray}
\begin{eqnarray}
y_0^m(j,\tilde{q};\xi) &=& \sum_{i_0=0}^{m} \frac{\left( \Delta_0^{-} \left( j,\tilde{q}\right) \right)_{i_0}\left( \Delta_0^{+} \left( j,\tilde{q}\right) \right)_{i_0}}{\left( 1 \right)_{i_0} \left( \gamma \right)_{i_0}} \eta ^{i_0} \nonumber\\
y_1^m(j,\tilde{q};\xi) &=& \left\{\sum_{i_0=0}^{m}\frac{\left( i_0 -j\right)\left( i_0+\beta \right) }{\left( i_0+2 \right) \left( i_0+1+ \gamma \right)} \frac{ \left( \Delta_0^{-} \left( j,\tilde{q}\right) \right)_{i_0}\left( \Delta_0^{+} \left( j,\tilde{q}\right) \right)_{i_0}}{\left( 1 \right)_{i_0} \left( \gamma \right)_{i_0}} \right. \nonumber\\
&&\times \left. \sum_{i_1 = i_0}^{m} \frac{ \left( \Delta_1^{-} \left( j,\tilde{q}\right) \right)_{i_1}\left( \Delta_1^{+} \left( j,\tilde{q}\right) \right)_{i_1} \left( 3 \right)_{i_0} \left( 2+ \gamma \right)_{i_0}}{\left( \Delta_1^{-} \left( j,\tilde{q}\right) \right)_{i_0}\left( \Delta_1^{+} \left( j,\tilde{q}\right) \right)_{i_0}\left( 3 \right)_{i_1} \left( 2+ \gamma \right)_{i_1}} \eta ^{i_1}\right\} z 
\nonumber\\
y_{\tau }^m(j,\tilde{q};\xi) &=& \left\{ \sum_{i_0=0}^{m} \frac{\left( i_0 -j\right)\left( i_0+\beta \right) }{\left( i_0+2 \right) \left( i_0+1+ \gamma \right)} \frac{ \left( \Delta_0^{-} \left( j,\tilde{q}\right) \right)_{i_0}\left( \Delta_0^{+} \left( j,\tilde{q}\right) \right)_{i_0}}{\left( 1 \right)_{i_0} \left( \gamma \right)_{i_0}} \right.\nonumber\\
&&\times \prod_{k=1}^{\tau -1} \left( \sum_{i_k = i_{k-1}}^{m} \frac{\left( i_k+ 2k-j\right)\left( i_k +2k+\beta \right)}{\left( i_k+2k+2 \right) \left( i_k+2k+1+ \gamma \right)} \right. \nonumber\\
&&\times  \left. \frac{ \left( \Delta_k^{-} \left( j,\tilde{q}\right) \right)_{i_k}\left( \Delta_k^{+} \left( j,\tilde{q}\right) \right)_{i_k} \left( 2k+1 \right)_{i_{k-1}} \left( 2k+ \gamma \right)_{i_{k-1}}}{\left( \Delta_k^{-} \left( j,\tilde{q}\right) \right)_{i_{k-1}}\left( \Delta_k^{+} \left( j,\tilde{q}\right) \right)_{i_{k-1}}\left( 2k+1 \right)_{i_k} \left( 2k+ \gamma \right)_{i_k}} \right) \nonumber\\
&&\times \left. \sum_{i_{\tau } = i_{\tau -1}}^{m}  \frac{ \left( \Delta_{\tau }^{-} \left( j,\tilde{q}\right) \right)_{i_{\tau }}\left( \Delta_{\tau }^{+} \left( j,\tilde{q}\right) \right)_{i_{\tau }} \left( 2\tau +1  \right)_{i_{\tau -1}} \left( 2\tau + \gamma \right)_{i_{\tau -1}}}{\left( \Delta_{\tau }^{-} \left( j,\tilde{q}\right) \right)_{i_{\tau -1}}\left( \Delta_{\tau }^{+} \left( j,\tilde{q}\right) \right)_{i_{\tau -1}}\left( 2\tau +1 \right)_{i_\tau } \left( 2\tau + \gamma \right)_{i_{\tau }}} \eta ^{i_{\tau }}\right\} z^{\tau } \nonumber
\end{eqnarray}
where
\begin{equation}
\begin{cases} \tau \geq 2 \cr
\xi = \frac{(1-a)x}{x-a} \cr
z = -\frac{1}{1-a}\xi^2 \cr
\eta = \frac{2-a}{1-a} \xi \cr
\tilde{q} =  -q+\gamma \beta \cr
\Delta_k^{\pm} \left( j,\tilde{q}\right) = \frac{\varphi +4(2-a)k \pm \sqrt{\varphi ^2-4(2-a)\tilde{q}}}{2(2-a)}  \cr
\varphi = -\alpha +\beta +\gamma +(1-a)(\alpha -1-j)
\end{cases}\nonumber
\end{equation}
\end{enumerate}
\subsection{ \footnotesize ${\displaystyle (1-x)^{1-\delta }\left(1-\frac{x}{a} \right)^{-\beta+\delta -1} Hl\left(1-a, -q+\gamma [(\delta -1)a+\beta -\delta +1]; -\alpha +\gamma +1, \beta -\delta+1, \gamma, 2-\delta; \frac{(1-a)x}{x-a} \right)}$ \normalsize}
\subsubsection{The first species complete polynomial}
\paragraph{The case of $\alpha $, $q$ = fixed values and $\gamma, \beta, \delta $ = free variables }
Replacing coefficients $a$, $q$, $\alpha $, $\beta $, $\delta $ and $x$ by $1-a$, $-q+\gamma [(\delta -1)a+\beta -\delta +1]$, $-\alpha +\gamma +1$, $\beta -\delta+1$, $2-\delta $ and $\frac{(1-a)x}{x-a}$ into (\ref{eq:50019})--(\ref{eq:50024c}). Multiply $(1-x)^{1-\delta }\left(1-\frac{x}{a} \right)^{-\beta+\delta -1}$ and the new (\ref{eq:50019}), (\ref{eq:50020b}), (\ref{eq:50021b}) and (\ref{eq:50022b}) together. 
\begin{enumerate} 
\item As $ \alpha =\gamma +1 $ and $q =\gamma [(\delta -1)a+\beta -\delta +1]-q_0^0$ where $q_0^0=0$,

The eigenfunction is given by
\begin{eqnarray}
&& (1-x)^{1-\delta }\left(1-\frac{x}{a} \right)^{-\beta+\delta -1} y(\xi ) \nonumber\\
 &=& (1-x)^{1-\delta }\left(1-\frac{x}{a} \right)^{-\beta+\delta -1} Hl\bigg( 1-a, 0; 0, \beta -\delta +1, \gamma , 2-\delta; \frac{(1-a)x}{x-a} \bigg) \nonumber\\
&=& (1-x)^{1-\delta }\left(1-\frac{x}{a} \right)^{-\beta+\delta -1}  \nonumber
\end{eqnarray}
\item As $ \alpha =\gamma +2$,

An algebraic equation of degree 2 for the determination of $q$ is given by
\begin{equation}
0 = (1-a)(\beta -\delta +1)\gamma 
+\prod_{l=0}^{1}\Big( \tilde{q} + l(\beta -2+l+a(\gamma -\delta +1+l))\Big) \nonumber
\end{equation}
The eigenvalue of $ q$ is written by $ \gamma [(\delta -1)a+\beta -\delta +1]-q_1^m$ where $m = 0,1 $; $q_{1}^0 < q_{1}^1$. Its eigenfunction is given by
\begin{eqnarray}
&& (1-x)^{1-\delta }\left(1-\frac{x}{a} \right)^{-\beta+\delta -1} y(\xi ) \nonumber\\
 &=& (1-x)^{1-\delta }\left(1-\frac{x}{a} \right)^{-\beta+\delta -1} Hl\bigg( 1-a, q_1^m; -1, \beta -\delta +1, \gamma , 2-\delta; \frac{(1-a)x}{x-a} \bigg) \nonumber\\
&=& (1-x)^{1-\delta }\left(1-\frac{x}{a} \right)^{-\beta+\delta -1} \left\{  1+\frac{ q_1^m}{(2-a)\gamma } \eta \right\} \nonumber  
\end{eqnarray}
\item As $ \alpha = \gamma +2N+3 $ where $N \in \mathbb{N}_{0}$,

An algebraic equation of degree $2N+3$ for the determination of $q$ is given by
\begin{equation}
0 = \sum_{r=0}^{N+1}\bar{c}\left( r, 2(N-r)+3; 2N+2,\tilde{q}\right)  \nonumber
\end{equation}
The eigenvalue of $q$ is written by $\gamma [(\delta -1)a+\beta -\delta +1]- q_{2N+2}^m$ where $m = 0,1,2,\cdots,2N+2 $; $q_{2N+2}^0 < q_{2N+2}^1 < \cdots < q_{2N+2}^{2N+2}$. Its eigenfunction is given by 
\begin{eqnarray} 
&& (1-x)^{1-\delta }\left(1-\frac{x}{a} \right)^{-\beta+\delta -1} y(\xi ) \nonumber\\
 &=& (1-x)^{1-\delta }\left(1-\frac{x}{a} \right)^{-\beta+\delta -1} Hl\bigg( 1-a, q_{2N+2}^m; -2N-2, \beta -\delta +1, \gamma , 2-\delta; \frac{(1-a)x}{x-a} \bigg) \nonumber\\
&=& (1-x)^{1-\delta }\left(1-\frac{x}{a} \right)^{-\beta+\delta -1} \sum_{r=0}^{N+1} y_{r}^{2(N+1-r)}\left( 2N+2, q_{2N+2}^m; \xi \right)  
\nonumber 
\end{eqnarray}
\item As $ \alpha =\gamma +2N+4 $ where $N \in \mathbb{N}_{0}$,

An algebraic equation of degree $2N+4$ for the determination of $q$ is given by
\begin{equation}  
0 = \sum_{r=0}^{N+2}\bar{c}\left( r, 2(N+2-r); 2N+3,\tilde{q}\right) \nonumber
\end{equation}
The eigenvalue of $q$ is written by $\gamma [(\delta -1)a+\beta -\delta +1]- q_{2N+3}^m$ where $m = 0,1,2,\cdots,2N+3 $; $q_{2N+3}^0 < q_{2N+3}^1 < \cdots < q_{2N+3}^{2N+3}$. Its eigenfunction is given by
\begin{eqnarray} 
&& (1-x)^{1-\delta }\left(1-\frac{x}{a} \right)^{-\beta+\delta -1} y(\xi ) \nonumber\\
 &=& (1-x)^{1-\delta }\left(1-\frac{x}{a} \right)^{-\beta+\delta -1} Hl\bigg( 1-a, q_{2N+3}^m; -2N-3, \beta -\delta +1, \gamma , 2-\delta; \frac{(1-a)x}{x-a} \bigg) \nonumber\\
&=& (1-x)^{1-\delta }\left(1-\frac{x}{a} \right)^{-\beta+\delta -1} \sum_{r=0}^{N+1} y_{r}^{2(N-r)+3} \left( 2N+3,q_{2N+3}^m;\xi \right) \nonumber
\end{eqnarray}
In the above,
\begin{eqnarray}
\bar{c}(0,n;j,\tilde{q})  &=& \frac{\left( \Delta_0^{-} \left( j,\tilde{q}\right) \right)_{n}\left( \Delta_0^{+} \left( j,\tilde{q}\right) \right)_{n}}{\left( 1 \right)_{n} \left( \gamma \right)_{n}} \left( \frac{2-a}{1-a} \right)^{n}\nonumber\\
\bar{c}(1,n;j,\tilde{q}) &=& \left( -\frac{1}{1-a}\right) \sum_{i_0=0}^{n}\frac{\left( i_0 -j\right)\left( i_0+1+\beta -\delta \right) }{\left( i_0+2 \right) \left( i_0+1+ \gamma \right)} \frac{ \left( \Delta_0^{-} \left( j,\tilde{q}\right) \right)_{i_0}\left( \Delta_0^{+} \left( j,\tilde{q}\right) \right)_{i_0}}{\left( 1 \right)_{i_0} \left( \gamma \right)_{i_0}} \nonumber\\
&&\times  \frac{ \left( \Delta_1^{-} \left( j,\tilde{q}\right) \right)_{n}\left( \Delta_1^{+} \left( j,\tilde{q}\right) \right)_{n} \left( 3 \right)_{i_0} \left( 2+ \gamma \right)_{i_0}}{\left( \Delta_1^{-} \left( j,\tilde{q}\right) \right)_{i_0}\left( \Delta_1^{+} \left( j,\tilde{q}\right) \right)_{i_0}\left( 3 \right)_{n} \left( 2+ \gamma \right)_{n}} \left(\frac{2-a}{1-a} \right)^{n }  
\nonumber\\
\bar{c}(\tau ,n;j,\tilde{q}) &=& \left( -\frac{1}{1-a}\right)^{\tau} \sum_{i_0=0}^{n}\frac{\left( i_0 -j\right)\left( i_0+1+\beta -\delta \right) }{\left( i_0+2 \right) \left( i_0+1+ \gamma \right)}  \frac{ \left( \Delta_0^{-} \left( j,\tilde{q}\right) \right)_{i_0}\left( \Delta_0^{+} \left( j,\tilde{q}\right) \right)_{i_0}}{\left( 1 \right)_{i_0} \left( \gamma \right)_{i_0}}  \nonumber\\
&&\times \prod_{k=1}^{\tau -1} \left( \sum_{i_k = i_{k-1}}^{n} \frac{\left( i_k+ 2k-j\right)\left( i_k +2k+1+\beta -\delta \right)}{\left( i_k+2k+2 \right) \left( i_k+2k+1+ \gamma \right)} \right. \nonumber\\
&&\times  \left. \frac{ \left( \Delta_k^{-} \left( j,\tilde{q}\right) \right)_{i_k}\left( \Delta_k^{+} \left( j,\tilde{q}\right) \right)_{i_k} \left( 2k+1 \right)_{i_{k-1}} \left( 2k+ \gamma \right)_{i_{k-1}}}{\left( \Delta_k^{-} \left( j,\tilde{q}\right) \right)_{i_{k-1}}\left( \Delta_k^{+} \left( j,\tilde{q}\right) \right)_{i_{k-1}}\left( 2k+1 \right)_{i_k} \left( 2k+ \gamma \right)_{i_k}} \right) \nonumber\\
&&\times \frac{ \left( \Delta_{\tau }^{-} \left( j,\tilde{q}\right) \right)_{n}\left( \Delta_{\tau }^{+} \left( j,\tilde{q}\right) \right)_{n} \left( 2\tau +1 \right)_{i_{\tau -1}} \left( 2\tau + \gamma \right)_{i_{\tau -1}}}{\left( \Delta_{\tau }^{-} \left( j,\tilde{q}\right) \right)_{i_{\tau -1}}\left( \Delta_{\tau }^{+} \left( j,\tilde{q}\right) \right)_{i_{\tau -1}}\left( 2\tau +1 \right)_{n} \left( 2\tau + \gamma  \right)_{n}} \left(\frac{2-a}{1-a}\right)^{n } \nonumber 
\end{eqnarray}
\begin{eqnarray}
y_0^m(j,\tilde{q};\xi) &=& \sum_{i_0=0}^{m} \frac{\left( \Delta_0^{-} \left( j,\tilde{q}\right) \right)_{i_0}\left( \Delta_0^{+} \left( j,\tilde{q}\right) \right)_{i_0}}{\left( 1 \right)_{i_0} \left( \gamma \right)_{i_0}} \eta ^{i_0} \nonumber\\
y_1^m(j,\tilde{q};\xi) &=& \left\{\sum_{i_0=0}^{m}\frac{\left( i_0 -j\right)\left( i_0+1+\beta -\delta \right) }{\left( i_0+2 \right) \left( i_0+1+ \gamma \right)} \frac{ \left( \Delta_0^{-} \left( j,\tilde{q}\right) \right)_{i_0}\left( \Delta_0^{+} \left( j,\tilde{q}\right) \right)_{i_0}}{\left( 1 \right)_{i_0} \left( \gamma \right)_{i_0}} \right. \nonumber\\
&&\times \left. \sum_{i_1 = i_0}^{m} \frac{ \left( \Delta_1^{-} \left( j,\tilde{q}\right) \right)_{i_1}\left( \Delta_1^{+} \left( j,\tilde{q}\right) \right)_{i_1} \left( 3 \right)_{i_0} \left( 2+ \gamma \right)_{i_0}}{\left( \Delta_1^{-} \left( j,\tilde{q}\right) \right)_{i_0}\left( \Delta_1^{+} \left( j,\tilde{q}\right) \right)_{i_0}\left( 3 \right)_{i_1} \left( 2+ \gamma \right)_{i_1}} \eta ^{i_1}\right\} z 
\nonumber\\
y_{\tau }^m(j,\tilde{q};\xi) &=& \left\{ \sum_{i_0=0}^{m} \frac{\left( i_0 -j\right)\left( i_0+1+\beta -\delta \right) }{\left( i_0+2 \right) \left( i_0+1+ \gamma \right)} \frac{ \left( \Delta_0^{-} \left( j,\tilde{q}\right) \right)_{i_0}\left( \Delta_0^{+} \left( j,\tilde{q}\right) \right)_{i_0}}{\left( 1 \right)_{i_0} \left( \gamma \right)_{i_0}} \right.\nonumber\\
&&\times \prod_{k=1}^{\tau -1} \left( \sum_{i_k = i_{k-1}}^{m} \frac{\left( i_k+ 2k-j\right)\left( i_k +2k+1+\beta -\delta \right)}{\left( i_k+2k+2 \right) \left( i_k+2k+1+ \gamma \right)} \right. \nonumber\\
&&\times  \left. \frac{ \left( \Delta_k^{-} \left( j,\tilde{q}\right) \right)_{i_k}\left( \Delta_k^{+} \left( j,\tilde{q}\right) \right)_{i_k} \left( 2k+1 \right)_{i_{k-1}} \left( 2k+ \gamma \right)_{i_{k-1}}}{\left( \Delta_k^{-} \left( j,\tilde{q}\right) \right)_{i_{k-1}}\left( \Delta_k^{+} \left( j,\tilde{q}\right) \right)_{i_{k-1}}\left( 2k+1 \right)_{i_k} \left( 2k+ \gamma \right)_{i_k}} \right) \nonumber\\
&&\times \left. \sum_{i_{\tau } = i_{\tau -1}}^{m}  \frac{ \left( \Delta_{\tau }^{-} \left( j,\tilde{q}\right) \right)_{i_{\tau }}\left( \Delta_{\tau }^{+} \left( j,\tilde{q}\right) \right)_{i_{\tau }} \left( 2\tau +1  \right)_{i_{\tau -1}} \left( 2\tau + \gamma \right)_{i_{\tau -1}}}{\left( \Delta_{\tau }^{-} \left( j,\tilde{q}\right) \right)_{i_{\tau -1}}\left( \Delta_{\tau }^{+} \left( j,\tilde{q}\right) \right)_{i_{\tau -1}}\left( 2\tau +1 \right)_{i_\tau } \left( 2\tau + \gamma \right)_{i_{\tau }}} \eta ^{i_{\tau }}\right\} z^{\tau } \nonumber
\end{eqnarray}
where
\begin{equation}
\begin{cases} \tau \geq 2 \cr
\xi = \frac{(1-a)x}{x-a} \cr
z = -\frac{1}{1-a}\xi^2 \cr
\eta = \frac{(2-a)}{1-a} \xi \cr
\tilde{q} = -q+\gamma [(\delta -1)a+\beta -\delta +1]\cr
\Delta_k^{\pm} \left( j,\tilde{q}\right) = \frac{\varphi +4(2-a)k \pm \sqrt{\varphi ^2-4(2-a)\tilde{q}}}{2(2-a)}  \cr
\varphi = \beta -1-j+a(\gamma -\delta +1)
\end{cases}\nonumber
\end{equation}
\end{enumerate}
\paragraph{The case of $\gamma $, $q$ = fixed values and $\alpha , \beta, \delta $ = free variables }
Replacing coefficients $a$, $q$, $\alpha $, $\beta $, $\delta $ and $x$ by $1-a$, $-q+\gamma [(\delta -1)a+\beta -\delta +1]$, $-\alpha +\gamma +1$, $\beta -\delta+1$, $2-\delta $ and $\frac{(1-a)x}{x-a}$ into (\ref{eq:50019})--(\ref{eq:50024c}). Replacing $\gamma $ by $\alpha -1-j$  into the new (\ref{eq:50023a})--(\ref{eq:50024c}). Multiply $(1-x)^{1-\delta }\left(1-\frac{x}{a} \right)^{-\beta+\delta -1}$ and the new (\ref{eq:50019}), (\ref{eq:50020b}), (\ref{eq:50021b}) and (\ref{eq:50022b}) together. 
\begin{enumerate} 
\item As $ \gamma =\alpha -1$ and $q=(\alpha -1) [(\delta -1)a+\beta -\delta +1]- q_0^0 $ where $q_0^0=0$,

The eigenfunction is given by
\begin{eqnarray}
&& (1-x)^{1-\delta }\left(1-\frac{x}{a} \right)^{-\beta+\delta -1} y(\xi ) \nonumber\\
 &=& (1-x)^{1-\delta }\left(1-\frac{x}{a} \right)^{-\beta+\delta -1} Hl\bigg( 1-a, 0; 0, \beta -\delta +1, \gamma , 2-\delta; \frac{(1-a)x}{x-a} \bigg) \nonumber\\
&=& (1-x)^{1-\delta }\left(1-\frac{x}{a} \right)^{-\beta+\delta -1} \nonumber 
\end{eqnarray}
\item As $ \gamma =\alpha -2$,

An algebraic equation of degree 2 for the determination of $q$ is given by
\begin{equation}
0 = (1-a)(\beta -\delta +1)\gamma 
+\prod_{l=0}^{1}\Big( -q+ (\alpha -2)[(\delta -1)a+\beta -\delta +1]+ l(\beta -2 +l+a(\gamma -\delta +1+l))\Big) \nonumber
\end{equation}
The eigenvalue of $q$ is written by $(\alpha -2)[(\delta -1)a+\beta -\delta +1] -q_1^m$ where $m = 0,1 $; $q_{1}^0 < q_{1}^1$. Its eigenfunction is given by
\begin{eqnarray}
&& (1-x)^{1-\delta }\left(1-\frac{x}{a} \right)^{-\beta+\delta -1} y(\xi ) \nonumber\\
 &=& (1-x)^{1-\delta }\left(1-\frac{x}{a} \right)^{-\beta+\delta -1} Hl\bigg( 1-a, q_1^m; -1, \beta -\delta +1, \gamma , 2-\delta; \frac{(1-a)x}{x-a} \bigg) \nonumber\\
&=& (1-x)^{1-\delta }\left(1-\frac{x}{a} \right)^{-\beta+\delta -1} \left\{  1+\frac{ q_1^m}{(2-a)(\alpha -2)} \eta \right\} \nonumber 
\end{eqnarray}
\item As $ \gamma =\alpha -2N-3 $ where $N \in \mathbb{N}_{0}$,

An algebraic equation of degree $2N+3$ for the determination of $q$ is given by
\begin{equation}
0 = \sum_{r=0}^{N+1}\bar{c}\left( r, 2(N-r)+3; 2N+2,-q+(\alpha -2N-3)[(\delta -1)a+\beta -\delta +1] \right)  \nonumber
\end{equation}
The eigenvalue of $q$ is written by $(\alpha -2N-3)[(\delta -1)a+\beta -\delta +1]-q_{2N+2}^m$ where $m = 0,1,2,\cdots,2N+2 $; $q_{2N+2}^0 < q_{2N+2}^1 < \cdots < q_{2N+2}^{2N+2}$. Its eigenfunction is given by 
\begin{eqnarray} 
&& (1-x)^{1-\delta }\left(1-\frac{x}{a} \right)^{-\beta+\delta -1} y(\xi ) \nonumber\\
 &=& (1-x)^{1-\delta }\left(1-\frac{x}{a} \right)^{-\beta+\delta -1} Hl\bigg( 1-a, q_{2N+2}^m; -2N-2, \beta -\delta +1, \gamma , 2-\delta; \frac{(1-a)x}{x-a} \bigg) \nonumber\\
&=& (1-x)^{1-\delta }\left(1-\frac{x}{a} \right)^{-\beta+\delta -1} \sum_{r=0}^{N+1} y_{r}^{2(N+1-r)}\left( 2N+2, q_{2N+2}^m; \xi \right)  
\nonumber
\end{eqnarray}
\item As $ \gamma =\alpha -2N-4 $ where $N \in \mathbb{N}_{0}$,

An algebraic equation of degree $2N+4$ for the determination of $q$ is given by
\begin{equation}  
0 = \sum_{r=0}^{N+2}\bar{c}\left( r, 2(N+2-r); 2N+3,-q+(\alpha -2N-4)[(\delta -1)a+\beta -\delta +1] \right) \nonumber
\end{equation}
The eigenvalue of $q$ is written by $(\alpha -2N-4)[(\delta -1)a+\beta -\delta +1]-q_{2N+3}^m$ where $m = 0,1,2,\cdots,2N+3 $; $q_{2N+3}^0 < q_{2N+3}^1 < \cdots < q_{2N+3}^{2N+3}$. Its eigenfunction is given by
\begin{eqnarray} 
&& (1-x)^{1-\delta }\left(1-\frac{x}{a} \right)^{-\beta+\delta -1} y(\xi ) \nonumber\\
 &=& (1-x)^{1-\delta }\left(1-\frac{x}{a} \right)^{-\beta+\delta -1} Hl\bigg( 1-a, q_{2N+3}^m; -2N-3, \beta -\delta +1, \gamma , 2-\delta; \frac{(1-a)x}{x-a} \bigg) \nonumber\\
&=& (1-x)^{1-\delta }\left(1-\frac{x}{a} \right)^{-\beta+\delta -1}  \sum_{r=0}^{N+1} y_{r}^{2(N-r)+3} \left( 2N+3,q_{2N+3}^m;\xi\right) \nonumber
\end{eqnarray}
In the above,
\begin{eqnarray}
\bar{c}(0,n;j,\tilde{q})  &=& \frac{\left( \Delta_0^{-} \left( j,\tilde{q}\right) \right)_{n}\left( \Delta_0^{+} \left( j,\tilde{q}\right) \right)_{n}}{\left( 1 \right)_{n} \left( \alpha -1-j \right)_{n}} \left( \frac{2-a}{1-a} \right)^{n}\nonumber\\
\bar{c}(1,n;j,\tilde{q}) &=& \left( -\frac{1}{1-a}\right) \sum_{i_0=0}^{n}\frac{\left( i_0 -j\right)\left( i_0+1+\beta -\delta \right) }{\left( i_0+2 \right) \left( i_0 -j+ \alpha \right)} \frac{ \left( \Delta_0^{-} \left( j,\tilde{q}\right) \right)_{i_0}\left( \Delta_0^{+} \left( j,\tilde{q}\right) \right)_{i_0}}{\left( 1 \right)_{i_0} \left( \alpha -1-j \right)_{i_0}} \nonumber\\
&&\times  \frac{ \left( \Delta_1^{-} \left( j,\tilde{q}\right) \right)_{n}\left( \Delta_1^{+} \left( j,\tilde{q}\right) \right)_{n} \left( 3 \right)_{i_0} \left( 1-j+ \alpha \right)_{i_0}}{\left( \Delta_1^{-} \left( j,\tilde{q}\right) \right)_{i_0}\left( \Delta_1^{+} \left( j,\tilde{q}\right) \right)_{i_0}\left( 3 \right)_{n} \left( 1-j+ \alpha \right)_{n}} \left(\frac{2-a}{1-a} \right)^{n }  
\nonumber\\
\bar{c}(\tau ,n;j,\tilde{q}) &=& \left( -\frac{1}{1-a}\right)^{\tau} \sum_{i_0=0}^{n}\frac{\left( i_0 -j\right)\left( i_0+1+\beta -\delta \right) }{\left( i_0+2 \right) \left( i_0 -j+ \alpha \right)} \frac{ \left( \Delta_0^{-} \left( j,\tilde{q}\right) \right)_{i_0}\left( \Delta_0^{+} \left( j,\tilde{q}\right) \right)_{i_0}}{\left( 1 \right)_{i_0} \left( \alpha -1-j \right)_{i_0}} \nonumber\\
&&\times \prod_{k=1}^{\tau -1} \left( \sum_{i_k = i_{k-1}}^{n} \frac{\left( i_k+ 2k-j\right)\left( i_k +2k+1+\beta -\delta \right)}{\left( i_k+2k+2 \right) \left( i_k+2k-j+ \alpha \right)} \right. \nonumber\\
&&\times  \left. \frac{ \left( \Delta_k^{-} \left( j,\tilde{q}\right) \right)_{i_k}\left( \Delta_k^{+} \left( j,\tilde{q}\right) \right)_{i_k} \left( 2k+1 \right)_{i_{k-1}} \left( 2k-1-j+ \alpha \right)_{i_{k-1}}}{\left( \Delta_k^{-} \left( j,\tilde{q}\right) \right)_{i_{k-1}}\left( \Delta_k^{+} \left( j,\tilde{q}\right) \right)_{i_{k-1}}\left( 2k+1 \right)_{i_k} \left( 2k-1-j+ \alpha \right)_{i_k}} \right) \nonumber\\
&&\times \frac{ \left( \Delta_{\tau }^{-} \left( j,\tilde{q}\right) \right)_{n}\left( \Delta_{\tau }^{+} \left( j,\tilde{q}\right) \right)_{n} \left( 2\tau +1 \right)_{i_{\tau -1}} \left( 2\tau -1-j+ \alpha \right)_{i_{\tau -1}}}{\left( \Delta_{\tau }^{-} \left( j,\tilde{q}\right) \right)_{i_{\tau -1}}\left( \Delta_{\tau }^{+} \left( j,\tilde{q}\right) \right)_{i_{\tau -1}}\left( 2\tau +1 \right)_{n} \left( 2\tau -1-j+ \alpha \right)_{n}} \left(\frac{2-a}{1-a}\right)^{n } \nonumber 
\end{eqnarray}
\begin{eqnarray}
y_0^m(j,\tilde{q};\xi) &=& \sum_{i_0=0}^{m} \frac{\left( \Delta_0^{-} \left( j,\tilde{q}\right) \right)_{i_0}\left( \Delta_0^{+} \left( j,\tilde{q}\right) \right)_{i_0}}{\left( 1 \right)_{i_0} \left( \alpha -1-j \right)_{i_0}} \eta ^{i_0} \nonumber\\
y_1^m(j,\tilde{q};\xi) &=& \left\{\sum_{i_0=0}^{m}\frac{\left( i_0 -j\right)\left( i_0+1+\beta -\delta  \right) }{\left( i_0+2 \right) \left( i_0 -j+ \alpha \right)} \frac{ \left( \Delta_0^{-} \left( j,\tilde{q}\right) \right)_{i_0}\left( \Delta_0^{+} \left( j,\tilde{q}\right) \right)_{i_0}}{\left( 1 \right)_{i_0} \left( \alpha -1-j \right)_{i_0}} \right. \nonumber\\
&&\times \left. \sum_{i_1 = i_0}^{m} \frac{ \left( \Delta_1^{-} \left( j,\tilde{q}\right) \right)_{i_1}\left( \Delta_1^{+} \left( j,\tilde{q}\right) \right)_{i_1} \left( 3 \right)_{i_0} \left( 1-j+ \alpha \right)_{i_0}}{\left( \Delta_1^{-} \left( j,\tilde{q}\right) \right)_{i_0}\left( \Delta_1^{+} \left( j,\tilde{q}\right) \right)_{i_0}\left( 3 \right)_{i_1} \left( 1-j+ \alpha \right)_{i_1}} \eta ^{i_1}\right\} z 
\nonumber\\
y_{\tau }^m(j,\tilde{q};\xi) &=& \left\{ \sum_{i_0=0}^{m} \frac{\left( i_0 -j\right)\left( i_0+1+\beta -\delta  \right) }{\left( i_0+2 \right) \left( i_0-j+ \alpha \right)} \frac{ \left( \Delta_0^{-} \left( j,\tilde{q}\right) \right)_{i_0}\left( \Delta_0^{+} \left( j,\tilde{q}\right) \right)_{i_0}}{\left( 1 \right)_{i_0} \left( \alpha -1-j \right)_{i_0}} \right.\nonumber\\
&&\times \prod_{k=1}^{\tau -1} \left( \sum_{i_k = i_{k-1}}^{m} \frac{\left( i_k+ 2k-j\right)\left( i_k +2k+1+\beta -\delta \right)}{\left( i_k+2k+2 \right) \left( i_k+2k-j+ \alpha \right)} \right. \nonumber\\
&&\times  \left. \frac{ \left( \Delta_k^{-} \left( j,\tilde{q}\right) \right)_{i_k}\left( \Delta_k^{+} \left( j,\tilde{q}\right) \right)_{i_k} \left( 2k+1 \right)_{i_{k-1}} \left( 2k-1-j+ \alpha \right)_{i_{k-1}}}{\left( \Delta_k^{-} \left( j,\tilde{q}\right) \right)_{i_{k-1}}\left( \Delta_k^{+} \left( j,\tilde{q}\right) \right)_{i_{k-1}}\left( 2k+1 \right)_{i_k} \left( 2k-1-j+ \alpha \right)_{i_k}} \right) \nonumber\\
&&\times \left. \sum_{i_{\tau } = i_{\tau -1}}^{m}  \frac{ \left( \Delta_{\tau }^{-} \left( j,\tilde{q}\right) \right)_{i_{\tau }}\left( \Delta_{\tau }^{+} \left( j,\tilde{q}\right) \right)_{i_{\tau }} \left( 2\tau +1  \right)_{i_{\tau -1}} \left( 2\tau -1-j+ \alpha \right)_{i_{\tau -1}}}{\left( \Delta_{\tau }^{-} \left( j,\tilde{q}\right) \right)_{i_{\tau -1}}\left( \Delta_{\tau }^{+} \left( j,\tilde{q}\right) \right)_{i_{\tau -1}}\left( 2\tau +1 \right)_{i_\tau } \left( 2\tau -1-j+\alpha \right)_{i_{\tau }}} \eta ^{i_{\tau }}\right\} z^{\tau } \nonumber
\end{eqnarray}
where
\begin{equation}
\begin{cases} \tau \geq 2 \cr
\xi = \frac{(1-a)x}{x-a} \cr
z = -\frac{1}{1-a}\xi^2 \cr
\eta = \frac{(2-a)}{1-a} \xi \cr
\tilde{q} = -q+(\alpha -1-j)[(\delta -1)a+\beta -\delta +1]\cr
\Delta_k^{\pm} \left( j,\tilde{q}\right) = \frac{\varphi +4(2-a)k \pm \sqrt{\varphi ^2-4(2-a)\tilde{q}}}{2(2-a)}  \cr
\varphi = \beta -1-j+(1-a)(\alpha -\delta -j)
\end{cases}\nonumber
\end{equation}
\end{enumerate}
\subsection{ \footnotesize ${\displaystyle x^{-\alpha } Hl\left(\frac{a-1}{a}, \frac{-q+\alpha (\delta a+\beta -\delta )}{a}; \alpha, \alpha -\gamma +1, \delta , \alpha -\beta +1; \frac{x-1}{x} \right)}$\normalsize}
\subsubsection{The first species complete polynomial}
Replacing coefficients $a$, $q$, $\alpha $, $\beta $, $\gamma $, $\delta $ and $x$ by $\frac{a-1}{a}$, $\frac{-q+\alpha (\delta a+\beta -\delta )}{a}$, $\alpha -\gamma +1$, $\alpha $, $\delta $, $\alpha -\beta +1$ and $\frac{x-1}{x}$ into (\ref{eq:50019})--(\ref{eq:50024c}). Multiply $x^{-\alpha }$ and the new (\ref{eq:50019}), (\ref{eq:50020b}), (\ref{eq:50021b}) and (\ref{eq:50022b}) together.\footnote{I treat $\alpha $, $\beta $ and $\delta$ as free variables and fixed values of $\gamma $ and $q$.}
\begin{enumerate} 
\item As $\gamma =\alpha +1 $ and $q=\alpha (\delta a+\beta -\delta )-a q_0^0 $  where $q_0^0=0$,

The eigenfunction is given by
\begin{eqnarray}
&& x^{-\alpha } y(\xi ) \nonumber\\
 &=& x^{-\alpha } Hl\left(\frac{a-1}{a}, 0; 0, \alpha, \delta , \alpha -\beta +1; \frac{x-1}{x} \right) \nonumber\\
&=& x^{-\alpha } \nonumber
\end{eqnarray}
\item As $\gamma = \alpha +2$,

An algebraic equation of degree 2 for the determination of $q$ is given by
\begin{equation}
0 = \frac{(a-1)\alpha (\alpha +2)}{a}
+\prod_{l=0}^{1}\Big( \tilde{q}+ l( \beta -2+l+a( \alpha -\beta +\delta +l))\Big) \nonumber
\end{equation}
The eigenvalue of $q$ is written by $\alpha (\delta a+\beta -\delta ) -a q_1^m$ where $m = 0,1 $; $q_{1}^0 < q_{1}^1$. Its eigenfunction is given by
\begin{eqnarray}
&& x^{-\alpha } y(\xi ) \nonumber\\
 &=& x^{-\alpha } Hl\left(\frac{a-1}{a}, q_1^m; -1, \alpha, \delta , \alpha -\beta +1; \frac{x-1}{x} \right) \nonumber\\
&=& x^{-\alpha } \left\{ 1+\frac{ a q_1^m}{(2a-1)\delta  } \eta \right\} \nonumber  
\end{eqnarray}
\item As $\gamma =\alpha  +2N+3 $ where $N \in \mathbb{N}_{0}$,

An algebraic equation of degree $2N+3$ for the determination of $q$ is given by
\begin{equation}
0 = \sum_{r=0}^{N+1}\bar{c}\left( r, 2(N-r)+3; 2N+2,\tilde{q}\right)  
\nonumber
\end{equation}
The eigenvalue of $q$ is written by $\alpha (\delta a+\beta -\delta ) -a q_{2N+2}^m$ where $m = 0,1,2,\cdots,2N+2 $; $q_{2N+2}^0 < q_{2N+2}^1 < \cdots < q_{2N+2}^{2N+2}$. Its eigenfunction is given by 
\begin{eqnarray} 
&& x^{-\alpha } y(\xi ) \nonumber\\
 &=& x^{-\alpha } Hl\left(\frac{a-1}{a}, q_{2N+2}^m; -2N-2, \alpha , \delta , \alpha -\beta +1; \frac{x-1}{x} \right) \nonumber\\
&=& x^{-\alpha }  \sum_{r=0}^{N+1} y_{r}^{2(N+1-r)}\left( 2N+2, q_{2N+2}^m; \xi \right)  
\nonumber
\end{eqnarray}
\item As $\gamma =\alpha  +2N+4$ where $N \in \mathbb{N}_{0}$,

An algebraic equation of degree $2N+4$ for the determination of $q$ is given by
\begin{equation}  
0 = \sum_{r=0}^{N+2}\bar{c}\left( r, 2(N+2-r); 2N+3,\tilde{q}\right) 
\nonumber
\end{equation}
The eigenvalue of $q$ is written by $\alpha (\delta a+\beta -\delta ) -a q_{2N+3}^m$ where $m = 0,1,2,\cdots,2N+3 $; $q_{2N+3}^0 < q_{2N+3}^1 < \cdots < q_{2N+3}^{2N+3}$. Its eigenfunction is given by
\begin{eqnarray} 
&& x^{-\alpha } y(\xi ) \nonumber\\
 &=& x^{-\alpha } Hl\left(\frac{a-1}{a}, q_{2N+3}^m; -2N-3, \alpha , \delta , \alpha -\beta +1; \frac{x-1}{x} \right) \nonumber\\
&=& x^{-\alpha }  \sum_{r=0}^{N+1} y_{r}^{2(N-r)+3} \left( 2N+3,q_{2N+3}^m;\xi\right) \nonumber
\end{eqnarray}
In the above,
\begin{eqnarray}
\bar{c}(0,n;j,\tilde{q})  &=& \frac{\left( \Delta_0^{-} \left( j,\tilde{q}\right) \right)_{n}\left( \Delta_0^{+} \left( j,\tilde{q}\right) \right)_{n}}{\left( 1 \right)_{n} \left( \delta  \right)_{n}} \left( \frac{2a-1}{a-1}\right)^{n}\nonumber\\
\bar{c}(1,n;j,\tilde{q}) &=& \left( -\frac{a}{a-1}\right) \sum_{i_0=0}^{n}\frac{\left( i_0 -j\right)\left( i_0+\alpha  \right) }{\left( i_0+2 \right) \left( i_0+1+ \delta \right)} \frac{ \left( \Delta_0^{-} \left( j,\tilde{q}\right) \right)_{i_0}\left( \Delta_0^{+} \left( j,\tilde{q}\right) \right)_{i_0}}{\left( 1 \right)_{i_0} \left( \delta \right)_{i_0}} \nonumber\\
&&\times  \frac{ \left( \Delta_1^{-} \left( j,\tilde{q}\right) \right)_{n}\left( \Delta_1^{+} \left( j,\tilde{q}\right) \right)_{n} \left( 3 \right)_{i_0} \left( 2+ \delta \right)_{i_0}}{\left( \Delta_1^{-} \left( j,\tilde{q}\right) \right)_{i_0}\left( \Delta_1^{+} \left( j,\tilde{q}\right) \right)_{i_0}\left( 3 \right)_{n} \left( 2+ \delta \right)_{n}} \left(\frac{2a-1}{a-1} \right)^{n }  
\nonumber\\
\bar{c}(\tau ,n;j,\tilde{q}) &=& \left( -\frac{a}{a-1}\right)^{\tau} \sum_{i_0=0}^{n}\frac{\left( i_0 -j\right)\left( i_0+\alpha \right) }{\left( i_0+2 \right) \left( i_0+1+ \delta \right)} \frac{ \left( \Delta_0^{-} \left( j,\tilde{q}\right) \right)_{i_0}\left( \Delta_0^{+} \left( j,\tilde{q}\right) \right)_{i_0}}{\left( 1 \right)_{i_0} \left( \delta \right)_{i_0}} \nonumber\\
&&\times \prod_{k=1}^{\tau -1} \left( \sum_{i_k = i_{k-1}}^{n} \frac{\left( i_k+ 2k-j\right)\left( i_k +2k+\alpha \right)}{\left( i_k+2k+2 \right) \left( i_k+2k+1+ \delta \right)} \right. \nonumber\\
&&\times  \left. \frac{ \left( \Delta_k^{-} \left( j,\tilde{q}\right) \right)_{i_k}\left( \Delta_k^{+} \left( j,\tilde{q}\right) \right)_{i_k} \left( 2k+1 \right)_{i_{k-1}} \left( 2k+ \delta \right)_{i_{k-1}}}{\left( \Delta_k^{-} \left( j,\tilde{q}\right) \right)_{i_{k-1}}\left( \Delta_k^{+} \left( j,\tilde{q}\right) \right)_{i_{k-1}}\left( 2k+1 \right)_{i_k} \left( 2k+ \delta \right)_{i_k}} \right) \nonumber\\
&&\times \frac{ \left( \Delta_{\tau }^{-} \left( j,\tilde{q}\right) \right)_{n}\left( \Delta_{\tau }^{+} \left( j,\tilde{q}\right) \right)_{n} \left( 2\tau +1 \right)_{i_{\tau -1}} \left( 2\tau + \delta \right)_{i_{\tau -1}}}{\left( \Delta_{\tau }^{-} \left( j,\tilde{q}\right) \right)_{i_{\tau -1}}\left( \Delta_{\tau }^{+} \left( j,\tilde{q}\right) \right)_{i_{\tau -1}}\left( 2\tau +1 \right)_{n} \left( 2\tau + \delta  \right)_{n}} \left(\frac{2a-1}{a-1}\right)^{n } \nonumber 
\end{eqnarray}
\begin{eqnarray}
y_0^m(j,\tilde{q};\xi) &=& \sum_{i_0=0}^{m} \frac{\left( \Delta_0^{-} \left( j,\tilde{q}\right) \right)_{i_0}\left( \Delta_0^{+} \left( j,\tilde{q}\right) \right)_{i_0}}{\left( 1 \right)_{i_0} \left( \delta \right)_{i_0}} \eta ^{i_0} \nonumber\\
y_1^m(j,\tilde{q};\xi) &=& \left\{\sum_{i_0=0}^{m}\frac{\left( i_0 -j\right)\left( i_0+\alpha \right) }{\left( i_0+2 \right) \left( i_0+1+ \delta \right)} \frac{ \left( \Delta_0^{-} \left( j,\tilde{q}\right) \right)_{i_0}\left( \Delta_0^{+} \left( j,\tilde{q}\right) \right)_{i_0}}{\left( 1 \right)_{i_0} \left( \delta \right)_{i_0}} \right. \nonumber\\
&&\times \left. \sum_{i_1 = i_0}^{m} \frac{ \left( \Delta_1^{-} \left( j,\tilde{q}\right) \right)_{i_1}\left( \Delta_1^{+} \left( j,\tilde{q}\right) \right)_{i_1} \left( 3 \right)_{i_0} \left( 2+ \delta \right)_{i_0}}{\left( \Delta_1^{-} \left( j,\tilde{q}\right) \right)_{i_0}\left( \Delta_1^{+} \left( j,\tilde{q}\right) \right)_{i_0}\left( 3 \right)_{i_1} \left( 2+ \delta \right)_{i_1}} \eta ^{i_1}\right\} z 
\nonumber\\
y_{\tau }^m(j,\tilde{q};\xi) &=& \left\{ \sum_{i_0=0}^{m} \frac{\left( i_0 -j\right)\left( i_0+\alpha \right) }{\left( i_0+2 \right) \left( i_0+1+ \delta \right)} \frac{ \left( \Delta_0^{-} \left( j,\tilde{q}\right) \right)_{i_0}\left( \Delta_0^{+} \left( j,\tilde{q}\right) \right)_{i_0}}{\left( 1 \right)_{i_0} \left( \delta \right)_{i_0}} \right.\nonumber\\
&&\times \prod_{k=1}^{\tau -1} \left( \sum_{i_k = i_{k-1}}^{m} \frac{\left( i_k+ 2k-j\right)\left( i_k +2k+\alpha \right)}{\left( i_k+2k+2 \right) \left( i_k+2k+1+ \delta \right)} \right. \nonumber\\
&&\times  \left. \frac{ \left( \Delta_k^{-} \left( j,\tilde{q}\right) \right)_{i_k}\left( \Delta_k^{+} \left( j,\tilde{q}\right) \right)_{i_k} \left( 2k+1 \right)_{i_{k-1}} \left( 2k+ \delta \right)_{i_{k-1}}}{\left( \Delta_k^{-} \left( j,\tilde{q}\right) \right)_{i_{k-1}}\left( \Delta_k^{+} \left( j,\tilde{q}\right) \right)_{i_{k-1}}\left( 2k+1 \right)_{i_k} \left( 2k+ \delta \right)_{i_k}} \right) \nonumber\\
&&\times \left. \sum_{i_{\tau } = i_{\tau -1}}^{m}  \frac{ \left( \Delta_{\tau }^{-} \left( j,\tilde{q}\right) \right)_{i_{\tau }}\left( \Delta_{\tau }^{+} \left( j,\tilde{q}\right) \right)_{i_{\tau }} \left( 2\tau +1  \right)_{i_{\tau -1}} \left( 2\tau + \delta \right)_{i_{\tau -1}}}{\left( \Delta_{\tau }^{-} \left( j,\tilde{q}\right) \right)_{i_{\tau -1}}\left( \Delta_{\tau }^{+} \left( j,\tilde{q}\right) \right)_{i_{\tau -1}}\left( 2\tau +1 \right)_{i_\tau } \left( 2\tau + \delta \right)_{i_{\tau }}} \eta ^{i_{\tau }}\right\} z^{\tau } \nonumber
\end{eqnarray}
where
\begin{equation}
\begin{cases} \tau \geq 2 \cr
\xi= \frac{x-1}{x} \cr
z = \frac{-a}{a-1}\xi^2 \cr
\eta = \frac{2a-1}{a-1} \xi \cr
\tilde{q} =  \frac{-q+\alpha (\delta a+\beta -\delta )}{a} \cr
\Delta_k^{\pm} \left( j,\tilde{q}\right) = \frac{a\varphi +4(2a-1)k \pm \sqrt{(a\varphi )^2-4(2a-1)\tilde{q}}}{2(2a-1)}  \cr
\varphi = \beta -1 -j+\frac{a-1}{a}(\alpha -\beta +\delta )
\end{cases}\nonumber
\end{equation}
\end{enumerate}
\subsection{ \footnotesize ${\displaystyle \left(\frac{x-a}{1-a} \right)^{-\alpha } Hl\left(a, q-(\beta -\delta )\alpha ; \alpha , -\beta+\gamma +\delta , \delta , \gamma; \frac{a(x-1)}{x-a} \right)}$\normalsize}
\subsubsection{The first species complete polynomial}
\paragraph{The case of $\beta $, $q$ = fixed values and $\alpha , \gamma, \delta $ = free variables }
Replacing coefficients $q$, $\alpha $, $\beta $, $\gamma $, $\delta $ and $x$ by $q-(\beta -\delta )\alpha $, $-\beta+\gamma +\delta $, $\alpha $, $\delta $,  $\gamma $ and $\frac{a(x-1)}{x-a}$ into (\ref{eq:50019})--(\ref{eq:50024c}). Multiply $\left(\frac{x-a}{1-a} \right)^{-\alpha }$ and the new (\ref{eq:50019}), (\ref{eq:50020b}), (\ref{eq:50021b}) and (\ref{eq:50022b}) together. 
\begin{enumerate} 
\item As $ \beta =\gamma +\delta $ and $q =\gamma \alpha +q_0^0 $ where $q_0^0=0$,

The eigenfunction is given by
\begin{eqnarray}
&& \left(\frac{x-a}{1-a} \right)^{-\alpha } y(\xi ) \nonumber\\
 &=& \left(\frac{x-a}{1-a} \right)^{-\alpha } Hl\left( a, 0; 0, \alpha , \delta , \gamma; \frac{a(x-1)}{x-a} \right) \nonumber\\
&=& \left(\frac{x-a}{1-a} \right)^{-\alpha } \nonumber
\end{eqnarray}
\item As $ \beta =\gamma +\delta +1 $,

An algebraic equation of degree 2 for the determination of $q$ is given by
\begin{equation}
0 = a\alpha \delta  
+\prod_{l=0}^{1}\Big( q-(\gamma +1)\alpha + l(\alpha  -\gamma  -1+l+a(\gamma +\delta -1+l))\Big) \nonumber
\end{equation}
The eigenvalue of $q$ is written by $(\gamma +1)\alpha +q_1^m$ where $m = 0,1 $; $q_{1}^0 < q_{1}^1$. Its eigenfunction is given by
\begin{eqnarray}
&& \left(\frac{x-a}{1-a} \right)^{-\alpha } y(\xi ) \nonumber\\
 &=& \left(\frac{x-a}{1-a} \right)^{-\alpha } Hl\left( a, q_1^m; -1, \alpha , \delta , \gamma; \frac{a(x-1)}{x-a} \right) \nonumber\\
&=& \left(\frac{x-a}{1-a} \right)^{-\alpha } \left\{ 1+\frac{ q_1^m}{(1+a)\delta } \eta \right\} \nonumber
\end{eqnarray}
\item As $\beta =\gamma +\delta +2N+2 $ where $N \in \mathbb{N}_{0}$,

An algebraic equation of degree $2N+3$ for the determination of $q$ is given by
\begin{equation}
0 = \sum_{r=0}^{N+1}\bar{c}\left( r, 2(N-r)+3; 2N+2,q-(\gamma +2N+2)\alpha\right)  \nonumber
\end{equation}
The eigenvalue of $q$ is written by $(\gamma +2N+2)\alpha + q_{2N+2}^m$ where $m = 0,1,2,\cdots,2N+2 $; $q_{2N+2}^0 < q_{2N+2}^1 < \cdots < q_{2N+2}^{2N+2}$. Its eigenfunction is given by 
\begin{eqnarray} 
&& \left(\frac{x-a}{1-a} \right)^{-\alpha } y(\xi ) \nonumber\\
 &=& \left(\frac{x-a}{1-a} \right)^{-\alpha } Hl\left( a, q_{2N+2}^m; -2N-2, \alpha , \delta , \gamma; \frac{a(x-1)}{x-a} \right) \nonumber\\
&=& \left(\frac{x-a}{1-a} \right)^{-\alpha } \sum_{r=0}^{N+1} y_{r}^{2(N+1-r)}\left( 2N+2, q_{2N+2}^m; \xi \right)  
\nonumber
\end{eqnarray}
\item As $ \beta = \gamma +\delta +2N+3 $ where $N \in \mathbb{N}_{0}$,

An algebraic equation of degree $2N+4$ for the determination of $q$ is given by
\begin{equation}  
0 = \sum_{r=0}^{N+2}\bar{c}\left( r, 2(N+2-r); 2N+3,q-(\gamma +2N+3)\alpha \right) \nonumber
\end{equation}
The eigenvalue of $q$ is written by $(\gamma +2N+3)\alpha + q_{2N+3}^m$ where $m = 0,1,2,\cdots,2N+3 $; $q_{2N+3}^0 < q_{2N+3}^1 < \cdots < q_{2N+3}^{2N+3}$. Its eigenfunction is given by
\begin{eqnarray} 
&& \left(\frac{x-a}{1-a} \right)^{-\alpha } y(\xi ) \nonumber\\
 &=& \left(\frac{x-a}{1-a} \right)^{-\alpha } Hl\left( a, q_{2N+3}^m; -2N-3, \alpha , \delta , \gamma; \frac{a(x-1)}{x-a} \right) \nonumber\\
&=& \left(\frac{x-a}{1-a} \right)^{-\alpha }  \sum_{r=0}^{N+1} y_{r}^{2(N-r)+3} \left( 2N+3,q_{2N+3}^m;\xi \right) \nonumber
\end{eqnarray}
In the above,
\begin{eqnarray}
\bar{c}(0,n;j,\tilde{q})  &=& \frac{\left( \Delta_0^{-} \left( j,\tilde{q}\right) \right)_{n}\left( \Delta_0^{+} \left( j,\tilde{q}\right) \right)_{n}}{\left( 1 \right)_{n} \left( \delta  \right)_{n}} \left( \frac{1+a}{a} \right)^{n}\nonumber\\
\bar{c}(1,n;j,\tilde{q}) &=& \left( -\frac{1}{a}\right) \sum_{i_0=0}^{n}\frac{\left( i_0 -j\right)\left( i_0+\alpha  \right) }{\left( i_0+2 \right) \left( i_0+1+ \delta \right)} \frac{ \left( \Delta_0^{-} \left( j,\tilde{q}\right) \right)_{i_0}\left( \Delta_0^{+} \left( j,\tilde{q}\right) \right)_{i_0}}{\left( 1 \right)_{i_0} \left( \delta \right)_{i_0}} \nonumber\\
&&\times  \frac{ \left( \Delta_1^{-} \left( j,\tilde{q}\right) \right)_{n}\left( \Delta_1^{+} \left( j,\tilde{q}\right) \right)_{n} \left( 3 \right)_{i_0} \left( 2+ \delta \right)_{i_0}}{\left( \Delta_1^{-} \left( j,\tilde{q}\right) \right)_{i_0}\left( \Delta_1^{+} \left( j,\tilde{q}\right) \right)_{i_0}\left( 3 \right)_{n} \left( 2+ \delta \right)_{n}} \left(\frac{1+a}{a} \right)^{n }  
\nonumber\\
\bar{c}(\tau ,n;j,\tilde{q}) &=& \left( -\frac{1}{a}\right)^{\tau} \sum_{i_0=0}^{n}\frac{\left( i_0 -j\right)\left( i_0+\alpha \right) }{\left( i_0+2 \right) \left( i_0+1+ \delta \right)}  \frac{ \left( \Delta_0^{-} \left( j,\tilde{q}\right) \right)_{i_0}\left( \Delta_0^{+} \left( j,\tilde{q}\right) \right)_{i_0}}{\left( 1 \right)_{i_0} \left( \delta \right)_{i_0}} \nonumber\\
&&\times \prod_{k=1}^{\tau -1} \left( \sum_{i_k = i_{k-1}}^{n} \frac{\left( i_k+ 2k-j\right)\left( i_k +2k+\alpha \right)}{\left( i_k+2k+2 \right) \left( i_k+2k+1+ \delta \right)} \right. \nonumber\\
&&\times  \left. \frac{ \left( \Delta_k^{-} \left( j,\tilde{q}\right) \right)_{i_k}\left( \Delta_k^{+} \left( j,\tilde{q}\right) \right)_{i_k} \left( 2k+1 \right)_{i_{k-1}} \left( 2k+ \delta \right)_{i_{k-1}}}{\left( \Delta_k^{-} \left( j,\tilde{q}\right) \right)_{i_{k-1}}\left( \Delta_k^{+} \left( j,\tilde{q}\right) \right)_{i_{k-1}}\left( 2k+1 \right)_{i_k} \left( 2k+ \delta \right)_{i_k}} \right) \nonumber\\
&&\times \frac{ \left( \Delta_{\tau }^{-} \left( j,\tilde{q}\right) \right)_{n}\left( \Delta_{\tau }^{+} \left( j,\tilde{q}\right) \right)_{n} \left( 2\tau +1 \right)_{i_{\tau -1}} \left( 2\tau + \delta \right)_{i_{\tau -1}}}{\left( \Delta_{\tau }^{-} \left( j,\tilde{q}\right) \right)_{i_{\tau -1}}\left( \Delta_{\tau }^{+} \left( j,\tilde{q}\right) \right)_{i_{\tau -1}}\left( 2\tau +1 \right)_{n} \left( 2\tau + \delta  \right)_{n}} \left(\frac{1+a}{a}\right)^{n } \nonumber 
\end{eqnarray}
\begin{eqnarray}
y_0^m(j,\tilde{q};\xi) &=& \sum_{i_0=0}^{m} \frac{\left( \Delta_0^{-} \left( j,\tilde{q}\right) \right)_{i_0}\left( \Delta_0^{+} \left( j,\tilde{q}\right) \right)_{i_0}}{\left( 1 \right)_{i_0} \left( \delta \right)_{i_0}} \eta ^{i_0} \nonumber\\
y_1^m(j,\tilde{q};\xi) &=& \left\{\sum_{i_0=0}^{m}\frac{\left( i_0 -j\right)\left( i_0+\alpha \right) }{\left( i_0+2 \right) \left( i_0+1+ \delta \right)} \frac{ \left( \Delta_0^{-} \left( j,\tilde{q}\right) \right)_{i_0}\left( \Delta_0^{+} \left( j,\tilde{q}\right) \right)_{i_0}}{\left( 1 \right)_{i_0} \left( \delta \right)_{i_0}} \right. \nonumber\\
&&\times \left. \sum_{i_1 = i_0}^{m} \frac{ \left( \Delta_1^{-} \left( j,\tilde{q}\right) \right)_{i_1}\left( \Delta_1^{+} \left( j,\tilde{q}\right) \right)_{i_1} \left( 3 \right)_{i_0} \left( 2+ \delta \right)_{i_0}}{\left( \Delta_1^{-} \left( j,\tilde{q}\right) \right)_{i_0}\left( \Delta_1^{+} \left( j,\tilde{q}\right) \right)_{i_0}\left( 3 \right)_{i_1} \left( 2+ \delta \right)_{i_1}} \eta ^{i_1}\right\} z 
\nonumber
\end{eqnarray}
\begin{eqnarray}
y_{\tau }^m(j,\tilde{q};\xi) &=& \left\{ \sum_{i_0=0}^{m} \frac{\left( i_0 -j\right)\left( i_0+\alpha \right) }{\left( i_0+2 \right) \left( i_0+1+ \delta \right)} \frac{ \left( \Delta_0^{-} \left( j,\tilde{q}\right) \right)_{i_0}\left( \Delta_0^{+} \left( j,\tilde{q}\right) \right)_{i_0}}{\left( 1 \right)_{i_0} \left( \delta \right)_{i_0}} \right.\nonumber\\
&&\times \prod_{k=1}^{\tau -1} \left( \sum_{i_k = i_{k-1}}^{m} \frac{\left( i_k+ 2k-j\right)\left( i_k +2k+\alpha \right)}{\left( i_k+2k+2 \right) \left( i_k+2k+1+ \delta \right)} \right. \nonumber\\
&&\times  \left. \frac{ \left( \Delta_k^{-} \left( j,\tilde{q}\right) \right)_{i_k}\left( \Delta_k^{+} \left( j,\tilde{q}\right) \right)_{i_k} \left( 2k+1 \right)_{i_{k-1}} \left( 2k+ \delta \right)_{i_{k-1}}}{\left( \Delta_k^{-} \left( j,\tilde{q}\right) \right)_{i_{k-1}}\left( \Delta_k^{+} \left( j,\tilde{q}\right) \right)_{i_{k-1}}\left( 2k+1 \right)_{i_k} \left( 2k+ \delta \right)_{i_k}} \right) \nonumber\\
&&\times \left. \sum_{i_{\tau } = i_{\tau -1}}^{m}  \frac{ \left( \Delta_{\tau }^{-} \left( j,\tilde{q}\right) \right)_{i_{\tau }}\left( \Delta_{\tau }^{+} \left( j,\tilde{q}\right) \right)_{i_{\tau }} \left( 2\tau +1  \right)_{i_{\tau -1}} \left( 2\tau + \delta \right)_{i_{\tau -1}}}{\left( \Delta_{\tau }^{-} \left( j,\tilde{q}\right) \right)_{i_{\tau -1}}\left( \Delta_{\tau }^{+} \left( j,\tilde{q}\right) \right)_{i_{\tau -1}}\left( 2\tau +1 \right)_{i_\tau } \left( 2\tau + \delta \right)_{i_{\tau }}} \eta ^{i_{\tau }}\right\} z^{\tau } \nonumber
\end{eqnarray}
where
\begin{equation}
\begin{cases} \tau \geq 2 \cr
\xi= \frac{a(x-1)}{x-a} \cr
z = -\frac{1}{a}\xi^2 \cr
\eta = \frac{(1+a)}{a} \xi \cr
\tilde{q} =  q-(\beta -\delta )\alpha \cr
\Delta_k^{\pm} \left( j,\tilde{q}\right) = \frac{\varphi +4(1+a)k \pm \sqrt{\varphi ^2-4(1+a)\tilde{q}}}{2(1+a)}  \cr
\varphi = \alpha -\gamma -j+a(\gamma +\delta -1)
\end{cases}\nonumber
\end{equation}
\end{enumerate}
\paragraph{The case of $\gamma $, $q$ = fixed values and $\alpha , \beta, \delta $ = free variables }
Replacing coefficients $q$, $\alpha $, $\beta $, $\gamma $, $\delta $ and $x$ by $q-(\beta -\delta )\alpha $, $-\beta+\gamma +\delta $, $\alpha $, $\delta $,  $\gamma $ and $\frac{a(x-1)}{x-a}$ into (\ref{eq:50019})--(\ref{eq:50024c}) with $\gamma \rightarrow \beta -\delta -j$ in $\varphi $. Multiply $\left(\frac{x-a}{1-a} \right)^{-\alpha }$ and the new (\ref{eq:50019}), (\ref{eq:50020b}), (\ref{eq:50021b}) and (\ref{eq:50022b}) together.
\begin{enumerate} 
\item As $ \gamma =\beta -\delta $ and $q=(\beta -\delta )\alpha + q_0^0$ where $q_0^0=0$,

The eigenfunction is given by
\begin{eqnarray}
&& \left(\frac{x-a}{1-a} \right)^{-\alpha } y(\xi ) \nonumber\\
 &=& \left(\frac{x-a}{1-a} \right)^{-\alpha } Hl\left(a, 0 ; 0, \alpha , \delta , \gamma; \frac{a(x-1)}{x-a} \right) \nonumber\\
&=& \left(\frac{x-a}{1-a} \right)^{-\alpha } \nonumber
\end{eqnarray}
\item As $ \gamma =\beta -\delta -1$,

An algebraic equation of degree 2 for the determination of $q$ is given by
\begin{equation}
0 = a\alpha \delta  
+\prod_{l=0}^{1}\Big( q-(\beta -\delta )\alpha + l(\alpha -\beta +\delta +l+a(\beta -2+l))\Big) \nonumber
\end{equation}
The eigenvalue of $q$ is written by $(\beta -\delta )\alpha +q_1^m$ where $m = 0,1 $; $q_{1}^0 < q_{1}^1$. Its eigenfunction is given by
\begin{eqnarray}
&& \left(\frac{x-a}{1-a} \right)^{-\alpha } y(\xi ) \nonumber\\
 &=& \left(\frac{x-a}{1-a} \right)^{-\alpha } Hl\left(a, q_1^m ; -1, \alpha , \delta , \gamma; \frac{a(x-1)}{x-a} \right) \nonumber\\
&=& \left(\frac{x-a}{1-a} \right)^{-\alpha } \left\{  1+\frac{ q_1^m}{(1+a)\delta } \eta \right\} \nonumber  
\end{eqnarray}
\item As $ \gamma =\beta -\delta -2N-2 $ where $N \in \mathbb{N}_{0}$,

An algebraic equation of degree $2N+3$ for the determination of $q$ is given by
\begin{equation}
0 =  \sum_{r=0}^{N+1}\bar{c}\left( r, 2(N-r)+3; 2N+2,\tilde{q}\right)  \nonumber
\end{equation}
The eigenvalue of $q$ is written by $(\beta -\delta )\alpha + q_{2N+2}^m$ where $m = 0,1,2,\cdots,2N+2 $; $q_{2N+2}^0 < q_{2N+2}^1 < \cdots < q_{2N+2}^{2N+2}$. Its eigenfunction is given by 
\begin{eqnarray} 
&& \left(\frac{x-a}{1-a} \right)^{-\alpha } y(\xi ) \nonumber\\
 &=& \left(\frac{x-a}{1-a} \right)^{-\alpha } Hl\left(a, q_{2N+2}^m ; -2N-2, \alpha , \delta , \gamma; \frac{a(x-1)}{x-a} \right) \nonumber\\
&=& \left(\frac{x-a}{1-a} \right)^{-\alpha } \sum_{r=0}^{N+1} y_{r}^{2(N+1-r)}\left( 2N+2, q_{2N+2}^m; \xi \right)
\nonumber
\end{eqnarray}
\item As $ \gamma =\beta -\delta -2N-3$ where $N \in \mathbb{N}_{0}$,

An algebraic equation of degree $2N+4$ for the determination of $q$ is given by
\begin{equation}  
0 = \sum_{r=0}^{N+2}\bar{c}\left( r, 2(N+2-r); 2N+3,\tilde{q}\right) \nonumber
\end{equation}
The eigenvalue of $q$ is written by $(\beta -\delta )\alpha + q_{2N+3}^m$ where $m = 0,1,2,\cdots,2N+3 $; $q_{2N+3}^0 < q_{2N+3}^1 < \cdots < q_{2N+3}^{2N+3}$. Its eigenfunction is given by
\begin{eqnarray} 
&& \left(\frac{x-a}{1-a} \right)^{-\alpha } y(\xi ) \nonumber\\
 &=& \left(\frac{x-a}{1-a} \right)^{-\alpha } Hl\left(a, q_{2N+3}^m ; -2N-3, \alpha , \delta , \gamma; \frac{a(x-1)}{x-a} \right) \nonumber\\
&=& \left(\frac{x-a}{1-a} \right)^{-\alpha } \sum_{r=0}^{N+1} y_{r}^{2(N-r)+3} \left( 2N+3,q_{2N+3}^m;\xi\right) \nonumber
\end{eqnarray}
In the above,
\begin{eqnarray}
\bar{c}(0,n;j,\tilde{q})  &=& \frac{\left( \Delta_0^{-} \left( j,\tilde{q}\right) \right)_{n}\left( \Delta_0^{+} \left( j,\tilde{q}\right) \right)_{n}}{\left( 1 \right)_{n} \left( \delta  \right)_{n}} \left( \frac{1+a}{a} \right)^{n}\nonumber\\
\bar{c}(1,n;j,\tilde{q}) &=& \left( -\frac{1}{a}\right) \sum_{i_0=0}^{n}\frac{\left( i_0 -j\right)\left( i_0+\alpha \right) }{\left( i_0+2 \right) \left( i_0+1+ \delta \right)} \frac{ \left( \Delta_0^{-} \left( j,\tilde{q}\right) \right)_{i_0}\left( \Delta_0^{+} \left( j,\tilde{q}\right) \right)_{i_0}}{\left( 1 \right)_{i_0} \left( \delta \right)_{i_0}} \nonumber\\
&&\times \frac{ \left( \Delta_1^{-} \left( j,\tilde{q}\right) \right)_{n}\left( \Delta_1^{+} \left( j,\tilde{q}\right) \right)_{n} \left( 3 \right)_{i_0} \left( 2+ \delta \right)_{i_0}}{\left( \Delta_1^{-} \left( j,\tilde{q}\right) \right)_{i_0}\left( \Delta_1^{+} \left( j,\tilde{q}\right) \right)_{i_0}\left( 3 \right)_{n} \left( 2+ \delta \right)_{n}} \left(\frac{1+a}{a} \right)^{n }  
\nonumber\\
\bar{c}(\tau ,n;j,\tilde{q}) &=& \left( -\frac{1}{a}\right)^{\tau} \sum_{i_0=0}^{n}\frac{\left( i_0 -j\right)\left( i_0+\alpha \right) }{\left( i_0+2 \right) \left( i_0+1+ \delta \right)} \frac{ \left( \Delta_0^{-} \left( j,\tilde{q}\right) \right)_{i_0}\left( \Delta_0^{+} \left( j,\tilde{q}\right) \right)_{i_0}}{\left( 1 \right)_{i_0} \left( \delta \right)_{i_0}}  \nonumber\\
&&\times \prod_{k=1}^{\tau -1} \left( \sum_{i_k = i_{k-1}}^{n} \frac{\left( i_k+ 2k-j\right)\left( i_k +2k+\alpha \right)}{\left( i_k+2k+2 \right) \left( i_k+2k+1+ \delta \right)} \right. \nonumber\\
&&\times  \left. \frac{ \left( \Delta_k^{-} \left( j,\tilde{q}\right) \right)_{i_k}\left( \Delta_k^{+} \left( j,\tilde{q}\right) \right)_{i_k} \left( 2k+1 \right)_{i_{k-1}} \left( 2k+ \delta \right)_{i_{k-1}}}{\left( \Delta_k^{-} \left( j,\tilde{q}\right) \right)_{i_{k-1}}\left( \Delta_k^{+} \left( j,\tilde{q}\right) \right)_{i_{k-1}}\left( 2k+1 \right)_{i_k} \left( 2k+ \delta \right)_{i_k}} \right) \nonumber\\
&&\times \frac{ \left( \Delta_{\tau }^{-} \left( j,\tilde{q}\right) \right)_{n}\left( \Delta_{\tau }^{+} \left( j,\tilde{q}\right) \right)_{n} \left( 2\tau +1 \right)_{i_{\tau -1}} \left( 2\tau + \delta \right)_{i_{\tau -1}}}{\left( \Delta_{\tau }^{-} \left( j,\tilde{q}\right) \right)_{i_{\tau -1}}\left( \Delta_{\tau }^{+} \left( j,\tilde{q}\right) \right)_{i_{\tau -1}}\left( 2\tau +1 \right)_{n} \left( 2\tau + \delta  \right)_{n}} \left(\frac{1+a}{a}\right)^{n } \nonumber 
\end{eqnarray}
\begin{eqnarray}
y_0^m(j,\tilde{q};\xi) &=& \sum_{i_0=0}^{m} \frac{\left( \Delta_0^{-} \left( j,\tilde{q}\right) \right)_{i_0}\left( \Delta_0^{+} \left( j,\tilde{q}\right) \right)_{i_0}}{\left( 1 \right)_{i_0} \left( \delta \right)_{i_0}} \eta ^{i_0} \nonumber\\
y_1^m(j,\tilde{q};\xi) &=& \left\{\sum_{i_0=0}^{m}\frac{\left( i_0 -j\right)\left( i_0+\alpha  \right) }{\left( i_0+2 \right) \left( i_0+1+ \delta \right)} \frac{ \left( \Delta_0^{-} \left( j,\tilde{q}\right) \right)_{i_0}\left( \Delta_0^{+} \left( j,\tilde{q}\right) \right)_{i_0}}{\left( 1 \right)_{i_0} \left( \delta \right)_{i_0}} \right. \nonumber\\
&&\times \left. \sum_{i_1 = i_0}^{m} \frac{ \left( \Delta_1^{-} \left( j,\tilde{q}\right) \right)_{i_1}\left( \Delta_1^{+} \left( j,\tilde{q}\right) \right)_{i_1} \left( 3 \right)_{i_0} \left( 2+ \delta \right)_{i_0}}{\left( \Delta_1^{-} \left( j,\tilde{q}\right) \right)_{i_0}\left( \Delta_1^{+} \left( j,\tilde{q}\right) \right)_{i_0}\left( 3 \right)_{i_1} \left( 2+ \delta \right)_{i_1}} \eta ^{i_1}\right\} z 
\nonumber\\
y_{\tau }^m(j,\tilde{q};\xi) &=& \left\{ \sum_{i_0=0}^{m} \frac{\left( i_0 -j\right)\left( i_0+\alpha \right) }{\left( i_0+2 \right) \left( i_0+1+ \delta \right)} \frac{ \left( \Delta_0^{-} \left( j,\tilde{q}\right) \right)_{i_0}\left( \Delta_0^{+} \left( j,\tilde{q}\right) \right)_{i_0}}{\left( 1 \right)_{i_0} \left( \delta \right)_{i_0}} \right.\nonumber\\
&&\times \prod_{k=1}^{\tau -1} \left( \sum_{i_k = i_{k-1}}^{m} \frac{\left( i_k+ 2k-j\right)\left( i_k +2k+\alpha \right)}{\left( i_k+2k+2 \right) \left( i_k+2k+1+ \delta \right)} \right. \nonumber\\
&&\times  \left. \frac{ \left( \Delta_k^{-} \left( j,\tilde{q}\right) \right)_{i_k}\left( \Delta_k^{+} \left( j,\tilde{q}\right) \right)_{i_k} \left( 2k+1 \right)_{i_{k-1}} \left( 2k+ \delta \right)_{i_{k-1}}}{\left( \Delta_k^{-} \left( j,\tilde{q}\right) \right)_{i_{k-1}}\left( \Delta_k^{+} \left( j,\tilde{q}\right) \right)_{i_{k-1}}\left( 2k+1 \right)_{i_k} \left( 2k+ \delta \right)_{i_k}} \right) \nonumber\\
&&\times \left. \sum_{i_{\tau } = i_{\tau -1}}^{m}  \frac{ \left( \Delta_{\tau }^{-} \left( j,\tilde{q}\right) \right)_{i_{\tau }}\left( \Delta_{\tau }^{+} \left( j,\tilde{q}\right) \right)_{i_{\tau }} \left( 2\tau +1  \right)_{i_{\tau -1}} \left( 2\tau + \delta \right)_{i_{\tau -1}}}{\left( \Delta_{\tau }^{-} \left( j,\tilde{q}\right) \right)_{i_{\tau -1}}\left( \Delta_{\tau }^{+} \left( j,\tilde{q}\right) \right)_{i_{\tau -1}}\left( 2\tau +1 \right)_{i_\tau } \left( 2\tau + \delta \right)_{i_{\tau }}} \eta ^{i_{\tau }}\right\} z^{\tau } \nonumber 
\end{eqnarray}
where
\begin{equation}
\begin{cases} \tau \geq 2 \cr
\xi= \frac{a(x-1)}{x-a} \cr
z = -\frac{1}{a}\xi^2 \cr
\eta = \frac{(1+a)}{a} \xi \cr
\tilde{q} =  q-(\beta -\delta )\alpha \cr
\Delta_k^{\pm} \left( j,\tilde{q}\right) = \frac{\varphi +4(1+a)k \pm \sqrt{\varphi ^2-4(1+a)\tilde{q}}}{2(1+a)}  \cr
\varphi = \alpha -\beta +\delta +a(\beta -1-j)
\end{cases}\nonumber
\end{equation}
\end{enumerate}
\paragraph{The case of $\delta $, $q$ = fixed values and $\alpha , \beta, \gamma $ = free variables } 
Replacing coefficients $q$, $\alpha $, $\beta $, $\gamma $, $\delta $ and $x$ by $q-(\beta -\delta )\alpha $, $-\beta+\gamma +\delta $, $\alpha $, $\delta $,  $\gamma $ and $\frac{a(x-1)}{x-a}$ into (\ref{eq:50019})--(\ref{eq:50024c}). Replacing $\delta $ by $\beta -\gamma -j$ into the new (\ref{eq:50023a})--(\ref{eq:50024c}). Multiply $\left(\frac{x-a}{1-a} \right)^{-\alpha }$ and the new (\ref{eq:50019}), (\ref{eq:50020b}), (\ref{eq:50021b}) and (\ref{eq:50022b}) together.
\begin{enumerate} 
\item As $ \delta =\beta -\gamma $ and $q = \gamma \alpha +q_0^0 $ where $q_0^0=0$,

The eigenfunction is given by
\begin{eqnarray}
&& \left(\frac{x-a}{1-a} \right)^{-\alpha } y(\xi ) \nonumber\\
 &=& \left(\frac{x-a}{1-a} \right)^{-\alpha } Hl\left( a, 0; 0, \alpha , \delta , \gamma; \frac{a(x-1)}{x-a} \right) \nonumber\\
&=& \left(\frac{x-a}{1-a} \right)^{-\alpha } \nonumber
\end{eqnarray}
\item As $ \delta  =\beta -\gamma -1$,

An algebraic equation of degree 2 for the determination of $q$ is given by
\begin{equation}
0 = a\alpha (\beta -\gamma -1) 
+\prod_{l=0}^{1}\Big( q-(\gamma +1)\alpha + l(\alpha -\gamma -1+l+a(\beta -2+l))\Big) \nonumber
\end{equation}
The eigenvalue of $q$ is written by $(\gamma +1)\alpha +q_1^m$ where $m = 0,1 $; $q_{1}^0 < q_{1}^1$. Its eigenfunction is given by
\begin{eqnarray}
&& \left(\frac{x-a}{1-a} \right)^{-\alpha } y(\xi ) \nonumber\\
 &=& \left(\frac{x-a}{1-a} \right)^{-\alpha } Hl\left( a, q_1^m; -1, \alpha , \delta , \gamma; \frac{a(x-1)}{x-a} \right) \nonumber\\
&=& \left(\frac{x-a}{1-a} \right)^{-\alpha } \left\{  1+\frac{ q_1^m}{(1+a)(\beta -\gamma -1)} \eta \right\} \nonumber  
\end{eqnarray}
\item As $ \delta  =\beta -\gamma -2N-2 $ where $N \in \mathbb{N}_{0}$,

An algebraic equation of degree $2N+3$ for the determination of $q$ is given by
\begin{equation}
0 =   \sum_{r=0}^{N+1}\bar{c}\left( r, 2(N-r)+3; 2N+2,q-(\gamma +2N+2)\alpha\right)  \nonumber
\end{equation}
The eigenvalue of $q$ is written by $(\gamma +2N+2)\alpha + q_{2N+2}^m$ where $m = 0,1,2,\cdots,2N+2 $; $q_{2N+2}^0 < q_{2N+2}^1 < \cdots < q_{2N+2}^{2N+2}$. Its eigenfunction is given by 
\begin{eqnarray} 
&& \left(\frac{x-a}{1-a} \right)^{-\alpha } y(\xi ) \nonumber\\
 &=& \left(\frac{x-a}{1-a} \right)^{-\alpha } Hl\left( a, q_{2N+2}^m; -2N-2, \alpha , \delta , \gamma; \frac{a(x-1)}{x-a} \right) \nonumber\\
&=& \left(\frac{x-a}{1-a} \right)^{-\alpha } \sum_{r=0}^{N+1} y_{r}^{2(N+1-r)}\left( 2N+2, q_{2N+2}^m; \xi \right)  
\nonumber
\end{eqnarray}
\item As $ \delta =\beta -\gamma -2N-3 $ where $N \in \mathbb{N}_{0}$,

An algebraic equation of degree $2N+4$ for the determination of $q$ is given by
\begin{equation}  
0 = \sum_{r=0}^{N+2}\bar{c}\left( r, 2(N+2-r); 2N+3,q-(\gamma +2N+3)\alpha \right) \nonumber
\end{equation}
The eigenvalue of $q$ is written by $(\gamma +2N+3)\alpha +q_{2N+3}^m$ where $m = 0,1,2,\cdots,2N+3 $; $q_{2N+3}^0 < q_{2N+3}^1 < \cdots < q_{2N+3}^{2N+3}$. Its eigenfunction is given by
\begin{eqnarray} 
&& \left(\frac{x-a}{1-a} \right)^{-\alpha } y(\xi ) \nonumber\\
 &=& \left(\frac{x-a}{1-a} \right)^{-\alpha } Hl\left( a, q_{2N+3}^m; -2N-3, \alpha , \delta , \gamma; \frac{a(x-1)}{x-a} \right) \nonumber\\
&=& \left(\frac{x-a}{1-a} \right)^{-\alpha }  \sum_{r=0}^{N+1} y_{r}^{2(N-r)+3} \left( 2N+3,q_{2N+3}^m;\xi \right) \nonumber
\end{eqnarray}
In the above,
\begin{eqnarray}
\bar{c}(0,n;j,\tilde{q})  &=& \frac{\left( \Delta_0^{-} \left( j,\tilde{q}\right) \right)_{n}\left( \Delta_0^{+} \left( j,\tilde{q}\right) \right)_{n}}{\left( 1 \right)_{n} \left( \beta -\gamma -j \right)_{n}} \left( \frac{1+a}{a} \right)^{n}\nonumber\\
\bar{c}(1,n;j,\tilde{q}) &=& \left( -\frac{1}{a}\right) \sum_{i_0=0}^{n}\frac{\left( i_0 -j\right)\left( i_0+\alpha  \right) }{\left( i_0+2 \right) \left( i_0+1 -j +\beta -\gamma \right)} \frac{ \left( \Delta_0^{-} \left( j,\tilde{q}\right) \right)_{i_0}\left( \Delta_0^{+} \left( j,\tilde{q}\right) \right)_{i_0}}{\left( 1 \right)_{i_0} \left( \beta -\gamma -j \right)_{i_0}} \nonumber\\
&&\times  \frac{ \left( \Delta_1^{-} \left( j,\tilde{q}\right) \right)_{n}\left( \Delta_1^{+} \left( j,\tilde{q}\right) \right)_{n} \left( 3 \right)_{i_0} \left( 2-j+ \beta -\gamma  \right)_{i_0}}{\left( \Delta_1^{-} \left( j,\tilde{q}\right) \right)_{i_0}\left( \Delta_1^{+} \left( j,\tilde{q}\right) \right)_{i_0}\left( 3 \right)_{n} \left( 2-j+ \beta -\gamma \right)_{n}} \left(\frac{1+a}{a} \right)^{n }  
\nonumber\\
\bar{c}(\tau ,n;j,\tilde{q}) &=& \left( -\frac{1}{a}\right)^{\tau} \sum_{i_0=0}^{n}\frac{\left( i_0 -j\right)\left( i_0+\alpha \right) }{\left( i_0+2 \right) \left( i_0+1-j+ \beta -\gamma \right)} \frac{ \left( \Delta_0^{-} \left( j,\tilde{q}\right) \right)_{i_0}\left( \Delta_0^{+} \left( j,\tilde{q}\right) \right)_{i_0}}{\left( 1 \right)_{i_0} \left( \beta -\gamma -j \right)_{i_0}}\nonumber\\
&&\times \prod_{k=1}^{\tau -1} \left( \sum_{i_k = i_{k-1}}^{n} \frac{\left( i_k+ 2k-j\right)\left( i_k +2k+\alpha \right)}{\left( i_k+2k+2 \right) \left( i_k+2k+1-j+ \beta -\gamma \right)} \right. \nonumber\\
&&\times \left. \frac{ \left( \Delta_k^{-} \left( j,\tilde{q}\right) \right)_{i_k}\left( \Delta_k^{+} \left( j,\tilde{q}\right) \right)_{i_k} \left( 2k+1 \right)_{i_{k-1}} \left( 2k-j+ \beta -\gamma \right)_{i_{k-1}}}{\left( \Delta_k^{-} \left( j,\tilde{q}\right) \right)_{i_{k-1}}\left( \Delta_k^{+} \left( j,\tilde{q}\right) \right)_{i_{k-1}}\left( 2k+1 \right)_{i_k} \left( 2k-j+ \beta -\gamma \right)_{i_k}} \right) \nonumber\\
&&\times \frac{ \left( \Delta_{\tau }^{-} \left( j,\tilde{q}\right) \right)_{n}\left( \Delta_{\tau }^{+} \left( j,\tilde{q}\right) \right)_{n} \left( 2\tau +1 \right)_{i_{\tau -1}} \left( 2\tau -j+ \beta -\gamma \right)_{i_{\tau -1}}}{\left( \Delta_{\tau }^{-} \left( j,\tilde{q}\right) \right)_{i_{\tau -1}}\left( \Delta_{\tau }^{+} \left( j,\tilde{q}\right) \right)_{i_{\tau -1}}\left( 2\tau +1 \right)_{n} \left( 2\tau -j+ \beta -\gamma  \right)_{n}} \left(\frac{1+a}{a}\right)^{n } \nonumber
\end{eqnarray}
\begin{eqnarray}
y_0^m(j,\tilde{q};\xi) &=& \sum_{i_0=0}^{m} \frac{\left( \Delta_0^{-} \left( j,\tilde{q}\right) \right)_{i_0}\left( \Delta_0^{+} \left( j,\tilde{q}\right) \right)_{i_0}}{\left( 1 \right)_{i_0} \left( \beta -\gamma -j \right)_{i_0}} \eta ^{i_0} \nonumber\\
y_1^m(j,\tilde{q};\xi) &=& \left\{\sum_{i_0=0}^{m}\frac{\left( i_0 -j\right)\left( i_0+\alpha \right) }{\left( i_0+2 \right) \left( i_0+1-j+ \beta -\gamma \right)} \frac{ \left( \Delta_0^{-} \left( j,\tilde{q}\right) \right)_{i_0}\left( \Delta_0^{+} \left( j,\tilde{q}\right) \right)_{i_0}}{\left( 1 \right)_{i_0} \left( \beta -\gamma -j \right)_{i_0}} \right. \nonumber\\
&&\times \left. \sum_{i_1 = i_0}^{m} \frac{ \left( \Delta_1^{-} \left( j,\tilde{q}\right) \right)_{i_1}\left( \Delta_1^{+} \left( j,\tilde{q}\right) \right)_{i_1} \left( 3 \right)_{i_0} \left( 2-j+ \beta -\gamma \right)_{i_0}}{\left( \Delta_1^{-} \left( j,\tilde{q}\right) \right)_{i_0}\left( \Delta_1^{+} \left( j,\tilde{q}\right) \right)_{i_0}\left( 3 \right)_{i_1} \left( 2-j+ \beta -\gamma \right)_{i_1}} \eta ^{i_1}\right\} z 
\nonumber\\
y_{\tau }^m(j,\tilde{q};\xi) &=& \left\{ \sum_{i_0=0}^{m} \frac{\left( i_0 -j\right)\left( i_0+\alpha \right) }{\left( i_0+2 \right) \left( i_0+1-j+ \beta -\gamma \right)} \frac{ \left( \Delta_0^{-} \left( j,\tilde{q}\right) \right)_{i_0}\left( \Delta_0^{+} \left( j,\tilde{q}\right) \right)_{i_0}}{\left( 1 \right)_{i_0} \left( \beta -\gamma -j \right)_{i_0}} \right.\nonumber\\
&&\times \prod_{k=1}^{\tau -1} \left( \sum_{i_k = i_{k-1}}^{m} \frac{\left( i_k+ 2k-j\right)\left( i_k +2k+\alpha \right)}{\left( i_k+2k+2 \right) \left( i_k+2k+1-j+ \beta -\gamma \right)} \right. \nonumber\\
&&\times \left. \frac{ \left( \Delta_k^{-} \left( j,\tilde{q}\right) \right)_{i_k}\left( \Delta_k^{+} \left( j,\tilde{q}\right) \right)_{i_k} \left( 2k+1 \right)_{i_{k-1}} \left( 2k-j+ \beta -\gamma \right)_{i_{k-1}}}{\left( \Delta_k^{-} \left( j,\tilde{q}\right) \right)_{i_{k-1}}\left( \Delta_k^{+} \left( j,\tilde{q}\right) \right)_{i_{k-1}}\left( 2k+1 \right)_{i_k} \left( 2k-j+ \beta -\gamma \right)_{i_k}} \right) \nonumber\\
&&\times \left. \sum_{i_{\tau } = i_{\tau -1}}^{m}  \frac{ \left( \Delta_{\tau }^{-} \left( j,\tilde{q}\right) \right)_{i_{\tau }}\left( \Delta_{\tau }^{+} \left( j,\tilde{q}\right) \right)_{i_{\tau }} \left( 2\tau +1  \right)_{i_{\tau -1}} \left( 2\tau -j+ \beta -\gamma \right)_{i_{\tau -1}}}{\left( \Delta_{\tau }^{-} \left( j,\tilde{q}\right) \right)_{i_{\tau -1}}\left( \Delta_{\tau }^{+} \left( j,\tilde{q}\right) \right)_{i_{\tau -1}}\left( 2\tau +1 \right)_{i_\tau } \left( 2\tau -j+ \beta -\gamma \right)_{i_{\tau }}} \eta ^{i_{\tau }}\right\} z^{\tau } \nonumber 
\end{eqnarray}
where
\begin{equation}
\begin{cases} \tau \geq 2 \cr
\xi= \frac{a(x-1)}{x-a} \cr
z = -\frac{1}{a}\xi^2 \cr
\eta = \frac{(1+a)}{a} \xi \cr
\tilde{q} =  q-( \gamma +j)\alpha \cr
\Delta_k^{\pm} \left( j,\tilde{q}\right) = \frac{\varphi +4(1+a)k \pm \sqrt{\varphi ^2-4(1+a)\tilde{q}}}{2(1+a)}  \cr
\varphi = \alpha  -\gamma -j +a( \beta -1-j)
\end{cases}\nonumber
\end{equation}
\end{enumerate}
\end{appendices}

\addcontentsline{toc}{section}{Bibliography}
\bibliographystyle{model1a-num-names}
\bibliography{<your-bib-database>}
 
\chapter{Complete polynomials of Confluent Heun equation about the regular singular point at zero}
\chaptermark{Complete polynomials of the CHE around $x=0$} 

The power series solution of a Confluent Heun equation (CHE) provides a 3-term recurrence relation between successive coefficients. 
In chapter 4 of Ref.\cite{6Choun2013}, by applying three term recurrence formula (3TRF), I construct power series solutions in closed forms and integral forms of the CHE for an infinite series and a polynomial of type 1. Indeed, generating functions for the Confluent Heun polynomial (CHP) of type 1 are constructed analytically.   

In chapter 5 of Ref.\cite{6Choun2013}, by applying reversible three term recurrence formula (R3TRF), I derive power series expansions in closed forms and integral representations of the CHE for an infinite series and a polynomial of type 2 including generating functions for the CHP of type 2.

In this chapter I construct Frobenius solutions of the CHE about the regular singular point at zero for a polynomial of type 3 by applying mathematical formulas of complete polynomials using 3TRF and R3TRF.

\section{Introduction}
A formal series solution of hypergeometric equation is composed of a 2-term recursion relation between successive coefficients. Hypergeometric equation generalizes any special funtions having a 2-term recursive relation in their power series such as Legendre, Laguerre, Kummer equations, etc.
In contrast, Heun equation has a 3-term recurrence relation between consecutive coefficients in its power series. 
Currently, the Heun equation is considered as the mother of all well-known special functions such as: Spheroidal Wave, Lame, Mathieu, hypergeometric type equations, etc.  According to Karl Heun \cite{6Heun1889,6Ronv1995}, Heun equation is a second-order linear ODE which has four regular singular points. 
\begin{equation}
\frac{d^2{y}}{d{x}^2} + \left(\frac{\gamma }{x} +\frac{\delta }{x-1} + \frac{\epsilon }{x-a}\right) \frac{d{y}}{d{x}} +  \frac{\alpha \beta x-q}{x(x-1)(x-a)} y = 0 \label{eq:6001}
\end{equation}
where $\epsilon = \alpha +\beta -\gamma -\delta +1$ for assuring the regularity of the point at $x=\infty $. It has four regular singular points which are 0, 1, $a$ and $\infty $ with exponents $\{ 0, 1-\gamma \}$, $\{ 0, 1-\delta \}$, $\{ 0, 1-\epsilon \}$ and $\{ \alpha, \beta \}$. 

Heun equation has four different confluent forms such as: (1) Confluent Heun (two regular and one irregular singularities), (2) Doubly Confluent Heun (two irregular singularities), (3) Biconfluent Heun (one regular and one irregular singularities), (4) Triconfluent Heun equations (one irregular singularity).
Like deriving of confluent hypergeometric equation from the hypergeometric equation, four confluent forms of Heun equation are obtained by merging two or more regular singularities to take an irregular singularity in Heun equation.

The non-symmetrical canonical form of the CHE is a second-order linear ordinary differential equation of the form \cite{6Deca1978,6Decar1978,6Ronv1995,6Slavy2000}
\begin{equation}
\frac{d^2{y}}{d{x}^2} + \left(\beta  +\frac{\gamma }{x} + \frac{\delta }{x-1}\right) \frac{d{y}}{d{x}} +  \frac{\alpha \beta x-q}{x(x-1)} y = 0 \label{eq:6002}
\end{equation}
(\ref{eq:6002}) has three singular points: two regular singular points which are 0 and 1 with exponents $\{0, 1-\gamma\}$ and $\{0, 1-\delta \}$, and one irregular singular point which is $\infty$ with an exponent $\alpha$. Its solution is denoted as $H_{c}(\alpha,\beta,\gamma,\delta,q;x)$.\footnote{Several authors denote as coefficients $4p$ and $\sigma $ instead of $\beta $ and $q$. And they define the solution of the CHE as $H_{c}^{(a)}(p,\alpha,\gamma,\delta,\sigma;x)$.} The CHE is a more general form than any other linear second ODEs such as Whittaker-Hill, spheroidal wave, Coulomb spheroidal and Mathieu's equations.
All possible local solutions of the CHE (Regge-Wheeler and Teukolsky equations) were constructed by Fiziev. \cite{6Fizi2009,6Fizi2010} The general solutions in series of the CHE around $x=\infty $ is not convergent by only asymptotic. \cite{6Fizi2010a,6Ronv1995} 

Assume that a formal series solution takes the form 
\begin{equation}
y(x)= \sum_{n=0}^{\infty } c_n x^{n+\lambda } \hspace{1cm}\mbox{where}\; \lambda =\mbox{indicial}\;\mbox{root}\label{eq:6003}
\end{equation}
Substituting (\ref{eq:6003}) into (\ref{eq:6002}) gives for the coefficients $c_n$ the recurrence relations
\begin{equation}
c_{n+1}=A_n \;c_n +B_n \;c_{n-1} \hspace{1cm};n\geq 1 \label{eq:6004}
\end{equation}
where,
\begin{subequations}
\begin{equation}
A_n =\frac{(n+\lambda )(n+\lambda -\beta +\gamma +\delta -1)-q}{(n+1+\lambda )(n+\gamma +\lambda  )}\label{eq:6005a}
\end{equation}
\begin{equation}
B_n = \frac{\beta (n+\lambda +\alpha -1)}{(n+1+\lambda )(n+\gamma +\lambda )} \label{eq:6005b}
\end{equation}
\end{subequations}
where $c_1= A_0 \;c_0$. Two indicial roots are given such as $\lambda = 0$ and $ 1-\gamma $.

There are 4 types of power series solutions in a 3-term recurrence relation of a linear ODE such as an infinite series and 3 types of polynomials: (1) a polynomial which makes $B_n$ term terminated; $A_n$ term is not terminated, (2) a polynomial which makes $A_n$ term terminated; $B_n$ term is not terminated, (3) a polynomial which makes $A_n$ and $B_n$ terms terminated at the same time, referred as `a complete polynomial.' 

Complete polynomials can be classified in two different types such as the first species complete polynomial and the second species complete polynomial. 
If a parameter of a numerator in $B_n$ term and a (spectral) parameter of a numerator in $A_n$ term are fixed constants, we should apply the first species complete polynomial. And if two parameters of a numerator in $B_n$ term and a parameter of a numerator in $A_n$ term are fixed constants, the second species complete polynomial is able to be utilized. 
The former has multi-valued roots of a parameter of a numerator in $A_n$ term, but the latter has only one fixed value of a parameter of a numerator in $A_n$ term for an eigenvalue.  
 
\begin{table}[h]
\begin{center}
\thispagestyle{plain}
\hspace*{-0.1\linewidth}\resizebox{1.2\linewidth}{!}
{
 \Tree[.{\Huge Confluent Heun differential equation about the regular singular point at zero} [.{\Huge 3TRF} [.{\Huge Infinite series} ]
              [.{\Huge Polynomials} [[.{\Huge Polynomial of type 1} ]
               [.{\Huge Polynomial of type 3} [.{\Huge $ \begin{array}{lcll}  1^{\mbox{st}}\;  \mbox{species}\\ \mbox{complete} \\ \mbox{polynomial} \end{array}$} ]  ]]]]                         
  [.{\Huge R3TRF} [.{\Huge Infinite series} ]
     [.{\Huge Polynomials} [[.{\Huge Polynomial of type 2} ]
       [.{\Huge  Polynomial of type 3} [.{\Huge $ \begin{array}{lcll}  1^{\mbox{st}}\;  \mbox{species} \\ \mbox{complete} \\ \mbox{polynomial} \end{array}$} ]  ]]]]]
}
\end{center}
\caption{Power series of the CHE about the regular singular point at zero}
\end{table}  

Table 7.1 informs us about all possible general solutions in series of the CHE about the regular singular point at zero. 
I construct power series solutions of the CHE around $x=0$ and their combined definite \& contour integral representations for an infinite series and a polynomial of type 1 by applying 3TRF in chapter 4 of Ref.\cite{6Choun2013}. 
The sequence $c_n$ combines into combinations of $A_n$ and $B_n$ terms in (\ref{eq:6004}). This is done by letting $A_n$ in the sequence $c_n$ is the leading term in the analytic function $y(x)$: I observe the term of sequence $c_n$ which includes zero term of $A_n's$, one term of $A_n's$, two terms of $A_n's$, three terms of $A_n's$, etc. 
 For a polynomial of type 1, I treat $\beta $, $\gamma $, $\delta $, $q$ as free variables and $\alpha $ as a fixed value. Also, generating functions for the CHP of type 1 are constructed analytically by applying generating functions for confluent hypergeometric (Kummer's) polynomials into each sub-integral form of the general integral representation for the CHP  of type 1.
  
I obtain Frobenius solutions of the CHE around $x=0$ for an infinite series and a polynomial of type 2 by applying R3TRF including their combined definite \& contour integrals in chapter 5 of Ref.\cite{6Choun2013}. This is done by letting $B_n$ in the sequence $c_n$ is the leading term in the analytic function $y(x)$: I observe the term of sequence $c_n$ which includes zero term of $B_n's$, one term of $B_n's$, two terms of $B_n's$, three terms of $B_n's$, etc. For a polynomial of type 2, I treat $\alpha $, $\beta $, $\gamma $, $\delta $ as free variables and $q$ as a fixed value. Generating functions for the CHP of type 2 are constructed analytically by applying the generating function for Jacobi polynomial using hypergeometric functions into each sub-integral form  of the general integral form for the CHP of type 2.
Infinite series of the CHE using 3TRF around $x=0$ are equivalent to infinite series of it using R3TRF. The former is that $A_n$ is the leading term in each sub-power series of the CHE. The latter is that $B_n$ is the leading term in each sub-power series of it. 
  
In general, the CHP has been hitherto defined as a type 3 polynomial which is equivalent to the first species complete polynomial. \cite{6Fizi2010a,6Fizi2009,6Fizi2010,6Fizi2010b} The CHP comes from the CHE around $x=0$ or $x=\infty $ that has a fixed integer value of $\alpha $, just as it has a fixed value of $q$. 
Various authors build an algebraic equation of the $(j+1)$th order for the determination of a parameter $q$ by using matrix formalism (the determinant of $(j+1)\times (j+1)$ matrices).  
And they leave a polynomial of the $j$th order for the CHP as solutions of recurrences because of complicated computations.\cite{6ElJa2009,6Fizi2009b,6Figu2005,6Hall2010,6Sesm2010}
 
In chapter 1, I construct the mathematical formulas of complete polynomials for the first and second species, by allowing $A_n$ as the leading term in each finite sub-power series of the general power series $y(x)$, designed as ``complete polynomials using 3-term recurrence formula (3TRF).''
In chapter 2, I construct the classical mathematical formulas in series of complete polynomials for the first and second species, by allowing $B_n$ as the leading term in each finite sub-power series of the general power series $y(x)$, denominated as ``complete polynomials using reversible 3-term recurrence formula (R3TRF).'' 

In this chapter, by applying complete polynomials using 3TRF and R3TRF, I construct the power series expansion in closed forms of the CHE around $x=0$ for a polynomial which makes $A_n$ and $B_n$ terms terminated.
Besides, I build an algebraic equation of the CHE for the determination of a parameter $q$ in the form of partial sums of the sequences $\{A_n\}$ and $\{B_n\}$ using 3TRF and R3TRF for the first species polynomials.\footnote{I do not include general series solutions of the CHE around $x=1$ in chapters 4 \& 5 of Ref.\cite{6Choun2013} and this chapter. Because since we construct power series solutions in compact forms of the CHE around $x=0$, we can obtain Frobenius solutions in closed forms of the CHE readily by a linear transformation interchanging the regular singular points at zero and one: $x\rightarrow 1-x$.}
    
\section{Power series about the regular singular point at zero}
\sectionmark{Power series about the regular singular point at $x=0 $}
For the first species complete polynomial of the CHE around $x=0$ in table 7.1, I treat $\beta $, $\gamma $, $\delta $ as free variables and $\alpha $, $q$ as fixed values. There is no such solution for the second species complete polynomial of the CHE. Because a parameter $\alpha $ of a numerator in $B_n$ term in (\ref{eq:6005b}) is only a fixed constant in order to make $B_n$ term terminated at a specific index summation $n$ for two independent solutions where $\lambda = 0$ and $ 1-\gamma $.  
\subsection{The first species complete polynomial of Confluent Heun equation using 3TRF}
For the first species complete polynomials using 3TRF and R3TRF, we need a condition which is given by
\begin{equation}
B_{j+1}= c_{j+1}=0\hspace{1cm}\mathrm{where}\;j\in \mathbb{N}_{0}  
 \label{eq:6006}
\end{equation}
(\ref{eq:6006}) gives successively $c_{j+2}=c_{j+3}=c_{j+4}=\cdots=0$. And $c_{j+1}=0$ is defined by a polynomial equation of degree $j+1$ for the determination of an accessory parameter in $A_n$ term. 
\begin{theorem}
In chapter 1, the general expression of a function $y(x)$ for the first species complete polynomial using 3-term recurrence formula and its algebraic equation for the determination of an accessory parameter in $A_n$ term are given by
\begin{enumerate} 
\item As $B_1=0$,
\begin{equation}
0 =\bar{c}(1,0) \label{eq:6007a}
\end{equation}
\begin{equation}
y(x) = y_{0}^{0}(x) \label{eq:6007b}
\end{equation}
\item As $B_{2N+2}=0$ where $N \in \mathbb{N}_{0}$,
\begin{equation}
0  = \sum_{r=0}^{N+1}\bar{c}\left( 2r, N+1-r\right) \label{eq:6008a}
\end{equation}
\begin{equation}
y(x)= \sum_{r=0}^{N} y_{2r}^{N-r}(x)+ \sum_{r=0}^{N} y_{2r+1}^{N-r}(x)  \label{eq:6008b}
\end{equation}
\item As $B_{2N+3}=0$ where $N \in \mathbb{N}_{0}$,
\begin{equation}
0  = \sum_{r=0}^{N+1}\bar{c}\left( 2r+1, N+1-r\right) \label{eq:6009a}
\end{equation}
\begin{equation}
y(x)= \sum_{r=0}^{N+1} y_{2r}^{N+1-r}(x)+ \sum_{r=0}^{N} y_{2r+1}^{N-r}(x)  \label{eq:6009b}
\end{equation}
In the above,
\begin{eqnarray}
\bar{c}(0,n)  &=& \prod _{i_0=0}^{n-1}B_{2i_0+1} \label{eq:60010a}\\
\bar{c}(1,n) &=&  \sum_{i_0=0}^{n} \left\{ A_{2i_0} \prod _{i_1=0}^{i_0-1}B_{2i_1+1} \prod _{i_2=i_0}^{n-1}B_{2i_2+2} \right\} 
\label{eq:60010b}\\
\bar{c}(\tau ,n) &=& \sum_{i_0=0}^{n} \left\{A_{2i_0}\prod _{i_1=0}^{i_0-1} B_{2i_1+1} 
\prod _{k=1}^{\tau -1} \left( \sum_{i_{2k}= i_{2(k-1)}}^{n} A_{2i_{2k}+k}\prod _{i_{2k+1}=i_{2(k-1)}}^{i_{2k}-1}B_{2i_{2k+1}+(k+1)}\right) \right. \nonumber\\ 
&&\times \left. \prod _{i_{2\tau}=i_{2(\tau -1)}}^{n-1} B_{2i_{2\tau }+(\tau +1)} \right\} 
\hspace{1cm}\label{eq:60010c} 
\end{eqnarray}
and
\begin{eqnarray}
y_0^m(x) &=& c_0 x^{\lambda } \sum_{i_0=0}^{m} \left\{ \prod _{i_1=0}^{i_0-1}B_{2i_1+1} \right\} x^{2i_0 } \label{eq:60011a}\\
y_1^m(x) &=& c_0 x^{\lambda } \sum_{i_0=0}^{m}\left\{ A_{2i_0} \prod _{i_1=0}^{i_0-1}B_{2i_1+1}  \sum_{i_2=i_0}^{m} \left\{ \prod _{i_3=i_0}^{i_2-1}B_{2i_3+2} \right\}\right\} x^{2i_2+1 } \label{eq:60011b}\\
y_{\tau }^m(x) &=& c_0 x^{\lambda } \sum_{i_0=0}^{m} \left\{A_{2i_0}\prod _{i_1=0}^{i_0-1} B_{2i_1+1} 
\prod _{k=1}^{\tau -1} \left( \sum_{i_{2k}= i_{2(k-1)}}^{m} A_{2i_{2k}+k}\prod _{i_{2k+1}=i_{2(k-1)}}^{i_{2k}-1}B_{2i_{2k+1}+(k+1)}\right) \right. \nonumber\\
&& \times  \left. \sum_{i_{2\tau} = i_{2(\tau -1)}}^{m} \left( \prod _{i_{2\tau +1}=i_{2(\tau -1)}}^{i_{2\tau}-1} B_{2i_{2\tau +1}+(\tau +1)} \right) \right\} x^{2i_{2\tau}+\tau }\hspace{1cm}\mathrm{where}\;\tau \geq 2
\label{eq:60011c} 
\end{eqnarray}
\end{enumerate}
\end{theorem}
Put $n= j+1$ in (\ref{eq:6005b}) and use the condition $B_{j+1}=0$ for $\alpha $.  
\begin{equation}
\alpha = -j-\lambda 
\label{eq:60012}
\end{equation}
Take (\ref{eq:60012}) into (\ref{eq:6005b}).
\begin{equation}
B_n = \frac{\beta (n-j -1)}{(n+1+\lambda )(n+\gamma +\lambda )} \label{eq:60013}
\end{equation}
Now the condition $c_{j+1}=0$ is clearly an algebraic equation in $q$ of degree $j+1$ and thus has $j+1$ zeros denoted them by $q_j^m$ eigenvalues where $m = 0,1,2, \cdots, j$. They can be arranged in the following order: $q_j^0 < q_j^1 < q_j^2 < \cdots < q_j^j$.

Substitute (\ref{eq:6005a}) and (\ref{eq:60013}) with $c_1= A_0 \;c_0$ into (\ref{eq:60010a})--(\ref{eq:60011c}).

As $B_{1}= c_{1}=0$, take the new (\ref{eq:60010b}) into (\ref{eq:6007a}) putting $j=0$. Substitute the new (\ref{eq:60011a}) into (\ref{eq:6007b}) putting $j=0$. 

As $B_{2N+2}= c_{2N+2}=0$, take the new (\ref{eq:60010a})--(\ref{eq:60010c}) into (\ref{eq:6008a}) putting $j=2N+1$. Substitute the new 
(\ref{eq:60011a})--(\ref{eq:60011c}) into (\ref{eq:6008b}) putting $j=2N+1$ and $q=q_{2N+1}^m$.

As $B_{2N+3}= c_{2N+3}=0$, take the new (\ref{eq:60010a})--(\ref{eq:60010c}) into (\ref{eq:6009a}) putting $j=2N+2$. Substitute the new 
(\ref{eq:60011a})--(\ref{eq:60011c}) into (\ref{eq:6009b}) putting $j=2N+2$ and $q=q_{2N+2}^m$.

After the replacement process, the general expression of power series of the CHE about $x=0$ for the first species complete polynomial using 3-term recurrence formula and its algebraic equation for the determination of an accessory parameter $q$ are given by
\begin{enumerate} 
\item As $\alpha =-\lambda $,

An algebraic equation of degree 1 for the determination of $q$ is given by
\begin{equation}
0= \bar{c}(1,0;0,q)= -q + \lambda (-\beta +\gamma + \delta -1+\lambda ) \label{eq:60014a}
\end{equation}
The eigenvalue of $q$ is written by $q_0^0$. Its eigenfunction is given by
\begin{equation}
y(x) = y_0^0\left( 0,q_0^0;x\right)= c_0 x^{\lambda } \label{eq:60014b}
\end{equation}
\item As $\alpha =-2N-1-\lambda $ where $N \in \mathbb{N}_{0}$,

An algebraic equation of degree $2N+2$ for the determination of $q$ is given by
\begin{equation}
0 = \sum_{r=0}^{N+1}\bar{c}\left( 2r, N+1-r; 2N+1,q\right)  \label{eq:60015a}
\end{equation}
The eigenvalue of $q$ is written by $q_{2N+1}^m$ where $m = 0,1,2,\cdots,2N+1 $; $q_{2N+1}^0 < q_{2N+1}^1 < \cdots < q_{2N+1}^{2N+1}$. Its eigenfunction is given by 
\begin{equation} 
y(x) = \sum_{r=0}^{N} y_{2r}^{N-r}\left( 2N+1,q_{2N+1}^m;x\right)+ \sum_{r=0}^{N} y_{2r+1}^{N-r}\left( 2N+1,q_{2N+1}^m;x\right)
\label{eq:60015b} 
\end{equation}
\item As $\alpha =-2N-2-\lambda $ where $N \in \mathbb{N}_{0}$,

An algebraic equation of degree $2N+3$ for the determination of $q$ is given by
\begin{equation}  
0 = \sum_{r=0}^{N+1}\bar{c}\left( 2r+1, N+1-r; 2N+2,q\right) \label{eq:60016a}
\end{equation}
The eigenvalue of $q$ is written by $q_{2N+2}^m$ where $m = 0,1,2,\cdots,2N+2 $; $q_{2N+2}^0 < q_{2N+2}^1 < \cdots < q_{2N+2}^{2N+2}$. Its eigenfunction is given by
\begin{equation} 
y(x) =  \sum_{r=0}^{N+1} y_{2r}^{N+1-r}\left( 2N+2,q_{2N+2}^m;x\right) + \sum_{r=0}^{N} y_{2r+1}^{N-r}\left( 2N+2,q_{2N+2}^m;x\right) \label{eq:60016b}
\end{equation}
In the above,
\begin{eqnarray}
\bar{c}(0,n;j,q)  &=& \frac{\left( -\frac{j}{2}\right)_{n}}{\left( 1+\frac{\lambda }{2}\right)_{n} \left(  \frac{1}{2}+ \frac{\gamma }{2}+\frac{\lambda }{2}\right)_{n}} \left( \frac{\beta }{2} \right)^{n}\label{eq:60017a}\\
\bar{c}(1,n;j,q) &=& \sum_{i_0=0}^{n}\frac{\left( i_0+\frac{\lambda }{2}\right)\left( i_0 +\Gamma _0 \right) -\frac{q}{4}}{\left( i_0+\frac{1}{2}+\frac{\lambda }{2}\right) \left( i_0+\frac{\gamma }{2}+\frac{\lambda }{2}\right)} \frac{\left( -\frac{j}{2}\right)_{i_0} }{\left( 1+\frac{\lambda }{2}\right)_{i_0} \left( \frac{1}{2}+ \frac{\gamma }{2}+ \frac{\lambda }{2}\right)_{i_0}} \nonumber\\
&&\times \frac{\left( \frac{1}{2}-\frac{j}{2} \right)_{n} \left( \frac{3}{2}+\frac{\lambda }{2}\right)_{i_0} \left( 1+\frac{\gamma }{2}+ \frac{\lambda }{2}\right)_{i_0}}{\left( \frac{1}{2}-\frac{j}{2}\right)_{i_0} \left( \frac{3}{2}+\frac{\lambda }{2}\right)_{n} \left( 1+\frac{\gamma }{2}+ \frac{\lambda }{2}\right)_{n}} \left( \frac{\beta }{2} \right)^{n }  
\label{eq:60017b}\\
\bar{c}(\tau ,n;j,q) &=& \sum_{i_0=0}^{n}\frac{\left( i_0+\frac{\lambda }{2}\right)\left( i_0 +\Gamma _0 \right) -\frac{q}{4}}{\left( i_0+\frac{1}{2}+\frac{\lambda }{2}\right) \left( i_0+\frac{\gamma }{2}+\frac{\lambda }{2}\right)} \frac{\left( -\frac{j}{2}\right)_{i_0} }{\left( 1+\frac{\lambda }{2}\right)_{i_0} \left( \frac{1}{2}+\frac{\gamma }{2}+  \frac{\lambda }{2}\right)_{i_0}}  \nonumber\\
&&\times   \prod_{k=1}^{\tau -1} \left( \sum_{i_k = i_{k-1}}^{n} \frac{\left( i_k+ \frac{k}{2}+\frac{\lambda }{2}\right)\left( i_k +\Gamma _k \right) -\frac{q}{4}}{\left( i_k+\frac{k}{2}+\frac{1}{2}+\frac{\lambda }{2}\right) \left( i_k+\frac{k}{2}+\frac{\gamma }{2}+\frac{\lambda }{2}\right)} \right. \nonumber\\
&&\times \left. \frac{\left( \frac{k}{2}-\frac{j}{2}\right)_{i_k} \left( \frac{k}{2}+1+ \frac{\lambda }{2}\right)_{i_{k-1}} \left( \frac{k}{2}+ \frac{1}{2}+\frac{\gamma }{2}+ \frac{\lambda }{2}\right)_{i_{k-1}}}{\left( \frac{k}{2}-\frac{j}{2}\right)_{i_{k-1}} \left( \frac{k}{2}+1+ \frac{\lambda }{2}\right)_{i_k} \left( \frac{k}{2}+\frac{1}{2}+ \frac{\gamma }{2}+  \frac{\lambda }{2}\right)_{i_k}} \right) \nonumber\\ 
&&\times \frac{\left( \frac{\tau }{2} -\frac{j}{2}\right)_{n} \left( \frac{\tau }{2}+1+\frac{\lambda }{2}\right)_{i_{\tau -1}} \left( \frac{\tau }{2}+\frac{1}{2}+\frac{\gamma }{2}+  \frac{\lambda }{2}\right)_{i_{\tau -1}}}{\left( \frac{\tau }{2}-\frac{j}{2}\right)_{i_{\tau -1}} \left( \frac{\tau }{2}+1+\frac{\lambda }{2}\right)_{n} \left( \frac{\tau }{2}+\frac{1}{2}+ \frac{\gamma }{2}+ \frac{\lambda }{2}\right)_{n}} \left( \frac{\beta }{2} \right)^{n }\hspace{1.5cm} \label{eq:60017c} 
\end{eqnarray}
\begin{eqnarray}
y_0^m(j,q;x) &=& c_0 x^{\lambda }  \sum_{i_0=0}^{m} \frac{\left( -\frac{j}{2}\right)_{i_0}}{\left( 1+\frac{\lambda }{2}\right)_{i_0} \left( \frac{1}{2}+ \frac{\gamma }{2}+ \frac{\lambda }{2}\right)_{i_0}} \eta ^{i_0} \label{eq:60018a}\\
y_1^m(j,q;x) &=& c_0 x^{\lambda } \left\{\sum_{i_0=0}^{m} \frac{\left( i_0+\frac{\lambda }{2}\right)\left(i_0 +\Gamma _0 \right) -\frac{q}{4}}{\left( i_0+\frac{1}{2}+\frac{\lambda }{2}\right) \left( i_0+\frac{\gamma }{2}+\frac{\lambda }{2}\right)} \frac{\left( -\frac{j}{2}\right)_{i_0} }{\left( 1+\frac{\lambda }{2}\right)_{i_0} \left( \frac{1}{2}+\frac{\gamma }{2}+ \frac{\lambda }{2}\right)_{i_0}} \right. \nonumber\\
&&\times \left. \sum_{i_1 = i_0}^{m} \frac{\left( \frac{1}{2}-\frac{j}{2} \right)_{i_1} \left( \frac{3}{2}+\frac{\lambda }{2}\right)_{i_0} \left( 1+\frac{\gamma }{2}+ \frac{\lambda }{2}\right)_{i_0}}{\left( \frac{1}{2}-\frac{j}{2} \right)_{i_0} \left( \frac{3}{2}+\frac{\lambda }{2}\right)_{i_1} \left( 1+\frac{\gamma }{2} + \frac{\lambda }{2}\right)_{i_1}} \eta ^{i_1}\right\}x 
\label{eq:60018b}
\end{eqnarray}
\begin{eqnarray}
y_{\tau }^m(j,q;x) &=& c_0 x^{\lambda } \left\{ \sum_{i_0=0}^{m} \frac{\left( i_0+\frac{\lambda }{2}\right)\left( i_0 +\Gamma _0 \right) -\frac{q}{4}}{\left( i_0+\frac{1}{2}+\frac{\lambda }{2}\right) \left( i_0+\frac{\gamma }{2}+\frac{\lambda }{2}\right)} \frac{\left( -\frac{j}{2}\right)_{i_0} }{\left( 1+\frac{\lambda }{2}\right)_{i_0} \left( \frac{1}{2}+\frac{\gamma }{2}+ \frac{\lambda }{2}\right)_{i_0}} \right.\nonumber\\
&&\times \prod_{k=1}^{\tau -1} \left( \sum_{i_k = i_{k-1}}^{m} \frac{\left( i_k+ \frac{k}{2}+\frac{\lambda }{2}\right)\left( i_k +\Gamma _k \right) -\frac{q}{4}}{\left( i_k+\frac{k}{2}+\frac{1}{2}+\frac{\lambda }{2}\right) \left( i_k+\frac{k}{2}+\frac{\gamma }{2}+\frac{\lambda }{2}\right)} \right.  \nonumber\\
&&\times  \left. \frac{\left( \frac{k}{2}-\frac{j}{2}\right)_{i_k} \left( \frac{k}{2}+1+ \frac{\lambda }{2}\right)_{i_{k-1}} \left(  \frac{k}{2}+\frac{1}{2}+\frac{\gamma }{2}+ \frac{\lambda }{2}\right)_{i_{k-1}}}{\left( \frac{k}{2}-\frac{j}{2}\right)_{i_{k-1}} \left( \frac{k}{2}+1+ \frac{\lambda }{2}\right)_{i_k} \left( \frac{k}{2}+\frac{1}{2}+\frac{\gamma }{2}+ \frac{\lambda }{2}\right)_{i_k}} \right) \nonumber\\
&&\times \left. \sum_{i_{\tau } = i_{\tau -1}}^{m} \frac{\left( \frac{\tau }{2}-\frac{j}{2}\right)_{i_{\tau }} \left( \frac{\tau }{2}+1+\frac{\lambda }{2}\right)_{i_{\tau -1}} \left( \frac{\tau }{2}+\frac{1}{2}+\frac{\gamma }{2}+ \frac{\lambda }{2}\right)_{i_{\tau -1}}}{\left( \frac{\tau }{2}-\frac{j}{2}\right)_{i_{\tau -1}} \left( \frac{\tau }{2}+1+\frac{\lambda }{2}\right)_{i_{\tau }} \left( \frac{\tau }{2}+\frac{1}{2}+\frac{\gamma }{2}+ \frac{\lambda }{2}\right)_{i_{\tau }}} \eta ^{i_{\tau }}\right\} x^{\tau } \hspace{1.5cm}\label{eq:60018c} 
\end{eqnarray}
where
\begin{equation}
\begin{cases} \tau \geq 2 \cr
\eta  = \frac{1}{2}\beta x^2 \cr
\Gamma _0 = \frac{1}{2}\left( -\beta + \gamma +\delta -1+\lambda \right) \cr
\Gamma _k = \frac{1}{2}\left( -\beta + \gamma +\delta -1+k+\lambda \right)
\end{cases}\nonumber
\end{equation}
\end{enumerate}
Put $c_0$= 1 as $\lambda =0$ for the first kind of independent solutions of the CHE and $\lambda = 1-\gamma $ for the second one in (\ref{eq:60014a})--(\ref{eq:60018c}). 
\begin{remark}
The power series expansion of Confluent Heun equation of the first kind for the first species complete polynomial using 3TRF about $x=0$ is given by
\begin{enumerate} 
\item As $\alpha =0$ and $q=q_0^0=0$,

The eigenfunction is given by
\begin{equation}
y(x) =\; _pH_c^{(a)}F_{0,0} \left( \alpha =0, \beta, \gamma, \delta, q=q_0^0=0; \eta = \frac{1}{2}\beta x^2 \right) =1 \label{eq:60019}
\end{equation}
\item As $\alpha =-2N-1$ where $N \in \mathbb{N}_{0}$,

An algebraic equation of degree $2N+2$ for the determination of $q$ is given by
\begin{equation}
0 = \sum_{r=0}^{N+1}\bar{c}\left( 2r, N+1-r; 2N+1,q\right)\label{eq:60020a}
\end{equation}
The eigenvalue of $q$ is written by $q_{2N+1}^m$ where $m = 0,1,2,\cdots,2N+1 $; $q_{2N+1}^0 < q_{2N+1}^1 < \cdots < q_{2N+1}^{2N+1}$. Its eigenfunction is given by
\begin{eqnarray} 
y(x) &=& _pH_c^{(a)}F_{2N+1,m} \left( \alpha =-2N-1, \beta, \gamma, \delta, q=q_{2N+1}^m; \eta = \frac{1}{2}\beta x^2 \right)\nonumber\\
&=& \sum_{r=0}^{N} y_{2r}^{N-r}\left( 2N+1,q_{2N+1}^m;x\right)+ \sum_{r=0}^{N} y_{2r+1}^{N-r}\left( 2N+1,q_{2N+1}^m;x\right)  
\label{eq:60020b}
\end{eqnarray}
\item As $\alpha =-2N-2$ where $N \in \mathbb{N}_{0}$,

An algebraic equation of degree $2N+3$ for the determination of $q$ is given by
\begin{eqnarray}
0  = \sum_{r=0}^{N+1}\bar{c}\left( 2r+1, N+1-r; 2N+2,q\right)\label{eq:60021a}
\end{eqnarray}
The eigenvalue of $q$ is written by $q_{2N+2}^m$ where $m = 0,1,2,\cdots,2N+2 $; $q_{2N+2}^0 < q_{2N+2}^1 < \cdots < q_{2N+2}^{2N+2}$. Its eigenfunction is given by
\begin{eqnarray} 
y(x) &=& _pH_c^{(a)}F_{2N+2,m} \left( \alpha =-2N-2, \beta, \gamma, \delta, q=q_{2N+2}^m; \eta = \frac{1}{2}\beta x^2 \right)\nonumber\\
&=& \sum_{r=0}^{N+1} y_{2r}^{N+1-r}\left( 2N+2,q_{2N+2}^m;x\right) + \sum_{r=0}^{N} y_{2r+1}^{N-r}\left( 2N+2,q_{2N+2}^m;x\right) 
\hspace{1.5cm} \label{eq:60021b}
\end{eqnarray}
In the above,
\begin{eqnarray}
\bar{c}(0,n;j,q)  &=& \frac{\left( -\frac{j}{2}\right)_{n}}{\left( 1 \right)_{n} \left(  \frac{1}{2}+ \frac{\gamma }{2} \right)_{n}} \left( \frac{\beta }{2} \right)^{n}\label{eq:60022a}\\
\bar{c}(1,n;j,q) &=& \sum_{i_0=0}^{n}\frac{ i_0 \left( i_0 +\Gamma _0 \right) -\frac{q}{4}}{\left( i_0+\frac{1}{2} \right) \left( i_0+\frac{\gamma }{2} \right)} \frac{\left( -\frac{j}{2}\right)_{i_0} }{\left( 1 \right)_{i_0} \left( \frac{1}{2}+ \frac{\gamma }{2} \right)_{i_0}}   \frac{\left( \frac{1}{2}-\frac{j}{2} \right)_{n} \left( \frac{3}{2} \right)_{i_0} \left( 1+\frac{\gamma }{2} \right)_{i_0}}{\left( \frac{1}{2}-\frac{j}{2}\right)_{i_0} \left( \frac{3}{2} \right)_{n} \left( 1+\frac{\gamma }{2} \right)_{n}} \left( \frac{\beta }{2} \right)^{n }  
\label{eq:60022b}\\
\bar{c}(\tau ,n;j,q) &=& \sum_{i_0=0}^{n}\frac{ i_0 \left( i_0 +\Gamma _0 \right) -\frac{q}{4}}{\left( i_0+\frac{1}{2} \right) \left( i_0+\frac{\gamma }{2} \right)} \frac{\left( -\frac{j}{2}\right)_{i_0} }{\left( 1 \right)_{i_0} \left( \frac{1}{2}+\frac{\gamma }{2} \right)_{i_0}}  \nonumber\\
&&\times  \prod_{k=1}^{\tau -1} \left( \sum_{i_k = i_{k-1}}^{n} \frac{\left( i_k+ \frac{k}{2} \right)\left( i_k +\Gamma _k \right) -\frac{q}{4}}{\left( i_k+\frac{k}{2}+\frac{1}{2} \right) \left( i_k+\frac{k}{2}+\frac{\gamma }{2} \right)} \right.   \left. \frac{\left( \frac{k}{2}-\frac{j}{2}\right)_{i_k} \left( \frac{k}{2}+1 \right)_{i_{k-1}} \left( \frac{k}{2}+ \frac{1}{2}+\frac{\gamma }{2} \right)_{i_{k-1}}}{\left( \frac{k}{2}-\frac{j}{2}\right)_{i_{k-1}} \left( \frac{k}{2}+1 \right)_{i_k} \left( \frac{k}{2}+\frac{1}{2}+ \frac{\gamma }{2} \right)_{i_k}} \right) \nonumber\\ 
&&\times \frac{\left( \frac{\tau }{2} -\frac{j}{2}\right)_{n} \left( \frac{\tau }{2}+1 \right)_{i_{\tau -1}} \left( \frac{\tau }{2}+\frac{1}{2}+\frac{\gamma }{2} \right)_{i_{\tau -1}}}{\left( \frac{\tau }{2}-\frac{j}{2}\right)_{i_{\tau -1}} \left( \frac{\tau }{2}+1 \right)_{n} \left( \frac{\tau }{2}+\frac{1}{2}+ \frac{\gamma }{2} \right)_{n}} \left( \frac{\beta }{2} \right)^{n } \label{eq:60022c} 
\end{eqnarray}
\begin{eqnarray}
y_0^m(j,q;x) &=& \sum_{i_0=0}^{m} \frac{\left( -\frac{j}{2}\right)_{i_0}}{\left( 1 \right)_{i_0} \left( \frac{1}{2}+ \frac{\gamma }{2} \right)_{i_0}} \eta ^{i_0} \label{eq:60023a}\\
y_1^m(j,q;x) &=& \left\{\sum_{i_0=0}^{m} \frac{ i_0 \left(i_0 +\Gamma _0 \right) -\frac{q}{4}}{\left( i_0+\frac{1}{2} \right) \left( i_0+\frac{\gamma }{2} \right)} \frac{\left( -\frac{j}{2}\right)_{i_0} }{\left( 1 \right)_{i_0} \left( \frac{1}{2}+\frac{\gamma }{2} \right)_{i_0}} \right. \nonumber\\
&&\times \left. \sum_{i_1 = i_0}^{m} \frac{\left( \frac{1}{2}-\frac{j}{2} \right)_{i_1} \left( \frac{3}{2} \right)_{i_0} \left( 1+\frac{\gamma }{2} \right)_{i_0}}{\left( \frac{1}{2}-\frac{j}{2} \right)_{i_0} \left( \frac{3}{2} \right)_{i_1} \left( 1+\frac{\gamma }{2} \right)_{i_1}} \eta ^{i_1}\right\}x 
\hspace{1.5cm}\label{eq:60023b}\\
y_{\tau }^m(j,q;x) &=& \left\{ \sum_{i_0=0}^{m} \frac{ i_0 \left( i_0 +\Gamma _0 \right) -\frac{q}{4}}{\left( i_0+\frac{1}{2} \right) \left( i_0+\frac{\gamma }{2} \right)} \frac{\left( -\frac{j}{2}\right)_{i_0} }{\left( 1 \right)_{i_0} \left( \frac{1}{2}+\frac{\gamma }{2} \right)_{i_0}} \right.\nonumber\\
&&\times \prod_{k=1}^{\tau -1} \left( \sum_{i_k = i_{k-1}}^{m} \frac{\left( i_k+ \frac{k}{2} \right)\left( i_k +\Gamma _k \right) -\frac{q}{4}}{\left( i_k+\frac{k}{2}+\frac{1}{2} \right) \left( i_k+\frac{k}{2}+\frac{\gamma }{2} \right)} \right. \nonumber\\
&&\times  \left. \frac{\left( \frac{k}{2}-\frac{j}{2}\right)_{i_k} \left( \frac{k}{2}+1 \right)_{i_{k-1}} \left(  \frac{k}{2}+\frac{1}{2}+\frac{\gamma }{2} \right)_{i_{k-1}}}{\left( \frac{k}{2}-\frac{j}{2}\right)_{i_{k-1}} \left( \frac{k}{2}+1 \right)_{i_k} \left( \frac{k}{2}+\frac{1}{2}+\frac{\gamma }{2} \right)_{i_k}} \right) \nonumber\\
&&\times  \left. \sum_{i_{\tau } = i_{\tau -1}}^{m} \frac{\left( \frac{\tau }{2}-\frac{j}{2}\right)_{i_{\tau }} \left( \frac{\tau }{2}+1 \right)_{i_{\tau -1}} \left( \frac{\tau }{2}+\frac{1}{2}+\frac{\gamma }{2} \right)_{i_{\tau -1}}}{\left( \frac{\tau }{2}-\frac{j}{2}\right)_{i_{\tau -1}} \left( \frac{\tau }{2}+1 \right)_{i_{\tau }} \left( \frac{\tau }{2}+\frac{1}{2}+\frac{\gamma }{2} \right)_{i_{\tau }}} \eta ^{i_{\tau }}\right\} x^{\tau } \hspace{1.5cm}\label{eq:60023c} 
\end{eqnarray}
where
\begin{equation}
\begin{cases} \tau \geq 2 \cr 
\Gamma _0 = \frac{1}{2}\left( -\beta + \gamma +\delta -1 \right) \cr
\Gamma _k = \frac{1}{2}\left( -\beta + \gamma +\delta -1+k \right)
\end{cases}\nonumber
\end{equation}
\end{enumerate}
\end{remark}
\begin{remark}
The power series expansion of Confluent Heun equation of the second kind for the first species complete polynomial using 3TRF about $x=0$ is given by
\begin{enumerate} 
\item As $\alpha = \gamma -1$ and $q=q_0^0= (\gamma -1)(\beta -\delta )$,

The eigenfunction is given by
\begin{eqnarray}
 y(x) &=&\; _pH_c^{(a)}S_{0,0} \left( \alpha =\gamma -1, \beta, \gamma, \delta, q=q_0^0=(\gamma -1)(\beta -\delta ); \eta = \frac{1}{2}\beta x^2 \right) \nonumber\\
&=& x^{1-\gamma } \label{eq:60024}
\end{eqnarray}
\item As $\alpha =\gamma -2N-2 $ where $N \in \mathbb{N}_{0}$, 

An algebraic equation of degree $2N+2$ for the determination of $q$ is given by
\begin{equation}
0 = \sum_{r=0}^{N+1}\bar{c}\left( 2r, N+1-r; 2N+1,q\right) \label{eq:60025a}
\end{equation}
The eigenvalue of $q$ is written by $q_{2N+1}^m$ where $m = 0,1,2,\cdots,2N+1 $; $q_{2N+1}^0 < q_{2N+1}^1 < \cdots < q_{2N+1}^{2N+1}$. Its eigenfunction is given by
\begin{eqnarray} 
y(x) &=& _pH_c^{(a)}S_{2N+1,m} \left( \alpha = \gamma -2N-2, \beta, \gamma, \delta, q=q_{2N+1}^m; \eta = \frac{1}{2}\beta x^2 \right) \nonumber\\
&=&  \sum_{r=0}^{N} y_{2r}^{N-r}\left( 2N+1,q_{2N+1}^m;x\right)+ \sum_{r=0}^{N} y_{2r+1}^{N-r}\left( 2N+1,q_{2N+1}^m;x\right)
\label{eq:60025b}
\end{eqnarray}
\item As $\alpha =\gamma -2N-3 $ where  $N \in \mathbb{N}_{0}$,

An algebraic equation of degree $2N+3$ for the determination of $q$ is given by
\begin{equation}
0 = \sum_{r=0}^{N+1}\bar{c}\left( 2r+1, N+1-r; 2N+2,q\right)  \label{eq:60026a}
\end{equation}
The eigenvalue of $q$ is written by $q_{2N+2}^m$ where $m = 0,1,2,\cdots,2N+2 $; $q_{2N+2}^0 < q_{2N+2}^1 < \cdots < q_{2N+2}^{2N+2}$. Its eigenfunction is given by
\begin{eqnarray} 
y(x) &=& _pH_c^{(a)}S_{2N+2,m} \left( \alpha = \gamma -2N-3, \beta, \gamma, \delta, q=q_{2N+2}^m; \eta = \frac{1}{2}\beta x^2 \right) \nonumber\\
&=& \sum_{r=0}^{N+1} y_{2r}^{N+1-r}\left( 2N+2,q_{2N+2}^m;x\right) + \sum_{r=0}^{N} y_{2r+1}^{N-r}\left( 2N+2,q_{2N+2}^m;x\right) 
\hspace{1.5cm}\label{eq:60026b}
\end{eqnarray}
In the above,
\begin{eqnarray}
\bar{c}(0,n;j,q)  &=& \frac{\left( -\frac{j}{2}\right)_{n}}{\left( \frac{3}{2}-\frac{\gamma }{2}\right)_{n} \left( 1\right)_{n}} \left( \frac{\beta }{2} \right)^{n}\label{eq:60027a}\\
\bar{c}(1,n;j,q) &=& \sum_{i_0=0}^{n}\frac{\left( i_0+\frac{1}{2}-\frac{\gamma }{2}\right)\left( i_0 +\Gamma _0 \right) -\frac{q}{4}}{\left( i_0+1-\frac{\gamma }{2}\right) \left( i_0 +\frac{1}{2}\right)} \frac{\left( -\frac{j}{2}\right)_{i_0} }{\left( \frac{3}{2}-\frac{\gamma }{2}\right)_{i_0} \left( 1\right)_{i_0}}  \nonumber\\
&&\times  \frac{\left( \frac{1}{2}-\frac{j}{2} \right)_{n} \left( 2-\frac{\gamma }{2}\right)_{i_0} \left( \frac{3}{2}\right)_{i_0}}{\left( \frac{1}{2}-\frac{j}{2}\right)_{i_0} \left( 2-\frac{\gamma }{2}\right)_{n} \left( \frac{3}{2}\right)_{n}} \left( \frac{\beta }{2} \right)^{n } \hspace{1.5cm}\label{eq:60027b}\\
\bar{c}(\tau ,n;j,q) &=& \sum_{i_0=0}^{n}\frac{\left( i_0+\frac{1}{2}-\frac{\gamma }{2}\right)\left( i_0 +\Gamma _0 \right) -\frac{q}{4}}{\left( i_0+1-\frac{\gamma }{2}\right) \left( i_0 +\frac{1}{2}\right)} \frac{\left( -\frac{j}{2}\right)_{i_0} }{\left( \frac{3}{2}-\frac{\gamma }{2}\right)_{i_0} \left( 1\right)_{i_0}}  \nonumber\\
&&\times  \prod_{k=1}^{\tau -1} \left( \sum_{i_k = i_{k-1}}^{n} \frac{\left( i_k+ \frac{k}{2}+\frac{1}{2}-\frac{\gamma }{2}\right)\left( i_k +\Gamma _k \right) -\frac{q}{4}}{\left( i_k+\frac{k}{2}+1-\frac{\gamma }{2}\right) \left( i_k+\frac{k}{2} +\frac{1}{2}\right)} \right. \nonumber\\
&&\times \left. \frac{\left( \frac{k}{2}-\frac{j}{2}\right)_{i_k} \left( \frac{k}{2} +\frac{3}{2}-\frac{\gamma }{2}\right)_{i_{k-1}} \left( \frac{k}{2}+ 1\right)_{i_{k-1}}}{\left( \frac{k}{2}-\frac{j}{2}\right)_{i_{k-1}} \left( \frac{k}{2}+\frac{3}{2}- \frac{\gamma }{2}\right)_{i_k} \left( \frac{k}{2}+1\right)_{i_k}} \right) \nonumber\\ 
&&\times \frac{\left( \frac{\tau }{2} -\frac{j}{2}\right)_{n} \left( \frac{\tau }{2} +\frac{3}{2}-\frac{\gamma }{2}\right)_{i_{\tau -1}} \left( \frac{\tau }{2}+1\right)_{i_{\tau -1}}}{\left( \frac{\tau }{2}-\frac{j}{2}\right)_{i_{\tau -1}} \left( \frac{\tau }{2} +\frac{3}{2}-\frac{\gamma }{2}\right)_{n} \left( \frac{\tau }{2}+1\right)_{n}} \left( \frac{\beta }{2} \right)^{n } \label{eq:60027c} 
\end{eqnarray}
\begin{eqnarray}
y_0^m(j,q;x) &=& x^{1-\gamma }  \sum_{i_0=0}^{m} \frac{\left( -\frac{j}{2}\right)_{i_0}}{\left( \frac{3}{2}-\frac{\gamma }{2}\right)_{i_0} \left( 1\right)_{i_0}} \eta ^{i_0} \label{eq:60028a}\\
y_1^m(j,q;x) &=& x^{1-\gamma } \left\{\sum_{i_0=0}^{m} \frac{\left( i_0+\frac{1}{2}-\frac{\gamma }{2}\right)\left(i_0 +\Gamma _0 \right) -\frac{q}{4}}{\left( i_0+1-\frac{\gamma }{2}\right) \left( i_0 +\frac{1}{2}\right)} \frac{\left( -\frac{j}{2}\right)_{i_0} }{\left( \frac{3}{2}-\frac{\gamma }{2}\right)_{i_0} \left( 1\right)_{i_0}} \right.  \nonumber\\
&&\times \left. \sum_{i_1 = i_0}^{m} \frac{\left( \frac{1}{2}-\frac{j}{2} \right)_{i_1} \left( 2-\frac{\gamma }{2}\right)_{i_0} \left( \frac{3}{2}\right)_{i_0}}{\left( \frac{1}{2}-\frac{j}{2} \right)_{i_0} \left( 2-\frac{\gamma }{2}\right)_{i_1} \left( \frac{3}{2}\right)_{i_1}} \eta ^{i_1}\right\}x 
\hspace{1.5cm}\label{eq:60028b}
\end{eqnarray}
\begin{eqnarray}
y_{\tau }^m(j,q;x) &=&  x^{1-\gamma } \left\{ \sum_{i_0=0}^{m} \frac{\left( i_0+\frac{1}{2}-\frac{\gamma }{2}\right)\left( i_0 +\Gamma _0 \right) -\frac{q}{4}}{\left( i_0+1-\frac{\gamma }{2}\right) \left( i_0 +\frac{1}{2}\right)} \frac{\left( -\frac{j}{2}\right)_{i_0} }{\left( \frac{3}{2}-\frac{\gamma }{2}\right)_{i_0} \left( 1\right)_{i_0}} \right.\nonumber\\
&&\times \prod_{k=1}^{\tau -1} \left( \sum_{i_k = i_{k-1}}^{m} \frac{\left( i_k+ \frac{k}{2}+\frac{1}{2}-\frac{\gamma }{2}\right)\left( i_k +\Gamma _k \right) -\frac{q}{4}}{\left( i_k+\frac{k}{2}+1-\frac{\gamma }{2}\right) \left( i_k+\frac{k}{2} +\frac{1}{2}\right)} \right.  \nonumber\\
&&\times  \left. \frac{\left( \frac{k}{2}-\frac{j}{2}\right)_{i_k} \left( \frac{k}{2} +\frac{3}{2}-\frac{\gamma }{2}\right)_{i_{k-1}} \left(  \frac{k}{2}+1\right)_{i_{k-1}}}{\left( \frac{k}{2}-\frac{j}{2}\right)_{i_{k-1}} \left( \frac{k}{2} +\frac{3}{2}-\frac{\gamma }{2}\right)_{i_k} \left( \frac{k}{2}+1\right)_{i_k}} \right) \nonumber\\
&&\times \left. \sum_{i_{\tau } = i_{\tau -1}}^{m} \frac{\left( \frac{\tau }{2}-\frac{j}{2}\right)_{i_{\tau }} \left( \frac{\tau }{2} +\frac{3}{2}-\frac{\gamma }{2}\right)_{i_{\tau -1}} \left( \frac{\tau }{2}+1\right)_{i_{\tau -1}}}{\left( \frac{\tau }{2}-\frac{j}{2}\right)_{i_{\tau -1}} \left( \frac{\tau }{2} +\frac{3}{2}-\frac{\gamma }{2}\right)_{i_{\tau }} \left( \frac{\tau }{2}+1\right)_{i_{\tau }}} \eta ^{i_{\tau }}\right\} x^{\tau } \hspace{1.5cm}\label{eq:60028c} 
\end{eqnarray}
where
\begin{equation}
\begin{cases} \tau \geq 2 \cr
\Gamma _0 = \frac{1}{2}\left( -\beta +\delta \right) \cr
\Gamma _k = \frac{1}{2}\left( -\beta +\delta +k \right)
\end{cases}\nonumber
\end{equation}
\end{enumerate}
\end{remark}
\subsection{The first species complete polynomial of Confluent Heun equation using R3TRF}
\begin{theorem}
In chapter 2, the general expression of a function $y(x)$ for the first species complete polynomial using reversible 3-term recurrence formula and its algebraic equation for the determination of an accessory parameter in $A_n$ term are given by
\begin{enumerate} 
\item As $B_1=0$,
\begin{equation}
0 =\bar{c}(0,1) \label{eq:60029a}
\end{equation}
\begin{equation}
y(x) = y_{0}^{0}(x) \label{eq:60029b}
\end{equation}
\item As $B_2=0$, 
\begin{equation}
0 = \bar{c}(0,2)+\bar{c}(1,0) \label{eq:60030a}
\end{equation}
\begin{equation}
y(x)= y_{0}^{1}(x) \label{eq:60030b}
\end{equation}
\item As $B_{2N+3}=0$ where $N \in \mathbb{N}_{0}$,
\begin{equation}
0  =  \sum_{r=0}^{N+1}\bar{c}\left( r, 2(N-r)+3\right) \label{eq:60031a}
\end{equation}
\begin{equation}
y(x)= \sum_{r=0}^{N+1} y_{r}^{2(N+1-r)}(x) \label{eq:60031b}
\end{equation}
\item As $B_{2N+4}=0$ where$N \in \mathbb{N}_{0}$,
\begin{equation}
0  = \sum_{r=0}^{N+2}\bar{c}\left( r, 2(N+2-r)\right) \label{eq:60032a}
\end{equation}
\begin{equation}
y(x)=  \sum_{r=0}^{N+1} y_{r}^{2(N-r)+3}(x) \label{eq:60032b}
\end{equation}
In the above,
\begin{eqnarray}
\bar{c}(0,n) &=& \prod _{i_0=0}^{n-1}A_{i_0} \label{eq:60033a}\\
\bar{c}(1,n) &=& \sum_{i_0=0}^{n} \left\{ B_{i_0+1} \prod _{i_1=0}^{i_0-1}A_{i_1} \prod _{i_2=i_0}^{n-1}A_{i_2+2} \right\} \label{eq:60033b}\\
\bar{c}(\tau ,n) &=& \sum_{i_0=0}^{n} \left\{B_{i_0+1}\prod _{i_1=0}^{i_0-1} A_{i_1} 
\prod _{k=1}^{\tau -1} \left( \sum_{i_{2k}= i_{2(k-1)}}^{n} B_{i_{2k}+(2k+1)}\prod _{i_{2k+1}=i_{2(k-1)}}^{i_{2k}-1}A_{i_{2k+1}+2k}\right) \right. \nonumber\\
&&\times \left. \prod _{i_{2\tau} = i_{2(\tau -1)}}^{n-1} A_{i_{2\tau }+ 2\tau} \right\} 
\hspace{1cm}\label{eq:60033c}
\end{eqnarray}
and
\begin{eqnarray}
y_0^m(x) &=& c_0 x^{\lambda} \sum_{i_0=0}^{m} \left\{ \prod _{i_1=0}^{i_0-1}A_{i_1} \right\} x^{i_0 } \label{eq:60034a}\\
y_1^m(x) &=& c_0 x^{\lambda} \sum_{i_0=0}^{m}\left\{ B_{i_0+1} \prod _{i_1=0}^{i_0-1}A_{i_1}  \sum_{i_2=i_0}^{m} \left\{ \prod _{i_3=i_0}^{i_2-1}A_{i_3+2} \right\}\right\} x^{i_2+2 } \label{eq:60034b}\\
y_{\tau }^m(x) &=& c_0 x^{\lambda} \sum_{i_0=0}^{m} \left\{B_{i_0+1}\prod _{i_1=0}^{i_0-1} A_{i_1} 
\prod _{k=1}^{\tau -1} \left( \sum_{i_{2k}= i_{2(k-1)}}^{m} B_{i_{2k}+(2k+1)}\prod _{i_{2k+1}=i_{2(k-1)}}^{i_{2k}-1}A_{i_{2k+1}+2k}\right) \right. \nonumber\\
&&\times \left. \sum_{i_{2\tau} = i_{2(\tau -1)}}^{m} \left( \prod _{i_{2\tau +1}=i_{2(\tau -1)}}^{i_{2\tau}-1} A_{i_{2\tau +1}+ 2\tau} \right) \right\} x^{i_{2\tau}+2\tau }\hspace{1cm}\mathrm{where}\;\tau \geq 2
\label{eq:60034c}
\end{eqnarray}
\end{enumerate}
\end{theorem}
(\ref{eq:6005a}) can be described as the different version such as
\begin{equation}
A_n =\frac{\left( n+ \Delta _{0}^{-}(q) \right) \left( n+ \Delta _{0}^{+}(q) \right)}{(n+1+\lambda )(n+\gamma +\lambda )}\label{eq:60035}
\end{equation}
where
\begin{equation}
\begin{cases} 
\Delta _{k}^{\pm}(q) = \frac{\varphi +2\lambda +4k \pm \sqrt{\varphi ^2 +4q}}{2}\cr
\varphi = -\beta +\gamma +\delta -1
\end{cases}
\nonumber
\end{equation}
According to (\ref{eq:6006}), $c_{j+1}=0$ is clearly an algebraic equation in $q$ of degree $j+1$ and thus has $j+1$ zeros denoted them by $q_j^m$ eigenvalues where $m = 0,1,2, \cdots, j$. They can be arranged in the following order: $q_j^0 < q_j^1 < q_j^2 < \cdots < q_j^j$.
 
Substitute (\ref{eq:60013}) and (\ref{eq:60035}) into (\ref{eq:60033a})--(\ref{eq:60034c}).

As $B_{1}= c_{1}=0$, take the new (\ref{eq:60033a}) into (\ref{eq:60029a}) putting $j=0$. Substitute the new (\ref{eq:60034a}) into (\ref{eq:60029b}) putting $j=0$.

As $B_{2}= c_{2}=0$, take the new (\ref{eq:60033a}) and (\ref{eq:60033b}) into (\ref{eq:60030a}) putting $j=1$. Substitute the new (\ref{eq:60034a}) into (\ref{eq:60030b}) putting $j=1$ and $q=q_1^m$. 

As $B_{2N+3}= c_{2N+3}=0$, take the new (\ref{eq:60033a})--(\ref{eq:60033c}) into (\ref{eq:60031a}) putting $j=2N+2$. Substitute the new 
(\ref{eq:60034a})--(\ref{eq:60034c}) into (\ref{eq:60031b}) putting $j=2N+2$ and $q=q_{2N+2}^m$.

As $B_{2N+4}= c_{2N+4}=0$, take the new (\ref{eq:60033a})--(\ref{eq:60033c}) into (\ref{eq:60032a}) putting $j=2N+3$. Substitute the new 
(\ref{eq:60034a})--(\ref{eq:60034c}) into (\ref{eq:60032b}) putting $j=2N+3$ and $q=q_{2N+3}^m$.

After the replacement process, the general expression of power series of the CHE about $x=0$ for the first species complete polynomial using reversible 3-term recurrence formula and its algebraic equation for the determination of an accessory parameter $q$ are given by
\begin{enumerate} 
\item As $\alpha = -\lambda $,

An algebraic equation of degree 1 for the determination of $q$ is given by
\begin{equation}
0= \bar{c}(0,1;0,q)= q - \lambda \left( -\beta +\gamma +\delta  -1+\lambda \right) \label{eq:60036a}
\end{equation}
The eigenvalue of $q$ is written by $q_0^0$. Its eigenfunction is given by
\begin{equation}
y(x) = y_0^0\left( 0,q_0^0;x\right)= c_0 x^{\lambda } \label{eq:60036b}  
\end{equation}
\item As $\alpha =-1-\lambda $,

An algebraic equation of degree 2 for the determination of $q$ is given by
\begin{eqnarray}
0 &=& \bar{c}(0,2;1,q)+\bar{c}(1,0;1,q) \nonumber\\
&=& -\beta (1 +\lambda )(\gamma +\lambda ) 
+\prod_{l=0}^{1}\Big( (\lambda +l)(-\beta +\gamma +\delta -1+l+\lambda ) -q\Big) \hspace{1cm}\label{eq:60037a}
\end{eqnarray}
The eigenvalue of $q$ is written by $q_1^m$ where $m = 0,1 $; $q_{1}^0 < q_{1}^1$. Its eigenfunction is given by
\begin{equation}
y(x) = y_{0}^{1}\left( 1,q_1^m;x\right)= c_0 x^{\lambda } \left\{ 1+\frac{\lambda (-\beta +\gamma +\delta -1+\lambda )-q_1^m}{(1+\lambda )(\gamma +\lambda )}x \right\} \label{eq:60037b}  
\end{equation}
\item As $\alpha =-2N-2-\lambda $ where $N \in \mathbb{N}_{0}$,

An algebraic equation of degree $2N+3$ for the determination of $q$ is given by
\begin{equation}
0 =  \sum_{r=0}^{N+1}\bar{c}\left( r, 2(N-r)+3; 2N+2,q\right)  \label{eq:60038a}
\end{equation}
The eigenvalue of $q$ is written by $q_{2N+2}^m$ where $m = 0,1,2,\cdots,2N+2 $; $q_{2N+2}^0 < q_{2N+2}^1 < \cdots < q_{2N+2}^{2N+2}$. Its eigenfunction is given by 
\begin{equation} 
y(x) = \sum_{r=0}^{N+1} y_{r}^{2(N+1-r)}\left( 2N+2, q_{2N+2}^m; x \right)  
\label{eq:60038b} 
\end{equation}
\item As $\alpha =-2N-3-\lambda $ where $N \in \mathbb{N}_{0}$,

An algebraic equation of degree $2N+4$ for the determination of $q$ is given by
\begin{equation}  
0 =  \sum_{r=0}^{N+2}\bar{c}\left( r, 2(N+2-r); 2N+3,q\right) \label{eq:60039a}
\end{equation}
The eigenvalue of $q$ is written by $q_{2N+3}^m$ where $m = 0,1,2,\cdots,2N+3 $; $q_{2N+3}^0 < q_{2N+3}^1 < \cdots < q_{2N+3}^{2N+3}$. Its eigenfunction is given by
\begin{equation} 
y(x) =  \sum_{r=0}^{N+1} y_{r}^{2(N-r)+3} \left( 2N+3,q_{2N+3}^m;x\right) \label{eq:60039b}
\end{equation}
In the above,
\begin{eqnarray}
\bar{c}(0,n;j,q)  &=& \frac{\left( \Delta_0^{-} \left( q\right) \right)_{n}\left( \Delta_0^{+} \left( q\right) \right)_{n}}{\left( 1+ \lambda \right)_{n} \left( \gamma + \lambda \right)_{n}} \label{eq:60040a}\\
\bar{c}(1,n;j,q) &=& \beta \sum_{i_0=0}^{n}\frac{\left( i_0 -j\right) }{\left( i_0+2+ \lambda \right) \left( i_0+1+ \gamma + \lambda \right)} \frac{ \left( \Delta_0^{-} \left( q\right) \right)_{i_0}\left( \Delta_0^{+} \left( q\right) \right)_{i_0}}{\left( 1+ \lambda \right)_{i_0} \left( \gamma + \lambda \right)_{i_0}} \nonumber\\
&&\times  \frac{ \left( \Delta_1^{-} \left( q\right) \right)_{n}\left( \Delta_1^{+} \left( q\right) \right)_{n} \left( 3 + \lambda \right)_{i_0} \left( 2+ \gamma + \lambda \right)_{i_0}}{\left( \Delta_1^{-} \left( q\right) \right)_{i_0}\left( \Delta_1^{+} \left(  q\right) \right)_{i_0}\left( 3 + \lambda \right)_{n} \left( 2+ \gamma + \lambda \right)_{n}} \label{eq:60040b}\\
\bar{c}(\tau ,n;j,q) &=& \beta ^{\tau} \sum_{i_0=0}^{n}\frac{\left( i_0 -j\right) }{\left( i_0+2+ \lambda \right) \left( i_0+1+ \gamma + \lambda \right)}  \frac{ \left( \Delta_0^{-} \left( q\right) \right)_{i_0}\left( \Delta_0^{+} \left( q\right) \right)_{i_0}}{\left( 1+ \lambda \right)_{i_0} \left( \gamma + \lambda \right)_{i_0}}  \nonumber\\
&&\times \prod_{k=1}^{\tau -1} \left( \sum_{i_k = i_{k-1}}^{n} \frac{\left( i_k+ 2k-j\right) }{\left( i_k+2k+2+ \lambda \right) \left( i_k+2k+1+ \gamma + \lambda \right)} \right. \nonumber\\
&&\times  \left. \frac{ \left( \Delta_k^{-} \left( q\right) \right)_{i_k}\left( \Delta_k^{+} \left( q\right) \right)_{i_k} \left( 2k+1 + \lambda \right)_{i_{k-1}} \left( 2k+ \gamma + \lambda \right)_{i_{k-1}}}{\left( \Delta_k^{-} \left( q\right) \right)_{i_{k-1}}\left( \Delta_k^{+} \left( q\right) \right)_{i_{k-1}}\left( 2k+1 + \lambda \right)_{i_k} \left( 2k+ \gamma + \lambda \right)_{i_k}} \right) \nonumber\\
&&\times \frac{ \left( \Delta_{\tau }^{-} \left( q\right) \right)_{n}\left( \Delta_{\tau }^{+} \left( q\right) \right)_{n} \left( 2\tau +1+ \lambda \right)_{i_{\tau -1}} \left( 2\tau + \gamma + \lambda \right)_{i_{\tau -1}}}{\left( \Delta_{\tau }^{-} \left( q\right) \right)_{i_{\tau -1}}\left( \Delta_{\tau }^{+} \left( q\right) \right)_{i_{\tau -1}}\left( 2\tau +1 + \lambda \right)_{n} \left( 2\tau + \gamma + \lambda \right)_{n}}  \hspace{1cm}\label{eq:60040c} 
\end{eqnarray}
\begin{eqnarray}
y_0^m(j,q;x) &=& c_0 x^{\lambda }  \sum_{i_0=0}^{m} \frac{\left( \Delta_0^{-} \left( q\right) \right)_{i_0}\left( \Delta_0^{+} \left( q\right) \right)_{i_0}}{\left( 1+ \lambda \right)_{i_0} \left( \gamma + \lambda \right)_{i_0}} x^{i_0} \label{eq:60041a}\\
y_1^m(j,q;x) &=& c_0 x^{\lambda } \left\{\sum_{i_0=0}^{m}\frac{\left( i_0 -j\right) }{\left( i_0+2+ \lambda \right) \left( i_0+1+ \gamma + \lambda \right)} \frac{ \left( \Delta_0^{-} \left( q\right) \right)_{i_0}\left( \Delta_0^{+} \left( q\right) \right)_{i_0}}{\left( 1+ \lambda \right)_{i_0} \left( \gamma + \lambda \right)_{i_0}} \right. \nonumber\\
&&\times \left. \sum_{i_1 = i_0}^{m} \frac{ \left( \Delta_1^{-} \left( q\right) \right)_{i_1}\left( \Delta_1^{+} \left( q\right) \right)_{i_1} \left( 3 + \lambda \right)_{i_0} \left( 2+ \gamma + \lambda \right)_{i_0}}{\left( \Delta_1^{-} \left( q\right) \right)_{i_0}\left( \Delta_1^{+} \left( q\right) \right)_{i_0}\left( 3 + \lambda \right)_{i_1} \left( 2+ \gamma + \lambda \right)_{i_1}} x^{i_1}\right\} \rho  
\hspace{1.5cm}\label{eq:60041b}
\end{eqnarray}
\begin{eqnarray}
y_{\tau }^m(j,q;x) &=& c_0 x^{\lambda } \left\{ \sum_{i_0=0}^{m} \frac{\left( i_0 -j\right) }{\left( i_0+2+ \lambda \right) \left( i_0+1+ \gamma + \lambda \right)} \frac{ \left( \Delta_0^{-} \left( q\right) \right)_{i_0}\left( \Delta_0^{+} \left( q\right) \right)_{i_0}}{\left( 1+ \lambda \right)_{i_0} \left( \gamma + \lambda \right)_{i_0}} \right.\nonumber\\
&&\times \prod_{k=1}^{\tau -1} \left( \sum_{i_k = i_{k-1}}^{m} \frac{\left( i_k+ 2k-j\right) }{\left( i_k+2k+2+ \lambda \right) \left( i_k+2k+1+ \gamma + \lambda \right)} \right. \nonumber\\
&&\times  \left. \frac{ \left( \Delta_k^{-} \left( q\right) \right)_{i_k}\left( \Delta_k^{+} \left( q\right) \right)_{i_k} \left( 2k+1 + \lambda \right)_{i_{k-1}} \left( 2k+ \gamma + \lambda \right)_{i_{k-1}}}{\left( \Delta_k^{-} \left( q\right) \right)_{i_{k-1}}\left( \Delta_k^{+} \left( q\right) \right)_{i_{k-1}}\left( 2k+1 + \lambda \right)_{i_k} \left( 2k+ \gamma + \lambda \right)_{i_k}} \right) \nonumber\\
&&\times \left. \sum_{i_{\tau } = i_{\tau -1}}^{m}  \frac{ \left( \Delta_{\tau }^{-} \left( q\right) \right)_{i_{\tau }}\left( \Delta_{\tau }^{+} \left( q\right) \right)_{i_{\tau }} \left( 2\tau +1 + \lambda \right)_{i_{\tau -1}} \left( 2\tau + \gamma + \lambda \right)_{i_{\tau -1}}}{\left( \Delta_{\tau }^{-} \left( q\right) \right)_{i_{\tau -1}}\left( \Delta_{\tau }^{+} \left( q\right) \right)_{i_{\tau -1}}\left( 2\tau +1 + \lambda \right)_{i_\tau } \left( 2\tau + \gamma + \lambda \right)_{i_{\tau }}} x^{i_{\tau }}\right\} \rho ^{\tau } \hspace{1.5cm}\label{eq:60041c} 
\end{eqnarray}
where
\begin{equation}
\begin{cases} \tau \geq 2 \cr
\rho = \beta x^2 \cr 
\Delta _{k}^{\pm}(q) = \frac{\varphi +2\lambda +4k \pm \sqrt{\varphi ^2 +4q}}{2}\cr
\varphi = -\beta +\gamma +\delta -1
\end{cases}\nonumber
\end{equation}
\end{enumerate}
Put $c_0$= 1 as $\lambda =0$ for the first kind of independent solutions of the CHE and $\lambda = 1-\gamma $ for the second one in (\ref{eq:60036a})--(\ref{eq:60041c}). 
\begin{remark}
The power series expansion of Confluent Heun equation of the first kind for the first species complete polynomial using R3TRF about $x=0$ is given by
\begin{enumerate} 
\item As $\alpha =0$ and $q=q_0^0=0$,

The eigenfunction is given by
\begin{equation}
y(x) =\; _pH_c^{(a)}F_{0,0}^R \left( \alpha =0, \beta, \gamma, \delta, q=q_0^0=0; \rho = \beta x^2 \right) =1\label{eq:60042}
\end{equation}
\item As $\alpha =-1$,

An algebraic equation of degree 2 for the determination of $q$ is given by
\begin{equation}
0 = -\beta \gamma 
+\prod_{l=0}^{1}\Big( l(-\beta +\gamma +\delta -1+l)-q\Big) \label{eq:60043a}
\end{equation}
The eigenvalue of $q$ is written by $q_1^m$ where $m = 0,1 $; $q_{1}^0 < q_{1}^1$. Its eigenfunction is given by
\begin{eqnarray}
y(x) &=& _pH_c^{(a)}F_{1,m}^R \left( \alpha =-1, \beta, \gamma, \delta, q=q_1^m; \rho = \beta x^2 \right)\nonumber\\
&=&  1-\frac{ q_1^m}{\gamma }x \label{eq:60043b}  
\end{eqnarray}
\item As $\alpha =-2N-2 $ where $N \in \mathbb{N}_{0}$,

An algebraic equation of degree $2N+3$ for the determination of $q$ is given by
\begin{equation}
0 = \sum_{r=0}^{N+1}\bar{c}\left( r, 2(N-r)+3; 2N+2,q\right)  \label{eq:60044a}
\end{equation}
The eigenvalue of $q$ is written by $q_{2N+2}^m$ where $m = 0,1,2,\cdots,2N+2 $; $q_{2N+2}^0 < q_{2N+2}^1 < \cdots < q_{2N+2}^{2N+2}$. Its eigenfunction is given by 
\begin{eqnarray} 
y(x) &=& _pH_c^{(a)}F_{2N+2,m}^R \left( \alpha =-2N-2, \beta, \gamma, \delta, q=q_{2N+2}^m ; \rho = \beta x^2 \right)\nonumber\\
&=& \sum_{r=0}^{N+1} y_{r}^{2(N+1-r)}\left( 2N+2, q_{2N+2}^m; x \right)  
\label{eq:60044b} 
\end{eqnarray}
\item As $\alpha =-2N-3 $ where $N \in \mathbb{N}_{0}$,

An algebraic equation of degree $2N+4$ for the determination of $q$ is given by
\begin{equation}  
0 = \sum_{r=0}^{N+2}\bar{c}\left( r, 2(N+2-r); 2N+3,q\right) \label{eq:60045a}
\end{equation}
The eigenvalue of $q$ is written by $q_{2N+3}^m$ where $m = 0,1,2,\cdots,2N+3 $; $q_{2N+3}^0 < q_{2N+3}^1 < \cdots < q_{2N+3}^{2N+3}$. Its eigenfunction is given by
\begin{eqnarray} 
y(x) &=& _pH_c^{(a)}F_{2N+3,m}^R \left( \alpha =-2N-3, \beta, \gamma, \delta, q=q_{2N+3}^m ; \rho = \beta x^2\right)\nonumber\\
&=& \sum_{r=0}^{N+1} y_{r}^{2(N-r)+3} \left( 2N+3,q_{2N+3}^m;x\right) \label{eq:60045b}
\end{eqnarray}
In the above,
\begin{eqnarray}
\bar{c}(0,n;j,q)  &=& \frac{\left( \Delta_0^{-} \left( q\right) \right)_{n}\left( \Delta_0^{+} \left( q\right) \right)_{n}}{\left( 1 \right)_{n} \left( \gamma \right)_{n}} \label{eq:60046a}\\
\bar{c}(1,n;j,q) &=& \beta \sum_{i_0=0}^{n}\frac{\left( i_0 -j\right) }{\left( i_0+2 \right) \left( i_0+1+ \gamma \right)} \frac{ \left( \Delta_0^{-} \left( q\right) \right)_{i_0}\left( \Delta_0^{+} \left( q\right) \right)_{i_0}}{\left( 1 \right)_{i_0} \left( \gamma \right)_{i_0}}  \nonumber\\
&&\times  \frac{ \left( \Delta_1^{-} \left( q\right) \right)_{n}\left( \Delta_1^{+} \left( q\right) \right)_{n} \left( 3\right)_{i_0} \left( 2+ \gamma \right)_{i_0}}{\left( \Delta_1^{-} \left( q\right) \right)_{i_0}\left( \Delta_1^{+} \left(  q\right) \right)_{i_0}\left( 3 \right)_{n} \left( 2+ \gamma \right)_{n}} \label{eq:60046b}\\
\bar{c}(\tau ,n;j,q) &=& \beta ^{\tau} \sum_{i_0=0}^{n}\frac{\left( i_0 -j\right) }{\left( i_0+2 \right) \left( i_0+1+ \gamma \right)}  \frac{ \left( \Delta_0^{-} \left( q\right) \right)_{i_0}\left( \Delta_0^{+} \left( q\right) \right)_{i_0}}{\left( 1\right)_{i_0} \left( \gamma \right)_{i_0}}  \nonumber\\
&&\times \prod_{k=1}^{\tau -1} \left( \sum_{i_k = i_{k-1}}^{n} \frac{\left( i_k+ 2k-j\right) }{\left( i_k+2k+2 \right) \left( i_k+2k+1+ \gamma \right)} \right. \nonumber\\
&&\times   \left. \frac{ \left( \Delta_k^{-} \left( q\right) \right)_{i_k}\left( \Delta_k^{+} \left( q\right) \right)_{i_k} \left( 2k+1 \right)_{i_{k-1}} \left( 2k+ \gamma \right)_{i_{k-1}}}{\left( \Delta_k^{-} \left( q\right) \right)_{i_{k-1}}\left( \Delta_k^{+} \left( q\right) \right)_{i_{k-1}}\left( 2k+1 \right)_{i_k} \left( 2k+ \gamma \right)_{i_k}} \right) \nonumber\\
&&\times \frac{ \left( \Delta_{\tau }^{-} \left( q\right) \right)_{n}\left( \Delta_{\tau }^{+} \left( q\right) \right)_{n} \left( 2\tau +1 \right)_{i_{\tau -1}} \left( 2\tau + \gamma \right)_{i_{\tau -1}}}{\left( \Delta_{\tau }^{-} \left( q\right) \right)_{i_{\tau -1}}\left( \Delta_{\tau }^{+} \left( q\right) \right)_{i_{\tau -1}}\left( 2\tau +1 \right)_{n} \left( 2\tau + \gamma \right)_{n}}  \label{eq:60046c} 
\end{eqnarray}
\begin{eqnarray}
y_0^m(j,q;x) &=& \sum_{i_0=0}^{m} \frac{\left( \Delta_0^{-} \left( q\right) \right)_{i_0}\left( \Delta_0^{+} \left( q\right) \right)_{i_0}}{\left( 1 \right)_{i_0} \left( \gamma \right)_{i_0}} x^{i_0} \label{eq:60047a}\\
y_1^m(j,q;x) &=& \left\{\sum_{i_0=0}^{m}\frac{\left( i_0 -j\right) }{\left( i_0+2 \right) \left( i_0+1+ \gamma \right)} \frac{ \left( \Delta_0^{-} \left( q\right) \right)_{i_0}\left( \Delta_0^{+} \left( q\right) \right)_{i_0}}{\left( 1 \right)_{i_0} \left( \gamma \right)_{i_0}} \right.  \nonumber\\
&&\times \left. \sum_{i_1 = i_0}^{m} \frac{ \left( \Delta_1^{-} \left( q\right) \right)_{i_1}\left( \Delta_1^{+} \left( q\right) \right)_{i_1} \left( 3 \right)_{i_0} \left( 2+ \gamma \right)_{i_0}}{\left( \Delta_1^{-} \left( q\right) \right)_{i_0}\left( \Delta_1^{+} \left( q\right) \right)_{i_0}\left( 3 \right)_{i_1} \left( 2+ \gamma \right)_{i_1}} x^{i_1}\right\} \rho  
\hspace{1cm}\label{eq:60047b}
\end{eqnarray}
\begin{eqnarray}
y_{\tau }^m(j,q;x) &=& \left\{ \sum_{i_0=0}^{m} \frac{\left( i_0 -j\right) }{\left( i_0+2 \right) \left( i_0+1+ \gamma \right)} \frac{ \left( \Delta_0^{-} \left( q\right) \right)_{i_0}\left( \Delta_0^{+} \left( q\right) \right)_{i_0}}{\left( 1 \right)_{i_0} \left( \gamma \right)_{i_0}} \right.\nonumber\\
&&\times \prod_{k=1}^{\tau -1} \left( \sum_{i_k = i_{k-1}}^{m} \frac{\left( i_k+ 2k-j\right) }{\left( i_k+2k+2 \right) \left( i_k+2k+1+ \gamma \right)} \right.  \nonumber\\
&&\times  \left. \frac{ \left( \Delta_k^{-} \left( q\right) \right)_{i_k}\left( \Delta_k^{+} \left( q\right) \right)_{i_k} \left( 2k+1 \right)_{i_{k-1}} \left( 2k+ \gamma \right)_{i_{k-1}}}{\left( \Delta_k^{-} \left( q\right) \right)_{i_{k-1}}\left( \Delta_k^{+} \left( q\right) \right)_{i_{k-1}}\left( 2k+1 + \lambda \right)_{i_k} \left( 2k+ \gamma + \lambda \right)_{i_k}} \right) \nonumber\\
&&\times \left. \sum_{i_{\tau } = i_{\tau -1}}^{m}  \frac{ \left( \Delta_{\tau }^{-} \left( q\right) \right)_{i_{\tau }}\left( \Delta_{\tau }^{+} \left( q\right) \right)_{i_{\tau }} \left( 2\tau +1 \right)_{i_{\tau -1}} \left( 2\tau + \gamma \right)_{i_{\tau -1}}}{\left( \Delta_{\tau }^{-} \left( q\right) \right)_{i_{\tau -1}}\left( \Delta_{\tau }^{+} \left( q\right) \right)_{i_{\tau -1}}\left( 2\tau +1 \right)_{i_\tau } \left( 2\tau + \gamma \right)_{i_{\tau }}} x^{i_{\tau }}\right\} \rho ^{\tau } \hspace{1.5cm}\label{eq:60047c} 
\end{eqnarray}
where
\begin{equation}
\begin{cases} \tau \geq 2 \cr 
\Delta _{k}^{\pm}(q) = \frac{\varphi +4k \pm \sqrt{\varphi ^2 +4q}}{2}\cr
\varphi = -\beta +\gamma +\delta -1
\end{cases}\nonumber
\end{equation}
\end{enumerate}
\end{remark}
\begin{remark}
The power series expansion of Confluent Heun equation of the second kind for the first species complete polynomial using R3TRF about $x=0$ is given by
\begin{enumerate} 
\item As $\alpha = \gamma -1$ and $q=q_0^0= (\gamma -1)(\beta -\delta )$,

The eigenfunction is given by
\begin{eqnarray}
y(x) &=& _pH_c^{(a)}S_{0,0}^R \left( \alpha =\gamma -1, \beta, \gamma, \delta, q=q_0^0=(\gamma -1)(\beta -\delta ); \rho = \beta x^2 \right) \nonumber\\
&=& x^{1-\gamma } \label{eq:60048}
\end{eqnarray}
\item As $\alpha =\gamma -2$,

An algebraic equation of degree 2 for the determination of $q$ is given by
\begin{equation}
0 = \beta ( \gamma -2)  
+\prod_{l=0}^{1}\Big( (\gamma -1-l)(\beta -\delta -l)-q \Big) \label{eq:60049a}
\end{equation}
The eigenvalue of $q$ is written by $q_1^m$ where $m = 0,1 $; $q_{1}^0 < q_{1}^1$. Its eigenfunction is given by
\begin{eqnarray}
y(x) &=& _pH_c^{(a)}S_{1,m}^R \left( \alpha =\gamma -2, \beta, \gamma, \delta, q=q_1^m; \rho = \beta x^2\right) \nonumber\\
&=& x^{1-\gamma } \left\{ 1+\frac{ (\gamma -1)(\beta -\delta )-q_1^m}{ (2-\gamma )}x \right\} \label{eq:60049b}  
\end{eqnarray}
\item As $\alpha =\gamma -2N-3  $ where $N \in \mathbb{N}_{0}$,

An algebraic equation of degree $2N+3$ for the determination of $q$ is given by
\begin{equation}
0 = \sum_{r=0}^{N+1}\bar{c}\left( r, 2(N-r)+3; 2N+2,q\right)  \label{eq:60050a}
\end{equation}
The eigenvalue of $q$ is written by $q_{2N+2}^m$ where $m = 0,1,2,\cdots,2N+2 $; $q_{2N+2}^0 < q_{2N+2}^1 < \cdots < q_{2N+2}^{2N+2}$. Its eigenfunction is given by 
\begin{eqnarray} 
y(x) &=& _pH_c^{(a)}S_{2N+2,m}^R \left( \alpha =\gamma -2N-3, \beta, \gamma, \delta, q=q_{2N+2}^m; \rho = \beta x^2 \right) \nonumber\\
&=& \sum_{r=0}^{N+1} y_{r}^{2(N+1-r)}\left( 2N+2, q_{2N+2}^m; x \right)  
\label{eq:60050b} 
\end{eqnarray}
\item As $\alpha =\gamma -2N-4 $ where $N \in \mathbb{N}_{0}$,

An algebraic equation of degree $2N+4$ for the determination of $q$ is given by
\begin{equation}  
0 = \sum_{r=0}^{N+2}\bar{c}\left( r, 2(N+2-r); 2N+3,q\right) \label{eq:60051a}
\end{equation}
The eigenvalue of $q$ is written by $q_{2N+3}^m$ where $m = 0,1,2,\cdots,2N+3 $; $q_{2N+3}^0 < q_{2N+3}^1 < \cdots < q_{2N+3}^{2N+3}$. Its eigenfunction is given by
\begin{eqnarray} 
y(x) &=& _pH_c^{(a)}S_{2N+3,m}^R \left( \alpha =\gamma -2N-4, \beta, \gamma, \delta, q=q_{2N+3}^m; \rho = \beta x^2\right) \nonumber\\
&=& \sum_{r=0}^{N+1} y_{r}^{2(N-r)+3} \left( 2N+3,q_{2N+3}^m;x\right) \label{eq:60051b}
\end{eqnarray}
In the above,
\begin{eqnarray}
\bar{c}(0,n;j,q)  &=& \frac{\left( \Delta_0^{-} \left( q\right) \right)_{n}\left( \Delta_0^{+} \left( q\right) \right)_{n}}{\left( 2-\gamma \right)_{n} \left( 1\right)_{n}} \label{eq:60052a}\\
\bar{c}(1,n;j,q) &=& \beta \sum_{i_0=0}^{n}\frac{\left( i_0 -j\right) }{\left( i_0+3-\gamma \right) \left( i_0+2 \right)} \frac{ \left( \Delta_0^{-} \left( q\right) \right)_{i_0}\left( \Delta_0^{+} \left( q\right) \right)_{i_0}}{\left( 2-\gamma \right)_{i_0} \left( 1\right)_{i_0}}  \nonumber\\
&&\times \frac{ \left( \Delta_1^{-} \left( q\right) \right)_{n}\left( \Delta_1^{+} \left( q\right) \right)_{n} \left( 4-\gamma \right)_{i_0} \left( 3\right)_{i_0}}{\left( \Delta_1^{-} \left( q\right) \right)_{i_0}\left( \Delta_1^{+} \left(  q\right) \right)_{i_0}\left( 4-\gamma \right)_{n} \left( 3\right)_{n}} \label{eq:60052b}\\
\bar{c}(\tau ,n;j,q) &=& \beta ^{\tau} \sum_{i_0=0}^{n}\frac{\left( i_0 -j\right) }{\left( i_0+3-\gamma \right) \left( i_0+2\right)}  \frac{ \left( \Delta_0^{-} \left( q\right) \right)_{i_0}\left( \Delta_0^{+} \left( q\right) \right)_{i_0}}{\left( 2-\gamma \right)_{i_0} \left( 1\right)_{i_0}}  \nonumber\\
&&\times \prod_{k=1}^{\tau -1} \left( \sum_{i_k = i_{k-1}}^{n} \frac{\left( i_k+ 2k-j\right) }{\left( i_k+2k+3-\gamma \right) \left( i_k+2k+2\right)} \right. \nonumber\\
&&\times \left. \frac{ \left( \Delta_k^{-} \left( q\right) \right)_{i_k}\left( \Delta_k^{+} \left( q\right) \right)_{i_k} \left( 2k+2-\gamma \right)_{i_{k-1}} \left( 2k+1 \right)_{i_{k-1}}}{\left( \Delta_k^{-} \left( q\right) \right)_{i_{k-1}}\left( \Delta_k^{+} \left( q\right) \right)_{i_{k-1}}\left( 2k+2-\gamma  \right)_{i_k} \left( 2k+1\right)_{i_k}} \right) \nonumber\\
&&\times \frac{ \left( \Delta_{\tau }^{-} \left( q\right) \right)_{n}\left( \Delta_{\tau }^{+} \left( q\right) \right)_{n} \left( 2\tau +2-\gamma  \right)_{i_{\tau -1}} \left( 2\tau +1 \right)_{i_{\tau -1}}}{\left( \Delta_{\tau }^{-} \left( q\right) \right)_{i_{\tau -1}}\left( \Delta_{\tau }^{+} \left( q\right) \right)_{i_{\tau -1}}\left( 2\tau +2-\gamma  \right)_{n} \left( 2\tau +1\right)_{n}}  \label{eq:60052c} 
\end{eqnarray}
\begin{eqnarray}
y_0^m(j,q;x) &=& x^{1-\gamma } \sum_{i_0=0}^{m} \frac{\left( \Delta_0^{-} \left( q\right) \right)_{i_0}\left( \Delta_0^{+} \left( q\right) \right)_{i_0}}{\left( 2-\gamma  \right)_{i_0} \left( 1\right)_{i_0}} x^{i_0} \label{eq:60053a}\\
y_1^m(j,q;x) &=& x^{1-\gamma } \left\{\sum_{i_0=0}^{m}\frac{\left( i_0 -j\right) }{\left( i_0+3-\gamma \right) \left( i_0+2 \right)} \frac{ \left( \Delta_0^{-} \left( q\right) \right)_{i_0}\left( \Delta_0^{+} \left( q\right) \right)_{i_0}}{\left( 2-\gamma \right)_{i_0} \left( 1\right)_{i_0}} \right. \nonumber\\
&&\times \left. \sum_{i_1 = i_0}^{m} \frac{ \left( \Delta_1^{-} \left( q\right) \right)_{i_1}\left( \Delta_1^{+} \left( q\right) \right)_{i_1} \left( 4-\gamma  \right)_{i_0} \left( 3\right)_{i_0}}{\left( \Delta_1^{-} \left( q\right) \right)_{i_0}\left( \Delta_1^{+} \left( q\right) \right)_{i_0}\left( 4-\gamma  \right)_{i_1} \left( 3\right)_{i_1}} x^{i_1}\right\} \rho  
\label{eq:60053b}
\end{eqnarray}
\begin{eqnarray}
y_{\tau }^m(j,q;x) &=& x^{1-\gamma } \left\{ \sum_{i_0=0}^{m} \frac{\left( i_0 -j\right) }{\left( i_0+3-\gamma  \right) \left( i_0+2\right)} \frac{ \left( \Delta_0^{-} \left( q\right) \right)_{i_0}\left( \Delta_0^{+} \left( q\right) \right)_{i_0}}{\left( 2-\gamma  \right)_{i_0} \left( 1\right)_{i_0}} \right.\nonumber\\
&&\times \prod_{k=1}^{\tau -1} \left( \sum_{i_k = i_{k-1}}^{m} \frac{\left( i_k+ 2k-j\right) }{\left( i_k+2k+3-\gamma \right) \left( i_k+2k+2\right)} \right.  \nonumber\\
&&\times \left. \frac{ \left( \Delta_k^{-} \left( q\right) \right)_{i_k}\left( \Delta_k^{+} \left( q\right) \right)_{i_k} \left( 2k+2-\gamma \right)_{i_{k-1}} \left( 2k+1\right)_{i_{k-1}}}{\left( \Delta_k^{-} \left( q\right) \right)_{i_{k-1}}\left( \Delta_k^{+} \left( q\right) \right)_{i_{k-1}}\left( 2k+2-\gamma \right)_{i_k} \left( 2k+1\right)_{i_k}} \right) \nonumber\\
&&\times \left. \sum_{i_{\tau } = i_{\tau -1}}^{m}  \frac{ \left( \Delta_{\tau }^{-} \left( q\right) \right)_{i_{\tau }}\left( \Delta_{\tau }^{+} \left( q\right) \right)_{i_{\tau }} \left( 2\tau +2-\gamma \right)_{i_{\tau -1}} \left( 2\tau +1\right)_{i_{\tau -1}}}{\left( \Delta_{\tau }^{-} \left( q\right) \right)_{i_{\tau -1}}\left( \Delta_{\tau }^{+} \left( q\right) \right)_{i_{\tau -1}}\left( 2\tau +2-\gamma \right)_{i_\tau } \left( 2\tau +1\right)_{i_{\tau }}} x^{i_{\tau }}\right\} \rho ^{\tau } \hspace{1.5cm}\label{eq:60053c} 
\end{eqnarray}
where
\begin{equation}
\begin{cases} \tau \geq 2 \cr
\Delta _{k}^{\pm}(q) = \frac{\varphi +2(2k+1-\gamma ) \pm \sqrt{\varphi ^2 +4q}}{2}\cr
\varphi = -\beta +\gamma +\delta -1
\end{cases}\nonumber
\end{equation}
\end{enumerate}
\end{remark}
\section{Summary}
In this chapter, by applying a mathematical expression of the first species complete polynomial using 3TRF; this is done by letting $A_n$ in sequences $c_n$ is the leading term in each finite sub-power series of the general series solution $y(x)$, I show formal series solutions in closed forms of the CHE for a polynomial of type 3 which makes $A_n$ and $B_n$ terms terminated at a specific index summation $n$.   
And the first species complete polynomials of the CHE around $x=0$ are constructed by applying a general summation formula of complete polynomials using R3TRF: this is done by letting $B_n$ in sequences $c_n$ is the leading term in each finite sub-power series of the general series solution $y(x)$.

It follows from the principle of the radius of convergence that $\gamma \ne 0,-1,-2,\cdots$ is required for the first kind of independent solutions of the CHE for first species complete polynomials using 3TRF and R3TRF. Similarly, it is also accompanied from the principle of analytic continuation that $\gamma \ne 2,3,4,\cdots$ is required for the second kind of independent solutions of the CHE for all cases.   

These two complete polynomial solutions are identical to each other analytically. Difference between two complete polynomials is that complete polynomial solutions of the CHE around $x=0$ by applying 3TRF contain the sum of two finite sub-power series of the general series solution of the CHE. In contrast, complete polynomial expressions of the CHE by applying R3TRF only consists of one finite sub-formal series. 

Two algebraic equations of the CHE around $x=0$ for the determination of a parameter $q$ in the form of partial sums of the sequences $\{A_n\}$ and $\{B_n\}$ using 3TRF and R3TRF are equivalent to each other analytically.
Dissimilarity for the mathematical structure is that $A_n$ is the leading term of each sequences $\bar{c}(l,n;j,q)$ where $l \in \mathbb{N}_{0}$ in a polynomial equation of $q$ for the first species complete polynomial using 3TRF. And $B_n$ is the leading term of each sequences $\bar{c}(l,n;j,q)$ in a algebraic equation of $q$ for the first species complete polynomial using R3TRF.

For the first species complete polynomials of the CHE around $x=0$ by applying 3TRF, the denominators and numerators in all $B_n$ terms of each finite sub-power series arise with Pochhammer symbols. For the first species complete polynomials of the CHE around $x=0$ by applying R3TRF, the denominators and numerators in all $A_n$ terms of each finite sub-power series arise with Pochhammer symbols.

In chapter 4 of Ref.\cite{6Choun2013}, combined definite and contour integrals for an infinite series and a type 1 polynomial of the CHE around $x=0$ are constructed by applying an integral representation of Confluent hypergeometric (Kummer's) functions into each sub-power series of the general series solution for the CHE. 
And generating functions for the CHP of type 1 are constructed by applying the generating function for Kummer's polynomials into each sub-integral of the general integral form of the CHP of type 1.

In chapter 5 of Ref.\cite{6Choun2013}, combined definite and contour integrals for an infinite series and a type 2 polynomial of the CHE around $x=0$ are constructed by applying an integral form of Gauss hypergeometric functions into each sub-power series of the general series solutions for the CHE. 
And generating functions for the CHP of type 2 are constructed by applying the generating function for Jacobi polynomial using hypergeometric functions into each sub-integral of the general integral form of the CHP of type 1.

By using similar methods, 3 analytic mathematical structures such as (1) integral forms, (2) generating functions and (3) orthogonal relations of the CHE around $x=0$ for the first species complete polynomials using 3TRF and R3TRF  will be constructed in the future series.

\addcontentsline{toc}{section}{Bibliography}
\bibliographystyle{model1a-num-names}
\bibliography{<your-bib-database>}

\chapter{Complete polynomials of Confluent Heun equation about the irregular singular point at infinity}
\chaptermark{Complete polynomials of the CHE around $x=\infty $} 


In chapter 7, I construct Frobenius solutions of the CHE about the regular singular point at zero for a polynomial of type 3 by applying mathematical formulas of complete polynomials using 3TRF and R3TRF.
In this chapter I derive power series solutions of the CHE about the irregular singular point at infinity for a polynomial of type 3 by applying general mathematical expressions of complete polynomials using 3TRF and R3TRF. 
 
\section{Introduction}
As I know, there are three types of the Confluent Heun equation (CHE): (1)Generalized spheroidal equation (GSE), (2) the generalized spheroidal wave equation (GSWE) in the Leaver version, (3) Non-symmetrical canonical form of the CHE.

The GSE was obtained by A.H. Wilson in 1928 since he looked into the wave equation for the ion of the hydrogen molecule $H_2^{+}$ with considering the protons as at rest and neglecting the mass of the electron. \cite{7Wils1928,7Wils1928a} 
The GSWE (the different version of the GSE), suggested by E.W. Leaver in 1986 \cite{7Leav1986}, arises in non-relativistic quantum scattering mechanics, Teukolsky's equations of Kerr black holes, general relativity, the gauge theories on thick brane words, hydrogen molecule ion in Stark effect, etc.  
These two differential equations have two regular singular points and one irregular singular point. 
 


The non-symmetrical canonical form of the CHE is a second-order linear ordinary differential equation of the form \cite{7Deca1978,7Decar1978,7Ronv1995,7Slavy2000}
\begin{equation}
\frac{d^2{y}}{d{x}^2} + \left(\beta  +\frac{\gamma }{x} + \frac{\delta }{x-1}\right) \frac{d{y}}{d{x}} +  \frac{\alpha \beta x-q}{x(x-1)} y = 0 \label{eq:7002} 
\end{equation}
(\ref{eq:7002}) has three singular points: two regular singular points which are 0 and 1 with exponents $\{0, 1-\gamma\}$ and $\{0, 1-\delta \}$, and one irregular singular point which is $\infty$ with an exponent $\alpha$. Its solution is denoted as $H_{c}(\alpha,\beta,\gamma,\delta,q;x)$.\footnote{Several authors denote as coefficients $4p$ and $\sigma $ instead of $\beta $ and $q$. And they define the solution of the CHE as $H_{c}^{(a)}(p,\alpha,\gamma,\delta,\sigma;x)$.} The CHE is a more general form than any other linear second ODEs such as Whittaker-Hill, spheroidal wave, Coulomb spheroidal and Mathieu's equations.

According to changing parameters, adding the new variables and combining regular singularities into the GSE and its Leaver version, we can convert to the non-symmetrical canonical form of the CHE. 
In chapters 4 \& 5 of Ref.\cite{7Choun2013}, chapter 7 and this chapter, (\ref{eq:7002}) is investigated for analytic solutions of the CHE rather than the GSE and its Leaver version because of more convenient numerical computations.  

Slavyanov \textit {et al.} showed asymptotic expansions of the CHE around $x=\infty $ with eigenvalues at large values of several parameters since he investigated phase shifts in the Coulomb two-center problem of quantum theory in the limit of large and small distances between the charges. \cite{7Abra1978,7Abra1979,7Slav1996}

In general, when we substitute a power series $y(x)= \sum_{n=0}^{\infty } c_n x^{n+\lambda }$ into any linear ODEs having a 3-term recurrence relation between consecutive coefficients in their power series expansions, the recurrence relation for frobenius solution starts to appear such as
\begin{equation}
c_{n+1}=A_n \;c_n +B_n \;c_{n-1} \hspace{1cm};n\geq 1
\label{eq:7003}
\end{equation}
where $c_1= A_0 \;c_0$ and $\lambda $ is an indicial root.

Formal series solutions for any ODEs having a 3-term recursive relation and their integral representations are hitherto unknown unfortunately because of their complex mathematical calculations.  
For example, hypergeometric differential equation provides a 2-term recursion relation between successive coefficients. Its power series solutions are classified into two types such as an infinite series and a polynomial called as `Jacobi polynomial.'
In contrast, there are $2^{3-1}$ possible formal series solutions of any ODEs having a 3-term recurrence relation in their power series dissimilarly such as an infinite series and 3 types of polynomials: (1) a polynomial which makes $B_n$ term terminated; $A_n$ term is not terminated, (2) a polynomial which makes $A_n$ term terminated; $B_n$ term is not terminated, (3) a polynomial which makes $A_n$ and $B_n$ terms terminated at the same time, referred as `a complete polynomial.'

The sequence $c_n$ consists of combinations of $A_n$ and $B_n$ terms in (\ref{eq:7003}). 
By allowing $A_n$ in the sequence $c_n$ is the leading term of each sub-power series $y_l(x)$ where $l\in \mathbb{N}_{0}$ in a function $y(x)= \sum_{n=0}^{\infty } y_n(x)$ \cite{7Chou2012}, the general summation formulas of the 3-term recurrence relation in a linear ODE are constructed for an infinite series and a polynomial of type 1,  designated as `three term recurrence formula (3TRF).' More precisely, a sub-power series $y_l(x)$ is obtained by observing the term of sequence $c_n$ which includes $l$ terms of $A_n's$.
  
By allowing $B_n$ in the sequence $c_n$ is the leading term of each sub-power series $y_l(x)$ in a function $y(x)= \sum_{n=0}^{\infty } y_n(x)$ in chapter 1 of Ref.\cite{7Choun2013}, the general summation formalism of the 3-term recurrence relation in a linear ODE are constructed for an infinite series and a polynomial of type 2,  denominated as `reversible three term recurrence formula (R3TRF)': a sub-power series $y_l(x)$ is obtained by observing the term of sequence $c_n$ which includes $l$ terms of $B_n's$.
   
Complete polynomials which make $A_n$ and $B_n$ terms terminated in the 3-term recurrence relation are divided into two different types such as (1) the first species complete polynomial where a parameter of a numerator in $B_n$ term and a (spectral) parameter of a numerator in $A_n$ term are fixed constants and (2) the second species complete polynomial where two parameters of a numerator in $B_n$ term and a parameter of a numerator in $A_n$ term are fixed constants.
With my definition, the first species complete polynomial has multi-valued roots of a parameter of a numerator in $A_n$ term, but the second species complete polynomial has only one fixed value of a parameter of a numerator in $A_n$ term for an eigenvalue. 

By allowing $A_n$ as the leading term in each finite sub-power series $y_{\tau }^m(x)$ where $\tau ,m \in \mathbb{N}_{0}$ of the general power series $y(x)$, the mathematical summation formulas of complete polynomials for the first and second species, designated as `complete polynomials using 3-term recurrence formula (3TRF),' are constructed in chapter 1: a finite sub-power series $y_{\tau }^m(x)$ is a polynomial that the lower bound of summation is zero and the upper one is $m$. It is obtained by observing the term of sequence $c_n$ which includes $\tau $ terms of $A_n's$. 
 
By allowing $B_n$ as the leading term in each finite sub-power series $y_{\tau }^m(x)$ of the general power series $y(x)$, the classical mathematical formulas in series of complete polynomials for the first and second species, designated as `complete polynomials using reversible 3-term recurrence formula (R3TRF),' are constructed in chapter 2: a finite sub-power series $y_{\tau }^m(x)$ is a polynomial which is obtained by observing the term of sequence $c_n$ which includes $\tau $ terms of $B_n's$. 

In this chapter I derive power series solutions in compact forms, called summation notation, of the CHE around  $x=\infty$ for two types of polynomials which make $A_n$ and $B_n$ terms terminated by applying complete polynomials using 3TRF and R3TRF.
 
\section{Power series about the irregular singular point at infinity}
\sectionmark{Power series about the irregular singular point at $x=\infty $} 
The general solution in series of the CHE about the singular point at infinity is not convergent by only asymptotic. \cite{7Fizi2010a,7Ronv1995} Putting $z=\frac{1}{x}$ into (\ref{eq:7002}),
\begin{equation}
z^2 \frac{d^2{y}}{d{z}^2} + z\left((2-\gamma ) -\frac{\beta }{z} + \frac{\delta }{z-1}\right) \frac{d{y}}{d{z}} +  \frac{q z-\alpha \beta }{z(z-1)} y = 0 \label{eq:70054}
\end{equation}
Assuming the formal solution of (\ref{eq:70054}) as
\begin{equation}
y(z)= \sum_{n=0}^{\infty } c_n z^{n+\lambda }  \label{eq:70055}
\end{equation}
We again obtain by substitution in (\ref{eq:70054}) the 3-term recurrence system
\begin{equation}
c_{n+1}=A_n \;c_n +B_n \;c_{n-1} \hspace{1cm};n\geq 1 \label{eq:70056}
\end{equation}
where,
\begin{subequations}
\begin{equation}
A_n =\frac{(n+\alpha )(n+\alpha +\beta -\gamma -\delta +1)-q}{\beta (n+1)}\label{eq:70057a}
\end{equation}
\begin{equation}
B_n = \frac{-(n+\alpha -1)(n+\alpha-\gamma )}{\beta (n+1)} \label{eq:70057b}
\end{equation}
\end{subequations}
where $c_1= A_0 \;c_0$. We only have an indicial root which is $\lambda =\alpha $.

\begin{table}[h]
\begin{center}
{
 \Tree[.{ Confluent Heun differential equation about the irregular singular point at infinity} [.{ 3TRF}  
               [ 
               [.{ Polynomial of type 3} [.{  $ \begin{array}{lcll}  1^{\mbox{st}}\;  \mbox{species}\\ \mbox{complete} \\ \mbox{polynomial} \end{array}$} ] [.{ $ \begin{array}{lcll}  2^{\mbox{nd}}\;  \mbox{species}\\ \mbox{complete} \\ \mbox{polynomial} \end{array}$} ]]]]                         
  [.{ R3TRF} 
     [ 
       [.{  Polynomial of type 3} [.{  $ \begin{array}{lcll}  1^{\mbox{st}}\;  \mbox{species} \\ \mbox{complete} \\ \mbox{polynomial} \end{array}$} ] [.{  $ \begin{array}{lcll}  2^{\mbox{nd}}\;  \mbox{species} \\ \mbox{complete} \\ \mbox{polynomial} \end{array}$} ]]]]] 
}
\end{center}
\caption{Power series of Confluent Heun equation about the irregular singular point at infinity}
\end{table}

Table 8.1 informs us about all possible general solutions in series of the CHE about the irregular singular point at infinity.
There are no general series solutions for an infinite series, the type 1 and type 2 polynomials because since the index of summation $n$ goes to infinity in $A_n$ and $B_n$ terms in (\ref{eq:70057a}) and (\ref{eq:70057b}), series solutions will be divergent. 
The Frobenius solutions of the CHE about the singular point at infinity are only valid for complete polynomials where $A_n$ and $B_n$ terms terminated. 

In this chapter, for the first species complete polynomial of the CHE around $x=\infty $, I treat $\beta $, $\gamma $, $\delta $ as free variables and $\alpha $, $q$ as fixed values. And for the second species complete polynomial of it, I treat $\beta $, $\delta $ as free variables and $\alpha $, $\gamma $, $q$ as fixed values.
\subsection{The first species complete polynomial of Confluent Heun equation using 3TRF}
For the first species complete polynomials using 3TRF and R3TRF, we need a condition which is given by
\begin{equation}
B_{j+1}= c_{j+1}=0\hspace{1cm}\mathrm{where}\;j\in \mathbb{N}_{0}  
 \label{eq:7006}
\end{equation}
(\ref{eq:7006}) gives successively $c_{j+2}=c_{j+3}=c_{j+4}=\cdots=0$. And $c_{j+1}=0$ is defined by a polynomial equation of degree $j+1$ for the determination of an accessory parameter in $A_n$ term. 
\begin{theorem}
In chapter 1, the general expression of a function $y(x)$ for the first species complete polynomial using 3-term recurrence formula and its algebraic equation for the determination of an accessory parameter in $A_n$ term are given by
\begin{enumerate} 
\item As $B_1=0$,
\begin{equation}
0 =\bar{c}(1,0) \label{eq:7007a}
\end{equation}
\begin{equation}
y(x) = y_{0}^{0}(x) \label{eq:7007b}
\end{equation}
\item As $B_{2N+2}=0$ where $N \in \mathbb{N}_{0}$,
\begin{equation}
0  = \sum_{r=0}^{N+1}\bar{c}\left( 2r, N+1-r\right) \label{eq:7008a}
\end{equation}
\begin{equation}
y(x)= \sum_{r=0}^{N} y_{2r}^{N-r}(x)+ \sum_{r=0}^{N} y_{2r+1}^{N-r}(x)  \label{eq:7008b}
\end{equation}
\item As $B_{2N+3}=0$ where $N \in \mathbb{N}_{0}$,
\begin{equation}
0  = \sum_{r=0}^{N+1}\bar{c}\left( 2r+1, N+1-r\right) \label{eq:7009a}
\end{equation}
\begin{equation}
y(x)= \sum_{r=0}^{N+1} y_{2r}^{N+1-r}(x)+ \sum_{r=0}^{N} y_{2r+1}^{N-r}(x)  \label{eq:7009b}
\end{equation}
In the above,
\begin{eqnarray}
\bar{c}(0,n)  &=& \prod _{i_0=0}^{n-1}B_{2i_0+1} \label{eq:70010a}\\
\bar{c}(1,n) &=&  \sum_{i_0=0}^{n} \left\{ A_{2i_0} \prod _{i_1=0}^{i_0-1}B_{2i_1+1} \prod _{i_2=i_0}^{n-1}B_{2i_2+2} \right\} 
\label{eq:70010b}\\
\bar{c}(\tau ,n) &=& \sum_{i_0=0}^{n} \left\{A_{2i_0}\prod _{i_1=0}^{i_0-1} B_{2i_1+1} 
\prod _{k=1}^{\tau -1} \left( \sum_{i_{2k}= i_{2(k-1)}}^{n} A_{2i_{2k}+k}\prod _{i_{2k+1}=i_{2(k-1)}}^{i_{2k}-1}B_{2i_{2k+1}+(k+1)}\right) \right. \nonumber\\
&&\times \left. \prod _{i_{2\tau}=i_{2(\tau -1)}}^{n-1} B_{2i_{2\tau }+(\tau +1)} \right\} 
\hspace{1cm}\label{eq:70010c} 
\end{eqnarray}
and
\begin{eqnarray}
y_0^m(x) &=& c_0 x^{\lambda } \sum_{i_0=0}^{m} \left\{ \prod _{i_1=0}^{i_0-1}B_{2i_1+1} \right\} x^{2i_0 } \label{eq:70011a}\\
y_1^m(x) &=& c_0 x^{\lambda } \sum_{i_0=0}^{m}\left\{ A_{2i_0} \prod _{i_1=0}^{i_0-1}B_{2i_1+1}  \sum_{i_2=i_0}^{m} \left\{ \prod _{i_3=i_0}^{i_2-1}B_{2i_3+2} \right\}\right\} x^{2i_2+1 } \label{eq:70011b}\\
y_{\tau }^m(x) &=& c_0 x^{\lambda } \sum_{i_0=0}^{m} \left\{A_{2i_0}\prod _{i_1=0}^{i_0-1} B_{2i_1+1} 
\prod _{k=1}^{\tau -1} \left( \sum_{i_{2k}= i_{2(k-1)}}^{m} A_{2i_{2k}+k}\prod _{i_{2k+1}=i_{2(k-1)}}^{i_{2k}-1}B_{2i_{2k+1}+(k+1)}\right) \right. \nonumber\\
&& \times  \left. \sum_{i_{2\tau} = i_{2(\tau -1)}}^{m} \left( \prod _{i_{2\tau +1}=i_{2(\tau -1)}}^{i_{2\tau}-1} B_{2i_{2\tau +1}+(\tau +1)} \right) \right\} x^{2i_{2\tau}+\tau }\hspace{1cm}\mathrm{where}\;\tau \geq 2
\label{eq:70011c} 
\end{eqnarray}
\end{enumerate}
\end{theorem} 
Put $n= j+1$ in (\ref{eq:70057b}) and use the condition $B_{j+1}=0$ for $\alpha $. We obtain two possible fixed values for $\alpha $ such as
\begin{subequations}
\begin{equation}
\alpha = -j \label{eq:70058a}
\end{equation}
\begin{equation}
\alpha = \gamma -j-1  \label{eq:70058b}
\end{equation}
\end{subequations}
\subsubsection{The case of $\alpha = -j $ }
Take (\ref{eq:70058a}) into (\ref{eq:70057a}) and (\ref{eq:70057b}).
\begin{subequations}
\begin{equation}
A_n =\frac{(n-j)(n-j+\beta -\gamma -\delta +1)-q}{\beta (n+1)}\label{eq:70059a}
\end{equation}
\begin{equation}
B_n = -\frac{ (n-j-1)(n-j-\gamma )}{\beta (n+1)} \label{eq:70059b}
\end{equation}
\end{subequations}
According to (\ref{eq:7006}), $c_{j+1}=0$ is clearly an algebraic equation in $q$ of degree $j+1$ and thus has $j+1$ zeros denoted them by $q_j^m$ eigenvalues where $m = 0,1,2, \cdots, j$. They can be arranged in the following order: $q_j^0 < q_j^1 < q_j^2 < \cdots < q_j^j$.

In (\ref{eq:7007b}), (\ref{eq:7008b}), (\ref{eq:7009b}) and (\ref{eq:70011a})--(\ref{eq:70011c}) replace $c_0$, $\lambda $ and $x$ by $1$, $\alpha$ and $z$. Substitute (\ref{eq:70059a}) and (\ref{eq:70059b}) into (\ref{eq:70010a})--(\ref{eq:70011c}).

As $B_{1}= c_{1}=0$, take the new (\ref{eq:70010b}) into (\ref{eq:7007a}) putting $j=0$. Substitute the new (\ref{eq:70011a}) with $\alpha =0$ into (\ref{eq:7007b}) putting $j=0$. 

As $B_{2N+2}= c_{2N+2}=0$, take the new (\ref{eq:70010a})--(\ref{eq:70010c}) into (\ref{eq:7008a}) putting $j=2N+1$. Substitute the new 
(\ref{eq:70011a})--(\ref{eq:70011c}) with $\alpha =-2N-1$ into (\ref{eq:7008b}) putting $j=2N+1$ and $q=q_{2N+1}^m$.

As $B_{2N+3}= c_{2N+3}=0$, take the new (\ref{eq:70010a})--(\ref{eq:70010c}) into (\ref{eq:7009a}) putting $j=2N+2$. Substitute the new 
(\ref{eq:70011a})--(\ref{eq:70011c}) with $\alpha =-2N-2$ into (\ref{eq:7009b}) putting $j=2N+2$ and $q=q_{2N+2}^m$.

After the replacement process, we obtain the independent solution of the CHE about $x=\infty $. The solution is as follows.
\begin{remark}
The power series expansion of Confluent Heun equation of the first kind for the first species complete polynomial using 3TRF about $x=\infty $ for $\alpha =$ non-positive integer is given by
\begin{enumerate} 
\item As $\alpha =0$ and $q=q_0^0=0$,

The eigenfunction is given by
\begin{equation}
y(z) =\; _pH_c^{(r,1)}F_{0,0} \left( \alpha =0, \beta, \gamma, \delta, q=q_0^0=0; z=\frac{1}{x}, \mu =-\frac{2}{\beta } z^2, \xi =\frac{2}{\beta } z\right) =1 \label{eq:70060}
\end{equation}
\item As $\alpha =-2N-1$ where $N \in \mathbb{N}_{0}$,

An algebraic equation of degree $2N+2$ for the determination of $q$ is given by
\begin{equation}
0 = \sum_{r=0}^{N+1}\bar{c}\left( 2r, N+1-r; 2N+1,q\right)\label{eq:70061a}
\end{equation}
The eigenvalue of $q$ is written by $q_{2N+1}^m$ where $m = 0,1,2,\cdots,2N+1 $; $q_{2N+1}^0 < q_{2N+1}^1 < \cdots < q_{2N+1}^{2N+1}$. Its eigenfunction is given by
\begin{eqnarray} 
y(z) &=& _pH_c^{(r,1)}F_{2N+1,m} \left( \alpha =-2N-1, \beta, \gamma, \delta, q=q_{2N+1}^m; z=\frac{1}{x}, \mu =-\frac{2}{\beta }z^2, \xi =\frac{2}{\beta } z \right)\nonumber\\
&=& z^{-2N-1} \left\{ \sum_{r=0}^{N} y_{2r}^{N-r}\left( 2N+1,q_{2N+1}^m;z\right)+ \sum_{r=0}^{N} y_{2r+1}^{N-r}\left( 2N+1,q_{2N+1}^m;z\right) \right\}  
\hspace{1.5cm}\label{eq:70061b}
\end{eqnarray}
\item As $\alpha =-2N-2$ where $N \in \mathbb{N}_{0}$,

An algebraic equation of degree $2N+3$ for the determination of $q$ is given by
\begin{eqnarray}
0  = \sum_{r=0}^{N+1}\bar{c}\left( 2r+1, N+1-r; 2N+2,q\right)\label{eq:70062a}
\end{eqnarray}
The eigenvalue of $q$ is written by $q_{2N+2}^m$ where $m = 0,1,2,\cdots,2N+2 $; $q_{2N+2}^0 < q_{2N+2}^1 < \cdots < q_{2N+2}^{2N+2}$. Its eigenfunction is given by
\begin{eqnarray} 
y(z) &=& _pH_c^{(r,1)}F_{2N+2,m} \left( \alpha =-2N-2, \beta, \gamma, \delta, q=q_{2N+2}^m; z=\frac{1}{x}, \mu =-\frac{2}{\beta } z^2, \xi =\frac{2}{\beta } z \right)\nonumber\\
&=& z^{-2N-2} \left\{ \sum_{r=0}^{N+1} y_{2r}^{N+1-r}\left( 2N+2,q_{2N+2}^m;z\right) + \sum_{r=0}^{N} y_{2r+1}^{N-r}\left( 2N+2,q_{2N+2}^m;z\right) \right\} 
\hspace{1.5cm}\label{eq:70062b}
\end{eqnarray}
In the above,
\begin{eqnarray}
\bar{c}(0,n;j,q)  &=& \frac{\left( -\frac{j}{2}\right)_{n}\left(  \frac{1}{2}-\frac{j}{2}-\frac{\gamma }{2} \right)_{n}}{\left( 1 \right)_{n}} \left( -\frac{2}{\beta } \right)^{n}\label{eq:70063a}\\
\bar{c}(1,n;j,q) &=& \left( \frac{2}{\beta } \right) \sum_{i_0=0}^{n}\frac{ \left( i_0 -\frac{j}{2}\right) \left( i_0 +\Pi _0^j \right) -\frac{q}{4}}{\left( i_0+\frac{1}{2} \right)} \frac{\left( -\frac{j}{2}\right)_{i_0} \left(  \frac{1}{2}-\frac{j}{2}-\frac{\gamma }{2} \right)_{i_0}}{\left( 1 \right)_{i_0}} \nonumber\\
&&\times  \frac{\left( \frac{1}{2}-\frac{j}{2} \right)_{n} \left(  1-\frac{j}{2}-\frac{\gamma }{2} \right)_n \left( \frac{3}{2} \right)_{i_0}}{\left( \frac{1}{2}-\frac{j}{2}\right)_{i_0}\left(  1-\frac{j}{2}-\frac{\gamma }{2} \right)_{i_0} \left( \frac{3}{2} \right)_{n}} \left( -\frac{2}{\beta } \right)^{n }  
\label{eq:70063b}\\
\bar{c}(\tau ,n;j,q) &=& \left( \frac{2}{\beta } \right)^{\tau } \sum_{i_0=0}^{n}\frac{ \left( i_0 -\frac{j}{2}\right) \left( i_0 +\Pi _0^j \right) -\frac{q}{4}}{\left( i_0+\frac{1}{2} \right)} \frac{\left( -\frac{j}{2}\right)_{i_0} \left(  \frac{1}{2}-\frac{j}{2}-\frac{\gamma }{2} \right)_{i_0}}{\left( 1 \right)_{i_0}}  \nonumber\\
&&\times  \prod_{k=1}^{\tau -1} \left( \sum_{i_k = i_{k-1}}^{n} \frac{\left( i_k+ \frac{k}{2}-\frac{j}{2} \right)\left( i_k +\Pi _k^j \right) -\frac{q}{4}}{\left( i_k+\frac{k}{2}+\frac{1}{2} \right)} \right.   \left. \frac{\left( \frac{k}{2}-\frac{j}{2}\right)_{i_k} \left( \frac{k}{2}+ \frac{1}{2}-\frac{j}{2}-\frac{\gamma }{2} \right)_{i_{k}}\left( \frac{k}{2}+1 \right)_{i_{k-1}} }{\left( \frac{k}{2}-\frac{j}{2}\right)_{i_{k-1}} \left( \frac{k}{2}+ \frac{1}{2}-\frac{j}{2}-\frac{\gamma }{2} \right)_{i_{k-1}}\left( \frac{k}{2}+1 \right)_{i_k}} \right) \nonumber\\ 
&&\times \frac{\left( \frac{\tau }{2} -\frac{j}{2}\right)_{n} \left( \frac{\tau }{2}+\frac{1}{2}-\frac{j}{2}-\frac{\gamma }{2} \right)_{n}\left( \frac{\tau }{2}+1 \right)_{i_{\tau -1}} }{\left( \frac{\tau }{2}-\frac{j}{2}\right)_{i_{\tau -1}} \left( \frac{\tau }{2}+\frac{1}{2}-\frac{j}{2}-\frac{\gamma }{2} \right)_{i_{\tau -1}}\left( \frac{\tau }{2}+1 \right)_{n} } \left( -\frac{\beta }{2} \right)^{n } \label{eq:70063c} 
\end{eqnarray}
\begin{eqnarray}
y_0^m(j,q;z) &=& \sum_{i_0=0}^{m} \frac{\left( -\frac{j}{2}\right)_{i_0}\left(  \frac{1}{2}-\frac{j}{2}-\frac{\gamma }{2} \right)_{i_0}}{\left( 1 \right)_{i_0}} \mu ^{i_0} \label{eq:70064a}\\
y_1^m(j,q;z) &=& \left\{\sum_{i_0=0}^{m} \frac{ \left( i_0 -\frac{j}{2}\right) \left( i_0 +\Pi _0^j \right) -\frac{q}{4}}{\left( i_0+\frac{1}{2} \right)} \frac{\left( -\frac{j}{2}\right)_{i_0} \left(  \frac{1}{2}-\frac{j}{2}-\frac{\gamma }{2} \right)_{i_0}}{\left( 1 \right)_{i_0}} \right. \nonumber\\
&&\times \left. \sum_{i_1 = i_0}^{m} \frac{\left( \frac{1}{2}-\frac{j}{2} \right)_{i_1} \left( 1-\frac{j}{2}-\frac{\gamma }{2} \right)_{i_1}\left( \frac{3}{2} \right)_{i_0} }{\left( \frac{1}{2}-\frac{j}{2} \right)_{i_0} \left( 1-\frac{j}{2}-\frac{\gamma }{2} \right)_{i_0}\left( \frac{3}{2} \right)_{i_1}} \mu ^{i_1}\right\}\xi 
\hspace{1.5cm}\label{eq:70064b}\\
y_{\tau }^m(j,q;z) &=& \left\{ \sum_{i_0=0}^{m} \frac{ \left( i_0 -\frac{j}{2}\right) \left( i_0 +\Pi _0^j \right) -\frac{q}{4}}{\left( i_0+\frac{1}{2} \right)} \frac{\left( -\frac{j}{2}\right)_{i_0} \left(  \frac{1}{2}-\frac{j}{2}-\frac{\gamma }{2} \right)_{i_0}}{\left( 1 \right)_{i_0}} \right.\nonumber\\
&&\times \prod_{k=1}^{\tau -1} \left( \sum_{i_k = i_{k-1}}^{m} \frac{\left( i_k+ \frac{k}{2} -\frac{j}{2}\right)\left( i_k +\Pi _k^j  \right) -\frac{q}{4}}{\left( i_k+\frac{k}{2}+\frac{1}{2} \right)} \right. \nonumber\\
&&\times  \left. \frac{\left( \frac{k}{2}-\frac{j}{2}\right)_{i_k} \left( \frac{k}{2}+ \frac{1}{2}-\frac{j}{2}-\frac{\gamma }{2} \right)_{i_{k}}\left( \frac{k}{2}+1 \right)_{i_{k-1}} }{\left( \frac{k}{2}-\frac{j}{2}\right)_{i_{k-1}} \left( \frac{k}{2}+ \frac{1}{2}-\frac{j}{2}-\frac{\gamma }{2} \right)_{i_{k-1}}\left( \frac{k}{2}+1 \right)_{i_k}} \right) \nonumber\\
&&\times  \left. \sum_{i_{\tau } = i_{\tau -1}}^{m} \frac{\left( \frac{\tau }{2} -\frac{j}{2}\right)_{i_{\tau}} \left( \frac{\tau }{2}+\frac{1}{2}-\frac{j}{2}-\frac{\gamma }{2} \right)_{i_{\tau}}\left( \frac{\tau }{2}+1 \right)_{i_{\tau -1}} }{\left( \frac{\tau }{2}-\frac{j}{2}\right)_{i_{\tau -1}} \left( \frac{\tau }{2}+\frac{1}{2}-\frac{j}{2}-\frac{\gamma }{2} \right)_{i_{\tau -1}}\left( \frac{\tau }{2}+1 \right)_{i_{\tau }} } \mu ^{i_{\tau }}\right\} \xi ^{\tau } \hspace{1.5cm}\label{eq:70064c} 
\end{eqnarray}
where
\begin{equation}
\begin{cases} \tau \geq 2 \cr 
\Pi _0^j = \frac{1}{2}\left( \beta -\gamma -\delta +1-j \right) \cr
\Pi _k^j = \frac{1}{2}\left( \beta -\gamma -\delta +1-j+k \right)
\end{cases}\nonumber
\end{equation}
\end{enumerate}
\end{remark}
\subsubsection{The case of $\alpha = \gamma -j-1 $}
Take (\ref{eq:70058b}) into (\ref{eq:70057a}) and (\ref{eq:70057b}).
\begin{subequations}
\begin{equation}
A_n =\frac{(n-j-1+\gamma )(n-j+\beta -\delta )-q}{\beta (n+1)}\label{eq:70065a}
\end{equation}
\begin{equation}
B_n = -\frac{ (n-j-1)(n-j-2+\gamma )}{\beta (n+1)} \label{eq:70065b}
\end{equation}
\end{subequations}
According to (\ref{eq:7006}), $c_{j+1}=0$ is clearly an algebraic equation in $q$ of degree $j+1$ and thus has $j+1$ zeros denoted them by $q_j^m$ eigenvalues where $m = 0,1,2, \cdots, j$. They can be arranged in the following order: $q_j^0 < q_j^1 < q_j^2 < \cdots < q_j^j$.
 
In (\ref{eq:7007b}), (\ref{eq:7008b}), (\ref{eq:7009b}) and (\ref{eq:70011a})--(\ref{eq:70011c}) replace $c_0$, $\lambda $ and $x$ by $1$, $\alpha$ and $z$. Substitute (\ref{eq:70065a}) and (\ref{eq:70065b}) into (\ref{eq:70010a})--(\ref{eq:70011c}).

As $B_{1}= c_{1}=0$, take the new (\ref{eq:70010b}) into (\ref{eq:7007a}) putting $j=0$. Substitute the new (\ref{eq:70011a})  with $\alpha =\gamma -1$ into (\ref{eq:7007b}) putting $j=0$. 

As $B_{2N+2}= c_{2N+2}=0$, take the new (\ref{eq:70010a})--(\ref{eq:70010c}) into (\ref{eq:7008a}) putting $j=2N+1$. Substitute the new 
(\ref{eq:70011a})--(\ref{eq:70011c}) with $\alpha =\gamma -2N-2$ into (\ref{eq:7008b}) putting $j=2N+1$ and $q=q_{2N+1}^m$.

As $B_{2N+3}= c_{2N+3}=0$, take the new (\ref{eq:70010a})--(\ref{eq:70010c}) into (\ref{eq:7009a}) putting $j=2N+2$. Substitute the new 
(\ref{eq:70011a})--(\ref{eq:70011c}) with $\alpha =\gamma -2N-3$ into (\ref{eq:7009b}) putting $j=2N+2$ and $q=q_{2N+2}^m$.

After the replacement process, we obtain the independent solution of the CHE about $x=\infty $. The solution is as follows.
\begin{remark}
The power series expansion of Confluent Heun equation of the first kind for the first species complete polynomial using 3TRF about $x=\infty $ for $\alpha =\gamma -1-j$ where $j\in \mathbb{N}_{0}$ is given by
\begin{enumerate} 
\item As $\alpha =\gamma -1$ and $q=q_0^0=(\gamma -1)(\beta -\delta )$,

The eigenfunction is given by
\begin{eqnarray}
y(z) &=& _pH_c^{(r,2)}F_{0,0} \left( \alpha =\gamma -1, \beta, \gamma, \delta, q=q_0^0=(\gamma -1)(\beta -\delta ); z=\frac{1}{x}, \mu =-\frac{2}{\beta } z^2, \xi =\frac{2}{\beta } z\right) \nonumber\\
&=& z^{\gamma -1} \label{eq:70066}
\end{eqnarray}
\item As $\alpha =\gamma -2N-2$ where $N \in \mathbb{N}_{0}$,

An algebraic equation of degree $2N+2$ for the determination of $q$ is given by
\begin{equation}
0 = \sum_{r=0}^{N+1}\bar{c}\left( 2r, N+1-r; 2N+1,q\right)\label{eq:70067a}
\end{equation}
The eigenvalue of $q$ is written by $q_{2N+1}^m$ where $m = 0,1,2,\cdots,2N+1 $; $q_{2N+1}^0 < q_{2N+1}^1 < \cdots < q_{2N+1}^{2N+1}$. Its eigenfunction is given by
\begin{eqnarray} 
y(z) &=& _pH_c^{(r,2)}F_{2N+1,m} \left( \alpha =\gamma -2N-2, \beta, \gamma, \delta, q=q_{2N+1}^m; z=\frac{1}{x}, \mu =-\frac{2}{\beta }z^2, \xi =\frac{2}{\beta } z \right)\nonumber\\
&=& z^{\gamma -2N-2} \left\{ \sum_{r=0}^{N} y_{2r}^{N-r}\left( 2N+1,q_{2N+1}^m;z\right)+ \sum_{r=0}^{N} y_{2r+1}^{N-r}\left( 2N+1,q_{2N+1}^m;z\right) \right\}   
\hspace{1.5cm}\label{eq:70067b}
\end{eqnarray}
\item As $\alpha =\gamma -2N-3$ where $N \in \mathbb{N}_{0}$,

An algebraic equation of degree $2N+3$ for the determination of $q$ is given by
\begin{eqnarray}
0  = \sum_{r=0}^{N+1}\bar{c}\left( 2r+1, N+1-r; 2N+2,q\right)\label{eq:70068a}
\end{eqnarray}
The eigenvalue of $q$ is written by $q_{2N+2}^m$ where $m = 0,1,2,\cdots,2N+2 $; $q_{2N+2}^0 < q_{2N+2}^1 < \cdots < q_{2N+2}^{2N+2}$. Its eigenfunction is given by
\begin{eqnarray} 
y(z) &=& _pH_c^{(r,2)}F_{2N+2,m} \left( \alpha =\gamma -2N-3, \beta, \gamma, \delta, q=q_{2N+2}^m; z=\frac{1}{x}, \mu =-\frac{2}{\beta } z^2, \xi =\frac{2}{\beta } z \right)\nonumber\\
&=& z^{\gamma -2N-3} \left\{ \sum_{r=0}^{N+1} y_{2r}^{N+1-r}\left( 2N+2,q_{2N+2}^m;z\right) + \sum_{r=0}^{N} y_{2r+1}^{N-r}\left( 2N+2,q_{2N+2}^m;z\right) \right\} 
\hspace{1.5cm}\label{eq:70068b}
\end{eqnarray}
In the above,
\begin{eqnarray}
\bar{c}(0,n;j,q)  &=& \frac{\left( -\frac{j}{2}\right)_{n}\left( \aleph^j \right)_{n}}{\left( 1 \right)_{n}} \left( -\frac{2}{\beta } \right)^{n}\label{eq:70069a}\\
\bar{c}(1,n;j,q) &=& \left( \frac{2}{\beta } \right) \sum_{i_0=0}^{n}\frac{ \left( i_0 + \aleph^j \right) \left( i_0 +\Upsilon^j \right) -\frac{q}{4}}{\left( i_0+\frac{1}{2} \right)} \nonumber\\
&&\times \frac{\left( -\frac{j}{2}\right)_{i_0} \left( \aleph^j \right)_{i_0}}{\left( 1 \right)_{i_0}}   \frac{\left( \frac{1}{2}-\frac{j}{2} \right)_{n} \left( \frac{1}{2} + \aleph ^j \right)_n \left( \frac{3}{2} \right)_{i_0}}{\left( \frac{1}{2}-\frac{j}{2}\right)_{i_0}\left( \frac{1}{2} + \aleph^j \right)_{i_0} \left( \frac{3}{2} \right)_{n}} \left( -\frac{2}{\beta } \right)^{n }  
\hspace{1.5cm}\label{eq:70069b}
\end{eqnarray}
\begin{eqnarray}
\bar{c}(\tau ,n;j,q) &=& \left( \frac{2}{\beta } \right)^{\tau } \sum_{i_0=0}^{n}\frac{ \left( i_0 +\aleph^j \right) \left( i_0 +\Upsilon^j \right) -\frac{q}{4}}{\left( i_0+\frac{1}{2} \right)} \frac{\left( -\frac{j}{2}\right)_{i_0} \left( \aleph^j \right)_{i_0}}{\left( 1 \right)_{i_0}}  \nonumber\\
&&\times  \prod_{k=1}^{\tau -1} \left( \sum_{i_k = i_{k-1}}^{n} \frac{\left( i_k+ \frac{k}{2}+\aleph^j \right)\left( i_k + \frac{k}{2}+\Upsilon ^j \right) -\frac{q}{4}}{\left( i_k+\frac{k}{2}+\frac{1}{2} \right)} \right. \nonumber\\
&&\times  \left. \frac{\left( \frac{k}{2}-\frac{j}{2}\right)_{i_k} \left( \frac{k}{2} +\aleph^j \right)_{i_{k}}\left( \frac{k}{2}+1 \right)_{i_{k-1}} }{\left( \frac{k}{2}-\frac{j}{2}\right)_{i_{k-1}} \left( \frac{k}{2}+\aleph^j \right)_{i_{k-1}}\left( \frac{k}{2}+1 \right)_{i_k}} \right) \nonumber\\ 
&&\times \frac{\left( \frac{\tau }{2} -\frac{j}{2}\right)_{n} \left( \frac{\tau }{2}+\aleph^j \right)_{n}\left( \frac{\tau }{2}+1 \right)_{i_{\tau -1}} }{\left( \frac{\tau }{2}-\frac{j}{2}\right)_{i_{\tau -1}} \left( \frac{\tau }{2}+\aleph^j \right)_{i_{\tau -1}}\left( \frac{\tau }{2}+1 \right)_{n} } \left( -\frac{\beta }{2} \right)^{n } \label{eq:70069c} 
\end{eqnarray}
\begin{eqnarray}
y_0^m(j,q;z) &=&  \sum_{i_0=0}^{m} \frac{\left( -\frac{j}{2}\right)_{i_0}\left( \aleph^j \right)_{i_0}}{\left( 1 \right)_{i_0}} \mu ^{i_0} \label{eq:70070a}\\
y_1^m(j,q;z) &=& \left\{\sum_{i_0=0}^{m}\frac{ \left( i_0 +\aleph^j \right) \left( i_0 +\Upsilon^j \right) -\frac{q}{4}}{\left( i_0+\frac{1}{2} \right)} \frac{\left( -\frac{j}{2}\right)_{i_0} \left( \aleph^j \right)_{i_0}}{\left( 1 \right)_{i_0}} \right. \nonumber\\
&&\times \left. \sum_{i_1 = i_0}^{m} \frac{\left( \frac{1}{2}-\frac{j}{2} \right)_{i_1} \left(  \frac{1}{2}+\aleph^j \right)_{i_1}\left( \frac{3}{2} \right)_{i_0} }{\left( \frac{1}{2}-\frac{j}{2} \right)_{i_0} \left( \frac{1}{2}+\aleph^j \right)_{i_0}\left( \frac{3}{2} \right)_{i_1}} \mu ^{i_1}\right\}\xi 
\hspace{1.5cm}\label{eq:70070b}\\
y_{\tau }^m(j,q;z) &=& \left\{ \sum_{i_0=0}^{m} \frac{ \left( i_0 +\aleph^j \right) \left( i_0 +\Upsilon^j \right) -\frac{q}{4}}{\left( i_0+\frac{1}{2} \right)} \frac{\left( -\frac{j}{2}\right)_{i_0} \left( \aleph^j \right)_{i_0}}{\left( 1 \right)_{i_0}} \right.\nonumber\\
&&\times \prod_{k=1}^{\tau -1} \left( \sum_{i_k = i_{k-1}}^{m} \frac{\left( i_k+ \frac{k}{2}+\aleph^j \right)\left( i_k +\frac{k}{2}+\Upsilon^j \right) -\frac{q}{4}}{\left( i_k+\frac{k}{2}+\frac{1}{2} \right)} \right. \nonumber\\
&&\times  \left. \frac{\left( \frac{k}{2}-\frac{j}{2}\right)_{i_k} \left( \frac{k}{2} +\aleph^j \right)_{i_{k}}\left( \frac{k}{2}+1 \right)_{i_{k-1}} }{\left( \frac{k}{2}-\frac{j}{2}\right)_{i_{k-1}} \left( \frac{k}{2}+\aleph^j \right)_{i_{k-1}}\left( \frac{k}{2}+1 \right)_{i_k}} \right) \nonumber\\
&&\times  \left. \sum_{i_{\tau } = i_{\tau -1}}^{m} \frac{\left( \frac{\tau }{2} -\frac{j}{2}\right)_{i_{\tau}} \left( \frac{\tau }{2}+\aleph^j \right)_{i_{\tau}}\left( \frac{\tau }{2}+1 \right)_{i_{\tau -1}} }{\left( \frac{\tau }{2}-\frac{j}{2}\right)_{i_{\tau -1}} \left( \frac{\tau }{2}+\aleph^j \right)_{i_{\tau -1}}\left( \frac{\tau }{2}+1 \right)_{i_{\tau }} } \mu ^{i_{\tau }}\right\} \xi ^{\tau } \label{eq:70070c} 
\end{eqnarray}
where
\begin{equation}
\begin{cases} \tau \geq 2 \cr 
\aleph^j = -\frac{1}{2}-\frac{j}{2}+\frac{\gamma }{2} \cr
\Upsilon^j = -\frac{j}{2}+\frac{\beta }{2}-\frac{\delta }{2}
\end{cases}\nonumber
\end{equation}
\end{enumerate}
\end{remark}
\subsection{The second species complete polynomial of Confluent Heun equation using 3TRF}

For the second species complete polynomials using 3TRF and R3TRF, we need a condition which is defined by
\begin{equation}
B_{j}=B_{j+1}= A_{j}=0\hspace{1cm}\mathrm{where}\;j \in \mathbb{N}_{0}    
 \label{eq:70071}
\end{equation}
\begin{theorem}
In chapter 1, the general expression of a function $y(x)$ for the second species complete polynomial using 3-term recurrence formula is given by
\begin{enumerate} 
\item As $B_1=A_0=0$,
\begin{equation}
y(x) = y_{0}^{0}(x) \label{eq:70072a}
\end{equation}
\item As $B_{2N+1}=B_{2N+2}=A_{2N+1}=0$  where $N \in \mathbb{N}_{0}$,
\begin{equation}
y(x)= \sum_{r=0}^{N} y_{2r}^{N-r}(x)+ \sum_{r=0}^{N} y_{2r+1}^{N-r}(x)  \label{eq:70072b}
\end{equation}
\item As $B_{2N+2}=B_{2N+3}=A_{2N+2}=0$  where $N \in \mathbb{N}_{0}$,
\begin{equation}
y(x)= \sum_{r=0}^{N+1} y_{2r}^{N+1-r}(x)+ \sum_{r=0}^{N} y_{2r+1}^{N-r}(x)  \label{eq:70072c}
\end{equation}
In the above,
\begin{eqnarray}
y_0^m(x) &=& c_0 x^{\lambda } \sum_{i_0=0}^{m} \left\{ \prod _{i_1=0}^{i_0-1}B_{2i_1+1} \right\} x^{2i_0 } \label{eq:70073a}\\
y_1^m(x) &=& c_0 x^{\lambda } \sum_{i_0=0}^{m}\left\{ A_{2i_0} \prod _{i_1=0}^{i_0-1}B_{2i_1+1}  \sum_{i_2=i_0}^{m} \left\{ \prod _{i_3=i_0}^{i_2-1}B_{2i_3+2} \right\}\right\} x^{2i_2+1 } \label{eq:70073b}\\
y_{\tau }^m(x) &=& c_0 x^{\lambda } \sum_{i_0=0}^{m} \left\{A_{2i_0}\prod _{i_1=0}^{i_0-1} B_{2i_1+1} 
\prod _{k=1}^{\tau -1} \left( \sum_{i_{2k}= i_{2(k-1)}}^{m} A_{2i_{2k}+k}\prod _{i_{2k+1}=i_{2(k-1)}}^{i_{2k}-1}B_{2i_{2k+1}+(k+1)}\right) \right. \nonumber\\
&& \times \left. \sum_{i_{2\tau} = i_{2(\tau -1)}}^{m} \left( \prod _{i_{2\tau +1}=i_{2(\tau -1)}}^{i_{2\tau}-1} B_{2i_{2\tau +1}+(\tau +1)} \right) \right\} x^{2i_{2\tau}+\tau }\hspace{1cm}\mathrm{where}\;\tau \geq 2
\label{eq:70073c}  
\end{eqnarray}
\end{enumerate}
\end{theorem}
Put $n= j+1$ in (\ref{eq:70057b}) and use the condition $B_{j+1}=0$ for $\alpha $. We obtain two possible fixed values for $\alpha $ such as
\begin{subequations}
\begin{equation}
\alpha = -j \label{eq:70074a}
\end{equation}
\begin{equation}
\alpha = \gamma -j-1  \label{eq:70074b}
\end{equation}
\end{subequations}  
Put $n= j$ with (\ref{eq:70074a}) in (\ref{eq:70057b}) and use the condition $B_{j}=0$ for $\gamma $.  
\begin{equation}
\gamma = 0
\label{eq:70075}
\end{equation}
Substitute (\ref{eq:70074a}) and (\ref{eq:70075}) into (\ref{eq:70057a}). Put $n= j$ in the new (\ref{eq:70057a}) and use the condition $A_{j}=0$ for $q$.  
\begin{equation}
q = 0
\label{eq:70076}
\end{equation}
Take (\ref{eq:70074a}), (\ref{eq:70075}) and (\ref{eq:70076}) into (\ref{eq:70057a}) and (\ref{eq:70057b}).
\begin{subequations}
\begin{equation}
A_n = \frac{(n-j)(n-j+1+\beta -\delta )}{\beta (n+1)} \label{eq:70077a}
\end{equation}
\begin{equation}
B_n = -\frac{(n-j)(n-j-1)}{\beta (n+1)} \label{eq:70077b}
\end{equation}
\end{subequations}
For the case of $\alpha = \gamma -j-1 $, putting $n= j$ with (\ref{eq:70074b}) in (\ref{eq:70057b}) and use the condition $B_{j}=0$ for $\gamma $, we obtain $\gamma =2$. 
Substitute $\alpha = 1-j $ and $\gamma =2$ into (\ref{eq:70057a}). Putting $n= j$ in the new (\ref{eq:70057a}) and use the condition $A_{j}=0$ for $q$, we derive $q= \beta -\delta $.
Take $\alpha = 1-j $, $\gamma =2$ and $q= \beta -\delta $ into (\ref{eq:70057a}) and (\ref{eq:70057b}). The new (\ref{eq:70057a}) is same as (\ref{eq:70077a}). And the new(\ref{eq:70057b}) is equivalent to (\ref{eq:70077b}).
Therefore, the second species complete polynomial $y(z)$ divided by $z^{\alpha }$ of the CHE about $x=\infty $ as $\alpha = -j$ and $\gamma =q= 0$ is equal to the independent solution $y(z)$ divided by $z^{\alpha }$ of the CHE as $\alpha = 1-j $, $\gamma =2$ and $q= \beta -\delta $.
  
In (\ref{eq:70072a})--(\ref{eq:70073c}) replace $c_0$, $\lambda $ and $x$ by $1$, $\alpha$ and $z$. 
Substitute (\ref{eq:70077a}) and (\ref{eq:70077b}) into (\ref{eq:70073a})--(\ref{eq:70073c}).

\underline{(1) The case of $\alpha = -j$ and $\gamma =q= 0$,} 
 
As $B_1=A_0=0$, substitute the new (\ref{eq:70073a}) with $\alpha =0$ into (\ref{eq:70072a}) putting $j=0$. 
As $B_{2N+1}=B_{2N+2}=A_{2N+1}=0$, substitute the new (\ref{eq:70073a})--(\ref{eq:70073c}) with $\alpha =-2N-1$ into (\ref{eq:70072b}) putting $j=2N+1$.
As $B_{2N+2}=B_{2N+3}=A_{2N+2}=0$, substitute the new (\ref{eq:70073a})--(\ref{eq:70073c})  with $\alpha =-2N-2$ into (\ref{eq:70072c}) putting $j=2N+2$.

\underline{(2) The case of $\alpha = 1-j $, $\gamma =2$ and $q= \beta -\delta $,} 
 
As $B_1=A_0=0$, substitute the new (\ref{eq:70073a}) with $\alpha =1$ into (\ref{eq:70072a}) putting $j=0$. 
As $B_{2N+1}=B_{2N+2}=A_{2N+1}=0$, substitute the new (\ref{eq:70073a})--(\ref{eq:70073c}) with $\alpha = -2N $ into (\ref{eq:70072b}) putting $j=2N+1$.
As $B_{2N+2}=B_{2N+3}=A_{2N+2}=0$, substitute the new (\ref{eq:70073a})--(\ref{eq:70073c})  with $\alpha = -2N-1$ into (\ref{eq:70072c}) putting $j=2N+2$.

After the replacement process, we obtain the independent solution of the CHE about $x=\infty $. The solution is as follows.

\begin{remark}
The power series expansion of Confluent Heun equation of the first kind for the second species complete polynomial using 3TRF about $x=\infty $ is given by
\begin{enumerate} 
\item As $\alpha =\gamma =q=0$, 

Its eigenfunction is given by
\begin{eqnarray}
y(z) &=& _pH_c^{(r,1)}F_0 \left( \alpha =0, \beta, \gamma =0, \delta, q=0; z=\frac{1}{x}, \mu =-\frac{2}{\beta } z^2, \xi =\frac{2}{\beta } z\right) \nonumber\\
 &=& 1 \label{eq:70078a}
\end{eqnarray}
\item As $\alpha =-2N-1$, $\gamma =q=0$ where $N \in \mathbb{N}_{0}$, 

Its eigenfunction is given by
\begin{eqnarray}
y(z) &=& _pH_c^{(r,1)}F_{2N+1} \left( \alpha =-2N-1, \beta, \gamma =0, \delta, q=0; z=\frac{1}{x}, \mu =-\frac{2}{\beta } z^2, \xi =\frac{2}{\beta } z\right) \nonumber\\
 &=& z^{-2N-1} \left\{ \sum_{r=0}^{N} y_{2r}^{N-r}\left( 2N+1;z\right) + \sum_{r=0}^{N} y_{2r+1}^{N-r}\left( 2N+1;z\right) \right\} \label{eq:70078b}
\end{eqnarray}
\item As $\alpha =-2N-2$, $\gamma =q=0$ where $N \in \mathbb{N}_{0}$, 

Its eigenfunction is given by
\begin{eqnarray}
y(z) &=& _pH_c^{(r,1)}F_{2N+2} \left( \alpha =-2N-2, \beta, \gamma =0, \delta, q=0; z=\frac{1}{x}, \mu =-\frac{2}{\beta } z^2, \xi =\frac{2}{\beta } z\right) \nonumber\\
 &=& z^{-2N-2} \left\{ \sum_{r=0}^{N+1} y_{2r}^{N+1-r}\left( 2N+2;z\right) + \sum_{r=0}^{N} y_{2r+1}^{N-r}\left( 2N+2;z\right) \right\} \label{eq:70078c}
\end{eqnarray}

\item As $\alpha =1, \gamma =2, q=\beta -\delta $,

Its eigenfunction is given by
\begin{eqnarray}
y(z) &=& _pH_c^{(r,2)}F_0  \left( \alpha =1, \beta, \gamma =2, \delta, q=\beta -\delta; z=\frac{1}{x}, \mu =-\frac{2}{\beta } z^2, \xi =\frac{2}{\beta } z\right) \nonumber\\
 &=& z \label{eq:70078a1}
\end{eqnarray}
\item As $ \alpha = -2N, \gamma =2, q=\beta -\delta $ where $N \in \mathbb{N}_{0}$,

Its eigenfunction is given by
\begin{eqnarray}
y(z) &=& _pH_c^{(r,2)}F_{2N+1}  \left( \alpha = -2N, \beta, \gamma =2, \delta, q=\beta -\delta; z=\frac{1}{x}, \mu =-\frac{2}{\beta } z^2, \xi =\frac{2}{\beta } z\right) \nonumber\\
 &=& z^{ -2N } \left\{ \sum_{r=0}^{N} y_{2r}^{N-r}\left( 2N+1;z\right) + \sum_{r=0}^{N} y_{2r+1}^{N-r}\left( 2N+1;z\right) \right\} \label{eq:70078b1}
\end{eqnarray}
\item As $ \alpha = -2N-1, \gamma =2, q=\beta -\delta $ where $N \in \mathbb{N}_{0}$,

Its eigenfunction is given by
\begin{eqnarray}
y(z) &=& _pH_c^{(r,2)}F_{2N+2}  \left( \alpha = -2N-1, \beta, \gamma =2, \delta, q=\beta -\delta; z=\frac{1}{x}, \mu =-\frac{2}{\beta } z^2, \xi =\frac{2}{\beta } z\right) \nonumber\\
 &=& z^{ -2N-1} \left\{ \sum_{r=0}^{N+1} y_{2r}^{N+1-r}\left( 2N+2;z\right) + \sum_{r=0}^{N} y_{2r+1}^{N-r}\left( 2N+2;z\right) \right\} \label{eq:70078c1}
\end{eqnarray}
In the above,
\begin{eqnarray}
y_0^m(j;z) &=& \sum_{i_0=0}^{m} \frac{\left( -\frac{j}{2}\right)_{i_0}\left(  \frac{1}{2}-\frac{j}{2}  \right)_{i_0}}{\left( 1 \right)_{i_0}} \mu ^{i_0} \label{eq:70079a}\\
y_1^m(j;z) &=& \left\{\sum_{i_0=0}^{m} \frac{ \left( i_0 -\frac{j}{2}\right) \left( i_0 +\mho_0^j \right)}{\left( i_0+\frac{1}{2} \right)} \frac{\left( -\frac{j}{2}\right)_{i_0} \left(  \frac{1}{2}-\frac{j}{2} \right)_{i_0}}{\left( 1 \right)_{i_0}} \right. \nonumber\\
&&\times \left. \sum_{i_1 = i_0}^{m} \frac{\left( \frac{1}{2}-\frac{j}{2} \right)_{i_1} \left( 1-\frac{j}{2} \right)_{i_1}\left( \frac{3}{2} \right)_{i_0} }{\left( \frac{1}{2}-\frac{j}{2} \right)_{i_0} \left( 1-\frac{j}{2} \right)_{i_0}\left( \frac{3}{2} \right)_{i_1}} \mu ^{i_1}\right\}\xi 
\hspace{1.5cm}\label{eq:70079b}\\
y_{\tau }^m(j;z) &=& \left\{ \sum_{i_0=0}^{m} \frac{ \left( i_0 -\frac{j}{2}\right) \left( i_0 +\mho_0^j \right) }{\left( i_0+\frac{1}{2} \right)} \frac{\left( -\frac{j}{2}\right)_{i_0} \left(  \frac{1}{2}-\frac{j}{2} \right)_{i_0}}{\left( 1 \right)_{i_0}} \right.\nonumber\\
&&\times \prod_{k=1}^{\tau -1} \left( \sum_{i_k = i_{k-1}}^{m} \frac{\left( i_k+ \frac{k}{2} -\frac{j}{2}\right)\left( i_k +\mho _k^j  \right) }{\left( i_k+\frac{k}{2}+\frac{1}{2} \right)} \right.   \left. \frac{\left( \frac{k}{2}-\frac{j}{2}\right)_{i_k} \left( \frac{k}{2}+ \frac{1}{2}-\frac{j}{2} \right)_{i_{k}}\left( \frac{k}{2}+1 \right)_{i_{k-1}} }{\left( \frac{k}{2}-\frac{j}{2}\right)_{i_{k-1}} \left( \frac{k}{2}+ \frac{1}{2}-\frac{j}{2} \right)_{i_{k-1}}\left( \frac{k}{2}+1 \right)_{i_k}} \right) \nonumber\\
&&\times  \left. \sum_{i_{\tau } = i_{\tau -1}}^{m} \frac{\left( \frac{\tau }{2} -\frac{j}{2}\right)_{i_{\tau}} \left( \frac{\tau }{2}+\frac{1}{2}-\frac{j}{2} \right)_{i_{\tau}}\left( \frac{\tau }{2}+1 \right)_{i_{\tau -1}} }{\left( \frac{\tau }{2}-\frac{j}{2}\right)_{i_{\tau -1}} \left( \frac{\tau }{2}+\frac{1}{2}-\frac{j}{2} \right)_{i_{\tau -1}}\left( \frac{\tau }{2}+1 \right)_{i_{\tau }} } \mu ^{i_{\tau }}\right\} \xi ^{\tau } \label{eq:70079c} 
\end{eqnarray}
where
\begin{equation}
\begin{cases} \tau \geq 2 \cr 
\mho_0^j = \frac{1}{2}\left( \beta -\delta +1-j \right) \cr
\mho_k^j = \frac{1}{2}\left( \beta -\delta +1-j+k \right)
\end{cases}\nonumber
\end{equation}
\end{enumerate}
\end{remark} 
\subsection{The first species complete polynomial of Confluent Heun equation using R3TRF}
\begin{theorem}
In chapter 2, the general expression of a function $y(x)$ for the first species complete polynomial using reversible 3-term recurrence formula and its algebraic equation for the determination of an accessory parameter in $A_n$ term are given by
\begin{enumerate} 
\item As $B_1=0$,
\begin{equation}
0 =\bar{c}(0,1) \label{eq:70029a}
\end{equation}
\begin{equation}
y(x) = y_{0}^{0}(x) \label{eq:70029b}
\end{equation}
\item As $B_2=0$, 
\begin{equation}
0 = \bar{c}(0,2)+\bar{c}(1,0) \label{eq:70030a}
\end{equation}
\begin{equation}
y(x)= y_{0}^{1}(x) \label{eq:70030b}
\end{equation}
\item As $B_{2N+3}=0$ where $N \in \mathbb{N}_{0}$,
\begin{equation}
0  = \sum_{r=0}^{N+1}\bar{c}\left( r, 2(N-r)+3\right) \label{eq:70031a}
\end{equation}
\begin{equation}
y(x)= \sum_{r=0}^{N+1} y_{r}^{2(N+1-r)}(x) \label{eq:70031b}
\end{equation}
\item As $B_{2N+4}=0$ where$N \in \mathbb{N}_{0}$,
\begin{equation}
0  = \sum_{r=0}^{N+2}\bar{c}\left( r, 2(N+2-r)\right) \label{eq:70032a}
\end{equation}
\begin{equation}
y(x)=  \sum_{r=0}^{N+1} y_{r}^{2(N-r)+3}(x) \label{eq:70032b}
\end{equation}
In the above,
\begin{eqnarray}
\bar{c}(0,n) &=& \prod _{i_0=0}^{n-1}A_{i_0} \label{eq:70033a}\\
\bar{c}(1,n) &=& \sum_{i_0=0}^{n} \left\{ B_{i_0+1} \prod _{i_1=0}^{i_0-1}A_{i_1} \prod _{i_2=i_0}^{n-1}A_{i_2+2} \right\} \label{eq:70033b}\\
\bar{c}(\tau ,n) &=& \sum_{i_0=0}^{n} \left\{B_{i_0+1}\prod _{i_1=0}^{i_0-1} A_{i_1} 
\prod _{k=1}^{\tau -1} \left( \sum_{i_{2k}= i_{2(k-1)}}^{n} B_{i_{2k}+(2k+1)}\prod _{i_{2k+1}=i_{2(k-1)}}^{i_{2k}-1}A_{i_{2k+1}+2k}\right)\right. \nonumber\\
&&\times \left. \prod _{i_{2\tau} = i_{2(\tau -1)}}^{n-1} A_{i_{2\tau }+ 2\tau} \right\} 
\hspace{1cm}\label{eq:70033c}
\end{eqnarray}
and
\begin{eqnarray}
y_0^m(x) &=& c_0 x^{\lambda} \sum_{i_0=0}^{m} \left\{ \prod _{i_1=0}^{i_0-1}A_{i_1} \right\} x^{i_0 } \label{eq:70034a}\\
y_1^m(x) &=& c_0 x^{\lambda} \sum_{i_0=0}^{m}\left\{ B_{i_0+1} \prod _{i_1=0}^{i_0-1}A_{i_1}  \sum_{i_2=i_0}^{m} \left\{ \prod _{i_3=i_0}^{i_2-1}A_{i_3+2} \right\}\right\} x^{i_2+2 } \label{eq:70034b}\\
y_{\tau }^m(x) &=& c_0 x^{\lambda} \sum_{i_0=0}^{m} \left\{B_{i_0+1}\prod _{i_1=0}^{i_0-1} A_{i_1} 
\prod _{k=1}^{\tau -1} \left( \sum_{i_{2k}= i_{2(k-1)}}^{m} B_{i_{2k}+(2k+1)}\prod _{i_{2k+1}=i_{2(k-1)}}^{i_{2k}-1}A_{i_{2k+1}+2k}\right) \right. \nonumber\\
&&\times \left. \sum_{i_{2\tau} = i_{2(\tau -1)}}^{m} \left( \prod _{i_{2\tau +1}=i_{2(\tau -1)}}^{i_{2\tau}-1} A_{i_{2\tau +1}+ 2\tau} \right) \right\} x^{i_{2\tau}+2\tau }\hspace{1cm}\mathrm{where}\;\tau \geq 2
\label{eq:70034c}
\end{eqnarray}
\end{enumerate}
\end{theorem}
For the first species complete polynomial of the CHE about $x=\infty $ in section 8.2.1, there are two possible fixed values of $\alpha $ such as $\alpha =-j$ and $\gamma -j-1$.  

According to (\ref{eq:7006}), $c_{j+1}=0$ is a polynomial equation of degree $j+1$ for the determination of the accessory parameter $q$ and thus has $j+1$ zeros denoted them by $q_j^m$ eigenvalues where $m = 0,1,2, \cdots, j$. They can be arranged in the following order: $q_j^0 < q_j^1 < q_j^2 < \cdots < q_j^j$.
\subsubsection{The case of $\alpha = -j $ }
In (\ref{eq:70029b}), (\ref{eq:70030b}), (\ref{eq:70031b}), (\ref{eq:70032b}) and (\ref{eq:70034a})--(\ref{eq:70034c}) replace $c_0$, $\lambda $ and $x$ by $1$, $\alpha$ and $z$. Substitute (\ref{eq:70059a}) and (\ref{eq:70059b}) into (\ref{eq:70033a})--(\ref{eq:70034c}).

As $B_{1}= c_{1}=0$, take the new (\ref{eq:70033a}) into (\ref{eq:70029a}) putting $j=0$. Substitute the new (\ref{eq:70034a}) with $\alpha =0$ into (\ref{eq:70029b}) putting $j=0$.

As $B_{2}= c_{2}=0$, take the new (\ref{eq:70033a}) and (\ref{eq:70033b}) into (\ref{eq:70030a}) putting $j=1$. Substitute the new (\ref{eq:70034a}) with $\alpha =-1$ into (\ref{eq:70030b}) putting $j=1$ and $q=q_1^m$. 

As $B_{2N+3}= c_{2N+3}=0$, take the new (\ref{eq:70033a})--(\ref{eq:70033c}) into (\ref{eq:70031a}) putting $j=2N+2$. Substitute the new 
(\ref{eq:70034a})--(\ref{eq:70034c}) with $\alpha =-2N-2$ into (\ref{eq:70031b}) putting $j=2N+2$ and $q=q_{2N+2}^m$.

As $B_{2N+4}= c_{2N+4}=0$, take the new (\ref{eq:70033a})--(\ref{eq:70033c}) into (\ref{eq:70032a}) putting $j=2N+3$. Substitute the new 
(\ref{eq:70034a})--(\ref{eq:70034c}) with $\alpha =-2N-3$ into (\ref{eq:70032b}) putting $j=2N+3$ and $q=q_{2N+3}^m$.

After the replacement process, we obtain the independent solution of the CHE about $x=\infty $. The solution is as follows.
\begin{remark}
The power series expansion of Confluent Heun equation of the first kind for the first species complete polynomial using R3TRF about $x=\infty $ for $\alpha =$ non-positive integer is given by
\begin{enumerate} 
\item As $\alpha =0$ and $q=q_0^0=0$,

The eigenfunction is given by
\begin{equation}
y(z) =\; _pH_c^{(r,1)}F_{0,0}^R \left( \alpha =0, \beta, \gamma, \delta, q=q_0^0=0; z=\frac{1}{x}, \sigma = \frac{1}{\beta }z, \varsigma =-\frac{1}{\beta } z^2\right) =1\label{eq:70080}
\end{equation}
\item As $\alpha =-1$,

An algebraic equation of degree 2 for the determination of $q$ is given by
\begin{equation}
0 = -\beta \gamma + q(\beta -\gamma -\delta +q)  \label{eq:70081a}
\end{equation}
The eigenvalue of $q$ is written by $q_1^m$ where $m = 0,1 $; $q_{1}^0 < q_{1}^1$. Its eigenfunction is given by
\begin{eqnarray}
y(z) &=& _pH_c^{(r,1)}F_{1,m}^R \left( \alpha =-1, \beta, \gamma, \delta, q=q_1^m; z=\frac{1}{x}, \sigma = \frac{1}{\beta }z, \varsigma =-\frac{1}{\beta } z^2\right)\nonumber\\
&=&  z^{-1}\left\{ 1 -(\beta -\gamma -\delta +q_1^m)\sigma \right\} \label{eq:70081b}  
\end{eqnarray}
\item As $\alpha =-2N-2 $ where $N \in \mathbb{N}_{0}$,

An algebraic equation of degree $2N+3$ for the determination of $q$ is given by
\begin{equation}
0 = \sum_{r=0}^{N+1}\bar{c}\left( r, 2(N-r)+3; 2N+2,q\right)  \label{eq:70082a}
\end{equation}
The eigenvalue of $q$ is written by $q_{2N+2}^m$ where $m = 0,1,2,\cdots,2N+2 $; $q_{2N+2}^0 < q_{2N+2}^1 < \cdots < q_{2N+2}^{2N+2}$. Its eigenfunction is given by 
\begin{eqnarray} 
y(z) &=& _pH_c^{(r,1)}F_{2N+2,m}^R \left( \alpha =-2N-2, \beta, \gamma, \delta, q=q_{2N+2}^m ; z=\frac{1}{x}, \sigma = \frac{1}{\beta }z, \varsigma =-\frac{1}{\beta } z^2 \right)\nonumber\\
&=& z^{-2N-2} \sum_{r=0}^{N+1} y_{r}^{2(N+1-r)}\left( 2N+2, q_{2N+2}^m; z \right) 
\label{eq:70082b} 
\end{eqnarray}
\item As $\alpha =-2N-3 $ where $N \in \mathbb{N}_{0}$,

An algebraic equation of degree $2N+4$ for the determination of $q$ is given by
\begin{equation}  
0 = \sum_{r=0}^{N+2}\bar{c}\left( r, 2(N+2-r); 2N+3,q\right) \label{eq:70083a}
\end{equation}
The eigenvalue of $q$ is written by $q_{2N+3}^m$ where $m = 0,1,2,\cdots,2N+3 $; $q_{2N+3}^0 < q_{2N+3}^1 < \cdots < q_{2N+3}^{2N+3}$. Its eigenfunction is given by
\begin{eqnarray} 
y(z) &=& _pH_c^{(r,1)}F_{2N+3,m}^R \left( \alpha =-2N-3, \beta, \gamma, \delta, q=q_{2N+3}^m ; z=\frac{1}{x}, \sigma = \frac{1}{\beta }z, \varsigma =-\frac{1}{\beta } z^2\right)\nonumber\\
&=& z^{-2N-3} \sum_{r=0}^{N+1} y_{r}^{2(N-r)+3} \left( 2N+3,q_{2N+3}^m;z\right) \label{eq:70083b}
\end{eqnarray}
In the above,
\begin{eqnarray}
\bar{c}(0,n;j,q)  &=& \frac{\left( \Lambda _0^{-} \left( j,q\right) \right)_{n}\left( \Lambda _0^{+} \left( j,q\right) \right)_{n}}{\left( 1 \right)_{n} }\left( \frac{1}{\beta }\right)^n \label{eq:70084a}\\
\bar{c}(1,n;j,q) &=& \left( -\frac{1}{\beta }\right) \sum_{i_0=0}^{n}\frac{\left( i_0 -j\right)\left( i_0-j+1- \gamma \right) }{\left( i_0+2 \right)} \frac{ \left( \Lambda _0^{-} \left( j,q\right) \right)_{i_0}\left( \Lambda _0^{+} \left( j,q\right) \right)_{i_0}}{\left( 1 \right)_{i_0}}  \nonumber\\
&&\times  \frac{ \left( \Lambda _1^{-} \left( j,q\right) \right)_{n}\left( \Lambda _1^{+} \left( j,q\right) \right)_{n} \left( 3\right)_{i_0} }{\left( \Lambda _1^{-} \left( j,q\right) \right)_{i_0}\left( \Lambda _1^{+} \left( j,q\right) \right)_{i_0}\left( 3 \right)_{n}} \left( \frac{1}{\beta }\right)^n \label{eq:70084b}\\
\bar{c}(\tau ,n;j,q) &=& \left( -\frac{1}{\beta }\right)^{\tau} \sum_{i_0=0}^{n} \frac{\left( i_0 -j\right)\left( i_0-j+1- \gamma \right) }{\left( i_0+2 \right)} \frac{ \left( \Lambda _0^{-} \left( j,q\right) \right)_{i_0}\left( \Lambda _0^{+} \left( j,q\right) \right)_{i_0}}{\left( 1 \right)_{i_0}}  \nonumber\\
&&\times \prod_{k=1}^{\tau -1} \left( \sum_{i_k = i_{k-1}}^{n} \frac{\left( i_k+ 2k-j\right)\left( i_k+2k-j+1- \gamma \right) }{\left( i_k+2k+2 \right)} \right.  \nonumber\\
&&\times \left. \frac{ \left( \Lambda _k^{-} \left( j,q\right) \right)_{i_k}\left( \Lambda _k^{+} \left( j,q\right) \right)_{i_k} \left( 2k+1 \right)_{i_{k-1}}}{\left( \Lambda _k^{-} \left( j,q\right) \right)_{i_{k-1}}\left( \Lambda _k^{+} \left( j,q\right) \right)_{i_{k-1}}\left( 2k+1 \right)_{i_k}} \right) \nonumber\\
&&\times \frac{ \left( \Lambda _{\tau }^{-} \left( j,q\right) \right)_{n}\left( \Lambda _{\tau }^{+} \left( j,q\right) \right)_{n} \left( 2\tau +1 \right)_{i_{\tau -1}}}{\left( \Lambda _{\tau }^{-} \left( j,q\right) \right)_{i_{\tau -1}}\left( \Lambda _{\tau }^{+} \left( j,q\right) \right)_{i_{\tau -1}}\left( 2\tau +1 \right)_{n}}\left( \frac{1}{\beta }\right)^n  \label{eq:70084c} 
\end{eqnarray}
\begin{eqnarray}
y_0^m(j,q;z) &=& \sum_{i_0=0}^{m} \frac{\left( \Lambda _0^{-} \left( j,q\right) \right)_{i_0}\left( \Lambda _0^{+} \left( j,q\right) \right)_{i_0}}{\left( 1 \right)_{i_0} } \sigma^{i_0} \label{eq:70085a}\\
y_1^m(j,q;z) &=& \left\{\sum_{i_0=0}^{m}\frac{\left( i_0 -j\right)\left( i_0-j+1- \gamma \right) }{\left( i_0+2 \right)} \frac{ \left( \Lambda _0^{-} \left( j,q\right) \right)_{i_0}\left( \Lambda _0^{+} \left( j,q\right) \right)_{i_0}}{\left( 1 \right)_{i_0}} \right.  \nonumber\\
&&\times \left. \sum_{i_1 = i_0}^{m} \frac{ \left( \Lambda _1^{-} \left( j,q\right) \right)_{i_1}\left( \Lambda _1^{+} \left( j,q\right) \right)_{i_1} \left( 3\right)_{i_0} }{\left( \Lambda _1^{-} \left( j,q\right) \right)_{i_0}\left( \Lambda _1^{+} \left( j,q\right) \right)_{i_0}\left( 3 \right)_{i_1}} \sigma^{i_1}\right\} \varsigma   
\hspace{1cm}\label{eq:70085b}\\
y_{\tau }^m(j,q;z) &=& \left\{ \sum_{i_0=0}^{m} \frac{\left( i_0 -j\right)\left( i_0-j+1- \gamma \right) }{\left( i_0+2 \right)} \frac{ \left( \Lambda _0^{-} \left( j,q\right) \right)_{i_0}\left( \Lambda _0^{+} \left( j,q\right) \right)_{i_0}}{\left( 1 \right)_{i_0}} \right.\nonumber\\
&&\times \prod_{k=1}^{\tau -1} \left( \sum_{i_k = i_{k-1}}^{m}  \frac{\left( i_k+ 2k-j\right)\left( i_k+2k-j+1- \gamma \right) }{\left( i_k+2k+2 \right)} \right.  \nonumber\\
&&\times \left. \frac{ \left( \Lambda _k^{-} \left( j,q\right) \right)_{i_k}\left( \Lambda _k^{+} \left( j,q\right) \right)_{i_k} \left( 2k+1 \right)_{i_{k-1}}}{\left( \Lambda _k^{-} \left( j,q\right) \right)_{i_{k-1}}\left( \Lambda _k^{+} \left( j,q\right) \right)_{i_{k-1}}\left( 2k+1 \right)_{i_k}} \right) \nonumber\\
&&\times \left. \sum_{i_{\tau } = i_{\tau -1}}^{m}  \frac{ \left( \Lambda _{\tau }^{-} \left( j,q\right) \right)_{i_{\tau }}\left( \Lambda _{\tau }^{+} \left( j,q\right) \right)_{i_{\tau }} \left( 2\tau +1 \right)_{i_{\tau -1}}}{\left( \Lambda _{\tau }^{-} \left( j,q\right) \right)_{i_{\tau -1}}\left( \Lambda _{\tau }^{+} \left( j,q\right) \right)_{i_{\tau -1}}\left( 2\tau +1 \right)_{i_{\tau }}} \sigma^{i_{\tau }}\right\} \varsigma ^{\tau } \hspace{1cm}\label{eq:70085c} 
\end{eqnarray}
where
\begin{equation}
\begin{cases} \tau \geq 2 \cr 
\Lambda _{k}^{\pm}(j,q) = \frac{\varpi -2j+4k \pm \sqrt{\varpi ^2 +4q}}{2}\cr
\varpi =  \beta -\gamma -\delta +1
\end{cases}\nonumber
\end{equation}
\end{enumerate}
\end{remark}
\subsubsection{The case of $\alpha =\gamma -j-1 $ }
In (\ref{eq:70029b}), (\ref{eq:70030b}), (\ref{eq:70031b}), (\ref{eq:70032b}) and (\ref{eq:70034a})--(\ref{eq:70034c}) replace $c_0$, $\lambda $ and $x$ by $1$, $\alpha$ and $z$. Substitute (\ref{eq:70065a}) and (\ref{eq:70065b}) into (\ref{eq:70033a})--(\ref{eq:70034c}).

As $B_{1}= c_{1}=0$, take the new (\ref{eq:70033a}) into (\ref{eq:70029a}) putting $j=0$. Substitute the new (\ref{eq:70034a}) with $\alpha =\gamma -1$ into (\ref{eq:70029b}) putting $j=0$.

As $B_{2}= c_{2}=0$, take the new (\ref{eq:70033a}) and (\ref{eq:70033b}) into (\ref{eq:70030a}) putting $j=1$. Substitute the new (\ref{eq:70034a}) with $\gamma -2$ into (\ref{eq:70030b}) putting $j=1$ and $q=q_1^m$. 

As $B_{2N+3}= c_{2N+3}=0$, take the new (\ref{eq:70033a})--(\ref{eq:70033c}) into (\ref{eq:70031a}) putting $j=2N+2$. Substitute the new 
(\ref{eq:70034a})--(\ref{eq:70034c}) with $\gamma -2N-3$ into (\ref{eq:70031b}) putting $j=2N+2$ and $q=q_{2N+2}^m$.

As $B_{2N+4}= c_{2N+4}=0$, take the new (\ref{eq:70033a})--(\ref{eq:70033c}) into (\ref{eq:70032a}) putting $j=2N+3$. Substitute the new 
(\ref{eq:70034a})--(\ref{eq:70034c}) with $\gamma -2N-4$ into (\ref{eq:70032b}) putting $j=2N+3$ and $q=q_{2N+3}^m$.

After the replacement process, we obtain the independent solution of the CHE about $x=\infty $. The solution is as follows.
\begin{remark}
The power series expansion of Confluent Heun equation of the first kind for the first species complete polynomial using R3TRF about $x=\infty $ for $\alpha =\gamma -1-j$ where $j\in \mathbb{N}_{0}$ is given by
\begin{enumerate} 
\item As $\alpha =\gamma -1$ and $q=q_0^0=(\gamma -1)(\beta -\delta )$,

The eigenfunction is given by
\begin{eqnarray}
y(z) &=& _pH_c^{(r,2)}F_{0,0}^R \left( \alpha =\gamma -1, \beta, \gamma, \delta, q=q_0^0=(\gamma -1)(\beta -\delta ); z=\frac{1}{x}, \sigma = \frac{1}{\beta }z, \varsigma =-\frac{1}{\beta } z^2\right) \nonumber\\
&=& z^{\gamma -1}\label{eq:70086}
\end{eqnarray}
\item As $\alpha =\gamma -2$,

An algebraic equation of degree 2 for the determination of $q$ is given by
\begin{equation}
0 = \beta (\gamma -2) -\prod_{l=0}^{1}\Big( q-(\beta -\delta -l)(\gamma -1-l) \Big) \label{eq:70087a}
\end{equation}
The eigenvalue of $q$ is written by $q_1^m$ where $m = 0,1 $; $q_{1}^0 < q_{1}^1$. Its eigenfunction is given by
\begin{eqnarray}
y(z) &=& _pH_c^{(r,2)}F_{1,m}^R \left( \alpha =\gamma -2, \beta, \gamma, \delta, q=q_1^m; z=\frac{1}{x}, \sigma = \frac{1}{\beta }z, \varsigma =-\frac{1}{\beta } z^2\right)\nonumber\\
&=&  z^{\gamma -2}\left\{ 1 +\left( (\beta -\delta -1)(\gamma -2)- q_1^m\right)\sigma \right\} \label{eq:70087b}  
\end{eqnarray}
\item As $\alpha =\gamma -2N-3 $ where $N \in \mathbb{N}_{0}$,

An algebraic equation of degree $2N+3$ for the determination of $q$ is given by
\begin{equation}
0 = \sum_{r=0}^{N+1}\bar{c}\left( r, 2(N-r)+3; 2N+2,q\right)  \label{eq:70088a}
\end{equation}
The eigenvalue of $q$ is written by $q_{2N+2}^m$ where $m = 0,1,2,\cdots,2N+2 $; $q_{2N+2}^0 < q_{2N+2}^1 < \cdots < q_{2N+2}^{2N+2}$. Its eigenfunction is given by 
\begin{eqnarray} 
y(z) &=& _pH_c^{(r,2)}F_{2N+2,m}^R \left( \alpha =\gamma -2N-3, \beta, \gamma, \delta, q=q_{2N+2}^m ; z=\frac{1}{x}, \sigma = \frac{1}{\beta }z, \varsigma =-\frac{1}{\beta } z^2 \right)\nonumber\\
&=& z^{\gamma -2N-3} \sum_{r=0}^{N+1} y_{r}^{2(N+1-r)}\left( 2N+2, q_{2N+2}^m; z \right) 
\label{eq:70088b} 
\end{eqnarray}
\item As $\alpha =\gamma -2N-4 $ where $N \in \mathbb{N}_{0}$,

An algebraic equation of degree $2N+4$ for the determination of $q$ is given by
\begin{equation}  
0 = \sum_{r=0}^{N+2}\bar{c}\left( r, 2(N+2-r); 2N+3,q\right) \label{eq:70089a}
\end{equation}
The eigenvalue of $q$ is written by $q_{2N+3}^m$ where $m = 0,1,2,\cdots,2N+3 $; $q_{2N+3}^0 < q_{2N+3}^1 < \cdots < q_{2N+3}^{2N+3}$. Its eigenfunction is given by
\begin{eqnarray} 
y(z) &=& _pH_c^{(r,2)}F_{2N+3,m}^R \left( \alpha =\gamma -2N-4, \beta, \gamma, \delta, q=q_{2N+3}^m ; z=\frac{1}{x}, \sigma = \frac{1}{\beta }z, \varsigma =-\frac{1}{\beta } z^2\right)\nonumber\\
&=& z^{\gamma -2N-4} \sum_{r=0}^{N+1} y_{r}^{2(N-r)+3} \left( 2N+3,q_{2N+3}^m;z\right) \label{eq:70089b}
\end{eqnarray}
In the above,
\begin{eqnarray}
\bar{c}(0,n;j,q)  &=& \frac{\left( \Psi _0^{-} \left( j,q\right) \right)_{n}\left( \Psi _0^{+} \left( j,q\right) \right)_{n}}{\left( 1 \right)_{n} }\left( \frac{1}{\beta }\right)^n \label{eq:70090a}\\
\bar{c}(1,n;j,q) &=& \left( -\frac{1}{\beta }\right) \sum_{i_0=0}^{n}\frac{\left( i_0 -j\right)\left( i_0-j-1+ \gamma \right) }{\left( i_0+2 \right)} \frac{ \left( \Psi _0^{-} \left( j,q\right) \right)_{i_0}\left( \Psi _0^{+} \left( j,q\right) \right)_{i_0}}{\left( 1 \right)_{i_0}}  \nonumber\\
&&\times  \frac{ \left( \Psi _1^{-} \left( j,q\right) \right)_{n}\left( \Psi _1^{+} \left( j,q\right) \right)_{n} \left( 3\right)_{i_0} }{\left( \Psi _1^{-} \left( j,q\right) \right)_{i_0}\left( \Psi _1^{+} \left( j,q\right) \right)_{i_0}\left( 3 \right)_{n}} \left( \frac{1}{\beta }\right)^n \label{eq:70090b}
\end{eqnarray}
\begin{eqnarray}
\bar{c}(\tau ,n;j,q) &=& \left( -\frac{1}{\beta }\right)^{\tau} \sum_{i_0=0}^{n} \frac{\left( i_0 -j\right)\left( i_0-j-1+ \gamma \right) }{\left( i_0+2 \right)} \frac{ \left( \Psi _0^{-} \left( j,q\right) \right)_{i_0}\left( \Psi _0^{+} \left( j,q\right) \right)_{i_0}}{\left( 1 \right)_{i_0}}  \nonumber\\
&&\times \prod_{k=1}^{\tau -1} \left( \sum_{i_k = i_{k-1}}^{n} \frac{\left( i_k+ 2k-j\right)\left( i_k+2k-j-1+ \gamma \right) }{\left( i_k+2k+2 \right)} \right. \nonumber\\
&&\times  \left. \frac{ \left( \Psi _k^{-} \left( j,q\right) \right)_{i_k}\left( \Psi _k^{+} \left( j,q\right) \right)_{i_k} \left( 2k+1 \right)_{i_{k-1}}}{\left( \Psi _k^{-} \left( j,q\right) \right)_{i_{k-1}}\left( \Psi _k^{+} \left( j,q\right) \right)_{i_{k-1}}\left( 2k+1 \right)_{i_k}} \right) \nonumber\\
&&\times \frac{ \left( \Psi _{\tau }^{-} \left( j,q\right) \right)_{n}\left( \Psi _{\tau }^{+} \left( j,q\right) \right)_{n} \left( 2\tau +1 \right)_{i_{\tau -1}}}{\left( \Psi _{\tau }^{-} \left( j,q\right) \right)_{i_{\tau -1}}\left( \Psi _{\tau }^{+} \left( j,q\right) \right)_{i_{\tau -1}}\left( 2\tau +1 \right)_{n}}\left( \frac{1}{\beta }\right)^n  \label{eq:70090c} 
\end{eqnarray}
\begin{eqnarray}
y_0^m(j,q;z) &=& \sum_{i_0=0}^{m} \frac{\left( \Psi _0^{-} \left( j,q\right) \right)_{i_0}\left( \Psi _0^{+} \left( j,q\right) \right)_{i_0}}{\left( 1 \right)_{i_0} } \sigma^{i_0} \label{eq:70091a}\\
y_1^m(j,q;z) &=& \left\{\sum_{i_0=0}^{m}\frac{\left( i_0 -j\right)\left( i_0-j-1+ \gamma \right) }{\left( i_0+2 \right)} \frac{ \left( \Psi _0^{-} \left( j,q\right) \right)_{i_0}\left( \Psi _0^{+} \left( j,q\right) \right)_{i_0}}{\left( 1 \right)_{i_0}} \right.  \nonumber\\
&&\times \left. \sum_{i_1 = i_0}^{m} \frac{ \left( \Psi _1^{-} \left( j,q\right) \right)_{i_1}\left( \Psi _1^{+} \left( j,q\right) \right)_{i_1} \left( 3\right)_{i_0} }{\left( \Psi _1^{-} \left( j,q\right) \right)_{i_0}\left( \Psi _1^{+} \left( j,q\right) \right)_{i_0}\left( 3 \right)_{i_1}} \sigma^{i_1}\right\} \varsigma   
\hspace{1cm}\label{eq:70091b}\\
y_{\tau }^m(j,q;z) &=& \left\{ \sum_{i_0=0}^{m} \frac{\left( i_0 -j\right)\left( i_0-j-1+ \gamma \right) }{\left( i_0+2 \right)} \frac{ \left( \Psi _0^{-} \left( j,q\right) \right)_{i_0}\left( \Psi _0^{+} \left( j,q\right) \right)_{i_0}}{\left( 1 \right)_{i_0}} \right.\nonumber\\
&&\times \prod_{k=1}^{\tau -1} \left( \sum_{i_k = i_{k-1}}^{m}  \frac{\left( i_k+ 2k-j\right)\left( i_k+2k-j-1+ \gamma \right) }{\left( i_k+2k+2 \right)} \right. \nonumber\\
&&\times  \left. \frac{ \left( \Psi _k^{-} \left( j,q\right) \right)_{i_k}\left( \Psi _k^{+} \left( j,q\right) \right)_{i_k} \left( 2k+1 \right)_{i_{k-1}}}{\left( \Psi _k^{-} \left( j,q\right) \right)_{i_{k-1}}\left( \Psi _k^{+} \left( j,q\right) \right)_{i_{k-1}}\left( 2k+1 \right)_{i_k}} \right) \nonumber\\
&&\times \left. \sum_{i_{\tau } = i_{\tau -1}}^{m}  \frac{ \left( \Psi _{\tau }^{-} \left( j,q\right) \right)_{i_{\tau }}\left( \Psi _{\tau }^{+} \left( j,q\right) \right)_{i_{\tau }} \left( 2\tau +1 \right)_{i_{\tau -1}}}{\left( \Psi _{\tau }^{-} \left( j,q\right) \right)_{i_{\tau -1}}\left( \Psi _{\tau }^{+} \left( j,q\right) \right)_{i_{\tau -1}}\left( 2\tau +1 \right)_{i_{\tau }}} \sigma^{i_{\tau }}\right\} \varsigma ^{\tau } \hspace{1cm}\label{eq:70091c} 
\end{eqnarray}
where
\begin{equation}
\begin{cases} \tau \geq 2 \cr 
\Psi _{k}^{\pm}(j,q) = \frac{\omega  -2j+4k \pm \sqrt{\varphi ^2 +4q}}{2}\cr
\omega  =  \beta +\gamma -\delta -1 \cr
\varphi  = -\beta +\gamma +\delta -1
\end{cases}\nonumber
\end{equation}
\end{enumerate}
\end{remark}
\subsection{The second species complete polynomial of Confluent Heun equation using R3TRF}
As previously mentioned, we need a condition for the second species complete polynomial such as $B_{j}=B_{j+1}= A_{j}=0$ where $j \in \mathbb{N}_{0}$.
\begin{theorem}
In chapter 2, the general expression of a function $y(x)$ for the second species complete polynomial using reversible  3-term recurrence formula is given by
\begin{enumerate} 
\item As $B_1=A_0=0$,
\begin{equation}
y(x) = y_{0}^{0}(x) \label{eq:70092a}
\end{equation}
\item As $B_1=B_2=A_1=0$, 
\begin{equation}
y(x)= y_{0}^{1}(x)  \label{eq:70092b}
\end{equation}
\item As $B_{2N+2}=B_{2N+3}=A_{2N+2}=0$ where $N \in \mathbb{N}_{0}$,
\begin{equation}
y(x)= \sum_{r=0}^{N+1} y_{r}^{2(N+1-r)}(x) \label{eq:70092c}
\end{equation}
\item As $B_{2N+3}=B_{2N+4}=A_{2N+3}=0$ where $N \in \mathbb{N}_{0}$,
\begin{equation}
y(x)= \sum_{r=0}^{N+1} y_{r}^{2(N-r)+3}(x) \label{eq:70092d}
\end{equation}
In the above,
\begin{eqnarray}
y_0^m(x) &=& c_0 x^{\lambda} \sum_{i_0=0}^{m} \left\{ \prod _{i_1=0}^{i_0-1}A_{i_1} \right\} x^{i_0 }\label{eq:70093a}\\
y_1^m(x) &=& c_0 x^{\lambda} \sum_{i_0=0}^{m}\left\{ B_{i_0+1} \prod _{i_1=0}^{i_0-1}A_{i_1}  \sum_{i_2=i_0}^{m} \left\{ \prod _{i_3=i_0}^{i_2-1}A_{i_3+2} \right\}\right\} x^{i_2+2 } \label{eq:70093b}\\
y_{\tau }^m(x) &=& c_0 x^{\lambda} \sum_{i_0=0}^{m} \left\{B_{i_0+1}\prod _{i_1=0}^{i_0-1} A_{i_1} 
\prod _{k=1}^{\tau -1} \left( \sum_{i_{2k}= i_{2(k-1)}}^{m} B_{i_{2k}+(2k+1)}\prod _{i_{2k+1}=i_{2(k-1)}}^{i_{2k}-1}A_{i_{2k+1}+2k}\right) \right. \nonumber\\
&& \times \left. \sum_{i_{2\tau} = i_{2(\tau -1)}}^{m} \left( \prod _{i_{2\tau +1}=i_{2(\tau -1)}}^{i_{2\tau}-1} A_{i_{2\tau +1}+ 2\tau} \right) \right\} x^{i_{2\tau}+2\tau }\hspace{1cm}\mathrm{where}\;\tau \geq 2
\label{eq:70093c}
\end{eqnarray} 
\end{enumerate}
\end{theorem}
For the second species complete polynomial of the CHE about $x=\infty $ in section 8.2.2, there are two possible fixed values of $\alpha $ such as $\alpha =-j$ and $\gamma -j-1$ for $B_{j+1}=0$.  
As $\alpha =-j$, we need a fixed value $\gamma =0$ for $B_{j}=0$ and a fixed constant $q=0$ for $A_{j}=0$.
As $\alpha =\gamma -j-1$, a fixed value $\gamma =2$ for $B_{j}=0$ and a fixed constant $q=\beta -\delta $ for $A_{j}=0$ are required.
A series solution $y(z)$ divided by $z^{\alpha }$ of the CHE about $x=\infty $ as $\alpha = -j$ and $\gamma =q= 0$ for the second species complete polynomial using R3TRF is equal to the independent solution $y(z)$ divided by $z^{\alpha }$ of the CHE as $\alpha = 1-j $, $\gamma =2$ and $q= \beta -\delta $.
  
In (\ref{eq:70092a})--(\ref{eq:70093c}) replace $c_0$, $\lambda $ and $x$ by $1$, $\alpha$ and $z$. 
Substitute (\ref{eq:70077a}) and (\ref{eq:70077b}) into (\ref{eq:70093a})--(\ref{eq:70093c}).

\underline{(1) The case of $\alpha = -j$ and $\gamma =q= 0$,} 

As $B_1=A_0=0$, substitute the new (\ref{eq:70093a}) with $\alpha =0$ into (\ref{eq:70092a}) putting $j=0$. 
As $B_1=B_2=A_1=0$, substitute the new (\ref{eq:70093a}) with $\alpha =-1$ into (\ref{eq:70092b}) putting $j=1$. 
As $B_{2N+2}=B_{2N+3}=A_{2N+2}=0$, substitute the new (\ref{eq:70093a})--(\ref{eq:70093c}) with $\alpha =-2N-2$ into (\ref{eq:70092c}) putting $j=2N+2$.
As $B_{2N+3}=B_{2N+4}=A_{2N+3}=0$, substitute the new (\ref{eq:70093a})--(\ref{eq:70093c}) with $\alpha =-2N-3$ into (\ref{eq:70092d}) putting $j=2N+3$.

\underline{(2) The case of $\alpha = 1-j $, $\gamma =2$ and $q= \beta -\delta $,} 
 
As $B_1=A_0=0$, substitute the new (\ref{eq:70093a}) with $\alpha =1$ into (\ref{eq:70092a}) putting $j=0$. 
As $B_1=B_2=A_1=0$, substitute the new (\ref{eq:70093a}) with $\alpha =0$ into (\ref{eq:70092b}) putting $j=1$.
As $B_{2N+2}=B_{2N+3}=A_{2N+2}=0$, substitute the new (\ref{eq:70093a})--(\ref{eq:70093c}) with $\alpha = -2N-1$ into (\ref{eq:70092c}) putting $j=2N+2$.
As $B_{2N+3}=B_{2N+4}=A_{2N+3}=0$, substitute the new (\ref{eq:70093a})--(\ref{eq:70093c}) with $\alpha = -2N-2$ into (\ref{eq:70092d}) putting $j=2N+3$.

After the replacement process, we obtain the independent solution of the CHE about $x=\infty $. The solution is as follows.
\begin{remark}
The power series expansion of Confluent Heun equation of the first kind for the second species complete polynomial using R3TRF about $x=\infty $ is given by
\begin{enumerate} 
\item As $\alpha =\gamma =q=0$, 

Its eigenfunction is given by
\begin{eqnarray}
y(z) &=& _pH_c^{(r,1)}F_0^R \left( \alpha =0, \beta, \gamma =0, \delta, q=0; z=\frac{1}{x}, \sigma = \frac{1}{\beta }z, \varsigma =-\frac{1}{\beta } z^2\right) \nonumber\\
&=& 1 \label{eq:70094a}
\end{eqnarray}
\item As $\alpha =-1$, $\gamma =q=0$, 

Its eigenfunction is given by
\begin{eqnarray}
y(z) &=& _pH_c^{(r,1)}F_1^R \left( \alpha =-1, \beta, \gamma =0, \delta, q=0; z=\frac{1}{x}, \sigma = \frac{1}{\beta }z, \varsigma =-\frac{1}{\beta } z^2\right) \nonumber\\
&=& z^{-1}\left\{ 1-(\beta -\delta )\sigma \right\} \label{eq:70094b}
\end{eqnarray}
\item As $\alpha =-2N-2$, $\gamma =q=0$ where $N \in \mathbb{N}_{0}$, 

Its eigenfunction is given by
\begin{eqnarray}
y(z) &=& _pH_c^{(r,1)}F_{2N+2}^R \left( \alpha =-2N-2, \beta, \gamma =0, \delta, q=0; z=\frac{1}{x}, \sigma = \frac{1}{\beta }z, \varsigma =-\frac{1}{\beta } z^2 \right) \nonumber\\
 &=& z^{-2N-2} \sum_{r=0}^{N+1} y_{2r}^{2(N+1-r)}\left( 2N+2;z\right)  \label{eq:70094c}
\end{eqnarray}
\item As $\alpha =-2N-3$, $\gamma =q=0$ where $N \in \mathbb{N}_{0}$, 

Its eigenfunction is given by
\begin{eqnarray}
y(z) &=& _pH_c^{(r,1)}F_{2N+3}^R \left( \alpha =-2N-3, \beta, \gamma =0, \delta, q=0; z=\frac{1}{x}, \sigma = \frac{1}{\beta }z, \varsigma =-\frac{1}{\beta } z^2 \right) \nonumber\\
 &=& z^{-2N-3} \sum_{r=0}^{N+1} y_{r}^{2(N-r)+3}\left( 2N+3;z\right)  \label{eq:70094d}
\end{eqnarray}
\item As $\alpha =1, \gamma =2, q=\beta -\delta $,  

Its eigenfunction is given by
\begin{eqnarray}
y(z) &=& _pH_c^{(r,2)}F_0^R \left( \alpha =1, \beta, \gamma =2, \delta, q=\beta -\delta ; z=\frac{1}{x}, \sigma = \frac{1}{\beta }z, \varsigma =-\frac{1}{\beta } z^2\right) \nonumber\\
&=& z  \label{eq:70094e}
\end{eqnarray}
\item As $\alpha =0, \gamma =2, q=\beta -\delta $,  

Its eigenfunction is given by
\begin{eqnarray}
y(z) &=& _pH_c^{(r,2)}F_1^R \left( \alpha =0, \beta, \gamma =2, \delta, q=\beta -\delta ; z=\frac{1}{x}, \sigma = \frac{1}{\beta }z, \varsigma =-\frac{1}{\beta } z^2\right) \nonumber\\
&=&  1-(\beta -\delta )\sigma  \label{eq:70094f}
\end{eqnarray}
\item As $ \alpha = -2N-1, \gamma =2, q=\beta -\delta $ where $N \in \mathbb{N}_{0}$,

Its eigenfunction is given by
\begin{eqnarray}
y(z) &=& _pH_c^{(r,2)}F_{2N+2}^R \left( \alpha = -2N-1, \beta, \gamma =2, \delta, q=\beta -\delta ; z=\frac{1}{x}, \sigma = \frac{1}{\beta }z, \varsigma =-\frac{1}{\beta } z^2 \right) \nonumber\\
 &=& z^{ -2N-1} \sum_{r=0}^{N+1} y_{2r}^{2(N+1-r)}\left( 2N+2;z\right) \label{eq:70094g}
\end{eqnarray}
\item As $ \alpha = -2N-2, \gamma =2, q=\beta -\delta $ where $N \in \mathbb{N}_{0}$,

Its eigenfunction is given by
\begin{eqnarray}
y(z) &=& _pH_c^{(r,2)}F_{2N+3}^R \left( \alpha =-2N-2, \beta, \gamma =2, \delta, q=\beta -\delta ; z=\frac{1}{x}, \sigma = \frac{1}{\beta }z, \varsigma =-\frac{1}{\beta } z^2 \right) \nonumber\\
 &=& z^{ -2N-2} \sum_{r=0}^{N+1} y_{r}^{2(N-r)+3}\left( 2N+3;z\right) \label{eq:70094h}
\end{eqnarray}
In the above,
\begin{eqnarray}
y_0^m(j;z) &=& \sum_{i_0=0}^{m} \frac{\left( -j\right)_{i_0}\left( \wp_{0}^j\right)_{i_0}}{\left( 1 \right)_{i_0}} \sigma^{i_0} \label{eq:70095a}\\
y_1^m(j;z) &=& \left\{\sum_{i_0=0}^{m} \frac{ \left( i_0 -j\right) \left( i_0 +1-j \right)}{\left( i_0+2\right)} \frac{\left( -j\right)_{i_0} \left( \wp_{0}^j\right)_{i_0}}{\left( 1 \right)_{i_0}} \right.  \left. \sum_{i_1 = i_0}^{m} \frac{\left( 2-j\right)_{i_1} \left( \wp_{1}^j\right)_{i_1}\left( 3\right)_{i_0} }{\left( 2-j\right)_{i_0} \left( \wp_{1}^j\right)_{i_0}\left( 3\right)_{i_1}} \sigma^{i_1}\right\}\varsigma 
\hspace{2cm}\label{eq:70095b}\\
y_{\tau }^m(j;z) &=& \left\{ \sum_{i_0=0}^{m} \frac{ \left( i_0 -j\right) \left( i_0 +1-j\right) }{\left( i_0+2\right)} \frac{\left( -j\right)_{i_0} \left( \wp_{0}^j \right)_{i_0}}{\left( 1 \right)_{i_0}} \right.\nonumber\\
&&\times \prod_{k=1}^{\tau -1} \left( \sum_{i_k = i_{k-1}}^{m} \frac{\left( i_k+2k-j\right)\left( i_k +2k+1-j\right) }{\left( i_k+2k+2\right)} \right.   \left. \frac{\left( 2k-j\right)_{i_k} \left( \wp_{k}^j\right)_{i_{k}}\left( 2k+1 \right)_{i_{k-1}} }{\left( 2k-j\right)_{i_{k-1}} \left( \wp_{k}^j\right)_{i_{k-1}}\left( 2k+1 \right)_{i_k}} \right) \nonumber\\
&&\times  \left. \sum_{i_{\tau } = i_{\tau -1}}^{m} \frac{\left( 2\tau -j\right)_{i_{\tau}} \left( \wp_{\tau }^j\right)_{i_{\tau}}\left( 2\tau +1 \right)_{i_{\tau -1}} }{\left( 2\tau -j \right)_{i_{\tau -1}} \left( \wp_{\tau }^j \right)_{i_{\tau -1}}\left( 2\tau +1 \right)_{i_{\tau }} } \sigma^{i_{\tau }}\right\} \varsigma^{\tau } \label{eq:70095c} 
\end{eqnarray}
where
\begin{equation}
\begin{cases} \tau \geq 2 \cr 
\wp_{k}^j = \beta -\delta +1+2k-j  
\end{cases}\nonumber
\end{equation}
\end{enumerate}
\end{remark}  
\section{Summary}

In chapters 4 \& 5 of Ref.\cite{7Choun2013} and chapter 7, all possible Frobenius solutions of the CHE about the regular singular point at zero are an infinite series, a polynomial of type 1, a polynomial of type 2 and the first species complete polynomial. 
For a polynomial of type 1, I treat $\beta $, $\gamma $, $\delta $, $q$ as free variables and $\alpha $ as a fixed value.
For a polynomial of type 2, I treat $\beta $, $\gamma $, $\delta $, $\alpha $ as free variables and $q$ as a fixed value.
For the first species complete polynomial, I treat $\beta $, $\gamma $, $\delta $ as free variables and $\alpha $, $q$ as fixed values. 
There is no such solution for the second species complete polynomial of the CHE around $x=0$. Because a parameter $\alpha $ of a numerator in $B_n$ term is only a fixed constant in order to make $B_n$ term terminated at a specific index summation $n$. 

By applying $A_n$ and $B_n$ terms, composed of different coefficients in a numerator and a denominator in a recursive relation for the CHE around $x=0$, into 3TRF such as the general summation formulas, an infinite series and a polynomial of type 1 of the CHE are constructed analytically: the denominators and numerators in all $B_n$ terms of each sub-power series arise with Pochhammer symbols in chapter 4 of Ref.\cite{7Choun2013}. 
Combined definite and contour integral forms of an infinite series and a polynomial of type 1 are derived by applying integrals of Confluent hypergeometric functions into the sub-power series of the general series solutions. 
Generating functions for a polynomial of type 1 are constructed analytically by applying generating functions for confluent hypergeometric polynomials into the sub-integral of general integral solutions. 

By applying $A_n$ and $B_n$ terms in a recurrence relation of the CHE around $x=0$ into R3TRF, an infinite series and a polynomial of type 2 of the CHE are constructed analytically: the denominators and numerators in all $A_n$ terms of each sub-power series arise with Pochhammer symbols in chapter 5 of Ref.\cite{7Choun2013}.
Combined definite and contour integral forms of an infinite series and a polynomial of type 2 are derived by applying integrals of Gauss hypergeometric functions into the sub-power series of the general series solutions. 
Generating functions for a polynomial of type 2 are constructed analytically by applying generating functions for Jacobi polynomial using hypergeometric functions into the sub-integral of general integral solutions.

By applying $A_n$ and $B_n$ terms in a recurrence relation  of the CHE around $x=0$ into the first species complete polynomials using 3TRF and R3TRF, the first species complete polynomials of the CHE are constructed analytically. For the first species complete polynomial using 3TRF, the denominators and numerators in all $B_n$ terms of each finite sub-power series arise with Pochhammer symbols in chapter 7. For the first species complete polynomial using R3TRF, the denominators and numerators in all $A_n$ terms of each finite sub-power series arise with Pochhammer symbols.
 
In this chapter, all possible power series solutions of the CHE about the irregular singular point at infinity are the first species complete polynomial and the second species complete polynomial. 
All possible local solutions of the CHE (Regge-Wheeler and Teukolsky equations) were constructed by Fiziev.\cite{7Fizi2009,7Fizi2010} The power series representation of the CHE about the singular point at infinity is not convergent by only asymptotic. \cite{7Ronv1995}  
In general Confluent Heun polynomial (CHP) has been hitherto defined as a type 3 polynomial which is equivalent to the first species complete polynomial. \cite{7Fizi2010a,7Fizi2009,7Fizi2010,7Fizi2010b}  
 However, there is an another type of a polynomial which makes $A_n$ and $B_n$ terms terminated, called as the second species complete polynomial.   
For the first species complete polynomial, I treat $\beta $, $\gamma $, $\delta $ as free variables and $\alpha $, $q$ as fixed values. 
For the second species complete polynomial, I treat $\beta $, $\delta $ as free variables and $\alpha $, $\gamma $, $q$ as fixed values. 
There are no general series solutions for an infinite series, a polynomial of type 1 and a polynomial of type 2 because  series solutions are not convergent any more since  $ { \displaystyle \lim_{n\gg 1} A_n = \lim_{n\gg 1} B_n \rightarrow \infty }$.  

By applying $A_n$ and $B_n$ terms in a recurrence relation of the CHE around $x=\infty $ into a general summation expression of the first species complete polynomial using 3TRF; this is done by letting $A_n$ in sequences $c_n$ is the leading term in each finite sub-power series of the general series solution $y(x)$, I show formal series solutions of the CHE for the first species complete polynomial.   
And the first species complete polynomials of the CHE around $x=\infty $ are constructed by applying a general summation formula of complete polynomials using R3TRF: this is done by letting $B_n$ in sequences $c_n$ is the leading term in each finite sub-power series of the general series solution $y(x)$.
 
By applying $A_n$ and $B_n$ terms in a recurrence relation of the CHE around $x=\infty $ into the second species complete polynomials using 3TRF and R3TRF, the second species complete polynomials of the CHE are constructed analytically. For the second species complete polynomial using 3TRF, the denominators and numerators in all $B_n$ terms of each finite sub-power series arise with Pochhammer symbols. For the second species complete polynomial using R3TRF, the denominators and numerators in all $A_n$ terms of each finite sub-power series arise with Pochhammer symbols.

In the future series, integral forms of the CHE around $x=\infty $ for the first and second species complete polynomials using 3TRF and R3TRF will be obtained by applying integrals of hypergeometric-type functions into the finite sub-power series of general series solutions. And generating functions for these complete polynomials will be constructed by applying generating functions for hypergeometric-type functions into the sub-integral of the general integral solutions.
 
\addcontentsline{toc}{section}{Bibliography}
\bibliographystyle{model1a-num-names}
\bibliography{<your-bib-database>}
\chapter{Complete polynomials of Grand Confluent Hypergeometric equation about the regular singular point at zero}
\chaptermark{Complete polynomials of the GCH around $x=0$} 
 
The classical method of solution in series of a Grand Confluent Hypergeometric (GCH) equation provides a 3-term recurrence relation between successive coefficients. 
By applying three term recurrence formula (3TRF) \cite{8Chou2012b}, I construct power series solutions in closed forms of the GCH equation about the regular singular point at zero for an infinite series and a polynomial of type 1 \cite{8Chou2012i}. And its combined definite and contour integrals involving hypergeometric-type functions are obtained including generating functions for the GCH polynomial of type 1 \cite{8Chou2012j}.

In chapter 6 of Ref.\cite{8Choun2013}, I apply reversible three term recurrence formula (R3TRF) to power series expansions in closed forms of the GCH equation about $x=0$ for an infinite series and a polynomial of type 2. And its representations for solutions as combined definite and contour integrals are constructed analytically including generating functions for the GCH polynomial of type 2. 
 
In this chapter I construct Frobenius solutions of the GCH equation about $x=0$ for a polynomial of type 3 by applying general summation formulas of complete polynomials using 3TRF and R3TRF. 
 
\section{Introduction}
Lichitenberg and collaborators\cite{8Lich1982} observed the semi-relativistic Hamiltonian (the ``Krolikowski'' type second order differential equation \cite{8Krol1980,8Krol1981,8Todo1971}) in order to calculate meson and baryon masses in 1982. From their analysis, G\"{u}rsey \textit{et al.} suggested the spin free Hamiltonian involving only scalar potential in the meson ($q-\bar{q}$) system \cite{81985,81988,81991,81990}. Since they neglected the mass of quark in their supersymmetric wave equation, they noticed that its differential equation is equivalent to confluent hypergeometric equation having a 2-term recurrence relation between successive coefficients in its classical formal series. 

By substituting a power series with unknown coefficients into their spin free wave equation including the mass of quark, I noticed that its recursive relation is composed of a 3-term between consecutive coefficients in a power series solution. Confluent hypergeometric equation is just the special case of this new type differential equation, designated as grand confluent hypergeometric (GCH) equation.

The canonical form of the GCH equation is a second-order linear ordinary differential equation of the form \cite{8Chou2012i}
\begin{equation}
x \frac{d^2{y}}{d{x}^2} + \left( \mu x^2 + \varepsilon x + \nu  \right) \frac{d{y}}{d{x}} + \left( \Omega x + \varepsilon \omega \right) y = 0
\label{eq:8002}
\end{equation}
where $\mu$, $\varepsilon$, $\nu $, $\Omega$ and $\omega$ are real or complex parameters.
GCH equation is of Fuchsian types with two singular points: one regular singular point which is zero with exponents $\{ 0,1-\nu \}$, and one irregular singular point which is infinity with an exponent $\frac{\Omega }{\mu }$.

Since a formal series with unknown coefficients is substituted into the hypergoemetric equation having three regular singular points, a 2-term recursion relation between successive coefficients in its power series solution starts to appear. Confluent types of the hypergeometric equation is derived when two or more singularities coalesce into an irregular singularity. This type includes equations of Legendre, Laguerre, Kummer, Bessel, Jacobi and etc, whose analytic solutions in compact forms are already constructed by great many scholars extensively including their definite or contour integrals.

For Heun equation having four regular singular points, the recurrence relation in its Frobenius solution involves 3 terms. The Heun equation generalizes all well-known equations of Spheroidal Wave, Lame, Mathieu, hypergeometric type and etc. 
Until now, its series solutions in which coefficients are given fully and clearly have been unknown because of its complex mathematical computations; their numerical calculations are still ambiguous. Of course, for definite or contour integrals of Heun equation, no analytic solutions have been constructed regrettably.    

According to Karl Heun \cite{8Heun1889,8Ronv1995}, the canonical form of general Heun's equation is taken as
\begin{equation}
\frac{d^2{y}}{d{x}^2} + \left(\frac{\gamma }{x} +\frac{\delta }{x-1} + \frac{\epsilon }{x-a}\right) \frac{d{y}}{d{x}} +  \frac{\alpha \beta x-q}{x(x-1)(x-a)} y = 0 \label{eq:8001}
\end{equation}
where $\epsilon = \alpha +\beta -\gamma -\delta +1$ for assuring the regularity of the point at $x=\infty $. It has four regular singular points which are 0, 1, $a$ and $\infty $ with exponents $\{ 0, 1-\gamma \}$, $\{ 0, 1-\delta \}$, $\{ 0, 1-\epsilon \}$ and $\{ \alpha, \beta \}$. 
 
Like deriving of confluent hypergeometric equation from the hypergeometric equation, 4 confluent types of Heun equation can be derived from merging two or more regular singularities to take an irregular singularity in Heun equation. These types include such as (1) Confluent Heun (two regular and one irregular singularities), (2) Doubly Confluent Heun (two irregular singularities), (3) Biconfluent Heun (one regular and one irregular singularities), (4) Triconfluent Heun equations (one irregular singularity).

The GCH equation generalizes Biconfluent Heun (BCH) equation having a regular singularity at $x=0$ and an irregular singularity at $\infty$ of rank 2. For the canonical form of the BCH equation\cite{8Ronv1995}, replace $\mu $, $\varepsilon $, $\nu $, $\Omega $ and $\omega $ by $-2$, $-\beta  $, $ 1+\alpha $, $\gamma -\alpha -2 $ and $ 1/2 (\delta /\beta +1+\alpha )$ in (\ref{eq:8002}). For DLFM version \cite{8NIST} or in Ref.\cite{8Slavy2000}, replace $\mu $ and $\omega $ by 1 and $-q/\varepsilon $ in (\ref{eq:8002}).  

The BCH equation, the special case of the GCH equation, is one of the center of attention in mathematics \cite{8Bato1977,8Belm2004,8Cheb2004,8Deca1978,8Decar1978,8Exto1989,8Haut1969,8Maro1967,8Tric1960,8Urwi1975} and modern physics \cite{8Cift2010,8Figu2012,8Fles1980,8Fles1982,8Leau1986,8Leau1990,8Mass1983,8Pons1988,8Ralk2002,8Truo2000}.
For instance, from a quantum mechanical point of view, the BCH equation arises in the radial Schr$\ddot{\mbox{o}}$dinger equation  of the new classes of solvable potentials for various harmonic oscillators \cite{8Chau1983,8Chau1984,8Fles1980,8Fles1982,8Leau1986,8Leau1990,8Lemi1969,8Mass1983,8Pons1988}, the relativistic quantum mechanics  
in a uniform magnetic field and scalar potentials in the cosmic string space-time \cite{8Figu2012}, quantum mechanical wave functions of two charged particles moving on a plane with a perpendicular uniform magnetic field \cite{8Ralk2002}, non-relativistic eigenfunctions for the Schr$\ddot{\mbox{o}}$dinger equation by two electrons interacting by a Coulomb potential and anharmonic oscillator potential \cite{8Caru2014}, etc.  
 
Assume that a formal series solution of the GCH equation takes the form 
\begin{equation}
y(x)= \sum_{n=0}^{\infty } c_n x^{n+\lambda } \hspace{1cm}\mbox{where}\; \lambda =\mbox{indicial}\;\mbox{root}\label{eq:8003}
\end{equation}
Substituting (\ref{eq:8003}) into (\ref{eq:8002}) gives for the coefficients $c_n$ the recurrence relations
\begin{equation}
c_{n+1}=A_n \;c_n +B_n \;c_{n-1} \hspace{1cm};n\geq 1 \label{eq:8004}
\end{equation}
where,
\begin{subequations}
\begin{equation}
A_n = - \frac{\varepsilon (n+\omega +\lambda ) }{ (n+1+\lambda )(n+\nu +\lambda)} \label{eq:8005a}
\end{equation}
\begin{equation}
B_n = - \frac{\mu \left( n-1 + \frac{\Omega }{\mu }+\lambda \right) }{ (n+1+\lambda )(n+\nu +\lambda)} \label{eq:8005b}
\end{equation}
\end{subequations}
where $c_1= A_0 \;c_0$. Two indicial roots are given such as $\lambda = 0$ and $ 1-\nu $.

There are 4 possible formal series solutions of a linear ODE having a 3-term recursive relation between successive coefficients such as an infinite series and 3 types of polynomials: (1) a polynomial which makes $B_n$ term terminated; $A_n$ term is not terminated, (2) a polynomial which makes $A_n$ term terminated; $B_n$ term is not terminated, (3) a polynomial which makes $A_n$ and $B_n$ terms terminated at the same time, referred as `a complete polynomial.' 

Complete polynomials have two different types such as the first species complete polynomial and the second species complete polynomial. The first species complete polynomial is applicable since a parameter of a numerator in $B_n$ term and a (spectral) parameter of a numerator in $A_n$ term are fixed constants. The second species complete polynomial is utilizable since two parameters of a numerator in $B_n$ term and a parameter of a numerator in $A_n$ term are fixed constants.
The former has multi-valued roots of a parameter of a numerator in $A_n$ term, but the latter has only one fixed value of a parameter of a numerator in $A_n$ term for an eigenvalue.  
 
\begin{table}[h]
\begin{center}
\thispagestyle{plain}
\hspace*{-0.1\linewidth}\resizebox{1.2\linewidth}{!}
{
 \Tree[.{\Huge The GCH differential equation about the regular singular point at zero} [.{\Huge 3TRF} [.{\Huge Infinite series} ]
              [.{\Huge Polynomials} [[.{\Huge Polynomial of type 1} ]
               [.{\Huge Polynomial of type 3} [.{\Huge $ \begin{array}{lcll}  1^{\mbox{st}}\;  \mbox{species}\\ \mbox{complete} \\ \mbox{polynomial} \end{array}$} ]  ]]]]                         
  [.{\Huge R3TRF} [.{\Huge Infinite series} ]
     [.{\Huge Polynomials} [[.{\Huge Polynomial of type 2} ]
       [.{\Huge  Polynomial of type 3} [.{\Huge $ \begin{array}{lcll}  1^{\mbox{st}}\;  \mbox{species} \\ \mbox{complete} \\ \mbox{polynomial} \end{array}$} ]  ]]]]]
}
\end{center}
\caption{Power series of the GCH equation about the regular singular point at zero}
\end{table}  
Table 9.1 tells us about all possible general power series solutions of the GCH equation about the regular singular point at zero. 
With my definition, there are 2 types of the general summation formulas in order to derive formal series solutions in closed forms of the GCH equation such as 3TRF and R3TRF. 

For $n=0,1,2,3,\cdots $, (\ref{eq:8004}) is expanded to combinations of $A_n$ and $B_n$ terms. I define that a sub-power series $y_l(x)$ where $l\in \mathbb{N}_{0}$ is constructed by observing the term of sequence $c_n$ which includes $l$ terms of $A_n's$. The general series solution is delineated by sums of each $y_l(x)$ such as $y(x)= \sum_{n=0}^{\infty } y_n(x)$. By allowing $A_n$ in the sequence $c_n$ is the leading term of each sub-power series $y_l(x)$, the general summation formulas of the 3-term recurrence relation in a linear ODE are constructed for an infinite series and a polynomial of type 1, denominated as `three term recurrence formula (3TRF).' \cite{8Chou2012b}
 
Frobenius solutions of the GCH equation around $x=0$ are expressed analytically for an infinite series and a polynomial of type 1 by applying 3TRF. 
For a polynomial of type 1, I treat $\mu $, $\varepsilon  $, $\nu $, $\omega $ as free variables and $\Omega $ as a fixed value. 
General integral representations of the GCH equation are constructed by applying the combined definite and contour integral of $_2F_2$ function into each of its sub-power series solutions: a $_1F_1$ function (the Kummer function of the first kind) recurs in each of sub-integral forms of the GCH equation.   
By applying generating functions for Kummer's polynomials of the first kind into each of sub-integrals of the general integral form for the GCH polynomial of type 1, generating functions of the GCH polynomial of type 1 are obtained analytically. 
 
In chapter 1 of Ref.\cite{8Choun2013}, the general solution in a formal series is described by sums of each $y_l(x)$. A sub-power series $y_l(x)$ is obtained by observing the term of sequence $c_n$ which includes $l$ terms of $B_n's$. 
By allowing $B_n$ in the sequence $c_n$ is the leading term of each sub-power series $y_l(x)$, the general summation expressions of the 3-term recurrence relation in a linear ODE are constructed for an infinite series and a polynomial of type 2, designated as `reversible three term recurrence formula (R3TRF).' 

In chapter 6 of Ref.\cite{8Choun2013}, power series solutions in closed forms of the GCH equation around $x=0$ are constructed for an infinite series and a polynomial of type 2 by applying R3TRF. For a polynomial of type 2, I treat $\mu $, $\varepsilon  $, $\nu $, $\Omega $ as free variables and $\omega $ as a fixed value. Integrals of the GCH equation are obtained by applying an integral of $_2F_2$ function into each of its sub-formal series. And generating functions of the GCH polynomial of type 2 are constructed by applying generating functions for Kummer's polynomials of the first kind into each of sub-integrals of the general integral form for the GCH polynomial of type 2.
 
Infinite series of the GCH equation around $x=0$ by applying 3TRF are tantamount to infinite series by applying R3TRF for power series solutions and integrals. The former is that $A_n$ is the leading term in each sub-power series of the GCH equation. The latter is that $B_n$ is the leading term in each sub-power series of it. 
  
In chapter 1, by allowing $A_n$ as the leading term in each of finite sub-power series of a formal series $y(x)$, general summation expressions for the first and second species complete polynomials are obtained in closed forms. I refer these mathematical formulas as ``complete polynomials using 3-term recurrence formula (3TRF).''
In chapter 2,  by allowing $B_n$ as the leading term in each finite sub-power series of the general power series $y(x)$, I construct general solutions in series, built in compact forms, for the first and second species complete polynomials. I designate these classical summation formulas as ``complete polynomials using reversible 3-term recurrence formula (R3TRF).''

In this chapter, by substituting (\ref{eq:8005a}) and (\ref{eq:8005b}) into complete polynomials using 3TRF and R3TRF, the construction of solutions in series of the GCH equation around $x=0$ is given for the first species complete polynomials.
Indeed, I show a polynomial equation of the GCH equation for the determination of a parameter $\omega $ in the combinational form of partial sums of the sequences $\{A_n\}$ and $\{B_n\}$ using 3TRF and R3TRF. 
    
\section{Power series about the regular singular point at zero}
The spectral polynomial (the first species complete polynomial) of the BCH equation around $x=0$, heretofore has been obtained by applying a power series with unknown coefficients \cite{8Bato1977,8Deca1978,8Haut1969,8Urwi1975}. 
The BCH polynomial comes from the BCH equation that has a fixed integer value of $\gamma -\alpha -2$, just as it has a fixed value of $\delta $. They left a polynomial equation of the $(j+1)$th order for the determination of a parameter $\delta $ as the determinant of $(j+1)\times (j+1)$ matrices \cite{8Caru2014,8Cift2010}. For the GCH polynomial in the canonical form,  $\gamma -\alpha -2$ and  $\delta $ correspond to $\Omega $ and $\varepsilon \left( \mu \omega + \nu \right)$.
They just left an analytic solution of the BCH polynomial as solutions of recurrences because of a 3-term recursive relation between successive coefficients in its solution in series. It is still unknown to find power series solutions in which the coefficients are given explicitly for a 3-term recursion relation.

For the first species complete polynomial of the GCH equation around $x=0$ in table 9.1, I treat $\varepsilon $, $\mu $, $\nu  $ as free variables and $\omega $, $\Omega $ as fixed values. There is no such solution for the second species complete polynomial of the GCH equation. Because a parameter $\Omega $ of a numerator in $B_n$ term in (\ref{eq:8005b}) is only a fixed constant in order to make $B_n$ term terminated at a specific index summation $n$ for two independent solutions where $\lambda = 0$ and $ 1-\nu $.  
\subsection{The first species complete polynomial of the GCH equation using 3TRF}

For the first species complete polynomials using 3TRF and R3TRF, we need a condition which is given by
\begin{equation}
B_{j+1}= c_{j+1}=0\hspace{1cm}\mathrm{where}\;j\in \mathbb{N}_{0}  
 \label{eq:8006}
\end{equation}
(\ref{eq:8006}) gives successively $c_{j+2}=c_{j+3}=c_{j+4}=\cdots=0$. And $c_{j+1}=0$ is defined by a polynomial equation of degree $j+1$ for the determination of an accessory parameter in $A_n$ term. 
\begin{theorem}
In chapter 1, the general summation expression of a function $y(x)$ for the first species complete polynomial using 3-term recurrence formula and its algebraic equation for the determination of an accessory parameter in $A_n$ term are given by
\begin{enumerate} 
\item As $B_1=0$,
\begin{equation}
0 =\bar{c}(1,0) \label{eq:8007a}
\end{equation}
\begin{equation}
y(x) = y_{0}^{0}(x) \label{eq:8007b}
\end{equation}
\item As $B_{2N+2}=0$ where $N \in \mathbb{N}_{0}$,
\begin{equation}
0  = \sum_{r=0}^{N+1}\bar{c}\left( 2r, N+1-r\right) \label{eq:8008a}
\end{equation}
\begin{equation}
y(x)= \sum_{r=0}^{N} y_{2r}^{N-r}(x)+ \sum_{r=0}^{N} y_{2r+1}^{N-r}(x)  \label{eq:8008b}
\end{equation}
\item As $B_{2N+3}=0$ where $N \in \mathbb{N}_{0}$,
\begin{equation}
0  = \sum_{r=0}^{N+1}\bar{c}\left( 2r+1, N+1-r\right) \label{eq:8009a}
\end{equation}
\begin{equation}
y(x)= \sum_{r=0}^{N+1} y_{2r}^{N+1-r}(x)+ \sum_{r=0}^{N} y_{2r+1}^{N-r}(x)  \label{eq:8009b}
\end{equation}
In the above,
\begin{eqnarray}
\bar{c}(0,n)  &=& \prod _{i_0=0}^{n-1}B_{2i_0+1} \label{eq:80010a}\\
\bar{c}(1,n) &=&  \sum_{i_0=0}^{n} \left\{ A_{2i_0} \prod _{i_1=0}^{i_0-1}B_{2i_1+1} \prod _{i_2=i_0}^{n-1}B_{2i_2+2} \right\} 
\label{eq:80010b}\\
\bar{c}(\tau ,n) &=& \sum_{i_0=0}^{n} \left\{A_{2i_0}\prod _{i_1=0}^{i_0-1} B_{2i_1+1} 
\prod _{k=1}^{\tau -1} \left( \sum_{i_{2k}= i_{2(k-1)}}^{n} A_{2i_{2k}+k}\prod _{i_{2k+1}=i_{2(k-1)}}^{i_{2k}-1}B_{2i_{2k+1}+(k+1)}\right) \right. \nonumber\\
&&\times \left. \prod _{i_{2\tau}=i_{2(\tau -1)}}^{n-1} B_{2i_{2\tau }+(\tau +1)} \right\} 
 \label{eq:80010c} 
\end{eqnarray}
and
\begin{eqnarray}
y_0^m(x) &=& c_0 x^{\lambda } \sum_{i_0=0}^{m} \left\{ \prod _{i_1=0}^{i_0-1}B_{2i_1+1} \right\} x^{2i_0 } \label{eq:80011a}\\
y_1^m(x) &=& c_0 x^{\lambda } \sum_{i_0=0}^{m}\left\{ A_{2i_0} \prod _{i_1=0}^{i_0-1}B_{2i_1+1}  \sum_{i_2=i_0}^{m} \left\{ \prod _{i_3=i_0}^{i_2-1}B_{2i_3+2} \right\}\right\} x^{2i_2+1 } \label{eq:80011b}\\
y_{\tau }^m(x) &=& c_0 x^{\lambda } \sum_{i_0=0}^{m} \left\{A_{2i_0}\prod _{i_1=0}^{i_0-1} B_{2i_1+1} 
\prod _{k=1}^{\tau -1} \left( \sum_{i_{2k}= i_{2(k-1)}}^{m} A_{2i_{2k}+k}\prod _{i_{2k+1}=i_{2(k-1)}}^{i_{2k}-1}B_{2i_{2k+1}+(k+1)}\right) \right. \nonumber\\
&& \times  \left. \sum_{i_{2\tau} = i_{2(\tau -1)}}^{m} \left( \prod _{i_{2\tau +1}=i_{2(\tau -1)}}^{i_{2\tau}-1} B_{2i_{2\tau +1}+(\tau +1)} \right) \right\} x^{2i_{2\tau}+\tau }\hspace{1cm}\mathrm{where}\;\tau \geq 2
\label{eq:80011c} 
\end{eqnarray}
\end{enumerate}
\end{theorem}
Put $n= j+1$ in (\ref{eq:8005b}) and use the condition $B_{j+1}=0$ for $\Omega  $.  
\begin{equation}
\Omega  = -\mu \left( j+\lambda \right)
\label{eq:80012}
\end{equation}
Take (\ref{eq:80012}) into (\ref{eq:8005b}).
\begin{equation}
B_n = -\frac{\mu (n-j-1)}{(n+1+\lambda )(n+\nu +\lambda )} \label{eq:80013}
\end{equation}
Now the condition $c_{j+1}=0$ is clearly an algebraic equation in $\omega $ of degree $j+1$ and thus has $j+1$ zeros denoted them by $\omega _j^m$ eigenvalues where $m = 0,1,2, \cdots, j$. They can be arranged in the following order: $\omega _j^0 < \omega _j^1 < \omega _j^2 < \cdots < \omega _j^j$.

Substitute (\ref{eq:8005a}) and (\ref{eq:80013}) into (\ref{eq:80010a})--(\ref{eq:80011c}).

As $B_{1}= c_{1}=0$, take the new (\ref{eq:80010b}) into (\ref{eq:8007a}) putting $j=0$. Substitute the new (\ref{eq:80011a}) into (\ref{eq:8007b}) putting $j=0$. 

As $B_{2N+2}= c_{2N+2}=0$, take the new (\ref{eq:80010a})--(\ref{eq:80010c}) into (\ref{eq:8008a}) putting $j=2N+1$. Substitute the new 
(\ref{eq:80011a})--(\ref{eq:80011c}) into (\ref{eq:8008b}) putting $j=2N+1$ and $\omega =\omega _{2N+1}^m$.

As $B_{2N+3}= c_{2N+3}=0$, take the new (\ref{eq:80010a})--(\ref{eq:80010c}) into (\ref{eq:8009a}) putting $j=2N+2$. Substitute the new 
(\ref{eq:80011a})--(\ref{eq:80011c}) into (\ref{eq:8009b}) putting $j=2N+2$ and $\omega =\omega _{2N+2}^m$.

After the replacement process, the general expression of power series of the GCH equation about $x=0$ for the first species complete polynomial using 3-term recurrence formula and its algebraic equation for the determination of an accessory parameter $\omega $ are given by
\begin{enumerate} 
\item As $\Omega =-\mu \lambda $,

An algebraic equation of degree 1 for the determination of $\omega $ is given by
\begin{equation}
0= \bar{c}(1,0;0,\omega )= \omega + \lambda \label{eq:80014a}
\end{equation}
The eigenvalue of $\omega $ is written by $\omega _0^0$. Its eigenfunction is given by
\begin{equation}
y(x) = y_0^0\left( 0,\omega _0^0;x\right)= c_0 x^{\lambda } \label{eq:80014b}
\end{equation}
\item As $\Omega  =-\mu \left( 2N+1+\lambda \right)$ where $N \in \mathbb{N}_{0}$,

An algebraic equation of degree $2N+2$ for the determination of $\omega $ is given by
\begin{equation}
0 = \sum_{r=0}^{N+1}\bar{c}\left( 2r, N+1-r; 2N+1,\omega \right)  \label{eq:80015a}
\end{equation}
The eigenvalue of $\omega $ is written by $\omega _{2N+1}^m$ where $m = 0,1,2,\cdots,2N+1 $; $\omega _{2N+1}^0 < \omega _{2N+1}^1 < \cdots < \omega _{2N+1}^{2N+1}$. Its eigenfunction is given by 
\begin{equation} 
y(x) = \sum_{r=0}^{N} y_{2r}^{N-r}\left( 2N+1,\omega _{2N+1}^m;x\right)+ \sum_{r=0}^{N} y_{2r+1}^{N-r}\left( 2N+1,\omega _{2N+1}^m;x\right)
\label{eq:80015b} 
\end{equation}
\item As $\Omega  =-\mu \left( 2N+2+\lambda \right) $ where $N \in \mathbb{N}_{0}$,

An algebraic equation of degree $2N+3$ for the determination of $\omega $ is given by
\begin{equation}  
0 = \sum_{r=0}^{N+1}\bar{c}\left( 2r+1, N+1-r; 2N+2,\omega \right) \label{eq:80016a}
\end{equation}
The eigenvalue of $\omega $ is written by $\omega _{2N+2}^m$ where $m = 0,1,2,\cdots,2N+2 $; $\omega _{2N+2}^0 < \omega _{2N+2}^1 < \cdots < \omega _{2N+2}^{2N+2}$. Its eigenfunction is given by
\begin{equation} 
y(x) =  \sum_{r=0}^{N+1} y_{2r}^{N+1-r}\left( 2N+2,\omega _{2N+2}^m;x\right) + \sum_{r=0}^{N} y_{2r+1}^{N-r}\left( 2N+2,\omega _{2N+2}^m;x\right) \label{eq:80016b}
\end{equation}
In the above,
\begin{eqnarray}
\bar{c}(0,n;j,\omega )  &=& \frac{\left( -\frac{j}{2}\right)_{n}}{\left( 1+\frac{\lambda }{2}\right)_{n} \left(  \frac{1}{2}+ \frac{\nu }{2}+\frac{\lambda }{2}\right)_{n}} \left( -\frac{1}{2}\mu \right)^{n}\label{eq:80017a}\\
\bar{c}(1,n;j,\omega ) &=& \left( -\frac{1}{2}\varepsilon \right) \sum_{i_0=0}^{n}\frac{\left( i_0+\frac{\omega }{2}+\frac{\lambda }{2}\right) }{\left( i_0+\frac{1}{2}+\frac{\lambda }{2}\right) \left( i_0+\frac{\nu }{2}+\frac{\lambda }{2}\right)} \frac{\left( -\frac{j}{2}\right)_{i_0} }{\left( 1+\frac{\lambda }{2}\right)_{i_0} \left( \frac{1}{2}+ \frac{\nu }{2}+ \frac{\lambda }{2}\right)_{i_0}} \nonumber\\
&&\times \frac{\left( \frac{1}{2}-\frac{j}{2} \right)_{n} \left( \frac{3}{2}+\frac{\lambda }{2}\right)_{i_0} \left( 1+\frac{\nu }{2}+ \frac{\lambda }{2}\right)_{i_0}}{\left( \frac{1}{2}-\frac{j}{2}\right)_{i_0} \left( \frac{3}{2}+\frac{\lambda }{2}\right)_{n} \left( 1+\frac{\nu }{2}+ \frac{\lambda }{2}\right)_{n}} \left( -\frac{1}{2}\mu \right)^{n }  
\label{eq:80017b}\\
\bar{c}(\tau ,n;j,\omega ) &=& \left( -\frac{1}{2}\varepsilon \right)^{\tau } \sum_{i_0=0}^{n}\frac{\left( i_0+\frac{\omega }{2}+\frac{\lambda }{2}\right)}{\left( i_0+\frac{1}{2}+\frac{\lambda }{2}\right) \left( i_0+\frac{\nu }{2}+\frac{\lambda }{2}\right)} \frac{\left( -\frac{j}{2}\right)_{i_0} }{\left( 1+\frac{\lambda }{2}\right)_{i_0} \left( \frac{1}{2}+\frac{\nu }{2}+  \frac{\lambda }{2}\right)_{i_0}}  \nonumber\\
&&\times \prod_{k=1}^{\tau -1} \left( \sum_{i_k = i_{k-1}}^{n} \frac{\left( i_k+ \frac{k}{2}+\frac{\omega }{2}+\frac{\lambda }{2}\right) }{\left( i_k+\frac{k}{2}+\frac{1}{2}+\frac{\lambda }{2}\right) \left( i_k+\frac{k}{2}+\frac{\nu }{2}+\frac{\lambda }{2}\right)} \right. \nonumber\\
&&\times \left. \frac{\left( \frac{k}{2}-\frac{j}{2}\right)_{i_k} \left( \frac{k}{2}+1+ \frac{\lambda }{2}\right)_{i_{k-1}} \left( \frac{k}{2}+ \frac{1}{2}+\frac{\nu }{2}+ \frac{\lambda }{2}\right)_{i_{k-1}}}{\left( \frac{k}{2}-\frac{j}{2}\right)_{i_{k-1}} \left( \frac{k}{2}+1+ \frac{\lambda }{2}\right)_{i_k} \left( \frac{k}{2}+\frac{1}{2}+ \frac{\nu }{2}+  \frac{\lambda }{2}\right)_{i_k}} \right) \nonumber\\ 
&&\times \frac{\left( \frac{\tau }{2} -\frac{j}{2}\right)_{n} \left( \frac{\tau }{2}+1+\frac{\lambda }{2}\right)_{i_{\tau -1}} \left( \frac{\tau }{2}+\frac{1}{2}+\frac{\nu }{2}+ \frac{\lambda }{2}\right)_{i_{\tau -1}}}{\left( \frac{\tau }{2}-\frac{j}{2}\right)_{i_{\tau -1}} \left( \frac{\tau }{2}+1+\frac{\lambda }{2}\right)_{n} \left( \frac{\tau }{2}+\frac{1}{2}+ \frac{\nu }{2}+ \frac{\lambda }{2}\right)_{n}} \left( -\frac{1}{2}\mu \right)^{n } \label{eq:80017c} 
\end{eqnarray}
\begin{eqnarray}
y_0^m(j,\omega ;x) &=& c_0 x^{\lambda }  \sum_{i_0=0}^{m} \frac{\left( -\frac{j}{2}\right)_{i_0}}{\left( 1+\frac{\lambda }{2}\right)_{i_0} \left( \frac{1}{2}+ \frac{\nu }{2}+ \frac{\lambda }{2}\right)_{i_0}} \tilde{\mu }^{i_0} \label{eq:80018a}\\
y_1^m(j,\omega ;x) &=& c_0 x^{\lambda } \left\{\sum_{i_0=0}^{m} \frac{\left( i_0+\frac{\omega }{2}+\frac{\lambda }{2}\right) }{\left( i_0+\frac{1}{2}+\frac{\lambda }{2}\right) \left( i_0+\frac{\nu }{2}+\frac{\lambda }{2}\right)} \frac{\left( -\frac{j}{2}\right)_{i_0} }{\left( 1+\frac{\lambda }{2}\right)_{i_0} \left( \frac{1}{2}+\frac{\nu }{2}+ \frac{\lambda }{2}\right)_{i_0}} \right. \nonumber\\
&&\times \left. \sum_{i_1 = i_0}^{m} \frac{\left( \frac{1}{2}-\frac{j}{2} \right)_{i_1} \left( \frac{3}{2}+\frac{\lambda }{2}\right)_{i_0} \left( 1+\frac{\nu }{2}+ \frac{\lambda }{2}\right)_{i_0}}{\left( \frac{1}{2}-\frac{j}{2} \right)_{i_0} \left( \frac{3}{2}+\frac{\lambda }{2}\right)_{i_1} \left( 1+\frac{\nu }{2} + \frac{\lambda }{2}\right)_{i_1}} \tilde{\mu }^{i_1}\right\} \tilde{\varepsilon } 
\label{eq:80018b}\\
y_{\tau }^m(j,\omega ;x) &=& c_0 x^{\lambda } \left\{ \sum_{i_0=0}^{m} \frac{\left( i_0+\frac{\omega }{2}+\frac{\lambda }{2}\right)}{\left( i_0+\frac{1}{2}+\frac{\lambda }{2}\right) \left( i_0+\frac{\nu }{2}+\frac{\lambda }{2}\right)} \frac{\left( -\frac{j}{2}\right)_{i_0} }{\left( 1+\frac{\lambda }{2}\right)_{i_0} \left( \frac{1}{2}+\frac{\nu }{2}+ \frac{\lambda }{2}\right)_{i_0}} \right.\nonumber\\
&&\times \prod_{k=1}^{\tau -1} \left( \sum_{i_k = i_{k-1}}^{m} \frac{\left( i_k+ \frac{k}{2}+\frac{\omega }{2}+\frac{\lambda }{2}\right) }{\left( i_k+\frac{k}{2}+\frac{1}{2}+\frac{\lambda }{2}\right) \left( i_k+\frac{k}{2}+\frac{\nu }{2}+\frac{\lambda }{2}\right)} \right. \nonumber\\
&&\times \left. \frac{\left( \frac{k}{2}-\frac{j}{2}\right)_{i_k} \left( \frac{k}{2}+1+ \frac{\lambda }{2}\right)_{i_{k-1}} \left(  \frac{k}{2}+\frac{1}{2}+\frac{\nu }{2}+ \frac{\lambda }{2}\right)_{i_{k-1}}}{\left( \frac{k}{2}-\frac{j}{2}\right)_{i_{k-1}} \left( \frac{k}{2}+1+ \frac{\lambda }{2}\right)_{i_k} \left( \frac{k}{2}+\frac{1}{2}+\frac{\nu }{2}+ \frac{\lambda }{2}\right)_{i_k}} \right) \nonumber\\
&&\times \left. \sum_{i_{\tau } = i_{\tau -1}}^{m} \frac{\left( \frac{\tau }{2}-\frac{j}{2}\right)_{i_{\tau }} \left( \frac{\tau }{2}+1+\frac{\lambda }{2}\right)_{i_{\tau -1}} \left( \frac{\tau }{2}+\frac{1}{2}+\frac{\nu }{2}+ \frac{\lambda }{2}\right)_{i_{\tau -1}}}{\left( \frac{\tau }{2}-\frac{j}{2}\right)_{i_{\tau -1}} \left( \frac{\tau }{2}+1+\frac{\lambda }{2}\right)_{i_{\tau }} \left( \frac{\tau }{2}+\frac{1}{2}+\frac{\nu }{2}+ \frac{\lambda }{2}\right)_{i_{\tau }}} \tilde{\mu }^{i_{\tau }}\right\} \tilde{\varepsilon }^{\tau } \hspace{1.5cm}  \label{eq:80018c} 
\end{eqnarray}
where
\begin{equation}
\begin{cases} \tau \geq 2 \cr
\tilde{\mu }  = -\frac{1}{2}\mu x^2 \cr
\tilde{\varepsilon } = -\frac{1}{2}\varepsilon x 
\end{cases}\nonumber
\end{equation}
\end{enumerate}
Put $c_0$= 1 as $\lambda =0$ for the first kind of independent solutions of the GCH equation and $\lambda = 1-\nu $ for the second one in (\ref{eq:80014a})--(\ref{eq:80018c}). 
\begin{remark}
The power series expansion of the GCH equation of the first kind for the first species complete polynomial using 3TRF about $x=0$ is given by
\begin{enumerate} 
\item As $\Omega =0$ and $\omega = \omega _0^0 =0$,

The eigenfunction is given by
\begin{equation}
y(x) = Q_pW_{0,0}\left( \mu ,\varepsilon ,\nu ,\Omega =0, \omega = \omega _0^0 =0; \tilde{\mu }  = -\frac{1}{2}\mu x^2; \tilde{\varepsilon } = -\frac{1}{2}\varepsilon x \right) =1 \label{eq:80019}
\end{equation}
\item As $\Omega  =-\mu \left( 2N+1 \right)$ where $N \in \mathbb{N}_{0}$,

An algebraic equation of degree $2N+2$ for the determination of $\omega $ is given by
\begin{equation}
0 = \sum_{r=0}^{N+1}\bar{c}\left( 2r, N+1-r; 2N+1,\omega \right)\label{eq:80020a}
\end{equation}
The eigenvalue of $\omega $ is written by $\omega _{2N+1}^m$ where $m = 0,1,2,\cdots,2N+1 $; $\omega _{2N+1}^0 < \omega _{2N+1}^1 < \cdots < \omega _{2N+1}^{2N+1}$. Its eigenfunction is given by
\begin{eqnarray} 
y(x) &=&  Q_pW_{2N+1,m} \left( \mu ,\varepsilon ,\nu ,\Omega =-\mu \left( 2N+1 \right), \omega = \omega _{2N+1}^m; \tilde{\mu }  = -\frac{1}{2}\mu x^2; \tilde{\varepsilon } = -\frac{1}{2}\varepsilon x \right)\nonumber\\
&=& \sum_{r=0}^{N} y_{2r}^{N-r}\left( 2N+1, \omega _{2N+1}^m;x\right)+ \sum_{r=0}^{N} y_{2r+1}^{N-r}\left( 2N+1, \omega _{2N+1}^m;x\right)  
\label{eq:80020b}
\end{eqnarray}
\item As $\Omega  =-\mu \left( 2N+2 \right)$ where $N \in \mathbb{N}_{0}$,

An algebraic equation of degree $2N+3$ for the determination of $\omega $ is given by
\begin{eqnarray}
0  = \sum_{r=0}^{N+1}\bar{c}\left( 2r+1, N+1-r; 2N+2,\omega \right)\label{eq:80021a}
\end{eqnarray}
The eigenvalue of $\omega $ is written by $\omega _{2N+2}^m$ where $m = 0,1,2,\cdots,2N+2 $; $\omega _{2N+2}^0 < \omega _{2N+2}^1 < \cdots < \omega _{2N+2}^{2N+2}$. Its eigenfunction is given by
\begin{eqnarray} 
y(x) &=& Q_pW_{2N+2,m} \left( \mu ,\varepsilon ,\nu ,\Omega =-\mu \left( 2N+2\right), \omega = \omega _{2N+2}^m; \tilde{\mu }  = -\frac{1}{2}\mu x^2; \tilde{\varepsilon } = -\frac{1}{2}\varepsilon x \right) \nonumber\\
&=& \sum_{r=0}^{N+1} y_{2r}^{N+1-r}\left( 2N+2, \omega _{2N+2}^m;x\right) + \sum_{r=0}^{N} y_{2r+1}^{N-r}\left( 2N+2, \omega _{2N+2}^m;x\right) 
\label{eq:80021b}
\end{eqnarray}
In the above,
\begin{eqnarray}
\bar{c}(0,n;j,\omega )  &=& \frac{\left( -\frac{j}{2}\right)_{n}}{\left( 1 \right)_{n} \left(  \frac{1}{2}+ \frac{\nu }{2} \right)_{n}} \left( -\frac{1}{2}\mu \right)^{n}\label{eq:80022a}\\
\bar{c}(1,n;j,\omega ) &=& \left( -\frac{1}{2}\varepsilon \right) \sum_{i_0=0}^{n}\frac{\left( i_0+\frac{\omega }{2} \right) }{\left( i_0+\frac{1}{2} \right) \left( i_0+\frac{\nu }{2} \right)} \frac{\left( -\frac{j}{2}\right)_{i_0} }{\left( 1 \right)_{i_0} \left( \frac{1}{2}+ \frac{\nu }{2} \right)_{i_0}}   \frac{\left( \frac{1}{2}-\frac{j}{2} \right)_{n} \left( \frac{3}{2} \right)_{i_0} \left( 1+\frac{\nu }{2} \right)_{i_0}}{\left( \frac{1}{2}-\frac{j}{2}\right)_{i_0} \left( \frac{3}{2} \right)_{n} \left( 1+\frac{\nu }{2} \right)_{n}} \left( -\frac{1}{2}\mu \right)^{n } \hspace{1.5cm} \label{eq:80022b}\\
\bar{c}(\tau ,n;j,\omega ) &=& \left( -\frac{1}{2}\varepsilon \right)^{\tau } \sum_{i_0=0}^{n}\frac{\left( i_0+\frac{\omega }{2} \right)}{\left( i_0+\frac{1}{2} \right) \left( i_0+\frac{\nu }{2} \right)} \frac{\left( -\frac{j}{2}\right)_{i_0} }{\left( 1 \right)_{i_0} \left( \frac{1}{2}+\frac{\nu }{2} \right)_{i_0}}  \nonumber\\
&&\times \prod_{k=1}^{\tau -1} \left( \sum_{i_k = i_{k-1}}^{n} \frac{\left( i_k+ \frac{k}{2}+\frac{\omega }{2} \right) }{\left( i_k+\frac{k}{2}+\frac{1}{2} \right) \left( i_k+\frac{k}{2}+\frac{\nu }{2} \right)} \right.  \left. \frac{\left( \frac{k}{2}-\frac{j}{2}\right)_{i_k} \left( \frac{k}{2}+1 \right)_{i_{k-1}} \left( \frac{k}{2}+ \frac{1}{2}+\frac{\nu }{2} \right)_{i_{k-1}}}{\left( \frac{k}{2}-\frac{j}{2}\right)_{i_{k-1}} \left( \frac{k}{2}+1 \right)_{i_k} \left( \frac{k}{2}+\frac{1}{2}+ \frac{\nu }{2} \right)_{i_k}} \right) \nonumber\\ 
&&\times \frac{\left( \frac{\tau }{2} -\frac{j}{2}\right)_{n} \left( \frac{\tau }{2}+1 \right)_{i_{\tau -1}} \left( \frac{\tau }{2}+\frac{1}{2}+\frac{\nu }{2} \right)_{i_{\tau -1}}}{\left( \frac{\tau }{2}-\frac{j}{2}\right)_{i_{\tau -1}} \left( \frac{\tau }{2}+1 \right)_{n} \left( \frac{\tau }{2}+\frac{1}{2}+ \frac{\nu }{2} \right)_{n}} \left( -\frac{1}{2}\mu \right)^{n } \label{eq:80022c} 
\end{eqnarray}
\begin{eqnarray}
y_0^m(j,\omega ;x) &=& \sum_{i_0=0}^{m} \frac{\left( -\frac{j}{2}\right)_{i_0}}{\left( 1 \right)_{i_0} \left( \frac{1}{2}+ \frac{\nu }{2} \right)_{i_0}} \tilde{\mu }^{i_0} \label{eq:80023a}\\
y_1^m(j,\omega ;x) &=& \left\{\sum_{i_0=0}^{m} \frac{\left( i_0+\frac{\omega }{2} \right) }{\left( i_0+\frac{1}{2} \right) \left( i_0+\frac{\nu }{2} \right)} \frac{\left( -\frac{j}{2}\right)_{i_0} }{\left( 1 \right)_{i_0} \left( \frac{1}{2}+\frac{\nu }{2} \right)_{i_0}} \right.   \left. \sum_{i_1 = i_0}^{m} \frac{\left( \frac{1}{2}-\frac{j}{2} \right)_{i_1} \left( \frac{3}{2} \right)_{i_0} \left( 1+\frac{\nu }{2} \right)_{i_0}}{\left( \frac{1}{2}-\frac{j}{2} \right)_{i_0} \left( \frac{3}{2} \right)_{i_1} \left( 1+\frac{\nu }{2}  \right)_{i_1}} \tilde{\mu }^{i_1}\right\} \tilde{\varepsilon } \hspace{1.5cm} \label{eq:80023b}\\
y_{\tau }^m(j,\omega ;x) &=& \left\{ \sum_{i_0=0}^{m} \frac{\left( i_0+\frac{\omega }{2} \right)}{\left( i_0+\frac{1}{2} \right) \left( i_0+\frac{\nu }{2} \right)} \frac{\left( -\frac{j}{2}\right)_{i_0} }{\left( 1 \right)_{i_0} \left( \frac{1}{2}+\frac{\nu }{2} \right)_{i_0}} \right.\nonumber\\
&&\times \prod_{k=1}^{\tau -1} \left( \sum_{i_k = i_{k-1}}^{m} \frac{\left( i_k+ \frac{k}{2}+\frac{\omega }{2} \right) }{\left( i_k+\frac{k}{2}+\frac{1}{2} \right) \left( i_k+\frac{k}{2}+\frac{\nu }{2} \right)} \right.  \left. \frac{\left( \frac{k}{2}-\frac{j}{2}\right)_{i_k} \left( \frac{k}{2}+1 \right)_{i_{k-1}} \left(  \frac{k}{2}+\frac{1}{2}+\frac{\nu }{2} \right)_{i_{k-1}}}{\left( \frac{k}{2}-\frac{j}{2}\right)_{i_{k-1}} \left( \frac{k}{2}+1 \right)_{i_k} \left( \frac{k}{2}+\frac{1}{2}+\frac{\nu }{2} \right)_{i_k}} \right) \nonumber\\
&&\times \left. \sum_{i_{\tau } = i_{\tau -1}}^{m} \frac{\left( \frac{\tau }{2}-\frac{j}{2}\right)_{i_{\tau }} \left( \frac{\tau }{2}+1 \right)_{i_{\tau -1}} \left( \frac{\tau }{2}+\frac{1}{2}+\frac{\nu }{2} \right)_{i_{\tau -1}}}{\left( \frac{\tau }{2}-\frac{j}{2}\right)_{i_{\tau -1}} \left( \frac{\tau }{2}+1 \right)_{i_{\tau }} \left( \frac{\tau }{2}+\frac{1}{2}+\frac{\nu }{2} \right)_{i_{\tau }}} \tilde{\mu }^{i_{\tau }}\right\} \tilde{\varepsilon }^{\tau }  \label{eq:80023c} 
\end{eqnarray}
\end{enumerate}
\end{remark}
\begin{remark}
The power series expansion of the GCH equation of the second kind for the first species complete polynomial using 3TRF about $x=0$ is given by
\begin{enumerate} 
\item As $\Omega = \mu \left( \nu -1 \right)$ and $\omega = \omega _0^0= \nu -1 $,

The eigenfunction is given by
\begin{eqnarray}
 y(x) &=& R_pW_{0,0}\left( \mu ,\varepsilon ,\nu ,\Omega = \mu \left( \nu -1 \right), \omega = \omega _0^0 =\nu -1; \tilde{\mu }  = -\frac{1}{2}\mu x^2; \tilde{\varepsilon } = -\frac{1}{2}\varepsilon x \right) \nonumber\\
&=& x^{1-\nu } \label{eq:80024}
\end{eqnarray}
\item As $\Omega = \mu \left( \nu -2N-2 \right) $ where $N \in \mathbb{N}_{0}$, 

An algebraic equation of degree $2N+2$ for the determination of $\omega $ is given by
\begin{equation}
0 = \sum_{r=0}^{N+1}\bar{c}\left( 2r, N+1-r; 2N+1,\omega \right) \label{eq:80025a}
\end{equation}
The eigenvalue of $\omega $ is written by $\omega _{2N+1}^m$ where $m = 0,1,2,\cdots,2N+1 $; $\omega _{2N+1}^0 < \omega _{2N+1}^1 < \cdots < \omega _{2N+1}^{2N+1}$. Its eigenfunction is given by
\begin{eqnarray} 
y(x) &=& R_pW_{2N+1,m} \left( \mu ,\varepsilon ,\nu ,\Omega = \mu \left( \nu -2N-2 \right), \omega = \omega _{2N+1}^m; \tilde{\mu }  = -\frac{1}{2}\mu x^2; \tilde{\varepsilon } = -\frac{1}{2}\varepsilon x \right) \nonumber\\
&=&  \sum_{r=0}^{N} y_{2r}^{N-r}\left( 2N+1,\omega _{2N+1}^m;x\right)+ \sum_{r=0}^{N} y_{2r+1}^{N-r}\left( 2N+1,\omega _{2N+1}^m;x\right)
\label{eq:80025b}
\end{eqnarray}
\item As $\Omega = \mu \left( \nu -2N-3 \right) $ where  $N \in \mathbb{N}_{0}$,

An algebraic equation of degree $2N+3$ for the determination of $\omega $ is given by
\begin{equation}
0 = \sum_{r=0}^{N+1}\bar{c}\left( 2r+1, N+1-r; 2N+2,\omega \right)  \label{eq:80026a}
\end{equation}
The eigenvalue of $\omega $ is written by $\omega _{2N+2}^m$ where $m = 0,1,2,\cdots,2N+2 $; $\omega _{2N+2}^0 < \omega _{2N+2}^1 < \cdots < \omega _{2N+2}^{2N+2}$. Its eigenfunction is given by
\begin{eqnarray} 
y(x) &=& R_pW_{2N+2,m} \left( \mu ,\varepsilon ,\nu ,\Omega = \mu \left( \nu -2N-3 \right), \omega = \omega _{2N+2}^m; \tilde{\mu }  = -\frac{1}{2}\mu x^2; \tilde{\varepsilon } = -\frac{1}{2}\varepsilon x \right) \nonumber\\
&=& \sum_{r=0}^{N+1} y_{2r}^{N+1-r}\left( 2N+2,\omega _{2N+2}^m;x\right) + \sum_{r=0}^{N} y_{2r+1}^{N-r}\left( 2N+2,\omega _{2N+2}^m;x\right) 
\label{eq:80026b}
\end{eqnarray}
In the above,
\begin{eqnarray}
\bar{c}(0,n;j,\omega )  &=& \frac{\left( -\frac{j}{2}\right)_{n}}{\left( \frac{3}{2}-\frac{\nu }{2}\right)_{n} \left( 1\right)_{n}} \left( -\frac{1}{2}\mu \right)^{n}\label{eq:80027a}\\
\bar{c}(1,n;j,\omega ) &=& \left( -\frac{1}{2}\varepsilon \right) \sum_{i_0=0}^{n}\frac{\left( i_0+\frac{1}{2}+\frac{\omega }{2}-\frac{\nu }{2}\right) }{\left( i_0+1-\frac{\nu }{2}\right) \left( i_0 +\frac{1}{2}\right)} \frac{\left( -\frac{j}{2}\right)_{i_0} }{\left( \frac{3}{2}-\frac{\nu }{2}\right)_{i_0} \left( 1\right)_{i_0}} \nonumber\\
&&\times \frac{\left( \frac{1}{2}-\frac{j}{2} \right)_{n} \left( 2-\frac{\nu  }{2}\right)_{i_0} \left( \frac{3}{2}\right)_{i_0}}{\left( \frac{1}{2}-\frac{j}{2}\right)_{i_0} \left( 2-\frac{\nu }{2}\right)_{n} \left( \frac{3}{2}\right)_{n}} \left( -\frac{1}{2}\mu \right)^{n } \label{eq:80027b}\\
\bar{c}(\tau ,n;j,\omega ) &=& \left( -\frac{1}{2}\varepsilon \right)^{\tau } \sum_{i_0=0}^{n}\frac{\left( i_0+\frac{1}{2}+\frac{\omega }{2}-\frac{\nu }{2}\right) }{\left( i_0+1-\frac{\nu }{2}\right) \left( i_0 +\frac{1}{2}\right)} \frac{\left( -\frac{j}{2}\right)_{i_0} }{\left( \frac{3}{2}-\frac{\nu }{2}\right)_{i_0} \left( 1\right)_{i_0}}  \nonumber\\
&&\times  \prod_{k=1}^{\tau -1} \left( \sum_{i_k = i_{k-1}}^{n} \frac{\left( i_k+ \frac{k}{2}+\frac{1}{2}+\frac{\omega }{2}-\frac{\nu }{2}\right) }{\left( i_k+\frac{k}{2}+1-\frac{\nu }{2}\right) \left( i_k+\frac{k}{2} +\frac{1}{2}\right)} \right.  \left. \frac{\left( \frac{k}{2}-\frac{j}{2}\right)_{i_k} \left( \frac{k}{2} +\frac{3}{2}-\frac{\nu }{2}\right)_{i_{k-1}} \left( \frac{k}{2}+ 1\right)_{i_{k-1}}}{\left( \frac{k}{2}-\frac{j}{2}\right)_{i_{k-1}} \left( \frac{k}{2}+\frac{3}{2}- \frac{\nu }{2}\right)_{i_k} \left( \frac{k}{2}+1\right)_{i_k}} \right) \nonumber\\ 
&&\times \frac{\left( \frac{\tau }{2} -\frac{j}{2}\right)_{n} \left( \frac{\tau }{2} +\frac{3}{2}-\frac{\nu }{2}\right)_{i_{\tau -1}} \left( \frac{\tau }{2}+1\right)_{i_{\tau -1}}}{\left( \frac{\tau }{2}-\frac{j}{2}\right)_{i_{\tau -1}} \left( \frac{\tau }{2} +\frac{3}{2}-\frac{\nu }{2}\right)_{n} \left( \frac{\tau }{2}+1\right)_{n}} \left( -\frac{1}{2}\mu \right)^{n } \label{eq:80027c} 
\end{eqnarray}
\begin{eqnarray}
y_0^m(j,\omega ;x) &=& x^{1-\nu }  \sum_{i_0=0}^{m} \frac{\left( -\frac{j}{2}\right)_{i_0}}{\left( \frac{3}{2}-\frac{\nu }{2}\right)_{i_0} \left( 1\right)_{i_0}} \tilde{\mu }^{i_0} \label{eq:80028a}\\
y_1^m(j,\omega ;x) &=& x^{1-\nu } \left\{\sum_{i_0=0}^{m} \frac{\left( i_0+\frac{1}{2}+\frac{\omega }{2}-\frac{\nu }{2}\right) }{\left( i_0+1-\frac{\nu }{2}\right) \left( i_0 +\frac{1}{2}\right)} \frac{\left( -\frac{j}{2}\right)_{i_0} }{\left( \frac{3}{2}-\frac{\nu }{2}\right)_{i_0} \left( 1\right)_{i_0}} \right. \nonumber\\
&&\times  \left. \sum_{i_1 = i_0}^{m} \frac{\left( \frac{1}{2}-\frac{j}{2} \right)_{i_1} \left( 2-\frac{\nu }{2}\right)_{i_0} \left( \frac{3}{2}\right)_{i_0}}{\left( \frac{1}{2}-\frac{j}{2} \right)_{i_0} \left( 2-\frac{\nu }{2}\right)_{i_1} \left( \frac{3}{2}\right)_{i_1}} \tilde{\mu }^{i_1}\right\} \tilde{\varepsilon } 
\hspace{1.5cm}\label{eq:80028b}\\
y_{\tau }^m(j,\omega ;x) &=&  x^{1-\nu } \left\{ \sum_{i_0=0}^{m} \frac{\left( i_0+\frac{1}{2}+\frac{\omega }{2}-\frac{\nu }{2}\right) }{\left( i_0+1-\frac{\nu }{2}\right) \left( i_0 +\frac{1}{2}\right)} \frac{\left( -\frac{j}{2}\right)_{i_0} }{\left( \frac{3}{2}-\frac{\nu }{2}\right)_{i_0} \left( 1\right)_{i_0}} \right.\nonumber\\
&&\times \prod_{k=1}^{\tau -1} \left( \sum_{i_k = i_{k-1}}^{m} \frac{\left( i_k+ \frac{k}{2}+\frac{1}{2}+\frac{\omega }{2}-\frac{\nu }{2}\right) }{\left( i_k+\frac{k}{2}+1-\frac{\nu }{2}\right) \left( i_k+\frac{k}{2} +\frac{1}{2}\right)} \right. \nonumber\\
&&\times \left. \frac{\left( \frac{k}{2}-\frac{j}{2}\right)_{i_k} \left( \frac{k}{2} +\frac{3}{2}-\frac{\nu }{2}\right)_{i_{k-1}} \left(  \frac{k}{2}+1\right)_{i_{k-1}}}{\left( \frac{k}{2}-\frac{j}{2}\right)_{i_{k-1}} \left( \frac{k}{2} +\frac{3}{2}-\frac{\nu }{2}\right)_{i_k} \left( \frac{k}{2}+1\right)_{i_k}} \right) \nonumber\\
&&\times \left. \sum_{i_{\tau } = i_{\tau -1}}^{m} \frac{\left( \frac{\tau }{2}-\frac{j}{2}\right)_{i_{\tau }} \left( \frac{\tau }{2} +\frac{3}{2}-\frac{\nu }{2}\right)_{i_{\tau -1}} \left( \frac{\tau }{2}+1\right)_{i_{\tau -1}}}{\left( \frac{\tau }{2}-\frac{j}{2}\right)_{i_{\tau -1}} \left( \frac{\tau }{2} +\frac{3}{2}-\frac{\nu }{2}\right)_{i_{\tau }} \left( \frac{\tau }{2}+1\right)_{i_{\tau }}} \tilde{\mu }^{i_{\tau }}\right\} \tilde{\varepsilon }^{\tau } \hspace{1.5cm}\label{eq:80028c} 
\end{eqnarray}
\end{enumerate}
\end{remark}
\subsection{The first species complete polynomial of the GCH equation using R3TRF}
\begin{theorem}
In chapter 2, the general summation expression of a function $y(x)$ for the first species complete polynomial using reversible 3-term recurrence formula and its algebraic equation for the determination of an accessory parameter in $A_n$ term are given by
\begin{enumerate} 
\item As $B_1=0$,
\begin{equation}
0 =\bar{c}(0,1) \label{eq:80029a}
\end{equation}
\begin{equation}
y(x) = y_{0}^{0}(x) \label{eq:80029b}
\end{equation}
\item As $B_2=0$, 
\begin{equation}
0 = \bar{c}(0,2)+\bar{c}(1,0) \label{eq:80030a}
\end{equation}
\begin{equation}
y(x)= y_{0}^{1}(x) \label{eq:80030b}
\end{equation}
\item As $B_{2N+3}=0$ where $N \in \mathbb{N}_{0}$,
\begin{equation}
0  = \sum_{r=0}^{N+1}\bar{c}\left( r, 2(N-r)+3\right) \label{eq:80031a}
\end{equation}
\begin{equation}
y(x)= \sum_{r=0}^{N+1} y_{r}^{2(N+1-r)}(x) \label{eq:80031b}
\end{equation}
\item As $B_{2N+4}=0$ where$N \in \mathbb{N}_{0}$,
\begin{equation}
0  =  \sum_{r=0}^{N+2}\bar{c}\left( r, 2(N+2-r)\right) \label{eq:80032a}
\end{equation}
\begin{equation}
y(x)=  \sum_{r=0}^{N+1} y_{r}^{2(N-r)+3}(x) \label{eq:80032b}
\end{equation}
In the above,
\begin{eqnarray}
\bar{c}(0,n) &=& \prod _{i_0=0}^{n-1}A_{i_0} \label{eq:80033a}\\
\bar{c}(1,n) &=& \sum_{i_0=0}^{n} \left\{ B_{i_0+1} \prod _{i_1=0}^{i_0-1}A_{i_1} \prod _{i_2=i_0}^{n-1}A_{i_2+2} \right\} \label{eq:80033b}\\
\bar{c}(\tau ,n) &=& \sum_{i_0=0}^{n} \left\{B_{i_0+1}\prod _{i_1=0}^{i_0-1} A_{i_1} 
\prod _{k=1}^{\tau -1} \left( \sum_{i_{2k}= i_{2(k-1)}}^{n} B_{i_{2k}+(2k+1)}\prod _{i_{2k+1}=i_{2(k-1)}}^{i_{2k}-1}A_{i_{2k+1}+2k}\right) \right. \nonumber\\
&&\times \left. \prod _{i_{2\tau} = i_{2(\tau -1)}}^{n-1} A_{i_{2\tau }+ 2\tau} \right\} \label{eq:80033c}
\end{eqnarray}
and
\begin{eqnarray}
y_0^m(x) &=& c_0 x^{\lambda} \sum_{i_0=0}^{m} \left\{ \prod _{i_1=0}^{i_0-1}A_{i_1} \right\} x^{i_0 } \label{eq:80034a}\\
y_1^m(x) &=& c_0 x^{\lambda} \sum_{i_0=0}^{m}\left\{ B_{i_0+1} \prod _{i_1=0}^{i_0-1}A_{i_1}  \sum_{i_2=i_0}^{m} \left\{ \prod _{i_3=i_0}^{i_2-1}A_{i_3+2} \right\}\right\} x^{i_2+2 } \label{eq:80034b}\\
y_{\tau }^m(x) &=& c_0 x^{\lambda} \sum_{i_0=0}^{m} \left\{B_{i_0+1}\prod _{i_1=0}^{i_0-1} A_{i_1} 
\prod _{k=1}^{\tau -1} \left( \sum_{i_{2k}= i_{2(k-1)}}^{m} B_{i_{2k}+(2k+1)}\prod _{i_{2k+1}=i_{2(k-1)}}^{i_{2k}-1}A_{i_{2k+1}+2k}\right) \right. \nonumber\\
&&\times \left. \sum_{i_{2\tau} = i_{2(\tau -1)}}^{m} \left( \prod _{i_{2\tau +1}=i_{2(\tau -1)}}^{i_{2\tau}-1} A_{i_{2\tau +1}+ 2\tau} \right) \right\} x^{i_{2\tau}+2\tau }\hspace{1cm}\mathrm{where}\;\tau \geq 2
\label{eq:80034c}
\end{eqnarray}
\end{enumerate}
\end{theorem}
According to (\ref{eq:8006}), $c_{j+1}=0$ is clearly an algebraic equation in $\omega $ of degree $j+1$ and thus has $j+1$ zeros denoted them by $\omega _j^m$ eigenvalues where $m = 0,1,2, \cdots, j$. They can be arranged in the following order: $\omega _j^0 < \omega _j^1 < \omega _j^2 < \cdots < \omega _j^j$.
 
Substitute (\ref{eq:8005a}) and (\ref{eq:80013}) into (\ref{eq:80033a})--(\ref{eq:80034c}).

As $B_{1}= c_{1}=0$, take the new (\ref{eq:80033a}) into (\ref{eq:80029a}) putting $j=0$. Substitute the new (\ref{eq:80034a}) into (\ref{eq:80029b}) putting $j=0$.

As $B_{2}= c_{2}=0$, take the new (\ref{eq:80033a}) and (\ref{eq:80033b}) into (\ref{eq:80030a}) putting $j=1$. Substitute the new (\ref{eq:80034a}) into (\ref{eq:80030b}) putting $j=1$ and $\omega = \omega _1^m$. 

As $B_{2N+3}= c_{2N+3}=0$, take the new (\ref{eq:80033a})--(\ref{eq:80033c}) into (\ref{eq:80031a}) putting $j=2N+2$. Substitute the new 
(\ref{eq:80034a})--(\ref{eq:80034c}) into (\ref{eq:80031b}) putting $j=2N+2$ and $\omega =\omega _{2N+2}^m$.

As $B_{2N+4}= c_{2N+4}=0$, take the new (\ref{eq:80033a})--(\ref{eq:80033c}) into (\ref{eq:80032a}) putting $j=2N+3$. Substitute the new 
(\ref{eq:80034a})--(\ref{eq:80034c}) into (\ref{eq:80032b}) putting $j=2N+3$ and $\omega =\omega _{2N+3}^m$.

After the replacement process, the general expression of power series of the GCH equation about $x=0$ for the first species complete polynomial using reversible 3-term recurrence formula and its algebraic equation for the determination of an accessory parameter $\omega $ are given by
\begin{enumerate} 
\item As $\Omega = -\mu \lambda $,

An algebraic equation of degree 1 for the determination of $\omega$ is given by
\begin{equation}
0= \bar{c}(0,1;0,\omega )= \omega + \lambda  \label{eq:80036a}
\end{equation}
The eigenvalue of $\omega $ is written by $\omega _0^0$. Its eigenfunction is given by
\begin{equation}
y(x) = y_0^0\left( 0,\omega _0^0;x\right)= c_0 x^{\lambda } \label{eq:80036b}  
\end{equation}
\item As $\Omega =-\mu \left( 1+\lambda \right)$,

An algebraic equation of degree 2 for the determination of $\omega $ is given by
\begin{eqnarray}
0 &=& \bar{c}(0,2;1,\omega )+\bar{c}(1,0;1,\omega ) \nonumber\\
&=& (\omega +\lambda )(\omega +1 +\lambda ) 
+ \frac{\mu }{\varepsilon ^2}(1+\lambda )(\nu +\lambda ) \label{eq:80037a}
\end{eqnarray}
The eigenvalue of $\omega $ is written by $\omega _1^m$ where $m = 0,1 $; $\omega _{1}^0 < \omega _{1}^1$. Its eigenfunction is given by
\begin{equation}
y(x) = y_{0}^{1}\left( 1,\omega _1^m;x\right)= c_0 x^{\lambda } \left\{ 1+\frac{\left( \omega _1^m + \lambda \right) }{(1+\lambda )(\nu +\lambda )}\Tilde{\Tilde{\varepsilon}} \right\} \label{eq:80037b}  
\end{equation}
\item As $\Omega =-\mu \left( 2N+2+\lambda \right) $ where $N \in \mathbb{N}_{0}$,

An algebraic equation of degree $2N+3$ for the determination of $\omega $ is given by
\begin{equation}
0 = \sum_{r=0}^{N+1}\bar{c}\left( r, 2(N-r)+3; 2N+2,\omega \right)  \label{eq:80038a}
\end{equation}
The eigenvalue of $\omega $ is written by $\omega _{2N+2}^m$ where $m = 0,1,2,\cdots,2N+2 $; $\omega _{2N+2}^0 < \omega _{2N+2}^1 < \cdots < \omega _{2N+2}^{2N+2}$. Its eigenfunction is given by 
\begin{equation} 
y(x) = \sum_{r=0}^{N+1} y_{r}^{2(N+1-r)}\left( 2N+2, \omega _{2N+2}^m; x \right)  
\label{eq:80038b} 
\end{equation}
\item As $\Omega =-\mu \left( 2N+3+\lambda \right) $ where $N \in \mathbb{N}_{0}$,

An algebraic equation of degree $2N+4$ for the determination of $\omega $ is given by
\begin{equation}  
0 = \sum_{r=0}^{N+2}\bar{c}\left( r, 2(N+2-r); 2N+3,\omega \right) \label{eq:80039a}
\end{equation}
The eigenvalue of $\omega $ is written by $\omega _{2N+3}^m$ where $m = 0,1,2,\cdots,2N+3 $; $\omega _{2N+3}^0 < \omega _{2N+3}^1 < \cdots < \omega _{2N+3}^{2N+3}$. Its eigenfunction is given by
\begin{equation} 
y(x) = \sum_{r=0}^{N+1} y_{r}^{2(N-r)+3} \left( 2N+3,\omega _{2N+3}^m;x\right) \label{eq:80039b}
\end{equation}
In the above,
\begin{eqnarray}
\bar{c}(0,n;j,\omega )  &=& \frac{\left( \omega +\lambda \right)_{n}}{\left( 1+ \lambda \right)_{n} \left( \nu +\lambda \right)_{n}} \left( -\varepsilon \right)^{n}\label{eq:80040a}\\
\bar{c}(1,n;j,\omega ) &=& \left( -\mu \right) \sum_{i_0=0}^{n}\frac{\left( i_0-j\right) }{\left( i_0+ 2 + \lambda \right) \left( i_0+1+ \nu + \lambda \right)} \frac{\left( \omega +\lambda \right)_{i_0} }{\left( 1+ \lambda \right)_{i_0} \left( \nu + \lambda \right)_{i_0}} \nonumber\\
&&\times \frac{\left( \omega + 2+\lambda  \right)_{n} \left( 3 + \lambda \right)_{i_0} \left( 2+ \nu + \lambda \right)_{i_0}}{\left( \omega +2+\lambda \right)_{i_0} \left( 3 + \lambda \right)_{n} \left( 2+ \nu + \lambda \right)_{n}} \left( -\varepsilon \right)^{n }  
\label{eq:80040b}\\
\bar{c}(\tau ,n;j,\omega ) &=& \left( -\mu \right)^{\tau } \sum_{i_0=0}^{n}\frac{\left( i_0-j\right)}{\left( i_0+2+ \lambda \right) \left( i_0+1+ \nu + \lambda \right)} \frac{\left( \omega +\lambda \right)_{i_0} }{\left( 1+ \lambda \right)_{i_0} \left( \nu + \lambda \right)_{i_0}}  \nonumber\\
&&\times \prod_{k=1}^{\tau -1} \left( \sum_{i_k = i_{k-1}}^{n} \frac{\left( i_k+ 2k-j\right) }{\left( i_k+2k+2+ \lambda \right) \left( i_k+2k+1+ \nu + \lambda \right)} \right. \nonumber\\
&&\times \left. \frac{\left( 2k+\omega +\lambda \right)_{i_k} \left( 2k+1+ \lambda \right)_{i_{k-1}} \left( 2k+\nu + \lambda \right)_{i_{k-1}}}{\left( 2k+\omega + \lambda \right)_{i_{k-1}} \left( 2k+1+ \lambda \right)_{i_k} \left( 2k+\nu + \lambda \right)_{i_k}} \right) \nonumber\\ 
&&\times \frac{\left( 2\tau +\omega +\lambda \right)_{n} \left( 2\tau +1+ \lambda \right)_{i_{\tau -1}} \left( 2\tau + \nu + \lambda  \right)_{i_{\tau -1}}}{\left( 2\tau +\omega +\lambda \right)_{i_{\tau -1}} \left( 2\tau +1+\lambda\right)_{n} \left( 2\tau +\nu  +\lambda \right)_{n}} \left( -\varepsilon \right)^{n } \hspace{1.5cm} \label{eq:80040c} 
\end{eqnarray}
\begin{eqnarray}
y_0^m(j,\omega ;x) &=& c_0 x^{\lambda }  \sum_{i_0=0}^{m} \frac{\left( \omega +\lambda \right)_{i_0}}{\left( 1+ \lambda \right)_{i_0} \left( \nu + \lambda \right)_{i_0}} \Tilde{\Tilde{\varepsilon}}^{i_0} \label{eq:80041a}\\
y_1^m(j,\omega ;x) &=& c_0 x^{\lambda } \left\{\sum_{i_0=0}^{m} \frac{\left( i_0-j \right) }{\left( i_0+2+ \lambda \right) \left( i_0+1+ \nu + \lambda \right)} \frac{\left( \omega +\lambda \right)_{i_0} }{\left( 1+ \lambda \right)_{i_0} \left( \nu + \lambda \right)_{i_0}} \right. \nonumber\\
&&\times \left. \sum_{i_1 = i_0}^{m} \frac{\left( \omega +2+ \lambda \right)_{i_1} \left( 3 + \lambda \right)_{i_0} \left( 2+ \nu  +  \lambda \right)_{i_0}}{\left( \omega +2+ \lambda \right)_{i_0} \left( 3 + \lambda \right)_{i_1} \left( 2+ \nu + \lambda \right)_{i_1}} \Tilde{\Tilde{\varepsilon}}^{i_1}\right\} \Tilde{\Tilde{\mu}}
\label{eq:80041b}\\
y_{\tau }^m(j,\omega ;x) &=& c_0 x^{\lambda } \left\{ \sum_{i_0=0}^{m} \frac{\left( i_0-j \right)}{\left( i_0+2+ \lambda \right) \left( i_0+1+ \nu + \lambda \right)} \frac{\left( \omega +\lambda \right)_{i_0} }{\left( 1+ \lambda \right)_{i_0} \left( \nu + \lambda \right)_{i_0}} \right.\nonumber\\
&&\times \prod_{k=1}^{\tau -1} \left( \sum_{i_k = i_{k-1}}^{m} \frac{\left( i_k+ 2k-j \right) }{\left( i_k+2k+2+ \lambda \right) \left( i_k+2k+1+ \nu + \lambda \right)} \right. \nonumber\\
&&\times \left. \frac{\left( 2k+\omega +\lambda \right)_{i_k} \left( 2k+1+\lambda \right)_{i_{k-1}} \left( 2k+\nu +\lambda \right)_{i_{k-1}}}{\left( 2k+\omega +\lambda \right)_{i_{k-1}} \left( 2k+1+\lambda \right)_{i_k} \left( 2k+\nu +\lambda \right)_{i_k}} \right) \nonumber\\
&&\times \left. \sum_{i_{\tau } = i_{\tau -1}}^{m} \frac{\left( 2\tau +\omega +\lambda \right)_{i_{\tau }} \left( 2\tau +1 +\lambda\right)_{i_{\tau -1}} \left( 2\tau +\nu +\lambda\right)_{i_{\tau -1}}}{\left( 2\tau +\omega +\lambda \right)_{i_{\tau -1}} \left( 2\tau +1+\lambda \right)_{i_{\tau }} \left( 2\tau +\nu +\lambda \right)_{i_{\tau }}} \Tilde{\Tilde{\varepsilon}}^{i_{\tau }}\right\} \Tilde{\Tilde{\mu}}^{\tau } \hspace{1.5cm} \label{eq:80041c} 
\end{eqnarray}
where
\begin{equation}
\begin{cases} \tau \geq 2 \cr
\Tilde{\Tilde{\varepsilon}} = -\varepsilon  x  \cr
\Tilde{\Tilde{\mu}} = -\mu  x^2 
\end{cases}\nonumber
\end{equation}
\end{enumerate}
Put $c_0$= 1 as $\lambda =0$ for the first kind of independent solutions of the CHE and $\lambda = 1-\nu $ for the second one in (\ref{eq:80036a})--(\ref{eq:80041c}). 
\begin{remark}
The power series expansion of the GCH equation of the first kind for the first species complete polynomial using R3TRF about $x=0$ is given by
\begin{enumerate} 
\item As $\Omega =0$ and $\omega =\omega _0^0=0$,

The eigenfunction is given by
\begin{equation}
y(x) = Q_pW_{0,0}^R \left( \mu ,\varepsilon ,\nu ,\Omega =0, \omega = \omega _0^0 =0; \Tilde{\Tilde{\mu}} = -\mu  x^2; \Tilde{\Tilde{\varepsilon}} = -\varepsilon  x  \right) =1 \label{eq:80042}
\end{equation}
\item As $\Omega =-\mu $,

An algebraic equation of degree 2 for the determination of $\omega $ is given by
\begin{equation}
0 = \omega (\omega +1 ) 
+ \frac{\mu }{\varepsilon ^2} \nu \label{eq:80043a}
\end{equation}
The eigenvalue of $\omega $ is written by $\omega _1^m$ where $m = 0,1 $; $\omega _{1}^0 < \omega _{1}^1$. Its eigenfunction is given by
\begin{eqnarray}
y(x) &=& Q_pW_{1,m}^R \left( \mu ,\varepsilon ,\nu ,\Omega =-\mu , \omega = \omega _1^m; \Tilde{\Tilde{\mu}} = -\mu  x^2; \Tilde{\Tilde{\varepsilon}} = -\varepsilon  x  \right) \nonumber\\
&=&  1+\frac{ \omega _1^m }{ \nu }\Tilde{\Tilde{\varepsilon}} \label{eq:80043b}  
\end{eqnarray}
\item As $\Omega =-\mu \left( 2N+2 \right) $ where $N \in \mathbb{N}_{0}$,

An algebraic equation of degree $2N+3$ for the determination of $\omega $ is given by
\begin{equation}
0 = \sum_{r=0}^{N+1}\bar{c}\left( r, 2(N-r)+3; 2N+2,\omega \right)  \label{eq:80044a}
\end{equation}
The eigenvalue of $\omega $ is written by $\omega _{2N+2}^m$ where $m = 0,1,2,\cdots,2N+2 $; $\omega _{2N+2}^0 < \omega _{2N+2}^1 < \cdots < \omega _{2N+2}^{2N+2}$. Its eigenfunction is given by 
\begin{eqnarray} 
y(x) &=&  Q_pW_{2N+2,m}^R  \left( \mu ,\varepsilon ,\nu ,\Omega =-\mu \left( 2N+2 \right), \omega = \omega _{2N+2}^m; \Tilde{\Tilde{\mu}} = -\mu  x^2; \Tilde{\Tilde{\varepsilon}} = -\varepsilon  x  \right)\nonumber\\
&=& \sum_{r=0}^{N+1} y_{r}^{2(N+1-r)}\left( 2N+2, \omega _{2N+2}^m; x \right)  
\label{eq:80044b} 
\end{eqnarray}
\item As $\Omega =-\mu \left( 2N+3 \right) $ where $N \in \mathbb{N}_{0}$,

An algebraic equation of degree $2N+4$ for the determination of $\omega $ is given by
\begin{equation}  
0 =  \sum_{r=0}^{N+2}\bar{c}\left( r, 2(N+2-r); 2N+3,\omega \right) \label{eq:80045a}
\end{equation}
The eigenvalue of $\omega $ is written by $\omega _{2N+3}^m$ where $m = 0,1,2,\cdots,2N+3 $; $\omega _{2N+3}^0 < \omega _{2N+3}^1 < \cdots < \omega _{2N+3}^{2N+3}$. Its eigenfunction is given by
\begin{eqnarray} 
y(x) &=& Q_pW_{2N+3,m}^R \left( \mu ,\varepsilon ,\nu ,\Omega =-\mu \left( 2N+3 \right), \omega = \omega _{2N+3}^m; \Tilde{\Tilde{\mu}} = -\mu  x^2; \Tilde{\Tilde{\varepsilon}} = -\varepsilon  x  \right)\nonumber\\
&=& \sum_{r=0}^{N+1} y_{r}^{2(N-r)+3} \left( 2N+3,\omega _{2N+3}^m;x\right) \label{eq:80045b}
\end{eqnarray}
In the above,
\begin{eqnarray}
\bar{c}(0,n;j,\omega )  &=& \frac{\left( \omega \right)_{n}}{\left( 1 \right)_{n} \left( \nu \right)_{n}} \left( -\varepsilon \right)^{n}\label{eq:80046a}\\
\bar{c}(1,n;j,\omega ) &=& \left( -\mu \right) \sum_{i_0=0}^{n}\frac{\left( i_0-j\right) }{\left( i_0+ 2 \right) \left( i_0+1+ \nu \right)} \frac{\left( \omega \right)_{i_0} }{\left( 1 \right)_{i_0} \left( \nu \right)_{i_0}}  \frac{\left( \omega + 2 \right)_{n} \left( 3 \right)_{i_0} \left( 2+ \nu \right)_{i_0}}{\left( \omega +2 \right)_{i_0} \left( 3 \right)_{n} \left( 2+ \nu \right)_{n}} \left( -\varepsilon \right)^{n }  
\hspace{1.5cm}\label{eq:80046b}\\
\bar{c}(\tau ,n;j,\omega ) &=& \left( -\mu \right)^{\tau } \sum_{i_0=0}^{n}\frac{\left( i_0-j\right)}{\left( i_0+2 \right) \left( i_0+1+ \nu \right)} \frac{\left( \omega \right)_{i_0} }{\left( 1 \right)_{i_0} \left( \nu \right)_{i_0}}  \nonumber\\
&&\times \prod_{k=1}^{\tau -1} \left( \sum_{i_k = i_{k-1}}^{n} \frac{\left( i_k+ 2k-j\right) }{\left( i_k+2k+2 \right) \left( i_k+2k+1+ \nu \right)} \right. \nonumber\\
&&\times \left. \frac{\left( 2k+\omega \right)_{i_k} \left( 2k+1 \right)_{i_{k-1}} \left( 2k+\nu \right)_{i_{k-1}}}{\left( 2k+\omega \right)_{i_{k-1}} \left( 2k+1 \right)_{i_k} \left( 2k+\nu \right)_{i_k}} \right) \nonumber\\ 
&&\times \frac{\left( 2\tau +\omega \right)_{n} \left( 2\tau +1 \right)_{i_{\tau -1}} \left( 2\tau + \nu \right)_{i_{\tau -1}}}{\left( 2\tau +\omega \right)_{i_{\tau -1}} \left( 2\tau +1 \right)_{n} \left( 2\tau +\nu \right)_{n}} \left( -\varepsilon \right)^{n } \label{eq:80046c} 
\end{eqnarray}
\begin{eqnarray}
y_0^m(j,\omega ;x) &=& \sum_{i_0=0}^{m} \frac{\left( \omega \right)_{i_0}}{\left( 1 \right)_{i_0} \left( \nu \right)_{i_0}} \Tilde{\Tilde{\varepsilon}}^{i_0} \label{eq:80047a}\\
y_1^m(j,\omega ;x) &=& \left\{\sum_{i_0=0}^{m} \frac{\left( i_0-j \right) }{\left( i_0+2 \right) \left( i_0+1+ \nu \right)} \frac{\left( \omega \right)_{i_0} }{\left( 1 \right)_{i_0} \left( \nu \right)_{i_0}} \right.   \left. \sum_{i_1 = i_0}^{m} \frac{\left( \omega +2 \right)_{i_1} \left( 3 \right)_{i_0} \left( 2+ \nu \right)_{i_0}}{\left( \omega +2 \right)_{i_0} \left( 3 \right)_{i_1} \left( 2+ \nu \right)_{i_1}} \Tilde{\Tilde{\varepsilon}}^{i_1}\right\} \Tilde{\Tilde{\mu}}
\label{eq:80047b}\\
y_{\tau }^m(j,\omega ;x) &=& \left\{ \sum_{i_0=0}^{m} \frac{\left( i_0-j \right)}{\left( i_0+2 \right) \left( i_0+1+ \nu  \right)} \frac{\left( \omega \right)_{i_0} }{\left( 1 \right)_{i_0} \left( \nu \right)_{i_0}} \right.\nonumber\\
&&\times \prod_{k=1}^{\tau -1} \left( \sum_{i_k = i_{k-1}}^{m} \frac{\left( i_k+ 2k-j \right) }{\left( i_k+2k+2 \right) \left( i_k+2k+1+ \nu \right)} \right.  \left. \frac{\left( 2k+\omega \right)_{i_k} \left( 2k+1 \right)_{i_{k-1}} \left( 2k+\nu  \right)_{i_{k-1}}}{\left( 2k+\omega \right)_{i_{k-1}} \left( 2k+1 \right)_{i_k} \left( 2k+\nu \right)_{i_k}} \right) \nonumber\\
&&\times \left. \sum_{i_{\tau } = i_{\tau -1}}^{m} \frac{\left( 2\tau +\omega \right)_{i_{\tau }} \left( 2\tau +1 \right)_{i_{\tau -1}} \left( 2\tau +\nu \right)_{i_{\tau -1}}}{\left( 2\tau +\omega \right)_{i_{\tau -1}} \left( 2\tau +1 \right)_{i_{\tau }} \left( 2\tau +\nu \right)_{i_{\tau }}} \Tilde{\Tilde{\varepsilon}}^{i_{\tau }}\right\} \Tilde{\Tilde{\mu}}^{\tau } \label{eq:80047c} 
\end{eqnarray}
\end{enumerate}
\end{remark}
\begin{remark}
The power series expansion of the GCH equation of the second kind for the first species complete polynomial using R3TRF about $x=0$ is given by
\begin{enumerate} 
\item As $\Omega = \mu \left( \nu -1 \right) $ and $\omega =\omega _0^0= \nu -1 $,

The eigenfunction is given by
\begin{eqnarray}
y(x) &=& R_pW_{0,0}^R \left( \mu ,\varepsilon ,\nu ,\Omega =\mu \left( \nu -1 \right), \omega = \omega _0^0 =\nu -1; \Tilde{\Tilde{\mu}} = -\mu  x^2; \Tilde{\Tilde{\varepsilon}} = -\varepsilon  x  \right) \nonumber\\
&=& x^{1-\nu } \label{eq:80048}
\end{eqnarray}
\item As $\Omega = \mu \left( \nu -2\right)$,

An algebraic equation of degree 2 for the determination of $\omega $ is given by
\begin{equation}
0 = (\omega -\nu +1)(\omega -\nu +2) 
+ \frac{\mu }{\varepsilon ^2}(2-\nu ) \label{eq:80049a}
\end{equation}
The eigenvalue of $\omega $ is written by $\omega _1^m$ where $m = 0,1 $; $\omega _{1}^0 < \omega _{1}^1$. Its eigenfunction is given by
\begin{eqnarray}
y(x) &=& R_pW_{1,m}^R \left( \mu ,\varepsilon ,\nu ,\Omega =\mu \left( \nu -2 \right), \omega = \omega _1^m; \Tilde{\Tilde{\mu}} = -\mu  x^2; \Tilde{\Tilde{\varepsilon}} = -\varepsilon  x  \right) \nonumber\\
&=& x^{1-\nu } \left\{ 1+\frac{\left( \omega _1^m + 1-\nu \right) }{ (2-\nu )}\Tilde{\Tilde{\varepsilon}} \right\} \label{eq:80049b}  
\end{eqnarray}
\item As $\Omega = \mu \left( \nu -2N-3\right) $ where $N \in \mathbb{N}_{0}$,

An algebraic equation of degree $2N+3$ for the determination of $\omega $ is given by
\begin{equation}
0 = \sum_{r=0}^{N+1}\bar{c}\left( r, 2(N-r)+3; 2N+2,\omega \right)  \label{eq:80050a}
\end{equation}
The eigenvalue of $\omega $ is written by $\omega _{2N+2}^m$ where $m = 0,1,2,\cdots,2N+2 $; $\omega _{2N+2}^0 < \omega _{2N+2}^1 < \cdots < \omega _{2N+2}^{2N+2}$. Its eigenfunction is given by 
\begin{eqnarray} 
y(x) &=& R_pW_{2N+2,m}^R \left( \mu ,\varepsilon ,\nu ,\Omega =\mu \left( \nu -2N-3 \right), \omega = \omega _{2N+2}^m; \Tilde{\Tilde{\mu}} = -\mu  x^2; \Tilde{\Tilde{\varepsilon}} = -\varepsilon  x  \right) \nonumber\\
&=& \sum_{r=0}^{N+1} y_{r}^{2(N+1-r)}\left( 2N+2, \omega _{2N+2}^m; x \right)  
\label{eq:80050b} 
\end{eqnarray}
\item As $\Omega = \mu \left( \nu -2N-4\right) $ where $N \in \mathbb{N}_{0}$,

An algebraic equation of degree $2N+4$ for the determination of $\omega $ is given by
\begin{equation}  
0 =  \sum_{r=0}^{N+2}\bar{c}\left( r, 2(N+2-r); 2N+3,\omega \right) \label{eq:80051a}
\end{equation}
The eigenvalue of $\omega $ is written by $\omega _{2N+3}^m$ where $m = 0,1,2,\cdots,2N+3 $; $\omega _{2N+3}^0 < \omega _{2N+3}^1 < \cdots < \omega _{2N+3}^{2N+3}$. Its eigenfunction is given by
\begin{eqnarray} 
y(x) &=& R_pW_{2N+3,m}^R \left( \mu ,\varepsilon ,\nu ,\Omega =\mu \left( \nu -2N-4 \right), \omega = \omega _{2N+3}^m; \Tilde{\Tilde{\mu}} = -\mu  x^2; \Tilde{\Tilde{\varepsilon}} = -\varepsilon  x  \right) \nonumber\\
&=& \sum_{r=0}^{N+1} y_{r}^{2(N-r)+3} \left( 2N+3,\omega _{2N+3}^m;x\right) \label{eq:80051b}
\end{eqnarray}
In the above,
\begin{eqnarray}
\bar{c}(0,n;j,\omega )  &=& \frac{\left( \omega -\nu +1 \right)_{n}}{\left( 2-\nu \right)_{n} \left( 1\right)_{n}} \left( -\varepsilon \right)^{n}\label{eq:80052a}\\
\bar{c}(1,n;j,\omega ) &=& \left( -\mu \right) \sum_{i_0=0}^{n}\frac{\left( i_0-j\right) }{\left( i_0+ 3-\nu \right) \left( i_0+2 \right)} \frac{\left( \omega -\nu +1 \right)_{i_0} }{\left( 2-\nu \right)_{i_0} \left( 1\right)_{i_0}} \nonumber\\
&&\times  \frac{\left( \omega -\nu +3 \right)_{n} \left( 4-\nu \right)_{i_0} \left( 3\right)_{i_0}}{\left( \omega -\nu +3 \right)_{i_0} \left( 4-\nu \right)_{n} \left( 3\right)_{n}} \left( -\varepsilon \right)^{n } \label{eq:80052b}\\
\bar{c}(\tau ,n;j,\omega ) &=& \left( -\mu \right)^{\tau } \sum_{i_0=0}^{n}\frac{\left( i_0-j\right)}{\left( i_0+3-\nu \right) \left( i_0+2 \right)} \frac{\left( \omega -\nu +1 \right)_{i_0} }{\left( 2-\nu \right)_{i_0} \left( 1\right)_{i_0}}  \nonumber\\
&&\times \prod_{k=1}^{\tau -1} \left( \sum_{i_k = i_{k-1}}^{n} \frac{\left( i_k+ 2k-j\right) }{\left( i_k+2k+3-\nu \right) \left( i_k+2k+2 \right)} \right. \nonumber\\
&&\times \left. \frac{\left( 2k+1+ \omega -\nu \right)_{i_k} \left( 2k+2-\nu \right)_{i_{k-1}} \left( 2k+1\right)_{i_{k-1}}}{\left( 2k+1+\omega -\nu \right)_{i_{k-1}} \left( 2k+2-\nu \right)_{i_k} \left( 2k+1\right)_{i_k}} \right) \nonumber\\ 
&&\times \frac{\left( 2\tau +1+\omega -\nu \right)_{n} \left( 2\tau +2-\nu \right)_{i_{\tau -1}} \left( 2\tau +1 \right)_{i_{\tau -1}}}{\left( 2\tau +1+\omega -\nu \right)_{i_{\tau -1}} \left( 2\tau +2-\nu \right)_{n} \left( 2\tau +1\right)_{n}} \left( -\varepsilon \right)^{n } \hspace{1.5cm} \label{eq:80052c} 
\end{eqnarray}
\begin{eqnarray}
y_0^m(j,\omega ;x) &=& x^{1-\nu }  \sum_{i_0=0}^{m} \frac{\left( \omega -\nu +1\right)_{i_0}}{\left( 2-\nu \right)_{i_0} \left( 1\right)_{i_0}} \Tilde{\Tilde{\varepsilon}}^{i_0} \label{eq:80053a}\\
y_1^m(j,\omega ;x) &=& x^{1-\nu } \left\{\sum_{i_0=0}^{m} \frac{\left( i_0-j \right) }{\left( i_0+3-\nu \right) \left( i_0+2 \right)} \frac{\left( \omega -\nu +1 \right)_{i_0} }{\left( 2-\nu \right)_{i_0} \left( 1\right)_{i_0}} \right. \nonumber\\
&&\times  \left. \sum_{i_1 = i_0}^{m} \frac{\left( \omega -\nu +3 \right)_{i_1} \left( 4-\nu \right)_{i_0} \left( 3\right)_{i_0}}{\left( \omega -\nu +3 \right)_{i_0} \left( 4-\nu \right)_{i_1} \left( 3 \right)_{i_1}} \Tilde{\Tilde{\varepsilon}}^{i_1}\right\} \Tilde{\Tilde{\mu}}
\label{eq:80053b}\\
y_{\tau }^m(j,\omega ;x) &=& x^{1-\nu } \left\{ \sum_{i_0=0}^{m} \frac{\left( i_0-j \right)}{\left( i_0+3-\nu \right) \left( i_0+2 \right)} \frac{\left( \omega -\nu +1 \right)_{i_0} }{\left( 2-\nu \right)_{i_0} \left( 1\right)_{i_0}} \right.\nonumber\\
&&\times \prod_{k=1}^{\tau -1} \left( \sum_{i_k = i_{k-1}}^{m} \frac{\left( i_k+ 2k-j \right) }{\left( i_k+2k+3-\nu \right) \left( i_k+2k+2 \right)} \right. \nonumber\\
&&\times \left. \frac{\left( 2k+1+\omega -\nu \right)_{i_k} \left( 2k+2-\nu \right)_{i_{k-1}} \left( 2k+1 \right)_{i_{k-1}}}{\left( 2k+1+\omega -\nu \right)_{i_{k-1}} \left( 2k+2-\nu \right)_{i_k} \left( 2k+1 \right)_{i_k}} \right) \nonumber\\
&&\times \left. \sum_{i_{\tau } = i_{\tau -1}}^{m} \frac{\left( 2\tau +1+\omega -\nu \right)_{i_{\tau }} \left( 2\tau +2-\nu \right)_{i_{\tau -1}} \left( 2\tau +1 \right)_{i_{\tau -1}}}{\left( 2\tau +1+\omega -\nu \right)_{i_{\tau -1}} \left( 2\tau +2-\nu \right)_{i_{\tau }} \left( 2\tau +1 \right)_{i_{\tau }}} \Tilde{\Tilde{\varepsilon}}^{i_{\tau }}\right\} \Tilde{\Tilde{\mu}}^{\tau } \hspace{1.5cm} \label{eq:80053c} 
\end{eqnarray}
\end{enumerate}
\end{remark}
\section{Summary}
Numerical solutions of the GCH equation for the first species complete polynomials using 3TRF and R3TRF are equivalent to each other. 
However, there are three significant differences of mathematical summation structures for these two complete polynomials:

(1) $A_n$ in sequences $c_n$ is the leading term in each of finite sub-power series of general series solutions of the GCH equation for the first species complete polynomials using 3TRF. But $B_n$ is the leading term in each of finite sub-power series of its general series solutions for the first species complete polynomials using R3TRF. 

(2) Summation solutions of the GCH equation for the first species complete polynomial using 3TRF contain the sum of two finite sub-power series of its general series solutions. In contrast, its solutions in series for the first species complete polynomial using R3TRF only consist of one finite sub-formal series of its general Frobenius solution. 

(3) As we observe solutions of the GCH equation for the first species complete polynomial using 3TRF, the denominators and numerators in all $B_n$ terms of each finite sub-power series arise with Pochhammer symbols. And for the first species complete polynomial using R3TRF, the denominators and numerators in all $A_n$ terms of each finite sub-power series arise with Pochhammer symbols.
 
Independent solutions of the GCH equation of the first kind for the first species complete polynomials using 3TRF and R3TRF require $\nu \ne 0,-1,-2,\cdots$ from the principle of the radius of convergence. And $\nu \ne 2,3,4,\cdots$ is required for its independent solutions of the second kind for the first species complete polynomials using 3TRF and R3TRF. For instance, if $\nu $ is a non-positive integer in the first kind of an GCH complete polynomial using 3TRF, its formal series solution diverges because several denominators turn to zero at the specific value of an index summation $n$. 

Numerical computations for the determination of a parameter $\omega $ in two polynomial equations of the GCH equation for the first species complete polynomials using 3TRF and R3TRF are equivalent to each other. 
Unsimilarity between two algebraic equations of the GCH equation is that $A_n$ is the leading term of each sequences $\bar{c}(l,n;j,\omega )$ where $l \in \mathbb{N}_{0}$ in a polynomial equation of $\omega $ for the first species complete polynomial using 3TRF. And $B_n$ is the leading term of each sequences $\bar{c}(l,n;j,\omega )$ in an algebraic equation of $\omega $ for the first species complete polynomial using R3TRF.
 
\addcontentsline{toc}{section}{Bibliography}
\bibliographystyle{model1a-num-names}
\bibliography{<your-bib-database>}
\chapter{Complete polynomials of Grand Confluent Hypergeometric equation about the irregular singular point at infinity}
\chaptermark{Complete polynomials of the GCH around $x=\infty $} 
 
In chapter 6 of Ref.\cite{9Choun2013}, by applying three term recurrence formula (3TRF) \cite{9Chou2012b}, I construct a formal series solution in the closed form of the GCH polynomial of type 1 about the irregular singular point at infinity. And its combined definite \& contour integral involving only $_1F_1$ functions and the generating function of it are derived analytically . 
 
In this chapter I show power series solutions of the GCH equation about $x=\infty $ for a polynomial of type 3 by applying general summation expressions of complete polynomials using 3TRF and reversible 3TRF (R3TRF).
 
\section{Introduction}
G\"{u}rsey and collaborators wrote the spin free semi-relativistic Hamiltonian for the quark-antiquark system neglecting small mass of quarks\cite{91985,91988,91991,91990}, suggested by Lichitenberg \textit{et al.}\cite{9Lich1982,9Krol1980,9Krol1981,9Todo1971}.
They noticed that their supersymmetric wave equation with ignoring the mass of quarks is essentially equal to confluent hypergeometric equation. By the method of a series solution, its recurrence relation for the coefficients consists of a 2-term and they showed the normalized spin free wave function involving only scalar potential neglecting the mass of quarks. 

Since I include the mass of quarks into their supersymmetric spin free Hamiltonian involving only scalar potential, I detect the recursion relation of this differential equation for the coefficients is composed of a 3-term and this ODE generalizes the confluent hypergeometric equation. I designate this new type equation as `grand confluent hypergeometric (GCH) equation.'  
 
We deal with the canonical form of the GCH equation\cite{9Chou2012i}
\begin{equation}
x \frac{d^2{y}}{d{x}^2} + \left( \mu x^2 + \varepsilon x + \nu  \right) \frac{d{y}}{d{x}} + \left( \Omega x + \varepsilon \omega \right) y = 0
\label{eq:9001}
\end{equation}
where $\mu$, $\varepsilon$, $\nu $, $\Omega$ and $\omega$ are real or complex parameters.
It is the proto-type of the Fuchsian equation with two singular points: one regular singular point which is zero with exponents $\{ 0,1-\nu \}$, and one irregular singular point which is infinity with an exponent $\frac{\Omega }{\mu }$.

Heun ordinary differential equation--Fuchsian equation with four singularities (named after Karl Heun, 1859--1929 \cite{9Heun1889}) is much more diversified than the hypergeomegtric-type equations. Any linear second-order differential equations with four singular points can be transformed to the Heun's equation by appropriate elementary transformations of the independent and dependent variables.
By the process of `confluence' such as deriving of confluent hypergeometric equation from the hypergeometric equation, we obtain four confluent forms by merging two or more regular singularities to take an irregular singularity in the Heun equation. 
Because of these reasons, several scholars treat the Heun equation as a successor of hypergeomegtric equations in $21^{th}$ century  \cite{9Fizi2012}. 
 
The canonical form of Heun's differential equation is written by \cite{9Heun1889,9Ronv1995,9Slavy2000}
\begin{equation}
\frac{d^2{y}}{d{x}^2} + \left(\frac{\gamma }{x} +\frac{\delta }{x-1} + \frac{\epsilon }{x-a}\right) \frac{d{y}}{d{x}} +  \frac{\alpha \beta x-q}{x(x-1)(x-a)} y = 0 \label{eq:9002}
\end{equation}
where $\epsilon = \alpha +\beta -\gamma -\delta +1$ for assuring the regularity of the point at $x=\infty $. It has four regular singular points which are 0, 1, $a$ and $\infty $ with exponents $\{ 0, 1-\gamma \}$, $\{ 0, 1-\delta \}$, $\{ 0, 1-\epsilon \}$ and $\{ \alpha, \beta \}$. 

Heun equation has four different confluent forms such as: (1) Confluent Heun (two regular and one irregular singularities), (2) Doubly Confluent Heun (two irregular singularities), (3) Biconfluent Heun (one regular and one irregular singularities), (4) Triconfluent Heun equations (one irregular singularity).
 
The canonical form of Biconfluent Heun (BCH) differential equation is given by\cite{9Ronv1995,9Maro1967}
\begin{equation}
x\frac{d^2{y}}{d{x}^2} + \left( -2x^2 -\beta x + 1+\alpha \right) \frac{d{y}}{d{x}} + \left( (\gamma -\alpha -2)x-\frac{1}{2}\left[ \delta +(1+\alpha )\beta \right] \right) y = 0 \label{eq:9003}
\end{equation}
where $\alpha $, $\beta $, $\gamma  $ and $\delta $ are real or complex parameters.
It is the Fuchsian-type equation with two singular points: one regular singular point which is zero with exponents $\{ 0,-\alpha \}$, and one irregular singular point which is infinity with an exponent $-\frac{1}{2}(\gamma -\alpha -2)$.

For DLFM version \cite{9NIST} or in Ref.\cite{9Slavy2000}, the BCH canonical equation is mentioned by
\begin{equation}
x\frac{d^2{y}}{d{x}^2} - \left( \gamma +\delta x +x^2\right) \frac{d{y}}{d{x}} + \left( \alpha x -q \right) y = 0 \label{eq:9004}
\end{equation}
It is also the Fuchsian-type equation with one regular singular point which is zero with exponents $\{ 0,1+\gamma \}$, and one irregular singular point which is infinity with an exponent $-\alpha $.

The GCH equation generalize Biconfluent Heun (BCH) equation having a regular singularity at $x=0$ and an irregular singularity at $\infty$ of rank 2. 
As we compare (\ref{eq:9001}) with (\ref{eq:9003}) and (\ref{eq:9004}), all coefficients on the above are correspondent to the following way.
\begin{equation}
\begin{array}{ll}
\left\{\begin{array}{c}
 \varepsilon  \longrightarrow  -\beta \\
\mu  \longrightarrow   -2 \\
\nu \longrightarrow  1+\alpha \\
\omega  \longrightarrow 1/2 (\delta /\beta +1+\alpha ) \\
\Omega  \longrightarrow  \gamma -\alpha -2  
\end{array}\right.  
\hspace{2cm} 
\left\{\begin{array}{c}
\varepsilon   \longrightarrow   -\delta  \\
\mu  \longrightarrow  -1 \\
\nu  \longrightarrow  -\gamma  \\
\omega  \longrightarrow  q/\delta  \\
\Omega  \longrightarrow \alpha 
\end{array}\right. 
\end{array}
\label{eq:9005}
\end{equation}
By the classical method of a series solution for the BCH equation, the special case of the GCH equation, assuming $y(x)=\sum_{n=0}^{\infty } c_n x^{n+\lambda }$, $c_0 \ne 0$, the recursive relation for the coefficients is governed by a 3-term instead of a 2-term. Consequently, no such solutions for its formal series in closed forms have been found in which coefficients are given explicitly yet. And there are no examples for the BCH equation in terms of definite or contour integrals involving simpler functions such as hypergeometric-type functions.

The BCH equation is applicable to diverse areas in mathematics \cite{9Bato1977,9Belm2004,9Cheb2004,9Deca1978,9Decar1978,9Exto1989,9Haut1969,9Maro1967,9Urwi1975}.  
Such solutions have recently been used in modern physics such as the radial Schr$\ddot{\mbox{o}}$dinger equation for various harmonic oscillators \cite{9Chau1983,9Chau1984,9Fles1980,9Fles1982,9Leau1986,9Leau1990,9Lemi1969,9Mass1983,9Pons1988}, the relativistic quantum mechanical model in a uniform magnetic field and scalar potentials \cite{9Figu2012}, quantum mechanical wave functions of two charged particles moving on a plane with a perpendicular uniform magnetic field \cite{9Ralk2002}, a two-electron quantum dot model \cite{9Caru2014}, the Dirac equation with the mixed scalar–vector–pseudoscalar Cornell potential \cite{9Hamz2013}, the two-dimensional Schr$\ddot{\mbox{o}}$dinger equation involving a parabolic potential for two interacting electrons in a uniform magnetic field \cite{9Kand2005}, the singular and the 2:1 anisotropic Dunkl oscillator models \cite{9Gene2013}, the Dirac equation with scalar, vector and tensor Cornell radial potentials \cite{9Cast2013} and etc. 
 
In general, by the method of a formal series solution into any linear ODEs, assuming $y(x)= \sum_{n=0}^{\infty } c_n x^{n+\lambda }$, $c_0 \ne 0$, the recurrence relation for the coefficients starts to appear. A 3-term recurrence relation is given by 
\begin{equation}
c_{n+1}=A_n \;c_n +B_n \;c_{n-1} \hspace{1cm};n\geq 1
\label{eq:9006}
\end{equation}
where $c_1= A_0 \;c_0$ and $\lambda $ is an indicial root. On the above, $A_n$ and $B_n$ are themselves polynomials of degree $n$: for the second-order ODEs, a numerator and a denominator of $A_n$ are usually equal or less than polynomials of degrees 2. Likewise, a numerator and a denominator of $B_n$ are also equal or less than polynomials of degrees 2

Any linear ODEs having a 2-term recurrence relation between successive coefficients, hypergeometric differential equation is one of examples, have two possible series solutions such as an infinite series and a polynomial.
In contrast, there are $2^{3-1}$ possible power series solutions of linear ODEs having a 3-term recursion relation between consecutive coefficients such as an infinite series and 3 types of polynomials: (1) a polynomial which makes $B_n$ term terminated; $A_n$ term is not terminated, (2) a polynomial which makes $A_n$ term terminated; $B_n$ term is not terminated, (3) a polynomial which makes $A_n$ and $B_n$ terms terminated at the same time, referred as `a complete polynomial.' 

The sequence $c_n$ is expanded to combinations of $A_n$ and $B_n$ terms in (\ref{eq:9006}).
I specify that a general series solution $y(x)$ is composed of the sum of sub-power series $y_l(x)$ where $l\in \mathbb{N}_{0}$
, referred as $y(x)= \sum_{n=0}^{\infty } y_n(x)$. And a $y_l(x)$ is constructed by observing the term of sequence $c_n$ which includes $l$ terms of $A_n's$.
By allowing $A_n$ in the sequence $c_n$ is the leading term of each of sub-power series solutions $y_l(x)$, 
power series solutions for an infinite series and a polynomial of type 1 are obtained by applying the general summation formulas of the 3-term recurrence relation in a linear ODE, referred as `three term recurrence formula (3TRF).' \cite{9Chou2012b}
 
A sub-power series $y_l(x)$ is obtained by observing the term of sequence $c_n$ which includes $l$ terms of $B_n's$ in chapter 1 of Ref.\cite{9Choun2013}.    
By allowing $B_n$ in the sequence $c_n$ is the leading term of each of sub-power series solutions $y_l(x)$ in a function $y(x)= \sum_{n=0}^{\infty } y_n(x)$, formal series solutions for an infinite series and a polynomial of type 2 are obtained by applying the general summation expressions of the 3-term recurrence relation in a linear ODE, denominated as `reversible three term recurrence formula (R3TRF).'

Numerical calculations of an infinite series by applying 3TRF is equivalent to its computations of an infinite series by applying R3TRF. Difference for mathematical summation structures between two infinite series solutions is that $A_n$ in sequences $c_n$ is the leading term in each of sub-power series solutions for an infinite series by applying 3TRF. And $B_n$ is the leading term in each of sub-formal series solutions for an infinite series by applying R3TRF. 
   
Complete polynomials are classified into two different types such as (1) the first species complete polynomial where a parameter of a numerator in $B_n$ term and a (spectral) parameter of a numerator in $A_n$ term are fixed values and (2) the second species complete polynomial where two parameters of a numerator in $B_n$ term and a parameter of a numerator in $A_n$ term are fixed values.
With my definition, the first species complete polynomial has multi-valued roots of a parameter in $A_n$ term, but the second  species complete polynomial has only one fixed value of a parameter in $A_n$ term for an eigenvalue. 
 
I define that a general series solution $y(x)$ is composed of the sum of finite sub-power series $y_{\tau }^m(x)$ where $\tau ,m \in \mathbb{N}_{0}$ in chapter 1.  And a finite sub-power series $y_{\tau }^m(x)$ is a polynomial that the lower bound of summation is zero and the upper one is $m$. $y_{\tau }^m(x)$ is obtained by observing the term of sequence $c_n$ which includes $\tau $ terms of $A_n's$. 
By allowing $A_n$ as the leading term in each finite sub-power series $y_{\tau }^m(x)$, I construct the mathematical summation formulas of the first and second species complete polynomials, designated as `complete polynomials using 3-term recurrence formula (3TRF).' 

In chapter 2, I suggest that a finite sub-power series $y_{\tau }^m(x)$ is a polynomial, obtained by observing the term of sequence $c_n$ which includes $\tau $ terms of $B_n's$.  
By allowing $B_n$ as the leading term in each finite sub-power series $y_{\tau }^m(x)$ of the general power series $y(x)$, I construct the classical mathematical expressions in series of the first and second species complete polynomials, denominated as `complete polynomials using reversible 3-term recurrence formula (R3TRF).' 

In this chapter I show general solutions in series in compact forms, called summation notation, of the GCH equation around  $x=\infty$ for polynomials of type 3 by applying complete polynomials using 3TRF and R3TRF.
 
\section{Power series about the irregular singular point at $x=\infty $}
Let $z=\frac{1}{x}$ in (\ref{eq:9001}) in order to obtain the analytic solution of the GCH equation about $x=\infty $.
\begin{equation}
z^4 \frac{d^2{y}}{d{z}^2} + \left( (2-\nu ) z^3 - \varepsilon z^2 - \mu z\right) \frac{d{y}}{d{z}} + \left( \Omega + \varepsilon \omega z\right) y = 0
\label{eq:9007}
\end{equation}
By the method of a series solution into (\ref{eq:9007}), assuming $y(z)= \sum_{n=0}^{\infty } c_n z^{n+\lambda }$, $c_0 \ne 0$, we find, at the singular point $z=0$, only one exponent $\lambda =\Omega /\mu $, and the recurrence relation between the coefficients: 
\begin{equation}
c_{n+1}=A_n \;c_n +B_n \;c_{n-1} \hspace{1cm};n\geq 1 \label{eq:9008}
\end{equation}
where,
\begin{subequations}
\begin{equation}
A_n = -\frac{\varepsilon }{\mu } \frac{\left( n-\omega +\frac{\Omega }{\mu } \right)}{\left( n+1 \right)} \label{eq:9009a}
\end{equation}
\begin{equation}
B_n = \frac{1}{\mu } \frac{\left( n-1 +\frac{\Omega }{\mu } \right)\left( n-\nu +\frac{\Omega }{\mu } \right)}{\left( n+1 \right)} \label{eq:9009b}
\end{equation}
\end{subequations}
where $c_1= A_0 \;c_0$.

All possible general solutions in series of the GCH equation about the irregular singular point at infinity are described in a following table.
\begin{table}[h]
\begin{center}
{
 \Tree[.{  The GCH differential equation about the irregular singular point at infinity} [.{ 3TRF} 
              [.{ Polynomials} [[.{ Polynomial of type 1} ]
               [.{ Polynomial of type 3} [.{ $ \begin{array}{lcll}  1^{\mbox{st}}\;  \mbox{species}\\ \mbox{complete} \\ \mbox{polynomial} \end{array}$} ] [.{ $ \begin{array}{lcll}  2^{\mbox{nd}}\;  \mbox{species}\\ \mbox{complete} \\ \mbox{polynomial} \end{array}$} ]]]]]                         
  [.{ R3TRF}        
       [.{ Polynomial of type 3} [.{ $ \begin{array}{lcll}  1^{\mbox{st}}\;  \mbox{species} \\ \mbox{complete} \\ \mbox{polynomial} \end{array}$} ] [.{ $ \begin{array}{lcll}  2^{\mbox{nd}}\;  \mbox{species} \\ \mbox{complete} \\ \mbox{polynomial} \end{array}$} ]] ]] 
}
\end{center}
\caption{Power series of the GCH equation about the irregular singular point at infinity}
\end{table}

As we see the table, there are no power series solutions for an infinite series and a polynomial of type 2 because of $\lim_{n \ge 1} A_n = -\frac{\varepsilon }{\mu }$ and $\lim_{n \ge 1} B_n \rightarrow \infty $. 
From the principle of the radius of convergence, $B_n$ terms must to be terminated at the specific index of summation $n$ in each of sub-power series in a general series solution. For an infinite series and a polynomial of type 2, power series solutions will be divergent, because $B_n$ terms are not convergent as the index of summation $n$ goes to infinity.  

There are two possible combinatorial cases for the first species complete polynomial of the GCH equation around $x=\infty $. For the first case, I treat $\omega $, $\Omega $ as fixed values and $\varepsilon $, $\mu $, $\nu $ as free variables. For the second one, I treat $\nu $, $\omega $ as fixed values and $\varepsilon $, $\mu $, $\Omega $ as free variables.   
And for the second species complete polynomial of the GCH equation, I treat $\nu $, $\omega $, $\Omega $ as fixed values and $\varepsilon $, $\mu $ as free variables.
 
\subsection{The first species complete polynomial of the GCH equation using 3TRF}
For the first species complete polynomials using 3TRF and R3TRF, we need a condition which is given by
\begin{equation}
B_{j+1}= c_{j+1}=0\hspace{1cm}\mathrm{where}\;j\in \mathbb{N}_{0}  
 \label{eq:90010}
\end{equation}
(\ref{eq:90010}) gives successively $c_{j+2}=c_{j+3}=c_{j+4}=\cdots=0$. And $c_{j+1}=0$ is defined by a polynomial equation of degree $j+1$ for the determination of an accessory parameter in $A_n$ term. 
\begin{theorem}
In chapter 1, the general expression of a function $y(x)$ for the first species complete polynomial using 3-term recurrence formula and its algebraic equation for the determination of an accessory parameter in $A_n$ term are given by
\begin{enumerate} 
\item As $B_1=0$,
\begin{equation}
0 =\bar{c}(1,0) \label{eq:90011a}
\end{equation}
\begin{equation}
y(x) = y_{0}^{0}(x) \label{eq:90011b}
\end{equation}
\item As $B_{2N+2}=0$ where $N \in \mathbb{N}_{0}$,
\begin{equation}
0  = \sum_{r=0}^{N+1}\bar{c}\left( 2r, N+1-r\right) \label{eq:90012a}
\end{equation}
\begin{equation}
y(x)= \sum_{r=0}^{N} y_{2r}^{N-r}(x)+ \sum_{r=0}^{N} y_{2r+1}^{N-r}(x)  \label{eq:90012b}
\end{equation}
\item As $B_{2N+3}=0$ where $N \in \mathbb{N}_{0}$,
\begin{equation}
0  = \sum_{r=0}^{N+1}\bar{c}\left( 2r+1, N+1-r\right) \label{eq:90013a}
\end{equation}
\begin{equation}
y(x)= \sum_{r=0}^{N+1} y_{2r}^{N+1-r}(x)+ \sum_{r=0}^{N} y_{2r+1}^{N-r}(x)  \label{eq:90013b}
\end{equation}
In the above,
\begin{eqnarray}
\bar{c}(0,n)  &=& \prod _{i_0=0}^{n-1}B_{2i_0+1} \label{eq:90014a}\\
\bar{c}(1,n) &=&  \sum_{i_0=0}^{n} \left\{ A_{2i_0} \prod _{i_1=0}^{i_0-1}B_{2i_1+1} \prod _{i_2=i_0}^{n-1}B_{2i_2+2} \right\} 
\label{eq:90014b}\\
\bar{c}(\tau ,n) &=& \sum_{i_0=0}^{n} \left\{A_{2i_0}\prod _{i_1=0}^{i_0-1} B_{2i_1+1} 
\prod _{k=1}^{\tau -1} \left( \sum_{i_{2k}= i_{2(k-1)}}^{n} A_{2i_{2k}+k}\prod _{i_{2k+1}=i_{2(k-1)}}^{i_{2k}-1}B_{2i_{2k+1}+(k+1)}\right)\right. \nonumber\\
&&\times \left. \prod _{i_{2\tau}=i_{2(\tau -1)}}^{n-1} B_{2i_{2\tau }+(\tau +1)} \right\} 
 \label{eq:90014c} 
\end{eqnarray}
and
\begin{eqnarray}
y_0^m(x) &=& c_0 x^{\lambda } \sum_{i_0=0}^{m} \left\{ \prod _{i_1=0}^{i_0-1}B_{2i_1+1} \right\} x^{2i_0 } \label{eq:90015a}\\
y_1^m(x) &=& c_0 x^{\lambda } \sum_{i_0=0}^{m}\left\{ A_{2i_0} \prod _{i_1=0}^{i_0-1}B_{2i_1+1}  \sum_{i_2=i_0}^{m} \left\{ \prod _{i_3=i_0}^{i_2-1}B_{2i_3+2} \right\}\right\} x^{2i_2+1 } \label{eq:90015b}\\
y_{\tau }^m(x) &=& c_0 x^{\lambda } \sum_{i_0=0}^{m} \left\{A_{2i_0}\prod _{i_1=0}^{i_0-1} B_{2i_1+1} 
\prod _{k=1}^{\tau -1} \left( \sum_{i_{2k}= i_{2(k-1)}}^{m} A_{2i_{2k}+k}\prod _{i_{2k+1}=i_{2(k-1)}}^{i_{2k}-1}B_{2i_{2k+1}+(k+1)}\right) \right. \nonumber\\
&& \times  \left. \sum_{i_{2\tau} = i_{2(\tau -1)}}^{m} \left( \prod _{i_{2\tau +1}=i_{2(\tau -1)}}^{i_{2\tau}-1} B_{2i_{2\tau +1}+(\tau +1)} \right) \right\} x^{2i_{2\tau}+\tau }\hspace{1cm}\mathrm{where}\;\tau \geq 2 \hspace{1.5cm}
\label{eq:90015c} 
\end{eqnarray}
\end{enumerate}
\end{theorem}
\subsubsection{The case of $\omega $, $\Omega $ as fixed values and $\varepsilon $, $\mu $, $\nu $ as free variables}
 
Put $n= j+1$ in (\ref{eq:9009b}) and use the condition $B_{j+1}=0$ for $\Omega $. We obtain two possible fixed values for $\Omega $ such as
\begin{subequations}
\begin{equation}
\Omega = -\mu j \label{eq:90016a}
\end{equation}
\begin{equation}
\Omega = \mu \left( \nu -j-1 \right)  \label{eq:90016b}
\end{equation}
\end{subequations}
\paragraph{The case of $ \Omega = -\mu j $}

Take (\ref{eq:90016a}) into (\ref{eq:9009a}) and (\ref{eq:9009b}).
\begin{subequations}
\begin{equation}
A_n =-\frac{\varepsilon }{\mu } \frac{\left( n-j-\omega \right)}{\left( n+1 \right)} \label{eq:90017a}
\end{equation}
\begin{equation}
B_n = \frac{1}{\mu } \frac{\left( n -j-1 \right) \left( n -j-\nu \right)}{\left( n+1 \right)}  \label{eq:90017b}
\end{equation}
\end{subequations}
According to (\ref{eq:90010}), $c_{j+1}=0$ is clearly an algebraic equation in $\omega $ of degree $j+1$ and thus has $j+1$ zeros denoted them by $\omega _j^m$ eigenvalues where $m = 0,1,2, \cdots, j$. They can be arranged in the following order: $\omega _j^0 < \omega _j^1 < \omega _j^2 < \cdots < \omega _j^j$.

In (\ref{eq:90011b}), (\ref{eq:90012b}), (\ref{eq:90013b}) and (\ref{eq:90015a})--(\ref{eq:90015c}) replace $c_0$, $\lambda $ and $x$ by $1$, $\Omega /\mu $ and $z$. Substitute (\ref{eq:90017a}) and (\ref{eq:90017b}) into (\ref{eq:90014a})--(\ref{eq:90015c}).

As $B_{1}= c_{1}=0$, take the new (\ref{eq:90014b}) into (\ref{eq:90011a}) putting $j=0$. Substitute the new (\ref{eq:90015a}) with $\Omega  =0$ into (\ref{eq:90011b}) putting $j=0$. 

As $B_{2N+2}= c_{2N+2}=0$, take the new (\ref{eq:90014a})--(\ref{eq:90014c}) into (\ref{eq:90012a}) putting $j=2N+1$. Substitute the new 
(\ref{eq:90015a})--(\ref{eq:90015c}) with $\Omega = -\mu \left( 2N+1 \right) $ into (\ref{eq:90012b}) putting $j=2N+1$ and $\omega = \omega q_{2N+1}^m$.

As $B_{2N+3}= c_{2N+3}=0$, take the new (\ref{eq:90014a})--(\ref{eq:90014c}) into (\ref{eq:90013a}) putting $j=2N+2$. Substitute the new 
(\ref{eq:90015a})--(\ref{eq:90015c}) with $\Omega = -\mu \left( 2N+2 \right)$ into (\ref{eq:90013b}) putting $j=2N+2$ and $\omega = \omega _{2N+2}^m$.

After the replacement process, we obtain the independent solution of the GCH equation about $x=\infty $. The solution is as follows.
\begin{remark}
The power series expansion of the GCH equation of the first kind for the first species complete polynomial using 3TRF about $x=\infty $ for $ \Omega = -\mu j$ where $j \in \mathbb{N}_{0}$ is given by
\begin{enumerate} 
\item As $\Omega =0$ and $\omega = \omega _0^0=0$,

The eigenfunction is given by
\begin{eqnarray}
y(z) &=& Q_p^{(i,1)}W_{0,0} \left( \mu ,\varepsilon ,\nu ,\Omega =0, \omega = \omega _0^0 =0; z=\frac{1}{x}, \sigma = \frac{2}{\mu }z^2; \xi = -\frac{\varepsilon }{\mu }z \right)\nonumber\\
&=& 1 \label{eq:90018}
\end{eqnarray}
\item As $\Omega = -\mu \left( 2N+1 \right)$ where $N \in \mathbb{N}_{0}$,

An algebraic equation of degree $2N+2$ for the determination of $\omega $ is given by
\begin{equation}
0 = \sum_{r=0}^{N+1}\bar{c}\left( 2r, N+1-r; 2N+1,\omega \right)\label{eq:90019a}
\end{equation}
The eigenvalue of $\omega $ is written by $\omega _{2N+1}^m$ where $m = 0,1,2,\cdots,2N+1 $; $\omega _{2N+1}^0 < \omega _{2N+1}^1 < \cdots < \omega _{2N+1}^{2N+1}$. Its eigenfunction is given by
\begin{eqnarray} 
y(z) &=& Q_p^{(i,1)}W_{2N+1,m} \left( \mu ,\varepsilon ,\nu ,\Omega =-\mu \left( 2N+1 \right), \omega = \omega _{2N+1}^m; z=\frac{1}{x}, \sigma = \frac{2}{\mu }z^2; \xi = -\frac{\varepsilon }{\mu }z \right) \nonumber\\
&=& z^{-2N-1} \left\{ \sum_{r=0}^{N} y_{2r}^{N-r}\left( 2N+1,\omega _{2N+1}^m;z\right)+ \sum_{r=0}^{N} y_{2r+1}^{N-r}\left( 2N+1,\omega _{2N+1}^m;z\right) \right\}  
 \label{eq:90019b}
\end{eqnarray}
\item As $\Omega = -\mu \left( 2N+2 \right)$ where $N \in \mathbb{N}_{0}$,

An algebraic equation of degree $2N+3$ for the determination of $\omega $ is given by
\begin{eqnarray}
0  = \sum_{r=0}^{N+1}\bar{c}\left( 2r+1, N+1-r; 2N+2,\omega \right)\label{eq:90020a}
\end{eqnarray}
The eigenvalue of $\omega $ is written by $\omega _{2N+2}^m$ where $m = 0,1,2,\cdots,2N+2 $; $\omega _{2N+2}^0 < \omega _{2N+2}^1 < \cdots < \omega _{2N+2}^{2N+2}$. Its eigenfunction is given by
\begin{eqnarray} 
y(z) &=& Q_p^{(i,1)}W_{2N+2,m} \left( \mu ,\varepsilon ,\nu ,\Omega =-\mu \left( 2N+2 \right), \omega = \omega _{2N+2}^m; z=\frac{1}{x}, \sigma = \frac{2}{\mu }z^2; \xi = -\frac{\varepsilon }{\mu }z \right) \nonumber\\
&=& z^{-2N-2} \left\{ \sum_{r=0}^{N+1} y_{2r}^{N+1-r}\left( 2N+2,\omega _{2N+2}^m;z\right) + \sum_{r=0}^{N} y_{2r+1}^{N-r}\left( 2N+2,\omega _{2N+2}^m;z\right) \right\} 
 \hspace{1.5cm}\label{eq:90020b}
\end{eqnarray}
In the above,
\begin{eqnarray}
\bar{c}(0,n;j,\omega )  &=& \frac{\left( -\frac{j}{2}\right)_{n}\left(  \frac{1}{2}-\frac{j}{2}-\frac{\nu }{2} \right)_{n}}{\left( 1 \right)_{n}} \left( \frac{2}{\mu } \right)^{n}\label{eq:90021a}\\
\bar{c}(1,n;j,\omega ) &=& \left( -\frac{\varepsilon }{\mu } \right) \sum_{i_0=0}^{n}\frac{ \left( i_0 -\frac{j}{2}-\frac{\omega }{2}\right) }{\left( i_0+\frac{1}{2} \right)} \frac{\left( -\frac{j}{2}\right)_{i_0} \left(  \frac{1}{2}-\frac{j}{2}-\frac{\nu }{2} \right)_{i_0}}{\left( 1 \right)_{i_0}} \nonumber\\
&&\times  \frac{\left( \frac{1}{2}-\frac{j}{2} \right)_{n} \left(  1-\frac{j}{2}-\frac{\nu }{2} \right)_n \left( \frac{3}{2} \right)_{i_0}}{\left( \frac{1}{2}-\frac{j}{2}\right)_{i_0}\left(  1-\frac{j}{2}-\frac{\nu }{2} \right)_{i_0} \left( \frac{3}{2} \right)_{n}} \left( \frac{2}{\mu } \right)^{n }  
\label{eq:90021b}\\
\bar{c}(\tau ,n;j,\omega ) &=& \left( -\frac{\varepsilon }{\mu } \right)^{\tau } \sum_{i_0=0}^{n}\frac{ \left( i_0 -\frac{j}{2}-\frac{\omega }{2}\right) }{\left( i_0+\frac{1}{2} \right)} \frac{\left( -\frac{j}{2}\right)_{i_0} \left(  \frac{1}{2}-\frac{j}{2}-\frac{\nu }{2} \right)_{i_0}}{\left( 1 \right)_{i_0}}  \nonumber\\
&&\times  \prod_{k=1}^{\tau -1} \left( \sum_{i_k = i_{k-1}}^{n} \frac{\left( i_k+ \frac{k}{2}-\frac{j}{2}-\frac{\omega }{2} \right) }{\left( i_k+\frac{k}{2}+\frac{1}{2} \right)} \right.   \left. \frac{\left( \frac{k}{2}-\frac{j}{2}\right)_{i_k} \left( \frac{k}{2}+ \frac{1}{2}-\frac{j}{2}-\frac{\nu }{2} \right)_{i_{k}}\left( \frac{k}{2}+1 \right)_{i_{k-1}} }{\left( \frac{k}{2}-\frac{j}{2}\right)_{i_{k-1}} \left( \frac{k}{2}+ \frac{1}{2}-\frac{j}{2}-\frac{\nu }{2} \right)_{i_{k-1}}\left( \frac{k}{2}+1 \right)_{i_k}} \right) \nonumber\\ 
&&\times \frac{\left( \frac{\tau }{2} -\frac{j}{2}\right)_{n} \left( \frac{\tau }{2}+\frac{1}{2}-\frac{j}{2}-\frac{\nu }{2} \right)_{n}\left( \frac{\tau }{2}+1 \right)_{i_{\tau -1}} }{\left( \frac{\tau }{2}-\frac{j}{2}\right)_{i_{\tau -1}} \left( \frac{\tau }{2}+\frac{1}{2}-\frac{j}{2}-\frac{\nu }{2} \right)_{i_{\tau -1}}\left( \frac{\tau }{2}+1 \right)_{n} } \left( \frac{2}{\mu } \right)^{n } \label{eq:90021c} 
\end{eqnarray}
\begin{eqnarray}
y_0^m(j,\omega ;z) &=& \sum_{i_0=0}^{m} \frac{\left( -\frac{j}{2}\right)_{i_0}\left(  \frac{1}{2}-\frac{j}{2}-\frac{\nu }{2} \right)_{i_0}}{\left( 1 \right)_{i_0}} \sigma ^{i_0} \label{eq:90022a}\\
y_1^m(j,\omega ;z) &=& \left\{\sum_{i_0=0}^{m} \frac{ \left( i_0 -\frac{j}{2}-\frac{\omega }{2}\right) }{\left( i_0+\frac{1}{2} \right)} \frac{\left( -\frac{j}{2}\right)_{i_0} \left(  \frac{1}{2}-\frac{j}{2}-\frac{\nu }{2} \right)_{i_0}}{\left( 1 \right)_{i_0}} \right. \nonumber\\
&&\times \left. \sum_{i_1 = i_0}^{m} \frac{\left( \frac{1}{2}-\frac{j}{2} \right)_{i_1} \left( 1-\frac{j}{2}-\frac{\nu }{2} \right)_{i_1}\left( \frac{3}{2} \right)_{i_0} }{\left( \frac{1}{2}-\frac{j}{2} \right)_{i_0} \left( 1-\frac{j}{2}-\frac{\nu }{2} \right)_{i_0}\left( \frac{3}{2} \right)_{i_1}} \sigma ^{i_1}\right\}\xi 
\label{eq:90022b}\\
y_{\tau }^m(j,\omega ;z) &=& \left\{ \sum_{i_0=0}^{m} \frac{ \left( i_0 -\frac{j}{2}-\frac{\omega }{2}\right) }{\left( i_0+\frac{1}{2} \right)} \frac{\left( -\frac{j}{2}\right)_{i_0} \left(  \frac{1}{2}-\frac{j}{2}-\frac{\nu }{2} \right)_{i_0}}{\left( 1 \right)_{i_0}} \right.\nonumber\\
&&\times \prod_{k=1}^{\tau -1} \left( \sum_{i_k = i_{k-1}}^{m} \frac{\left( i_k+ \frac{k}{2} -\frac{j}{2}-\frac{\omega }{2}\right) }{\left( i_k+\frac{k}{2}+\frac{1}{2} \right)} \right.   \left. \frac{\left( \frac{k}{2}-\frac{j}{2}\right)_{i_k} \left( \frac{k}{2}+ \frac{1}{2}-\frac{j}{2}-\frac{\nu }{2} \right)_{i_{k}}\left( \frac{k}{2}+1 \right)_{i_{k-1}} }{\left( \frac{k}{2}-\frac{j}{2}\right)_{i_{k-1}} \left( \frac{k}{2}+ \frac{1}{2}-\frac{j}{2}-\frac{\nu }{2} \right)_{i_{k-1}}\left( \frac{k}{2}+1 \right)_{i_k}} \right) \nonumber\\
&&\times  \left. \sum_{i_{\tau } = i_{\tau -1}}^{m} \frac{\left( \frac{\tau }{2} -\frac{j}{2}\right)_{i_{\tau}} \left( \frac{\tau }{2}+\frac{1}{2}-\frac{j}{2}-\frac{\nu }{2} \right)_{i_{\tau}}\left( \frac{\tau }{2}+1 \right)_{i_{\tau -1}} }{\left( \frac{\tau }{2}-\frac{j}{2}\right)_{i_{\tau -1}} \left( \frac{\tau }{2}+\frac{1}{2}-\frac{j}{2}-\frac{\nu }{2} \right)_{i_{\tau -1}}\left( \frac{\tau }{2}+1 \right)_{i_{\tau }} } \sigma ^{i_{\tau }}\right\} \xi ^{\tau } \label{eq:90022c} 
\end{eqnarray}
\end{enumerate}
\end{remark}
\paragraph{The case of $ \Omega = \mu \left( \nu -j-1 \right) $}

Take (\ref{eq:90016b}) into (\ref{eq:9009a}) and (\ref{eq:9009b}).
\begin{subequations}
\begin{equation}
A_n =-\frac{\varepsilon }{\mu } \frac{\left( n-j-1-\omega +\nu \right)}{\left( n+1 \right)} \label{eq:90023a}
\end{equation}
\begin{equation}
B_n = \frac{1}{\mu } \frac{\left( n-j-1 \right) \left( n-j-2+\nu \right)}{\left( n+1 \right)}  \label{eq:90023b}
\end{equation}
\end{subequations}
According to (\ref{eq:90010}), $c_{j+1}=0$ is clearly an algebraic equation in $\omega $ of degree $j+1$ and thus has $j+1$ zeros denoted them by $\omega _j^m$ eigenvalues where $m = 0,1,2, \cdots, j$. They can be arranged in the following order: $\omega _j^0 < \omega _j^1 < \omega _j^2 < \cdots < \omega _j^j$.

In (\ref{eq:90011b}), (\ref{eq:90012b}), (\ref{eq:90013b}) and (\ref{eq:90015a})--(\ref{eq:90015c}) replace $c_0$, $\lambda $ and $x$ by $1$, $\Omega /\mu $ and $z$. Substitute (\ref{eq:90023a}) and (\ref{eq:90023b}) into (\ref{eq:90014a})--(\ref{eq:90015c}).

As $B_{1}= c_{1}=0$, take the new (\ref{eq:90014b}) into (\ref{eq:90011a}) putting $j=0$. Substitute the new (\ref{eq:90015a}) with $\Omega = \mu \left( \nu -1 \right)$ into (\ref{eq:90011b}) putting $j=0$. 

As $B_{2N+2}= c_{2N+2}=0$, take the new (\ref{eq:90014a})--(\ref{eq:90014c}) into (\ref{eq:90012a}) putting $j=2N+1$. Substitute the new 
(\ref{eq:90015a})--(\ref{eq:90015c}) with $\Omega = \mu \left( \nu -2N-2 \right) $ into (\ref{eq:90012b}) putting $j=2N+1$ and $\omega = \omega _{2N+1}^m$.

As $B_{2N+3}= c_{2N+3}=0$, take the new (\ref{eq:90014a})--(\ref{eq:90014c}) into (\ref{eq:90013a}) putting $j=2N+2$. Substitute the new 
(\ref{eq:90015a})--(\ref{eq:90015c}) with $\Omega = \mu \left( \nu -2N-3 \right) $ into (\ref{eq:90013b}) putting $j=2N+2$ and $\omega = \omega _{2N+2}^m$.

After the replacement process, we obtain the independent solution of the GCH equation about $x=\infty $. The solution is as follows.
\begin{remark}
The power series expansion of the GCH equation of the first kind for the first species complete polynomial using 3TRF about $x=\infty $ for $\Omega = \mu \left( \nu -j-1 \right)$ where $j\in \mathbb{N}_{0}$ is given by
\begin{enumerate} 
\item As $\Omega =\mu \left( \nu -1 \right)$ and $\omega =\omega _0^0= \nu -1 $,

The eigenfunction is given by
\begin{eqnarray}
y(z) &=& Q_p^{(i,2)}W_{0,0} \left( \mu ,\varepsilon ,\nu ,\Omega =\mu \left( \nu -1 \right), \omega = \omega _0^0 = \nu -1; z=\frac{1}{x}, \sigma = \frac{2}{\mu }z^2; \xi = -\frac{\varepsilon }{\mu }z \right) \nonumber\\
&=& z^{\nu -1} \label{eq:90024}
\end{eqnarray}
\item As $\Omega =\mu \left( \nu -2N-2 \right)$ where $N \in \mathbb{N}_{0}$,

An algebraic equation of degree $2N+2$ for the determination of $\omega $ is given by
\begin{equation}
0 = \sum_{r=0}^{N+1}\bar{c}\left( 2r, N+1-r; 2N+1,\omega \right)\label{eq:90025a}
\end{equation}
The eigenvalue of $\omega $ is written by $\omega _{2N+1}^m$ where $m = 0,1,2,\cdots,2N+1 $; $\omega _{2N+1}^0 < \omega _{2N+1}^1 < \cdots < \omega _{2N+1}^{2N+1}$. Its eigenfunction is given by
\begin{eqnarray} 
y(z) &=& Q_p^{(i,2)}W_{2N+1,m} \left( \mu ,\varepsilon ,\nu ,\Omega =\mu \left( \nu -2N-2 \right), \omega = \omega _{2N+1}^m; z=\frac{1}{x}, \sigma = \frac{2}{\mu }z^2; \xi = -\frac{\varepsilon }{\mu }z \right) \nonumber\\
&=& z^{\nu -2N-2} \left\{ \sum_{r=0}^{N} y_{2r}^{N-r}\left( 2N+1,\omega _{2N+1}^m;z\right)+ \sum_{r=0}^{N} y_{2r+1}^{N-r}\left( 2N+1,\omega _{2N+1}^m;z\right) \right\}   
 \label{eq:90025b}
\end{eqnarray}
\item As $\Omega =\mu \left( \nu -2N-3 \right)$ where $N \in \mathbb{N}_{0}$,

An algebraic equation of degree $2N+3$ for the determination of $\omega $ is given by
\begin{eqnarray}
0  = \sum_{r=0}^{N+1}\bar{c}\left( 2r+1, N+1-r; 2N+2,\omega \right)\label{eq:90026a}
\end{eqnarray}
The eigenvalue of $\omega $ is written by $\omega _{2N+2}^m$ where $m = 0,1,2,\cdots,2N+2 $; $\omega _{2N+2}^0 < \omega _{2N+2}^1 < \cdots < \omega _{2N+2}^{2N+2}$. Its eigenfunction is given by
\begin{eqnarray} 
y(z) &=& Q_p^{(i,2)}W_{2N+2,m} \left( \mu ,\varepsilon ,\nu ,\Omega =\mu \left( \nu -2N-3 \right), \omega = \omega _{2N+2}^m; z=\frac{1}{x}, \sigma = \frac{2}{\mu }z^2; \xi = -\frac{\varepsilon }{\mu }z \right) \nonumber\\
&=& z^{\nu -2N-3} \left\{ \sum_{r=0}^{N+1} y_{2r}^{N+1-r}\left( 2N+2,\omega _{2N+2}^m;z\right) + \sum_{r=0}^{N} y_{2r+1}^{N-r}\left( 2N+2,\omega _{2N+2}^m;z\right) \right\} 
 \label{eq:90026b}
\end{eqnarray}
In the above,
\begin{eqnarray}
\bar{c}(0,n;j,\omega )  &=& \frac{\left( -\frac{j}{2}\right)_{n}\left(  -\frac{1}{2}-\frac{j}{2}+\frac{\nu }{2} \right)_{n}}{\left( 1 \right)_{n}} \left( \frac{2}{\mu } \right)^{n}\label{eq:90027a}\\
\bar{c}(1,n;j,\omega ) &=& \left( -\frac{\varepsilon }{\mu } \right) \sum_{i_0=0}^{n}\frac{ \left( i_0 -\frac{1}{2}-\frac{j}{2}-\frac{\omega }{2}+\frac{\nu }{2}\right) }{\left( i_0+\frac{1}{2} \right)} \frac{\left( -\frac{j}{2}\right)_{i_0} \left( -\frac{1}{2}-\frac{j}{2}+\frac{\nu }{2} \right)_{i_0}}{\left( 1 \right)_{i_0}} \nonumber\\
&&\times  \frac{\left( \frac{1}{2}-\frac{j}{2} \right)_{n} \left( -\frac{j}{2}+\frac{\nu }{2} \right)_n \left( \frac{3}{2} \right)_{i_0}}{\left( \frac{1}{2}-\frac{j}{2}\right)_{i_0}\left( -\frac{j}{2}+\frac{\nu }{2} \right)_{i_0} \left( \frac{3}{2} \right)_{n}} \left( \frac{2}{\mu } \right)^{n }  
\label{eq:90027b}\\
\bar{c}(\tau ,n;j,\omega ) &=& \left( -\frac{\varepsilon }{\mu } \right)^{\tau } \sum_{i_0=0}^{n}\frac{ \left( i_0 -\frac{1}{2}-\frac{j}{2}-\frac{\omega }{2}+\frac{\nu }{2}\right) }{\left( i_0+\frac{1}{2} \right)} \frac{\left( -\frac{j}{2}\right)_{i_0} \left( -\frac{1}{2}-\frac{j}{2}+\frac{\nu }{2} \right)_{i_0}}{\left( 1 \right)_{i_0}}  \nonumber\\
&&\times  \prod_{k=1}^{\tau -1} \left( \sum_{i_k = i_{k-1}}^{n} \frac{\left( i_k+ \frac{k}{2}-\frac{1}{2}-\frac{j}{2}-\frac{\omega }{2}+\frac{\nu }{2} \right) }{\left( i_k+\frac{k}{2}+\frac{1}{2} \right)} \right. \nonumber\\
&&\times   \left. \frac{\left( \frac{k}{2}-\frac{j}{2}\right)_{i_k} \left( \frac{k}{2}- \frac{1}{2}-\frac{j}{2}+\frac{\nu }{2} \right)_{i_{k}}\left( \frac{k}{2}+1 \right)_{i_{k-1}} }{\left( \frac{k}{2}-\frac{j}{2}\right)_{i_{k-1}} \left( \frac{k}{2}-\frac{1}{2}-\frac{j}{2}+\frac{\nu }{2} \right)_{i_{k-1}}\left( \frac{k}{2}+1 \right)_{i_k}} \right) \nonumber\\ 
&&\times \frac{\left( \frac{\tau }{2} -\frac{j}{2}\right)_{n} \left( \frac{\tau }{2}-\frac{1}{2}-\frac{j}{2}+\frac{\nu }{2} \right)_{n}\left( \frac{\tau }{2}+1 \right)_{i_{\tau -1}} }{\left( \frac{\tau }{2}-\frac{j}{2}\right)_{i_{\tau -1}} \left( \frac{\tau }{2}-\frac{1}{2}-\frac{j}{2}+\frac{\nu }{2} \right)_{i_{\tau -1}}\left( \frac{\tau }{2}+1 \right)_{n} } \left( \frac{2}{\mu } \right)^{n } \label{eq:90027c} 
\end{eqnarray}
\begin{eqnarray}
y_0^m(j,\omega ;z) &=& \sum_{i_0=0}^{m} \frac{\left( -\frac{j}{2}\right)_{i_0}\left( -\frac{1}{2}-\frac{j}{2}+\frac{\nu }{2} \right)_{i_0}}{\left( 1 \right)_{i_0}} \sigma ^{i_0} \label{eq:90028a}\\
y_1^m(j,\omega ;z) &=& \left\{\sum_{i_0=0}^{m} \frac{ \left( i_0 -\frac{1}{2}-\frac{j}{2}-\frac{\omega }{2}+\frac{\nu }{2}\right) }{\left( i_0+\frac{1}{2} \right)} \frac{\left( -\frac{j}{2}\right)_{i_0} \left( -\frac{1}{2}-\frac{j}{2}+\frac{\nu }{2} \right)_{i_0}}{\left( 1 \right)_{i_0}} \right. \nonumber\\
&&\times \left. \sum_{i_1 = i_0}^{m} \frac{\left( \frac{1}{2}-\frac{j}{2} \right)_{i_1} \left( -\frac{j}{2}+\frac{\nu }{2} \right)_{i_1}\left( \frac{3}{2} \right)_{i_0} }{\left( \frac{1}{2}-\frac{j}{2} \right)_{i_0} \left( -\frac{j}{2}+\frac{\nu }{2} \right)_{i_0}\left( \frac{3}{2} \right)_{i_1}} \sigma ^{i_1}\right\}\xi 
\label{eq:90028b}\\
y_{\tau }^m(j,\omega ;z) &=& \left\{ \sum_{i_0=0}^{m} \frac{ \left( i_0 -\frac{1}{2}-\frac{j}{2}-\frac{\omega }{2}+\frac{\nu }{2}\right) }{\left( i_0+\frac{1}{2} \right)} \frac{\left( -\frac{j}{2}\right)_{i_0} \left( -\frac{1}{2}-\frac{j}{2}+\frac{\nu }{2} \right)_{i_0}}{\left( 1 \right)_{i_0}} \right.\nonumber\\
&&\times \prod_{k=1}^{\tau -1} \left( \sum_{i_k = i_{k-1}}^{m} \frac{\left( i_k+ \frac{k}{2} -\frac{1}{2}-\frac{j}{2}-\frac{\omega }{2}+\frac{\nu }{2}\right) }{\left( i_k+\frac{k}{2}+\frac{1}{2} \right)} \right.   \left. \frac{\left( \frac{k}{2}-\frac{j}{2}\right)_{i_k} \left( \frac{k}{2}- \frac{1}{2}-\frac{j}{2}+\frac{\nu }{2} \right)_{i_{k}}\left( \frac{k}{2}+1 \right)_{i_{k-1}} }{\left( \frac{k}{2}-\frac{j}{2}\right)_{i_{k-1}} \left( \frac{k}{2}-\frac{1}{2}-\frac{j}{2}+\frac{\nu }{2} \right)_{i_{k-1}}\left( \frac{k}{2}+1 \right)_{i_k}} \right) \nonumber\\
&&\times  \left. \sum_{i_{\tau } = i_{\tau -1}}^{m} \frac{\left( \frac{\tau }{2} -\frac{j}{2}\right)_{i_{\tau}} \left( \frac{\tau }{2}-\frac{1}{2}-\frac{j}{2}+\frac{\nu }{2} \right)_{i_{\tau}}\left( \frac{\tau }{2}+1 \right)_{i_{\tau -1}} }{\left( \frac{\tau }{2}-\frac{j}{2}\right)_{i_{\tau -1}} \left( \frac{\tau }{2}-\frac{1}{2}-\frac{j}{2}+\frac{\nu }{2} \right)_{i_{\tau -1}}\left( \frac{\tau }{2}+1 \right)_{i_{\tau }} } \sigma ^{i_{\tau }}\right\} \xi ^{\tau } \label{eq:90028c} 
\end{eqnarray}
\end{enumerate}
\end{remark}
\subsubsection{The case of $\nu $, $\omega $ as fixed values and $\varepsilon $, $\mu $, $\Omega $ as free variables}

Put $n= j+1$ in (\ref{eq:9009b}) and use the condition $B_{j+1}=0$ for a fixed value $\nu $.  
\begin{equation}
\nu = \frac{\Omega }{\mu } +j+1 \label{eq:90029}
\end{equation}
 Take (\ref{eq:90029}) into (\ref{eq:9009b}).
\begin{equation}
B_n = \frac{1}{\mu } \frac{\left( n-j-1 \right) \left( n-1 +\frac{\Omega }{\mu } \right)}{\left( n+1 \right)} \label{eq:90030}
\end{equation}

According to (\ref{eq:90010}), $c_{j+1}=0$ is clearly an algebraic equation in $\omega $ of degree $j+1$ and thus has $j+1$ zeros denoted them by $\omega _j^m$ eigenvalues where $m = 0,1,2, \cdots, j$. They can be arranged in the following order: $\omega _j^0 < \omega _j^1 < \omega _j^2 < \cdots < \omega _j^j$.

In (\ref{eq:90011b}), (\ref{eq:90012b}), (\ref{eq:90013b}) and (\ref{eq:90015a})--(\ref{eq:90015c}) replace $c_0$, $\lambda $ and $x$ by $1$, $\Omega /\mu $ and $z$. Substitute (\ref{eq:9009a}) and (\ref{eq:90030}) into (\ref{eq:90014a})--(\ref{eq:90015c}).

As $B_{1}= c_{1}=0$, take the new (\ref{eq:90014b}) into (\ref{eq:90011a}) putting $j=0$. Substitute the new (\ref{eq:90015a}) with $\nu = \frac{\Omega }{\mu } +1 $ into (\ref{eq:90011b}) putting $j=0$. 

As $B_{2N+2}= c_{2N+2}=0$, take the new (\ref{eq:90014a})--(\ref{eq:90014c}) into (\ref{eq:90012a}) putting $j=2N+1$. Substitute the new 
(\ref{eq:90015a})--(\ref{eq:90015c}) with $\nu = \frac{\Omega }{\mu } +2N+2 $ into (\ref{eq:90012b}) putting $j=2N+1$ and $\omega = \omega _{2N+1}^m$.

As $B_{2N+3}= c_{2N+3}=0$, take the new (\ref{eq:90014a})--(\ref{eq:90014c}) into (\ref{eq:90013a}) putting $j=2N+2$. Substitute the new 
(\ref{eq:90015a})--(\ref{eq:90015c}) with $\nu = \frac{\Omega }{\mu } +2N+3 $ into (\ref{eq:90013b}) putting $j=2N+2$ and $\omega = \omega _{2N+2}^m$.

After the replacement process, we obtain the independent solution of the GCH equation about $x=\infty $. The solution is as follows.
\begin{remark}
The power series expansion of the GCH equation of the first kind for the first species complete polynomial using 3TRF about $x=\infty $ for $\nu = \frac{\Omega }{\mu } +j+1 $ where $j\in \mathbb{N}_{0}$ is given by
\begin{enumerate} 
\item As $\nu = \frac{\Omega }{\mu } +1 $ and $\omega =\omega _0^0= \frac{\Omega }{\mu } $,

The eigenfunction is given by
\begin{eqnarray}
y(z) &=& Q_p^{(i,3)}W_{0,0} \left( \mu ,\varepsilon ,\nu = \frac{\Omega }{\mu } +1, \Omega , \omega = \omega _0^0 = \frac{\Omega }{\mu } ; z=\frac{1}{x}, \sigma = \frac{2}{\mu }z^2; \xi = -\frac{\varepsilon }{\mu }z \right) \nonumber\\
&=& z^{\frac{\Omega }{\mu }} \label{eq:90031}
\end{eqnarray}
\item As $\nu = \frac{\Omega }{\mu } +2N+2 $ where $N \in \mathbb{N}_{0}$,

An algebraic equation of degree $2N+2$ for the determination of $\omega $ is given by
\begin{equation}
0 = \sum_{r=0}^{N+1}\bar{c}\left( 2r, N+1-r; 2N+1,\omega \right)\label{eq:90032a}
\end{equation}
The eigenvalue of $\omega $ is written by $\omega _{2N+1}^m$ where $m = 0,1,2,\cdots,2N+1 $; $\omega _{2N+1}^0 < \omega _{2N+1}^1 < \cdots < \omega _{2N+1}^{2N+1}$. Its eigenfunction is given by
\begin{eqnarray} 
y(z) &=& Q_p^{(i,3)}W_{2N+1,m} \left( \mu ,\varepsilon ,\nu = \frac{\Omega }{\mu } +2N+2, \Omega , \omega = \omega _{2N+1}^m ; z=\frac{1}{x}, \sigma = \frac{2}{\mu }z^2; \xi = -\frac{\varepsilon }{\mu }z \right) \nonumber\\
&=& \sum_{r=0}^{N} y_{2r}^{N-r}\left( 2N+1,\omega _{2N+1}^m;z\right)+ \sum_{r=0}^{N} y_{2r+1}^{N-r}\left( 2N+1,\omega _{2N+1}^m;z\right) 
 \label{eq:90032b}
\end{eqnarray}
\item As $\nu = \frac{\Omega }{\mu } +2N+3 $ where $N \in \mathbb{N}_{0}$,

An algebraic equation of degree $2N+3$ for the determination of $\omega $ is given by
\begin{eqnarray}
0  = \sum_{r=0}^{N+1}\bar{c}\left( 2r+1, N+1-r; 2N+2,\omega \right)\label{eq:90033a}
\end{eqnarray}
The eigenvalue of $\omega $ is written by $\omega _{2N+2}^m$ where $m = 0,1,2,\cdots,2N+2 $; $\omega _{2N+2}^0 < \omega _{2N+2}^1 < \cdots < \omega _{2N+2}^{2N+2}$. Its eigenfunction is given by
\begin{eqnarray} 
y(z) &=& Q_p^{(i,3)}W_{2N+2,m} \left( \mu ,\varepsilon ,\nu = \frac{\Omega }{\mu } +2N+3, \Omega , \omega = \omega _{2N+2}^m ; z=\frac{1}{x}, \sigma = \frac{2}{\mu }z^2; \xi = -\frac{\varepsilon }{\mu }z \right) \nonumber\\
&=& \sum_{r=0}^{N+1} y_{2r}^{N+1-r}\left( 2N+2,\omega _{2N+2}^m;z\right) + \sum_{r=0}^{N} y_{2r+1}^{N-r}\left( 2N+2,\omega _{2N+2}^m;z\right)  
 \label{eq:90033b}
\end{eqnarray}
In the above,
\begin{eqnarray}
\bar{c}(0,n;j,\omega )  &=& \frac{\left( -\frac{j}{2}\right)_{n}\left( \frac{\Omega }{2\mu }  \right)_{n}}{\left( 1 \right)_{n}} \left( \frac{2}{\mu } \right)^{n}\label{eq:90034a}\\
\bar{c}(1,n;j,\omega ) &=& \left( -\frac{\varepsilon }{\mu } \right) \sum_{i_0=0}^{n}\frac{ \left( i_0 +\frac{\Omega }{2\mu }-\frac{\omega }{2} \right) }{\left( i_0+\frac{1}{2} \right)} \frac{\left( -\frac{j}{2}\right)_{i_0} \left( \frac{\Omega }{2\mu } \right)_{i_0}}{\left( 1 \right)_{i_0}} \nonumber\\
&&\times  \frac{\left( \frac{1}{2}-\frac{j}{2} \right)_{n} \left( \frac{1}{2}+\frac{\Omega }{2\mu } \right)_n \left( \frac{3}{2} \right)_{i_0}}{\left( \frac{1}{2}-\frac{j}{2}\right)_{i_0}\left( \frac{1}{2}+\frac{\Omega }{2\mu } \right)_{i_0} \left( \frac{3}{2} \right)_{n}} \left( \frac{2}{\mu } \right)^{n }  
\label{eq:90034b}\\
\bar{c}(\tau ,n;j,\omega ) &=& \left( -\frac{\varepsilon }{\mu } \right)^{\tau } \sum_{i_0=0}^{n}\frac{ \left( i_0 +\frac{\Omega }{2\mu }-\frac{\omega }{2} \right) }{\left( i_0+\frac{1}{2} \right)} \frac{\left( -\frac{j}{2}\right)_{i_0} \left( \frac{\Omega }{2\mu } \right)_{i_0}}{\left( 1 \right)_{i_0}}  \nonumber\\
&&\times  \prod_{k=1}^{\tau -1} \left( \sum_{i_k = i_{k-1}}^{n} \frac{\left( i_k+ \frac{k}{2}+\frac{\Omega }{2\mu }-\frac{\omega }{2} \right) }{\left( i_k+\frac{k}{2}+\frac{1}{2} \right)} \right.   \left. \frac{\left( \frac{k}{2}-\frac{j}{2}\right)_{i_k} \left( \frac{k}{2}+\frac{\Omega }{2\mu } \right)_{i_{k}}\left( \frac{k}{2}+1 \right)_{i_{k-1}} }{\left( \frac{k}{2}-\frac{j}{2}\right)_{i_{k-1}} \left( \frac{k}{2}+\frac{\Omega }{2\mu } \right)_{i_{k-1}}\left( \frac{k}{2}+1 \right)_{i_k}} \right) \nonumber\\ 
&&\times \frac{\left( \frac{\tau }{2} -\frac{j}{2}\right)_{n} \left( \frac{\tau }{2}+\frac{\Omega }{2\mu } \right)_{n}\left( \frac{\tau }{2}+1 \right)_{i_{\tau -1}} }{\left( \frac{\tau }{2}-\frac{j}{2}\right)_{i_{\tau -1}} \left( \frac{\tau }{2}+\frac{\Omega }{2\mu } \right)_{i_{\tau -1}}\left( \frac{\tau }{2}+1 \right)_{n} } \left( \frac{2}{\mu } \right)^{n } \label{eq:90034c} 
\end{eqnarray}
\begin{eqnarray}
y_0^m(j,\omega ;z) &=& z^{\frac{\Omega }{\mu }} \sum_{i_0=0}^{m} \frac{\left( -\frac{j}{2}\right)_{i_0}\left( \frac{\Omega }{2\mu } \right)_{i_0}}{\left( 1 \right)_{i_0}} \sigma ^{i_0} \label{eq:90035a}\\
y_1^m(j,\omega ;z) &=& z^{\frac{\Omega }{\mu }} \left\{\sum_{i_0=0}^{m} \frac{ \left( i_0 +\frac{\Omega }{2\mu }-\frac{\omega }{2} \right) }{\left( i_0+\frac{1}{2} \right)} \frac{\left( -\frac{j}{2}\right)_{i_0} \left( \frac{\Omega }{2\mu } \right)_{i_0}}{\left( 1 \right)_{i_0}} \right. \nonumber\\
&&\times \left. \sum_{i_1 = i_0}^{m} \frac{\left( \frac{1}{2}-\frac{j}{2} \right)_{i_1} \left( \frac{1}{2}+\frac{\Omega }{2\mu } \right)_{i_1}\left( \frac{3}{2} \right)_{i_0} }{\left( \frac{1}{2}-\frac{j}{2} \right)_{i_0} \left( \frac{1}{2}+\frac{\Omega }{2\mu } \right)_{i_0}\left( \frac{3}{2} \right)_{i_1}} \sigma ^{i_1}\right\}\xi 
\label{eq:90035b}\\
y_{\tau }^m(j,\omega ;z) &=& z^{\frac{\Omega }{\mu }} \left\{ \sum_{i_0=0}^{m} \frac{ \left( i_0 +\frac{\Omega }{2\mu }-\frac{\omega }{2} \right) }{\left( i_0+\frac{1}{2} \right)} \frac{\left( -\frac{j}{2}\right)_{i_0} \left( \frac{\Omega }{2\mu } \right)_{i_0}}{\left( 1 \right)_{i_0}} \right.\nonumber\\
&&\times \prod_{k=1}^{\tau -1} \left( \sum_{i_k = i_{k-1}}^{m} \frac{\left( i_k+ \frac{k}{2} +\frac{\Omega }{2\mu }-\frac{\omega }{2} \right) }{\left( i_k+\frac{k}{2}+\frac{1}{2} \right)} \right.   \left. \frac{\left( \frac{k}{2}-\frac{j}{2}\right)_{i_k} \left( \frac{k}{2}+\frac{\Omega }{2\mu } \right)_{i_{k}}\left( \frac{k}{2}+1 \right)_{i_{k-1}} }{\left( \frac{k}{2}-\frac{j}{2}\right)_{i_{k-1}} \left( \frac{k}{2}+\frac{\Omega }{2\mu } \right)_{i_{k-1}}\left( \frac{k}{2}+1 \right)_{i_k}} \right) \nonumber\\
&&\times  \left. \sum_{i_{\tau } = i_{\tau -1}}^{m} \frac{\left( \frac{\tau }{2} -\frac{j}{2}\right)_{i_{\tau}} \left( \frac{\tau }{2}+\frac{\Omega }{2\mu } \right)_{i_{\tau}}\left( \frac{\tau }{2}+1 \right)_{i_{\tau -1}} }{\left( \frac{\tau }{2}-\frac{j}{2}\right)_{i_{\tau -1}} \left( \frac{\tau }{2}+\frac{\Omega }{2\mu } \right)_{i_{\tau -1}}\left( \frac{\tau }{2}+1 \right)_{i_{\tau }} } \sigma ^{i_{\tau }}\right\} \xi ^{\tau } \label{eq:90035c} 
\end{eqnarray}
\end{enumerate}
\end{remark}
\subsection{The second species complete polynomial of the GCH equation using 3TRF}

For the second species complete polynomials using 3TRF and R3TRF, we need a condition which is defined by
\begin{equation}
B_{j}=B_{j+1}= A_{j}=0\hspace{1cm}\mathrm{where}\;j \in \mathbb{N}_{0}    
 \label{eq:90036}
\end{equation}
\begin{theorem}
In chapter 1, the general expression of a function $y(x)$ for the second species complete polynomial using 3-term recurrence formula is given by
\begin{enumerate} 
\item As $B_1=A_0=0$,
\begin{equation}
y(x) = y_{0}^{0}(x) \label{eq:90037a}
\end{equation}
\item As $B_{2N+1}=B_{2N+2}=A_{2N+1}=0$  where $N \in \mathbb{N}_{0}$,
\begin{equation}
y(x)= \sum_{r=0}^{N} y_{2r}^{N-r}(x)+ \sum_{r=0}^{N} y_{2r+1}^{N-r}(x)  \label{eq:90037b}
\end{equation}
\item As $B_{2N+2}=B_{2N+3}=A_{2N+2}=0$  where $N \in \mathbb{N}_{0}$,
\begin{equation}
y(x)= \sum_{r=0}^{N+1} y_{2r}^{N+1-r}(x)+ \sum_{r=0}^{N} y_{2r+1}^{N-r}(x)  \label{eq:90037c}
\end{equation}
In the above,
\begin{eqnarray}
y_0^m(x) &=& c_0 x^{\lambda } \sum_{i_0=0}^{m} \left\{ \prod _{i_1=0}^{i_0-1}B_{2i_1+1} \right\} x^{2i_0 } \label{eq:90038a}\\
y_1^m(x) &=& c_0 x^{\lambda } \sum_{i_0=0}^{m}\left\{ A_{2i_0} \prod _{i_1=0}^{i_0-1}B_{2i_1+1}  \sum_{i_2=i_0}^{m} \left\{ \prod _{i_3=i_0}^{i_2-1}B_{2i_3+2} \right\}\right\} x^{2i_2+1 } \label{eq:90038b}\\
y_{\tau }^m(x) &=& c_0 x^{\lambda } \sum_{i_0=0}^{m} \left\{A_{2i_0}\prod _{i_1=0}^{i_0-1} B_{2i_1+1} 
\prod _{k=1}^{\tau -1} \left( \sum_{i_{2k}= i_{2(k-1)}}^{m} A_{2i_{2k}+k}\prod _{i_{2k+1}=i_{2(k-1)}}^{i_{2k}-1}B_{2i_{2k+1}+(k+1)}\right) \right. \nonumber\\
&& \times \left. \sum_{i_{2\tau} = i_{2(\tau -1)}}^{m} \left( \prod _{i_{2\tau +1}=i_{2(\tau -1)}}^{i_{2\tau}-1} B_{2i_{2\tau +1}+(\tau +1)} \right) \right\} x^{2i_{2\tau}+\tau }\hspace{1cm}\mathrm{where}\;\tau \geq 2 \hspace{1.5cm}
\label{eq:90038c}  
\end{eqnarray}
\end{enumerate}
\end{theorem}
In (\ref{eq:90016a}) and (\ref{eq:90016b}), for a condition $B_{j+1}=0$, two possible fixed values for $\Omega  $ are $\Omega = -\mu j$ and $\Omega = \mu \left( \nu -j-1 \right)$.
Put $n= j$ with (\ref{eq:90016a}) in (\ref{eq:9009b}) and use the condition $B_{j}=0$ for a fixed value $\nu $.  
\begin{equation}
\nu = 0
\label{eq:90039}
\end{equation}
Substitute (\ref{eq:90016a}) and (\ref{eq:90039}) into (\ref{eq:9009a}). Put $n= j$ in the new (\ref{eq:9009a}) and use the condition $A_{j}=0$ for a fixed value $\omega $.  
\begin{equation}
\omega  = 0
\label{eq:90040}
\end{equation}
Take (\ref{eq:90016a}), (\ref{eq:90039}) and (\ref{eq:90040}) into (\ref{eq:9009a}) and (\ref{eq:9009b}).
\begin{subequations}
\begin{equation}
A_n = -\frac{\varepsilon }{\mu } \frac{(n-j) }{ (n+1)} \label{eq:90041a}
\end{equation}
\begin{equation}
B_n = \frac{1}{\mu } \frac{(n-1-j)(n-j)}{ (n+1)} \label{eq:90041b}
\end{equation}
\end{subequations}
For the case of $\Omega = \mu \left( \nu -j-1 \right) $, putting $n= j$ with (\ref{eq:90016b}) in (\ref{eq:9009b}) and use the condition $B_{j}=0$ for a fixed value $\nu $, we obtain $\nu =2$. 
Substitute $\Omega = \mu \left( 1-j \right) $ and $\nu =2$ into (\ref{eq:9009a}). Putting $n= j$ in the new (\ref{eq:9009a}) and use the condition $A_{j}=0$ for a fixed value $\omega $, we derive $\omega =1 $.
Take $\Omega = \mu \left( 1-j \right) $, $\nu =2$ and $\omega = 1 $ into (\ref{eq:9009a}) and (\ref{eq:9009b}). The new (\ref{eq:9009a}) is same as (\ref{eq:90041a}). And the new(\ref{eq:9009b}) is equivalent to (\ref{eq:90041b}).
Therefore, the second species complete polynomial $y(z)$ divided by $z^{\Omega /\mu  }$ of the GCH equation about $x=\infty $ as $\Omega = -\mu j$ and $\nu =\omega = 0$ is equal to the independent solution $y(z)$ divided by $z^{\Omega /\mu  }$ of the GCH equation as $\Omega = \mu \left( 1-j \right) $, $\nu =2$ and $\omega =1 $.
  
In (\ref{eq:90037a})--(\ref{eq:90038c}) replace $c_0$, $\lambda $ and $x$ by $1$, $\Omega /\mu $ and $z$. 
Substitute (\ref{eq:90041a}) and (\ref{eq:90041b}) into (\ref{eq:90038a})--(\ref{eq:90038c}).

\underline{(1) The case of $\Omega = -\mu j$ and $\nu =\omega = 0$,} 
 
As $B_1=A_0=0$, substitute the new (\ref{eq:90038a}) with $\Omega  =0$ into (\ref{eq:90037a}) putting $j=0$. 
As $B_{2N+1}=B_{2N+2}=A_{2N+1}=0$, substitute the new (\ref{eq:90038a})--(\ref{eq:90038c}) with $\Omega =-\mu \left( 2N+1 \right)$ into (\ref{eq:90037b}) putting $j=2N+1$.
As $B_{2N+2}=B_{2N+3}=A_{2N+2}=0$, substitute the new (\ref{eq:90038a})--(\ref{eq:90038c})  with $\Omega =-\mu \left( 2N+2 \right)$ into (\ref{eq:90037c}) putting $j=2N+2$.

\underline{(2) The case of $\Omega = \mu \left( 1-j \right) $, $\nu =2$ and $\omega =1 $,} 
 
As $B_1=A_0=0$, substitute the new (\ref{eq:90038a}) with $\Omega = \mu $ into (\ref{eq:90037a}) putting $j=0$. 
As $B_{2N+1}=B_{2N+2}=A_{2N+1}=0$, substitute the new (\ref{eq:90038a})--(\ref{eq:90038c}) with $\Omega = -2\mu N $ into (\ref{eq:90037b}) putting $j=2N+1$.
As $B_{2N+2}=B_{2N+3}=A_{2N+2}=0$, substitute the new (\ref{eq:90038a})--(\ref{eq:90038c})  with $\Omega = -\mu \left( 2N+1 \right)$ into (\ref{eq:90037c}) putting $j=2N+2$.

After the replacement process, we obtain the independent solution of the GCH equation about $x=\infty $. The solution is as follows.
\begin{remark}
The power series expansion of the GCH equation of the first kind for the second species complete polynomial using 3TRF about $x=\infty $ is given by
\begin{enumerate} 
\item As $\Omega  =\nu =\omega =0$, 

Its eigenfunction is given by
\begin{eqnarray}
y(z) &=& Q_p^{(i,1)}W_0 \left( \mu ,\varepsilon ,\nu =0,\Omega =0, \omega = 0; z=\frac{1}{x}, \sigma = \frac{2}{\mu }z^2; \xi = -\frac{\varepsilon }{\mu }z \right) \nonumber\\
&=& 1 \label{eq:90042a}
\end{eqnarray}
\item As $\Omega =-\mu \left( 2N+1 \right)$, $\nu =\omega =0$ where $N \in \mathbb{N}_{0}$, 

Its eigenfunction is given by
\begin{eqnarray}
y(z) &=& Q_p^{(i,1)}W_{2N+1} \left( \mu ,\varepsilon ,\nu =0,\Omega =-\mu \left( 2N+1 \right), \omega = 0; z=\frac{1}{x}, \sigma = \frac{2}{\mu }z^2; \xi = -\frac{\varepsilon }{\mu }z \right) \nonumber\\
 &=& z^{-2N-1} \left\{ \sum_{r=0}^{N} y_{2r}^{N-r}\left( 2N+1;z\right) + \sum_{r=0}^{N} y_{2r+1}^{N-r}\left( 2N+1;z\right) \right\} \label{eq:90042b}
\end{eqnarray}
\item As $\Omega =-\mu \left( 2N+2 \right)$, $\nu =\omega =0$ where $N \in \mathbb{N}_{0}$, 

Its eigenfunction is given by
\begin{eqnarray}
y(z) &=& Q_p^{(i,1)}W_{2N+2} \left( \mu ,\varepsilon ,\nu =0,\Omega =-\mu \left( 2N+2 \right), \omega = 0; z=\frac{1}{x}, \sigma = \frac{2}{\mu }z^2; \xi = -\frac{\varepsilon }{\mu }z \right) \nonumber\\
 &=& z^{-2N-2} \left\{ \sum_{r=0}^{N+1} y_{2r}^{N+1-r}\left( 2N+2;z\right) + \sum_{r=0}^{N} y_{2r+1}^{N-r}\left( 2N+2;z\right) \right\} \label{eq:90042c}
\end{eqnarray}

\item As $\Omega = \mu, \nu =2, \omega =1 $,  

Its eigenfunction is given by
\begin{eqnarray}
y(z) &=& Q_p^{(i,2)}W_0 \left( \mu ,\varepsilon ,\nu =2,\Omega =\mu, \omega = 1; z=\frac{1}{x}, \sigma = \frac{2}{\mu }z^2; \xi = -\frac{\varepsilon }{\mu }z \right) \nonumber\\
 &=& z \label{eq:90042d}
\end{eqnarray}
\item As $\Omega = -2\mu N , \nu =2, \omega =1 $ where $N \in \mathbb{N}_{0}$,

Its eigenfunction is given by
\begin{eqnarray}
y(z) &=& Q_p^{(i,2)}W_{2N+1} \left( \mu ,\varepsilon ,\nu =2,\Omega =-2\mu N, \omega = 1; z=\frac{1}{x}, \sigma = \frac{2}{\mu }z^2; \xi = -\frac{\varepsilon }{\mu }z \right) \nonumber\\
 &=& z^{ -2N } \left\{ \sum_{r=0}^{N} y_{2r}^{N-r}\left( 2N+1;z\right) + \sum_{r=0}^{N} y_{2r+1}^{N-r}\left( 2N+1;z\right) \right\} \label{eq:90042e}
\end{eqnarray}
\item As  $\Omega = -\mu \left( 2N+1 \right) , \nu =2, \omega =1 $ where $N \in \mathbb{N}_{0}$,

Its eigenfunction is given by
\begin{eqnarray}
y(z) &=& Q_p^{(i,2)}W_{2N+2}  \left( \mu ,\varepsilon ,\nu =2,\Omega =-\mu \left( 2N+1 \right), \omega = 1; z=\frac{1}{x}, \sigma = \frac{2}{\mu }z^2; \xi = -\frac{\varepsilon }{\mu }z \right) \nonumber\\
 &=& z^{ -2N-1} \left\{ \sum_{r=0}^{N+1} y_{2r}^{N+1-r}\left( 2N+2;z\right) + \sum_{r=0}^{N} y_{2r+1}^{N-r}\left( 2N+2;z\right) \right\} \label{eq:90042f}
\end{eqnarray}
In the above,
\begin{eqnarray}
y_0^m(j;z) &=& \sum_{i_0=0}^{m} \frac{\left( -\frac{j}{2}\right)_{i_0}\left(  \frac{1}{2}-\frac{j}{2}  \right)_{i_0}}{\left( 1 \right)_{i_0}} \sigma ^{i_0} \label{eq:90043a}\\
y_1^m(j;z) &=& \left\{\sum_{i_0=0}^{m} \frac{ \left( i_0 -\frac{j}{2}\right) }{\left( i_0+\frac{1}{2} \right)} \frac{\left( -\frac{j}{2}\right)_{i_0} \left(  \frac{1}{2}-\frac{j}{2} \right)_{i_0}}{\left( 1 \right)_{i_0}} \right.  \left. \sum_{i_1 = i_0}^{m} \frac{\left( \frac{1}{2}-\frac{j}{2} \right)_{i_1} \left( 1-\frac{j}{2} \right)_{i_1}\left( \frac{3}{2} \right)_{i_0} }{\left( \frac{1}{2}-\frac{j}{2} \right)_{i_0} \left( 1-\frac{j}{2} \right)_{i_0}\left( \frac{3}{2} \right)_{i_1}} \sigma ^{i_1}\right\}\xi 
 \hspace{1.5cm}\label{eq:90043b}\\
y_{\tau }^m(j;z) &=& \left\{ \sum_{i_0=0}^{m} \frac{ \left( i_0 -\frac{j}{2}\right) }{\left( i_0+\frac{1}{2} \right)} \frac{\left( -\frac{j}{2}\right)_{i_0} \left(  \frac{1}{2}-\frac{j}{2} \right)_{i_0}}{\left( 1 \right)_{i_0}} \right.\nonumber\\
&&\times \prod_{k=1}^{\tau -1} \left( \sum_{i_k = i_{k-1}}^{m} \frac{\left( i_k+ \frac{k}{2} -\frac{j}{2}\right) }{\left( i_k+\frac{k}{2}+\frac{1}{2} \right)} \right.   \left. \frac{\left( \frac{k}{2}-\frac{j}{2}\right)_{i_k} \left( \frac{k}{2}+ \frac{1}{2}-\frac{j}{2} \right)_{i_{k}}\left( \frac{k}{2}+1 \right)_{i_{k-1}} }{\left( \frac{k}{2}-\frac{j}{2}\right)_{i_{k-1}} \left( \frac{k}{2}+ \frac{1}{2}-\frac{j}{2} \right)_{i_{k-1}}\left( \frac{k}{2}+1 \right)_{i_k}} \right) \nonumber\\
&&\times  \left. \sum_{i_{\tau } = i_{\tau -1}}^{m} \frac{\left( \frac{\tau }{2} -\frac{j}{2}\right)_{i_{\tau}} \left( \frac{\tau }{2}+\frac{1}{2}-\frac{j}{2} \right)_{i_{\tau}}\left( \frac{\tau }{2}+1 \right)_{i_{\tau -1}} }{\left( \frac{\tau }{2}-\frac{j}{2}\right)_{i_{\tau -1}} \left( \frac{\tau }{2}+\frac{1}{2}-\frac{j}{2} \right)_{i_{\tau -1}}\left( \frac{\tau }{2}+1 \right)_{i_{\tau }} } \sigma ^{i_{\tau }}\right\} \xi ^{\tau } \label{eq:90043c} 
\end{eqnarray}
\end{enumerate}
\end{remark} 
\subsection{The first species complete polynomial of the GCH equation using R3TRF}
\begin{theorem}
In chapter 2, the general expression of a function $y(x)$ for the first species complete polynomial using reversible 3-term recurrence formula and its algebraic equation for the determination of an accessory parameter in $A_n$ term are given by
\begin{enumerate} 
\item As $B_1=0$,
\begin{equation}
0 =\bar{c}(0,1) \label{eq:90044a}
\end{equation}
\begin{equation}
y(x) = y_{0}^{0}(x) \label{eq:90044b}
\end{equation}
\item As $B_2=0$, 
\begin{equation}
0 = \bar{c}(0,2)+\bar{c}(1,0) \label{eq:90045a}
\end{equation}
\begin{equation}
y(x)= y_{0}^{1}(x) \label{eq:90045b}
\end{equation}
\item As $B_{2N+3}=0$ where $N \in \mathbb{N}_{0}$,
\begin{equation}
0  = \sum_{r=0}^{N+1}\bar{c}\left( r, 2(N-r)+3\right) \label{eq:90046a}
\end{equation}
\begin{equation}
y(x)= \sum_{r=0}^{N+1} y_{r}^{2(N+1-r)}(x) \label{eq:90046b}
\end{equation}
\item As $B_{2N+4}=0$ where$N \in \mathbb{N}_{0}$,
\begin{equation}
0  = \sum_{r=0}^{N+2}\bar{c}\left( r, 2(N+2-r)\right) \label{eq:90047a}
\end{equation}
\begin{equation}
y(x)=  \sum_{r=0}^{N+1} y_{r}^{2(N-r)+3}(x) \label{eq:90047b}
\end{equation}
In the above,
\begin{eqnarray}
\bar{c}(0,n) &=& \prod _{i_0=0}^{n-1}A_{i_0} \label{eq:90048a}\\
\bar{c}(1,n) &=& \sum_{i_0=0}^{n} \left\{ B_{i_0+1} \prod _{i_1=0}^{i_0-1}A_{i_1} \prod _{i_2=i_0}^{n-1}A_{i_2+2} \right\} \label{eq:90048b}\\
\bar{c}(\tau ,n) &=& \sum_{i_0=0}^{n} \left\{B_{i_0+1}\prod _{i_1=0}^{i_0-1} A_{i_1} 
\prod _{k=1}^{\tau -1} \left( \sum_{i_{2k}= i_{2(k-1)}}^{n} B_{i_{2k}+(2k+1)}\prod _{i_{2k+1}=i_{2(k-1)}}^{i_{2k}-1}A_{i_{2k+1}+2k}\right) \right. \nonumber\\
&&\times \left. \prod _{i_{2\tau} = i_{2(\tau -1)}}^{n-1} A_{i_{2\tau }+ 2\tau} \right\} 
 \label{eq:90048c}
\end{eqnarray}
and
\begin{eqnarray}
y_0^m(x) &=& c_0 x^{\lambda} \sum_{i_0=0}^{m} \left\{ \prod _{i_1=0}^{i_0-1}A_{i_1} \right\} x^{i_0 } \label{eq:90049a}\\
y_1^m(x) &=& c_0 x^{\lambda} \sum_{i_0=0}^{m}\left\{ B_{i_0+1} \prod _{i_1=0}^{i_0-1}A_{i_1}  \sum_{i_2=i_0}^{m} \left\{ \prod _{i_3=i_0}^{i_2-1}A_{i_3+2} \right\}\right\} x^{i_2+2 } \label{eq:90049b}\\
y_{\tau }^m(x) &=& c_0 x^{\lambda} \sum_{i_0=0}^{m} \left\{B_{i_0+1}\prod _{i_1=0}^{i_0-1} A_{i_1} 
\prod _{k=1}^{\tau -1} \left( \sum_{i_{2k}= i_{2(k-1)}}^{m} B_{i_{2k}+(2k+1)}\prod _{i_{2k+1}=i_{2(k-1)}}^{i_{2k}-1}A_{i_{2k+1}+2k}\right) \right. \nonumber\\
&&\times \left. \sum_{i_{2\tau} = i_{2(\tau -1)}}^{m} \left( \prod _{i_{2\tau +1}=i_{2(\tau -1)}}^{i_{2\tau}-1} A_{i_{2\tau +1}+ 2\tau} \right) \right\} x^{i_{2\tau}+2\tau }\hspace{1cm}\mathrm{where}\;\tau \geq 2
\label{eq:90049c}
\end{eqnarray}
\end{enumerate}
\end{theorem}
\subsubsection{The case of $\omega $, $\Omega $ as fixed values and $\varepsilon $, $\mu $, $\nu $ as free variables}

For the first species complete polynomial of the GCH equation about $x=\infty $ in section 10.2.1, there are two possible fixed values of $\Omega $ such as $\Omega = -\mu j$ and $ \mu \left( \nu -j-1 \right)$.
According to (\ref{eq:90010}), $c_{j+1}=0$ is a polynomial equation of degree $j+1$ for the determination of the accessory parameter $\omega $ and thus has $j+1$ zeros denoted them by $\omega _j^m$ eigenvalues where $m = 0,1,2, \cdots, j$. They can be arranged in the following order: $\omega _j^0 < \omega _j^1 < \omega _j^2 < \cdots < \omega _j^j$.

\paragraph{The case of $ \Omega = -\mu j $}
In (\ref{eq:90044b}), (\ref{eq:90045b}), (\ref{eq:90046b}), (\ref{eq:90047b}) and (\ref{eq:90049a})--(\ref{eq:90049c}) replace $c_0$, $\lambda $ and $x$ by $1$, $\Omega /\mu $ and $z$. Substitute (\ref{eq:90017a}) and (\ref{eq:90017b}) into (\ref{eq:90048a})--(\ref{eq:90049c}).

As $B_{1}= c_{1}=0$, take the new (\ref{eq:90048a}) into (\ref{eq:90044a}) putting $j=0$. Substitute the new (\ref{eq:90049a}) with $\Omega  =0$ into (\ref{eq:90044b}) putting $j=0$.

As $B_{2}= c_{2}=0$, take the new (\ref{eq:90048a}) and (\ref{eq:90048b}) into (\ref{eq:90045a}) putting $j=1$. Substitute the new (\ref{eq:90049a}) with $\Omega  =-\mu $ into (\ref{eq:90045b}) putting $j=1$ and $\omega =\omega _1^m$. 

As $B_{2N+3}= c_{2N+3}=0$, take the new (\ref{eq:90048a})--(\ref{eq:90048c}) into (\ref{eq:90046a}) putting $j=2N+2$. Substitute the new 
(\ref{eq:90049a})--(\ref{eq:90049c}) with $\Omega  =-\mu \left( 2N+2 \right) $ into (\ref{eq:90046b}) putting $j=2N+2$ and $\omega =\omega _{2N+2}^m$.

As $B_{2N+4}= c_{2N+4}=0$, take the new (\ref{eq:90048a})--(\ref{eq:90048c}) into (\ref{eq:90047a}) putting $j=2N+3$. Substitute the new 
(\ref{eq:90049a})--(\ref{eq:90049c}) with $\Omega  =-\mu \left( 2N+3 \right) $ into (\ref{eq:90047b}) putting $j=2N+3$ and $\omega =\omega _{2N+3}^m$.

After the replacement process, we obtain the independent solution of the GCH equation about $x=\infty $. The solution is as follows.
\begin{remark}
The power series expansion of the GCH equation of the first kind for the first species complete polynomial using R3TRF about $x=\infty $ for $\Omega = -\mu j$ where $j \in \mathbb{N}_{0}$ is given by
\begin{enumerate} 
\item As $\Omega =0$ and $\omega =\omega _0^0=0$,

The eigenfunction is given by
\begin{eqnarray}
y(z) &=& Q_p^{(i,1)}W_{0,0}^R \left( \mu ,\varepsilon ,\nu ,\Omega =0, \omega = \omega _0^0 =0; z=\frac{1}{x}, \vartheta = \frac{1}{\mu }z^2; \xi = -\frac{\varepsilon }{\mu }z \right) \nonumber\\
&=& 1\label{eq:90050}
\end{eqnarray}
\item As $\Omega =-\mu $,

An algebraic equation of degree 2 for the determination of $\omega $ is given by
\begin{equation}
0 = \mu \nu  + \varepsilon ^2 \omega \left( \omega +1\right)  \label{eq:90051a}
\end{equation}
The eigenvalue of $\omega $ is written by $\omega _1^m$ where $m = 0,1 $; $\omega _{1}^0 < \omega _{1}^1$. Its eigenfunction is given by
\begin{eqnarray}
y(z) &=& Q_p^{(i,1)}W_{1,m}^R \left( \mu ,\varepsilon ,\nu ,\Omega =-\mu, \omega = \omega _1^m; z=\frac{1}{x}, \vartheta = \frac{1}{\mu }z^2; \xi = -\frac{\varepsilon }{\mu }z \right)\nonumber\\
&=&  z^{-1}\left\{ 1 -\left( 1+\omega _1^m \right)\xi \right\} \label{eq:90051b}  
\end{eqnarray}
\item As $\Omega =-\mu \left( 2N+2 \right) $ where $N \in \mathbb{N}_{0}$,

An algebraic equation of degree $2N+3$ for the determination of $\omega $ is given by
\begin{equation}
0 = \sum_{r=0}^{N+1}\bar{c}\left( r, 2(N-r)+3; 2N+2,\omega \right)  \label{eq:90052a}
\end{equation}
The eigenvalue of $\omega $ is written by $\omega _{2N+2}^m$ where $m = 0,1,2,\cdots,2N+2 $; $\omega _{2N+2}^0 < \omega _{2N+2}^1 < \cdots < \omega _{2N+2}^{2N+2}$. Its eigenfunction is given by 
\begin{eqnarray} 
y(z) &=& Q_p^{(i,1)}W_{2N+2,m}^R \left( \mu ,\varepsilon ,\nu ,\Omega =-\mu \left( 2N+2 \right), \omega = \omega _{2N+2}^m; z=\frac{1}{x}, \vartheta = \frac{1}{\mu }z^2; \xi = -\frac{\varepsilon }{\mu }z \right)\nonumber\\
&=& z^{-2N-2} \sum_{r=0}^{N+1} y_{r}^{2(N+1-r)}\left( 2N+2, \omega _{2N+2}^m; z \right) 
\label{eq:90052b} 
\end{eqnarray}
\item As $\Omega =-\mu \left( 2N+3 \right) $ where $N \in \mathbb{N}_{0}$,

An algebraic equation of degree $2N+4$ for the determination of $\omega $ is given by
\begin{equation}  
0 = \sum_{r=0}^{N+2}\bar{c}\left( r, 2(N+2-r); 2N+3,\omega \right) \label{eq:90053a}
\end{equation}
The eigenvalue of $\omega $ is written by $\omega _{2N+3}^m$ where $m = 0,1,2,\cdots,2N+3 $; $\omega _{2N+3}^0 < \omega _{2N+3}^1 < \cdots < \omega _{2N+3}^{2N+3}$. Its eigenfunction is given by
\begin{eqnarray} 
y(z) &=& Q_p^{(i,1)}W_{2N+3,m}^R \left( \mu ,\varepsilon ,\nu ,\Omega =-\mu \left( 2N+3 \right), \omega = \omega _{2N+3}^m; z=\frac{1}{x}, \vartheta = \frac{1}{\mu }z^2; \xi = -\frac{\varepsilon }{\mu }z \right)\nonumber\\
&=& z^{-2N-3} \sum_{r=0}^{N+1} y_{r}^{2(N-r)+3} \left( 2N+3,\omega _{2N+3}^m;z\right) \label{eq:90053b}
\end{eqnarray}
In the above,
\begin{eqnarray}
\bar{c}(0,n;j,\omega )  &=& \frac{\left( -j-\omega  \right)_{n} }{\left( 1 \right)_{n} }\left( -\frac{\varepsilon }{\mu  }\right)^n \label{eq:90054a}\\
\bar{c}(1,n;j,\omega ) &=& \left( \frac{1}{\mu }\right) \sum_{i_0=0}^{n}\frac{\left( i_0 -j\right)\left( i_0+1-j-\nu \right) }{\left( i_0+2 \right)} \frac{ \left( -j-\omega  \right)_{i_0}}{\left( 1 \right)_{i_0}}  \nonumber\\
&&\times  \frac{ \left( 2-j-\omega  \right)_{n} \left( 3\right)_{i_0} }{\left( 2-j-\omega \right)_{i_0} \left( 3 \right)_{n}} \left(  -\frac{\varepsilon }{\mu } \right)^n \label{eq:90054b}\\
\bar{c}(\tau ,n;j,\omega ) &=& \left( \frac{1}{\mu }\right)^{\tau} \sum_{i_0=0}^{n} \frac{\left( i_0 -j\right)\left( i_0+1-j-\nu \right) }{\left( i_0+2 \right)} \frac{ \left( -j-\omega  \right)_{i_0}}{\left( 1 \right)_{i_0}}  \nonumber\\
&&\times \prod_{k=1}^{\tau -1} \left( \sum_{i_k = i_{k-1}}^{n} \frac{\left( i_k+ 2k-j\right)\left( i_k+2k+1-j-\nu \right) }{\left( i_k+2k+2 \right)} \right.   \left. \frac{ \left( 2k-j-\omega \right)_{i_k} \left( 2k+1 \right)_{i_{k-1}}}{\left( 2k-j-\omega  \right)_{i_{k-1}} \left( 2k+1 \right)_{i_k}} \right) \nonumber\\
&&\times \frac{ \left( 2\tau -j-\omega \right)_{n} \left( 2\tau +1 \right)_{i_{\tau -1}}}{\left( 2\tau -j-\omega  \right)_{i_{\tau -1}} \left( 2\tau +1 \right)_{n}}\left( -\frac{\varepsilon }{\mu } \right)^n  \label{eq:90054c} 
\end{eqnarray}
\begin{eqnarray}
y_0^m(j,\omega ;z) &=& \sum_{i_0=0}^{m} \frac{\left( -j-\omega \right)_{i_0} }{\left( 1 \right)_{i_0} } \xi ^{i_0} \label{eq:90055a}\\
y_1^m(j,\omega ;z) &=& \left\{\sum_{i_0=0}^{m}\frac{\left( i_0 -j\right)\left( i_0+1-j-\nu \right) }{\left( i_0+2 \right)} \frac{ \left( -j-\omega \right)_{i_0} }{\left( 1 \right)_{i_0}} \right.  \nonumber\\
&&\times \left. \sum_{i_1 = i_0}^{m} \frac{ \left( 2-j-\omega \right)_{i_1} \left( 3\right)_{i_0} }{\left( 2-j-\omega  \right)_{i_0} \left( 3 \right)_{i_1}} \xi ^{i_1}\right\} \vartheta    
 \label{eq:90055b}\\
y_{\tau }^m(j,\omega ;z) &=& \left\{ \sum_{i_0=0}^{m} \frac{\left( i_0 -j\right)\left( i_0+1-j-\nu \right) }{\left( i_0+2 \right)} \frac{ \left( -j-\omega \right)_{i_0} }{\left( 1 \right)_{i_0}} \right.\nonumber\\
&&\times \prod_{k=1}^{\tau -1} \left( \sum_{i_k = i_{k-1}}^{m}  \frac{\left( i_k+ 2k-j\right)\left( i_k+2k+1-j-\nu \right) }{\left( i_k+2k+2 \right)} \right.   \left. \frac{ \left( 2k-j-\omega  \right)_{i_k} \left( 2k+1 \right)_{i_{k-1}}}{\left( 2k-j-\omega \right)_{i_{k-1}} \left( 2k+1 \right)_{i_k}} \right) \nonumber\\
&&\times \left. \sum_{i_{\tau } = i_{\tau -1}}^{m}  \frac{ \left(  2\tau -j-\omega  \right)_{i_{\tau }} \left( 2\tau +1 \right)_{i_{\tau -1}}}{\left( 2\tau -j-\omega \right)_{i_{\tau -1}} \left( 2\tau +1 \right)_{i_{\tau }}} \xi^{i_{\tau }}\right\} \vartheta ^{\tau } \label{eq:90055c} 
\end{eqnarray}
\end{enumerate}
\end{remark}
\paragraph{The case of $ \Omega = \mu \left( \nu -j-1 \right) $}

In (\ref{eq:90044b}), (\ref{eq:90045b}), (\ref{eq:90046b}), (\ref{eq:90047b}) and (\ref{eq:90049a})--(\ref{eq:90049c}) replace $c_0$, $\lambda $ and $x$ by $1$, $\Omega /\mu $ and $z$. Substitute (\ref{eq:90023a}) and (\ref{eq:90023b}) into (\ref{eq:90048a})--(\ref{eq:90049c}).

As $B_{1}= c_{1}=0$, take the new (\ref{eq:90048a}) into (\ref{eq:90044a}) putting $j=0$. Substitute the new (\ref{eq:90049a}) with $\Omega  = \mu \left( \nu -1 \right) $ into (\ref{eq:90044b}) putting $j=0$.

As $B_{2}= c_{2}=0$, take the new (\ref{eq:90048a}) and (\ref{eq:90048b}) into (\ref{eq:90045a}) putting $j=1$. Substitute the new (\ref{eq:90049a}) with $\Omega  = \mu \left( \nu -2 \right) $ into (\ref{eq:90045b}) putting $j=1$ and $\omega =\omega _1^m$. 

As $B_{2N+3}= c_{2N+3}=0$, take the new (\ref{eq:90048a})--(\ref{eq:90048c}) into (\ref{eq:90046a}) putting $j=2N+2$. Substitute the new 
(\ref{eq:90049a})--(\ref{eq:90049c}) with $\Omega  = \mu \left( \nu -2N-3 \right) $ into (\ref{eq:90046b}) putting $j=2N+2$ and $\omega =\omega _{2N+2}^m$.

As $B_{2N+4}= c_{2N+4}=0$, take the new (\ref{eq:90048a})--(\ref{eq:90048c}) into (\ref{eq:90047a}) putting $j=2N+3$. Substitute the new 
(\ref{eq:90049a})--(\ref{eq:90049c}) with $\Omega  = \mu \left( \nu -2N-4 \right) $ into (\ref{eq:90047b}) putting $j=2N+3$ and $\omega =\omega _{2N+3}^m$.

After the replacement process, we obtain the independent solution of the GCH equation about $x=\infty $. The solution is as follows.
\begin{remark}
The power series expansion of the GCH equation of the first kind for the first species complete polynomial using R3TRF about $x=\infty $ for $\Omega = \mu \left( \nu -j-1 \right)$ where $j \in \mathbb{N}_{0}$ is given by
\begin{enumerate} 
\item As $\Omega =\mu \left( \nu -1 \right)$ and $\omega =\omega _0^0=\nu -1$,

The eigenfunction is given by
\begin{eqnarray}
y(z) &=& Q_p^{(i,2)}W_{0,0}^R \left( \mu ,\varepsilon ,\nu ,\Omega =\mu \left( \nu -1 \right), \omega = \omega _0^0 =\nu -1; z=\frac{1}{x}, \vartheta = \frac{1}{\mu }z^2; \xi = -\frac{\varepsilon }{\mu }z \right) \nonumber\\
&=& z^{\nu -1}\label{eq:90056}
\end{eqnarray}
\item As $\Omega =\mu \left( \nu -2 \right) $,

An algebraic equation of degree 2 for the determination of $\omega $ is given by
\begin{equation}
0 = \mu \left( 2-\nu \right) + \varepsilon ^2 \left( \omega +1-\nu \right) \left( \omega +2-\nu \right)  \label{eq:90057a}
\end{equation}
The eigenvalue of $\omega $ is written by $\omega _1^m$ where $m = 0,1 $; $\omega _{1}^0 < \omega _{1}^1$. Its eigenfunction is given by
\begin{eqnarray}
y(z) &=& Q_p^{(i,2)}W_{1,m}^R \left( \mu ,\varepsilon ,\nu ,\Omega =\mu \left( \nu -2 \right), \omega = \omega _1^m; z=\frac{1}{x}, \vartheta = \frac{1}{\mu }z^2; \xi = -\frac{\varepsilon }{\mu }z \right)\nonumber\\
&=&  z^{\nu -2}\left\{ 1 -\left( 2-\nu +\omega _1^m \right)\xi \right\} \label{eq:90057b}  
\end{eqnarray}
\item As $\Omega =\mu \left( \nu -2N-3 \right) $ where $N \in \mathbb{N}_{0}$,

An algebraic equation of degree $2N+3$ for the determination of $\omega $ is given by
\begin{equation}
0 = \sum_{r=0}^{N+1}\bar{c}\left( r, 2(N-r)+3; 2N+2,\omega \right)  \label{eq:90058a}
\end{equation}
The eigenvalue of $\omega $ is written by $\omega _{2N+2}^m$ where $m = 0,1,2,\cdots,2N+2 $; $\omega _{2N+2}^0 < \omega _{2N+2}^1 < \cdots < \omega _{2N+2}^{2N+2}$. Its eigenfunction is given by 
\begin{eqnarray} 
y(z) &=& Q_p^{(i,2)}W_{2N+2,m}^R \left( \mu ,\varepsilon ,\nu ,\Omega =\mu \left( \nu -2N-3 \right), \omega = \omega _{2N+2}^m; z=\frac{1}{x}, \vartheta = \frac{1}{\mu }z^2; \xi = -\frac{\varepsilon }{\mu }z \right)\nonumber\\
&=& z^{\nu -2N-3} \sum_{r=0}^{N+1} y_{r}^{2(N+1-r)}\left( 2N+2, \omega _{2N+2}^m; z \right) 
\label{eq:90058b} 
\end{eqnarray}
\item As $\Omega =\mu \left( \nu -2N-4 \right) $ where $N \in \mathbb{N}_{0}$,

An algebraic equation of degree $2N+4$ for the determination of $\omega $ is given by
\begin{equation}  
0 = \sum_{r=0}^{N+2}\bar{c}\left( r, 2(N+2-r); 2N+3,\omega \right) \label{eq:90059a}
\end{equation}
The eigenvalue of $\omega $ is written by $\omega _{2N+3}^m$ where $m = 0,1,2,\cdots,2N+3 $; $\omega _{2N+3}^0 < \omega _{2N+3}^1 < \cdots < \omega _{2N+3}^{2N+3}$. Its eigenfunction is given by
\begin{eqnarray} 
y(z) &=& Q_p^{(i,2)}W_{2N+3,m}^R \left( \mu ,\varepsilon ,\nu ,\Omega =\mu \left( \nu -2N-4 \right), \omega = \omega _{2N+3}^m; z=\frac{1}{x}, \vartheta = \frac{1}{\mu }z^2; \xi = -\frac{\varepsilon }{\mu }z \right)\nonumber\\
&=& z^{\nu -2N-4} \sum_{r=0}^{N+1} y_{r}^{2(N-r)+3} \left( 2N+3,\omega _{2N+3}^m;z\right) \label{eq:90059b}
\end{eqnarray}
In the above,
\begin{eqnarray}
\bar{c}(0,n;j,\omega )  &=& \frac{\left( -1-j+\nu -\omega  \right)_{n} }{\left( 1 \right)_{n} }\left( -\frac{\varepsilon }{\mu  }\right)^n \label{eq:90060a}\\
\bar{c}(1,n;j,\omega ) &=& \left( \frac{1}{\mu }\right) \sum_{i_0=0}^{n}\frac{\left( i_0 -j\right)\left( i_0-1-j+\nu \right) }{\left( i_0+2 \right)} \frac{ \left( -1-j+\nu -\omega  \right)_{i_0}}{\left( 1 \right)_{i_0}}  \nonumber\\
&&\times  \frac{ \left( 1-j+\nu -\omega \right)_{n} \left( 3\right)_{i_0} }{\left( 1-j+\nu -\omega \right)_{i_0} \left( 3 \right)_{n}} \left(  -\frac{\varepsilon }{\mu } \right)^n \label{eq:90060b}\\
\bar{c}(\tau ,n;j,\omega ) &=& \left( \frac{1}{\mu }\right)^{\tau} \sum_{i_0=0}^{n} \frac{\left( i_0 -j\right)\left( i_0-1-j+\nu \right) }{\left( i_0+2 \right)} \frac{ \left( -1-j+\nu -\omega  \right)_{i_0}}{\left( 1 \right)_{i_0}}  \nonumber\\
&&\times \prod_{k=1}^{\tau -1} \left( \sum_{i_k = i_{k-1}}^{n} \frac{\left( i_k+ 2k-j\right)\left( i_k+2k-1-j+\nu \right) }{\left( i_k+2k+2 \right)} \right. \left. \frac{ \left( 2k-1-j+\nu -\omega \right)_{i_k} \left( 2k+1 \right)_{i_{k-1}}}{\left( 2k-1-j+\nu -\omega \right)_{i_{k-1}} \left( 2k+1 \right)_{i_k}} \right) \nonumber\\
&&\times \frac{ \left( 2\tau -1-j+\nu -\omega \right)_{n} \left( 2\tau +1 \right)_{i_{\tau -1}}}{\left( 2\tau -1-j+\nu -\omega  \right)_{i_{\tau -1}} \left( 2\tau +1 \right)_{n}}\left( -\frac{\varepsilon }{\mu } \right)^n  \label{eq:90060c} 
\end{eqnarray}
\begin{eqnarray}
y_0^m(j,\omega ;z) &=& \sum_{i_0=0}^{m} \frac{\left( -1-j+\nu -\omega \right)_{i_0} }{\left( 1 \right)_{i_0} } \xi ^{i_0} \label{eq:90061a}\\
y_1^m(j,\omega ;z) &=& \left\{\sum_{i_0=0}^{m}\frac{\left( i_0 -j\right)\left( i_0-1-j+\nu \right) }{\left( i_0+2 \right)} \frac{ \left( -1-j+\nu -\omega \right)_{i_0} }{\left( 1 \right)_{i_0}} \right.  \nonumber\\
&&\times \left. \sum_{i_1 = i_0}^{m} \frac{ \left( 1-j+\nu -\omega \right)_{i_1} \left( 3\right)_{i_0} }{\left( 1-j+\nu -\omega \right)_{i_0} \left( 3 \right)_{i_1}} \xi ^{i_1}\right\} \vartheta    
 \label{eq:90061b}\\
y_{\tau }^m(j,\omega ;z) &=& \left\{ \sum_{i_0=0}^{m} \frac{\left( i_0 -j\right)\left( i_0-1-j+\nu \right) }{\left( i_0+2 \right)} \frac{ \left( -1-j+\nu -\omega \right)_{i_0} }{\left( 1 \right)_{i_0}} \right.\nonumber\\
&&\times \prod_{k=1}^{\tau -1} \left( \sum_{i_k = i_{k-1}}^{m}  \frac{\left( i_k+ 2k-j\right)\left( i_k+2k-1-j+\nu \right) }{\left( i_k+2k+2 \right)} \right. \left. \frac{ \left( 2k-1-j+\nu -\omega \right)_{i_k} \left( 2k+1 \right)_{i_{k-1}}}{\left( 2k-1-j+\nu -\omega  \right)_{i_{k-1}} \left( 2k+1 \right)_{i_k}} \right) \nonumber\\
&&\times \left. \sum_{i_{\tau } = i_{\tau -1}}^{m}  \frac{ \left(  2\tau -1-j+\nu -\omega  \right)_{i_{\tau }} \left( 2\tau +1 \right)_{i_{\tau -1}}}{\left( 2\tau -1-j+\nu -\omega \right)_{i_{\tau -1}} \left( 2\tau +1 \right)_{i_{\tau }}} \xi^{i_{\tau }}\right\} \vartheta ^{\tau } \label{eq:90061c} 
\end{eqnarray}
\end{enumerate}
\end{remark}
\subsubsection{The case of $\nu $, $\omega $ as fixed values and $\varepsilon $, $\mu $, $\Omega $ as free variables}
For the first species complete polynomial of the GCH equation about $x=\infty $ in section 10.2.1, there is a fixed values of $\nu $ such as $\Omega /\mu +j+1$.
According to (\ref{eq:90010}), $c_{j+1}=0$ is a polynomial equation of degree $j+1$ for the determination of the accessory parameter $\omega $ and thus has $j+1$ zeros denoted them by $\omega _j^m$ eigenvalues where $m = 0,1,2, \cdots, j$. They can be arranged in the following order: $\omega _j^0 < \omega _j^1 < \omega _j^2 < \cdots < \omega _j^j$.

In (\ref{eq:90044b}), (\ref{eq:90045b}), (\ref{eq:90046b}), (\ref{eq:90047b}) and (\ref{eq:90049a})--(\ref{eq:90049c}) replace $c_0$, $\lambda $ and $x$ by $1$, $\Omega /\mu $ and $z$. Substitute (\ref{eq:9009a}) and (\ref{eq:90030}) into (\ref{eq:90048a})--(\ref{eq:90049c}).

As $B_{1}= c_{1}=0$, take the new (\ref{eq:90048a}) into (\ref{eq:90044a}) putting $j=0$. Substitute the new (\ref{eq:90049a}) with $\nu = \Omega /\mu +1 $ into (\ref{eq:90044b}) putting $j=0$.

As $B_{2}= c_{2}=0$, take the new (\ref{eq:90048a}) and (\ref{eq:90048b}) into (\ref{eq:90045a}) putting $j=1$. Substitute the new (\ref{eq:90049a}) with $\nu = \Omega /\mu +2 $ into (\ref{eq:90045b}) putting $j=1$ and $\omega =\omega _1^m$. 

As $B_{2N+3}= c_{2N+3}=0$, take the new (\ref{eq:90048a})--(\ref{eq:90048c}) into (\ref{eq:90046a}) putting $j=2N+2$. Substitute the new 
(\ref{eq:90049a})--(\ref{eq:90049c}) with $\nu = \Omega /\mu +2N+3  $ into (\ref{eq:90046b}) putting $j=2N+2$ and $\omega =\omega _{2N+2}^m$.

As $B_{2N+4}= c_{2N+4}=0$, take the new (\ref{eq:90048a})--(\ref{eq:90048c}) into (\ref{eq:90047a}) putting $j=2N+3$. Substitute the new 
(\ref{eq:90049a})--(\ref{eq:90049c}) with $\nu = \Omega /\mu +2N+4 $ into (\ref{eq:90047b}) putting $j=2N+3$ and $\omega =\omega _{2N+3}^m$.

After the replacement process, we obtain the independent solution of the GCH equation about $x=\infty $. The solution is as follows.
\begin{remark}
The power series expansion of the GCH equation of the first kind for the first species complete polynomial using R3TRF about $x=\infty $ for $\nu  =\frac{\Omega }{\mu }+j+1$ where $j\in \mathbb{N}_{0}$ is given by
\begin{enumerate} 
\item As $\nu  =\frac{\Omega }{\mu } +1 $ and $\omega =\omega _0^0= \frac{\Omega }{\mu } $,

The eigenfunction is given by
\begin{eqnarray}
y(z) &=& Q_p^{(i,3)}W_{0,0}^R \left( \mu ,\varepsilon ,\nu = \frac{\Omega }{\mu } +1, \Omega , \omega = \omega _0^0 = \frac{\Omega }{\mu } ; z=\frac{1}{x}, \vartheta = \frac{1}{\mu }z^2; \xi = -\frac{\varepsilon }{\mu }z \right) \nonumber\\
&=& z^{\frac{\Omega }{\mu }} \label{eq:90062}
\end{eqnarray}
\item As $\nu  =\frac{\Omega }{\mu } +2 $ where $N \in \mathbb{N}_{0}$,

An algebraic equation of degree $2N+2$ for the determination of $\omega $ is given by
\begin{equation}
0 = \left( \omega -\frac{\Omega }{\mu }\right)\left( \omega -1-\frac{\Omega }{\mu }\right) -\frac{\Omega }{2\varepsilon ^2}  \label{eq:90063a}
\end{equation}
The eigenvalue of $\omega $ is written by $\omega _1^m$ where $m = 0,1 $; $\omega _1^0 < \omega _1^1$. Its eigenfunction is given by
\begin{eqnarray} 
y(z) &=& Q_p^{(i,3)}W_{1,m}^R \left( \mu ,\varepsilon ,\nu = \frac{\Omega }{\mu } +2, \Omega , \omega = \omega _1^m ; z=\frac{1}{x}, \vartheta = \frac{1}{\mu }z^2; \xi = -\frac{\varepsilon }{\mu }z \right) \nonumber\\
&=& z^{\frac{\Omega }{\mu }} \left\{ 1+ \left( \frac{\Omega }{\mu } -\omega _1^m \right) \xi\right\}
 \label{eq:90063b}
\end{eqnarray}
\item As $\nu = \frac{\Omega }{\mu } +2N+3 $ where $N \in \mathbb{N}_{0}$,

An algebraic equation of degree $2N+3$ for the determination of $\omega $ is given by
\begin{eqnarray}
0  = \sum_{r=0}^{N+1}\bar{c}\left( r, 2(N-r)+3; 2N+2,\omega \right)\label{eq:90064a}
\end{eqnarray}
The eigenvalue of $\omega $ is written by $\omega _{2N+2}^m$ where $m = 0,1,2,\cdots,2N+2 $; $\omega _{2N+2}^0 < \omega _{2N+2}^1 < \cdots < \omega _{2N+2}^{2N+2}$. Its eigenfunction is given by
\begin{eqnarray} 
y(z) &=& Q_p^{(i,3)}W_{2N+2,m}^R \left( \mu ,\varepsilon ,\nu = \frac{\Omega }{\mu } +2N+3, \Omega , \omega = \omega _{2N+2}^m ; z=\frac{1}{x}, \vartheta = \frac{1}{\mu }z^2; \xi = -\frac{\varepsilon }{\mu }z \right) \nonumber\\
&=& \sum_{r=0}^{N+1} y_{r}^{2(N+1-r)}\left( 2N+2,\omega _{2N+2}^m;z\right) \label{eq:90064b}
\end{eqnarray}
\item As $\nu = \frac{\Omega }{\mu } +2N+4 $ where $N \in \mathbb{N}_{0}$,

An algebraic equation of degree $2N+4$ for the determination of $\omega $ is given by
\begin{eqnarray}
0  = \sum_{r=0}^{N+2}\bar{c}\left( r, 2(N+2-r); 2N+3,\omega \right)\label{eq:90065a}
\end{eqnarray}
The eigenvalue of $\omega $ is written by $\omega _{2N+3}^m$ where $m = 0,1,2,\cdots,2N+3 $; $\omega _{2N+3}^0 < \omega _{2N+3}^1 < \cdots < \omega _{2N+3}^{2N+3}$. Its eigenfunction is given by
\begin{eqnarray} 
y(z) &=& Q_p^{(i,3)}W_{2N+3,m}^R \left( \mu ,\varepsilon ,\nu = \frac{\Omega }{\mu } +2N+4, \Omega , \omega = \omega _{2N+3}^m ; z=\frac{1}{x}, \vartheta = \frac{1}{\mu }z^2; \xi = -\frac{\varepsilon }{\mu }z \right) \nonumber\\
&=& \sum_{r=0}^{N+1} y_{r}^{2(N-r)+3}\left( 2N+3,\omega _{2N+3}^m;z\right) \label{eq:90065b}
\end{eqnarray}
In the above,
\begin{eqnarray}
\bar{c}(0,n;j,\omega )  &=& \frac{\left( \frac{\Omega }{\mu }-\omega \right)_{n} }{\left( 1 \right)_{n}} \left( -\frac{\varepsilon }{\mu } \right)^{n}\label{eq:90066a}\\
\bar{c}(1,n;j,\omega ) &=& \left( \frac{1}{\mu } \right) \sum_{i_0=0}^{n}\frac{ \left( i_0 -j \right) \left( i_0 +\frac{\Omega }{\mu } \right) }{\left( i_0+2 \right)} \frac{\left( \frac{\Omega }{\mu }-\omega \right)_{i_0} }{\left( 1 \right)_{i_0}}   \frac{\left( 2+ \frac{\Omega }{\mu }-\omega \right)_{n} \left( 3 \right)_{i_0}}{\left( 2+ \frac{\Omega }{\mu }-\omega \right)_{i_0} \left( 3 \right)_{n}} \left( -\frac{\varepsilon }{\mu } \right)^{n }  
\label{eq:90066b}\\
\bar{c}(\tau ,n;j,\omega ) &=& \left( \frac{1}{\mu } \right)^{\tau } \sum_{i_0=0}^{n} \frac{ \left( i_0 -j \right) \left( i_0 +\frac{\Omega }{\mu } \right) }{\left( i_0+2 \right)} \frac{\left( \frac{\Omega }{\mu }-\omega \right)_{i_0} }{\left( 1 \right)_{i_0}}  \nonumber\\
&&\times  \prod_{k=1}^{\tau -1} \left( \sum_{i_k = i_{k-1}}^{n} \frac{\left( i_k+ 2k-j \right)\left( i_k+ 2k +\frac{\Omega }{ \mu } \right) }{\left( i_k+2k+2 \right)} \right.   \left. \frac{ \left( 2k+\frac{\Omega }{ \mu } -\omega  \right)_{i_{k}}\left( 2k+1 \right)_{i_{k-1}} }{\left( 2k+\frac{\Omega }{ \mu } -\omega \right)_{i_{k-1}} \left( \frac{k}{2}+1 \right)_{i_k}} \right) \nonumber\\ 
&&\times \frac{ \left( 2\tau +\frac{\Omega }{ \mu }-\omega  \right)_{n}\left( 2 \tau +1 \right)_{i_{\tau -1}} }{\left( 2\tau +\frac{\Omega }{ \mu }-\omega \right)_{i_{\tau -1}}  \left( 2 \tau +1 \right)_{n} } \left( -\frac{\varepsilon }{\mu } \right)^{n } \label{eq:90066c} 
\end{eqnarray}
\begin{eqnarray}
y_0^m(j,\omega ;z) &=& z^{\frac{\Omega }{\mu }} \sum_{i_0=0}^{m} \frac{\left( \frac{\Omega }{\mu }-\omega \right)_{i_0} }{\left( 1 \right)_{i_0}}  \xi ^{i_0} \label{eq:90067a}\\
y_1^m(j,\omega ;z) &=& z^{\frac{\Omega }{\mu }} \left\{\sum_{i_0=0}^{m} \frac{ \left( i_0 -j \right) \left( i_0 +\frac{\Omega }{\mu } \right) }{\left( i_0+2 \right)} \frac{\left( \frac{\Omega }{\mu }-\omega \right)_{i_0} }{\left( 1 \right)_{i_0}} \right.   \left. \sum_{i_1 = i_0}^{m} \frac{\left( 2+ \frac{\Omega }{\mu }-\omega  \right)_{i_1}  \left( 3 \right)_{i_0} }{\left( 2+ \frac{\Omega }{\mu }-\omega  \right)_{i_0}  \left( 3 \right)_{i_1}} \xi ^{i_1}\right\}\vartheta 
\hspace{1.8cm} \label{eq:90067b}\\
y_{\tau }^m(j,\omega ;z) &=& z^{\frac{\Omega }{\mu }} \left\{ \sum_{i_0=0}^{m}  \frac{ \left( i_0 -j \right) \left( i_0 +\frac{\Omega }{\mu } \right) }{\left( i_0+2 \right)} \frac{\left( \frac{\Omega }{\mu }-\omega \right)_{i_0} }{\left( 1 \right)_{i_0}} \right.\nonumber\\
&&\times \prod_{k=1}^{\tau -1} \left( \sum_{i_k = i_{k-1}}^{m} \frac{\left( i_k+ 2k-j \right)\left( i_k+ 2k +\frac{\Omega }{ \mu } \right) }{\left( i_k+2k+2 \right)} \right.   \left. \frac{ \left( 2k+\frac{\Omega }{ \mu } -\omega  \right)_{i_{k}}\left( 2k+1 \right)_{i_{k-1}} }{\left( 2k+\frac{\Omega }{ \mu } -\omega \right)_{i_{k-1}} \left( \frac{k}{2}+1 \right)_{i_k}} \right) \nonumber\\
&&\times  \left. \sum_{i_{\tau } = i_{\tau -1}}^{m} \frac{ \left( 2 \tau +\frac{\Omega }{ \mu }-\omega  \right)_{i_{\tau}}\left( 2 \tau +1 \right)_{i_{\tau -1}} }{\left( 2 \tau +\frac{\Omega }{ \mu }-\omega  \right)_{i_{\tau -1}}  \left( 2 \tau +1 \right)_{i_{\tau }} } \xi ^{i_{\tau }}\right\} \vartheta ^{\tau } \label{eq:90067c} 
\end{eqnarray}
\end{enumerate}
\end{remark}
\subsection{The second species complete polynomial of the GCH equation using R3TRF}
As previously mentioned, we need a condition for the second species complete polynomial such as $B_{j}=B_{j+1}= A_{j}=0$ where $j \in \mathbb{N}_{0}$.
\begin{theorem}
In chapter 2, the general expression of a function $y(x)$ for the second species complete polynomial using reversible  3-term recurrence formula is given by
\begin{enumerate} 
\item As $B_1=A_0=0$,
\begin{equation}
y(x) = y_{0}^{0}(x) \label{eq:90068a}
\end{equation}
\item As $B_1=B_2=A_1=0$, 
\begin{equation}
y(x)= y_{0}^{1}(x)  \label{eq:90068b}
\end{equation}
\item As $B_{2N+2}=B_{2N+3}=A_{2N+2}=0$ where $N \in \mathbb{N}_{0}$,
\begin{equation}
y(x)= \sum_{r=0}^{N+1} y_{r}^{2(N+1-r)}(x) \label{eq:90068c}
\end{equation}
\item As $B_{2N+3}=B_{2N+4}=A_{2N+3}=0$ where $N \in \mathbb{N}_{0}$,
\begin{equation}
y(x)= \sum_{r=0}^{N+1} y_{r}^{2(N-r)+3}(x) \label{eq:90068d}
\end{equation}
In the above,
\begin{eqnarray}
y_0^m(x) &=& c_0 x^{\lambda} \sum_{i_0=0}^{m} \left\{ \prod _{i_1=0}^{i_0-1}A_{i_1} \right\} x^{i_0 }\label{eq:90069a}\\
y_1^m(x) &=& c_0 x^{\lambda} \sum_{i_0=0}^{m}\left\{ B_{i_0+1} \prod _{i_1=0}^{i_0-1}A_{i_1}  \sum_{i_2=i_0}^{m} \left\{ \prod _{i_3=i_0}^{i_2-1}A_{i_3+2} \right\}\right\} x^{i_2+2 } \label{eq:90069b}\\
y_{\tau }^m(x) &=& c_0 x^{\lambda} \sum_{i_0=0}^{m} \left\{B_{i_0+1}\prod _{i_1=0}^{i_0-1} A_{i_1} 
\prod _{k=1}^{\tau -1} \left( \sum_{i_{2k}= i_{2(k-1)}}^{m} B_{i_{2k}+(2k+1)}\prod _{i_{2k+1}=i_{2(k-1)}}^{i_{2k}-1}A_{i_{2k+1}+2k}\right) \right. \nonumber\\
&& \times \left. \sum_{i_{2\tau} = i_{2(\tau -1)}}^{m} \left( \prod _{i_{2\tau +1}=i_{2(\tau -1)}}^{i_{2\tau}-1} A_{i_{2\tau +1}+ 2\tau} \right) \right\} x^{i_{2\tau}+2\tau }\hspace{1cm}\mathrm{where}\;\tau \geq 2
\label{eq:90069c}
\end{eqnarray} 
\end{enumerate}
\end{theorem}
For the second species complete polynomial of the GCH equation about $x=\infty $ in section 10.2.2, there are two possible fixed values of $\Omega $ such as $\Omega  =-\mu j$ and $\mu \left( \nu -j-1 \right)$ for $B_{j+1}=0$.  
As $\Omega =-\mu j$, we need a fixed value $\nu =0$ for $B_{j}=0$ and a fixed constant $\omega =0$ for $A_{j}=0$.
As $\Omega =\mu \left( \nu -j-1 \right)$, a fixed value $\nu =2$ for $B_{j}=0$ and a fixed constant $\omega =1$ for $A_{j}=0$ are required.
A series solution $y(z)$ divided by $z^{\frac{\Omega }{\mu } }$ of the GCH equation about $x=\infty $ as $\Omega =-\mu j$ and $\nu =\omega = 0$ for the second species complete polynomial using R3TRF is equal to an independent solution $y(z)$ divided by $z^{\frac{\Omega }{\mu } }$ of the CHE as $\Omega = \mu \left( 1-j \right)$, $\nu =2$ and $\omega =1 $.
  
In (\ref{eq:90068a})--(\ref{eq:90069c}) replace $c_0$, $\lambda $ and $x$ by $1$, $\Omega /\mu $ and $z$. 
Substitute (\ref{eq:90041a}) and (\ref{eq:90041b}) into (\ref{eq:90069a})--(\ref{eq:90069c}).

\underline{(1) The case of $\Omega  =-\mu j$ and $\nu =\omega = 0$,} 

As $B_1=A_0=0$, substitute the new (\ref{eq:90069a}) with $\Omega =0$ into (\ref{eq:90068a}) putting $j=0$. 
As $B_1=B_2=A_1=0$, substitute the new (\ref{eq:90069a}) with $\Omega  =-\mu$ into (\ref{eq:90068b}) putting $j=1$. 
As $B_{2N+2}=B_{2N+3}=A_{2N+2}=0$, substitute the new (\ref{eq:90069a})--(\ref{eq:90069c}) with $\Omega  =-\mu \left( 2N+2 \right)$ into (\ref{eq:90068c}) putting $j=2N+2$.
As $B_{2N+3}=B_{2N+4}=A_{2N+3}=0$, substitute the new (\ref{eq:90069a})--(\ref{eq:90069c}) with $\Omega  =-\mu \left( 2N+3 \right) $ into (\ref{eq:90068d}) putting $j=2N+3$.

\underline{(2) The case of $\Omega = \mu \left( 1-j \right) $, $\nu =2$ and $\omega = 1 $,} 
 
As $B_1=A_0=0$, substitute the new (\ref{eq:90069a}) with $\Omega =\mu $ into (\ref{eq:90068a}) putting $j=0$. 
As $B_1=B_2=A_1=0$, substitute the new (\ref{eq:90069a}) with $\Omega  =0$ into (\ref{eq:90068b}) putting $j=1$. 
As $B_{2N+2}=B_{2N+3}=A_{2N+2}=0$, substitute the new (\ref{eq:90069a})--(\ref{eq:90069c}) with $\Omega  =-\mu \left( 2N+1 \right)$ into (\ref{eq:90068c}) putting $j=2N+2$.
As $B_{2N+3}=B_{2N+4}=A_{2N+3}=0$, substitute the new (\ref{eq:90069a})--(\ref{eq:90069c}) with $\Omega  =-\mu \left( 2N+2 \right) $ into (\ref{eq:90068d}) putting $j=2N+3$.

After the replacement process, we obtain the independent solution of the GCH equation about $x=\infty $. The solution is as follows.
\begin{remark}
The power series expansion of the GCH equation of the first kind for the second species complete polynomial using R3TRF about $x=\infty $ is given by
\begin{enumerate} 
\item As $\Omega  =\nu =\omega =0$, 

Its eigenfunction is given by
\begin{eqnarray}
y(z) &=& Q_p^{(i,1)}W_0^R \left( \mu ,\varepsilon ,\nu =0,\Omega =0, \omega = 0; z=\frac{1}{x}, \vartheta = \frac{1}{\mu }z^2; \xi = -\frac{\varepsilon }{\mu }z \right) \nonumber\\
&=& 1 \label{eq:90070a}
\end{eqnarray}
\item As $\Omega =-\mu $, $\nu =\omega =0$, 

Its eigenfunction is given by
\begin{eqnarray}
y(z) &=& Q_p^{(i,1)}W_1^R \left( \mu ,\varepsilon ,\nu =0,\Omega =-\mu, \omega = 0; z=\frac{1}{x}, \vartheta = \frac{1}{\mu }z^2; \xi = -\frac{\varepsilon }{\mu }z \right) \nonumber\\
 &=& z^{-1} \left\{ 1-\xi \right\} \label{eq:90070b}
\end{eqnarray}
\item As $\Omega =-\mu \left( 2N+2 \right)$, $\nu =\omega =0$ where $N \in \mathbb{N}_{0}$, 

Its eigenfunction is given by
\begin{eqnarray}
y(z) &=& Q_p^{(i,1)}W_{2N+2}^R \left( \mu ,\varepsilon ,\nu =0,\Omega =-\mu \left( 2N+2 \right), \omega = 0; z=\frac{1}{x}, \vartheta = \frac{1}{\mu }z^2; \xi = -\frac{\varepsilon }{\mu }z \right) \nonumber\\
 &=& z^{-2N-2} \sum_{r=0}^{N+1} y_{r}^{2(N+1-r)}\left( 2N+2;z\right) \label{eq:90070c}
\end{eqnarray}

\item As $\Omega =-\mu \left( 2N+3 \right)$, $\nu =\omega =0$ where $N \in \mathbb{N}_{0}$, 

Its eigenfunction is given by
\begin{eqnarray}
y(z) &=& Q_p^{(i,1)}W_{2N+3}^R \left( \mu ,\varepsilon ,\nu =0,\Omega =-\mu \left( 2N+3 \right), \omega = 0; z=\frac{1}{x}, \vartheta = \frac{1}{\mu }z^2; \xi = -\frac{\varepsilon }{\mu }z \right) \nonumber\\
 &=& z^{-2N-3} \sum_{r=0}^{N+1} y_{r}^{2(N-r)+3}\left( 2N+3;z\right) \label{eq:90070d}
\end{eqnarray}
\item As $\Omega =\mu $, $\nu =2$, $\omega =1$, 

Its eigenfunction is given by
\begin{eqnarray}
y(z) &=& Q_p^{(i,2)}W_0^R \left( \mu ,\varepsilon ,\nu =2,\Omega =\mu , \omega = 1; z=\frac{1}{x}, \vartheta = \frac{1}{\mu }z^2; \xi = -\frac{\varepsilon }{\mu }z \right) \nonumber\\
&=& z \label{eq:90070e}
\end{eqnarray}
\item As $\Omega =0$, $\nu =2$, $\omega =1$, 

Its eigenfunction is given by
\begin{eqnarray}
y(z) &=& Q_p^{(i,2)}W_1^R \left( \mu ,\varepsilon ,\nu =2,\Omega =0, \omega = 1; z=\frac{1}{x}, \vartheta = \frac{1}{\mu }z^2; \xi = -\frac{\varepsilon }{\mu }z \right) \nonumber\\
 &=& 1-\xi  \label{eq:90070f}
\end{eqnarray}
\item As $\Omega =-\mu \left( 2N+1\right)$, $\nu =2$, $\omega =1$ where $N \in \mathbb{N}_{0}$, 

Its eigenfunction is given by
\begin{eqnarray}
y(z) &=& Q_p^{(i,2)}W_{2N+2}^R \left( \mu ,\varepsilon ,\nu =2,\Omega =-\mu \left( 2N+1 \right), \omega = 1; z=\frac{1}{x}, \vartheta = \frac{1}{\mu }z^2; \xi = -\frac{\varepsilon }{\mu }z \right) \nonumber\\
 &=& z^{-2N-1} \sum_{r=0}^{N+1} y_{r}^{2(N+1-r)}\left( 2N+2;z\right) \label{eq:90070g}
\end{eqnarray}

\item As $\Omega =-\mu \left( 2N+2 \right)$, $\nu =2$, $\omega =1$ where $N \in \mathbb{N}_{0}$, 

Its eigenfunction is given by
\begin{eqnarray}
y(z) &=& Q_p^{(i,2)}W_{2N+3}^R \left( \mu ,\varepsilon ,\nu =2,\Omega =-\mu \left( 2N+2 \right), \omega = 1; z=\frac{1}{x}, \vartheta = \frac{1}{\mu }z^2; \xi = -\frac{\varepsilon }{\mu }z \right) \nonumber\\
 &=& z^{-2N-2} \sum_{r=0}^{N+1} y_{r}^{2(N-r)+3}\left( 2N+3;z\right) \label{eq:90070h}
\end{eqnarray}
In the above,
\begin{eqnarray}
y_0^m(j;z) &=& \sum_{i_0=0}^{m} \frac{\left( -j\right)_{i_0} }{\left( 1 \right)_{i_0}} \xi ^{i_0} \label{eq:90071a}\\
y_1^m(j;z) &=& \left\{\sum_{i_0=0}^{m} \frac{ \left( i_0 -j\right)\left( i_0+1-j \right) }{\left( i_0+2 \right)} \frac{\left( -j\right)_{i_0} }{\left( 1 \right)_{i_0}} \right.  \left. \sum_{i_1 = i_0}^{m} \frac{\left( 2-j\right)_{i_1}  \left( 3 \right)_{i_0} }{\left( 2-j \right)_{i_0} \left( 3 \right)_{i_1}} \xi ^{i_1}\right\}\vartheta 
 \label{eq:90071b}\\
y_{\tau }^m(j;z) &=& \left\{ \sum_{i_0=0}^{m} \frac{ \left( i_0 -j\right)\left( i_0+1-j \right) }{\left( i_0+2 \right)} \frac{\left( -j\right)_{i_0} }{\left( 1 \right)_{i_0}} \right.\nonumber\\
&&\times \prod_{k=1}^{\tau -1} \left( \sum_{i_k = i_{k-1}}^{m} \frac{\left( i_k+ 2k -j\right) \left( i_k+ 2k +1-j\right)}{\left( i_k+2k+2\right)} \right. \left. \frac{\left( 2k-j\right)_{i_k} \left( 2k+1 \right)_{i_{k-1}} }{\left( 2k-j\right)_{i_{k-1}}  \left( 2k+1 \right)_{i_k}} \right) \nonumber\\
&&\times  \left. \sum_{i_{\tau } = i_{\tau -1}}^{m} \frac{\left( 2\tau -j\right)_{i_{\tau}}  \left( 2 \tau +1 \right)_{i_{\tau -1}} }{\left( 2 \tau -j\right)_{i_{\tau -1}} \left( 2 \tau +1 \right)_{i_{\tau }} } \xi ^{i_{\tau }}\right\} \vartheta ^{\tau } \label{eq:90071c} 
\end{eqnarray}
\end{enumerate}
\end{remark}  
\section{Summary}

In chapter 6 of Ref.\cite{9Choun2013} and this chapter, all possible general solutions in series of the GCH equation about the irregular singular point at infinity are a polynomial of type 1, the first species complete polynomial and the second species complete polynomial. 
For a polynomial of type 1,  I treat $\nu $ as a fixed value and $\varepsilon $, $\mu $, $\omega $, $\Omega $ as free variables.
For the first species complete polynomial, I treat $\omega $, $\Omega $ as fixed values and $\varepsilon $, $\mu $, $\nu $ as free variables.
For an another combinatorial case of the first species complete polynomial, I treat $\nu $, $\omega $ as fixed values and $\varepsilon $, $\mu $, $\Omega $ as free variables.   
For the second species complete polynomial, I treat $\nu $, $\omega $, $\Omega $ as fixed values and $\varepsilon $, $\mu $ as free variables.
There is no such analytic solutions in series for an infinite series and a polynomial of type 2  because of $\lim_{n \ge 1} B_n \rightarrow \infty $ in (\ref{eq:9009b}). 

By applying (\ref{eq:9009a}) and (\ref{eq:9009b}) with $\nu = 2 \nu _i+i +1+ \Omega /\mu $ where $i,\nu \in \mathbb{N}_{0}$, into 3TRF such as the general summation formulas in closed forms, a polynomial of type 1 of the GCH equation around $x=\infty $ is constructed in chapter 6 of Ref.\cite{9Choun2013}. As we observe each of sub-power series solutions $y_l(x)$, obtained by observing the term of sequence $c_n$ which includes $l$ terms of $A_n's$, in a general solution in series $y(x)= \sum_{n=0}^{\infty } y_n(x)$, the denominators and numerators in all $B_n$ terms of $y_l(x)$ arise with Pochhammer symbols. 
An integral solution of the GCH equation (each sub-integral $y_l(x)$ is composed of $2l$ terms of definite integrals and $l$ terms of contour integrals) for a polynomial of type 1 is derived by applying the contour integral of Tricomi's function into the sub-power series of a general series solution. 
The generating function for a polynomial of type 1 is constructed analytically by applying generating functions for Tricomi's polynomials into sub-integrals of a general integral solution.
 
By putting $A_n$ and $B_n$ terms, composed of various coefficients of a numerator and a denominator
in a recursive relation of the GCH equation around $x=\infty $, into a general summation formula of the first species complete polynomial using 3TRF, I show its power series solutions for the first species complete polynomial in this chapter. Similarly, another series solutions of the GCH equation for a polynomial of type 3 are constructed by applying a general summation expression of the first species complete polynomials using R3TRF. 

By putting (\ref{eq:9009a}) and (\ref{eq:9009b}) with $\Omega = -\mu j$, $\nu =\omega = 0$ or $\Omega = \mu \left( 1-j \right) $, $\nu =2$, $\omega =1 $  where $j \in \mathbb{N}_{0}$ into a general summation formula of the second species complete polynomial using 3TRF, formal series solutions of the GCH equation around $x=\infty $ for the second species complete polynomial are obtained. By using similar summation techniques, another formal solutions of the GCH equation for a polynomial of type 3 is constructed by applying a general summation expression of the second species complete polynomials using R3TRF.

Even though numerical computations of the GCH equation for complete polynomials (either the first species complete polynomials or the second species complete polynomials) using 3TRF and R3TRF are equivalent to each other, classical summation structures between these two polynomials of a type 3 have several crucial differences such as:
(1) $A_n$ term in each of finite sub-polynomials $y_{\tau }^m(x)$, obtained by observing the term of sequence $c_n$ which includes $\tau $ terms of $A_n's$, of general series solutions is the leading term for complete polynomials using 3TRF. For complete polynomials using R3TRF, $B_n$ term in each $y_{\tau }^m(x)$ of general series solutions is the leading term.
(2) Summation solutions for complete polynomials using 3TRF contain the sum of two partial sums of the sequence $\{y_{\tau }^m(x)\}$ for its general series solutions. 
And, for complete polynomials using R3TRF, there is only one partial sum of the sequence $\{y_{\tau }^m(x)\}$ for its general Frobenius solutions. 
(3) The denominators and numerators in all $B_n$ terms of each $y_{\tau }^m(x)$ arise with Pochhammer symbols in series solutions of the GCH equation for complete polynomials using 3TRF. And for complete polynomials using R3TRF, the denominators and numerators in all $A_n$ terms of each $y_{\tau }^m(x)$ arise with Pochhammer symbols.   
(4) The first species complete polynomials using 3TRF and R3TRF have multi-valued roots of a parameter $\omega $, but the second species complete polynomials using 3TRF and R3TRF have only two fixed values of a parameter $\omega $ such as $ \omega = 0$ or 1.   

In the future papers, I will show how to construct solutions as combined definite \& contour integrals of the GCH equation around $x=\infty $ for the first and second species complete polynomials using 3TRF and R3TRF by applying integrals of hypergeometric-type into finite sub-power series of general series solutions.
Moreover, I will also represent how to build generating functions for these polynomials of type 3 by applying generating functions for hypergeometric-type polynomials into sub-integrals of general integrals.
 
\addcontentsline{toc}{section}{Bibliography}
\bibliographystyle{model1a-num-names}
\bibliography{<your-bib-database>}
\chapter{Complete polynomials of Lam\'{e} equation}
\chaptermark{Complete polynomials of Lam\'{e} equation} 
 
Power series solutions of Lam\'{e} equation either using an algebraic form or Weierstrass's form provide a 3-term recurrence relation between successive coefficients. 
By applying three term recurrence formula (3TRF) \cite{10Chou2012b}, I construct power series solutions in closed forms of Lam\'{e} equation in the algebraic form for an infinite series and a polynomial of type 1. And its combined definite and contour integrals involving only $_2F_1$ functions are obtained including generating functions for Lam\'{e} polynomials of type 1 \cite{10Chou2012f,10Chou2012h}.  

In chapter 8 of Ref.\cite{10Choun2013}, I apply reversible three term recurrence formula (R3TRF) to power series expansions in closed forms of Lam\'{e} equation in the algebraic form for an infinite series and a polynomial of type 2. And its representations for solutions as combined definite and contour integrals are constructed analytically including generating functions for Lam\'{e} polynomials of type 2.

In this chapter I construct Frobenius solutions of Lam\'{e} equation for a polynomial of type 3 by applying general summation formulas of complete polynomials using 3TRF and R3TRF.
\section{Introduction}
Actually, the nature is nonlinear systems and non-symmetry geometrically. Many great scholars linearize those system with better numerical computations for the purpose of simplification. The sphere is a perfect symmetrical object geometrically, any points at its surface are always equidistant from its center. In contrast, an ellipsoid is a non-symmetrical and imperfect one, a surface whose plane sections are all ellipses or circles. 
For example, even though Romain and Jean-Pierre showed that Laplace's ellipsoidal harmonics were represented in calculations of Earth's gravitational potential\cite{10Romai2001}, Laplace's spherical harmonics are preferred rather than ellipsoid harmonics because of its complex mathematical calculations: A recurrence relation of Laplace's equation in spherical coordinates consists of a 2-term between successive coefficients; but a recursion relation for an ellipsoidal geometry is composed of a 3-term.

In 1837, Gabriel Lam\'{e} introduced a second-order ordinary differential equation having four regular singular points in the method of separation of variables applied to Laplace's equation in certain systems of curvilinear coordinates such as `elliptic coordinates' \cite{10Lame1837}. This equation has been named as `Lam\'{e} equation' or `ellipsoidal harmonic equation' by various scholars\cite{10Erde1955}. 
Separation of variables in the three-dimensional Laplace equation including Laplace's ellipsoid harmonics is studied by many authors 
\cite{10Levi1949,10Boch1894,10Moon1952a,10Moon1952b,10Moon1953}.
In 1934, Eisenhart determined various canonical forms for euclidean 3-space and higher order in more general differential equations. For euclidean 4-space, he showed the St\"{a}ckel forms for the constant Riemannian metric of a space $V_3$ \cite{10Eise1934}.
Lam\'{e} equation is appeared in modern physics such as boundary value problems in ellipsoidal geometry, chaotic Hamiltonian systems, the theory of Bose-Einstein condensates, scalar Dirichlet scattering problems of a plane wave by a general ellipsoid, group theory to quantum mechanical bound states and band structures of the Lam\'{e} Hamiltonian, etc \cite{10Brac2001,10Bron2001,10Qian2003,10Slee1967,10Gurs1983,10Iach2000}.    

Generally, two different forms of Lam\'{e} equation are given such as the algebraic form and Weierstrass's form. The Jacobian (Weierstrass's) form is preferred rather than the algebraic form for mathematical calculations. Because the algebraic form has three singularity parameters such as $a$, $b$ and $c$. And the Jacobian form only has two singularity parameters which are 1 and the modulus of the elliptic function $sn\; z$.
For the general mathematical properties of analytic solutions of Lame equation, the algebraic form is more convenient. Its asymptotic expansions and boundary conditions of it are more various because of three different singularity parameters. 

By substituting a power series with unknown coefficients into Lam\'{e} equation, its recursive relation is composed of a 3-term between consecutive coefficients in a power series solution \cite{10Hobs1931,10Whit1927}. In contrast, the recurrence relation of a hypergeoemtric equation has a 2-term.
For functions of hypergeometric type such as Gauss hypergeometric, Kummer functions, and etc, their assorted definite or contour integrals are already known analytically. For functions of Lam\'{e} type represented either in the algebraic form or Weierstrass's form, no such solutions for integral representations have been constructed \cite{10Erde1955,10Hobs1931,10Whit1927}. 
Instead, solutions can be represented to satisfy relatively a Fredholm integral equation of the first kind with a degenerate kernel, and such a Lam\'{e} function describe integral relationships in terms of an another analytic solution \cite{10Ince1922,10Ince1940a,10Ince1940b,10Lamb1934,10Erde1942,10Whit1915,10Whit1927}.
Several mathematicians developed integral equations of Lam\'{e} equation by using simple kernels involving Legendre functions of the Jacobian elliptic function or ellipsoidal harmonic functions, and various forms of its integral equation were found \cite{10Erde1955,10Arsc1964a,10Shai1980,10Slee1968a,10Volk1982,10Volk1983,10Volk1984,10Wang1989,10Whit1927}. 
Many great scholars just left analytic solutions of Lam\'{e} equation as solutions of recurrences because of a 3-term recurrence relation between successive coefficients in its power series expansion. 3 or more terms in a recursive relation of any linear ODEs create complicated mathematical computations.

According to Gabriel Lam\'{e} \cite{10Lame1837,10Stie1885,10Zwil1997}, the algebraic form of Lam\'{e} equation is taken as 
\begin{equation}
\frac{d^2{y}}{d{x}^2} + \frac{1}{2}\left(\frac{1}{x-a} +\frac{1}{x-b} + \frac{1}{x-c}\right) \frac{d{y}}{d{x}} +  \frac{-\alpha (\alpha +1) x+q}{4 (x-a)(x-b)(x-c)} y = 0\label{eq:10001}
\end{equation}
Lam\'{e} equation has four regular singular points: $a$, $b$, $c$ and $\infty $; the exponents at the first three are all $0$ and $\frac{1}{2}$, and those at infinity are $-\frac{1}{2}\alpha $ and $\frac{1}{2}(\alpha +1)$ \cite{10Moon1961,10Boch1894}. 
The Heun differential equation generalizes Lam\'{e} equation represented either in the algebraic form or Weierstrass's form.

According to Karl Heun \cite{10Heun1889,10Ronv1995}, the canonical form of general Heun's equation is taken as
\begin{equation}
\frac{d^2{y}}{d{x}^2} + \left(\frac{\gamma }{x} +\frac{\delta }{x-1} + \frac{\epsilon }{x-a}\right) \frac{d{y}}{d{x}} +  \frac{\alpha \beta x-q}{x(x-1)(x-a)} y = 0 \label{eq:10002}
\end{equation}
where the parameters should satisfy the condition such as $\epsilon = \alpha +\beta -\gamma -\delta +1$ in order to have the regular singular point at $x =\infty $.
The parameters play different roles: $a \ne 0 $ is the singularity parameter, $\alpha $, $\beta $, $\gamma $, $\delta $, $\epsilon $ are exponent parameters, $q$ is the accessory parameter which in many physical applications appears as a spectral parameter. Also, $\alpha $ and $\beta $ are identical to each other. The total number of free parameters is six. It has four regular singular points which are 0, 1, $a$ and $\infty $ with exponents $\{ 0, 1-\gamma \}$, $\{ 0, 1-\delta \}$, $\{ 0, 1-\epsilon \}$ and $\{ \alpha, \beta \}$.

Let $z=x-a$ in (\ref{eq:10001}).
\begin{equation}
\frac{d^2{y}}{d{z}^2} + \frac{1}{2}\left(\frac{1}{z} +\frac{1}{z-(b-a)} + \frac{1}{z-(c-a)}\right) \frac{d{y}}{d{z}} +  \frac{-\frac{1}{4}\alpha (\alpha +1)z- \frac{1}{4}\left(- q+ \alpha (\alpha +1)a\right)}{ z(z-(b-a))(z-(c-a))} y = 0\label{eq:10003}
\end{equation}
As we compare (\ref{eq:10002}) with (\ref{eq:10003}), all coefficients on the above are correspondent to the following way.
\begin{equation}
\begin{split}
& \gamma ,\delta ,\epsilon  \longleftrightarrow   \frac{1}{2} \\ & 1\longleftrightarrow  b-a \\ & a\longleftrightarrow  c-a \\ & \alpha  \longleftrightarrow \frac{1}{2}(\alpha +1) \\
& \beta   \longleftrightarrow -\frac{1}{2}\alpha \\
& q \longleftrightarrow  \frac{1}{4}\left( - q+ \alpha (\alpha +1)a\right) \\ & x \longleftrightarrow z
\end{split}\label{eq:10004}   
\end{equation}
Assume that the solution of (\ref{eq:10003}) is
\begin{equation}
y(z)= \sum_{n=0}^{\infty } c_n z^{n+\lambda }  \label{eq:10005}
\end{equation}
where $\lambda $ is an indicial root. Plug (\ref{eq:10005}) into (\ref{eq:10003}).
\begin{equation}
c_{n+1}=A_n \;c_n +B_n \;c_{n-1} \hspace{1cm};n\geq 1\label{eq:10006}
\end{equation}
where,
\begin{subequations}
\begin{eqnarray}
A_n = -\frac{(2a-b-c)}{(a-b)(a-c)} \frac{(n +\lambda -\varphi)(n +\lambda + \varphi)}{(n+1+\lambda )(n+\frac{1}{2}+\lambda )} \label{eq:10007a}
\end{eqnarray}
\begin{equation}
B_n =\frac{-1}{(a-b)(a-c)}\frac{\left( n+\frac{\alpha }{2}-\frac{1}{2}+\lambda \right) \left( n-\frac{\alpha }{2}-1+\lambda\right)}{(n+1+\lambda )(n+\frac{1}{2}+\lambda )}\label{eq:10007b}
\end{equation}
and
\begin{equation}
\varphi = \sqrt{\frac{\alpha (\alpha +1)a-q}{4(2a-b-c)}} \nonumber
\end{equation}
\end{subequations}
where $c_1= A_0 \;c_0$. Two indicial roots are given such as $\lambda = 0$ and $ \frac{1}{2}$.
 
There are 4 different possible power series solutions of a linear ODE having a 3-term recursive relation between successive coefficients such as an infinite series and 3 types of polynomials: (1) a polynomial which makes $B_n$ term terminated; $A_n$ term is not terminated, (2) a polynomial which makes $A_n$ term terminated; $B_n$ term is not terminated, (3) a polynomial which makes $A_n$ and $B_n$ terms terminated at the same time, referred as `a complete polynomial.' 

There are two different types of complete polynomials such as the first species complete polynomial and the second species complete polynomial. 
The first species complete polynomial is applicable since a parameter of a numerator in $B_n$ term and a (spectral) parameter of a numerator in $A_n$ term are fixed constants. The second species complete polynomial is employed since two parameters of a numerator in $B_n$ term and a parameter of a numerator in $A_n$ term are fixed constants.
The former has multi-valued roots of a parameter of a numerator in $A_n$ term, but the latter has only one fixed value of a parameter of a numerator in $A_n$ term for an eigenvalue. 

\begin{table}[h]
\begin{center}
\thispagestyle{plain}
\hspace*{-0.1\linewidth}\resizebox{1.3\linewidth}{!}
{
 \Tree[.{\Huge Lam\'{e} differential equation in the algebraic form or Weierstrass's form} [.{\Huge 3TRF} [.{\Huge Infinite series} ]
              [.{\Huge Polynomials} [[.{\Huge Polynomial of type 1} ]
               [.{\Huge Polynomial of type 3} [.{\Huge $ \begin{array}{lcll}  1^{\mbox{st}}\;  \mbox{species}\\ \mbox{complete} \\ \mbox{polynomial} \end{array}$} ]  ]]]]                         
  [.{\Huge R3TRF} [.{\Huge Infinite series} ]
     [.{\Huge Polynomials} [[.{\Huge Polynomial of type 2} ]
       [.{\Huge  Polynomial of type 3} [.{\Huge $ \begin{array}{lcll}  1^{\mbox{st}}\;  \mbox{species} \\ \mbox{complete} \\ \mbox{polynomial} \end{array}$} ]  ]]]]]
}
\end{center}
\caption{Power series of Lam\'{e} equation}
\end{table}  
Table 11.1 informs us about all possible general formal series solutions of Lam\'{e} equation; it is available at one of any regular singular points. According to my definition, 2 types of the general summation formulas, named as 3TRF and R3TRF, are applicable in order to obtain power series solutions of any ODEs having a 3-term recurrence relation between successive coefficients.

Traditionally, a power series considering a Maclaurin series in one variable is defined as an infinite series of the form $y(x) =\sum_{n=0}^{\infty } c_n x^n$; $c_n$ represents the coefficient of the $n^{th}$ term. We can convert any analytic functions as a formal solution in series. For instance, by applying a power series into hypergeometric type differential equations, we obtain their analytic solutions directly including definite or contour integrals. However, if we substitute a power series into any ODEs having a 3-term in a recurrence relation, such as Heun, Lam\'{e}, Mathieu and Confluent Heun equations, there are no ways to find series solutions in which coefficients are given explicitly \cite{10Arsc1983}. 

In general, for a 3-term case, since we substitute $\sum_{n=0}^{\infty } c_n x^n$ into one of any linear ODEs, a recursive relation for its Frobenius solution is given by (\ref{eq:10006}) where $c_1= A_0 \;c_0$. In order to describe general solutions in series for the 3-term case, I suggest that a power series $y(x)$ is an infinite sum of sub-power series.
More precisely, for $n=0,1,2,3,\cdots $, (\ref{eq:10006}) is expanded to combinations of $A_n$ and $B_n$ terms. 
A sub-power series $y_l(x)$ where $l\in \mathbb{N}_{0}$ is designated by observing the term of sequence $c_n$ which includes $l$ terms of $A_n's$.  
The formal solution in series is represented by sums of each $y_l(x)$ such as $y(x)= \sum_{n=0}^{\infty } y_n(x)$.
By allowing $A_n$ in the sequence $c_n$ is the leading term of each sub-power series $y_l(x)$, the general summation formulas of a 3-term recurrence relation in a linear ODE are constructed for an infinite series and a polynomial of type 1, referred as `three term recurrence formula (3TRF).' \cite{10Chou2012b} 

Maclaurin series solutions of Lam\'{e} equation in the algebraic form around $x=a$ are derived analytically for an infinite series and a polynomial of type 1 by substitute (\ref{eq:10007a}) and (\ref{eq:10007b}) into 3TRF \cite{10Chou2012f}. Its asymptotic series in the compact form for an infinite series is constructed with the radius of convergence $z=x-a$.
For a polynomial of type 1, I treat an exponent parameter $q$ as a free variable and $\alpha $ as a fixed value.
Solutions as general integrals of Lam\'{e} equation are built by applying the contour integral of a $_2F_1$ function into each of  sub-power series $y_l(z)$. 
And as observing its general forms, a $_2F_1$ function recurs in each of $y_m(x)$ functions; $y_m(z)$ means the sub-integral form contains $m$ term of $A_n's$. 
Each sub-integral $y_m(z)$ where $m=0,1,2,\cdots$ is composed of $2m$ terms of definite integrals and $m$ terms of contour integrals.   
Formal series solutions and integrals of Weierstrass's form of Lam\'{e} equation about $ \xi =0$ using 3TRF, for an infinite series and a polynomial of type 1, are obtained by changing all coefficients in analytic solutions of its equation in the algebraic form  \cite{10Chou2012g}.
Generating functions of Lam\'{e} equation in the Jacobian form for a polynomial of type 1 are found by applying the generating function for a Jacobi polynomial using a $_2F_1$ function into each of sub-integrals $y_m(\xi )$ in general integrals of Lam\'{e} polynomials \cite{10Chou2012h}.

In chapter 1 of Ref.\cite{10Choun2013}, a sub-power series $y_l(x)$ is found by observing the term of sequence $c_n$ which includes $l$ terms of $B_n's$ instead of $A_n's$. 
The general series solution is also described by sums of each $y_l(x)$, denoted as $y(x)= \sum_{n=0}^{\infty } y_n(x)$.  
By allowing $B_n$ in the sequence $c_n$ is the leading term of each sub-power series $y_l(x)$, general summation formulas of the 3-term recurrence relation in a linear ODE are constructed for an infinite series and a polynomial of type 2, designated as `reversible three term recurrence formula (R3TRF).' 

In chapter 8 of Ref.\cite{10Choun2013}, power series solutions of Lam\'{e} equation in the algebraic form around $z=0$ are derived  for an infinite series and a polynomial of type 2 by substitute (\ref{eq:10007a}) and (\ref{eq:10007b}) into R3TRF.  
For a polynomial of type 2, I treat an exponent parameter $\alpha $ as a free variable and $q$ as a fixed value.
General integrals of Lam\'{e} equation are obtained by applying the contour integral of a $_2F_1$ function into each of sub-power series $y_l(z)$; one interesting observation resulting from the calculations is the fact that a $_2F_1$ function recurs in each of sub-integrals $y_m(z)$ where $m$ term of $B_n's$ are contained.
And generating functions of Lam\'{e} polynomial of type 2 are found by applying the generating function for a Jacobi polynomial using a $_2F_1$ function into each of sub-integrals $y_m(z)$ in general integrals of Lam\'{e} polynomials.
(1) Power series solutions, (2) integrals and (3) generating functions of Lam\'{e} equation in Weierstrass's form using R3TRF, for an infinite series and a polynomial of type 2, are obtained by changing all coefficients in analytic solutions of its equation in the algebraic form in chapter 9 of Ref.\cite{10Choun2013}. 
 
For numerical computation, infinite series of Lam\'{e} equation around $x=a$ by applying 3TRF are equivalent to infinite series by applying R3TRF for solutions in series and integrals. 
There are mathematical difference in summation structures between two infinite series solutions of Lam\'{e} equation:
$A_n$ in sequences $c_n$ is the leading term in each of sub-power series and sub-integrals for infinite series by applying 3TRF; and $B_n$ is the leading term in each of sub-formal series and sub-integrals for infinite series by applying R3TRF.

In chapter 1, I suggest that a general solution in series $y(x)$ is a finite sum of finite sub-power series $y_{\tau }^m(x)$ where $\tau ,m \in \mathbb{N}_{0}$. $y_{\tau }^m(x)$, where the lower bound of summation is zero and the upper one is $m$, is designated by observing the term of sequence $c_n$ which includes $\tau $ terms of $A_n's$.  
By allowing $A_n$ as the leading term in each finite sub-power series $y_{\tau }^m(x)$ of a general solution in series $y(x)$, 
general summation formulas for the first and second species complete polynomials are constructed analytically. I refer these mathematical summation expressions as ``complete polynomials using 3-term recurrence formula (3TRF).'' 

In chapter 2, I define that a general solution in series $y(x)$ is a finite sum of finite sub-power series $y_{\tau }^m(x)$, and a finite sub-power series $y_{\tau }^m(x)$ is obtained by observing the term of sequence $c_n$ which includes $\tau $ terms of $B_n's$. 
By allowing $B_n$ as the leading term in $y_{\tau }^m(x)$ of the general power series $y(x)$, I construct the classical mathematical solutions in series, built in compact forms, for the first and second species complete polynomials. I denominate these classical summation formulas as ``complete polynomials using reversible 3-term recurrence formula (R3TRF).'' 

In this chapter, by substituting (\ref{eq:10007a}) and (\ref{eq:10007b}) into complete polynomials using 3TRF and R3TRF, 

(1) I construct power series solutions of Lam\'{e} equation in the algebraic form around $x=a$ for the first species complete polynomials. 

(2) I derive an algebraic equation of Lam\'{e} equation for the determination of a parameter $q$ in the combinational form of partial sums of the sequences $\{A_n\}$ and $\{B_n\}$ using 3TRF and R3TRF. 
 
(3) By changing all coefficients of formal series solutions, represented in closed forms, of Lam\'{e} equation in the algebraic form, I obtain its Frobenius solutions in Weierstrass's form for the first species complete polynomials including its polynomial equation for the determination of a parameter $h$.

(4) Four examples of 192 local solutions of the Heun equation (Maier, 2007 \cite{10Maie2007}) are provided in the appendix.  For each example, I show how to convert a local solution of Heun equation for complete polynomials using 3TRF and R3TRF to an analytic solution of Lam\'{e} equation in Weierstrass's form. 
\section{Lam\'{e} equation in the algebraic form}
  
For the first species complete polynomial of Lam\'{e} equation in the algebraic form around $x=a$ in table 11.1, I treat parameters $\alpha $ and $q$ as fixed values. There is no such solution in series for the second species complete polynomial of Lam\'{e} equation. Because a parameter $\alpha $ of a numerator in $B_n$ term in (\ref{eq:10007b}) is only a fixed constant in order to make $B_n$ term terminated at a specific index summation $n$ for Frobenius solutions for a polynomial of type 3.
\subsection{The first species complete polynomial of Lam\'{e} equation using 3TRF}

For the first species complete polynomials using 3TRF and R3TRF, we need a condition which is given by
\begin{equation}
B_{j+1}= c_{j+1}=0\hspace{1cm}\mathrm{where}\;j\in \mathbb{N}_{0}  
 \label{eq:10008}
\end{equation}
(\ref{eq:10008}) gives successively $c_{j+2}=c_{j+3}=c_{j+4}=\cdots=0$. And $c_{j+1}=0$ is defined by a polynomial equation of degree $j+1$ for the determination of an accessory parameter in $A_n$ term. 
\begin{theorem}
In chapter 1, the general summation expression of a function $y(x)$ for the first species complete polynomial using 3-term recurrence formula and its algebraic equation for the determination of an accessory parameter in $A_n$ term are given by
\begin{enumerate} 
\item As $B_1=0$,
\begin{equation}
0 =\bar{c}(1,0) \label{eq:10009a}
\end{equation}
\begin{equation}
y(x) = y_{0}^{0}(x) \label{eq:10009b}
\end{equation}
\item As $B_{2N+2}=0$ where $N \in \mathbb{N}_{0}$,
\begin{equation}
0  = \sum_{r=0}^{N+1}\bar{c}\left( 2r, N+1-r\right) \label{eq:100010a}
\end{equation}
\begin{equation}
y(x)= \sum_{r=0}^{N} y_{2r}^{N-r}(x)+ \sum_{r=0}^{N} y_{2r+1}^{N-r}(x)  \label{eq:100010b}
\end{equation}
\item As $B_{2N+3}=0$ where $N \in \mathbb{N}_{0}$,
\begin{equation}
0  = \sum_{r=0}^{N+1}\bar{c}\left( 2r+1, N+1-r\right) \label{eq:100011a}
\end{equation}
\begin{equation}
y(x)= \sum_{r=0}^{N+1} y_{2r}^{N+1-r}(x)+ \sum_{r=0}^{N} y_{2r+1}^{N-r}(x)  \label{eq:100011b}
\end{equation}
In the above,
\begin{eqnarray}
\bar{c}(0,n)  &=& \prod _{i_0=0}^{n-1}B_{2i_0+1} \label{eq:100012a}\\
\bar{c}(1,n) &=&  \sum_{i_0=0}^{n} \left\{ A_{2i_0} \prod _{i_1=0}^{i_0-1}B_{2i_1+1} \prod _{i_2=i_0}^{n-1}B_{2i_2+2} \right\} 
\label{eq:100012b}\\
\bar{c}(\tau ,n) &=& \sum_{i_0=0}^{n} \left\{A_{2i_0}\prod _{i_1=0}^{i_0-1} B_{2i_1+1} 
\prod _{k=1}^{\tau -1} \left( \sum_{i_{2k}= i_{2(k-1)}}^{n} A_{2i_{2k}+k}\prod _{i_{2k+1}=i_{2(k-1)}}^{i_{2k}-1}B_{2i_{2k+1}+(k+1)}\right) \right. \nonumber\\
&& \times  \left. \prod _{i_{2\tau}=i_{2(\tau -1)}}^{n-1} B_{2i_{2\tau }+(\tau +1)} \right\} 
\hspace{1cm}\label{eq:100012c} 
\end{eqnarray}
and
\begin{eqnarray}
y_0^m(x) &=& c_0 x^{\lambda } \sum_{i_0=0}^{m} \left\{ \prod _{i_1=0}^{i_0-1}B_{2i_1+1} \right\} x^{2i_0 } \label{eq:100013a}\\
y_1^m(x) &=& c_0 x^{\lambda } \sum_{i_0=0}^{m}\left\{ A_{2i_0} \prod _{i_1=0}^{i_0-1}B_{2i_1+1}  \sum_{i_2=i_0}^{m} \left\{ \prod _{i_3=i_0}^{i_2-1}B_{2i_3+2} \right\}\right\} x^{2i_2+1 } \label{eq:100013b}\\
y_{\tau }^m(x) &=& c_0 x^{\lambda } \sum_{i_0=0}^{m} \left\{A_{2i_0}\prod _{i_1=0}^{i_0-1} B_{2i_1+1} 
\prod _{k=1}^{\tau -1} \left( \sum_{i_{2k}= i_{2(k-1)}}^{m} A_{2i_{2k}+k}\prod _{i_{2k+1}=i_{2(k-1)}}^{i_{2k}-1}B_{2i_{2k+1}+(k+1)}\right) \right. \nonumber\\
&& \times  \left. \sum_{i_{2\tau} = i_{2(\tau -1)}}^{m} \left( \prod _{i_{2\tau +1}=i_{2(\tau -1)}}^{i_{2\tau}-1} B_{2i_{2\tau +1}+(\tau +1)} \right) \right\} x^{2i_{2\tau}+\tau }\;\;\;\mathrm{where}\;\tau \geq 2
\label{eq:100013c} 
\end{eqnarray}
\end{enumerate}
\end{theorem} 
Put $n= j+1$ in (\ref{eq:10007b}) and use the condition $B_{j+1}=0$ for $\alpha $. We obtain two possible fixed values for $\alpha $ such as 
\begin{subequations}
\begin{equation}
\alpha  = -2\left( j+\frac{1}{2}+\lambda \right) \label{eq:100014a}
\end{equation}
\begin{equation}
\alpha = 2\left( j+\lambda \right) \label{eq:100014b}
\end{equation}
\end{subequations}
Fro the case of $ \alpha  = -2\left( j+\frac{1}{2}+\lambda \right) $, take (\ref{eq:100014a}) into (\ref{eq:10007a}) and (\ref{eq:10007b}).
\begin{subequations}
\begin{eqnarray}
A_n = -\frac{(2a-b-c)}{(a-b)(a-c)} \frac{\left( n+ \lambda \right)^2 - 4\varpi _j(q) }{(n+1+\lambda )(n+\frac{1}{2}+\lambda )} \label{eq:100015a}
\end{eqnarray}
\begin{equation}
B_n =\frac{-1}{(a-b)(a-c)}\frac{\left( n-j-1 \right) \left( n+j-\frac{1}{2} +2\lambda \right)}{(n+1+\lambda )(n+\frac{1}{2}+\lambda )}\label{eq:100015b}
\end{equation}
and
\begin{equation}
\varpi _j(q) = \frac{ \left( j+\lambda \right)\left( j+\frac{1}{2}+\lambda\right) a-\frac{q}{4}}{ 4(2a-b-c)} \nonumber
\end{equation}
\end{subequations}
Fro the case of $ \alpha  = 2\left( j+\lambda \right) $, take (\ref{eq:100014b}) into (\ref{eq:10007a}) and (\ref{eq:10007b}). Solutions for the new $A_n$ and $B_n$ terms are also equivalent to (\ref{eq:100015a}) and (\ref{eq:100015b}).

According to (\ref{eq:10008}), $c_{j+1}=0$ is clearly an algebraic equation in $q $ of degree $j+1$ and thus has $j+1$ zeros denoted them by $q_j^m$ eigenvalues where $m = 0,1,2, \cdots, j$. They can be arranged in the following order: $q_j^0 < q_j^1 < q _j^2 < \cdots < q_j^j$. 

Substitute (\ref{eq:100015a}) and (\ref{eq:100015b}) into (\ref{eq:100012a})--(\ref{eq:100013c}). Replace $x$ by $z$ in (\ref{eq:10009b}), (\ref{eq:100010b}), (\ref{eq:100011b}) and (\ref{eq:100013a})--(\ref{eq:100013c}).

As $B_{1}= c_{1}=0$, take the new (\ref{eq:100012b}) into (\ref{eq:10009a}) putting $j=0$. Substitute the new (\ref{eq:100013a}) into (\ref{eq:10009b}) putting $j=0$. 

As $B_{2N+2}= c_{2N+2}=0$, take the new (\ref{eq:100012a})--(\ref{eq:100012c}) into (\ref{eq:100010a}) putting $j=2N+1$. Substitute the new 
(\ref{eq:100013a})--(\ref{eq:100013c}) into (\ref{eq:100010b}) putting $j=2N+1$ and $q =q_{2N+1}^m$.

As $B_{2N+3}= c_{2N+3}=0$, take the new (\ref{eq:100012a})--(\ref{eq:100012c}) into (\ref{eq:100011a}) putting $j=2N+2$. Substitute the new 
(\ref{eq:100013a})--(\ref{eq:100013c}) into (\ref{eq:100011b}) putting $j=2N+2$ and $q =q_{2N+2}^m$.

After the replacement process, the general expression of power series of Lam\'{e} equation about $z=0$ for the first species complete polynomial using 3-term recurrence formula and its algebraic equation for the determination of an accessory parameter $q$ are given by

\begin{enumerate} 
\item As $\alpha = -2\left( \frac{1}{2}+\lambda \right) $ or $2\lambda $,

An algebraic equation of degree 1 for the determination of $q$ is given by
\begin{equation}
0= \bar{c}(1,0;0,q)= q+4\lambda \left( (a-b-c)\lambda -\frac{a}{2}\right) \label{eq:100016a}
\end{equation}
The eigenvalue of $q$ is written by $q_0^0$. Its eigenfunction is given by
\begin{equation}
y(z) = y_0^0\left( 0,q_0^0;z\right)= c_0 z^{\lambda } \label{eq:100016b}
\end{equation}
\item As $\alpha  =-2\left( 2N+\frac{3}{2}+\lambda \right)$ or $2\left( 2N+1+\lambda \right) $ where $N \in \mathbb{N}_{0}$,

An algebraic equation of degree $2N+2$ for the determination of $q$ is given by
\begin{equation}
0 = \sum_{r=0}^{N+1}\bar{c}\left( 2r, N+1-r; 2N+1,q \right)  \label{eq:100017a}
\end{equation}
The eigenvalue of $q$ is written by $q_{2N+1}^m$ where $m = 0,1,2,\cdots,2N+1 $; $q_{2N+1}^0 < q_{2N+1}^1 < \cdots < q_{2N+1}^{2N+1}$. Its eigenfunction is given by 
\begin{equation} 
y(z) = \sum_{r=0}^{N} y_{2r}^{N-r}\left( 2N+1,q_{2N+1}^m;z\right)+ \sum_{r=0}^{N} y_{2r+1}^{N-r}\left( 2N+1,q_{2N+1}^m;z\right)
\label{eq:100017b} 
\end{equation}
\item As  $\alpha  =-2\left( 2N+\frac{5}{2}+\lambda \right)$ or $2\left( 2N+2+\lambda \right) $ where $N \in \mathbb{N}_{0}$,

An algebraic equation of degree $2N+3$ for the determination of $q$ is given by
\begin{equation}  
0 = \sum_{r=0}^{N+1}\bar{c}\left( 2r+1, N+1-r; 2N+2,q \right) \label{eq:100018a}
\end{equation}
The eigenvalue of $q$ is written by $q_{2N+2}^m$ where $m = 0,1,2,\cdots,2N+2 $; $q_{2N+2}^0 < q_{2N+2}^1 < \cdots < q_{2N+2}^{2N+2}$. Its eigenfunction is given by
\begin{equation} 
y(z) =  \sum_{r=0}^{N+1} y_{2r}^{N+1-r}\left( 2N+2,q_{2N+2}^m;z\right) + \sum_{r=0}^{N} y_{2r+1}^{N-r}\left( 2N+2,q_{2N+2}^m;z\right) \label{eq:100018b}
\end{equation}
In the above,
\begin{eqnarray}
\bar{c}(0,n;j,q)  &=& \frac{\left( -\frac{j}{2}\right)_{n}\left( \frac{1}{4}+\frac{j}{2}+\lambda \right)_{n}}{ \left(  1+\frac{\lambda }{2} \right)_{n}\left( \frac{3}{4}+\frac{\lambda }{2} \right)_{n}} \overline{\eta}^{n}\label{eq:100019a}\\
\bar{c}(1,n;j,q) &=& \overline{\mu} \sum_{i_0=0}^{n}\frac{\left( i_0 +\frac{\lambda }{2}\right)^2 -\varpi _j(q) }{\left( i_0+\frac{1}{2}+\frac{\lambda }{2}\right) \left( i_0+\frac{1}{4}+\frac{\lambda }{2}\right)} \frac{\left( -\frac{j}{2}\right)_{i_0}\left( \frac{1}{4}+\frac{j}{2}+\lambda \right)_{i_0} }{\left( 1+\frac{\lambda }{2}\right)_{i_0} \left( \frac{3}{4} + \frac{\lambda }{2}\right)_{i_0}} \nonumber\\
&&\times \frac{\left( \frac{1}{2}-\frac{j}{2} \right)_{n} \left( \frac{3}{4}+\frac{j}{2}+ \lambda \right)_{n} \left( \frac{3}{2}+ \frac{\lambda }{2}\right)_{i_0}\left( \frac{5}{4}+ \frac{\lambda }{2}\right)_{i_0}}{\left( \frac{1}{2}-\frac{j}{2} \right)_{i_0} \left( \frac{3}{4}+\frac{j}{2}+ \lambda \right)_{i_0} \left( \frac{3}{2}+ \frac{\lambda }{2}\right)_n \left( \frac{5}{4}+ \frac{\lambda }{2}\right)_n} \overline{\eta}^{n }  
\label{eq:100019b}\\
\bar{c}(\tau ,n;j,q) &=& \overline{\mu}^{\tau } \sum_{i_0=0}^{n} \frac{\left( i_0 +\frac{\lambda }{2}\right)^2 -\varpi _j(q) }{\left( i_0+\frac{1}{2}+\frac{\lambda }{2}\right) \left( i_0+\frac{1}{4}+\frac{\lambda }{2}\right)} \frac{\left( -\frac{j}{2}\right)_{i_0}\left( \frac{1}{4}+\frac{j}{2}+\lambda \right)_{i_0} }{\left( 1+\frac{\lambda }{2}\right)_{i_0} \left( \frac{3}{4} + \frac{\lambda }{2}\right)_{i_0}}  \nonumber\\
&&\times \prod_{k=1}^{\tau -1} \left( \sum_{i_k = i_{k-1}}^{n} \frac{\left( i_k+ \frac{k}{2} +\frac{\lambda }{2}\right)^2 -\varpi _j(q)}{\left( i_k+\frac{k}{2}+\frac{1}{2}+\frac{\lambda }{2}\right) \left( i_k+\frac{k}{2}+\frac{1}{4}+\frac{\lambda }{2}\right)} \right. \nonumber\\
&&\times \left. \frac{\left( \frac{k}{2}-\frac{j}{2}\right)_{i_k} \left( \frac{k}{2}+\frac{1}{4}+\frac{j}{2}+\lambda \right)_{i_k}\left( \frac{k}{2}+1+ \frac{\lambda }{2}\right)_{i_{k-1}} \left( \frac{k}{2}+ \frac{3}{4} + \frac{\lambda }{2}\right)_{i_{k-1}}}{\left( \frac{k}{2}-\frac{j}{2}\right)_{i_{k-1}} \left( \frac{k}{2}+\frac{1}{4}+\frac{j}{2}+\lambda \right)_{i_{k-1}}\left( \frac{k}{2}+1+ \frac{\lambda }{2}\right)_{i_k} \left( \frac{k}{2}+ \frac{3}{4} + \frac{\lambda }{2}\right)_{i_k}} \right)  \nonumber \\ 
&&\times \frac{\left( \frac{\tau }{2} -\frac{j}{2}\right)_{n}\left( \frac{\tau }{2}+\frac{1}{4} +\frac{j}{2}+\lambda \right)_{n} \left( \frac{\tau }{2}+1+\frac{\lambda }{2}\right)_{i_{\tau -1}} \left( \frac{\tau }{2}+\frac{3}{4} + \frac{\lambda }{2}\right)_{i_{\tau -1}}}{\left( \frac{\tau }{2} -\frac{j}{2}\right)_{i_{\tau -1}}\left( \frac{\tau }{2}+\frac{1}{4} +\frac{j}{2}+\lambda \right)_{i_{\tau -1}} \left( \frac{\tau }{2}+1+\frac{\lambda }{2}\right)_n \left( \frac{\tau }{2}+\frac{3}{4} + \frac{\lambda }{2}\right)_n} \overline{\eta}^{n }\hspace{1.5cm} \label{eq:100019c}
\end{eqnarray}
\begin{eqnarray}
y_0^m(j,q ;z) &=& c_0 z^{\lambda }  \sum_{i_0=0}^{m} \frac{\left( -\frac{j}{2}\right)_{i_0}\left( \frac{1}{4}+\frac{j}{2}+\lambda \right)_{i_0} }{\left( 1+\frac{\lambda }{2}\right)_{i_0} \left( \frac{3}{4} + \frac{\lambda }{2}\right)_{i_0}} \eta ^{i_0} \label{eq:100020a}\\
y_1^m(j,q ;z) &=& c_0 z^{\lambda } \left\{\sum_{i_0=0}^{m} \frac{\left( i_0 +\frac{\lambda }{2}\right)^2 -\varpi _j(q) }{\left( i_0+\frac{1}{2}+\frac{\lambda }{2}\right) \left( i_0+\frac{1}{4}+\frac{\lambda }{2}\right)} \frac{\left( -\frac{j}{2}\right)_{i_0}\left( \frac{1}{4}+\frac{j}{2}+\lambda \right)_{i_0} }{\left( 1+\frac{\lambda }{2}\right)_{i_0} \left( \frac{3}{4} + \frac{\lambda }{2}\right)_{i_0}} \right. \nonumber\\
&&\times \left. \sum_{i_1 = i_0}^{m} \frac{\left( \frac{1}{2}-\frac{j}{2} \right)_{i_1} \left( \frac{3}{4}+\frac{j}{2}+ \lambda \right)_{i_1} \left( \frac{3}{2}+ \frac{\lambda }{2}\right)_{i_0}\left( \frac{5}{4}+ \frac{\lambda }{2}\right)_{i_0}}{\left( \frac{1}{2}-\frac{j}{2} \right)_{i_0} \left( \frac{3}{4}+\frac{j}{2}+ \lambda \right)_{i_0} \left( \frac{3}{2}+ \frac{\lambda }{2}\right)_{i_1} \left( \frac{5}{4}+ \frac{\lambda }{2}\right)_{i_1}} \eta ^{i_1}\right\} \mu 
\label{eq:100020b}\\
y_{\tau }^m(j,q ;z) &=& c_0 z^{\lambda } \left\{ \sum_{i_0=0}^{m} \frac{\left( i_0 +\frac{\lambda }{2}\right)^2 -\varpi _j(q) }{\left( i_0+\frac{1}{2}+\frac{\lambda }{2}\right) \left( i_0+\frac{1}{4}+\frac{\lambda }{2}\right)} \frac{\left( -\frac{j}{2}\right)_{i_0}\left( \frac{1}{4}+\frac{j}{2}+\lambda \right)_{i_0} }{\left( 1+\frac{\lambda }{2}\right)_{i_0} \left( \frac{3}{4} + \frac{\lambda }{2}\right)_{i_0}} \right.\nonumber\\
&&\times \prod_{k=1}^{\tau -1} \left( \sum_{i_k = i_{k-1}}^{m} \frac{\left( i_k+ \frac{k}{2} +\frac{\lambda }{2}\right)^2 -\varpi _j(q)}{\left( i_k+\frac{k}{2}+\frac{1}{2}+\frac{\lambda }{2}\right) \left( i_k+\frac{k}{2}+\frac{1}{4}+\frac{\lambda }{2}\right)} \right. \nonumber\\
&&\times \left. \frac{\left( \frac{k}{2}-\frac{j}{2}\right)_{i_k} \left( \frac{k}{2}+\frac{1}{4}+\frac{j}{2}+\lambda \right)_{i_k}\left( \frac{k}{2}+1+ \frac{\lambda }{2}\right)_{i_{k-1}} \left( \frac{k}{2}+ \frac{3}{4} + \frac{\lambda }{2}\right)_{i_{k-1}}}{\left( \frac{k}{2}-\frac{j}{2}\right)_{i_{k-1}} \left( \frac{k}{2}+\frac{1}{4}+\frac{j}{2}+\lambda \right)_{i_{k-1}}\left( \frac{k}{2}+1+ \frac{\lambda }{2}\right)_{i_k} \left( \frac{k}{2}+ \frac{3}{4} + \frac{\lambda }{2}\right)_{i_k}} \right) \label{eq:100020c}\\
&&\times \left. \sum_{i_{\tau } = i_{\tau -1}}^{m} \frac{\left( \frac{\tau }{2} -\frac{j}{2}\right)_{i_{\tau }}\left( \frac{\tau }{2}+\frac{1}{4} +\frac{j}{2}+\lambda \right)_{i_{\tau }} \left( \frac{\tau }{2}+1+\frac{\lambda }{2}\right)_{i_{\tau -1}} \left( \frac{\tau }{2}+\frac{3}{4} + \frac{\lambda }{2}\right)_{i_{\tau -1}}}{\left( \frac{\tau }{2} -\frac{j}{2}\right)_{i_{\tau -1}}\left( \frac{\tau }{2}+\frac{1}{4} +\frac{j}{2}+\lambda \right)_{i_{\tau -1}} \left( \frac{\tau }{2}+1+\frac{\lambda }{2}\right)_{i_{\tau }} \left( \frac{\tau }{2}+\frac{3}{4} + \frac{\lambda }{2}\right)_{i_{\tau }}} \eta ^{i_{\tau }}\right\} \mu ^{\tau }  \nonumber 
\end{eqnarray}
where
\begin{equation}
\begin{cases} \tau \geq 2 \cr
\eta = \frac{-z^2}{(a-b)(a-c)}  \cr
\mu =  \frac{-(2a-b-c)}{(a-b)(a-c)} z \cr
\overline{\eta}= \frac{-1}{(a-b)(a-c)}  \cr
\overline{\mu}= \frac{-(2a-b-c)}{(a-b)(a-c)}  \cr
\varpi _j(q) = \frac{ \left( j+\lambda \right)\left( j+\frac{1}{2}+\lambda\right) a-\frac{q}{4}}{ 4(2a-b-c)}
\end{cases}\nonumber
\end{equation}
\end{enumerate}
Put $c_0$= 1 as $\lambda =0$ for the first kind of independent solutions of Lam\'{e} equation and $\lambda = \frac{1}{2}$ for the second one in (\ref{eq:100016a})--(\ref{eq:100020c}). 
\begin{remark}
The power series expansion of Lam\'{e} equation in the algebraic form of the first kind for the first species complete polynomial using 3TRF about $z=0$ is given by
\begin{enumerate} 
\item As $\alpha = -1 $ or 0 with $q=q_0^0=0 $,
 
The eigenfunction is given by
\begin{eqnarray}
y(z) &=& L_p^{(a)}F_{0,0}\left( a,b,c,q=q_0^0=0, \alpha = -1\; \mbox{or}\; 0; z=x-a,\eta = \frac{-z^2}{(a-b)(a-c)}\right. \nonumber\\
&&,\left. \mu = \frac{-(2a-b-c)}{(a-b)(a-c)} z \right) \nonumber\\
&=& 1 \label{eq:100021a}
\end{eqnarray}
\item As $\alpha  =-2\left( 2N+\frac{3}{2} \right)$ or $2\left( 2N+1 \right) $ where $N \in \mathbb{N}_{0}$,

An algebraic equation of degree $2N+2$ for the determination of $q$ is given by
\begin{equation}
0 = \sum_{r=0}^{N+1}\bar{c}\left( 2r, N+1-r; 2N+1,q \right)  \label{eq:100022a}
\end{equation}
The eigenvalue of $q$ is written by $q_{2N+1}^m$ where $m = 0,1,2,\cdots,2N+1 $; $q_{2N+1}^0 < q_{2N+1}^1 < \cdots < q_{2N+1}^{2N+1}$. Its eigenfunction is given by 
\begin{eqnarray} 
y(z) &=& L_p^{(a)}F_{2N+1,m}\left( a,b,c,q=q_{2N+1}^m, \alpha = -2\left( 2N+\frac{3}{2} \right)\; \mbox{or}\; 2\left( 2N+1 \right); z=x-a \right. \nonumber\\
&&,\left. \eta = \frac{-z^2}{(a-b)(a-c)},\mu = \frac{-(2a-b-c)}{(a-b)(a-c)} z \right) \nonumber\\
&=& \sum_{r=0}^{N} y_{2r}^{N-r}\left( 2N+1,q_{2N+1}^m;z\right)+ \sum_{r=0}^{N} y_{2r+1}^{N-r}\left( 2N+1,q_{2N+1}^m;z\right)
\label{eq:100022b} 
\end{eqnarray}
\item As  $\alpha  =-2\left( 2N+\frac{5}{2} \right)$ or $2\left( 2N+2 \right) $ where $N \in \mathbb{N}_{0}$,

An algebraic equation of degree $2N+3$ for the determination of $q$ is given by
\begin{equation}  
0 = \sum_{r=0}^{N+1}\bar{c}\left( 2r+1, N+1-r; 2N+2,q \right) \label{eq:100023a}
\end{equation}
The eigenvalue of $q$ is written by $q_{2N+2}^m$ where $m = 0,1,2,\cdots,2N+2 $; $q_{2N+2}^0 < q_{2N+2}^1 < \cdots < q_{2N+2}^{2N+2}$. Its eigenfunction is given by
\begin{eqnarray} 
y(z) &=& L_p^{(a)}F_{2N+2,m}\left( a,b,c,q=q_{2N+2}^m, \alpha = -2\left( 2N+\frac{5}{2} \right)\; \mbox{or}\; 2\left( 2N+2 \right); z=x-a \right. \nonumber\\
&&,\left. \eta = \frac{-z^2}{(a-b)(a-c)},\mu = \frac{-(2a-b-c)}{(a-b)(a-c)} z \right) \nonumber\\
&=& \sum_{r=0}^{N+1} y_{2r}^{N+1-r}\left( 2N+2,q_{2N+2}^m;z\right) + \sum_{r=0}^{N} y_{2r+1}^{N-r}\left( 2N+2,q_{2N+2}^m;z\right) \label{eq:100023b}
\end{eqnarray}
In the above,
\begin{align}
\bar{c}(0,n;j,q)  &=  \frac{\left( -\frac{j}{2}\right)_{n}\left( \frac{1}{4}+\frac{j}{2} \right)_{n}}{ \left(  1 \right)_{n}\left( \frac{3}{4} \right)_{n}} \overline{\eta}^{n}\label{eq:100024a}\\
\bar{c}(1,n;j,q) &=  \overline{\mu} \sum_{i_0=0}^{n}\frac{ i_0^2 -\varpi _j(q) }{\left( i_0+\frac{1}{2} \right) \left( i_0+\frac{1}{4} \right)} \frac{\left( -\frac{j}{2}\right)_{i_0}\left( \frac{1}{4}+\frac{j}{2}  \right)_{i_0} }{\left( 1 \right)_{i_0} \left( \frac{3}{4} \right)_{i_0}}   \frac{\left( \frac{1}{2}-\frac{j}{2} \right)_{n} \left( \frac{3}{4}+\frac{j}{2} \right)_{n} \left( \frac{3}{2} \right)_{i_0}\left( \frac{5}{4} \right)_{i_0}}{\left( \frac{1}{2}-\frac{j}{2} \right)_{i_0} \left( \frac{3}{4}+\frac{j}{2} \right)_{i_0} \left( \frac{3}{2} \right)_n \left( \frac{5}{4} \right)_n} \overline{\eta}^{n }   
\label{eq:100024b}\\
\bar{c}(\tau ,n;j,q) &=  \overline{\mu}^{\tau } \sum_{i_0=0}^{n} \frac{ i_0^2 -\varpi _j(q) }{\left( i_0+\frac{1}{2} \right) \left( i_0+\frac{1}{4} \right)} \frac{\left( -\frac{j}{2}\right)_{i_0}\left( \frac{1}{4}+\frac{j}{2} \right)_{i_0} }{\left( 1 \right)_{i_0} \left( \frac{3}{4} \right)_{i_0}}  \nonumber\\
 &\times \prod_{k=1}^{\tau -1} \left( \sum_{i_k = i_{k-1}}^{n} \frac{\left( i_k+ \frac{k}{2} \right)^2 -\varpi _j(q)}{\left( i_k+\frac{k}{2}+\frac{1}{2} \right) \left( i_k+\frac{k}{2}+\frac{1}{4} \right)} \right.   \left. \frac{\left( \frac{k}{2}-\frac{j}{2}\right)_{i_k} \left( \frac{k}{2}+\frac{1}{4}+\frac{j}{2} \right)_{i_k}\left( \frac{k}{2}+1 \right)_{i_{k-1}} \left( \frac{k}{2}+ \frac{3}{4} \right)_{i_{k-1}}}{\left( \frac{k}{2}-\frac{j}{2}\right)_{i_{k-1}} \left( \frac{k}{2}+\frac{1}{4}+\frac{j}{2} \right)_{i_{k-1}}\left( \frac{k}{2}+1 \right)_{i_k} \left( \frac{k}{2}+ \frac{3}{4} \right)_{i_k}} \right) \nonumber\\ 
 &\times \frac{\left( \frac{\tau }{2} -\frac{j}{2}\right)_{n}\left( \frac{\tau }{2}+\frac{1}{4} +\frac{j}{2} \right)_{n} \left( \frac{\tau }{2}+1 \right)_{i_{\tau -1}} \left( \frac{\tau }{2}+\frac{3}{4} \right)_{i_{\tau -1}}}{\left( \frac{\tau }{2} -\frac{j}{2}\right)_{i_{\tau -1}}\left( \frac{\tau }{2}+\frac{1}{4} +\frac{j}{2} \right)_{i_{\tau -1}} \left( \frac{\tau }{2}+1 \right)_n \left( \frac{\tau }{2}+\frac{3}{4} \right)_n} \overline{\eta}^{n } \label{eq:100024c} 
\end{align}
\begin{align}
y_0^m(j,q ;z) &=  \sum_{i_0=0}^{m} \frac{\left( -\frac{j}{2}\right)_{i_0}\left( \frac{1}{4}+\frac{j}{2} \right)_{i_0} }{\left( 1 \right)_{i_0} \left( \frac{3}{4} \right)_{i_0}} \eta ^{i_0} \label{eq:100025a}\\
y_1^m(j,q ;z) &=  \left\{\sum_{i_0=0}^{m} \frac{ i_0 ^2 -\varpi _j(q) }{\left( i_0+\frac{1}{2} \right) \left( i_0+\frac{1}{4} \right)} \frac{\left( -\frac{j}{2}\right)_{i_0}\left( \frac{1}{4}+\frac{j}{2} \right)_{i_0} }{\left( 1 \right)_{i_0} \left( \frac{3}{4} \right)_{i_0}} \right.   \left. \sum_{i_1 = i_0}^{m} \frac{\left( \frac{1}{2}-\frac{j}{2} \right)_{i_1} \left( \frac{3}{4}+\frac{j}{2} \right)_{i_1} \left( \frac{3}{2} \right)_{i_0}\left( \frac{5}{4} \right)_{i_0}}{\left( \frac{1}{2}-\frac{j}{2} \right)_{i_0} \left( \frac{3}{4}+\frac{j}{2} \right)_{i_0} \left( \frac{3}{2} \right)_{i_1} \left( \frac{5}{4} \right)_{i_1}} \eta ^{i_1}\right\} \mu  \label{eq:100025b}\\
y_{\tau }^m(j,q ;z) &= \left\{ \sum_{i_0=0}^{m} \frac{ i_0 ^2 -\varpi _j(q) }{\left( i_0+\frac{1}{2} \right) \left( i_0+\frac{1}{4} \right)} \frac{\left( -\frac{j}{2}\right)_{i_0}\left( \frac{1}{4}+\frac{j}{2} \right)_{i_0} }{\left( 1 \right)_{i_0} \left( \frac{3}{4} \right)_{i_0}} \right.\nonumber\\
 &\times \prod_{k=1}^{\tau -1} \left( \sum_{i_k = i_{k-1}}^{m} \frac{\left( i_k+ \frac{k}{2} \right)^2 -\varpi _j(q)}{\left( i_k+\frac{k}{2}+\frac{1}{2} \right) \left( i_k+\frac{k}{2}+\frac{1}{4} \right)} \right.  \left. \frac{\left( \frac{k}{2}-\frac{j}{2}\right)_{i_k} \left( \frac{k}{2}+\frac{1}{4}+\frac{j}{2} \right)_{i_k}\left( \frac{k}{2}+1 \right)_{i_{k-1}} \left( \frac{k}{2}+ \frac{3}{4} \right)_{i_{k-1}}}{\left( \frac{k}{2}-\frac{j}{2}\right)_{i_{k-1}} \left( \frac{k}{2}+\frac{1}{4}+\frac{j}{2} \right)_{i_{k-1}}\left( \frac{k}{2}+1 \right)_{i_k} \left( \frac{k}{2}+ \frac{3}{4} \right)_{i_k}} \right) \nonumber\\
 &\times \left. \sum_{i_{\tau } = i_{\tau -1}}^{m} \frac{\left( \frac{\tau }{2} -\frac{j}{2}\right)_{i_{\tau }}\left( \frac{\tau }{2}+\frac{1}{4} +\frac{j}{2} \right)_{i_{\tau }} \left( \frac{\tau }{2}+1 \right)_{i_{\tau -1}} \left( \frac{\tau }{2}+\frac{3}{4}  \right)_{i_{\tau -1}}}{\left( \frac{\tau }{2} -\frac{j}{2}\right)_{i_{\tau -1}}\left( \frac{\tau }{2}+\frac{1}{4} +\frac{j}{2} \right)_{i_{\tau -1}} \left( \frac{\tau }{2}+1 \right)_{i_{\tau }} \left( \frac{\tau }{2}+\frac{3}{4} \right)_{i_{\tau }}} \eta ^{i_{\tau }}\right\} \mu ^{\tau }  \label{eq:100025c}
\end{align}
where
\begin{equation}
\begin{cases} \tau \geq 2 \cr
\overline{\eta}= \frac{-1}{(a-b)(a-c)}  \cr
\overline{\mu}= \frac{-(2a-b-c)}{(a-b)(a-c)}  \cr
\varpi _j(q) = \frac{ j\left( j+\frac{1}{2} \right) a-\frac{q}{4}}{ 4(2a-b-c)}
\end{cases}\nonumber
\end{equation}
\end{enumerate}
\end{remark}
\begin{remark}
The power series expansion of Lam\'{e} equation in the algebraic form of the second kind for the first species complete polynomial using 3TRF about $z=0$ is given by
\begin{enumerate} 
\item As $\alpha = -2 $ or 1 with $q=q_0^0=b+c $,
 
The eigenfunction is given by
\begin{eqnarray}
y(z) &=& L_p^{(a)}S_{0,0}\left( a,b,c,q=q_0^0=b+c, \alpha = -2\; \mbox{or}\; 1; z=x-a,\eta = \frac{-z^2}{(a-b)(a-c)}\right.\nonumber\\ 
&&,\left. \mu = \frac{-(2a-b-c)}{(a-b)(a-c)} z \right) \nonumber\\
&=& z^{\frac{1}{2}} \label{eq:100026a}
\end{eqnarray}
\item As $\alpha  =-2\left( 2N+2 \right)$ or $2\left( 2N+\frac{3}{2} \right) $ where $N \in \mathbb{N}_{0}$,

An algebraic equation of degree $2N+2$ for the determination of $q$ is given by
\begin{equation}
0 = \sum_{r=0}^{N+1}\bar{c}\left( 2r, N+1-r; 2N+1,q \right)  \label{eq:100027a}
\end{equation}
The eigenvalue of $q$ is written by $q_{2N+1}^m$ where $m = 0,1,2,\cdots,2N+1 $; $q_{2N+1}^0 < q_{2N+1}^1 < \cdots < q_{2N+1}^{2N+1}$. Its eigenfunction is given by 
\begin{eqnarray} 
y(z) &=& L_p^{(a)}S_{2N+1,m}\left( a,b,c,q=q_{2N+1}^m, \alpha = -2\left( 2N+2 \right)\; \mbox{or}\; 2\left( 2N+\frac{3}{2} \right); z=x-a \right. \nonumber\\
&&,\left. \eta = \frac{-z^2}{(a-b)(a-c)},\mu = \frac{-(2a-b-c)}{(a-b)(a-c)} z \right) \nonumber\\
&=& \sum_{r=0}^{N} y_{2r}^{N-r}\left( 2N+1,q_{2N+1}^m;z\right)+ \sum_{r=0}^{N} y_{2r+1}^{N-r}\left( 2N+1,q_{2N+1}^m;z\right)
\label{eq:100027b} 
\end{eqnarray}
\item As  $\alpha  =-2\left( 2N+3 \right)$ or $2\left( 2N+\frac{5}{2} \right) $ where $N \in \mathbb{N}_{0}$,

An algebraic equation of degree $2N+3$ for the determination of $q$ is given by
\begin{equation}  
0 = \sum_{r=0}^{N+1}\bar{c}\left( 2r+1, N+1-r; 2N+2,q \right) \label{eq:100028a}
\end{equation}
The eigenvalue of $q$ is written by $q_{2N+2}^m$ where $m = 0,1,2,\cdots,2N+2 $; $q_{2N+2}^0 < q_{2N+2}^1 < \cdots < q_{2N+2}^{2N+2}$. Its eigenfunction is given by
\begin{eqnarray} 
y(z) &=& L_p^{(a)}S_{2N+2,m}\left( a,b,c,q=q_{2N+2}^m, \alpha = -2\left( 2N+3 \right)\; \mbox{or}\; 2\left( 2N+\frac{5}{2} \right); z=x-a \right. \nonumber\\
&&,\left. \eta = \frac{-z^2}{(a-b)(a-c)},\mu = \frac{-(2a-b-c)}{(a-b)(a-c)} z \right) \nonumber\\
&=& \sum_{r=0}^{N+1} y_{2r}^{N+1-r}\left( 2N+2,q_{2N+2}^m;z\right) + \sum_{r=0}^{N} y_{2r+1}^{N-r}\left( 2N+2,q_{2N+2}^m;z\right) \label{eq:100028b}
\end{eqnarray}
In the above,
\begin{align}
\bar{c}(0,n;j,q)  &=  \frac{\left( -\frac{j}{2}\right)_{n}\left( \frac{3}{4}+\frac{j}{2} \right)_{n}}{ \left(  \frac{5}{4} \right)_{n}\left( 1 \right)_{n}} \overline{\eta}^{n}\label{eq:100029a}\\
\bar{c}(1,n;j,q) &=  \overline{\mu} \sum_{i_0=0}^{n}\frac{\left( i_0 +\frac{1}{4}\right)^2 -\varpi _j(q) }{\left( i_0+\frac{3}{4} \right) \left( i_0+\frac{1}{2} \right)} \frac{\left( -\frac{j}{2}\right)_{i_0}\left( \frac{3}{4}+\frac{j}{2} \right)_{i_0} }{\left( \frac{5}{4}\right)_{i_0} \left( 1\right)_{i_0}}   \frac{\left( \frac{1}{2}-\frac{j}{2} \right)_{n} \left( \frac{5}{4}+\frac{j}{2} \right)_{n} \left( \frac{7}{4} \right)_{i_0}\left( \frac{3}{2} \right)_{i_0}}{\left( \frac{1}{2}-\frac{j}{2} \right)_{i_0} \left( \frac{5}{4}+\frac{j}{2}  \right)_{i_0} \left( \frac{7}{4} \right)_n \left( \frac{3}{2} \right)_n} \overline{\eta}^{n }  
 \label{eq:100029b}\\
\bar{c}(\tau ,n;j,q) &=  \overline{\mu}^{\tau } \sum_{i_0=0}^{n} \frac{\left( i_0 +\frac{1}{4}\right)^2 -\varpi _j(q) }{\left( i_0+\frac{3}{4} \right) \left( i_0+\frac{1}{2} \right)} \frac{\left( -\frac{j}{2}\right)_{i_0}\left( \frac{3}{4}+\frac{j}{2} \right)_{i_0} }{\left( \frac{5}{4}\right)_{i_0} \left( 1\right)_{i_0}}  \nonumber\\
 &\times \prod_{k=1}^{\tau -1} \left( \sum_{i_k = i_{k-1}}^{n} \frac{\left( i_k+ \frac{k}{2} +\frac{1}{4}\right)^2 -\varpi _j(q)}{\left( i_k+\frac{k}{2}+\frac{3}{4} \right) \left( i_k+\frac{k}{2}+\frac{1}{2} \right)} \right.   \left. \frac{\left( \frac{k}{2}-\frac{j}{2}\right)_{i_k} \left( \frac{k}{2}+\frac{3}{4}+\frac{j}{2} \right)_{i_k}\left( \frac{k}{2} + \frac{5}{4}\right)_{i_{k-1}} \left( \frac{k}{2}+ 1\right)_{i_{k-1}}}{\left( \frac{k}{2}-\frac{j}{2}\right)_{i_{k-1}} \left( \frac{k}{2}+\frac{3}{4}+\frac{j}{2}  \right)_{i_{k-1}}\left( \frac{k}{2} + \frac{5}{4}\right)_{i_k} \left( \frac{k}{2}+ 1\right)_{i_k}} \right) \nonumber \\ 
 &\times \frac{\left( \frac{\tau }{2} -\frac{j}{2}\right)_{n}\left( \frac{\tau }{2}+\frac{3}{4} +\frac{j}{2}  \right)_{n} \left( \frac{\tau }{2} +\frac{5}{4}\right)_{i_{\tau -1}} \left( \frac{\tau }{2}+1\right)_{i_{\tau -1}}}{\left( \frac{\tau }{2} -\frac{j}{2}\right)_{i_{\tau -1}}\left( \frac{\tau }{2}+\frac{3}{4} +\frac{j}{2} \right)_{i_{\tau -1}} \left( \frac{\tau }{2} +\frac{5}{4}\right)_n \left( \frac{\tau }{2}+1\right)_n} \overline{\eta}^{n } \label{eq:100029c}
\end{align}
\begin{align}
y_0^m(j,q ;z) &= z^{\frac{1}{2}}  \sum_{i_0=0}^{m} \frac{\left( -\frac{j}{2}\right)_{i_0}\left( \frac{3}{4}+\frac{j}{2} \right)_{i_0} }{\left( \frac{5}{4}\right)_{i_0} \left( 1\right)_{i_0}} \eta ^{i_0} \label{eq:100030a}\\
y_1^m(j,q ;z) &= z^{\frac{1}{2}} \left\{\sum_{i_0=0}^{m} \frac{\left( i_0 +\frac{1}{4}\right)^2 -\varpi _j(q) }{\left( i_0+\frac{3}{4} \right) \left( i_0+\frac{1}{2} \right)} \frac{\left( -\frac{j}{2}\right)_{i_0}\left( \frac{3}{4}+\frac{j}{2} \right)_{i_0} }{\left( \frac{5}{4}\right)_{i_0} \left( 1\right)_{i_0}} \right. \nonumber\\ 
 &\times  \left. \sum_{i_1 = i_0}^{m} \frac{\left( \frac{1}{2}-\frac{j}{2} \right)_{i_1} \left( \frac{5}{4}+\frac{j}{2} \right)_{i_1} \left( \frac{7}{4} \right)_{i_0}\left( \frac{3}{2} \right)_{i_0}}{\left( \frac{1}{2}-\frac{j}{2} \right)_{i_0} \left( \frac{5}{4}+\frac{j}{2} \right)_{i_0} \left( \frac{7}{4} \right)_{i_1} \left( \frac{3}{2} \right)_{i_1}} \eta ^{i_1}\right\} \mu \hspace{1cm} 
\label{eq:100030b}\\
y_{\tau }^m(j,q ;z) &= z^{\frac{1}{2}} \left\{ \sum_{i_0=0}^{m} \frac{\left( i_0 +\frac{1}{4}\right) ^2 -\varpi _j(q) }{\left( i_0+\frac{3}{4} \right) \left( i_0+\frac{1}{2} \right)} \frac{\left( -\frac{j}{2}\right)_{i_0}\left( \frac{3}{4}+\frac{j}{2} \right)_{i_0} }{\left( \frac{5}{4}\right)_{i_0} \left( 1\right)_{i_0}} \right.\nonumber\\
 &\times \prod_{k=1}^{\tau -1} \left( \sum_{i_k = i_{k-1}}^{m} \frac{\left( i_k+ \frac{k}{2} +\frac{1}{4}\right) ^2 -\varpi _j(q)}{\left( i_k+\frac{k}{2}+\frac{3}{4} \right) \left( i_k+\frac{k}{2}+\frac{1}{2} \right)} \right.\nonumber\\ 
 &\times \left. \frac{\left( \frac{k}{2}-\frac{j}{2}\right)_{i_k} \left( \frac{k}{2}+\frac{3}{4}+\frac{j}{2} \right)_{i_k}\left( \frac{k}{2} + \frac{5}{4}\right)_{i_{k-1}} \left( \frac{k}{2}+ 1 \right)_{i_{k-1}}}{\left( \frac{k}{2}-\frac{j}{2}\right)_{i_{k-1}} \left( \frac{k}{2}+\frac{3}{4}+\frac{j}{2} \right)_{i_{k-1}}\left( \frac{k}{2} + \frac{5}{4}\right)_{i_k} \left( \frac{k}{2}+ 1\right)_{i_k}} \right) \nonumber \\
 &\times \left. \sum_{i_{\tau } = i_{\tau -1}}^{m} \frac{\left( \frac{\tau }{2} -\frac{j}{2}\right)_{i_{\tau }}\left( \frac{\tau }{2}+\frac{3}{4} +\frac{j}{2} \right)_{i_{\tau }} \left( \frac{\tau }{2} +\frac{5}{4}\right)_{i_{\tau -1}} \left( \frac{\tau }{2}+1\right)_{i_{\tau -1}}}{\left( \frac{\tau }{2} -\frac{j}{2}\right)_{i_{\tau -1}}\left( \frac{\tau }{2}+\frac{3}{4} +\frac{j}{2} \right)_{i_{\tau -1}} \left( \frac{\tau }{2} +\frac{5}{4}\right)_{i_{\tau }} \left( \frac{\tau }{2}+1\right)_{i_{\tau }}} \eta ^{i_{\tau }}\right\} \mu ^{\tau } \hspace{1.5cm} \label{eq:100030c}
\end{align}
where
\begin{equation}
\begin{cases} \tau \geq 2 \cr
\overline{\eta}= \frac{-1}{(a-b)(a-c)}  \cr
\overline{\mu}= \frac{-(2a-b-c)}{(a-b)(a-c)}  \cr
\varpi _j(q) = \frac{ \left( j+\frac{1}{2} \right)\left( j+1 \right) a-\frac{q}{4}}{ 4(2a-b-c)}
\end{cases}\nonumber
\end{equation}
\end{enumerate}
\end{remark}
\subsection{The first species complete polynomial of Lam\'{e} equation using R3TRF} 
\begin{theorem}
In chapter 2, the general summation expression of a function $y(x)$ for the first species complete polynomial using reversible 3-term recurrence formula and its algebraic equation for the determination of an accessory parameter in $A_n$ term are given by
\begin{enumerate} 
\item As $B_1=0$,
\begin{equation}
0 =\bar{c}(0,1) \label{eq:100031a}
\end{equation}
\begin{equation}
y(x) = y_{0}^{0}(x) \label{eq:100031b}
\end{equation}
\item As $B_2=0$, 
\begin{equation}
0 = \bar{c}(0,2)+\bar{c}(1,0) \label{eq:100032a}
\end{equation}
\begin{equation}
y(x)= y_{0}^{1}(x) \label{eq:100032b}
\end{equation}
\item As $B_{2N+3}=0$ where $N \in \mathbb{N}_{0}$,
\begin{equation}
0  = \sum_{r=0}^{N+1}\bar{c}\left( r, 2(N-r)+3\right) \label{eq:100033a}
\end{equation}
\begin{equation}
y(x)= \sum_{r=0}^{N+1} y_{r}^{2(N+1-r)}(x) \label{eq:100033b}
\end{equation}
\item As $B_{2N+4}=0$ where$N \in \mathbb{N}_{0}$,
\begin{equation}
0  =  \sum_{r=0}^{N+2}\bar{c}\left( r, 2(N+2-r)\right) \label{eq:100034a}
\end{equation}
\begin{equation}
y(x)=  \sum_{r=0}^{N+1} y_{r}^{2(N-r)+3}(x) \label{eq:100034b}
\end{equation}
In the above,
\begin{eqnarray}
\bar{c}(0,n) &=& \prod _{i_0=0}^{n-1}A_{i_0} \label{eq:100035a}\\
\bar{c}(1,n) &=& \sum_{i_0=0}^{n} \left\{ B_{i_0+1} \prod _{i_1=0}^{i_0-1}A_{i_1} \prod _{i_2=i_0}^{n-1}A_{i_2+2} \right\} \label{eq:100035b}\\
\bar{c}(\tau ,n) &=& \sum_{i_0=0}^{n} \left\{B_{i_0+1}\prod _{i_1=0}^{i_0-1} A_{i_1} 
\prod _{k=1}^{\tau -1} \left( \sum_{i_{2k}= i_{2(k-1)}}^{n} B_{i_{2k}+(2k+1)}\prod _{i_{2k+1}=i_{2(k-1)}}^{i_{2k}-1}A_{i_{2k+1}+2k}\right) \right. \nonumber\\
&&\times \left. \prod _{i_{2\tau} = i_{2(\tau -1)}}^{n-1} A_{i_{2\tau }+ 2\tau} \right\} 
\hspace{1cm}\label{eq:100035c}
\end{eqnarray}
and
\begin{eqnarray}
y_0^m(x) &=& c_0 x^{\lambda} \sum_{i_0=0}^{m} \left\{ \prod _{i_1=0}^{i_0-1}A_{i_1} \right\} x^{i_0 } \label{eq:100036a}\\
y_1^m(x) &=& c_0 x^{\lambda} \sum_{i_0=0}^{m}\left\{ B_{i_0+1} \prod _{i_1=0}^{i_0-1}A_{i_1}  \sum_{i_2=i_0}^{m} \left\{ \prod _{i_3=i_0}^{i_2-1}A_{i_3+2} \right\}\right\} x^{i_2+2 } \label{eq:100036b}\\
y_{\tau }^m(x) &=& c_0 x^{\lambda} \sum_{i_0=0}^{m} \left\{B_{i_0+1}\prod _{i_1=0}^{i_0-1} A_{i_1} 
\prod _{k=1}^{\tau -1} \left( \sum_{i_{2k}= i_{2(k-1)}}^{m} B_{i_{2k}+(2k+1)}\prod _{i_{2k+1}=i_{2(k-1)}}^{i_{2k}-1}A_{i_{2k+1}+2k}\right) \right. \nonumber\\
&&\times \left. \sum_{i_{2\tau} = i_{2(\tau -1)}}^{m} \left( \prod _{i_{2\tau +1}=i_{2(\tau -1)}}^{i_{2\tau}-1} A_{i_{2\tau +1}+ 2\tau} \right) \right\} x^{i_{2\tau}+2\tau }\hspace{1cm}\mathrm{where}\;\tau \geq 2
\label{eq:100036c}
\end{eqnarray}
\end{enumerate}
\end{theorem} 
According to (\ref{eq:10008}), $c_{j+1}=0$ is clearly an algebraic equation in $q$ of degree $j+1$ and thus has $j+1$ zeros denoted them by $q_j^m$ eigenvalues where $m = 0,1,2, \cdots, j$. They can be arranged in the following order: $q_j^0 < q_j^1 < q_j^2 < \cdots < q_j^j$.
 
Substitute (\ref{eq:100015a}) and (\ref{eq:100015b}) into (\ref{eq:100035a})--(\ref{eq:100036c}). Replace $x$ by $z$ in (\ref{eq:100031b}), (\ref{eq:100032b}), (\ref{eq:100033b}), (\ref{eq:100034b}) and (\ref{eq:100036a})--(\ref{eq:100036c}).

As $B_{1}= c_{1}=0$, take the new (\ref{eq:100035a}) into (\ref{eq:100031a}) putting $j=0$. Substitute the new (\ref{eq:100036a}) into (\ref{eq:100031b}) putting $j=0$.

As $B_{2}= c_{2}=0$, take the new (\ref{eq:100035a}) and (\ref{eq:100035b}) into (\ref{eq:100032a}) putting $j=1$. Substitute the new (\ref{eq:100036a}) into (\ref{eq:100032b}) putting $j=1$ and $q=q_1^m$. 

As $B_{2N+3}= c_{2N+3}=0$, take the new (\ref{eq:100035a})--(\ref{eq:100035c}) into (\ref{eq:100033a}) putting $j=2N+2$. Substitute the new 
(\ref{eq:100036a})--(\ref{eq:100036c}) into (\ref{eq:100033b}) putting $j=2N+2$ and $q =q_{2N+2}^m$.

As $B_{2N+4}= c_{2N+4}=0$, take the new (\ref{eq:100035a})--(\ref{eq:100035c}) into (\ref{eq:100034a}) putting $j=2N+3$. Substitute the new 
(\ref{eq:100036a})--(\ref{eq:100036c}) into (\ref{eq:100034b}) putting $j=2N+3$ and $q =q_{2N+3}^m$.

After the replacement process, the general expression of power series of Lam\'{e} equation about $z=0$ for the first species complete polynomial using reversible 3-term recurrence formula and its algebraic equation for the determination of an accessory parameter $q$ are given by 
\begin{enumerate} 
\item As $\alpha = -2\left( \frac{1}{2}+\lambda \right) $ or $2\lambda $,

An algebraic equation of degree 1 for the determination of $q$ is given by
\begin{equation}
0= \bar{c}(0,1;0,q)= q+4\lambda \left( (a-b-c)\lambda -\frac{a}{2}\right) \label{eq:100037a}
\end{equation}
The eigenvalue of $q$ is written by $q_0^0$. Its eigenfunction is given by
\begin{equation}
y(z) = y_0^0\left( 0,q_0^0;z\right)= c_0 z^{\lambda } \label{eq:100037b}
\end{equation}
\item As $\alpha = -2\left( \frac{3}{2}+\lambda \right) $ or $2\left( 1+\lambda \right) $,

An algebraic equation of degree 2 for the determination of $q$ is given by
\begin{eqnarray}
0 &=& \bar{c}(0,2;1,q)+\bar{c}(1,0;1,q) \nonumber\\
&=& \left( 1+\lambda \right)\left( \frac{1}{2}+\lambda \right)\left( \frac{3}{2}+2\lambda \right)\frac{\left( a-b\right)\left( a-c\right)}{\left( 2a-b-c\right)^2} +\prod_{l=0}^{1} \left( \left( l+\lambda \right)^2 -\frac{\left(  1+\lambda \right)\left( \frac{3}{2}+\lambda  \right) a-\frac{q}{4}}{(2a-b-c)} \right) \hspace{1.5cm}\label{eq:100038a}
\end{eqnarray}
The eigenvalue of $q$ is written by $q_1^m$ where $m = 0,1 $; $q_{1}^0 < q_{1}^1$. Its eigenfunction is given by
\begin{equation}
y(z) = y_0^1\left( 0,q_1^m;z\right)= c_0 z^{\lambda }\left( 1+\frac{\lambda ^2 -\frac{\left( 1+\lambda \right) \left( \frac{3}{2}+\lambda \right) a-\frac{q_1^m}{4}}{(2a-b-c)}}{\left( 1+\lambda \right) \left( \frac{1}{2}+\lambda \right)} \mu \right) \label{eq:100038b}
\end{equation}
\item As $\alpha = -2\left( 2N +\frac{5}{2}+\lambda \right) $ or $2\left( 2N+2+\lambda \right) $ where $N \in \mathbb{N}_{0}$,

An algebraic equation of degree $2N+3$ for the determination of $q$ is given by
\begin{equation}
0 = \sum_{r=0}^{N+1}\bar{c}\left( r, 2(N-r)+3; 2N+2,q \right)  \label{eq:100039a}
\end{equation}
The eigenvalue of $q$ is written by $q_{2N+2}^m$ where $m = 0,1,2,\cdots,2N+2 $; $q_{2N+2}^0 < q_{2N+2}^1 < \cdots < q_{2N+2}^{2N+2}$. Its eigenfunction is given by 
\begin{equation} 
y(z) = \sum_{r=0}^{N+1} y_{r}^{2(N+1-r)}\left( 2N+2, q_{2N+2}^m; z\right)  
\label{eq:100039b} 
\end{equation}
\item As $\alpha = -2\left( 2N +\frac{7}{2}+\lambda \right) $ or $2\left( 2N+3+\lambda \right) $ where $N \in \mathbb{N}_{0}$,

An algebraic equation of degree $2N+4$ for the determination of $q$ is given by
\begin{equation}  
0 = \sum_{r=0}^{N+2}\bar{c}\left( r, 2(N+2-r); 2N+3,q \right) \label{eq:100040a}
\end{equation}
The eigenvalue of $q$ is written by $q_{2N+3}^m$ where $m = 0,1,2,\cdots,2N+3 $; $q_{2N+3}^0 < q_{2N+3}^1 < \cdots < q_{2N+3}^{2N+3}$. Its eigenfunction is given by
\begin{equation} 
y(z) = \sum_{r=0}^{N+1} y_{r}^{2(N-r)+3} \left( 2N+3, q_{2N+3}^m; z\right) \label{eq:100040b}
\end{equation}
In the above,
\begin{align}
\bar{c}(0,n;j,q)  &=  \frac{\left( -\varphi _j(q) +\lambda \right)_{n}\left( \varphi _j(q) +\lambda \right)_{n}}{\left( 1+ \lambda \right)_{n} \left( \frac{1}{2} +\lambda \right)_{n}} \overline{\mu}^{n}\label{eq:100041a}\\
\bar{c}(1,n;j,q) &=  \overline{\eta} \sum_{i_0=0}^{n}\frac{\left( i_0-j\right) \left( i_0+\frac{1}{2}+j+2\lambda \right)}{\left( i_0+ 2 + \lambda \right) \left( i_0+ \frac{3}{2} + \lambda \right)} \frac{\left( -\varphi _j(q) +\lambda \right)_{i_0}\left( \varphi _j(q) +\lambda \right)_{i_0}}{\left( 1+ \lambda \right)_{i_0} \left( \frac{1}{2} +\lambda \right)_{i_0}} \nonumber\\
 &\times \frac{\left(  2 -\varphi _j(q)+\lambda \right)_{n}\left(  2 +\varphi _j(q)+\lambda \right)_{n} \left( 3 + \lambda \right)_{i_0} \left( \frac{5}{2} + \lambda \right)_{i_0}}{\left(  2 -\varphi _j(q)+\lambda \right)_{i_0}\left(  2 +\varphi _j(q)+\lambda \right)_{i_0} \left( 3 + \lambda \right)_{n} \left( \frac{5}{2}+\lambda \right)_{n}}\overline{\mu}^n  
 \label{eq:100041b}\\
\bar{c}(\tau ,n;j,q) &=  \overline{\eta}^{\tau } \sum_{i_0=0}^{n} \frac{\left( i_0-j\right) \left( i_0+\frac{1}{2}+j+2\lambda \right)}{\left( i_0+ 2 + \lambda \right) \left( i_0+ \frac{3}{2} + \lambda \right)} \frac{\left( -\varphi _j(q) +\lambda \right)_{i_0}\left( \varphi _j(q) +\lambda \right)_{i_0}}{\left( 1+ \lambda \right)_{i_0} \left( \frac{1}{2} +\lambda \right)_{i_0}}  \nonumber\\
 &\times \prod_{k=1}^{\tau -1} \left( \sum_{i_k = i_{k-1}}^{n} \frac{\left( i_k+ 2k-j\right)\left( i_k+ 2k+\frac{1}{2}+j+2\lambda \right) }{\left( i_k+2k+2+ \lambda \right) \left( i_k+2k+\frac{3}{2} + \lambda \right)} \right. \nonumber\\
 &\times \left. \frac{\left( 2k-\varphi _j(q) +\lambda \right)_{i_k} \left( 2k+\varphi _j(q) +\lambda \right)_{i_k}\left( 2k+1+ \lambda \right)_{i_{k-1}} \left( 2k+\frac{1}{2} + \lambda \right)_{i_{k-1}}}{\left( 2k-\varphi _j(q) +\lambda \right)_{i_{k-1}} \left( 2k+\varphi _j(q) +\lambda \right)_{i_{k-1}}\left( 2k+1+ \lambda \right)_{i_k} \left( 2k+\frac{1}{2} + \lambda \right)_{i_k}} \right) \nonumber\\ 
 &\times \frac{\left( 2\tau -\varphi _j(q) +\lambda \right)_n \left( 2\tau +\varphi _j(q) +\lambda \right)_n \left( 2\tau +1+ \lambda \right)_{i_{\tau -1}} \left( 2\tau +\frac{1}{2} + \lambda \right)_{i_{\tau -1}}}{\left( 2\tau -\varphi _j(q) +\lambda \right)_{i_{\tau -1}} \left( 2\tau +\varphi _j(q) +\lambda \right)_{i_{\tau -1}}\left( 2\tau +1+ \lambda \right)_n \left( 2\tau +\frac{1}{2} + \lambda \right)_n} \overline{\mu}^{n }  \label{eq:100041c}
\end{align}
\begin{align}
y_0^m(j,q ;z) &=  c_0 z^{\lambda }  \sum_{i_0=0}^{m} \frac{\left( -\varphi _j(q) +\lambda \right)_{i_0}\left( \varphi _j(q) +\lambda \right)_{i_0}}{\left( 1+ \lambda \right)_{i_0} \left( \frac{1}{2} +\lambda \right)_{i_0}} \mu ^{i_0} \label{eq:100042a}\\
y_1^m(j,q ;z) &=  c_0 z^{\lambda } \left\{\sum_{i_0=0}^{m} \frac{\left( i_0-j\right) \left( i_0+\frac{1}{2}+j+2\lambda \right)}{\left( i_0+ 2 + \lambda \right) \left( i_0+ \frac{3}{2} + \lambda \right)} \frac{\left( -\varphi _j(q) +\lambda \right)_{i_0}\left( \varphi _j(q) +\lambda \right)_{i_0}}{\left( 1+ \lambda \right)_{i_0} \left( \frac{1}{2} +\lambda \right)_{i_0}} \right. \nonumber\\
 &\times \left. \sum_{i_1 = i_0}^{m} \frac{\left( 2 -\varphi _j(q)+\lambda \right)_{i_1}\left(  2 +\varphi _j(q)+\lambda \right)_{i_1} \left( 3 + \lambda \right)_{i_0} \left( \frac{5}{2} + \lambda \right)_{i_0}}{\left( 2 -\varphi _j(q)+\lambda \right)_{i_0}\left( 2 +\varphi _j(q)+\lambda \right)_{i_0} \left( 3 + \lambda \right)_{i_1} \left( \frac{5}{2}+\lambda \right)_{i_1}} \mu ^{i_1}\right\} \eta 
\label{eq:100042b}
\end{align}
\begin{align}
  y_{\tau }^m(j,q ;z) &= c_0 z^{\lambda } \left\{ \sum_{i_0=0}^{m} \frac{\left( i_0-j\right) \left( i_0+\frac{1}{2}+j+2\lambda \right)}{\left( i_0+ 2 + \lambda \right) \left( i_0+ \frac{3}{2} + \lambda \right)} \frac{\left( -\varphi _j(q) +\lambda \right)_{i_0}\left( \varphi _j(q) +\lambda \right)_{i_0}}{\left( 1+ \lambda \right)_{i_0} \left( \frac{1}{2} +\lambda \right)_{i_0}} \right.\nonumber\\
 &\times \prod_{k=1}^{\tau -1} \left( \sum_{i_k = i_{k-1}}^{m} \frac{\left( i_k+ 2k-j\right)\left( i_k+ 2k+\frac{1}{2}+j+2\lambda \right) }{\left( i_k+2k+2+ \lambda \right) \left( i_k+2k+\frac{3}{2} + \lambda \right)} \right. \nonumber\\
 &\times \left. \frac{\left( 2k-\varphi _j(q) +\lambda \right)_{i_k} \left( 2k+\varphi _j(q) +\lambda \right)_{i_k}\left( 2k+1+ \lambda \right)_{i_{k-1}} \left( 2k+\frac{1}{2} + \lambda \right)_{i_{k-1}}}{\left( 2k-\varphi _j(q) +\lambda \right)_{i_{k-1}} \left( 2k+\varphi _j(q) +\lambda \right)_{i_{k-1}}\left( 2k+1+ \lambda \right)_{i_k} \left( 2k+\frac{1}{2} + \lambda \right)_{i_k}} \right)  \nonumber\\
 &\times \left. \sum_{i_{\tau } = i_{\tau -1}}^{m}  \frac{\left( 2\tau -\varphi _j(q) +\lambda \right)_{i_{\tau }} \left( 2\tau +\varphi _j(q) +\lambda \right)_{i_{\tau }} \left( 2\tau +1+ \lambda \right)_{i_{\tau -1}} \left( 2\tau +\frac{1}{2} + \lambda \right)_{i_{\tau -1}}}{\left( 2\tau -\varphi _j(q) +\lambda \right)_{i_{\tau -1}} \left( 2\tau +\varphi _j(q) +\lambda \right)_{i_{\tau -1}}\left( 2\tau +1+ \lambda \right)_{i_{\tau }} \left( 2\tau +\frac{1}{2} + \lambda \right)_{i_{\tau }}} \mu ^{i_{\tau }}\right\} \eta ^{\tau }  \label{eq:100042c}
\end{align}
where
\begin{equation}
\begin{cases} \tau \geq 2 \cr
\eta = \frac{-z^2}{(a-b)(a-c)}  \cr
\mu =  \frac{-(2a-b-c)}{(a-b)(a-c)} z \cr
\overline{\eta}= \frac{-1}{(a-b)(a-c)}  \cr
\overline{\mu}= \frac{-(2a-b-c)}{(a-b)(a-c)}  \cr
\varphi _j(q) = \sqrt{ \frac{ \left( j+\lambda \right)\left( j+\frac{1}{2}+\lambda\right) a-\frac{q}{4}}{ (2a-b-c)}}
\end{cases}\nonumber
\end{equation}
\end{enumerate} 
Put $c_0$= 1 as $\lambda =0$ for the first kind of independent solutions of Lam\'{e} equation and $\lambda =\frac{1}{2} $ for the second one in (\ref{eq:100037a})--(\ref{eq:100042c}). 
\begin{remark}
The power series expansion of Lam\'{e} equation in the algebraic form of the first kind for the first species complete polynomial using R3TRF about $z=0$ is given by
\begin{enumerate} 
\item As $\alpha = -1 $ or 0 with $q=q_0^0=0 $,
 
The eigenfunction is given by
\begin{eqnarray}
y(z) &=& L_p^{(a)}F_{0,0}^R\left( a,b,c,q=q_0^0=0, \alpha = -1\; \mbox{or}\; 0; z=x-a,\eta = \frac{-z^2}{(a-b)(a-c)}\right. \nonumber\\
&&,\left. \mu = \frac{-(2a-b-c)}{(a-b)(a-c)} z \right) \nonumber\\
&=& 1 \label{eq:100043}
\end{eqnarray}
\item As $\alpha = -3 $ or $2 $,

An algebraic equation of degree 2 for the determination of $q$ is given by
\begin{equation}
0 =\frac{3\left( a-b\right)\left( a-c\right)}{4\left( 2a-b-c\right)^2} +\prod_{l=0}^{1} \left( l^2 -\frac{ 6a-q }{4(2a-b-c)} \right) \label{eq:100044a}
\end{equation}
The eigenvalue of $q$ is written by $q_1^m$ where $m = 0,1 $; $q_{1}^0 < q_{1}^1$. Its eigenfunction is given by
\begin{eqnarray}
y(z) &=& L_p^{(a)}F_{1,m}^R\left( a,b,c,q=q_1^m, \alpha = -3\; \mbox{or}\; 2; z=x-a,\eta = \frac{-z^2}{(a-b)(a-c)}\right.\nonumber\\
&&,\left. \mu = \frac{-(2a-b-c)}{(a-b)(a-c)} z \right) \nonumber\\
&=& 1- \frac{ 6a- q_1^m }{2(2a-b-c)} \mu \label{eq:100044b}
\end{eqnarray}
\item As $\alpha = -2\left( 2N +\frac{5}{2} \right) $ or $2\left( 2N+2 \right) $ where $N \in \mathbb{N}_{0}$,

An algebraic equation of degree $2N+3$ for the determination of $q$ is given by
\begin{equation}
0 = \sum_{r=0}^{N+1}\bar{c}\left( r, 2(N-r)+3; 2N+2,q \right)  \label{eq:100045a}
\end{equation}
The eigenvalue of $q$ is written by $q_{2N+2}^m$ where $m = 0,1,2,\cdots,2N+2 $; $q_{2N+2}^0 < q_{2N+2}^1 < \cdots < q_{2N+2}^{2N+2}$. Its eigenfunction is given by 
\begin{eqnarray} 
y(z) &=& L_p^{(a)}F_{2N+2,m}^R\left( a,b,c,q=q_{2N+2}^m, \alpha = -2\left( 2N+\frac{5}{2} \right)\; \mbox{or}\; 2\left( 2N+2 \right); z=x-a \right. \nonumber\\
&&,\left. \eta = \frac{-z^2}{(a-b)(a-c)},\mu = \frac{-(2a-b-c)}{(a-b)(a-c)} z \right) \nonumber\\
&=& \sum_{r=0}^{N+1} y_{r}^{2(N+1-r)}\left( 2N+2, q_{2N+2}^m; z\right)  
\label{eq:100045b} 
\end{eqnarray}
\item As $\alpha = -2\left( 2N +\frac{7}{2} \right) $ or $2\left( 2N+3 \right) $ where $N \in \mathbb{N}_{0}$,

An algebraic equation of degree $2N+4$ for the determination of $q$ is given by
\begin{equation}  
0 = \sum_{r=0}^{N+2}\bar{c}\left( r, 2(N+2-r); 2N+3,q \right) \label{eq:100046a}
\end{equation}
The eigenvalue of $q$ is written by $q_{2N+3}^m$ where $m = 0,1,2,\cdots,2N+3 $; $q_{2N+3}^0 < q_{2N+3}^1 < \cdots < q_{2N+3}^{2N+3}$. Its eigenfunction is given by
\begin{eqnarray} 
y(z) &=& L_p^{(a)}F_{2N+3,m}^R\left( a,b,c,q=q_{2N+3}^m, \alpha = -2\left( 2N+\frac{7}{2} \right)\; \mbox{or}\; 2\left( 2N+3 \right); z=x-a \right. \nonumber\\
&&,\left. \eta = \frac{-z^2}{(a-b)(a-c)},\mu = \frac{-(2a-b-c)}{(a-b)(a-c)} z \right) \nonumber\\
&=& \sum_{r=0}^{N+1} y_{r}^{2(N-r)+3} \left( 2N+3, q_{2N+3}^m; z\right) \label{eq:100046b}
\end{eqnarray}
In the above,
\begin{eqnarray}
\bar{c}(0,n;j,q)  &=& \frac{\left( -\varphi _j(q) \right)_{n}\left( \varphi _j(q) \right)_{n}}{\left( 1\right)_{n} \left( \frac{1}{2} \right)_{n}} \overline{\mu}^{n}\label{eq:100047a}\\
\bar{c}(1,n;j,q) &=& \overline{\eta} \sum_{i_0=0}^{n}\frac{\left( i_0-j\right) \left( i_0+\frac{1}{2}+j\right)}{\left( i_0+2\right) \left( i_0+ \frac{3}{2} \right)} \frac{\left( -\varphi _j(q) \right)_{i_0}\left( \varphi _j(q) \right)_{i_0}}{\left( 1 \right)_{i_0} \left( \frac{1}{2} \right)_{i_0}}  \nonumber\\
&&\times  \frac{\left(  2 -\varphi _j(q) \right)_{n}\left(  2 +\varphi _j(q) \right)_{n} \left( 3 \right)_{i_0} \left( \frac{5}{2} \right)_{i_0}}{\left(  2 -\varphi _j(q) \right)_{i_0}\left(  2 +\varphi _j(q) \right)_{i_0} \left( 3 \right)_{n} \left( \frac{5}{2} \right)_{n}}\overline{\mu}^n  
 \label{eq:100047b}\\
\bar{c}(\tau ,n;j,q) &=& \overline{\eta}^{\tau } \sum_{i_0=0}^{n} \frac{\left( i_0-j\right) \left( i_0+\frac{1}{2}+j \right)}{\left( i_0+2 \right) \left( i_0+ \frac{3}{2} \right)} \frac{\left( -\varphi _j(q) \right)_{i_0}\left( \varphi _j(q) \right)_{i_0}}{\left( 1\right)_{i_0} \left( \frac{1}{2} \right)_{i_0}}  \nonumber\\
&&\times \prod_{k=1}^{\tau -1} \left( \sum_{i_k = i_{k-1}}^{n} \frac{\left( i_k+ 2k-j\right)\left( i_k+ 2k+\frac{1}{2}+j \right) }{\left( i_k+2k+2 \right) \left( i_k+2k+\frac{3}{2} \right)} \right. \nonumber\\
&&\times \left. \frac{\left( 2k-\varphi _j(q) \right)_{i_k} \left( 2k+\varphi _j(q) \right)_{i_k}\left( 2k+1 \right)_{i_{k-1}} \left( 2k+\frac{1}{2} \right)_{i_{k-1}}}{\left( 2k-\varphi _j(q) \right)_{i_{k-1}} \left( 2k+\varphi _j(q) \right)_{i_{k-1}}\left( 2k+1 \right)_{i_k} \left( 2k+\frac{1}{2} \right)_{i_k}} \right) \nonumber\\ 
&&\times \frac{\left( 2\tau -\varphi _j(q) \right)_n \left( 2\tau +\varphi _j(q) \right)_n \left( 2\tau +1 \right)_{i_{\tau -1}} \left( 2\tau +\frac{1}{2} \right)_{i_{\tau -1}}}{\left( 2\tau -\varphi _j(q) \right)_{i_{\tau -1}} \left( 2\tau +\varphi _j(q) \right)_{i_{\tau -1}}\left( 2\tau +1 \right)_n \left( 2\tau +\frac{1}{2}  \right)_n} \overline{\mu}^{n } \hspace{1.5cm}\label{eq:100047c}
\end{eqnarray}
\begin{eqnarray}
y_0^m(j,q ;z) &=& \sum_{i_0=0}^{m} \frac{\left( -\varphi _j(q) \right)_{i_0}\left( \varphi _j(q) \right)_{i_0}}{\left( 1 \right)_{i_0} \left( \frac{1}{2} \right)_{i_0}} \mu ^{i_0} \label{eq:100048a}\\
y_1^m(j,q ;z) &=&  \left\{\sum_{i_0=0}^{m} \frac{\left( i_0-j\right) \left( i_0+\frac{1}{2}+j \right)}{\left( i_0+ 2 \right) \left( i_0+ \frac{3}{2} \right)} \frac{\left( -\varphi _j(q) \right)_{i_0}\left( \varphi _j(q)  \right)_{i_0}}{\left( 1 \right)_{i_0} \left( \frac{1}{2}  \right)_{i_0}} \right. \nonumber\\
&&\times \left. \sum_{i_1 = i_0}^{m} \frac{\left( 2 -\varphi _j(q) \right)_{i_1}\left(  2 +\varphi _j(q)  \right)_{i_1} \left( 3 \right)_{i_0} \left( \frac{5}{2} \right)_{i_0}}{\left( 2 -\varphi _j(q) \right)_{i_0}\left( 2 +\varphi _j(q) \right)_{i_0} \left( 3 \right)_{i_1} \left( \frac{5}{2} \right)_{i_1}} \mu ^{i_1}\right\} \eta 
\label{eq:100048b}\\
y_{\tau }^m(j,q ;z) &=& \left\{ \sum_{i_0=0}^{m} \frac{\left( i_0-j\right) \left( i_0+\frac{1}{2}+j \right)}{\left( i_0+ 2 \right) \left( i_0+ \frac{3}{2} \right)} \frac{\left( -\varphi _j(q) \right)_{i_0}\left( \varphi _j(q) \right)_{i_0}}{\left( 1 \right)_{i_0} \left( \frac{1}{2} \right)_{i_0}} \right.\nonumber\\
&&\times \prod_{k=1}^{\tau -1} \left( \sum_{i_k = i_{k-1}}^{m} \frac{\left( i_k+ 2k-j\right)\left( i_k+ 2k+\frac{1}{2}+j \right) }{\left( i_k+2k+2 \right) \left( i_k+2k+\frac{3}{2} \right)} \right. \nonumber\\
&&\times \left. \frac{\left( 2k-\varphi _j(q) \right)_{i_k} \left( 2k+\varphi _j(q) \right)_{i_k}\left( 2k+1 \right)_{i_{k-1}} \left( 2k+\frac{1}{2} \right)_{i_{k-1}}}{\left( 2k-\varphi _j(q) \right)_{i_{k-1}} \left( 2k+\varphi _j(q) \right)_{i_{k-1}}\left( 2k+1 \right)_{i_k} \left( 2k+\frac{1}{2} \right)_{i_k}} \right) \label{eq:100048c}\\
&&\times \left. \sum_{i_{\tau } = i_{\tau -1}}^{m}  \frac{\left( 2\tau -\varphi _j(q) \right)_{i_{\tau }} \left( 2\tau +\varphi _j(q) \right)_{i_{\tau }} \left( 2\tau +1 \right)_{i_{\tau -1}} \left( 2\tau +\frac{1}{2} \right)_{i_{\tau -1}}}{\left( 2\tau -\varphi _j(q) \right)_{i_{\tau -1}} \left( 2\tau +\varphi _j(q) \right)_{i_{\tau -1}}\left( 2\tau +1 \right)_{i_{\tau }} \left( 2\tau +\frac{1}{2} \right)_{i_{\tau }}} \mu ^{i_{\tau }}\right\} \eta ^{\tau }   \nonumber
\end{eqnarray}
where
\begin{equation}
\begin{cases} \tau \geq 2 \cr
\overline{\eta}= \frac{-1}{(a-b)(a-c)}  \cr
\overline{\mu}= \frac{-(2a-b-c)}{(a-b)(a-c)}  \cr
\varphi _j(q) = \sqrt{\frac{ j \left( j+\frac{1}{2} \right) a-\frac{q}{4}}{ (2a-b-c)}}
\end{cases}\nonumber
\end{equation}
\end{enumerate}
\end{remark} 
\begin{remark}
The power series expansion of Lam\'{e} equation in the algebraic form of the second kind for the first species complete polynomial using R3TRF about $z=0$ is given by
\begin{enumerate} 
\item As $\alpha = -2 $ or 1 with $q=q_0^0=b+c $,
 
The eigenfunction is given by
\begin{eqnarray}
y(z) &=& L_p^{(a)}S_{0,0}^R\left( a,b,c,q=q_0^0=b+c, \alpha = -2\; \mbox{or}\; 1; z=x-a,\eta = \frac{-z^2}{(a-b)(a-c)}\right. \nonumber\\
&&,\left.\mu = \frac{-(2a-b-c)}{(a-b)(a-c)} z \right) \nonumber\\
&=& z^{\frac{1}{2}} \label{eq:100049}
\end{eqnarray}
\item As $\alpha = -4 $ or $3 $,

An algebraic equation of degree 2 for the determination of $q$ is given by
\begin{equation}
0 = \frac{15\left( a-b\right)\left( a-c\right)}{4\left( 2a-b-c\right)^2} +\prod_{l=0}^{1} \left( \left( l+\frac{1}{2} \right)^2 -\frac{12a-q}{4(2a-b-c)} \right)  \label{eq:100050a}
\end{equation}
The eigenvalue of $q$ is written by $q_1^m$ where $m = 0,1 $; $q_{1}^0 < q_{1}^1$. Its eigenfunction is given by
\begin{eqnarray}
y(z) &=& L_p^{(a)}S_{1,m}^R\left( a,b,c,q=q_1^m, \alpha = -4\; \mbox{or}\; 3; z=x-a,\eta = \frac{-z^2}{(a-b)(a-c)}\right. \nonumber\\
&&,\left. \mu = \frac{-(2a-b-c)}{(a-b)(a-c)} z \right) \nonumber\\
&=& z^{\frac{1}{2}}\left\{ 1+ \frac{1}{6} \left( 1 -\frac{12a- q_1^m }{(2a-b-c)} \right) \mu \right\} \label{eq:100050b}
\end{eqnarray}
\item As $\alpha = -2\left( 2N +3 \right) $ or $2\left( 2N+\frac{5}{2} \right) $ where $N \in \mathbb{N}_{0}$,

An algebraic equation of degree $2N+3$ for the determination of $q$ is given by
\begin{equation}
0 = \sum_{r=0}^{N+1}\bar{c}\left( r, 2(N-r)+3; 2N+2,q \right)  \label{eq:100051a}
\end{equation}
The eigenvalue of $q$ is written by $q_{2N+2}^m$ where $m = 0,1,2,\cdots,2N+2 $; $q_{2N+2}^0 < q_{2N+2}^1 < \cdots < q_{2N+2}^{2N+2}$. Its eigenfunction is given by 
\begin{eqnarray} 
y(z) &=& L_p^{(a)}S_{2N+2,m}^R\left( a,b,c,q=q_{2N+2}^m, \alpha = -2\left( 2N+3 \right)\; \mbox{or}\; 2\left( 2N+\frac{5}{2} \right); z=x-a \right. \nonumber\\
&&,\left. \eta = \frac{-z^2}{(a-b)(a-c)},\mu = \frac{-(2a-b-c)}{(a-b)(a-c)} z \right) \nonumber\\
&=& \sum_{r=0}^{N+1} y_{r}^{2(N+1-r)}\left( 2N+2, q_{2N+2}^m; z\right)  
\label{eq:100051b} 
\end{eqnarray}
\item As $\alpha = -2\left( 2N +4 \right) $ or $2\left( 2N+\frac{7}{2} \right) $ where $N \in \mathbb{N}_{0}$,

An algebraic equation of degree $2N+4$ for the determination of $q$ is given by
\begin{equation}  
0 = \sum_{r=0}^{N+2}\bar{c}\left( r, 2(N+2-r); 2N+3,q \right) \label{eq:100052a}
\end{equation}
The eigenvalue of $q$ is written by $q_{2N+3}^m$ where $m = 0,1,2,\cdots,2N+3 $; $q_{2N+3}^0 < q_{2N+3}^1 < \cdots < q_{2N+3}^{2N+3}$. Its eigenfunction is given by
\begin{eqnarray} 
y(z) &=& L_p^{(a)}S_{2N+3,m}^R\left( a,b,c,q=q_{2N+3}^m, \alpha = -2\left( 2N+4 \right)\; \mbox{or}\; 2\left( 2N+\frac{7}{2} \right); z=x-a \right. \nonumber\\
&&,\left. \eta = \frac{-z^2}{(a-b)(a-c)},\mu = \frac{-(2a-b-c)}{(a-b)(a-c)} z \right) \nonumber\\
&=& \sum_{r=0}^{N+1} y_{r}^{2(N-r)+3} \left( 2N+3, q_{2N+3}^m; z\right) \label{eq:100052b}
\end{eqnarray}
In the above,
\begin{eqnarray}
\bar{c}(0,n;j,q)  &=& \frac{\left( \frac{1}{2}-\varphi _j(q) \right)_{n}\left( \frac{1}{2}+\varphi _j(q) \right)_{n}}{\left( \frac{3}{2}\right)_{n} \left( 1 \right)_{n}} \overline{\mu}^{n}\label{eq:100053a}\\
\bar{c}(1,n;j,q) &=& \overline{\eta} \sum_{i_0=0}^{n}\frac{\left( i_0-j\right) \left( i_0+\frac{3}{2}+j\right)}{\left( i_0+ \frac{5}{2} \right) \left( i_0+ 2 \right)} \frac{\left( \frac{1}{2}-\varphi _j(q) \right)_{i_0}\left( \frac{1}{2}+\varphi _j(q) \right)_{i_0}}{\left( \frac{3}{2} \right)_{i_0} \left( 1\right)_{i_0}} \nonumber\\
&&\times \frac{\left(  \frac{5}{2}-\varphi _j(q) \right)_{n}\left(  \frac{5}{2}+\varphi _j(q) \right)_{n} \left( \frac{7}{2} \right)_{i_0} \left( 3 \right)_{i_0}}{\left(  \frac{5}{2}-\varphi _j(q)  \right)_{i_0}\left( \frac{5}{2}+\varphi _j(q) \right)_{i_0} \left( \frac{7}{2} \right)_{n} \left( 3\right)_{n}}\overline{\mu}^n  
 \label{eq:100053b}\\
\bar{c}(\tau ,n;j,q) &=& \overline{\eta}^{\tau } \sum_{i_0=0}^{n} \frac{\left( i_0-j\right) \left( i_0+\frac{3}{2}+j  \right)}{\left( i_0+ \frac{5}{2} \right) \left( i_0+ 2\right)} \frac{\left( \frac{1}{2}-\varphi _j(q) \right)_{i_0}\left( \frac{1}{2}+\varphi _j(q) \right)_{i_0}}{\left( \frac{3}{2} \right)_{i_0} \left( 1 \right)_{i_0}}  \nonumber\\
&&\times \prod_{k=1}^{\tau -1} \left( \sum_{i_k = i_{k-1}}^{n} \frac{\left( i_k+ 2k-j\right)\left( i_k+ 2k+\frac{3}{2}+j \right) }{\left( i_k+2k+\frac{5}{2} \right) \left( i_k+2k+2 \right)} \right. \nonumber\\
&&\times \left. \frac{\left( 2k+\frac{1}{2}-\varphi _j(q)  \right)_{i_k} \left( 2k+\frac{1}{2}+\varphi _j(q)  \right)_{i_k}\left( 2k+\frac{3}{2} \right)_{i_{k-1}} \left( 2k+1 \right)_{i_{k-1}}}{\left( 2k+\frac{1}{2}-\varphi _j(q) \right)_{i_{k-1}} \left( 2k+\frac{1}{2}+\varphi _j(q) \right)_{i_{k-1}}\left( 2k+\frac{3}{2} \right)_{i_k} \left( 2k+1 \right)_{i_k}} \right) \nonumber\\ 
&&\times \frac{\left( 2\tau +\frac{1}{2}-\varphi _j(q) \right)_n \left( 2\tau +\frac{1}{2}+\varphi _j(q) \right)_n \left( 2\tau +\frac{3}{2} \right)_{i_{\tau -1}} \left( 2\tau +1 \right)_{i_{\tau -1}}}{\left( 2\tau +\frac{1}{2}-\varphi _j(q) \right)_{i_{\tau -1}} \left( 2\tau +\frac{1}{2}+\varphi _j(q) \right)_{i_{\tau -1}}\left( 2\tau +\frac{3}{2} \right)_n \left( 2\tau +1 \right)_n} \overline{\mu}^{n } \hspace{2cm}\label{eq:100053c}
\end{eqnarray}
\begin{eqnarray}
y_0^m(j,q ;z) &=& z^{\frac{1}{2}}  \sum_{i_0=0}^{m} \frac{\left( \frac{1}{2}-\varphi _j(q) \right)_{i_0}\left( \frac{1}{2}+ \varphi _j(q) \right)_{i_0}}{\left( \frac{3}{2} \right)_{i_0} \left( 1\right)_{i_0}} \mu ^{i_0} \label{eq:100054a}\\
y_1^m(j,q ;z) &=& z^{\frac{1}{2}} \left\{\sum_{i_0=0}^{m} \frac{\left( i_0-j\right) \left( i_0+\frac{3}{2}+j \right)}{\left( i_0+ \frac{5}{2} \right) \left( i_0+ 2 \right)} \frac{\left( \frac{1}{2}-\varphi _j(q) \right)_{i_0}\left( \frac{1}{2}+ \varphi _j(q) \right)_{i_0}}{\left( \frac{3}{2} \right)_{i_0} \left( 1 \right)_{i_0}} \right. \nonumber\\
&&\times \left. \sum_{i_1 = i_0}^{m} \frac{\left( \frac{5}{2} -\varphi _j(q)  \right)_{i_1}\left( \frac{5}{2} +\varphi _j(q) \right)_{i_1} \left( \frac{7}{2} \right)_{i_0} \left( 3 \right)_{i_0}}{\left( \frac{5}{2} -\varphi _j(q) \right)_{i_0}\left( \frac{5}{2}+\varphi _j(q) \right)_{i_0} \left( \frac{7}{2} \right)_{i_1} \left( 3 \right)_{i_1}} \mu ^{i_1}\right\} \eta 
\label{eq:100054b}\\
y_{\tau }^m(j,q ;z) &=& z^{\frac{1}{2}} \left\{ \sum_{i_0=0}^{m} \frac{\left( i_0-j\right) \left( i_0+\frac{3}{2}+j \right)}{\left( i_0+ \frac{5}{2} \right) \left( i_0+ 2 \right)} \frac{\left( \frac{1}{2}-\varphi _j(q) \right)_{i_0}\left( \frac{1}{2}+ \varphi _j(q)  \right)_{i_0}}{\left( \frac{3}{2} \right)_{i_0} \left( 1 \right)_{i_0}} \right.\nonumber\\
&&\times \prod_{k=1}^{\tau -1} \left( \sum_{i_k = i_{k-1}}^{m} \frac{\left( i_k+ 2k-j\right)\left( i_k+ 2k+\frac{3}{2}+j \right) }{\left( i_k+2k+\frac{5}{2} \right) \left( i_k+2k+2 \right)} \right. \nonumber\\
&&\times \left. \frac{\left( 2k+\frac{1}{2}-\varphi _j(q) \right)_{i_k} \left( 2k+\frac{1}{2}+\varphi _j(q) \right)_{i_k}\left( 2k+\frac{3}{2} \right)_{i_{k-1}} \left( 2k+1 \right)_{i_{k-1}}}{\left( 2k+\frac{1}{2}-\varphi _j(q) \right)_{i_{k-1}} \left( 2k+\frac{1}{2}+\varphi _j(q) \right)_{i_{k-1}}\left( 2k+\frac{3}{2} \right)_{i_k} \left( 2k+1 \right)_{i_k}} \right) \label{eq:100054c}\\
&&\times \left. \sum_{i_{\tau } = i_{\tau -1}}^{m}  \frac{\left( 2\tau +\frac{1}{2}-\varphi _j(q)  \right)_{i_{\tau }} \left( 2\tau +\frac{1}{2}+\varphi _j(q)  \right)_{i_{\tau }} \left( 2\tau +\frac{3}{2} \right)_{i_{\tau -1}} \left( 2\tau +1 \right)_{i_{\tau -1}}}{\left( 2\tau +\frac{1}{2}-\varphi _j(q)  \right)_{i_{\tau -1}} \left( 2\tau +\frac{1}{2}+\varphi _j(q)  \right)_{i_{\tau -1}}\left( 2\tau +\frac{3}{2} \right)_{i_{\tau }} \left( 2\tau +1 \right)_{i_{\tau }}} \mu ^{i_{\tau }}\right\} \eta ^{\tau }   \nonumber
\end{eqnarray}
where
\begin{equation}
\begin{cases} \tau \geq 2 \cr
\overline{\eta}= \frac{-1}{(a-b)(a-c)}  \cr
\overline{\mu}= \frac{-(2a-b-c)}{(a-b)(a-c)}  \cr
\varphi _j(q) = \sqrt{ \frac{ \left( j+\frac{1}{2}\right)\left( j+1\right) a-\frac{q}{4}}{ (2a-b-c)}}
\end{cases}\nonumber
\end{equation}
\end{enumerate} 
\end{remark} 
\section{Lam\'{e} equation in Weierstrass's form}
Weierstrass's form of Lam\'{e} equation is a second-order linear ordinary differential equation (ODE) of the form
\begin{equation}
\frac{d^2{y}}{d{z}^2} = \{ \alpha (\alpha +1)\rho^2\;sn^2(z,\rho )-h\} y(z)\label{eq:100055}
\end{equation}
where the Jacobian elliptic function $sn^2(z,\rho )$ has modulus $\rho$. And, in classical treatment \cite{10Whit1927}, $\rho$, $\alpha $  and $h$ are real parameters such that $0<\rho <1$ and $\alpha \geq -\frac{1}{2}$.
The Jacobian form has been investigated by Hermite with Floquet's theorem, and he established the solutions of Lam\'{e} equation are doubly periodic of the second kind for every value of $h$ \cite{10Herm1877}.   

Let $sn^2(z,\rho)=\xi $ as an independent variable in (\ref{eq:100055}), Lam\'{e} equation is written by
\begin{equation}
\frac{d^2{y}}{d{\xi }^2} + \frac{1}{2}\left(\frac{1}{\xi } +\frac{1}{\xi -1} + \frac{1}{\xi -\rho ^{-2}}\right) \frac{d{y}}{d{\xi }} +  \frac{-\alpha (\alpha +1) \xi +h\rho ^{-2}}{4 \xi (\xi -1)(\xi -\rho ^{-2})} y(\xi ) = 0\label{eq:100056}
\end{equation}
This is an equation of Fuchsian type with the four regular singularities: $\xi=0, 1, \rho ^{-2}, \infty $. The first three, namely $0, 1, \rho ^{-2}$, have the property that the corresponding exponents are $0, \frac{1}{2}$ which is equivalent to the case of Lam\'{e} equation in the algebraic form.

As we compare (\ref{eq:100056}) with (\ref{eq:10001}), all coefficients on the above are correspondent to the following way.
\begin{equation}
\begin{split}
& a \longrightarrow   0 \\ & b \longrightarrow  1 \\ & c \longrightarrow  \rho ^{-2} \\
& q \longrightarrow  h \rho ^{-2} \\ & x \longrightarrow \xi = sn^2(z,\rho ) 
\end{split}\label{eq:100057}   
\end{equation}
\subsection{The first species complete polynomial of Lam\'{e} equation using 3TRF}
Substitute (\ref{eq:100057}) into (\ref{eq:100016a})--(\ref{eq:100020c}) with replacing a spectral parameter $q_j^m$ by $h_j^m$ where $j,m \in \mathbb{N}_{0}$ and $m=0,1,2,\cdots, j$. 
Take $c_0$= 1 as $\lambda =0$ for the first kind of independent solutions of Lam\'{e} equation and $\lambda = \frac{1}{2}$ for the second one in the new (\ref{eq:100016a})--(\ref{eq:100020c}).
\begin{remark}
The power series expansion of Lam\'{e} equation in Weierstrass's form of the first kind for the first species complete polynomial using 3TRF about $\xi =0$ is given by
\begin{enumerate} 
\item As $\alpha = -1 $ or 0 with $h=h_0^0=0 $,
 
The eigenfunction is given by
\begin{eqnarray}
y(\xi ) &=& L_p^{(w)}F_{0,0}\left( \rho ,h=h_0^0=0, \alpha = -1\; \mbox{or}\; 0; \xi = sn^2(z,\rho ),\eta = -\rho ^2 \xi ^2,\mu = ( 1+\rho ^2 ) \xi \right) \nonumber\\
&=& 1 \label{eq:100058a}
\end{eqnarray}
\item As $\alpha  =-2\left( 2N+\frac{3}{2} \right)$ or $2\left( 2N+1 \right) $ where $N \in \mathbb{N}_{0}$,

An algebraic equation of degree $2N+2$ for the determination of $h$ is given by
\begin{equation}
0 = \sum_{r=0}^{N+1}\bar{c}\left( 2r, N+1-r; 2N+1,h \right)  \label{eq:100059a}
\end{equation}
The eigenvalue of $h$ is written by $h_{2N+1}^m$ where $m = 0,1,2,\cdots,2N+1 $; $h_{2N+1}^0 < h_{2N+1}^1 < \cdots < h_{2N+1}^{2N+1}$. Its eigenfunction is given by 
\begin{eqnarray} 
y(\xi ) &=& L_p^{(w)}F_{2N+1,m}\bigg( \rho, h=h_{2N+1}^m, \alpha = -2\left( 2N+\frac{3}{2} \right)\; \mbox{or}\; 2\left( 2N+1 \right); \xi = sn^2(z,\rho )\nonumber\\
&&,\eta = -\rho ^2 \xi ^2,\mu = ( 1+\rho ^2 ) \xi \bigg) \nonumber\\
&=& \sum_{r=0}^{N} y_{2r}^{N-r}\left( 2N+1,h_{2N+1}^m;\xi \right)+ \sum_{r=0}^{N} y_{2r+1}^{N-r}\left( 2N+1,h_{2N+1}^m;\xi \right)
\label{eq:100059b} 
\end{eqnarray}
\item As  $\alpha  =-2\left( 2N+\frac{5}{2} \right)$ or $2\left( 2N+2 \right) $ where $N \in \mathbb{N}_{0}$,

An algebraic equation of degree $2N+3$ for the determination of $h$ is given by
\begin{equation}  
0 = \sum_{r=0}^{N+1}\bar{c}\left( 2r+1, N+1-r; 2N+2,h \right) \label{eq:100060a}
\end{equation}
The eigenvalue of $h$ is written by $h_{2N+2}^m$ where $m = 0,1,2,\cdots,2N+2 $; $h_{2N+2}^0 < h_{2N+2}^1 < \cdots < h_{2N+2}^{2N+2}$. Its eigenfunction is given by
\begin{eqnarray} 
y(\xi ) &=& L_p^{(w)}F_{2N+2,m}\bigg( \rho ,h=h_{2N+2}^m, \alpha = -2\left( 2N+\frac{5}{2} \right)\; \mbox{or}\; 2\left( 2N+2 \right); \xi = sn^2(z,\rho )\nonumber\\
&&,\eta = -\rho ^2 \xi ^2,\mu = ( 1+\rho ^2 ) \xi \bigg) \nonumber\\
&=& \sum_{r=0}^{N+1} y_{2r}^{N+1-r}\left( 2N+2,h_{2N+2}^m;\xi \right) + \sum_{r=0}^{N} y_{2r+1}^{N-r}\left( 2N+2,h_{2N+2}^m;\xi \right) \label{eq:100060b}
\end{eqnarray}
In the above,
\begin{eqnarray}
\bar{c}(0,n;j,h)  &=& \frac{\left( -\frac{j}{2}\right)_{n}\left( \frac{1}{4}+\frac{j}{2} \right)_{n}}{ \left(  1 \right)_{n}\left( \frac{3}{4} \right)_{n}} \left( -\rho ^2\right)^{n}\label{eq:100061a}\\
\bar{c}(1,n;j,h) &=& \left( 1+\rho ^2\right) \sum_{i_0=0}^{n}\frac{ i_0^2 -\varpi(h) }{\left( i_0+\frac{1}{2} \right) \left( i_0+\frac{1}{4} \right)} \frac{\left( -\frac{j}{2}\right)_{i_0}\left( \frac{1}{4}+\frac{j}{2}  \right)_{i_0} }{\left( 1 \right)_{i_0} \left( \frac{3}{4} \right)_{i_0}} \nonumber\\
&&\times  \frac{\left( \frac{1}{2}-\frac{j}{2} \right)_{n} \left( \frac{3}{4}+\frac{j}{2} \right)_{n} \left( \frac{3}{2} \right)_{i_0}\left( \frac{5}{4} \right)_{i_0}}{\left( \frac{1}{2}-\frac{j}{2} \right)_{i_0} \left( \frac{3}{4}+\frac{j}{2} \right)_{i_0} \left( \frac{3}{2} \right)_n \left( \frac{5}{4} \right)_n} \left( -\rho ^2\right)^{n } \hspace{1cm}
\label{eq:100061b}\\
\bar{c}(\tau ,n;j,h) &=& \left( 1+\rho ^2\right)^{\tau } \sum_{i_0=0}^{n} \frac{ i_0^2 -\varpi(h) }{\left( i_0+\frac{1}{2} \right) \left( i_0+\frac{1}{4} \right)} \frac{\left( -\frac{j}{2}\right)_{i_0}\left( \frac{1}{4}+\frac{j}{2} \right)_{i_0} }{\left( 1 \right)_{i_0} \left( \frac{3}{4} \right)_{i_0}}  \nonumber\\
&&\times \prod_{k=1}^{\tau -1} \left( \sum_{i_k = i_{k-1}}^{n} \frac{\left( i_k+ \frac{k}{2} \right)^2 -\varpi(h)}{\left( i_k+\frac{k}{2}+\frac{1}{2} \right) \left( i_k+\frac{k}{2}+\frac{1}{4} \right)} \right. \nonumber\\
&&\times  \left. \frac{\left( \frac{k}{2}-\frac{j}{2}\right)_{i_k} \left( \frac{k}{2}+\frac{1}{4}+\frac{j}{2} \right)_{i_k}\left( \frac{k}{2}+1 \right)_{i_{k-1}} \left( \frac{k}{2}+ \frac{3}{4} \right)_{i_{k-1}}}{\left( \frac{k}{2}-\frac{j}{2}\right)_{i_{k-1}} \left( \frac{k}{2}+\frac{1}{4}+\frac{j}{2} \right)_{i_{k-1}}\left( \frac{k}{2}+1 \right)_{i_k} \left( \frac{k}{2}+ \frac{3}{4} \right)_{i_k}} \right) \nonumber\\ 
&&\times \frac{\left( \frac{\tau }{2} -\frac{j}{2}\right)_{n}\left( \frac{\tau }{2}+\frac{1}{4} +\frac{j}{2} \right)_{n} \left( \frac{\tau }{2}+1 \right)_{i_{\tau -1}} \left( \frac{\tau }{2}+\frac{3}{4} \right)_{i_{\tau -1}}}{\left( \frac{\tau }{2} -\frac{j}{2}\right)_{i_{\tau -1}}\left( \frac{\tau }{2}+\frac{1}{4} +\frac{j}{2} \right)_{i_{\tau -1}} \left( \frac{\tau }{2}+1 \right)_n \left( \frac{\tau }{2}+\frac{3}{4} \right)_n} \left( -\rho ^2\right)^{n } \hspace{1.5cm}\label{eq:100061c} 
\end{eqnarray}
\begin{eqnarray}
y_0^m(j,h ;\xi ) &=& \sum_{i_0=0}^{m} \frac{\left( -\frac{j}{2}\right)_{i_0}\left( \frac{1}{4}+\frac{j}{2} \right)_{i_0} }{\left( 1 \right)_{i_0} \left( \frac{3}{4} \right)_{i_0}} \eta ^{i_0} \label{eq:100062a}\\
y_1^m(j,h ;\xi ) &=& \left\{\sum_{i_0=0}^{m} \frac{ i_0 ^2 -\varpi(h) }{\left( i_0+\frac{1}{2} \right) \left( i_0+\frac{1}{4} \right)} \frac{\left( -\frac{j}{2}\right)_{i_0}\left( \frac{1}{4}+\frac{j}{2} \right)_{i_0} }{\left( 1 \right)_{i_0} \left( \frac{3}{4} \right)_{i_0}} \right. \nonumber\\
&&\times   \left. \sum_{i_1 = i_0}^{m} \frac{\left( \frac{1}{2}-\frac{j}{2} \right)_{i_1} \left( \frac{3}{4}+\frac{j}{2} \right)_{i_1} \left( \frac{3}{2} \right)_{i_0}\left( \frac{5}{4} \right)_{i_0}}{\left( \frac{1}{2}-\frac{j}{2} \right)_{i_0} \left( \frac{3}{4}+\frac{j}{2} \right)_{i_0} \left( \frac{3}{2} \right)_{i_1} \left( \frac{5}{4} \right)_{i_1}} \eta ^{i_1}\right\} \mu \hspace{1cm}\label{eq:100062b}
\end{eqnarray}
\begin{eqnarray}
y_{\tau }^m(j,h ;\xi ) &=& \left\{ \sum_{i_0=0}^{m} \frac{ i_0 ^2 -\varpi(h) }{\left( i_0+\frac{1}{2} \right) \left( i_0+\frac{1}{4} \right)} \frac{\left( -\frac{j}{2}\right)_{i_0}\left( \frac{1}{4}+\frac{j}{2} \right)_{i_0} }{\left( 1 \right)_{i_0} \left( \frac{3}{4} \right)_{i_0}} \right.\nonumber\\
&&\times \prod_{k=1}^{\tau -1} \left( \sum_{i_k = i_{k-1}}^{m} \frac{\left( i_k+ \frac{k}{2} \right)^2 -\varpi(h)}{\left( i_k+\frac{k}{2}+\frac{1}{2} \right) \left( i_k+\frac{k}{2}+\frac{1}{4} \right)} \right. \nonumber\\
&&\times   \left. \frac{\left( \frac{k}{2}-\frac{j}{2}\right)_{i_k} \left( \frac{k}{2}+\frac{1}{4}+\frac{j}{2} \right)_{i_k}\left( \frac{k}{2}+1 \right)_{i_{k-1}} \left( \frac{k}{2}+ \frac{3}{4} \right)_{i_{k-1}}}{\left( \frac{k}{2}-\frac{j}{2}\right)_{i_{k-1}} \left( \frac{k}{2}+\frac{1}{4}+\frac{j}{2} \right)_{i_{k-1}}\left( \frac{k}{2}+1 \right)_{i_k} \left( \frac{k}{2}+ \frac{3}{4} \right)_{i_k}} \right) \nonumber\\
&&\times \left. \sum_{i_{\tau } = i_{\tau -1}}^{m} \frac{\left( \frac{\tau }{2} -\frac{j}{2}\right)_{i_{\tau }}\left( \frac{\tau }{2}+\frac{1}{4} +\frac{j}{2} \right)_{i_{\tau }} \left( \frac{\tau }{2}+1 \right)_{i_{\tau -1}} \left( \frac{\tau }{2}+\frac{3}{4}  \right)_{i_{\tau -1}}}{\left( \frac{\tau }{2} -\frac{j}{2}\right)_{i_{\tau -1}}\left( \frac{\tau }{2}+\frac{1}{4} +\frac{j}{2} \right)_{i_{\tau -1}} \left( \frac{\tau }{2}+1 \right)_{i_{\tau }} \left( \frac{\tau }{2}+\frac{3}{4} \right)_{i_{\tau }}} \eta ^{i_{\tau }}\right\} \mu ^{\tau }\hspace{1.5cm} \label{eq:100062c}
\end{eqnarray}
where
\begin{equation}
\begin{cases} \tau \geq 2 \cr
\varpi(h) = \frac{h}{4\left( 1+\rho ^2\right)}
\end{cases}\nonumber
\end{equation}
\end{enumerate}
\end{remark}
\begin{remark}
The power series expansion of Lam\'{e} equation in Weierstrass's form of the second kind for the first species complete polynomial using 3TRF about $\xi =0$ is given by
\begin{enumerate} 
\item As $\alpha = -2 $ or 1 with $h=h_0^0 =1+\rho ^2 $,
 
The eigenfunction is given by
\begin{eqnarray}
y(\xi ) &=& L_p^{(w)}S_{0,0}\left( \rho ,h=h_0^0 =1+\rho ^2, \alpha = -2\; \mbox{or}\; 1; \xi = sn^2(z,\rho ),\eta = -\rho ^2 \xi ^2,\mu = ( 1+\rho ^2 ) \xi \right) \nonumber\\
&=& \xi ^{\frac{1}{2}} \label{eq:100063a}
\end{eqnarray}
\item As $\alpha  =-2\left( 2N+2 \right)$ or $2\left( 2N+\frac{3}{2} \right) $ where $N \in \mathbb{N}_{0}$,

An algebraic equation of degree $2N+2$ for the determination of $h$ is given by
\begin{equation}
0 = \sum_{r=0}^{N+1}\bar{c}\left( 2r, N+1-r; 2N+1,h \right)  \label{eq:100064a}
\end{equation}
The eigenvalue of $h$ is written by $h_{2N+1}^m$ where $m = 0,1,2,\cdots,2N+1 $; $h_{2N+1}^0 < h_{2N+1}^1 < \cdots < h_{2N+1}^{2N+1}$. Its eigenfunction is given by 
\begin{eqnarray} 
y(\xi ) &=& L_p^{(w)}S_{2N+1,m}\bigg( \rho ,h=h_{2N+1}^m, \alpha = -2\left( 2N+2 \right)\; \mbox{or}\; 2\left( 2N+\frac{3}{2} \right); \xi = sn^2(z,\rho )\nonumber\\
&&,\eta = -\rho ^2 \xi ^2,\mu = ( 1+\rho ^2 ) \xi \bigg) \nonumber\\
&=& \sum_{r=0}^{N} y_{2r}^{N-r}\left( 2N+1,h_{2N+1}^m;\xi \right)+ \sum_{r=0}^{N} y_{2r+1}^{N-r}\left( 2N+1,h_{2N+1}^m;\xi \right)
\label{eq:100064b} 
\end{eqnarray}
\item As  $\alpha  =-2\left( 2N+3 \right)$ or $2\left( 2N+\frac{5}{2} \right) $ where $N \in \mathbb{N}_{0}$,

An algebraic equation of degree $2N+3$ for the determination of $h$ is given by
\begin{equation}  
0 = \sum_{r=0}^{N+1}\bar{c}\left( 2r+1, N+1-r; 2N+2,h \right) \label{eq:100065a}
\end{equation}
The eigenvalue of $h$ is written by $h_{2N+2}^m$ where $m = 0,1,2,\cdots,2N+2 $; $h_{2N+2}^0 < h_{2N+2}^1 < \cdots < h_{2N+2}^{2N+2}$. Its eigenfunction is given by
\begin{eqnarray} 
y(\xi ) &=& L_p^{(w)}S_{2N+2,m}\bigg( \rho ,h=h_{2N+2}^m, \alpha = -2\left( 2N+3 \right)\; \mbox{or}\; 2\left( 2N+\frac{5}{2} \right); \xi = sn^2(z,\rho )\nonumber\\
&&,\eta = -\rho ^2 \xi ^2,\mu = ( 1+\rho ^2 ) \xi \bigg) \nonumber\\
&=& \sum_{r=0}^{N+1} y_{2r}^{N+1-r}\left( 2N+2,h_{2N+2}^m;\xi \right) + \sum_{r=0}^{N} y_{2r+1}^{N-r}\left( 2N+2,h_{2N+2}^m;\xi \right) \label{eq:100065b}
\end{eqnarray}
In the above,
\begin{eqnarray}
\bar{c}(0,n;j,h)  &=& \frac{\left( -\frac{j}{2}\right)_{n}\left( \frac{3}{4}+\frac{j}{2} \right)_{n}}{ \left(  \frac{5}{4} \right)_{n}\left( 1 \right)_{n}} \left( -\rho ^2\right)^{n}\label{eq:100066a}\\
\bar{c}(1,n;j,h) &=& \left( 1+\rho ^2\right) \sum_{i_0=0}^{n}\frac{\left( i_0 +\frac{1}{4}\right)^2 -\varpi(h) }{\left( i_0+\frac{3}{4} \right) \left( i_0+\frac{1}{2} \right)} \frac{\left( -\frac{j}{2}\right)_{i_0}\left( \frac{3}{4}+\frac{j}{2} \right)_{i_0} }{\left( \frac{5}{4}\right)_{i_0} \left( 1\right)_{i_0}} \nonumber\\
&&\times \frac{\left( \frac{1}{2}-\frac{j}{2} \right)_{n} \left( \frac{5}{4}+\frac{j}{2} \right)_{n} \left( \frac{7}{4} \right)_{i_0}\left( \frac{3}{2} \right)_{i_0}}{\left( \frac{1}{2}-\frac{j}{2} \right)_{i_0} \left( \frac{5}{4}+\frac{j}{2}  \right)_{i_0} \left( \frac{7}{4} \right)_n \left( \frac{3}{2} \right)_n} \left( -\rho ^2\right)^{n }\hspace{1cm}\label{eq:100066b}\\
\bar{c}(\tau ,n;j,h) &=& \left( 1+\rho ^2\right)^{\tau } \sum_{i_0=0}^{n} \frac{\left( i_0 +\frac{1}{4}\right)^2 -\varpi(h) }{\left( i_0+\frac{3}{4} \right) \left( i_0+\frac{1}{2} \right)} \frac{\left( -\frac{j}{2}\right)_{i_0}\left( \frac{3}{4}+\frac{j}{2} \right)_{i_0} }{\left( \frac{5}{4}\right)_{i_0} \left( 1\right)_{i_0}}  \nonumber\\
&&\times \prod_{k=1}^{\tau -1} \left( \sum_{i_k = i_{k-1}}^{n} \frac{\left( i_k+ \frac{k}{2} +\frac{1}{4}\right)^2 -\varpi(h)}{\left( i_k+\frac{k}{2}+\frac{3}{4} \right) \left( i_k+\frac{k}{2}+\frac{1}{2} \right)} \right. \nonumber\\
&&\times \left. \frac{\left( \frac{k}{2}-\frac{j}{2}\right)_{i_k} \left( \frac{k}{2}+\frac{3}{4}+\frac{j}{2} \right)_{i_k}\left( \frac{k}{2} + \frac{5}{4}\right)_{i_{k-1}} \left( \frac{k}{2}+ 1\right)_{i_{k-1}}}{\left( \frac{k}{2}-\frac{j}{2}\right)_{i_{k-1}} \left( \frac{k}{2}+\frac{3}{4}+\frac{j}{2}  \right)_{i_{k-1}}\left( \frac{k}{2} + \frac{5}{4}\right)_{i_k} \left( \frac{k}{2}+ 1\right)_{i_k}} \right) \nonumber \\ 
&&\times \frac{\left( \frac{\tau }{2} -\frac{j}{2}\right)_{n}\left( \frac{\tau }{2}+\frac{3}{4} +\frac{j}{2}  \right)_{n} \left( \frac{\tau }{2} +\frac{5}{4}\right)_{i_{\tau -1}} \left( \frac{\tau }{2}+1\right)_{i_{\tau -1}}}{\left( \frac{\tau }{2} -\frac{j}{2}\right)_{i_{\tau -1}}\left( \frac{\tau }{2}+\frac{3}{4} +\frac{j}{2} \right)_{i_{\tau -1}} \left( \frac{\tau }{2} +\frac{5}{4}\right)_n \left( \frac{\tau }{2}+1\right)_n} \left( -\rho ^2\right)^{n }\hspace{1.5cm} \label{eq:100066c}
\end{eqnarray}
\begin{eqnarray}
y_0^m(j,h ;\xi ) &=& \xi ^{\frac{1}{2}}  \sum_{i_0=0}^{m} \frac{\left( -\frac{j}{2}\right)_{i_0}\left( \frac{3}{4}+\frac{j}{2} \right)_{i_0} }{\left( \frac{5}{4}\right)_{i_0} \left( 1\right)_{i_0}} \eta ^{i_0} \label{eq:100067a}\\
y_1^m(j,h ;\xi ) &=& \xi ^{\frac{1}{2}} \left\{\sum_{i_0=0}^{m} \frac{\left( i_0 +\frac{1}{4}\right)^2 -\varpi(h) }{\left( i_0+\frac{3}{4} \right) \left( i_0+\frac{1}{2} \right)} \frac{\left( -\frac{j}{2}\right)_{i_0}\left( \frac{3}{4}+\frac{j}{2} \right)_{i_0} }{\left( \frac{5}{4}\right)_{i_0} \left( 1\right)_{i_0}} \right.  \nonumber\\
&&\times \left. \sum_{i_1 = i_0}^{m} \frac{\left( \frac{1}{2}-\frac{j}{2} \right)_{i_1} \left( \frac{5}{4}+\frac{j}{2} \right)_{i_1} \left( \frac{7}{4} \right)_{i_0}\left( \frac{3}{2} \right)_{i_0}}{\left( \frac{1}{2}-\frac{j}{2} \right)_{i_0} \left( \frac{5}{4}+\frac{j}{2} \right)_{i_0} \left( \frac{7}{4} \right)_{i_1} \left( \frac{3}{2} \right)_{i_1}} \eta ^{i_1}\right\} \mu \hspace{1cm} 
\label{eq:100067b}
\end{eqnarray}
\begin{eqnarray}
y_{\tau }^m(j,h ;\xi ) &=& \xi ^{\frac{1}{2}} \left\{ \sum_{i_0=0}^{m} \frac{\left( i_0 +\frac{1}{4}\right) ^2 -\varpi(h)}{\left( i_0+\frac{3}{4} \right) \left( i_0+\frac{1}{2} \right)} \frac{\left( -\frac{j}{2}\right)_{i_0}\left( \frac{3}{4}+\frac{j}{2} \right)_{i_0} }{\left( \frac{5}{4}\right)_{i_0} \left( 1\right)_{i_0}} \right.\nonumber\\
&&\times \prod_{k=1}^{\tau -1} \left( \sum_{i_k = i_{k-1}}^{m} \frac{\left( i_k+ \frac{k}{2} +\frac{1}{4}\right) ^2 -\varpi(h)}{\left( i_k+\frac{k}{2}+\frac{3}{4} \right) \left( i_k+\frac{k}{2}+\frac{1}{2} \right)} \right.\nonumber\\
&&\times \left. \frac{\left( \frac{k}{2}-\frac{j}{2}\right)_{i_k} \left( \frac{k}{2}+\frac{3}{4}+\frac{j}{2} \right)_{i_k}\left( \frac{k}{2} + \frac{5}{4}\right)_{i_{k-1}} \left( \frac{k}{2}+ 1 \right)_{i_{k-1}}}{\left( \frac{k}{2}-\frac{j}{2}\right)_{i_{k-1}} \left( \frac{k}{2}+\frac{3}{4}+\frac{j}{2} \right)_{i_{k-1}}\left( \frac{k}{2} + \frac{5}{4}\right)_{i_k} \left( \frac{k}{2}+ 1\right)_{i_k}} \right) \nonumber \\
&&\times \left. \sum_{i_{\tau } = i_{\tau -1}}^{m} \frac{\left( \frac{\tau }{2} -\frac{j}{2}\right)_{i_{\tau }}\left( \frac{\tau }{2}+\frac{3}{4} +\frac{j}{2} \right)_{i_{\tau }} \left( \frac{\tau }{2} +\frac{5}{4}\right)_{i_{\tau -1}} \left( \frac{\tau }{2}+1\right)_{i_{\tau -1}}}{\left( \frac{\tau }{2} -\frac{j}{2}\right)_{i_{\tau -1}}\left( \frac{\tau }{2}+\frac{3}{4} +\frac{j}{2} \right)_{i_{\tau -1}} \left( \frac{\tau }{2} +\frac{5}{4}\right)_{i_{\tau }} \left( \frac{\tau }{2}+1\right)_{i_{\tau }}} \eta ^{i_{\tau }}\right\} \mu ^{\tau } \hspace{1.5cm} \label{eq:100067c}
\end{eqnarray}
where
\begin{equation}
\begin{cases} \tau \geq 2 \cr
\varpi(h) = \frac{h}{4\left( 1+\rho ^2\right)}
\end{cases}\nonumber
\end{equation}
\end{enumerate}
\end{remark}
If we take $\alpha \geq -\frac{1}{2}$ for solutions of Lam\'{e} equation in Weierstrass's form, $\alpha  =-2\left( j+\frac{1}{2}+\lambda \right)$ with $\lambda =0$ or $1/2$ where $j \in \mathbb{N}_{0}$ is not available any more for the first species complete polynomial. In this chapter, I suggest that $\alpha  \in \mathbb{N}_{0}$ is allowed for formal series solutions of the first species complete polynomial.
\subsection{The first species complete polynomial of Lam\'{e} equation using R3TRF} 
Substitute (\ref{eq:100057}) into (\ref{eq:100037a})--(\ref{eq:100042c}) with replacing a spectral parameter $q_j^m$ by $h_j^m$ where $j,m \in \mathbb{N}_{0}$ and $m=0,1,2,\cdots, j$. 
Take $c_0$= 1 as $\lambda =0$ for the first kind of independent solutions of Lam\'{e} equation and $\lambda = \frac{1}{2}$ for the second one in the new (\ref{eq:100037a})--(\ref{eq:100042c}). 
\begin{remark}
The power series expansion of Lam\'{e} equation in Weierstrass's form of the first kind for the first species complete polynomial using R3TRF about $\xi =0$ is given by
\begin{enumerate}
\item As $\alpha = -1 $ or 0 with $h=h_0^0=0 $,
 
The eigenfunction is given by
\begin{eqnarray}
y(\xi ) &=& L_p^{(w)}F_{0,0}^R\left( \rho ,h=h_0^0=0, \alpha = -1\; \mbox{or}\; 0; \xi = sn^2(z,\rho ),\eta = -\rho ^2 \xi ^2,\mu = ( 1+\rho ^2 ) \xi \right) \nonumber\\
&=& 1 \label{eq:100068}
\end{eqnarray}
\item As $\alpha = -3 $ or $2 $,

An algebraic equation of degree 2 for the determination of $h$ is given by
\begin{equation}
0 =\frac{3}{4\left( \rho + \rho ^{-1}\right)^2} +\prod_{l=0}^{1} \left( l^2 -\frac{h}{4(1+\rho ^2 )} \right) \label{eq:100069a}
\end{equation}
The eigenvalue of $h$ is written by $h_1^m$ where $m = 0,1 $; $h_{1}^0 < h_{1}^1$. Its eigenfunction is given by
\begin{eqnarray}
y(\xi ) &=& L_p^{(w)}F_{1,m}^R\left( \rho ,h=h_1^m, \alpha = -3\; \mbox{or}\; 2; \xi = sn^2(z,\rho ),\eta = -\rho ^2 \xi ^2,\mu = ( 1+\rho ^2 ) \xi \right) \nonumber\\
&=& 1- \frac{ h_1^m }{2(1+\rho ^2)} \mu \label{eq:100069b}
\end{eqnarray}
\item As $\alpha = -2\left( 2N +\frac{5}{2} \right) $ or $2\left( 2N+2 \right) $ where $N \in \mathbb{N}_{0}$,

An algebraic equation of degree $2N+3$ for the determination of $h$ is given by
\begin{equation}
0 = \sum_{r=0}^{N+1}\bar{c}\left( r, 2(N-r)+3; 2N+2,h \right)  \label{eq:100070a}
\end{equation}
The eigenvalue of $h$ is written by $h_{2N+2}^m$ where $m = 0,1,2,\cdots,2N+2 $; $h_{2N+2}^0 < h_{2N+2}^1 < \cdots < h_{2N+2}^{2N+2}$. Its eigenfunction is given by 
\begin{eqnarray} 
y(\xi ) &=& L_p^{(w)}F_{2N+2,m}^R\bigg( \rho ,h=h_{2N+2}^m, \alpha = -2\left( 2N+\frac{5}{2} \right)\; \mbox{or}\; 2\left( 2N+2 \right); \xi = sn^2(z,\rho )\nonumber\\
&&,\eta = -\rho ^2 \xi ^2,\mu = ( 1+\rho ^2 ) \xi  \bigg) \nonumber\\
&=& \sum_{r=0}^{N+1} y_{r}^{2(N+1-r)}\left( 2N+2, h_{2N+2}^m; \xi \right)  
\label{eq:100070b} 
\end{eqnarray}
\item As $\alpha = -2\left( 2N +\frac{7}{2} \right) $ or $2\left( 2N+3 \right) $ where $N \in \mathbb{N}_{0}$,

An algebraic equation of degree $2N+4$ for the determination of $h$ is given by
\begin{equation}  
0 = \sum_{r=0}^{N+2}\bar{c}\left( r, 2(N+2-r); 2N+3,h \right) \label{eq:100071a}
\end{equation}
The eigenvalue of $h$ is written by $h_{2N+3}^m$ where $m = 0,1,2,\cdots,2N+3 $; $h_{2N+3}^0 < h_{2N+3}^1 < \cdots < h_{2N+3}^{2N+3}$. Its eigenfunction is given by
\begin{eqnarray} 
y(\xi ) &=& L_p^{(w)}F_{2N+3,m}^R\bigg( \rho ,h=h_{2N+3}^m, \alpha = -2\left( 2N+\frac{7}{2} \right)\; \mbox{or}\; 2\left( 2N+3 \right); \xi = sn^2(z,\rho )\nonumber\\
&&,\eta = -\rho ^2 \xi ^2,\mu = ( 1+\rho ^2 ) \xi \bigg) \nonumber\\
&=& \sum_{r=0}^{N+1} y_{r}^{2(N-r)+3} \left( 2N+3, h_{2N+3}^m; \xi \right) \label{eq:100071b}
\end{eqnarray}
In the above,
\begin{eqnarray}
\bar{c}(0,n;j,h)  &=& \frac{\left( -\varphi (h) \right)_{n}\left( \varphi (h) \right)_{n}}{\left( 1\right)_{n} \left( \frac{1}{2} \right)_{n}} \left( 1+\rho ^2 \right)^{n}\label{eq:100072a}\\
\bar{c}(1,n;j,h) &=& \left( -\rho ^2\right) \sum_{i_0=0}^{n}\frac{\left( i_0-j\right) \left( i_0+\frac{1}{2}+j\right)}{\left( i_0+2\right) \left( i_0+ \frac{3}{2} \right)} \frac{\left( -\varphi (h) \right)_{i_0}\left( \varphi (h) \right)_{i_0}}{\left( 1 \right)_{i_0} \left( \frac{1}{2} \right)_{i_0}} \nonumber\\
&&\times \frac{\left(  2 -\varphi (h) \right)_{n}\left(  2 +\varphi (h) \right)_{n} \left( 3 \right)_{i_0} \left( \frac{5}{2} \right)_{i_0}}{\left(  2 -\varphi (h) \right)_{i_0}\left(  2 +\varphi (h) \right)_{i_0} \left( 3 \right)_{n} \left( \frac{5}{2} \right)_{n}}\left( 1+\rho ^2 \right)^n  
 \label{eq:100072b}\\
\bar{c}(\tau ,n;j,h) &=& \left( -\rho ^2\right)^{\tau } \sum_{i_0=0}^{n} \frac{\left( i_0-j\right) \left( i_0+\frac{1}{2}+j \right)}{\left( i_0+2 \right) \left( i_0+ \frac{3}{2} \right)} \frac{\left( -\varphi (h) \right)_{i_0}\left( \varphi (h) \right)_{i_0}}{\left( 1\right)_{i_0} \left( \frac{1}{2} \right)_{i_0}}  \nonumber\\
&&\times \prod_{k=1}^{\tau -1} \left( \sum_{i_k = i_{k-1}}^{n} \frac{\left( i_k+ 2k-j\right)\left( i_k+ 2k+\frac{1}{2}+j \right) }{\left( i_k+2k+2 \right) \left( i_k+2k+\frac{3}{2} \right)} \right. \nonumber\\
&&\times \left. \frac{\left( 2k-\varphi (h) \right)_{i_k} \left( 2k+\varphi (h) \right)_{i_k}\left( 2k+1 \right)_{i_{k-1}} \left( 2k+\frac{1}{2} \right)_{i_{k-1}}}{\left( 2k-\varphi (h) \right)_{i_{k-1}} \left( 2k+\varphi (h) \right)_{i_{k-1}}\left( 2k+1 \right)_{i_k} \left( 2k+\frac{1}{2} \right)_{i_k}} \right) \nonumber\\ 
&&\times \frac{\left( 2\tau -\varphi (h) \right)_n \left( 2\tau +\varphi (h) \right)_n \left( 2\tau +1 \right)_{i_{\tau -1}} \left( 2\tau +\frac{1}{2} \right)_{i_{\tau -1}}}{\left( 2\tau -\varphi (h) \right)_{i_{\tau -1}} \left( 2\tau +\varphi (h) \right)_{i_{\tau -1}}\left( 2\tau +1 \right)_n \left( 2\tau +\frac{1}{2}  \right)_n} \left( 1+\rho ^2 \right)  \hspace{2cm}\label{eq:100072c}
\end{eqnarray}
\begin{eqnarray}
y_0^m(j,h ;\xi ) &=& \sum_{i_0=0}^{m} \frac{\left( -\varphi (h) \right)_{i_0}\left( \varphi (h) \right)_{i_0}}{\left( 1 \right)_{i_0} \left( \frac{1}{2} \right)_{i_0}} \mu ^{i_0} \label{eq:100073a}\\
y_1^m(j,h ;\xi ) &=&  \left\{\sum_{i_0=0}^{m} \frac{\left( i_0-j\right) \left( i_0+\frac{1}{2}+j \right)}{\left( i_0+ 2 \right) \left( i_0+ \frac{3}{2} \right)} \frac{\left( -\varphi (h) \right)_{i_0}\left( \varphi (h)  \right)_{i_0}}{\left( 1 \right)_{i_0} \left( \frac{1}{2}  \right)_{i_0}} \right. \nonumber\\
&&\times \left. \sum_{i_1 = i_0}^{m} \frac{\left( 2 -\varphi (h) \right)_{i_1}\left(  2 +\varphi (h)  \right)_{i_1} \left( 3 \right)_{i_0} \left( \frac{5}{2} \right)_{i_0}}{\left( 2 -\varphi (h) \right)_{i_0}\left( 2 +\varphi (h) \right)_{i_0} \left( 3 \right)_{i_1} \left( \frac{5}{2} \right)_{i_1}} \mu ^{i_1}\right\} \eta 
\label{eq:100073b}
\end{eqnarray}
\begin{align}
y_{\tau }^m(j,h ;\xi ) &=  \left\{ \sum_{i_0=0}^{m} \frac{\left( i_0-j\right) \left( i_0+\frac{1}{2}+j \right)}{\left( i_0+ 2 \right) \left( i_0+ \frac{3}{2} \right)} \frac{\left( -\varphi (h) \right)_{i_0}\left( \varphi (h) \right)_{i_0}}{\left( 1 \right)_{i_0} \left( \frac{1}{2} \right)_{i_0}} \right.\nonumber\\
 &\times \prod_{k=1}^{\tau -1} \left( \sum_{i_k = i_{k-1}}^{m} \frac{\left( i_k+ 2k-j\right)\left( i_k+ 2k+\frac{1}{2}+j \right) }{\left( i_k+2k+2 \right) \left( i_k+2k+\frac{3}{2} \right)} \right. \nonumber\\
 &\times \left. \frac{\left( 2k-\varphi (h) \right)_{i_k} \left( 2k+\varphi (h) \right)_{i_k}\left( 2k+1 \right)_{i_{k-1}} \left( 2k+\frac{1}{2} \right)_{i_{k-1}}}{\left( 2k-\varphi (h) \right)_{i_{k-1}} \left( 2k+\varphi (h) \right)_{i_{k-1}}\left( 2k+1 \right)_{i_k} \left( 2k+\frac{1}{2} \right)_{i_k}} \right)  \nonumber\\
 &\times \left. \sum_{i_{\tau } = i_{\tau -1}}^{m}  \frac{\left( 2\tau -\varphi (h) \right)_{i_{\tau }} \left( 2\tau +\varphi (h) \right)_{i_{\tau }} \left( 2\tau +1 \right)_{i_{\tau -1}} \left( 2\tau +\frac{1}{2} \right)_{i_{\tau -1}}}{\left( 2\tau -\varphi (h) \right)_{i_{\tau -1}} \left( 2\tau +\varphi (h) \right)_{i_{\tau -1}}\left( 2\tau +1 \right)_{i_{\tau }} \left( 2\tau +\frac{1}{2} \right)_{i_{\tau }}} \mu ^{i_{\tau }}\right\} \eta ^{\tau }  \label{eq:100073c}
\end{align}
where
\begin{equation}
\begin{cases} \tau \geq 2 \cr
\varphi (h) = \sqrt{\frac{h}{ 4(1+\rho ^2)}}
\end{cases}\nonumber
\end{equation}
\end{enumerate}
\end{remark} 
\begin{remark}
The power series expansion of Lam\'{e} equation in Weierstrass's form of the second kind for the first species complete polynomial using R3TRF about $\xi =0$ is given by
\begin{enumerate} 
\item As $\alpha = -2 $ or 1 with $h=h_0^0= 1+ \rho ^2 $,
 
The eigenfunction is given by
\begin{eqnarray}
y(\xi ) &=& L_p^{(w)}S_{0,0}^R\left( \rho ,h=h_0^0=1+ \rho ^2, \alpha = -2\; \mbox{or}\; 1; \xi = sn^2(z,\rho ),\eta = -\rho ^2 \xi ^2,\mu = ( 1+\rho ^2 ) \xi \right) \nonumber\\
&=& \xi^{\frac{1}{2}} \label{eq:100074}
\end{eqnarray}
\item As $\alpha = -4 $ or $3 $,

An algebraic equation of degree 2 for the determination of $h$ is given by
\begin{equation}
0 = \frac{15}{4\left( \rho + \rho ^{-1}\right)^2} +\prod_{l=0}^{1} \left( \left( l+\frac{1}{2} \right)^2 -\frac{h}{4(1+ \rho ^2)} \right)  \label{eq:100075a}
\end{equation}
The eigenvalue of $h$ is written by $h_1^m$ where $m = 0,1 $; $h_{1}^0 < h_{1}^1$. Its eigenfunction is given by
\begin{eqnarray}
y(\xi ) &=& L_p^{(w)}S_{1,m}^R\left( \rho ,h=h_1^m, \alpha = -4\; \mbox{or}\; 3; \xi = sn^2(z,\rho ),\eta = -\rho ^2 \xi ^2,\mu = ( 1+\rho ^2 ) \xi \right) \nonumber\\
&=& \xi^{\frac{1}{2}}\left\{ 1+ \frac{1}{6} \left( 1 -\frac{ h_1^m }{(1+\rho ^2)} \right) \mu \right\} \label{eq:100075b}
\end{eqnarray}
\item As $\alpha = -2\left( 2N +3 \right) $ or $2\left( 2N+\frac{5}{2} \right) $ where $N \in \mathbb{N}_{0}$,

An algebraic equation of degree $2N+3$ for the determination of $h$ is given by
\begin{equation}
0 = \sum_{r=0}^{N+1}\bar{c}\left( r, 2(N-r)+3; 2N+2,h \right)  \label{eq:100076a}
\end{equation}
The eigenvalue of $h$ is written by $h_{2N+2}^m$ where $m = 0,1,2,\cdots,2N+2 $; $h_{2N+2}^0 < h_{2N+2}^1 < \cdots < h_{2N+2}^{2N+2}$. Its eigenfunction is given by 
\begin{eqnarray} 
y(\xi ) &=& L_p^{(w)}S_{2N+2,m}^R\bigg( \rho ,h=h_{2N+2}^m, \alpha = -2\left( 2N+3 \right)\; \mbox{or}\; 2\left( 2N+\frac{5}{2} \right); \xi = sn^2(z,\rho )\nonumber\\
&&,\eta = -\rho ^2 \xi ^2,\mu = ( 1+\rho ^2 ) \xi \bigg) \nonumber\\
&=& \sum_{r=0}^{N+1} y_{r}^{2(N+1-r)}\left( 2N+2, h_{2N+2}^m; \xi \right)  
\label{eq:100076b} 
\end{eqnarray}
\item As $\alpha = -2\left( 2N +4 \right) $ or $2\left( 2N+\frac{7}{2} \right) $ where $N \in \mathbb{N}_{0}$,

An algebraic equation of degree $2N+4$ for the determination of $h$ is given by
\begin{equation}  
0 = \sum_{r=0}^{N+2}\bar{c}\left( r, 2(N+2-r); 2N+3,h \right) \label{eq:100077a}
\end{equation}
The eigenvalue of $h$ is written by $h_{2N+3}^m$ where $m = 0,1,2,\cdots,2N+3 $; $h_{2N+3}^0 < h_{2N+3}^1 < \cdots < h_{2N+3}^{2N+3}$. Its eigenfunction is given by
\begin{eqnarray} 
y(\xi ) &=& L_p^{(w)}S_{2N+3,m}^R\bigg( \rho ,h=h_{2N+3}^m, \alpha = -2\left( 2N+4 \right)\; \mbox{or}\; 2\left( 2N+\frac{7}{2} \right); \xi = sn^2(z,\rho )\nonumber\\
&&,\eta = -\rho ^2 \xi ^2,\mu = ( 1+\rho ^2 ) \xi  \bigg) \nonumber\\
&=& \sum_{r=0}^{N+1} y_{r}^{2(N-r)+3} \left( 2N+3, h_{2N+3}^m; \xi \right) \label{eq:100077b}
\end{eqnarray}
In the above,
\begin{eqnarray}
\bar{c}(0,n;j,h)  &=& \frac{\left( \frac{1}{2}-\varphi (h) \right)_{n}\left( \frac{1}{2}+\varphi (h) \right)_{n}}{\left( \frac{3}{2}\right)_{n} \left( 1 \right)_{n}} \left( 1+\rho ^2 \right)^{n}\label{eq:100078a}\\
\bar{c}(1,n;j,h) &=& \left( -\rho ^2\right) \sum_{i_0=0}^{n}\frac{\left( i_0-j\right) \left( i_0+\frac{3}{2}+j\right)}{\left( i_0+ \frac{5}{2} \right) \left( i_0+ 2 \right)} \frac{\left( \frac{1}{2}-\varphi (h) \right)_{i_0}\left( \frac{1}{2}+\varphi (h) \right)_{i_0}}{\left( \frac{3}{2} \right)_{i_0} \left( 1\right)_{i_0}} \nonumber\\
&&\times \frac{\left(  \frac{5}{2}-\varphi (h) \right)_{n}\left(  \frac{5}{2}+\varphi (h) \right)_{n} \left( \frac{7}{2} \right)_{i_0} \left( 3 \right)_{i_0}}{\left(  \frac{5}{2}-\varphi (h)  \right)_{i_0}\left( \frac{5}{2}+\varphi (h) \right)_{i_0} \left( \frac{7}{2} \right)_{n} \left( 3\right)_{n}}\left( 1+\rho ^2 \right)^n  
 \label{eq:100078b}\\
\bar{c}(\tau ,n;j,h) &=& \left( -\rho ^2\right)^{\tau } \sum_{i_0=0}^{n} \frac{\left( i_0-j\right) \left( i_0+\frac{3}{2}+j  \right)}{\left( i_0+ \frac{5}{2} \right) \left( i_0+ 2\right)} \frac{\left( \frac{1}{2}-\varphi (h) \right)_{i_0}\left( \frac{1}{2}+\varphi (h) \right)_{i_0}}{\left( \frac{3}{2} \right)_{i_0} \left( 1 \right)_{i_0}}  \nonumber\\
&&\times \prod_{k=1}^{\tau -1} \left( \sum_{i_k = i_{k-1}}^{n} \frac{\left( i_k+ 2k-j\right)\left( i_k+ 2k+\frac{3}{2}+j \right) }{\left( i_k+2k+\frac{5}{2} \right) \left( i_k+2k+2 \right)} \right. \nonumber\\
&&\times \left. \frac{\left( 2k+\frac{1}{2}-\varphi (h)  \right)_{i_k} \left( 2k+\frac{1}{2}+\varphi (h) \right)_{i_k}\left( 2k+\frac{3}{2} \right)_{i_{k-1}} \left( 2k+1 \right)_{i_{k-1}}}{\left( 2k+\frac{1}{2}-\varphi (h) \right)_{i_{k-1}} \left( 2k+\frac{1}{2}+\varphi (h) \right)_{i_{k-1}}\left( 2k+\frac{3}{2} \right)_{i_k} \left( 2k+1 \right)_{i_k}} \right) \nonumber\\ 
&&\times \frac{\left( 2\tau +\frac{1}{2}-\varphi (h) \right)_n \left( 2\tau +\frac{1}{2}+\varphi (h) \right)_n \left( 2\tau +\frac{3}{2} \right)_{i_{\tau -1}} \left( 2\tau +1 \right)_{i_{\tau -1}}}{\left( 2\tau +\frac{1}{2}-\varphi (h) \right)_{i_{\tau -1}} \left( 2\tau +\frac{1}{2}+\varphi (h) \right)_{i_{\tau -1}}\left( 2\tau +\frac{3}{2} \right)_n \left( 2\tau +1 \right)_n} \left( 1+\rho ^2 \right)^{n } \hspace{1.5cm}\label{eq:100078c}
\end{eqnarray}
\begin{eqnarray}
y_0^m(j,h ;\xi ) &=& \xi^{\frac{1}{2}}  \sum_{i_0=0}^{m} \frac{\left( \frac{1}{2}-\varphi (h) \right)_{i_0}\left( \frac{1}{2}+ \varphi (h) \right)_{i_0}}{\left( \frac{3}{2} \right)_{i_0} \left( 1\right)_{i_0}} \mu ^{i_0} \label{eq:100079a}\\
y_1^m(j,h ;\xi ) &=& \xi^{\frac{1}{2}} \left\{\sum_{i_0=0}^{m} \frac{\left( i_0-j\right) \left( i_0+\frac{3}{2}+j \right)}{\left( i_0+ \frac{5}{2} \right) \left( i_0+ 2 \right)} \frac{\left( \frac{1}{2}-\varphi (h) \right)_{i_0}\left( \frac{1}{2}+ \varphi (h) \right)_{i_0}}{\left( \frac{3}{2} \right)_{i_0} \left( 1 \right)_{i_0}} \right. \nonumber\\
&&\times \left. \sum_{i_1 = i_0}^{m} \frac{\left( \frac{5}{2} -\varphi (h)  \right)_{i_1}\left( \frac{5}{2} +\varphi (h) \right)_{i_1} \left( \frac{7}{2} \right)_{i_0} \left( 3 \right)_{i_0}}{\left( \frac{5}{2} -\varphi (h) \right)_{i_0}\left( \frac{5}{2}+\varphi (h) \right)_{i_0} \left( \frac{7}{2} \right)_{i_1} \left( 3 \right)_{i_1}} \mu ^{i_1}\right\} \eta 
\label{eq:100079b}
\end{eqnarray}
\begin{align}
y_{\tau }^m(j,h ;\xi ) &=  \xi^{\frac{1}{2}} \left\{ \sum_{i_0=0}^{m} \frac{\left( i_0-j\right) \left( i_0+\frac{3}{2}+j \right)}{\left( i_0+ \frac{5}{2} \right) \left( i_0+ 2 \right)} \frac{\left( \frac{1}{2}-\varphi (h) \right)_{i_0}\left( \frac{1}{2}+ \varphi (h)  \right)_{i_0}}{\left( \frac{3}{2} \right)_{i_0} \left( 1 \right)_{i_0}} \right.\nonumber\\
 &\times \prod_{k=1}^{\tau -1} \left( \sum_{i_k = i_{k-1}}^{m} \frac{\left( i_k+ 2k-j\right)\left( i_k+ 2k+\frac{3}{2}+j \right) }{\left( i_k+2k+\frac{5}{2} \right) \left( i_k+2k+2 \right)} \right. \nonumber\\
 &\times \left. \frac{\left( 2k+\frac{1}{2}-\varphi (h) \right)_{i_k} \left( 2k+\frac{1}{2}+\varphi (h) \right)_{i_k}\left( 2k+\frac{3}{2} \right)_{i_{k-1}} \left( 2k+1 \right)_{i_{k-1}}}{\left( 2k+\frac{1}{2}-\varphi (h) \right)_{i_{k-1}} \left( 2k+\frac{1}{2}+\varphi (h) \right)_{i_{k-1}}\left( 2k+\frac{3}{2} \right)_{i_k} \left( 2k+1 \right)_{i_k}} \right) \nonumber\\
 &\times \left. \sum_{i_{\tau } = i_{\tau -1}}^{m}  \frac{\left( 2\tau +\frac{1}{2}-\varphi (h)  \right)_{i_{\tau }} \left( 2\tau +\frac{1}{2}+\varphi (h)  \right)_{i_{\tau }} \left( 2\tau +\frac{3}{2} \right)_{i_{\tau -1}} \left( 2\tau +1 \right)_{i_{\tau -1}}}{\left( 2\tau +\frac{1}{2}-\varphi (h)  \right)_{i_{\tau -1}} \left( 2\tau +\frac{1}{2}+\varphi (h)  \right)_{i_{\tau -1}}\left( 2\tau +\frac{3}{2} \right)_{i_{\tau }} \left( 2\tau +1 \right)_{i_{\tau }}} \mu ^{i_{\tau }}\right\} \eta ^{\tau }   \label{eq:100079c}
\end{align}
where
\begin{equation}
\begin{cases} \tau \geq 2 \cr
\varphi (h) = \sqrt{\frac{h}{ 4(1+\rho ^2)}}
\end{cases}\nonumber
\end{equation}
\end{enumerate} 
\end{remark}    
\section{Summary}

All possible formal solutions in series of Lam\'{e} equation in the algebraic form around $x=a$ are derived analytically such as a polynomial of type 1, a polynomial of type 2, the first species complete polynomial and an infinite series in Ref.\cite{10Chou2012f}, chapter 8 of Ref.\cite{10Choun2013} and this chapter. Asymptotic expansions in closed forms of Lam\'{e} equation are obtained including the radius of convergence \cite{10Chou2012f}. This boundary condition is quiet important from a mathematical point if views; if the condition of convergence is not required for general solutions in series for an infinite series, we can not obtain appropriate and explicit numerical values using machine calculations.\footnote{We can derive asymptotic series solutions and boundary conditions for polynomials of type 1 and type 2 by applying a theorem in chapter 3.} 

For a polynomial of type 1 of Lam\'{e} equation, I treat an exponent parameter $q$ as a free variable and $\alpha $ as a fixed value \cite{10Chou2012f}. Its power series solution in closed forms are constructed by applying 3TRF. I suggest that $\alpha $ is $ -2 (2\alpha _i+i +\lambda )-1$ or $ 2 (2\alpha _i+i +\lambda )$ where $i,\alpha _i \in \mathbb{N}_{0}$; $\lambda =0$ for the first independent solution and $\lambda =1/2$ for the second one.

For a polynomial of type 2, I treat a parameter $\alpha $ as a free variable and $q$ as a fixed value in chapter 8 of Ref.\cite{10Choun2013}. Its formal series solutions in compact forms are constructed by applying R3TRF. I suggest that $q $ is $ \alpha (\alpha +1)a- 4(2a-b-c)(q_j+2j+\lambda )^2$ where $j,q_j\in \mathbb{N}_{0}$.

For the first species complete polynomials, I treat parameters $\alpha $ and $q$ as fixed values in this chapter. Its summation solutions in series are obtained by applying the first species complete polynomial using either 3TRF or R3TRF. I allow that $\alpha $ is $ -2\left( j+1/2+\lambda \right) $ or $  2\left( j+\lambda \right) $. $c_{j+1}=0$ is a algebraic equation of degree $j+1$ for the determination of a spectral parameter $q$ and has $j+1$ zeros denoted them by $q_j^m$ eigenvalues where $m = 0,1,2, \cdots, j$. 

An infinite series solution of Lam\'{e} equation can be constructed by applying either 3TRF or R3TRF in Ref.\cite{10Chou2012f} and chapter 8 of Ref.\cite{10Choun2013}. Two infinite series solutions using different general summation formulas are equivalent to each other for numerical calculations.   
For the difference in the aspect of mathematical structures between these two series solutions, $A_n$ in sequences $c_n$ is the leading term in each of sub-power series of a general series solution for an infinite series by applying 3TRF:  And $B_n$ is the leading term in each of sub-formal series of a solution in series for an infinite series by applying R3TRF.

Also, in the aspect of numerical computations, summation solutions in closed forms of Lam\'{e} equation for the first species complete polynomials using 3TRF and R3TRF are tantamount to each other in this chapter. And there are several difference in general summation structures for these polynomials using 3TRF and R3TRF: 
For power series solutions of Lam\'{e} equation for the first species complete polynomial using 3TRF, (1) $A_n$ is the leading term in each of finite sub-power series in a general series solution, (2) Frobenius solutions contain the sum of two finite sub-power series in its general solution and (3) the denominators and numerators in all $B_n$ terms of each finite sub-power series arise with Pochhammer symbols;
for formal solutions in series of Lam\'{e} equation for the first species complete polynomial using R3TRF, (1) $B_n$ is the leading term in each of finite sub-power series in a general series solution, (2) formal series solutions only consist of one finite sub-formal series in its general power series solution and (3) the denominators and numerators in all $A_n$ terms of each finite sub-power series arise with Pochhammer symbols.
 
Similarly, for numerical values, algebraic equations $c_{j+1}=0$ of Lam\'{e} equation for the determination of a parameter $q$ for the first species complete polynomials using 3TRF and R3TRF are equivalent to each other in this chapter. Difference in sequences structures between two algebraic equations is that $A_n$ is the leading term of each sequences $\bar{c}(l,n;j,q )$ where $l \in \mathbb{N}_{0}$ in a polynomial equation of $q$ for the first species complete polynomial using 3TRF: And $B_n$ is the leading term of each sequences $\bar{c}(l,n;j,q )$ for the first species complete polynomial using R3TRF. 
 
In Ref.\cite{10Chou2012g}, chapter 9 of Ref.\cite{10Choun2013} and this chapter, all available power series solutions of Lam\'{e} equation in Weierstrass's form around $\xi =0$ are constructed by changing coefficients $ a \rightarrow 0 $ $b \rightarrow  1$, $ c \rightarrow  \rho ^{-2}$, $ q \rightarrow  h \rho ^{-2} $ and $ x \rightarrow \xi = sn^2(z,\rho )$ in formal series solutions of Lam\'{e} equation in the algebraic form around $x=a$ for a polynomial of type 1, a polynomial of type 2, the first species complete polynomial and an infinite series. 
  
\begin{appendices}
\section*{Appendix. Conversion from 4 out of 192 local solutions of Heun equation to 4 local solutions of Lam\'{e} equation in Weierstrass's form for the first species complete polynomials}
\addtocontents{toc}{\protect\setcounter{tocdepth}{1}}
A machine-generated list of 192 (isomorphic to the Coxeter group of the Coxeter diagram $D_4$) local solutions of the Heun equation was obtained by Robert S. Maier(2007) \cite{10Maie2007}.

Lame equation in Weierstrass's form is a special case of Heun's equation.
As we compare (\ref{eq:10002}) with (\ref{eq:100056}), all coefficients in Heun's equation are correspondent to the following way.
\begin{equation}
\begin{split}
& \gamma ,\delta ,\epsilon  \longleftrightarrow   \frac{1}{2} \\ & a\longleftrightarrow  \rho ^{-2} \\ & \alpha  \longleftrightarrow \frac{1}{2}(\alpha +1) \\
& \beta   \longleftrightarrow -\frac{1}{2} \alpha \\
& q \longleftrightarrow  -\frac{1}{4}h \rho ^{-2} \\ & x \longleftrightarrow \xi = sn^2(z,\rho ) 
\end{split}\label{eq:100080}   
\end{equation} 
 
In appendix of Ref.\cite{10Chou2012g}, by changing all coefficients and independent variables of nine examples out of 192 local solutions of Heun function in Table 2 \cite{10Maie2007} into a formal series solution and an integral of Heun equation for the first kind by applying 3TRF in Ref.\cite{10Chou2012c,10Chou2012d}, 9 local solutions of Lam\'{e} equation in Weierstrass's form by applying 3TRF are constructed by substituting (\ref{eq:100080}): Formal solutions in series of Lam\'{e} equation and their combined definite \& contour integrals are obtained for an infinite series and a polynomial of type 1.

In appendix of chapters 9 of Ref.\cite{10Choun2013}, I change all coefficients and independent variables of the previous nine examples of 192 local solutions of Heun function in Table 2 \cite{10Maie2007} into (1) a Frobenius solution, (2) an integral and (3) a generating function of Heun equation for the first kind by applying R3TRF in chapters 2 and 3 of Ref.\cite{10Choun2013}.  And for an infinite series and a polynomial of type 2, power series solutions and integrals for 9 local solutions of Lam\'{e} equation in Weierstrass's form are constructed analytically by taking (\ref{eq:100080}), including their generating functions for Lam\'{e} polynomial of type 2.    

In this appendix, by changing all coefficients and independent variables of four examples out of 192 local solutions of Heun function \cite{10Maie2007} into power series solutions of Heun equation for the first species complete polynomials using 3TRF and R3TRF in chapter 5 and 6, I construct 4 local solutions in series of Lam\'{e} equation in Weierstrass's form for the first species complete polynomials by applying (\ref{eq:100080}).\footnote{In this appendix, I treat $\alpha $ and $h$ as fixed values to construct polynomials of type 2 for 4 out of 192 local solutions of Heun function \cite{10Maie2007}. And an independent variable $sn^2(z,\rho )$ is denoted by $\xi$.}
\section{The first species complete polynomials using 3TRF}

In chapter 5, the power series expansion of Heun equation of the first kind for the first species complete polynomial using 3TRF about $x=0$ is given by
\begin{enumerate} 
\item As $\alpha =0$ and $q=q_0^0=0$,

The eigenfunction is given by
\begin{equation}
y(x) = H_pF_{0,0} \left( a, q=q_0^0=0; \alpha =0, \beta, \gamma, \delta ; \eta = \frac{(1+a)}{a} x ; z= -\frac{1}{a} x^2 \right) =1 \label{eq:100081}
\end{equation}
\item As $\alpha =-2N-1$ where $N \in \mathbb{N}_{0}$,

An algebraic equation of degree $2N+2$ for the determination of $q$ is given by
\begin{equation}
0 = \sum_{r=0}^{N+1}\bar{c}\left( 2r, N+1-r; 2N+1,q\right)\label{eq:100082a}
\end{equation}
The eigenvalue of $q$ is written by $q_{2N+1}^m$ where $m = 0,1,2,\cdots,2N+1 $; $q_{2N+1}^0 < q_{2N+1}^1 < \cdots < q_{2N+1}^{2N+1}$. Its eigenfunction is given by
\begin{eqnarray} 
y(x) &=& H_pF_{2N+1,m} \left( a, q=q_{2N+1}^m; \alpha =-2N-1, \beta, \gamma, \delta ; \eta = \frac{(1+a)}{a} x ; z= -\frac{1}{a} x^2 \right)\nonumber\\
&=& \sum_{r=0}^{N} y_{2r}^{N-r}\left( 2N+1,q_{2N+1}^m;x\right)+ \sum_{r=0}^{N} y_{2r+1}^{N-r}\left( 2N+1,q_{2N+1}^m;x\right)  
\label{eq:100082b}
\end{eqnarray}
\item As $\alpha =-2N-2$ where $N \in \mathbb{N}_{0}$,

An algebraic equation of degree $2N+3$ for the determination of $q$ is given by
\begin{eqnarray}
0  = \sum_{r=0}^{N+1}\bar{c}\left( 2r+1, N+1-r; 2N+2,q\right)\label{eq:100083a}
\end{eqnarray}
The eigenvalue of $q$ is written by $q_{2N+2}^m$ where $m = 0,1,2,\cdots,2N+2 $; $q_{2N+2}^0 < q_{2N+2}^1 < \cdots < q_{2N+2}^{2N+2}$. Its eigenfunction is given by
\begin{eqnarray} 
y(x) &=& H_pF_{2N+2,m} \left( a, q=q_{2N+2}^m; \alpha =-2N-2, \beta, \gamma, \delta ; \eta = \frac{(1+a)}{a} x ; z= -\frac{1}{a} x^2 \right)\nonumber\\
&=& \sum_{r=0}^{N+1} y_{2r}^{N+1-r}\left( 2N+2,q_{2N+2}^m;x\right) + \sum_{r=0}^{N} y_{2r+1}^{N-r}\left( 2N+2,q_{2N+2}^m;x\right) 
\label{eq:100083b}
\end{eqnarray}
In the above,
\begin{eqnarray}
\bar{c}(0,n;j,q)  &=& \frac{\left( -\frac{j}{2}\right)_{n} \left( \frac{\beta }{2} \right)_{n}}{\left( 1 \right)_{n} \left( \frac{\gamma }{2}+ \frac{1}{2} \right)_{n}} \left(-\frac{1}{a} \right)^{n}\label{eq:100084a}\\
\bar{c}(1,n;j,q) &=& \left( \frac{1+a}{a}\right) \sum_{i_0=0}^{n}\frac{ i_0 \left(i_0 +\Gamma _0^{(F)} \right)+\frac{q}{4(1+a)}}{\left( i_0+\frac{1}{2} \right) \left( i_0+\frac{\gamma }{2} \right)} \frac{\left( -\frac{j}{2}\right)_{i_0} \left( \frac{\beta }{2} \right)_{i_0}}{\left( 1 \right)_{i_0} \left( \frac{\gamma }{2}+ \frac{1}{2} \right)_{i_0}}\nonumber\\
&\times&  \frac{\left( -\frac{j}{2}+\frac{1}{2} \right)_{n} \left( \frac{\beta }{2}+ \frac{1}{2} \right)_{n}\left( \frac{3}{2} \right)_{i_0} \left( \frac{\gamma }{2}+1 \right)_{i_0}}{\left( -\frac{j}{2}+\frac{1}{2}\right)_{i_0} \left( \frac{\beta }{2}+ \frac{1}{2} \right)_{i_0}\left( \frac{3}{2} \right)_{n} \left( \frac{\gamma }{2}+1 \right)_{n}} \left(-\frac{1}{a} \right)^{n}  
\label{eq:100084b}\\
\bar{c}(\tau ,n;j,q) &=& \left( \frac{1+a}{a}\right)^{\tau} \sum_{i_0=0}^{n}\frac{ i_0 \left( i_0 +\Gamma _0^{(F)} \right)+\frac{q}{4(1+a)}}{\left( i_0+\frac{1}{2} \right) \left( i_0+\frac{\gamma }{2} \right)} \frac{\left( -\frac{j}{2}\right)_{i_0} \left( \frac{\beta }{2} \right)_{i_0}}{\left( 1 \right)_{i_0} \left( \frac{\gamma }{2}+ \frac{1}{2} \right)_{i_0}} \nonumber\\
&\times& \prod_{k=1}^{\tau -1} \left( \sum_{i_k = i_{k-1}}^{n} \frac{\left( i_k+ \frac{k}{2} \right)\left(i_k +\Gamma _k^{(F)} \right)+\frac{q}{4(1+a)}}{\left( i_k+\frac{k}{2}+\frac{1}{2} \right) \left( i_k+\frac{k}{2}+\frac{\gamma }{2} \right)} \right.\nonumber\\
&\times&   \left. \frac{\left( -\frac{j}{2}+\frac{k}{2}\right)_{i_k} \left( \frac{\beta }{2}+ \frac{k}{2} \right)_{i_k}\left( \frac{k}{2}+1 \right)_{i_{k-1}} \left( \frac{\gamma }{2}+ \frac{1}{2}+ \frac{k}{2} \right)_{i_{k-1}}}{\left( -\frac{j}{2}+\frac{k}{2}\right)_{i_{k-1}} \left( \frac{\beta }{2}+ \frac{k}{2} \right)_{i_{k-1}}\left( \frac{k}{2} +1\right)_{i_k} \left( \frac{\gamma }{2}+ \frac{1}{2}+ \frac{k}{2} \right)_{i_k}} \right) \nonumber\\
&\times& \frac{\left( -\frac{j}{2}+\frac{\tau }{2}\right)_{n} \left( \frac{\beta }{2}+ \frac{\tau }{2} \right)_{n}\left( \frac{\tau }{2} +1\right)_{i_{\tau -1}} \left( \frac{\gamma }{2}+ \frac{1}{2}+\frac{\tau }{2} \right)_{i_{\tau -1}}}{\left( -\frac{j}{2}+\frac{\tau }{2}\right)_{i_{\tau -1}} \left( \frac{\beta }{2}+\frac{\tau }{2} \right)_{i_{\tau -1}}\left( \frac{\tau }{2}+1 \right)_{n} \left( \frac{\gamma }{2}+ \frac{1}{2}+\frac{\tau }{2} \right)_{n}} \left(-\frac{1}{a} \right)^{n } \label{eq:100084c}
\end{eqnarray}
\begin{eqnarray}
y_0^m(j,q;x) &=& \sum_{i_0=0}^{m} \frac{\left( -\frac{j}{2}\right)_{i_0} \left( \frac{\beta }{2} \right)_{i_0}}{\left( 1 \right)_{i_0} \left( \frac{\gamma }{2}+ \frac{1}{2} \right)_{i_0}} z^{i_0} \label{eq:100085a}\\
y_1^m(j,q;x) &=& \left\{\sum_{i_0=0}^{m} \frac{ i_0 \left(i_0 +\Gamma _0^{(F)} \right)+\frac{q}{4(1+a)}}{\left( i_0+\frac{1}{2} \right) \left( i_0+\frac{\gamma }{2} \right)} \frac{\left( -\frac{j}{2}\right)_{i_0} \left( \frac{\beta }{2} \right)_{i_0}}{\left( 1 \right)_{i_0} \left( \frac{\gamma }{2}+ \frac{1}{2} \right)_{i_0}} \right. \nonumber\\
&\times& \left. \sum_{i_1 = i_0}^{m} \frac{\left( -\frac{j}{2}+\frac{1}{2} \right)_{i_1} \left( \frac{\beta }{2}+ \frac{1}{2} \right)_{i_1}\left( \frac{3}{2} \right)_{i_0} \left( \frac{\gamma }{2}+1 \right)_{i_0}}{\left(-\frac{j}{2}+\frac{1}{2} \right)_{i_0} \left( \frac{\beta }{2}+ \frac{1}{2} \right)_{i_0}\left( \frac{3}{2} \right)_{i_1} \left( \frac{\gamma }{2} +1 \right)_{i_1}} z^{i_1}\right\} \eta 
\label{eq:100085b}\\
y_{\tau }^m(j,q;x) &=& \left\{ \sum_{i_0=0}^{m} \frac{ i_0 \left(i_0 +\Gamma _0^{(F)} \right)+\frac{q}{4(1+a)}}{\left( i_0+\frac{1}{2} \right) \left( i_0+\frac{\gamma }{2} \right)} \frac{\left( -\frac{j}{2}\right)_{i_0} \left( \frac{\beta }{2} \right)_{i_0}}{\left( 1 \right)_{i_0} \left( \frac{\gamma }{2}+ \frac{1}{2} \right)_{i_0}} \right.\nonumber\\
&\times& \prod_{k=1}^{\tau -1} \left( \sum_{i_k = i_{k-1}}^{m} \frac{\left( i_k+ \frac{k}{2} \right)\left(i_k +\Gamma _k^{(F)} \right)+\frac{q}{4(1+a)}}{\left( i_k+\frac{k}{2}+\frac{1}{2} \right) \left( i_k+\frac{k}{2}+\frac{\gamma }{2} \right)} \right. \nonumber\\
&\times&  \left. \frac{\left( -\frac{j}{2}+\frac{k}{2}\right)_{i_k} \left( \frac{\beta }{2}+ \frac{k}{2} \right)_{i_k}\left( \frac{k}{2}+1 \right)_{i_{k-1}} \left( \frac{\gamma }{2}+ \frac{1}{2}+ \frac{k}{2} \right)_{i_{k-1}}}{\left( -\frac{j}{2}+\frac{k}{2}\right)_{i_{k-1}} \left( \frac{\beta }{2}+ \frac{k}{2} \right)_{i_{k-1}}\left( \frac{k}{2} +1\right)_{i_k} \left( \frac{\gamma }{2}+ \frac{1}{2}+ \frac{k}{2} \right)_{i_k}} \right) \nonumber\\
&\times& \left. \sum_{i_{\tau } = i_{\tau -1}}^{m} \frac{\left( -\frac{j}{2}+\frac{\tau }{2}\right)_{i_{\tau }} \left( \frac{\beta }{2}+ \frac{\tau }{2} \right)_{i_{\tau }}\left( \frac{\tau }{2}+1 \right)_{i_{\tau -1}} \left( \frac{\gamma }{2}+ \frac{1}{2}+\frac{\tau }{2} \right)_{i_{\tau -1}}}{\left( -\frac{j}{2}+\frac{\tau }{2}\right)_{i_{\tau -1}} \left( \frac{\beta }{2}+ \frac{\tau }{2} \right)_{i_{\tau -1}}\left( \frac{\tau }{2}+1 \right)_{i_{\tau }} \left( \frac{\gamma }{2}+ \frac{1}{2}+\frac{\tau }{2} \right)_{i_{\tau }}} z^{i_{\tau }}\right\} \eta ^{\tau }\hspace{1.5cm}\label{eq:100085c}
\end{eqnarray}
where
\begin{equation}
\begin{cases} \tau \geq 2 \cr
\Gamma _0^{(F)} = \frac{1}{2(1+a)}\left( \beta -\delta -j +a\left( \gamma +\delta -1 \right)\right) \cr
\Gamma _k^{(F)} = \frac{1}{2(1+a)}\left( \beta -\delta -j+k +a\left( \gamma +\delta -1+k \right)\right)
\end{cases}\nonumber
\end{equation}
\end{enumerate}
\subsection{ ${\displaystyle (1-x)^{1-\delta } Hl(a, q - (\delta  - 1)\gamma a; \alpha - \delta  + 1, \beta - \delta + 1, \gamma ,2 - \delta ; x)}$ }
Replace coefficients $q$, $\alpha$, $\beta$ and $\delta$ by $q - (\delta - 1)\gamma a $, $\alpha - \delta  + 1 $, $\beta - \delta + 1$ and $2 - \delta$ into (\ref{eq:100081})--(\ref{eq:100085c}). Multiply (\ref{eq:100081}), (\ref{eq:100082b}) and (\ref{eq:100083b}) by $(1-x)^{1-\delta }$. Substitute (\ref{eq:100080}) into the new (\ref{eq:100081})--(\ref{eq:100085c}) with replacing both new $\alpha $ and $\beta $ by $-j$ where $-j \in \mathbb{N}_{0}$ in order to make $A_n$ and $B_n$ terms terminated at specific index summation $n$: There are two possible values for a coefficient $\alpha $ in (\ref{eq:100081}) such as $\alpha =-2$ or 1; two fixed values for $\alpha $ in (\ref{eq:100082a}) and (\ref{eq:100082b}) are given by $\alpha =-2\left( 2N+2\right)$ or $2\left( 2N+3/2\right)$; two fixed constants for $\alpha $ in (\ref{eq:100083a}) and (\ref{eq:100083b}) are described by $\alpha =-2\left( 2N+3\right)$ or $2\left( 2N+5/2\right)$. Replace $\alpha$ by $-2(j+1)$ in the new (\ref{eq:100084a})--(\ref{eq:100085c}).
\begin{enumerate} 
\item As $\alpha = -2 $ or 1 with $h=h_0^0=1 $,
 
The eigenfunction is given by
\begin{equation}
(1-\xi )^{\frac{1}{2}} y(\xi ) = (1-\xi )^{\frac{1}{2}} Hl\left( \rho ^{-2}, 0; 0,0,\frac{1}{2}, \frac{3}{2};\xi \right) = (1-\xi )^{\frac{1}{2}} \nonumber
\end{equation}
\item As $\alpha  =-2\left( 2N+2 \right)$ or $2\left( 2N+\frac{3}{2} \right) $ where $N \in \mathbb{N}_{0}$,

An algebraic equation of degree $2N+2$ for the determination of $h$ is given by
\begin{equation}
0 = \sum_{r=0}^{N+1}\bar{c}\left( 2r, N+1-r; 2N+1,h \right)  \nonumber
\end{equation}
The eigenvalue of $h$ is written by $h_{2N+1}^m$ where $m = 0,1,2,\cdots,2N+1 $; $h_{2N+1}^0 < h_{2N+1}^1 < \cdots < h_{2N+1}^{2N+1}$. Its eigenfunction is given by 
\begin{eqnarray} 
(1-\xi )^{\frac{1}{2}} y(\xi ) &=& (1-\xi )^{\frac{1}{2}} Hl\left( \rho ^{-2}, \frac{\rho ^{-2}}{4}\left( 1-h_{2N+1}^m \right); -\left( 2N+1\right), -\left( 2N+1\right),\frac{1}{2}, \frac{3}{2};\xi \right) \nonumber\\
&=& (1-\xi )^{\frac{1}{2}} \left\{ \sum_{r=0}^{N} y_{2r}^{N-r}\left( 2N+1,h_{2N+1}^m;\xi \right)+ \sum_{r=0}^{N} y_{2r+1}^{N-r}\left( 2N+1,h_{2N+1}^m;\xi \right) \right\}
\nonumber
\end{eqnarray}
\item As  $\alpha  =-2\left( 2N+3 \right)$ or $2\left( 2N+\frac{5}{2} \right) $ where $N \in \mathbb{N}_{0}$,

An algebraic equation of degree $2N+3$ for the determination of $h$ is given by
\begin{equation}  
0 = \sum_{r=0}^{N+1}\bar{c}\left( 2r+1, N+1-r; 2N+2,h \right) \nonumber
\end{equation}
The eigenvalue of $h$ is written by $h_{2N+2}^m$ where $m = 0,1,2,\cdots,2N+2 $; $h_{2N+2}^0 < h_{2N+2}^1 < \cdots < h_{2N+2}^{2N+2}$. Its eigenfunction is given by
\begin{eqnarray} 
(1-\xi )^{\frac{1}{2}} y(\xi ) &=& (1-\xi )^{\frac{1}{2}} Hl\left( \rho ^{-2}, \frac{\rho ^{-2}}{4}\left( 1-h_{2N+2}^m \right); -\left( 2N+2\right), -\left( 2N+2\right),\frac{1}{2}, \frac{3}{2};\xi \right) \nonumber\\
&=& (1-\xi )^{\frac{1}{2}} \left\{ \sum_{r=0}^{N+1} y_{2r}^{N+1-r}\left( 2N+2,h_{2N+2}^m;\xi \right) + \sum_{r=0}^{N} y_{2r+1}^{N-r}\left( 2N+2,h_{2N+2}^m;\xi \right) \right\} \nonumber
\end{eqnarray}
In the above,
\begin{eqnarray}
\bar{c}(0,n;j,h)  &=& \frac{\left( -\frac{j}{2}\right)_{n}\left( \frac{3}{4}+\frac{j}{2} \right)_{n}}{ \left(  1 \right)_{n}\left( \frac{3}{4} \right)_{n}} \left( -\rho ^2\right)^{n}\nonumber\\
\bar{c}(1,n;j,h) &=& \left( 1+\rho ^2\right) \sum_{i_0=0}^{n}\frac{ i_0 \left( i_0 +\Gamma _0^{(F)} \right) +Q(h)}{\left( i_0+\frac{1}{2} \right) \left( i_0+\frac{1}{4} \right)} \frac{\left( -\frac{j}{2}\right)_{i_0}\left( \frac{3}{4}+\frac{j}{2}  \right)_{i_0} }{\left( 1 \right)_{i_0} \left( \frac{3}{4} \right)_{i_0}} \nonumber\\
&\times& \frac{\left( \frac{1}{2}-\frac{j}{2} \right)_{n} \left( \frac{5}{4}+\frac{j}{2} \right)_{n} \left( \frac{3}{2} \right)_{i_0}\left( \frac{5}{4} \right)_{i_0}}{\left( \frac{1}{2}-\frac{j}{2} \right)_{i_0} \left( \frac{5}{4}+\frac{j}{2} \right)_{i_0} \left( \frac{3}{2} \right)_n \left( \frac{5}{4} \right)_n} \left( -\rho ^2\right)^{n } \nonumber\\
\bar{c}(\tau ,n;j,h) &=& \left( 1+\rho ^2\right)^{\tau } \sum_{i_0=0}^{n} \frac{ i_0 \left( i_0 +\Gamma _0^{(F)} \right) +Q(h)}{\left( i_0+\frac{1}{2} \right) \left( i_0+\frac{1}{4} \right)} \frac{\left( -\frac{j}{2}\right)_{i_0}\left( \frac{3}{4}+\frac{j}{2} \right)_{i_0} }{\left( 1 \right)_{i_0} \left( \frac{3}{4} \right)_{i_0}}  \nonumber\\
&&\times \prod_{k=1}^{\tau -1} \left( \sum_{i_k = i_{k-1}}^{n} \frac{\left( i_k+ \frac{k}{2} \right) \left( i_k +\Gamma _k^{(F)} \right) +Q(h)}{\left( i_k+\frac{k}{2}+\frac{1}{2} \right) \left( i_k+\frac{k}{2}+\frac{1}{4} \right)} \right. \nonumber\\
&\times& \left. \frac{\left( \frac{k}{2}-\frac{j}{2}\right)_{i_k} \left( \frac{k}{2}+\frac{3}{4}+\frac{j}{2} \right)_{i_k}\left( \frac{k}{2}+1 \right)_{i_{k-1}} \left( \frac{k}{2}+ \frac{3}{4} \right)_{i_{k-1}}}{\left( \frac{k}{2}-\frac{j}{2}\right)_{i_{k-1}} \left( \frac{k}{2}+\frac{3}{4}+\frac{j}{2} \right)_{i_{k-1}}\left( \frac{k}{2}+1 \right)_{i_k} \left( \frac{k}{2}+ \frac{3}{4} \right)_{i_k}} \right) \nonumber\\ 
&&\times \frac{\left( \frac{\tau }{2} -\frac{j}{2}\right)_{n}\left( \frac{\tau }{2}+\frac{3}{4} +\frac{j}{2} \right)_{n} \left( \frac{\tau }{2}+1 \right)_{i_{\tau -1}} \left( \frac{\tau }{2}+\frac{3}{4} \right)_{i_{\tau -1}}}{\left( \frac{\tau }{2} -\frac{j}{2}\right)_{i_{\tau -1}}\left( \frac{\tau }{2}+\frac{3}{4} +\frac{j}{2} \right)_{i_{\tau -1}} \left( \frac{\tau }{2}+1 \right)_n \left( \frac{\tau }{2}+\frac{3}{4} \right)_n} \left( -\rho ^2\right)^{n } \nonumber
\end{eqnarray}
\begin{eqnarray}
y_0^m(j,h ;\xi ) &=& \sum_{i_0=0}^{m} \frac{\left( -\frac{j}{2}\right)_{i_0}\left( \frac{3}{4}+\frac{j}{2} \right)_{i_0} }{\left( 1 \right)_{i_0} \left( \frac{3}{4} \right)_{i_0}} z^{i_0} \nonumber\\
y_1^m(j,h ;\xi ) &=& \left\{\sum_{i_0=0}^{m} \frac{ i_0 \left( i_0 +\Gamma _0^{(F)} \right) +Q(h)}{\left( i_0+\frac{1}{2} \right) \left( i_0+\frac{1}{4} \right)} \frac{\left( -\frac{j}{2}\right)_{i_0}\left( \frac{3}{4}+\frac{j}{2} \right)_{i_0} }{\left( 1 \right)_{i_0} \left( \frac{3}{4} \right)_{i_0}} \right.  \nonumber\\
&\times&  \left. \sum_{i_1 = i_0}^{m} \frac{\left( \frac{1}{2}-\frac{j}{2} \right)_{i_1} \left( \frac{5}{4}+\frac{j}{2} \right)_{i_1} \left( \frac{3}{2} \right)_{i_0}\left( \frac{5}{4} \right)_{i_0}}{\left( \frac{1}{2}-\frac{j}{2} \right)_{i_0} \left( \frac{5}{4}+\frac{j}{2} \right)_{i_0} \left( \frac{3}{2} \right)_{i_1} \left( \frac{5}{4} \right)_{i_1}} z^{i_1}\right\} \eta \nonumber
\end{eqnarray}
\begin{eqnarray}
y_{\tau }^m(j,h ;\xi ) &=& \left\{ \sum_{i_0=0}^{m} \frac{ i_0 \left( i_0 +\Gamma _0^{(F)} \right) +Q(h)}{\left( i_0+\frac{1}{2} \right) \left( i_0+\frac{1}{4} \right)} \frac{\left( -\frac{j}{2}\right)_{i_0}\left( \frac{3}{4}+\frac{j}{2} \right)_{i_0} }{\left( 1 \right)_{i_0} \left( \frac{3}{4} \right)_{i_0}} \right.\nonumber\\
&&\times \prod_{k=1}^{\tau -1} \left( \sum_{i_k = i_{k-1}}^{m} \frac{\left( i_k+ \frac{k}{2} \right) \left( i_k +\Gamma _k^{(F)} \right) +Q(h)}{\left( i_k+\frac{k}{2}+\frac{1}{2} \right) \left( i_k+\frac{k}{2}+\frac{1}{4} \right)} \right. \nonumber\\
&\times&   \left. \frac{\left( \frac{k}{2}-\frac{j}{2}\right)_{i_k} \left( \frac{k}{2}+\frac{3}{4}+\frac{j}{2} \right)_{i_k}\left( \frac{k}{2}+1 \right)_{i_{k-1}} \left( \frac{k}{2}+ \frac{3}{4} \right)_{i_{k-1}}}{\left( \frac{k}{2}-\frac{j}{2}\right)_{i_{k-1}} \left( \frac{k}{2}+\frac{3}{4}+\frac{j}{2} \right)_{i_{k-1}}\left( \frac{k}{2}+1 \right)_{i_k} \left( \frac{k}{2}+ \frac{3}{4} \right)_{i_k}} \right) \nonumber\\
&&\times \left. \sum_{i_{\tau } = i_{\tau -1}}^{m} \frac{\left( \frac{\tau }{2} -\frac{j}{2}\right)_{i_{\tau }}\left( \frac{\tau }{2}+\frac{3}{4} +\frac{j}{2} \right)_{i_{\tau }} \left( \frac{\tau }{2}+1 \right)_{i_{\tau -1}} \left( \frac{\tau }{2}+\frac{3}{4}  \right)_{i_{\tau -1}}}{\left( \frac{\tau }{2} -\frac{j}{2}\right)_{i_{\tau -1}}\left( \frac{\tau }{2}+\frac{3}{4} +\frac{j}{2} \right)_{i_{\tau -1}} \left( \frac{\tau }{2}+1 \right)_{i_{\tau }} \left( \frac{\tau }{2}+\frac{3}{4} \right)_{i_{\tau }}} z^{i_{\tau }}\right\} \eta ^{\tau } \nonumber
\end{eqnarray}
where
\begin{equation}
\begin{cases} \tau \geq 2 \cr
z = -\rho ^2 \xi^2 \cr
\eta = \left( 1+\rho ^2 \right) \xi \cr
\Gamma _0^{(F)} = \frac{1}{2\left( 1+\rho ^2 \right)}  \cr
\Gamma _k^{(F)} = \frac{k}{2} +\frac{1}{2\left( 1+\rho ^2 \right)} \cr
Q(h) = \frac{1-h}{16\left( 1+\rho ^2 \right)}
\end{cases}\nonumber
\end{equation}
\end{enumerate}
\subsection{\footnotesize ${\displaystyle x^{1-\gamma } (1-x)^{1-\delta }Hl(a, q-(\gamma +\delta -2)a -(\gamma -1)(\alpha +\beta -\gamma -\delta +1), \alpha - \gamma -\delta +2, \beta - \gamma -\delta +2, 2-\gamma, 2 - \delta ; x)}$ \normalsize}
Replace coefficients $q$, $\alpha$, $\beta$, $\gamma $ and $\delta$ by $q-(\gamma +\delta -2)a-(\gamma -1)(\alpha +\beta -\gamma -\delta +1)$, $\alpha - \gamma -\delta +2 $, $\beta - \gamma -\delta +2$, $2 -\gamma $ and $2 - \delta$ into (\ref{eq:100081})--(\ref{eq:100085c}). Multiply (\ref{eq:100081}), (\ref{eq:100082b}) and (\ref{eq:100083b}) by $x^{1-\gamma } (1-x)^{1-\delta }$. Substitute (\ref{eq:100080}) into the new (\ref{eq:100081})--(\ref{eq:100085c}) with replacing both new $\alpha $ and $\beta $ by $-j$ where $-j \in \mathbb{N}_{0}$ in order to make $A_n$ and $B_n$ terms terminated at specific index summation $n$: There are two possible values for a coefficient $\alpha $ in (\ref{eq:100081}) such as $\alpha =-3$ or 2; two fixed values for $\alpha $ in (\ref{eq:100082a}) and (\ref{eq:100082b}) are given by $\alpha =-2\left( 2N+5/2\right)$ or $2\left( 2N+2\right)$; two fixed constants for $\alpha $ in (\ref{eq:100083a}) and (\ref{eq:100083b}) are described by $\alpha =-2\left( 2N+7/2\right)$ or $2\left( 2N+3\right)$. Replace $\alpha$ by $-2(j+3/2)$ in the new (\ref{eq:100084a})--(\ref{eq:100085c}).
\begin{enumerate} 
\item As $\alpha = -3 $ or 2 with $h=h_0^0= 4+ \rho ^2$,
 
The eigenfunction is given by
\begin{equation}
\xi ^{\frac{1}{2}}(1-\xi )^{\frac{1}{2}} y(\xi ) = \xi ^{\frac{1}{2}}(1-\xi )^{\frac{1}{2}} Hl\left( \rho ^{-2}, 0; 0,0,\frac{3}{2}, \frac{3}{2};\xi \right) = \xi ^{\frac{1}{2}} (1-\xi )^{\frac{1}{2}} \nonumber
\end{equation}
\item As $\alpha  =-2\left( 2N+\frac{5}{2} \right)$ or $2\left( 2N+2 \right) $ where $N \in \mathbb{N}_{0}$,

An algebraic equation of degree $2N+2$ for the determination of $h$ is given by
\begin{equation}
0 = \sum_{r=0}^{N+1}\bar{c}\left( 2r, N+1-r; 2N+1,h \right)  \nonumber
\end{equation}
The eigenvalue of $h$ is written by $h_{2N+1}^m$ where $m = 0,1,2,\cdots,2N+1 $; $h_{2N+1}^0 < h_{2N+1}^1 < \cdots < h_{2N+1}^{2N+1}$. Its eigenfunction is given by 
\begin{eqnarray} 
\xi ^{\frac{1}{2}} (1-\xi )^{\frac{1}{2}} y(\xi ) &=& \xi ^{\frac{1}{2}} (1-\xi )^{\frac{1}{2}} Hl\left( \rho ^{-2}, \frac{1}{4}\left( 1-(h_{2N+1}^m-4)\rho ^{-2} \right); -\left( 2N+1\right), -\left( 2N+1\right),\frac{3}{2}, \frac{3}{2};\xi \right) \nonumber\\
&=& \xi ^{\frac{1}{2}} (1-\xi )^{\frac{1}{2}} \left\{ \sum_{r=0}^{N} y_{2r}^{N-r}\left( 2N+1,h_{2N+1}^m;\xi \right) + \sum_{r=0}^{N} y_{2r+1}^{N-r}\left( 2N+1,h_{2N+1}^m;\xi \right) \right\}
\nonumber
\end{eqnarray}
\item As  $\alpha  =-2\left( 2N+\frac{7}{2} \right)$ or $2\left( 2N+3 \right) $ where $N \in \mathbb{N}_{0}$,

An algebraic equation of degree $2N+3$ for the determination of $h$ is given by
\begin{equation}  
0 = \sum_{r=0}^{N+1}\bar{c}\left( 2r+1, N+1-r; 2N+2,h \right) \nonumber
\end{equation}
The eigenvalue of $h$ is written by $h_{2N+2}^m$ where $m = 0,1,2,\cdots,2N+2 $; $h_{2N+2}^0 < h_{2N+2}^1 < \cdots < h_{2N+2}^{2N+2}$. Its eigenfunction is given by
\begin{eqnarray} 
\xi ^{\frac{1}{2}} (1-\xi )^{\frac{1}{2}} y(\xi ) &=&\xi ^{\frac{1}{2}} (1-\xi )^{\frac{1}{2}} Hl\left( \rho ^{-2}, \frac{1}{4}\left( 1-(h_{2N+2}^m-4)\rho ^{-2} \right); -\left( 2N+2\right), -\left( 2N+2\right),\frac{3}{2}, \frac{3}{2};\xi \right) \nonumber\\
&=&\xi ^{\frac{1}{2}} (1-\xi )^{\frac{1}{2}} \left\{ \sum_{r=0}^{N+1} y_{2r}^{N+1-r}\left( 2N+2,h_{2N+2}^m;\xi \right) + \sum_{r=0}^{N} y_{2r+1}^{N-r}\left( 2N+2,h_{2N+2}^m;\xi \right) \right\} \nonumber
\end{eqnarray}
In the above,
\begin{eqnarray}
\bar{c}(0,n;j,h)  &=& \frac{\left( -\frac{j}{2}\right)_{n}\left( \frac{5}{4}+\frac{j}{2} \right)_{n}}{ \left(  1 \right)_{n}\left( \frac{5}{4} \right)_{n}} \left( -\rho ^2\right)^{n}\nonumber\\
\bar{c}(1,n;j,h) &=& \left( 1+\rho ^2\right) \sum_{i_0=0}^{n}\frac{ i_0 \left( i_0 +\Gamma _0^{(F)} \right) +Q(h)}{\left( i_0+\frac{1}{2} \right) \left( i_0+\frac{3}{4} \right)} \frac{\left( -\frac{j}{2}\right)_{i_0}\left( \frac{5}{4}+\frac{j}{2}  \right)_{i_0} }{\left( 1 \right)_{i_0} \left( \frac{5}{4} \right)_{i_0}} \nonumber\\
&&\times\frac{\left( \frac{1}{2}-\frac{j}{2} \right)_{n} \left( \frac{7}{4}+\frac{j}{2} \right)_{n} \left( \frac{3}{2} \right)_{i_0}\left( \frac{7}{4} \right)_{i_0}}{\left( \frac{1}{2}-\frac{j}{2} \right)_{i_0} \left( \frac{7}{4}+\frac{j}{2} \right)_{i_0} \left( \frac{3}{2} \right)_n \left( \frac{7}{4} \right)_n} \left( -\rho ^2\right)^{n } \nonumber\\
\bar{c}(\tau ,n;j,h) &=& \left( 1+\rho ^2\right)^{\tau } \sum_{i_0=0}^{n} \frac{ i_0 \left( i_0 +\Gamma _0^{(F)} \right) +Q(h)}{\left( i_0+\frac{1}{2} \right) \left( i_0+\frac{3}{4} \right)} \frac{\left( -\frac{j}{2}\right)_{i_0}\left( \frac{5}{4}+\frac{j}{2} \right)_{i_0} }{\left( 1 \right)_{i_0} \left( \frac{5}{4} \right)_{i_0}}  \nonumber\\
&&\times \prod_{k=1}^{\tau -1} \left( \sum_{i_k = i_{k-1}}^{n} \frac{\left( i_k+ \frac{k}{2} \right) \left( i_k +\Gamma _k^{(F)} \right) +Q(h)}{\left( i_k+\frac{k}{2}+\frac{1}{2} \right) \left( i_k+\frac{k}{2}+\frac{3}{4} \right)} \right. \nonumber\\
&&\times\left. \frac{\left( \frac{k}{2}-\frac{j}{2}\right)_{i_k} \left( \frac{k}{2}+\frac{5}{4}+\frac{j}{2} \right)_{i_k}\left( \frac{k}{2}+1 \right)_{i_{k-1}} \left( \frac{k}{2}+ \frac{5}{4} \right)_{i_{k-1}}}{\left( \frac{k}{2}-\frac{j}{2}\right)_{i_{k-1}} \left( \frac{k}{2}+\frac{5}{4}+\frac{j}{2} \right)_{i_{k-1}}\left( \frac{k}{2}+1 \right)_{i_k} \left( \frac{k}{2}+ \frac{5}{4} \right)_{i_k}} \right) \nonumber\\ 
&&\times \frac{\left( \frac{\tau }{2} -\frac{j}{2}\right)_{n}\left( \frac{\tau }{2}+\frac{5}{4} +\frac{j}{2} \right)_{n} \left( \frac{\tau }{2}+1 \right)_{i_{\tau -1}} \left( \frac{\tau }{2}+\frac{5}{4} \right)_{i_{\tau -1}}}{\left( \frac{\tau }{2} -\frac{j}{2}\right)_{i_{\tau -1}}\left( \frac{\tau }{2}+\frac{5}{4} +\frac{j}{2} \right)_{i_{\tau -1}} \left( \frac{\tau }{2}+1 \right)_n \left( \frac{\tau }{2}+\frac{5}{4} \right)_n} \left( -\rho ^2\right)^{n } \nonumber
\end{eqnarray}
\begin{eqnarray}
y_0^m(j,h ;\xi ) &=& \sum_{i_0=0}^{m} \frac{\left( -\frac{j}{2}\right)_{i_0}\left( \frac{5}{4}+\frac{j}{2} \right)_{i_0} }{\left( 1 \right)_{i_0} \left( \frac{5}{4} \right)_{i_0}} z^{i_0} \nonumber\\
y_1^m(j,h ;\xi ) &=& \left\{\sum_{i_0=0}^{m} \frac{ i_0 \left( i_0 +\Gamma _0^{(F)} \right) +Q(h)}{\left( i_0+\frac{1}{2} \right) \left( i_0+\frac{3}{4} \right)} \frac{\left( -\frac{j}{2}\right)_{i_0}\left( \frac{5}{4}+\frac{j}{2} \right)_{i_0} }{\left( 1 \right)_{i_0} \left( \frac{5}{4} \right)_{i_0}} \right.  \nonumber\\
&&\times \left. \sum_{i_1 = i_0}^{m} \frac{\left( \frac{1}{2}-\frac{j}{2} \right)_{i_1} \left( \frac{7}{4}+\frac{j}{2} \right)_{i_1} \left( \frac{3}{2} \right)_{i_0}\left( \frac{7}{4} \right)_{i_0}}{\left( \frac{1}{2}-\frac{j}{2} \right)_{i_0} \left( \frac{7}{4}+\frac{j}{2} \right)_{i_0} \left( \frac{3}{2} \right)_{i_1} \left( \frac{7}{4} \right)_{i_1}} z^{i_1}\right\} \eta \nonumber
\end{eqnarray}
\begin{eqnarray}
y_{\tau }^m(j,h ;\xi ) &=& \left\{ \sum_{i_0=0}^{m} \frac{ i_0 \left( i_0 +\Gamma _0^{(F)} \right) +Q(h)}{\left( i_0+\frac{1}{2} \right) \left( i_0+\frac{3}{4} \right)} \frac{\left( -\frac{j}{2}\right)_{i_0}\left( \frac{5}{4}+\frac{j}{2} \right)_{i_0} }{\left( 1 \right)_{i_0} \left( \frac{5}{4} \right)_{i_0}} \right.\nonumber\\
&&\times \prod_{k=1}^{\tau -1} \left( \sum_{i_k = i_{k-1}}^{m} \frac{\left( i_k+ \frac{k}{2} \right) \left( i_k +\Gamma _k^{(F)} \right) +Q(h)}{\left( i_k+\frac{k}{2}+\frac{1}{2} \right) \left( i_k+\frac{k}{2}+\frac{3}{4} \right)} \right. \nonumber\\
&&\times  \left. \frac{\left( \frac{k}{2}-\frac{j}{2}\right)_{i_k} \left( \frac{k}{2}+\frac{5}{4}+\frac{j}{2} \right)_{i_k}\left( \frac{k}{2}+1 \right)_{i_{k-1}} \left( \frac{k}{2}+ \frac{5}{4} \right)_{i_{k-1}}}{\left( \frac{k}{2}-\frac{j}{2}\right)_{i_{k-1}} \left( \frac{k}{2}+\frac{5}{4}+\frac{j}{2} \right)_{i_{k-1}}\left( \frac{k}{2}+1 \right)_{i_k} \left( \frac{k}{2}+ \frac{5}{4} \right)_{i_k}} \right) \nonumber\\
&&\times \left. \sum_{i_{\tau } = i_{\tau -1}}^{m} \frac{\left( \frac{\tau }{2} -\frac{j}{2}\right)_{i_{\tau }}\left( \frac{\tau }{2}+\frac{5}{4} +\frac{j}{2} \right)_{i_{\tau }} \left( \frac{\tau }{2}+1 \right)_{i_{\tau -1}} \left( \frac{\tau }{2}+\frac{5}{4}  \right)_{i_{\tau -1}}}{\left( \frac{\tau }{2} -\frac{j}{2}\right)_{i_{\tau -1}}\left( \frac{\tau }{2}+\frac{5}{4} +\frac{j}{2} \right)_{i_{\tau -1}} \left( \frac{\tau }{2}+1 \right)_{i_{\tau }} \left( \frac{\tau }{2}+\frac{5}{4} \right)_{i_{\tau }}} z^{i_{\tau }}\right\} \eta ^{\tau } \nonumber
\end{eqnarray}
where
\begin{equation}
\begin{cases} \tau \geq 2 \cr
z = -\rho ^2 \xi^2 \cr
\eta = \left( 1+\rho ^2 \right) \xi \cr
\Gamma _0^{(F)} = \frac{\left( 2+\rho ^2 \right)}{2\left( 1+\rho ^2 \right)}  \cr
\Gamma _k^{(F)} = \frac{k}{2} +\frac{\left( 2+\rho ^2 \right)}{2\left( 1+\rho ^2 \right)} \cr
Q(h) = \frac{4+\rho ^2 -h}{16\left( 1+\rho ^2 \right)}
\end{cases}\nonumber
\end{equation}
\end{enumerate}
\subsection{ ${\displaystyle  Hl(1-a,-q+\alpha \beta; \alpha,\beta, \delta, \gamma; 1-x)}$} 
Replace coefficients $a$, $q$, $\gamma $, $\delta$ and $x$ by $1-a$, $-q+ \alpha \beta $, $\delta $, $\gamma $ and $1-x$ into (\ref{eq:100081})--(\ref{eq:100085c}). Substitute (\ref{eq:100080}) into the new (\ref{eq:100081})--(\ref{eq:100085c}) with replacing the new $\alpha $, $\beta $ and $q$ by $-j$, $-j$ and $1/4h\rho ^{-2}-j(j+1/2)$ where $-j \in \mathbb{N}_{0}$ in order to make $A_n$ and $B_n$ terms terminated at specific index summation $n$: There are two possible values for a coefficient $\alpha $ in (\ref{eq:100081}) such as $\alpha =-1$ or 0; two fixed values for $\alpha $ in (\ref{eq:100082a}) and (\ref{eq:100082b}) are given by $\alpha =-2\left( 2N+3/2\right)$ or $2\left( 2N+1\right)$; two fixed constants for $\alpha $ in (\ref{eq:100083a}) and (\ref{eq:100083b}) are described by $\alpha =-2\left( 2N+5/2\right)$ or $2\left( 2N+2\right)$. Replace $\alpha$ by $-2(j+1/2)$ in the new (\ref{eq:100084a})--(\ref{eq:100085c}).
\begin{enumerate} 
\item As $\alpha = -1$ or 0 with $h=h_0^0= 0$,
 
The eigenfunction is given by
\begin{equation}
 y(\sigma ) = Hl\left( 1-\rho ^{-2}, 0; 0,0,\frac{1}{2}, \frac{1}{2};\sigma \right) = 1 \nonumber
\end{equation}
\item As $\alpha  =-2\left( 2N+\frac{3}{2} \right)$ or $2\left( 2N+1 \right) $ where $N \in \mathbb{N}_{0}$,

An algebraic equation of degree $2N+2$ for the determination of $h$ is given by
\begin{equation}
0 = \sum_{r=0}^{N+1}\bar{c}\left( 2r, N+1-r; 2N+1,h \right)  \nonumber
\end{equation}
The eigenvalue of $h$ is written by $h_{2N+1}^m$ where $m = 0,1,2,\cdots,2N+1 $; $h_{2N+1}^0 < h_{2N+1}^1 < \cdots < h_{2N+1}^{2N+1}$. Its eigenfunction is given by 
\begin{eqnarray} 
 y(\sigma ) &=& Hl\left( 1-\rho ^{-2}, \frac{\rho ^{-2}}{4}h_{2N+1}^m-\left( 2N+1\right)\left( 2N+\frac{3}{2}\right); -\left( 2N+1\right), -\left( 2N+1\right),\frac{1}{2}, \frac{1}{2};\sigma \right) \nonumber\\
&=& \sum_{r=0}^{N} y_{2r}^{N-r}\left( 2N+1,h_{2N+1}^m;\sigma \right) + \sum_{r=0}^{N} y_{2r+1}^{N-r}\left( 2N+1,h_{2N+1}^m;\sigma \right) 
\nonumber
\end{eqnarray}
\item As  $\alpha  =-2\left( 2N+\frac{5}{2} \right)$ or $2\left( 2N+2 \right) $ where $N \in \mathbb{N}_{0}$,

An algebraic equation of degree $2N+3$ for the determination of $h$ is given by
\begin{equation}  
0 = \sum_{r=0}^{N+1}\bar{c}\left( 2r+1, N+1-r; 2N+2,h \right) \nonumber
\end{equation}
The eigenvalue of $h$ is written by $h_{2N+2}^m$ where $m = 0,1,2,\cdots,2N+2 $; $h_{2N+2}^0 < h_{2N+2}^1 < \cdots < h_{2N+2}^{2N+2}$. Its eigenfunction is given by
\begin{eqnarray} 
  y(\sigma ) &=&  Hl\left( 1-\rho ^{-2}, \frac{\rho ^{-2}}{4}h_{2N+2}^m-\left( 2N+2\right)\left( 2N+\frac{5}{2}\right); -\left( 2N+2\right), -\left( 2N+2\right),\frac{1}{2}, \frac{1}{2};\sigma \right) \nonumber\\
&=&  \sum_{r=0}^{N+1} y_{2r}^{N+1-r}\left( 2N+2,h_{2N+2}^m;\sigma \right) + \sum_{r=0}^{N} y_{2r+1}^{N-r}\left( 2N+2,h_{2N+2}^m;\sigma \right)  \nonumber
\end{eqnarray}
In the above,
\begin{eqnarray}
\bar{c}(0,n;j,h)  &=& \frac{\left( -\frac{j}{2}\right)_{n}\left( \frac{1}{4}+\frac{j}{2} \right)_{n}}{ \left( 1\right)_{n}\left( \frac{3}{4} \right)_{n}} \left( \frac{-1}{1-\rho ^{-2}}\right)^{n}\nonumber\\
\bar{c}(1,n;j,h) &=& \left( \frac{2-\rho ^{-2}}{1-\rho ^{-2}}\right) \sum_{i_0=0}^{n}\frac{ i_0^2 +Q(j,h)}{\left( i_0+\frac{1}{2} \right) \left( i_0+\frac{1}{4} \right)} \frac{\left( -\frac{j}{2}\right)_{i_0}\left( \frac{1}{4}+\frac{j}{2}  \right)_{i_0} }{\left( 1 \right)_{i_0} \left( \frac{3}{4} \right)_{i_0}} \nonumber\\
&&\times \frac{\left( \frac{1}{2}-\frac{j}{2} \right)_{n} \left( \frac{3}{4}+\frac{j}{2} \right)_{n} \left( \frac{3}{2} \right)_{i_0}\left( \frac{5}{4} \right)_{i_0}}{\left( \frac{1}{2}-\frac{j}{2} \right)_{i_0} \left( \frac{3}{4}+\frac{j}{2} \right)_{i_0} \left( \frac{3}{2} \right)_n \left( \frac{5}{4} \right)_n} \left( \frac{-1}{1-\rho ^{-2}}\right)^{n } \nonumber\\
\bar{c}(\tau ,n;j,h) &=& \left( \frac{2-\rho ^{-2}}{1-\rho ^{-2}}\right)^{\tau } \sum_{i_0=0}^{n} \frac{ i_0^2 +Q(j,h)}{\left( i_0+\frac{1}{2} \right) \left( i_0+\frac{1}{4} \right)} \frac{\left( -\frac{j}{2}\right)_{i_0}\left( \frac{1}{4}+\frac{j}{2} \right)_{i_0} }{\left( 1 \right)_{i_0} \left( \frac{3}{4} \right)_{i_0}}  \nonumber\\
&&\times \prod_{k=1}^{\tau -1} \left( \sum_{i_k = i_{k-1}}^{n} \frac{\left( i_k+ \frac{k}{2} \right)^2 +Q(j,h)}{\left( i_k+\frac{k}{2}+\frac{1}{2} \right) \left( i_k+\frac{k}{2}+\frac{1}{4} \right)} \right.\nonumber\\
&&\times  \left. \frac{\left( \frac{k}{2}-\frac{j}{2}\right)_{i_k} \left( \frac{k}{2}+\frac{1}{4}+\frac{j}{2} \right)_{i_k}\left( \frac{k}{2}+1 \right)_{i_{k-1}} \left( \frac{k}{2}+ \frac{3}{4} \right)_{i_{k-1}}}{\left( \frac{k}{2}-\frac{j}{2}\right)_{i_{k-1}} \left( \frac{k}{2}+\frac{1}{4}+\frac{j}{2} \right)_{i_{k-1}}\left( \frac{k}{2}+1 \right)_{i_k} \left( \frac{k}{2}+ \frac{3}{4} \right)_{i_k}} \right) \nonumber\\ 
&&\times \frac{\left( \frac{\tau }{2} -\frac{j}{2}\right)_{n}\left( \frac{\tau }{2}+\frac{1}{4} +\frac{j}{2} \right)_{n} \left( \frac{\tau }{2}+1 \right)_{i_{\tau -1}} \left( \frac{\tau }{2}+\frac{3}{4} \right)_{i_{\tau -1}}}{\left( \frac{\tau }{2} -\frac{j}{2}\right)_{i_{\tau -1}}\left( \frac{\tau }{2}+\frac{1}{4} +\frac{j}{2} \right)_{i_{\tau -1}} \left( \frac{\tau }{2}+1 \right)_n \left( \frac{\tau }{2}+\frac{3}{4} \right)_n} \left( \frac{-1}{1-\rho ^{-2}}\right)^{n } \nonumber
\end{eqnarray}
\begin{eqnarray}
y_0^m(j,h ;\sigma ) &=& \sum_{i_0=0}^{m} \frac{\left( -\frac{j}{2}\right)_{i_0}\left( \frac{1}{4}+\frac{j}{2} \right)_{i_0} }{\left( 1 \right)_{i_0} \left( \frac{3}{4} \right)_{i_0}} z^{i_0} \nonumber\\
y_1^m(j,h ;\sigma ) &=& \left\{\sum_{i_0=0}^{m} \frac{ i_0^2 +Q(j,h)}{\left( i_0+\frac{1}{2} \right) \left( i_0+\frac{1}{4} \right)} \frac{\left( -\frac{j}{2}\right)_{i_0}\left( \frac{1}{4}+\frac{j}{2} \right)_{i_0} }{\left( 1 \right)_{i_0} \left( \frac{3}{4} \right)_{i_0}} \right.  \nonumber\\
&&\times  \left. \sum_{i_1 = i_0}^{m} \frac{\left( \frac{1}{2}-\frac{j}{2} \right)_{i_1} \left( \frac{3}{4}+\frac{j}{2} \right)_{i_1} \left( \frac{3}{2} \right)_{i_0}\left( \frac{5}{4} \right)_{i_0}}{\left( \frac{1}{2}-\frac{j}{2} \right)_{i_0} \left( \frac{3}{4}+\frac{j}{2} \right)_{i_0} \left( \frac{3}{2} \right)_{i_1} \left( \frac{5}{4} \right)_{i_1}} z^{i_1}\right\} \eta \nonumber
\end{eqnarray}
\begin{eqnarray}
y_{\tau }^m(j,h ;\sigma ) &=& \left\{ \sum_{i_0=0}^{m} \frac{ i_0^2 +Q(j,h)}{\left( i_0+\frac{1}{2} \right) \left( i_0+\frac{1}{4} \right)} \frac{\left( -\frac{j}{2}\right)_{i_0}\left( \frac{1}{4}+\frac{j}{2} \right)_{i_0} }{\left( 1 \right)_{i_0} \left( \frac{3}{4} \right)_{i_0}} \right.\nonumber\\
&&\times \prod_{k=1}^{\tau -1} \left( \sum_{i_k = i_{k-1}}^{m} \frac{\left( i_k+ \frac{k}{2} \right)^2 +Q(j,h)}{\left( i_k+\frac{k}{2}+\frac{1}{2} \right) \left( i_k+\frac{k}{2}+\frac{1}{4} \right)} \right. \nonumber\\
&&\times   \left. \frac{\left( \frac{k}{2}-\frac{j}{2}\right)_{i_k} \left( \frac{k}{2}+\frac{1}{4}+\frac{j}{2} \right)_{i_k}\left( \frac{k}{2}+1 \right)_{i_{k-1}} \left( \frac{k}{2}+ \frac{3}{4} \right)_{i_{k-1}}}{\left( \frac{k}{2}-\frac{j}{2}\right)_{i_{k-1}} \left( \frac{k}{2}+\frac{1}{4}+\frac{j}{2} \right)_{i_{k-1}}\left( \frac{k}{2}+1 \right)_{i_k} \left( \frac{k}{2}+ \frac{3}{4} \right)_{i_k}} \right) \nonumber\\
&&\times \left. \sum_{i_{\tau } = i_{\tau -1}}^{m} \frac{\left( \frac{\tau }{2} -\frac{j}{2}\right)_{i_{\tau }}\left( \frac{\tau }{2}+\frac{1}{4} +\frac{j}{2} \right)_{i_{\tau }} \left( \frac{\tau }{2}+1 \right)_{i_{\tau -1}} \left( \frac{\tau }{2}+\frac{3}{4}  \right)_{i_{\tau -1}}}{\left( \frac{\tau }{2} -\frac{j}{2}\right)_{i_{\tau -1}}\left( \frac{\tau }{2}+\frac{1}{4} +\frac{j}{2} \right)_{i_{\tau -1}} \left( \frac{\tau }{2}+1 \right)_{i_{\tau }} \left( \frac{\tau }{2}+\frac{3}{4} \right)_{i_{\tau }}} z^{i_{\tau }}\right\} \eta ^{\tau } \nonumber
\end{eqnarray}
where
\begin{equation}
\begin{cases} \tau \geq 2 \cr
\sigma = 1-\xi \cr
z = \frac{-1}{1-\rho ^{-2}} \sigma^2 \cr
\eta =  \frac{2-\rho ^{-2}}{1-\rho ^{-2}} \sigma \cr
Q(j,h) = \frac{h\rho ^{-2}-4j\left( j+\frac{1}{2}\right) }{16\left( 2-\rho ^{-2} \right)}
\end{cases}\nonumber
\end{equation}
\end{enumerate}
\subsection{\footnotesize ${\displaystyle (1-x)^{1-\delta } Hl(1-a,-q+(\delta -1)\gamma a+(\alpha -\delta +1)(\beta -\delta +1); \alpha-\delta +1,\beta-\delta +1, 2-\delta, \gamma; 1-x)}$ \normalsize}
Replace coefficients $a$, $q$, $\alpha$, $\beta$, $\gamma $, $\delta$ and $x$ by $1-a$, $-q+( \delta -1)\gamma a+(\alpha -\delta +1)( \beta - \delta +1)$, $\alpha -\delta +1 $, $\beta -\delta +1$, $2 -\delta $, $\gamma $ and $1-x$ into (\ref{eq:100081})--(\ref{eq:100085c}). Multiply (\ref{eq:100081}), (\ref{eq:100082b}) and (\ref{eq:100083b}) by $ (1-x)^{1-\delta }$. Substitute (\ref{eq:100080}) into the new (\ref{eq:100081})--(\ref{eq:100085c}) with replacing the new $\alpha $, $\beta $ and $q$ by $-j$, $-j$ and $1/4(h-1)\rho ^{-2}-j(j+3/2)$ where $-j \in \mathbb{N}_{0}$ in order to make $A_n$ and $B_n$ terms terminated at specific index summation $n$: There are two possible values for a coefficient $\alpha $ in (\ref{eq:100081}) such as $\alpha =-2$ or 1; two fixed values for $\alpha $ in (\ref{eq:100082a}) and (\ref{eq:100082b}) are given by $\alpha =-2\left( 2N+2\right)$ or $2\left( 2N+3/2\right)$; two fixed constants for $\alpha $ in (\ref{eq:100083a}) and (\ref{eq:100083b}) are described by $\alpha =-2\left( 2N+3\right)$ or $2\left( 2N+5/2\right)$. Replace $\alpha$ by $-2(j+1)$ in the new (\ref{eq:100084a})--(\ref{eq:100085c}).
\begin{enumerate} 
\item As $\alpha = -2 $ or 1 with $h=h_0^0= 1$,
 
The eigenfunction is given by
\begin{equation}
\sigma ^{\frac{1}{2}} y(\sigma ) = \sigma ^{\frac{1}{2}} Hl\left( 1-\rho ^{-2}, 0; 0,0,\frac{3}{2}, \frac{1}{2};\sigma \right) = \sigma ^{\frac{1}{2}} \nonumber
\end{equation}
\item As $\alpha  =-2\left( 2N+2 \right)$ or $2\left( 2N+\frac{3}{2} \right) $ where $N \in \mathbb{N}_{0}$,

An algebraic equation of degree $2N+2$ for the determination of $h$ is given by
\begin{equation}
0 = \sum_{r=0}^{N+1}\bar{c}\left( 2r, N+1-r; 2N+1,h \right)  \nonumber
\end{equation}
The eigenvalue of $h$ is written by $h_{2N+1}^m$ where $m = 0,1,2,\cdots,2N+1 $; $h_{2N+1}^0 < h_{2N+1}^1 < \cdots < h_{2N+1}^{2N+1}$. Its eigenfunction is given by 
\begin{eqnarray} 
\sigma ^{\frac{1}{2}} y(\sigma ) &=& \sigma ^{\frac{1}{2}} Hl\bigg( 1-\rho ^{-2}, \frac{\rho ^{-2}}{4}\left(  h_{2N+1}^m -1\right)-\left( 2N+1\right)\left( 2N+\frac{5}{2}\right); -\left( 2N+1\right)\nonumber\\
&&, -\left( 2N+1\right),\frac{3}{2}, \frac{1}{2};\sigma \bigg) \nonumber\\
&=& \sigma ^{\frac{1}{2}} \left\{ \sum_{r=0}^{N} y_{2r}^{N-r}\left( 2N+1,h_{2N+1}^m;\sigma \right) + \sum_{r=0}^{N} y_{2r+1}^{N-r}\left( 2N+1,h_{2N+1}^m;\sigma \right) \right\}
\nonumber
\end{eqnarray}
\item As  $\alpha  =-2\left( 2N+3 \right)$ or $2\left( 2N+\frac{5}{2} \right) $ where $N \in \mathbb{N}_{0}$,

An algebraic equation of degree $2N+3$ for the determination of $h$ is given by
\begin{equation}  
0 = \sum_{r=0}^{N+1}\bar{c}\left( 2r+1, N+1-r; 2N+2,h \right) \nonumber
\end{equation}
The eigenvalue of $h$ is written by $h_{2N+2}^m$ where $m = 0,1,2,\cdots,2N+2 $; $h_{2N+2}^0 < h_{2N+2}^1 < \cdots < h_{2N+2}^{2N+2}$. Its eigenfunction is given by
\begin{eqnarray} 
\sigma ^{\frac{1}{2}} y(\sigma ) &=&\sigma ^{\frac{1}{2}} Hl\bigg( 1-\rho ^{-2}, \frac{\rho ^{-2}}{4}\left(  h_{2N+2}^m -1\right) -\left( 2N+2\right)\left( 2N+\frac{7}{2}\right); -\left( 2N+2\right)\nonumber\\
&&, -\left( 2N+2\right),\frac{3}{2}, \frac{1}{2};\sigma \bigg) \nonumber\\
&=&\sigma ^{\frac{1}{2}} \left\{ \sum_{r=0}^{N+1} y_{2r}^{N+1-r}\left( 2N+2,h_{2N+2}^m;\sigma \right) + \sum_{r=0}^{N} y_{2r+1}^{N-r}\left( 2N+2,h_{2N+2}^m;\sigma \right) \right\} \nonumber
\end{eqnarray}
In the above,
\begin{eqnarray}
\bar{c}(0,n;j,h)  &=& \frac{\left( -\frac{j}{2}\right)_{n}\left( \frac{3}{4}+\frac{j}{2} \right)_{n}}{ \left( 1\right)_{n}\left( \frac{5}{4} \right)_{n}} \left( \frac{-1}{1-\rho ^{-2}}\right)^{n}\nonumber\\
\bar{c}(1,n;j,h) &=& \left( \frac{2-\rho ^{-2}}{1-\rho ^{-2}}\right) \sum_{i_0=0}^{n}\frac{ i_0 \left( i_0 +\frac{1}{2} \right) +Q(j,h)}{\left( i_0+\frac{1}{2} \right) \left( i_0+\frac{3}{4} \right)} \frac{\left( -\frac{j}{2}\right)_{i_0}\left( \frac{3}{4}+\frac{j}{2}  \right)_{i_0} }{\left( 1 \right)_{i_0} \left( \frac{5}{4} \right)_{i_0}}\nonumber\\
&&\times \frac{\left( \frac{1}{2}-\frac{j}{2} \right)_{n} \left( \frac{5}{4}+\frac{j}{2} \right)_{n} \left( \frac{3}{2} \right)_{i_0}\left( \frac{7}{4} \right)_{i_0}}{\left( \frac{1}{2}-\frac{j}{2} \right)_{i_0} \left( \frac{5}{4}+\frac{j}{2} \right)_{i_0} \left( \frac{3}{2} \right)_n \left( \frac{7}{4} \right)_n} \left( \frac{-1}{1-\rho ^{-2}}\right)^{n } \nonumber\\
\bar{c}(\tau ,n;j,h) &=& \left( \frac{2-\rho ^{-2}}{1-\rho ^{-2}}\right)^{\tau } \sum_{i_0=0}^{n} \frac{ i_0 \left( i_0 +\frac{1}{2} \right) +Q(j,h)}{\left( i_0+\frac{1}{2} \right) \left( i_0+\frac{3}{4} \right)} \frac{\left( -\frac{j}{2}\right)_{i_0}\left( \frac{3}{4}+\frac{j}{2} \right)_{i_0} }{\left( 1 \right)_{i_0} \left( \frac{5}{4} \right)_{i_0}}  \nonumber\\
&&\times \prod_{k=1}^{\tau -1} \left( \sum_{i_k = i_{k-1}}^{n} \frac{\left( i_k+ \frac{k}{2} \right) \left( i_k +\frac{k}{2}+\frac{1}{2} \right) +Q(j,h)}{\left( i_k+\frac{k}{2}+\frac{1}{2} \right) \left( i_k+\frac{k}{2}+\frac{3}{4} \right)} \right. \nonumber\\
&&\times \left. \frac{\left( \frac{k}{2}-\frac{j}{2}\right)_{i_k} \left( \frac{k}{2}+\frac{3}{4}+\frac{j}{2} \right)_{i_k}\left( \frac{k}{2}+1 \right)_{i_{k-1}} \left( \frac{k}{2}+ \frac{5}{4} \right)_{i_{k-1}}}{\left( \frac{k}{2}-\frac{j}{2}\right)_{i_{k-1}} \left( \frac{k}{2}+\frac{3}{4}+\frac{j}{2} \right)_{i_{k-1}}\left( \frac{k}{2}+1 \right)_{i_k} \left( \frac{k}{2}+ \frac{5}{4} \right)_{i_k}} \right) \nonumber\\ 
&&\times \frac{\left( \frac{\tau }{2} -\frac{j}{2}\right)_{n}\left( \frac{\tau }{2}+\frac{3}{4} +\frac{j}{2} \right)_{n} \left( \frac{\tau }{2}+1 \right)_{i_{\tau -1}} \left( \frac{\tau }{2}+\frac{5}{4} \right)_{i_{\tau -1}}}{\left( \frac{\tau }{2} -\frac{j}{2}\right)_{i_{\tau -1}}\left( \frac{\tau }{2}+\frac{3}{4} +\frac{j}{2} \right)_{i_{\tau -1}} \left( \frac{\tau }{2}+1 \right)_n \left( \frac{\tau }{2}+\frac{5}{4} \right)_n} \left( \frac{-1}{1-\rho ^{-2}}\right)^{n } \nonumber
\end{eqnarray}
\begin{eqnarray}
y_0^m(j,h ;\sigma ) &=& \sum_{i_0=0}^{m} \frac{\left( -\frac{j}{2}\right)_{i_0}\left( \frac{3}{4}+\frac{j}{2} \right)_{i_0} }{\left( 1 \right)_{i_0} \left( \frac{5}{4} \right)_{i_0}} z^{i_0} \nonumber\\
y_1^m(j,h ;\sigma ) &=& \left\{\sum_{i_0=0}^{m} \frac{ i_0 \left( i_0 +\frac{1}{2} \right) +Q(j,h)}{\left( i_0+\frac{1}{2} \right) \left( i_0+\frac{3}{4} \right)} \frac{\left( -\frac{j}{2}\right)_{i_0}\left( \frac{3}{4}+\frac{j}{2} \right)_{i_0} }{\left( 1 \right)_{i_0} \left( \frac{5}{4} \right)_{i_0}} \right.  \nonumber\\
&&\times  \left. \sum_{i_1 = i_0}^{m} \frac{\left( \frac{1}{2}-\frac{j}{2} \right)_{i_1} \left( \frac{5}{4}+\frac{j}{2} \right)_{i_1} \left( \frac{3}{2} \right)_{i_0}\left( \frac{7}{4} \right)_{i_0}}{\left( \frac{1}{2}-\frac{j}{2} \right)_{i_0} \left( \frac{5}{4}+\frac{j}{2} \right)_{i_0} \left( \frac{3}{2} \right)_{i_1} \left( \frac{7}{4} \right)_{i_1}} z^{i_1}\right\} \eta \nonumber
\end{eqnarray}
\begin{eqnarray}
y_{\tau }^m(j,h ;\sigma ) &=& \left\{ \sum_{i_0=0}^{m} \frac{ i_0 \left( i_0 +\frac{1}{2} \right) +Q(j,h)}{\left( i_0+\frac{1}{2} \right) \left( i_0+\frac{3}{4} \right)} \frac{\left( -\frac{j}{2}\right)_{i_0}\left( \frac{3}{4}+\frac{j}{2} \right)_{i_0} }{\left( 1 \right)_{i_0} \left( \frac{5}{4} \right)_{i_0}} \right.\nonumber\\
&&\times \prod_{k=1}^{\tau -1} \left( \sum_{i_k = i_{k-1}}^{m} \frac{\left( i_k+ \frac{k}{2} \right) \left( i_k +\frac{k}{2}+\frac{1}{2} \right) +Q(j,h)}{\left( i_k+\frac{k}{2}+\frac{1}{2} \right) \left( i_k+\frac{k}{2}+\frac{3}{4} \right)} \right. \nonumber\\
&&\times   \left. \frac{\left( \frac{k}{2}-\frac{j}{2}\right)_{i_k} \left( \frac{k}{2}+\frac{3}{4}+\frac{j}{2} \right)_{i_k}\left( \frac{k}{2}+1 \right)_{i_{k-1}} \left( \frac{k}{2}+ \frac{5}{4} \right)_{i_{k-1}}}{\left( \frac{k}{2}-\frac{j}{2}\right)_{i_{k-1}} \left( \frac{k}{2}+\frac{3}{4}+\frac{j}{2} \right)_{i_{k-1}}\left( \frac{k}{2}+1 \right)_{i_k} \left( \frac{k}{2}+ \frac{5}{4} \right)_{i_k}} \right) \nonumber\\
&&\times \left. \sum_{i_{\tau } = i_{\tau -1}}^{m} \frac{\left( \frac{\tau }{2} -\frac{j}{2}\right)_{i_{\tau }}\left( \frac{\tau }{2}+\frac{3}{4} +\frac{j}{2} \right)_{i_{\tau }} \left( \frac{\tau }{2}+1 \right)_{i_{\tau -1}} \left( \frac{\tau }{2}+\frac{5}{4}  \right)_{i_{\tau -1}}}{\left( \frac{\tau }{2} -\frac{j}{2}\right)_{i_{\tau -1}}\left( \frac{\tau }{2}+\frac{3}{4} +\frac{j}{2} \right)_{i_{\tau -1}} \left( \frac{\tau }{2}+1 \right)_{i_{\tau }} \left( \frac{\tau }{2}+\frac{5}{4} \right)_{i_{\tau }}} z^{i_{\tau }}\right\} \eta ^{\tau } \nonumber
\end{eqnarray}
where
\begin{equation}
\begin{cases} \tau \geq 2 \cr
\sigma = 1-\xi \cr
z = \frac{-1}{1-\rho ^{-2}} \sigma^2 \cr
\eta = \frac{2-\rho ^{-2}}{1-\rho ^{-2}} \sigma \cr
Q(j,h) = \frac{(h-1)\rho ^{-2}-4j\left( j+\frac{3}{2}\right)}{16(2-\rho ^{-2})}
\end{cases}\nonumber
\end{equation}
\end{enumerate}
\section{The first species complete polynomials using R3TRF} 
In chapter 6, the power series expansion of Heun equation of the first kind for the first species complete polynomial using R3TRF about $x=0$ is given by
\begin{enumerate} 
\item As $\alpha =0$ and $q=q_0^0=0$,

The eigenfunction is given by
\begin{equation}
y(x) = H_pF_{0,0}^{R} \left( a, q=q_0^0=0; \alpha =0, \beta, \gamma, \delta ; \eta = \frac{(1+a)}{a} x ; z= -\frac{1}{a} x^2 \right) =1 \label{eq:100086}
\end{equation}
\item As $\alpha =-1$,

An algebraic equation of degree 2 for the determination of $q$ is given by
\begin{equation}
0 = a\beta \gamma 
+\prod_{l=0}^{1}\Big( q+ l(\beta -\delta -1+l+a(\gamma +\delta -1+l))\Big) \label{eq:100087a}
\end{equation}
The eigenvalue of $q$ is written by $q_1^m$ where $m = 0,1 $; $q_{1}^0 < q_{1}^1$. Its eigenfunction is given by
\begin{eqnarray}
y(x) &=& H_pF_{1,m}^R \left( a, q=q_1^m; \alpha = -1, \beta, \gamma, \delta ; \eta = \frac{(1+a)}{a} x ; z= -\frac{1}{a} x^2 \right)\nonumber\\
&=&  1+\frac{ q_1^m}{(1+a)\gamma } \eta \label{eq:100087b}  
\end{eqnarray}
\item As $\alpha =-2N-2 $ where $N \in \mathbb{N}_{0}$,

An algebraic equation of degree $2N+3$ for the determination of $q$ is given by
\begin{equation}
0 = \sum_{r=0}^{N+1}\bar{c}\left( r, 2(N-r)+3; 2N+2,q\right)  \label{eq:100088a}
\end{equation}
The eigenvalue of $q$ is written by $q_{2N+2}^m$ where $m = 0,1,2,\cdots,2N+2 $; $q_{2N+2}^0 < q_{2N+2}^1 < \cdots < q_{2N+2}^{2N+2}$. Its eigenfunction is given by 
\begin{eqnarray} 
y(x) &=& H_pF_{2N+2,m}^R \left( a, q=q_{2N+2}^m; \alpha =-2N-2, \beta, \gamma, \delta ; \eta = \frac{(1+a)}{a} x ; z= -\frac{1}{a} x^2 \right)\nonumber\\
&=& \sum_{r=0}^{N+1} y_{r}^{2(N+1-r)}\left( 2N+2, q_{2N+2}^m; x \right)  
\label{eq:100088b} 
\end{eqnarray}
\item As $\alpha =-2N-3 $ where $N \in \mathbb{N}_{0}$,

An algebraic equation of degree $2N+4$ for the determination of $q$ is given by
\begin{equation}  
0 = \sum_{r=0}^{N+2}\bar{c}\left( r, 2(N+2-r); 2N+3,q\right) \label{eq:100089a}
\end{equation}
The eigenvalue of $q$ is written by $q_{2N+3}^m$ where $m = 0,1,2,\cdots,2N+3 $; $q_{2N+3}^0 < q_{2N+3}^1 < \cdots < q_{2N+3}^{2N+3}$. Its eigenfunction is given by
\begin{eqnarray} 
y(x) &=& H_pF_{2N+3,m}^R \left( a, q=q_{2N+3}^m; \alpha =-2N-3, \beta, \gamma, \delta ; \eta = \frac{(1+a)}{a} x ; z= -\frac{1}{a} x^2 \right)\nonumber\\
&=&   \sum_{r=0}^{N+1} y_{r}^{2(N-r)+3} \left( 2N+3,q_{2N+3}^m;x\right) \label{eq:100089b}
\end{eqnarray}
In the above,
\begin{eqnarray}
\bar{c}(0,n;j,q)  &=& \frac{\left( \Delta_0^{-} \left( j,q\right) \right)_{n}\left( \Delta_0^{+} \left( j,q\right) \right)_{n}}{\left( 1 \right)_{n} \left( \gamma \right)_{n}} \left( \frac{1+a}{a} \right)^{n}\label{eq:100090a}\\
\bar{c}(1,n;j,q) &=& \left( -\frac{1}{a}\right) \sum_{i_0=0}^{n}\frac{\left( i_0 -j\right)\left( i_0+\beta \right) }{\left( i_0+2 \right) \left( i_0+1+ \gamma \right)} \frac{ \left( \Delta_0^{-} \left( j,q\right) \right)_{i_0}\left( \Delta_0^{+} \left( j,q\right) \right)_{i_0}}{\left( 1 \right)_{i_0} \left( \gamma \right)_{i_0}} \nonumber\\
&\times&  \frac{ \left( \Delta_1^{-} \left( j,q\right) \right)_{n}\left( \Delta_1^{+} \left( j,q\right) \right)_{n} \left( 3 \right)_{i_0} \left( 2+ \gamma \right)_{i_0}}{\left( \Delta_1^{-} \left( j,q\right) \right)_{i_0}\left( \Delta_1^{+} \left( j,q\right) \right)_{i_0}\left( 3 \right)_{n} \left( 2+ \gamma \right)_{n}} \left(\frac{1+a}{a} \right)^{n }  
\label{eq:100090b}\\
\bar{c}(\tau ,n;j,q) &=& \left( -\frac{1}{a}\right)^{\tau} \sum_{i_0=0}^{n}\frac{\left( i_0 -j\right)\left( i_0+\beta \right) }{\left( i_0+2 \right) \left( i_0+1+ \gamma \right)} \frac{ \left( \Delta_0^{-} \left( j,q\right) \right)_{i_0}\left( \Delta_0^{+} \left( j,q\right) \right)_{i_0}}{\left( 1 \right)_{i_0} \left( \gamma \right)_{i_0}}  \nonumber\\
&\times& \prod_{k=1}^{\tau -1} \left( \sum_{i_k = i_{k-1}}^{n} \frac{\left( i_k+ 2k-j\right)\left( i_k +2k+\beta \right)}{\left( i_k+2k+2 \right) \left( i_k+2k+1+ \gamma \right)} \right. \nonumber\\
&\times&\left. \frac{ \left( \Delta_k^{-} \left( j,q\right) \right)_{i_k}\left( \Delta_k^{+} \left( j,q\right) \right)_{i_k} \left( 2k+1 \right)_{i_{k-1}} \left( 2k+ \gamma \right)_{i_{k-1}}}{\left( \Delta_k^{-} \left( j,q\right) \right)_{i_{k-1}}\left( \Delta_k^{+} \left( j,q\right) \right)_{i_{k-1}}\left( 2k+1 \right)_{i_k} \left( 2k+ \gamma \right)_{i_k}} \right) \nonumber\\
&\times& \frac{ \left( \Delta_{\tau }^{-} \left( j,q\right) \right)_{n}\left( \Delta_{\tau }^{+} \left( j,q\right) \right)_{n} \left( 2\tau +1 \right)_{i_{\tau -1}} \left( 2\tau + \gamma \right)_{i_{\tau -1}}}{\left( \Delta_{\tau }^{-} \left( j,q\right) \right)_{i_{\tau -1}}\left( \Delta_{\tau }^{+} \left( j,q\right) \right)_{i_{\tau -1}}\left( 2\tau +1 \right)_{n} \left( 2\tau + \gamma  \right)_{n}} \left(\frac{1+a}{a}\right)^{n } \hspace{1.5cm}\label{eq:100090c} 
\end{eqnarray}
\begin{eqnarray}
y_0^m(j,q;x) &=& \sum_{i_0=0}^{m} \frac{\left( \Delta_0^{-} \left( j,q\right) \right)_{i_0}\left( \Delta_0^{+} \left( j,q\right) \right)_{i_0}}{\left( 1 \right)_{i_0} \left( \gamma \right)_{i_0}} \eta ^{i_0} \label{eq:100091a}\\
y_1^m(j,q;x) &=& \left\{\sum_{i_0=0}^{m}\frac{\left( i_0 -j\right)\left( i_0+\beta \right) }{\left( i_0+2 \right) \left( i_0+1+ \gamma \right)} \frac{ \left( \Delta_0^{-} \left( j,q\right) \right)_{i_0}\left( \Delta_0^{+} \left( j,q\right) \right)_{i_0}}{\left( 1 \right)_{i_0} \left( \gamma \right)_{i_0}} \right. \nonumber\\
&\times& \left. \sum_{i_1 = i_0}^{m} \frac{ \left( \Delta_1^{-} \left( j,q\right) \right)_{i_1}\left( \Delta_1^{+} \left( j,q\right) \right)_{i_1} \left( 3 \right)_{i_0} \left( 2+ \gamma \right)_{i_0}}{\left( \Delta_1^{-} \left( j,q\right) \right)_{i_0}\left( \Delta_1^{+} \left( j,q\right) \right)_{i_0}\left( 3 \right)_{i_1} \left( 2+ \gamma \right)_{i_1}} \eta ^{i_1}\right\} z 
\label{eq:100091b}
\end{eqnarray}
\begin{eqnarray}
y_{\tau }^m(j,q;x) &=& \left\{ \sum_{i_0=0}^{m} \frac{\left( i_0 -j\right)\left( i_0+\beta \right) }{\left( i_0+2 \right) \left( i_0+1+ \gamma \right)} \frac{ \left( \Delta_0^{-} \left( j,q\right) \right)_{i_0}\left( \Delta_0^{+} \left( j,q\right) \right)_{i_0}}{\left( 1 \right)_{i_0} \left( \gamma \right)_{i_0}} \right.\nonumber\\
&\times& \prod_{k=1}^{\tau -1} \left( \sum_{i_k = i_{k-1}}^{m} \frac{\left( i_k+ 2k-j\right)\left( i_k +2k+\beta \right)}{\left( i_k+2k+2 \right) \left( i_k+2k+1+ \gamma \right)} \right.\nonumber\\
&\times& \left. \frac{ \left( \Delta_k^{-} \left( j,q\right) \right)_{i_k}\left( \Delta_k^{+} \left( j,q\right) \right)_{i_k} \left( 2k+1 \right)_{i_{k-1}} \left( 2k+ \gamma \right)_{i_{k-1}}}{\left( \Delta_k^{-} \left( j,q\right) \right)_{i_{k-1}}\left( \Delta_k^{+} \left( j,q\right) \right)_{i_{k-1}}\left( 2k+1 \right)_{i_k} \left( 2k+ \gamma \right)_{i_k}} \right) \nonumber\\
&\times & \left. \sum_{i_{\tau } = i_{\tau -1}}^{m}  \frac{ \left( \Delta_{\tau }^{-} \left( j,q\right) \right)_{i_{\tau }}\left( \Delta_{\tau }^{+} \left( j,q\right) \right)_{i_{\tau }} \left( 2\tau +1  \right)_{i_{\tau -1}} \left( 2\tau + \gamma \right)_{i_{\tau -1}}}{\left( \Delta_{\tau }^{-} \left( j,q\right) \right)_{i_{\tau -1}}\left( \Delta_{\tau }^{+} \left( j,q\right) \right)_{i_{\tau -1}}\left( 2\tau +1 \right)_{i_\tau } \left( 2\tau + \gamma \right)_{i_{\tau }}} \eta ^{i_{\tau }}\right\} z^{\tau } \hspace{1.5cm}\label{eq:100091c} 
\end{eqnarray}
where
\begin{equation}
\begin{cases} \tau \geq 2 \cr
\Delta_k^{\pm} \left( j,q\right) = \frac{\varphi +4(1+a)k \pm \sqrt{\varphi ^2-4(1+a)q}}{2(1+a)}  \cr
\varphi = \beta -\delta -j+a(\gamma +\delta -1)
\end{cases}\nonumber
\end{equation}
\end{enumerate} 
\subsection{ ${\displaystyle (1-x)^{1-\delta } Hl(a, q - (\delta  - 1)\gamma a; \alpha - \delta  + 1, \beta - \delta + 1, \gamma ,2 - \delta ; x)}$ }
Replace coefficients $q$, $\alpha$, $\beta$ and $\delta$ by $q - (\delta - 1)\gamma a $, $\alpha - \delta  + 1 $, $\beta - \delta + 1$ and $2 - \delta$ into (\ref{eq:100086})--(\ref{eq:100091c}). Multiply (\ref{eq:100086}), (\ref{eq:100087b}), (\ref{eq:100088b}) and (\ref{eq:100089b}) by $(1-x)^{1-\delta }$. Substitute (\ref{eq:100080}) into the new (\ref{eq:100086})--(\ref{eq:100091c}) with replacing both new $\alpha $ and $\beta $ by $-j$ where $-j \in \mathbb{N}_{0}$ in order to make $A_n$ and $B_n$ terms terminated at specific index summation $n$: There are two possible values for a coefficient $\alpha $ in (\ref{eq:100086}) such as $\alpha =-2$ or 1; the two fixed values for $\alpha $ in (\ref{eq:100087a}) and (\ref{eq:100087b}) are give by $\alpha =-4$ or 3; two fixed constants for $\alpha $ in (\ref{eq:100088a}) and (\ref{eq:100088b}) are given by $\alpha =-2\left( 2N+3\right)$ or $2\left( 2N+5/2\right)$; two fixed values for $\alpha $ in (\ref{eq:100089a}) and (\ref{eq:100089b}) are described by $\alpha =-2\left( 2N+4\right)$ or $2\left( 2N+7/2\right)$. Replace $\alpha$ by $-2(j+1)$ in the new (\ref{eq:100090a})--(\ref{eq:100091c}).
\begin{enumerate} 
\item As $\alpha = -2 $ or 1 with $h=h_0^0=1 $,
 
The eigenfunction is given by
\begin{equation}
(1-\xi )^{\frac{1}{2}} y(\xi ) = (1-\xi )^{\frac{1}{2}} Hl\left( \rho ^{-2}, 0; 0,0,\frac{1}{2}, \frac{3}{2};\xi \right) = (1-\xi )^{\frac{1}{2}} \nonumber
\end{equation}
\item As $\alpha = -4 $ or 3,

An algebraic equation of degree 2 for the determination of $h$ is given by
\begin{equation}
0 =  \frac{5}{4}\rho ^{-2}+ \prod_{l=0}^{1}\left( \frac{\rho ^{-2}}{4}(1-h) +l\left( l+(l+1)\rho ^{-2}\right)\right)\nonumber
\end{equation}
The eigenvalue of $h$ is written by $h_1^m$ where $m = 0,1 $; $h_{1}^0 < h_{1}^1$. Its eigenfunction is given by
\begin{eqnarray}
(1-\xi )^{\frac{1}{2}} y(\xi ) &=& (1-\xi )^{\frac{1}{2}} Hl\left( \rho ^{-2}, \frac{\rho ^{-2}}{4}(1-h_1^m); -1,-1,\frac{1}{2}, \frac{3}{2};\xi \right) \nonumber\\
&=&  (1-\xi )^{\frac{1}{2}} \left\{ 1+\frac{ 1-h_1^m}{2( 1+\rho ^2 )} \eta \right\} \nonumber
\end{eqnarray}
\item As $\alpha  =-2\left( 2N+3 \right)$ or $2\left( 2N+\frac{5}{2} \right) $ where $N \in \mathbb{N}_{0}$,

An algebraic equation of degree $2N+3$ for the determination of $h$ is given by
\begin{equation}
0 = \sum_{r=0}^{N+1}\bar{c}\left( r, 2(N-r)+3; 2N+2,h\right)  \nonumber
\end{equation}
The eigenvalue of $h$ is written by $h_{2N+2}^m$ where $m = 0,1,2,\cdots,2N+2 $; $h_{2N+2}^0 < h_{2N+2}^1 < \cdots < h_{2N+2}^{2N+2}$. Its eigenfunction is given by 
\begin{eqnarray} 
(1-\xi )^{\frac{1}{2}} y(\xi ) &=& (1-\xi )^{\frac{1}{2}} Hl\left( \rho ^{-2}, \frac{\rho ^{-2}}{4}(1-h_{2N+2}^m); -\left( 2N+2\right), -\left( 2N+2\right),\frac{1}{2}, \frac{3}{2};\xi \right) \nonumber\\
&=& (1-\xi )^{\frac{1}{2}} \sum_{r=0}^{N+1} y_{r}^{2(N+1-r)}\left( 2N+2, h_{2N+2}^m; \xi \right)  
\nonumber
\end{eqnarray}
\item As $\alpha  =-2\left( 2N+4 \right)$ or $2\left( 2N+\frac{7}{2} \right) $ where $N \in \mathbb{N}_{0}$,

An algebraic equation of degree $2N+4$ for the determination of $h$ is given by
\begin{equation}  
0 = \sum_{r=0}^{N+2}\bar{c}\left( r, 2(N+2-r); 2N+3,h\right) \nonumber
\end{equation}
The eigenvalue of $h$ is written by $h_{2N+3}^m$ where $m = 0,1,2,\cdots,2N+3 $; $h_{2N+3}^0 < h_{2N+3}^1 < \cdots < h_{2N+3}^{2N+3}$. Its eigenfunction is given by
\begin{eqnarray} 
(1-\xi )^{\frac{1}{2}} y(\xi ) &=&(1-\xi )^{\frac{1}{2}} Hl\left( \rho ^{-2}, \frac{\rho ^{-2}}{4}(1-h_{2N+3}^m); -\left( 2N+3\right), -\left( 2N+3\right),\frac{1}{2}, \frac{3}{2};\xi \right)\nonumber\\
&=&  (1-\xi )^{\frac{1}{2}} \sum_{r=0}^{N+1} y_{r}^{2(N-r)+3} \left( 2N+3,h_{2N+3}^m; \xi \right) \nonumber
\end{eqnarray}
In the above,
\begin{eqnarray}
\bar{c}(0,n;j,h)  &=& \frac{\left( \Delta_0^{-} \left( h\right) \right)_{n}\left( \Delta_0^{+} \left( h\right) \right)_{n}}{\left( 1 \right)_{n} \left( \frac{1}{2} \right)_{n}} \left( 1+\rho ^2 \right)^{n}\nonumber\\
\bar{c}(1,n;j,h) &=& \left( -\rho ^2\right) \sum_{i_0=0}^{n}\frac{\left( i_0 -j\right)\left( i_0+j+\frac{3}{2} \right) }{\left( i_0+2 \right) \left( i_0+\frac{3}{2} \right)} \frac{ \left( \Delta_0^{-} \left( h\right) \right)_{i_0}\left( \Delta_0^{+} \left( h\right) \right)_{i_0}}{\left( 1 \right)_{i_0} \left( \frac{1}{2} \right)_{i_0}} \nonumber\\
&\times&  \frac{ \left( \Delta_1^{-} \left( h\right) \right)_{n}\left( \Delta_1^{+} \left( h\right) \right)_{n} \left( 3 \right)_{i_0} \left( \frac{5}{2} \right)_{i_0}}{\left( \Delta_1^{-} \left( h\right) \right)_{i_0}\left( \Delta_1^{+} \left( h\right) \right)_{i_0}\left( 3 \right)_{n} \left( \frac{5}{2} \right)_{n}} \left( 1+\rho ^2 \right)^{n }  
\nonumber\\
\bar{c}(\tau ,n;j,h) &=& \left( -\rho ^2\right)^{\tau} \sum_{i_0=0}^{n}\frac{\left( i_0 -j\right)\left( i_0+j+\frac{3}{2} \right) }{\left( i_0+2 \right) \left( i_0+\frac{3}{2} \right)} \frac{ \left( \Delta_0^{-} \left( h\right) \right)_{i_0}\left( \Delta_0^{+} \left( h\right) \right)_{i_0}}{\left( 1 \right)_{i_0} \left( \frac{1}{2} \right)_{i_0}}  \nonumber\\
&\times& \prod_{k=1}^{\tau -1} \left( \sum_{i_k = i_{k-1}}^{n} \frac{\left( i_k+ 2k-j\right)\left( i_k +2k+j+\frac{3}{2} \right)}{\left( i_k+2k+2 \right) \left( i_k+2k+\frac{3}{2} \right)} \right. \nonumber\\
&\times&\left. \frac{ \left( \Delta_k^{-} \left( h\right) \right)_{i_k}\left( \Delta_k^{+} \left( h\right) \right)_{i_k} \left( 2k+1 \right)_{i_{k-1}} \left( 2k+ \frac{1}{2} \right)_{i_{k-1}}}{\left( \Delta_k^{-} \left( h\right) \right)_{i_{k-1}}\left( \Delta_k^{+} \left( h\right) \right)_{i_{k-1}}\left( 2k+1 \right)_{i_k} \left( 2k+ \frac{1}{2} \right)_{i_k}} \right) \nonumber\\
&\times& \frac{ \left( \Delta_{\tau }^{-} \left( h\right) \right)_{n}\left( \Delta_{\tau }^{+} \left( h\right) \right)_{n} \left( 2\tau +1 \right)_{i_{\tau -1}} \left( 2\tau + \frac{1}{2} \right)_{i_{\tau -1}}}{\left( \Delta_{\tau }^{-} \left( h\right) \right)_{i_{\tau -1}}\left( \Delta_{\tau }^{+} \left( h\right) \right)_{i_{\tau -1}}\left( 2\tau +1 \right)_{n} \left( 2\tau + \frac{1}{2} \right)_{n}} \left( 1+\rho ^2 \right)^{n } \nonumber
\end{eqnarray}
\begin{eqnarray}
y_0^m(j,h;\xi ) &=& \sum_{i_0=0}^{m} \frac{\left( \Delta_0^{-} \left( h\right) \right)_{i_0}\left( \Delta_0^{+} \left( h\right) \right)_{i_0}}{\left( 1 \right)_{i_0} \left( \frac{1}{2} \right)_{i_0}} \eta ^{i_0} \nonumber\\
y_1^m(j,h;\xi ) &=& \left\{\sum_{i_0=0}^{m}\frac{\left( i_0 -j\right)\left( i_0+j+\frac{3}{2} \right) }{\left( i_0+2 \right) \left( i_0+\frac{3}{2} \right)} \frac{ \left( \Delta_0^{-} \left( h\right) \right)_{i_0}\left( \Delta_0^{+} \left( h\right) \right)_{i_0}}{\left( 1 \right)_{i_0} \left( \frac{1}{2} \right)_{i_0}} \right. \nonumber\\
&\times& \left. \sum_{i_1 = i_0}^{m} \frac{ \left( \Delta_1^{-} \left( h\right) \right)_{i_1}\left( \Delta_1^{+} \left( h\right) \right)_{i_1} \left( 3 \right)_{i_0} \left( \frac{5}{2} \right)_{i_0}}{\left( \Delta_1^{-} \left( h\right) \right)_{i_0}\left( \Delta_1^{+} \left( h\right) \right)_{i_0}\left( 3 \right)_{i_1} \left( \frac{5}{2} \right)_{i_1}} \eta ^{i_1}\right\} z 
\nonumber
\end{eqnarray}
\begin{eqnarray}
y_{\tau }^m(j,h;\xi ) &=& \left\{ \sum_{i_0=0}^{m} \frac{\left( i_0 -j\right)\left( i_0+j+\frac{3}{2} \right) }{\left( i_0+2 \right) \left( i_0+\frac{3}{2} \right)} \frac{ \left( \Delta_0^{-} \left( h\right) \right)_{i_0}\left( \Delta_0^{+} \left( h\right) \right)_{i_0}}{\left( 1 \right)_{i_0} \left( \frac{1}{2} \right)_{i_0}} \right.\nonumber\\
&\times& \prod_{k=1}^{\tau -1} \left( \sum_{i_k = i_{k-1}}^{m} \frac{\left( i_k+ 2k-j\right)\left( i_k +2k+j+\frac{3}{2}\right)}{\left( i_k+2k+2 \right) \left( i_k+2k+\frac{3}{2} \right)} \right.\nonumber\\
&\times& \left. \frac{ \left( \Delta_k^{-} \left( h\right) \right)_{i_k}\left( \Delta_k^{+} \left( h\right) \right)_{i_k} \left( 2k+1 \right)_{i_{k-1}} \left( 2k+ \frac{1}{2} \right)_{i_{k-1}}}{\left( \Delta_k^{-} \left( h\right) \right)_{i_{k-1}}\left( \Delta_k^{+} \left( h\right) \right)_{i_{k-1}}\left( 2k+1 \right)_{i_k} \left( 2k+ \frac{1}{2} \right)_{i_k}} \right) \nonumber\\
&\times & \left. \sum_{i_{\tau } = i_{\tau -1}}^{m}  \frac{ \left( \Delta_{\tau }^{-} \left( h\right) \right)_{i_{\tau }}\left( \Delta_{\tau }^{+} \left( h\right) \right)_{i_{\tau }} \left( 2\tau +1  \right)_{i_{\tau -1}} \left( 2\tau + \frac{1}{2} \right)_{i_{\tau -1}}}{\left( \Delta_{\tau }^{-} \left( h\right) \right)_{i_{\tau -1}}\left( \Delta_{\tau }^{+} \left( h\right) \right)_{i_{\tau -1}}\left( 2\tau +1 \right)_{i_\tau } \left( 2\tau + \frac{1}{2} \right)_{i_{\tau }}} \eta ^{i_{\tau }}\right\} z^{\tau } \nonumber
\end{eqnarray}
where
\begin{equation}
\begin{cases} \tau \geq 2 \cr
\Delta_k^{\pm} \left( h\right) = 2k +\frac{1 \pm \sqrt{(1+\rho ^2)h-\rho ^2}}{2(1+\rho ^2)}  \cr
z= -\rho ^2\xi ^2\cr
\eta = (1+\rho ^2) \xi
\end{cases}\nonumber
\end{equation}
\end{enumerate} 
\subsection{\footnotesize ${\displaystyle x^{1-\gamma } (1-x)^{1-\delta }Hl(a, q-(\gamma +\delta -2)a -(\gamma -1)(\alpha +\beta -\gamma -\delta +1), \alpha - \gamma -\delta +2, \beta - \gamma -\delta +2, 2-\gamma, 2 - \delta ; x)}$ \normalsize}
Replace coefficients $q$, $\alpha$, $\beta$, $\gamma $ and $\delta$ by $q-(\gamma +\delta -2)a-(\gamma -1)(\alpha +\beta -\gamma -\delta +1) $, $\alpha - \gamma -\delta +2 $, $\beta - \gamma -\delta +2$, $2-\gamma $ and $2 - \delta$ into (\ref{eq:100086})--(\ref{eq:100091c}). Multiply (\ref{eq:100086}), (\ref{eq:100087b}), (\ref{eq:100088b}) and (\ref{eq:100089b}) by $x^{1-\gamma } (1-x)^{1-\delta }$. Substitute (\ref{eq:100080}) into the new (\ref{eq:100086})--(\ref{eq:100091c}) with replacing both new $\alpha $ and $\beta $ by $-j$ where $-j \in \mathbb{N}_{0}$ in order to make $A_n$ and $B_n$ terms terminated at specific index summation $n$: There are two possible values for a coefficient $\alpha $ in (\ref{eq:100086}) such as $\alpha =-3$ or 2; the two fixed values for $\alpha $ in (\ref{eq:100087a}) and (\ref{eq:100087b}) are give by $\alpha =-5$ or 4; two fixed constants for $\alpha $ in (\ref{eq:100088a}) and (\ref{eq:100088b}) are given by $\alpha =-2\left( 2N+7/2\right)$ or $2\left( 2N+3\right)$; two fixed values for $\alpha $ in (\ref{eq:100089a}) and (\ref{eq:100089b}) are described by $\alpha =-2\left( 2N+9/2\right)$ or $2\left( 2N+4\right)$. Replace $\alpha$ by $-2(j+3/2)$ in the new (\ref{eq:100090a})--(\ref{eq:100091c}).
\begin{enumerate}
\item As $\alpha = -3 $ or 2 with $h=h_0^0=4+ \rho ^2 $,
 
The eigenfunction is given by
\begin{equation}
\xi^{\frac{1}{2}}(1-\xi )^{\frac{1}{2}} y(\xi ) = \xi^{\frac{1}{2}}(1-\xi )^{\frac{1}{2}} Hl\left( \rho ^{-2}, 0; 0,0,\frac{3}{2}, \frac{3}{2};\xi \right) = \xi^{\frac{1}{2}}(1-\xi )^{\frac{1}{2}} \nonumber
\end{equation}
\item As $\alpha = -5 $ or 4,

An algebraic equation of degree 2 for the determination of $h$ is given by
\begin{equation}
0 =  \frac{21}{4}\rho ^{-2}+ \prod_{l=0}^{1}\left( \frac{1}{4}\left( (4-h)\rho ^{-2} +1\right) +l\left( l+1+(l+2)\rho ^{-2}\right)\right)\nonumber
\end{equation}
The eigenvalue of $h$ is written by $h_1^m$ where $m = 0,1 $; $h_{1}^0 < h_{1}^1$. Its eigenfunction is given by
\begin{eqnarray}
\xi^{\frac{1}{2}}(1-\xi )^{\frac{1}{2}} y(\xi ) &=& \xi^{\frac{1}{2}}(1-\xi )^{\frac{1}{2}} Hl\left( \rho ^{-2}, \frac{1}{4}\left( (4-h_1^m)\rho ^{-2} +1\right); -1,-1,\frac{3}{2}, \frac{3}{2};\xi \right) \nonumber\\
&=& \xi^{\frac{1}{2}} (1-\xi )^{\frac{1}{2}} \left\{ 1+\frac{ 4-h_1^m+\rho ^2}{6(1+\rho ^2)} \eta \right\} \nonumber
\end{eqnarray}
\item As $\alpha  =-2\left( 2N+\frac{7}{2} \right)$ or $2\left( 2N+3 \right) $ where $N \in \mathbb{N}_{0}$,

An algebraic equation of degree $2N+3$ for the determination of $h$ is given by
\begin{equation}
0 = \sum_{r=0}^{N+1}\bar{c}\left( r, 2(N-r)+3; 2N+2,h\right)  \nonumber
\end{equation}
The eigenvalue of $h$ is written by $h_{2N+2}^m$ where $m = 0,1,2,\cdots,2N+2 $; $h_{2N+2}^0 < h_{2N+2}^1 < \cdots < h_{2N+2}^{2N+2}$. Its eigenfunction is given by 
\begin{eqnarray} 
\xi^{\frac{1}{2}}(1-\xi )^{\frac{1}{2}} y(\xi ) &=& \xi^{\frac{1}{2}}(1-\xi )^{\frac{1}{2}} Hl\left( \rho ^{-2}, \frac{1}{4}\left( (4- h_{2N+2}^m)\rho ^{-2} +1\right); -\left( 2N+2\right), -\left( 2N+2\right),\frac{3}{2}, \frac{3}{2};\xi \right) \nonumber\\
&=& \xi^{\frac{1}{2}}(1-\xi )^{\frac{1}{2}} \sum_{r=0}^{N+1} y_{r}^{2(N+1-r)}\left( 2N+2, h_{2N+2}^m; \xi \right)  
\nonumber
\end{eqnarray}
\item As $\alpha  =-2\left( 2N+\frac{9}{2} \right)$ or $2\left( 2N+4 \right) $ where $N \in \mathbb{N}_{0}$,

An algebraic equation of degree $2N+4$ for the determination of $h$ is given by
\begin{equation}  
0 = \sum_{r=0}^{N+2}\bar{c}\left( r, 2(N+2-r); 2N+3,h\right) \nonumber
\end{equation}
The eigenvalue of $h$ is written by $h_{2N+3}^m$ where $m = 0,1,2,\cdots,2N+3 $; $h_{2N+3}^0 < h_{2N+3}^1 < \cdots < h_{2N+3}^{2N+3}$. Its eigenfunction is given by
\begin{eqnarray} 
\xi^{\frac{1}{2}}(1-\xi )^{\frac{1}{2}} y(\xi ) &=& \xi^{\frac{1}{2}}(1-\xi )^{\frac{1}{2}} Hl\left( \rho ^{-2}, \frac{1}{4}\left( (4- h_{2N+3}^m)\rho ^{-2} +1\right); -\left( 2N+3\right), -\left( 2N+3\right),\frac{3}{2}, \frac{3}{2};\xi \right)\nonumber\\
&=&  \xi^{\frac{1}{2}}(1-\xi )^{\frac{1}{2}} \sum_{r=0}^{N+1} y_{r}^{2(N-r)+3} \left( 2N+3,h_{2N+3}^m; \xi \right) \nonumber
\end{eqnarray}
In the above,
\begin{eqnarray}
\bar{c}(0,n;j,h)  &=& \frac{\left( \Delta_0^{-} \left( h\right) \right)_{n}\left( \Delta_0^{+} \left( h\right) \right)_{n}}{\left( 1 \right)_{n} \left( \frac{3}{2} \right)_{n}} \left( 1+\rho ^2 \right)^{n}\nonumber\\
\bar{c}(1,n;j,h) &=& \left( -\rho ^2\right) \sum_{i_0=0}^{n}\frac{\left( i_0 -j\right)\left( i_0+j+\frac{5}{2} \right) }{\left( i_0+2 \right) \left( i_0+\frac{5}{2} \right)} \frac{ \left( \Delta_0^{-} \left( h\right) \right)_{i_0}\left( \Delta_0^{+} \left( h\right) \right)_{i_0}}{\left( 1 \right)_{i_0} \left( \frac{3}{2} \right)_{i_0}} \nonumber\\
&\times&  \frac{ \left( \Delta_1^{-} \left( h\right) \right)_{n}\left( \Delta_1^{+} \left( h\right) \right)_{n} \left( 3 \right)_{i_0} \left( \frac{7}{2} \right)_{i_0}}{\left( \Delta_1^{-} \left( h\right) \right)_{i_0}\left( \Delta_1^{+} \left( h\right) \right)_{i_0}\left( 3 \right)_{n} \left( \frac{7}{2} \right)_{n}} \left( 1+\rho ^2 \right)^{n }  
\nonumber\\
\bar{c}(\tau ,n;j,h) &=& \left( -\rho ^2\right)^{\tau} \sum_{i_0=0}^{n}\frac{\left( i_0 -j\right)\left( i_0+j+\frac{5}{2} \right) }{\left( i_0+2 \right) \left( i_0+\frac{5}{2} \right)} \frac{ \left( \Delta_0^{-} \left( h\right) \right)_{i_0}\left( \Delta_0^{+} \left( h\right) \right)_{i_0}}{\left( 1 \right)_{i_0} \left( \frac{3}{2} \right)_{i_0}}  \nonumber\\
&\times& \prod_{k=1}^{\tau -1} \left( \sum_{i_k = i_{k-1}}^{n} \frac{\left( i_k+ 2k-j\right)\left( i_k +2k+j+\frac{5}{2} \right)}{\left( i_k+2k+2 \right) \left( i_k+2k+\frac{5}{2} \right)} \right. \nonumber\\
&\times&\left. \frac{ \left( \Delta_k^{-} \left( h\right) \right)_{i_k}\left( \Delta_k^{+} \left( h\right) \right)_{i_k} \left( 2k+1 \right)_{i_{k-1}} \left( 2k+ \frac{3}{2} \right)_{i_{k-1}}}{\left( \Delta_k^{-} \left( h\right) \right)_{i_{k-1}}\left( \Delta_k^{+} \left( h\right) \right)_{i_{k-1}}\left( 2k+1 \right)_{i_k} \left( 2k+ \frac{3}{2} \right)_{i_k}} \right) \nonumber\\
&\times& \frac{ \left( \Delta_{\tau }^{-} \left( h\right) \right)_{n}\left( \Delta_{\tau }^{+} \left( h\right) \right)_{n} \left( 2\tau +1 \right)_{i_{\tau -1}} \left( 2\tau + \frac{3}{2} \right)_{i_{\tau -1}}}{\left( \Delta_{\tau }^{-} \left( h\right) \right)_{i_{\tau -1}}\left( \Delta_{\tau }^{+} \left( h\right) \right)_{i_{\tau -1}}\left( 2\tau +1 \right)_{n} \left( 2\tau + \frac{3}{2} \right)_{n}} \left( 1+\rho ^2 \right)^{n } \nonumber
\end{eqnarray}
\begin{eqnarray}
y_0^m(j,h;\xi ) &=& \sum_{i_0=0}^{m} \frac{\left( \Delta_0^{-} \left( h\right) \right)_{i_0}\left( \Delta_0^{+} \left( h\right) \right)_{i_0}}{\left( 1 \right)_{i_0} \left( \frac{3}{2} \right)_{i_0}} \eta ^{i_0} \nonumber\\
y_1^m(j,h;\xi ) &=& \left\{\sum_{i_0=0}^{m}\frac{\left( i_0 -j\right)\left( i_0+j+\frac{5}{2} \right) }{\left( i_0+2 \right) \left( i_0+\frac{5}{2} \right)} \frac{ \left( \Delta_0^{-} \left( h\right) \right)_{i_0}\left( \Delta_0^{+} \left( h\right) \right)_{i_0}}{\left( 1 \right)_{i_0} \left( \frac{3}{2} \right)_{i_0}} \right. \nonumber\\
&\times& \left. \sum_{i_1 = i_0}^{m} \frac{ \left( \Delta_1^{-} \left( h\right) \right)_{i_1}\left( \Delta_1^{+} \left( h\right) \right)_{i_1} \left( 3 \right)_{i_0} \left( \frac{7}{2} \right)_{i_0}}{\left( \Delta_1^{-} \left( h\right) \right)_{i_0}\left( \Delta_1^{+} \left( h\right) \right)_{i_0}\left( 3 \right)_{i_1} \left( \frac{7}{2} \right)_{i_1}} \eta ^{i_1}\right\} z 
\nonumber
\end{eqnarray}
\begin{eqnarray}
y_{\tau }^m(j,h;\xi ) &=& \left\{ \sum_{i_0=0}^{m} \frac{\left( i_0 -j\right)\left( i_0+j+\frac{5}{2} \right) }{\left( i_0+2 \right) \left( i_0+\frac{5}{2} \right)} \frac{ \left( \Delta_0^{-} \left( h\right) \right)_{i_0}\left( \Delta_0^{+} \left( h\right) \right)_{i_0}}{\left( 1 \right)_{i_0} \left( \frac{3}{2} \right)_{i_0}} \right.\nonumber\\
&\times& \prod_{k=1}^{\tau -1} \left( \sum_{i_k = i_{k-1}}^{m} \frac{\left( i_k+ 2k-j\right)\left( i_k +2k+j+\frac{5}{2}\right)}{\left( i_k+2k+2 \right) \left( i_k+2k+\frac{5}{2} \right)} \right. \nonumber\\
&\times&\left. \frac{ \left( \Delta_k^{-} \left( h\right) \right)_{i_k}\left( \Delta_k^{+} \left( h\right) \right)_{i_k} \left( 2k+1 \right)_{i_{k-1}} \left( 2k+ \frac{3}{2} \right)_{i_{k-1}}}{\left( \Delta_k^{-} \left( h\right) \right)_{i_{k-1}}\left( \Delta_k^{+} \left( h\right) \right)_{i_{k-1}}\left( 2k+1 \right)_{i_k} \left( 2k+ \frac{3}{2} \right)_{i_k}} \right) \nonumber\\
&\times & \left. \sum_{i_{\tau } = i_{\tau -1}}^{m}  \frac{ \left( \Delta_{\tau }^{-} \left( h\right) \right)_{i_{\tau }}\left( \Delta_{\tau }^{+} \left( h\right) \right)_{i_{\tau }} \left( 2\tau +1  \right)_{i_{\tau -1}} \left( 2\tau + \frac{3}{2} \right)_{i_{\tau -1}}}{\left( \Delta_{\tau }^{-} \left( h\right) \right)_{i_{\tau -1}}\left( \Delta_{\tau }^{+} \left( h\right) \right)_{i_{\tau -1}}\left( 2\tau +1 \right)_{i_\tau } \left( 2\tau + \frac{3}{2} \right)_{i_{\tau }}} \eta ^{i_{\tau }}\right\} z^{\tau } \nonumber
\end{eqnarray}
where
\begin{equation}
\begin{cases} \tau \geq 2 \cr
\Delta_k^{\pm} \left( h\right) = 2k +\frac{2+\rho ^2 \pm \sqrt{(1+\rho ^2)h-\rho ^2}}{2(1+\rho ^2)}  \cr
z= -\rho ^2\xi ^2\cr
\eta = (1+\rho ^2) \xi
\end{cases}\nonumber
\end{equation}
\end{enumerate}   
\subsection{ ${\displaystyle  Hl(1-a,-q+\alpha \beta; \alpha,\beta, \delta, \gamma; 1-x)}$}  
 Replace coefficients $a$, $q$, $\gamma $, $\delta$ and $x$ by $1-a $, $-q+\alpha \beta $, $ \delta $, $ \gamma $ and $1-x$ into (\ref{eq:100086})--(\ref{eq:100091c}). Substitute (\ref{eq:100080}) into the new (\ref{eq:100086})--(\ref{eq:100091c}) with replacing the new $\alpha $, $\beta $ and $q$ by $-j$, $-j$ and $1/4h\rho ^{-2}-j(j+1/2)$ where $-j \in \mathbb{N}_{0}$ in order to make $A_n$ and $B_n$ terms terminated at specific index summation $n$: There are two possible values for a coefficient $\alpha $ in (\ref{eq:100086}) such as $\alpha =-1$ or 0; the two fixed values for $\alpha $ in (\ref{eq:100087a}) and (\ref{eq:100087b}) are give by $\alpha =-3$ or 2; two fixed constants for $\alpha $ in (\ref{eq:100088a}) and (\ref{eq:100088b}) are given by $\alpha =-2\left( 2N+5/2\right)$ or $2\left( 2N+2\right)$; two fixed values for $\alpha $ in (\ref{eq:100089a}) and (\ref{eq:100089b}) are described by $\alpha =-2\left( 2N+7/2\right)$ or $2\left( 2N+3\right)$. Replace $\alpha$ by $-2(j+1/2)$ in the new (\ref{eq:100090a})--(\ref{eq:100091c}). 
\begin{enumerate} 
\item As $\alpha = -1 $ or 0 with $h=h_0^0=0 $,
 
The eigenfunction is given by
\begin{equation}
 y(\sigma ) = Hl\left( 1-\rho ^{-2}, 0; 0,0,\frac{1}{2}, \frac{1}{2};\sigma \right) = 1 \nonumber
\end{equation}
\item As $\alpha = -3 $ or 2,

An algebraic equation of degree 2 for the determination of $h$ is given by
\begin{equation}
0 =  \frac{3}{4}\left( 1-\rho ^{-2} \right) + \prod_{l=0}^{1}\left( \frac{\rho ^{-2}}{4}h -\frac{3}{2} + l^2\left( 2-\rho ^{-2}\right)  \right)\nonumber
\end{equation}
The eigenvalue of $h$ is written by $h_1^m$ where $m = 0,1 $; $h_{1}^0 < h_{1}^1$. Its eigenfunction is given by
\begin{eqnarray}
 y(\sigma ) &=& Hl\left( 1-\rho ^{-2}, \frac{\rho ^{-2}}{4} h_1^m-\frac{3}{2}; -1,-1,\frac{1}{2}, \frac{1}{2};\sigma \right) \nonumber\\
&=& 1+\frac{ h_1^m \rho ^{-2}-6}{2(2-\rho ^{-2})} \eta \nonumber
\end{eqnarray}
\item As $\alpha  =-2\left( 2N+\frac{5}{2} \right)$ or $2\left( 2N+2 \right) $ where $N \in \mathbb{N}_{0}$,

An algebraic equation of degree $2N+3$ for the determination of $h$ is given by
\begin{equation}
0 = \sum_{r=0}^{N+1}\bar{c}\left( r, 2(N-r)+3; 2N+2,h\right)  \nonumber
\end{equation}
The eigenvalue of $h$ is written by $h_{2N+2}^m$ where $m = 0,1,2,\cdots,2N+2 $; $h_{2N+2}^0 < h_{2N+2}^1 < \cdots < h_{2N+2}^{2N+2}$. Its eigenfunction is given by 
\begin{eqnarray} 
 y(\sigma ) &=& Hl\left( 1-\rho ^{-2}, \frac{\rho ^{-2}}{4} h_{2N+2}^m -\left( 2N+2\right) \left( 2N+\frac{5}{2} \right); -\left( 2N+2\right), -\left( 2N+2\right),\frac{1}{2}, \frac{1}{2};\sigma \right) \nonumber\\
&=& \sum_{r=0}^{N+1} y_{r}^{2(N+1-r)}\left( 2N+2, h_{2N+2}^m; \sigma \right)  
\nonumber
\end{eqnarray}
\item As $\alpha  =-2\left( 2N+\frac{7}{2} \right)$ or $2\left( 2N+3 \right) $ where $N \in \mathbb{N}_{0}$,

An algebraic equation of degree $2N+4$ for the determination of $h$ is given by
\begin{equation}  
0 = \sum_{r=0}^{N+2}\bar{c}\left( r, 2(N+2-r); 2N+3,h\right) \nonumber
\end{equation}
The eigenvalue of $h$ is written by $h_{2N+3}^m$ where $m = 0,1,2,\cdots,2N+3 $; $h_{2N+3}^0 < h_{2N+3}^1 < \cdots < h_{2N+3}^{2N+3}$. Its eigenfunction is given by
\begin{eqnarray} 
 y(\sigma ) &=& Hl\left( 1-\rho ^{-2}, \frac{\rho ^{-2}}{4} h_{2N+3}^m -\left( 2N+3\right) \left( 2N+\frac{7}{2} \right); -\left( 2N+3\right), -\left( 2N+3\right),\frac{1}{2}, \frac{1}{2};\sigma \right)\nonumber\\
&=& \sum_{r=0}^{N+1} y_{r}^{2(N-r)+3} \left( 2N+3,h_{2N+3}^m; \sigma \right) \nonumber
\end{eqnarray}
In the above,
\begin{eqnarray}
\bar{c}(0,n;j,h)  &=& \frac{\left( \Delta_0^{-} \left( j,h\right) \right)_{n}\left( \Delta_0^{+} \left( j,h\right) \right)_{n}}{\left( 1 \right)_{n} \left( \frac{1}{2} \right)_{n}} \left( \frac{2-\rho ^{-2}}{1-\rho ^{-2}} \right)^{n}\nonumber\\
\bar{c}(1,n;j,h) &=& \left( \frac{-1}{1-\rho ^{-2}} \right) \sum_{i_0=0}^{n}\frac{\left( i_0 -j\right)\left( i_0+j+\frac{1}{2} \right) }{\left( i_0+2 \right) \left( i_0+\frac{3}{2} \right)} \frac{ \left( \Delta_0^{-} \left( j,h\right) \right)_{i_0}\left( \Delta_0^{+} \left( j,h\right) \right)_{i_0}}{\left( 1 \right)_{i_0} \left( \frac{1}{2} \right)_{i_0}} \nonumber\\
&\times&  \frac{ \left( \Delta_1^{-} \left( j,h\right) \right)_{n}\left( \Delta_1^{+} \left( j,h\right) \right)_{n} \left( 3 \right)_{i_0} \left( \frac{5}{2} \right)_{i_0}}{\left( \Delta_1^{-} \left( j,h\right) \right)_{i_0}\left( \Delta_1^{+} \left( j,h\right) \right)_{i_0}\left( 3 \right)_{n} \left( \frac{5}{2} \right)_{n}} \left( \frac{2-\rho ^{-2}}{1-\rho ^{-2}} \right)^{n }  
\nonumber\\
\bar{c}(\tau ,n;j,h) &=& \left( \frac{-1}{1-\rho ^{-2}} \right)^{\tau} \sum_{i_0=0}^{n}\frac{\left( i_0 -j\right)\left( i_0+j+\frac{1}{2} \right) }{\left( i_0+2 \right) \left( i_0+\frac{3}{2} \right)} \frac{ \left( \Delta_0^{-} \left( j,h\right) \right)_{i_0}\left( \Delta_0^{+} \left( j,h\right) \right)_{i_0}}{\left( 1 \right)_{i_0} \left( \frac{1}{2} \right)_{i_0}}  \nonumber\\
&\times& \prod_{k=1}^{\tau -1} \left( \sum_{i_k = i_{k-1}}^{n} \frac{\left( i_k+ 2k-j\right)\left( i_k +2k+j+\frac{1}{2} \right)}{\left( i_k+2k+2 \right) \left( i_k+2k+\frac{3}{2} \right)} \right. \nonumber\\
&\times&\left. \frac{ \left( \Delta_k^{-} \left( j,h\right) \right)_{i_k}\left( \Delta_k^{+} \left( j,h\right) \right)_{i_k} \left( 2k+1 \right)_{i_{k-1}} \left( 2k+ \frac{1}{2} \right)_{i_{k-1}}}{\left( \Delta_k^{-} \left( j,h\right) \right)_{i_{k-1}}\left( \Delta_k^{+} \left( j,h\right) \right)_{i_{k-1}}\left( 2k+1 \right)_{i_k} \left( 2k+ \frac{1}{2} \right)_{i_k}} \right) \nonumber\\
&\times& \frac{ \left( \Delta_{\tau }^{-} \left( j,h\right) \right)_{n}\left( \Delta_{\tau }^{+} \left( j,h\right) \right)_{n} \left( 2\tau +1 \right)_{i_{\tau -1}} \left( 2\tau + \frac{1}{2} \right)_{i_{\tau -1}}}{\left( \Delta_{\tau }^{-} \left( j,h\right) \right)_{i_{\tau -1}}\left( \Delta_{\tau }^{+} \left( j,h\right) \right)_{i_{\tau -1}}\left( 2\tau +1 \right)_{n} \left( 2\tau + \frac{1}{2} \right)_{n}} \left( \frac{2-\rho ^{-2}}{1-\rho ^{-2}} \right)^{n } \nonumber
\end{eqnarray}
\begin{eqnarray}
y_0^m(j,h;\sigma ) &=& \sum_{i_0=0}^{m} \frac{\left( \Delta_0^{-} \left( j,h\right) \right)_{i_0}\left( \Delta_0^{+} \left( j,h\right) \right)_{i_0}}{\left( 1 \right)_{i_0} \left( \frac{1}{2} \right)_{i_0}} \eta ^{i_0} \nonumber\\
y_1^m(j,h;\sigma ) &=& \left\{\sum_{i_0=0}^{m}\frac{\left( i_0 -j\right)\left( i_0+j+\frac{1}{2} \right) }{\left( i_0+2 \right) \left( i_0+\frac{3}{2} \right)} \frac{ \left( \Delta_0^{-} \left( j,h\right) \right)_{i_0}\left( \Delta_0^{+} \left( j,h\right) \right)_{i_0}}{\left( 1 \right)_{i_0} \left( \frac{1}{2} \right)_{i_0}} \right. \nonumber\\
&\times& \left. \sum_{i_1 = i_0}^{m} \frac{ \left( \Delta_1^{-} \left( j,h\right) \right)_{i_1}\left( \Delta_1^{+} \left( j,h\right) \right)_{i_1} \left( 3 \right)_{i_0} \left( \frac{5}{2} \right)_{i_0}}{\left( \Delta_1^{-} \left( j,h\right) \right)_{i_0}\left( \Delta_1^{+} \left( j,h\right) \right)_{i_0}\left( 3 \right)_{i_1} \left( \frac{5}{2} \right)_{i_1}} \eta ^{i_1}\right\} z 
\nonumber
\end{eqnarray}
\begin{eqnarray}
y_{\tau }^m(j,h;\sigma ) &=& \left\{ \sum_{i_0=0}^{m} \frac{\left( i_0 -j\right)\left( i_0+j+\frac{1}{2} \right) }{\left( i_0+2 \right) \left( i_0+\frac{3}{2} \right)} \frac{ \left( \Delta_0^{-} \left( j,h\right) \right)_{i_0}\left( \Delta_0^{+} \left( j,h\right) \right)_{i_0}}{\left( 1 \right)_{i_0} \left( \frac{1}{2} \right)_{i_0}} \right.\nonumber\\
&\times& \prod_{k=1}^{\tau -1} \left( \sum_{i_k = i_{k-1}}^{m} \frac{\left( i_k+ 2k-j\right)\left( i_k +2k+j+\frac{1}{2}\right)}{\left( i_k+2k+2 \right) \left( i_k+2k+\frac{3}{2} \right)} \right. \nonumber\\
&\times&\left. \frac{ \left( \Delta_k^{-} \left( j,h\right) \right)_{i_k}\left( \Delta_k^{+} \left( j,h\right) \right)_{i_k} \left( 2k+1 \right)_{i_{k-1}} \left( 2k+ \frac{1}{2} \right)_{i_{k-1}}}{\left( \Delta_k^{-} \left( j,h\right) \right)_{i_{k-1}}\left( \Delta_k^{+} \left( j,h\right) \right)_{i_{k-1}}\left( 2k+1 \right)_{i_k} \left( 2k+ \frac{1}{2} \right)_{i_k}} \right) \nonumber\\
&\times & \left. \sum_{i_{\tau } = i_{\tau -1}}^{m}  \frac{ \left( \Delta_{\tau }^{-} \left( j,h\right) \right)_{i_{\tau }}\left( \Delta_{\tau }^{+} \left( j,h\right) \right)_{i_{\tau }} \left( 2\tau +1  \right)_{i_{\tau -1}} \left( 2\tau + \frac{1}{2} \right)_{i_{\tau -1}}}{\left( \Delta_{\tau }^{-} \left( j,h\right) \right)_{i_{\tau -1}}\left( \Delta_{\tau }^{+} \left( j,h\right) \right)_{i_{\tau -1}}\left( 2\tau +1 \right)_{i_\tau } \left( 2\tau + \frac{1}{2} \right)_{i_{\tau }}} \eta ^{i_{\tau }}\right\} z^{\tau } \nonumber
\end{eqnarray}
where
\begin{equation}
\begin{cases} \tau \geq 2 \cr
\Delta_k^{\pm} \left( j,h\right) = 2k \pm \sqrt{ \frac{j\left( j+\frac{1}{2}\right) -\frac{1}{4}h\rho ^{-2} }{2-\rho ^{-2}}}  \cr
\sigma = 1-\xi \cr
z=  \frac{-1}{1-\rho ^{-2}} \sigma ^2\cr
\eta = \frac{2-\rho ^{-2}}{1-\rho ^{-2}} \sigma
\end{cases}\nonumber
\end{equation}
\end{enumerate}  
\subsection{\footnotesize ${\displaystyle (1-x)^{1-\delta } Hl(1-a,-q+(\delta -1)\gamma a+(\alpha -\delta +1)(\beta -\delta +1); \alpha-\delta +1,\beta-\delta +1, 2-\delta, \gamma; 1-x)}$ \normalsize}  
Replace coefficients $a$, $q$, $\alpha$, $\beta$, $\gamma $, $\delta$ and $x$ by $1-a$, $-q+(\delta -1)\gamma a+(\alpha -\delta +1)(\beta -\delta +1)$, $\alpha-\delta +1$, $\beta-\delta +1 $, $2 - \delta$, $\gamma $ and $1-x$ into (\ref{eq:100086})--(\ref{eq:100091c}). Multiply (\ref{eq:100086}), (\ref{eq:100087b}), (\ref{eq:100088b}) and (\ref{eq:100089b}) by $(1-x)^{1-\delta }$. Substitute (\ref{eq:100080}) into the new (\ref{eq:100086})--(\ref{eq:100091c}) with replacing the new $\alpha $, $\beta $ and $q$ by $-j$, $-j$ and  $1/4(h-1)\rho ^{-2}-j(j+3/2)$ by $-j$ where $-j \in \mathbb{N}_{0}$ in order to make $A_n$ and $B_n$ terms terminated at specific index summation $n$: There are two possible values for a coefficient $\alpha $ in (\ref{eq:100086}) such as $\alpha =-2$ or 1; the two fixed values for $\alpha $ in (\ref{eq:100087a}) and (\ref{eq:100087b}) are give by $\alpha =-4$ or 3; two fixed constants for $\alpha $ in (\ref{eq:100088a}) and (\ref{eq:100088b}) are given by $\alpha =-2\left( 2N+3\right)$ or $2\left( 2N+5/2\right)$; two fixed values for $\alpha $ in (\ref{eq:100089a}) and (\ref{eq:100089b}) are described by $\alpha =-2\left( 2N+4\right)$ or $2\left( 2N+7/2\right)$. Replace $\alpha$ by $-2(j+1)$ in the new (\ref{eq:100090a})--(\ref{eq:100091c}).
\begin{enumerate} 
\item As $\alpha = -2 $ or 1 with $h=h_0^0=1 $,
 
The eigenfunction is given by
\begin{equation}
\sigma^{\frac{1}{2}} y(\sigma ) = \sigma^{\frac{1}{2}} Hl\left( 1-\rho ^{-2}, 0; 0,0,\frac{3}{2}, \frac{1}{2};\sigma \right) = \sigma^{\frac{1}{2}} \nonumber
\end{equation}
\item As $\alpha = -4 $ or 3,

An algebraic equation of degree 2 for the determination of $h$ is given by
\begin{equation}
0 =  \frac{15}{4}\left( 1-\rho ^{-2} \right) + \prod_{l=0}^{1}\left( \frac{\rho ^{-2}}{4} (h-1) -\frac{5}{2}+l(l+1)\left( 2-\rho ^{-2}\right)  \right)\nonumber
\end{equation}
The eigenvalue of $h$ is written by $h_1^m$ where $m = 0,1 $; $h_{1}^0 < h_{1}^1$. Its eigenfunction is given by
\begin{eqnarray}
\sigma^{\frac{1}{2}} y(\sigma ) &=& \sigma^{\frac{1}{2}} Hl\left( 1-\rho ^{-2}, \frac{\rho ^{-2}}{4} (h_1^m-1); -1,-1,\frac{3}{2}, \frac{1}{2};\sigma \right) \nonumber\\
&=& \sigma^{\frac{1}{2}} \left\{ 1+\frac{ (h_1^m-1)\rho ^{-2} -10}{6(2-\rho ^{-2})} \eta \right\} \nonumber
\end{eqnarray}
\item As $\alpha  =-2\left( 2N+3 \right)$ or $2\left( 2N+\frac{5}{2} \right) $ where $N \in \mathbb{N}_{0}$,

An algebraic equation of degree $2N+3$ for the determination of $h$ is given by
\begin{equation}
0 = \sum_{r=0}^{N+1}\bar{c}\left( r, 2(N-r)+3; 2N+2,h\right)  \nonumber
\end{equation}
The eigenvalue of $h$ is written by $h_{2N+2}^m$ where $m = 0,1,2,\cdots,2N+2 $; $h_{2N+2}^0 < h_{2N+2}^1 < \cdots < h_{2N+2}^{2N+2}$. Its eigenfunction is given by 
\begin{eqnarray} 
\sigma^{\frac{1}{2}} y(\sigma ) &=& \sigma^{\frac{1}{2}} Hl\left( 1-\rho ^{-2}, \frac{\rho ^{-2}}{4}\left( h_{2N+2}^m -1\right)-\left( 2N+2\right) \left( 2N+\frac{7}{2}\right); -\left( 2N+2\right), -\left( 2N+2\right),\frac{3}{2}, \frac{1}{2};\sigma \right) \nonumber\\
&=& \sigma^{\frac{1}{2}} \sum_{r=0}^{N+1} y_{r}^{2(N+1-r)}\left( 2N+2, h_{2N+2}^m; \sigma \right)  
\nonumber
\end{eqnarray}
\item As $\alpha  =-2\left( 2N+4 \right)$ or $2\left( 2N+\frac{7}{2} \right) $ where $N \in \mathbb{N}_{0}$,

An algebraic equation of degree $2N+4$ for the determination of $h$ is given by
\begin{equation}  
0 = \sum_{r=0}^{N+2}\bar{c}\left( r, 2(N+2-r); 2N+3,h\right) \nonumber
\end{equation}
The eigenvalue of $h$ is written by $h_{2N+3}^m$ where $m = 0,1,2,\cdots,2N+3 $; $h_{2N+3}^0 < h_{2N+3}^1 < \cdots < h_{2N+3}^{2N+3}$. Its eigenfunction is given by
\begin{eqnarray} 
\sigma^{\frac{1}{2}} y(\sigma ) &=& \sigma^{\frac{1}{2}} Hl\left( 1-\rho ^{-2}, \frac{\rho ^{-2}}{4}\left( h_{2N+3}^m -1\right)-\left( 2N+3\right) \left( 2N+\frac{9}{2}\right); -\left( 2N+3\right), -\left( 2N+3\right),\frac{3}{2}, \frac{1}{2};\sigma \right)\nonumber\\
&=&  \sigma^{\frac{1}{2}} \sum_{r=0}^{N+1} y_{r}^{2(N-r)+3} \left( 2N+3,h_{2N+3}^m; \sigma \right) \nonumber
\end{eqnarray}
In the above,
\begin{eqnarray}
\bar{c}(0,n;j,h)  &=& \frac{\left( \Delta_0^{-} \left( j,h\right) \right)_{n}\left( \Delta_0^{+} \left( j,h\right) \right)_{n}}{\left( 1 \right)_{n} \left( \frac{3}{2} \right)_{n}} \left( \frac{2-\rho ^{-2}}{1-\rho ^{-2}} \right)^{n}\nonumber\\
\bar{c}(1,n;j,h) &=& \left( \frac{-1}{1-\rho ^{-2}}\right) \sum_{i_0=0}^{n}\frac{\left( i_0 -j\right)\left( i_0+j+\frac{3}{2} \right) }{\left( i_0+2 \right) \left( i_0+\frac{5}{2} \right)} \frac{ \left( \Delta_0^{-} \left( j,h\right) \right)_{i_0}\left( \Delta_0^{+} \left( j,h\right) \right)_{i_0}}{\left( 1 \right)_{i_0} \left( \frac{3}{2} \right)_{i_0}} \nonumber\\
&\times&  \frac{ \left( \Delta_1^{-} \left( j,h\right) \right)_{n}\left( \Delta_1^{+} \left( j,h\right) \right)_{n} \left( 3 \right)_{i_0} \left( \frac{7}{2} \right)_{i_0}}{\left( \Delta_1^{-} \left( j,h\right) \right)_{i_0}\left( \Delta_1^{+} \left( j,h\right) \right)_{i_0}\left( 3 \right)_{n} \left( \frac{7}{2} \right)_{n}} \left(\frac{2-\rho ^{-2}}{1-\rho ^{-2}} \right)^{n }  
\nonumber\\
\bar{c}(\tau ,n;j,h) &=& \left( \frac{-1}{1-\rho ^{-2}}\right)^{\tau} \sum_{i_0=0}^{n}\frac{\left( i_0 -j\right)\left( i_0+j+\frac{3}{2} \right) }{\left( i_0+2 \right) \left( i_0+\frac{5}{2} \right)} \frac{ \left( \Delta_0^{-} \left( j,h\right) \right)_{i_0}\left( \Delta_0^{+} \left( j,h\right) \right)_{i_0}}{\left( 1 \right)_{i_0} \left( \frac{3}{2} \right)_{i_0}}  \nonumber\\
&\times& \prod_{k=1}^{\tau -1} \left( \sum_{i_k = i_{k-1}}^{n} \frac{\left( i_k+ 2k-j\right)\left( i_k +2k+j+\frac{3}{2} \right)}{\left( i_k+2k+2 \right) \left( i_k+2k+\frac{5}{2} \right)} \right. \nonumber\\
&\times&\left. \frac{ \left( \Delta_k^{-} \left( j,h\right) \right)_{i_k}\left( \Delta_k^{+} \left( j,h\right) \right)_{i_k} \left( 2k+1 \right)_{i_{k-1}} \left( 2k+ \frac{3}{2} \right)_{i_{k-1}}}{\left( \Delta_k^{-} \left( j,h\right) \right)_{i_{k-1}}\left( \Delta_k^{+} \left( j,h\right) \right)_{i_{k-1}}\left( 2k+1 \right)_{i_k} \left( 2k+ \frac{3}{2} \right)_{i_k}} \right) \nonumber\\
&\times& \frac{ \left( \Delta_{\tau }^{-} \left( j,h\right) \right)_{n}\left( \Delta_{\tau }^{+} \left( j,h\right) \right)_{n} \left( 2\tau +1 \right)_{i_{\tau -1}} \left( 2\tau + \frac{3}{2} \right)_{i_{\tau -1}}}{\left( \Delta_{\tau }^{-} \left( j,h\right) \right)_{i_{\tau -1}}\left( \Delta_{\tau }^{+} \left( j,h\right) \right)_{i_{\tau -1}}\left( 2\tau +1 \right)_{n} \left( 2\tau + \frac{3}{2} \right)_{n}} \left( \frac{2-\rho ^{-2}}{1-\rho ^{-2}} \right)^{n } \nonumber
\end{eqnarray}
\begin{eqnarray}
y_0^m(j,h;\sigma ) &=& \sum_{i_0=0}^{m} \frac{\left( \Delta_0^{-} \left( j,h\right) \right)_{i_0}\left( \Delta_0^{+} \left( j,h\right) \right)_{i_0}}{\left( 1 \right)_{i_0} \left( \frac{3}{2} \right)_{i_0}} \eta ^{i_0} \nonumber\\
y_1^m(j,h;\sigma ) &=& \left\{\sum_{i_0=0}^{m}\frac{\left( i_0 -j\right)\left( i_0+j+\frac{3}{2} \right) }{\left( i_0+2 \right) \left( i_0+\frac{5}{2} \right)} \frac{ \left( \Delta_0^{-} \left( j,h\right) \right)_{i_0}\left( \Delta_0^{+} \left( j,h\right) \right)_{i_0}}{\left( 1 \right)_{i_0} \left( \frac{3}{2} \right)_{i_0}} \right. \nonumber\\
&\times& \left. \sum_{i_1 = i_0}^{m} \frac{ \left( \Delta_1^{-} \left( j,h\right) \right)_{i_1}\left( \Delta_1^{+} \left( j,h\right) \right)_{i_1} \left( 3 \right)_{i_0} \left( \frac{7}{2} \right)_{i_0}}{\left( \Delta_1^{-} \left( j,h\right) \right)_{i_0}\left( \Delta_1^{+} \left( j,h\right) \right)_{i_0}\left( 3 \right)_{i_1} \left( \frac{7}{2} \right)_{i_1}} \eta ^{i_1}\right\} z 
\nonumber
\end{eqnarray}
\begin{eqnarray}
y_{\tau }^m(j,h;\sigma ) &=& \left\{ \sum_{i_0=0}^{m} \frac{\left( i_0 -j\right)\left( i_0+j+\frac{3}{2} \right) }{\left( i_0+2 \right) \left( i_0+\frac{5}{2} \right)} \frac{ \left( \Delta_0^{-} \left( j,h\right) \right)_{i_0}\left( \Delta_0^{+} \left( j,h\right) \right)_{i_0}}{\left( 1 \right)_{i_0} \left( \frac{3}{2} \right)_{i_0}} \right.\nonumber\\
&\times& \prod_{k=1}^{\tau -1} \left( \sum_{i_k = i_{k-1}}^{m} \frac{\left( i_k+ 2k-j\right)\left( i_k +2k+j+\frac{3}{2}\right)}{\left( i_k+2k+2 \right) \left( i_k+2k+\frac{5}{2} \right)} \right. \nonumber\\
&\times&\left. \frac{ \left( \Delta_k^{-} \left( j,h\right) \right)_{i_k}\left( \Delta_k^{+} \left( j,h\right) \right)_{i_k} \left( 2k+1 \right)_{i_{k-1}} \left( 2k+ \frac{3}{2} \right)_{i_{k-1}}}{\left( \Delta_k^{-} \left( j,h\right) \right)_{i_{k-1}}\left( \Delta_k^{+} \left( j,h\right) \right)_{i_{k-1}}\left( 2k+1 \right)_{i_k} \left( 2k+ \frac{3}{2} \right)_{i_k}} \right) \nonumber\\
&\times & \left. \sum_{i_{\tau } = i_{\tau -1}}^{m}  \frac{ \left( \Delta_{\tau }^{-} \left( j,h\right) \right)_{i_{\tau }}\left( \Delta_{\tau }^{+} \left( j,h\right) \right)_{i_{\tau }} \left( 2\tau +1  \right)_{i_{\tau -1}} \left( 2\tau + \frac{3}{2} \right)_{i_{\tau -1}}}{\left( \Delta_{\tau }^{-} \left( j,h\right) \right)_{i_{\tau -1}}\left( \Delta_{\tau }^{+} \left( j,h\right) \right)_{i_{\tau -1}}\left( 2\tau +1 \right)_{i_\tau } \left( 2\tau + \frac{3}{2} \right)_{i_{\tau }}} \eta ^{i_{\tau }}\right\} z^{\tau } \nonumber
\end{eqnarray}
where
\begin{equation}
\begin{cases} \tau \geq 2 \cr
\Delta_k^{\pm} \left( j,h\right) = 2k +\frac{1}{2}\left\{ 1 \pm \sqrt{1-\frac{(h-1)\rho ^{-2}-4j\left( j+\frac{3}{2}\right)}{2-\rho ^{-2}}} \right\}  \cr
\sigma = 1-\xi \cr
z= \frac{-1}{1-\rho ^{-2}}\sigma ^2\cr
\eta = \frac{2-\rho ^{-2}}{1-\rho ^{-2}} \sigma
\end{cases}\nonumber
\end{equation}
\end{enumerate}  
\end{appendices} 

\addcontentsline{toc}{section}{Bibliography}
\bibliographystyle{model1a-num-names}
\bibliography{<your-bib-database>}
 
\chapter{Double Confluent Heun functions using reversible three term recurrence formula}
\chaptermark{DCH functions using R3TRF} 
 
Double Confluent Heun (DCH) equation is one of 4 confluent forms of Heun's differential equation \cite{11Heun1889,11Ronv1995}; this has irregular singularities at the origin and infinity, each of rank 1.
Power series solutions of the DCH equation (DCHE) provide a 3-term recurrence relation between successive coefficients. 

In this chapter I apply reversible three term recurrence formula (R3TRF) in chapter 1 of Ref.\cite{11Choun2013} to power series expansions of the DCHE with two irregular singularities for a polynomial of type 2. Their combined definite and contour integrals are derived analytically including generating functions for DCH polynomials of type 2.
 
\section{Introduction}
Heun's equation is a second-order linear ordinary differential equation of the form \cite{11Heun1889,11Ronv1995}
\begin{equation}
\frac{d^2{y}}{d{x}^2} + \left(\frac{\gamma }{x} +\frac{\delta }{x-1} + \frac{\epsilon }{x-a}\right) \frac{d{y}}{d{x}} +  \frac{\alpha \beta x-q}{x(x-1)(x-a)} y = 0 \label{eq:11001}
\end{equation}
with the condition $\epsilon = \alpha +\beta -\gamma -\delta +1$. The parameters play different roles: $a \ne 0 $ is the singularity parameter, $\alpha $, $\beta $, $\gamma $, $\delta $, $\epsilon $ are exponent parameters, $q$ is the accessory parameter. Also, $\alpha $ and $\beta $ are identical to each other. The total number of free parameters is six. It has four regular singular points which are 0, 1, $a$ and $\infty $ with exponents $\{ 0, 1-\gamma \}$, $\{ 0, 1-\delta \}$, $\{ 0, 1-\epsilon \}$ and $\{ \alpha, \beta \}$.

The DCH equation (DCHE), one of confluent forms of Heun equation, is applicable to areas such as Dirac equations in the Nutku's helicoid spacetime \cite{11Birk2007,11Hort2007} and the rotating Bertotti-Robinson spacetime \cite{11AlBa2008}, the Schr$\ddot{\mbox{o}}$dinger equation for inverse fourth and sixth-power potentials \cite{11Figu2005,11Fran1971,11Klei1968}, Klein-Gordon equation in curved spacetimes \cite{11Beze2014,11Birr1982,11Cost1989,11Nove1979}, scattering theory in non-relativistic quantum mechanics, Teukolsky's equations in general relativity, etc.

In general, as far as I know, there are three types of the DCHE such as (1) non-symmetrical canonical form of the DCHE, (2) canonical form of the general DCHE and (3) generalized spheroidal equation in the Leaver version as a regular singular point $x_0 \rightarrow 0$.

\subsection{Non-symmetrical canonical form of the DCHE}
Heun equation has the four kinds of confluent forms: (1) Confluent Heun (two regular and one irregular singularities), (2) Double Confluent Heun (two irregular singularities), (3) Biconfluent Heun (one regular and one irregular singularities), (4) Triconfluent Heun equations (one irregular singularity).
We can derive these four confluent forms from Heun equation by combining two or more regular singularities to take form an irregular singularity. Its process, converting Heun equation to other confluent forms, is similar to deriving of confluent hypergeometric equation from the hypergeometric equation.

For confluent Heun equation (CHE), divide (\ref{eq:11001}) by $a$. Let $a \rightarrow \infty $, and simultaneously saying $\beta, \epsilon, q\rightarrow \infty $ in such a way that $\beta /a, \epsilon /a \rightarrow -\beta $ and $q/a \rightarrow -q$ in the new (\ref{eq:11001}).
\begin{equation}
\frac{d^2{y}}{d{x}^2} + \left(\beta  +\frac{\gamma }{x} + \frac{\delta }{x-1}\right) \frac{d{y}}{d{x}} +  \frac{\alpha \beta x-q}{x(x-1)} y = 0 \label{eq:11002}
\end{equation}
(\ref{eq:11002}) is defined as the non-symmetrical canonical form of the CHE \cite{11Deca1978,11Decar1978,11Ronv1995}. It has three singular points: two regular singular points which are 0 and 1 with exponents $\{0, 1-\gamma\}$ and $\{0, 1-\delta \}$, and one irregular singular point which is $\infty$ with an exponent $\alpha$.   

We are able to derive the DCHE by changing coefficients and combining two regular singularities. Let us allow $x\rightarrow x/\epsilon $, $\beta \rightarrow \beta \epsilon $, $\gamma \rightarrow \gamma -\delta /\epsilon $ and $\delta \rightarrow \delta /\epsilon $ in (\ref{eq:11002}). Assuming $\epsilon \rightarrow 0$ in the new (\ref{eq:11002})
\begin{equation}
\frac{d^2{y}}{d{x}^2} + \left(\beta  +\frac{\gamma }{x} + \frac{\delta }{x^2}\right) \frac{d{y}}{d{x}} +  \frac{\alpha \beta x-q}{x^2} y = 0 \label{eq:11003}
\end{equation}
(\ref{eq:11003}) is the non-symmetrical canonical form of the DCHE which has irregular singularities at $x=0$ and $\infty $, each of rank 1. Its solution is denoted as $H_d(\alpha,\beta,\gamma,\delta,q;x)$.
 For DLFM version \cite{11NIST} or in Ref.\cite{11Slavy2000}, replace $\beta $ by 1 in (\ref{eq:11003}). The reason, why the parameter $\beta $ is included in (\ref{eq:11003}) instead of the unity, is that we can obtain an equation of the Whittaker-Ince limit of the DCHE by putting $\beta \rightarrow 0$, $\alpha \rightarrow \infty $, such that $\alpha \beta \rightarrow \alpha $ in (\ref{eq:11003}).
\subsection{Canonical form of the general DCHE}
The canonical form of the general DCHE is given by \cite{11Ronv1995}
\begin{equation}
 \frac{d^2{y}}{d{x}^2} + \left( \alpha _1 + \frac{1}{x}+ \frac{\alpha _{-1}}{x^2}\right) \frac{d{y}}{d{x}} + \left( \frac{\left( B_1 +\frac{\alpha _1}{2}\right)}{x} +\frac{\left( B_0 +\frac{\alpha _1\alpha _{-1}}{2}\right)}{x^2} +\frac{\left( B_{-1} -\frac{ \alpha _{-1}}{2}\right)}{x^3}\right) y = 0 \label{eq:11004}
\end{equation}
where $\alpha _1, \alpha _{-1}, B_{-1}, B_0, B_1$ are arbitrary complex parameter. If $\alpha _1, \alpha _{-1}\rightarrow \alpha $, $B_1 \rightarrow \alpha \beta _1$, $B_{-1} \rightarrow \alpha \beta _{-1}$ and $B_0\rightarrow -\gamma $, it turns to be the symmetric canonical form of the DCHE. We transform (\ref{eq:11004}) by
\begin{equation}
y(x) = x^{\frac{1}{2}-\frac{B_{-1}}{\alpha _{-1}}}\; \tilde{y}(x) \nonumber
\end{equation}
And we obtain the second-order linear ODE such as
\begin{equation}
 \frac{d^2{\tilde{y}}}{d{x}^2} + \left( \alpha _1 + \frac{2\left( 1- \frac{B_{-1}}{\alpha _{-1}}\right)}{x}+ \frac{\alpha _{-1}}{x^2}\right) \frac{d{\tilde{y}}}{d{x}} + \frac{\left( B_1 + \alpha _1 - \frac{\alpha _1 B_{-1}}{\alpha _{-1}}\right) x + \left( \frac{1}{2}-\frac{B_{-1}}{\alpha _{-1}} \right)^2 + B_0 +\frac{\alpha _1 \alpha _{-1}}{2}}{x^2} \tilde{y} = 0 \label{eq:11005}
\end{equation}
As we compare all coefficients in (\ref{eq:11005}) with (\ref{eq:11003}), we find that
\begin{equation}
\begin{cases} \alpha \rightarrow  \frac{B_1}{\alpha _1} -\frac{B_{-1}}{\alpha _{-1}} +1  \cr
\beta \rightarrow  \alpha _1  \cr
\gamma \rightarrow  2 \left( 1- \frac{B_{-1}}{\alpha _{-1}}\right) \cr
\delta \rightarrow \alpha _{-1}  \cr
q \rightarrow -\frac{\alpha _1 \alpha _{-1}}{2} -B_0 -\left( \frac{1}{2}-\frac{B_{-1}}{\alpha _{-1}} \right)^2   
\end{cases}\nonumber  
\end{equation}
Then the solution $y(x)$ for a canonical form of the general DCHE is described by the non-symmetrical canonical form of the DCHE.
\begin{align}
&y\left( \alpha _1, \alpha _{-1}, B_{-1}, B_0, B_1; x\right) \nonumber\\ 
&= x^{\frac{1}{2}-\frac{B_{-1}}{\alpha _{-1}}}\; H_d\left( \frac{B_1}{\alpha _1} -\frac{B_{-1}}{\alpha _{-1}} +1, \alpha _1, 2 \left( 1- \frac{B_{-1}}{\alpha _{-1}}\right), \alpha _{-1}, -\frac{\alpha _1 \alpha _{-1}}{2} -B_0 -\left( \frac{1}{2}-\frac{B_{-1}}{\alpha _{-1}} \right)^2; x\right) \label{eq:11006}
\end{align}
\subsection{Generalized spheroidal equation in the Leaver version}
The generalized spheroidal wave equation (GSWE), one of standard forms of the CHE, in the Leaver version reads \cite{11Leav1986}
\begin{equation}
x(x-x_0)\frac{d^2{y}}{d{x}^2} + (B_1+B_2 x)\frac{d{y}}{d{x}}+ \left( B_3 -2\eta \omega (x-x_0)+\omega ^2x(x-x_0)\right) y(x) = 0 \hspace{.5cm}\mbox{where}\; \omega \ne0 \label{eq:11007}
\end{equation}
where $B_1, B_2, B_3, \eta $ and $\omega $ are arbitrary constant. It is one of standard from of the CHE having three singular points; two regular singular points which are 0 and $x_0$ with exponents $\left\{ 0, 1+\frac{B_1}{x_0}\right\}$ and $\left\{ 0, 1-B_2-\frac{B_1}{x_0}\right\}$, and one irregular singular point which is $\infty $.
Taking $x_0\rightarrow 0$ in (\ref{eq:11007}), the GSWE in the Leaver version gives the following the DCHE with five parameters
\begin{equation}
x^2\frac{d^2{y}}{d{x}^2} + (B_1+B_2 x)\frac{d{y}}{d{x}}+ \left( B_3 -2\eta \omega x+\omega ^2x^2\right) y(x) = 0 \hspace{.5cm}\mbox{where}\; B_1 \ne0, \omega \ne0 \label{eq:11008}
\end{equation}
To permit the non-symmetrical canonical form of the DCHE limit, the first step is given by the substitution
\begin{equation}
y(x) = \exp(\pm i\omega x) \tilde{y}(x) \nonumber
\end{equation}
which convert (\ref{eq:11008}) into
\begin{equation}
 \frac{d^2{\tilde{y}}}{d{x}^2} + \left( \pm 2i\omega + \frac{B_2}{x} + \frac{B_1}{x^2} \right) \frac{d{\tilde{y}}}{d{x}} + \frac{\omega \left( -2\eta \pm i B_2\right) x+ B_3 \pm i\omega B_1}{x^2} \tilde{y} = 0 \label{eq:11009}
\end{equation}
As we compare all coefficients in (\ref{eq:11009}) with (\ref{eq:11003}), we find that
\begin{equation}
\begin{cases} \alpha \rightarrow \frac{B_2}{2} \pm i\eta \cr
\beta \rightarrow \pm 2i\omega \cr
\gamma \rightarrow B_2\cr
\delta \rightarrow B_1\cr
q \rightarrow  -\left( B_3 \pm i\omega B_1\right)  
\end{cases}\nonumber  
\end{equation}
Then the solution $y(x)$ for the GSWE in the Leaver version as $x_0 =0$ is converted to the non-symmetrical canonical form of the DCHE.
\begin{align}
&y\left( x_0 =0; B_1, B_2, B_3, \omega , \eta ; x\right) \nonumber\\ 
&= \exp(\pm i\omega x)\; H_d\left( \frac{B_2}{2} \pm i\eta, \pm 2i\omega, B_2, B_1, -\left( B_3 \pm i\omega B_1\right); x\right) \label{eq:110010}
\end{align}
\section{A 3-term recurrence relation of the DCHE}
Generally, by substituting a power series with unknown coefficients into linear ODEs, the recurrence relation of coefficients starts to appear. There can be between two and infinity term in the recursion relation of a Frobenius solution. 

For instead, a hypergeometric equation is a Fuchsian differential equation, having three regular singular points at 0, 1 and $\infty $ with exponents $\{ 0,1-\gamma \}$, $\{ 0,\gamma -\alpha -\beta \}$ and $\{ \alpha ,\beta \}$. 
\begin{equation}
x(1-x)\frac{d^2{y}}{d{x}^2} + \left( \gamma -(\alpha +\beta +1)x\right)\frac{d{y}}{d{x}}-\alpha \beta y =0\nonumber
\end{equation}
All well-known special functions such as Bessel, Legendre, Laguerre, Kummer functions, etc are just the special case of a hypergeometric function.
By putting a function $y(x)=\sum_{n=0}^{\infty }c_n x^{n+\lambda }$ into the above equation, a 2-term recursion relation between successive coefficients arise such as
\begin{equation}
c_{n+1}= A_n \;c_n  \hspace{1cm};n\geq 0
\nonumber
\end{equation}
where,
\begin{equation}
A_n = \frac{(n+\alpha +\lambda )(n+\beta +\lambda)}{(n+1+\lambda )(n+\gamma +\lambda)}
\nonumber
\end{equation}
where $ c_0 \ne 0$ and $\lambda $ is an indicial root.
Hitherto we have described tremendous differential equations and explained physical phenomena with hypergeometric-type functions, having a 2-term recursion relation between successive coefficients, for convenient computations; the formal series solutions in closed forms have been analyzed including its definite or contour integrals. 

However, since modern (particle) physics come into existence, around the $21^{th}$ century, we can not depict and make an understanding of the nature with linear ODEs, especially second-order linear ODEs, having a 2-term recursion relation in general solutions in series any more. 
The most complex and obscure problems in physics such as QCD, SUSY, general relativity, string theory, etc require more than a 3-term in a recurrence relation for analytic solutions in series \cite{11Hortacsu:2011rr}. 
A 3-term recursive relation between successive coefficients in linear ODEs make difficulties for mathematical computations and its numerical analysis. Its methods of proof for general summation schemes seem obscure including definite or contour integrals \cite{11Erde1955,11Hobs1931,11Whit1914,11Whit1952}.
  
Heun's equation including its confluent forms, except triconfluent Heun one (its recursion relation consists of a 4-term), is one of examples for a 3-term recurrence relation. Heun equation is considered as the mother of all well-known linear ODEs such as spheroidal wave, Lame, Mathieu, and hypergeometric-type equations.
Even though its solutions are treated as the most outstanding analytic functions in among every special functions, the general summation solutions of it including integrals have been unknown until now. Unlike formal series solutions of a hypergeoemtric equation, the solutions in series of Heun equation is just left as solutions of recurrence because a 3-term in the recurrence relation.  

By the method of a power series solution into linear ODEs, assuming $y(x)= \sum_{n=0}^{\infty } c_n x^{n+\lambda }$, the recurrence relation for the coefficients starts to appear. For a 3-term case, its recurrence relation is given by
\begin{equation}
c_{n+1}=A_n \;c_n +B_n \;c_{n-1} \hspace{1cm};n\geq 1
\label{ysc:1}
\end{equation}
where $c_1= A_0 \;c_0$ and $ c_0 \ne 0$ . 

Any linear ODEs having a 2-term recurrence relation have two general solutions in series such as an infinite series and a polynomial.
Contradistinctively, $2^{3-1}$ possible power series solutions for a 3-term recursion relation arise such as an infinite series and 3 types of polynomials: (1) a polynomial which makes $B_n$ term terminated; $A_n$ term is not terminated, (2) a polynomial which makes $A_n$ term terminated; $B_n$ term is not terminated, (3) a polynomial which makes $A_n$ and $B_n$ terms terminated at the same time, referred as `a complete polynomial.'  

In this chapter, by putting $A_n$ and $B_n$ terms in a recurrence relation of a non-symmetrical canonical form of the DCHE around $x=0$ and $x=\infty $ into general summation formulas such as R3TRF in chapter 1 of Ref.\cite{11Choun2013}, power series solutions of the DCHE are constructed for a polynomial of type 2. The combined definite and contour integrals of DCH polynomials of type 2 are derived analytically including their generating functions.\footnote{(\ref{eq:11003}) is preferred as an analytic solution of the DCHE rather than a canonical form of the general DCHE and the GSWE in the Leaver version because of more convenient for mathematical computations and the application of the Laplace transform.}
\section[The DCHE with a irregular singular point at $x=0$]{The DCHE with a irregular singular point at the origin} 
\subsection{Power series for a polynomial of type 2}
Assume that an solution in series of (\ref{eq:11003}) is written by
\begin{equation}
y(x)= \sum_{n=0}^{\infty } c_n x^{n+\lambda }  \label{eq:110011}
\end{equation}
where $c_0 \ne 0$. Substituting (\ref{eq:110011}) into (\ref{eq:11003}) gives for the coefficients $c_n$ the recurrence relations
\begin{equation}
c_{n+1}=A_n \;c_n +B_n \;c_{n-1} \hspace{1cm};n\geq 1 \label{eq:110012}
\end{equation}
where,
\begin{subequations}
\begin{align}
 A_n &= -\frac{n\left( n-1+\gamma \right) -q}{\delta \left( n+1\right)} \label{eq:110013a}\\ 
&= -\frac{\left( n+\frac{\gamma -1-\sqrt{\left( \gamma -1\right)^2 +4q}}{2}\right) \left( n+\frac{\gamma -1+\sqrt{\left( \gamma -1\right)^2 +4q}}{2}\right)}{\delta \left( n+1\right)} \label{eq:110013b}\\
B_n &= -\frac{\beta \left( n-1+\alpha \right)}{\delta \left( n+1\right)}  \label{eq:110013c}\\
c_1 &= A_0 \;c_0 \label{eq:110013d}
\end{align}
\end{subequations} 
with $\delta \ne 0$. We only have one indicial root such as $\lambda = 0$.

Now let's test for convergence of the analytic function $y(x)$. As $n\gg 1$ (for sufficiently large), (\ref{eq:110012})--(\ref{eq:110013d}) are given by
\begin{subequations}
\begin{equation}
c_{n+1}=A\;c_n +B\;c_{n-1} \hspace{1cm}\mbox{as}\;\; n\geq 1 \label{eq:110014}
\end{equation}
where

\begin{tabular}{p{8cm}p{6cm}}
{\begin{align}
&\lim_{n\gg 1} A_n = A=  -\frac{n}{\delta } \rightarrow \infty \label{eq:110015a}\ 
\end{align}}
&
{\begin{align}
&\lim_{n\gg 1} B_n = B= -\frac{\beta }{\delta } \label{eq:110015b}\ 
\end{align} }
\end{tabular}
\end{subequations}

by letting $c_1\sim  A c_0$.\footnote{We only have the sense of curiosity about an asymptotic series as  $n\gg 1$ for given $x$. Actually, $c_1 =  A_0 c_0$. But for a huge value of an index $n$, I treat the coefficient $c_1$ as $ A c_0$ for simple computations.} 
There are no general solutions in series for a polynomial of type 1 and an infinite series. Because the $y(x)$ is divergent as $n\gg 1$ in (\ref{eq:110015a}). So there are only two types of formal series solutions of the DCHE about $x=0$ such as a polynomial of type 2 and a complete polynomial. For a polynomial of type 2, I treat $\alpha $, $\beta $, $\gamma $, $\delta $ as free variables and $q$ as a fixed value. For a complete polynomial, I treat $\beta $, $\gamma $, $\delta $ as free variables and $\alpha $, $q$ as fixed values. 

For a polynomial of type 2, (\ref{eq:110015b}) is only available for an asymptotic behavior of the minimum $y(x)$; (\ref{eq:110015a}) is negligible for the minimum value of a $y(x)$ because $A_n$ term is truncated at the specific index summation $n$. 
Substituting (\ref{eq:110015b}) into (\ref{eq:110014}) with $A=0$ gives the recurrence relations for the coefficients $c_n$. For $n=0,1,2,\cdots$, we can classify $c_n$ as to even and odd terms such as
\begin{equation}
\begin{tabular}{  l  l }
  \vspace{2 mm}
   $c_0$ &\hspace{1cm}  $c_1$ \\
  \vspace{2 mm}
   $c_2 = -\frac{\beta }{\delta } c_0 $  &\hspace{1cm}  $c_3 = -\frac{\beta }{\delta } c_1$  \\
  \vspace{2 mm}
  $c_4 = \left(-\frac{\beta }{\delta }\right)^2 c_0 $ &\hspace{1cm}  $c_5 = \left(-\frac{\beta }{\delta }\right)^2 c_1$\\
  \vspace{2 mm}
  $c_6 = \left(-\frac{\beta }{\delta }\right)^3 c_0 $ &\hspace{1cm}  $c_7 = \left(-\frac{\beta }{\delta }\right)^3 c_1$\\
 \hspace{2 mm} \large{\vdots} & \hspace{1.5 cm}\large{\vdots} \\
 \vspace{2 mm}
  $c_{2n} = \left(-\frac{\beta }{\delta }\right)^n c_0 $ &\hspace{1cm}  $c_{2n+1} = \left(-\frac{\beta }{\delta }\right)^n c_1$\\
\end{tabular}
\label{eq:110016}
\end{equation}
I suggest $c_1\sim  A c_0 =0$ in (\ref{eq:110016}). Because $A_n$ term is negligible for the minimum $y(x)$ since $A_n$ term is terminated at the specific eigenvalues. Put the coefficients $c_{2n}$ on the above into a power series $\sum_{n=0}^{\infty }c_{2n}x^{2n}$, letting $c_0=1$ for simplicity.
\begin{equation}
\mbox{min}\left( \lim_{n\gg 1}y(x) \right) = \frac{1}{1+\frac{\beta }{\delta }x^2} \label{eq:110017}
\end{equation}
In (\ref{eq:110017}), a polynomial of type 2 requires $\left| \beta /\delta x^2\right|<1$ for the convergence of the radius.  
 
In chapter 1 of Ref.\cite{11Choun2013}, the general expression of power series of $y(x)$ for polynomial of type 2 is defined by
\begin{eqnarray}
y(x) &=& \sum_{n=0}^{\infty } y_{n}(x) = y_0(x)+ y_1(x)+ y_2(x)+y_3(x)+\cdots \nonumber\\
&=&  c_0 \Bigg\{ \sum_{i_0=0}^{\alpha _0} \left( \prod _{i_1=0}^{i_0-1}A_{i_1} \right) x^{i_0+\lambda }
+ \sum_{i_0=0}^{\alpha _0}\left\{ B_{i_0+1} \prod _{i_1=0}^{i_0-1}A_{i_1}  \sum_{i_2=i_0}^{\alpha _1} \left( \prod _{i_3=i_0}^{i_2-1}A_{i_3+2} \right)\right\} x^{i_2+2+\lambda }  \nonumber\\
&& + \sum_{N=2}^{\infty } \Bigg\{ \sum_{i_0=0}^{\alpha _0} \Bigg\{B_{i_0+1}\prod _{i_1=0}^{i_0-1} A_{i_1} 
\prod _{k=1}^{N-1} \Bigg( \sum_{i_{2k}= i_{2(k-1)}}^{\alpha _k} B_{i_{2k}+2k+1}\prod _{i_{2k+1}=i_{2(k-1)}}^{i_{2k}-1}A_{i_{2k+1}+2k}\Bigg)\nonumber\\
&& \times  \sum_{i_{2N} = i_{2(N-1)}}^{\alpha _N} \Bigg( \prod _{i_{2N+1}=i_{2(N-1)}}^{i_{2N}-1} A_{i_{2N+1}+2N} \Bigg) \Bigg\} \Bigg\} x^{i_{2N}+2N+\lambda }\Bigg\}  \label{eq:110018}
\end{eqnarray}
In the above, $\alpha _i\leq \alpha _j$ only if $i\leq j$ where $i,j,\alpha _i, \alpha _j \in \mathbb{N}_{0}$.

For a polynomial, we need a condition which is:
\begin{equation}
A_{\alpha _i+ 2i}=0 \hspace{1cm} \mathrm{where}\;i,\alpha _i =0,1,2,\cdots
\label{eq:110019}
\end{equation}
In the above, $ \alpha _i$ is an eigenvalue that makes $A_n$ term terminated at certain value of index $n$. (\ref{eq:110019}) makes each sub-power series $y_i(x)$ where $i=0,1,2,\cdots$ as a polynomial in (\ref{eq:110018}).

Replace $\alpha _i$ by $q_i$ and put $n= q_i+ 2i$ in (\ref{eq:110013b}) with the condition $A_{q_i+ 2i}=0$.  
\begin{equation}
\pm \sqrt{\left( \gamma -1\right)^2 +4q}= -2\left( q_i +2i\right) -\gamma +1
\nonumber
\end{equation}
\subsubsection{The case of $ \sqrt{\left( \gamma -1\right)^2 +4q}= 2\left( q_i +2i\right) +\gamma -1$ where $i,q_i \in \mathbb{N}_{0}$}

In (\ref{eq:110013b}) replace $\sqrt{\left( \gamma -1\right)^2 +4q}$ by $2\left( q_i +2i\right) +\gamma -1$. In (\ref{eq:110018}) replace an index $\alpha _i$ by $q_i$.  Take the new (\ref{eq:110013b}) and (\ref{eq:110013c}) in the new (\ref{eq:110018}) with $c_0=1$ and $\lambda =0$.
After the replacement process, 
\begin{remark}
The power series expansion of the DCHE of the first kind for a polynomial of type 2 about $x=0$ as $q=\left( q_j+2j \right)\left( q_j+2j-1+\gamma \right) $ where $j,q_j \in \mathbb{N}_{0}$ is
\begin{eqnarray}
 y(x)&=& \sum_{n=0}^{\infty } y_{n}(x) = y_0(x)+ y_1(x)+ y_2(x)+y_3(x)+\cdots \nonumber\\ 
&=& H_d^{(o)}F_{q_{j}}^R\left( \alpha ,\beta ,\gamma ,\delta ,q= \left( q_j+2j \right)\left( q_j+2j-1+\gamma \right); \eta =-\frac{1}{\delta }x, \mu =-\frac{\beta }{\delta }x^2\right)\nonumber\\
&=& \sum_{i_0=0}^{q_0} \frac{\left( -q_0\right)_{i_0} \left( q_0-1+\gamma \right)_{i_0}}{(1)_{i_0}} \eta ^{i_0} \nonumber\\
&&+ \left\{ \sum_{i_0=0}^{q_0}\frac{ (i_0+ \alpha ) }{ (i_0+ 2) }\frac{\left( -q_0\right)_{i_0} \left( q_0-1+\gamma \right)_{i_0} }{(1)_{i_0} } \sum_{i_1=i_0}^{q_1} \frac{\left( -q_1\right)_{i_1}\left( q_1+3+\gamma \right)_{i_1}(3)_{i_0} }{\left( -q_1\right)_{i_0}\left( q_1+3+\gamma \right)_{i_0}(3)_{i_1}} \eta ^{i_1}\right\} \mu \nonumber\\
&&+ \sum_{n=2}^{\infty } \left\{ \sum_{i_0=0}^{q_0}\frac{ (i_0+ \alpha ) }{ (i_0+ 2) }\frac{\left( -q_0\right)_{i_0} \left( q_0-1+\gamma \right)_{i_0} }{(1)_{i_0} }\right.\nonumber\\
&&\times \prod _{k=1}^{n-1} \left\{ \sum_{i_k=i_{k-1}}^{q_k} \frac{ (i_k+ 2k +\alpha )}{ (i_k+2k+2) } \frac{\left( -q_k\right)_{i_k}\left( q_k+4k-1+\gamma \right)_{i_k}(2k+1)_{i_{k-1}} }{\left( -q_k\right)_{i_{k-1}}\left( q_k+4k-1+\gamma \right)_{i_{k-1}}(2k+1)_{i_k}}\right\} \nonumber\\
&&\times \left.\sum_{i_n= i_{n-1}}^{q_n} \frac{\left( -q_n\right)_{i_n}\left( q_n+4n-1+\gamma \right)_{i_n}(2n+1)_{i_{n-1}} }{\left( -q_n\right)_{i_{k-1}}\left( q_n+4n-1+\gamma \right)_{i_{n-1}}(2n+1)_{i_n}} \eta ^{i_n} \right\} \mu ^n  \label{eq:110020}
\end{eqnarray}
\end{remark} 
\subsubsection{The case of $ \sqrt{\left( \gamma -1\right)^2 +4q}= -2\left( q_i +2i\right) -\gamma +1$}
In (\ref{eq:110013b}) replace $\sqrt{\left( \gamma -1\right)^2 +4q}$ by $-2\left( q_i +2i\right) -\gamma +1$. In (\ref{eq:110018}) replace  index $\alpha _i$ by $q_i$.  Take the new (\ref{eq:110013b}) and (\ref{eq:110013c}) in (\ref{eq:110018}) with $c_0=1$ and $\lambda =0$.
After the replacement process, its solution is equivalent to (\ref{eq:110020}).

For the minimum value of the DCHE for a polynomial of type 2 about $x=0$, put $q_0=q_1=q_2=\cdots=0$ in (\ref{eq:110020}).
\begin{eqnarray}
 y(x) &=& H_d^{(o)}F_{0}^R\left( \alpha ,\beta ,\gamma ,\delta ,q= 2j\left( 2j-1+\gamma \right); \eta =-\frac{1}{\delta }x, \mu =-\frac{\beta }{\delta }x^2\right)\nonumber\\
&=&  \left( 1-\mu \right)^{-\frac{\alpha }{2}} \hspace{1cm}\mbox{where}\;\;|\mu|<1 \label{eq:110021}
\end{eqnarray} 
\subsection{Integral representation for a polynomial of type 2}

There is a generalized hypergeometric function which is given by
\begin{eqnarray}
L_l &=& \sum_{i_l= i_{l-1}}^{q_l} \frac{(-q_l)_{i_l}(q_l+4l-1+\gamma )_{i_l}(2l+1)_{i_{l-1}}}{(-q_l)_{i_{l-1}}(q_l+4l-1+\gamma )_{i_{l-1}}(2l+1)_{i_l}} \eta ^{i_l}\nonumber\\
&=& (i_{l-1}+2l) 
\sum_{j=0}^{\infty } \frac{B\left( i_{l-1}+2l,j+1\right) (i_{l-1}-q_l)_j\; (i_{l-1}+q_l+4l-1+\gamma )_j}{(1)_j} \eta ^{j+i_{l-1}} 
\hspace{1.5cm}\label{eq:110022}
\end{eqnarray}
By using integral form of beta function,
\begin{equation}
B\left( i_{l-1}+2l,j+1\right) = \int_{0}^{1} dt_l\;t_l^{i_{l-1}+2l-1} (1-t_l)^j
\label{eq:110023}
\end{equation}
Substitute (\ref{eq:110023}) into (\ref{eq:110022}), and divide $(i_{l-1}+2l)$ into $L_l$.
\begin{eqnarray}
G_l &=& \frac{1}{(i_{l-1}+2l)}
 \sum_{i_l= i_{l-1}}^{q_l} \frac{(-q_l)_{i_l}(q_l+4l-1+\gamma )_{i_l}(2l+1)_{i_{l-1}}}{(-q_l)_{i_{l-1}}(q_l+4l-1+\gamma )_{i_{l-1}}(2l+1)_{i_l}} \eta ^{i_l} \nonumber\\
&=&  \int_{0}^{1} dt_l\;t_l^{2l-1} \left(\eta t_l \right)^{i_{l-1}} \sum_{j=0}^{\infty } \frac{(i_{l-1}-q_l)_j\; (i_{l-1}+q_l+4l-1+\gamma )_j}{(1)_j} \left(\eta (1-t_l) \right)^j
\hspace{1.5cm} \label{eq:110024}
\end{eqnarray}
There is a Kummer's formula such as 
\begin{equation}
M(a,b,z)= \sum_{j= 0}^{\infty } \frac{(a)_j}{(b)_j}\frac{z^j}{j!} = \; _1F_1 (a;b;z)
\nonumber
\end{equation}
Tricomi's function is defined by
\begin{equation}
U(a,b,z)=  \frac{\Gamma (1-b)}{\Gamma (a-b+1)} M(a,b,z) +\frac{\Gamma (b-1)}{\Gamma (a)}z^{1-b} M(a-b+1,2-b,z)
\label{eq:110025}
\end{equation}
The contour integral form of (\ref{eq:110025}) is given by\cite{11NIST1}
\begin{equation}
U(a,b,z)=  e^{-a\pi i}\frac{\Gamma (1-a)}{2\pi i} \int_{\infty }^{(0+)} dv_l\; e^{-z v_l} v_l^{a-1} (1+ v_l)^{b-a-1}\;\mbox{where}\;a\ne 1,2,3, \cdots, \;\;\; |\mbox{ph} z|< \frac{1}{2}\pi
\label{eq:110026}
\end{equation}
Also (\ref{eq:110025}) is written by \cite{11NIST1}
\begin{equation}
U(a,b,z)=  z^{-a}  \sum_{j= 0}^{\infty } \frac{(a)_j (a-b+1)_j}{(1)_j} (-z^{-1})^j = z^{-a}\; _2F_0 (a,a-b+1;-;-z^{-1})
\label{eq:110027}
\end{equation}
Replace $a$, $b$ and $z$ by $i_{l-1} -q_l$, $-2q_l -4l+2-\gamma $ and $\frac{-1}{\eta (1-t_l)}$ into (\ref{eq:110027}).
\begin{align}
&\sum_{j= 0}^{\infty } \frac{(i_{l-1} - q_l )_j (i_{l-1}+q_l+4l-1+\gamma )_j}{(1)_j} (\eta (1-t_l))^j \nonumber\\
&= \left( \frac{-1}{\eta (1-t_l)}\right)^{i_{l-1}-q_l} U\left( i_{l-1} -q_l, -2q_l -4l+2-\gamma ,\frac{-1}{\eta (1-t_l)} \right)
\label{eq:110028}
\end{align}
Replace $a$, $b$ and $z$ by $i_{l-1} -q_l$, $-2q_l -4l+2-\gamma $ and $\frac{-1}{\eta (1-t_l)}$ into (\ref{eq:110026}).
Take the new (\ref{eq:110026}) into (\ref{eq:110028}).
\begin{eqnarray}
&& \sum_{j= 0}^{\infty } \frac{(i_{l-1} - q_l )_j (i_{l-1}+q_l+4l-1+\gamma )_j}{(1)_j} (\eta (1-t_l))^j \nonumber\\ 
&&= \frac{\Gamma (q_l -i_{l-1}+1)}{2\pi i} \int_{\infty }^{(0+)} dv_l\; \exp\left(\frac{ v_l}{\eta (1-v_l)}\right) v_l^{-1} (1+v_l)^{-4l+1-\gamma } \left( \frac{\eta (1-t_l)}{v_l (1+v_l)}\right)^{q_l}\nonumber\\
&&\times \left( \frac{v_l}{\eta (1-t_l)(1+v_l)}\right)^{i_{l-1}}
\hspace{1cm}\label{eq:110029}
\end{eqnarray}
There is the definition of the Gamma function $\Gamma (z)$ such as
\begin{equation}
\Gamma (z) = \int_{0}^{\infty} du_l\; e^{-u_l} u_l^{z-1} \;\mbox{where}\; Re(z)>0
\label{eq:110030}
\end{equation}
Put $z=q_l -i_{l-1}+1$ in (\ref{eq:110030}).
\begin{equation}
\Gamma (q_l -i_{l-1}+1) = \int_{0}^{\infty} du_l\; e^{-u_l} u_l^{q_l -i_{l-1}} \label{eq:110031}
\end{equation}
Put (\ref{eq:110031}) in (\ref{eq:110029}). Take the new (\ref{eq:110029}) into (\ref{eq:110024}).
\begin{eqnarray}
G_l &=& \frac{1}{(i_{l-1}+2l)} \sum_{i_l= i_{l-1}}^{q_l} \frac{(-q_l)_{i_l}(q_l+4l-1+\gamma )_{i_l}(2l+1)_{i_{l-1}}}{(-q_l)_{i_{l-1}}(q_l+4l-1+\gamma )_{i_{l-1}}(2l+1)_{i_l}} \eta ^{i_l} \nonumber\\
&=&  \int_{0}^{1} dt_l\;t_l^{2l-1} \int_{0}^{\infty} du_l\; \exp\left( -u_l\right)
\frac{1}{2\pi i} \int_{\infty }^{(0+)} dv_l\; \exp\left(\frac{ v_l}{\eta (1-t_l)}\right) \frac{1}{v_l(1+v_l)^{4l-1+\gamma }} \nonumber\\
&&\times \left( \frac{\eta u_l (1-t_l)}{v_l(1+v_l)}\right)^{q_l} \left( \frac{t_l v_l}{u_l (1-t_l)(1+v_l)}\right)^{i_{l-1}}
 \label{eq:110032}
\end{eqnarray}
Substitute (\ref{eq:110032}) into (\ref{eq:110020}) where $l=1,2,3,\cdots$: Apply $G_1$ into the second summation of sub-power series $y_1(x)$; apply $G_2$ into the third summation and $G_1$ into the second summation of sub-power series $y_2(x)$; apply $G_3$ into the forth summation, $G_2$ into the third summation and $G_1$ into the second summation of sub-power series $y_3(x)$, etc.\footnote{$y_1(x)$ means the sub-power series in (\ref{eq:110020}) contains one term of $B_n's$, $y_2(x)$ means the sub-power series in (\ref{eq:110020}) contains two terms of $B_n's$, $y_3(x)$ means the sub-power series in (\ref{eq:110020}) contains three terms of $B_n's$, etc.}
\begin{theorem}
The general representation in the form of integral of the DCH polynomial of type 2 about $x=0$ is given by 
\begin{eqnarray}
y(x)&=& \sum_{n=0}^{\infty } y_n(x) = y_0(x)+ y_1(x)+ y_2(x)+y_3(x)+\cdots \nonumber\\
&=& \sum_{i_0=0}^{q_0 }\frac{(-q_0)_{i_0} (q_0-1+\gamma )_{i_0}}{(1)_{i_0}}  \eta ^{i_0} 
+ \sum_{n=1}^{\infty } \left\{\prod _{k=0}^{n-1} \Bigg\{ \int_{0}^{1} dt_{n-k}\;t_{n-k}^{2(n-k)-1} \int_{0}^{\infty } du_{n-k}\;\exp\left( -u_{n-k}\right)  \right. \nonumber\\
&&\times \frac{1}{2\pi i} \int_{\infty }^{(0+)} dv_{n-k}\; \frac{1}{v_{n-k}\left( 1+v_{n-k}\right)^{4(n-k)-1+\gamma }} \exp \left( \frac{v_{n-k}}{w_{n-k+1,n}(1-t_{n-k})}\right) \nonumber\\
&&\times  \left( \frac{w_{n-k+1,n} u_{n-k} (1-t_{n-k})}{v_{n-k}\left( 1+v_{n-k}\right)}\right)^{q_{n-k}} w_{n-k,n}^{-2(n-k-1)-\alpha } \left(  w_{n-k,n} \partial _{w_{n-k,n}} \right) w_{n-k,n}^{2(n-k-1)+\alpha }\Bigg\}  \nonumber\\
&&\times \left. \sum_{i_0=0}^{q_0}\frac{(-q_0)_{i_0} (q_0-1+\gamma )_{i_0}}{(1)_{i_0}}  w_{1,n}^{i_0}\right\} \mu ^n  
\label{eq:110033}
\end{eqnarray}
where
\begin{equation} w_{i,j}=
\begin{cases} \displaystyle {\frac{t_i v_i}{u_i (1-t_i)(1+v_i)} }\;\;\mbox{where}\; i\leq j\cr
\eta \;\;\mbox{only}\;\mbox{if}\; i>j
\end{cases}
\nonumber
\end{equation}
In the above, the first sub-integral form contains one term of $B_n's$, the second one contains two terms of $B_n$'s, the third one contains three terms of $B_n$'s, etc.
\end{theorem} 
\begin{proof} 
In (\ref{eq:110020}) the power series expansions of sub-summation $y_0(x) $, $y_1(x)$, $y_2(x)$ and $y_3(x)$ of the DCH polynomial of type 2 about $x=0$ are given by
\begin{equation}
 y(x)= \sum_{n=0}^{\infty } y_n(x) = y_0(x)+ y_1(x)+ y_2(x)+y_3(x)+\cdots  \label{eq:110034}
\end{equation}
where
\begin{subequations}
\begin{equation}
 y_0(x)= \sum_{i_0=0}^{q_0} \frac{\left( -q_0\right)_{i_0} \left( q_0-1+\gamma \right)_{i_0}}{(1)_{i_0}} \eta ^{i_0} \label{eq:110035a}
\end{equation}
\begin{equation}
 y_1(x)=\left\{ \sum_{i_0=0}^{q_0}\frac{ (i_0+ \alpha ) }{ (i_0+ 2) }\frac{\left( -q_0\right)_{i_0} \left( q_0-1+\gamma \right)_{i_0} }{(1)_{i_0} } \right.   \left. \sum_{i_1=i_0}^{q_1} \frac{\left( -q_1\right)_{i_1}\left( q_1+3+\gamma \right)_{i_1}(3)_{i_0} }{\left( -q_1\right)_{i_0}\left( q_1+3+\gamma \right)_{i_0}(3)_{i_1}} \eta ^{i_1}\right\} \mu  \label{eq:110035b}
\end{equation}
\begin{eqnarray}
 y_2(x) &=& \left\{ \sum_{i_0=0}^{q_0}\frac{ (i_0+ \alpha ) }{ (i_0+ 2) }\frac{\left( -q_0\right)_{i_0} \left( q_0-1+\gamma \right)_{i_0} }{(1)_{i_0} } \right. \sum_{i_1=i_0}^{q_1} \frac{ (i_1+2+ \alpha ) }{ (i_1+ 4) } \frac{\left( -q_1\right)_{i_1}\left( q_1+3+\gamma \right)_{i_1}(3)_{i_0} }{\left( -q_1\right)_{i_0}\left( q_1+3+\gamma \right)_{i_0}(3)_{i_1}}  \nonumber\\
&&\times \left. \sum_{i_2= i_1}^{q_2} \frac{\left( -q_2\right)_{i_2}\left( q_2+7+\gamma \right)_{i_2}(5)_{i_1} }{\left( -q_2\right)_{i_1}\left( q_2+7+\gamma \right)_{i_1}(5)_{i_2}} \eta ^{i_2} \right\} \mu ^2 
\label{eq:110035c}
\end{eqnarray}
\begin{eqnarray}
 y_3(x)&=& \left\{ \sum_{i_0=0}^{q_0}\frac{ (i_0+ \alpha ) }{ (i_0+ 2) }\frac{\left( -q_0\right)_{i_0} \left( q_0-1+\gamma \right)_{i_0} }{(1)_{i_0} } \right. \sum_{i_1=i_0}^{q_1} \frac{ (i_1+2+ \alpha ) }{ (i_1+ 4) } \frac{\left( -q_1\right)_{i_1}\left( q_1+3+\gamma \right)_{i_1}(3)_{i_0} }{\left( -q_1\right)_{i_0}\left( q_1+3+\gamma \right)_{i_0}(3)_{i_1}}  \nonumber\\
&&\times \sum_{i_2=i_1}^{q_2}  \frac{ (i_2+4+ \alpha ) }{ (i_2+ 6) } \frac{\left( -q_2\right)_{i_2}\left( q_2+7+\gamma \right)_{i_2}(5)_{i_1} }{\left( -q_2\right)_{i_1}\left( q_2+7+\gamma \right)_{i_1}(5)_{i_2}}\nonumber\\
&&\times \left. \sum_{i_3=i_2}^{q_3} \frac{\left( -q_3\right)_{i_3}\left( q_3+11+\gamma \right)_{i_3}(7)_{i_2} }{\left( -q_3\right)_{i_2}\left( q_3+11+\gamma \right)_{i_2}(7)_{i_3}} \eta^{i_3} \right\} \mu ^3 
\label{eq:110035d}
\end{eqnarray}
\end{subequations}
Put $l=1$ in (\ref{eq:110032}). Take the new (\ref{eq:110032}) into (\ref{eq:110035b}).
\begin{eqnarray}
 y_1(x)&=& \int_{0}^{1} dt_1\;t_1 \int_{0}^{\infty } du_1\;e^{-u_1} 
\frac{1}{2\pi i}  \int_{\infty }^{(0+)} dv_1\; \exp\left(\frac{ v_1}{\eta (1-t_1)}\right) \frac{1}{v_1\left( 1+v_1\right)^{3+\gamma }} \nonumber\\
&&\times \left( \frac{\eta u_1 (1-t_1)}{v_1\left( 1+v_1\right)}\right)^{q_1} \left\{ \sum_{i_0=0}^{q_0} \left( i_0+\alpha \right) \frac{(-q_0)_{i_0}\left( q_0-1+\gamma \right)_{i_0}}{(1)_{i_0}} \left(\frac{t_1 v_1}{u_1 (1-t_1)(1+v_1)}\right)^{i_0} \right\} \mu \nonumber\\
&=&  \int_{0}^{1} dt_1\;t_1 \int_{0}^{\infty } du_1\;e^{-u_1} 
\frac{1}{2\pi i}  \int_{\infty }^{(0+)} dv_1\; \exp\left(\frac{ v_1}{\eta (1-t_1)}\right) \frac{1}{v_1\left( 1+v_1\right)^{3+\gamma }} \nonumber\\
&&\times \left( \frac{\eta u_1 (1-t_1)}{v_1\left( 1+v_1\right)}\right)^{q_1} w_{1,1}^{-\alpha } \left( w_{1,1}\partial_{w_{1,1}} \right) w_{1,1}^{\alpha } \left\{ \sum_{i_0=0}^{q_0} \frac{(-q_0)_{i_0}\left( q_0-1+\gamma \right)_{i_0}}{ (1)_{i_0}} w_{1,1}^{i_0} \right\} \mu \hspace{1.5cm}\label{eq:110036}\\
&& \mathrm{where}\hspace{.5cm} w_{1,1}= \frac{t_1 v_1}{u_1 (1-t_1)(1+v_1)}  \nonumber
\end{eqnarray}
Put $l=2$ in (\ref{eq:110032}). Take the new (\ref{eq:110032}) into (\ref{eq:110035c}). 
\begin{eqnarray}
 y_2(x) &=& \int_{0}^{1} dt_2\;t_2^3 \int_{0}^{\infty } du_2\;e^{-u_2} 
\frac{1}{2\pi i}  \int_{\infty }^{(0+)} dv_2\; \exp\left(\frac{v_2}{\eta (1-t_2)}\right) \frac{1}{v_2\left( 1+v_2\right)^{7+\gamma }}\nonumber\\
&&\times \left( \frac{\eta u_2 (1-t_2)}{v_2 \left( 1+v_2\right)}\right)^{q_2}  w_{2,2}^{-\left( 2+\alpha \right)} \left( w_{2,2}\partial_{w_{2,2}} \right) w_{2,2}^{2+\alpha } \nonumber\\
&&\times \left\{\sum_{i_0=0}^{q_0}\frac{ (i_0+ \alpha ) }{ (i_0+ 2) }\frac{\left( -q_0\right)_{i_0} \left( q_0-1+\gamma \right)_{i_0} }{(1)_{i_0} } \right.   \left. \sum_{i_1=i_0}^{q_1} \frac{\left( -q_1\right)_{i_1}\left( q_1+3+\gamma \right)_{i_1}(3)_{i_0} }{\left( -q_1\right)_{i_0}\left( q_1+3+\gamma \right)_{i_0}(3)_{i_1}} w_{2,2}^{i_1} \right\} \mu ^2 \hspace{1.5cm}\label{eq:110037}\\
&& \mathrm{where}\hspace{.5cm} w_{2,2}=  \frac{t_2 v_2}{u_2 (1-t_2)(1+v_2)} \nonumber
\end{eqnarray}
Put $l=1$ and $\eta = w_{2,2}$ in (\ref{eq:110032}). Take the new (\ref{eq:110032}) into (\ref{eq:110037}).
\begin{eqnarray}
 y_2(z)&=& \int_{0}^{1} dt_2\;t_2^3 \int_{0}^{\infty } du_2\;e^{-u_2} 
\frac{1}{2\pi i}  \int_{\infty }^{(0+)} dv_2\; \exp\left(\frac{v_2}{\eta (1-t_2)}\right) \frac{1}{v_2\left( 1+v_2\right)^{7+\gamma }}\nonumber\\
&&\times \left( \frac{\eta u_2 (1-t_2)}{v_2 \left( 1+v_2\right)}\right)^{q_2}  w_{2,2}^{-\left( 2+\alpha \right)} \left( w_{2,2}\partial_{w_{2,2}} \right) w_{2,2}^{2+\alpha } \nonumber\\
&&\times \int_{0}^{1} dt_1\;t_1 \int_{0}^{\infty } du_1\;e^{-u_1} 
\frac{1}{2\pi i}  \int_{\infty }^{(0+)} dv_1\; \exp\left(\frac{v_1}{w_{2,2} (1-t_1)}\right) \frac{1}{v_1\left( 1+v_1\right)^{3+\gamma }}\nonumber\\
&&\times \left( \frac{w_{2,2} u_1 (1-t_1)}{v_1 \left( 1+v_1\right)}\right)^{q_1}  w_{1,2}^{-\alpha } \left( w_{1,2}\partial_{w_{1,2}} \right) w_{1,2}^{\alpha } \left\{ \sum_{i_0=0}^{q_0} \frac{\left( -q_0\right)_{i_0} \left( q_0-1+\gamma \right)_{i_0} }{(1)_{i_0} } w_{1,2}^{i_0} \right\} \mu ^2  \hspace{1.5cm} \label{eq:110038}\\
&& \mathrm{where}\hspace{.5cm} w_{1,2}= \frac{t_1 v_1}{u_1 (1-t_1)(1+v_1)} \nonumber
\end{eqnarray}
By using similar process for the previous cases of integral forms of $y_1(x)$ and $y_2(x)$, the integral form of a sub-power series $y_3(x)$ is obtained by
\begin{eqnarray}
 y_3(x)&=&  \int_{0}^{1} dt_3\;t_3^5 \int_{0}^{\infty } du_3\;e^{-u_3} 
\frac{1}{2\pi i}  \int_{\infty }^{(0+)} dv_3\; \exp\left(\frac{ v_3}{\eta (1-t_3)}\right) \frac{1}{v_3\left( 1+v_3\right)^{11+\gamma }} \nonumber\\
&&\times \left( \frac{\eta u_3 (1-t_3)}{v_3 \left( 1+v_3\right)}\right)^{q_3}  w_{3,3}^{-\left( 4+\alpha \right)} \left( w_{3,3}\partial_{w_{3,3}} \right) w_{3,3}^{4+\alpha } \nonumber\\
&&\times \int_{0}^{1} dt_2\;t_2^3 \int_{0}^{\infty } du_2\;e^{-u_2} 
\frac{1}{2\pi i}  \int_{\infty }^{(0+)} dv_2\; \exp\left(\frac{v_2}{w_{3,3}(1-t_2)}\right) \frac{1}{v_2\left( 1+v_2\right)^{7+\gamma }}\nonumber\\
&&\times \left( \frac{w_{3,3} u_2 (1-t_2)}{v_2 \left( 1+v_2\right)}\right)^{q_2}  w_{2,3}^{-\left( 2+\alpha \right)} \left( w_{2,3}\partial_{w_{2,3}} \right) w_{2,3}^{2+\alpha } \nonumber\\
&&\times \int_{0}^{1} dt_1\;t_1 \int_{0}^{\infty } du_1\;e^{-u_1} 
\frac{1}{2\pi i}  \int_{\infty }^{(0+)} dv_1\; \exp\left(\frac{v_1}{w_{2,3} (1-t_1)}\right) \frac{1}{v_1\left( 1+v_1\right)^{3+\gamma }}\nonumber\\
&&\times \left( \frac{w_{2,3} u_1 (1-t_1)}{v_1 \left( 1+v_1\right)}\right)^{q_1}  w_{1,3}^{-\alpha } \left( w_{1,3}\partial_{w_{1,3}} \right) w_{1,3}^{\alpha } \left\{ \sum_{i_0=0}^{q_0} \frac{\left( -q_0\right)_{i_0} \left( q_0-1+\gamma \right)_{i_0} }{(1)_{i_0} } w_{1,3}^{i_0} \right\} \mu ^3 \hspace{1.5cm} \label{eq:110039} 
\end{eqnarray}
where
\begin{equation}
\begin{cases} \displaystyle{w_{3,3} =\frac{t_3 v_3}{u_3 (1-t_3)(1+v_3)}} \cr
\displaystyle{w_{2,3} = \frac{t_2 v_2}{u_2 (1-t_2)(1+v_2)}} \cr
\displaystyle{w_{1,3}= \frac{t_1 v_1}{u_1 (1-t_1)(1+v_1)}}
\end{cases}
\nonumber
\end{equation}
By repeating this process for all higher terms of integral forms of sub-summation $y_m(x)$ terms where $m \geq 4$, we obtain every integral forms of $y_m(x)$ terms. 
Since we substitute (\ref{eq:110035a}), (\ref{eq:110036}), (\ref{eq:110038}), (\ref{eq:110039}) and including all integral forms of $y_m(x)$ terms where $m \geq 4$ into (\ref{eq:110034}), we obtain (\ref{eq:110033}).
\qed
\end{proof}   
\begin{remark}
The integral representation of the DCHE of the first kind for a polynomial of type 2 about $x=0$ as $q = \left( q_j +2j \right)\left( q_j +2j-1+\gamma \right)$ where $j, q_j \in \mathbb{N}_{0}$ is
\begin{eqnarray}
y(x) &=& H_d^{(o)}F_{q_{j}}^R\left( \alpha ,\beta ,\gamma ,\delta ,q= \left( q_j+2j \right)\left( q_j+2j-1+\gamma \right); \eta =-\frac{1}{\delta }x, \mu =-\frac{\beta }{\delta }x^2\right)\nonumber\\
&=& \left(-\eta \right)^{q_0} U\left( -q_0, -2q_0+2-\gamma ,-\eta ^{-1}\right) \nonumber\\ 
&&+ \sum_{n=1}^{\infty } \Bigg\{\prod _{k=0}^{n-1} \Bigg\{ \int_{0}^{1} dt_{n-k}\;t_{n-k}^{2(n-k)-1} \int_{0}^{\infty } du_{n-k}\;\exp\left( -u_{n-k}\right)  \nonumber\\
&&\times \frac{1}{2\pi i} \int_{\infty }^{(0+)} dv_{n-k}\; \frac{1}{v_{n-k}\left( 1+v_{n-k}\right)^{4(n-k)-1+\gamma }} \exp \left( \frac{v_{n-k}}{w_{n-k+1,n}(1-t_{n-k})}\right) \nonumber\\
&&\times \left( \frac{w_{n-k+1,n} u_{n-k} (1-t_{n-k})}{v_{n-k}\left( 1+v_{n-k}\right)}\right)^{q_{n-k}} w_{n-k,n}^{-2(n-k-1)-\alpha } \left(  w_{n-k,n} \partial _{w_{n-k,n}} \right) w_{n-k,n}^{2(n-k-1)+\alpha }\Bigg\} \nonumber\\
&&\times \left( -w_{1,n}\right)^{q_0} U\left( -q_0, -2q_0+2-\gamma , -w_{1,n}^{-1}\right)\Bigg\} \mu ^n  
\label{eq:110040}
\end{eqnarray}
\end{remark}
\begin{proof}
Replace $a$, $b$ and $z$ by $ -q_0$, $-2q_0+2-\gamma $ and $ -\eta^{-1} $ into (\ref{eq:110027}).
\begin{equation}
\sum_{i_0=0}^{q_0} \frac{\left( -q_0\right)_{i_0} \left( q_0-1+\gamma \right)_{i_0}}{(1)_{i_0}} \eta ^{i_0}  = \left(-\eta \right)^{q_0} U\left( -q_0, -2q_0+2-\gamma ,-\eta ^{-1}\right)  \label{eq:110041}
\end{equation} 
Replace $a$, $b$ and $z$ by $ -q_0$, $-2q_0+2-\gamma $ and $ -w_{1,n}^{-1} $ into (\ref{eq:110027}).
\begin{equation}
\sum_{i_0=0}^{q_0} \frac{\left( -q_0\right)_{i_0} \left( q_0-1+\gamma \right)_{i_0}}{(1)_{i_0}} w_{1,n}^{i_0}  = \left( -w_{1,n} \right)^{q_0} U\left( -q_0, -2q_0+2-\gamma ,-w_{1,n}^{-1}\right)  \label{eq:110042}
\end{equation} 
We obtain (\ref{eq:110040}) by substituting (\ref{eq:110041}) and (\ref{eq:110042}) into (\ref{eq:110033}).
\qed
\end{proof} 
\subsection{Generating function for a polynomial of type 2}
Let's investigate the generating function for the DCH polynomial of type 2 of the first kind about $x=0$. 
\begin{definition}
I define that
\begin{equation}
\begin{cases}
\displaystyle { s_{a,b}} = \begin{cases} \displaystyle {  s_a\cdot s_{a+1}\cdot s_{a+2}\cdots s_{b-2}\cdot s_{b-1}\cdot s_b}\;\;\mbox{where}\;a<b \cr
s_a \;\;\mbox{where}\;a=b\end{cases}
\cr
\cr
\displaystyle {  \widetilde{w}_{i,j}= \begin{cases} \displaystyle { \frac{t_i}{u_i(1-t_i)}\frac{-1+\sqrt{1+4\widetilde{w}_{i+1,j} u_i(1-t_i)s_i}}{1+\sqrt{1+4\widetilde{w}_{i+1,j} u_i(1-t_i)s_i}}}\;\;\mbox{where}\;i<j \cr
\displaystyle { \frac{t_i}{u_i(1-t_i)}\frac{-1+\sqrt{1+4\eta u_i(1-t_i)s_{i,\infty }}}{1+\sqrt{1+4\eta u_i(1-t_i)s_{i,\infty }}}}  \;\;\mbox{where}\;i=j\end{cases}}
\end{cases}\label{eq:110041}
\end{equation}
where $a,b,i,j \in \mathbb{N}_{0}$, $0\leq a\leq b\leq \infty $ and $1\leq i\leq j\leq \infty $.
\end{definition}
And we have
\begin{equation}
\sum_{q_i = q_j}^{\infty } s_i^{q_i} = \frac{s_i^{q_j}}{(1-s_i)}\label{eq:110042}
\end{equation}
Acting the summation operator $ \sum_{q_0 =0}^{\infty } \frac{s_0^{q_0}}{ q_0 !}  \prod _{n=1}^{\infty } \left\{ \sum_{ q_n = q_{n-1}}^{\infty } s_n^{q_n }\right\} $ on (\ref{eq:110033}) where $|s_i|<1$ as $i=0,1,2,\cdots$ by using (\ref{eq:110041}) and (\ref{eq:110042}),
\begin{theorem}
The general expression of the generating function for the DCH polynomial of type 2 about $x=0$ is given by
\begin{eqnarray}
&& \sum_{q_0 =0}^{\infty } \frac{s_0^{q_0}}{ q_0 !}  \prod _{n=1}^{\infty } \left\{ \sum_{ q_n = q_{n-1}}^{\infty } s_n^{q_n }\right\} y(x) \nonumber\\
&&= \prod_{k=1}^{\infty } \frac{1}{(1-s_{k,\infty })} \mathbf{\Upsilon}(s_{0,\infty } ;\eta )  \nonumber\\
&&+ \left\{ \prod_{k=2}^{\infty } \frac{1}{(1-s_{k,\infty })} \int_{0}^{1} dt_1\;t_1 \int_{0}^{\infty } du_1\; \exp\left( -u_1 +\frac{-1+\sqrt{1+4\eta u_1(1-t_1)s_{1,\infty }}}{2\eta (1-t_1)}\right) \right. \nonumber\\
&&\times \left. \left( \frac{2}{1+\sqrt{1+4\eta u_1(1-t_1)s_{1,\infty }}}\right)^{2+\gamma } \frac{\widetilde{w}_{1,1}^{-\alpha } \left(  \widetilde{w}_{1,1} \partial _{\widetilde{w}_{1,1}} \right) \widetilde{w}_{1,1}^{\alpha }}{\sqrt{1+4\eta u_1(1-t_1)s_{1,\infty }}} \mathbf{\Upsilon}(s_0; \widetilde{w}_{1,1}) \right\} \mu  \nonumber\\
&&+ \sum_{n=2}^{\infty } \Bigg\{ \frac{(1-s_{n,\infty })}{\prod_{k=n}^{\infty }(1-s_{k,\infty })}  \int_{0}^{1} dt_n\;t_n^{2n-1} \int_{0}^{\infty } du_n\; \exp\left( -u_n +\frac{-1+\sqrt{1+4\eta u_n(1-t_n)s_{n,\infty }}}{2\eta (1-t_n)}\right) \nonumber\\
&&\times  \left( \frac{2}{1+\sqrt{1+4\eta u_n(1-t_n)s_{n,\infty }}}\right)^{4n-2+\gamma } \frac{\widetilde{w}_{n,n}^{-2(n-1)-\alpha } \left(  \widetilde{w}_{n,n} \partial _{\widetilde{w}_{n,n}} \right) \widetilde{w}_{n,n}^{2(n-1)+\alpha }}{\sqrt{1+4\eta u_n(1-t_n)s_{n,\infty }}}  \nonumber\\
&&\times  \prod_{j=1}^{n-1} \Bigg\{ \int_{0}^{1} dt_{n-j}\;t_{n-j}^{2(n-j)-1} \int_{0}^{\infty } du_{n-j}\; \exp\left( -u_{n-j} +\frac{-1+\sqrt{1+4\widetilde{w}_{n-j+1,n} u_{n-j}(1-t_{n-j})s_{n-j}}}{2\widetilde{w}_{n-j+1,n} (1-t_{n-j})} \right) \nonumber\\
&&\times \left( \frac{2}{1+\sqrt{1+4\widetilde{w}_{n-j+1,n} u_{n-j}(1-t_{n-j})s_{n-j}}}\right)^{4(n-j)-2+\gamma }\nonumber\\ 
&&\times \frac{\widetilde{w}_{n-j,n}^{-2(n-j-1)-\alpha } \left(  \widetilde{w}_{n-j,n} \partial _{\widetilde{w}_{n-j,n}} \right) \widetilde{w}_{n-j,n}^{2(n-j-1)+\alpha }}{\sqrt{1+4\widetilde{w}_{n-j+1,n} u_{n-j}(1-t_{n-j})s_{n-j}}} \Bigg\}
 \mathbf{\Upsilon}( s_0 ;\widetilde{w}_{1,n}) \Bigg\} \mu ^n \label{eq:110043}
\end{eqnarray}
where
\begin{equation}
\begin{cases} 
{ \displaystyle \mathbf{\Upsilon}( s_{0,\infty } ;\eta )= \sum_{q_0=0}^{\infty } \frac{s_{0,\infty }^{q_0}}{q_0!}  \left\{ \sum_{i_0=0}^{q_0} \frac{(-q_0)_{i_0} (q_0-1+\gamma )_{i_0}}{(1)_{i_0}} \eta^{i_0} \right\} }
\cr
{ \displaystyle \mathbf{\Upsilon}( s_0;\widetilde{w}_{1,1}) =  \sum_{q_0=0}^{\infty } \frac{s_0^{q_0}}{q_0!}  \left\{ \sum_{i_0=0}^{q_0} \frac{(-q_0)_{i_0} (q_0-1+\gamma )_{i_0}}{(1)_{i_0}} \widetilde{w}_{1,1}^{i_0} \right\}} \cr
{ \displaystyle \mathbf{\Upsilon}( s_0 ;\widetilde{w}_{1,n}) = \sum_{q_0=0}^{\infty } \frac{s_0^{q_0}}{q_0!}  \left\{ \sum_{i_0=0}^{q_0} \frac{(-q_0)_{i_0} (q_0-1+\gamma )_{i_0}}{(1)_{i_0}} \widetilde{w}_{1,n}^{i_0} \right\}}
\end{cases}\nonumber  
\end{equation}
\end{theorem}
\begin{proof} 
Acting the summation operator  $ \sum_{q_0 =0}^{\infty } \frac{s_0^{q_0}}{ q_0 !}  \prod _{n=1}^{\infty } \left\{ \sum_{ q_n = q_{n-1}}^{\infty } s_n^{q_n }\right\} $ on the form of integral of the type 2 DCH polynomial $y(x)$,
\begin{eqnarray}
&& \sum_{q_0 =0}^{\infty } \frac{s_0^{q_0}}{ q_0 !}  \prod _{n=1}^{\infty } \left\{ \sum_{ q_n = q_{n-1}}^{\infty } s_n^{q_n }\right\} y(x) \nonumber\\
&&=  \sum_{q_0 =0}^{\infty } \frac{s_0^{q_0}}{ q_0 !}  \prod _{n=1}^{\infty } \left\{ \sum_{ q_n = q_{n-1}}^{\infty } s_n^{q_n }\right\} \Big( y_0(x)+y_1(x)+y_2(x)+y_3(x)+ \cdots \Big) \label{eq:110044}
\end{eqnarray}
Acting the summation operator $ \sum_{q_0 =0}^{\infty } \frac{s_0^{q_0}}{ q_0 !}  \prod _{n=1}^{\infty } \left\{ \sum_{ q_n = q_{n-1}}^{\infty } s_n^{q_n }\right\} $ on (\ref{eq:110035a}) by using (\ref{eq:110041}) and (\ref{eq:110042}),
\begin{eqnarray}
&&\sum_{q_0 =0}^{\infty } \frac{s_0^{q_0}}{ q_0 !}  \prod _{n=1}^{\infty } \left\{ \sum_{ q_n = q_{n-1}}^{\infty } s_n^{q_n }\right\} y_0(x) \nonumber\\
&&= \prod_{k=1}^{\infty } \frac{1}{(1-s_{k,\infty })} \sum_{q_0 =0}^{\infty } \frac{s_{0,\infty }^{q_0}}{q_0!} \left\{ \sum_{i_0=0}^{q_0} \frac{(-q_0)_{i_0} (q_0-1+\gamma )_{i_0}}{ (1)_{i_0}} \eta ^{i_0} \right\} \label{eq:110045}
\end{eqnarray}
Acting the summation operator $ \sum_{q_0 =0}^{\infty } \frac{s_0^{q_0}}{ q_0 !}  \prod _{n=1}^{\infty } \left\{ \sum_{ q_n = q_{n-1}}^{\infty } s_n^{q_n }\right\} $ on (\ref{eq:110036}) by using (\ref{eq:110041}) and (\ref{eq:110042}),
\begin{eqnarray}
&&\sum_{q_0 =0}^{\infty } \frac{s_0^{q_0}}{ q_0 !}  \prod _{n=1}^{\infty } \left\{ \sum_{ q_n = q_{n-1}}^{\infty } s_n^{q_n }\right\} y_1(x) \nonumber\\
&&=  \prod_{k=2}^{\infty } \frac{1}{(1-s_{k,\infty })}  \int_{0}^{1} dt_1\;t_1 \int_{0}^{\infty } du_1\;e^{-u_1} 
\frac{1}{2\pi i}  \int_{\infty }^{(0+)} dv_1\; \exp\left(\frac{ v_1}{\eta (1-t_1)}\right) \frac{1}{v_1(1+v_1)^{3+\gamma }} \nonumber\\
&&\times \sum_{ q_1 = q_0}^{\infty } \left( \frac{\eta u_1 (1-t_1)s_{1,\infty }}{v_1(1+v_1)}\right)^{q_1} w_{1,1}^{-\alpha } \left( w_{1,1}\partial_{w_{1,1}} \right) w_{1,1}^{\alpha } \nonumber\\
&&\times \sum_{q_0 =0}^{\infty } \frac{s_0^{q_0}}{ q_0 !} \left\{ \sum_{i_0=0}^{q_0} \frac{(-q_0)_{i_0}(q_0-1+\gamma )_{i_0}}{(1)_{i_0}} w_{1,1}^{i_0} \right\} \mu \hspace{1.5cm} \label{eq:110046}
\end{eqnarray}
Replace $q_i$, $q_j$ and $s_i$ by $q_1$, $q_0$ and $ \displaystyle{\frac{\eta u_1 (1-t_1)s_{1,\infty }}{v_1(1+v_1)}} $ in (\ref{eq:110042}). Take the new (\ref{eq:110042}) into (\ref{eq:110046}).
\begin{eqnarray}
&&\sum_{q_0 =0}^{\infty } \frac{s_0^{q_0}}{ q_0 !}  \prod _{n=1}^{\infty } \left\{ \sum_{ q_n = q_{n-1}}^{\infty } s_n^{q_n }\right\} y_1(x) \nonumber\\
&&= \prod_{k=2}^{\infty } \frac{1}{(1-s_{k,\infty })}  \int_{0}^{1} dt_1\;t_1 \int_{0}^{\infty } du_1\;e^{-u_1} 
 \nonumber\\
&&\times \frac{1}{2\pi i}  \int_{\infty }^{(0+)} dv_1\; \exp\left(\frac{ v_1}{\eta (1-t_1)}\right) \frac{1}{(1+v_1)^{2+\gamma }} \frac{ w_{1,1}^{-\alpha } \left( w_{1,1}\partial_{w_{1,1}} \right) w_{1,1}^{\alpha }}{v_1^2+v_1-\eta u_1(1-t_1)s_{1,\infty }} \nonumber\\
&&\times \sum_{q_0 =0}^{\infty } \left( \frac{\eta u_1 (1-t_1)s_{0,\infty }}{v_1(1+v_1)}\right)^{q_0}\frac{1}{ q_0 !} \left\{ \sum_{i_0=0}^{q_0} \frac{(-q_0)_{i_0}(q_0-1+\gamma )_{i_0}}{(1)_{i_0}} w_{1,1}^{i_0} \right\} \mu  \label{eq:110047}
\end{eqnarray}
By using Cauchy's integral formula, the contour integrand has poles at $v_1 = \frac{-1-\sqrt{1+4\eta u_1(1-t_1)s_{1,\infty }}}{2}$ 
or $\frac{-1+\sqrt{1+4\eta u_1(1-t_1)s_{1,\infty }}}{2}$ and $\frac{-1+\sqrt{1+4\eta u_1(1-t_1)s_{1,\infty }}}{2}$ is only inside the unit circle. As we compute the residue there in (\ref{eq:110047}) we obtain
\begin{eqnarray}
&&\sum_{q_0 =0}^{\infty } \frac{s_0^{q_0}}{ q_0 !}  \prod _{n=1}^{\infty } \left\{ \sum_{ q_n = q_{n-1}}^{\infty } s_n^{q_n }\right\} y_1(x) \nonumber\\
&&= \left\{ \prod_{k=2}^{\infty } \frac{1}{(1-s_{k,\infty })} \int_{0}^{1} dt_1\;t_1 \int_{0}^{\infty } du_1\; \exp\left( -u_1 +\frac{-1+\sqrt{1+4\eta u_1(1-t_1)s_{1,\infty }}}{2\eta (1-t_1)}\right) \right. \nonumber\\
&&\times \left( \frac{2}{1+\sqrt{1+4\eta u_1(1-t_1)s_{1,\infty }}}\right)^{2+\gamma } \frac{\widetilde{w}_{1,1}^{-\alpha } \left(  \widetilde{w}_{1,1} \partial _{\widetilde{w}_{1,1}} \right) \widetilde{w}_{1,1}^{\alpha }}{\sqrt{1+4\eta u_1(1-t_1)s_{1,\infty }}} \nonumber\\
&&\times \left. \sum_{q_0=0}^{\infty } \frac{s_0^{q_0}}{q_0!}  \left\{ \sum_{i_0=0}^{q_0} \frac{(-q_0)_{i_0} (q_0-1+\gamma )_{i_0}}{(1)_{i_0}} \widetilde{w}_{1,1}^{i_0} \right\} \right\} \mu  \label{eq:110048}
\end{eqnarray}
where
\begin{equation}
\widetilde{w}_{1,1} =  \frac{t_1 v_1}{u_1 (1-t_1)(1+v_1)}\Bigg|_{ v_1 = \frac{1}{2}\left( -1+\sqrt{1+4\eta u_1(1-t_1)s_{1,\infty }}\right)} = \frac{t_1}{u_1(1-t_1)}\frac{-1+\sqrt{1+4\eta u_1(1-t_1)s_{1,\infty }}}{1+\sqrt{1+4\eta u_1(1-t_1)s_{1,\infty }}} \nonumber
\end{equation}
Acting the summation operator $ \sum_{q_0 =0}^{\infty } \frac{s_0^{q_0}}{ q_0 !}  \prod _{n=1}^{\infty } \left\{ \sum_{ q_n = q_{n-1}}^{\infty } s_n^{q_n }\right\} $ on (\ref{eq:110038}) by using (\ref{eq:110041}) and (\ref{eq:110042}),
\begin{eqnarray}
&&\sum_{q_0 =0}^{\infty } \frac{s_0^{q_0}}{ q_0 !}  \prod _{n=1}^{\infty } \left\{ \sum_{ q_n = q_{n-1}}^{\infty } s_n^{q_n }\right\} y_2(x) \nonumber\\
&&=  \prod_{k=3}^{\infty } \frac{1}{(1-s_{k,\infty })}  \int_{0}^{1} dt_2\;t_2^{3} \int_{0}^{\infty } du_2\;e^{-u_2} 
\frac{1}{2\pi i}  \int_{\infty }^{(0+)} dv_2\; \exp\left(\frac{ v_2}{\eta (1-t_2)}\right) \frac{1}{v_2(1+v_2)^{7+\gamma }} \nonumber\\
&&\times \sum_{ q_2=q_1}^{\infty } \left( \frac{\eta u_2 (1-t_2)s_{2,\infty }}{v_2(1+v_2)}\right)^{q_2} w_{2,2}^{-2-\alpha } \left( w_{2,2}\partial_{w_{2,2}} \right) w_{2,2}^{2+\alpha } \nonumber\\
&&\times  \int_{0}^{1} dt_1\;t_1 \int_{0}^{\infty } du_1\;e^{-u_1} 
\frac{1}{2\pi i}  \int_{\infty }^{(0+)} dv_1\; \exp\left(\frac{ v_1}{w_{2,2} (1-t_1)}\right) \frac{1}{v_1(1+v_1)^{3+\gamma }} \nonumber\\
&&\times \sum_{q_1 =q_0}^{\infty } \left( \frac{ w_{2,2} u_1 (1-t_1)s_1}{v_1(1+v_1)}\right)^{q_1} w_{1,2}^{-\alpha } \left( w_{1,2}\partial_{w_{1,2}} \right) w_{1,2}^{\alpha } \nonumber\\
&&\times \sum_{q_0 =0}^{\infty } \frac{s_0^{q_0}}{q_0!} \left\{ \sum_{i_0=0}^{q_0} \frac{(-q_0)_{i_0} (q_0-1+\gamma )_{i_0}}{ (1)_{i_0}} w_{1,2}^{i_0} \right\} \mu ^2  \label{eq:110049}
\end{eqnarray}
Replace $q_i$, $q_j$ and $s_i$ by $q_2$, $q_1$ and ${ \displaystyle \frac{\eta u_2 (1-t_2)s_{2,\infty }}{v_2(1+v_2)}}$ in (\ref{eq:110042}). Take the new (\ref{eq:110042}) into (\ref{eq:110049}).
\begin{eqnarray}
&&\sum_{q_0 =0}^{\infty } \frac{s_0^{q_0}}{ q_0 !}  \prod _{n=1}^{\infty } \left\{ \sum_{ q_n = q_{n-1}}^{\infty } s_n^{q_n }\right\} y_2(x) \nonumber\\
&&=  \prod_{k=3}^{\infty } \frac{1}{(1-s_{k,\infty })}  \int_{0}^{1} dt_2\;t_2^3 \int_{0}^{\infty } du_2\;e^{-u_2} 
 \nonumber\\
&&\times \frac{1}{2\pi i}  \int_{\infty }^{(0+)} dv_2\; \exp\left(\frac{ v_2}{\eta (1-t_2)}\right) \frac{1}{(1+v_2)^{ 6+\gamma }}\frac{w_{2,2}^{-2-\alpha } \left( w_{2,2}\partial_{w_{2,2}} \right) w_{2,2}^{2+\alpha }}{v_2^2+v_2-\eta u_2(1-t_2)s_{2,\infty }}  \nonumber\\
&&\times \int_{0}^{1} dt_1\;t_1 \int_{0}^{\infty } du_1\;e^{-u_1} 
\frac{1}{2\pi i}  \int_{\infty }^{(0+)} dv_1\; \exp\left(\frac{ v_1}{w_{2,2} (1-t_1)}\right) \frac{1}{v_1(1+v_1)^{3+\gamma }} \nonumber\\ 
&&\times \sum_{ q_1=q_0}^{\infty } \left( \frac{\eta u_2(1-t_2)s_{1,\infty }}{v_2(1+v_2)}\frac{ w_{2,2} u_1 (1-t_1)}{v_1(1+v_1)}\right)^{q_1} w_{1,2}^{-\alpha } \left( w_{1,2}\partial_{w_{1,2}} \right) w_{1,2}^{\alpha } \nonumber\\
&&\times  \sum_{q_0 =0}^{\infty } \frac{s_0^{q_0}}{q_0!} \left\{ \sum_{i_0=0}^{q_0} \frac{(-q_0)_{i_0} (q_0-1+\gamma )_{i_0}}{(1)_{i_0}} w_{1,2}^{i_0} \right\} \mu ^2 \label{eq:110050}
\end{eqnarray}
By using Cauchy's integral formula, the contour integrand has poles at $v_2 = \frac{-1-\sqrt{1+4\eta u_2(1-t_2)s_{2,\infty }}}{2}$ 
or $\frac{-1+\sqrt{1+4\eta u_2(1-t_2)s_{2,\infty }}}{2}$ and $\frac{-1+\sqrt{1+4\eta u_2(1-t_2)s_{2,\infty }}}{2}$ is only inside the unit circle. As we compute the residue there in (\ref{eq:110050}) we obtain
\begin{eqnarray}
&&\sum_{q_0 =0}^{\infty } \frac{s_0^{q_0}}{ q_0 !}  \prod _{n=1}^{\infty } \left\{ \sum_{ q_n = q_{n-1}}^{\infty } s_n^{q_n }\right\} y_2(x) \nonumber\\
&&=  \prod_{k=3}^{\infty } \frac{1}{(1-s_{k,\infty })}  \int_{0}^{1} dt_2\;t_2^3 \int_{0}^{\infty } du_2\; \exp\left( -u_2 +\frac{-1+\sqrt{1+4\eta u_2(1-t_2)s_{2,\infty }}}{2\eta (1-t_2)}\right) \nonumber\\
&&\times \left( \frac{2}{1+\sqrt{1+4\eta u_2(1-t_2)s_{2,\infty }}}\right)^{6+\gamma } \frac{\widetilde{w}_{2,2}^{-2-\alpha } \left(  \widetilde{w}_{2,2} \partial _{\widetilde{w}_{2,2}} \right) \widetilde{w}_{2,2}^{2+\alpha }}{\sqrt{1+4\eta u_2(1-t_2)s_{2,\infty }}} \nonumber\\
&&\times \int_{0}^{1} dt_1\;t_1 \int_{0}^{\infty } du_1\;e^{-u_1} 
\frac{1}{2\pi i}  \int_{\infty }^{(0+)} dv_1\; \exp\left(\frac{ v_1}{\widetilde{w}_{2,2} (1-t_1)}\right) \frac{1}{v_1(1+v_1)^{3+\gamma }} \nonumber\\ 
&&\times \sum_{ q_1=q_0}^{\infty } \left( \frac{ \widetilde{w}_{2,2} u_1 (1-t_1)s_1}{v_1(1+v_1)}\right)^{q_1} w_{1,2}^{-\alpha } \left( w_{1,2}\partial_{w_{1,2}} \right) w_{1,2}^{\alpha } \nonumber\\
&&\times  \sum_{q_0 =0}^{\infty } \frac{s_0^{q_0}}{q_0!} \left\{ \sum_{i_0=0}^{q_0} \frac{(-q_0)_{i_0} (q_0-1+\gamma )_{i_0}}{(1)_{i_0}} w_{1,2}^{i_0} \right\} \mu ^2 \label{eq:110051}
\end{eqnarray}
where
\begin{equation}
\widetilde{w}_{2,2} =  \frac{t_2 v_2}{u_2 (1-t_2)(1+v_2)}\Bigg|_{ v_2 = \frac{1}{2}\left( -1+\sqrt{1+4\eta u_2(1-t_2)s_{2,\infty }}\right)} = \frac{t_2}{u_2(1-t_2)}\frac{-1+\sqrt{1+4\eta u_2(1-t_2)s_{2,\infty }}}{1+\sqrt{1+4\eta u_2(1-t_2)s_{2,\infty }}} \nonumber
\end{equation}
Replace $q_i$, $q_j$ and $s_i$ by $q_1$, $q_0$ and ${ \displaystyle \frac{\widetilde{w}_{2,2} u_1 (1-t_1)s_1}{v_1(1+v_1)}}$ in (\ref{eq:110042}). Take the new (\ref{eq:110042}) into (\ref{eq:110051}).
\begin{eqnarray}
&&\sum_{q_0 =0}^{\infty } \frac{s_0^{q_0}}{ q_0 !}  \prod _{n=1}^{\infty } \left\{ \sum_{ q_n = q_{n-1}}^{\infty } s_n^{q_n }\right\} y_2(x) \nonumber\\
&&=  \prod_{k=3}^{\infty } \frac{1}{(1-s_{k,\infty })}  \int_{0}^{1} dt_2\;t_2^3 \int_{0}^{\infty } du_2\; \exp\left( -u_2 +\frac{-1+\sqrt{1+4\eta u_2(1-t_2)s_{2,\infty }}}{2\eta (1-t_2)}\right) \nonumber\\
&&\times \left( \frac{2}{1+\sqrt{1+4\eta u_2(1-t_2)s_{2,\infty }}}\right)^{6+\gamma } \frac{\widetilde{w}_{2,2}^{-2-\alpha } \left(  \widetilde{w}_{2,2} \partial _{\widetilde{w}_{2,2}} \right) \widetilde{w}_{2,2}^{2+\alpha }}{\sqrt{1+4\eta u_2(1-t_2)s_{2,\infty }}} \nonumber\\
&&\times \int_{0}^{1} dt_1\;t_1 \int_{0}^{\infty } du_1\;e^{-u_1} 
\frac{1}{2\pi i} \int_{\infty }^{(0+)} dv_1\; \exp\left(\frac{ v_1}{\widetilde{w}_{2,2} (1-t_1)}\right) \nonumber\\
&&\times \frac{1}{(1+v_1)^{2+\gamma }} \frac{w_{1,2}^{-\alpha } \left( w_{1,2}\partial_{w_{1,2}} \right) w_{1,2}^{\alpha }}{v_1^2+v_1-\widetilde{w}_{2,2}u_1(1-t_1)s_1} \nonumber\\
&&\times  \sum_{q_0 =0}^{\infty } \left( \frac{\widetilde{w}_{2,2} u_1 (1-t_1)s_{0,1}}{v_1(1+v_1)}\right)^{q_0}  \frac{1}{q_0!} \left\{ \sum_{i_0=0}^{q_0} \frac{(-q_0)_{i_0} (q_0-1+\gamma )_{i_0}}{(1)_{i_0}} w_{1,2}^{i_0} \right\} \mu ^2 \label{eq:110052}
\end{eqnarray}
By using Cauchy's integral formula, the contour integrand has poles at $v_1 = \frac{-1-\sqrt{1+4\widetilde{w}_{2,2} u_1(1-t_1)s_1}}{2}$ 
or $\frac{-1+\sqrt{1+4\widetilde{w}_{2,2} u_1(1-t_1)s_1}}{2}$ and $\frac{-1+\sqrt{1+4\widetilde{w}_{2,2} u_1(1-t_1)s_1}}{2}$ is only inside the unit circle. As we compute the residue there in (\ref{eq:110052}) we obtain
\begin{eqnarray}
&&\sum_{q_0 =0}^{\infty } \frac{s_0^{q_0}}{ q_0 !}  \prod _{n=1}^{\infty } \left\{ \sum_{ q_n = q_{n-1}}^{\infty } s_n^{q_n }\right\} y_2(x) \nonumber\\
&&=  \prod_{k=3}^{\infty } \frac{1}{(1-s_{k,\infty })}  \int_{0}^{1} dt_2\;t_2^3 \int_{0}^{\infty } du_2\; \exp\left( -u_2 +\frac{-1+\sqrt{1+4\eta u_2(1-t_2)s_{2,\infty }}}{2\eta (1-t_2)}\right) \nonumber\\
&&\times \left( \frac{2}{1+\sqrt{1+4\eta u_2(1-t_2)s_{2,\infty }}}\right)^{6+\gamma } \frac{\widetilde{w}_{2,2}^{-2-\alpha } \left(  \widetilde{w}_{2,2} \partial _{\widetilde{w}_{2,2}} \right) \widetilde{w}_{2,2}^{2+\alpha }}{\sqrt{1+4\eta u_2(1-t_2)s_{2,\infty }}} \nonumber\\
&&\times \int_{0}^{1} dt_1\;t_1 \int_{0}^{\infty } du_1\; \exp\left( -u_1 +\frac{-1+\sqrt{1+4\widetilde{w}_{2,2} u_1(1-t_1)s_1}}{2\widetilde{w}_{2,2} (1-t_1)}\right) \nonumber\\
&&\times \left( \frac{2}{1+\sqrt{1+4\widetilde{w}_{2,2} u_1(1-t_1)s_1}}\right)^{2+\gamma } \frac{\widetilde{w}_{1,2}^{-\alpha } \left(  \widetilde{w}_{1,2} \partial _{\widetilde{w}_{1,2}} \right) \widetilde{w}_{1,2}^{\alpha }}{\sqrt{1+4\widetilde{w}_{2,2} u_1(1-t_1)s_1}} \nonumber\\
&&\times  \sum_{q_0 =0}^{\infty } \frac{s_0^{q_0}}{q_0!} \left\{ \sum_{i_0=0}^{q_0} \frac{(-q_0)_{i_0} (q_0-1+\gamma )_{i_0}}{(1)_{i_0}} \widetilde{w}_{1,2}^{i_0} \right\} \mu ^2 \label{eq:110053}
\end{eqnarray}
where
\begin{equation}
\widetilde{w}_{1,2} = \frac{t_1 v_1}{u_1 (1-t_1)(1+v_1)}\Bigg|_{ v_1 = \frac{1}{2}\left( -1+\sqrt{1+4\widetilde{w}_{2,2} u_1(1-t_1)s_1}\right) }=\frac{t_1}{u_1(1-t_1)}\frac{-1+\sqrt{1+4\widetilde{w}_{2,2} u_1(1-t_1)s_1}}{1+\sqrt{1+4\widetilde{w}_{2,2} u_1(1-t_1)s_1}} \nonumber
\end{equation}
Acting the summation operator $ \sum_{q_0 =0}^{\infty } \frac{s_0^{q_0}}{ q_0 !}  \prod _{n=1}^{\infty } \left\{ \sum_{ q_n = q_{n-1}}^{\infty } s_n^{q_n }\right\} $ on (\ref{eq:110039}) by using (\ref{eq:110041}) and (\ref{eq:110042}),
\begin{eqnarray}
&&\sum_{q_0 =0}^{\infty } \frac{s_0^{q_0}}{ q_0 !}  \prod _{n=1}^{\infty } \left\{ \sum_{ q_n = q_{n-1}}^{\infty } s_n^{q_n }\right\} y_3(x) \nonumber\\
&&=  \prod_{k=4}^{\infty } \frac{1}{(1-s_{k,\infty })}  \int_{0}^{1} dt_3\;t_3^5 \int_{0}^{\infty } du_3\; \exp\left( -u_3 +\frac{-1+\sqrt{1+4\eta u_3(1-t_3)s_{3,\infty }}}{2\eta (1-t_3)}\right) \nonumber\\
&&\times \left( \frac{2}{1+\sqrt{1+4\eta u_3(1-t_3)s_{3,\infty }}}\right)^{10+\gamma } \frac{\widetilde{w}_{3,3}^{-4-\alpha } \left(  \widetilde{w}_{3,3} \partial _{\widetilde{w}_{3,3}} \right) \widetilde{w}_{3,3}^{4+\alpha }}{\sqrt{1+4\eta u_3(1-t_3)s_{3,\infty }}} \nonumber\\
&&\times \int_{0}^{1} dt_2\;t_2^3 \int_{0}^{\infty } du_2\; \exp\left( -u_2 +\frac{-1+\sqrt{1+4\widetilde{w}_{3,3} u_2(1-t_2)s_2}}{2\widetilde{w}_{3,3} (1-t_2)}\right) \nonumber\\
&&\times \left( \frac{2}{1+\sqrt{1+4\widetilde{w}_{3,3} u_2(1-t_2)s_2}}\right)^{6+\gamma } \frac{\widetilde{w}_{2,3}^{-2-\alpha } \left(  \widetilde{w}_{2,3} \partial _{\widetilde{w}_{2,3}} \right) \widetilde{w}_{2,3}^{2+\alpha }}{\sqrt{1+4\widetilde{w}_{3,3} u_2(1-t_2)s_2}} \nonumber\\
&&\times \int_{0}^{1} dt_1\;t_1 \int_{0}^{\infty } du_1\; \exp\left( -u_1 +\frac{-1+\sqrt{1+4\widetilde{w}_{2,3} u_1(1-t_1)s_1}}{2\widetilde{w}_{2,3} (1-t_1)}\right) \nonumber\\
&&\times \left( \frac{2}{1+\sqrt{1+4\widetilde{w}_{2,3} u_1(1-t_1)s_1}}\right)^{2+\gamma } \frac{\widetilde{w}_{1,3}^{-\alpha } \left( \widetilde{w}_{1,3} \partial _{\widetilde{w}_{1,3}} \right) \widetilde{w}_{1,3}^{\alpha }}{\sqrt{1+4\widetilde{w}_{2,3} u_1(1-t_1)s_1}} \nonumber\\
&&\times \sum_{q_0 =0}^{\infty }  \frac{s_0^{q_0}}{q_0!} \left\{ \sum_{i_0=0}^{q_0} \frac{(-q_0)_{i_0}(q_0-1+\gamma )_{i_0}}{(1)_{i_0}} \widetilde{w}_{1,3}^{i_0} \right\} \mu ^3 \label{eq:110054}
\end{eqnarray}
where
\begin{equation}
\begin{cases} \displaystyle{\widetilde{w}_{3,3}=\frac{t_3 v_3}{u_3 (1-t_3)(1+v_3)}\Bigg|_{v_3 = \frac{1}{2}\left( -1+\sqrt{1+4\eta u_3(1-t_3)s_{3,\infty }}\right)}=\frac{t_3}{u_3(1-t_3)}\frac{-1+\sqrt{1+4\eta u_3(1-t_3)s_{3,\infty }}}{1+\sqrt{1+4\eta u_3(1-t_3)s_{3,\infty }}}} \cr
\displaystyle{\widetilde{w}_{2,3} =  \frac{t_2 v_2}{u_2 (1-t_2)(1+v_2)}\Bigg|_{ v_2 = \frac{1}{2}\left( -1+\sqrt{1+4\widetilde{w}_{3,3} u_2(1-t_2)s_2}\right)}=\frac{t_2}{u_2(1-t_2)}\frac{-1+\sqrt{1+4\widetilde{w}_{3,3} u_2(1-t_2)s_2}}{1+\sqrt{1+4\widetilde{w}_{3,3} u_2(1-t_2)s_2}}}  \cr
\displaystyle{\widetilde{w}_{1,3} =  \frac{t_1 v_1}{u_1 (1-t_1)(1+v_1)}\Bigg|_{v_1 = \frac{1}{2} \left(-1+\sqrt{1+4\widetilde{w}_{2,3} u_1(1-t_1)s_1}\right)}=\frac{t_1}{u_1(1-t_1)}\frac{-1+\sqrt{1+4\widetilde{w}_{2,3} u_1(1-t_1)s_1}}{1+\sqrt{1+4\widetilde{w}_{2,3} u_1(1-t_1)s_1}}}  
\end{cases}\nonumber  
\end{equation}
By repeating this process for all higher terms of integral forms of sub-summation $y_m(x)$ terms where $m > 3$, we obtain every  $  \sum_{q_0 =0}^{\infty } \frac{s_0^{q_0}}{ q_0 !}  \prod _{n=1}^{\infty } \left\{ \sum_{ q_n = q_{n-1}}^{\infty } s_n^{q_n }\right\} y_m(x) $ terms. 
Since we substitute (\ref{eq:110045}), (\ref{eq:110048}), (\ref{eq:110053}), (\ref{eq:110054}) and including all $ \sum_{q_0 =0}^{\infty } \frac{s_0^{q_0}}{ q_0 !}  \prod _{n=1}^{\infty } \left\{ \sum_{ q_n = q_{n-1}}^{\infty } s_n^{q_n }\right\} y_m(x) $ terms where $m > 3$ into (\ref{eq:110044}), we obtain (\ref{eq:110043})
\qed
\end{proof}
\begin{remark}
The generating function for the DCH polynomial of type 2 of the first kind about $x=0 $ as  $q = \left( q_j +2j \right)\left( q_j +2j-1+\gamma \right)$ where $j,q_j \in \mathbb{N}_{0}$ is
\begin{eqnarray}
&&\sum_{q_0 =0}^{\infty } \frac{s_0^{q_0}}{q_0!}  \prod _{n=1}^{\infty } \left\{ \sum_{q_n =q_{n-1}}^{\infty } s_n^{q_n }\right\} H_d^{(o)}F_{q_{j}}^R\left( \alpha ,\beta ,\gamma ,\delta ,q= \left( q_j+2j \right)\left( q_j+2j-1+\gamma \right); \eta =-\frac{1}{\delta }x, \mu =-\frac{\beta }{\delta }x^2\right) \nonumber\\
&&=  \prod_{k=1}^{\infty } \frac{1}{(1-s_{k,\infty })} \mathbf{A} \left( s_{0,\infty } ;\eta \right) \nonumber\\
&&+ \left\{ \prod_{k=2}^{\infty } \frac{1}{(1-s_{k,\infty })} \int_{0}^{1} dt_1\;t_1 \int_{0}^{\infty } du_1\; \overleftrightarrow {\mathbf{\Gamma}}_1 \left(s_{1,\infty };t_1,u_1,\eta \right)   \widetilde{w}_{1,1}^{-\alpha } \left(  \widetilde{w}_{1,1} \partial _{\widetilde{w}_{1,1}} \right) \widetilde{w}_{1,1}^{ \alpha }\; \mathbf{A} \left(s_{0} ;\widetilde{w}_{1,1} \right) \right\}\mu \nonumber\\
&&+ \sum_{n=2}^{\infty } \Bigg\{ \frac{(1-s_{n,\infty })}{\prod_{k=n}^{\infty }(1-s_{k,\infty })} \int_{0}^{1} dt_n\;t_n^{2n-1} \int_{0}^{\infty } du_n\; \overleftrightarrow {\mathbf{\Gamma}}_n \left(s_{n,\infty };t_n,u_n,\eta \right)  \widetilde{w}_{n,n}^{-2(n-1)-\alpha } \left(  \widetilde{w}_{n,n} \partial _{\widetilde{w}_{n,n}} \right) \widetilde{w}_{n,n}^{2(n-1)+\alpha } \nonumber\\
&&\times  \prod_{j=1}^{n-1} \Bigg\{ \int_{0}^{1} dt_{n-j}\;t_{n-j}^{2(n-j)-1} \int_{0}^{\infty } du_{n-j}\; \overleftrightarrow {\mathbf{\Gamma}}_{n-j} \left(s_{n-j};t_{n-j},u_{n-j},\widetilde{w}_{n-j+1,n} \right) \nonumber\\
&&\times \widetilde{w}_{n-j,n}^{-2(n-j-1)-\alpha } \left(  \widetilde{w}_{n-j,n} \partial _{\widetilde{w}_{n-j,n}} \right) \widetilde{w}_{n-j,n}^{2(n-j-1)+\alpha } \Bigg\}\mathbf{A} \left(s_{0} ;\widetilde{w}_{1,n} \right) \Bigg\} \mu ^n   \label{eq:110055}
\end{eqnarray}
where
\begin{equation}
\begin{cases} 
{ \displaystyle \overleftrightarrow {\mathbf{\Gamma}}_1 \left(s_{1,\infty };t_1,u_1,\eta \right) }\cr
{ \displaystyle = \exp\left( -u_1 +\frac{-1+\sqrt{1+4\eta u_1(1-t_1)s_{1,\infty }}}{2\eta (1-t_1)}\right)  \frac{\left( \frac{2}{1+\sqrt{1+4\eta u_1(1-t_1)s_{1,\infty }}}\right)^{2+\gamma }}{\sqrt{ 1+4\eta u_1(1-t_1)s_{1,\infty }}}}\cr
{ \displaystyle  \overleftrightarrow {\mathbf{\Gamma}}_n \left(s_{n,\infty };t_n,u_n,\eta \right) }\cr
{ \displaystyle = \exp\left( -u_n +\frac{-1+\sqrt{1+4\eta u_n(1-t_n)s_{n,\infty }}}{2\eta (1-t_n)}\right)  \frac{\left( \frac{2}{1+\sqrt{1+4\eta u_n(1-t_n)s_{n,\infty }}}\right)^{4n-2+\gamma }}{\sqrt{1+4\eta u_n(1-t_n)s_{n,\infty }}}}\cr
{ \displaystyle \overleftrightarrow {\mathbf{\Gamma}}_{n-j} \left(s_{n-j};t_{n-j},u_{n-j},\widetilde{w}_{n-j+1,n} \right) }\cr
{ \displaystyle = \exp \left( -u_{n-j} +\frac{-1+\sqrt{1+4\widetilde{w}_{n-j+1,n} u_{n-j}(1-t_{n-j})s_{n-j}}}{2\widetilde{w}_{n-j+1,n} (1-t_{n-j})} \right) }\cr 
{ \displaystyle \hspace{.5cm}\times  \frac{\left( \frac{2}{1+\sqrt{1+4\widetilde{w}_{n-j+1,n} u_{n-j}(1-t_{n-j})s_{n-j}}}\right)^{4(n-j)-2+\gamma }}{\sqrt{1+4\widetilde{w}_{n-j+1,n} u_{n-j}(1-t_{n-j})s_{n-j}}}}
\end{cases}\nonumber  
\end{equation}
and
\begin{equation}
\begin{cases} 
{ \displaystyle \mathbf{A} \left( s_{0,\infty } ;\eta \right)  = \exp\left( \frac{-1+\sqrt{1+4\eta  s_{0,\infty }}}{2\eta }\right)  \frac{\left( \frac{2}{1+\sqrt{1+4\eta s_{0,\infty }}}\right)^{-2+\gamma }}{\sqrt{1+4\eta s_{0,\infty }}}}\cr
{ \displaystyle  \mathbf{A} \left(s_{0} ;\widetilde{w}_{1,1} \right)   = \exp\left( \frac{-1+\sqrt{1+4\widetilde{w}_{1,1}s_0}}{2\widetilde{w}_{1,1}}\right)  \frac{\left( \frac{2}{1+\sqrt{1+4\widetilde{w}_{1,1} s_0}}\right)^{-2+\gamma }}{\sqrt{1+4\widetilde{w}_{1,1}s_0}}}\cr
{ \displaystyle \mathbf{A} \left( s_{0} ;\widetilde{w}_{1,n} \right) = \exp \left( \frac{-1+\sqrt{1+4\widetilde{w}_{1,n}s_0}}{2\widetilde{w}_{1,n}} \right) \frac{\left( \frac{2}{1+\sqrt{1+4\widetilde{w}_{1,n}s_0}}\right)^{-2+\gamma }}{\sqrt{1+4\widetilde{w}_{1,n}s_0}}}
\end{cases}\nonumber  
\end{equation}
\end{remark}
\begin{proof}
Replace $a$, $b$, $j$ and $z$ by $-q_0$, $-2q_0+2-\gamma $, $i_0$ and $-\eta ^{-1}$ in (\ref{eq:110027}). Acting the summation operator $\displaystyle{ \sum_{q_0 =0}^{\infty } \frac{s_{0,\infty }^{q_0}}{q_0!}}$ on the new (\ref{eq:110027}) 
\begin{equation}
\sum_{q_0 =0}^{\infty } \frac{s_{0,\infty }^{q_0}}{q_0!} \sum_{i_0=0}^{q_0} \frac{(-q_0)_{i_0} (q_0-1+\gamma )_{i_0}}{(1)_{i_0}} \eta ^{i_0} = \sum_{q_0 =0}^{\infty } \frac{(-\eta s_{0,\infty } )^{q_0}}{q_0!} U\left( -q_0, -2q_0 +2-\gamma ,-\eta ^{-1}\right)   
\label{eq:110056}
\end{equation}
Replace $a$, $b$, $v_l$ and $z$ by $-q_0$, $-2q_0+2-\gamma $, $v$ and $-\eta ^{-1}$ in (\ref{eq:110026}).
\begin{equation}
U\left( -q_0, -2q_0+2-\gamma ,-\eta ^{-1}\right) =  \frac{q_0!}{2\pi i} \int_{\infty }^{(0+)} dv\;\frac{\exp\left(\frac{v}{\eta }\right)}{v(1+v)^{-1+\gamma }}\left( \frac{-1}{v(1+v)}\right)^{q_0}
\label{eq:110057}
\end{equation}
Put (\ref{eq:110057}) in (\ref{eq:110056}).
\begin{eqnarray}
\sum_{q_0 =0}^{\infty } \frac{s_{0,\infty }^{q_0}}{q_0!} \sum_{i_0=0}^{q_0} \frac{(-q_0)_{i_0} (q_0-1+\gamma )_{i_0}}{(1)_{i_0}} \eta ^{i_0} &=& \frac{1}{2\pi i} \int_{\infty }^{(0+)} dv\;\frac{\exp\left(\frac{v}{\eta }\right)}{v(1+v)^{-1+\gamma }} \sum_{q_0 =0}^{\infty } \left(\frac{\eta s_{0,\infty }}{v(1+v)}\right)^{q_0} \nonumber\\
&=& \frac{1}{2\pi i} \int_{\infty }^{(0+)} dv\; \exp\left(\frac{v}{\eta }\right) \frac{(1+v)^{2-\gamma }}{v^2+v-\eta s_{0,\infty }} 
\hspace{1.5cm}\label{eq:110058}
\end{eqnarray}
By using Cauchy's integral formula, the contour integrand has poles at $v = \frac{-1-\sqrt{1+4\eta s_{0,\infty }}}{2}$ 
or $\frac{-1+\sqrt{1+4\eta s_{0,\infty }}}{2}$ and $\frac{-1+\sqrt{1+4\eta s_{0,\infty }}}{2}$ is only inside the unit circle. As we compute the residue there in (\ref{eq:110058}) we obtain
\begin{align}
&\sum_{q_0 =0}^{\infty } \frac{s_{0,\infty }^{q_0}}{q_0!} \sum_{i_0=0}^{q_0}\frac{(-q_0)_{i_0} (q_0-1+\gamma )_{i_0}}{(1)_{i_0}} \eta ^{i_0} \nonumber\\
&= \exp\left( \frac{-1+\sqrt{1+4\eta  s_{0,\infty }}}{2\eta }\right)  \frac{\left( \frac{2}{1+\sqrt{1+4\eta s_{0,\infty }}}\right)^{-2+\gamma }}{\sqrt{1+4\eta s_{0,\infty }}} 
\label{eq:110059}
\end{align}
Replace $s_{0,\infty }$ and $\eta $  by $s_0$ and $\widetilde{w}_{1,1}$ in (\ref{eq:110059}).
\begin{align}
&\sum_{q_0 =0}^{\infty } \frac{s_0^{q_0}}{q_0!} \sum_{i_0=0}^{q_0} \frac{(-q_0)_{i_0} (q_0-1+\gamma )_{i_0}}{(1)_{i_0}} \widetilde{w}_{1,1}^{i_0} \nonumber\\
&= \exp\left( \frac{-1+\sqrt{1+4\widetilde{w}_{1,1}s_0}}{2\widetilde{w}_{1,1} }\right)  \frac{\left( \frac{2}{1+\sqrt{1+4\widetilde{w}_{1,1} s_0}}\right)^{-2+\gamma }}{\sqrt{1+4\widetilde{w}_{1,1} s_0}} 
\label{eq:110060}
\end{align}
Replace $s_{0,\infty }$ and $\eta $  by $s_0$ and $\widetilde{w}_{1,n}$ in (\ref{eq:110059}).
\begin{align}
&\sum_{q_0 =0}^{\infty } \frac{s_0^{q_0}}{q_0!} \sum_{i_0=0}^{q_0} \frac{(-q_0)_{i_0} (q_0-1+\gamma )_{i_0}}{(1)_{i_0}} \widetilde{w}_{1,n}^{i_0} \nonumber\\
&= \exp\left( \frac{-1+\sqrt{1+4\widetilde{w}_{1,n}s_0}}{2\widetilde{w}_{1,n} }\right) \frac{\left( \frac{2}{1+\sqrt{1+4\widetilde{w}_{1,n} s_0}}\right)^{-2+\gamma }}{\sqrt{1+4\widetilde{w}_{1,n} s_0}}  
\label{eq:110061}
\end{align}
We obtain (\ref{eq:110055}) by substituting (\ref{eq:110059}), (\ref{eq:110060}) and (\ref{eq:110061}) into (\ref{eq:110043}).
\qed
\end{proof}
\section[The DCHE with a irregular singular point at $x=\infty $]{The DCHE with a irregular singular point at infinity} 
\subsection{Power series for a polynomial of type 2}
Let $z=\frac{1}{x}$ in (\ref{eq:11003}) in order to obtain the Frobenius solution of the DCHE about $x=\infty $.
\begin{equation}
\frac{d^2{y}}{d{z}^2} + \left( -\delta +\frac{2-\gamma }{z} -\frac{\beta }{z^2}\right) \frac{d{y}}{d{z}} +  \frac{-qz+\alpha \beta }{z^3} y = 0 \label{eq:110062}
\end{equation}
We take the power series
\begin{equation}
y(z)= \sum_{n=0}^{\infty } c_n z^{n+\lambda }  \label{eq:110063}
\end{equation}
we obtain by substitution in (\ref{eq:110062}) a 3-term recurrence system between successive coefficients $c_n$:
\begin{equation}
c_{n+1}=A_n \;c_n +B_n \;c_{n-1} \hspace{1cm};n\geq 1 \label{eq:110064}
\end{equation}
where,
\begin{subequations}
\begin{align}
 A_n &= \frac{\left( n+\alpha \right)\left( n+1+\alpha -\gamma \right) -q}{\beta \left( n+1\right)} \label{eq:110065a}\\ 
&= \frac{\left( n+\frac{2\alpha -\gamma +1-\sqrt{\left( \gamma -1\right)^2 +4q}}{2}\right) \left( n+\frac{2\alpha -\gamma +1+\sqrt{\left( \gamma -1\right)^2 +4q}}{2}\right)}{\beta \left( n+1\right)} \label{eq:110065b}\\
B_n &= -\frac{\delta \left( n-1+\alpha \right)}{\beta \left( n+1\right)}  \label{eq:110065c}\\
c_1 &= A_0 \;c_0 \label{eq:110065d}
\end{align}
\end{subequations} 
with $\beta  \ne 0$. We only have one indicial root such as $\lambda = \alpha $.

Let us allow to test for convergence of the analytic function $y(z)$. As $n\gg 1$, a recursion system of (\ref{eq:110064})--(\ref{eq:110065d}) is approximate to
\begin{subequations}
\begin{equation}
c_{n+1}=A\;c_n +B\;c_{n-1} \hspace{1cm}\mbox{as}\;\; n\geq 1 \label{eq:110066}
\end{equation}
where

\begin{tabular}{p{7.5cm}p{7.5cm}}
{\begin{align}
&\lim_{n\gg 1} A_n = A= \frac{n}{\beta } \rightarrow \infty \label{eq:110067a}\ 
\end{align}}
&
{\begin{align}
&\lim_{n\gg 1} B_n = B= -\frac{\delta }{\beta } \label{eq:110067b}\ 
\end{align} }
\end{tabular}
\end{subequations}

by letting $c_1\sim  A c_0$ for simple calculations of a convergence series for $y(z)$. 
A $y(z)$ series is not convergent any more as $n\gg 1$ in (\ref{eq:110067a}); it makes that formal series solutions for a polynomial of type 1 and an infinite series do not exist. Therefore, two types of power series solutions of the DCHE about $x=\infty $ are only available such as a polynomial of type 2 and a complete polynomial. 
For a polynomial of type 2, I treat $\alpha $, $\beta $, $\gamma $, $\delta $ as free variables and $q$ as a fixed value. For a complete polynomial, I treat $\beta $, $\gamma $, $\delta $ as free variables and $\alpha $, $q$ as fixed values. 

Like the case of the DCHE about $x=0$, for a polynomial of type 2, (\ref{eq:110067a}) is negligible for the minimum value of a $y(z)$ because $A_n$ term is terminated at the specific index summation $n$. (\ref{eq:110067b}) is only available for an asymptotic behavior of the minimum $y(x)$. 
Substituting (\ref{eq:110067b}) into (\ref{eq:110066}) with $A=0$ gives a recurrence system. 
And the coefficients $c_n$ are classified as even and odd terms for $n=0,1,2,\cdots$.
\begin{equation}
\begin{tabular}{  l  l }
  \vspace{2 mm}
   $c_0$ &\hspace{1cm}  $c_1$ \\
  \vspace{2 mm}
   $c_2 = -\frac{\delta}{\beta} c_0 $  &\hspace{1cm}  $c_3 = -\frac{\delta}{\beta} c_1$  \\
  \vspace{2 mm}
  $c_4 = \left(-\frac{\delta}{\beta}\right)^2 c_0 $ &\hspace{1cm}  $c_5 = \left(-\frac{\delta}{\beta}\right)^2 c_1$\\
  \vspace{2 mm}
  $c_6 = \left(-\frac{\delta}{\beta}\right)^3 c_0 $ &\hspace{1cm}  $c_7 = \left(-\frac{\delta}{\beta}\right)^3 c_1$\\
 \hspace{2 mm} \large{\vdots} & \hspace{1.5 cm}\large{\vdots} \\
 \vspace{2 mm}
  $c_{2n} = \left(-\frac{\delta}{\beta}\right)^n c_0 $ &\hspace{1cm}  $c_{2n+1} = \left(-\frac{\delta}{\beta}\right)^n c_1$\\
\end{tabular}
\label{eq:110068}
\end{equation}
$A_n$ term is negligible for the minimum $y(z)$ since $A_n$ term is terminated at the specific summation index $n$. Then we are allowed to $c_1\sim  A c_0 =0$ in (\ref{eq:110068}).
Put the coefficients $c_{2n}$ on the above into an asymptotic series $\sum_{n=0}^{\infty }c_{2n}x^{2n}$, putting $c_0=1$ for simplicity.
\begin{equation}
\mbox{min}\left( \lim_{n\gg 1}y(z) \right) = \frac{1}{1+\frac{\delta}{\beta}z^2} \label{eq:110069}
\end{equation}
A polynomial of type 2 requires $\left| \delta/\beta z^2\right|<1$ for the convergence of the radius. 
An asymptotic series (\ref{eq:110069}) is constructed as similarly as a series solution (\ref{eq:110017}).

We obtain an eigenvalue $q$ for a polynomial of type 2 by replacing $\alpha _i$ by $q_i$ and put $n= q_i+ 2i$ in (\ref{eq:110065b}) with a boundary condition $A_{q_i+ 2i}=0$.  
\begin{equation}
\pm \sqrt{\left( \gamma -1\right)^2 +4q}= 2\left( q_i +2i+ \alpha \right) -\gamma +1
\nonumber
\end{equation} 
\subsubsection{The case of $ \sqrt{\left( \gamma -1\right)^2 +4q}= 2\left( q_i +2i+ \alpha \right) -\gamma +1$ where $i,q_i \in \mathbb{N}_{0}$}

In (\ref{eq:110065b}) replace $\sqrt{\left( \gamma -1\right)^2 +4q}$ by $2\left( q_i +2i+ \alpha \right) -\gamma +1$. In (\ref{eq:110018}) replace a variable $x$ and an index $\alpha _i$ by $z$ and $q_i$.  Take the new (\ref{eq:110065b}) and (\ref{eq:110065c}) in the new (\ref{eq:110018}) with $c_0=1$ and $\lambda =\alpha $.
After the replacement process, 
\begin{remark}
The power series expansion of the DCHE of the first kind for a polynomial of type 2 about $x=\infty $ as $q=\left( q_j+2j +\alpha \right)\left( q_j+2j+1+\alpha -\gamma \right) $ where $j,q_j \in \mathbb{N}_{0}$ is
\begin{eqnarray}
 y(z)&=& \sum_{n=0}^{\infty } y_{n}(z) = y_0(z)+ y_1(z)+ y_2(z)+y_3(z)+\cdots \nonumber\\ 
&=& H_d^{(i)}F_{q_{j}}^R\left( \alpha ,\beta ,\gamma ,\delta ,q= \left( q_j+2j +\alpha \right)\left( q_j+2j+1+\alpha -\gamma \right); z=\frac{1}{x}, \xi =\frac{1}{\beta }z, \rho =-\frac{\delta }{\beta }z^2\right)\nonumber\\
&=& z^{\alpha }\left\{ \sum_{i_0=0}^{q_0} \frac{\left( -q_0\right)_{i_0} \left( q_0+1+2\alpha -\gamma \right)_{i_0}}{(1)_{i_0}} \xi ^{i_0} \right. \nonumber\\
&&+ \left\{ \sum_{i_0=0}^{q_0}\frac{ (i_0+ \alpha ) }{ (i_0+ 2) }\frac{\left( -q_0\right)_{i_0} \left( q_0+1+2\alpha -\gamma \right)_{i_0} }{(1)_{i_0} } \right.   \left. \sum_{i_1=i_0}^{q_1} \frac{\left( -q_1\right)_{i_1}\left( q_1 +5+2\alpha -\gamma \right)_{i_1}(3)_{i_0} }{\left( -q_1\right)_{i_0}\left( q_1 +5+2\alpha -\gamma \right)_{i_0}(3)_{i_1}} \xi ^{i_1}\right\} \rho \nonumber\\
&&+ \sum_{n=2}^{\infty } \left\{ \sum_{i_0=0}^{q_0}\frac{ (i_0+ \alpha ) }{ (i_0+ 2) }\frac{\left( -q_0\right)_{i_0} \left( q_0+1+2\alpha -\gamma \right)_{i_0} }{(1)_{i_0} }\right.\nonumber\\
&&\times \prod _{k=1}^{n-1} \left\{ \sum_{i_k=i_{k-1}}^{q_k} \frac{ (i_k+ 2k +\alpha )}{ (i_k+2k+2) } \frac{\left( -q_k\right)_{i_k}\left( q_k+4k +1+2\alpha -\gamma \right)_{i_k}(2k+1)_{i_{k-1}} }{\left( -q_k\right)_{i_{k-1}}\left( q_k+4k +1+2\alpha -\gamma \right)_{i_{k-1}}(2k+1)_{i_k}}\right\} \nonumber\\
&&\times \left.\left.\sum_{i_n= i_{n-1}}^{q_n} \frac{\left( -q_n\right)_{i_n}\left( q_n+4n +1+2\alpha -\gamma  \right)_{i_n}(2n+1)_{i_{n-1}} }{\left( -q_n\right)_{i_{k-1}}\left( q_n+4n +1+2\alpha -\gamma \right)_{i_{n-1}}(2n+1)_{i_n}} \xi ^{i_n} \right\} \rho ^n \right\} \label{eq:110070}
\end{eqnarray}
\end{remark} 
\subsubsection{The case of $ \sqrt{\left( \gamma -1\right)^2 +4q}= -2\left( q_i +2i+ \alpha \right) +\gamma -1$}
In (\ref{eq:110065b}) replace $\sqrt{\left( \gamma -1\right)^2 +4q}$ by $ -2\left( q_i +2i+ \alpha \right) +\gamma -1$. In (\ref{eq:110018}) replace a variable $x$ and an index $\alpha _i$ by $z$ and $q_i$.  Take the new (\ref{eq:110065b}) and (\ref{eq:110065c}) in the new (\ref{eq:110018}) with $c_0=1$ and $\lambda =\alpha $.
After the replacement process, its solution is equivalent to (\ref{eq:110070}). 

For the minimum value of the DCHE for a polynomial of type 2 about $x=\infty $, put $q_0=q_1=q_2=\cdots=0$ in (\ref{eq:110070}).
\begin{eqnarray}
 y(z) &=& H_d^{(i)}F_{q_{j}}^R\left( \alpha ,\beta ,\gamma ,\delta ,q= \left( 2j +\alpha \right)\left( 2j+1+\alpha -\gamma \right); z=\frac{1}{x}, \xi =\frac{1}{\beta }z, \rho =-\frac{\delta }{\beta }z^2\right)\nonumber\\
&=& z^{\alpha } \left( 1-\rho \right)^{-\frac{\alpha }{2}} \hspace{1cm}\mbox{where}\;\;|\rho |<1 \label{eq:110071}
\end{eqnarray} 
\subsection{Integral representation for a polynomial of type 2}
General summation structures between (\ref{eq:110020}) and (\ref{eq:110070}) are similar to each other.
As we compare all coefficients and variables in (\ref{eq:110070}) with (\ref{eq:110020}), we find that 
\begin{equation}
\begin{cases} x \rightarrow  z=\frac{1}{x}  \cr
y(x) \rightarrow  z^{-\alpha } y(z)  \cr
q= \left( q_j+2j \right)\left( q_j+2j-1+\gamma \right) \rightarrow  q= \left( q_j+2j +\alpha \right)\left( q_j+2j+1+\alpha -\gamma \right) \cr
\gamma  \rightarrow 2+2\alpha -\gamma  \cr
\eta =-\frac{1}{\delta }x  \rightarrow \xi =\frac{1}{\beta }z\cr
\mu =-\frac{\beta }{\delta }x^2 \rightarrow \rho =-\frac{\delta }{\beta }z^2   
\end{cases}\label{eq:110072}
\end{equation}
We obtain an integral form of $y(z)$ by substituting (\ref{eq:110072}) in (\ref{eq:110040}) such as  
\begin{remark}
The integral representation of the DCHE of the first kind for a polynomial of type 2 about $x=\infty $ as $q= \left( q_j+2j +\alpha \right)\left( q_j+2j+1+\alpha -\gamma \right)$ where $j, q_j \in \mathbb{N}_{0}$ is
\begin{eqnarray}
y(z)&=& H_d^{(i)}F_{q_{j}}^R\left( \alpha ,\beta ,\gamma ,\delta ,q= \left( q_j+2j +\alpha \right)\left( q_j+2j+1+\alpha -\gamma \right); z=\frac{1}{x}, \xi =\frac{1}{\beta }z, \rho =-\frac{\delta }{\beta }z^2\right)\nonumber\\
&=& z^{\alpha }\Bigg\{ \left(-\xi \right)^{q_0} U\left( -q_0, -2q_0-2\alpha +\gamma ,-\xi ^{-1}\right)  
+ \sum_{n=1}^{\infty } \Bigg\{\prod _{k=0}^{n-1} \Bigg\{ \int_{0}^{1} dt_{n-k}\;t_{n-k}^{2(n-k)-1} \int_{0}^{\infty } du_{n-k}\;\exp\left( -u_{n-k}\right)  \nonumber\\
&&\times \frac{1}{2\pi i} \int_{\infty }^{(0+)} dv_{n-k}\; \frac{1}{v_{n-k}\left( 1+v_{n-k}\right)^{4(n-k)+1+2\alpha -\gamma }} \exp \left( \frac{v_{n-k}}{w_{n-k+1,n}(1-t_{n-k})}\right)  \nonumber\\
&&\times \left( \frac{w_{n-k+1,n} u_{n-k} (1-t_{n-k})}{v_{n-k}\left( 1+v_{n-k}\right)}\right)^{q_{n-k}} w_{n-k,n}^{-2(n-k-1)-\alpha } \left(  w_{n-k,n} \partial _{w_{n-k,n}} \right) w_{n-k,n}^{2(n-k-1)+\alpha }\Bigg\}  \nonumber\\
&&\times \left( -w_{1,n}\right)^{q_0} U\left( -q_0, -2q_0-2\alpha +\gamma , -w_{1,n}^{-1}\right)\Bigg\} \rho ^n \Bigg\}  
\label{eq:110073}
\end{eqnarray}
where
\begin{equation} w_{i,j}=
\begin{cases} \displaystyle {\frac{t_i v_i}{u_i (1-t_i)(1+v_i)} }\;\;\mbox{where}\; i\leq j\cr
\xi \;\;\mbox{only}\;\mbox{if}\; i>j
\end{cases}
\nonumber
\end{equation}
\end{remark} 
\subsection{Generating function for a polynomial of type 2}
We find the generating function of $y(z)$ by substituting (\ref{eq:110072}) in (\ref{eq:110041}) and (\ref{eq:110055}) such as
\begin{remark}
The generating function for the DCH polynomial of type 2 of the first kind about $x=\infty $ as  $q= \left( q_j+2j +\alpha \right)\left( q_j+2j+1+\alpha -\gamma \right)$ where $j,q_j \in \mathbb{N}_{0}$ is
\begin{eqnarray}
&&\sum_{q_0 =0}^{\infty } \frac{s_0^{q_0}}{q_0!}  \prod _{n=1}^{\infty } \left\{ \sum_{q_n =q_{n-1}}^{\infty } s_n^{q_n }\right\} H_d^{(i)}F_{q_{j}}^R\left( \alpha ,\beta ,\gamma ,\delta ,q= \left( q_j+2j +\alpha \right)\left( q_j+2j+1+\alpha -\gamma \right); z=\frac{1}{x}\right. \nonumber\\
&&\hspace{1cm},\left. \xi =\frac{1}{\beta }z, \rho =-\frac{\delta }{\beta }z^2\right) \nonumber\\
&&= z^{\alpha }\Bigg\{ \prod_{k=1}^{\infty } \frac{1}{(1-s_{k,\infty })} \mathbf{A} \left( s_{0,\infty } ;\xi \right) \nonumber\\
&&+ \left\{ \prod_{k=2}^{\infty } \frac{1}{(1-s_{k,\infty })} \int_{0}^{1} dt_1\;t_1 \int_{0}^{\infty } du_1\; \overleftrightarrow {\mathbf{\Gamma}}_1 \left(s_{1,\infty };t_1,u_1,\xi \right) \widetilde{w}_{1,1}^{-\alpha } \left(  \widetilde{w}_{1,1} \partial _{\widetilde{w}_{1,1}} \right) \widetilde{w}_{1,1}^{ \alpha }\; \mathbf{A} \left(s_{0} ;\widetilde{w}_{1,1} \right) \right\}\rho \nonumber\\
&&+ \sum_{n=2}^{\infty } \Bigg\{ \frac{(1-s_{n,\infty })}{\prod_{k=n}^{\infty }(1-s_{k,\infty })} \int_{0}^{1} dt_n\;t_n^{2n-1} \int_{0}^{\infty } du_n\; \overleftrightarrow {\mathbf{\Gamma}}_n \left(s_{n,\infty };t_n,u_n,\xi \right) \nonumber\\
&&\times \widetilde{w}_{n,n}^{-2(n-1)-\alpha } \left(  \widetilde{w}_{n,n} \partial _{\widetilde{w}_{n,n}} \right) \widetilde{w}_{n,n}^{2(n-1)+\alpha } \nonumber\\
&&\times  \prod_{j=1}^{n-1} \Bigg\{ \int_{0}^{1} dt_{n-j}\;t_{n-j}^{2(n-j)-1} \int_{0}^{\infty } du_{n-j}\; \overleftrightarrow {\mathbf{\Gamma}}_{n-j} \left(s_{n-j};t_{n-j},u_{n-j},\widetilde{w}_{n-j+1,n} \right) \nonumber\\
&&\times \widetilde{w}_{n-j,n}^{-2(n-j-1)-\alpha } \left(  \widetilde{w}_{n-j,n} \partial _{\widetilde{w}_{n-j,n}} \right) \widetilde{w}_{n-j,n}^{2(n-j-1)+\alpha } \Bigg\}\mathbf{A} \left(s_{0} ;\widetilde{w}_{1,n} \right) \Bigg\} \rho ^n \Bigg\} \hspace{1cm} \label{eq:110074}
\end{eqnarray}
where
\begin{equation}
\begin{cases} 
{ \displaystyle \overleftrightarrow {\mathbf{\Gamma}}_1 \left(s_{1,\infty };t_1,u_1,\xi \right) }\cr
{ \displaystyle = \exp\left( -u_1 +\frac{-1+\sqrt{1+4\xi u_1(1-t_1)s_{1,\infty }}}{2\xi (1-t_1)}\right) \frac{\left( \frac{2}{1+\sqrt{1+4\xi u_1(1-t_1)s_{1,\infty }}}\right)^{4+2\alpha -\gamma }}{\sqrt{1+4\xi u_1(1-t_1)s_{1,\infty }}}}\cr
{ \displaystyle  \overleftrightarrow {\mathbf{\Gamma}}_n \left(s_{n,\infty };t_n,u_n,\xi \right) }\cr
{ \displaystyle = \exp\left( -u_n +\frac{-1+\sqrt{1+4\xi u_n(1-t_n)s_{n,\infty }}}{2\xi (1-t_n)}\right)  \frac{\left( \frac{2}{1+\sqrt{1+4\xi u_n(1-t_n)s_{n,\infty }}}\right)^{4n+2\alpha -\gamma }}{\sqrt{1+4\xi u_n(1-t_n)s_{n,\infty }}}}\cr
{ \displaystyle \overleftrightarrow {\mathbf{\Gamma}}_{n-j} \left(s_{n-j};t_{n-j},u_{n-j},\widetilde{w}_{n-j+1,n} \right) }\cr
{ \displaystyle = \exp \left( -u_{n-j} +\frac{-1+\sqrt{1+4\widetilde{w}_{n-j+1,n} u_{n-j}(1-t_{n-j})s_{n-j}}}{2\widetilde{w}_{n-j+1,n} (1-t_{n-j})} \right) }\cr 
{ \displaystyle \hspace{.5cm}\times \frac{\left( \frac{2}{1+\sqrt{1+4\widetilde{w}_{n-j+1,n} u_{n-j}(1-t_{n-j})s_{n-j}}}\right)^{4(n-j)+2\alpha -\gamma }}{\sqrt{1+4\widetilde{w}_{n-j+1,n} u_{n-j}(1-t_{n-j})s_{n-j}}}}
\end{cases}\nonumber  
\end{equation}
and
\begin{equation}
\begin{cases} 
{ \displaystyle \mathbf{A} \left( s_{0,\infty } ;\xi \right)  = \exp\left( \frac{-1+\sqrt{1+4\xi  s_{0,\infty }}}{2\xi }\right)  \frac{\left( \frac{2}{1+\sqrt{1+4\xi s_{0,\infty }}}\right)^{2\alpha -\gamma }}{\sqrt{1+4\xi s_{0,\infty }}} }\cr
{ \displaystyle  \mathbf{A} \left(s_{0} ;\widetilde{w}_{1,1} \right) = \exp\left( \frac{-1+\sqrt{1+4\widetilde{w}_{1,1}s_0}}{2\widetilde{w}_{1,1}}\right)  \frac{\left( \frac{2}{1+\sqrt{1+4\widetilde{w}_{1,1} s_0}}\right)^{2\alpha -\gamma }}{\sqrt{1+4\widetilde{w}_{1,1}s_0}}}\cr
{ \displaystyle \mathbf{A} \left(s_{0} ;\widetilde{w}_{1,n} \right) = \exp \left( \frac{-1+\sqrt{1+4\widetilde{w}_{1,n}s_0}}{2\widetilde{w}_{1,n}} \right) \frac{\left( \frac{2}{1+\sqrt{1+4\widetilde{w}_{1,n}s_0}}\right)^{2\alpha -\gamma }}{\sqrt{1+4\widetilde{w}_{1,n}s_0}}} 
\end{cases}\nonumber  
\end{equation}
with a definition such as
\begin{equation}
\displaystyle {  \widetilde{w}_{i,j}= \begin{cases} \displaystyle { \frac{t_i}{u_i(1-t_i)}\frac{-1+\sqrt{1+4\widetilde{w}_{i+1,j} u_i(1-t_i)s_i}}{1+\sqrt{1+4\widetilde{w}_{i+1,j} u_i(1-t_i)s_i}}}\;\;\mbox{where}\;i<j\;\;\mbox{and}\;i,j\in \mathbb{N}_{0} \cr
\displaystyle { \frac{t_i}{u_i(1-t_i)}\frac{-1+\sqrt{1+4\xi u_i(1-t_i)s_{i,\infty }}}{1+\sqrt{1+4\xi u_i(1-t_i)s_{i,\infty }}}}  \;\;\mbox{where}\;\;i=j\end{cases}}\nonumber
\end{equation}
\end{remark}
\section{Summary}

In chapter 1 of Ref.\cite{11Choun2013}, for $n=0,1,2,3,\cdots $, (\ref{ysc:1}) is expanded to combinations of $A_n$ and $B_n$ terms. I define that a sub-power series $y_l(x)$ where $l\in \mathbb{N}_{0}$ is constructed by observing the term of sequence $c_n$ which includes $l$ terms of $B_n's$. 
The formal series solution is described by sums of each $y_l(x)$ such as $y(x)= \sum_{n=0}^{\infty } y_n(x)$. By allowing $B_n$ in the sequence $c_n$ is the leading term of each sub-power series $y_l(x)$, the general summation formulas of the 3-term recurrence relation in a linear ODE are constructed for an infinite series and a polynomial of type 2, denominated as `reversible three term recurrence formula (R3TRF).' 
 
By applying $A_n$ and $B_n$ terms in the recurrence relations for the DCHE around the origin and infinity into R3TRF, Frobenius series solutions in closed forms for a polynomial of type 2 are constructed in use for hypergeometric-type equations. There are no such series solutions of the DCHE for an infinite series and a polynomial of type 1; because formal series solutions of the DCHE around $x=0$ and $x=\infty $ are divergent as $n\gg 1$ in $A_n$ term of a 3-term recursive relation.

For a polynomial of type 2, it discusses special parameter values $q$ for which sub-power series of the general formal series solutions are truncated to polynomials: For the type 2 DCH polynomial around $x=0$, I suggest that a parameter $q$ is equivalent to $ ( q_j+2j )( q_j+2j-1+\gamma ) $ where $j,q_j \in \mathbb{N}_{0}$; for its polynomial around $x=\infty $, I treat $q$ as $( q_j+2j +\alpha )( q_j+2j+1+\alpha -\gamma )$.

In particular, type 2 polynomial of the DCHE around the origin requires $\left| \beta /\delta x^2\right|<1$ for the convergence of the radius. Similarly, its type 2 polynomial around infinity necessitates $\left| \delta/\beta z^2\right|<1$. Otherwise, general series solutions for a polynomial of type 2 will not convergent any more with respect to an independent variable.

For summation series solution of the DCHE around $x=0$ and $x=\infty $, the denominators and numerators in all $A_n$ terms of sub-power series arise with Pochhammer symbols. 
Their representation in terms of integrals for a polynomial of type 2, akin to those of the hypergeometric functions, is constructed by applying the contour integral of Tricomi's functions (Kummer's function of the second kind) into sub-power series of their general series solutions; each sub-integral $y_l(x)$ is composed of $2l$ terms of definite integrals and $l$ terms of contour integrals.
Their generating functions for a polynomial of type 2 are obtained analytically by applying generating functions for Tricomi's polynomials into sub-integrals of general integral solutions. 

The normalized wave function of hydrogen-like atoms and expectation values of its physical quantities such as position and momentum are obtained by applying the generating function of associated Laguerre polynomials.
likewise, we are able to find orthogonal relations and recursion relations of the DCH polynomial of type 2 from their generating functions including physical expectation values.
   
\addcontentsline{toc}{section}{Bibliography}
\bibliographystyle{model1a-num-names}
\bibliography{<your-bib-database>}

\chapter{Complete polynomials of Double Confluent Heun equation}
\chaptermark{Complete polynomials of the DCHE} 

Reversible three term recurrence formula (R3TRF) in chapter 1 of Ref.\cite{12Choun2013} is applied into the DCH equation (DCHE) with irregular singular points at the origin and infinity for a polynomial which makes the series solutions are truncated with special parameter values.  
And their power series solutions in closed forms are derived for a polynomial of type 2 in chapter 12.
Also their combined definite and contour integrals involving only Tricomi's functions are constructed including generating functions for DCH polynomials of type 2.
 
In this chapter I construct Frobenius solutions of the DCH equation with two singular points for a polynomial of type 3 by applying general summation formulas of complete polynomials using 3TRF and R3TRF.

\section{Introduction}
In 1889, Karl Heun suggested a second linear ODE having four regular singular points such as \cite{12Heun1889,12Ronv1995}
\begin{equation}
\frac{d^2{y}}{d{x}^2} + \left(\frac{\gamma }{x} +\frac{\delta }{x-1} + \frac{\epsilon }{x-a}\right) \frac{d{y}}{d{x}} +  \frac{\alpha \beta x-q}{x(x-1)(x-a)} y = 0 \label{eq:12001}
\end{equation}
where $\epsilon = \alpha +\beta -\gamma -\delta +1$ for assuring the regularity of the point at $x=\infty $. Its equation, called Heun's differential equation, has four regular singular points which are 0, 1, $a$ and $\infty $ with exponents $\{ 0, 1-\gamma \}$, $\{ 0, 1-\delta \}$, $\{ 0, 1-\epsilon \}$ and $\{ \alpha, \beta \}$.  

There are 4 different types of confluent forms of Heun equation such as (1) confluent Heun (two regular and one irregular singularities), (2) double confluent Heun (two irregular singularities), (3) biconfluent Heun (one regular and one irregular singularities), (4) triconfluent Heun equations (one irregular singularity).
Like deriving of confluent hypergeometric equation from a hypergeometric equation, 4 confluent types of Heun equation can be derived from merging two or more regular singularities to take an irregular singularity in Heun equation.
 
The non-symmetrical canonical form of the confluent Heun equation (CHE) \cite{12Deca1978,12Decar1978,12Ronv1995} reads
\begin{equation}
\frac{d^2{y}}{d{x}^2} + \left(\beta  +\frac{\gamma }{x} + \frac{\delta }{x-1}\right) \frac{d{y}}{d{x}} +  \frac{\alpha \beta x-q}{x(x-1)} y = 0 \label{eq:12002}
\end{equation}
It has three singular points: two regular singular points which are 0 and 1 with exponents $\{0, 1-\gamma\}$ and $\{0, 1-\delta \}$, and one irregular singular point which is $\infty$ with an exponent $\alpha$.   

In general, there are three types of the DCH equation (DCHE): (1) non-symmetrical canonical form of the DCHE, (2) canonical form of the general DCHE, (3) generalized spheroidal equation in the Leaver version, by letting a regular singular point as zero \cite{12Leav1986}.
We obtain the DCHE by changing coefficients and combining two regular singularities in the CHE. Let us allow $x\rightarrow x/\epsilon $, $\beta \rightarrow \beta \epsilon $, $\gamma \rightarrow \gamma -\delta /\epsilon $ and $\delta \rightarrow \delta /\epsilon $ in (\ref{eq:12002}). Assuming $\epsilon \rightarrow 0$ in the new (\ref{eq:12002})
\begin{equation}
\frac{d^2{y}}{d{x}^2} + \left(\beta  +\frac{\gamma }{x} + \frac{\delta }{x^2}\right) \frac{d{y}}{d{x}} +  \frac{\alpha \beta x-q}{x^2} y = 0 \label{eq:12003}
\end{equation}
(\ref{eq:12003}) is the non-symmetrical canonical form of the DCHE which has irregular singularities at $x=0$ and $\infty $, each of rank 1. Its general solution is denoted as $H_d(\alpha,\beta,\gamma,\delta,q;x)$.
 For DLFM version \cite{12NIST} or in Ref.\cite{12Slavy2000}, replace $\beta $ by 1 in (\ref{eq:12003}). The reason, why the parameter $\beta $ is included in (\ref{eq:12003}) instead of the unity, is that we can obtain an equation of the Whittaker-Ince limit of the DCHE by putting $\beta \rightarrow 0$, $\alpha \rightarrow \infty $, such that $\alpha \beta \rightarrow \alpha $ in (\ref{eq:12003}). In this chapter, (\ref{eq:12003}) is selected as an analytic solution of the DCHE because of more convenient for mathematical computations and the application of the Laplace transform.

Recently The DCH functions started to appear in theoretical modern physics. For examples the DCHE appears in diverse areas such as Dirac equations in the Nutku's helicoid spacetime \cite{12Birk2007,12Hort2007} and the rotating electromagnetic spacetime \cite{12AlBa2008}, the Schr$\ddot{\mbox{o}}$dinger equation for inverse fourth and sixth-power potentials (asymmetric double-Morse potentials) \cite{12Buhr1974,12Figu2005,12Fran1971,12Klei1968}, the time-dependence of Klein-Gordon and Dirac test-fields in curved spacetimes \cite{12Beze2014,12Birr1982,12Cost1989,12Nove1979}, quantum scattering problems of non-relativistic electrons from intermolecular forces \cite{12Babe1935}, Teukolsky's equations in general relativity, linearized perturbation theories on the backgrounds of Schwarzschild and Kerr black holes, etc.

Since a power series with unknown coefficients is put into Heun-type equations except triconfluent Heun equation, the recurrence relation between consecutive coefficients has a 3-term. Currently, it is impossible to obtain Frobenius solutions for a 3-term recursive relation in which coefficients are given explicitly. Indeed, their representation in terms of definite or contour integrals are still ambiguous including even their numerical computation. 
In contrast to Heun-type equations, formal series in compact forms of hypergeometric-type equations are already built by great many mathematicians extensively including their integral forms. Because their recurrence relation between successive coefficients can always be made to yield a 2-term, which is easily handled. 
\section[4 power series solutions in a 3-term recursive relation]{4 power series solutions in linear ODEs having a 3-term recursive relation}
By substituting a power series with unknown coefficients such as $y(x)=\sum_{n=0}^{\infty }c_nx^{n+\lambda }$ into linear ODEs, we obtain a 3-term recurrence system such as
\begin{equation}
c_{n+1}=A_n \;c_n +B_n \;c_{n-1} \hspace{1cm};n\geq 1
\label{eq:12004}
\end{equation}
where $\lambda =$ an indicial root, $c_1= A_0 \;c_0$ and $ c_0 \ne 0$ . On the above, $A_n$ and $B_n$ are themselves polynomials of degree $m$: for the second-order ODEs, a numerator and a denominator of $A_n$ are usually equal or less than polynomials of degrees 2. 
 
In linear ODEs consisting of a 3-term recursive relation between successive coefficient, there are 4 types power series solutions such as an infinite series and 3 possible polynomials which makes the series solutions are truncated with fixed parameter values such as (1) a polynomial which makes $B_n$ term terminated; $A_n$ term is not terminated, (2) a polynomial which makes $A_n$ term terminated; $B_n$ term is not terminated, (3) a polynomial which makes $A_n$ and $B_n$ terms terminated at the same time, referred as `a complete polynomial.'    
 
In general, an infinite series and 3 different polynomials of linear ODEs having a 3-term recursive relation are constructed by two possible general summation formulas such as a three term recurrence formula (3TRF) and a reversible three term recurrence formula (R3TRF). 

For $n=0,1,2,3,\cdots $ in (\ref{eq:12004}), the coefficients $c_{n+1}$ of the series are expanded to combinations of $A_n$ and $B_n$ terms. First of all, I suggest that a sub-power series $y_l(x)$ where $l\in \mathbb{N}_{0}$ is composed by observing the term of sequence $c_n$ which includes $l$ terms of $A_n's$. And the power series solution is defined by sums of each $y_l(x)$ such as $y(x)= \sum_{n=0}^{\infty } y_n(x)$. The Frobenius solutions in linear ODEs having a 3-term recurrence relation between successive coefficients are constructed, since $A_n$ in the sequence $c_n$ is allowed as the leading term of each sub-power series $y_l(x)$ for an infinite series and a polynomial of type 1. In particular, a polynomial of type 1 is determined that sub-power series solutions $y_l(x)$ are truncated to polynomials with fixed parameter values in a numerator of $B_n$ term. This general summation formulas is designated as `three term recurrence formula (3TRF).' \cite{12Chou2012b} 

In chapter 1 of Ref.\cite{12Choun2013}, by contrast with 3TRF, I define that a sub-power series $y_l(x)$ is obtained by observing the term of sequence $c_n$ which includes $l$ terms of $B_n's$. And the general solution in a formal series is described by sums of each $y_l(x)$. By allowing $B_n$ in the sequence $c_n$ is the leading term of each sub-power series $y_l(x)$, the general summation expressions of the 3-term recurrence relation in linear ODEs are constructed for an infinite series and a polynomial of type 2, designated as `reversible three term recurrence formula (R3TRF).'
A polynomial of type 2 is specified that sub-power series solutions $y_l(x)$ are truncated to polynomials with fixed parameter values in a numerator of $A_n$ term, compared with a polynomial of type 1.

With my definition, there are 2 types of complete polynomials such as the first species complete polynomial and the second species complete polynomial. 
The first species complete polynomial is obeyed since a parameter of a numerator in $B_n$ term and a (spectral) parameter of a numerator in $A_n$ term are fixed constants. The second species complete polynomial is applicable since two parameters of a numerator in $B_n$ term and a parameter of a numerator in $A_n$ term are fixed constants.
The former has multi-valued roots of a parameter of a numerator in $A_n$ term, but the latter has only one fixed parameter value of a numerator in $A_n$ term.   

In chapter 1, representation in terms of general summation formulas for the first and second species complete polynomials are constructed by allowing $A_n$ as a leading term in each of finite sub-power series of a general series solution $y(x)$. Their mathematical expressions are referred as `complete polynomials using 3-term recurrence formula (3TRF).''

In comparison with complete polynomials using 3TRF, general series solutions in closed forms for the first and second species complete polynomials are obtained by allowing $B_n$ as the leading term in each finite sub-power series of the general power series $y(x)$ in chapter 2. These classical summation formulas are designated as ``complete polynomials using reversible 3-term recurrence formula (R3TRF).''
\section{All possible formal series solutions in the DCHE} 
By the method of Frobenius, all potential formal series solutions of the DCHE are described in the following table.
\begin{table}[h]
\begin{center}
{
 \Tree[.{ The DCHE with irregular singular points at the origin and infinity} [.{ 3TRF} 
              [[[.{ Polynomial of type 3} [.{ $ \begin{array}{lcll}  1^{\mbox{st}}\;  \mbox{species}\\ \mbox{complete} \\ \mbox{polynomial} \end{array}$} ]  ]]]]                         
  [.{ R3TRF} 
     [ [[.{ Polynomial of type 2} ]
       [.{ Polynomial of type 3} [.{ $ \begin{array}{lcll}  1^{\mbox{st}}\;  \mbox{species} \\ \mbox{complete} \\ \mbox{polynomial} \end{array}$} ]  ]]]]]
}
\end{center}
\caption{Power series of the DCHE with irregular singular points at the origin and infinity}
\end{table}  

The non-symmetrical canonical form of the DCHE has 4 parameters such as $\alpha $, $\beta $, $\gamma $, $\delta $ and $q$.  
In chapter 12, the parameters play different roles for a polynomial of type 2 in the DCHE at the origin and infinity: $\alpha $, $\beta $, $\gamma $ and $\delta $ are treated as free variables; $q$ is considered as a fixed value.
Conspicuously, type 2 polynomials for the DCHE require $\left| \beta /\delta x^2\right|<1$ at the origin and $\left| \delta/\beta z^2\right|<1$ at infinity for the radius of convergence.
A polynomial of type 1 and an infinite series for the DCHE with both irregular singular points do not exist prominently because $A_n$ term is divergent as $n\gg 1$. 
Representation in terms of integrals of the DCHE at $x=0$ and $x=\infty $ are derived by applying the contour integral of Tricomi's functions into sub-power series of their general series solutions. 
By applying generating functions for Kummer's function of the second kind into each of sub-integrals of general integrals for the DCH polynomials of type 2, their generating functions with both irregular singular points are found analytically.

The first species complete polynomial of the DCHE with both irregular singular points requires $\beta $, $\gamma $, $\delta$ as free variables and $\alpha$, $q$ as fixed values. 
The second species complete polynomial of the DCHE around $x=0$ and $x=\infty $ are not available because $\alpha $ of a numerator in $B_n$ term in its differential equation is only a fixed parameter value for which $B_n$ term are truncated at a specific index summation $n$.

Spectral polynomials (the first species complete polynomial) of the symmetric canonical form of the DCHE, the special case of canonical form of the general DCHE, heretofore have been found by applying a Laurent series with a characteristic exponent $\nu $ \cite{12Meix1954,12Pham1970,12Ronv1995}. 

Ishkhanyan \textit {et al.} show that power series solutions of the non-symmetrical canonical form of the DCHE in terms of the Kummer functions of the first kind. They find that several Kummer functions lead to expansions the coefficients of which in general obey a 3-term recurrence relation. They present spectral polynomials of the DCHE for the special cases as recursion relations \cite{12Ishk2014}.

Power series solutions of coulomb wave functions including the DCHE (generalized spheroidal equation in the Leaver version) 
have been obtained from its Laurent series \cite{12Figu2002}. 
Figueiredo \textit {et al.} examined the Whittaker-Ince limits of the DCHE having polynomial solutions, and showed the  possibility having finite series in a fixed set for the DCHE \cite{12ElJa2012}.  

They just left analytic solutions of the DCH spectral polynomials as solutions of recurrences because of a 3-term recursive relation between successive coefficients in its Laurent series. And they left an algebraic equation of the $(j+1)$th order for the determination of a special parameter as the determinant of $(j+1)\times (j+1)$ matrices. Currently, there are still no general series solutions of the DCHE in which the coefficients are given explicitly.
 
In this chapter, by substituting $A_n$ and $B_n$ terms of the DCHE with both irregular singular points into complete polynomials using 3TRF and R3TRF, representation in terms of their power series solutions in closed forms are given for the first species complete polynomials. 
Furthermore, I shew algebraic equations of the DCHE for the determination of a parameter $q$ in the combinational form of partial sums of the sequences $\{A_n\}$ and $\{B_n\}$ using 3TRF and R3TRF.     
\section[The DCHE with a irregular singular point at $x=0$]{The DCHE with a irregular singular point at the origin}
We take the power series 
\begin{equation}
y(x)= \sum_{n=0}^{\infty } c_n x^{n+\lambda } \;\;\;\mbox{where}\;\;c_0 \ne 0 \label{eq:12005}
\end{equation}
where $\lambda $ as an indicial root. Substitute (\ref{eq:12005}) into (\ref{eq:12003}), and equalize coefficients of consecutive power of $x$ with zero. We observe the recurrence system for determining the coefficients $c_n$ of the series:
\begin{equation}
c_{n+1}=A_n \;c_n +B_n \;c_{n-1} \hspace{1cm};n\geq 1 \label{eq:12006}
\end{equation}
where,
\begin{subequations}
\begin{align}
 A_n &= -\frac{n\left( n-1+\gamma \right) -q}{\delta \left( n+1\right)} \label{eq:12007a}\\ 
&= -\frac{\left( n+\frac{\gamma -1-\sqrt{\left( \gamma -1\right)^2 +4q}}{2}\right) \left( n+\frac{\gamma -1+\sqrt{\left( \gamma -1\right)^2 +4q}}{2}\right)}{\delta \left( n+1\right)} \label{eq:12007b}\\
B_n &= -\frac{\beta \left( n-1+\alpha \right)}{\delta \left( n+1\right)}  \label{eq:12007c}\\
c_1 &= A_0 \;c_0 \label{eq:12007d}
\end{align}
\end{subequations} 
with $\delta \ne 0$. We only have one indicial root such as $\lambda = 0$.
\subsection{The first species complete polynomial using 3TRF}
For the first species complete polynomials using 3TRF and R3TRF, we need a condition such as
\begin{equation}
B_{j+1}= c_{j+1}=0\hspace{1cm}\mathrm{where}\;j\in \mathbb{N}_{0}  
 \label{eq:12008}
\end{equation}
(\ref{eq:12008}) gives successively $c_{j+2}=c_{j+3}=c_{j+4}=\cdots=0$. And $c_{j+1}=0$ is defined by a polynomial equation of degree $j+1$ for the determination of an accessory parameter in $A_n$ term. 
\begin{theorem}
In chapter 1, the general summation expression of a function $y(x)$ for the first species complete polynomial using 3-term recurrence formula and its algebraic equation for the determination of an accessory parameter in $A_n$ term are given by
\begin{enumerate} 
\item As $B_1=0$,
\begin{equation}
0 =\bar{c}(1,0) \label{eq:12009a}
\end{equation}
\begin{equation}
y(x) = y_{0}^{0}(x) \label{eq:12009b}
\end{equation}
\item As $B_{2N+2}=0$ where $N \in \mathbb{N}_{0}$,
\begin{equation}
0  = \sum_{r=0}^{N+1}\bar{c}\left( 2r, N+1-r\right) \label{eq:120010a}
\end{equation}
\begin{equation}
y(x)= \sum_{r=0}^{N} y_{2r}^{N-r}(x)+ \sum_{r=0}^{N} y_{2r+1}^{N-r}(x)  \label{eq:120010b}
\end{equation}
\item As $B_{2N+3}=0$ where $N \in \mathbb{N}_{0}$,
\begin{equation}
0  = \sum_{r=0}^{N+1}\bar{c}\left( 2r+1, N+1-r\right) \label{eq:120011a}
\end{equation}
\begin{equation}
y(x)= \sum_{r=0}^{N+1} y_{2r}^{N+1-r}(x)+ \sum_{r=0}^{N} y_{2r+1}^{N-r}(x)  \label{eq:120011b}
\end{equation}
In the above,
\begin{eqnarray}
\bar{c}(0,n)  &=& \prod _{i_0=0}^{n-1}B_{2i_0+1} \label{eq:120012a}\\
\bar{c}(1,n) &=&  \sum_{i_0=0}^{n} \left\{ A_{2i_0} \prod _{i_1=0}^{i_0-1}B_{2i_1+1} \prod _{i_2=i_0}^{n-1}B_{2i_2+2} \right\} 
\label{eq:120012b}\\
\bar{c}(\tau ,n) &=& \sum_{i_0=0}^{n} \left\{A_{2i_0}\prod _{i_1=0}^{i_0-1} B_{2i_1+1} 
\prod _{k=1}^{\tau -1} \left( \sum_{i_{2k}= i_{2(k-1)}}^{n} A_{2i_{2k}+k}\prod _{i_{2k+1}=i_{2(k-1)}}^{i_{2k}-1}B_{2i_{2k+1}+(k+1)}\right) \right.\nonumber\\
&&\times \left.\prod _{i_{2\tau}=i_{2(\tau -1)}}^{n-1} B_{2i_{2\tau }+(\tau +1)} \right\} 
\hspace{1cm}\label{eq:120012c} 
\end{eqnarray}
and
\begin{eqnarray}
y_0^m(x) &=& c_0 x^{\lambda } \sum_{i_0=0}^{m} \left\{ \prod _{i_1=0}^{i_0-1}B_{2i_1+1} \right\} x^{2i_0 } \label{eq:120013a}\\
y_1^m(x) &=& c_0 x^{\lambda } \sum_{i_0=0}^{m}\left\{ A_{2i_0} \prod _{i_1=0}^{i_0-1}B_{2i_1+1}  \sum_{i_2=i_0}^{m} \left\{ \prod _{i_3=i_0}^{i_2-1}B_{2i_3+2} \right\}\right\} x^{2i_2+1 } \label{eq:120013b}\\
y_{\tau }^m(x) &=& c_0 x^{\lambda } \sum_{i_0=0}^{m} \left\{A_{2i_0}\prod _{i_1=0}^{i_0-1} B_{2i_1+1} 
\prod _{k=1}^{\tau -1} \left( \sum_{i_{2k}= i_{2(k-1)}}^{m} A_{2i_{2k}+k}\prod _{i_{2k+1}=i_{2(k-1)}}^{i_{2k}-1}B_{2i_{2k+1}+(k+1)}\right) \right. \nonumber\\
&& \times  \left. \sum_{i_{2\tau} = i_{2(\tau -1)}}^{m} \left( \prod _{i_{2\tau +1}=i_{2(\tau -1)}}^{i_{2\tau}-1} B_{2i_{2\tau +1}+(\tau +1)} \right) \right\} x^{2i_{2\tau}+\tau }\hspace{1cm}\mathrm{where}\;\tau \geq 2
\label{eq:120013c} 
\end{eqnarray}
\end{enumerate}
\end{theorem}
Put $n= j+1$ in (\ref{eq:12007c}) and use the condition $B_{j+1}=0$ for $\alpha $.  
\begin{equation}
\alpha  = -j
\label{eq:120014}
\end{equation}
Take (\ref{eq:120014}) into (\ref{eq:12007c}).
\begin{equation}
B_n = -\frac{\beta (n-1-j)}{\delta (n+1)} \label{eq:120015}
\end{equation}
Now the condition $c_{j+1}=0$ is clearly an algebraic equation in $q$ of degree $j+1$ and thus has $j+1$ zeros denoted them by $q_j^m$ eigenvalues where $m = 0,1,2, \cdots, j$. They can be arranged in the following order: $q_j^0 < q_j^1 < q_j^2 < \cdots < q_j^j$.

Substitute (\ref{eq:12007a}) and (\ref{eq:120015}) into (\ref{eq:120012a})--(\ref{eq:120013c}) by letting $c_0=1$ and $\lambda =0$.

As $B_{1}= c_{1}=0$, take the new (\ref{eq:120012b}) into (\ref{eq:12009a}) putting $j=0$. Substitute the new (\ref{eq:120013a}) into (\ref{eq:12009b}) putting $j=0$. 

As $B_{2N+2}= c_{2N+2}=0$, take the new (\ref{eq:120012a})--(\ref{eq:120012c}) into (\ref{eq:120010a}) putting $j=2N+1$. Substitute the new 
(\ref{eq:120013a})--(\ref{eq:120013c}) into (\ref{eq:120010b}) putting $j=2N+1$ and $q =q_{2N+1}^m$.

As $B_{2N+3}= c_{2N+3}=0$, take the new (\ref{eq:120012a})--(\ref{eq:120012c}) into (\ref{eq:120011a}) putting $j=2N+2$. Substitute the new 
(\ref{eq:120013a})--(\ref{eq:120013c}) into (\ref{eq:120011b}) putting $j=2N+2$ and $q =q_{2N+2}^m$.
After the replacement process, we obtain an independent solution of the DCHE. The solution is as follows.
\begin{remark}
The power series expansion of the DCHE of the first kind about $x=0$ for the first species complete polynomial using 3TRF as $\alpha = -j$ where $j \in \mathbb{N}_{0}$ and its algebraic equation for the determination of an accessory parameter $q$ are given by
\begin{enumerate} 
\item As $\alpha =0$ and $q=q_0^0 =0$,

The eigenfunction is given by
\begin{equation}
y(x) = H_d^{(o)}F_{0,0}\left( \alpha =0,\beta ,\gamma ,\delta ,q=q_0^0 =0; \tilde{\eta } =-\frac{2}{\delta }x, \mu =-\frac{\beta }{\delta }x^2\right) =1 \label{eq:120016}
\end{equation}
\item As $\alpha = -\left( 2N+1 \right)$ where $N \in \mathbb{N}_{0}$,

An algebraic equation of degree $2N+2$ for the determination of $q$ is given by
\begin{equation}
0 = \sum_{r=0}^{N+1}\bar{c}\left( 2r, N+1-r; 2N+1,q\right)\label{eq:120017a}
\end{equation}
The eigenvalue of $q$ is written by $q_{2N+1}^m$ where $m = 0,1,2,\cdots,2N+1 $; $q_{2N+1}^0 < q_{2N+1}^1 < \cdots < q_{2N+1}^{2N+1}$. Its eigenfunction is given by
\begin{eqnarray} 
y(x) &=&  H_d^{(o)}F_{2N+1,m}\left( \alpha =-\left( 2N+1 \right),\beta ,\gamma ,\delta ,q=q_{2N+1}^m; \tilde{\eta } =-\frac{2}{\delta }x, \mu =-\frac{\beta }{\delta }x^2\right)\nonumber\\
&=& \sum_{r=0}^{N} y_{2r}^{N-r}\left( 2N+1, q_{2N+1}^m;x\right)+ \sum_{r=0}^{N} y_{2r+1}^{N-r}\left( 2N+1, q _{2N+1}^m;x\right)  
\label{eq:120017b}
\end{eqnarray}
\item As $\alpha = -\left( 2N+2\right)$ where $N \in \mathbb{N}_{0}$,

An algebraic equation of degree $2N+3$ for the determination of $q$ is given by
\begin{eqnarray}
0  = \sum_{r=0}^{N+1}\bar{c}\left( 2r+1, N+1-r; 2N+2,q\right)\label{eq:120018a}
\end{eqnarray}
The eigenvalue of $q$ is written by $q_{2N+2}^m$ where $m = 0,1,2,\cdots,2N+2 $; $q_{2N+2}^0 < q_{2N+2}^1 < \cdots < q_{2N+2}^{2N+2}$. Its eigenfunction is given by
\begin{eqnarray} 
y(x) &=& H_d^{(o)}F_{2N+2,m}\left( \alpha =-\left( 2N+2\right),\beta ,\gamma ,\delta ,q=q_{2N+2}^m; \tilde{\eta } =-\frac{2}{\delta }x, \mu =-\frac{\beta }{\delta }x^2\right)\nonumber\\
&=& \sum_{r=0}^{N+1} y_{2r}^{N+1-r}\left( 2N+2, q_{2N+2}^m;x\right) + \sum_{r=0}^{N} y_{2r+1}^{N-r}\left( 2N+2, q_{2N+2}^m;x\right) 
\label{eq:120018b}
\end{eqnarray}
In the above,
\begin{eqnarray}
\bar{c}(0,n;j,q)  &=& \frac{\left( -\frac{j}{2}\right)_{n}}{\left( 1 \right)_{n}} \left( -\frac{\beta }{\delta } \right)^{n}\label{eq:120019a}\\
\bar{c}(1,n;j,q) &=& \left( -\frac{2}{\delta } \right) \sum_{i_0=0}^{n}\frac{i_0 \left( i_0-\frac{1}{2}+\frac{\gamma }{2} \right)-\frac{q}{4} }{\left( i_0+\frac{1}{2} \right)} \frac{\left( -\frac{j}{2}\right)_{i_0} }{\left( 1\right)_{i_0}}   \frac{\left( \frac{1}{2}-\frac{j}{2} \right)_{n} \left( \frac{3}{2} \right)_{i_0}}{\left( \frac{1}{2}-\frac{j}{2}\right)_{i_0} \left( \frac{3}{2} \right)_{n}} \left( -\frac{\beta }{\delta } \right)^{n } \label{eq:120019b}\\
\bar{c}(\tau ,n;j,q) &=& \left( -\frac{2}{\delta } \right)^{\tau } \sum_{i_0=0}^{n}\frac{i_0 \left( i_0-\frac{1}{2}+\frac{\gamma }{2} \right)-\frac{q}{4}}{\left( i_0+\frac{1}{2} \right)} \frac{\left( -\frac{j}{2}\right)_{i_0} }{\left( 1\right)_{i_0}} \nonumber\\
&&\times \prod_{k=1}^{\tau -1} \left( \sum_{i_k = i_{k-1}}^{n} \frac{\left( i_k+ \frac{k}{2} \right)\left( i_k+ \frac{k}{2}-\frac{1}{2}+\frac{\gamma }{2} \right)-\frac{q}{4} }{\left( i_k+\frac{k}{2}+\frac{1}{2} \right)} \right.  \left. \frac{\left( \frac{k}{2}-\frac{j}{2}\right)_{i_k} \left( \frac{k}{2}+1 \right)_{i_{k-1}}}{\left( \frac{k}{2}-\frac{j}{2}\right)_{i_{k-1}} \left( \frac{k}{2}+1 \right)_{i_k}} \right) \nonumber\\ 
&&\times \frac{\left( \frac{\tau }{2} -\frac{j}{2}\right)_{n} \left( \frac{\tau }{2}+1 \right)_{i_{\tau -1}}}{\left( \frac{\tau }{2}-\frac{j}{2}\right)_{i_{\tau -1}} \left( \frac{\tau }{2}+1 \right)_{n}} \left( -\frac{\beta }{\delta } \right)^{n } \label{eq:120019c} 
\end{eqnarray}
\begin{eqnarray}
y_0^m(j,q;x) &=& \sum_{i_0=0}^{m} \frac{\left( -\frac{j}{2}\right)_{n}}{\left( 1 \right)_{n}} \mu ^{i_0} \label{eq:120020a}\\
y_1^m(j,q;x) &=& \left\{\sum_{i_0=0}^{m} \frac{i_0 \left( i_0-\frac{1}{2}+\frac{\gamma }{2} \right)-\frac{q}{4} }{\left( i_0+\frac{1}{2} \right)} \frac{\left( -\frac{j}{2}\right)_{i_0} }{\left( 1\right)_{i_0}} \right.   \left. \sum_{i_1 = i_0}^{m} \frac{\left( \frac{1}{2}-\frac{j}{2} \right)_{i_1} \left( \frac{3}{2} \right)_{i_0}}{\left( \frac{1}{2}-\frac{j}{2} \right)_{i_0} \left( \frac{3}{2} \right)_{i_1}} \mu ^{i_1}\right\} \tilde{\eta } \label{eq:120020b}\\
y_{\tau }^m(j,q;x) &=& \left\{ \sum_{i_0=0}^{m} \frac{i_0 \left( i_0-\frac{1}{2}+\frac{\gamma }{2} \right)-\frac{q}{4} }{\left( i_0+\frac{1}{2} \right)} \frac{\left( -\frac{j}{2}\right)_{i_0} }{\left( 1\right)_{i_0}} \right.\nonumber\\
&&\times \prod_{k=1}^{\tau -1} \left( \sum_{i_k = i_{k-1}}^{m} \frac{\left( i_k+ \frac{k}{2} \right)\left( i_k+ \frac{k}{2}-\frac{1}{2}+\frac{\gamma }{2} \right)-\frac{q}{4} }{\left( i_k+\frac{k}{2}+\frac{1}{2} \right)} \right.  \left. \frac{\left( \frac{k}{2}-\frac{j}{2}\right)_{i_k} \left( \frac{k}{2}+1 \right)_{i_{k-1}}}{\left( \frac{k}{2}-\frac{j}{2}\right)_{i_{k-1}} \left( \frac{k}{2}+1 \right)_{i_k}} \right) \nonumber\\
&&\times \left. \sum_{i_{\tau } = i_{\tau -1}}^{m} \frac{\left( \frac{\tau }{2}-\frac{j}{2}\right)_{i_{\tau }} \left( \frac{\tau }{2}+1 \right)_{i_{\tau -1}}}{\left( \frac{\tau }{2}-\frac{j}{2}\right)_{i_{\tau -1}} \left( \frac{\tau }{2}+1 \right)_{i_{\tau }} } \mu ^{i_{\tau }}\right\} \tilde{\eta }^{\tau }  \label{eq:120020c} 
\end{eqnarray}
where $\tau \geq 2$.
\end{enumerate}
\end{remark}
\subsection{The first species complete polynomial using R3TRF}
\begin{theorem}
In chapter 2, the general summation expression of a function $y(x)$ for the first species complete polynomial using reversible 3-term recurrence formula and its algebraic equation for the determination of an accessory parameter in $A_n$ term are given by
\begin{enumerate} 
\item As $B_1=0$,
\begin{equation}
0 =\bar{c}(0,1) \label{eq:120021a}
\end{equation}
\begin{equation}
y(x) = y_{0}^{0}(x) \label{eq:120021b}
\end{equation}
\item As $B_2=0$, 
\begin{equation}
0 = \bar{c}(0,2)+\bar{c}(1,0) \label{eq:120022a}
\end{equation}
\begin{equation}
y(x)= y_{0}^{1}(x) \label{eq:120022b}
\end{equation}
\item As $B_{2N+3}=0$ where $N \in \mathbb{N}_{0}$,
\begin{equation}
0  = \sum_{r=0}^{N+1}\bar{c}\left( r, 2(N-r)+3\right) \label{eq:120023a}
\end{equation}
\begin{equation}
y(x)= \sum_{r=0}^{N+1} y_{r}^{2(N+1-r)}(x) \label{eq:120023b}
\end{equation}
\item As $B_{2N+4}=0$ where$N \in \mathbb{N}_{0}$,
\begin{equation}
0  =  \sum_{r=0}^{N+2}\bar{c}\left( r, 2(N+2-r)\right) \label{eq:120024a}
\end{equation}
\begin{equation}
y(x)=  \sum_{r=0}^{N+1} y_{r}^{2(N-r)+3}(x) \label{eq:120024b}
\end{equation}
In the above,
\begin{eqnarray}
\bar{c}(0,n) &=& \prod _{i_0=0}^{n-1}A_{i_0} \label{eq:120025a}\\
\bar{c}(1,n) &=& \sum_{i_0=0}^{n} \left\{ B_{i_0+1} \prod _{i_1=0}^{i_0-1}A_{i_1} \prod _{i_2=i_0}^{n-1}A_{i_2+2} \right\} \label{eq:120025b}\\
\bar{c}(\tau ,n) &=& \sum_{i_0=0}^{n} \left\{B_{i_0+1}\prod _{i_1=0}^{i_0-1} A_{i_1} 
\prod _{k=1}^{\tau -1} \left( \sum_{i_{2k}= i_{2(k-1)}}^{n} B_{i_{2k}+(2k+1)}\prod _{i_{2k+1}=i_{2(k-1)}}^{i_{2k}-1}A_{i_{2k+1}+2k}\right) \right.\nonumber\\
&&\times \left. \prod _{i_{2\tau} = i_{2(\tau -1)}}^{n-1} A_{i_{2\tau }+ 2\tau} \right\} 
\hspace{1cm}\label{eq:120025c}
\end{eqnarray}
and
\begin{eqnarray}
y_0^m(x) &=& c_0 x^{\lambda} \sum_{i_0=0}^{m} \left\{ \prod _{i_1=0}^{i_0-1}A_{i_1} \right\} x^{i_0 } \label{eq:120026a}\\
y_1^m(x) &=& c_0 x^{\lambda} \sum_{i_0=0}^{m}\left\{ B_{i_0+1} \prod _{i_1=0}^{i_0-1}A_{i_1}  \sum_{i_2=i_0}^{m} \left\{ \prod _{i_3=i_0}^{i_2-1}A_{i_3+2} \right\}\right\} x^{i_2+2 } \label{eq:120026b}\\
y_{\tau }^m(x) &=& c_0 x^{\lambda} \sum_{i_0=0}^{m} \left\{B_{i_0+1}\prod _{i_1=0}^{i_0-1} A_{i_1} 
\prod _{k=1}^{\tau -1} \left( \sum_{i_{2k}= i_{2(k-1)}}^{m} B_{i_{2k}+(2k+1)}\prod _{i_{2k+1}=i_{2(k-1)}}^{i_{2k}-1}A_{i_{2k+1}+2k}\right) \right. \nonumber\\
&&\times \left. \sum_{i_{2\tau} = i_{2(\tau -1)}}^{m} \left( \prod _{i_{2\tau +1}=i_{2(\tau -1)}}^{i_{2\tau}-1} A_{i_{2\tau +1}+ 2\tau} \right) \right\} x^{i_{2\tau}+2\tau }\hspace{1cm}\mathrm{where}\;\tau \geq 2
\label{eq:120026c}
\end{eqnarray}
\end{enumerate}
\end{theorem}
According to (\ref{eq:12008}), $c_{j+1}=0$ is clearly an algebraic equation in $q$ of degree $j+1$ and thus has $j+1$ zeros denoted them by $q_j^m$ eigenvalues where $m = 0,1,2, \cdots, j$. They can be arranged in the following order: $q_j^0 < q_j^1 < q_j^2 < \cdots < q_j^j$.
 
Substitute (\ref{eq:12007b}) and (\ref{eq:120015}) into (\ref{eq:120025a})--(\ref{eq:120026c}) by letting $c_0=1$ and $\lambda =0$.

As $B_{1}= c_{1}=0$, take the new (\ref{eq:120025a}) into (\ref{eq:120021a}) putting $j=0$. Substitute the new (\ref{eq:120026a}) into (\ref{eq:120021b}) putting $j=0$.

As $B_{2}= c_{2}=0$, take the new (\ref{eq:120025a}) and (\ref{eq:120025b}) into (\ref{eq:120022a}) putting $j=1$. Substitute the new (\ref{eq:120026a}) into (\ref{eq:120022b}) putting $j=1$ and $q=q_1^m$. 

As $B_{2N+3}= c_{2N+3}=0$, take the new (\ref{eq:120025a})--(\ref{eq:120025c}) into (\ref{eq:120023a}) putting $j=2N+2$. Substitute the new 
(\ref{eq:120026a})--(\ref{eq:120026c}) into (\ref{eq:120023b}) putting $j=2N+2$ and $q=q_{2N+2}^m$.

As $B_{2N+4}= c_{2N+4}=0$, take the new (\ref{eq:120025a})--(\ref{eq:120025c}) into (\ref{eq:120024a}) putting $j=2N+3$. Substitute the new 
(\ref{eq:120026a})--(\ref{eq:120026c}) into (\ref{eq:120024b}) putting $j=2N+3$ and $q=q_{2N+3}^m$.
After the replacement process, we obtain an independent solution of the DCHE. The solution is as follows.
\begin{remark}
The power series expansion of the DCHE of the first kind about $x=0$ for the first species complete polynomial using R3TRF as $\alpha = -j$ where $j \in \mathbb{N}_{0}$ and its algebraic equation for the determination of an accessory parameter $q$ are given by
\begin{enumerate} 
\item As $\alpha =0$ and $q=q_0^0=0$,

The eigenfunction is given by
\begin{equation}
y(x) = H_d^{(o)}F_{0,0}^R\left( \alpha =0,\beta ,\gamma ,\delta ,q=q_0^0 =0; \eta =-\frac{1}{\delta }x, \mu =-\frac{\beta }{\delta }x^2\right) =1 \label{eq:120027}
\end{equation}
\item As $\alpha =-1$,

An algebraic equation of degree 2 for the determination of $q$ is given by
\begin{equation}
0 = q(q-\gamma ) +\beta \delta  \label{eq:120028a}
\end{equation}
The eigenvalue of $q$ is written by $q_1^m$ where $m = 0,1 $; $q_{1}^0 < q_{1}^1$. Its eigenfunction is given by
\begin{eqnarray}
y(x) &=& H_d^{(o)}F_{1,m}^R\left( \alpha =-1,\beta ,\gamma ,\delta ,q=q_1^m; \eta =-\frac{1}{\delta }x, \mu =-\frac{\beta }{\delta }x^2\right) \nonumber\\
&=&  1-q_1^m \eta \label{eq:120028b}  
\end{eqnarray}
\item As $\alpha =-\left( 2N+2 \right) $ where $N \in \mathbb{N}_{0}$,

An algebraic equation of degree $2N+3$ for the determination of $q$ is given by
\begin{equation}
0 = \sum_{r=0}^{N+1}\bar{c}\left( r, 2(N-r)+3; 2N+2,q\right)  \label{eq:120029a}
\end{equation}
The eigenvalue of $q$ is written by $q_{2N+2}^m$ where $m = 0,1,2,\cdots,2N+2 $; $q_{2N+2}^0 < q_{2N+2}^1 < \cdots < q_{2N+2}^{2N+2}$. Its eigenfunction is given by 
\begin{eqnarray} 
y(x) &=&  H_d^{(o)}F_{2N+2,m}^R\left( \alpha =-\left( 2N+2 \right),\beta ,\gamma ,\delta ,q=q_{2N+2}^m; \eta =-\frac{1}{\delta }x, \mu =-\frac{\beta }{\delta }x^2\right) \nonumber\\
&=& \sum_{r=0}^{N+1} y_{r}^{2(N+1-r)}\left( 2N+2, q_{2N+2}^m; x \right)  
\label{eq:120029b} 
\end{eqnarray}
\item As $\Omega =-\left( 2N+3 \right) $ where $N \in \mathbb{N}_{0}$,

An algebraic equation of degree $2N+4$ for the determination of $q$ is given by
\begin{equation}  
0 =  \sum_{r=0}^{N+2}\bar{c}\left( r, 2(N+2-r); 2N+3,q\right) \label{eq:120030a}
\end{equation}
The eigenvalue of $q$ is written by $q_{2N+3}^m$ where $m = 0,1,2,\cdots,2N+3 $; $q_{2N+3}^0 < q_{2N+3}^1 < \cdots < q_{2N+3}^{2N+3}$. Its eigenfunction is given by
\begin{eqnarray} 
y(x) &=& H_d^{(o)}F_{2N+3,m}^R\left( \alpha =-\left( 2N+3 \right),\beta ,\gamma ,\delta ,q=q_{2N+3}^m; \eta =-\frac{1}{\delta }x, \mu =-\frac{\beta }{\delta }x^2\right)\nonumber\\
&=& \sum_{r=0}^{N+1} y_{r}^{2(N-r)+3} \left( 2N+3,q_{2N+3}^m;x\right) \label{eq:120030b}
\end{eqnarray}
In the above,
\begin{eqnarray}
\bar{c}(0,n;j,q)  &=& \frac{\left( \Delta_0^{-} \left( q\right) \right)_{n}\left( \Delta_0^{+} \left( q\right) \right)_{n}}{\left( 1 \right)_{n}} \left( -\frac{1}{\delta }\right)^n\label{eq:120031a}\\
\bar{c}(1,n;j,q) &=& \left( -\frac{\beta }{\delta }\right) \sum_{i_0=0}^{n}\frac{\left( i_0 -j\right) }{\left( i_0+2 \right)} \frac{ \left( \Delta_0^{-} \left( q\right) \right)_{i_0}\left( \Delta_0^{+} \left( q\right) \right)_{i_0}}{\left( 1 \right)_{i_0}}     \frac{ \left( \Delta_1^{-} \left( q\right) \right)_{n}\left( \Delta_1^{+} \left( q\right) \right)_{n} \left( 3\right)_{i_0}}{\left( \Delta_1^{-} \left( q\right) \right)_{i_0}\left( \Delta_1^{+} \left(  q\right) \right)_{i_0}\left( 3 \right)_{n} }\left( -\frac{1}{\delta }\right)^n \hspace{1.5cm}\label{eq:120031b}\\
\bar{c}(\tau ,n;j,q) &=& \left( -\frac{\beta }{\delta }\right)^{\tau} \sum_{i_0=0}^{n}\frac{\left( i_0 -j\right) }{\left( i_0+2 \right)} \frac{ \left( \Delta_0^{-} \left( q\right) \right)_{i_0}\left( \Delta_0^{+} \left( q\right) \right)_{i_0}}{\left( 1\right)_{i_0}}  \nonumber\\
&&\times \prod_{k=1}^{\tau -1} \left( \sum_{i_k = i_{k-1}}^{n} \frac{\left( i_k+ 2k-j\right) }{\left( i_k+2k+2 \right)} \right.  \left. \frac{ \left( \Delta_k^{-} \left( q\right) \right)_{i_k}\left( \Delta_k^{+} \left( q\right) \right)_{i_k} \left( 2k+1 \right)_{i_{k-1}}}{\left( \Delta_k^{-} \left( q\right) \right)_{i_{k-1}}\left( \Delta_k^{+} \left( q\right) \right)_{i_{k-1}}\left( 2k+1 \right)_{i_k}} \right) \nonumber\\
&&\times \frac{ \left( \Delta_{\tau }^{-} \left( q\right) \right)_{n}\left( \Delta_{\tau }^{+} \left( q\right) \right)_{n} \left( 2\tau +1 \right)_{i_{\tau -1}}}{\left( \Delta_{\tau }^{-} \left( q\right) \right)_{i_{\tau -1}}\left( \Delta_{\tau }^{+} \left( q\right) \right)_{i_{\tau -1}}\left( 2\tau +1 \right)_{n}} \left( -\frac{1}{\delta }\right)^n \label{eq:120031c} 
\end{eqnarray}
\begin{eqnarray}
y_0^m(j,q;x) &=& \sum_{i_0=0}^{m} \frac{\left( \Delta_0^{-} \left( q\right) \right)_{i_0}\left( \Delta_0^{+} \left( q\right) \right)_{i_0}}{\left( 1 \right)_{i_0}} \eta ^{i_0} \label{eq:120032a}\\
y_1^m(j,q;x) &=& \left\{\sum_{i_0=0}^{m}\frac{\left( i_0 -j\right) }{\left( i_0+2 \right)} \frac{ \left( \Delta_0^{-} \left( q\right) \right)_{i_0}\left( \Delta_0^{+} \left( q\right) \right)_{i_0}}{\left( 1 \right)_{i_0}} \right. \left. \sum_{i_1 = i_0}^{m} \frac{ \left( \Delta_1^{-} \left( q\right) \right)_{i_1}\left( \Delta_1^{+} \left( q\right) \right)_{i_1} \left( 3 \right)_{i_0}}{\left( \Delta_1^{-} \left( q\right) \right)_{i_0}\left( \Delta_1^{+} \left( q\right) \right)_{i_0}\left( 3 \right)_{i_1}} \eta ^{i_1}\right\} \mu 
\hspace{1.5cm} \label{eq:120032b}\\
y_{\tau }^m(j,q;x) &=& \left\{ \sum_{i_0=0}^{m} \frac{\left( i_0 -j\right) }{\left( i_0+2 \right)} \frac{ \left( \Delta_0^{-} \left( q\right) \right)_{i_0}\left( \Delta_0^{+} \left( q\right) \right)_{i_0}}{\left( 1 \right)_{i_0}} \right.\nonumber\\
&&\times \prod_{k=1}^{\tau -1} \left( \sum_{i_k = i_{k-1}}^{m} \frac{\left( i_k+ 2k-j\right) }{\left( i_k+2k+2 \right)} \right.   \left. \frac{ \left( \Delta_k^{-} \left( q\right) \right)_{i_k}\left( \Delta_k^{+} \left( q\right) \right)_{i_k} \left( 2k+1 \right)_{i_{k-1}}}{\left( \Delta_k^{-} \left( q\right) \right)_{i_{k-1}}\left( \Delta_k^{+} \left( q\right) \right)_{i_{k-1}}\left( 2k+1 \right)_{i_k}} \right) \nonumber\\
&&\times \left. \sum_{i_{\tau } = i_{\tau -1}}^{m}  \frac{ \left( \Delta_{\tau }^{-} \left( q\right) \right)_{i_{\tau }}\left( \Delta_{\tau }^{+} \left( q\right) \right)_{i_{\tau }} \left( 2\tau +1 \right)_{i_{\tau -1}}}{\left( \Delta_{\tau }^{-} \left( q\right) \right)_{i_{\tau -1}}\left( \Delta_{\tau }^{+} \left( q\right) \right)_{i_{\tau -1}}\left( 2\tau +1 \right)_{i_\tau } } \eta ^{i_{\tau }}\right\} \mu ^{\tau } \label{eq:120032c} 
\end{eqnarray}
where
\begin{equation}
\begin{cases} \tau \geq 2 \cr
\Delta _{k}^{\pm}(q) = \frac{\gamma -1+4k \pm \sqrt{(\gamma -1)^2 +4q}}{2}
\end{cases}\nonumber
\end{equation}
\end{enumerate}
\end{remark}
\section[The DCHE with a irregular singular point at $x=\infty $]{The DCHE with a irregular singular point at infinity}
Let $z=\frac{1}{x}$ in (\ref{eq:12003}). Then the DCHE about $x=\infty $ takes the form
\begin{equation}
\frac{d^2{y}}{d{z}^2} + \left( -\delta +\frac{2-\gamma }{z} -\frac{\beta }{z^2}\right) \frac{d{y}}{d{z}} +  \frac{-qz+\alpha \beta }{z^3} y = 0 \label{eq:120033}
\end{equation}
assuming a solution on the form
\begin{equation}
y(z)= \sum_{n=0}^{\infty } c_n z^{n+\lambda }  \label{eq:120034}
\end{equation}
we obtain by substitution in (\ref{eq:120033}) a 3-term recurrence relation among the consecutive coefficients $c_n$:
\begin{equation}
c_{n+1}=A_n \;c_n +B_n \;c_{n-1} \hspace{1cm};n\geq 1 \label{eq:120035}
\end{equation}
where,
\begin{subequations}
\begin{align}
 A_n &= \frac{\left( n+\alpha \right)\left( n+1+\alpha -\gamma \right) -q}{\beta \left( n+1\right)} \label{eq:120036a}\\ 
&= \frac{\left( n+\frac{2\alpha -\gamma +1-\sqrt{\left( \gamma -1\right)^2 +4q}}{2}\right) \left( n+\frac{2\alpha -\gamma +1+\sqrt{\left( \gamma -1\right)^2 +4q}}{2}\right)}{\beta \left( n+1\right)} \label{eq:120036b}\\
B_n &= -\frac{\delta \left( n-1+\alpha \right)}{\beta \left( n+1\right)}  \label{eq:120036c}\\
c_1 &= A_0 \;c_0 \label{eq:120036d}
\end{align}
\end{subequations} 
with $\beta  \ne 0$. We only have one indicial root such as $\lambda = \alpha $.
\subsection{The first species complete polynomial using 3TRF}
Put $n= j+1$ in (\ref{eq:120036c}) and use the condition $B_{j+1}=0$ for $\alpha $.  
\begin{equation}
\alpha  = -j
\label{eq:120037}
\end{equation}
Take (\ref{eq:120037}) into (\ref{eq:120036a})--(\ref{eq:120036c}).
\begin{subequations}
\begin{align}
 A_n &= \frac{\left( n-j \right)\left( n+1-j-\gamma \right) -q}{\beta \left( n+1\right)} \label{eq:120038a}\\ 
&= \frac{\left( n+\frac{-\gamma +1-2j-\sqrt{\left( \gamma -1\right)^2 +4q}}{2}\right) \left( n+\frac{-\gamma +1-2j+\sqrt{\left( \gamma -1\right)^2 +4q}}{2}\right)}{\beta \left( n+1\right)} \label{eq:120038b}\\
B_n &= -\frac{\delta \left( n-1-j\right)}{\beta \left( n+1\right)}  \label{eq:120038c}
\end{align}
\end{subequations} 
The condition $c_{j+1}=0$ is clearly an algebraic equation in $q$ of degree $j+1$ and thus has $j+1$ zeros denoted them by $q_j^m$ eigenvalues where $m = 0,1,2, \cdots, j$. They can be arranged in the following order: $q_j^0 < q_j^1 < q_j^2 < \cdots < q_j^j$.

Substitute (\ref{eq:120038a}) and (\ref{eq:120038c}) into (\ref{eq:120012a})--(\ref{eq:120013c}) with replacing $x$, $c_0$ and $\lambda $ by $z$, $1$ and $\alpha $.

As $B_{1}= c_{1}=0$, take the new (\ref{eq:120012b}) into (\ref{eq:12009a}) putting $j=0$. Substitute the new (\ref{eq:120013a}) into (\ref{eq:12009b}) putting $j=0$ and $x=z$. 

As $B_{2N+2}= c_{2N+2}=0$, take the new (\ref{eq:120012a})--(\ref{eq:120012c}) into (\ref{eq:120010a}) putting $j=2N+1$. Substitute the new 
(\ref{eq:120013a})--(\ref{eq:120013c}) into (\ref{eq:120010b}) putting $j=2N+1$, $q =q_{2N+1}^m$ and $x=z$.

As $B_{2N+3}= c_{2N+3}=0$, take the new (\ref{eq:120012a})--(\ref{eq:120012c}) into (\ref{eq:120011a}) putting $j=2N+2$. Substitute the new 
(\ref{eq:120013a})--(\ref{eq:120013c}) into (\ref{eq:120011b}) putting $j=2N+2$, $q =q_{2N+2}^m$ and $x=z$.
After the replacement process, we obtain an independent solution of the DCHE. The solution is as follows.
\begin{remark}
The power series expansion of the DCHE of the first kind about $x=\infty $ for the first species complete polynomial using 3TRF as $\alpha = -j$ where $j \in \mathbb{N}_{0}$ and its algebraic equation for the determination of an accessory parameter $q$ are given by
\begin{enumerate} 
\item As $\alpha =0$ and $q=q_0^0 =0$,

The eigenfunction is given by
\begin{equation}
y(z) = H_d^{(i)}F_{0,0}\left( \alpha =0,\beta ,\gamma ,\delta ,q=q_0^0 =0; z=\frac{1}{x}, \tilde{\varepsilon } = \frac{2}{\beta }z, \rho =-\frac{\delta }{\beta }z^2\right) =z^{\alpha } \label{eq:120039}
\end{equation}
\item As $\alpha = -\left( 2N+1 \right)$ where $N \in \mathbb{N}_{0}$,

An algebraic equation of degree $2N+2$ for the determination of $q$ is given by
\begin{equation}
0 = \sum_{r=0}^{N+1}\bar{c}\left( 2r, N+1-r; 2N+1,q\right)\label{eq:120040a}
\end{equation}
The eigenvalue of $q$ is written by $q_{2N+1}^m$ where $m = 0,1,2,\cdots,2N+1 $; $q_{2N+1}^0 < q_{2N+1}^1 < \cdots < q_{2N+1}^{2N+1}$. Its eigenfunction is given by
\begin{eqnarray} 
y(z) &=&  H_d^{(i)}F_{2N+1,m}\left( \alpha =-\left( 2N+1 \right),\beta ,\gamma ,\delta ,q=q_{2N+1}^m; z=\frac{1}{x}, \tilde{\varepsilon } = \frac{2}{\beta }z, \rho =-\frac{\delta }{\beta }z^2\right)\nonumber\\
&=& z^{\alpha } \left\{ \sum_{r=0}^{N} y_{2r}^{N-r}\left( 2N+1, q_{2N+1}^m;z\right)+ \sum_{r=0}^{N} y_{2r+1}^{N-r}\left( 2N+1, q _{2N+1}^m;z\right)\right\}  
\label{eq:120040b}
\end{eqnarray}
\item As $\alpha = -\left( 2N+2\right)$ where $N \in \mathbb{N}_{0}$,

An algebraic equation of degree $2N+3$ for the determination of $q$ is given by
\begin{eqnarray}
0  = \sum_{r=0}^{N+1}\bar{c}\left( 2r+1, N+1-r; 2N+2,q\right)\label{eq:120041a}
\end{eqnarray}
The eigenvalue of $q$ is written by $q_{2N+2}^m$ where $m = 0,1,2,\cdots,2N+2 $; $q_{2N+2}^0 < q_{2N+2}^1 < \cdots < q_{2N+2}^{2N+2}$. Its eigenfunction is given by
\begin{eqnarray} 
y(z) &=&  H_d^{(i)}F_{2N+2,m}\left( \alpha =-\left( 2N+2\right),\beta ,\gamma ,\delta ,q=q_{2N+2}^m; z=\frac{1}{x}, \tilde{\varepsilon } = \frac{2}{\beta }z, \rho =-\frac{\delta }{\beta }z^2\right)\nonumber\\
&=& z^{\alpha } \left\{ \sum_{r=0}^{N+1} y_{2r}^{N+1-r}\left( 2N+2, q_{2N+2}^m;z\right) + \sum_{r=0}^{N} y_{2r+1}^{N-r}\left( 2N+2, q_{2N+2}^m;z\right) \right\} 
\label{eq:120041b}
\end{eqnarray}
In the above,
\begin{eqnarray}
\bar{c}(0,n;j,q)  &=& \frac{\left( -\frac{j}{2}\right)_{n}}{\left( 1 \right)_{n}} \left( -\frac{\delta }{\beta } \right)^{n}\label{eq:120042a}\\
\bar{c}(1,n;j,q) &=& \left( \frac{2}{\beta } \right) \sum_{i_0=0}^{n}\frac{\left( i_0 -\frac{j}{2} \right) \left( i_0+\frac{1}{2}-\frac{j}{2}-\frac{\gamma }{2} \right)-\frac{q}{4} }{\left( i_0+\frac{1}{2} \right)} \frac{\left( -\frac{j}{2}\right)_{i_0} }{\left( 1\right)_{i_0}}   \frac{\left( \frac{1}{2}-\frac{j}{2} \right)_{n} \left( \frac{3}{2} \right)_{i_0}}{\left( \frac{1}{2}-\frac{j}{2}\right)_{i_0} \left( \frac{3}{2} \right)_{n}} \left( -\frac{\delta }{\beta } \right)^{n } \hspace{1.5cm}\label{eq:120042b}\\
\bar{c}(\tau ,n;j,q) &=& \left( \frac{2}{\beta } \right)^{\tau } \sum_{i_0=0}^{n}\frac{\left( i_0 -\frac{j}{2} \right) \left( i_0+\frac{1}{2}-\frac{j}{2}-\frac{\gamma }{2} \right)-\frac{q}{4}}{\left( i_0+\frac{1}{2} \right)} \frac{\left( -\frac{j}{2}\right)_{i_0} }{\left( 1\right)_{i_0}} \nonumber\\
&&\times \prod_{k=1}^{\tau -1} \left( \sum_{i_k = i_{k-1}}^{n} \frac{\left( i_k+ \frac{k}{2}-\frac{j}{2} \right)\left( i_k+ \frac{k}{2}+\frac{1}{2}-\frac{j}{2}-\frac{\gamma }{2} \right)-\frac{q}{4} }{\left( i_k+\frac{k}{2}+\frac{1}{2} \right)} \right.  \left. \frac{\left( \frac{k}{2}-\frac{j}{2}\right)_{i_k} \left( \frac{k}{2}+1 \right)_{i_{k-1}}}{\left( \frac{k}{2}-\frac{j}{2}\right)_{i_{k-1}} \left( \frac{k}{2}+1 \right)_{i_k}} \right) \nonumber\\ 
&&\times \frac{\left( \frac{\tau }{2} -\frac{j}{2}\right)_{n} \left( \frac{\tau }{2}+1 \right)_{i_{\tau -1}}}{\left( \frac{\tau }{2}-\frac{j}{2}\right)_{i_{\tau -1}} \left( \frac{\tau }{2}+1 \right)_{n}} \left( -\frac{\delta }{\beta }\right)^{n } \label{eq:120042c} 
\end{eqnarray}
\begin{eqnarray}
y_0^m(j,q;z) &=& \sum_{i_0=0}^{m} \frac{\left( -\frac{j}{2}\right)_{n}}{\left( 1 \right)_{n}} \rho ^{i_0} \label{eq:120043a}\\
y_1^m(j,q;z) &=& \left\{\sum_{i_0=0}^{m} \frac{\left( i_0 -\frac{j}{2} \right) \left( i_0+\frac{1}{2}-\frac{j}{2}-\frac{\gamma }{2} \right)-\frac{q}{4} }{\left( i_0+\frac{1}{2} \right)} \frac{\left( -\frac{j}{2}\right)_{i_0} }{\left( 1\right)_{i_0}} \right.   \left. \sum_{i_1 = i_0}^{m} \frac{\left( \frac{1}{2}-\frac{j}{2} \right)_{i_1} \left( \frac{3}{2} \right)_{i_0}}{\left( \frac{1}{2}-\frac{j}{2} \right)_{i_0} \left( \frac{3}{2} \right)_{i_1}} \rho ^{i_1}\right\} \tilde{\varepsilon } \hspace{1.5cm}\label{eq:120043b}\\
y_{\tau }^m(j,q;z) &=& \left\{ \sum_{i_0=0}^{m} \frac{\left( i_0 -\frac{j}{2} \right) \left( i_0+\frac{1}{2}-\frac{j}{2}-\frac{\gamma }{2} \right)-\frac{q}{4} }{\left( i_0+\frac{1}{2} \right)} \frac{\left( -\frac{j}{2}\right)_{i_0} }{\left( 1\right)_{i_0}} \right.\nonumber\\
&&\times \prod_{k=1}^{\tau -1} \left( \sum_{i_k = i_{k-1}}^{m} \frac{\left( i_k+ \frac{k}{2}-\frac{j}{2} \right)\left( i_k+ \frac{k}{2}+\frac{1}{2}-\frac{j}{2}-\frac{\gamma }{2} \right)-\frac{q}{4} }{\left( i_k+\frac{k}{2}+\frac{1}{2} \right)} \right.  \left. \frac{\left( \frac{k}{2}-\frac{j}{2}\right)_{i_k} \left( \frac{k}{2}+1 \right)_{i_{k-1}}}{\left( \frac{k}{2}-\frac{j}{2}\right)_{i_{k-1}} \left( \frac{k}{2}+1 \right)_{i_k}} \right) \nonumber\\
&&\times \left. \sum_{i_{\tau } = i_{\tau -1}}^{m} \frac{\left( \frac{\tau }{2}-\frac{j}{2}\right)_{i_{\tau }} \left( \frac{\tau }{2}+1 \right)_{i_{\tau -1}}}{\left( \frac{\tau }{2}-\frac{j}{2}\right)_{i_{\tau -1}} \left( \frac{\tau }{2}+1 \right)_{i_{\tau }} } \rho ^{i_{\tau }}\right\} \tilde{\varepsilon }^{\tau }  \label{eq:120043c} 
\end{eqnarray}
where $\tau \geq 2$.
\end{enumerate}
\end{remark}
\subsection{The first species complete polynomial using R3TRF}
According to (\ref{eq:12008}), $c_{j+1}=0$ is clearly an algebraic equation in $q$ of degree $j+1$ and thus has $j+1$ zeros denoted them by $q_j^m$ eigenvalues where $m = 0,1,2, \cdots, j$. They can be arranged in the following order: $q_j^0 < q_j^1 < q_j^2 < \cdots < q_j^j$.
 
Substitute (\ref{eq:120038b}) and (\ref{eq:120038c}) into (\ref{eq:120025a})--(\ref{eq:120026c}) with $c_0=1$, $\lambda =\alpha $ and $x=z$.

As $B_{1}= c_{1}=0$, take the new (\ref{eq:120025a}) into (\ref{eq:120021a}) putting $j=0$. Substitute the new (\ref{eq:120026a}) into (\ref{eq:120021b}) putting $j=0$ and $x=z$.

As $B_{2}= c_{2}=0$, take the new (\ref{eq:120025a}) and (\ref{eq:120025b}) into (\ref{eq:120022a}) putting $j=1$. Substitute the new (\ref{eq:120026a}) into (\ref{eq:120022b}) putting $j=1$, $q=q_1^m$ and $x=z$. 

As $B_{2N+3}= c_{2N+3}=0$, take the new (\ref{eq:120025a})--(\ref{eq:120025c}) into (\ref{eq:120023a}) putting $j=2N+2$. Substitute the new 
(\ref{eq:120026a})--(\ref{eq:120026c}) into (\ref{eq:120023b}) putting $j=2N+2$, $q=q_{2N+2}^m$ and $x=z$.

As $B_{2N+4}= c_{2N+4}=0$, take the new (\ref{eq:120025a})--(\ref{eq:120025c}) into (\ref{eq:120024a}) putting $j=2N+3$. Substitute the new 
(\ref{eq:120026a})--(\ref{eq:120026c}) into (\ref{eq:120024b}) putting $j=2N+3$, $q=q_{2N+3}^m$ and $x=z$.
After the replacement process, we obtain an independent solution of the DCHE. The solution is as follows.
\begin{remark}
The power series expansion of the DCHE of the first kind about $x=\infty $ for the first species complete polynomial using R3TRF as $\alpha = -j$ where $j \in \mathbb{N}_{0}$ and its algebraic equation for the determination of an accessory parameter $q$ are given by
\begin{enumerate} 
\item As $\alpha =0$ and $q=q_0^0=0$,

The eigenfunction is given by
\begin{equation}
y(z) = H_d^{(i)}F_{0,0}^R\left( \alpha =0,\beta ,\gamma ,\delta ,q=q_0^0 =0; z=\frac{1}{x}, \xi = \frac{1}{\beta }z, \rho =-\frac{\delta }{\beta }z^2\right) =z^{\alpha } \label{eq:120044}
\end{equation}
\item As $\alpha =-1$,

An algebraic equation of degree 2 for the determination of $q$ is given by
\begin{equation}
0 = q(q-\gamma ) +\beta \delta  \label{eq:120045a}
\end{equation}
The eigenvalue of $q$ is written by $q_1^m$ where $m = 0,1 $; $q_{1}^0 < q_{1}^1$. Its eigenfunction is given by
\begin{eqnarray}
y(z) &=& H_d^{(i)}F_{1,m}^R\left( \alpha =-1,\beta ,\gamma ,\delta ,q=q_1^m; z=\frac{1}{x}, \xi = \frac{1}{\beta }z, \rho =-\frac{\delta }{\beta }z^2\right) \nonumber\\
&=& z^{\alpha }\left\{ 1+(\gamma -q_1^m )\xi \right\} \label{eq:120045b}  
\end{eqnarray}
\item As $\alpha =-\left( 2N+2 \right) $ where $N \in \mathbb{N}_{0}$,

An algebraic equation of degree $2N+3$ for the determination of $q$ is given by
\begin{equation}
0 = \sum_{r=0}^{N+1}\bar{c}\left( r, 2(N-r)+3; 2N+2,q\right)  \label{eq:120046a}
\end{equation}
The eigenvalue of $q$ is written by $q_{2N+2}^m$ where $m = 0,1,2,\cdots,2N+2 $; $q_{2N+2}^0 < q_{2N+2}^1 < \cdots < q_{2N+2}^{2N+2}$. Its eigenfunction is given by 
\begin{eqnarray} 
y(z) &=&  H_d^{(i)}F_{2N+2,m}^R\left( \alpha =-\left( 2N+2 \right),\beta ,\gamma ,\delta ,q=q_{2N+2}^m; z=\frac{1}{x}, \xi = \frac{1}{\beta }z, \rho =-\frac{\delta }{\beta }z^2\right) \nonumber\\
&=& z^{\alpha } \sum_{r=0}^{N+1} y_{r}^{2(N+1-r)}\left( 2N+2, q_{2N+2}^m;z\right)  
\label{eq:120046b} 
\end{eqnarray}
\item As $\alpha =-\left( 2N+3 \right) $ where $N \in \mathbb{N}_{0}$,

An algebraic equation of degree $2N+4$ for the determination of $q$ is given by
\begin{equation}  
0 =  \sum_{r=0}^{N+2}\bar{c}\left( r, 2(N+2-r); 2N+3,q\right) \label{eq:120047a}
\end{equation}
The eigenvalue of $q$ is written by $q_{2N+3}^m$ where $m = 0,1,2,\cdots,2N+3 $; $q_{2N+3}^0 < q_{2N+3}^1 < \cdots < q_{2N+3}^{2N+3}$. Its eigenfunction is given by
\begin{eqnarray} 
y(x) &=& H_d^{(i)}F_{2N+3,m}^R\left( \alpha =-\left( 2N+3 \right),\beta ,\gamma ,\delta ,q=q_{2N+3}^m; z=\frac{1}{x}, \xi = \frac{1}{\beta }z, \rho =-\frac{\delta }{\beta }z^2\right)\nonumber\\
&=& z^{\alpha } \sum_{r=0}^{N+1} y_{r}^{2(N-r)+3} \left( 2N+3,q_{2N+3}^m;z\right) \label{eq:120047b}
\end{eqnarray}
In the above,
\begin{eqnarray}
\bar{c}(0,n;j,q)  &=& \frac{\left( \Delta_0^{-} \left( j,q\right) \right)_{n}\left( \Delta_0^{+} \left( j,q\right) \right)_{n}}{\left( 1 \right)_{n}} \left( \frac{1}{\beta }\right)^n\label{eq:120048a}\\
\bar{c}(1,n;j,q) &=& \left( -\frac{\delta }{\beta }\right) \sum_{i_0=0}^{n}\frac{\left( i_0 -j\right) }{\left( i_0+2 \right)} \frac{ \left( \Delta_0^{-} \left( j,q\right) \right)_{i_0}\left( \Delta_0^{+} \left( j,q\right) \right)_{i_0}}{\left( 1 \right)_{i_0}} \nonumber\\
&&\times \frac{ \left( \Delta_1^{-} \left( j,q\right) \right)_{n}\left( \Delta_1^{+} \left( j,q\right) \right)_{n} \left( 3\right)_{i_0}}{\left( \Delta_1^{-} \left( j,q\right) \right)_{i_0}\left( \Delta_1^{+} \left( j,q\right) \right)_{i_0}\left( 3 \right)_{n} }\left( \frac{1}{\beta }\right)^n  \label{eq:120048b}\\
\bar{c}(\tau ,n;j,q) &=& \left( -\frac{\delta }{\beta }\right)^{\tau} \sum_{i_0=0}^{n}\frac{\left( i_0 -j\right) }{\left( i_0+2 \right)} \frac{ \left( \Delta_0^{-} \left( j,q\right) \right)_{i_0}\left( \Delta_0^{+} \left( j,q\right) \right)_{i_0}}{\left( 1\right)_{i_0}}  \nonumber\\
&&\times \prod_{k=1}^{\tau -1} \left( \sum_{i_k = i_{k-1}}^{n} \frac{\left( i_k+ 2k-j\right) }{\left( i_k+2k+2 \right)} \right.  \left. \frac{ \left( \Delta_k^{-} \left( j,q\right) \right)_{i_k}\left( \Delta_k^{+} \left( j,q\right) \right)_{i_k} \left( 2k+1 \right)_{i_{k-1}}}{\left( \Delta_k^{-} \left( j,q\right) \right)_{i_{k-1}}\left( \Delta_k^{+} \left( j,q\right) \right)_{i_{k-1}}\left( 2k+1 \right)_{i_k}} \right) \nonumber\\
&&\times \frac{ \left( \Delta_{\tau }^{-} \left( j,q\right) \right)_{n}\left( \Delta_{\tau }^{+} \left( j,q\right) \right)_{n} \left( 2\tau +1 \right)_{i_{\tau -1}}}{\left( \Delta_{\tau }^{-} \left( j,q\right) \right)_{i_{\tau -1}}\left( \Delta_{\tau }^{+} \left( j,q\right) \right)_{i_{\tau -1}}\left( 2\tau +1 \right)_{n}} \left( \frac{1}{\beta }\right)^n \label{eq:120048c} 
\end{eqnarray}
\begin{eqnarray}
y_0^m(j,q;z) &=& \sum_{i_0=0}^{m} \frac{\left( \Delta_0^{-} \left( j,q\right) \right)_{i_0}\left( \Delta_0^{+} \left( j,q\right) \right)_{i_0}}{\left( 1 \right)_{i_0}} \xi ^{i_0} \label{eq:120049a}\\
y_1^m(j,q;z) &=& \left\{\sum_{i_0=0}^{m}\frac{\left( i_0 -j\right) }{\left( i_0+2 \right)} \frac{ \left( \Delta_0^{-} \left( j,q\right) \right)_{i_0}\left( \Delta_0^{+} \left( j,q\right) \right)_{i_0}}{\left( 1 \right)_{i_0}} \right.\nonumber\\
&&\times \left. \sum_{i_1 = i_0}^{m} \frac{ \left( \Delta_1^{-} \left( j,q\right) \right)_{i_1}\left( \Delta_1^{+} \left( j,q\right) \right)_{i_1} \left( 3 \right)_{i_0}}{\left( \Delta_1^{-} \left( j,q\right) \right)_{i_0}\left( \Delta_1^{+} \left( j,q\right) \right)_{i_0}\left( 3 \right)_{i_1}} \xi ^{i_1}\right\} \rho 
 \label{eq:120049b}\\
y_{\tau }^m(j,q;z) &=& \left\{ \sum_{i_0=0}^{m} \frac{\left( i_0 -j\right) }{\left( i_0+2 \right)} \frac{ \left( \Delta_0^{-} \left( j,q\right) \right)_{i_0}\left( \Delta_0^{+} \left( j,q\right) \right)_{i_0}}{\left( 1 \right)_{i_0}} \right.\nonumber\\
&&\times \prod_{k=1}^{\tau -1} \left( \sum_{i_k = i_{k-1}}^{m} \frac{\left( i_k+ 2k-j\right) }{\left( i_k+2k+2 \right)} \right.   \left. \frac{ \left( \Delta_k^{-} \left( j,q\right) \right)_{i_k}\left( \Delta_k^{+} \left( j,q\right) \right)_{i_k} \left( 2k+1 \right)_{i_{k-1}}}{\left( \Delta_k^{-} \left( j,q\right) \right)_{i_{k-1}}\left( \Delta_k^{+} \left( j,q\right) \right)_{i_{k-1}}\left( 2k+1 \right)_{i_k}} \right) \nonumber\\
&&\times \left. \sum_{i_{\tau } = i_{\tau -1}}^{m}  \frac{ \left( \Delta_{\tau }^{-} \left( j,q\right) \right)_{i_{\tau }}\left( \Delta_{\tau }^{+} \left( j,q\right) \right)_{i_{\tau }} \left( 2\tau +1 \right)_{i_{\tau -1}}}{\left( \Delta_{\tau }^{-} \left( j,q\right) \right)_{i_{\tau -1}}\left( \Delta_{\tau }^{+} \left( j,q\right) \right)_{i_{\tau -1}}\left( 2\tau +1 \right)_{i_\tau } } \xi ^{i_{\tau }}\right\} \rho ^{\tau } \label{eq:120049c} 
\end{eqnarray}
where
\begin{equation}
\begin{cases} \tau \geq 2 \cr
\Delta _{k}^{\pm}\left( j,q\right) = \frac{-\gamma +1-2j+4k \pm \sqrt{(\gamma -1)^2 +4q}}{2}
\end{cases}\nonumber
\end{equation}
\end{enumerate}
\end{remark}
\section{Summary}

There are only two possible power series solutions in the DCHE with two irregular singular points such as a polynomial of type 2 and the first species complete polynomial: (1) coefficients $\alpha $, $\beta $, $\gamma $, $\delta $ are treated as free variables and a parameter $q$ as a fixed value for a polynomial of type 2; (2) $\beta $, $\gamma $, $\delta $ as free variables and 
$\alpha $, $q$ as fixed values for the first species complete polynomial.

In particular, specific parameter values $q$ for a polynomial of type 2 make that sub-power series of general power series solutions of the DCHE with both singular points are truncated to polynomials for a polynomial of type 2: a coefficient $q$ has infinite eigenvalues. In contrast, the first species complete polynomial has multi-valued roots of a coefficient $q$, referred as a spectral parameter. 
In the quantum point of view, a polynomial of type 2 is applicable since an eigenvalue is located at $A_n$ term in the recurrence system. If one of two eigenvalues at $A_n$ term and an another value at $B_n$ term are placed in the recursive relation between successive coefficients, the first species complete polynomial is utilized.
 
Numerical calculations of the DCHE with two singular points for the first species complete polynomial using 3TRF and R3TRF are tantamount to each other. However, general power series structures between two type polynomials have several important difference such as:
  
(1) $A_n$ is the leading term in each of finite sub-formal series of general Frobenius solutions for the first species complete polynomials using 3TRF. On the contrary, $B_n$ is the leading one in each of finite sub-power series for the first species complete polynomials using R3TRF.

(2) Denominators and numerators in all $B_n$ terms of each finite sub-power series arise with Pochhammer symbols for the first species complete polynomial using 3TRF. And the denominators and numerators in all $A_n$ terms of each finite sub-formal series arise with Pochhammer symbols for the first species complete polynomial using R3TRF.

(3) Classical summation structure consist of the sum of two finite sub-power series solutions of general formal series solutions for the first species complete polynomial using 3TRF. But general formal series solutions for the first species complete polynomial using R3TRF are only composed of one finite sub-power series solutions.  

In chapter 12, contour integrals of Tricomi's functions are applied into each sub-power series of general series solutions in the DCHE with two singular points for a polynomial of type 2, therethrough their integral representation are constructed analytically. 
Similarly, in the future papers, integrals of hypergeometric-type functions are employed in finite sub-power series of general formal series solutions of the DCHE for the first species complete polynomials using 3TRF and R3TRF, and their representation in terms of integral for this type 3 polynomial are given explicitly.
  
\addcontentsline{toc}{section}{Bibliography}
\bibliographystyle{model1a-num-names}
\bibliography{<your-bib-database>}

\end{document}